\pdfoutput=1
\documentclass[english]{shinybook}
\makeindex[title=Index,intoc,columns=2]
\indexsetup{headers={\indexname}{\indexname}}

\newcommand{\term}[2][]{\emph{#2}\if\relax\detokenize{#1}\relax\index{#2}\else\index{#1}\fi}
\newcommand{\noemphterm}[2][]{#2\if\relax\detokenize{#1}\relax\index{#2}\else\index{#1}\fi}

\usepackage{asymptote}
\begin{asydef}
    texpreamble("\usepackage[T1]{fontenc}\usepackage{libertine}\usepackage[libertine,smallerops,liby,vvarbb]{newtxmath}\newcommand{\Rbar}{\overline{\mathbb{R}}}
\newcommand{\N}{\mathbb{N}}
\newcommand{\Q}{\mathbb{Q}}
\newcommand{\R}{\mathbb{R}}
\newcommand{\1}{\mathbb{1}}
\newcommand{\calA}{\mathcal{A}}
\newcommand{\calI}{\mathcal{I}}
\newcommand{\calR}{\mathcal{R}}

\renewcommand{\implies}{\Rightarrow}
\newcommand{\setto}{\rightrightarrows}
\newcommand{\equivalent}{\Leftrightarrow}

\newcommand{\eps}{\varepsilon}
\renewcommand{\epsilon}{\varepsilon}
\renewcommand{\phi}{\varphi}
\newcommand{\Id}{\mathrm{Id}}
\newcommand{\prox}{\mathrm{prox}}
\newcommand{\proj}{\mathrm{proj}}
\newcommand{\soft}{\mathrm{soft}}

\newcommand{\defeq}{\coloneq}
\newcommand{\setof}[2]{\left\{#1 \;\middle|\; #2\right\}}
\newcommand{\abs}[1]{|#1|}
\newcommand{\norm}[1]{\|#1\|}
\newcommand{\adaptnorm}[1]{\left\|#1\right\|}
\newcommand{\innerprod}[2]{( #1 \,|\, #2)}
\newcommand{\dual}[1]{\langle #1 \rangle}
\newcommand{\iprod}[2]{\langle #1, #2\rangle}
\newcommand{\dualprod}[2]{\langle #1, #2\rangle}

\def\weaki{\protect{\u{\i}}}

\DeclareMathOperator{\dom}{\mathrm{dom}}
\DeclareMathOperator{\kernel}{\mathrm{ker}}
\DeclareMathOperator{\range}{\mathrm{ran}}
\DeclareMathOperator{\graph}{\mathrm{graph}}
\DeclareMathOperator{\epi}{\mathrm{epi}}
\DeclareMathOperator{\sub}{\mathrm{sub}}
\DeclareMathOperator{\sign}{\mathrm{sign}}
\DeclareMathOperator{\interior}{\mathrm{int}}
\DeclareMathOperator*{\lip}{\mathrm{lip}}
\DeclareMathOperator{\conv}{\mathrm{co}}
\DeclareMathOperator{\closure}{\mathrm{cl}}
\DeclareMathOperator{\ri}{\mathrm{ri}}
\DeclareMathOperator{\BD}{\mathrm{bd}}

\DeclareMathOperator{\dist}{\mathrm{dist}}
\DeclareMathOperator*{\argmin}{\mathrm{arg\,min}}
\DeclareMathOperator*{\argmax}{\mathrm{arg\,max}}
\DeclareMathOperator{\reg}{\mathrm{reg}}
\DeclareMathOperator{\subreg}{\mathrm{subreg}}
\DeclareMathOperator{\calm}{\mathrm{calm}}

\def\ddd{\,d}
\newcommand{\inv}[1]{#1^{-1}}
\newcommand{\pinv}[1]{#1^{\dagger}}
\newcommand{\grad}{\nabla}
\newcommand{\B}{\mathbb{B}}
\newcommand{\OB}{\mathbb{O}}
\newcommand{\frechetNormal}{\widehat N}
\newcommand{\frechetCod}{\widehat D^*}
\newcommand{\clarkeTangent}{\widehat T}
\newcommand{\clarkeGD}{\widehat D}
\newcommand{\coderivative}{D^*}
\newcommand{\polar}[1]{#1^{\circ}}
\newcommand{\bipolar}[1]{#1^{\circ\circ}}
\newcommand{\upperadj}[1]{#1^{\circ+}}
\newcommand{\loweradj}[1]{#1^{\circ-}}
\newcommand{\subdiff}{\partial}
\DeclareRobustCommand{\downto}{{{\mathchoice%
            {\rotatebox[origin=c]{-20}{$\to$}}%
            {\rotatebox[origin=c]{-20}{$\to$}}%
            {\rotatebox[origin=c]{-20}{\scalebox{0.75}{$\to$}}}%
            {\rotatebox[origin=c]{-20}{\scalebox{0.6}{$\to$}}}%
}}}
\DeclareRobustCommand{\upto}{{{\mathchoice%
            {\rotatebox[origin=c]{20}{$\to$}}%
            {\rotatebox[origin=c]{20}{$\to$}}%
            {\rotatebox[origin=c]{20}{\scalebox{0.75}{$\to$}}}%
            {\rotatebox[origin=c]{20}{\scalebox{0.6}{$\to$}}}%
}}}
\newcommand{\freevar}{{\,\boldsymbol\cdot\,}}
\newcommand{\Union}{\bigcup}
\newcommand{\Isect}{\bigcap}
\newcommand{\union}{\cup}
\newcommand{\isect}{\cap}

\newcommand{\weakto}{\mathrel{\rightharpoonup}}
\makeatletter
\def \weaktostar@sym{\setbox0=\hbox{$\rightharpoonup$}\rlap{\hbox 
to\wd0{\hss\raise1ex\hbox{$\scriptscriptstyle{*\,}$}\hss}}\box0}
\def \weaktostar    {\mathrel{\weaktostar@sym}}
\makeatother

\def\linear{\mathbb{L}}

\def\opt#1{\bar #1}
\def\realopt#1{\widehat #1}
\def\this#1{#1^k}
\def\nexxt#1{#1^{k+1}}
\def\overnext#1{\bar #1^{k+1}}
\def\prev#1{#1^{k-1}}

\def\optu{{\opt{u}}}
\def\optw{{\opt{w}}}
\def\optx{{\opt{x}}}
\def\opty{{\opt{y}}}
\def\optz{{\opt{z}}}

\def\realoptu{{\realopt{u}}}
\def\realoptw{{\realopt{w}}}
\def\realoptx{{\realopt{x}}}
\def\realopty{{\realopt{y}}}
\def\realoptz{{\realopt{z}}}

\def\nextu{\nexxt{u}}
\def\nextw{\nexxt{w}}
\def\nextx{\nexxt{x}}
\def\nexty{\nexxt{y}}
\def\nextz{\nexxt{z}}

\def\thisu{\this{u}}
\def\thisw{\this{w}}
\def\thisx{\this{x}}
\def\thisy{\this{y}}
\def\thisz{\this{z}}

\def\prevx{\prev{x}}

\def\overnextx{\overnext{x}}

\def\E{\mathbb{E}}

\newcommand{\rinv}[1]{#1^{-R}}
\newcommand{\linv}[1]{#1^{-L}}
\newcommand{\rinvstar}[1]{#1^{-R,*}}
\newcommand{\linvstar}[1]{#1^{-L,*}}
\newcommand{\starlinv}[1]{#1^{*,-L}}

\def\Happrox{\widetilde H}
\def\Hany{H}
\def\Hsaddle{H}
\def\tauTest{\phi}
\def\Tau{T}
\def\TauTest{\Phi}
\def\SigmaTest{\Psi}
\def\sigmaTest{\psi}
\def\gap{\mathcal{G}}
\def\GammaLift#1{\Xi_{#1}}
\def\Space{U}

\newcommand{\Test}{Z}
\newcommand{\Precond}{M}
\newcommand{\Step}{W}
\newcommand{\GenGap}{\mathcal{V}}

\def\kgrad#1{\grad K(#1)}
\def\kgradconj#1{[\grad K(#1)]^*}

\def\dir#1{\Delta #1}
\def\alt#1{\tilde #1}

\def\minval#1{#1_{\mathrm{min}}}

\def\kgrad#1{\grad K(#1)}
\def\kgradconj#1{[\grad K(#1)]^*}
\def\neighu{\mathcal{U}}
\def\neighx{\mathcal{X}}
\def\metricRhoY{\rho_y}

\newcommand{\mathvar}[1]{\textup{#1}}
\def\NL{\mathvar{NL}}
\def\LIN{\mathvar{L}}
\def\Ynl{Y_{\NL}}
\def\Ylin{Y_{\LIN}}
\def\Pnl{P_{\NL}}

\def\weaklim{\mathop{\operatorname{w-\kern.07em lim}\,}}
\def\weaklimsup{\mathop{\operatorname{w-\kern.07em lim\,sup}\,}}
\def\weakliminf{\mathop{\operatorname{w-\kern.07em lim\,inf}\,}}
\def\weakstarlimsup{\mathop{\operatorname{w-\!\ast\!-\kern.07em lim\,sup}\,}}
\def\weakstarliminf{\mathop{\operatorname{w-\!\ast\!-\kern.07em lim\,inf}\,}}

\def\lebesgue{\mathcal{L}}

\newcommand{\Uad}{U_{\mathrm{ad}}}
\newcommand{\Yad}{Y_{\mathrm{ad}}}
");
    import graph;
    import math;
    unitsize(50,50);
    pair pt(real f(real), real x) { return (x, f(x)); };
    pair orthog(pair tan, real sign=1) { return sign*(tan.y, -tan.x); }
    DefaultHead.size=new real(pen p=currentpen) {return 2mm;};
    pen dotted = linetype(new real[] {0,3})+linewidth(1.1);
    pen shortdashed = linetype(new real[] {4,4});
    pen lightfill = rgb(0.09,0.09,0.44)+opacity(.15);
    pen darkfill = rgb(0.09,0.09,0.44)+opacity(.3);
    pen primalline = rgb(0.09,0.09,0.44);
    pen dualfill = rgb(0.5,0.5,0)+opacity(.33);
    pen dualline = rgb(0.5,0.5,0);
    pen tertfill = rgb(0.72,0.45,0.2)+opacity(.33);
    pen tertline = rgb(0.72,0.45,0.2);
    pen violation = dualfill;
\end{asydef}

\usepackage{tikz,pgfplots,pgfplotstable}
\usepgfplotslibrary{colorbrewer,fillbetween}
\pgfplotsset{compat=newest}
\pgfplotsset{plot coordinates/math parser=false}
\usetikzlibrary{plotmarks}
\pgfplotsset{
    table/search path={numerics/},
}
\pgfplotsset{every tick label/.append style={font=\scriptsize}}

\pgfplotsset{
    tick label style = {/pgf/number format/fixed, /pgf/number format/set thousands separator={\,}},
    legend style = {
        inner sep=0pt, outer sep=10pt, legend cell align=left, align=left,
        draw=none, fill=none, font=\scriptsize
    },
    width = 0.95\linewidth,
    height = 0.3\linewidth,
    scaled x ticks=false,
    xminorticks=true,
    unbounded coords=jump,
    axis x line*=bottom,
    axis y line*=left,
    yminorticks=true,
}

\definecolor{primal}{rgb}{0.09,0.09,0.44}
\definecolor{dual}{rgb}{0.5,0.5,0}
\definecolor{tert}{rgb}{0.72,0.45,0.2}

\definecolor{cb1}{HTML}{1b9e77}
\definecolor{cb2}{HTML}{d95f02}
\definecolor{cb3}{HTML}{7570b3}
\definecolor{cb4}{HTML}{e7298a}
\definecolor{cb5}{HTML}{66a61e}

\pgfplotsset{
    fb/.style                 = { color = cb3, line width = 1pt },
    projgrad/.style           = { color = cb3, line width = 1.25pt, dotted },
    fista/.style              = { color = cb3, line width = 0.75pt, dashed },
    fbrelax/.style            = { color = cb3, line width = 1.25pt, dotted },
    pdps/.style               = { color = cb1, line width = 1pt },
    pdps accel/.style         = { color = cb1, line width = 0.75pt, dashed },
    pdps forward/.style       = { color = cb2, line width = 1pt },
    pdps forward accel/.style = { color = cb2, line width = 0.75pt, dashed  },
    ssn/.style                = { color = cb2, line width = 1pt },
    pdes/.style               = { color = cb3, line width = 1pt },
}
\pgfplotscreateplotcyclelist{primaldual}{
    {cb1},
    {cb2},
    {cb3},
    {cb4},
    {cb5},
}

\newcommand{\findminmax}[5]{
    \def\tempa{#3}
    \def\tempb{true}
    \def\firstelem{1}
    \ifx\tempa\tempb
    \pgfplotstablegetelem{\firstelem}{#2}\of{#1}%
    \pgfmathparse{\pgfplotsretval}
    \global\let#5\pgfmathresult
    \pgfplotstablegetelem{\firstelem}{#2}\of{#1}%
    \pgfmathparse{\pgfplotsretval}
    \global\let#4\pgfmathresult
    \fi
    \pgfplotstablegetrowsof#1
    \pgfmathparse{\pgfplotsretval-1}
    \foreach \i in {\firstelem,...,\pgfmathresult}{
        \pgfplotstablegetelem{\i}{#2}\of{#1}%
        \pgfmathparse{max(#5, \pgfplotsretval)}
        \global\let#5\pgfmathresult
        \pgfplotstablegetelem{\i}{#2}\of{#1}%
        \pgfmathparse{min(#4, \pgfplotsretval)}
        \global\let#4\pgfmathresult
    }
}

\newcommand{\SetMinMax}[3]{%
    \pgfplotstableread{#1}{#3}
    \findminmax{#3}{#2}{true}{\ExperimentMinValue}{\ExperimentMaxValue}
}

\newcommand{\UpdMinMax}[3]{
    \pgfplotstableread{#1}{#3}
    \findminmax{#3}{#2}{false}{\ExperimentMinValue}{\ExperimentMaxValue}
}

\pgfplotsset{
    fixedylog/.code 2 args={
        \pgfmathsetmacro{\base}{10^(-#1)}
        \pgfmathsetmacro{\logscale}{#1}
        \gdef\scaledlog##1{%
            (log10(\base + ((##1-\ExperimentMinValue)/(\ExperimentMaxValue-\ExperimentMinValue))) + \logscale)%
        }
        \pgfmathsetmacro{\logmax}{\scaledlog{\ExperimentMaxValue}}
        \pgfplotsset{
            y filter/.expression={\scaledlog{y}},
            ytick={0, 0.25*\logmax, 0.50*\logmax, 0.75*\logmax, \logmax},
            /pgf/number format/.cd,#2,
            scaled y ticks = manual:{$\pgfmathprintnumber{\ExperimentMinValue}\,+$}%
            {\pgfmathparse{10^((##1) - \logscale) - \base) * (\ExperimentMaxValue - \ExperimentMinValue)}}
        }
    }
}

\pgfplotsset{ignore legend/.style={every axis legend/.code={\let\addlegendentry\relax}}}

\usepackage{booktabs}

\usepackage[most]{tcolorbox}
\tcbset{%
    colback=dual!7, colframe=dual!7, sharp corners=all,
    top=0.3em, bottom=0.3em, left=0.3em, right=0.3em,
}
\newtcolorbox{algeqbox}{%
    ams nodisplayskip,
    left=-0.3em, right=-0.3em,
}
\newtcolorbox{algeqbox*}{}
\newtcolorbox{algenumbox}{}
\tcolorboxenvironment{example}{%
    colback=primal!5, colframe=primal!5,
    breakable, pad after break=0.3em
}
\usepackage{tablefootnote}
\makeatletter
\AfterEndEnvironment{tcolorbox}{%
    \tfn@tablefootnoteprintout%
    \gdef\tfn@fnt{0}%
}
\makeatother

\title{Introduction to Nonsmooth Analysis and Optimization}
\author{%
    \begin{minipage}[t]{0.4\linewidth}
        \centering
        Christian Clason\\[1ex]
        \small\email{c.clason@uni-graz.at}\\
        \url{https://imsc.uni-graz.at/clason}\\
        \orcid{0000-0002-9948-8426}
    \end{minipage}
    \and
    \begin{minipage}[t]{0.4\linewidth}
        \centering
        Tuomo Valkonen\\[1ex]
        \small\email{tuomov@iki.fi}\\
        \url{https://tuomov.iki.fi}\\
        \orcid{0000-0001-6683-3572}
    \end{minipage}
    \vspace*{1cm}
}
\date{%
    \today
    \\
    {\small\sffamily\textsc{arxiv}:\,\href{https://arxiv.org/abs/2001.00216}{\nolinkurl{2001.00216v7}}}
}
\publishers{%
    \centering
    \begin{tikzpicture}[domain=-1:1, scale=8, opacity=0.8, overlay, yshift=-10pt]
        \draw[color=structure!20,fill=structure!20,decoration={zigzag}] (0,0.5) -- (0.9,1.4) decorate {-- (0.9,-0.4)};
        \draw[color=structure!20,fill=structure!20,decoration={zigzag}] (0,0.5) -- (-0.9,1.4) decorate {-- (-0.9,-0.4)};
        \node[color=structure,opacity=0.3,scale=5] at (0.75,0.7) {$\partial f(\bar u)$};
        \draw[color=structure,thick] plot (\x,0.5);
        \node[color=structure,opacity=0.4,scale=2.5] at (-0.66,0.4) {\footnotesize $f(\bar u)+\langle 0,u\rangle$ };
        \draw[thick] plot (\x,{abs(\x)+0.5});
        \node[opacity=0.3,scale=4] at (-0.75,1.2) {$f(u)$};
        \draw[color=structure,thin,dashed] plot (\x,\x+0.5);
        \draw[color=structure,thin,dashed] plot (\x,-\x+0.5);
        \node[left,color=structure,opacity=0.3,scale=2] at (1,-0.25) {\footnotesize $f(\bar u) + \langle f'(\bar u;-1),u\rangle$};
        \draw[fill] (0,0.5) circle(0.33pt) node[opacity=0.3,below, yshift=-5pt]{$\bar u$};
    \end{tikzpicture}
}

\hypersetup{
    pdftitle={Nonsmooth Analysis and Optimization},
    pdfauthor={Christian Clason and Tuomo Valkonen},
    pdfkeywords={Convex analysis, proximal point methods, Clarke subdifferential, semismooth Newton methods, set-valued analysis, metric regularity, inverse problems, optimal control}
}

\begin{document}

\maketitle

\frontmatter

\tableofcontents

\chapter*{Preface}
\markboth{Preface}{Preface}
\addcontentsline{toc}{chapter}{Preface}

One of the major applications of classical analysis is {optimization} or the search for minima (or maxima) of a given function; this search may be motivated by a function directly representing an outcome of which lower values are desirable (say, the total cost of an economic production plan) or by a minimizing property indirectly being an essential characterization of a point of interest (say, a physical state as the minimizer of an energy functional, or the solution of an inverse or imaging problem as the minimizer of a regularization functional).
In particular, analytical concepts are crucial in every stage of the treatment of optimization problems: continuity properties for showing existence of solutions (that the minimal value is actually attained at a feasible point), first derivatives for intrinsic characterizations of solutions (via {Fermat principles} or {optimality conditions}) and for the numerical solution via {steepest descent} or {gradient methods}, and second derivatives for the numerical solution via {Newton methods} and for deriving stability results, e.g., with respect to computational errors (via implicit function theorems).

However, there are many practically relevant functions that are \emph{not} differentiable, such as the absolute value or maximum function. The goal of nonsmooth analysis is therefore to find generalized derivative concepts that on the one hand allow the above sketched approach for such functions and on the other hand admit a sufficiently rich calculus to give \emph{explicit} derivatives for a sufficiently large class of functions.
In this book, we specifically aim at treating problems of the form
\begin{equation}\label{eq:intro:prob}
    \tag{P}
    \min_{x\in C} \frac1p\norm{S(x)-z}_Y^p + \frac{\alpha}{q}\norm{x}_X^q
\end{equation}
for a closed convex {constraint} or feasible set $C\subset X$, a (possibly nonlinear but differentiable) operator $S:X\to Y$, $\alpha\geq 0$ and $p,q\in[1,\infty)$ (in particular, $p=1$ and/or $q=1$). Such problems are ubiquitous in inverse problems, imaging, and optimal control of differential equations.
Hence, we consider optimization in {infinite-dimensional} function spaces; i.e., we are looking for functions as minimizers. The main benefit (beyond the frequently cleaner notation) is that the developed algorithms become {discretization independent}: they can be applied to any (reasonable) finite-dimensional approximation, and the details -- in particular, the fineness -- of the approximation do not influence the convergence behavior of the algorithms. Naturally, such benefits do not come without  additional mathematical challenges. Many results in finite-dimensional optimization exploit one or more of the following properties:
\begin{enumerate}
    \item all norms are equivalent;
    \item all bounded sequences contain a convergent subsequence;
    \item derivative-like objects live in the same space as the minimizer is sought in.
\end{enumerate}
These allow defining iterations like the classical gradient method in a straightforward way and -- if the iterates can be shown to stay bounded -- showing that the iteration converges at least up to a subsequence. Unfortunately, these useful properties in general no longer hold in infinite-dimensional spaces, which makes it necessary to carefully keep track of the different spaces, norms, and modes of convergence that are natural to the specific problem and to make sure that they match up. We thus need to combine the classical approaches from nonlinear optimization with tools from functional analysis and the calculus of variations.
Throughout, a special role will be played by integral functionals and superposition operators that act pointwise on functions, since these allow transferring the often more explicit finite-dimensional calculus to the infinite-dimensional setting.

\bigskip

Nonsmooth analysis and optimization in finite dimensions has a long history; we refer here only to the classical textbooks \cite{makela1992nonsmooth,lemarechal1993-1,lemarechal1993-2,Rockafellar:1998,Boyd:2004,Ruszczynski:2006a} as well as the recent \cite{bagirov2014nonsmooth,Beck:2017,Nesterov:2018,CuiPang:2021,RoysetWets:2022}. There also exists a large body of literature on specific nonsmooth optimization problems, in particular ones involving variational inequalities and equilibrium constraints; see, e.g., \cite{outrata1998nonsmooth,facchinei2003finite-1,facchinei2003finite-2}. In contrast, the infinite-dimensional setting is still being actively developed, with monographs and textbooks focusing on either theory \cite{Clarke:1990a,Zalinescu:2002,Mordukhovich:2006,Schirotzek:2007,Barbu:2012,Penot:2013,Clarke:2013,ioffe2017variational,mordukhovich2018variational,Dontchev:2021} or algorithms \cite{Kunisch:2008a,Ulbrich:2011}. Two exceptions are \cite{Bauschke} and \cite{Peypouquet:2015}, the former containing an impressively comprehensive and integrated treatment of convex analysis and proximal point methods in Hilbert spaces, and the latter giving an equally impressively concise introduction to these topics in normed vector spaces. As this book neared completion, \cite{BauschkeMoursi:2023} was published, which serves as a very gentle introduction to convex optimization and first-order methods in Hilbert spaces as treated in \cite{Bauschke}. On the other end of the spectrum, the recent monograph \cite{Mordukhovich:2024} presents the state of the art of second-order variational analysis and generalized Newton methods in infinite dimensions.
The aim of this book is thus to draw together results scattered throughout the literature in order to give a unified presentation of theory -- both convex and nonconvex -- and algorithms -- both first- and second-order -- in Banach spaces that is suitable for an advanced class on mathematical optimization and can serve as a gateway to the research literature (including the monographs cited above). In order to do this, we focus on optimization of nonsmooth functionals rather than nonsmooth constraints; in particular, we do not treat optimization with complementarity or equilibrium constraints, which still see significant active development in infinite dimensions.
We also restrict the treatment to the two classes of
\begin{enumerate}[label=\roman*)]
    \item convex functions and
    \item locally Lipschitz continuous functions,
\end{enumerate}
which together cover a wide spectrum of applications. In particular, the first class will lead us to generalized gradient methods, while the second class is the basis for generalized Newton methods. These methods are chosen since they have become increasingly popular in recent years and fit particularly well within the integrated approach of this book. On the other hand, this focus leads us to omit other, more classical, methods and in particular subgradient and bundle methods, which have very recently seen developments in Hilbert spaces. Here, too, we can only refer to the research literature as well as to the classical books cited above for finite-dimensional treatments.
Regarding generalized derivatives of set-valued mappings required for the mentioned stability results, we similarly do not aim for a (possibly fuzzy) general theory and instead restrict ourselves to situations where a regularity condition (one out of the veritable zoo of conditions found in the literature) holds that allows deriving exact results that still apply to problems of the form \eqref{eq:intro:prob}. The general theory can be found in, e.g., \cite{aubin1990sva,Rockafellar:1998,mordukhovich2018variational,Mordukhovich:2006,Mordukhovich:2024}.

\bigskip

The book is intended for students and researchers with a solid background in analysis and linear algebra and an interest in the mathematical foundations of nonsmooth optimization.
Since we deal with infinite-dimensional spaces, some knowledge of functional analysis is assumed, but the necessary background will be summarized in \cref{chap:functan}. Similarly, \cref{chap:variation} collects needed fundamental results from the calculus of variations, including the direct method for existence of minimizers and the related notion of lower semicontinuity as well as differential calculus in Banach spaces, where the results on pointwise superposition operators on Lebesgue spaces require elementary (Lebesgue) measure and integration theory. Basic familiarity with classical nonlinear optimization is helpful but not necessary.

In \cref{part:convex} we then start our study of {convex} optimization problems. After introducing convex functionals and their basic properties in \cref{chap:convex}, we define our first generalized derivative in \cref{chap:subdiff}: the {convex subdifferential}, which is no longer a single unique derivative but consists of a {set} of equally admissible subderivatives. Nevertheless, we obtain a useful corresponding Fermat principle as well as calculus rules. A particularly useful calculus rule in convex optimization is {Fenchel duality}, which assigns to any optimization problem a {dual problem} that can help treating the original {primal} problem; this is the content of \cref{chap:fenchel}. We change our viewpoint in \cref{chap:monotone} slightly to study the subdifferential as a set-valued {monotone operator}, which leads us to the corresponding {resolvent} or {proximal point mapping}, which will later become the basis of all algorithms. The following \cref{chap:smoothness} discusses the relation between convexity and smoothness of primal and dual problem and introduces the {Moreau--Yosida regularization}, which has better properties in both regards that can be used to accelerate the convergence of algorithms. We turn to these in \cref{chap:proximal}, where we start by deriving a number of popular first-order methods including {forward-backward splitting} and {primal-dual proximal splitting} (also known as the \emph{Chambolle--Pock method}). Their convergence under rather general assumptions is then shown in \cref{chap:convergence}. If additional convexity properties hold, we can even show convergence rates for the iterates using a general {testing approach}; this is carried out in \cref{chap:testing}. Otherwise we either have to restrict ourselves to more abstract criticality measures as in \cref{chap:gap} or modify the algorithms to include {over-relaxation} or {inertia} as in \cref{chap:meta}.
One philosophy we here wish to pass to the reader is that the development of optimization methods consists, firstly, in suitable {reformulation} of the problem; secondly, in the {preconditioning} of the raw optimality conditions; and, thirdly, in {testing} with appropriate operators whether this yields fast convergence.

We leave the convex world in \cref{part:nonconvex}. For locally Lipschitz continuous functions, we introduce the {Clarke subdifferential} in \cref{chap:clarke} and derive calculus rules. Not only is this useful for obtaining a Fermat principle for problems of the form \eqref{eq:intro:prob}, it is also the basis for defining a further generalized derivative that can be used in place of the Hessian in a generalized Newton method. This {Newton derivative} and the corresponding {semismooth Newton method} is studied in \cref{chap:semismooth}. We also derive and analyze a variant of the primal-dual proximal splitting method suitable for \eqref{eq:intro:prob} in \cref{chap:nlpdps}. We end this part in \cref{chap:limiting,chap:epsilon} with a short outlook to further subdifferential concepts that can lead to sharper optimality conditions but in general admit a weaker calculus; we will treat some of these in detail in the next part.

To derive stability properties of minimization problems, we need to study the sensitivity of subdifferentials to perturbations and hence generalized derivative concepts for set-valued mappings; this is the goal of \cref{part:setvalued}. The construction of the generalized derivatives is geometric, based on {tangent} and {normal cones} introduced in \cref{chap:cones}. From these, we obtain {Fréchet} and {limiting (co)derivatives} in \cref{chap:graphical} and derive calculus rules for them in \crefrange{chap:gderiv}{chap:colimiting}.
In particular, we show how to lift the (more extensive) finite-dimensional theory to the special case of pointwise-defined sets and mappings on Lebesgue spaces in \cref{chap:pointcones,chap:superposition}.
We then address second-order conditions for nonsmooth nonconvex optimization problems in \cref{chap:secondorder}.
In \cref{chap:regularity}, we use these derivatives to characterize Lipschitz-like properties of set-valued mappings, which then are used to obtain the desired stability properties in \cref{chap:stability}. We also show in \cref{chap:fasterconvergence} that these regularity properties imply faster convergence of first-order methods.

Finally, \cref{part:applications} illustrates how these results apply to concrete optimization problems arising in inverse problems and mathematical imaging (\cref{chap:sparse,chap:l1fit,chap:tv}) and in optimal control (\cref{chap:control,chap:discretecontrol}), where we freely admit that the selection of examples is subjective and driven by the authors' research interests.
These chapters are accompanied by Julia implementations \cite{nonsmoothbook-codes} of the discussed algorithms, which can be used to recreate the presented numerical results.

\bigskip

This book can serve as a textbook for several different classes:
\begin{enumerate}
    \item an introductory course on convex optimization based on \crefrange{chap:convex}{chap:testing} (excluding \cref{sec:convex:continuity} and results on superposition operators) and adding \cref{chap:gap,chap:meta,chap:nlpdps} as time permits;
    \item an intermediate course on nonsmooth optimization based on \crefrange{chap:convex}{chap:convergence}  (including \cref{sec:convex:continuity} and results on superposition operators) together with \cref{chap:clarke,chap:semismooth,chap:limiting};
    \item an intermediate course on nonsmooth analysis based on \crefrange{chap:convex}{chap:monotone} together with \cref{chap:clarke} and \crefrange{chap:limiting}{chap:superposition}, adding \crefrange{chap:gderiv}{chap:secondorder} as time permits;
    \item an advanced course on set-valued analysis based on \crefrange{chap:limiting}{chap:fasterconvergence}.
\end{enumerate}

\bigskip

This book is based in part on such graduate lectures given by the first author in 2014 (in slightly different form) and 2016--2017 at the University of Duisburg-Essen and by the second author at the University of Cambridge in 2015 and at the MODEMAT Research Center in Mathematical Modeling and Optimization in Quito in 2020. Shorter seminars were also delivered at the University of Jyväskylä and MODEMAT in 2017.
\Cref{part:setvalued} of the book was also used in a course on variational analysis in 2019.
Parts of the book were also taught by both authors at the Winter School \enquote{Modern Methods in Nonsmooth Optimization} organized by Christian Kanzow and Daniel Wachsmuth at the University of Würzburg in February 2018, for which the notes were further adapted and extended.
As such, much (but not all) of this material is classical. In particular, \crefrange{chap:convex}{chap:smoothness} as well as \cref{chap:clarke} are based on \cite{Barbu:2012,Brokate,Schirotzek:2007,Attouch,Bauschke,Clarke:2013}, \cref{chap:semismooth} is based on \cite{Ulbrich:2002a,Kunisch:2008a,Schiela:2008a}, \cref{chap:limiting} is extracted from \cite{Mordukhovich:2006}, and \crefrange{chap:cones}{chap:colimiting} are adapted from \cite{Rockafellar:1998,Mordukhovich:2006}.
Parts of \cref{chap:epsilon} are adapted from \cite{ioffe2017variational}, and other parts are original work.
On the other hand, \crefrange{chap:proximal}{chap:meta} as well as \cref{chap:nlpdps,chap:superposition,chap:fasterconvergence} are adapted from \cite{tuomov-proxtest,tuomov-subreg,tuomov-nlpdhgm-redo}, \cite{tuomov-pdex2stability}, and \cite{tuomov-subreg}, respectively.

Finally, we would like to express our gratitude to Sebastian Angerhausen, Ronny Bergmann, Alberto Domínguez Corella, Edison Felipe Guerra Urgiles, Andreas Habring, Jyrki Jauhiainen, Bjørn Jensen, Fernando Jimenez Torres, Heikki von Koch, Anton Schiela, Ensio Suonperä, Stefan Ulbrich, Diego Vargas Jaramillo, Daniel Wachsmuth, Gerd Wachsmuth, and an anonymous reviewer for carefully reading parts of the manuscript, finding mistakes and arguments that could be expressed more clearly, or making other helpful suggestions. We would also like to thank Julia Cochrane, Elizabeth Greenspan, Kiara Hicks, Cheryl Hufnagle, Rose Kolassiba, David Riegelhaupt, and Doug Smock at SIAM for their expert work on the manuscript and for making the publication process a pleasure. All remaining errors are of course our own.

\enlargethispage{1cm}
\bigskip

\begin{flushright}\noindent
    \itshape{Essen/Graz and Quito/Helsinki, October 2025}
\end{flushright}

\mainmatter

\part{Background}\label{part:background}

\chapter{Functional analysis}\label{chap:functan}

Functional analysis is the study of infinite-dimensional vector spaces and of the operators acting between them, and has since its foundations in the beginning of the 20th century grown into the \emph{lingua franca} of modern applied mathematics.
In this chapter we collect the basic concepts and results (and, more importantly, fix notations) from linear functional analysis that will be used throughout the rest of the book. For details and proofs, the reader is referred to the standard literature, e.g., \cite{Alt:2016,Brezis:2010a,rynne2008functional}, or to \cite{Clason}.

\section{Normed vector spaces}\label{sec:functan:normed}

In the following, $X$ will denote a real vector space.
A mapping $\norm{\cdot}:X\to \R^+\defeq[0,\infty)$ is called a \term{norm} (on $X$), if for all $x\in X$ there holds
\begin{enumerate}
    \item $\norm{\lambda x} = |\lambda| \norm{x}$ for all $\lambda\in\R$;
    \item $\norm{x+y} \leq \norm{x} + \norm{y}$ for all $y\in X$;
    \item $\norm{x} = 0$ if and only if $x = 0\in X$.
\end{enumerate}
\begin{example}\label{ex:functan:norm}
    \begin{enumerate}
        \item\label{ex:functan:norm:i}
            The following mappings define norms on $X = \R^N$:
            \begin{equation*}
                \begin{aligned}
                    \norm{x}_p &= \left(\sum_{i=1}^N |x_i|^p\right)^{1/p},\qquad1\leq p<\infty,\\
                    \norm{x}_\infty &= \max_{i=1,\dots,N} |x_i|.
                \end{aligned}
            \end{equation*}
        \item \label{ex:functan:norm:ii}
            The following mappings define norms on $X = \ell^p$ (the space of real-valued sequences for which these terms are finite):
            \begin{equation*}
                \begin{aligned}
                    \norm{x}_p &= \left(\sum_{i=1}^\infty |x_i|^p\right)^{1/p},\qquad 1\leq p<\infty,\\
                    \norm{x}_\infty &= \sup_{i=1,\dots,\infty} |x_i|.
                \end{aligned}
            \end{equation*}
        \item \label{ex:functan:norm:iii}
            The following mappings define norms on $X = L^p(\Omega)$ (the space of real-valued measurable functions on the domain $\Omega\subset \R^d$ for which these terms are finite):
            \begin{equation*}
                \begin{aligned}
                    \norm{u}_{L^p} &= \left(\int_\Omega |u(x)|^p\right)^{1/p},\qquad 1\leq p<\infty,\\
                    \norm{u}_{L^\infty} &= \mathop\mathrm{ess\,\sup}_{x\in \Omega} |u(x)|,
                \end{aligned}
            \end{equation*}
            where $\mathrm{ess\,\sup}$ stands for the essential supremum; for details on these definitions, see, e.g., \cite{Alt:2016}.
        \item \label{ex:functan:norm:iv}
            The following mapping defines a norm on $X = C(\overline\Omega)$ (the space of continuous functions on $\overline\Omega$):
            \begin{equation*}
                \norm{u}_C = \sup_{x\in \overline\Omega} |u(x)|.
            \end{equation*}
            An analogous norm is defined on $X=C_0(\Omega)$ (the space of continuous functions on $\Omega$ with compact support), if the supremum is taken only over $x\in\Omega$.
    \end{enumerate}
\end{example}
If $\norm{\cdot}$ is a norm on $X$, the tuple $(X,\norm{\cdot})$ is called a \term[space!normed]{normed vector space}, and one frequently denotes this by writing $\norm{\cdot}_X$. If the norm is canonical (as in \cref{ex:functan:norm}\,\ref{ex:functan:norm:ii}--\ref{ex:functan:norm:iv}), it is often omitted, and one speaks simply of ``the normed vector space $X$''.

Two norms $\norm{\cdot}_1$, $\norm{\cdot}_2$ are called \term[norm!equivalent]{equivalent} on $X$, if there are constants $c_1,c_2 >0$ such that
\begin{equation*}
    c_1 \norm{x}_2 \leq \norm{x}_1 \leq c_2 \norm{x}_2 \qquad\text{for all } x\in X.
\end{equation*}
If $X$ is finite-dimensional, all norms on $X$ are equivalent. However, the corresponding constants $c_1$ and $c_2$ may depend on the dimension $N$ of $X$; avoiding such dimension-dependent constants is one of the main reasons to consider optimization in infinite-dimensional spaces.

If $(X,\norm{\cdot}_X)$ and $(Y,\norm{\cdot}_Y)$ are normed vector spaces with $X\subset Y$, we call $X$ \term[space!continuously embedded]{continuously embedded} in $Y$, denoted by $X\hookrightarrow Y$, if there exists a $C>0$ with
\begin{equation*}
    \norm{x}_Y \leq C\norm{x}_X \qquad\text{for all } x\in X.
\end{equation*}
For example, if $\Omega\subset \R^d$ is a bounded domain, $L^q(\Omega)\hookrightarrow L^p(\Omega)$ for every $1\leq p\leq q\leq\infty$.

\bigskip

A norm directly induces a notion of convergence, the so-called \term[convergence!strong]{strong convergence}.
A sequence $\{x_n\}_{n\in\N}\subset X$ \emph{converges} (\emph{strongly} in $X$) to a $x\in X$, denoted by $x_n\to x$, if
\begin{equation*}
    \lim_{n\to\infty} \norm{x_n-x}_X = 0.
\end{equation*}

\clearpage

A set $U\subset X$ is called
\begin{itemize}
    \item \term[set!closed]{closed}, if for every convergent sequence $\{x_n\}_{n\in\N}\subset U$ the limit $x\in X$ is an element of $U$ as well;
    \item \term[set!compact]{compact}, if every sequence $\{x_n\}_{n\in\N}\subset U$ contains a convergent subsequence $\{x_{n_k}\}_{k\in\N}$ with limit $x\in U$.
\end{itemize}
A mapping $F:X\to Y$ is \term[mapping!continuous]{continuous} if and only if $x_n\to x$ implies $F(x_n)\to F(x)$. If $x_n\to x$ and $F(x_n)\to y$ imply that $F(x) = y$ (i.e., $\graph F\subset X\times Y$ is a closed set), we say that $F$ has \term[graph!closed]{closed graph}.

Further we define for later use for $x\in X$ and $r>0$
\begin{itemize}
    \item the \term[ball!open]{open ball} $\OB(x,r) \defeq \setof{z\in X}{\norm{x-z}_X< r}$ and
    \item the \term[ball!closed]{closed ball} $\B(x,r) \defeq \setof{z\in X}{\norm{x-z}_X\leq r}$.
\end{itemize}
The closed ball around $0\in X$ with radius $1$ is also referred to as the \term[ball!unit]{unit ball} $\B_X$.
A set $U\subset X$ is called
\begin{itemize}
    \item \term[set!open]{open}, if for all $x\in U$ there exists an $r>0$ with $\OB(x,r)\subset U$ (i.e., all $x\in U$ are \term[point!interior]{interior points} of $U$);
    \item \term[set!bounded]{bounded}, if it is contained in $\B(0,r)$ for some $r>0$;
    \item \term[set!convex]{convex}, if for any $x,y\in U$ and $\lambda\in[0,1]$ also $\lambda x + (1-\lambda)y\in U$.
\end{itemize}
In normed vector spaces it always holds that the complement of an open set is closed and vice versa (i.e., the closed sets in the sense of topology are exactly the (sequentially) closed sets as defined above). The definition of a norm directly implies that both open and closed balls are convex.

For arbitrary $U$, we denote by $\closure U$ the \term{closure} of $U$, defined as the smallest closed set that contains $U$ (which coincides with the set of all limit points of convergent sequences in $U$); we write $\interior U$ for the \term{interior} of $U$, which is the largest open set contained in $U$; and we write $\BD U\defeq \closure U \setminus \interior U$ for the \term{boundary} of $U$.
Finally, we write $\conv U$ for the \term[hull, convex]{convex hull} of $U$, defined as the smallest convex set that contains $U$.

A normed vector space $X$ is called \term[space]{complete} if every Cauchy sequence in $X$ is convergent; in this case, $X$ is called a \term[space!Banach]{Banach space}.
All spaces in \cref{ex:functan:norm} are Banach spaces.
Convex subsets of Banach spaces have the following useful property which derives from the Baire theorem.
\begin{lemma}[core--int\protect\footnote{e.g., \cite[Lemma 5.2]{Clason}}]\index{lemma!core--int}
    \label{lem:functan:coreint}
    Let $X$ be a Banach space and $U\subset X$ be closed and convex. Then
    \begin{equation*}
        \interior U = \setof{x\in U}{\text{for all } h\in X \text{ there is a } \delta>0      \text{ with } x+ t h \in U \text{ for all } t\in [0,\delta]}.
    \end{equation*}
\end{lemma}
The set on the right-hand side is called \term[interior!algebraic]{algebraic interior} or \term[core|see{interior, algebraic}]{core}, which explains the name of the lemma. Note that the inclusion \enquote{$\subset$} always holds in normed vector spaces due to the definition of interior points via open balls.

\bigskip

We now consider mappings between normed vector spaces. In the following, let $(X,\norm{\cdot}_X)$ and $(Y,\norm{\cdot}_Y)$ be normed vector spaces, $U\subset X$, and $F: U\to Y$ be a mapping. We denote by
\begin{itemize}
    \item $\kernel F  \defeq \setof{x\in U}{F(x) = 0}$ the \term{kernel} or \term[space!null]{null space} of $F$;
    \item $\range F  \defeq \setof{F(x)\in Y}{x\in U}$ the \term[range!of an operator]{range} of $F$;
    \item $\graph F  \defeq \setof{(x,y)\in X\times Y}{y=F(x)}$ the \term[graph!of an operator]{graph} of $F$.
\end{itemize}
We call $F:U\to Y$
\begin{itemize}
    \item \term[mapping!continuous]{continuous} at $x\in U$, if for all $\eps>0$ there exists a $\delta >0$ with
        \begin{equation*}
            \norm{F(x)-F(z)}_Y \leq \eps\qquad \text{for all } z\in \OB(x,\delta) \cap U;
        \end{equation*}
    \item \term[mapping!continuous!Lipschitz]{Lipschitz continuous}, if there exists an $L>0$ (called \term[constant, Lipschitz]{Lipschitz constant} or \term[factor, Lipschitz|see{constant, Lipschitz}]{Lipschitz factor}) with
        \begin{equation*}
            \norm{F(x_1)-F(x_2)}_Y \leq L \norm{x_1-x_2}_X \qquad\text{for all } x_1,x_2 \in U.
        \end{equation*}
    \item \term[mapping!continuous!locally Lipschitz]{locally Lipschitz continuous at} $x\in U$, if there exists a $\delta>0$ and a $L=L(x,\delta)>0$ with
        \begin{equation*}
            \norm{F(x)-F(\tilde x)}_Y \leq L \norm{x-\tilde x}_X \qquad\text{for all } \tilde x\in \OB(x,\delta) \cap U;
        \end{equation*}
    \item \emph{locally Lipschitz continuous near} $x\in U$, if there exists a $\delta>0$ and a $L=L(x,\delta)>0$ with
        \begin{equation*}
            \norm{F(x_1)-F(x_2)}_Y \leq L \norm{x_1-x_2}_X \qquad\text{for all } x_1,x_2 \in \OB(x,\delta)\cap U.
        \end{equation*}
        We will refer to the $\OB(x,\delta)$ as the \term[neighborhood, Lipschitz]{Lipschitz neighborhood} of $x$ (for $F$).
        If $F$ is locally Lipschitz continuous near every $x\in U$, we call $F$ \emph{locally Lipschitz continuous on} $U$.
\end{itemize}

If $T:X\to Y$ is linear, continuity is equivalent to the existence of a constant $C>0$ with
\begin{equation*}
    \norm{Tx}_Y \leq C\norm{x}_X \qquad\text{for all }x\in X.
\end{equation*}
For this reason, continuous linear mappings are called \term[operator!bounded linear]{bounded}; one speaks of a bounded linear \emph{operator}.
The space $\linear(X;Y)$ of bounded linear operators is itself a normed vector space if endowed with the \term[norm!operator]{operator norm}
\begin{equation*}
    \norm{T}_{\linear(X;Y)}     = \sup_{x\in X\setminus\{0\}}\frac{\norm{Tx}_Y}{\norm{x}_X}
    = \sup_{\norm{x}_X= 1} \norm{Tx}_Y =   \sup_{\norm{x}_X\leq 1} \norm{Tx}_Y,
\end{equation*}
which is equal to the smallest possible constant $C$ in the definition of boundedness. (In the pathological case $X=\{0\}$ -- which we exclude from now on as it is irrelevant for optimization -- only the last definition is applicable.)
If $(Y,\norm{\cdot}_Y)$ is a Banach space, then so is $(\linear(X;Y),\norm{\cdot}_{\linear(X;Y)})$.

Finally, if $T\in \linear(X;Y)$ is bijective, the inverse $T^{-1}:Y\to X$ is continuous if and only if there exists a $c>0$ with
\begin{equation*}
    c\norm{x}_X \leq \norm{Tx}_Y \qquad\text{for all }x\in X.
\end{equation*}
In this case, $\norm{T^{-1}}_{\linear(Y;X)} = c^{-1}$ for the largest possible choice of $c$.

\section{Dual spaces, separation, and weak convergence}\label{sec:functan:dual}

Of particular importance to us is the special case $\linear(X;Y)$ for $Y=\R$, the space of \term[functional!bounded linear]{bounded linear functionals} on $X$. In this case, $X^*\defeq\linear(X;\R)$ is called the \term[space!dual]{dual space} (or just \emph{dual}) of $X$. For $x^*\in X^*$ and $x\in X$, we set
\begin{equation*}
    \dual{x^*,x}_X \defeq x^*(x) \in\R.
\end{equation*}
This \term[pairing, duality]{duality pairing} indicates that we can also interpret it as $x$ acting on $x^*$, which will become important later. The definition of the operator norm immediately implies that
\begin{equation}\label{eq:functan:cs_banach}
    \dual{x^*,x}_X \leq \norm{x^*}_{X^*}\norm{x}_X\qquad\text{for all }x\in X,x^*\in X^*.
\end{equation}

In many cases, the dual of a Banach space can be identified with another known Banach space.
\begin{example}\label{ex:functan:dual}
    \begin{enumerate}
        \item \label{ex:functan:dual:i}
            $(\R^N,\norm{\cdot}_p)^* \cong (\R^N,\norm{\cdot}_q)$ with $p^{-1}+q^{-1} = 1$, where we set $0^{-1}=\infty$ and $\infty^{-1} = 0$. The duality pairing is given by
            \begin{equation*}
                \dual{x^*,x}_{p} = \sum_{i=1}^N x^*_i x_i.
            \end{equation*}
        \item \label{ex:functan:dual:ii}
            $(\ell^p)^* \cong (\ell^q)$ for $1< p < \infty$. The duality pairing is given by
            \begin{equation*}
                \dual{x^*,x}_{p} = \sum_{i=1}^\infty x^*_i x_i.
            \end{equation*}
            Furthermore, $(\ell^1)^*=\ell^\infty$, but $(\ell^\infty)^*$ is not a sequence space.
        \item \label{ex:functan:dual:iii}
            Analogously, $L^p(\Omega)^* \cong L^q(\Omega)$ with $p^{-1}+q^{-1}=1$ for $1< p < \infty$. The duality pairing is given by
            \begin{equation*}
                \dual{u^*,u}_{p} = \int_\Omega u^*(x)u(x)\,dx.
            \end{equation*}
            Furthermore, $L^1(\Omega)^*\cong L^\infty(\Omega)$, but $L^\infty(\Omega)^*$ is not a function space.
        \item \label{ex:functan:dual:iv}
            $C_0(\Omega)^*\cong \mathcal{M}(\Omega)$, the space of \term[measure, Radon]{Radon measures}; it contains among others the Lebesgue measure as well as Dirac measures $\delta_x$ for $x\in\Omega$, defined via $\delta_x(u) = u(x)$ for $u\in C_0(\Omega)$. The duality pairing is given by
            \begin{equation*}
                \dual{u^*,u}_{C} = \int_\Omega u(x)\,du^*.
            \end{equation*}
    \end{enumerate}
\end{example}

A central result on dual spaces is the \emph{Hahn--Banach theorem}, which comes in a variety of analytic and  geometric forms; we will specifically need the following.
\begin{theorem}[Hahn--Banach, analytic\protect\footnote{e.g., \cite[Theorem 8.3]{Clason}}]\label{thm:functan:hb-extension}\index{theorem!Hahn--Banach!analytic}
    Let $X$ be a normed vector space and $x\in X\setminus\{0\}$. Then there exists a \term[functional!norming]{norming functional} $x^*\in X^*$ with
    \begin{equation*}
        \norm{x^*}_{X^*} = 1 \qquad\text{and}\qquad \dual{x^*,x}_X = \norm{x}_X.
    \end{equation*}
\end{theorem}
\begin{theorem}[Hahn--Banach, geometric\protect\footnote{e.g., \cite[Theorems 1.6, 1.7]{Brezis:2010a}}]\label{thm:functan:hb-separation}\index{theorem!Hahn--Banach!geometric}
    Let $X$ be a normed vector space and $A,B\subset X$ be convex, nonempty, and disjoint.
    \begin{enumerate}
        \item \label{item:functan:hb-separation:i}
            If $A$ is open, there exists an $x^*\in X^*$ and a $\lambda\in \R$ with
            \begin{equation*}
                \dual{x^*,x_1}_X < \lambda \leq \dual{x^*,x_2}_X \qquad\text{for all }x_1\in A, x_2\in B.
            \end{equation*}
        \item \label{item:functan:hb-separation:ii}
            If $A$ is closed and $B$ is compact, there exists an $x^*\in X^*$ and a $\lambda\in \R$ with
            \begin{equation*}
                \dual{x^*,x_1}_X \leq \lambda < \dual{x^*,x_2}_X \qquad\text{for all }x_1\in A, x_2\in B.
            \end{equation*}
    \end{enumerate}
\end{theorem}
Particularly the geometric version -- also referred to as \term[theorem!separation]{separation theorems} -- is of crucial importance in convex analysis. We will also require their following variant, which is known as \term[theorem!Eidelheit]{Eidelheit's theorem}.
\begin{corollary}\label{lem:convex:eidelheit}
    Let $X$ be a normed vector space and $A,B\subset X$ be convex and nonempty. If the interior $\interior A$ of $A$ is nonempty and disjoint with $B$, there exists an $x^*\in X^*\setminus\{0\}$ and a $\lambda\in \R$ with
    \begin{equation*}
        \dual{x^*,x_1}_X \leq \lambda \leq \dual{x^*,x_2}_X \qquad\text{for all }x_1\in A, x_2\in B.
    \end{equation*}
\end{corollary}
\begin{proof}
    \Cref{thm:functan:hb-separation}\,\ref{item:functan:hb-separation:i} yields the existence of $x^*$ and $\lambda$ satisfying the claim for all $x_1 \in \interior A$; this inequality is even strict, which also implies $x^*\neq 0$.
    It thus remains to show that the first inequality also holds for the remaining $x_1 \in A \setminus \interior A$.
    Since $\interior A$ is nonempty, there exists an $x_0\in \interior A$, i.e., there is an $r>0$ with $\OB(x_0,r)\subset A$. The convexity of $A$ then implies that $t\tilde x + (1-t)x_1\in A$ for all $\tilde x \in \OB(x_0,r)$ and $t\in [0,1]$. Hence,
    \begin{equation*}
        t\OB(x_0,r) + (1-t) x = \OB(tx_0 + (1-t) x_1, tr) \subset A,
    \end{equation*}
    and in particular $z(t) \defeq t x_0 + (1-t) x_1 \in \interior A$ for all $t\in (0,1)$.

    We can thus find a sequence $\{z_n\}_{n\in\N}\subset \interior A$ (e.g., $z_n = z(n^{-1})$) with $z_n\to x_1$. Due to the continuity of $x^*\in X^*=\linear(X;\R)$ we can thus pass to the limit $n\to\infty$ and obtain
    \begin{equation*}
        \dual{x^*,x_1}_X = \lim_{n\to \infty} \dual{x^*,z_n}_X \leq \lambda.
        \qedhere
    \end{equation*}
\end{proof}

This can be used to characterize a normed vector space by its dual.
For example, a direct consequence of \cref{thm:functan:hb-extension} is that the norm on a Banach space can be expressed as an operator norm.
\begin{corollary}\label{cor:functan:norm_dual}
    Let $X$ be a Banach space. Then for all $x\in X$,
    \begin{equation*}
        \norm{x}_X = \sup_{\norm{x^*}_{X^*}\leq 1} \dual{x^*,x}_X,
    \end{equation*}
    and the supremum is attained.
\end{corollary}
\begin{proof}
    This follows from
    \begin{equation*}
        \dual{x^*,x}_X \leq \norm{x^*}_{X^*}\norm{x}_X \leq \norm{x}_X \qquad\text{for all }x^*\in B_{X^*}, \ x\in X,
    \end{equation*}
    with equality for the norming functional from \cref{thm:functan:hb-extension}.
\end{proof}
A vector $x\in X$ can therefore be considered as a linear and, by \eqref{eq:functan:cs_banach}, bounded functional on $X^*$, i.e., as an element of the \term[space!bidual]{bidual} $X^{**}\defeq(X^*)^*$.
The embedding $X\hookrightarrow X^{**}$ is realized by the \term[injection, canonical]{canonical injection}
\begin{equation}\label{eq:functan:canonical-injection}
    J:X\to X^{**},\qquad \dual{Jx,x^*}_{X^{*}} \defeq \dual{x^*,x}_X\quad\text{for all }x^*\in X^*.
\end{equation}
Clearly, $J$ is linear; \cref{cor:functan:norm_dual} furthermore implies that $\norm{Jx}_{X^{**}} = \norm{x}_X$.
If the canonical injection is surjective and we can thus identify $X^{**}$ with $X$, the space $X$ is called \term[space!reflexive]{reflexive}.
All finite-dimensional spaces are reflexive, as are \cref{ex:functan:norm}\,\ref{ex:functan:norm:ii} and \ref{ex:functan:norm:iii} for $1<p<\infty$; however, $\ell^1,\ell^\infty$ as well as $L^1(\Omega),L^\infty(\Omega)$ and $C(\overline\Omega)$ are not reflexive.
In general, a normed vector space is reflexive if and only if its dual space is reflexive.

The following consequence of the separation \cref{thm:functan:hb-separation} will be of crucial importance in \cref{part:setvalued}. For a set $A\subset X$, we define the \term[cone!polar]{polar cone}
\begin{align*}
    \polar A &\defeq \setof{ x^* \in X^*}{\dual{x^*,x}_X \le 0 \text{ for all } x \in A},
    \intertext{cf.~\cref{fig:polar}. Similarly, we define for $B\subset X^*$ the \term[cone!prepolar]{prepolar cone}}
     B_\circ &\defeq \setof{ x \in X}{\dual{x^*,x}_X \le 0 \text{ for all } x^* \in B}.
\end{align*}

\begin{figure}
    \centering
    \begin{asy}
        pair p1=(.6, 1.2);
        pair p2=(1.2, .1);
        pair O=(0, 0);
        pair n1=orthog(p1, -1);
        pair n2=orthog(p2, 1);

        path p=(.6,.5)..p1..(1.1,1.2)..p2..controls(.2, .1)..cycle;
        draw(p);
        fill(p, lightfill);
        label("$A$", (p1+p2)/2);

        draw(O--1.5*p1);
        draw(O--1.7*p2);

        fill(O--1.05*n1--1.1*n2--cycle, darkfill);
        draw(O--.9*n1, primalline + linewidth(1.1), Arrow);
        draw(O--.9*n2, primalline + linewidth(1.1), Arrow);
        label("$A^\circ$", (-0.25,-0.2));

        dot(O);
        label("$0$", (0.16,0.12));

        draw(.16*n1--.16*n1+.16*p1--.16*p1, dashed);
        draw(.18*n2--.18*n2+.18*p2--.18*p2, dashed);
    \end{asy}
    \caption{The polar cone $\polar A$ is bounded by the half-lines at right angles to the smallest cone containing $A$.}
    \label{fig:polar}
\end{figure}
The \term[cone!bipolar]{bipolar cone} of $A\subset X$ is then defined as
\begin{equation*}
    \bipolar A \defeq (\polar A)_\circ \subset X.
\end{equation*}
(If $X$ is reflexive, $\bipolar A=\polar{(\polar A)}$.)
For the following statement about polar cones, recall that a set $C\subset X$ is called a \term{cone} if  $x\in C$ and $\lambda> 0$ implies that $\lambda x\in C$ (such that (pre-, bi-)polar cones are indeed cones).
\begin{theorem}[bipolar]\index{theorem!bipolar}
    \label{lemma:functan:polar-inclusion}
    Let $X$ be a normed vector space and $A \subset X$. Then
    \begin{enumerate}
        \item \label{item:functan:polar-inclusion:i}
            $\polar A$ is closed and convex;

        \item \label{item:functan:polar-inclusion:ii}
            $A \subset \bipolar A$;

        \item \label{item:functan:polar-inclusion:iii}
            if $A \subset B$, then $\polar B \subset \polar A$;

        \item \label{item:functan:polar-inclusion:iv}
            if $C$ is a nonempty, closed, and convex cone, then $C=\bipolar C$.
    \end{enumerate}
\end{theorem}
\begin{proof}
    \emph{\ref{item:functan:polar-inclusion:i}:} This follows directly from the definition and the continuity of the duality pairing.

    \emph{\ref{item:functan:polar-inclusion:ii}:}
    Let $x \in A$ be arbitrary. Then by definition of the polar cone, every $x^*\in\polar A$ satisfies
    \begin{equation*}
        \dual{x^*,x}_X \le 0,
    \end{equation*}
    i.e., $x\in (\polar A)_\circ = \bipolar A$.

    \emph{\ref{item:functan:polar-inclusion:iii}:} This is immediate from the definition.

    \emph{\ref{item:functan:polar-inclusion:iv}:}
    By \ref{item:functan:polar-inclusion:ii}, we only need to prove $\bipolar C \subset C$ which we do by contradiction. First, observe that since $C$ is a nonempty and closed cone, we always have $0\in C$.
    Assume therefore that there exists $x \in \bipolar C \setminus \{0\}$ with $x \not\in C$.
    Applying \cref{thm:functan:hb-separation}\,\ref{item:functan:hb-separation:ii} to the nonempty (due to \ref{item:functan:polar-inclusion:ii}) closed, and convex set $\bipolar C$ and the disjoint compact convex set $\{x\}$, we obtain $x^* \in X^* \setminus \{0\}$ and $\lambda \in \R$ such that
    \begin{equation}
        \label{eq:cones:polarity-hb}
        \dual{x^*,\tilde x}_X
        \le
        \lambda
        <
        \dual{x^*,x}_X
        \quad
        \text{for all }\tilde x\in C.
    \end{equation}
    Since $C$ is a cone, the first inequality must also hold for $t\tilde x\in C$ for every $t>0$. This implies that
    \begin{equation*}
        \dual{x^*,\tilde x}_X \le t^{-1} \lambda\to 0 \quad \text{for }t\to \infty,
    \end{equation*}
    i.e., $\dual{x^*,\tilde x}_X \le 0$ for all $\tilde x\in C$ must hold, which implies that $x^*\in \polar C$.
    On the other hand, if $\lambda<0$, we obtain by the same argument that
    \begin{equation*}
        \dual{x^*,\tilde x}_X \le t^{-1} \lambda\to -\infty \quad \text{for }t\to 0,
    \end{equation*}
    which cannot hold. Hence, we can take $\lambda=0$ in \eqref{eq:cones:polarity-hb}. Together, we obtain from $x\in \bipolar C$ the contradiction
    \begin{equation*}
        0 < \dual{x^*,x}_X \leq 0.
        \qedhere
    \end{equation*}
\end{proof}

\bigskip

The duality pairing induces further notions of convergence.
\begin{enumerate}
    \item A sequence $\{x_n\}_{n\in\N}\subset X$ \term[convergence!weak]{converges weakly} (in $X$) to $x\in X$, denoted by  $x_n\weakto x$, if
        \begin{equation*}
            \dual{x^*,x_n}_X\to  \dual{x^*,x}_X \qquad\text{for all }x^*\in X^*.
        \end{equation*}
    \item A sequence $\{x^*_n\}_{n\in\N}\subset X^*$ \term[convergence!weak-$*$]{converges weakly-$*$} (in $X^*$) to $x^*\in X^*$, denoted by $x^*_n\weaktostar x^*$, if
        \begin{equation*}
            \dual{x_n^*,x}_X\to  \dual{x^*,x}_X \qquad\text{for all }x\in X.
        \end{equation*}
\end{enumerate}
Weak convergence generalizes the concept of componentwise convergence in $\R^N$, which -- as can be seen from the proof of the Heine--Borel theorem -- is the appropriate concept in the context of compactness.
Strong convergence in $X$ implies weak convergence by continuity of the duality pairing; in the same way, strong convergence in $X^*$ implies weak-$*$ convergence. If $X$ is reflexive, weak and weak-$*$ convergence (both in $X=X^{**}$) coincide. In finite-dimensional spaces, all these convergence notions coincide.

Weakly convergent sequences are always bounded; if $X$ is a Banach space, so are weakly-$*$ convergent sequences.
If $x_n\to x$ and $x_n^*\weaktostar x^*$ or $x_n\weakto x$ and $x_n^*\to x^*$, then  $\dual{x_n^*,x_n}_X\to \dual{x^*,x}_X$. However, the duality pairing of weak(-$*$) convergent sequences does not converge in general.

As for strong convergence, one defines weak(-$*$) continuity and closedness of mappings as well as weak(-$*$) sequential closedness and compactness of sets. The last property is of fundamental importance in optimization; its characterization is therefore a central result of this chapter.
\begin{theorem}[Eberlein--\u{S}mulyan\protect\footnote{e.g., \cite[Theorem 11.8]{Clason}}]\label{thm:ebsmul}\index{theorem!Eberlein--\u{S}mulyan}
    If $X$ is a reflexive Banach space, then $\B_X$ is weak sequentially compact.
\end{theorem}
Hence in a reflexive space, all bounded sequences contain a \emph{weakly} (but in general not strongly) convergent subsequence. Note that weak closedness is a \emph{stronger} claim than closedness, since the property has to hold for more sequences. For convex sets, however, both concepts coincide.
\begin{lemma}\label{lem:convex_closed}
    Let $X$ be a normed vector space and $U\subset X$ be convex. Then $U$ is weak sequentially closed if and only if $U$ is closed.
\end{lemma}
\begin{proof}
    Weak sequentially closed sets are always closed since a convergent sequence is also weakly convergent.
    Let now $U\subset X$ be convex closed and nonempty (otherwise nothing has to be shown) and consider a sequence $\{x_n\}_{n\in\N}\subset U$ with $x_n\weakto x\in X$. Assume that $x\in X\setminus U$. Then the sets $U$ and $\{x\}$ satisfy the premise of \cref{thm:functan:hb-separation}\,\ref{item:functan:hb-separation:ii}; we thus find an $x^*\in X^*$ and a $\lambda\in\R$ with
    \begin{equation*}
        \dual{x^*,x_n}_X \leq \lambda < \dual{x^*,x}_X\quad \text{for all } n\in\N.
    \end{equation*}
    Passing to the limit $n\to\infty$ in the first inequality yields the contradiction
    \begin{equation*}
        \dual{x^*,x}_X < \dual{x^*,x}_X.
        \qedhere
    \end{equation*}
\end{proof}

If $X$ is not reflexive (e.g., $X=L^\infty(\Omega)$), we have to turn to weak-$*$ convergence.
\begin{theorem}[Banach--Alaoglu\protect\footnote{e.g., \cite[Theorem 11.6]{Clason}}]\label{thm:banachal}\index{theorem!Banach--Alaoglu}
    If $X$ is a separable normed vector space (i.e., contains a countable dense subset), then $\B_{X^*}$ is weak-$*$ sequentially compact.
\end{theorem}
By the Weierstraß Approximation theorem, both $C(\overline\Omega)$ and $L^p(\Omega)$ for $1\leq p<\infty$ are separable; also, $\ell^p$ is separable for $1\leq p <\infty$.
Hence, bounded and weak-$*$ sequentially closed balls in  $\ell^\infty$, $L^\infty(\Omega)$, and $\mathcal{M}(\Omega)$ are weak-$*$ sequentially compact.

Finally, we will also need the following \enquote{weak-$*$} separation theorem, whose proof is analogous to the proof of \cref{thm:functan:hb-separation} (using the fact that the linear weakly-$*$ continuous functionals are exactly those of the form $x^*\mapsto \dual{x^*,x}_X$ for some $x\in X$); see also \cite[Theorem~3.4(b)]{Rudin:1991}.
\begin{theorem}\label{thm:clarke:hb}
    Let $X$ be a normed vector space and $A\subset X^*$ be a nonempty, convex, and weak-$*$ sequentially closed subset and $x^*\in X^*\setminus A$. Then there exist an $x\in X$ and a $\lambda\in\R$ with
    \begin{equation*}
        \dual{z^*,x}_X \leq \lambda < \dual{x^*,x}_X \qquad\text{for all }z^*\in A.
    \end{equation*}
\end{theorem}
Note, however, that arbitrary closed convex sets in nonreflexive spaces do \emph{not} have to be weak-$*$ sequentially closed; this only holds for specific sets such as the unit ball (since compact sets are \emph{a fortiori} closed) and polar cones (for which this follows directly from the definition).

\bigskip

Since a normed vector space is characterized by its dual, this is also the case for linear operators acting on this space.
For any $T\in \linear(X;Y)$, the \term[operator!adjoint]{adjoint operator}\index{adjoint!of linear operator} $T^*\in \linear(Y^*;X^*)$ is defined via
\begin{equation*}
    \dual{T^*y^*,x}_X = \dual{y^*,Tx}_Y\qquad\text{for all } x\in X, y^*\in Y^*.
\end{equation*}
It always holds that $\norm{T^*}_{\linear(Y^*;X^*)} = \norm{T}_{\linear(X;Y)}$.
Furthermore, the continuity of $T$ implies that $T^*$ is weakly-$*$ continuous (and $T$ weakly continuous).

\section{Hilbert spaces}

Especially strong duality properties hold in Hilbert spaces.
A mapping $\innerprod{\cdot}{\cdot}:X\times X\to \R$ on a vector space $X$ over $\R$ is called \term[product!inner]{inner product}, if\index{product!scalar|see{product, inner}}
\begin{enumerate}
    \item $\innerprod{\alpha x+\beta y}{z} = \alpha\innerprod{x}{z} + \beta \innerprod{y}{z}$ for all $x,y,z\in X$ and $\alpha,\beta\in \R$;
    \item $\innerprod{x}{y} = \innerprod{y}{x}$ for all $x,y\in X$;
    \item $\innerprod{x}{x} \geq 0$ for all $x\in X$ with equality if and only if $x=0$.
\end{enumerate}
An inner product induces a norm
\begin{equation*}
    \norm{x}_X \defeq \sqrt{\innerprod{x}{x}_X},
\end{equation*}
which satisfies the \term[inequality!Cauchy--Schwarz]{Cauchy--Schwarz inequality}
\begin{equation*}
    \innerprod{x}{y}_X \leq \norm{x}_X\norm{y}_X.
\end{equation*}
If $X$ is complete with respect to the induced norm (i.e., if $(X,\norm{\cdot}_X)$ is a Banach space),
then $X$ is called a \term[space!Hilbert]{Hilbert space}; if the inner product is canonical, it is frequently omitted, and the Hilbert space is simply denoted by $X$. The spaces in \cref{ex:functan:dual}\,\ref{ex:functan:dual:i}--\ref{ex:functan:dual:iii} for $p=2(=q)$ are all Hilbert spaces, where the inner product coincides with the duality pairing and induces the canonical norm.

Directly from the definition of the induced norm we obtain the \term[expansion, binomial]{binomial expansion}
\begin{equation*}
    \norm{x+y}_X^2 = \norm{x}_X^2 + 2\innerprod{x}{y}_X + \norm{y}_X^2,
\end{equation*}
which in turn can be used to verify the \term[identity!three-point]{three-point identity}
\begin{equation}
    \label{eq:hilbert:three-point-identity}
    \innerprod{x-y}{x-z}_X
    = \frac{1}{2}\norm{x-y}_X^2
    - \frac{1}{2}\norm{y-z}_X^2
    + \frac{1}{2}\norm{x-z}_X^2 \quad\text{for all }x,y,z\in X.
\end{equation}
(This can be seen as a generalization of the classical Pythagorean theorem in plane geometry.)

\bigskip

The relevant point in our context is that the dual of a Hilbert space $X$ can be identified with $X$ itself.
\begin{theorem}[Fréchet--Riesz\protect\footnote{e.g., \cite[Theorem 16.1]{Clason}}]\label{thm:frechetriesz}\index{theorem!Fréchet--Riesz}
    Let $X$ be a Hilbert space. Then for each $x^*\in X^*$ there exists a unique $z_{x^*}\in X$ with $\norm{x^*}_{X^*} = \norm{z_{x^*}}_X$ and
    \begin{equation*}
        \dual{x^*,x}_X = \innerprod{x}{z_{x^*}}_X \qquad\text{for all } x\in X.
    \end{equation*}
\end{theorem}
The element $z_{x^*}$ is called \term[representation, Riesz]{Riesz representation} of $x^*$.
The (linear) mapping $J_X:X^*\to X$, $x^*\mapsto z_{x^*}$, is called \term[isomorphism, Riesz]{Riesz isomorphism}, and can be used to show that every Hilbert space is reflexive.

\Cref{thm:frechetriesz} allows to use the inner product instead of the duality pairing in Hilbert spaces. For example, a sequence $\{x_n\}_{n\in\N}\subset X$ converges weakly to $x\in X$ if and only if
\begin{equation*}
    \innerprod{x_n}{z}_X \to \innerprod{x}{z}_X\qquad\text{for all } z\in X.
\end{equation*}
This implies that if $x_n\weakto x$ and in addition $\norm{x_n}_X \to \norm{x}_X$ (in which case we say that $x_n$ \term[convergence!strict]{strictly converges} to $x$),
\begin{equation}
\norm{x_n-x}_X^2 = \norm{x_n}_X^2 - 2 \innerprod{x_n}{x}_X +\norm{x}_X^2 \to 0,
\end{equation}
i.e., $x_n\to x$. A normed vector space in which strict convergence implies strong convergence is said to have the \term[property!Radon--Riesz]{Radon--Riesz property}.

Similar statements hold for linear operators on Hilbert spaces. For a linear operator $T\in \linear(X;Y)$ between Hilbert spaces $X$ and $Y$, the \term[operator!adjoint!Hilbert space]{Hilbert space adjoint operator} $T^\star\in \linear(Y;X)$ is defined via
\begin{equation*}
    \innerprod{T^\star y}{x}_X = \innerprod{Tx}{y}_Y \qquad\text{for all } x\in X, y\in Y.
\end{equation*}
If $T^\star = T$, the operator $T$ is called \term[operator!self-adjoint]{self-adjoint}. A self-adjoint operator is called \term[operator!positive definite]{positive definite}, if there exists a $c>0$ such that
\begin{equation*}
    \innerprod{Tx}{x}_X \geq c \norm{x}_X^2\quad\text{for all } x\in X.
\end{equation*}
In this case, $T$ has a bounded inverse $T^{-1}$ with $\norm{T^{-1}}_{\linear(X;X)} \leq c^{-1}$. We will also use the notation $S\geq T$ for two operators $S,T:X\to X$ if
\begin{equation*}
    \innerprod{Sx}{x}_X \geq \innerprod{Tx}{x}_X \quad\text{for all }x\in X.
\end{equation*}
Hence $T$ is positive definite if and only if $T\geq c \Id$ for some $c>0$; if $T\geq 0$, we say that $T$ is merely \term[operator!positive semi-definite]{positive semi-definite}.

The Hilbert space adjoint is related to the (Banach space) adjoint via $T^\star = J_X T^* J_Y^{-1}$.
If the context is obvious, we will not distinguish the two in notation. Similarly, we will also -- by a moderate abuse of notation -- use angled brackets to denote inner products in Hilbert spaces except where we need to refer to both at the same time (which will rarely be the case, and the danger of confusing inner products with elements of a product space is much greater).

\chapter{Calculus of variations}\label{chap:variation}

We first consider the question of the existence of solutions to optimization problems of the form
\begin{equation*}
    \min_{x\in U} F(x)
\end{equation*}
for a (nonlinear) functional $F:U\to\R$ and a subset $U$ of a Banach space $X$.
Answering such questions is one of the goals of the \term{calculus of variations}.
Here we will only treat those aspects of direct relevance to nonsmooth optimization; for a full treatment, the reader is referred to, e.g., \cite{Attouch,Rindler:2018,Beck:2024}.

Note that we don't require $F$ to be defined on all of $X$; this is important for example when $F$ involves the solution of a nonlinear partial differential equation which may only exist if $x$ is sufficiently small.
For the purposes of existence of a minimizer, however, we do not need to distinguish whether $U$ represents such a domain of definition or an additional constraint in the optimization problem. In both cases, we can get rid of the constraint by extending $F$ to all of $X$ with the value $\infty$ by setting
\begin{equation*}
    \overline F:X\to\Rbar\defeq\R\cup \{\infty\},\qquad
    \overline F(x) =
    \begin{cases}
        F(x) & \text{if }x\in U,\\
        \infty & \text{if }x\in X\setminus U.
    \end{cases}
\end{equation*}
We extend the usual arithmetic on $\R$ to $\Rbar$ by letting $t<\infty$ and $t+\infty = \infty$ for all $t\in\R$; we let $0\cdot\infty=0$, but
subtraction and multiplication of negative numbers with $\infty$ and in particular $F(x) = -\infty$ is not allowed. Thus if there is any $x\in U$ at all, a minimizer $\bar x$ of $\overline F$ necessarily must lie in $U$ and coincide with a minimizer of $F$ over $U$.

\section{The direct method}

Our goal now is to find conditions under which a functional $F:X\to\Rbar$ attains a (real-valued) minimum over $X$. First, there clearly must exist a point with finite value.
We call the set on which $F$ is finite the \term[domain!effective]{effective domain}
\begin{equation*}
    \dom F \defeq \setof{x\in X}{F(x)<\infty}.
\end{equation*}
If $\dom F \neq \emptyset$, the functional $F$ is called \term[functional!proper]{proper}.\footnote{Note that some references include the value $-\infty$ in their definition of $\Rbar$; in this case a functional $F$ is called proper if $\dom F\neq \emptyset$ and additionally $F>-\infty$ everywhere. In either case, functionals taking on the value $-\infty$ are excluded from the discussion.}

Next, we require a form of continuity to prevent the function from \enquote{jumping over} possible minima.
We call $F$ \term[functional!lower semicontinuous]{lower semicontinuous} in $x\in X$ if
\begin{equation*}
    F(x) \leq \liminf_{n\to\infty} F(x_n)\qquad\text{for every } \{x_n\}_{n\in\N}\subset X \text{ with  }x_n\to x,
\end{equation*}
see \cref{fig:variation:lsc}, where $F_2$ is an example for a function that is not lower semicontinuous and does not attain a minimum.
Analogously, we define \emph{weakly(-$*$) lower semicontinuous} functionals via weakly(-$*$) convergent sequences.

Finally, we need to prevent the function from having a \enquote{minimum at infinity}. Here we use the following property: We call $F$ \term[functional!coercive]{coercive} if for every sequence $\{x_n\}_{n\in\N}\subset X$ with $\norm{x_n}_X\to\infty$ we also have $F(x_n)\to \infty$.
\begin{figure}
    \centering
    \begin{subfigure}[t]{0.495\textwidth}
        \centering
        \begin{asy}
            unitsize(75,75);
            draw((-1.2,0)..(1.2,0),linewidth(0.5),Arrow);
            draw((0,-0.2)..(0,1.2),linewidth(0.5),Arrow);

            real getx(int n){ return -1/n ;}

            real x = 0;
            real F1(real x) {return 0.25 + x^2/2;}
            real F2(real x) {return 0.5 + x^2;}

            path F1g = graph(F1,-1,0);
            path F2g = graph(F2,0,0.75);
            draw(F1g, primalline + linewidth(1.5));
            draw(F2g, primalline + linewidth(1.5));

            dot((0,0));
            dot((0,F1(0)));
            dot((0,F2(0)), primalline, filltype=UnFill);

            label("$x$",(x,0),SE);
            label("$F_1(x)$", (x,F1(x)),E);

            for(int n=2; n<25; ++n){
                real xn = getx(n);
                dot((xn,0));
                dot((xn,F1(xn)));
            }
            real xn = getx(4);
            label("$x_n$",(xn,0),S);
            label("$F_1(x_n)$",(xn,F1(xn)),N);
        \end{asy}
        \caption{$F_1$ is lower semicontinuous at $x$}\label{fig:variation:lsc1}
    \end{subfigure}
    \hfill
    \begin{subfigure}[t]{0.495\textwidth}
        \centering
        \begin{asy}
            unitsize(75,75);
            draw((-1.2,0)..(1.2,0),linewidth(0.5),Arrow);
            draw((0,-0.2)..(0,1.2),linewidth(0.5),Arrow);

            real getx(int n){ return -1/n ;}

            real x = 0;
            real F1(real x) {return 0.25 + x^2/2;}
            real F2(real x) {return 0.5 + x^2;}

            path F1g = graph(F1,-1,0);
            path F2g = graph(F2,0,0.75);
            draw(F1g, primalline + linewidth(1.5));
            draw(F2g, primalline + linewidth(1.5));

            dot((0,0));
            dot((0,F1(0)), primalline, filltype=UnFill);
            dot((0,F2(0)));

            label("$x$",(x,0),SE);
            label("$F_2(x)$", (x,F2(x)),W);

            for(int n=2; n<25; ++n){
                real xn = getx(n);
                dot((xn,0));
                dot((xn,F1(xn)));
            }
            real xn = getx(4);
            label("$x_n$",(xn,0),S);
            label("$F_2(x_n)$",(xn,F1(xn)),S);
        \end{asy}
        \caption{$F_2$ is \emph{not} lower semicontinuous at $x$}\label{fig:variation:lsc2}
    \end{subfigure}
    \caption{Illustration of lower semicontinuity: two functions $F_1,F_2:\R\to\R$ and a sequence $\{x_n\}_{n\in\N}$ realizing their (identical) limes inferior.}
    \label{fig:variation:lsc}
\end{figure}

We now have everything at hand to prove the central existence result in the calculus of variations. The strategy for its proof is known as the \term[method!direct]{direct method}.\footnote{This strategy is applied so often in the literature that one usually just writes \enquote{Existence of a minimizer follows from the direct method.} or even just \enquote{Existence follows from standard arguments.}
The basic idea goes back to Hilbert; the version based on lower semicontinuity which we use here is due to \href{https://mathshistory.st-andrews.ac.uk/Biographies/Tonelli/}{Leonida Tonelli} (1885--1946), who through it had a lasting influence on the modern calculus of variations.}
\begin{theorem}\label{thm:variation:existence}
    Let $X$ be a reflexive Banach space and $F:X\to\Rbar$ be proper, coercive, and weakly lower semicontinuous. Then the minimization problem
    \begin{equation*}
        \min_{x\in X} F(x)
    \end{equation*}
    has a solution $\bar x\in \dom F$.
\end{theorem}
\begin{proof}
    The proof can be separated into three steps.
    \begin{enumerate}
        \item\label{thm:variation:existence:i}
            \emph{Pick a minimizing sequence.}

            Since $F$ is proper, there exists an $M\defeq\inf_{x\in X} F(x)<\infty$ (although $M=-\infty$ is not excluded so far). We can thus find a sequence $\{y_n\}_{n\in\N}\subset \range F\setminus\{\infty\}\subset \R$ with $y_n\to M$, i.e., there exists a sequence $\{x_n\}_{n\in\N}\subset X$ with
            \begin{equation*}
                F(x_n) \to M = \inf_{x\in X} F(x).
            \end{equation*}
            Such a sequence is called \term[sequence!minimizing]{minimizing sequence}. Note that from the convergence of $\{F(x_n)\}_{n\in\N}$ we cannot conclude the convergence of $\{x_n\}_{n\in\N}$ (yet).

        \item\label{thm:variation:existence:ii}
            \emph{Show that the minimizing sequence contains a weakly convergent subsequence.}

            Assume to the contrary that $\{x_n\}_{n\in\N}$ is unbounded, i.e., that $\norm{x_n}_X\to\infty$ for $n\to \infty$. The coercivity of $F$ then implies that $F(x_n)\to \infty$ as well, in contradiction to $F(x_n)\to M<\infty$ by definition of the minimizing sequence.
            Hence, the sequence is bounded, i.e., there is an $R>0$ with $\norm{x_n}_X\leq R$ for all $n\in\N$.
            In particular, $\{x_n\}_{n\in\N} \subset \B(0,R)$. The Eberlein--\u{S}mulyan theorem (\cref{thm:ebsmul}) therefore implies the existence of a weakly converging subsequence $\{x_{n_k}\}_{k\in\N}$ with limit $\bar x\in X$. (This limit is a candidate for the minimizer.)

        \item\label{thm:variation:existence:iii}
            \emph{Show that its limit is a minimizer.}

            From the definition of the minimizing sequence, we also have $F(x_{n_k})\to M$ for $k\to\infty$. Together with the weak lower semicontinuity of $F$ and the definition of the infimum we thus obtain
            \begin{equation*}
                \inf_{x\in X} F(x) \leq F(\bar x) \leq \liminf_{k\to\infty} F(x_{n_k}) = M = \inf_{x\in X} F(x)<\infty.
            \end{equation*}
            This implies that $\bar x\in \dom F$ and that $\inf_{x\in X} F(x)=F(\bar x)> -\infty$. Hence, the infimum is attained in $\bar x$ which is therefore the desired minimizer.
            \qedhere
    \end{enumerate}
\end{proof}
\begin{remark}\label{rem:variation:existence-dual}
    If $X$ is not reflexive but the dual of a separable Banach space, we can argue analogously for weakly-$*$ lower semicontinuous functionals using the Banach--Alaoglu theorem (\cref{thm:banachal}).
\end{remark}

Note how the topology on $X$ used in the proof is restricted in step \ref{thm:variation:existence:ii} and \ref{thm:variation:existence:iii}: Step \ref{thm:variation:existence:ii} profits from a coarse topology (in which more sequences are convergent), while step \ref{thm:variation:existence:iii} profits from a fine topology (the fewer sequences are convergent, the easier it is to satisfy the $\liminf$ conditions). Since in the cases of interest to us no more than boundedness of a minimizing sequence can be expected, we cannot use a finer than the weak topology. We thus have to ask whether a sufficiently large class of (interesting) functionals are weakly lower semicontinuous.

A first example is the class of bounded linear functionals: For any $x^*\in X^*$, the functional
\begin{equation*}
    F:X\to\R, \qquad x\mapsto \dual{x^*,x}_X,
\end{equation*}
is weakly continuous by definition of weak convergence and hence \emph{a fortiori} weakly lower semicontinuous.
Another advantage of (weak) lower semicontinuity is that it is preserved under certain operations.
\begin{lemma}\label{lem:variation:wlsc}
    Let $X$ and $Y$ be Banach spaces and $F:X\to\Rbar$ be weakly(-$*$) lower semicontinuous. Then the following functionals are weakly(-$*$) lower semicontinuous as well:
    \begin{enumerate}
        \item\label{lem:variation:wlsc:i}
            $\alpha F$ for all $\alpha \geq 0$;
        \item \label{lem:variation:wlsc:ii}
            $F+G$ for $G:X\to\Rbar$ weakly(-$*$) lower semicontinuous;
        \item \label{lem:variation:wlsc:iii}
            $\phi\circ F$ for $\phi:\Rbar\to\Rbar$ lower semicontinuous and monotonically increasing.
        \item \label{lem:variation:wlsc:iv}
            $F\circ \Phi$ for $\Phi:Y\to X$ weakly(-$*$) continuous, i.e., $y_n\weakto^{(*)} y$ implies $\Phi(y_n)\weakto^{(*)} \Phi(y)$;
        \item \label{lem:variation:wlsc:v}
            $x\mapsto \sup_{i\in I} F_i(x)$ with $F_i:X\to\Rbar$ weakly(-$*$) lower semicontinuous for all $i\in I$ and an arbitrary set $I$.
    \end{enumerate}
\end{lemma}
Note that \ref{lem:variation:wlsc:v} does \emph{not} hold for continuous functions.
\begin{proof}
    We only show the claim for the case of weak lower semicontinuity; the statements for weak-$*$ lower semicontinuity follow by the same arguments.

    Statements \ref{lem:variation:wlsc:i} and \ref{lem:variation:wlsc:ii} follow directly from the properties of the limes inferior.

    For statement \ref{lem:variation:wlsc:iii}, it first follows from the monotonicity of $\phi$ and the weak lower semicontinuity of $F$ that $x_n\weakto x$ implies
    \begin{equation*}
        \phi(F(x)) \leq \phi(\liminf_{n\to\infty} F(x_n)).
    \end{equation*}
    It remains to show that the right-hand side can be bounded by $\liminf_{n\to\infty} \phi(F(x_n))$. For that purpose, we consider the subsequence $\{x_{n_k}\}_{k\in\N}$ for which $\liminf_{n\to\infty} \phi(F(x_n)) = \lim_{k\to\infty} \phi(F(x_{n_k}))$. By passing to a further subsequence which we index by $k'$, we can also obtain that $\liminf_{k\to \infty} F(x_{n_k}) = \lim_{k'\to\infty} F(x_{n_{k'}})$. Since the $\liminf$ restricted to a subsequence can never be smaller than that of the full sequence, the monotonicity of $\phi$ together with its lower semicontinuity now implies that
    \begin{equation*}
        \phi(\liminf_{n\to\infty} F(x_n)) \leq \phi(\lim_{k'\to\infty} F(x_{n_{k'}}))
        \leq \liminf_{k'\to\infty} \phi(F(x_{n_{k'}}))
        = \liminf_{n\to\infty} \phi(F(x_{n})),
    \end{equation*}
    where we have used in the last step that a subsequence of the convergent sequence $\{\phi(F(x_{n_k}))\}_{k\in\N}$ has the same limit, which by construction realizes the $\liminf$.

    Statement \ref{lem:variation:wlsc:iv} follows directly from the weak continuity of $\Phi$, as $y_n\weakto y$ implies that $x_n \defeq \Phi(y_n) \weakto \Phi(y) =: x$, and the lower semicontinuity of $F$ yields
    \begin{equation*}
        F(\Phi(y)) \leq \liminf_{n\to\infty} F(\Phi(y_n)).
    \end{equation*}

    Finally, let $\{x_n\}_{n\in\N}$ be a weakly converging sequence with limit $x\in X$. Then the definition of the supremum implies that
    \begin{equation*}
        F_j(x) \leq \liminf_{n\to\infty} F_j(x_n) \leq \liminf_{n\to\infty} \sup_{i\in I} F_i(x_n) \qquad\text{for all }j\in I.
    \end{equation*}
    Taking the supremum over all $j\in I$ on both sides yields statement \ref{lem:variation:wlsc:v}.
\end{proof}

\begin{corollary}\label{cor:variation:norm}
    If $X$ is a Banach space, then the norm $\norm{\cdot}_X$ is proper, coercive, and weakly lower semicontinuous. Similarly, the dual norm $\norm{\cdot}_{X^*}$ is proper, coercive, and weakly-$*$ lower semicontinuous.
\end{corollary}
\begin{proof}
    Coercivity and $\dom \norm{\cdot}_X = X$ follow directly from the definition. Weak lower semicontinuity follows from \cref{lem:variation:wlsc}\,\ref{lem:variation:wlsc:v} and \cref{cor:functan:norm_dual} since
    \begin{equation*}
        \norm{x}_X = \sup_{\norm{x^*}_{X^*}\leq 1} |\dual{x^*,x}_X|.
    \end{equation*}
    The claim for $\norm{\cdot}_{X^*}$ follows analogously using the definition of the operator norm in place of \cref{cor:functan:norm_dual}.
\end{proof}

Another frequently occurring functional is the \term[function!indicator]{indicator function}\footnote{not to be confused with the \term[function!characteristic|infn]{characteristic function} $\1_U$ with $\1_U(x) = 1$ for $x\in U$ and $0$ else} of a set $U\subset X$, defined as
\begin{equation*}
    \delta_U(x) = \begin{cases} 0 & x\in U,\\ \infty& x \in X\setminus U.\end{cases}
\end{equation*}
The purpose of this definition is of course to write the minimization of a functional $F:X\to\R$ (i.e., defined on all of $X$) under the additional constraint $x\in U$ as the minimization of $\overline F\defeq F+\delta_U$ over $X$. The following result is therefore important for showing the existence of such a constrained minimizer.
\begin{lemma}\label{lem:variation:indicator}
    Let $X$ be a Banach space and $U\subset X$. Then $\delta_U:X\to \Rbar$ is
    \begin{enumerate}
        \item \label{lem:variation:indicator:i}
            proper if $U$ is nonempty;
        \item \label{lem:variation:indicator:ii}
            weakly lower semicontinuous if $U$ is convex and closed;
        \item \label{lem:variation:indicator:iii}
            coercive if $U$ is bounded.
    \end{enumerate}
\end{lemma}
\begin{proof}
    Statement \ref{lem:variation:indicator:i} is clear. For \ref{lem:variation:indicator:ii}, consider a weakly converging sequence $\{x_n\}_{n\in\N}\subset X$ with limit $x\in X$. If $x\in U$, then $\delta_U\geq 0$ immediately yields
    \begin{equation*}
        \delta_U(x) = 0 \leq \liminf_{n\to\infty} \delta_U(x_n).
    \end{equation*}
    Let now $x\notin U$. Since $U$ is convex and closed and hence by \cref{lem:convex_closed} also weakly closed, there must be a $N\in\N$ with $x_n\notin U$ for all $n\geq N$ (otherwise we could -- by passing to a subsequence if necessary -- construct a sequence with $x_{n}\weakto x\in U$, in contradiction to the assumption). Thus, $\delta_U(x_n) = \infty$ for all $n\geq N$, and therefore
    \begin{equation*}
        \delta_U(x) = \infty = \liminf_{n\to\infty} \delta_U(x_n).
    \end{equation*}

    For \ref{lem:variation:indicator:iii}, let $U$ be bounded, i.e., there exists an $M>0$ with $U\subset \B(0,M)$.
    If $\norm{x_n}_X\to\infty$, then there exists an $N\in\N$ with $\norm{x_n}_X >M$ for all $n\geq N$, and thus $x_n\notin \B(0,M)\supset U$ and hence $\delta_U(x_n)=\infty$ for all $n\geq N$. This implies that $\delta_U(x_n) \to \infty$ as well.
\end{proof}

\section{Differential calculus in normed vector spaces}

To characterize minimizers of functionals on infinite-dimensional spaces using the Fermat principle, we transfer the classical derivative concepts to normed vector spaces.

Let $X$ and $Y$ be normed vector spaces, $F:X\to Y$ be a mapping, and $x,h\in X$ be given.
\begin{itemize}
    \item If the one-sided limit
        \begin{equation*}
            F'(x;h) \defeq \lim_{t\downto 0} \frac{F(x+th)-F(x)}{t}\in Y
        \end{equation*}
        (where $t\downto 0$ denotes the limit for arbitrary positive decreasing null sequences)
        exists, it is called the \term[derivative!directional]{directional derivative} of $F$ in $x$ in direction $h$.
    \item If $F'(x;h)$ exists for all $h\in X$ and
        \begin{equation*}
            DF(x):X\to Y,\qquad h\mapsto F'(x;h)
        \end{equation*}
        defines a bounded linear operator, we call $F$ \term[mapping!differentiable!Gateaux]{Gateaux differentiable} (at $x$) and $DF(x)\in \linear(X;Y)$ its \term[derivative!Gateaux]{Gateaux derivative}.
    \item If additionally
        \begin{equation*}
            \lim_{\norm{h}_X\to 0} \frac{\norm{F(x+h) - F(x) - DF(x)h}_Y}{\norm{h}_X} = 0,
        \end{equation*}
        then $F$ is called \term[mapping!differentiable!Fréchet]{Fréchet differentiable} (at $x$) and $F'(x)\defeq DF(x)\in \linear(X;Y)$ its \term[derivative!Fréchet]{Fréchet derivative}.
    \item If additionally the mapping $F':X\to \linear(X; Y)$ is (Lipschitz) continuous, we call $F$ \term[mapping!differentiable!continuously]{(Lipschitz) continuously differentiable}.
\end{itemize}
The difference between Gateaux and Fréchet differentiability lies in the approximation error of $F$ near $x$ by $F(x)+DF(x)h$: while it only has to be bounded in $\norm{h}_X$ -- i.e., linear in $\norm{h}_X$ -- for a Gateaux differentiable function, it has to be superlinear in $\norm{h}_X$ if $F$ is Fréchet differentiable. (For a \emph{fixed} direction $h$, this is of course also the case for Gateaux differentiable functions; Fréchet differentiability thus additionally requires a uniformity in $h$.) We also point out that continuous differentiability always entails Fréchet differentiability.

\begin{remark}\label{rem:strictlydiff}
    Sometimes a weaker notion than continuous differentiability is used. A mapping $F:X\to Y$ is called \term[mapping!differentiable!strictly]{strictly differentiable} in $x$ if
    \begin{equation}
        \lim_{\substack{y\to x\\\norm{h}_X\to 0}} \frac{\norm{F(y+h)-F(y)-F'(x)h}_Y}{\norm{h}_X} = 0.
    \end{equation}
    The benefit of this definition over that of continuous differentiability is that the limit process is now in the function $F$ rather than the derivative $F'$; strict differentiability can therefore hold if every neighborhood of $x$ contains points where $F$ is not differentiable. However, if $F$ is differentiable everywhere in a neighborhood of $x$, then $F$ is strictly differentiable if and only if $F'$ is continuous; see \cite[Proposition 1D.7]{dontchev2014implicit}. Although many results of \crefrange{chap:clarke}{chap:colimiting} actually hold under the weaker assumption of strict differentiability, we will therefore work only with the more standard notion of continuous differentiability.
\end{remark}

If $F$ is Gateaux differentiable, the Gateaux derivative can be computed via
\begin{equation*}
    DF(x) h = \left(\tfrac{d}{dt}F(x+th)\right)\Big|_{t=0}.
\end{equation*}
Bounded linear operators $F\in \linear(X;Y)$ are obviously Fréchet differentiable with derivative $F'(x) = F \in \linear(X;Y)$ for all $x\in X$.
Derivatives of other functions can be obtained through the usual calculus rules, whose proofs in normed vector spaces are exactly as in $\R^N$. As an example, we prove a chain rule.
\begin{theorem}\label{thm:frechet_chain}
    Let $X$, $Y$, and $Z$ be normed vector spaces, and let $F:X\to Y$ be Fréchet differentiable at $x\in X$ and $G:Y\to Z$ be Fréchet differentiable at $y\defeq F(x)\in Y$. Then $G\circ F$ is Fréchet differentiable at $x$ and
    \begin{equation*}
        (G\circ F)'(x) = G'(F(x)) F'(x).
    \end{equation*}
\end{theorem}
\begin{proof}
    For $h\in X$ with $x+h\in\dom F$ we have
    \begin{equation*}
        (G\circ F)(x+h ) - (G\circ F)(x) = G(F(x+h))-G(F(x)) = G(y+g)  - G(y)
    \end{equation*}
    with $g \defeq F(x+h)-F(x)$. The Fréchet differentiability of $G$ thus implies that
    \begin{equation*}
        \norm{(G\circ F)(x+h ) - (G\circ F)(x) - G'(y)g }_Z =  r_1(\norm{g}_Y)
    \end{equation*}
    with $r_1(t)/t \to 0$ for $t\to 0$. The Fréchet differentiability of $F$ further implies
    \begin{equation*}
        \norm{g - F'(x)h}_Y = r_2(\norm{h}_X)
    \end{equation*}
    with $r_2(t)/t \to 0$ for $t\to 0$. In particular,
    \begin{equation}\label{eq:frechet_chain:est}
        \norm{g}_Y \leq \norm{F'(x)h}_Y + r_2(\norm{h}_X).
    \end{equation}
    Hence, with $c\defeq\norm{G'(F(x))}_{\linear(Y;Z)}$ we have
    \begin{equation*}
        \norm{(G\circ F)(x+h) - (G\circ F)(x) -  G'(F(x)) F'(x)h}_Z \leq  r_1(\norm{g}_Y) +  c\, r_2(\norm{h}_X).
    \end{equation*}
    If $\norm{h}_X\to 0$, we obtain from \eqref{eq:frechet_chain:est} and $F'(x)\in \linear(X;Y)$ that $\norm{g}_Y\to 0$ as well, and the claim follows.
\end{proof}
A similar rule for Gateaux derivatives does not hold, however.

Of special importance in \cref{part:setvalued} will be the following inverse function theorem, whose proof can be found, e.g., in \cite[Theorem XIV.1.2]{Lang:1993}.
\begin{theorem}[inverse function theorem]\label{thm:inversefunctiontheorem}\index{theorem!inverse function}
    Let $X,Y$ be Banach spaces and $F:X\to Y$ be continuously differentiable on a neighborhood of $x\in X$. If  $F'(x):X\to Y$ is bijective, then there exists a neighborhood $U$ of $x$ and a continuously differentiable mapping $F^{-1}:F(U)\to U$ such that $F^{-1}(F(x)) = x$ for all $x\in U$ and $F(F^{-1}(y))=y$ for all $y\in F(U)$.
\end{theorem}

\bigskip

Of particular relevance in optimization is of course the special case $F:X\to\R$, where $DF(x)\in \linear(X;\R)=X^*$ (if the Gateaux derivative exists). Following the usual notation from \cref{sec:functan:dual}, we will then write $F'(x; h)=\dual{DF(x),h}_X$ for the directional derivative in direction $h\in X$. Our first result is the classical Fermat principle characterizing minimizers of differentiable functions.
\begin{theorem}[Fermat principle]\label{thm:variation:fermat}\index{principle!Fermat!smooth}
    Let $F:X\to\R$ be Gateaux differentiable and $\bar x \in X$ be a local minimizer of $F$. Then $DF(\bar x) = 0$, i.e.,
    \begin{equation*}
        \dual{DF(\bar x), h}_X = 0 \qquad\text{for all }h\in X.
    \end{equation*}
\end{theorem}
\begin{proof}
    Let $h\in X$ be arbitrary. Since $\bar x$ is a local minimizer, the core--int \cref{lem:functan:coreint} implies that there exists an $\epsilon >0$ such that $F(\bar x) \leq F(\bar x+th)$ for all $t\in (0,\eps)$, i.e.,
    \begin{equation}
        0\leq \frac{F(\bar x+th) - F(\bar x)}{t} \to F'(\bar x; h) = \dual{DF(\bar x),h}_X \quad\text{for }t\to 0,
    \end{equation}
    where we have used the Gateaux differentiability and hence directional differentiability of $F$. Since the right-hand side is linear in $h$, the same argument for $-h$ yields $\dual{DF(\bar x), h}_X\leq 0$ and therefore the claim.
\end{proof}

We will also need the following version of the mean value theorem.
\begin{theorem}\label{thm:frechet:mean}
    Let $F:X\to \R$ be continuously differentiable. Then for all $x,h\in X$,
    \begin{equation*}
        F(x+h) - F(x) =\int_0^1 \dual{F'(x+th),h}_X\,dt.
    \end{equation*}
\end{theorem}
\begin{proof}
    Consider the scalar function
    \begin{equation*}
        f:[0,1]\to \R,\qquad t\mapsto F(x+th),
    \end{equation*}
    which by the definition of the directional derivative and of Gateaux differentiability satisfies
    \[
        f'(t) = F'(x + th; h) = \dual{F'(x + th),h}_X.
    \]
    By the assumption on $F$ and the continuity of the duality pairing, $f'$ is continuous on $[0,1]$, and hence we can apply the fundamental theorem of calculus in $\R$ to obtain that
    \begin{equation*}
        F(x+h)-F(x) = f(1)-f(0) = \int_0^1f'(t)\,dt =  \int_0^1 \dual{F'(x+th),h}_X\,dt.
        \qedhere
    \end{equation*}
\end{proof}

As in classical analysis, this result is useful for relating local and pointwise properties of smooth functions. A typical example is the following lemma.
\begin{lemma}\label{lem:variation:c1-lipschitz}
    Let $F:X\to Y$ be continuously differentiable at $x\in X$. Then $F$ is locally Lipschitz continuous near $x$.
\end{lemma}
\begin{proof}
    By assumption, there exists a neighborhood $U\subset X$ of $x$ such that $F':U\to \linear(X; Y)$ exists and is continuous at $x\in U$. Hence there exists a $\delta>0$ with $\norm{F'(z)-F'(x)}_{\linear(X;Y)} \leq 1$ and hence $\norm{F'(z)}_{\linear(X;Y)} \leq 1 +\norm{F'(x)}_{\linear(X;Y)}$ for all $z\in \B(x,\delta)\subset U$.
    For any $x_1,x_2\in \B(x,\delta)$ we also have $x_2 + t(x_1-x_2)\in \B(x,\delta)$ for all $t\in [0,1]$ (since balls in normed vector spaces are convex), and hence \cref{thm:frechet:mean} implies that
    \begin{equation*}
        \begin{aligned}[b]
            \norm{F(x_1)-F(x_2)}_Y &\leq \int_0^1 \norm{F'(x_2 + t(x_1-x_2))}_{\linear(X;Y)}\norm{x_1-x_2}_X \,dt  \\
            &\leq (1 +\norm{F'(x)}_{\linear(X;Y)})\norm{x_1-x_2}_X,
        \end{aligned}
    \end{equation*}
    and thus local Lipschitz continuity near $x$ with constant $L=1+\norm{F'(x)}_{\linear(X; Y)}$.
\end{proof}
\begin{remark}
    \label{rem:variation:frechet-lipschitz-at}
    If $F$ is merely Fréchet differentiable at $x\in X$ but not continuously differentiable near $x$, we can still deduce locally Lipschitz continuity \emph{at} (but not necessarily \emph{near}) $x$ for any factor $L>\norm{F'(x)}_{\linear(X; Y)}$ since
    by the definition of Fréchet differentiability and the inverse triangle inequality, for any $\epsilon>0$ there exists a $\delta>0$ such that
    \[
        \norm{F(x+h) - F(x)}_Y \le \norm{F'(x)h}_Y + \epsilon\norm{h}_X\leq (\norm{F'(x)}_{\linear(X; Y)} + \eps)\norm{h}_X
    \]
    for every $\norm{h}_X \le \delta$.
\end{remark}

Note that since the Gateaux derivative of $F:X\to\R$ is an element of $X^*$, it cannot be added to elements in $X$ (as required for, e.g., a steepest descent method). However, in Hilbert spaces (and in particular in $\R^N$), we can use the Fréchet--Riesz \cref{thm:frechetriesz} to identify $DF(x)\in X^*$ with an element $\nabla F(x) \in X$, called  the \term{gradient} of $F$ at $x$, in a canonical way via
\begin{equation*}
    \dual{DF(x),h}_X = \innerprod{\nabla F(x)}{h}_X \qquad\text{for all } h\in X.
\end{equation*}
We illustrate this with a simple example.
\begin{example}\label{ex:variation:gradient1}
    Let $F(x) = \frac12\norm{x}_X^2=\frac12\innerprod{x}{x}_X$. Then we have for all $x,h\in X$ that
    \begin{equation*}
        F'(x;h) = \lim_{t\downto 0} \frac{\frac12\innerprod{x+th}{x+th}_X - \frac12\innerprod{x}{x}_X}{t} = \innerprod{x}{h}_X = \dual{DF(x),h}_X,
    \end{equation*}
    since the inner product is linear in $h$ for fixed $x$.
    Hence, the squared norm is Gateaux differentiable at every $x\in X$ with derivative $DF(x) = h\mapsto \innerprod{x}{h}_X\in X^*$; it is even Fréchet differentiable since
    \begin{equation*}
        \lim_{\norm{h}_X\to 0} \frac{\left|\frac12\norm{x+h}_X^2 - \frac12\norm{x}_X^2 - (x,h)_X\right|}{\norm{h}_X} = \lim_{\norm{h}_X\to 0}\frac12 \norm{h}_X = 0.
    \end{equation*}
    The gradient $\nabla F(x) \in X$ by definition is given by
    \begin{equation*}
        \innerprod{\nabla F(x)}{h}_X = \dual{DF(x),h}_X = \innerprod{x}{h}_X \qquad\text{for all } h\in X,
    \end{equation*}
    i.e., $\nabla F(x) = x$.
\end{example}
The following example demonstrates how the gradient (in contrast to the derivative) depends on the inner product on $X$ -- which may be different from the inner product inducing the squared norm.
\begin{example}\label{ex:variation:gradient2}
    Let $M\in\linear(X; X)$ be self-adjoint and positive definite (and thus continuously invertible). Then $\innerprod{x}{y}_Z\defeq\innerprod{Mx}{y}_X$ also defines an inner product on the vector space $X$ and induces an (equivalent) norm $\norm{x}_Z \defeq \innerprod{x}{x}_Z^{1/2}$ on $X$. Hence $(X,\innerprod{\cdot}{\cdot}_Z)$ is a Hilbert space as well, which we will denote by $Z$. Consider now the functional $\tilde F:Z\to\R$ with $\tilde F(x) \defeq \frac12\norm{x}_X^2$ (which is well-defined since $\norm{\cdot}_X$ is also an equivalent norm on $Z$). Then, the derivative $D\tilde{F}(x)\in Z^*$ is still given by $\dual{D \tilde F(x),h}_Z = \innerprod{x}{h}_X$ for all $h\in Z$ (or, equivalently, for all $h\in X$ since we defined $Z$ via the same vector space). However, $\nabla \tilde F(x)\in Z$ is now characterized by
    \begin{equation*}
        \innerprod{x}{h}_X = \dual{D\tilde{F}(x),h}_Z = \innerprod{\nabla \tilde F(x)}{h}_{Z} = \innerprod{M\nabla \tilde F(x)}{h}_X \qquad\text{for all } h\in Z,
    \end{equation*}
    i.e., $\nabla \tilde F(x) = M^{-1} x\neq \nabla F(x)$.
\end{example}
(The situation is even more delicate if $M$ is only positive definite on a subspace, as in the case of $X=L^2(\Omega)$ and $Z=H^1(\Omega)$.)

\bigskip

Higher-order derivatives can be obtained by iterating these definitions; in particular, a mapping $F:X\to Y$ is called \term[mapping!differentiable!twice]{twice (Fréchet) differentiable} if it has a \term[derivative!second-order]{second-order (Fréchet) derivative}
\begin{equation*}
    F'':X\to \linear(X; \linear(X; Y)).
\end{equation*}
If the mapping $x\mapsto F''(x)$ is continuous, then $F$ is called \emph{twice continuously differentiable}. As Fréchet differentiable mappings are continuous, any twice differentiable mapping is \emph{a fortiori} continuously differentiable.

Again, we are most interested in the special case where $F:X \to \R$ is twice differentiable. In this case, for every $x,h\in X$ it holds that $F''(x)h\in \linear(X; \R) = X^*$; hence $F''(x)$ defines a \emph{quadratic form}
\begin{equation}\label{eq:variation:quadratic}
    F''(x): X\times X \to \R,\qquad (h_1,h_2) \mapsto \dualprod{F''(x)h_1}{h_2}_X,
\end{equation}
which in Hilbert spaces can be identified with the \term{Hessian} $\grad^2 F(x)\in \linear(X; X)$.

\section{Superposition operators}\label{sec:superposition}

A special class of operators on function spaces arises from pointwise application of a real-valued function, e.g., $u(x)\mapsto\sin(u(x))$. We thus consider for $f:\Omega\times \R\to\R$ with $\Omega\subset\R^d$ open and bounded as well as $p,q\in[1,\infty]$ the corresponding \term[operator!superposition]{superposition} or \term[operator!Nemytskii]{Nemytskii operator}
\begin{equation}\label{eq:superposition}
    F:L^p(\Omega)\to L^q(\Omega),\qquad [F(u)](x) = f(x,u(x)) \quad\text{for almost every }x\in \Omega.
\end{equation}
For this operator to be well-defined requires certain restrictions on $f$. We call $f:\Omega\times\R\to\R$ a \term[function!Carathéodory]{Carathéodory function} if
\begin{enumerate}
    \item for all $z\in\R$, the mapping $x\mapsto f(x,z)$ is measurable;
    \item for almost every $x\in\Omega$, the mapping $z\mapsto f(x,z)$ is continuous.
\end{enumerate}
We additionally require the following growth condition: For given $p,q\in[1,\infty)$ there exist $a\in L^q(\Omega)$ and $b\in L^\infty(\Omega)$ with
\begin{equation}\label{eq:superpos:growth}
    |f(x,z)| \leq a(x) +b(x)|z|^{p/q}.
\end{equation}
Under these conditions, $F$ is well-defined and even continuous.
\begin{theorem}\label{thm:superpos:continuous}
    If the Carathéodory function $f:\Omega\times \R\to\R$ satisfies the growth condition \eqref{eq:superpos:growth} for $p,q\in[1,\infty)$, then the superposition operator  $F:L^p(\Omega)\to L^q(\Omega)$ defined via \eqref{eq:superposition} is continuous.
\end{theorem}
\begin{proof}
    We sketch the essential steps; a complete proof can be found in, e.g., \cite[Theorems 3.1, 3.7]{Appell:1990}. First, one shows for given $u\in L^p(\Omega)$ the measurability of $F(u)$ using the Carathéodory properties. It then follows from \eqref{eq:superpos:growth} and the triangle inequality that
    \begin{equation*}
        \norm{F(u)}_{L^q} \leq \norm{a}_{L^q} + \norm{b}_{L^\infty}\norm{|u|^{p/q}}_{L^q}  = \norm{a}_{L^q} + \norm{b}_{L^\infty}\norm{u}^{p/q}_{L^p}<\infty,
    \end{equation*}
    i.e., $F(u)\in L^q(\Omega)$.

    To show continuity, we consider a sequence $\{u_n\}_{n\in\N}\subset L^p(\Omega)$ with $u_n\to u\in L^p(\Omega)$. Then there exists a subsequence, again denoted by $\{u_n\}_{n\in \N}$, that converges pointwise almost everywhere in $\Omega$, as well as a $v\in L^p(\Omega)$ with $|u_n(x)| \leq |v(x)| + |u_1(x)|=:g(x)$ for all $n\in \N$ and almost every $x\in\Omega$ (see, e.g., \cite[Lemma 3.22 as well as (3-14) in the proof of Theorem 3.17]{Alt:2016}).
    The continuity of $z\mapsto f(x,z)$ then implies $F(u_n)\to F(u)$ pointwise almost everywhere as well as
    \begin{equation*}
        |[F(u_n)](x)| \leq a(x) + b(x)|u_n(x)|^{p/q} \leq  a(x) + b(x)|g(x)|^{p/q}\quad\text{for almost every }x\in\Omega.
    \end{equation*}
    Since $g\in L^p(\Omega)$, the right-hand side defines a function in $L^q(\Omega)$, and we can apply Lebesgue's dominated convergence theorem to deduce that $F(u_n)\to F(u)$ in $L^q(\Omega)$.
    As this argument can be applied to any subsequence, the whole sequence must converge to $F(u)$, which yields the claimed continuity.
\end{proof}
In fact, the growth condition \eqref{eq:superpos:growth} is also necessary for continuity; see \cite[Theorem 3.2]{Appell:1990}.
In addition, it is straightforward to show that for $p=q=\infty$, the growth condition \eqref{eq:superpos:growth} (with $p/q \defeq 0$ in this case) implies that $F$ is even locally Lipschitz continuous.

Similarly, one would like to show that differentiability of $f$ implies differentiability of the corresponding superposition operator $F$, ideally with ``pointwise'' derivative $[F'(u)h](x) = f'(u(x))h(x)$ (compare \cref{ex:functan:dual}\,\ref{ex:functan:dual:iii}).
However, this does not hold in general; for example, the superposition operator defined by $f(x,z)=\sin(z)$ is \emph{not} differentiable at $u=0$ for  $1\leq p=q<\infty$. The reason is that for a Fréchet differentiable superposition operator $F:L^p(\Omega)\to L^q(\Omega)$ and a direction $h\in L^p(\Omega)$, the pointwise(!) product $F'(u)h$ has to be in $L^q(\Omega)$.
This leads to additional conditions on the superposition operator $F'$ defined by $f'$, which is known as \term[discrepancy, two-norm]{two-norm discrepancy}.
\begin{theorem}\label{thm:superpos:differentiable}
    Let $f:\Omega\times\R\to\R$ be a Carathéodory function that satisfies the growth condition \eqref{eq:superpos:growth} for $1\leq q<p<\infty$. If the partial derivative $f'_z$ is a Carathéodory function as well and satisfies \eqref{eq:superpos:growth} for $p'=p-q$, the superposition operator $F:L^p(\Omega)\to L^q(\Omega)$ is continuously differentiable, and its derivative in $u\in L^p(\Omega)$ in direction $h\in L^p(\Omega)$ is given by
    \begin{equation*}
        [F'(u)h](x) = f_z'(x,u(x))h(x) \qquad\text{for almost every }x\in\Omega.
    \end{equation*}
\end{theorem}
\begin{proof}
    \Cref{thm:superpos:continuous} yields that for $r \defeq \frac{pq}{p-q}$ (i.e., $\frac{r}{p} = \frac{q}{p'}$), the superposition operator
    \begin{equation*}
        G:L^p(\Omega)\to L^r(\Omega),\qquad [G(u)](x)=f'_z(x,u(x))\quad\text{for almost every }x\in\Omega,
    \end{equation*}
    is well-defined and continuous.
    The Hölder inequality further implies that for any $u\in L^p(\Omega)$,
    \begin{equation}\label{eq:superpos:hoelder}
        \norm{G(u)h}_{L^q} \leq \norm{G(u)}_{L^r}\norm{h}_{L^p}\qquad\text{for all }h\in L^p(\Omega),
    \end{equation}
    i.e., the pointwise multiplication $h\mapsto G(u)h$ defines a bounded linear operator $DF(u):L^p(\Omega)\to L^q(\Omega)$.

    Let now $h\in L^p(\Omega)$ be arbitrary. Since $z\mapsto f(x,z)$ is continuously differentiable by assumption, the classical mean value theorem together with the properties of the integral (in particular, monotonicity, Jensen's inequality on $[0,1]$, and Fubini's theorem) and \eqref{eq:superpos:hoelder} implies that
    \begin{equation*}
        \begin{multlined}[c][0.9\displaywidth]
            \norm{F(u+h)-F(u)-DF(u)h}_{L^q}\\
            \begin{aligned}[t]
                &= \left(\int_\Omega |f(x,u(x)+h(x)) - f(x,u(x)) - f'_z(x,u(x))h(x)|^q\,dx\right)^{\frac1q}\\
                &= \left(\int_\Omega \left|\int_0^1 f'_z(x,u(x)+th(x))h(x)\,dt  - f'_z(x,u(x))h(x)\right|^q\,dx\right)^{\frac1q}\\
                &\leq \left(\int_0^1\int_\Omega \left| \left(f'_z(x,u(x)+th(x)) - f'_z(x,u(x))\right)h(x)\right|^q\,dx\,dt\right)^{\frac1q}\\
                &= \int_0^1 \norm{(G(u+th)-G(u))h}_{L^q}\,dt\\
                &\leq \int_0^1 \norm{G(u+th)-G(u)}_{L^r}\,dt\ \norm{h}_{L^p}.
            \end{aligned}
        \end{multlined}
    \end{equation*}
    Due to the continuity of $G:L^p(\Omega)\to L^r(\Omega)$, the integrand tends to zero uniformly in $[0,1]$ for $\norm{h}_{L^p}\to 0$, and hence $F$ is by definition Fréchet differentiable with derivative $F'(u) = DF(u)$ (whose continuity we have already shown).
\end{proof}
In fact, this result is sharp: except for the case $p=q=\infty$, no superposition operator is differentiable from $L^p(\Omega)$ to $L^p(\Omega)$ (unless it is affine-linear); see, e.g., \cite[Theorem 3.12]{Appell:1990}.

\section{Variational principles}\label{sec:variation:ekeland}

As the example $f(t) = 1/t$ on $\{t\in \R:t\geq 1\}$ shows, the coercivity requirement in \cref{thm:variation:existence} is necessary to obtain minimizers even if the functional is bounded from below. However, sometimes one does not need an exact minimizer and is satisfied with \enquote{almost minimizers}. \term[principle!variational]{Variational principles} state that such almost minimizers can be obtained as minimizers of a perturbed functional and even give a precise relation between the size of the perturbation needed in terms of the desired distance from the infimum.

The most well-known variational principle is \term[principle!variational!Ekeland]{Ekeland's variational principle}, which holds in general complete metric spaces but which we here state in Banach spaces for the sake of notation. In the statement of the following theorem, note that we do not assume the functional to be \emph{weakly} lower semicontinuous.
\begin{theorem}[Ekeland's variational principle]\label{thm:variation:ekeland}
    Let $X$ be a Banach space and $F:X\to\Rbar$ be proper, lower semicontinuous, and bounded from below. Let $\eps>0$ and $z_\eps\in X$ be such that
    \begin{equation*}
        F(z_\eps) < \inf_{x\in X} F(x) + \eps.
    \end{equation*}
    Then for any $\lambda>0$, there exists an $x_\lambda \in X$ with
    \begin{enumerate}
        \item $\norm{x_\lambda-z_\eps}_X \leq \lambda$,
        \item $F(x_\lambda) +\frac\eps\lambda \norm{x_\lambda-z_\eps}_X \leq F(z_\eps)$,
        \item $F(x_\lambda) < F(x) + \frac\eps\lambda\norm{x-x_\lambda}_X$ for all $x\in X\setminus\{x_\lambda\}$.
    \end{enumerate}
\end{theorem}
\begin{proof}
    The proof proceeds similarly to that of \cref{thm:variation:existence}: We construct an \enquote{almost minimizing} sequence, show that it converges, and verify that the limit has the desired properties. Here we proceed inductively. First, set $x_0 \defeq z_\eps$. For given $x_n$, define now
    \begin{equation*}
        S_n \defeq \setof{x\in X}{ F(x) +\frac\eps\lambda \norm{x-x_n}_X \leq F(x_n)}.
    \end{equation*}
    Since $x_n\in S_n$, this set is nonempty. We can then choose $x_{n+1}\in S_n$ such that
    \begin{equation}\label{eq:variation:ekeland:induction}
        F(x_{n+1}) \leq \frac12 F(x_n) + \frac12 \inf_{x\in S_n} F(x).
    \end{equation}
    To see why this is possible, assume that $F(x_n) > \inf_{x\in S_n} F(x)$ (otherwise we can simply take $x_{n+1}=x_n\in S_n$ to satisfy the inequality with equality). Then the properties of the infimum imply that for $\eps:=\frac12(F(x_n)-\inf_{x\in S_n} F(x))>0$ we can find $x_{n+1}\in S_n$ such that $F(x_{n+1}) \leq \inf_{x\in S_n} F(x) + \eps$, from which \eqref{eq:variation:ekeland:induction} follows.
    By construction, the sequence $\{F(x_n)\}_{n\in\N}$ is thus decreasing as well as bounded from below and therefore convergent. Using the triangle inequality, the fact that $x_{n+1}\in S_n$, and the telescoping sum, we also obtain that for any $m\geq n\in \N$,
    \begin{equation}\label{eq:variation:ekeland:cauchy}
        \frac\eps\lambda \norm{x_n-x_m}_X \leq \sum_{j=n}^{m-1} \frac\eps\lambda  \norm{x_j-x_{j+1}}_X \leq F(x_n)-F(x_m).
    \end{equation}
    Hence, $\{x_n\}_{n\in\N}$ is a Cauchy sequence since $\{F(x_n)\}_{n\in\N}$ is one and hence converges to some $x_\lambda \in X$ since $X$ is complete.

    We now show that this limit has the claimed properties. We begin with (ii), for which we use the fact that both $F$ and the norm in $X$ are lower semicontinuous and hence obtain from \eqref{eq:variation:ekeland:cauchy} by taking $m\to\infty$ that
    \begin{equation}\label{eq:variation:ekeland:value}
        \frac\eps\lambda\norm{x_n -x_\lambda}_X + F(x_\lambda)
        \leq \limsup_{m\to \infty}
        \frac\eps\lambda\norm{x_n -x_m}_X + F(x_m)
        \leq F(x_n)  \quad\text{for any }n\geq 0.
    \end{equation}
    Choosing in particular $n=0$ such that $x_0 = z_\eps$ yields (ii).

    Furthermore, by definition of $z_\eps$, this implies that
    \begin{equation*}
        \frac\eps\lambda\norm{z_\eps-x_\lambda}_X \leq  F(z_\eps) - F(x_\lambda) \leq F(z_\eps) - \inf_{x\in X} F(x) < \eps
    \end{equation*}
    and hence (i) and (ii).

    Assume now that (iii) does not hold, i.e., that there exists an $x\in X\setminus\{x_\lambda\}$ such that
    \begin{equation}\label{eq:variation:ekeland:counter}
        F(x) \leq F(x_\lambda) - \frac\eps\lambda \norm{x-x_\lambda}_X <F(x_\lambda).
    \end{equation}
    Estimating $F(x_\lambda)$ using \eqref{eq:variation:ekeland:value} and then using the \enquote{productive zero} (i.e., adding and subtracting the same term) together with the triangle inequality, we obtain from the first inequality that for all $n\in \N$,
    \begin{equation*}
        F(x) \leq F(x_n) - \frac\eps\lambda\norm{x_n-x_\lambda}_X - \frac\eps\lambda \norm{x-x_\lambda}_X
        \leq F(x_n) - \frac\eps\lambda\norm{x_n-x}_X.
    \end{equation*}
    Hence, $x\in S_n$ for all $n\in \N$. From \eqref{eq:variation:ekeland:induction}, we then deduce that
    \begin{equation*}
        2 F(x_{n+1}) - F(x_n) \leq F(x) \quad\text{for all } n\in\N.
    \end{equation*}
    The convergence of $\{F(x_n)\}_{n\in\N}$ together with \eqref{eq:variation:ekeland:counter} and the lower semicontinuity of $F$ thus yields the contradiction
    \begin{equation*}
        \lim_{n\to\infty} F(x_n) \leq F(x) < F(x_\lambda) \leq \lim_{n\to\infty} F(x_n).
        \qedhere
    \end{equation*}
\end{proof}

Ekeland's variational principle has the disadvantage that even for differentiable $F$, the perturbed function that is minimized by $x_\lambda$ is inherently nonsmooth. This is different for \term[principle!variational!smooth]{smooth variational principles} such as the following one due to Borwein and Preiss \cite{borwein1987smooth}.

\begin{theorem}[Borwein--Preiss variational principle]\label{thm:variation:borweinpreiss}\index{principle!variational!Borwein--Preiss}
    Let $X$ be a Banach space and $F:X\to\Rbar$ be proper, lower semicontinuous, and bounded from below. Let $\eps>0$ and $z_\eps\in X$ be such that
    \begin{equation*}
        F(z_\eps) < \inf_{x\in X} F(x) + \eps.
    \end{equation*}
    Then for any $\lambda>0$ and $p\geq 1$, there exists
    \begin{itemize}
        \item a sequence $\{x_n\}_{n\in\N_0}\subset X$ with $x_0 = z_\eps$ converging strongly to some $x_\lambda\in X$ and
        \item a sequence $\{\mu_n\}_{n\in\N_0}\subset (0,\infty)$ with $\sum_{n=0}^\infty \mu_n = 1$
    \end{itemize}
    such that
    \begin{enumerate}
        \item\label{it:variation:borweinpreiss:1} $\norm{x_\lambda-x_n}_X \leq \lambda$ for all $n\in \N\cup\{0\}$,
        \item\label{it:variation:borweinpreiss:2} $F(x_\lambda) +\frac\eps{\lambda^p} \sum_{n=0}^\infty \mu_n \norm{x_\lambda-x_n}_X^p \leq F(z_\eps)$,
        \item\label{it:variation:borweinpreiss:3} $F(x_\lambda) +\frac\eps{\lambda^p} \sum_{n=0}^\infty \mu_n \norm{x_\lambda-x_n}_X^p \leq
            F(x) +\frac\eps{\lambda^p} \sum_{n=0}^\infty \mu_n \norm{x-x_n}_X^p$ for all $x\in X$.
    \end{enumerate}
\end{theorem}
\begin{proof}
    We proceed similarly to the proof of \cref{thm:variation:ekeland} by induction. First, we choose constants $\gamma,\eta,\mu,\theta>0$ such that
    \begin{itemize}
        \item $F(z_\eps) - \inf_{x\in X} F(x) < \eta < \gamma < \eps$,
        \item $\mu < 1-\frac\gamma\eps$,
        \item $\textstyle \theta <  \mu \left(1-\left(\frac\eta\gamma\right)^{1/p}\right)^p$.
    \end{itemize}
    Let now $x_0 \defeq z_\eps$ and $F_0 \defeq F$ and set $\delta \defeq  (1-\mu) \frac\eps{\lambda^p}>0$. We then define
    \begin{equation*}
        F_1(x) \defeq  F_0(x) + \delta \mu \norm{x-x_0}_X^p \qquad \text{for all }x\in X.
    \end{equation*}
    By construction, we then have
    \begin{equation*}
        \inf_{x\in X} F_1(x) \leq F_1(x_0) = F_0(x_0),
    \end{equation*}
    and thus we can find, by the same argument as for \eqref{eq:variation:ekeland:induction}, an $x_1\in X$ with
    \begin{equation*}
        F_1(x_1) \leq \theta F_0(x_0) + (1-\theta) \inf_{x\in X} F_1(x).
    \end{equation*}
    Continuing in this manner, we obtain sequences $\{x_n\}_{n\in\N}$ and $\{F_n\}_{n\in\N}$ with
    \begin{equation}\label{eq:variation:bp_rec1}
        F_{n+1}(x) = F_n(x) + \delta \mu^n \norm{x-x_n}_X^p
    \end{equation}
    and
    \begin{equation}\label{eq:variation:bp_rec2}
        F_{n+1}(x_{n+1}) \leq \theta F_n(x_n) + (1-\theta) \inf_{x\in X} F(x).
    \end{equation}
    Set now $s_n\defeq \inf_{x\in X} F_n(x)$ and $a_n \defeq F_n(x_n)$. Then \eqref{eq:variation:bp_rec1} implies that $\{s_n\}_{n\geq 0}$ is monotonically increasing, while \eqref{eq:variation:bp_rec2} implies that $\{a_n\}_{n\geq 0}$ is monotonically decreasing. We thus have
    \begin{equation}\label{eq:variation:bp_est1}
        s_n \leq s_{n+1} \leq a_{n+1} \leq \theta a_n + (1-\theta) s_{n+1} \leq a_n,
    \end{equation}
    which can be rearranged to show for all $n\geq 0$ that
    \begin{equation}\label{eq:variation:bp_est2}
        a_{n+1} - s_{n+1} \leq \theta a_n + (1-\theta) s_{n+1} - s_{n+1} = \theta(a_n - s_{n+1}) \leq \theta (a_n - s_n) \leq \theta^n(a_0-s_0).
    \end{equation}
    This together with the monotonicity of the two sequences and the boundedness of $F$ from below shows that $\lim_{n\to\infty} a_n = \lim_{n\to\infty} s_n \in \R$.
    We now use \eqref{eq:variation:bp_rec1} in \eqref{eq:variation:bp_est1} to obtain that
    \begin{equation*}
        a_n \geq a_{n+1} = F_n(x_n) + \delta\mu^n \norm{x_{n+1}-x_n}_X^p \geq s_n + \delta\mu^n \norm{x_{n+1}-x_n}_X^p,
    \end{equation*}
    which together with \eqref{eq:variation:bp_est2} and the choice of $\eta$ yields
    \begin{equation*}
        \delta \mu^n\norm{x_{n+1}-x_n}_X^p \leq a_n - s_n \leq \theta^n (a_0-s_0) < \eta \theta^n.
    \end{equation*}
    The choice of $\theta$ and $\mu$ now ensures that $0<\frac\theta\mu<1$, which implies that
    \begin{equation}\label{eq:variation:bp_bound}
        \begin{aligned}[t]
            \norm{x_m - x_n}_X &\leq \sum_{k=n}^{m-n-1}\norm{x_{k+1}-x_k}_X
            \leq \left(\frac\eta\delta\right)^{1/p} \sum_{k=n}^{m-n-1} \left(\frac\theta\mu\right)^{k/p}\\
            &\leq \left(\frac\eta\delta\right)^{1/p} \left(\frac\theta\mu\right)^{n/p} \left(1-\left(\frac\theta\mu\right)^{1/p}\right)^{-1}\quad\text{for all }m,n\geq 0
        \end{aligned}
    \end{equation}
    using the partial geometric series
    \begin{equation*}
        \sum_{k=n}^{m-n-1} \alpha^k = \sum_{k=0}^{m-n-1} \alpha^k - \sum_{k=0}^{n-1} \alpha^k
        = \frac{1-\alpha^{m-n}}{1-\alpha} - \frac{1-\alpha^n}{1-\alpha} < \frac{\alpha^n}{1-\alpha}
    \end{equation*}
    valid for any $\alpha\in (0,1)$.
    Hence $\{x_n\}_n\in\N$ is a Cauchy sequence which therefore converges to some $x_\lambda \in X$.
    Setting $\mu_n \defeq \mu^n (1-\mu)>0$, we also have $\sum_{n=0}^\infty \mu_n = 1$ by the choice of $\mu<1$.
    Furthermore, the definition of $\mu_n$ and $\delta$ implies for all $x\in X$ that
    \begin{equation}\label{eq:variation:bp_limit}
        F(x) + \frac\eps{\lambda^p}\sum_{k=0}^\infty\mu_k\norm{x-x_k}_X^p =
        \lim_{n\to \infty} F(x) + \sum_{k=0}^n \delta \mu^k \norm{x-x_k}_X^p = \lim_{n\to\infty} F_n(x).
    \end{equation}

    It remains to verify the claims on $x_\lambda$. First, \eqref{eq:variation:bp_bound} together with the choice of $\theta$ and $\delta$ implies for all $n,m\geq 0$ that
    \begin{equation*}
        \norm{x_m - x_n}_X \leq \left(\frac\eta\delta\right)^{1/p} \left(\frac\eta\gamma\right)^{-1/p} = \left(\frac\gamma\delta\right)^{1/p} < \left(\frac\eps\delta\right)^{1/p} (1-\mu)^{1/p} = \lambda.
    \end{equation*}
    Letting $m\to \infty$ for fixed $n\in\N\cup\{0\}$ now shows \ref{it:variation:borweinpreiss:1}.

    Second, by \eqref{eq:variation:bp_rec1} and the definition of $\delta$, we have
    \begin{equation*}
        F(x_n) + \frac\eps{\lambda^p}\sum_{k=0}^\infty\mu_k\norm{x_n-x_k}_X^p
        =
        F_n(x_n) + \frac\eps{\lambda^p}\sum_{k=n+1}^\infty\mu_k\norm{x_n-x_k}_X^p
        \leq
        a_n + \eps\sum_{k=n+1}^\infty \mu_k,
    \end{equation*}
    where the inequality follows from \ref{it:variation:borweinpreiss:1}.
    The lower semicontinuity of $F$ and of the norm thus yield
    \begin{equation}\label{eq:variation:bp_claim2}
        F(x_\lambda) +\frac\eps{\lambda^p}\sum_{k=0}^\infty\mu_k\norm{x_\lambda-x_k}_X^p
        \leq \lim_{n\to\infty} a_n \leq a_0 = F(z_\eps)
    \end{equation}
    since $\{a_n\}_{n\geq 0}$ is monotonically decreasing. This shows \ref{it:variation:borweinpreiss:2}.

    Finally, \eqref{eq:variation:bp_limit} and the definition of $s_n$ imply for all $x\in X$ that
    \begin{equation*}
        F(x) + \frac\eps{\lambda^p}\sum_{k=0}^\infty\mu_k\norm{x-x_k}_X^p = \lim_{n\to \infty} F_n(x) \geq \lim_{n\to \infty} s_n = \lim_{n\to\infty} a_n,
    \end{equation*}
    which together with \eqref{eq:variation:bp_claim2} yields \ref{it:variation:borweinpreiss:3}.
\end{proof}
The Borwein--Preiss variational principle therefore guarantees a smooth perturbation if, e.g., $X$ is a Hilbert space and $p=2$.
Further smooth variational principles that allow for more general smooth perturbations such as the \term[principle!variational!Deville--Godefroy--Zizler]{Deville--Godefroy--Zizler variational principle} can be found in, e.g., \cite{BorweinZhu:2005,Schirotzek:2007}.

\section{Variational convergence}\label{sec:variation:convergence}

Often one is not only interested in minimizing a single functional but a sequence of functionals. These could involve smooth approximations of nonsmooth functions, penalizations of constraints, finite-dimensional discretizations, or sample average approximations of stochastic functionals. The central question is then under which conditions minimizers of the approximating functionals converge to a minimizer of the original functional. A simple example shows that pointwise convergence of the functionals is not enough.
\begin{example}
    Consider for $n\in\N$
    \begin{equation*}
        F_n :\R\to\R,\qquad F_n(x) = \frac1n|x-n| - 1.
    \end{equation*}
    Then we obviously have
    \begin{equation*}
        \lim_{n\to\infty} F_n(x) = 0\qquad\text{for all }x\in\R
    \end{equation*}
    and hence $F_n\to 0=:F$ pointwise. But
    \begin{equation*}
        \lim_{n\to\infty} \inf_{x\in\R} F_n(x) = \lim_{n\to\infty} F_n(n) = -1 \neq 0 = \inf_{x\in \R} F(x).
    \end{equation*}
\end{example}

This shows that the convergence has to be uniform in some appropriate sense. This leads to the notion of \term[convergence!$\Gamma$-]{$\Gamma$-convergence}: A sequence $\{F_n\}_{n\in\N}$ of functionals $F_n:X\to \Rbar$ \emph{$\Gamma$-converges} to some $F:X\to \Rbar$ if for all $x\in X$,
\begin{enumerate}
    \item\label{it:variation:gamma-liminf} for all sequences $\{x_n\}_{n\in\N}\subset X$ with $x_n\to x$,
        \begin{equation*}
            \liminf_{n\to\infty} F_n(x_n) \geq F(x);
        \end{equation*}
    \item\label{it:variation:gamma-recovery} there exists a \term[sequence!recovery]{recovery sequence} $\{x_n\}_{n\in\N}\subset X$ such that $x_n\to x$ and
        \begin{equation*}
            \lim_{n\to\infty} F_n(x_n) = F(x).
        \end{equation*}
\end{enumerate}
In this case, we call $F$ the \term[limit!$\Gamma$-]{$\Gamma$-limit} of the $F_n$ and write $F_n\to_\Gamma F$. In view of (i), it suffices to show that the recovery sequence satisfies
\begin{enumerate}
    \item[(ii$'$)] $x_n\to x$ and
        \begin{equation*}
            \limsup_{n\to\infty} F_n(x_n) \leq F(x).
        \end{equation*}
\end{enumerate}
Analogously, we can define \emph{weak $\Gamma$-convergence} using weakly convergent sequences $x_n\weakto x$, denoted by $F_n\weakto_\Gamma F$.
\begin{remark}
    In the context of variational analysis, $\Gamma$-convergence is also referred to as \term[convergence!epigraphical]{epigraphical convergence}, since $F_n \to_\Gamma F$ if and only if $\epi F_n$ converges to $\epi F$ in the sense that the inner and outer limits (as later defined in \cref{sec:monotone:basic}) of $\{\epi F_n\}_{n\in\N}$ coincide and are equal to $\epi F$; see, e.g., \cite[Chapter 7]{Rockafellar:1998}.
\end{remark}

The $\Gamma$-limit has the following useful property.
\begin{lemma}\label{lem:variation:gamma-lsc}
    If $F_n\to_\Gamma F$, then $F$ is lower semicontinuous.
\end{lemma}
\begin{proof}
    Let $x\in X$ be arbitrary and $\{x_n\}_{n\in\N}\subset X$ with $x_n\to x$.
    By assumption \ref{it:variation:gamma-recovery}, we can then find for every $x_n$ a recovery sequence $\{x_{n,k}\}_{k\in\N}$ such that $x_{n,k}\to x_n$ and $F_k(x_{n,k})\to F(x_n)$.
    Hence for every $n\in\N$ we can choose $k(n)\in\N$ such that
    \begin{equation*}
        \norm{x_{n,k(n)}-x_n}_X\leq \frac1n\quad\text{ and }\quad
        |F_{k(n)}(x_{n,k(n)}) - F(x_n)|\leq \frac1n.
    \end{equation*}
    Without loss of generality, we can assume that the sequence $\{k(n)\}_{n\in\N}$ is strictly increasing.
    We now define the \enquote{step sequence} $\{x_m\}_{m\in \N}$ via
    \begin{equation*}
        x_m \defeq x_{n,k(n)}\qquad\text{for } k(n)\leq m <k(n+1).
    \end{equation*}
    Then we also have $\norm{x_{m}-x_n}_X\leq \frac1n$  and hence $x_m\to x$ since $x_n\to x$.
    We can thus apply the $\liminf$-inequality \ref{it:variation:gamma-liminf} to deduce that
    \begin{equation*}
        F(x) \leq \liminf_{m\to\infty} F_m(x_m) \leq \liminf_{n\to\infty} F_{k(n)}(x_{n,k(n)})
        = \liminf_{n\to\infty} F(u_n),
    \end{equation*}
    where we have passed to the subsequence of steps (whose $\liminf$ cannot be smaller) in the second step and used the choice of $k(n)$ in the last step.
    Hence $F$ is lower semicontinuous.
\end{proof}
This implies that even a constant sequence may not $\Gamma$-converge (choose any functional that is not lower semicontinuous)! Similarly, the sum of two $\Gamma$-converging sequences may not $\Gamma$-converge to the sum of their $\Gamma$-limits. However, it is easy to show that $\Gamma$-convergence is stable under \emph{continuous} perturbations.
\begin{lemma}\label{lem:variation:gamma-perturbation}
    Let $F_n:X\to\Rbar$, $n\in\N$, and $F:X\to\Rbar$. If $F_n\to_\Gamma F$ and $G:X\to \Rbar$ is continuous, then $F_n+G\to_\Gamma F+G$.
\end{lemma}
\begin{proof}
    We simply verify the defining properties.

    \emph{\ref{it:variation:gamma-liminf}:} For any sequence $\{x_n\}_{n\in\N}$ with $x_n\to x$, we have by assumption that
    \begin{equation*}
        \liminf_{n\to\infty}(F_n + G)(x_n) = \liminf_{n\to\infty} F_n(x_n) + \lim_{n\to\infty} G(x_n) \geq F(x)+ G(x).
    \end{equation*}
    \emph{\ref{it:variation:gamma-recovery}:} For $x\in X$, let $\{x_n\}_{n\in\N}$ be a recovery sequence with respect to $\{F_n\}_{n\in\N}$. Then it follows from $x_n\to x$ and the continuity of $G$ that
    \begin{equation*}
        \lim_{n\to\infty}(F_n +G)(x_n) = \lim_{n\to\infty} F_n(x_n) + \lim_{n\to\infty} G(x_n) = F(x) + G(x).
        \qedhere
    \end{equation*}
\end{proof}

We now come to the central result linking $\Gamma$-convergence of functionals to convergence of (almost) minimizers.
\begin{theorem}\label{thm:variation:gamma-min}
    Let $F_n:X\to\Rbar$, $n\in\N$, and $F:X\to\Rbar$. If $F_n\to_\Gamma F$, then
    \begin{equation*}
        \limsup_{n\to\infty} \inf_{x\in X} F_n(x) \leq \inf_{x\in X} F(x).
    \end{equation*}
    Furthermore, if there exists a null sequence $\{\eps_n\}_{n\in\N}$ and a sequence $\{x_n\}_{n\in\N}\subset X$ such that $x_n\to \bar x\in X$ and
    \begin{equation*}
        F_n(x_n) \leq \inf_{x\in X} F_n(x) + \eps_n \qquad\text{for all }n\in\N,
    \end{equation*}
    then $\bar x\in X$ is a minimizer of $F$ and
    \begin{equation*}
        \lim_{n\to\infty}\inf_{x\in X} F_n(x) = \min_{x\in X} F(x).
    \end{equation*}
\end{theorem}
\begin{proof}
    Let $x\in X$ be arbitrary and $\{x_n\}_{n\in\N}$ be a recovery sequence for $x$. Then the $\liminf$ inequality \ref{it:variation:gamma-liminf} immediately implies
    \begin{equation*}
        F(x) \geq \limsup_{n\to\infty} F_n (x_n) \geq \limsup_{n\to\infty} \inf_{x\in X} F_n (x).
    \end{equation*}
    Taking the infimum over all $x\in X$ then yields the first claim.

    Let now $\{x_n\}_{n\in\N}$ be a sequence of $\eps_n$-minimizers converging to $\bar x\in X$. Combining the definition of the $x_n$ with the $\liminf$ inequality \ref{it:variation:gamma-liminf} then yields
    \begin{equation*}
        \liminf_{n\to \infty} \inf_{x\in X} F_n(x) =
        \liminf_{n\to \infty} (\inf_{x\in X} F_n(x)+\eps_n) \geq
        \liminf_{n\to \infty}  F_n(x_n) \geq
        F(\bar x) \geq \inf_{x\in X} F(x).
    \end{equation*}
    Together with the first claim, this implies the claimed equality as well as $F(\bar x) = \inf_{x\in X} F(x)$.
\end{proof}
The same result holds for weak $\Gamma$-convergence if the sequence of almost minimizers converges weakly.

This begs the question under which conditions such a sequence of almost minimizers converges. This requires some sort of compactness assumption. In analogy to the direct method (\cref{thm:variation:existence}), we assume here that the sequence $\{F_n\}_{n\in\N}$ is \term[sequence!equicoercive]{equicoercive}, meaning that there exists a coercive functional $G:X\to \R$ such that $F_n(x)\geq G(x)$ for all $x\in X$.
\begin{theorem}\label{thm:variation:gamma-convergence}
    Let $X$ be a reflexive Banach space and $F:X\to \Rbar$ be proper. Let $\{F_n\}_{n\in\N}$ be equicoercive such that $F_n\weakto_\Gamma F$ and $\{\eps_n\}_{n\in\N}$ be a null sequence. Then
    every sequence $\{x_n\}_{n\in\N}$ of $\eps_n$-minimizers of $F_n$ has a weak accumulation point, and
    every weak accumulation point of such a sequence is a minimizer of $F$.
    Furthermore,
    \begin{equation*}
        \lim_{n\to\infty}\inf_{x\in X}F_n(x) = \min_{x\in X} F(x).
    \end{equation*}

    Conversely, every minimizer of $F$ is the limit of a sequence of $\eps_n$-minimizers of $F_n$.
\end{theorem}
\begin{proof}
    Since $F$ is proper, there exists $z\in X$ such that $F(z)<\infty$. Let now $\{z_n\}_{n\in\N}$ be a recovery sequence for $z$ and $\{x_n\}_{n\in\N}$ be a sequence of $\eps_n$-minimizers. Then we have that
    \begin{equation*}
        F_n(x_n) = \inf_{x\in X} F_n(x_n) + \eps_n \leq F_n(z_n) + \eps_n \to F(z)<\infty.
    \end{equation*}
    Hence $\{F_n(x_n)\}_{n\in\N}$ is bounded from above. By equicoercivity, this implies that $\{x_n\}_{n\in\N}$ is bounded, since otherwise we would have $F_n(x_n) \geq G(x_n) \to \infty$. Hence there exists a weakly converging subsequence. Since it is easy to see that every subsequence of a $\Gamma$-converging sequence $\Gamma$-converges to the same limit, we can apply \cref{thm:variation:gamma-min} to every weakly converging subsequence to obtain the first claim.

    Choose now $\{x_n\}_{n\in\N}$ such that
    \begin{equation*}
        F_n(x_n) \leq \inf_{x\in X} F_n(x) + \frac1n,
    \end{equation*}
    i.e., $\{x_n\}_{n\in\N}$ is a sequence of $\frac1n$-minimizers of $F_n$.
    Then we can apply the first claim to any subsequence to obtain a further subsequence -- which we do not relabel -- such that $F_{n}(x_n)\to\min_{x\in X} F(x)$. A subsequence--subsequence argument then shows that the whole sequence converges, i.e., that
    \begin{equation*}
        \lim_{n\to\infty} \inf F_n(x) = \lim_{n\to\infty} F_n(x_n) = \min_{x\in X} F(x).
    \end{equation*}

    Let now $\bar x\in X$ be a minimizer of $F$ and $\{x_n\}_{n\in\N}$ be a recovery sequence for $\bar x$. By the previous claim, we then have
    \begin{equation*}
        \lim_{n\to\infty} F_n(x_n) = F(\bar x) = \min_{x\in X} F(x) = \lim_{n\to\infty} \inf_{x\in X} F_n(x).
        \qedhere
    \end{equation*}
\end{proof}

\begin{remark}
    It is possible to define $\Gamma$-convergence in abstract metric and even topological spaces. For details, further properties, and applications to problems in the calculus of variations, we refer to \cite{Attouch:1984,DalMaso:1993,Braides:2002}.
\end{remark}

\part{Convex analysis}\label{part:convex}

\chapter{Convex functions}\label{chap:convex}

Now that we know from the direct method of the calculus of variations when a functional $F:X\to\Rbar\defeq \R\cup \{\infty\}$ admits a minimizer $\bar x\in X$, our next goal is to characterize such minimizers using optimality conditions, i.e., without comparing its function value to that at every other point. If $F$ is differentiable at $\bar x$, the classical optimality condition is by Fermat's principle, $F'(\bar x)=0$, and we can use calculus rules to evaluate this derivative in order to make this condition as explicit as possible. We wish to extend this as far as possible to nonsmooth $F$, i.e., not classically (Fréchet or Gateaux) differentiable. Clearly, \emph{not} being differentiable is not much to work with, so we have to assume other properties instead. One possibility is to replace the \emph{local} property of differentiability with the \emph{global} property of convexity. As we will see in this and the following chapters, this property will allow us to recover a satisfying calculus for a class of relevant nonsmooth functionals. We begin in this chapter with deriving several fundamental properties of convex functions relevant for optimization, while the corresponding Fermat principle and calculus rules are the topic of the next \cref{chap:subdiff}.

Throughout this and the following chapters, $X$ will be a normed vector space unless noted otherwise.

\section{Definition and basic properties}

A functional $F:X\to\Rbar$ is called \term[functional!convex]{convex} if for all $x,y\in X$ and $\lambda\in [0,1]$, it holds that
\begin{equation}\label{eq:convex:def}
    F(\lambda x + (1-\lambda)y)\leq \lambda F(x) + (1-\lambda)F(y)
\end{equation}
(where the function value $\infty$ is allowed on both sides).
If for all $x,y\in\dom F$ with $x\neq y$ and all $\lambda \in (0,1)$ we even have
\begin{equation*}
    F(\lambda x + (1-\lambda)y)< \lambda F(x) + (1-\lambda)F(y),
\end{equation*}
we call $F$ \term[functional!convex!strictly]{strictly convex}.

\begin{figure}
    \centering
    \begin{subfigure}[t]{0.325\textwidth}
        \centering
        \begin{asy}
            path g=(0,0)..(1,-0.7)..(2.2,0.2)..(1,0.8)..cycle;
            fill(g, lightfill);
            draw(g);
            pair x=(0.1,0.1);
            pair y=(1.1,-0.5);
            draw(x--y, defaultpen+linewidth(1.1));
            dot(x);
            dot(y);
            label("$x$", x, NE);
            label("$y$", y, NE);
        \end{asy}
        \caption{a convex set $C$}
    \end{subfigure}
    \hfil
    \begin{subfigure}[t]{0.325\textwidth}
        \centering
        \begin{asy}
            real f(real x){ return x^2; };
            path g=graph(f, -1, 1);
            fill(g--(1, 0.5+f(1))--(-1, 0.3+f(-1))--cycle, lightfill);
            label("$F$", midpoint(g), S);
            label("$\epi F$", midpoint(g), 5*N);
            draw(g, defaultpen);
            path l=pt(f, -.1)--pt(f, .7);
            draw(l, defaultpen+linewidth(1.1));
            dot(point(l, 0));
            dot(point(l, 1));
        \end{asy}
        \caption{the (convex) epigraph of a convex function}
    \end{subfigure}
    \caption{Illustration of a convex set and of the characterization of a convex function in terms of the convexity of its epigraph: all line segments between two points of corresponding set are completely contained in that set.}
    \label{fig:convex:illustration}
\end{figure}

As illustrated in \cref{fig:convex:illustration}, an alternative characterization of the convexity of a functional $F:X\to\Rbar$ is based on its \term{epigraph}
\begin{equation*}
    \epi F \defeq\setof{(x,t)\in X\times \R}{F(x)\leq t}.
\end{equation*}
\begin{lemma}\label{lem:convex:epi}
    Let $F:X\to\Rbar$. Then $\epi F$ is
    \begin{enumerate}
        \item \label{lem:convex:epi:i}
            nonempty if and only if $F$ is proper;
        \item \label{lem:convex:epi:ii}
            convex if and only if $F$ is convex;
        \item \label{lem:convex:epi:iii}
            (weakly) closed if and only if $F$ is (weakly) lower semicontinuous.%
            \footnote{For that reason, some authors use the term \term[functional!closed|infn]{closed} also to refer to lower semicontinuous functionals. We will stick with the latter, much less ambiguous, term throughout the following.}
    \end{enumerate}
\end{lemma}
\begin{proof}
    Statement \ref{lem:convex:epi:i} follows directly from the definition: $F$ is proper if and only if there exists an $x\in \dom F$, i.e., $(x,F(x))\in \epi F$.

    For \ref{lem:convex:epi:ii}, let $F$ be convex and $(x,r),(y,s)\in \epi F$ be given. For any $\lambda\in[0,1]$, the definition \eqref{eq:convex:def} then implies that
    \begin{equation*}
        F(\lambda x + (1-\lambda)y)\leq \lambda F(x) + (1-\lambda)F(y) \leq \lambda r + (1-\lambda)s,
    \end{equation*}
    i.e., that
    \begin{equation*}
        \lambda(x,r) + (1-\lambda)(y,s) = (\lambda x + (1-\lambda)y,\lambda r + (1-\lambda)s) \in \epi F,
    \end{equation*}
    and hence $\epi F$ is convex.
    Let conversely $\epi F$ be convex and $x,y\in X$ be arbitrary, where we can assume that $F(x)<\infty$ and $F(y)<\infty$ (otherwise \eqref{eq:convex:def} is trivially satisfied). We  clearly have $(x,F(x)),(y,F(y))\in\epi F$. The convexity of $\epi F$ then implies for all  $\lambda\in[0,1]$ that
    \begin{equation*}
        (\lambda x + (1-\lambda)y,\lambda F(x) + (1-\lambda)F(y)) = \lambda(x,F(x)) + (1-\lambda)(y,F(y)) \in \epi F,
    \end{equation*}
    and hence by definition of $\epi F$ that \eqref{eq:convex:def} holds.

    Finally, we show \ref{lem:convex:epi:iii}: Let first $F$ be lower semicontinuous, and let $\{(x_n,t_n)\}_{n\in\N}\subset \epi F$ be an arbitrary sequence with $(x_n,t_n)\to(x,t)\in X\times \R$. Then we have that
    \begin{equation*}
        F(x)\leq \liminf_{n\to\infty} F(x_n) \leq  \lim_{n\to\infty} t_n = t,
    \end{equation*}
    i.e., $(x,t)\in \epi F$. Let conversely $\epi F$ be closed and assume that $F$ is proper (otherwise the claim holds trivially) and not lower semicontinuous. Then there exists a sequence $\{x_n\}_{n\in\N}\subset X$ with $x_n\to x\in X$ and
    \begin{equation*}
        F(x) > \liminf_{n\to\infty} F(x_n) =: M \in [-\infty,\infty).
    \end{equation*}
    We now distinguish two cases.
    \begin{enumerate}[label=\alph*)]
        \item $x\in \dom F$: In this case, we can select a subsequence, again denoted by  $\{x_n\}_{n\in\N}$, such that there exists an $\eps>0$ with $F(x_n) \leq F(x)-\eps$ and thus $(x_n,F(x)-\eps)\in \epi F$ for all $n\in \N$.
            From $x_n\to x$ and the closedness of $\epi F$, we deduce that $(x,F(x)-\eps)\in \epi F$ and hence $F(x) \leq F(x) - \eps$, contradicting $\eps>0$.

        \item $x\not\in \dom F$: In this case, we can argue similarly using $F(x_n) \leq M+\eps$ for $M>-\infty$ or $F(x_n) \leq \eps$ for $M=-\infty$ to obtain a contradiction with $F(x)=\infty$.
    \end{enumerate}
    The equivalence of weak lower semicontinuity and weak closedness follows in exactly the same way.
\end{proof}
Note that $(x,t)\in \epi F$ implies that $x\in \dom F$; hence the effective domain of a proper, convex, and lower semicontinuous functional is always nonempty and convex as well.
Also, together with \cref{lem:convex_closed} we immediately obtain
\begin{corollary}\label{cor:convex:uhs}
    Let $F:X\to\Rbar$ be convex. Then $F$ is weakly lower semicontinuous if and only $F$ is lower semicontinuous.
\end{corollary}

Also useful for the study of a functional $F:X\to\Rbar$ are the corresponding \term[set!sublevel]{sublevel sets}
\begin{equation*}
    \sub_t F \defeq \setof{x\in X}{F(x)\leq t},\qquad t \in\R,
\end{equation*}
for which one shows as in \cref{lem:convex:epi} the following properties.
\begin{lemma}\label{lem:convex:sublevel}
    Let $F:X\to\Rbar$.
    \begin{enumerate}
        \item If $F$ is convex, $\sub_t F$ is convex for all $t\in\R$ (but the converse does not hold).
        \item $F$ is (weakly) lower semicontinuous if and only if $\sub_t F$ is (weakly) closed  for all $t\in\R$.
    \end{enumerate}
\end{lemma}

Directly from the definition we obtain the convexity of
\begin{enumerate}
    \item \term[functional!continuous affine]{continuous affine functionals} of the form $x\mapsto \dual{x^*,x}_X - \alpha$ for fixed $x^*\in X^*$ and $\alpha\in\R$;
    \item the norm $\norm{\cdot}_X$ in a normed vector space $X$;
    \item the indicator function $\delta_C$ for a convex set $C$.
\end{enumerate}
If $X$ is a Hilbert space, $F(x) = \norm{x}_X^2$ is even strictly convex: For $x,y\in X$ with $x\neq y$ and any $\lambda \in (0,1)$,
\begin{equation*}
    \begin{aligned}
        \norm{\lambda x+(1-\lambda )y}_X^2 &= \iprod{\lambda x+(1-\lambda )y}{\lambda x+(1-\lambda )y}_X \\
                                           &= \lambda ^2\iprod{x}{x}_X + 2\lambda (1-\lambda )\iprod{x}{y}_X+(1-\lambda )^2\iprod{y}{y}_X\\
                                           &=\lambda \Big(\lambda \iprod{x}{x}_X-(1-\lambda )\iprod{x-y}{x}_X+(1-\lambda )\iprod{y}{y}_X\Big)\\
        \MoveEqLeft[-1]
        +(1-\lambda )\Big(\lambda \iprod{x}{x}_X+\lambda \iprod{x-y}{y}_X+(1-\lambda )\iprod{y}{y}_X\Big)\\
        &=(\lambda +(1-\lambda ))\Big(\lambda \iprod{x}{x}_X+(1-\lambda )\iprod{y}{y}_X\Big)-\lambda (1-\lambda )\iprod{x-y}{x-y}_X\\
        &=\lambda \norm{x}_X^2 +(1-\lambda )\norm{y}_X^2-\lambda (1-\lambda )\norm{x-y}_X^2\\
        &<\lambda \norm{x}_X^2 +(1-\lambda )\norm{y}_X^2.
    \end{aligned}
\end{equation*}

Further examples can be constructed as in \cref{lem:variation:wlsc} through the following operations.
\begin{lemma}\label{lem:convex:func}
    Let $X$ and $Y$ be normed vector spaces and let $F:X\to\Rbar$ be convex. Then the following functionals are convex as well:
    \begin{enumerate}
        \item \label{lem:convex:func:i}
            $\alpha F$ for all $\alpha \geq 0$;
        \item \label{lem:convex:func:ii}
            $F+G$ for $G:X\to\Rbar$ convex (strictly if $F$ \emph{or} $G$ is strictly convex);
        \item \label{lem:convex:func:iii}
            $\phi\circ F$ for $\phi:\Rbar\to\Rbar$ convex and increasing;
        \item \label{lem:convex:func:iv}
            $F\circ K$ for $K:Y\to X$ linear;
        \item \label{lem:convex:func:v}
            $x\mapsto \sup_{i\in I} F_i(x)$ with $F_i:X\to\Rbar$ convex for an arbitrary set $I$.
    \end{enumerate}
\end{lemma}

\Cref{lem:convex:func}\,\ref{lem:convex:func:v} in particular implies that the pointwise supremum of continuous affine functionals is always convex. In fact, any convex functional can be written in this way. To show this, we define for a proper functional $F:X\to\Rbar$ the \term[envelope!convex]{convex envelope}
\begin{equation*}
    F^\Gamma(x)\defeq \sup\setof{a(x)}{a \text{ continuous affine with } a(\tilde x)\leq F(\tilde x)\text{ for all }\tilde x\in X}.
\end{equation*}
Note that $F^\Gamma$ could take the value $-\infty$ without further assumptions on $F$.
\begin{lemma}\label{lem:convex:gamma}
    Let $F:X\to\Rbar$ be proper. Then $F$ is convex and lower semicontinuous if and only if $F=F^\Gamma$.
\end{lemma}
\begin{proof}
    Since affine functionals are convex and continuous, \cref{lem:convex:func}\,\ref{lem:convex:func:v} and \cref{lem:variation:wlsc}\,\ref{lem:variation:wlsc:v} imply that $F^\Gamma$ is always convex and lower semicontinuous. Hence if $F=F^\Gamma$, the same obviously holds for $F$.

    Conversely, let $F:X\to\Rbar$ be proper, convex, and lower semicontinuous.
    It is clear from the definition of $F^\Gamma$ as a pointwise supremum that $F^\Gamma\leq F$ always holds. Assume therefore that $F^\Gamma\neq F$. Then there exists an $x_0\in X$ and a $\lambda\in\R$ with
    \begin{equation*}
        F^\Gamma(x_0) < \lambda < F(x_0).
    \end{equation*}
    We now use the Hahn--Banach separation theorem to construct a continuous affine functional $a$ with $a\leq F$ but $a(x_0)> \lambda>F^\Gamma(x_0)$, which would contradict the definition of $F^\Gamma$.
    Since $F$ is proper, convex, and lower semicontinuous, $\epi F$ is nonempty, convex, and closed by \cref{lem:convex:epi}. Furthermore, $\{(x_0,\lambda)\}$ is compact and, as $\lambda < F(x_0)$, disjoint with $\epi F$. \Cref{thm:functan:hb-separation}\,\ref{item:functan:hb-separation:ii} hence yields a $z^*\in (X\times \R)^*$ and an $\alpha\in \R$ with
    \begin{equation*}
        \dual{z^*,(x,t)}_{X\times \R} \leq \alpha <  \dual{z^*,(x_0,\lambda)}_{X\times \R} \qquad\text{for all }(x,t)\in\epi F.
    \end{equation*}
    We now define an $x^*\in X^*$ via $\dual{x^*,x}_X = \dual{z^*,(x,0)}_{X\times \R}$ for all $x\in X$ and set $s\defeq\dual{z^*,(0,1)}_{X\times \R}\in \R$. Then $\dual{z^*,(x,t)}_{X\times \R} =\dual{x^*,x}_X + st$ and hence
    \begin{equation}\label{eq:convex:func1}
        \dual{x^*,x}_X + st  \leq \alpha < \dual{x^*,x_0}_X +s\lambda  \qquad\text{for all } (x,t) \in\epi F.
    \end{equation}
    Now for $(x,t)\in \epi F$ we also have $(x,t')\in \epi F$ for all $t'>t$, and the first inequality in \eqref{eq:convex:func1} implies that for all sufficiently large $t>0$,
    \begin{equation*}
        s \leq \frac{\alpha - \dual{x^*,x}_X}{t} \to 0 \qquad\text{for }t\to\infty.
    \end{equation*}
    Hence $s\leq 0$. We continue with a case distinction.
    \begin{enumerate}
        \item \label{lem:convex:gamma:i} $s<0$:
            We set
            \begin{equation*}
                a:X\to\R,\qquad x\mapsto  \frac{\alpha - \dual{x^*,x}_X}{s},
            \end{equation*}
            which is continuous affine. Furthermore, using the productive zero in the first inequality in \eqref{eq:convex:func1} for $(x,F(x))\in \epi F$ implies (noting $s<0$!) that
            \begin{equation*}
                a(x) = \tfrac1s \left(\alpha - \dual{x^*,x}_X -sF(x)\right)+F(x) \leq F(x).
            \end{equation*}
            (For $x\notin \dom F$ this holds trivially.) But the second inequality in \eqref{eq:convex:func1} implies that
            \begin{equation*}
                a(x_0) = \tfrac1s\left(\alpha - \dual{x^*,x_0}_X\right) > \lambda.
            \end{equation*}
        \item \label{lem:convex:gamma:ii} $s=0$:
            Then $\dual{x^*,x}_X\leq \alpha<\dual{x^*,x_0}_X$ for all $x\in\dom F$, which can only hold for $x_0\notin \dom F$. But $F$ is proper, and hence we can find a $y_0\in \dom F$, for which we can construct as in case \ref{lem:convex:gamma:i} by separating $\epi F$ and $(y_0,\mu)$ for sufficiently small $\mu$ a continuous affine functional $a_0:X\to\R$ with $a_0 \leq F$ pointwise. For $\rho>0$ we now set
            \begin{equation*}
                a_\rho:X\to\R, \qquad x\mapsto a_0(x) + \rho\left( \dual{x^*,x}_X-\alpha\right),
            \end{equation*}
            which is continuous affine as well. Since $\dual{x^*,x}_X\leq \alpha$, we also have that $a_\rho(x) \leq a_0(x) \leq F(x)$ for all $x\in\dom F$ and any $\rho >0$. But due to $\dual{x^*,x_0}_X>\alpha$, we can choose $\rho>0$ with $a_\rho(x_0)>\lambda$.
    \end{enumerate}
    In both cases, the definition of $F^\Gamma$ as a supremum implies that $F^\Gamma(x_0)> \lambda$ as well, contradicting the assumption $F^\Gamma(x_0) <\lambda$.
\end{proof}

\begin{remark}\label{rem:convex:weak-star-hull}
    Using the weak-$*$ Hahn--Banach \cref{thm:clarke:hb} in place of \cref{thm:functan:hb-separation}, the same proof shows that a proper functional $F:X^*\to \Rbar$ is convex and weakly-$*$ lower semicontinuous if and only if $F=F_\Gamma$ for
    \begin{equation*}
        F_\Gamma(x^*) \defeq \sup \setof{\dual{x^*,x}_X + \alpha}{x\in X, \alpha\in\R, \dual{\tilde x^*,x}_X + \alpha \leq F(\tilde x^*) \text{ for all }\tilde x^* \in X^*}.
    \end{equation*}
    (Note that a convex and weakly lower semicontinuous functional need not be weakly-$*$ lower semicontinuous, since convex and closed sets need not be weakly-$*$ closed.)
\end{remark}

A particularly useful class of convex functionals in the calculus of variations arises from integral functionals with convex integrands defined through superposition operators.
\begin{lemma}\label{lem:lebesgue:lsc}
    Let $f:\R\to\Rbar$ be proper, convex, and lower semicontinuous. If $\Omega\subset\R^d$ is bounded and $1\leq p\leq \infty$, this also holds for
    \begin{equation*}
        F:L^p(\Omega)\to\Rbar,\qquad u\mapsto
        \begin{cases}
            \int_\Omega f(u(x))\,dx &\text{if }f\circ u \in L^1(\Omega),\\
            \infty &\text{else.}
        \end{cases}
    \end{equation*}
\end{lemma}
\begin{proof}
    First, \cref{lem:convex:gamma} implies that there exist $a,\alpha \in\R$ such that
    \begin{equation}\label{eq:lebesgue:lsc_bound}
        f(t) \geq at - \alpha \qquad\text{for all }t\in\R.
    \end{equation}
    Since $\Omega$ is bounded and hence $L^p(\Omega)\subset L^1(\Omega)$ for any $p\geq 1$, this implies that
    \begin{equation*}
        F(u) \geq \int_\Omega a u(x) -\alpha   \,dx \in \R \qquad \text{for any }u\in L^p(\Omega).
    \end{equation*}
    In particular, $F(u) > -\infty$ for all $u \in L^p(\Omega)$.
    Since $f$ is proper, there is a $t_0\in \dom f$. Hence (using again that $\Omega$ is bounded) the constant function $u_0 \equiv t_0 \in \dom F$ satisfies $F(u_0) < \infty$. This shows that $F$ is proper.

    To show convexity, we take $u,v\in \dom F$ (since otherwise \eqref{eq:convex:def} is trivially satisfied) and $\lambda\in[0,1]$ arbitrary.
    The convexity of $f$ now implies that
    \begin{equation*}
        f(\lambda u(x) + (1-\lambda)v(x)) \leq \lambda f(u(x))+(1-\lambda )f(v(x))\quad\text{for almost every }x\in\Omega.
    \end{equation*}
    Since $u,v\in \dom F$ and $L^1(\Omega)$ is a vector space, $\lambda f(u(x)) + (1-\lambda)f(v(x)) \in L^1(\Omega)$ as well. Similarly, the left-hand side is bounded from below by $a(\lambda u(x) + (1-\lambda)v(x))-\alpha\in L^1(\Omega)$ by \eqref{eq:lebesgue:lsc_bound}. We can thus integrate the inequality over $\Omega$ to obtain the convexity of $F$.

    To show lower semicontinuity, we use \cref{lem:convex:epi}. Let $\{(u_n,t_n)\}_{n\in\N}\subset \epi F$ with $u_n\to u$ in $L^p(\Omega)$ and $t_n\to t$ in $\R$.
    Then there exists a subsequence $\{u_{n_k}\}_{k\in\N}$ with $u_{n_k}(x)\to u(x)$ almost everywhere. Hence, the lower semicontinuity of $f$ together with Fatou's lemma implies that
    \begin{equation*}
        \begin{aligned}
            \int_\Omega f(u(x))-(au(x)-\alpha) \,dx
            &\leq \int_\Omega \liminf_{k\to\infty}(f(u_{n_k}(x))-(au_{n_k}(x)-\alpha ))\,dx\\
            &\leq  \liminf_{k\to\infty}\int_\Omega f(u_{n_k}(x))-(au_{n_k}(x)-\alpha )\,dx\\
            &= \liminf_{k\to\infty}\int_\Omega f(u_{n_k}(x))\,dx-\int_\Omega au(x)-\alpha\,dx
        \end{aligned}
    \end{equation*}
    as the integrands are nonnegative due to \eqref{eq:lebesgue:lsc_bound}.
    Since $(u_{n_k},t_{n_k})\in \epi F$, this yields
    \begin{equation*}
        F(u) =  \int_\Omega f(u(x))\,dx \leq \liminf_{k\to\infty} \int_\Omega f(u_{n_k}(x))\,dx = \liminf_{k\to\infty} F(u_{n_k}) \leq \lim_{k\to\infty} t_{n_k} = t,
    \end{equation*}
    i.e., $(u,t)\in \epi F$. Hence $\epi F$ is closed, and the lower semicontinuity of $F$ follows from  \cref{lem:convex:epi}\,\ref{lem:convex:epi:iii}.
\end{proof}

\section{Existence of minimizers}

After all this preparation, we can quickly prove the main result on existence of solutions to convex minimization problems.
\begin{theorem}\label{thm:convex:existence}
    Let $X$ be a reflexive Banach space and let
    \begin{enumerate}
        \item \label{thm:convex:existence:i}
            $U\subset X$ be nonempty, convex, and closed;
        \item \label{thm:convex:existence:ii}
            $F:X\to\Rbar$ be proper, convex, and lower semicontinuous with $\dom F \cap U \neq \emptyset$;
        \item \label{thm:convex:existence:iii}
            $U$ be bounded or $F$ be coercive.
    \end{enumerate}
    Then the problem
    \begin{equation*}
        \min_{x\in U} F(x)
    \end{equation*}
    admits a solution $\bar x\in U\cap \dom F$. If $F$ is strictly convex, the solution is unique.
\end{theorem}
\begin{proof}
    We consider the extended functional $\overline F = F + \delta_U:X\to\Rbar$.
    Assumption \ref{thm:convex:existence:i} together with \cref{lem:variation:indicator} implies that $\delta_U$ is proper, convex, and weakly lower semicontinuous. From \ref{thm:convex:existence:i} and \ref{thm:convex:existence:ii} we obtain an $x_0\in U$ with $\overline F(x_0)<\infty$ and hence that $\overline F$ is proper, convex, and (by \cref{cor:convex:uhs}) weakly lower semicontinuous. Finally, using that $F(x)>-\infty$ and $\delta_U(x) \geq 0$ for all $x\in X$, it follows from \ref{thm:convex:existence:iii} that
    \begin{itemize}
        \item if $U$ is bounded: $\overline F(x) = \infty$ for all $x\in\R\setminus U$;
        \item if $F$ is coercive: $\overline F(x) \geq F(x) \to \infty$ for $\norm{x}_X\to \infty$;
    \end{itemize}
    and hence that $\overline F$ is coercive.
    We can thus apply \cref{thm:variation:existence} to obtain the existence of a minimizer $\bar x \in \dom \overline F =  U\cap \dom F$ of $\overline F$ with
    \begin{equation*}
        F(\bar x) = \overline F(\bar x) \leq \overline F(x) = F(x) \qquad\text{for all }x\in U,
    \end{equation*}
    i.e., $\bar x$ is the claimed solution.

    Let now $F$ be strictly convex, and let $\bar x$ and $\bar x'\in U$ be two different minimizers, i.e., $F(\bar x) = F(\bar x') = \min_{x\in U}F(x)$ and $\bar x\neq \bar x'$.
    Then by the convexity of $U$ we have for all $\lambda\in (0,1)$ that
    \begin{equation*}
        x_\lambda\defeq\lambda \bar x + (1-\lambda) \bar x' \in U,
    \end{equation*}
    while the strict convexity of $F$ implies that
    \begin{equation*}
        F(x_\lambda) < \lambda F(\bar x) + (1-\lambda) F(\bar x') = F(\bar x).
    \end{equation*}
    But this is a contradiction to $F(\bar x)\leq F(x)$ for all $x\in U$.
\end{proof}

Note that for a sum of two convex functionals to be coercive, it is in general not sufficient that only one of them is. Functionals for which this is the case -- such as the indicator function of a bounded set -- are called \term[functional!supercoercive]{supercoercive}; another example which will be helpful later is the squared norm.
\begin{lemma}\label{lem:convex:supercoercive}
    Let $F:X\to\Rbar$ be proper, convex, and lower semicontinuous, and $x_0\in X$ be given. Then the functional
    \begin{equation*}
        J:X\to\Rbar,\qquad x\mapsto F(x) + \frac12\norm{x-x_0}_X^2
    \end{equation*}
    is coercive.
\end{lemma}
\begin{proof}
    Since $F$ is proper, convex, and lower semicontinuous, it follows from \cref{lem:convex:gamma} that $F$ is bounded from below by a continuous affine functional, i.e., there exists an $x^*\in X^*$ and an $\alpha\in\R$ with $F(x)\geq \dual{x^*,x}_X -\alpha$ for all $x\in X$. Together with the reverse triangle inequality and \eqref{eq:functan:cs_banach}, we obtain that
    \begin{equation*}
        \begin{aligned}
            J(x) &\geq \dual{x^*,x}_X - \alpha + \tfrac12\left(\norm{x}_X-\norm{x_0}_X\right)^2\\
                 &\geq -\norm{x^*}_{X^*}\norm{x}_X -\alpha + \tfrac12\norm{x}_X^2 -\norm{x}_X\norm{x_0}_X\\
                 &= \norm{x}_X\left(\tfrac12\norm{x}_X -\norm{x^*}_{X^*} - \norm{x_0}_X\right) - \alpha.
        \end{aligned}
    \end{equation*}
    Since $x^*$ and $x_0$ are fixed, the term in parentheses is positive for $\norm{x}_X$ sufficiently large, and hence $J(x)\to \infty$ for $\norm{x}_X\to \infty$ as claimed.
\end{proof}

\section{Continuity properties}\label{sec:convex:continuity}

To close this chapter, we show the following remarkable result: \emph{Any (locally) bounded convex functional is (locally) Lipschitz continuous.}
Besides being of use in later chapters, this result illustrates the beauty of convex analysis: an algebraic but global property (convexity) connects two topological but local properties (neighborhood and continuity).
Here we consider of course the strong topology in a normed vector space.
\begin{lemma}\label{thm:convex:cont_bounded}
    Let $X$ be a normed vector space, $F:X\to\Rbar$ be convex, and $x\in X$. If there is a $\rho>0$ such that $F$ is bounded from above on $\OB(x,\rho)$, then $F$ is locally Lipschitz continuous near~$x$.
\end{lemma}
\begin{proof}
    By assumption, there exists an $M\in\R$ with $F(y)\leq M$ for all $y\in \OB(x,\rho)$.
    We first show that $F$ is locally bounded from below as well. Let $y\in \OB(x,\rho)$ be arbitrary. Since $\norm{x-y}_X<\rho$, we also have that $z\defeq2x-y = x-(y-x)\in \OB(x,\rho)$, and the convexity of $F$ implies that
    $F(x) = F\left(\tfrac12 y + \tfrac12 z\right) \leq \tfrac12 F(y)+\tfrac12 F(z)$
    and hence that
    \begin{equation*}
        -F(y) \leq F(z) - 2 F(x) \leq M - 2F(x) =:m,
    \end{equation*}
    i.e., $-m\leq F(y) \leq M$ for all $y\in \OB(x,\rho)$.

    We now show that this implies Lipschitz continuity on $\OB(x,\frac\rho2)$. Let $y_1,y_2\in \OB(x,\frac\rho2)$ with $y_1\neq y_2$ and set
    \begin{equation*}
        z\defeq y_1 + \frac\rho2 \frac{y_1-y_2}{\norm{y_1-y_2}_X} \in \OB(x,\rho),
    \end{equation*}
    which holds because $\norm{z-x}_X \leq \norm{y_1-x}_X +\frac\rho2 < \rho$. By construction, we thus have that
    \begin{equation*}
        y_1 = \lambda z + (1-\lambda) y_2 \quad\text{for}\quad
        \lambda\defeq\frac{\norm{y_1-y_2}_X}{\norm{y_1-y_2}_X + \tfrac\rho2} \in(0,1),
    \end{equation*}
    and the convexity of $F$ now implies that $F(y_1) \leq \lambda F(z) + (1-\lambda)F(y_2)$.
    Together with the definition of $\lambda$ as well as $F(z)\leq M$ and $-F(y_2)\leq m=M-2F(x)$, this yields the estimate
    \begin{equation*}
        \begin{aligned}
            F(y_1) -F(y_2) \leq \lambda(F(z)-F(y_2)) &\leq \lambda(2M-2F(x))\\
                                                     &= \frac{2(M-F(x))}{\norm{y_1-y_2}_X + \frac\rho2} \norm{y_1-y_2}_X\\
                                                     &\leq \frac{2(M-F(x))}{\frac\rho2} \norm{y_1-y_2}_X.
        \end{aligned}
    \end{equation*}
    Exchanging the roles of $y_1$ and $y_2$, we obtain that
    \begin{equation*}
    |F(y_1)-F(y_2)| \leq \frac{4}{\rho}(M-F(x))\norm{y_1-y_2}_X \quad\text{for all }y_1,y_2\in \OB\left(x,\frac\rho2\right)
    \end{equation*}
    and hence the local Lipschitz continuity with constant $L(x,\rho/2)\defeq \frac4\rho(M-F(x))$.
\end{proof}
This result can be extended by showing that convex functions are bounded everywhere in the interior (again a topological concept!) of their effective domain.
As an intermediary step, we first consider the scalar case.%
\footnote{With a bit more effort, one can show that the claim holds for $F:\R^N\to\Rbar$ with arbitrary $N\in\N$; see, e.g., \cite[Corollary 1.4.2]{Schirotzek:2007}.}
\begin{lemma}\label{cor:convex:cont_r}
    If $f:\R\to\Rbar$ is convex, then $f$ is locally bounded from above on $\interior(\dom f)$.
\end{lemma}
\begin{proof}
    Let $x\in \interior(\dom f)$, i.e., there exist $a,b\in\R$ with $x\in(a,b)\subset \dom f$; by possibly shrinking the interval we can even assume that $[a,b]\subset \dom f$. Let now $z\in(a,b)$. Since intervals are convex, there exists a $\lambda\in(0,1)$ with $z=\lambda a+(1-\lambda)b$. By convexity, we thus have
    \begin{equation*}
        f(z) \leq \lambda f(a) + (1-\lambda)f(b) \leq \max\{|f(a)|,|f(b)|\} < \infty.
    \end{equation*}
    Hence $f$ is locally bounded from above in $x$.
\end{proof}

The proof of the general case requires further assumptions on $X$ and $F$.
\begin{theorem}\label{thm:convex:bounded}
    Let $X$ be a Banach space. If $F:X\to\Rbar$ is proper, convex, and lower semicontinuous, then $F$ is locally bounded from above on $\interior(\dom F)$.
\end{theorem}
\begin{proof}
    We first show the claim for the case $x=0\in\interior(\dom F)$, which implies in particular that $M\defeq |F(0)|$ is finite. Consider now for arbitrary $h\in X$ the mapping
    \begin{equation*}
        f:\R \to\Rbar,\qquad t\mapsto F(th).
    \end{equation*}
    It is straightforward to verify that $f$ is convex and satisfies $0\in\interior(\dom f)$.
    By \cref{thm:convex:cont_bounded,cor:convex:cont_r}, $f$ is thus locally Lipschitz continuous near $0$; hence in particular $|f(t)-f(0)|\leq L t \leq 1$ for sufficiently small $t>0$. The reverse triangle inequality therefore yields a $\delta>0$ with
    \begin{equation*}
        F(0+t h) \leq |F(0+th)| = |f(t)| \leq |f(0)|+1 = M+1 \qquad\text{for all } t\in [0,\delta].
    \end{equation*}
    Hence $0$ lies in the algebraic interior of the sublevel set $\sub_{M+1} F$, which is convex and closed (since we assumed $F$ to be lower semicontinuous) by \cref{lem:convex:sublevel}. The core--int \cref{lem:functan:coreint} thus yields that $0\in\interior(\sub_{M+1} F)$, i.e., there exists a $\rho>0$ with $F(z)\leq M+1$ for all $z\in \OB(0,\rho)$.

    For the general case $x\in \interior(\dom F)$, consider
    \begin{equation*}
        \tilde F:X\to\Rbar,\qquad y\mapsto F(y+x).
    \end{equation*}
    Again, it is straightforward to verify convexity and lower semicontinuity of  $\tilde F$ and that $0\in \interior(\dom \tilde F)$.
    It follows from what we have shown so far that $\tilde F$ is locally bounded from above on $\OB(0,\rho)$, which immediately implies that $F$ is locally bounded from above on $\OB(x,\rho)$.
\end{proof}
Together with \cref{thm:convex:cont_bounded}, we thus obtain the desired result.
\begin{theorem}\label{thm:convex:cont}
    Let $X$ be a Banach space. If $F:X\to\Rbar$ is proper, convex, and lower semicontinuous, then $F$ is locally Lipschitz continuous on $\interior(\dom F)$.
\end{theorem}
We shall have several more occasions to observe the unreasonably nice behavior of convex lower semicontinuous functions on the interior of their effective domain.

\chapter{Convex subdifferentials}
\label{chap:subdiff}

For convex functionals, we can use the general properties from the previous chapter to obtain explicit optimality conditions. We do this by first deriving a Fermat principle in terms of a generalized derivative that can be used to characterize global minimizers of nonsmooth functionals. The remainder of the chapter is then devoted to the explicit characterization of this generalized derivative specifically for convex lower semicontinuous functionals; first directly for elementary examples, then for more complicated functions by deriving calculus rules like a sum and a chain rule.

\section{Definition and basic properties}

The motivation for our notion of generalized derivative is geometric: The classical derivative $f'(t)$ of a scalar function $f:\R\to\R$ at $t$ can be interpreted as the slope of the tangent line to $f$ at $t$. If the function is not differentiable, the tangent line -- if it exists at all -- need no longer be unique. The idea is then to define as the generalized derivative the \emph{set of all} tangent slopes.
In multiple dimensions, tangent lines generalize to \emph{supporting hyperplanes}; we thus define in a normed vector space $X$ the \term[subdifferential!convex]{(convex) subdifferential} of $F:X\to\Rbar$ at $x\in \dom F$ as
\begin{equation}\label{eq:convex:def_subdiff}
    \partial F(x) \defeq  \setof{x^*\in X^*}{\dual{x^*,\tilde x - x}_X \leq F(\tilde x) - F(x)\quad\text{for all }\tilde x\in X}.
\end{equation}
(Note that $\tilde x\notin \dom F$ is allowed since in this case the inequality is trivially satisfied.) For $x\notin\dom F$, we set $\partial F(x) = \emptyset$.
An element $x^*\in \partial F(x)$ is called a \term{subderivative}. (Following the terminology for classical derivatives, we reserve the more common term \term{subgradient} for its Riesz representation $z_{x^*} \in X$ when $X$ is a Hilbert space.)

The following example shows that the subdifferential can also be empty for $x\in \dom F$, even if $F$ is convex.
\begin{example}\label{ex:subdiff:empty}
    We take $X=\R$ (and hence $X^*\cong X = \R$) and consider
    \begin{equation*}
        F(x) =
        \begin{cases}
            -\sqrt{x} & \text{if } x\geq 0,\\
            \infty & \text{if } x<0.
        \end{cases}
    \end{equation*}
    Since \eqref{eq:convex:def} is trivially satisfied if $x$ or $y$ is negative, we can assume $x,y\geq 0$ so that we are allowed to take the square of both sides of \eqref{eq:convex:def}. A straightforward algebraic manipulation then shows that this is equivalent to $\lambda(\lambda-1)(\sqrt{x}-\sqrt{y})^2\geq 0$, which holds for any $x,y\geq 0$ and $\lambda\in[0,1]$. Hence $F$ is convex.

    However, for $x=0$, any $x^*\in \partial F(0)$ by definition must in particular satisfy
    \begin{equation*}
        x^*\cdot \tilde x \leq -\sqrt{\tilde x} \qquad\text{for all }\tilde x \ge 0.
    \end{equation*}
    Taking now $\tilde x>0$ arbitrary, we can divide by it on both sides and let $\tilde x\to 0$ to obtain
    \begin{equation*}
        x^* \leq -\left(\sqrt{\tilde x}\right)^{-1} \to -\infty.
    \end{equation*}
    This is impossible for $x^*\in \R\cong X^*$. Hence, $\partial F(0)$ is empty.
\end{example}
In fact, it will become clear that the nonexistence of tangent lines is much more problematic than the nonuniqueness.
However, we will later show that for proper, convex, and lower semicontinuous functionals, $\partial F(x)$ is nonempty (and bounded) for all $x\in\interior(\dom F)$; see \cref{cor:convex:nonempty}.
Furthermore, it follows directly from the definition that for all $x\in X$, the set $\partial F(x)$ is convex and weakly-$*$ closed.

The definition immediately yields a Fermat principle.
\begin{theorem}[Fermat principle]\label{thm:convex:fermat}\index{principle!Fermat!convex}
    Let $F:X\to\Rbar$ and $\bar x \in \dom F$. Then the following statements are equivalent:
    \begin{enumerate}
        \item $\displaystyle  0\in \partial F(\bar x)$;
        \item $\displaystyle F(\bar x) = \min_{x\in X} F(x)$.
    \end{enumerate}
\end{theorem}
\begin{proof}
    This is a direct consequence of the definitions: $0\in\partial F(\bar x)$ if and only if
    \begin{equation*}
        0= \dual{0,\tilde x-\bar x}_X \leq F(\tilde x) - F(\bar x) \qquad\text{for all } \tilde x\in X,
    \end{equation*}
    i.e., $F(\bar x)\leq F(\tilde x)$ for all $\tilde x\in X$.%
    \footnote{Note that convexity of $F$ is not required for \cref{thm:convex:fermat}. The condition $0\in \partial F(\bar x)$ therefore characterizes the global(!) minimizers of \emph{any} function $F$. However, nonconvex functionals can also have local minimizers, for which the subdifferential inclusion is not satisfied.
    In fact, (convex) subdifferentials of nonconvex functionals are usually empty. (And conversely, one can show that $\partial F(x)\neq \emptyset$ for all $x\in\dom F$ implies that $F$ is convex.) This leads to problems in particular for the proof of calculus rules, for which we will indeed have to assume convexity.}
\end{proof}
This matches the geometrical intuition: If $X=\R\cong X^*$, the affine function $\tilde F(\tilde x) \defeq F(x) + x^*(\tilde x - x)$ with $x^*\in \partial F(x)$ describes a tangent line at $(x,F(x))$ with slope $x^*$; the condition $x^*=0\in \partial F(\tilde x)$ thus means that $F$ has a horizontal tangent line in $\bar x$.
(Conversely, the function from \cref{ex:subdiff:empty} only has a vertical tangent line in $x=0$, which corresponds to an infinite slope that is not an element of $\R\cong X^*$.)

\bigskip

Not surprisingly, the convex subdifferential behaves more nicely for convex functions. The key property is an alternative characterization using directional derivatives, which exist (at least in the extended real-valued sense) for any convex function.
\begin{lemma}\label{lem:convex:direct}
    Let $F:X\to\Rbar$ be convex and let $x\in \dom F$ and $h\in X$ be given. Then
    \begin{enumerate}
        \item \label{lem:convex:direct:i}
            the function
            \begin{equation*}
                \phi:(0,\infty) \to \Rbar,\qquad t\mapsto \frac{F(x+th)-F(x)}{t},
            \end{equation*}
            is increasing;
        \item \label{lem:convex:direct:ii}
            there exists a limit $F'(x;h)=\lim_{t\downto 0}\phi(t)\in[-\infty,\infty]$, which satisfies
            \begin{equation*}
                F'(x;h) \leq F(x+h) - F(x);
            \end{equation*}
        \item \label{lem:convex:direct:iii}
            if $x\in \interior(\dom F)$, the limit $F'(x;h)$ is finite.
    \end{enumerate}
\end{lemma}
\begin{proof}
    \emph{\ref{lem:convex:direct:i}:} Inserting the definition and sorting terms shows that for all $0<s\leq t$, the condition $\phi(s)\leq \phi(t)$ is equivalent to
    \begin{equation*}
        F(x+sh) \leq \frac{s}{t} F(x+th) + \left(1-\frac{s}{t} \right)F(x),
    \end{equation*}
    which follows from the convexity of $F$ since $x+sh = \frac{s}{t}(x+th)+(1-\frac{s}{t})x$.

    \emph{\ref{lem:convex:direct:ii}:} The claim immediately follows from \ref{lem:convex:direct:i} since
    \begin{equation*}
        F'(x;h) = \lim_{t\downto 0} \phi(t) = \inf_{t>0} \phi(t) \leq \phi(1) = F(x+h)-F(x)\in \Rbar.
    \end{equation*}

    \emph{\ref{lem:convex:direct:iii}:} Since $\interior(\dom F)$ is contained in the algebraic interior of $\dom F$, there exists an $\eps>0$ such that $x +t h\in \dom F$ for all $t\in[-\eps,\eps]$. Proceeding as in \ref{lem:convex:direct:i}, we obtain that $\phi(s)\leq \phi(t)$ for all $s<t<0$ as well. From $x=\frac12(x+th)+\frac12(x-th)$ for $t>0$, we also obtain that
    \begin{equation*}
        \phi(-t) = \frac{F(x-th)-F(x)}{-t} \leq \frac{F(x+th)-F(x)}{t} = \phi(t)
    \end{equation*}
    and hence that $\phi$ is increasing on all $\R\setminus\{0\}$. As in \ref{lem:convex:direct:ii}, the choice of $\eps$ now implies that
    \begin{equation*}
        -\infty < \phi(-\eps) \leq F'(x;h) \leq \phi(\eps) < \infty.
        \qedhere
    \end{equation*}
\end{proof}

\begin{lemma}\label{lem:convex:equiv}
    Let $F:X\to\Rbar$ be convex and $x\in \dom F$. Then
    \begin{equation*}
        \subdiff F(x) = \setof{x^* \in X^*}{\dual{x^*,h}_X \leq F'(x; h) \text{ for all } h \in X}.
    \end{equation*}
\end{lemma}
\begin{proof}
    Since any $\tilde x\in X$ can be written as $\tilde x = x+h$ for some $h\in X$ and vice versa, it suffices to show that for any $x^*\in X^*$, the following statements are equivalent:
    \begin{enumerate}
        \item \label{lem:convex:equiv:i}
            $\dual{x^*,h}_X \leq F'(x;h)$ \qquad for all $h\in X$;
        \item \label{lem:convex:equiv:ii}
            $\dual{x^*,h}_X \leq F(x+h) - F(x)$ \quad for all $h\in X$.
    \end{enumerate}
    If $x^*\in X^*$ satisfies $\dual{x^*,h}_X\leq F'(x;h)$ for all $h\in X$,
    we immediately obtain from \cref{lem:convex:direct}\,\ref{lem:convex:direct:ii} that
    \begin{equation*}
        \dual{x^*,h}_X \leq F'(x;h) \leq F(x+h) - F(x)\qquad \text{for all }h\in X.
    \end{equation*}
    Setting $\tilde x = x+h\in X$ then yields $x^*\in \partial F(x)$.

    Conversely, if $\dual{x^*,h}\leq F(x+h)-F(x)$ holds for all $h\in X$, it also holds for $th$ for all $h\in X$ and $t>0$. Dividing by $t$ and passing to the limit (which exists by \cref{lem:convex:direct}\,\ref{lem:convex:direct:ii}) then yields that
    \begin{equation*}
        \dual{x^*,h}_X \leq \lim_{t\downto 0} \frac{F(x+th)-F(x)}{t} = F'(x;h).
        \qedhere
    \end{equation*}
\end{proof}

\section{Fundamental examples}

We now look at some examples. First, the characterization via the directional derivative indicates that the subdifferential is indeed a generalization of the Gateaux derivative.
\begin{theorem}\label{thm:convex:gateaux}
    Let $F:X\to\Rbar$ be convex. If $F$ is Gateaux differentiable at $x$, then $\partial F(x) = \{DF(x)\}$.
\end{theorem}
\begin{proof}
    By definition of the Gateaux derivative, we have that
    \begin{equation*}
        \dual{DF(x),h}_X = DF(x)h =  F'(x;h) \quad\text{for all } h\in X.
    \end{equation*}
    \Cref{lem:convex:equiv} now immediately yields $DF(x)\in \partial F(x)$.
    Conversely, $x^*\in \partial F(x)$ again by \cref{lem:convex:equiv} implies that
    \begin{equation*}
        \dual{x^*,h}_X \leq F'(x;h) =\dual{DF(x),h}_X \quad\text{for all } h\in X.
    \end{equation*}
    Taking the supremum over all $h$ with $\norm{h}_X\leq 1$ now yields that $\norm{x^* - DF(x)}_{X^*} \leq 0$, i.e., $x^* = DF(x)$.
\end{proof}
The converse holds as well: If $x\in \interior(\dom F)$ and $\partial F(x)$ is a singleton, then $F$ is Gateaux differentiable; see \cref{thm:convex:singleton}.

Of course, we also want to compute subdifferentials of functionals that are not differentiable. The canonical example is the norm $\norm{\cdot}_X$ on a normed vector space, which even for $X=\R$ is not differentiable at $x=0$.
\begin{theorem}\label{thm:subdifferential:norm}
    For any $x\in X$,
    \begin{equation*}\label{eq:subdifferential:norm}
        \partial(\norm{\cdot}_X)(x) =
        \begin{cases}
            \setof{x^*\in X^*}{\dual{x^*,x}_X = \norm{x}_X \text{ and } \norm{x^*}_{X^*} = 1} &\text{if } x\neq 0,\\
            \B_{X^*} &\text{if } x = 0.
        \end{cases}
    \end{equation*}
\end{theorem}
\begin{proof}
    For $x=0$, we have $x^*\in\partial(\norm{\cdot}_X)(x)$ by definition if and only if
    \begin{equation*}
        \dual{x^*,\tilde x}_X \leq \norm{\tilde x}_X\qquad\text{for all }\tilde x\in X\setminus\{0\}
    \end{equation*}
    (since the inequality is trivial for $\tilde x=0$), which by the definition of the operator norm holds if and only if $\norm{x^*}_{X^*} \leq 1$.

    Let now $x\neq 0$ and consider $x^*\in\partial(\norm{\cdot}_X)(x)$. Inserting first $\tilde x=0$ and then $\tilde x=2x$ into the definition \eqref{eq:convex:def_subdiff} yields the sequence of inequalities
    \begin{equation*}
        \norm{x}_X \leq \dual{x^*,x}_X = \dual{x^*,2x -x } \leq \norm{2x}_X - \norm{x}_X = \norm{x}_X,
    \end{equation*}
    which imply that $\dual{x^*,x}_X=\norm{x}_X$. Similarly, we have for all $\tilde x \in X$ that
    \begin{equation*}
        \dual{x^*,\tilde x}_X = \dual{x^*,(\tilde x+x)- x}_X
        \leq \norm{\tilde x+x}_X - \norm{x}_X \leq \norm{\tilde x}_X.
    \end{equation*}
    As in the case $x=0$, this implies that $\norm{x^*}_{X^*}\leq 1$. For $\tilde x = x/\norm{x}_{X}$ we further have that
    \begin{equation*}
        \dual{x^*,\tilde x}_X = \norm{x}_X^{-1}\dual{x^*,x}_X = \norm{x}_X^{-1} \norm{x}_X = 1.
    \end{equation*}
    Hence, $\norm{x^*}_{X^*}=1$ is in fact attained.

    Conversely, let $x^*\in X^*$ with $\dual{x^*,x}_X =\norm{x}_X$ and $\norm{x^*}_{X^*} = 1$. Then we obtain for all $\tilde x \in X$ from \eqref{eq:functan:cs_banach} the relation
    \begin{equation*}
        \dual{x^*,\tilde x- x}_X = \dual{x^*,\tilde x}_X -  \dual{x^*, x}_X\leq \norm{\tilde x}_X - \norm{x}_X,
    \end{equation*}
    and hence $x^*\in\partial(\norm{\cdot}_X)(x)$ by definition.
\end{proof}
\begin{example}\label{ex:convex:subdiff_abs}
    In particular, we obtain for $X=\R$ the subdifferential of the absolute value function as%
    \tablefootnote{Note that this set-valued definition of $\sign(t)$ differs from the usual (single-valued) one, in particular for $t=0$; to make this distinction clear, one often refers to \eqref{eq:convex:subdiff_abs} as the \term[sign|infn]{sign in the sense of convex analysis}. Throughout this book, we will always use the sign in this sense.}
    \begin{equation}\label{eq:convex:subdiff_abs}
        \partial(|\cdot|)(t) = \sign(t) \defeq
        \begin{cases}
            \{1\} & \text{if }t>0,\\
            \{{-}1\} & \text{if }t<0,\\
            [-1,1] & \text{if }t=0,
        \end{cases}
    \end{equation}
    cf.~\cref{fig:subdiff:abs}.
\end{example}
\begin{figure}[t]
    \centering
    \begin{subfigure}[t]{0.495\textwidth}
        \centering
        \begin{asy}
            draw((-1.2,0)..(1.2,0),linewidth(0.5),Arrow); label("$x$",(1.2,0),S);
            draw((0,-1.2)..(0,1.2),linewidth(0.5),Arrow); label("$\partial F(x)$",(0,1.2),NE);
            dot((0,-1)); label("$-1$",(0,-1),E);
            dot((0,0)); label("$0$",(0,0),SE);
            dot((0,1)); label("$1$",(0,1),W);

            draw((0,-1)..(-1.2,-1),primalline+linewidth(1.5));
            draw((0,-1)..(0,1),primalline+linewidth(1.5));
            draw((0,1)..(1.2,1),primalline+linewidth(1.5));
        \end{asy}
            \caption{$F(x)=|x|$}\label{fig:subdiff:abs}
    \end{subfigure}
    \hfill
    \begin{subfigure}[t]{0.495\textwidth}
        \centering
        \begin{asy}
            draw((-1.2,0)..(1.2,0),linewidth(0.5),Arrow); label("$x$",(1.2,0),S);
            draw((0,-1.2)..(0,1.2),linewidth(0.5),Arrow); label("$\partial F(x)$",(0,1.2),NE);
            dot((0,-1)); label("$-1$",(0,-1),E);
            dot((-1,0)); label("$-1$",(-1,0),N);
            dot((0,0)); label("$0$",(0,0),SE);
            dot((1,0)); label("$1$",(1,0),S);

            draw((-1,0)..(-1,-1.2),primalline+linewidth(1.5));
            draw((-1,0)..(1,0),primalline+linewidth(1.5));
            draw((1,0)..(1,1.2),primalline+linewidth(1.5));
        \end{asy}
        \caption{$F(x)=\delta_{[-1,1]}(x)$}\label{fig:subdiff:ind}
    \end{subfigure}
    \caption{Illustration of $\graph \partial F$ for two different functions $F:\R\to\Rbar$.}
    \label{fig:subdiff}
\end{figure}
We can also give a more explicit characterization of the subdifferential of the indicator functional of a set $C\subset X$.
\begin{lemma}\label{lem:convex:normalcone}
    For any $C\subset X$,
    \begin{equation*}
        \partial\delta_C(x) = \setof{x^*\in X^*}{\dual{x^*,\tilde x - x}_X \leq 0
        \quad\text{for all }\tilde x \in C}.
    \end{equation*}
\end{lemma}
\begin{proof}
    For any $x\in C=\dom\delta_C$, we have that
    \begin{equation*}
        \begin{aligned}
            x^*\in\partial\delta_C(x) &\equivalent \dual{x^*,\tilde x - x}_X \leq \delta_C(\tilde x) && \text{for all }\tilde x \in X\\
            &\equivalent \dual{x^*,\tilde x - x}_X \leq 0 && \text{for all }\tilde x \in C,
        \end{aligned}
    \end{equation*}
    since the first inequality is trivially satisfied for all $\tilde x\notin C$.
\end{proof}
The set $N_C(x) \defeq \partial\delta_C(x)$ is also called the (convex) \term[cone!normal!convex]{normal cone} to $C$ at $x$ (which may be empty if $C$ is not convex).
We illustrate such sets in \cref{fig:subdiff:indicator}.
Depending on the set $C$, this can be made even more explicit.
\begin{example}\label{ex:convex:subdiff_ind}
    Let $X=\R$ and $C=[-1,1]$, and let $t\in C$. Then we have  $x^* \in \partial\delta_{[-1,1]}(t)$ if and only if $x^*(\tilde t -t)\leq 0$ for all $\tilde t \in [-1,1]$. We proceed by distinguishing three cases.
    \begin{enumerate}[label={Case }\arabic*:, align=left]
        \item  $t=1$. Then $\tilde t - t \in [-2,0]$, and hence the product is nonpositive if and only if $x^*\geq 0$.
        \item  $t=-1$. Then $\tilde t - t \in [0,2]$, and hence the product is nonpositive if and only if $x^*\leq 0$.
        \item  $t\in(-1,1)$. Then $\tilde t -t$ can be positive as well as negative, and hence only $x^* = 0$ is possible.
    \end{enumerate}
    We thus obtain that
    \begin{equation}\label{eq:subdiff:indicator}
        \partial\delta_{[-1,1]}(t) = \begin{cases}
            [0,\infty) & \text{if }t=1,\\
            (-\infty,0] & \text{if }t=-1,\\
            \{0\} & \text{if }t\in (-1,1),\\
            \emptyset & \text{if }t \in \R\setminus[-1,1],
        \end{cases}
    \end{equation}
    cf.~\cref{fig:subdiff:ind}.
    Readers familiar with (non)linear optimization will recognize these as the \term[condition!complementarity]{complementarity conditions} for Lagrange multipliers corresponding to the inequalities $-1\leq t \leq 1$.
\end{example}
Conversely, subdifferentials of functionals can be obtained from normal cones to corresponding epigraphs (which for convex functionals are convex sets by \cref{lem:convex:epi}). This relation will be the basis for defining further subdifferentials for more general classes of mappings in \cref{part:setvalued}.
We illustrate this result for the absolute value function of \cref{ex:convex:subdiff_abs} in \cref{fig:subdiff:epi-abs}.

\begin{figure}
    \centering
    \begin{asy}
        path p=(0, 0)..(0.3, 1)..(2,1)..controls (3, .5)..(1, .1)--cycle;
        fill(p, lightfill);
        draw(p);
        label("$C$", (.8, .6));

        int t1=0;
        pair x1=point(p, t1);
        pair n1=orthog(dir(p, t1, 1), -1);
        pair n2=orthog(dir(p, t1, -1), -1);
        fill(x1--(x1+1.2*n1)--(x1+1.1*n2)--cycle, darkfill);
        draw(x1--(x1+n1), primalline + linewidth(1.1), Arrow);
        draw(x1--(x1+n2), primalline + linewidth(1.1), Arrow);
        label("$N_C(x_1)$", x1+(n1+n2)/2);
        dot(x1);

        real t2=2.8;
        pair x2=point(p, t2);
        pair n=orthog(dir(p, t2), -1);
        draw(x2--(x2+n), primalline + linewidth(1.1), Arrow);
        label("$N_C(x_2)$", x2+n, S);
        dot(x2);

    \end{asy}
    \caption{Normal cones of a convex set $C$ at two points $x_1$ and $x_2$.}
    \label{fig:subdiff:indicator}
\end{figure}

\begin{lemma}\label{lem:convex:subdiff_epi}
    Let $F:X\to \Rbar$ and $x\in \dom F$. Then $x^*\in \partial F(x)$ if and only if $(x^*,-1) \in N_{\epi F}(x,F(x))$.
\end{lemma}
\begin{proof}
    By definition of the normal cone, $(x^*,-1)\in N_{\epi F}(x,F(x))$ is equivalent to
    \begin{equation}\label{eq:convex:subdiff_epi}
        \dual{x^*,\tilde x-x}_X - (\tilde t- F(x)) \leq 0 \quad \text{for all }(\tilde x,\tilde t)\in \epi F,
    \end{equation}
    i.e., for all $\tilde x\in X$ and $\tilde t \geq F(\tilde x)$. Taking $\tilde t = F(\tilde x)$ and rearranging, this yields that $x^*\in\partial F(x)$.

    Conversely, from $x^*\in\partial F(x)$ we immediately obtain that
    \begin{equation*}
        \dual{x^*,\tilde x-x}_X \leq F(\tilde x) - F(x) \leq \tilde t - F(x) \quad \text{for all }\tilde x\in X, \tilde t\geq F(\tilde x),
    \end{equation*}
    i.e., \eqref{eq:convex:subdiff_epi} and thus $(x^*,-1)\in N_{\epi F}(x, F(x))$.
\end{proof}

\begin{figure}
    \centering
    \begin{subfigure}[t]{.4\textwidth}
    \centering
    \begin{asy}
        real f(real x){ return abs(x); };
        real gf(real x){ return x/abs(x); }
        path g=graph(f, -1, 1);
        fill(g--(1, 0.5+f(1))--(-1, 0.3+f(-1))--cycle, lightfill);
        label("$f$", point(g, .5), W);
        draw(g, defaultpen);
        real x0=.5; pair xy0=pt(f, x0);
        real grad=gf(x0); pair gv=.7*(1, grad); pair nv=.7*(grad, -1);
        draw(xy0--(xy0+nv), primalline + linewidth(1.1), Arrow);
        label("$(g, -1)$", xy0+nv, E);
        dot(xy0);
        dot(.7*(1, -1), invisible);
    \end{asy}
    \caption{$\subdiff f(x)=\{\sign x\}$ at $x \ne 0$}
    \end{subfigure}
    \begin{subfigure}[t]{.4\textwidth}
    \centering
    \begin{asy}
        real f(real x){ return abs(x); };
        path g=graph(f, -1, 1);
        fill(g--(1, 0.5+f(1))--(-1, 0.3+f(-1))--cycle, lightfill);
        label("$f$", point(g, .25), W);
        draw(g, defaultpen);
        real x0=0; pair xy0=pt(f, x0);
        real gradl=-1; pair gvl=.7*(1, gradl); pair nvl=.7*(gradl, -1);
        real gradu=1; pair gvu=.7*(1, gradu); pair nvu=.7*(gradu, -1);
        fill(xy0--(xy0+1.2*nvl)--(xy0+1.1*nvu)--cycle, darkfill);
        draw(xy0--(xy0+nvl), primalline + linewidth(1.1), Arrow);
        draw(xy0--(xy0+nvu), primalline + linewidth(1.1), Arrow);
        label("$(g, -1)$", xy0+(nvl+nvu)/2);
        dot(xy0);
    \end{asy}
    \caption{$\subdiff f(x)=[-1, 1]$ at $x=0$}
    \end{subfigure}
    \caption{Subdifferentials of $f(x)=\abs{x}$ in terms of the normal cone of the epigraph.}
    \label{fig:subdiff:epi-abs}
\end{figure}

The following result furnishes a crucial link between finite- and infinite-dimensional convex optimization. We again assume (as we will from now on) that $\Omega\subset \R^d$ is open and bounded.
\begin{theorem}\label{thm:lebesgue:subdiff}
    Let $f:\R\to\Rbar$ be proper, convex, and lower semicontinuous, and let $F:L^p(\Omega)\to\Rbar$ with $1\leq p <\infty$ be as in \cref{lem:lebesgue:lsc}. Then we have for all $u\in\dom F$ with $q\defeq\frac{p}{p-1}$ that
    \begin{equation*}
        \partial F(u)  = \setof{u^*\in L^q(\Omega)}{u^*(x) \in\partial f(u(x))\quad\text{for almost every } x\in\Omega}.
    \end{equation*}
\end{theorem}
\begin{proof}
    Let $u, \tilde u\in \dom F$, i.e., $f\circ u, f \circ \tilde u\in L^1(\Omega)$ (otherwise there is nothing to show), and let $u^*\in L^q(\Omega)$ be arbitrary.
    If $u^*(x)\in \partial f(u(x))$ almost everywhere, we can integrate over all $x\in\Omega$ to obtain
    \begin{equation*}
        F(\tilde u)-F(u) =
        \int_\Omega f(\tilde u(x))-f(u(x))\,dx \geq
        \int_\Omega u^*(x)(\tilde u(x)-u(x))\,dx =
        \dual{u^*,\tilde u-u}_{L^p},
    \end{equation*}
    i.e., $u^*\in\partial F(u)$.

    Conversely, let $u^*\in\partial F(u)$. Then by definition it holds that
    \begin{equation*}
        \int_\Omega u^*(x)(\tilde u(x)-u(x))\,dx \leq \int_\Omega f(\tilde u(x))-f(u(x))\,dx \quad\text{for all }\tilde u\in L^p(\Omega).
    \end{equation*}
    Let now $t\in\R$ be arbitrary and let $A\subset \Omega$ be an arbitrary measurable set. Setting
    \begin{equation*}
        \tilde u(x) \defeq \begin{cases} t & \text{if }x\in A,\\ u(x) & \text{if }x\notin A, \end{cases}
    \end{equation*}
    the above inequality implies due to $\tilde u\in L^p(\Omega)$ that
    \begin{equation*}
        \int_A u^*(x)(t-u(x))\,dx \leq \int_A f(t)-f(u(x))\,dx.
    \end{equation*}
    Since $A$ was arbitrary, it must hold that
    \begin{equation*}
        u^*(x)(t-u(x))\leq f(t)-f(u(x)) \qquad\text{for almost every }x\in\Omega.
    \end{equation*}
    Furthermore, since $t\in\R$ was arbitrary, we obtain that $u^*(x) \in\partial f(u(x))$ for almost every $x\in\Omega$.
\end{proof}
\begin{remark}\label{rem:lebesgue:subdiff}
    A similar representation can be shown for vector-valued and spatially-dependent integrands $f:\Omega\times \R\to \R^m$ under stronger assumptions; see, e.g., \cite[Corollary 3F]{Rockafellar:1976}.
\end{remark}

A similar proof shows that for $F:\R^N\to\Rbar$ with $F(x) = \sum_{i=1}^N f_i(x_i)$ and $f_i:\R\to\Rbar$ convex, we have for any $x\in \dom F$ that
\begin{equation*}
    \partial F(x) = \setof{x^*\in\R^N}{x^*_i\in \partial f_i(x_i),\quad 1\leq i \leq N}.
\end{equation*}
Together with the above examples, this yields componentwise expressions for the subdifferential of the norm $\norm{\cdot}_1$ as well as of the indicator functional of the unit ball with respect to the supremum norm in $\R^N$.

\section{Calculus rules}

As for classical derivatives, one rarely obtains subdifferentials from the fundamental definition but rather by applying calculus rules. It stands to reason that these are more difficult to derive the weaker the derivative concept is (i.e., the more functionals are differentiable in that sense).
For convex subdifferentials, the following two rules still follow directly from the definition.
\begin{lemma}\label{lem:convex:subdiff_calc}
    Let $F:X\to\Rbar$ and $x\in \dom F$. Then,
    \begin{enumerate}
        \item \label{lem:convex:subdiff_calc:i}
            $\partial(\lambda F)(x) = \lambda(\partial F(x))\defeq \setof{\lambda x^*}{x^*\in\partial F(x)}$ for $\lambda \geq 0$;
        \item \label{lem:convex:subdiff_calc:ii}
            $\partial F(\cdot + x_0)(x) = \partial F(x+x_0)$ for $x_0\in X$ with $x+x_0\in\dom F$.
    \end{enumerate}
\end{lemma}

Already the sum rule is considerably more delicate and requires additional assumptions.
\begin{theorem}[sum rule]\label{thm:subdiff:sum}
    Let $X$ be a Banach space, $F,G:X\to\Rbar$ be convex and lower semicontinuous, and $x\in \dom F\cap \dom G$. Then
    \begin{equation*}
        \partial F(x) + \partial G(x) \subset \partial (F+G)(x),
    \end{equation*}
    with equality if there exists an $x_0 \in \interior(\dom F)\cap \dom G$.
\end{theorem}
\begin{proof}
    The inclusion follows directly from adding the definitions of the two subdifferentials.
    Let therefore $x\in \dom F\cap \dom G$ and $x^*\in \partial(F+G)(x)$, i.e., satisfying
    \begin{equation}\label{eq:convex:sum:fpg}
        \dual{x^*,\tilde x - x}_X \leq(F(\tilde x) + G(\tilde x)) - (F(x) + G(x)) \quad\text{for all }\tilde x \in X.
    \end{equation}
    Our goal is now to use (as in the proof of \cref{lem:convex:gamma}) the characterization of convex functionals via their epigraph together with the Hahn--Banach separation theorem to construct a bounded linear functional $y^*\in \partial G(x)\subset X^*$ with $x^*-y^*\in \partial F(x)$, i.e.,
    \begin{equation*}
        \begin{aligned}
            F(\tilde x) - F(x) - \dual{x^*,\tilde x-x}_X &\geq \dual{y^*,x- \tilde x}_X \quad\text{for all }\tilde x \in \dom F,\\
            G(x) - G(\tilde x) &\leq\dual{y^*,x- \tilde x}_X \quad\text{for all }\tilde x \in \dom G.
        \end{aligned}
    \end{equation*}
    For that purpose, we define the sets
    \begin{align*}
        C_1 &\defeq \setof{(\tilde x,t-(F(x) -\dual{x^*,x}_X))}{\tilde x\in \dom F,\ t \geq F(\tilde x) - \dual{x^*,\tilde x}_X },\\
        C_2 &\defeq \setof{(\tilde x,G(x)-t)}{\tilde x \in \dom G,\ t\geq G(\tilde x)},
    \end{align*}
    i.e.,
    \begin{equation*}
        C_1 = \epi (F-x^*) - (0,F(x) - \dual{x^*,x}_X),\qquad
        C_2 = -(\epi G - (0,G(x))),
    \end{equation*}
    cf.~\cref{fig:subdiff:sum1}.
    To apply \cref{lem:convex:eidelheit} to these sets, we have to verify its conditions.
    \begin{enumerate}
        \item Since $x\in \dom F\cap \dom G$, both $C_1$ and $C_2$ are nonempty.
            Furthermore, since $F$ and $G$ are convex, it is straightforward (if tedious) to verify from the definition that $C_1$ and $C_2$ are convex.
        \item The critical point is of course the nonemptiness of $\interior C_1$, for which we argue as follows. Since $x_0\in \interior(\dom F)$, we know from \cref{thm:convex:bounded} that $F$ is bounded in an open neighborhood $U\subset \interior(\dom F)$ of $x_0$. We can thus find an open interval $I\subset \R$ such that $U\times I \subset C_1$. Since $U\times I$ is open by the definition of the product topology on $X\times \R$, any $(x_0,\alpha)$ with $\alpha\in I$ is an interior point of $C_1$.
        \item It remains to show that $\interior C_1\cap C_2=\emptyset$.
            Assume there exists a $(\tilde x,\alpha) \in \interior C_1 \cap C_2$. But then the definitions of these sets and of the product topology imply that
            \begin{equation*}
                F(\tilde x) - F(x) - \dual{x^*,\tilde x-x}_X < \alpha
                \leq G(x) - G(\tilde x),
            \end{equation*}
            contradicting \eqref{eq:convex:sum:fpg}. Hence $\interior C_1$ and $C_2$ are disjoint.
    \end{enumerate}
    We can thus apply \cref{lem:convex:eidelheit} to separate $C_1$ from $C_2$ by a hyperplane, i.e., to obtain a pair $(z^*,s)\in (X^*\times \R)\setminus\{(0,0)\} \cong (X\times \R)^*\setminus\{(0,0)\}$ and a $\lambda\in\R$ with
    \begin{subequations}
        \begin{align}
            \label{eq:convex:sum:hb1}
            \dual{z^*,\tilde x}_X + s(t-(F(x) -\dual{x^*,x}_X)) &\leq \lambda, \quad\tilde{x}\in \dom F, t\geq F(\tilde x) - \dual{x^*,\tilde x}_X,\\
            \label{eq:convex:sum:hb2}
            \dual{z^*,\tilde x}_X + s(G(x)-t) &\geq \lambda, \quad \tilde{x}\in \dom G, t\geq G(\tilde x).
        \end{align}
    \end{subequations}

    We now show that $s<0$. To show that $s=0$ is not possible, we apply \cref{thm:functan:hb-separation}\,\cref{item:functan:hb-separation:i} to separate $\interior C_1$ from $C_2$. Due to the structure of the product topology and since $x_0 \in \interior \dom F$, this shows that \eqref{eq:convex:sum:hb1} holds \emph{strictly} for $\tilde x=x_0$ and $t$ large enough. For $s=0$, we can thus combine \eqref{eq:convex:sum:hb1} and \eqref{eq:convex:sum:hb2} to obtain the contradiction
    \begin{equation*}
        \dual{z^*,x_0}_X < \lambda \leq \dual{z^*,x_0}_X.
    \end{equation*}
    If $s>0$, choosing $t>F(x) - \dual{x^*,x}_X$ makes the term in parentheses in \eqref{eq:convex:sum:hb1} strictly positive, and taking $t\to \infty$ with fixed $\tilde x$ leads to a contradiction to the boundedness by $\lambda$.

    Hence $s<0$, and \eqref{eq:convex:sum:hb1} with $t=F(\tilde x) - \dual{x^*,\tilde x}_X$ and \eqref{eq:convex:sum:hb2} with $t=G(\tilde x)$ imply that
    \begin{align*}
        F(\tilde x) - F(x) - \dual{x^*,\tilde x - x}_X &\geq s^{-1}(\lambda -\dual{z^*,\tilde x}_X) \quad\text{for all } \tilde x\in \dom F,\\
        G(x) - G(\tilde x) &\leq s^{-1}(\lambda -\dual{z^*,\tilde x}_X) \quad \text{for all }\tilde x\in \dom G.
    \end{align*}
    Taking $\tilde x = x\in \dom F\cap \dom G$ in both inequalities immediately yields that
    $\lambda = \dual{z^*,x}_X$. Hence, $y^* = s^{-1}z^*\in X^*$ is the desired functional with $(x^* - y^*)\in \partial F(x)$ and $y^* \in \partial G(x)$, i.e., $x^* \in \partial F(x) + \partial G(x)$.
\end{proof}
\begin{figure}
    \centering
    \begin{minipage}[t]{0.5\textwidth}
        \centering
        \begin{asy}
            unitsize(75,75);
            real xi = 1/2;

            real Fxi(real x) {return 1/2*x^2-xi*x;}
            real mG(real x) {return -abs(x);}
            real zeta(real x) {return -xi*x;}

            path Fp = graph(Fxi,-1,1);
            path Gp = graph(mG,-1,1);
            path Zp = graph(zeta,-1,1);

            draw(Fp);
            label("$F(x)-x^*\cdot x$",(-1,1),N);
            fill(Fp--(1,1)--cycle,lightfill);
            label("$C_1$",(0.5,0.5),S);

            draw(Gp);
            label("$-G(x)$",(-1,-1),W);
            fill(Gp--cycle,lightfill);
            label("$C_2$",(-0.25,-0.75),N);

            draw(Zp,dualline+shortdashed+linewidth(1.5));
            label("$-y^*\cdot x$",(0.8,-0.5),SE);

            dot((0,0));
            draw((0,-1.2)..(0,1.2),linewidth(0.5),Arrow); label("$t$",(0,1.2),E);
            draw((-1.2,0)..(1.2,0),linewidth(0.5),Arrow); label("$x$",(1.2,0),S);
        \end{asy}
        \caption{Illustration of the proof of \cref{thm:subdiff:sum} for $F(x) = \frac12|x|^2$, $G(x) = |x|$, and $x^*=\frac12 \in \partial(F+G)(0)$. The dashed line is the separating hyperplane $\{(x,t)\mid z^*\cdot x + st = \lambda\}$, i.e., $\lambda =0$, $z^*=-1$, $s=-2$ and hence $y^* = \frac12\in \partial G(0)$.}\label{fig:subdiff:sum1}
    \end{minipage}
    \hfill
    \begin{minipage}[t]{0.45\textwidth}
        \centering
        \begin{asy}
            unitsize(75,75);
            real Fxi(real x) {return -sqrt(x)-x/2;}
            path Fp = graph(Fxi,0.000001,1);

            path Gp = (-1,0)--(0,0)--(0,-1.5);
            path Zp = (0,-1.5)--(0,0.75);

            draw(Fp);
            fill((0,0.75)--Fp--(1,0.75)--cycle,lightfill);
            label("$C_1$",(0.5,0.5),S);

            draw(Gp);
            fill(Gp--(-1,-1.5)--cycle,lightfill);
            label("$C_2$",(-0.5,-1),N);

            draw(Zp,dualline+shortdashed+linewidth(1.5));

            dot((0,0));
        \end{asy}
        \caption{Illustration of the situation in \cref{ex:subdiff:sum}. Here the dashed separating hyperplane corresponds to the vertical line $\{(x,t)\mid x=0\}$ (i.e., $z^*=1$ and $s=0$), and hence $y^* \notin\R$.}\label{fig:subdiff:sum2}
    \end{minipage}
\end{figure}

The following example demonstrates that the inclusion is strict in general (although naturally the situation in infinite-dimensional vector spaces is nowhere near as obvious).

\begin{example}\label{ex:subdiff:sum}
    We take again $X=\R$ and $F:X\to\Rbar$ from \cref{ex:subdiff:empty}, i.e.,
    \begin{equation*}
        F(x) =
        \begin{cases}
            -\sqrt{x} & \text{if } x\geq 0,\\
            \infty & \text{if } x<0,
        \end{cases}
    \end{equation*}
    as well as $G(x) = \delta_{(-\infty,0]}(x)$. Both $F$ and $G$ are convex, and $0\in \dom F\cap \dom G$. In fact, $(F+G)(x) = \delta_{\{0\}}(x)$ and hence it is straightforward to verify that $\partial(F+G)(0)=\R$.

    On the other hand, we know from \cref{ex:subdiff:empty} and the argument leading to \eqref{eq:subdiff:indicator} that
    \begin{equation*}
        \partial F(0) = \emptyset,\qquad \partial G(0) = [0,\infty),
    \end{equation*}
    and hence that
    \begin{equation*}
        \partial F(0)+\partial G(0) = \emptyset \subsetneq \R = \partial (F+G)(0).
    \end{equation*}
    (As $F$ admits only a vertical tangent line at $x=0$, this example corresponds to the situation where $s=0$ in \eqref{eq:convex:sum:hb1}, cf.~\cref{fig:subdiff:sum2}.)
\end{example}

\begin{remark}
    \label{rem:convex:attouch-brezis}
    There exist alternative conditions that guarantee that the sum rule holds with equality. For example, if $X$ is a Banach space and $F$ and $G$ are in addition lower semicontinuous, this holds under the \term[condition!Attouch--Brézis]{Attouch--Brézis condition} that
    \begin{equation*}
        \bigcup_{\lambda\geq 0} \lambda \left(\dom F - \dom G\right) =:Z \text{ is a closed subspace of }X,
    \end{equation*}
    see \cite{AttouchBrezis}.
    (Note that this condition is not satisfied in \cref{ex:subdiff:sum} either, since in this case $Z=-\dom G= [0,\infty)$ which is closed but not a subspace.)

    It is not difficult to see that the condition $x_0 \in \interior(\dom F)\cap \dom G$ in the statement of \cref{lem:convex:subdiff_calc} implies the Attouch--Brézis condition. In fact, the latter allows us to generalize the condition to $x_0 \in \ri(\dom F)\cap \dom G$ where $\ri A$ for a set $A$ denotes the \term[interior!relative]{relative interior}: the interior of $A$ with respect to the smallest closed affine set that contains $A$. As an example, $\ri\{c\}=\{c\}$ for a point $c \in X$.
\end{remark}
By induction, we obtain from this sum rules for an arbitrary (finite) number of functionals (where $x_0$ has to be an interior point of all but one effective domain). A chain rule for linear operators can be proved similarly.
\begin{theorem}[chain rule]\label{thm:convex:chain}
    Let $X,Y$ be Banach spaces, $K\in \linear(X;Y)$, $F:Y\to\Rbar$ be proper, convex, and lower semicontinuous, and $x\in \dom (F\circ K)$. Then,
    \begin{equation*}
        \partial (F\circ K)(x) \supset K^*\partial F(Kx) \defeq \setof{K^*y^*}{y^* \in \partial F(Kx)}
    \end{equation*}
    with equality if there exists an $x_0\in X$ with $Kx_0\in \interior(\dom F)$.
\end{theorem}
\begin{proof}
    The inclusion is again a direct consequence of the definition: If $y^* \in \partial F(Kx)\subset Y^*$, we in particular have for all $\tilde y = K\tilde x\in Y$ with $\tilde x\in X$ that
    \begin{equation*}
        F(K\tilde x) - F(Kx) \geq \dual{y^*,K\tilde x - Kx}_Y = \dual{K^*y^*,\tilde x-x}_X,
    \end{equation*}
    i.e., $x^* \defeq  K^*y^*\in \partial (F\circ K) \subset X^*$.

    To show the claimed equality under the additional assumption, let $x\in\dom (F\circ K)$ and $x^*\in \partial(F\circ K)(x)$, i.e.,
    \begin{equation*}\label{eq:convex:chain:subdiff}
        F(Kx) + \dual{x^*,\tilde x-x}_X \leq F(K\tilde x) \quad\text{for all }\tilde x\in X.
    \end{equation*}
    We now construct a $y^*\in\partial F(Kx)$ with $x^* = K^*y^*$ by applying the sum rule to%
    \footnote{This technique of \enquote{lifting}\index{lifting|infn} a problem to a product space in order to separate operators is also useful in many other contexts.}
    \begin{equation*}
        H:X\times Y\to\Rbar, \qquad  (x,y) \mapsto F(y) + \delta_{\graph K}(x,y).
    \end{equation*}
    Since $K$ is linear and continuous, $\graph K$ is convex and closed, and hence $\delta_{\graph K}$ is convex and lower semicontinuous. Furthermore, $Kx\in\dom F$ by assumption and thus $(x,Kx) \in \dom H$.

    We begin by showing that $x^*\in\partial (F\circ K)(x)$ if and only if $(x^*,0)\in\partial H(x,Kx)$. First, let $(x^*, 0)\in \partial H(x,Kx)$. Then we have for all $\tilde x\in X,\tilde y\in Y$ that
    \begin{equation*}
        \dual{x^*,\tilde x- x}_X + \dual{0,\tilde y -Kx}_Y \leq F(\tilde y) - F(Kx) + \delta_{\graph K}(\tilde x,\tilde y) - \delta_{\graph K}(x,Kx).
    \end{equation*}
    In particular, this holds for all $\tilde y \in \range(K)=\setof{K\tilde x}{\tilde x\in X}$. By $\delta_{\graph K}(\tilde x,K\tilde x) = 0$ we thus obtain that
    \begin{equation*}
        \dual{x^*,\tilde x- x}_X \leq F(K \tilde x) - F(Kx) \quad\text{for all }\tilde x\in X,
    \end{equation*}
    i.e., $x^* \in \partial(F\circ K)(x)$. Conversely, let $x^*\in \partial(F\circ K)(x)$. Since $\delta_{\graph K}(x,Kx) = 0$ and $\delta_{\graph K}(\tilde x,\tilde y)\geq 0$, it then follows for all $\tilde x \in X$ and $\tilde y \in Y$ that
    \begin{equation*}
        \begin{aligned}
            \dual{x^*,\tilde x- x}_X + \dual{0,\tilde y -Kx}_Y
            &= \dual{x^*,\tilde x- x}_X\\
            &\leq F(K \tilde x) - F(Kx) + \delta_{\graph K}(\tilde x,\tilde y)  - \delta_{\graph K}(x,Kx)\\
            &= F( \tilde y) - F(Kx) + \delta_{\graph K}(\tilde x,\tilde y) - \delta_{\graph K}(x,Kx),
        \end{aligned}
    \end{equation*}
    where we have used that the last equality holds trivially as $\infty=\infty$ for $\tilde y \neq K\tilde x$. Hence, $(x^*,0)\in \partial H(x,Kx)$.

    We now consider $\tilde F:X\times Y\to \Rbar$, $(x,y)\mapsto F(y)$, and  $(x_0,Kx_0)\in \graph K=\dom\delta_{\graph K}$.
    Since $Kx_0\in \interior(\dom F) \subset Y$ by assumption, $(x_0,Kx_0)\in \interior(\dom \tilde F)=X\times \interior(\dom F)\subset X\times Y$ as well.
    We can thus apply \cref{thm:subdiff:sum} to obtain
    \begin{equation*}
        (x^*,0) \in \partial H(x,Kx) = \partial \tilde F(x,Kx) + \partial \delta_{\graph K}(x,Kx),
    \end{equation*}
    i.e., $(x^*,0)=(x_1^*,y_1^*)+(x_2^*,y_2^*)$ for some $(x_1^*,y_1^*)\in \partial \tilde F(x,Kx)$ and $(x_2^*,y_2^*)\in\partial \delta_{\graph K}(x,Kx)$.

    Finally, we \enquote{collapse} these subdifferentials back to the individual spaces to obtain the desired characterization.
    First, we have $(x_1^*,y_1^*)\in \partial \tilde F(x,Kx)$ if and only if
    \begin{equation*}
        \dual{x_1^*,\tilde x -x}_X + \dual{y_1^*,\tilde y-Kx}_Y \leq F(\tilde y)- F(Kx)\quad\text{for all }\tilde x\in X,\tilde y \in Y.
    \end{equation*}
    Fixing in turn $\tilde x = x$ and $\tilde y=Kx$ implies that $y_1^*\in\partial F(Kx)$  and $x_1^* =0$, respectively.
    Second, $(x_2^*,y_2^*)\in \partial \delta_{\graph K}(x,Kx)$ if and only if
    \begin{equation*}
        \dual{x_2^*,\tilde x -x}_X + \dual{y_2^*,\tilde y-Kx}_Y \leq 0 \quad\text{for all }(\tilde x,\tilde y)\in \graph K,
    \end{equation*}
    i.e., for all $\tilde x\in X$ and $\tilde y = K\tilde x$. Therefore,
    \begin{equation*}
        \dual{x_2^*+K^*y_2^*,\tilde x - x}_X \leq 0 \quad\text{for all }\tilde x\in X
    \end{equation*}
    and hence $x_2^*=-K^*y_2^*\in X^*$. Together we obtain
    \begin{equation*}
        (x^*,0) = (0,y_1^*) + (-K^*y_2^*,y_2^*),
    \end{equation*}
    which implies that $y_1^* = -y_2^*$ and thus that $x^* = -K^*y_2^* = K^*y_1^*$ with $y_1^*\in\partial F(Kx)$ as claimed.
\end{proof}
The condition for equality in particular holds if $K$ is surjective and $\dom F$ has nonempty interior. Again, the inequality can be strict.
\begin{example}
    Here we take $X=Y=\R$ and again $F:X\to\Rbar$ from \cref{ex:subdiff:empty,ex:subdiff:sum} as well as
    \begin{equation*}
        K:\R\to\R,\qquad Kx = 0.
    \end{equation*}
    Clearly, $(F\circ K)(x) = 0$ for all $x\in \R$ and hence $\partial(F\circ K)(x) = \{0\}$ by \cref{thm:convex:gateaux}. On the other hand, $\partial F(0)=\emptyset$ by \cref{ex:subdiff:empty} and hence
    \begin{equation*}
        K^*\partial F(Kx) = K^* \partial F(0) = \emptyset \subsetneq \{0\}.
    \end{equation*}
    (Note the problem: $K^*$ is far from surjective, and $\range K\cap \interior(\dom F)=\emptyset$.)
\end{example}

We can also obtain a chain rule when the \emph{inner} mapping is nondifferentiable.

\begin{theorem}
    \label{thm:convex:increasing-post}
    Let $F: X \to \R$ be convex and $\phi: \R \to \R$ be convex, increasing, and differentiable. Then $\phi \circ F$ is convex, and for all $x\in X$,
    \[
        \subdiff(\phi \circ F)(x)
        = \phi'(F(x))\subdiff F(x) = \setof{ \phi'(F(x))x^*}{x^* \in \subdiff F(x)}.
    \]
\end{theorem}

\begin{proof}
    First, the convexity of $\phi \circ F$ follows from \cref{lem:convex:func}\,\ref{lem:convex:func:iii}.
    To calculate the subdifferential, we fix $x \in X$ and observe from \cref{thm:convex:cont} that $\phi$ is Lipschitz continuous with some constant $L$ near $F(x) \in \interior (\dom \phi) = \R$.
    Furthermore, by \cref{lem:convex:direct}\,\ref{lem:convex:direct:iii} the directional derivative $F'(x; h)$ exists and is finite in $x\in X=\interior (\dom F)$ for every $h\in X$.
    Thus, for any $h \in X$,
    \begin{equation*}
        \begin{aligned}[t]
            (\phi \circ F)'(x; h)
            &= \lim_{t \downto 0} \frac{(\phi \circ F)(x+th)-(\phi \circ F)(x)}{t}
            \\
            & = \lim_{t \downto 0} \frac{\phi(F(x+th))-\phi(F(x)+tF'(x; h))}{t}
            \\
            \MoveEqLeft[-1]
            +\lim_{t \downto 0} \frac{\phi(F(x)+tF'(x; h))-\phi(F(x))}{t}
            \\
            & \le \lim_{t \downto 0} L\left| \frac{F(x+th)-F(x)}{t}-F'(x;h)\right| + \phi'(F(x); F'(x; h))
            \\
            & = \phi'(F(x); F'(x; h)),
        \end{aligned}
    \end{equation*}
    where we have again used the directional differentiability of $F$ in the last step.
    Similarly, we prove the converse inequality using $\phi(t_1)-\phi(t_2) \geq -L|t_1-t_2|$ for all $t_1,t_2$ sufficiently close to $F(x)$. Hence
    \begin{equation*}
        (\phi \circ F)'(x; h) = \phi'(F(x); F'(x; h))=\phi'(F(x))F'(x;h)
    \end{equation*}
    by the differentiability of $\phi$.

    Now \cref{lem:convex:equiv} yields that
    \begin{equation*}
            \partial(\phi \circ F)(x)
             =
            \setof{ z^* \in X^*}{\dualprod{z^*}{h}_X \le \phi'(F(x))F'(x; h) \text{ for all } h \in X}.
    \end{equation*}
    Since $\phi:\R\to\R$ is monotone and differentiable, $\phi'(F(x))\geq 0$. Hence if $\phi'(F(x))>0$, we can set $x^*\defeq \phi'(F(x))^{-1} z^*\in X^*$; otherwise $z^*=0$ is the only element of $\partial(\phi\circ F)(x)$.
    In either case, we can write
    \begin{equation*}
        \partial(\phi \circ F)(x)
        =\setof{\phi'(F(x))x^*}{\dualprod{x^*}{h}_X \le F'(x; h) \text{ for all } h \in X}
    \end{equation*}
    so that the claim follows by \cref{lem:convex:equiv}.
\end{proof}

\begin{remark}
    The differentiability assumption on $\phi$ in \cref{thm:convex:increasing-post} is not necessary, but the proof is otherwise much more involved and demands the support functional machinery of \cref{sec:clarke:support}. See also \cite[Section D.4.3]{Hiriart:2001} for a version with set-valued $F$ in finite dimensions.
\end{remark}

\bigskip

The Fermat principle together with the sum rule yields the following characterization of minimizers of convex functionals under convex constraints.
\begin{corollary}
    Let $U\subset X$ be nonempty, convex, and closed, and let $F:X\to\Rbar$ be proper, convex, and lower semicontinuous. If there exists an $x_0 \in \interior U \cap \dom F$, then $\bar x\in U$ solves
    \begin{equation*}
        \min_{x\in U} F(x)
    \end{equation*}
    if and only if $0 \in \partial F(\bar x) + N_U(\bar x)$ or, in other words, if there exists an $x^*\in X^*$ with
    \begin{equation}\label{eq:convex:kkt}
        \left\{\begin{aligned}
                &x^*\in \partial F(\bar x),\\
                &\dual{x^*,\tilde x - \bar x}_X \geq 0 \quad\text{for all } \tilde x \in U.
        \end{aligned}\right.
    \end{equation}
\end{corollary}
\begin{proof}
    Due to the assumptions on $F$ and $U$, we can apply \cref{thm:convex:fermat} to $J\defeq F+\delta_U$. Furthermore, since $x_0 \in \interior U = \interior(\dom \delta_U)$, we can also apply \cref{thm:subdiff:sum}. Hence $F$ has a minimum in $\bar x$ if and only if
    \begin{equation*}
        0\in \partial J(\bar x) = \partial F(\bar x) + \partial \delta_U(\bar x).
    \end{equation*}
    Together with the characterization of subdifferentials of indicator functionals as normal cones, this yields \eqref{eq:convex:kkt}.
\end{proof}
If $F:X\to\R$ is Gateaux differentiable (and hence finite-valued), \eqref{eq:convex:kkt} coincide with the classical \term[condition!Karush--Kuhn--Tucker]{Karush--Kuhn--Tucker conditions}; the existence of an interior point $x_0\in \interior U$ is related to a \term[condition!Slater]{Slater condition} in nonlinear optimization needed to show existence of the Lagrange multiplier $x^*$ for inequality constraints.

\chapter{Fenchel duality}\label{chap:fenchel}

One of the main tools in convex optimization is \term[duality!Fenchel]{duality}: Any convex optimization problem can be related to a \term[problem!dual]{dual problem}, and the joint study of both problems yields additional information about the solution.
Our main objective in this chapter, the \emph{Fenchel--Rockafellar duality theorem}, will be our main tool for deriving explicit optimality conditions as well as numerical algorithms for convex minimization problems that can be expressed as the sum of (simple) functionals.

\section{Fenchel conjugates}

Let $X$ be a normed vector space and $F:X\to\Rbar$ be proper but not necessarily convex. We then define the \term[conjugate!Fenchel]{Fenchel conjugate} (or \term[conjugate!convex]{convex conjugate}) of $F$ as
\begin{equation*}
    F^*:X^*\to \Rbar,\qquad F^*(x^*) = \sup_{x\in X}\,\left\{ \dual{x^*,x}_X - F(x)\right\}.
\end{equation*}
(Since $\dom F = \emptyset$ is excluded, we have that $F^*(x^*)>-\infty$ for all $x^*\in X^*$, and hence the definition is meaningful.)
An alternative interpretation is that $F^*(x^*)$ is the (negative of the) affine part of the tangent to $F$ (in the point $x$ at which the supremum is attained) with slope $x^*$, see \cref{fig:convex:fenchel}.
\Cref{lem:convex:func}\,\ref{lem:convex:func:v} and \cref{lem:variation:wlsc}\,\ref{lem:variation:wlsc:v} immediately imply that $F^*$ is always convex and weakly-$*$ lower semicontinuous (as long as $F$ is indeed proper). If $F$ is bounded from below by an affine functional (which is always the case if $F$ is proper, convex, and lower semicontinuous by \cref{lem:convex:gamma}), then $F^*$ is proper as well.
Finally, the definition directly yields the \term[inequality!Fenchel--Young]{Fenchel--Young inequality}
\begin{equation}\label{eq:convex:fenchel-young}
    \dual{x^*,x}_X \leq F(x) + F^*(x^*)\qquad\text{for all }x\in X, x^*\in X^*.
\end{equation}

\begin{figure}
    \centering
    \begin{subfigure}[t]{0.495\textwidth}
        \centering
        \begin{asy}
            unitsize(75,75);
            draw((-1.2,0)..(1.2,0),linewidth(0.5),Arrow); label("$x$",(1.2,0),S);
            draw((0,-1.2)..(0,1.2),linewidth(0.5),Arrow); label("$t$",(0,1.2),E);

            real xi = 1.0;
            real F(real x) {return x^2/2 + abs(x-0.5);}
            real Fs(real x) {return xi*x - F(x);}

            path Fsg = graph(Fs,-0.25,1.1);
            draw(Fsg, primalline+linewidth(1.5));
            label("$x^*\cdot x - F(x)$",(-0.3,-1.1),W);

            dot((0.5,Fs(0.5)));
            dot((0,Fs(0.5)));
            draw((0,Fs(0.5))--(0.5,Fs(0.5)),linewidth(0.5)+dotted);
            label("$F^*(x^*)$",(0,Fs(0.5)),W);
        \end{asy}
        \caption{$F^*(x^*)$ as maximizer of $x^*\cdot x-F(x)$}\label{fig:convex:fenchel:def}
    \end{subfigure}
    \hfill
    \begin{subfigure}[t]{0.495\textwidth}
        \centering
        \begin{asy}
            unitsize(75,75);
            draw((-1.2,0)..(1.2,0),linewidth(0.5),Arrow); label("$x$",(1.2,0),S);
            draw((0,-1.2)..(0,1.2),linewidth(0.5),Arrow); label("$t$",(0,1.2),E);

            real F(real x) {return x^2/2 + abs(x-0.5);}
            path Fg = graph(F,-0.25,1.1);
            draw(Fg, primalline+linewidth(1.5));
            label("$F(x)$",(-0.25,0.8),W);
            dot((0.5,F(0.5)));

            real xi = 1.0;
            real Fs(real x) {return xi*x - F(x);}
            real Fxi(real x) {return -Fs(0.5) + xi*x;}
            path Fxig = graph(Fxi,-0.54,1.1);
            draw(Fxig, dualline+shortdashed+linewidth(1.5));
            label("$-F^*(x^*)+x^*\cdot x$",(-0.33,-1),W);
            label("$(\bar x, F(\bar x))$",(0.5,F(0.5)),ESE);

            dot((0,-Fs(0.5)));
            label("$-F^*(x^*)$",(0,-Fs(0.5)),W);
        \end{asy}
        \caption{Alternative interpretation: $-F^*(x^*)$ as offset for tangent to $F$ with given slope $x^*$. Note that in this case, $x^* \in \partial F(\bar x)$ and $-F^*(x^*)+x^*\cdot \bar x =F(\bar x)$.}
        \label{fig:convex:fenchel:alt}
    \end{subfigure}
    \caption{Geometrical illustration of the Fenchel conjugate.}
    \label{fig:convex:fenchel}
\end{figure}
Another simple consequence of the definition is the following order-reversing property of Fenchel conjugation.
\begin{lemma}\label{lem:convex:fenchel-monotone}
    If $F:X\to\Rbar$ and $G:X\to\Rbar$ satisfy $F\leq G$, then $F^*\geq G^*$.
\end{lemma}

If $X$ is not reflexive, we can similarly define for (weakly-$*$ lower semicontinuous) $F:X^*\to \R$ the \term[preconjugate, Fenchel]{Fenchel preconjugate}
\begin{equation*}
    F_*:X\to \Rbar,\qquad F_*(x) = \sup_{x^*\in X^*}\,\left\{ \dual{x^*,x}_X - F(x^*) \right\}.
\end{equation*}
The point of this convention is that even in nonreflexive spaces, the \term[biconjugate, Fenchel]{biconjugate}
\begin{equation*}
    F^{**}: X\to \Rbar, \qquad F^{**}(x) = (F^*)_* (x)
\end{equation*}
is again defined on $X$ (rather than $X^{**}\supset X$). For reflexive spaces, of course, we have $F^{**}=(F^*)^*$.
Intuitively, $F^{**}$ is the convex envelope of $F$, which by \cref{lem:convex:gamma} coincides with $F$ itself if $F$ is convex.
\begin{theorem}[Fenchel--Moreau--Rockafellar]\label{thm:convex:moreau}\index{theorem!Fenchel--Moreau--Rockafellar}
    Let $F:X\to\Rbar$ be proper. Then,
    \begin{enumerate}
        \item \label{thm:convex:moreau:i}
            $F^{**}\leq F$;
        \item \label{thm:convex:moreau:ii}
            $F^{**} = F^{\Gamma}$;
        \item \label{thm:convex:moreau:iii}
            $F^{**} = F$ if and only if $F$ is convex and lower semicontinuous.
    \end{enumerate}
\end{theorem}
\begin{proof}
    For \ref{thm:convex:moreau:i}, we take the supremum over all $x^*\in X^*$ in the Fenchel--Young inequality \eqref{eq:convex:fenchel-young} and obtain that
    \begin{equation*}
        F(x) \geq \sup_{x^*\in X^*} \left\{\dual{x^*,x}_X -F^*(x^*)\right\} = F^{**}(x).
    \end{equation*}

    For \ref{thm:convex:moreau:ii}, we first note that $F^{**}$ is convex and lower semicontinuous by definition as a Fenchel conjugate as well as proper by \ref{thm:convex:moreau:i}. Hence, \cref{lem:convex:gamma} yields that
    \begin{equation*}
        F^{**}(x) = (F^{**})^\Gamma(x) = \sup\setof{a(x)}{a:X\to\R\text{ continuous affine with } a\leq F^{**}}.
    \end{equation*}
    We now show that we can replace $F^{**}$ with $F$ on the right-hand side. For this, let $a(x)=\dual{x^*,x}_X-\alpha$ with arbitrary $x^*\in X^*$ and $\alpha\in\R$.
    If $a\leq F^{**}$, then \ref{thm:convex:moreau:i} implies that $a\leq F$. Conversely, if $a\leq F$, we have that  $\dual{x^*,x}_X-F(x)\leq \alpha$ for all $x\in X$, and taking the supremum over all $x\in X$ yields that $\alpha\geq F^{*}(x^*)$. By definition of $F^{**}$, we thus obtain that
    \begin{equation*}
        a(x) = \dual{x^*,x}_X - \alpha \leq \dual{x^*,x}_X - F^*(x^*) \leq F^{**}(x) \quad\text{for all }x\in X,
    \end{equation*}
    i.e., $a\leq F^{**}$.

    Statement \ref{thm:convex:moreau:iii} now directly follows from \ref{thm:convex:moreau:ii} and \cref{lem:convex:gamma}.
\end{proof}
\begin{remark}\label{rem:fenchel:weak-star}
    Continuing from \cref{rem:convex:weak-star-hull}, we can adapt the proof of \cref{thm:convex:moreau} to proper functionals $F:X^*\to \Rbar$ to show that $F=(F_*)^*$ if and only if $F$ is convex and weakly-$*$ lower semicontinuous.
\end{remark}

We again consider some relevant examples.
\begin{example}\label{ex:convex:fenchel}~
    \begin{enumerate}
        \item \label{ex:convex:fenchel:ii} Let $\B_X$ be the unit ball in the normed vector space $X$ and take $F=\delta_{\B_X}$. Then we have for any $x^*\in X^*$ that
            \begin{equation*}
                (\delta_{\B_X})^*(x^*) = \sup_{x\in X}\left\{ \dual{x^*,x}_X - \delta_{\B_X}(x) \right\}= \sup_{\norm{x}_{X}\leq 1}\left\{  \dual{x^*,x}_X \right\} = \norm{x^*}_{X^*}.
            \end{equation*}
            Similarly, one shows using the definition of the Fenchel preconjugate and \cref{cor:functan:norm_dual} that $(\delta_{\B_{X^*}})_*(x) = \norm{x}_X$.

        \item \label{ex:convex:fenchel:iii} Let $X$ be a normed vector space and take $F(x) = \norm{x}_X$. We now distinguish two cases for a given $x^*\in X^*$.
            \begin{enumerate}[label={Case} \arabic*:, align=left]
                \item $\norm{x^*}_{X^*} \leq 1$. Then it follows from \eqref{eq:functan:cs_banach} that $\dual{x^*,x}_X- \norm{x}_X\leq 0$ for all $x\in X$. Furthermore, $\dual{x^*,0} = 0 = \norm{0}_X$, which implies that
                    \begin{equation*}
                        F^*(x^*) = \sup_{x\in X}\left\{ \dual{x^*,x}_X -\norm{x}_X \right\}= 0.
                    \end{equation*}
                \item $\norm{x^*}_{X^*} >1$. Then by definition of the dual norm, there exists an $x_0\in X$ with $\dual{x^*,x_0}_X > \norm{x_0}_X$. Hence, taking $t\to\infty$ in
                    \begin{equation*}
                        0<t(\dual{x^*,x_0}_X - \norm{x_0}_X) = \dual{x^*,tx_0}_X - \norm{tx_0}_X \leq F^*(x^*)
                    \end{equation*}
                    yields $F^*(x^*) = \infty$.
            \end{enumerate}
            Together we obtain that $F^* = \delta_{\B_{X^*}}$.
            As above, a similar argument shows that $(\norm{\cdot}_{X^*})_* = \delta_{\B_X}$.
    \end{enumerate}
\end{example}

We can generalize \cref{ex:convex:fenchel}\,\ref{ex:convex:fenchel:iii} to powers of norms.
\begin{lemma}\label{lem:convex:power-conjugate}
    Let $X$ be a normed vector space and $F(x) \defeq \frac1p\norm{x}_X^p$ for $p\in (1,\infty)$. Then $F^*(x^*) = \frac1q\norm{x^*}_{X^*}^q$ for $q\defeq\frac{p}{p-1}$.
\end{lemma}
\begin{proof}
    We first consider the scalar function $\phi(t) \defeq\frac1p|t|^p$
    and compute the Fenchel conjugate $\phi^*(s)$ for $s\in \R$. By the choice of $p$ and $q$, we then can write $\frac1q = 1-\frac1p$ as well as $|s|^q = \sign(s)s |s|^{1/(p-1)} = |\sign(s)|s|^{1/(p-1)}|^p$ for any $s\in \R$ and therefore obtain
    \begin{equation*}
        \frac1q |s|^q = \left(\sign(s) |s|^{1/(p-1)}\right)s - \frac1p \left|\sign(s)|s|^{1/(p-1)}\right|^p \leq \sup_{t\in \R} \left\{ts - \frac1p|t|^p\right\} \leq \frac1q |s|^q,
    \end{equation*}
    where we have used the classical Young inequality  $ts \leq \frac1p |t|^p +\frac1q|s|^q$ in the last step.
    This shows that $\phi^*(s) = \frac1q |s|^q$.\footnote{Which is how the Fenchel--Young inequality got its name.}

    We now write using the definition of the norm in $X^*$ that
    \begin{equation*}
        \begin{aligned}
            F^*(x^*) &= \sup_{x\in X} \left\{\dual{x^*,x}_X - \frac1p\norm{x}_X^p\right\}
            = \sup_{t\geq 0} \left\{\sup_{\norm{x}_X = 1} \left\{\dual{x^*,tx}_X - \frac1p\norm{tx}_X^p\right\}\right\}\\
            &= \sup_{t\geq 0} \left\{t\norm{x^*}_{X^*} - \frac1p |t|^p\right\}
            = \frac1q \left|\norm{x^*}_{X^*}\right|^q
        \end{aligned}
    \end{equation*}
    since $\phi$ is even and the supremum over all $t\in \R$ is thus attained for $t\geq 0$.
\end{proof}

As for convex subdifferentials, Fenchel conjugates of integral functionals can be computed pointwise.
\begin{theorem}\label{thm:lebesgue:fenchel}
    Let $f:\R\to\Rbar$ be measurable, proper and lower semicontinuous, and let $F:L^p(\Omega)\to\Rbar$ with $1\leq p <\infty$ be defined as in \cref{lem:lebesgue:lsc}. Then we have for $q=\frac{p}{p-1}$ that
    \begin{equation*}
        F^*:L^q(\Omega) \to\Rbar, \qquad F^*(u^*) =
        \int_\Omega f^*(u^*(x))\,dx .
    \end{equation*}
\end{theorem}
\begin{proof}
    We argue similarly as in the proof of \cref{thm:lebesgue:subdiff}, with some changes that are needed since measurability of $f\circ u$ does not immediately imply that of $f^*\circ u^*$. Let $u^*\in L^q(\Omega)$ be arbitrary and consider for all $x\in \Omega$ the functions
    \begin{align*}
        \phi(x) &\defeq \sup_{t\in\R} \left\{tu^*(x) - f(t)\right\} = f^*(u^*(x)),
        \shortintertext{as well as for $n\in \N$}
        \phi_n(x) &\defeq \sup_{|t|\leq n} \left\{tu^*(x) - f(t)\right\} \leq f^*(u^*(x)).
    \end{align*}
    By a measurable selection theorem (\cite[Theorem VIII.1.2]{Ekeland:1999a}), the pointwise supremum in the definition of $\phi_n$ is attained at some $t^*_x$ for almost every $x\in \Omega$ and defines a measurable mapping $x\mapsto u_n(x)\defeq t^*_x$ with $\norm{u_n}_{L^\infty}\leq n$. This also implies that $\phi_n = u_n \cdot u^*-f\circ u_n$ is measurable.
    Furthermore, by assumption there exists a $t_0\in \dom f$, and hence
    $u_0 \defeq  t_0u^*(x) - f(t_0)$ is measurable and satisfies $u_0 \leq \phi_n(x)$ for all $n\geq |t_0|$.
    Finally, by construction, $\phi_n(x)$ is monotonically increasing and converges to $\phi(x)$ for all $x\in \Omega$. The sequence $\{\phi_n-u_0\}_{n\in\N}$ of functions is thus measurable and nonnegative, and the monotone convergence theorem yields that
    \begin{equation*}
        \int_\Omega \phi(x) - u_0(x) \,dx = \int_\Omega \sup_{n\in\N} \phi_n(x)-u_0(x)\,dx =
        \sup_{n\in\N} \int_\Omega \phi_n(x)-u_0(x)\,dx.
    \end{equation*}
    Hence the pointwise limit $\phi=f^*\circ u^*$ is measurable as well.

    The measurable selection theorem also yields that
    \begin{equation*}
        \begin{aligned}[t]
            \int_\Omega f^*(u^*(x))\,dx &= \sup_{n\in\N} \int_\Omega \sup_{|t|\leq n}\left\{ tu^*(x)-f(t)\right\}\,dx \\
            &= \sup_{n\in\N} \int_\Omega u^*(x)u_n(x)-f(u_n(x))\,dx \\
            &\leq \sup_{u\in L^p(\Omega)} \int_\Omega u^*(x)u(x)-f(u(x))\,dx = F^*(u^*),
        \end{aligned}
    \end{equation*}
    since $u_n \in L^\infty(\Omega)\subset L^p(\Omega)$ for all $n\in\N$.

    For the converse inequality, we can now proceed as in the proof of \cref{thm:lebesgue:subdiff}. For any $u\in L^p(\Omega)$ and $u^*\in L^q(\Omega)$, we have by the Fenchel--Young inequality \eqref{eq:convex:fenchel-young} applied to $f$ and $f^*$ that
    \begin{equation*}
        f(u(x)) + f^*(u^*(x)) \geq u^*(x)u(x)\quad \text{for almost every }x\in \Omega.
    \end{equation*}
    Since both sides are measurable, this implies that
    \begin{equation*}
        \int_\Omega f^*(u^*(x))\,dx \geq \int_\Omega u^*(x)u(x) - f(u(x))\,dx,
    \end{equation*}
    and taking the supremum over all $u\in L^p(\Omega)$ yields the claim.
\end{proof}
\begin{remark}
    A similar representation can be shown for vector-valued and spatially-dependent integrands $f:\Omega\times \R\to \R^m$ under stronger assumptions; see, e.g., \cite[Corollary 3C]{Rockafellar:1976}.
\end{remark}

\bigskip

Fenchel conjugates satisfy a number of useful calculus rules, which follow directly from the properties of the supremum.
\begin{lemma}\label{lem:convex:fenchel_calc}
    Let $F:X\to\Rbar$ be proper. Then,
    \begin{enumerate}
        \item \label{lem:convex:fenchel_calc:i}
            $(\alpha F)^* = \alpha F^*\circ (\alpha^{-1} \Id)$ for any $\alpha >0$;
        \item \label{lem:convex:fenchel_calc:ii}
            $(F(\cdot + x_0) + \dual{x_0^*,\cdot}_X)^* = F^*(\cdot - x_0^*) - \dual{\cdot - x_0^*,x_0}_X$ for all $x_0\in X$, $x_0^*\in X^*$;
        \item \label{lem:convex:fenchel_calc:iii}
            $(F\circ K)^* = F^* \circ K^{-*}$ for continuously invertible $K\in \linear(Y;X)$, where $K^{-*}\defeq(K^{-1})^*$.
    \end{enumerate}
\end{lemma}
\begin{proof}
    \emph{\ref{lem:convex:fenchel_calc:i}:} For any $\alpha >0$, we have that
    \begin{equation*}
        (\alpha F)^*(x^*) =  \sup_{x\in X} \left\{\alpha\dual{\alpha^{-1}x^*,x}_X - \alpha F(x)\right\}
        = \alpha \sup_{x\in X} \left\{\dual{\alpha^{-1}x^*,x}_X - F(x)\right\}
        = \alpha F^* (\alpha^{-1} x^*).
    \end{equation*}

    \emph{\ref{lem:convex:fenchel_calc:ii}:} Since $\setof{x+x_0}{x\in X}=X$, we have that
    \begin{equation*}
        \begin{aligned}
            (F(\cdot + x_0) + \dual{x_0^*,\cdot}_X)^*(x^*) &= \sup_{x\in X}\left\{ \dual{x^*,x}_X - F(x +x_0) - \dual{x_0^*,x}_X\right\}\\
            &= \sup_{x\in X} \left\{ \dual{x^*-x_0^*,x+x_0}_X - F(x +x_0)\right\} - \dual{x^*-x_0^*,x_0}_X\\
            &= \sup_{\tilde x=x+x_0,x\in X} \left\{ \dual{x^*-x^*_0,\tilde x}_X - F(\tilde x)\right\} - \dual{x^*-x_0^*,x_0}_X\\
            &= F^*(x^*-x_0^*) - \dual{x^*-x_0^*,x_0}_X.
        \end{aligned}
    \end{equation*}

    \emph{\ref{lem:convex:fenchel_calc:iii}:} Since $X=\range K$, we have that
    \begin{equation*}
        \begin{aligned}[b]
            (F\circ K)^*(y^*) &= \sup_{y\in Y}\left\{ \dual{y^*,K^{-1}Ky}_Y - F(Ky)\right\}\\
            & = \sup_{x=Ky,y\in Y} \left\{\dual{K^{-*}y^*,x}_X - F(x)\right\} = F^*(K^{-*}y^*).
        \end{aligned}
        \qedhere
    \end{equation*}
\end{proof}

There are some obvious similarities between the definitions of the Fenchel conjugate and of the subdifferential, which yield the following very useful property that plays the role of a \enquote{convex inverse function theorem}\index{theorem!inverse function!convex}. (See also \cref{fig:convex:fenchel:alt} and compare \cref{fig:subdiff:abs,fig:subdiff:ind}.)
\begin{lemma}[Fenchel--Young]\label{lem:convex:fenchel-young}\index{lemma!Fenchel--Young}
    Let $F:X\to\Rbar$ be proper, convex, and lower semicontinuous. Then the following statements are equivalent for any $x\in X$ and $x^*\in X^*$:
    \begin{enumerate}
        \item \label{lem:convex:fenchel-young:i}
            $\dual{x^*,x}_X = F(x) + F^*(x^*)$;
        \item \label{lem:convex:fenchel-young:ii}
            $x^*\in\partial F(x)$;
        \item \label{lem:convex:fenchel-young:iii}
            $x\in\partial F^*(x^*)$.
    \end{enumerate}
\end{lemma}
\begin{proof}
    If \ref{lem:convex:fenchel-young:i} holds, the definition of $F^*$ as a supremum immediately implies that
    \begin{equation}\label{eq:convex:fy1}
        \dual{x^*,x}_X - F(x) = F^*(x^*) \geq  \dual{x^*,\tilde x}_X - F(\tilde x) \qquad\text{for all }\tilde x\in X,
    \end{equation}
    which again by definition is equivalent to \ref{lem:convex:fenchel-young:ii}. Conversely, taking the supremum over all $\tilde x\in X$ in \eqref{eq:convex:fy1} yields
    \begin{equation*}
        \dual{x^*,x}_X \geq F(x) + F^*(x^*),
    \end{equation*}
    which together with the Fenchel--Young inequality \eqref{eq:convex:fenchel-young} leads to \ref{lem:convex:fenchel-young:i}.

    Similarly, \ref{lem:convex:fenchel-young:i} in combination with \cref{thm:convex:moreau} implies that
    \begin{equation*}
        \dual{x^*,x}_X - F^*(x^*) = F(x) = F^{**}(x) \geq \dual{\tilde x^*,x}-F^*(\tilde x^*) \qquad\text{for all }\tilde x^*\in X^*,
    \end{equation*}
    yielding as above the equivalence of \ref{lem:convex:fenchel-young:i} and \ref{lem:convex:fenchel-young:iii}.
\end{proof}
\begin{remark}\label{rem:convex:f-nonconvex}
    If $F$ is not convex, the above proof shows that we still have the equivalence \ref{lem:convex:fenchel-young:i}\,$\Leftrightarrow$\,\ref{lem:convex:fenchel-young:ii}. Furthermore since always $F^{**}\leq F$ by \cref{thm:convex:moreau}\,\ref{thm:convex:moreau:i}, it still holds that \ref{lem:convex:fenchel-young:i}\,$\implies$\,\ref{lem:convex:fenchel-young:iii}. However, we can only conclude from \ref{lem:convex:fenchel-young:iii} that \ref{lem:convex:fenchel-young:i} and \ref{lem:convex:fenchel-young:ii} hold for $F^{**}\neq F$ in place of $F$. Applying \cref{lem:convex:fenchel-young} to nonconvex functionals therefore inevitably introduces a \term{convexification} (by replacing the nonconvex $F$ with its convex envelope $F^{**}$).
\end{remark}
\begin{remark}\label{rem:convex:fy-preconjugate}
    Recall that $\partial F^*(x^*) \subset X^{**}$. Therefore, if $X$ is not reflexive, $x\in\partial F^*(x^*)$ in \ref{lem:convex:fenchel-young:iii} has to be understood via the canonical injection $J: X \hookrightarrow X^{**}$ as $Jx\in\partial F^*(x^*)$, i.e., as
    \begin{equation*}
        \dual{Jx,\tilde x^*-x^*}_{X^*} = \dual{\tilde x^* -x^*,x}_X \leq F^*(\tilde x^*) - F^*(x^*) \quad\text{for all $\tilde x^* \in X$.}
    \end{equation*}
    Using \ref{lem:convex:fenchel-young:iii} to conclude equality in \ref{lem:convex:fenchel-young:i} or, equivalently, the subdifferential inclusion \ref{lem:convex:fenchel-young:ii} therefore requires the additional condition that $x\in X\hookrightarrow X^{**}$. Conversely, if \ref{lem:convex:fenchel-young:i} or \ref{lem:convex:fenchel-young:ii} hold, \ref{lem:convex:fenchel-young:iii} also guarantees that the subderivative $x$ is an element of $\partial F^*(x^*)\cap X$, which is a stronger claim
    (see \cite{Gerd:2022} for a counterexample).

    Similar statements apply to (weakly-$*$ lower semicontinuous) $F:X^*\to\Rbar$ and $F_*:X\to\Rbar$.
\end{remark}

\section{Duality of optimization problems}

\Cref{lem:convex:fenchel-young} can be used to replace the subdifferential of a (complicated) norm with that of a (simpler) conjugate indicator functional (or vice versa).
For example, given a problem of the form
\begin{equation}\label{eq:convex:primal}
    \inf_{x\in X} F(x) + G(Kx)
\end{equation}
for $F:X\to\Rbar$ and $G:Y\to\Rbar$ proper, convex, and lower semicontinuous, and $K\in \linear(X;Y)$, we can use \cref{thm:convex:moreau} to replace $G$ with the definition of $G^{**}$ and obtain the \term[problem!saddle-point]{saddle-point problem}
\begin{equation}\label{eq:convex:saddle}
    \adjustlimits\inf_{x\in X}\sup_{y^*\in Y^*} F(x) + \dual{y^*,Kx}_Y - G^*(y^*).
\end{equation}
If(!) we were now able to exchange $\inf$ and $\sup$, we could write (with $\inf F = -\sup (-F)$)
\begin{equation*}
    \begin{aligned}
        \adjustlimits\inf_{x\in X}\sup_{y^*\in Y^*} F(x) + \dual{y^*,Kx}_Y - G^*(y^*)
        &=\adjustlimits\sup_{y^*\in Y^*} \inf_{x\in X} F(x) + \dual{y^*,Kx}_Y - G^*(y^*) \\
        &=\sup_{y^*\in Y^*} -\left\{\sup_{x\in X} -F(x) + \dual{-K^*y^*,x}_X\right\} - G^*(y^*).
    \end{aligned}
\end{equation*}
From the definition of $F^*$, we thus obtain the \term[problem!dual]{dual problem}
\begin{equation}\label{eq:convex:dual}
    \sup_{y^*\in Y^*} -G^*(y^*)-F^*(-K^*y^*).
\end{equation}
As a side effect, we have shifted the operator $K$ from $G$ to $F^*$ without having to invert it.

The following theorem uses in an elegant way the Fermat principle, the sum and chain rules, and the Fenchel--Young equality to derive sufficient conditions for the exchangeability.
\begin{theorem}[Fenchel--Rockafellar]\label{thm:convex:fenchel}\index{theorem!Fenchel--Rockafellar}
    Let $X$ and $Y$ be Banach spaces, $F:X\to\Rbar$ and $G:Y\to \Rbar$ be proper, convex, and lower semicontinuous, and $K\in \linear(X;Y)$. Assume furthermore that
    \begin{enumerate}
        \item \label{thm:convex:fenchel:i}
            the \term[problem!primal]{primal problem} \eqref{eq:convex:primal} admits a solution $\bar x \in X$;
        \item \label{thm:convex:fenchel:ii}
            there exists an $x_0\in \dom (G\circ K) \cap \dom F$ with $Kx_0\in \interior(\dom G)$ .
    \end{enumerate}
    Then the dual problem \eqref{eq:convex:dual} admits a solution $\bar {y}^*\in Y^*$ and
    \begin{equation}\label{eq:convex:fenchel:equal}
        \min_{x\in X} F(x) + G(Kx) = \max_{y^*\in Y^*} -G^*(y^*) - F^*(-K^*y^*).
    \end{equation}
    Furthermore, $\bar x$ and $\bar{y}^*$ are solutions to \eqref{eq:convex:primal} and \eqref{eq:convex:dual}, respectively, if and only if
    \begin{equation}\label{eq:convex:extremal}
        \left\{
            \begin{aligned}
                \bar{y}^* &\in \partial G(K\bar x),\\
                -K^*\bar{y}^* &\in \partial F(\bar x).
            \end{aligned}
        \right.
    \end{equation}
\end{theorem}
\begin{proof}
Let first $\bar x\in X$ be a solution to \eqref{eq:convex:primal}. By assumption \ref{thm:convex:fenchel:ii}, \cref{thm:subdiff:sum} (noting that since $K$ is a bounded linear operator, $Kx_0\in \interior(\dom G)$ implies $x_0\in \interior (\dom G\circ K)$), and \cref{thm:convex:chain} are applicable; \cref{thm:convex:fermat} thus implies that
    \begin{equation*}
        0\in  \partial (F + G\circ K)(\bar x) = K^*\partial G(K\bar x) + \partial F(\bar x)
    \end{equation*}
    and thus the existence of a $\bar{y}^*\in \partial G(K\bar x)$ with $-K^*\bar{y}^* \in \partial F(\bar x)$, i.e., satisfying \eqref{eq:convex:extremal}.

    Conversely, let \eqref{eq:convex:extremal} hold for $\bar x\in X$ and $\bar y^*\in Y^*$. Then again by \cref{thm:convex:fermat,thm:subdiff:sum,thm:convex:chain}, $\bar x$ is a solution to \eqref{eq:convex:primal}. Furthermore, \eqref{eq:convex:extremal} together with \cref{lem:convex:fenchel-young} imply equality in the Fenchel--Young inequalities for $F$ and $G$, i.e.,
    \begin{equation}\label{eq:convex:fenchel0}
        \left\{
            \begin{aligned}
                \dual{\bar y^*,K\bar x}_Y &=G(K\bar x) + G^*(\bar y^*),\\
                \dual{-K^*\bar y^*,\bar x}_X &=F(\bar x) + F^*(-K^*\bar y^*).
            \end{aligned}
        \right.
    \end{equation}
    Adding both equations and rearranging now yields
    \begin{equation}\label{eq:convex:fenchel1}
        F(\bar x) + G(K\bar x) =  -F^*(-K^*\bar y^*) - G^*(\bar y^*).
    \end{equation}
    It remains to show that $\bar y^*$ is a solution to \eqref{eq:convex:dual}.
    For this purpose, we introduce
    \begin{equation}
        \label{eq:convex:fenchel:L-fenchel}
        L:X\times Y^*\to\Rbar,\qquad L(x,y^*) =  F(x) + \dual{y^*,Kx}_Y - G^*(y^*).
    \end{equation}
    For all  $\tilde x\in X$ and $\tilde y^*\in Y^*$, we always have that
    \begin{equation}
        \label{eq:convex:fenchel:L-ineq-always}
        \sup_{y^*\in Y^*} L(\tilde x,y^*) \geq L(\tilde x,\tilde y^*) \geq \inf_{x\in X}  L(x,\tilde y^*),
    \end{equation}
    and hence (taking the infimum over all $\tilde x$ in the first and the supremum over all $\tilde y^*$ in the second inequality) that
    \begin{equation}
        \label{eq:convex:fenchel:infsup-supinf}
        \adjustlimits\inf_{x\in X} \sup_{y^*\in Y^*} L(x,y^*) \geq
        \adjustlimits\sup_{y^*\in Y^*} \inf_{x\in X} L(x,y^*).
    \end{equation}
    We thus obtain that
    \begin{equation}\label{eq:convex:fenchel2}
        \begin{aligned}[t]
            F(\bar x) + G(K\bar x)  &=
            \adjustlimits\inf_{x\in X}\sup_{y^*\in Y^*} F(x) + \dual{y^*,Kx}_Y - G^*(y^*)\\
            &\geq
            \adjustlimits\sup_{y^*\in Y^*} \inf_{x\in X} F(x) + \dual{y^*,Kx}_Y - G^*(y^*)\\
            & = \sup_{y^*\in Y^*} -G^*(y^*) - F^*(-K^*y^*)
        \end{aligned}
    \end{equation}
    (i.e., \term[duality!weak]{weak duality} holds merely under assumption \ref{thm:convex:fenchel:i}).
    Combining this with \eqref{eq:convex:fenchel1} yields that
    \begin{equation*}
        -G^*(\bar{y}^*) - F^*(-K^*\bar{y}^*) =  F(\bar x) + G(K\bar x)  \geq \sup_{y^*\in Y^*} -G^*(y^*) - F^*(-K^*y^*),
    \end{equation*}
    i.e., $\bar{y}^*$ is a solution to \eqref{eq:convex:dual}, which in particular shows the claimed existence of a solution.

    Since all solutions to \eqref{eq:convex:dual} have by definition the same (maximal) functional value, \eqref{eq:convex:fenchel1} also implies \eqref{eq:convex:fenchel:equal}.

    Finally, if $\bar x\in X$ and $\bar y^*\in Y^*$ are solutions to \eqref{eq:convex:primal} and \eqref{eq:convex:dual}, respectively, the just derived strong duality \eqref{eq:convex:fenchel:equal} conversely implies that \eqref{eq:convex:fenchel1} holds. Together with the productive zero, we obtain from this that
    \begin{equation*}
        0 = \left[ G(K\bar x) + G^*(\bar y^*) -\dual{\bar y^*,K\bar x}_X \right]+
        \left[F(\bar x) + F^*(-K^*\bar y^*) - \dual{-K^*\bar y^*,\bar x}_Y\right].
    \end{equation*}
    Since both brackets have to be nonnegative due to the Fenchel--Young inequality, they each have to be zero. We therefore deduce that \eqref{eq:convex:fenchel0} holds, and hence \cref{lem:convex:fenchel-young} implies~\eqref{eq:convex:extremal}.
\end{proof}

\begin{remark}\label{rem:fenchel:predual}
    If $X$ is the dual of a separable Banach space $X_*$, it is possible to derive a similar duality result with the (weakly-$*$ lower semicontinuous) preconjugate $F_*:X_*\to\Rbar$ in place of $F^*:X^*\to\Rbar$ under the additional assumption that $\range K^* \subset X_* \subsetneq X^*$ (using \cref{rem:convex:fy-preconjugate} in \eqref{eq:convex:fenchel0}). If $X_*$ is a \enquote{nicer} space than $X^*$ (e.g., for $X=\mathcal{M}(\Omega)$, the space of bounded Radon measures on a domain $\Omega\subset \R^d$ with $X_*=C_0(\Omega)$, the space of continuous functions with compact support), the \term[problem!predual]{predual problem}
    \begin{equation*}
        \sup_{y^*\in Y^*} -G^*(y^*) - F_*(-K^*y^*)
    \end{equation*}
    may be easier to treat than the dual problem \eqref{eq:convex:dual}. This is the basis of the ``preduality trick'' used in, e.g., \cite{Hintermuller:2004a,Clason:2010a}.
\end{remark}

\begin{remark}
    \label{rem:fenchel:attouch-brezis}
    The condition \ref{thm:convex:fenchel:ii} was only used to guarantee equality in the sum and chain rules \cref{thm:convex:chain,thm:subdiff:sum} applied to $F+G \circ K$. Since these rules hold under the weaker condition of \cref{rem:convex:attouch-brezis} (recall that the chain rule was proved by reduction to the sum rule), \cref{thm:convex:fenchel,cor:convex:fenchel-saddle} hold under this weaker condition as well.
\end{remark}

The relations \eqref{eq:convex:extremal} are referred to as \term[condition!Fenchel extremality]{Fenchel extremality conditions}; we can use \cref{lem:convex:fenchel-young} to generate further, equivalent, optimality conditions by inverting one or the other subdifferential inclusion.
We will later exploit this to derive implementable algorithms for solving optimization problems of the form \eqref{eq:convex:primal}. Furthermore, \cref{thm:convex:fenchel} characterizes the subderivative $\bar y^*$ produced by the sum and chain rules as solution to a convex minimization problem, which may be useful. For example, if either $F^*$ or $G^*$ is strongly convex, this subderivative will be unique, which has beneficial consequences for the stability and the convergence of algorithms for the computation of solutions to \eqref{eq:convex:extremal}.

For the analysis of such algorithms, it will sometimes be more convenient to apply the consequences of \cref{thm:convex:fenchel} in the form of the saddle-point problem \eqref{eq:convex:saddle}.
For a general mapping $L: X \times Y^* \to \Rbar$, we call
$(\tilde x, \tilde y^*)$ a \term[point!saddle]{saddle point} of $L$ if
\begin{equation}
    \label{eq:convex:fenchel:L-saddle}
    \sup_{y^*\in Y^*} L(\tilde x,y^*) \leq L(\tilde x,\tilde y^*) \leq \inf_{x\in X}  L(x,\tilde y^*).
\end{equation}
(Note that the converse inequality \eqref{eq:convex:fenchel:L-ineq-always} always holds.)
\begin{corollary}\label{cor:convex:fenchel-saddle}
    Assume that the conditions of \cref{thm:convex:fenchel} hold. Then any $(\bar x,\bar y^*)\in X\times Y^*$ satisfying \eqref{eq:convex:extremal} is a saddle point to
    \begin{equation*}
        L(x,y^*) \defeq F(x) + \dual{y^*,Kx}_Y - G^*(y^*).
    \end{equation*}
\end{corollary}
\begin{proof}
    Under the assumptions, which imply strong duality, the inequality in \eqref{eq:convex:fenchel2} holds as an equality, i.e.,
    \[
        \begin{aligned}
        \adjustlimits\inf_{x\in X}\sup_{y^*\in Y^*}  L(x, y^*)
        &
        =
        \adjustlimits\inf_{x\in X}\sup_{y^*\in Y^*} F(x) + \dual{y^*,Kx}_Y - G^*(y^*)
        \\
        &
        =
        \adjustlimits\sup_{y^*\in Y^*} \inf_{x\in X} F(x) + \dual{y^*,Kx}_Y - G^*(y^*)
        =
        \adjustlimits\sup_{y^*\in Y^*} \inf_{x\in X} L(x, y^*).
        \end{aligned}
    \]
    But by the definition of the Fenchel conjugate, we have that
    \[
        \adjustlimits\inf_{x\in X}\sup_{y^*\in Y^*}  L(x, y^*)
        = \adjustlimits\inf_{x\in X}F(x) + G(Kx)
        = F(\bar x) + G(K\bar x)
        = \sup_{y^*\in Y^*}  L(\bar x, y^*).
    \]
    We similarly have that
    \[
        \adjustlimits\sup_{y^*\in Y^*} \inf_{x\in X} L(x, y^*)
        = \adjustlimits\sup_{y^*\in Y^*}{-G^*(y^*) - F^*(-K^*y^*)}
        =  -G^*(\bar y^*) - F^*(-K^*\bar y^*)
        = \inf_{x \in X}  L(x, \bar y^*).
    \]
    Together, we obtain that
    \[
        \sup_{y^*\in Y^*}  L(\bar x, y^*) =  \inf_{x \in X}  L(x, \bar y^*)\leq L(\bar x,\bar y^*) \leq \sup_{y^*\in Y^*} L(\bar x,y^*),
    \]
    where the last two inequalities follow from \eqref{eq:convex:fenchel:L-ineq-always} (which always holds) for $\tilde x=\bar x$ and $\tilde y^* = \bar y^*$. Hence these inequalities hold with equality, which shows that
    \eqref{eq:convex:fenchel:L-saddle} holds for all $(x, y^*) \in X \times Y^*$, i.e., $(\bar x,\bar y^*)$ is a saddle point.
\end{proof}

With the notation $u=(x, y)$, let us define the  \term[gap!duality]{(Fenchel--Rockafellar) duality gap}
\begin{equation}\label{eq:convex:fenchel:duality-gap}
    \bar \gap(u) \defeq F(x)+G(Kx) + G^*(y^*) + F^*(-K^*y^*).
\end{equation}
By \cref{thm:convex:fenchel}, we have $\bar\gap \ge 0$ and $\bar\gap(\opt u)=0$ if and only if $\opt u$ is a saddle point.

On the other hand, for any saddle point $\opt u=(\opt x, \opt y^*)$ of a Lagrangian $L:X\times Y^*\to \Rbar$, we can also define the \term[gap!duality!Lagrangian]{Lagrangian duality gap}
\[
    \gap_L(u; \opt u)
    \defeq L(x, \opt y^*)-L(\opt x, y^*).
\]
For $L$ defined in \eqref{eq:convex:fenchel:L-fenchel}, we always have by
\cref{cor:convex:fenchel-saddle}, the definition \eqref{eq:convex:fenchel:L-saddle} of the saddle point, and the definition of the convex conjugate that
\begin{equation}
    \label{eq:convex:fenchel:lagrangian-duality-gap-bound}
    0 \le \gap_L(u; \opt u) \le \bar \gap(u).
\end{equation}
However, $\gap_L(u; \opt u)=0$ does not necessarily imply that $u$ is a saddle point.
(This is only the case if, e.g., $L$ is strictly convex in $x$ or strictly concave in $y$, i.e., if either $F$ or $G^*$ is strictly convex.)
Nevertheless, we will see in later chapters that for iterates of optimization algorithms, it is possible to show convergence of their Lagrangian duality gap, while this is in general more difficult for their Fenchel--Rockafellar duality gap.

\chapter{Monotone operators and proximal points}\label{chap:monotone}

Any minimizer $\bar x\in X$ of a convex functional $F:X\to \Rbar$ satisfies by \cref{thm:convex:fermat} the Fermat principle $0\in\partial F(\bar x)$.
To use this to characterize $\bar x$, and, later, to derive implementable algorithms for its iterative computation, we now study the mapping $x\mapsto \partial F(x)$ in more detail.

\section{Basic properties of set-valued mappings}
\label{sec:monotone:basic}

We start with some basic concepts.
For two normed vector spaces $X$ and $Y$ we consider a \term[mapping!set-valued]{set-valued mapping} $A:X\to\mathcal{P}(Y)$, also denoted by $A:X \setto Y$, and define
\begin{itemize}
    \item its \term[domain!of a set-valued mapping]{domain of definition} $\dom A = \setof{x\in X}{A(x)\neq \emptyset}$;
    \item its \term[range!of a set-valued mapping]{range} $\range A = \bigcup_{x\in X} A(x)$;
    \item its \term[graph!of a set-valued mapping]{graph} $\graph A = \setof{(x,y)\in X\times Y}{y\in A(x)}$;
    \item its \term[inverse!of a set-valued mapping]{inverse} $A^{-1}:Y\setto X$ via $A^{-1}(y) = \setof{x\in X}{y\in A(x)}$ for all $y\in Y$.
\end{itemize}
(Note that $A^{-1}(y) = \emptyset$ is allowed by definition; hence for set-valued mappings, the inverse always exists.) Similarly, we will say that $A:X\setto Y$ is \term[mapping!surjective]{surjective} if $\range A = Y$.

For $A,B:X\setto Y$, $C:Y\setto Z$, and $\lambda\in\R$ we further define
\begin{itemize}
    \item $\lambda A:X\setto Y$ via $(\lambda A)(x) = \setof{\lambda y}{y\in A(x)}$;
    \item $A+B:X\setto Y$ via $(A+B)(x) = \setof{y+z}{y\in A(x), z\in B(x)}$;
    \item $C\circ A:X\setto Z$ via $(C\circ A)(x) = \setof{z}{\text{there is }y\in A(x) \text{ with } z\in C(y)}$.
\end{itemize}

\bigskip

Of particular importance not only in the following but also in \cref{part:setvalued} is the continuity of set-valued mappings. We first introduce  notions of convergence of sets.
So let $\{X_n\}_{n\in\N}$ be a sequence of subsets of $X$. We define
\begin{enumerate}
    \item the \term[limit!outer]{outer limit} as the set
        \begin{equation*}
            \limsup_{n \to \infty} X_n \defeq \setof{x \in X}{\text{there exists } \{n_k\}_{k\in\N}\text{ with } x_{n_k}\in X_{n_k} \text{ and } \lim_{k\to\infty} x_{n_k} = x},
        \end{equation*}
    \item the \term[limit!inner]{inner limit} as the set
        \begin{equation*}
            \liminf_{n \to \infty} X_n \defeq \setof{x\in X}{\text{there exist } x_{n}\in X_{n} \text{ with } \lim_{n\to\infty} x_{n} = x}.
        \end{equation*}
\end{enumerate}
Correspondingly, we define the \term[limit!outer!weak]{weak outer limit} and the \term[limit!inner!weak]{weak inner limit}, denoted by $\weaklimsup_{n\to\infty} X_n$ and $\weakliminf_{n\to\infty}X_n$, respectively, using weakly converging (sub)\-se\-quences. Similarly, for a dual space $X^*$, we define the \term[limit!outer!weak-$*$]{weak-$*$ outer limit} $\weakstarlimsup_{n\to\infty}X^*_n$ and the \term[limit!inner!weak-$*$]{weak-$*$ inner limit} $\weakstarliminf_{n\to\infty}X^*_n$.

The outer limit consists of all points approximable through \emph{some} subsequence of the sets $X_n$, while the inner limit has to be approximable through \emph{every} subsequence.
The vast difference between inner and outer limits is illustrated by the following extreme example.

\begin{example}\label{ex:inner-outer-difference}
    Let $X=\R$ and $\{X_n\}_{n\in\N}$, $X_n\subset [0,1]$, be given as
    \begin{equation*}
        X_n \defeq
        \begin{cases}
            [0,\tfrac13) & \text{if $n=3k-2$ for some $k\in\N$},\\
            [\tfrac13,\tfrac23) & \text{if $n=3k-1$ for some $k\in\N$},\\
            [\tfrac23, 1] & \text{if $n=3k$ for some $k\in\N$},
        \end{cases}
    \end{equation*}
    see \cref{fig:inner-outer-difference}. Then,
    \begin{align*}
        \limsup_{n\to\infty} X_n &= [0,1],
        \intertext{since for any $x\in [0,1]$, we can find a subsequence of $\{X_n\}_{n\in\N}$ (by selecting subsequences with, e.g., $n=3k-2$ for $k\in \N$ if $x<\frac13$) that contain $x$. On the other hand,}
        \liminf_{n\to\infty} X_n &= \emptyset,
    \end{align*}
    since for any $x \in [0, 1]$, there will be a subsequence of $X_n$ (again, selecting only subsequences with, e.g., $n=3k$ for $k\in \N$ if $x<\frac13$) that will not contain points arbitrarily close to $x$.
\end{example}
\begin{figure}
    \centering
    \begin{asy}
        int m=700;
        real log2 (real x){ return log(x)/log(2); }
        real getx(int i){ return 2.5*log2(log2((i+10)/5)); }

        for(int i=0; i<m; ++i){
            real x=getx(i);
            if(i%3==0){
                draw((x, 0)--(x, 0.332), primalline);
            }else if(i%3==1){
                draw((x, 0.333)--(x, 0.665), primalline);
            }else{
                draw((x, 0.665)--(x, 1), primalline);
            }
        }
        real x=getx(m);
        draw((x,0)--(x,1), primalline);

        label("$X_1$", (getx(0),0),S);
        label("$X_2$", (getx(1),0),S);
        label("$X_3$", (getx(2),1),N);
        label("$X_4$", (getx(3),0),S);
        label("$X_5$", (getx(4),1),N);
        label("$\phantom{X_4}\ldots$", (getx(4),0),S);
        label("$\phantom{X_3}\ldots$", (getx(5),1),N);

        real x=1.1*getx(m);

        draw((x, 0)--(x, 1), Bars);
        label("$0$", (x, 0), S);
        label("$1$", (x, 1), N);
    \end{asy}
    \caption{Illustration of \cref{ex:inner-outer-difference} with $\limsup_{n \to \infty} X_n=[0, 1]$ while $\liminf_{n \to \infty} X_n=\emptyset$.}
    \label{fig:inner-outer-difference}
\end{figure}

\begin{lemma}
    \label{lemma:limsup-setlimit}
    Let $\{X_n\}_{n\in\N}$, $X_n\subset X$. Then $\limsup_{n\to\infty} X_n$ and $\liminf_{n\to\infty} X_n$
    are (possibly empty) closed sets.
\end{lemma}
\begin{proof}
    Let $X_\infty \defeq \limsup_{n \to \infty} X_n$.
    If $X_\infty$ is empty, there is nothing to prove.
    So suppose that $\{x_k\}_{k \in \N} \subset X_\infty$ converges to some $\hat x \in X$.
    Since each $x_k$ is an element of an outer limit, for each $k\in \N$ there exist infinite subsets $N_k \subset \N$ and subsequences $\{x_{k,n}\}_{n\in N_k}$ with $x_{k,n} \in X_n$ for all $n\in\N$ and $\lim_{N_k \ni n \to \infty} x_{k,n}=x_k$. We can therefore find for each $k \in \N$ an index $n_k \in N_k$ such that $\norm{x_k-x_{k,n_k}}_X \le 1/k$. Together, this implies that
    \[
        \norm{\hat x - x_{k,n_k}}_X \le \norm{\hat x - x_k}_X + \norm{x_k-x_{k,n_k}}_X \to 0
    \]
    as $k\to\infty$ and hence that $X_{n_k} \ni x_{k,n_k} \to \hat x$. This shows that $\hat x \in X_\infty$.

    Let then $X_\infty \defeq \liminf_{n \to \infty} X_n$.
    If $X_\infty$ is empty, there is nothing to prove.
    So suppose that $\{x_k\}_{k \in \N} \subset X_\infty$ converges to some $\hat x \in X$. We want to proceed similarly to the proof of the first claim, but since the inner limit requires taking a sequence where an element is taken from each $X_n$, we cannot simply \enquote{skip over} sets when constructing the diagonal sequence but must instead \enquote{delay} the convergence $x_k\to \hat x$ suitably.
    First, since each $x_k$ is an element of an inner limit, for each $k\in\N$ there exists a (full) sequence $\{x_{k,n}\}_{n\in\N}$ with $x_{k,n} \in X_n$ and $\lim_{n \to \infty} x_{k,n}=x_k$. This allows us to choose a sequence $\{k_n\}_{n\in\N}$ via
    \begin{equation*}
        k_1 = 1,\qquad
        k_n =
        \begin{cases}
            k_{n-1}+1 &\text{if } \norm{x_{k_{n-1}+1}-x_{k_{n-1}+1,n}}_X \leq \frac1{k_{n-1}+1},\\
            k_{n-1}   &\text{otherwise.}
        \end{cases}
    \end{equation*}
    Then $k_n\to\infty$ as $n\to \infty$, since otherwise $k_n=k^*$ would eventually remain constant, which by construction implies that $\norm{x_{k^*}-x_{k^*,n}}_X > \frac1{k^*}$ for all $n$ large enough. But this contradicts the convergence $x_{k,n}\to x_{k}$ for all $k\in\N$.
    Hence $k_n\to\infty$ and therefore $x_{k_n}\to \hat x$ as $n\to\infty$.
    Together, this yields
    \begin{equation*}
        \norm{\hat x - x_{k_n,n}}_X \le \norm{\hat x - x_{k_n}}_X + \norm{x_{k_n}-x_{k_n,n}}_X \to 0
    \end{equation*}
    as $n\to \infty$ and hence that $X_n \ni x_{k_n,n}\to \hat x$. This shows that $\hat x\in X_\infty$.
\end{proof}

With these definitions, we can define limits and continuity of set-valued mappings. Specifically, for $A: X \setto Y$, and a subset $C \subset X$, we define the inner and outer limits (relative to $C$, if $C \ne X$) as
\begin{align*}
    \limsup_{C \ni \tilde x \to x} A(\tilde x) &\defeq \Union_{C \ni x_n \to x} \limsup_{n \to \infty} A(x_n),
    \shortintertext{and}
    \liminf_{C \ni \tilde x \to x} A(\tilde x) &\defeq \Isect_{C \ni x_n \to x} \liminf_{n \to \infty} A(x_n).
\end{align*}
If $C=X$, we drop $C$ from the notations. Analogously, we define weak-to-strong, strong-to-weak, and weak-to-weak limits by replacing $x_n\to x$ by $x_n\weakto x$ and/or the outer/inner limit by the weak outer/inner limit.

\begin{corollary}
    \label{lemma:limsup}
    Let $A:X\setto Y$ and $x\in X$.
    Then  $\limsup_{\tilde x \to x} A(\tilde x)$ and $\liminf_{\tilde x \to x} A(\tilde x)$ are (possibly empty) closed sets.
\end{corollary}

\begin{proof}
    The proof of the closedness of the outer limit is analogous to \cref{lemma:limsup-setlimit},
    while the proof of the closedness of the inner limit is a consequence of \cref{lemma:limsup-setlimit} and of the fact that the intersections of closed sets are closed.
\end{proof}

Let then $A: X \setto Y$ be a set-valued mapping. We say that
\begin{enumerate}
    \item $A$ is \term[mapping!semicontinuous!outer]{outer semicontinuous at $x\in X$} if $\limsup_{C \ni \tilde x \to x} A(\tilde x) \subset A(x)$ with $C=X$.
    \item $A$ is \term[mapping!semicontinuous!inner]{inner semicontinuous at $x\in X$} if $\liminf_{C \ni \tilde x \to x} A(\tilde x) \supset A(x)$ with $C=X$.
    \item $A$ is \emph{outer/inner semicontinuous} if it is outer/inner semicontinuous  at all $x \in X$.
    \item $A$ is\term[mapping!continuous]{continuous} (at $x$) if it is both outer and inner semicontinuous (at $x$).
\end{enumerate}
We say that these properties are \enquote{relative to $C$} when we restrict $\tilde x \in C$ for some $C \subset X$.
These concepts are illustrated in Figure \ref{fig:semicontinuity}.
\begin{figure}
    \centering
    \begin{asy}
        real x1=-0.5;
        real x2=0.6;
        real fup(real x){ return 1+0.1*x^3; }
        real fdown0(real x) { return x^2*0.2+0.3; };
        real fdown1(real x) { return -0.1-0.3*x^3; };
        real fdown2(real x) { return 0; };

        path pup=graph(fup, -1, 1);
        path pdown0=graph(fdown0, -1, x1);
        path pdown1=graph(fdown1, x1, x2);
        path pdown2=graph(fdown2, x2, 1);
        path pupext=graph(fup, -1.1, 1.15);
        path pdownext0=graph(fdown0, -1.15, x1);
        path pdownext2=graph(fdown2, x2, 1.15);

        fill(pupext--reverse(pdownext0--pdown1--pdownext2)--cycle, lightfill);
        draw(pup);
        draw(pdown0);
        draw(pdown1--pdown2);

        dot(point(pdown0, length(pdown0)));
        dot(point(pdown1, length(pdown1)));

        pair up0=pt(fup, 0);
        pair down0=pt(fdown1, 0);

        label("$A$", (up0+down0)/2, S);

        label("$x_1$", (x1, -.3), S);
        label("$x_2$", (x2, -.3), S);
    \end{asy}
    \caption{Illustration of outer and inner semicontinuity. The black line indicates the bounds on the boundary of $\graph F$ that belong to the graph. The set-valued mapping $A$ is not outer semicontinuous at $x_1$, because $A(x_1)$ does not include all limits from the right. It is outer semicontinuous at the \enquote{discontinuous} point $x_2$, as $A(x_2)$ includes all limits from both sides. The mapping $A$ is not inner semicontinuous at $x_2$, because at this point, $A(x)$ cannot be approximated from both sides. It is inner semicontinuous at every other point $x$, including $x_1$, as at these points $A(x)$ can be approximated from both sides.}
    \label{fig:semicontinuity}
\end{figure}

Just like lower semicontinuity of functionals, the outer semicontinuity of set-valued mappings can be interpreted as a closedness property and will be crucial. The following lemma is stated for strong-to-strong outer semicontinuity, but corresponding statements hold (with identical proof) for weak-to-strong, strong-to-weak, and weak-to-weak outer semicontinuity as well.
\begin{lemma}
    \label{ex:monotone:outersemi-closed}
    A set-valued mapping $A: X \setto Y$ is outer semicontinuous if and only if $\graph A \subset X\times Y$ is sequentially closed, i.e., $x_n\to x$ and $A(x_n)\ni y_n\to y$ imply that $y\in A(x)$.
\end{lemma}
\begin{proof}
    Let $x_n \to x$ and $y_n \in A(x_n)$, and suppose also $y_n \to y$.
    Then if $\graph A$ is closed, $(x, y) \in \graph A$ and hence $y \in A(x)$. Since this holds for arbitrary sequences $\{x_n\}_{n\in\N}$, $A$ is outer semicontinuous.

    If, on the other hand, $A$ is outer semicontinuous, and $(x_n, y_n) \in \graph A$ converge to $(x, y)\in X\times Y$‚ then $y \in A(x)$ and hence $(x, y) \in \graph A$.
    Since this holds for arbitrary sequences $\{(x_n, y_n)\}_{n\in\N}$, $\graph A$ is closed.
\end{proof}

\section{Monotone operators}
\label{sec:monotone:monotone}

For the codomain $Y=X^*$ (as in the case of $x\mapsto\partial F(x)$), additional properties become important.
A set-valued mapping $A:X\setto X^*$ is called \term[mapping!monotone]{monotone} if
\begin{equation}\label{eq:monoton:def}
    \dual{x^*_1-x^*_2,x_1-x_2}_X \geq 0 \quad\text{for all}\quad (x_1,x_1^*),(x_2,x_2^*) \in \graph A.
\end{equation}
Straight from the definition, we obtain the monotonicity of the following mappings.
\begin{example}\label{ex:monotone}
    \begin{enumerate}
        \item If $A:X\setto X^*$ is monotone and $\lambda \geq 0$, then $\lambda A$ is monotone as well.
        \item If $A,B:X\setto X^*$ are monotone, then $A+B$ is monotone as well.
        \item\label{ex:monotone:subdiff} If $F:X\to\Rbar$ is proper, then $\partial F:X\setto X^*$, $x\mapsto \partial F(x)$, is monotone since for any $x_1,x_2\in X$ with $x^*_1 \in \partial F(x_1)$ and $x^*_2\in\partial F(x_2)$, we have by definition that
            \begin{align*}
                &\dual{x_1^*,\tilde x - x_1}_X\leq F(\tilde x)-F(x_1)\qquad\text{for all}\quad \tilde x\in X,\\
                &\dual{x_2^*,\tilde x - x_2}_X\leq F(\tilde x)-F(x_2)\qquad\text{for all }\quad \tilde x\in X.
            \end{align*}
            Adding the first inequality for $\tilde x=x_2$ and the second for $\tilde x=x_1$ and rearranging the result yields \eqref{eq:monoton:def}.
    \end{enumerate}
\end{example}
(\Cref{ex:monotone}\,\ref{ex:monotone:subdiff} generalizes the well-known fact that if $f:\R\to\R$ is convex and differentiable, its derivative $f'$ is monotonically increasing.)

In fact, we will need the following, stronger, property, which guarantees that $A$ is outer semicontinuous: A monotone operator $A: X\setto X^*$ is called \term[mapping!monotone!maximally]{maximally monotone} if there does not exist another monotone operator $\tilde A: X \setto X^*$ such that $\graph A \subsetneq \graph \tilde A$. In other words, $A$ is maximal monotone if for any $x\in X$ and $x^*\in X^*$ the condition
\begin{equation}\label{eq:monoton:max_char}
    \dual{x^*-\tilde x^*,x-\tilde x}_X \geq 0 \qquad \text{for all }(\tilde x , \tilde x^*) \in \graph A
\end{equation}
implies that $x^*\in A(x)$. (In other words, \eqref{eq:monoton:max_char} holds if \emph{and only if} $(x,x^*)\in\graph A$.)
For fixed $x\in X$ and $x^*\in X^*$, the condition claims that if $A$ is monotone, then so is the extension
\begin{equation*}
    \tilde A:X\setto X^*,\qquad
    \tilde x \mapsto
    \begin{cases}
        A(x) \cup \{x^*\}&  \text{if } \tilde x =  x,\\
        A(\tilde x) & \text{if }\tilde x\neq  x.
    \end{cases}
\end{equation*}
For $A$ to be maximally monotone means that this is not a true extension, i.e., $\tilde A=A$.
\begin{example}
The operator
\begin{equation*}
    A:\R\setto\R,\qquad t \mapsto
    \begin{cases} \{1\} & \text{if }t> 0,\\
        \{0\} & \text{if } t= 0,\\
    \{-1\} &\text{if } t< 0,\end{cases}
\end{equation*}
is monotone but not maximally monotone, since $A$ is a proper subset of the monotone operator defined by $\tilde A(t) = \sign(t)=\partial (|\cdot|)(t)$ from \cref{ex:convex:subdiff_abs}.
\end{example}

Several useful properties follow directly from the definition.
\begin{lemma}\label{lem:monoton:scalar_maxmon}
    If $A:X\setto X^*$ is maximally monotone, then so is $\lambda A$ for all $\lambda >0$.
\end{lemma}
\begin{proof}
    Let $x\in X$ and $x^*\in X^*$, and assume that
    \begin{equation*}
        0 \leq \dual{x^*-\tilde x^*,x-\tilde x}_X = \lambda\dual{\lambda^{-1}x^*-\lambda^{-1}\tilde x^*,x-\tilde x}_X \quad \text{for all }(\tilde x,  \tilde x^*) \in \graph \lambda A.
    \end{equation*}
    Since $\tilde x^*\in \lambda A(\tilde x)$ if and only if $\lambda^{-1}\tilde x^* \in A(\tilde x)$ and $A$ is maximally monotone, this implies that $\lambda^{-1}x^* \in A(x)$, i.e., $x^* \in (\lambda A)(x)$. Hence, $\lambda A$ is maximally monotone.
\end{proof}
\begin{lemma}%
    \label{lemma:monotone:convex}
    If $A: X \setto X^*$ is maximally monotone, then $A(x)$ is convex for all $x \in X$.
\end{lemma}
\begin{proof}
    Assume that $A(x)$ is not convex, i.e., $x_\lambda^* \defeq \lambda x^* + (1-\lambda) \alt x^* \notin A(x)$ for some $x^*, \alt x^* \in A(x)$ and $\lambda \in (0, 1)$. We then show that $A$ is not maximal. To see this, we define $\alt A$ via
    \begin{equation*}
        \alt A(y)\defeq
        \begin{cases}
            A(y) & y \ne x, \\
            A(x) \union \{x_\lambda^*\}, & y = x,
        \end{cases}
    \end{equation*}
    and show that $\alt A$ is monotone. By the definition of $\alt A$, it suffices to show for all $y \in X$ and $y^* \in A(y)$ that
    \begin{equation*}
        \dualprod{x_\lambda^*-y^*}{x-y}_{X} \ge 0.
    \end{equation*}
    But this follows directly from the definition of $x_\lambda^*$ and the monotonicity of $A$.
\end{proof}
\begin{lemma}%
    \label{lemma:monotone:inverse}
    Let $X$ be a reflexive Banach space. If $A:X\setto X^*$ is maximally monotone, then so is $A^{-1}:X^*\setto X^{**}\simeq X$.
\end{lemma}
\begin{proof}
    First, recall that the inverse $A^{-1}:X^*\setto X$ always exists as a set-valued mapping and can be identified with a set-valued mapping from $X^*$ to $X^{**}$ with the aid of the canonical injection $J:X\to X^{**}$ from \eqref{eq:functan:canonical-injection}, i.e.,
    \begin{equation*}
        A^{-1}(x^*) \defeq \setof{Jx\in X^{**}}{x^*\in A(x)} \qquad\text{for all }x^*\in X^*
    \end{equation*}
    From this and the definition \eqref{eq:functan:canonical-injection}, it is clear that $A^{-1}$ is monotone if and only if $A$ is.

    Let now $x^*\in X^*$ and $x^{**}\in X^{**}$ be given, and assume that
    \begin{equation}\label{eq:monotone:inverse1}
        \dual{x^{**} - \tilde x^{**},x^*-\tilde x^*}_{X^*} \ge 0 \qquad\text{for all }(\tilde x^*,\tilde x^{**})\in \graph A^{-1}.
    \end{equation}
    Since $X$ is reflexive, $J$ is surjective such that there exists an $x\in X$ with $x^{**}=Jx$. Similarly, we can write $\tilde x^{**} = J\tilde x$ for some $\tilde x\in X$ with $\tilde x^*\in A(\tilde x)$. By definition of the duality pairing, \eqref{eq:monotone:inverse1} is thus equivalent to
    \begin{equation*}
        \dual{x^{*} - \tilde x^{*},x-\tilde x}_{X} \geq 0
    \end{equation*}
    for all $\tilde x\in X$ and $\tilde x^* \in A(\tilde x)$. But since $A$ is maximally monotone, this implies that $x^*\in A(x)$ and hence $x^{**}=Jx\in A^{-1}(x)$.
\end{proof}

We now come to the outer semicontinuity.
\begin{lemma}\label{cor:monoton:closed}
    Let $A:X\setto X^*$ be maximally monotone. Then $A$ is both weak-to-strong and strong-to-weak-$*$ outer semicontinuous. In particular, $A(x)$ is closed for all $x\in X$.
\end{lemma}
\begin{proof}
    Let $x\in X$ and $x^*\in X^*$ and consider sequences $\{x_n\}_{n\in\N}\subset X$ with $x_n\weakto x$ and $\{x_n^*\}_{n\in\N}\subset X^*$ with $x_n^*\in A(x_n)$ and $x_n^*\to x^*$ (or $x_n\to x$ and $x^*_n\weaktostar x^*$).
    For arbitrary $\tilde x\in X$ and $\tilde x^*\in A(\tilde x)$, the monotonicity of $A$ implies that
    \begin{equation*}
        0\leq \dual{x^*_n - \tilde x^*,x_n - \tilde x}_X \to \dual{x^* - \tilde x^*,x - \tilde x}_X
    \end{equation*}
    since the duality pairing of strongly and weakly (or weakly-$*$ and strongly) converging sequences is convergent. Since $A$ is maximally monotone, we obtain that $x^*\in A(x)$ and hence $A$ is weak-to-strong (or strong-to-weak-$*$) outer semicontinuous by \cref{ex:monotone:outersemi-closed}.
\end{proof}

Since the pairing of weakly and weakly-$*$ convergent sequences does not converge in general, weak-to-weak-$*$ outer semicontinuity requires additional assumptions on the two sequences. Although we will not need to make use of it, the following notion can prove useful in other contexts. We call a set-valued mapping $A:X\setto X^*$ \term[mapping!semicontinuous!BCP outer]{BCP outer semicontinuous} (for \emph{Brezis--Crandall--Pazy}), if for any sequences $\{x_n\}_{n\in \N}\subset X$ and $\{x_n^*\}_{n\in\N}\subset X^*$ with
\begin{enumerate}
    \item \label{lemma:BCP-outer-semicontinuous:i}
        $x_n\weakto x$ and $A(x_n)\ni x_n^*\weaktostar x^*$,
    \item \label{lemma:BCP-outer-semicontinuous:ii}
        $\displaystyle \limsup_{n\to\infty}\,\dual{x_n^*-x^*,x_n-x}_X\leq 0$,
\end{enumerate}
we have $x^*\in A(x)$. The following result from \cite[Lemma 1.2]{BrezisCrandallPazy} (hence the name) shows that maximally monotone operators are BCP outer semicontinuous.
\begin{lemma}
    \label{lemma:max-mon-BCP-outer-semicontinuous}
    Let $X$ be a Banach space and let $A: X \setto X^*$ be maximally monotone.
    Then $A$ is BCP outer semicontinuous.
\end{lemma}
\begin{proof}
    First, the monotonicity of $A$ and assumption \ref{lemma:BCP-outer-semicontinuous:ii} imply that
    \begin{equation}
        \label{eq:monotone:bcp-max-mono-0}
        0\leq \liminf_{n\to\infty}\,\dual{x_n^*-x^*,x_n-x}_X \leq
        \limsup_{n\to\infty}\,\dual{x_n^*-x^*,x_n-x}_X \leq 0.
    \end{equation}
    Furthermore, from assumption \ref{lemma:BCP-outer-semicontinuous:i} and the fact that $X$ is a Banach space, it follows that $\{x_n\}_{n\in\N}$ and $\{x_n^*\}_{n\in\N}$ and hence also $\{\dual{x_n^*,x_n}_X\}_{n\in\N}$ are bounded. Thus there exists a subsequence such that $\dual{x_{n_k}^*,x_{n_k}}_X\to L$ for some $L\in \R$. Passing to the limit, and using \eqref{eq:monotone:bcp-max-mono-0}, we obtain that
    \begin{equation*}
        \begin{aligned}
            0 &= \lim_{k\to\infty} \dual{x_{n_k}^*-x^*,x_{n_k}-x}_X\\
              &= \lim_{k\to\infty} \dual{x_{n_k}^*,x_{n_k}}_X - \lim_{k\to\infty} \dual{x_{n_k}^*,x}_X -
            \lim_{k\to\infty} \dual{x^*,x_{n_k}}_X + \dual{x^*,x}_X \\
            &= L - \dual{x^*,x}_X.
        \end{aligned}
    \end{equation*}
    Since the limit does not depend on the subsequence, we have that $\dual{x_n^*,x_n}_X \to \dual{x^*,x}_X$.

    Let now $\tilde x\in X$ and $\tilde x^*\in A(\tilde x)$ be arbitrary. Using again the monotonicity of $A$ and assumption \ref{lemma:BCP-outer-semicontinuous:i} together with the first claim yields
    \begin{equation*}
        \begin{aligned}
        0&\leq \liminf_{n\to\infty}\,\dual{x_n^*-\tilde x^*,x_n-\tilde x}_X\\
             &\leq  \lim_{n\to\infty}\dual{x_n^*,x_n}_X - \lim_{n\to\infty} \dual{x_{n}^*,\tilde x}_X -
            \lim_{n\to\infty} \dual{\tilde x^*,x_{n}}_X + \dual{\tilde x^*,\tilde x}_X\\
            &= \dual{x^*-\tilde x^*,x-\tilde x}_X
        \end{aligned}
    \end{equation*}
    and hence that $x^*\in A(x)$ by the maximal monotonicity of $A$.
\end{proof}

The usefulness of BCP outer semicontinuity arises from the fact that it also implies weak-to-strong outer semicontinuity under slightly weaker conditions on $A$.
\begin{lemma}\label{lemma:bpr-strong-dual-corollary}
    Suppose $A: X \setto X^*$ is monotone (but not necessarily maximally monotone) and BCP outer semicontinuous. Then $A$ is also weak-to-strong outer semicontinuous.
\end{lemma}
\begin{proof}
    Let $x_n\weakto x$ and $x_n^*\to x^*$ with $x_n^*\in A(x_n)$ for all $n\in \N$.
    This implies that $x_n^*\weaktostar x^*$ as well and that $\{x_n\}_{n\in \N}$ is bounded.
    We thus have for some $C>0$ that
    \begin{equation*}
        \limsup_{n \to \infty}\,
        \dualprod{x_n^*-x^*}{x_n-x}_X
        \le
        C \limsup_{n \to \infty} \norm{x_n^*-x^*}_{X^*}
        =
        0.
    \end{equation*}
    Hence, condition \ref{lemma:BCP-outer-semicontinuous:ii} is satisfied, and the BCP outer semicontinuity yields $x^*\in A(x)$.
\end{proof}

We now show that convex subdifferentials are maximally monotone. Although this result (known as \term[theorem!Rockafellar]{Rockafellar's theorem}, see \cite{Rockafellar:1970}) holds in arbitrary Banach spaces, the proof (adapted here from \cite{Simons:2009}) greatly simplifies in reflexive Banach spaces.
\begin{theorem}\label{thm:monoton:subdiff}
    Let $X$ be a reflexive Banach space and $F:X\to\Rbar$ be proper, convex, and lower semicontinuous. Then $\partial F:X\setto X^*$ is maximally monotone.
\end{theorem}
\begin{proof}
    First, we already know from \cref{ex:monotone}\,\ref{ex:monotone:subdiff} that $\partial F$ is monotone.
    Let now $x\in X$ and $x^*\in X^*$ be given such that
    \begin{equation}\label{eq:monoton:subdiff1}
        \dual{x^*-\tilde x^*,x-\tilde x}_X \geq 0 \qquad \text{for all }\tilde x\in X, \tilde x^*\in\partial F(\tilde x).
    \end{equation}
    We consider
    \begin{equation*}
        J:X\to\Rbar,\qquad z\mapsto F(z+x) - \dual{x^*,z}_X + \frac12\norm{z}_X^2,
    \end{equation*}
    which is proper, convex and lower semicontinuous by the assumptions on $F$. Furthermore, $J$ is coercive by \cref{lem:convex:supercoercive}.
    \Cref{thm:convex:existence} thus yields a $\bar z\in X$ with $J(\bar z) = \min_{z\in X} J(z)$. By \cref{thm:convex:fermat,thm:subdiff:sum,thm:convex:gateaux,lem:convex:subdiff_calc}\,\ref{lem:convex:subdiff_calc:ii} then
    \begin{equation}
        \label{eq:monoton:subdiff:z-oc}
        0\in \partial F(\bar z +x) - \{x^*\} + \partial j(\bar z),
    \end{equation}
    where we have introduced $j(z)\defeq\frac12\norm{z}_X^2$. In other words, there exists a $z^*\in\partial j(\bar z)$ such that $x^*-z^*\in\partial F(\bar z + x)$.
    Combining \cref{lem:convex:power-conjugate} for $p=q=2$ and \cref{lem:convex:fenchel-young}, we furthermore have that $z^*\in \partial j(\bar z)$ if and only if
    \begin{equation}\label{eq:monoton:fenchel1}
        \dual{z^*,\bar z}_X = \frac12\norm{\bar z}_X^2 + \frac12\norm{z^*}_{X^*}^2.
    \end{equation}
    Applying now \eqref{eq:monoton:subdiff1} for $\tilde x = \bar z+x$ and $\tilde x^* = x^*-z^* \in \subdiff F(\tilde x)$, we obtain using \eqref{eq:monoton:fenchel1} that
    \begin{equation*}
        0\leq \dual{x^*- x^* +z^*,x-\bar z -x}_X = -\dual{z^*,\bar z}_X = -\frac12 \norm{z^*}_{X^*}^2-\frac12\norm{\bar z}_X^2,
    \end{equation*}
    implying that both $\bar z=0$ and $z^*=0$.
    Hence by \eqref{eq:monoton:subdiff:z-oc} we conclude that $x^*\in \partial F(x)$, which shows that $\partial F$ is maximally monotone.
\end{proof}

The argument in the preceding proof can be modified to give a characterization of maximal monotonicity for general monotone operators; this is known as \term[theorem!Minty]{Minty's theorem} and is a central result in the theory of monotone operators. We again make use of the \term[mapping!duality]{duality mapping} $\partial j:X\setto X^*$ for $j(x)=\frac12\norm{x}_X^2$.
\begin{theorem}[Minty]\label{thm:monoton:max_surj}
    Let $X$ be a reflexive Banach space and $A:X\setto X^*$ be monotone with $\graph A \neq \emptyset$. If $A$ is maximally monotone, then $\partial j+A$ is surjective.
\end{theorem}
\begin{proof}
    We proceed similarly as in the proof of \cref{thm:monoton:subdiff} by constructing a functional $F_A$ which plays the same role for $A$ as $F$ does for $\partial F$.
    Specifically, we define for a maximally monotone operator $A:X\setto X^*$ with $\graph A \neq \emptyset$ the \term[functional!Fitzpatrick]{Fitzpatrick functional}
    \begin{equation}\label{eq:monoton:fitzpatrick1}
        F_A:X\times X^*\to (-\infty,\infty],\qquad (x,x^*)\mapsto \sup_{(z,z^*)\in \graph A} \left(\dualprod{x^*}{z}_X + \dualprod{z^*}{x}_X - \dualprod{z^*}{z}_X\right),
    \end{equation}
    which can be written equivalently as
    \begin{equation}\label{eq:monoton:fitzpatrick2}
        F_A(x,x^*) = \dualprod{x^*}{x}_X - \inf_{(z,z^*)\in\graph A} \dualprod{x^*-z^*}{x-z}_X.
    \end{equation}
    Each characterization implies useful properties.
    \begin{enumerate}
        \item \label{thm:monoton:max_surj:i}
            By maximal monotonicity of $A$, we have by definition that $\dualprod{x^*-z^*}{x-z}_X\geq 0$ for all $(z,z^*)\in\graph A$ if and only if $(x,x^*)\in\graph A$. In particular, for all $(x,x^*)\notin\graph A$ there exists $(z,z^*)\in\graph A$ with $\dualprod{x^*-z^*}{x-z}_X< 0$, and therefore $\inf_{(z,z^*)\in\graph A}\dualprod{x^*-z^*}{x-z}_X< 0$ for all $(x,x^*)\notin\graph A$. Furthermore, for $(x,x^*)\in\graph A$ the infimum is attained in $(z,z^*)=(x,x^*)$). Hence \eqref{eq:monoton:fitzpatrick2} implies that $F_A(x,x^*)\geq \dualprod{x^*}{x}_X$ with equality for $(x,x^*)\in\graph A$. Since $\graph A\neq \emptyset$, this shows that $F_A$ is proper.

        \item \label{thm:monoton:max_surj:ii}
            On the other hand, the definition \eqref{eq:monoton:fitzpatrick1} yields that
            \begin{equation*}
                F_A = (G_A)^* \qquad\text{for}\qquad G_A(z^*,z) = \dualprod{z^*}{z}_X + \delta_{\graph A^{-1}}(z^*,z)
            \end{equation*}
            (since $(z,z^*)\in\graph A$ if and only if $(z^*,z)\in\graph A^{-1}$).
            Furthermore, since $\graph A\neq \emptyset$ was assumed, $F_A$ is the Fenchel conjugate of a proper functional and therefore convex and lower semicontinuous.
    \end{enumerate}

    As a first step, we show that $0\in\range (\partial j+A)$. We set $\Xi\defeq X\times X^*$ as well as $\xi\defeq(x,x^*)\in\Xi$ and consider the functional
    \begin{equation*}
        J_A:\Xi\to\Rbar,\qquad \xi \mapsto F_A(\xi) + \frac12\norm{\xi}_\Xi^2.
    \end{equation*}
    We first note that property \ref{thm:monoton:max_surj:i} implies for all $\xi\in \Xi$ that
    \begin{equation}\label{eq:monoton:fitzpatrick3}
        \begin{aligned}[t]
            J_A(\xi) = F_A(\xi) + \frac12\norm{\xi}_\Xi^2 &= F_A(x,x^*) + \frac12\norm{x}_X^2 +\frac12 \norm{x^*}_{X^*}^2 \\
            &\geq \dualprod{x^*}{x}_X +\frac12 \norm{x}_X^2 +\frac12 \norm{x^*}_{X^*}^2 \\
            &\geq 0,
        \end{aligned}
    \end{equation}
    where the last inequality follows from the Fenchel--Young inequality for $j$ applied to $(x,-x^*)$.
    Furthermore, $J_A$ is proper, convex, lower semicontinuous, and (by \cref{lem:convex:supercoercive}) coercive.
    \Cref{thm:convex:existence} thus yields a $\bar \xi\defeq(\bar x,\bar x^*)\in \Xi$ with $J_A(\bar \xi) = \min_{\xi\in \Xi} J_A(\xi)$, which by \cref{thm:convex:fermat,thm:subdiff:sum,thm:convex:gateaux} satisfies that
    \begin{equation*}
        0\in \partial J_A(\bar \xi)  = \partial\left(\frac12\norm{\bar \xi}_X^2\right) + \partial F_A(\bar \xi),
    \end{equation*}
    i.e., there exists a $\bar \xi^* =(\bar w^*,\bar w)\in \Xi^* \simeq X^*\times X$ (since $X$ is reflexive) such that $\bar \xi^*\in \partial F_A(\bar \xi)$ and $-\bar \xi^* \in \partial(\frac12\norm{\bar \xi}_X^2)$.

    By definition of the subdifferential, we thus have for all $\xi\in \Xi$ that
    \begin{equation*}
        \begin{aligned}
            F_A(\xi) \geq F_A(\bar \xi) + \dualprod{\bar \xi^*}{\xi-\bar \xi}_\Xi
            &= J_A(\bar \xi) + \frac12\norm{\bar \xi^*}_{\Xi^*}^2 + \dualprod{\bar \xi^*}{\xi}_\Xi
            \geq \frac12\norm{\bar \xi^*}_\Xi^2 + \dualprod{\bar \xi^*}{\xi}_\Xi,
        \end{aligned}
    \end{equation*}
    where the second step uses again the Fenchel--Young inequality holding with equality for $(\bar \xi,-\bar\xi^*)$, and the
    last step follows from \eqref{eq:monoton:fitzpatrick3}. Property \ref{thm:monoton:max_surj:i} then implies for all $(x,x^*)\in\graph A$ that
    \begin{equation*}
        \dualprod{x^*}{x}_X =F_A(x,x^*)\geq \frac12\norm{\bar w^*}_{X^*}^2 + \frac12\norm{\bar w}_X^2+ \dualprod{\bar w^*}{x}_X + \dualprod{x^*}{\bar w}_X.
    \end{equation*}
    Adding $\dualprod{\bar w^*}{\bar w}_X$ on both sides and rearranging yields
    \begin{equation}
        \label{eq:monoton:fitzpatrick4}
        \dualprod{x^*-\bar w^*}{x-\bar w}_X \geq \dualprod{\bar w^*}{\bar w}_X +\frac12 \norm{\bar w^*}_{X^*}^2 + \frac12 \norm{\bar w}_X^2 \geq 0,
    \end{equation}
    again by the Fenchel--Young inequality.  The maximal monotonicity of $A$ thus yields that $\bar w^* \in A (\bar w)$, i.e., $(\bar w,\bar w^*)\in\graph A$. Inserting this for $(x,x^*)$ in \eqref{eq:monoton:fitzpatrick4} then shows that
    \begin{equation*}
        \dualprod{\bar w^*}{\bar w}_X +\frac12 \norm{\bar w^*}_{X^*}^2 + \frac12 \norm{\bar w}_X^2 = 0.
    \end{equation*}
    Hence the Fenchel--Young inequality for $\partial j$ holds with equality at $(\bar w,-\bar w^*)$, implying $-\bar w^* \in \partial j(\bar w)$.
    Together, we obtain that $0=-\bar w^* + \bar w^*\in (\partial j + A)(\bar w)$.

    Finally, let $z^*\in X^*$ be arbitrary and set $B:X\setto X^*$, $x\mapsto \{-z^*\}+A(x)$.
    Using the definition, it is straightforward to verify that $B$ is maximally monotone with $\graph B \neq \emptyset$ as well. As we have just shown, there now exists a $\bar x\in X$ with $0\in (\partial j + B)(\bar x) = \partial j(\bar x) + \{-z^*\} + A(\bar x)$, i.e., $z^*\in (\partial j + A)(\bar x)$. Hence $\partial j+A$ is surjective.
\end{proof}

In Hilbert spaces, we can easily show that the converse implication holds as well since in this case, we can identify $X^*$ with $X^*$ using the Riesz isomorphism such that the duality pairing can be replaced with the inner product.
In particular, we can then identify the set $\partial F(x)\subset X^*$ of subderivatives with the corresponding set in $X$ of subgradients (i.e., their Riesz representations) such that in particular $\partial j = \Id$.
\begin{lemma}
    \label{lem:monoton:surj_max}
    Let $X$ be a Hilbert space and let $A:X\setto X$ be monotone. If $\Id+A$ is surjective, then $A$ is maximally monotone.
\end{lemma}
\begin{proof}
    Let $x\in X$ and $x^*\in X$ be such that
    \begin{equation}\label{eq:monoton:max_surj1}
        \iprod{x^*-\tilde x^*}{x-\tilde x}_X \geq 0 \qquad \text{for all }(\tilde x,  \tilde x^*) \in \graph A.
    \end{equation}
    Since $\Id+A$ is surjective, there exist for $x  + x^*\in X$ a $z \in X$ and a $z^* \in A(z)$ such that
    \begin{equation}\label{eq:monoton:max_surj2}
        x+ x^* = z + z^* \in (\Id + A)z.
    \end{equation}
    Inserting $(\tilde x,\tilde x^*)=(z,z^*)$ into \eqref{eq:monoton:max_surj1}
    then yields that
    \begin{equation*}
        0\leq \iprod{x^* - z^*}{x-z}_X = \iprod{z-x}{x-z}_X = -\norm{x-z}_X^2\leq 0,
    \end{equation*}
    i.e., $x=z$. From \eqref{eq:monoton:max_surj2} we further obtain that $x^* = z^*\in A(z) = A(x)$ and hence that $A$ is maximally monotone.
\end{proof}

\begin{remark}
    The converse implication also holds if $X$ is a reflexive Banach space, although the proof in this case requires heavy machinery from the geometry of Banach spaces (in particular, Lindenstrauss's Theorem); see \cite[Theorems III.2.9, II.1.8, V.3.11]{Cioranescu}.
\end{remark}

\section{Resolvents and proximal points}\label{sec:monotone:resolvents}

The proof of \cref{thm:monoton:subdiff} is based on associating to any $x^*\in \partial F(x)$ an element $\bar z\in X$ as the minimizer of a suitable functional. If $X$ is a Hilbert space, this functional is even strictly convex and hence the minimizer $\bar z$ is unique. This property can be exploited to define a new \emph{single-valued} mapping that is more useful for algorithms than the set-valued subdifferential mapping.
For this purpose, we restrict the discussion in the remainder of this chapter to Hilbert spaces (but see \cref{rem:monotone:banach} below), where we again identify $X^*$ with $X$ and subderivatives with subgradients using the Riesz isomorphism.
By the same token, we will also use the same notation for inner products as for duality pairings to avoid the danger of confusing pairs of elements $(x,x^*)\in\graph \partial F$ with their inner product.

We can then define for a maximally monotone operator $A:X\setto X$ with $\graph A\neq \emptyset$ the \term{resolvent}
\begin{equation*}
    \calR_A : X\setto X,\qquad \calR_A(x) = (\Id + A)^{-1}x,
\end{equation*}
as well as for a proper, convex, and lower semicontinuous functional $F:X\to\Rbar$ the \term[mapping!proximal point]{proximal point mapping}
\begin{equation}
    \label{eq:proximal:proximal}
    \prox_F:X\to X,\qquad \prox_F(x) = \argmin_{z\in X}~ \frac12\norm{z-x}_X^2 + F(z).
\end{equation}
Since a similar argument as in the proof of \cref{thm:monoton:subdiff} shows that $w\in\calR_{\partial F}(x)$ is equivalent to the necessary and sufficient conditions for the \term[point!proximal]{proximal point} $w$ to be a minimizer of the strictly convex functional in \eqref{eq:proximal:proximal}, we have that
\begin{equation}
    \label{eq:proximal:resolvent}
    \prox_F = (\Id +\partial F)^{-1} = \calR_{\partial F}.
\end{equation}

Resolvents of monotone and, in particular, maximally monotone operators have useful properties.
\begin{lemma}\label{lem:proximal:firmly-nonexpansive}
    If $A:X\setto X$ is monotone, $\calR_A$ is \term[operator!nonexpansive!firmly]{firmly nonexpansive}, i.e.,
    \begin{equation}\label{eq:proximal:firmly-nonexpansive}
        \norm{z_1-z_2}_X^2 \leq \iprod{x_1-x_2}{z_1-z_2}_X\qquad\text{for all }(x_1,z_1),(x_2,z_2)\in \graph \calR_A
    \end{equation}
    or equivalently,
    \begin{multline}\label{eq:proximal:firmly-nonexpansive2}
        \norm{z_1-z_2}_X^2 + \norm{(x_1-z_1)-(x_2-z_2)}_X^2 \leq \norm{x_1-x_2}_X^2\\
        \text{for all }(x_1,z_1),(x_2,z_2)\in \graph \calR_A.
    \end{multline}
\end{lemma}
\begin{proof}
    Let $x_1,x_2\in \dom \calR_A$ as well as $z_1\in \calR_A(x_1)$ and $z_2\in \calR_A(x_2)$. By definition of the resolvent, this implies that $x_1-z_1\in A(z_1)$ and $x_2-z_2\in A(z_2)$. By the monotonicity of $A$, we thus have
    \begin{equation*}
        0\leq \iprod{(x_1-z_1)-(x_2-z_2)}{z_1-z_2}_X,
    \end{equation*}
    which after rearranging yields \eqref{eq:proximal:firmly-nonexpansive}. The equivalence of \eqref{eq:proximal:firmly-nonexpansive} and \eqref{eq:proximal:firmly-nonexpansive2} is straightforward to verify using binomial expansion.
\end{proof}

\begin{corollary}\label{lem:proximal:lipschitz}
    Let $A:X \setto X$ be maximally monotone with $\graph A\neq \emptyset$. Then $\calR_A:X\to X$ is single-valued and Lipschitz continuous with constant $L=1$.
\end{corollary}
\begin{proof}
    Under the stated assumptions, $\Id+A$ is surjective by \cref{thm:monoton:max_surj}, which implies that $\calR_A(x) \neq \emptyset$ for all $x\in X$, i.e., $\dom \calR_A=X$. Let now $x\in X$ and $z_1,z_2\in \calR_A(x)$. Since $A$ is monotone, $\calR_A$ is nonexpansive by \cref{lem:proximal:firmly-nonexpansive}, which yields both single-valuedness of $\calR_A$ (by taking $x_1=x_2=x$ implies $z_1=z_2$) and its Lipschitz continuity (by applying the Cauchy--Schwarz inequality).
\end{proof}
In particular, by \cref{thm:monoton:subdiff}, this holds for the proximal point mapping $\prox_F:X\to X$ of a proper, convex, and lower semicontinuous functional $F:X\to\R$.
\begin{remark}
    Conversely, it can be shown that every nonexpansive mapping $T:X\to X$ that satisfies $T(x) \in \subdiff G(x)$ for all $x \in X$ for some proper, convex, and lower semicontinuous functional $G:X\to \Rbar$ is the proximal mapping of some proper, convex, and lower semicontinuous functional $F:X\to\Rbar$; see \cite{moreau1965proximite,gribonval2019characterization}.
\end{remark}

Lipschitz continuous mappings with constant $L=1$ are also called \term[operator!nonexpansive]{nonexpansive}. Such mappings furnish a useful class of maximally monotone mappings.
\begin{lemma}\label{lem:monotone:nonexpansive}
 Let $T: X\to X$ be nonexpansive. Then $\Id+\alpha T$ is maximally monotone for every $\alpha\in [-1,1]$.
\end{lemma}
\begin{proof}
    Let $\alpha\in [-1,1]$ be arbitrary. We first show that $\Id+\alpha T$ is monotone: Since $T$ is single-valued, we have for every $x,\tilde x\in X$ that
    \begin{equation*}
        \begin{aligned}
            \iprod{x+\alpha T(x)-(\tilde x + \alpha T(\tilde x))}{x-\tilde x}_X &=
            \norm{x-\tilde x}_X^2 + \alpha \iprod{T(x)-T(\tilde x)}{x-\tilde x}_X \\
            &\geq \norm{x-\tilde x}_X \left(\norm{x-\tilde x}_X-|\alpha|\norm{T(x)-T(\tilde x)}_X\right)\\
            &\geq 0
        \end{aligned}
    \end{equation*}
    since $|\alpha|\leq 1$ and $\norm{T(x)-T(\tilde x)}_X\leq \norm{x-\tilde x}_X$.

    To show maximal monotonicity, let $x\in X$ and $z\in X$ be such that
    \begin{equation*}
        \iprod{z-(\tilde x+\alpha T(\tilde x))}{x-\tilde x}_X \geq 0 \qquad\text{for all } \tilde x\in X.
    \end{equation*}
    Let now $t>0$ be arbitrary and consider $\tilde x = x_t \defeq x+t(z-(x+\alpha T(x)))$ such that $z-(x+\alpha T(x)) = -\frac1t(x-x_t)$ and hence
    \begin{equation}\label{eq:monotone:nonexpansive-1}
        \iprod{z-(x_t + \alpha T(x_t))}{z-(x+\alpha T(x))}_X =
        -\frac1t\iprod{z-(x_t + \alpha T(x_t))}{x-x_t}_X \leq 0
    \end{equation}
    for all $t>0$. Letting now $t\downto 0$, we have that $x_t\to x$ and hence by the Lipschitz continuity of $T$ that $x_t+\alpha T(x_t) \to x+\alpha T$. Passing to the limit in \eqref{eq:monotone:nonexpansive-1} thus yields that
    \begin{equation*}
        \norm{z-(x+\alpha T(x))}_X^2 = \iprod{z-(x + \alpha T(x))}{z-(x+\alpha T(x))}_X \leq 0
    \end{equation*}
    and hence that $z=x+\alpha T(x)$, showing the claimed maximal monotonicity.
\end{proof}

A related concept that is sometimes used is the following. A mapping $T: X \to X$ is called \term[operator!$\alpha$-averaged]{$\alpha$-averaged} for some $\alpha \in (0, 1)$, if there exists a nonexpansive mapping $J:X\to X$ such that $T=(1-\alpha)\Id+\alpha J$. This allows the following characterization of all firmly nonexpansive mappings.
\begin{lemma}
    \label{lemma:proximal:averaged}
    A mapping $T: X \to X$ is firmly nonexpansive if and only if $T$ is $(1/2)$-averaged.
\end{lemma}
\begin{proof}
    If $T$ is $(1/2)$-averaged, then $T=\frac12(\Id+J)$ for some nonexpansive mapping $J$.
    Inserting the definition and using the nonexpansivity of $J$ then implies that
    \begin{equation*}
        \begin{aligned}
        \norm{T(x)-T(y)}_X^2
        & =\frac{1}{4}\left(\norm{J(x)-J(y)}_X^2 + 2\iprod{J(x)-J(y)}{x-y}_X + \norm{x-y}_X^2\right)
        \\
        & \le \frac{1}{2}\left(\iprod{J(x)-J(y)}{x-y}_X+\norm{x-y}_X^2\right)
        \\
        &
        = \iprod{T(x)-T(y)}{x-y}_X,
        \end{aligned}
    \end{equation*}
    i.e., that $T$ is firmly nonexpansive.

    Conversely, let $T$ be firmly nonexpansive. To show that $T$ is $(1/2)$-averaged, it suffices to show that $J \defeq 2T-\Id$ is nonexpansive. This follows again by inserting the definition and using the firm nonexpansivity of $T$ to estimate
    \begin{equation*}
        \begin{aligned}[b]
            \norm{J(x)-J(y)}_X^2
            & = 4\norm{T(x)-T(y)}_X^2 - 4\iprod{T(x)-T(y)}{x-y}_X + \norm{x-y}_X^2
            \\
            & \le \norm{x-y}_X^2.
        \end{aligned}
        \qedhere
    \end{equation*}
\end{proof}

\bigskip

We now return to proximal point mappings.
The following result is central in optimization since it provides a \emph{single-valued} albeit implicit characterization of subgradients.
\begin{lemma}\label{lem:proximal:subdiff}
    Let $F:X\to \Rbar$ be proper, convex, and lower semicontinuous, and $x,x^*\in X$. Then for any $\gamma>0$,
    \begin{equation*}
        x^*\in \partial  F(x) \quad\equivalent\quad x = \prox_{\gamma F}(x + \gamma x^*).
    \end{equation*}
\end{lemma}
\begin{proof}
    Multiplying both sides of the subdifferential inclusion by $\gamma>0$ and adding $x$ yields that
    \begin{equation*}
        \begin{aligned}[b]
            x^* \in \partial F(x) &\Leftrightarrow x + \gamma x^* \in (\Id + \gamma\partial F)(x)\\
                                  &\Leftrightarrow x \in (\Id + \gamma\partial F)^{-1}(x+\gamma x^*)\\
                                  &\Leftrightarrow x = \prox_{\gamma F}(x+\gamma x^*),
        \end{aligned}
    \end{equation*}
    where in the last step we have used that $\gamma\partial F=\partial(\gamma F)$ by \cref{lem:convex:subdiff_calc}\,\ref{lem:convex:subdiff_calc:i} and hence that $\prox_{\gamma F} = \calR_{\partial(\gamma F)} = \calR_{\gamma\partial F}$.
\end{proof}

By applying \cref{lem:proximal:subdiff} to the Fermat principle $0\in\partial F(\bar x)$, we obtain the following fixed-point characterization of minimizers of $F$.
\begin{corollary}\label{thm:proximal:fermat}
    Let $F:X\to \Rbar$ be proper, convex and lower semicontinuous, and $\gamma>0$ be arbitrary. Then $\bar x\in \dom F$ is a minimizer of $F$ if and only if
    \begin{equation*}
        \bar x = \prox_{\gamma F}(\bar x).
    \end{equation*}
\end{corollary}
This simple result should not be underestimated: It allows replacing (explicit) set inclusions in optimality conditions by equivalent (implicit) Lipschitz continuous equations, which (as we will show in following chapters) can be solved by fixed-point iteration or Newton-type methods.

We can also derive a generalization of the orthogonal decomposition of vector spaces.
\begin{theorem}[Moreau decomposition]\label{thm:proximal:moreau}\index{theorem!Moreau}
    Let $F:X\to \Rbar$ be proper, convex, and lower semicontinuous. Then we have for all $x\in X$ that
    \begin{equation*}
        x = \prox_F(x) + \prox_{F^*}(x).
    \end{equation*}
\end{theorem}
\begin{proof}
    Setting $w = \prox_{F}(x)$, \cref{lem:proximal:subdiff,lem:convex:fenchel-young} for $\gamma=1$ imply that
    \begin{equation*}
        \begin{aligned}[b]
            w = \prox_F(x) = \prox_F(w+(x-w)) &\Leftrightarrow x-w \in \partial F(w)\\
                                              &\Leftrightarrow w \in \partial F^*(x-w)\\
                                              &\Leftrightarrow x-w = \prox_{F^*}((x-w)+w)
            =\prox_{F^*}(x).
        \end{aligned}
        \qedhere
    \end{equation*}
\end{proof}

The following calculus rules will prove useful.
\begin{lemma}\label{lem:proximal:calculus}
    Let $F:X\to\Rbar$ be proper, convex, and lower semicontinuous. Then,
    \begin{enumerate}
        \item \label{lem:proximal:calculus:i}
            for $\lambda \neq 0$ and $z\in X$ we have with $H(x) \defeq F(\lambda x + z)$ that
            \begin{equation*}
                \prox_H(x) = \lambda^{-1} (\prox_{\lambda^2 F}(\lambda x + z) -z);
            \end{equation*}
        \item \label{lem:proximal:calculus:ii}\label{item:proximal:calculus:conjugate}  for $\gamma>0$ we have that
            \begin{equation*}
                \prox_{\gamma F^*}(x) = x - \gamma\, \prox_{\gamma^{-1} F}(\gamma^{-1}x);
            \end{equation*}
        \item \label{lem:proximal:calculus:iii}
            for proper, convex, lower semicontinuous $G:Y\to\Rbar$ and $\gamma>0$ we have with $H(x,y) \defeq F(x)+G(y)$ that
            \begin{equation*}
                \prox_{\gamma H}(x,y) = \begin{pmatrix}\prox_{\gamma F}(x)\\\prox_{\gamma G}(y)
                \end{pmatrix}.
            \end{equation*}
    \end{enumerate}
\end{lemma}
\begin{proof}
    \emph{\ref{lem:proximal:calculus:i}:} By definition,
    \begin{equation*}
        \prox_H(x) = \argmin_{w\in X} \frac12\norm{w-x}^2_X + F(\lambda w + z) =: \bar w.
    \end{equation*}
    Now note that since $X$ is a vector space,
    \begin{equation*}
        \min_{w\in X} \frac12\norm{w-x}^2_X + F(\lambda w + z) = \min_{v\in X} \frac12\norm{\lambda^{-1}(v-z) -x}^2_X + F(v),
    \end{equation*}
    and the respective minimizers $\bar w$ and $\bar v$ are related by $\bar v=\lambda \bar w+z$. The claim then follows from
    \begin{equation*}
        \begin{aligned}
            \bar v &= \argmin_{v\in X}~ \frac12\norm{\lambda^{-1}(v-z) -x}^2_X + F(v)\\
                   &= \argmin_{v\in X}~ \frac1{2\lambda^2}\norm{v -(\lambda x+z)}^2_X + F(v)\\
                   &= \argmin_{v\in X}~ \frac1{2}\norm{v -(\lambda x+z)}^2_X + \lambda^2F(v)\\
                   &= \prox_{\lambda^2 F}(\lambda x+z).
        \end{aligned}
    \end{equation*}

    \emph{\ref{lem:proximal:calculus:ii}:} \Cref{thm:proximal:moreau}, \cref{lem:convex:fenchel_calc}\,\ref{lem:convex:fenchel_calc:i}, and \ref{lem:proximal:calculus:i} for $\lambda=\gamma^{-1}$ and $z=0$ together imply that
    \begin{equation*}
        \begin{aligned}
            \prox_{\gamma F}(x) &= x- \prox_{(\gamma F)^*}(x)\\
                                &= x - \prox_{\gamma F^*\circ (\gamma^{-1}\Id)}(x)\\
                                &= x - \gamma \,\prox_{\gamma( \gamma^{-2} F^*)}(\gamma^{-1}x).
        \end{aligned}
    \end{equation*}
    Applying this to $F^*$ and using that $F^{**} = F$ by \cref{thm:convex:moreau}\,\ref{thm:convex:moreau:iii} now yields the claim.

    \emph{\ref{lem:proximal:calculus:iii}:} By definition of the norm on the product space $X\times Y$, we have that
    \begin{equation*}
        \begin{aligned}
            \prox_{\gamma H}(x,y) &= \argmin_{(u,v)\in X\times Y}  \frac12 \norm{(u,v) - (x,y)}^2_{X\times Y} + \gamma H(u,v)\\
                                  &= \argmin_{u\in X, v\in Y}  \left(\frac12 \norm{u - x}^2_{X} + \gamma F(u) \right)
            +\left(\frac12 \norm{v - y}^2_{Y} + \gamma G(v)\right).
        \end{aligned}
    \end{equation*}
    Since there are no mixed terms in $u$ and $v$, the two terms in parentheses can be minimized separately. Hence, $\prox_{\gamma H}(x,y) = (\bar u, \bar v)$ for
    \begin{equation*}
        \begin{aligned}[b]
            \bar u =  \argmin_{u\in X}  \frac12 \norm{u - x}^2_{X} + \gamma F(u)  = \prox_{\gamma F(x)},\\
            \bar v =  \argmin_{v\in Y}  \frac12 \norm{v - y}^2_{Y} + \gamma G(v)  = \prox_{\gamma G(x)}.
        \end{aligned}
        \qedhere
    \end{equation*}
\end{proof}

Computing proximal points is difficult in general since evaluating $\prox_F$ by its definition entails minimizing $F$. In some cases, however, it is possible to give an explicit formula for $\prox_F$.
\begin{example}\label{ex:proximal:reell}
    We first consider scalar functions $f:\R\to\Rbar$.
    \begin{enumerate}
        \item \label{ex:proximal:reell:i}
            $f(t) = \frac12|t|^2$. Since $f$ is differentiable, we can set the derivative of $\frac12(s-t)^2+\frac\gamma2s^2$ to zero and solve for $s$ to obtain $\prox_{\gamma f}(t) = (1+\gamma)^{-1}t$.

        \item \label{ex:proximal:reell:ii}
            $f(t) = |t|$. By \cref{ex:convex:subdiff_abs}, we have that $\partial f(t) = \sign(t)$; hence
            $s \defeq\prox_{\gamma f}(t)=(\Id + \gamma\sign)^{-1}(t)$
            if and only if $t\in \{s\}+\gamma\sign(s)$. Let $t$ be given and assume this holds for some $\bar s$. We now proceed by case distinction.
            \begin{enumerate}[label={Case }\arabic*:, align=left]
                \item $\bar s>0$. This implies that $t = \bar s+\gamma$, i.e., $\bar s = t-\gamma$, and hence that $t>\gamma$.
                \item $\bar s<0$. This implies that $t = \bar s-\gamma$, i.e., $\bar s = t+\gamma$, and hence that $t<-\gamma$.
                \item $\bar s=0$. This implies that $t\in \gamma[-1,1]= [-\gamma,\gamma]$.
            \end{enumerate}
            Since this yields a complete and disjoint case distinction for $t$, we can conclude that
            \begin{equation*}
                \prox_{\gamma f}(t) = \begin{cases} t-\gamma & \text{if }t>\gamma,\\
                    0 &\text{if } t\in [-\gamma,\gamma],\\
                    t+\gamma &\text{if } t < -\gamma.
                \end{cases}
            \end{equation*}
            This mapping is also known as the \term[operator!soft-shrinkage]{soft-shrinkage} or \term[operator!soft-thresholding|see{operator, soft-shrinkage}]{soft-thresholding} operator.

        \item \label{ex:proximal:reell:iii}
          $f(t) = \delta_{[-1,1]}(t)$. We can proceed here in the same way as in \ref{ex:proximal:reell:ii}, but for the sake of variety we instead use \cref{lem:proximal:calculus}\,\ref{lem:proximal:calculus:ii} to compute the proximal point mapping from that of $f^*(t) = |t|$ (see \cref{ex:convex:fenchel}\,\ref{ex:convex:fenchel:iii}) via
            \begin{equation*}
                \begin{aligned}
                    \prox_{\gamma f}(t) &= t-\gamma\,\prox_{\gamma^{-1} f^*}(\gamma^{-1} t) \\
                                        &=\begin{cases}
                    t - \gamma(\gamma^{-1} t - \gamma^{-1})  &\text{if } \gamma^{-1}t > \gamma^{-1},\\
                    t - 0  &\text{if } \gamma^{-1}t\in [-\gamma^{-1},\gamma^{-1}],\\
                    t - \gamma(\gamma^{-1} t + \gamma^{-1})  &\text{if } \gamma^{-1}t <- \gamma^{-1}
                \end{cases}\\
                &=\begin{cases}
                \phantom{-}1&\text{if } t>1, \\
                \phantom{-}t   &\text{if } t\in [-1,1],\\
                -1&\text{if } t<-1.
            \end{cases}
        \end{aligned}
    \end{equation*}
    For every $\gamma>0$, the proximal point of $t$ is thus its projection onto $[-1,1]$.
    \end{enumerate}
\end{example}
\begin{example}\label{ex:proximal:rn}
    We can generalize \cref{ex:proximal:reell} to $X=\R^N$ (endowed with the Euclidean inner product) by applying \cref{lem:proximal:calculus}\,\ref{lem:proximal:calculus:iii} $N$ times. We thus obtain componentwise
    \begin{enumerate}
        \item \label{ex:proximal:rn:i}
            for $F(x)=\frac12\norm{x}_2^2=\sum_{i=1}^N \tfrac12 x_i^2$ that
            \begin{equation*}
                [\prox_{\gamma F}(x)]_i = \left(\frac1{1+\gamma}\right) x_i,\quad 1\leq i\leq N;
            \end{equation*}
        \item \label{ex:proximal:rn:ii}
            for $F(x) = \norm{x}_1=\sum_{i=1}^N|x_i|$ that
            \begin{equation*}
                [\prox_{\gamma F}(x)]_i = (|x_i|-\gamma)^+\sign(x_i),\quad 1\leq i\leq N;
            \end{equation*}
        \item \label{ex:proximal:rn:iii}
            for $F(x) = \delta_{\B_\infty}(x)=\sum_{i=1}^N \delta_{[-1,1]}(x_i)$ that
            \begin{equation*}
                [\prox_{\gamma F}(x)]_i = x_i-(x_i-1)^+-(x_i+1)^-
                = \frac{x_i}{\max\{1,|x_i|\}},\qquad 1\leq i\leq N.
            \end{equation*}
    \end{enumerate}
    Here we have used the convenient notation $(t)^+ \defeq \max\{t,0\}$ and $(t)^-\defeq\min\{t,0\}$.
\end{example}
Many more examples of projection operators and proximal mappings can be found in \cite{Cegielski:2012}, \cite[\S\,6.5]{Boyd:2014}, \cite{Beck:2017}, as well as at \url{https://www.proximity-operator.net}.

Since the subdifferential of convex integral functionals can be evaluated pointwise by \cref{thm:lebesgue:subdiff}, the same holds for the definition \eqref{eq:proximal:resolvent} of the proximal point mapping.
\begin{corollary}\label{lem:lebesgue:proximal}
    Let $f:\R\to\Rbar$ be proper, convex, and lower semicontinuous, and $F:L^2(\Omega)\to\Rbar$ be defined by superposition as in \cref{lem:lebesgue:lsc}. Then we have for all $\gamma>0$ and $u\in L^2(\Omega)$ that
    \begin{equation*}
        [\prox_{\gamma F}(u)](x) = \prox_{\gamma f}(u(x))\qquad\text{for almost every }x\in \Omega.
    \end{equation*}
\end{corollary}
\begin{example}\label{ex:proximal:hilbert}
    Let $X$ be a Hilbert space. Similarly to \cref{ex:proximal:reell} one can show that
    \begin{enumerate}
        \item \label{ex:proximal:hilbert:i}
            for $F=\frac12\norm{\cdot}_X^2 = \frac12\iprod{\cdot}{\cdot}_X$, that
            \begin{equation*}
                \prox_{\gamma F}(x) = \left(\frac{1}{1+\gamma}\right)x;
            \end{equation*}
        \item  \label{ex:proximal:hilbert:ii}
            for $F=\norm{\cdot}_X$, using a case distinction as in \cref{thm:subdifferential:norm}, that
            \begin{equation*}
                \prox_{\gamma F}(x) = \left(1-\frac\gamma{\norm{x}_X}\right)^+ x;
            \end{equation*}
        \item  \label{ex:proximal:hilbert:iii}
            for $F=\delta_C$ with $C\subset X$ nonempty, convex, and closed, that by definition
            \begin{equation*}
                \prox_{\gamma F}(x) = \proj_C(x) \defeq \argmin_{z\in C} \norm{z-x}_X
            \end{equation*}
            the \term[projection!metric]{metric projection} of $x$ onto $C$; the proximal point mapping thus generalizes the concept projection onto convex sets. Explicit or at least constructive formulas for the projection onto different classes of sets can be found in \cite[Chapter 4.1]{Cegielski:2012}.
    \end{enumerate}
\end{example}

\begin{remark}\label{rem:monotone:banach}
    The results of this section can be extended to (reflexive) Banach spaces if the identity is replaced by the duality mapping $\partial j:X\setto X^*$ for $j(x)=\frac12\norm{x}_X^2$. If the norm is differentiable (which is the case if the unit ball of $X^*$ is \emph{strictly} convex as for, e.g., $X=L^p(\Omega)$ with $p\in(1,\infty)$), the duality mapping is in fact single-valued \cite[Theorem 2.16]{Cioranescu}, and hence the corresponding resolvent $(\partial j + A)^{-1}$ is well-defined.
    However, the proximal mapping need no longer be Lipschitz continuous, although the definition can be modified to obtain uniform continuity; see \cite{Bacak:2017}. Similarly, the Moreau decomposition (\cref{thm:proximal:moreau}) needs to be modified appropriately; see \cite{Combettes2013}.
    The main difficulty from our point of view, however, lies in the evaluation of the proximal mapping, which then rarely admits a closed form even for simple functionals.
\end{remark}

\chapter{Smoothness and convexity}\label{chap:smoothness}

Before we turn to algorithms for the solution of nonsmooth optimization problems, we derive consequences of convexity for \emph{differentiable} functionals that will be useful in proving convergence of splitting methods for functionals involving a smooth component. In particular, we will show that Lipschitz continuous differentiability is linked via Fenchel duality to strong convexity.

\section{Smoothness}

We now derive useful consequences of Lipschitz differentiability and their relation to convexity.
Recall from \cref{thm:convex:gateaux} that for $F:X\to\Rbar$ convex and Gateaux differentiable, $\partial F(x)=\{DF(x)\}$ (which can be identified with $\{\nabla F(x)\}\subset X$ in Hilbert spaces).
\begin{lemma}
    \label{lemma:smoothness}
    Let $X$ be a Banach space and let $F: X \to \R$ be Gateaux differentiable.
    Consider for $L>0$ the properties:
    \begin{enumerate}
        \item\label{item:smoothness:grad-smoothness}
            The property
            \begin{equation}
                \label{eq:smoothness-grad-smoothness}
                F(y) \le F(x) + \dualprod{D F(y)}{y-x}_{X}-\frac{1}{2L}\norm{D F(x)-D F(y)}_{X^*}^2
                \quad \text{for all }x, y \in X.
            \end{equation}
        \item\label{item:smoothness:coco}
            The \term{co-coercivity} of $D F$ with factor $\inv L$:
            \begin{equation}
                \label{eq:smoothness-coco}
                \inv L \norm{D F(x)-D F(y)}_{X^*}^2
                \le
                \dualprod{D F(x)-D F(y)}{x-y}_{X}
                \quad\text{for all }x, y \in X.
            \end{equation}
        \item\label{item:smoothness:lipschitz}
            Lipschitz continuity of $D F$ with factor $L$:
            \begin{equation}
                \label{eq:smoothness:lipschitz}
                \norm{DF(x)-DF(y)}_{X^*} \le L \norm{x-y}_X
                \quad\text{for all }x, y \in X.
            \end{equation}
        \item\label{item:smoothness:smoothness}
            The \term{smoothness} (also known as \term[lemma!descent]{descent lemma}) of $F$ with factor $L$:
            \begin{equation}
                \label{eq:smoothness:smoothness}
                F(x+h) \le F(x) + \dualprod{D F(x)}{h}_{X}+\frac{L}{2}\norm{h}_X^2
                \quad\text{for all }x, h \in X.
            \end{equation}
        \item\label{item:smoothness:uniform}
            The \term[smoothness!uniform]{uniform smoothness} of $F$ with factor $L$:
            \begin{multline}
                \label{eq:smoothness-uniform}
                F(\lambda x + (1-\lambda)y) + \lambda(1-\lambda)\frac{L}{2}\norm{x-y}_X^2 \\
                \geq \lambda F(x) + (1-\lambda)F(y)
                \quad\text{for all }x,y \in X, \lambda\in[0,1].
            \end{multline}
    \end{enumerate}
    Then \ref{item:smoothness:grad-smoothness} $\implies$ \ref{item:smoothness:coco} $\implies$ \ref{item:smoothness:lipschitz} $\implies$ \ref{item:smoothness:smoothness} $\equivalent$ \ref{item:smoothness:uniform}.
    If $F$ is convex and $X$ is reflexive, then all the properties are equivalent.
\end{lemma}

\begin{proof}
    \emph{\ref{item:smoothness:grad-smoothness} $\implies$ \ref{item:smoothness:coco}:}
    Summing the estimate \eqref{eq:smoothness-grad-smoothness} with the same estimate with $x$ and $y$ exchanged yields \eqref{eq:smoothness-coco}.

    \ref{item:smoothness:coco} $\implies$ \ref{item:smoothness:lipschitz}: This follows immediately from \eqref{eq:functan:cs_banach}.

    \emph{\ref{item:smoothness:lipschitz} $\implies$ \ref{item:smoothness:smoothness}:}
    When $DF$ is Lipschitz continuous, $F$ is continuously differentiable with $F'=DF$.
    We can thus apply the mean value theorem (\cref{thm:frechet:mean}) to obtain that
    \begin{equation*}
        \begin{aligned}
            F(x+h)  - F(x) - \dualprod{D F(x)}{h}_{X}
            &=\int_0^1 \dualprod{D F(x+th)}{h}_{X} \ddd t - \dualprod{D F(x)}{h}_{X}
            \\
            &
            =\int_0^1 \dualprod{D F(x+th)-D F(x)}{h}_{X} \ddd t\\
            &\le \int_0^1 t \ddd t \cdot L\norm{h}_X^2
            = \frac{L}{2} \norm{h}_X^2.
        \end{aligned}
    \end{equation*}

    \emph{\ref{item:smoothness:smoothness} $\implies$ \ref{item:smoothness:uniform}:}
    Set $x_\lambda\defeq\lambda x+ (1-\lambda)y$. Multiplying \eqref{eq:smoothness:smoothness} first for $x=x_\lambda$ and $h=x-x_\lambda=(1-\lambda)(x-y)$ with $\lambda$ and then for $x=x_\lambda$ and $h=y-x_\lambda=\lambda(y-x)$ with $1-\lambda$ and adding the results yields \eqref{eq:smoothness-uniform}.

    \emph{\ref{item:smoothness:uniform} $\implies$ \ref{item:smoothness:smoothness}:}
    This follows by taking $x+h$ in place of $x$ and $x$ in place of $y$ in \eqref{eq:smoothness-uniform}, dividing by  $\lambda>0$, and taking the limit $\lambda \to 0$.

    \emph{\ref{item:smoothness:smoothness} $\implies$ \ref{item:smoothness:grad-smoothness} when $F$ is convex and $X$ is reflexive:}
    Since $F$ is convex, we have from \cref{thm:convex:gateaux} that
    \begin{equation*}
        \dualprod{DF(y)}{(x+h)-y}_X \leq F(x+h)-F(y).
    \end{equation*}
    Combining this with \eqref{eq:smoothness:smoothness} yields
    \begin{equation}
        \label{eq:smoothness:uniform-to-first-property:0}
        \begin{aligned}[t]
            F(y) &\leq F(x) + \dualprod{DF(x)}{h}_X - \dualprod{DF(y)}{(x+h)-y}_X + \frac{L}{2}\norm{h}_X^2\\
            &= F(x) + \dualprod{DF(y)}{y-x}_X + \dualprod{DF(x)-DF(y)}{h}_X +\frac{L}{2} \norm{h}_X^2.
        \end{aligned}
    \end{equation}
    Let $z^* \defeq -\inv L(DF(x)-DF(y))$.
    Since $X$ is reflexive, the analytic Hahn--Banach \cref{thm:functan:hb-extension} yields (after multiplication by $\norm{z^*}_{X^*}$) an $h \in X$ such that
    \begin{equation*}
        \norm{h}_X = \norm{z^*}_{X^*}
        \quad\text{and}\quad
        \dualprod{z^*}{h}_X=\norm{z^*}_{X^*}^2.
    \end{equation*}
    Consequently, continuing from \eqref{eq:smoothness:uniform-to-first-property:0},
    \begin{equation*}
        \begin{aligned}[t]
            F(y)
            &
            \le F(x) + \dualprod{D F(y)}{y-x}_{X} - L\dualprod{z^*}{h}_X +\frac{L}{2} \norm{h}_X^2
            \\
            &
            = F(x) + \dualprod{D F(y)}{y-x}_{X} - \frac{L}{2} \norm{z^*}_{X^*}^2
            \\
            &
            = F(x) + \dualprod{D F(y)}{y-x}_{X}-\frac{1}{2L}\norm{D F(x)-D F(y)}_{X^*}^2.
        \end{aligned}
    \end{equation*}
    This proves \eqref{eq:smoothness-grad-smoothness}.
\end{proof}
Due to \cref{lemma:smoothness}\,\ref{item:smoothness:smoothness}, it is customary to refer to functionals with Lipschitz continuous derivative with factor $L$ as \term[functional!$L$-smooth]{$L$-smooth}.

The next ``smoothness three-point corollary'' will be valuable for the study of splitting methods that involve a smooth component function.

\begin{corollary}
    \label{cor:smoothness:three-point}
    Let $X$ be a reflexive Banach space and let $F:X\to\R$ be convex and Gateaux differentiable. Then for any $L>0$, the following properties are equivalent:
    \begin{enumerate}[label=(\roman*)]
        \item\label{itm:smoothness:three-point:twopoint}
            $F$ has $\inv L$-co-coercive derivative (or any of the equivalent properties of \cref{lemma:smoothness}).
        \item\label{itm:smoothness:three-point:smoothness}
            The \term[smoothness!three-point]{three-point smoothness}
            \begin{equation}
                \label{eq:smoothness:three-point:smoothness}
                \dualprod{D F(z)}{x-\realoptx}_X
                \ge
                F(x)-F(\realoptx) -  \frac{L}{2}\norm{x-z}_X^2
                \quad\text{for all }\realoptx, z, x \in X.
            \end{equation}
        \item\label{itm:smoothness:three-point:monotonicity}
            The \term[monotonicity!three-point]{three-point monotonicity}
            \begin{equation}
                \label{eq:smoothness:three-point:monotonicity}
                \dualprod{D F(z)-D F(\realoptx)}{x-\realoptx}_X \ge
                -\frac{L}{4}\norm{x-z}_X^2
                \quad\text{for all }\realoptx, z, x \in X.
            \end{equation}
    \end{enumerate}
\end{corollary}

\begin{proof}
    \emph{\cref{itm:smoothness:three-point:twopoint} $\implies$ \cref{itm:smoothness:three-point:smoothness}:}
    Since $\grad F$ is $\inv L$-co-coercive, using \cref{lemma:smoothness}, we have the $L$-smoothness
    \begin{equation*}
        \label{eq:smoothness:corollary-proof:smoothness}
        F(z)-F(x) \ge \dualprod{D F(z)}{z-x}_X - \frac{L}{2}\norm{x-z}_X^2.
    \end{equation*}
    By convexity $F(\realoptx)-F(z) \ge \dualprod{D F(z)}{\realoptx-z}_X$.
    Summing up, we obtain \eqref{eq:smoothness:three-point:smoothness}.

    \emph{\cref{itm:smoothness:three-point:twopoint} $\implies$ \cref{itm:smoothness:three-point:monotonicity}:}
    By assumption we have the co-coercivity
    \begin{equation*}
        \dualprod{D F(z)-D F(\realoptx)}{z-\realoptx}_X
        \ge \inv L \norm{D F(z)-D F(\realoptx)}_{X^*}^2.
    \end{equation*}
    Thus, using \eqref{eq:functan:cs_banach} and Young's inequality in the form
    $ab \le \frac{1}{2\alpha}a^2 + \frac{\alpha}{2}b^2$ for $a,b\in \R$ and $\alpha>0$, we obtain
    \begin{equation*}
        \begin{aligned}[b]
            \dualprod{D F(z)-D F(\realoptx)}{x-\realoptx}_X
            &
            =\dualprod{D F(z)-D F(\realoptx)}{z-\realoptx}_X
            +\dualprod{D F(z)-D F(\realoptx)}{x-z}_X
            \\
            &
            \ge
            \inv L \norm{D F(z)-D F(\realoptx)}_{X^*}^2
            -\norm{D F(z)-D F(\realoptx)}_{X^*}\norm{x-z}_X
            \\
            &
            \ge -\frac{L}{4}\norm{x-z}_X^2.
        \end{aligned}
    \end{equation*}
    This is \eqref{eq:smoothness:three-point:monotonicity}

    \emph{\cref{itm:smoothness:three-point:monotonicity} $\implies$ \cref{itm:smoothness:three-point:twopoint}:}
    For the reverse implications, we assume that \eqref{eq:smoothness:three-point:monotonicity} holds and set $z^* \defeq -2\frac1L(DF(z)-DF(\realoptx))$.
    By the assumed reflexivity, we can again apply the analytic Hahn--Banach theorem (\cref{thm:functan:hb-extension}) to obtain an $h \in X$ such that
    \begin{equation*}
        \norm{h}_X = \norm{z^*}_{X^*}
        \quad\text{and}\quad
        \dualprod{z^*}{h}_X=\norm{z^*}_{X^*}^2.
    \end{equation*}
    With $x=z+h$, \eqref{eq:smoothness:three-point:monotonicity} gives
    \begin{equation*}
        \begin{aligned}[t]
            \dualprod{D F(z)-D F(\realoptx)}{z-\realoptx}_X
            &
            \ge
            -\dualprod{D F(z)-D F(\realoptx)}{h}_X
            -\frac{L}{4}\norm{h}_X^2
            \\
            &
            =\frac{L}{2}\dualprod{z^*}{h}_X-\frac{L}{4}\norm{z^*}_{X^*}^2
            \\
            &
            =\frac{L}{4}\norm{z^*}_{X^*}^2
            =\frac{1}{L}\norm{D F(z)-D F(\realoptx)}_{X^*}^2.
        \end{aligned}
    \end{equation*}
    This is the $\inv L$-co-coercivity \eqref{eq:smoothness-coco}.

    \emph{\cref{itm:smoothness:three-point:smoothness} $\implies$ \cref{itm:smoothness:three-point:twopoint}:}
    We take $z=\realoptx$ in \eqref{eq:smoothness:three-point:smoothness} to obtain \eqref{eq:smoothness:smoothness}, and finish with \cref{lemma:smoothness}.
\end{proof}

\section{Strong convexity}

The central notion in this chapter (and later for obtaining higher convergence rates for first-order algorithms) is the following \enquote{quantitative} version of convexity. We say that $F:X\to\Rbar$ is \term[functional!convex!strongly]{strongly convex} with the factor $\gamma>0$ if for all $x,y\in X$ and $\lambda\in [0,1]$,
\begin{equation}\label{eq:smoothness:strong-convexity}
    F(\lambda x+ (1-\lambda)y) + \lambda(1-\lambda)\frac{\gamma}{2}\norm{x-y}_X^2 \leq \lambda F(x) + (1-\lambda)F(y).
\end{equation}
Obviously, strong convexity implies strict convexity, so strongly convex functions have a unique minimizer. If $X$ is a Hilbert space, it is straightforward if tedious to verify by expanding the squared norm that \eqref{eq:smoothness:strong-convexity} is equivalent to $F-\frac\gamma2\norm{\cdot}_X^2$ being convex.

We have the following important duality result that was first shown in \cite{Aze:1995}.

\begin{theorem}\label{thm:smoothness:dual}
    Let $X$ be a Banach space and let $F:X\to\Rbar$ be proper and convex.
    Then the following are true:
    \begin{enumerate}
        \item \label{thm:smoothness:dual:i}
            If $F$ is strongly convex with factor $\gamma$, then $F^*$ is uniformly smooth with factor $\gamma^{-1}$.
        \item \label{thm:smoothness:dual:ii}
            If $F$ is uniformly smooth with factor $L$, then $F^*$ is strongly convex with factor $L^{-1}$.
        \item \label{thm:smoothness:dual:iii}
            If $F$ is lower semicontinuous, then $F$ is uniformly smooth with factor $L$ if and only if $F^*$ is strongly convex with factor $\inv L$.
    \end{enumerate}
\end{theorem}

\begin{proof}
    \emph{\ref{thm:smoothness:dual:i}:} Let $x^*,y^*\in X^*$ and $\alpha_x,\alpha_y\in \R$ with $\alpha_x < F^*(x^*)$ and $\alpha_y < F^*(y^*)$. From the definition of the Fenchel conjugate, there exist $x,y\in X$ such that
    \begin{equation*}
        \alpha_x < \dualprod{x^*}{x}_X-F(x),\qquad
        \alpha_y < \dualprod{y^*}{y}_X-F(y).
    \end{equation*}
    Multiplying the first inequality with $\lambda\in [0,1]$, the second with $(1-\lambda)$, and using the Fenchel--Young inequality \eqref{eq:convex:fenchel-young} in the form
    \begin{equation*}
        0\leq F(x_\lambda) + F^*(x^*_\lambda) - \dualprod{x^*_\lambda}{x_\lambda}_X
    \end{equation*}
    for $x^*_\lambda\defeq\lambda x^* + (1-\lambda)y^*$ and $x_\lambda \defeq \lambda x + (1-\lambda) y$ then yields
    \begin{equation*}
        \begin{aligned}
            \lambda \alpha_x + (1-\lambda)\alpha_y
            &\leq  F(x_\lambda) + F^*(x^*_\lambda) -\lambda F(x) - (1-\lambda)F(y) + \lambda(1-\lambda)\dualprod{x^*-y^*}{x-y}_X\\
            &\leq F^*(x^*_\lambda) + \lambda(1-\lambda)\left(\dualprod{x^*-y^*}{x-y}_X - \frac{\gamma}2\norm{x-y}_{X}^2\right)\\
            &\leq F^*(x^*_\lambda) + \lambda(1-\lambda)\sup_{z\in X}\left\{\dualprod{x^*-y^*}{z}_X - \frac{\gamma}2\norm{z}_{X}^2\right\}\\
            &= F^*(x^*_\lambda) + \lambda(1-\lambda)\frac1{2\gamma}\norm{x^*-y^*}_{X^*}^2,
        \end{aligned}
    \end{equation*}
    where we have used the definition \eqref{eq:smoothness:strong-convexity} of strong convexity in the second inequality and \cref{lem:convex:power-conjugate} together with \cref{lem:convex:fenchel_calc}\,\ref{lem:convex:fenchel_calc:i} in the final equality. Letting now $\alpha_x\to F^*(x^*)$ and $\alpha_y\to F^*(y^*)$, we obtain \eqref{eq:smoothness-uniform} for $F^*$ with $L\defeq\gamma^{-1}$.

    \emph{\ref{thm:smoothness:dual:ii}:}  Let $x^*,y^*\in X^*$ and $\lambda \in [0,1]$. Set again $x^*_\lambda\defeq\lambda x^*+(1-\lambda)y^*$. Then we obtain from the definition of the Fenchel conjugate and \eqref{eq:smoothness-uniform} that for any $x,y\in X$,
    \begin{equation*}
        \begin{aligned}
            \lambda F^*(x^*) + (1-\lambda) F^*(y^*)
            &\geq \lambda \left[\dualprod{x^*}{x+(1-\lambda)y}_X - F(x+(1-\lambda)y)\right]\\
            \MoveEqLeft[-1]
            +(1-\lambda) \left[\dualprod{y^*}{x-\lambda y}_X - F(x-\lambda y)\right]\\
            &\geq \lambda \dualprod{x^*}{x+(1-\lambda)y}_X + (1-\lambda)\dualprod{y^*}{x-\lambda y}_X \\
            \MoveEqLeft[-1] - F(x) - \lambda(1-\lambda)\frac{L}2\norm{y}_X^2\\
            &= \dualprod{x^*_\lambda}{x}_X - F(x) + \lambda(1-\lambda)\left(\dualprod{y^*-x^*}{y}_X - \frac{L}2\norm{y}_X^2\right).
        \end{aligned}
    \end{equation*}
    Taking now the supremum over all $x,y\in X$ and using again \cref{lem:convex:power-conjugate} together with \cref{lem:convex:fenchel_calc}\,\ref{lem:convex:fenchel_calc:i}, we obtain the strong convexity \eqref{eq:smoothness:strong-convexity} with $\gamma \defeq L^{-1}$.

    \emph{\ref{thm:smoothness:dual:iii}:}
    One direction of the claim is clear from \ref{thm:smoothness:dual:ii}.
    For the other direction, if $F^*$ is strongly convex with factor $\inv L$, then its preconjugate $(F^*)_*$ is uniformly smooth with factor $L$ by a proof completely analogous to \ref{thm:smoothness:dual:i}. Then we use \cref{thm:convex:moreau} to see that $F=F^{**} \defeq (F^*)_*$ under the lower semicontinuity assumption.
\end{proof}

Just as convexity of $F$ implies monotonicity of $\partial F$, strong convexity has the following consequences.
\begin{lemma}\label{thm:smoothness:strong-convexity}
    Let $X$ be a Banach space and $F:X\to\Rbar$. Consider the properties:
    \begin{enumerate}
        \item \label{item:strong-convexity}
            $F$ is strongly convex with factor $\gamma>0$.
        \item \label{item:strong-subdifferentiability}
            $F$ is \term[functional!subdifferentiable!strongly]{strongly subdifferentiable} with factor $\gamma$:
            \begin{equation}\label{eq:strong-subdifferentiability}
                F(y)-F(x) \ge \dualprod{x^*}{y-x}_X + \frac{\gamma}{2}\norm{y-x}_X^2
                \quad \text{for all } x, y \in X;\, x^* \in \subdiff F(x).
            \end{equation}
        \item \label{item:strong-monotonicity}
            $\partial F$ is \term[mapping!monotone!strongly]{strongly monotone} with factor $\gamma$:
            \begin{equation}\label{eq:strong-monotonicity}
                \dualprod{y^*-x^*}{y-x}_X\geq \gamma\norm{y-x}_X^2 \quad \text{for all } x,y\in X; \, x^*\in\partial F(x), \,y^*\in \partial F(y).
            \end{equation}
    \end{enumerate}
    Then  \ref{item:strong-convexity} $\implies$ \ref{item:strong-subdifferentiability} $\implies$  \ref{item:strong-monotonicity}. If $X$ is reflexive and $F$ is proper, convex, and lower semicontinuous, then also \ref{item:strong-monotonicity} $\implies$ \ref{item:strong-convexity}.
\end{lemma}

\begin{proof}
    \emph{\ref{item:strong-convexity} $\implies$ \ref{item:strong-subdifferentiability}:} Let $x,y\in X$ and $\lambda\in (0,1)$ be arbitrary. Dividing \eqref{eq:smoothness:strong-convexity} by $\lambda$ and rearranging yields
    \begin{equation*}
        \frac{F(y+\lambda(x-y))-F(y)}{\lambda} \leq F(x) - F(y) - (1-\lambda)\frac{\gamma}{2}\norm{x-y}_X^2.
    \end{equation*}
    Since strongly convex functions are also convex, we can apply \cref{lem:convex:direct}\,\ref{lem:convex:direct:ii} to pass to the limit $\lambda\to 0$ on both sides to obtain
    \begin{equation*}
        F'(y,x-y) \leq F(x)-F(y) - \frac\gamma2 \norm{x-y}_X^2.
    \end{equation*}
    Using \cref{lem:convex:equiv} for $h=x-y$, we thus obtain that for any $y^*\in\partial F(y)$,
    \begin{equation*}
        \dualprod{y^*}{x-y}_X \leq F'(y,x-y) \leq F(x)-F(y) - \frac\gamma2 \norm{x-y}_X^2.
    \end{equation*}
    Exchanging the roles of $x$ and $y$ and rearranging yields \eqref{eq:strong-subdifferentiability}.

    \emph{\ref{item:strong-subdifferentiability} $\implies$ \ref{item:strong-monotonicity}:} Adding \eqref{eq:strong-subdifferentiability} with the same inequality with $x$ and $y$ exchanged immediately yields \eqref{eq:strong-monotonicity}.

    \emph{\ref{item:strong-monotonicity} $\implies$ \ref{item:strong-convexity}:}
    Suppose first that $\subdiff F$ is surjective. Then $\dom \subdiff F^*=X^*$.
    Using the duality between $\subdiff F$ and $\subdiff F^*$ in \cref{lem:convex:fenchel-young}, we rewrite \eqref{eq:strong-monotonicity} as
    \begin{equation}
        \label{eq:smoothness:strong-convexity-dual}
        \dualprod{y^*-x^*}{y-x}_X\geq \gamma\norm{y-x}_X^2 \quad \text{for all } x^*,y^*\in X^*; \, x \in \partial F^*(x^*), \,y \in \partial F^*(y^*).
    \end{equation}
    Taking $y^*=x^*$, this implies that $x=y$, i.e., $\subdiff F^*(x^*)$ is a singleton for all $x^* \in X^*$. Here we use that $\dom \subdiff F^*=X^*$ to avoid the possibility that $\subdiff F^*(x^*)=\emptyset$.
    By \cref{thm:convex:gateaux} it follows that $F^*$ is Gateaux differentiable. Thus \eqref{eq:smoothness:strong-convexity-dual} describes the co-coercivity \eqref{eq:smoothness-coco} of $DF^*$ with factor $\gamma$. By \cref{lemma:smoothness} it follows that $F^*$ is uniformly smooth with factor $\inv\gamma$.
    Consequently, by \cref{thm:smoothness:dual} $F$ is strongly convex with factor $\gamma$.

    If $\subdiff F$ is not surjective, we replace $F$ by $F + \epsilon j$ for the duality mapping $j(x) \defeq \tfrac{1}{2}\norm{x}_X^2$ and some $\epsilon>0$. By \cref{thm:monoton:subdiff} and Minty's theorem (\cref{thm:monoton:max_surj}) now $\subdiff(F+\epsilon j)$ is surjective. It also remains strongly monotone with factor $\gamma$ as $\subdiff j$ is monotone.
    Now, by the above reasoning, $F+\epsilon j$ is strongly convex with factor $\gamma$. Since $\epsilon > 0$ was arbitrary, we deduce from the defining \eqref{eq:smoothness:strong-convexity} that $F$ is strongly convex with factor $\gamma$.
\end{proof}

Note that the factor $\gamma$ enters into the strong monotonicity \eqref{eq:strong-monotonicity} directly rather than as $\frac\gamma2$ as in the strong subdifferentiability \eqref{eq:strong-subdifferentiability} (and strong convexity).

We can also derive a stronger, quantitative, version of the fact that for convex functions, points that satisfy the Fermat principle are minimizers.
\begin{lemma}\label{thm:polyak-lojasewicz}
    Let $X$ be a Banach space and let $F:X\to\Rbar$ be strongly convex with factor $\gamma>0$. Assume that $F$ admits a minimum $M\defeq\min_{x\in X} F(x)$. Then the \term[inequality!Polyak--Łojasewicz]{Polyak--Łojasewicz inequality} holds:
    \begin{equation}\label{eq:polyak-lojasewicz}
        F(x) - M \leq \frac1{2\gamma} \norm{x^*}_{X^*}^2 \quad\text{for all } x\in X, \, x^*\in\partial F(x).
    \end{equation}
\end{lemma}
\begin{proof}
    Let $x\in X$ and $x^*\in\partial F(x)$ be arbitrary. Then from \cref{thm:smoothness:strong-convexity}\,\ref{item:strong-subdifferentiability} we have that
    \begin{equation*}
        -F(x) + \dualprod{x^*}{x-y}_X - \frac{\gamma}{2}\norm{x-y}_X^2 \geq -F(y).
    \end{equation*}
    Taking the supremum over all $y\in X$, noting that this is equivalent to taking the supremum over all $x-y\in X$, and inserting the Fenchel conjugate of the squared norm from \cref{lem:convex:power-conjugate} together with \cref{lem:convex:fenchel_calc}\,\ref{lem:convex:fenchel_calc:i}, we obtain
    \begin{equation*}
        -F(x) + \frac{1}{2\gamma}\norm{x^*}_{X^*}^2 \geq \sup_{y\in X} -F(y) = -\min_{y\in X} F(y)
    \end{equation*}
    and hence, after rearranging, \eqref{eq:polyak-lojasewicz}.
\end{proof}

Comparing the consequences of strong convexity in \cref{thm:smoothness:strong-convexity} and those of uniform smoothness in \cref{lemma:smoothness}, we can already see a certain duality between them: While the former give lower bounds, the latter give upper bounds and vice versa. A simple example is the following
\begin{corollary}
    If $X$ is a Banach space and $F:X\to \R$ is strongly convex with factor $\gamma$ and uniformly smooth with factor $L$, then
    \begin{equation}
        \gamma\norm{x-y}_X^2 \leq \dualprod{DF(x)-DF(y)}{x-y}_X \leq L \norm{x-y}_X^2 \quad\text{for all }x,y\in X.
    \end{equation}
\end{corollary}
\begin{proof}
    The first inequality follows from \cref{thm:smoothness:strong-convexity}\,\ref{item:strong-monotonicity}, while the second follows from \eqref{eq:functan:cs_banach} together with \cref{lemma:smoothness}\,\ref{item:smoothness:lipschitz}.
\end{proof}

The estimates of \cref{cor:smoothness:three-point} can be improved if $F$ is in addition strongly convex.

\begin{corollary}
    \label{cor:smoothness:three-point:sc}
    Let $X$ be a Banach space and let $F:X\to\R$ be strongly convex with factor $\gamma>0$ as well as Lipschitz differentiable with constant $L>0$. Then for any $\alpha>0$,
    \begin{equation}
        \label{eq:smoothness:three-point:smoothness-sc}
        \dualprod{D F(z)}{x-\realoptx}_X
        \ge
        F(x)-F(\realoptx) + \frac{\gamma-\alpha L}{2}\norm{x-\realoptx}_X^2
        -\frac{L}{2\alpha}\norm{x-z}_X^2
        \quad\text{for all }\realoptx, z, x \in X,
    \end{equation}
    as well as
    \begin{equation}
        \label{eq:smoothness:three-point:monotonicity-sc}
        \dualprod{D F(z)-D F(\realoptx)}{x-\realoptx}_X \ge
        (\gamma-\alpha L)\norm{x-\realoptx}_X^2
        -\frac{L}{4\alpha}\norm{x-z}_X^2
        \quad\text{for all }\realoptx, z, x \in X.
    \end{equation}
\end{corollary}
\begin{proof}
    Using the strong subdifferentiability from \cref{thm:smoothness:strong-convexity}\,\ref{item:strong-subdifferentiability}, the Lipschitz continuity of $DF$, \eqref{eq:functan:cs_banach}, and Young's inequality, we obtain
    \begin{equation*}
        \begin{aligned}
            \dualprod{D F(z)}{x-\realoptx}_X
            &
            =\dualprod{D F(x)}{x-\realoptx}_X
            +\dualprod{D F(z)-D F(x)}{x-\realoptx}_X
            \\
            & \ge F(x)-F(\realoptx) + \frac{\gamma}{2}\norm{x-\realoptx}_X^2
            - \frac{\alpha L}{2}\norm{x-\realoptx}_X^2 -\frac{1}{2\alpha L}\norm{D F(z)-D F(x)}_{X^*}^2
            \\
            & \ge F(x)-F(\realoptx) + \frac{\gamma}{2}\norm{x-\realoptx}_X^2- \frac{\alpha L}{2}\norm{x-\realoptx}_X^2
            -\frac{L}{2\alpha}\norm{x-z}_X^2.
        \end{aligned}
    \end{equation*}

    For \eqref{eq:smoothness:three-point:monotonicity-sc}, we can use the strong monotonicity of $DF$ from \cref{thm:smoothness:strong-convexity}\,\ref{item:strong-monotonicity} to estimate analogously
    \begin{equation*}
        \begin{aligned}[b]
            \dualprod{D F(z)-D F(\realoptx)}{x-\realoptx}_X
            &
            =\dualprod{D F(x)-D F(\realoptx)}{x-\realoptx}_X
            +\dualprod{D F(z)-D F(x)}{x-\realoptx}_X
            \\
            & \ge \gamma\norm{x-\realoptx}_X^2 - \alpha L\norm{x-\realoptx}_X^2 -\frac{L}{4\alpha}\norm{x-z}_X^2.
        \end{aligned}
        \qedhere
    \end{equation*}
\end{proof}

\section{Moreau--Yosida regularization}\label{sec:moreau-yosida}

We now look at another way to reformulate optimality conditions using proximal point mappings. Although these are no longer equivalent reformulations, they will serve as a link to the Newton-type methods which will be introduced in \cref{chap:semismooth}.

We again assume that $X$ is a Hilbert space and identify $X^*$ with $X$ via the Riesz isomorphism. Let $A:X\setto X$ be a maximally monotone operator with $\graph A \neq \emptyset$ and $\gamma >0$.
Then we define the \term[approximation, Yosida]{Yosida approximation} of $A$ as
\begin{equation*}
    A_\gamma \defeq \frac1\gamma\left(\Id - \calR_{\gamma A}\right).
\end{equation*}
In particular, the Yosida approximation of the subdifferential of a proper, convex, and lower semicontinuous functional $F:X\to\Rbar$ is given by
\begin{equation}\label{eq:proximal:yosida}
    (\partial F)_\gamma \defeq \frac1\gamma\left(\Id - \prox_{\gamma F}\right),
\end{equation}
which by \cref{lem:proximal:lipschitz,thm:monoton:subdiff} is always Lipschitz continuous with constant $L=\gamma^{-1}$.

An alternative point of view is the following. For a proper, convex, and lower semicontinuous functional $F:X\to\Rbar$ and $\gamma>0$, we define the \term[envelope!Moreau]{Moreau envelope}\footnote{not to be confused with the \emph{convex} envelope $F^\Gamma$!}
\begin{equation}\label{eq:proximal:moreau-envelope}
    F_\gamma :X\to\R,\qquad x\mapsto \inf_{z\in X}~ \frac{1}{2\gamma}\norm{z-x}_X^2 + F(z),
\end{equation}
see \cref{fig:smoothness:moreau-yosida}.
Comparing this with the definition \eqref{eq:proximal:proximal} of the proximal point mapping of $F$, we see that
\begin{equation}\label{eq:proximal:moreau-decomposition}
    F_\gamma(x) = \frac{1}{2\gamma}\norm{\prox_{\gamma F}(x)-x}_X^2 + F(\prox_{\gamma F}(x)).
\end{equation}
(Note that multiplying a functional by $\gamma>0$ does not change its minimizers.)
Hence $F_\gamma$ is indeed well-defined on $X$ and single-valued. Furthermore, we can deduce from \eqref{eq:proximal:moreau-decomposition} that $F_\gamma$ is convex as well.
\begin{lemma}\label{lem:moreau:convex}
    Let $F:X\to\Rbar$ be proper, convex, and lower semicontinuous, and $\gamma>0$. Then $F_\gamma$ is convex.
\end{lemma}
\begin{proof}
    We first show that for any convex $G:X\to\Rbar$, the mapping
    \begin{equation*}
        H:X\times X \to \Rbar,\qquad (x,z)\mapsto F(z)+G(z-x)
    \end{equation*}
    is convex as well. Indeed, for any $(x_1,z_1),(x_2,z_2)\in X\times X$ and $\lambda\in[0,1]$, the convexity of $F$ and $G$ implies that
    \begin{equation*}
        \begin{aligned}
            H(\lambda(x_1,z_1)+(1-\lambda)(x_2,z_2)) &= F\left(\lambda z_1 + (1-\lambda )z_2\right)+G\left(\lambda(z_1-x_1) + (1-\lambda)(z_2-x_2)\right)\\
            &\leq \lambda\left(F(z_1)+G(z_1-x_1)\right)+(1-\lambda)\left(F(z_2)+G(z_2-x_2)\right)\\
            &= \lambda H(x_1,z_1) + (1-\lambda)H(x_2,z_2).
        \end{aligned}
    \end{equation*}
    Let now $x_1,x_2\in X$ and $\lambda\in[0,1]$. Since $F_\gamma(x) = \inf_{z\in X} H(x,z)$ for $G(y)\defeq\frac1{2\gamma}\norm{y}_X^2$, there exist two minimizing sequences $\{z^1_n\}_{n\in\N},\{z^2_n\}_{n\in\N}\subset X$ with
    \begin{equation*}
        H(x_1,z^1_n)\to F_\gamma(x_1),\qquad H(x_2,z^2_n)\to F_\gamma(x_2).
    \end{equation*}
    From the properties of the infimum together with the convexity of $H$, we thus obtain for all $n\in\N$ that
    \begin{equation*}
        \begin{aligned}
            F_\gamma(\lambda x_1 + (1-\lambda)x_2) & \leq H(\lambda (x_1,z^1_n) + (1-\lambda) (x_2,z^2_n) ) \\
            &\leq \lambda H(x_1,z^1_n) + (1-\lambda)H(x_2,z^2_n),
        \end{aligned}
    \end{equation*}
    and passing to the limit $n\to\infty$ yields the desired convexity.
\end{proof}
We will also show later that Moreau--Yosida regularization preserves (global!) Lipschitz continuity.

The next theorem links the two concepts of Moreau envelope and of Yosida approximation and hence justifies the term \term[regularization!Moreau--Yosida]{Moreau--Yosida regularization}.
\begin{theorem}\label{thm:moreau-yosida}
    Let $F:X\to\Rbar$ be proper, convex, and lower semicontinuous, and $\gamma>0$.
    Then $F_\gamma$ is Fréchet differentiable with
    \begin{equation*}
        \nabla (F_\gamma) = (\partial F)_\gamma.
    \end{equation*}
\end{theorem}
\begin{proof}
    Let $x,y\in X$ be arbitrary and set $x^* = \prox_{\gamma F}(x)$ and $y^* = \prox_{\gamma F}(y)$. We first show that
    \begin{equation}\label{eq:proximal:my_frechet1}
        \frac1\gamma \iprod{y^*-x^*}{x-x^*}_X  \leq F(y^*) - F(x^*).
    \end{equation}
    (Note that for proper $F$, the definition of proximal points as minimizers necessarily implies that $x^*,y^*\in\dom F$.)
    To this purpose, consider for $t\in (0,1)$ the point $x^*_t \defeq ty^* + (1-t)x^*$.
    Using the minimizing property of the proximal point $x^*$ together with the convexity of $F$ and completing the square, we obtain that
    \begin{equation*}
        \begin{aligned}
            F(x^*) &\leq F(x^*_t) + \frac1{2\gamma}\norm{x^*_t - x}_X^2 - \frac1{2\gamma}\norm{x^* - x}_X^2\\
            &\leq t F(y^*) + (1-t)F(x^*) - \frac{t}{\gamma} \iprod{x-x^*}{y^*-x^*}_X +  \frac{t^2}{2\gamma}\norm{x^*-y^*}_X^2.
        \end{aligned}
    \end{equation*}
    Rearranging the terms, dividing by $t>0$ and passing to the limit $t\to 0$ then yields \eqref{eq:proximal:my_frechet1}.
    Combining this with \eqref{eq:proximal:moreau-decomposition} implies that
    \begin{equation*}
        \begin{aligned}
            F_\gamma(y) - F_\gamma(x) &= F(y^*) - F(x^*) + \frac1{2\gamma}\left(\norm{y-y^*}_X^2 - \norm{x-x^*}_X^2\right)\\
            &\geq \frac1{2\gamma}\left(2\iprod{y^*-x^*}{x-x^*}_X + \norm{y-y^*}_X^2 - \norm{x-x^*}_X^2\right)\\
            &= \frac1{2\gamma}\left(2\iprod{y-x}{x-x^*}_X + \norm{y-y^*-x+x^*}_X^2\right)\\
            &\geq \frac1\gamma \iprod{y-x}{x-x^*}_X.
        \end{aligned}
    \end{equation*}
    By exchanging the roles of $x^*$ and $y^*$ in \eqref{eq:proximal:my_frechet1} and repeating the above calculations, we obtain that
    \begin{equation*}
        F_\gamma(y) - F_\gamma(x) \leq \frac1\gamma \iprod{y-x}{y-y^*}_X.
    \end{equation*}
    Together, these two inequalities yield that
    \begin{equation*}
        \begin{aligned}
            0 &\leq F_\gamma(y) - F_\gamma(x) - \frac1\gamma \iprod{y-x}{x-x^*}_X\\
            &\leq \frac1\gamma \iprod{y-x}{(y-y^*)-(x-x^*)}_X\\
            &\leq \frac1\gamma\left(\norm{y-x}_X^2 - \norm{y^*-x^*}_X^2\right)\\
            &\leq \frac1\gamma \norm{y-x}_X^2,
        \end{aligned}
    \end{equation*}
    where the next-to-last inequality follows from the firm nonexpansivity of proximal point mappings (\cref{lem:proximal:firmly-nonexpansive}).

    If we now set $y=x+h$ for arbitrary $h\in X$, we obtain that
    \begin{equation*}
        0\leq \frac{F_\gamma(x+h) - F_\gamma(x) - \iprod{\gamma^{-1}(x-x^*)}{h}_X}{\norm{h}_X} \leq \frac 1\gamma \norm{h}_X \to 0 \qquad\text{for } h\to 0,
    \end{equation*}
    i.e., $F_\gamma$ is Fréchet differentiable with gradient $\frac1\gamma(x-x^*)  =  (\partial F)_\gamma(x)$.
\end{proof}
Since $F_\gamma$ is convex by \cref{lem:moreau:convex}, this result together with \cref{thm:convex:gateaux} yields the catchy relation $\partial(F_\gamma) = (\partial F)_\gamma$.

\begin{example}\label{ex:moreau}
    We consider again $X=\R^N$.
    \begin{enumerate}
        \item\label{it:moreau:norm-squared} For $F(x) = \frac12\norm{x}_2^2$, \cref{ex:proximal:rn}\,\ref{ex:proximal:rn:ii} yields $\prox_{\gamma F}(x) = \frac{1}{1+\gamma} x$. Inserting this into the definition
            of the Yosida approximation and the Moreau envelope and simplifying yields that
            \begin{equation*}
                (\partial F)_\gamma (x) = \frac1\gamma\left(x-\frac1{1+\gamma}\right)x = \frac1{1+\gamma} x
            \end{equation*}
            and
            \begin{equation*}
                F_\gamma (x) = \frac1{2\gamma} \norm{\tfrac1{1+\gamma}x}_2^2 + \frac12 \norm{\tfrac{1}{1+\gamma}x}_2^2 = \frac1{2(1+\gamma)} \norm{x}_2^2.
            \end{equation*}
            (Unsurprisingly, the Moreau envelope of a quadratic function remains quadratic and is simply scaled.)
        \item\label{it:moreau:abs} For $F(x)=\|x\|_1$, we have from \cref{ex:proximal:rn}\,\ref{ex:proximal:rn:ii} that the proximal point mapping is given by the componentwise soft-shrinkage operator. Inserting this into the definition yields that
            \begin{equation*}
                \left[(\partial \|\cdot\|_1)_\gamma(x)\right]_i =
                \begin{cases}
                    \frac1\gamma(x_i-(x_i-\gamma)) = 1 & \text{if }x_i>\gamma,\\
                    \frac1\gamma x_i & \text{if }x_i\in[-\gamma,\gamma],\\
                    \frac1\gamma(x_i-(x_i+\gamma)) = -1& \text{if }x_i<-\gamma.
                \end{cases}
            \end{equation*}
            Comparing this to the corresponding subdifferential \eqref{eq:convex:subdiff_abs}, we see that the set-valued case at the point $x_i=0$ has been replaced by a linear function on a small interval.

            Similarly, inserting the definition of the proximal point into \eqref{eq:proximal:moreau-decomposition} shows that
            \begin{equation*}
                F_\gamma(x) = \sum_{i=1}^N f_\gamma(x_i)
                \, \text{ for }\,
                f_\gamma(t) =
                \begin{cases}
                    \frac1{2\gamma}|t-(t-\gamma)|^2 + |t-\gamma| = t - \frac\gamma2 & \text{if }t>\gamma,\\
                    \frac1{2\gamma}|t|^2 & \text{if }t \in [-\gamma,\gamma],\\
                    \frac1{2\gamma}|t-(t+\gamma)|^2 + |t+\gamma| = -t - \frac\gamma2 &\text{if } t<-\gamma.
                \end{cases}
            \end{equation*}
            For small values, the absolute value is thus replaced by a quadratic function (which removes the nondifferentiability at $0$). This modification is well-known under the name \term[norm!Huber]{Huber norm}; see \cref{fig:smoothness:moreau-yosida:abs}.

        \item\label{it:moreau:indicator} For $F(x)=\delta_{\B_\infty}(x)$, we have from \cref{ex:proximal:rn}\,\ref{ex:proximal:rn:iii} that the proximal mapping is given by the componentwise projection onto $[-1,1]$ and hence that
            \begin{equation*}
                \left[(\partial \delta_{\B_\infty})_\gamma(x)\right]_i = \frac1\gamma\Big(x_i-\big(x_i - (x_i-1)^+ - (x_i+1)^-\big)\Big) = \frac1\gamma (x_i-1)^+ + \frac1\gamma (x_i+1)^-.
            \end{equation*}
            Similarly, inserting this and using that $\prox_{\gamma F}(x)\in \B_\infty$ and $\iprod{(x+1)^-}{(x-1)^+}_X=0$ yields that
            \begin{equation*}
                (\delta_{\B_\infty})_\gamma(x) = \frac{1}{2\gamma}\norm{(x-1)^+}_2^2 +  \frac{1}{2\gamma}\norm{(x+1)^-}_2^2,
            \end{equation*}
            which corresponds to the classical penalty functional for the inequality constraints $x-1\leq 0$ and $x+1\geq 0$ in nonlinear optimization; see \cref{fig:smoothness:moreau-yosida:indicator}.
    \end{enumerate}
\end{example}
\begin{figure}[t]
    \centering
    \begin{subfigure}[t]{0.45\textwidth}\centering
        \begin{asy}
            unitsize(70,70);

            real f(real x){ return abs(x); };
            real gamma=0.2;
            real q(real x, real base){ return (x-base)^2/(2*gamma); };
            path g=graph(f, -1.5, 1.5);
            draw(g, primalline+linewidth(1.5)+opacity(0.75));

            typedef real realop(real);

            realop quadratic(real base) {
                real g(real x){
                    return f(base)+q(x, base);
                }
                return g;
            }

            realop moreau(real base) {
                real g(real x){
                    return f(x)+q(x, base);
                }
                return g;
            }

            real prox(real x) {
                if(x>gamma){
                    return x-gamma;
                }else if(x<-gamma){
                    return x+gamma;
                }else{
                    return 0;
                }
            }

            real regval(real x) {
                if(x>gamma){
                    return x-gamma/2;
                }else if(x<-gamma){
                    return -x-gamma/2;
                }else{
                    return x^2/(2*gamma);
                }
            }

            void drawat(real base){
                real iv=2*gamma;
                real proxbase=prox(base);
                real val=regval(base);

                path g1q=graph(quadratic(base), base-iv, base+iv);
                draw(g1q, defaultpen+dotted);

                path g1m=graph(moreau(base), proxbase-iv, proxbase+iv);
                draw(g1m, defaultpen+dashed);

                pair minat=(base, val);
                dot(minat);
                draw((proxbase, val)--minat, defaultpen+linewidth(0.5)+opacity(0.75));
                draw((base, f(base))--minat, defaultpen+linewidth(0.5)+opacity(0.75));
            }

            drawat(1);
            drawat(0);

            path mor=graph(regval, -1.5, 1.5);
            draw(mor, dualline+linewidth(2)+opacity(0.75));
        \end{asy}
        \caption{$f(t) = \abs{t}$}
        \label{fig:smoothness:moreau-yosida:abs}
    \end{subfigure}
    \begin{subfigure}[t]{0.45\textwidth}\centering
        \begin{asy}
            unitsize(70,70);

            real f(real x){
                if(x<=-1 || x>=1){
                    return infinity;
                }else{
                    return 0;
                }
            }
            real gamma=0.2;
            real q(real x, real base){ return (x-base)^2/(2*gamma); };
            draw((-1,1.5)--(-1,0)--(1,0)--(1,1.5), primalline+linewidth(1.5)+opacity(0.75));

            typedef real realop(real);

            realop quadratic(real base) {
                real g(real x){
                    return q(x, base);
                }
                return g;
            }

            realop moreau(real base) {
                real g(real x){
                    return q(x, base);
                }
                return g;
            }

            real prox(real x) {
                if(x>1){
                    return 1;
                }else if(x<-1){
                    return -1;
                }else{
                    return x;
                }
            }

            real regval(real x) {
                return (max(0, x-1)^2+min(0, x+1)^2)/(2*gamma);
            }

            void drawat(real base){
                real iv=2*gamma;
                real proxbase=prox(base);
                real val=regval(base);

                path g1q=graph(quadratic(base), base-iv, base+iv);
                draw(g1q, defaultpen+dotted);

                path g1m=graph(moreau(base), max(-1, proxbase-iv), min(1, proxbase+iv));
                draw(g1m, defaultpen+dashed);

                pair minat=(base, val);
                dot(minat);
                draw((proxbase, val)--minat, defaultpen+linewidth(0.5)+opacity(0.75));
                draw((base, f(base))--minat, defaultpen+linewidth(0.5)+opacity(0.75));
            }

            drawat(0.99);
            drawat(0);
            drawat(-0.75);

            path mor=graph(regval, -1.5, 1.5);
            draw(mor, dualline+linewidth(1.5)+opacity(0.75));
        \end{asy}
        \caption{$f(t) = \delta_{[-1,1]}(t)$}
        \label{fig:smoothness:moreau-yosida:indicator}
    \end{subfigure}
    \caption{Illustration of the Moreau--Yosida regularization (solid green line) of $F$ (solid blue line).
        The dotted line indicates the quadratic function $z \mapsto \frac{1}{2\gamma}\norm{x-z}_X^2$, while the dashed line is $z \mapsto F(z)+ \frac{1}{2\gamma}\norm{x-z}_X^2$.
        The dots and the horizontal and vertical lines (nontrivial only in the second point of (\subref{fig:smoothness:moreau-yosida:abs})) emanating from the dots indicate the pair $(x, F_\gamma(x))$ and how it relates to the minimization of the shifted quadratic functional.
        (In (\subref{fig:smoothness:moreau-yosida:indicator}) the two lines are overlaid within $[-1,1]$, as only the domain of definition of the two functions is different.)
    }
    \label{fig:smoothness:moreau-yosida}
\end{figure}
Using \cref{lem:lebesgue:proximal}, analogous characterizations can be derived for the squared $L^2$-norm, the $L^1$-norm, and the indicator functional of the $L^\infty$-ball on $L^2(\Omega)$.

\begin{remark}
    \label{rem:monotone:banach-moreau-yosida}
    Continuing from \cref{rem:monotone:banach}, Moreau--Yosida regularization can also be defined in reflexive Banach spaces; we refer to \cite{BrezisCrandallPazy} for details. Again, the main issue is the practical evaluation of $F_\gamma$ and $(\partial F)_\gamma$ if the duality mapping is no longer the identity.
\end{remark}

\bigskip

By \cref{thm:moreau-yosida}, $F_\gamma$ is Fréchet differentiable with Lipschitz continuous gradient with factor $\gamma^{-1}$. From \cref{thm:smoothness:dual}, we thus know that $F_\gamma^*$ is strongly convex with factor $\gamma$, which in Hilbert spaces is equivalent to $F_\gamma^* - \frac\gamma2\norm{\cdot}_X^2$ being convex. In fact, this can be made even more explicit.

\begin{theorem}\label{thm:moreau:conjugate}
    Let $F:X\to\Rbar$ be proper, convex, and lower semicontinuous. Then we have for all $\gamma>0$ that
    \begin{equation*}
        (F_\gamma)^* = F^* + \frac{\gamma}{2} \norm{\cdot}_X^2.
    \end{equation*}
\end{theorem}
\begin{proof}
    We obtain directly from the definition of the Fenchel conjugate in Hilbert spaces and of the Moreau envelope that
    \begin{equation*}
        \begin{aligned}
            (F_\gamma)^*(x^*) &= \sup_{x\in X} \left\{\iprod{x^*}{x}_X - \inf_{z\in X} \left[\tfrac{1}{2\gamma}\norm{x-z}_X^2 + F(z)\right]\right\}\\
            &= \sup_{x\in X} \left\{\iprod{x^*}{x}_X + \sup_{z\in X} \left\{-\tfrac{1}{2\gamma}\norm{x-z}_X^2 - F(z)\right\}\right\}\\
            &= \sup_{z\in X} \left\{\iprod{x^*}{z}_X - F(z)+\sup_{x\in X} \left\{\iprod{x^*}{x-z}_X - \tfrac{1}{2\gamma}\norm{x-z}_X^2\right\}\right\}\\
            &= F^*(x^*)  + \left(\tfrac{1}{2\gamma}\norm{\cdot}_X^2\right)^*(x^*),
        \end{aligned}
    \end{equation*}
    since for any given $z\in X$, the inner supremum is always taken over the full space $X$.
    The claim now follows from \cref{lem:convex:power-conjugate} with $p=2$ (using again the fact that we have identified $X^*$ with $X$) and \cref{lem:convex:fenchel_calc}\,\ref{lem:convex:fenchel_calc:i}.
\end{proof}

From this, we obtain the following order relation for the Moreau envelope; cf. \cref{fig:smoothness:moreau-yosida}.
\begin{corollary}\label{lem:moreau:order}
    Let $F:X\to\Rbar$ be proper, convex, and lower semicontinuous. Then for all $\gamma_1\geq \gamma_2\geq 0$,
    \begin{equation*}
        F_{\gamma_1} \leq F_{\gamma_2} \leq F.
    \end{equation*}
\end{corollary}
\begin{proof}
    By \cref{thm:moreau:conjugate}, we have for all $x\in X$ that
    \begin{equation*}
        (F_{\gamma_1})^*(x) \geq (F_{\gamma_2})^*(x) \geq F^*(x).
    \end{equation*}
    Since $F_\gamma$ is convex by \cref{lem:moreau:convex} and smooth and hence \emph{a fortiori} proper and lower semicontinuous by \cref{thm:moreau-yosida}, we can combine this result with \cref{lem:convex:fenchel-monotone,thm:convex:moreau} to deduce that for all $x\in X$,
    \begin{equation*}
        F_{\gamma_1}(x) = (F_{\gamma_1})^{**}(x) \leq (F_{\gamma_2})^{**}(x) \leq F^{**} (x) = F(x).
        \qedhere
    \end{equation*}
\end{proof}

\Cref{thm:moreau:conjugate} also yields a Moreau decomposition of the envelope; cf.~\cref{lem:proximal:calculus}\,\ref{lem:proximal:calculus:ii}.
\begin{corollary}
    Let $F:X\to\Rbar$ be proper, convex, and lower semicontinuous. Then for all $x\in X$ and $\gamma>0$,
    \begin{equation*}
        \frac{1}{2\gamma}\norm{x}_X^2 = F_\gamma(x) + (F^*)_{\gamma^{-1}}(\gamma^{-1}x).
    \end{equation*}
\end{corollary}
\begin{proof}
    By definition of the Moreau envelope, we have that
    \begin{equation*}
        \begin{aligned}[b]
            F_\gamma(x) &= \inf_{z\in X} F(z)  + \frac{1}{2\gamma}\norm{x-z}_X^2 \\
            &= \inf_{z\in X} F(z)  + \frac{1}{2\gamma}\norm{x}_X^2 - \frac1\gamma\iprod{x}{z}_X + \frac1{2\gamma}\norm{z}_X^2\\
            &= \frac{1}{2\gamma}\norm{x}_X^2 - \sup_{z\in X} \left\{ \iprod{\gamma^{-1}x}{z}_X - F(z) - \frac{1}{2\gamma}\norm{z}_X^2 \right\}\\
            &= \frac{1}{2\gamma}\norm{x}_X^2 - \left(F + \frac1{2\gamma}\norm{\cdot}_X^2\right)^*(\gamma^{-1}x).
        \end{aligned}
    \end{equation*}
    The claim now follows since $F$ (by assumption) and $F_\gamma$ (by \cref{lem:moreau:convex,thm:moreau-yosida}) are convex and lower semicontinuous, and hence \cref{thm:convex:moreau}\,\ref{thm:convex:moreau:iii} together with \cref{thm:moreau:conjugate} implies
    \begin{equation*}
        \left(F + \frac1{2\gamma}\norm{\cdot}_X^2\right)^* =
        \left(F^{**} + \frac1{2\gamma}\norm{\cdot}_X^2\right)^* =
        \left((F^*)_{\gamma^{-1}}\right)^{**} =
        (F^*)_{\gamma^{-1}}.
        \qedhere
    \end{equation*}
\end{proof}
Taking the derivative of this identity and using \cref{thm:moreau-yosida} together with \eqref{eq:proximal:yosida}, we again obtain \cref{lem:proximal:calculus}\,\ref{lem:proximal:calculus:ii}.

With the help of \cref{thm:moreau:conjugate}, we can also show the converse of \cref{thm:moreau-yosida}: every smooth function can be obtained through Moreau--Yosida regularization.
\begin{corollary}
    \label{cor:moreau:smooth-are}
    Let $F: X \to \Rbar$ be convex and $L$-smooth.
    Then for all $x\in X$,
    \begin{align*}
        F(x)=(G^*)_{\inv L}(x)
        \quad\text{and}\quad
        \grad F(x)=\prox_{LG}(L x)
    \end{align*}
    for
    \begin{equation*}
        G:X\to \Rbar,\qquad G(x) = F^*(x) - \frac{1}{2L}\norm{x}_X^2.
    \end{equation*}
\end{corollary}
\begin{proof}
    Since $F$ is convex and $L$-smooth and $X$ is a Hilbert space, \cref{lemma:smoothness,thm:smoothness:dual} yield that $F^*$ is strongly convex with factor $\inv L$ and thus that $G$ is convex. Furthermore, as a Fenchel conjugate of a proper convex functional, $F^*$ and thus $G$ are proper and lower semicontinuous.
    \Cref{thm:moreau:conjugate,thm:convex:moreau} now imply that for all $x \in X$,
    \begin{equation*}
        (G^*)_{\inv L}(x) = (G^*)_{\inv L}^{**}(x)
        = \left(G + \frac{1}{2L}\norm{\cdot}_X^2\right)^{*}(x) = F^{**}(x)=F(x).
    \end{equation*}

    Furthermore, by \cref{lem:convex:subdiff_calc,thm:subdiff:sum,thm:convex:gateaux}, we have that
    \begin{equation*}
        \subdiff G(z)=\subdiff F^*(z)-\{\inv L z\} \qquad\text{for all } z \in X.
    \end{equation*}
    By the definition of the proximal mapping, this is equivalent to $z = \prox_{LG}{Lx}$ for any $x\in \subdiff F^*(z)$. But by \cref{lem:convex:fenchel-young}, $x\in \subdiff F^*(z)$ holds if and only if $z\in \partial F(x) = \{\grad F(x)\}$, and combining these two yields the first expression for the gradient.
\end{proof}
\begin{remark}[conversion between gradients and proximal mappings]
    To motivate the relevance of the previous result to optimization, recall that according to \cref{cor:moreau:smooth-are}, solving $\min_x F(x)$ for an $L$-smooth function $F$ is equivalent to solving
    \begin{equation*}
        \min_{x, \alt x \in X} G^*(x) + \frac{1}{2L}\norm{x-\alt x}_X^2.
    \end{equation*}
    Observe that $G^*$ may be nonsmooth. Suppose we apply an algorithm for the latter that makes use of the proximal mapping of $G^*$ (such as the splitting methods that will be discussed in the following chapters).
    Then using the Moreau decomposition of \cref{lem:proximal:calculus}\,\ref{lem:proximal:calculus:ii} with \cref{cor:moreau:smooth-are}, we see that
    \begin{equation*}
        \prox_{\inv L G^*}(x) = x - \inv L \grad F(x).
    \end{equation*}
    Therefore, this can still be done purely in terms of the gradient evaluations of $F$.
\end{remark}

\bigskip

\Cref{thm:moreau-yosida} allows approximating a (nonsmooth) convex optimization problem by a family of \emph{smooth} optimization problems. This begs the question of the relation of minimizers of the smooth approximation $F_\gamma$ to minimizers of the original functional $F$, which we will answer using the variational convergence theory of \cref{sec:variation:convergence}. Recall that a sequence $\{F_n\}_{n\in\N}$ of functionals $F_n:X\to\Rbar$ is said to $\Gamma$-converge weakly to $F$, written $F_n\weakto_\Gamma F$, if
\begin{enumerate}
    \item\label{item:moreau:gamma-liminf} for all sequences $\{x_n\}_{n\in\N}\subset X$ with $x_n\weakto x$,
        \begin{equation*}
            \liminf_{n\to\infty} F_n(x_n) \geq F(x);
        \end{equation*}
    \item\label{item:moreau:gamma-recovery} there exists a recovery sequence $\{x_n\}_{n\in\N}\subset X$ such that $x_n\weakto x$ and
        \begin{equation*}
            \limsup_{n\to\infty} F_n(x_n)\leq F(x).
        \end{equation*}
\end{enumerate}
\begin{theorem}\label{thm:moreau:gamma}
    Let $F:X\to\Rbar$ be proper, convex, lower semicontinuous, and bounded from below. Then $F_\gamma \weakto_\Gamma F$ as $\gamma\to 0$.
\end{theorem}
\begin{proof}
    Let $\{\gamma_n\}_{n\in\N}\subset (0,\infty)$ be an arbitrary null sequence and set $F_n:=F_{\gamma_n}$ for all $n\in \N$. We now verify the two defining properties.

    \emph{\ref{item:moreau:gamma-liminf}:} Let $x\in X$ and $\{x_n\}_{n\in\N} \subset X$ be such that $x_n\weakto x$. Furthermore, let $K\subset \N$ denote the subsequence realizing the $\liminf$ on the left-hand side, i.e.
    \begin{equation*}
        \liminf_{n\to \infty} F_n (x_n)  = \lim_{K\ni n\to \infty} F_n(x_n)< \infty,
    \end{equation*}
    where we can assume the limit to be finite since otherwise \ref{item:moreau:gamma-liminf} holds trivially. Set now $z_n:=\prox_{\gamma_n F}(x_n)\in \dom F$ for all $n\in K$ such that by \eqref{eq:proximal:moreau-decomposition},
    \begin{equation*}
        F_n(x_n) = \frac{1}{2\gamma_n}\norm{z_n - x_n}_X^2 + F(z_n)<\infty \qquad\text{for all } n\in K.
    \end{equation*}
    Since $F$ is bounded from below by some $\bar F\in\R$, this implies that for all $n\in K$,
    \begin{equation*}
        \bar F \leq \frac1{2\gamma_n} \norm{z_n-x_n}_X^2 + \bar F\leq \frac1{2\gamma_n}\norm{z_n-x_n}_X^2 + F(z_n).
    \end{equation*}
    Passing to the limit $K\ni n\to \infty$ and using the fact that $\gamma_n\to 0$ and that the right-hand side converges and hence is bounded along this subsequence yields that $\norm{z_n - x_n}_X\to 0$. This implies that for all $z\in X$,
    \begin{equation*}
        \begin{aligned}
            \left|\iprod{z_n - x}{z}_X\right| &\leq \left|\iprod{z_n - x_n}{z}_X\right| + \left|\iprod{x_n -x}{ z}_X\right| \\
            &\leq \norm{z_n-x_n}_X\norm{z}_X + \left|\iprod{x_n-x}{z}_X\right| \to 0 \quad\text{as }K\ni n\to\infty
        \end{aligned}
    \end{equation*}
    since $x_n\weakto x$ by assumption. Hence $z_n\weakto x$ as $K\ni n\to\infty$ as well. We can now use the lower semicontinuity of $F$ to obtain
    \begin{equation*}
        \begin{aligned}
            \liminf_{n\to \infty} F_n (x_n)  = \lim_{K\ni n\to \infty} F_n(x_n) &= \lim_{K\ni n\to\infty} \frac1{2\gamma_n}\norm{z_n-x_n}_X^2 + F(z_n) \\
            &\geq \liminf_{K\ni n\to\infty} F(z_n) \geq F(x).
        \end{aligned}
    \end{equation*}

    \emph{\ref{item:moreau:gamma-recovery}:} For any $x\in X$, we can simply choose the constant sequence $x_n=x$ for all $n\in \N$ as the recovery sequence, since in this case \cref{lem:moreau:order} shows that
    \begin{equation*}
        F_n(x_n) = F_{\gamma_n}(x_n) \leq F(x_n) = F(x)\qquad\text{for all }n\in\N,
    \end{equation*}
    and taking the $\limsup$ as $n\to\infty$ yields the claim.
\end{proof}
In fact, we can show even more: For convex functionals, the Moreau envelope is exact in minimizers.
\begin{lemma}
    Let $F:X\to\Rbar$ be proper, convex, and lower semicontinuous. Then for every $\gamma>0$, a point $\bar x\in X$ is a minimizer of $F_\gamma$ if and only if $\bar x$ is a minimizer of $F$,
    and in this case $F_\gamma(\bar x) = F(\bar x)$ holds.
\end{lemma}
\begin{proof}
    Since $F$ by assumption and $F_\gamma$ by \cref{lem:moreau:convex,thm:moreau-yosida} are convex and lower semicontinuous, we can combine \cref{lem:convex:fenchel-young,thm:moreau:conjugate,thm:subdiff:sum,thm:convex:gateaux} to deduce that
    \begin{equation*}
        \begin{aligned}
            0\in \partial F_\gamma (\bar x) &\equivalent \bar x \in \partial (F_\gamma)^*(0) = (\partial F^* + \gamma \Id)(0) = \partial F^*(0)\\
            &\equivalent 0\in \partial F(\bar x).
        \end{aligned}
    \end{equation*}
    The first claim now follows from \cref{thm:convex:fermat}.

    For the second claim, let $\bar x$ be a minimizer of $F$. Then we have for any $\gamma>0$ and all $z\in X$ that
    \begin{equation*}
        F(\bar x) \leq F(z) \leq F(z) + \frac1{2\gamma} \norm{z-\bar x}_X^2,
    \end{equation*}
    and taking the infimum over all $z\in X$ yields $F(\bar x)\leq F_\gamma(\bar x)$. Together with $F_\gamma(\bar x)\leq F(\bar x)$ from \cref{lem:moreau:order}, this yields the second claim.
\end{proof}
This is hardly surprising (and does not violate any Law of Conservation of Difficulties) since evaluating the Moreau envelope itself involves minimization of $F$. This approach is therefore more useful when considering composite functionals of the form $F+G$ where $F$ has a proximal mapping that is simple to evaluate and $G$ is weakly continuous (e.g., $G(x) = \norm{Kx-y}_Y^2$ for some compact operator $K:X\to Y$, which is a typical situation in applications; cf. \cref{part:applications}) since $\Gamma$-convergence is stable under continuous perturbations.
\begin{theorem}\label{thm:moreau:convergence}
    Let $F:X\to\Rbar$ be proper, convex, lower semicontinuous, and bounded from below, and let $G:X\to\Rbar$ be weakly continuous and bounded from below. Assume further that either $F$ or $G$ is coercive. Then the family $\{x_\gamma\}_{\gamma>0}$ of minimizers $x_\gamma \in X$ of $F_\gamma + G$ has a weak accumulation point, and every weak accumulation point is a minimizer of $F+G$.
\end{theorem}
\begin{proof}
    Let $\{\gamma_n\}_{n\in\N}\subset (0,\infty)$ be an arbitrary decreasing null sequence and set $F_n:=F_{\gamma_n}$ for all $n\in \N$. Then we have by \cref{thm:moreau:gamma} that $F_n\weakto_\Gamma F$ and hence by \cref{lem:variation:gamma-perturbation} that $F_n+G \weakto_\Gamma F+G$.
    To apply \cref{thm:variation:gamma-convergence}, it thus remains to show that $\{F_n+G\}_{n\in\N}$ is equicoercive. To that end, we first note that by \cref{lem:moreau:order}, $F_n+G\geq F_1+G$ and hence it suffices to show that $F_\gamma+G$ is coercive for an arbitrary $\gamma>0$.

    Let therefore $\{x_n\}_{n\in\N}\subset X$ be such that $\norm{x}_X\to \infty$ as $n\to\infty$, and set $z_n :=\prox_{\gamma F}{x_n}$ such that
    \begin{equation*}
        F_\gamma(x_n) = \frac1{2\gamma_1}\norm{z_n-x_n}_X^2 + F(z_n) \qquad\text{for all }n\in \N.
    \end{equation*}
    We now distinguish the two cases:
    \begin{enumerate}
        \item \emph{$G$ is coercive:} Then we have immediately that
            \begin{equation*}
                F_\gamma(x_n) + G(x_n) \geq F(z_n) + G(x_n) \geq \bar F + G(x_n) \to \infty
            \end{equation*}
            using the assumed lower bound $\bar F$ of $F$.
        \item \emph{$F$ is coercive:} Then we have that
            \begin{equation*}
                F_\gamma(x_n) + G(x_n) \geq \frac1{2\gamma}\norm{z_n-x_n}_X^2 + F(z_n) + \bar G
            \end{equation*}
            using the assumed lower bound $\bar G$ of $G$. We now distinguish further: If $\{F(z_n)\}_{n\in\N}$ is bounded, we can deduce from the coercivity of $F$ that$ \{z_n\}_{n\in\N}$ is bounded as well and hence that
            \begin{equation*}
                F_\gamma(x_n) + G(x_n) \geq \frac1{2\gamma}(\norm{x_n}_X-\norm{z_n}_X)^2 + \bar F + \bar G \to \infty.
            \end{equation*}
            Otherwise we have $F(z_n)\to\infty$ and thus can conclude as in the first case that
            \begin{equation*}
                F_\gamma(x_n) + G(x_n) \geq F(z_n) + \bar G \to \infty.
            \end{equation*}
    \end{enumerate}

    Finally, \cref{thm:moreau-yosida} and the assumptions on $G$ imply that $F_\gamma+G$ is proper, weakly lower semicontinuous, and coercive for every $\gamma>0$, and hence \cref{thm:variation:existence} ensures the existence of the claimed family $\{x_\gamma\}_{\gamma>0}$ of minimizers.
\end{proof}
An alternative approach is \term[regularization!Moreau--Yosida!dual]{dual Moreau--Yosida regularization}:
For proper, convex, and lower semicontinuous functionals $F,G:X\to\Rbar$, every minimizer $\bar x\in X$ of $F+G$ satisfies the Fermat principle $0\in \partial (F+G)(\bar x)$, which using the sum rule (\cref{thm:subdiff:sum}, under the standard regularity condition) is equivalent to the existence of a $\bar p\in X$ satisfying
\begin{equation*}
    \left\{
        \begin{aligned}
            -\bar p &\in \partial F(\bar x),\\
            \bar p &\in \partial G(\bar x),
        \end{aligned}
    \right.
\end{equation*}
where the last relation can by \cref{lem:convex:fenchel-young} equivalently be written as $\bar x\in \partial G^*(\bar p)$.
Assuming we can easily deal with $F$ (e.g., because $F$ is differentiable), we can replace $\partial G^*$ with its Yosida approximation $(\partial G^*)_\gamma$ to obtain the regularized relation
\begin{equation*}
    x_\gamma = (\partial G^*)_\gamma (p_\gamma) = \frac1\gamma (p_\gamma -\prox_{\gamma G^*}(p_\gamma))
\end{equation*}
for $-p_\gamma\in \partial F(x_\gamma)$.
This is now an \emph{explicit} and even Lipschitz continuous relation.
Although $x_\gamma$ is no longer a minimizer of $F+G$, the convexity of $G_\gamma$ together with \cref{thm:moreau:conjugate} implies that $x_\gamma = (\partial G^*)_\gamma (p_\gamma) = \partial(G^*_\gamma)(p_\gamma)$ is equivalent to
\begin{equation*}
    p_\gamma\in \partial (G^*_\gamma)^*(x_\gamma) = \partial\left(G^{**} + \tfrac{\gamma}2\norm{\cdot}_X^2\right)(x_\gamma) = \partial\left(G + \tfrac{\gamma}2\norm{\cdot}_X^2\right)(x_\gamma),
\end{equation*}
i.e., $x_\gamma$ is the (unique due to the strict convexity of the squared norm) minimizer of the functional $F+G+\tfrac{\gamma}2\norm{\cdot}_X^2$.
Hence, the regularization of $\partial G^*$ has not made the original problem smooth but merely (more) strongly convex.
The equivalence can also be used to show (similarly to the proof of \cref{thm:variation:existence}) that $x_\gamma \weakto \bar x$ for $\gamma\to 0$.

We will see applications of both approaches in \cref{chap:control,chap:discretecontrol}.

\index{algorithm|see{method}}
\index{iteration|see{method}}
\chapter{Proximal point and splitting methods}\label{chap:proximal}

We now turn to the development of algorithms for computing minimizers of functionals $J:X\to\Rbar$ of the form
\begin{equation*}
    J(x) \defeq F(x) + G(x)
\end{equation*}
for $F,G:X\to\Rbar$ convex but not necessarily differentiable.
A natural approach (in Hilbert spaces) would be to take the steepest descent method and simply replace the gradient with an arbitrary subgradient, which leads to the \term[method!subgradient]{subgradient method}
\begin{algeqbox*}
    \begin{equation*}
        x^{k+1} = x^k - \tau_k \xi^k, \qquad \xi^k\in \partial J(x^k).
    \end{equation*}
\end{algeqbox*}
However, this iteration will in general not converge: Even in finite dimensions, arbitrary subgradients need not be descent directions -- this can only be guaranteed for the subgradient of minimal norm, and the minimal norm subgradient of $J$ cannot be computed easily from those of $F$ and $G$. Hence convergence can only be shown for finite-valued $J$ under rather strong a priori assumptions on the choice of $\xi^k$ and $\tau_k$; we refer to \cite[Chapter 7.1]{Ruszczynski:2006a} for a detailed treatment in finite-dimensional vector spaces and \cite{Alber:1998} for a convergence result in Hilbert spaces.

We thus follow a different approach and look for a root $\realoptx$ of the set-valued mapping $x\mapsto \partial J(x)$ (which coincides with the minimizer $\bar x$ of $J$ if $J$ is convex). In this chapter, we only derive methods, postponing proofs of convergence -- in various different senses -- to \cref{chap:convergence,chap:testing,chap:gap}.
For the reasons mentioned in the beginning of \cref{sec:monotone:resolvents}, we will assume in this and the following chapters that $X$ (as well as all further occurring spaces) is a Hilbert space so that we can identify $X^*\cong X$.

\section{Proximal point method}
\label{sec:proximal}

We have seen in \cref{thm:proximal:fermat} that a root $\realoptx$ of $\partial J:X\setto X$ can be characterized as a fixed point of $\prox_{\tau J}$ for any $\tau>0$. This suggests a fixed-point iteration: Choose $x^0\in X$ and for an appropriate sequence $\{\tau_k\}_{k\in \N}$ of step sizes set
\begin{algeqbox}
    \begin{equation}\label{eq:proximal:ppa}
        x^{k+1} \defeq \prox_{\tau_k J}(x^k).
    \end{equation}
\end{algeqbox}
This iteration naturally generalizes to finding a root $\realoptx \in \inv A(0)$ of a set-valued (usually monotone) operator $A: X \setto X$ as
\begin{algeqbox}
    \begin{equation}
        \nextx \defeq \calR_{\tau_k A}(\thisx).
    \end{equation}
\end{algeqbox}
This is the \term[method!proximal point]{proximal point method}, which is the basic building block for all methods in this chapter.
Using the definition of the resolvent, this can also be written in implicit form as
\begin{equation}
    \label{eq:proximal:ppa-monotone}
    0 \in \tau_k A(\nextx) + (\nextx-\thisx),
\end{equation}
which will be useful for the analysis of the method.

If $A$ is maximally monotone (in particular if $A=\subdiff J$), \cref{lem:proximal:firmly-nonexpansive} shows that the iteration mapping $x \mapsto \calR_{\tau_k A}(x)$ is firmly nonexpansive. Mere (nonfirm) nonexpansivity already implies that
\begin{equation*}
    \norm{\nextx-\realoptx}_X = \norm{\calR_{\tau_k A}(\thisx)-\realoptx}_X \le \norm{\thisx-\realoptx}_X.
\end{equation*}
In other words, the method does not escape from a fixed point. Either a more refined analysis based on firm nonexpansivity of the iteration mapping or a more direct analysis based on the maximal monotonicity of $A$ can be used to further show that the iterates $\{\thisx\}_{k \in \N}$ indeed converge to a fixed point $\realoptx$ for an initial iterate $x^0$. The latter will be the topic of \cref{chap:convergence}.

A practical issue is that the steps \eqref{eq:proximal:ppa} of the basic proximal point method are typically just as difficult as the original problem, so the method is not feasible for problems that demand an iterative method for their solution in the first place. However, the proximal step does form an important building block of several more practical \emph{splitting} methods for problems of the form $J=F+G$, which we derive in the following by additional clever manipulations.

\begin{remark}
    The proximal point algorithm can be traced back to
    Krasnosel'ski\weaki{}
    \cite{krasnoselskii1955remarks} and Mann \cite{Mann1953mean} (as a special case of the \term[method!Krasnoselskii--Mann]{Krasnoselskii--Mann iteration}); it was also studied in \cite{martinet1970regularisation}.
    The formulation considered here was proposed in \cite{rockafellar1976proximal}.
\end{remark}

\section{Explicit splitting: forward-backward splitting}
\label{sec:splitting:explicit}

As we have noted, the proximal point method is not feasible for most functionals of the form  $J(x) = F(x)+G(x)$, since the evaluation of $\prox_J$ is not significantly easier than solving the original minimization problem -- even if $\prox_F$ and $\prox_G$ have a closed-form expression. (Such functionals are called \term[functional!prox-simple]{prox-simple}).
We thus proceed differently: instead of applying the proximal point reformulation directly to $0\in\partial J(\realoptx)$, we first apply the subdifferential sum rule (\cref{thm:subdiff:sum}) to deduce the existence of $\realopt{p}\in X$ with
\begin{equation}\label{eq:splitting:optsys}
    \left\{\begin{aligned}
            \realopt{p} &\in \partial F(\realoptx),\\
            -\realopt{p} &\in \partial G(\realoptx).
    \end{aligned}\right.
\end{equation}
We can now replace one or both of these subdifferential inclusions by a proximal point reformulation that only involves $F$ or $G$.

Explicit splitting methods\index{method!explicit splitting} -- also known as \term[method!forward-backward splitting]{forward-backward splitting} -- are based on applying \cref{lem:proximal:subdiff} only to, e.g., the second inclusion in \eqref{eq:splitting:optsys} to obtain
\begin{equation}\label{eq:splitting:optsys_exp}
    \left\{\begin{aligned}
            \realopt{p} &\in \partial F(\realoptx),\\
            \realoptx &= \prox_{\tau G}(\realoptx -\tau \realopt{p}).
    \end{aligned}\right.
\end{equation}
The corresponding fixed-point iteration then consists in
\begin{enumerate}[label=\arabic*.]
    \item choosing $p^k \in \partial F(x^k)$ (with minimal norm);
    \item setting $x^{k+1} =  \prox_{\tau_k G}(x^k -\tau_k p^k)$.
\end{enumerate}
Again, computing a subgradient with minimal norm can be complicated in general. It is, however, easy if $F$ is additionally differentiable since in this case $\partial F(x) = \{\nabla F(x)\}$ by \cref{thm:convex:gateaux}. This leads to the \term[method!proximal gradient]{proximal gradient} or \term[method!forward-backward splitting]{forward-backward splitting method}
\begin{algeqbox}
    \begin{equation}\label{eq:splitting:fb}
        x^{k+1} \defeq \prox_{\tau_k G}(x^k - \tau_k \nabla F(x^k)).
    \end{equation}
\end{algeqbox}
(The special case $G=\delta_C$ -- i.e., $\prox_{\tau G}(x) = \proj_C(x)$ -- is also known as the \term[method!projected gradient]{projected gradient method}).
Similarly to the proximal point method, this method can be written in implicit form as
\begin{equation}
    \label{eq:splitting:fb:implicit}
    0 \in \tau_k[\subdiff G(\nextx) + \grad F(\thisx)] + (\nextx-\thisx).
\end{equation}

Based on this, we will see in \cref{chap:convergence} that the iterates $\{\thisx\}_{k\in\N}$ converge weakly if $\tau_k L < 2$ for $L$ the Lipschitz factor of $\grad F$.
However, the need to know $L$ is one drawback of the explicit splitting method, which can to some extent be circumvented by performing a \emph{line search}: testing for various choices of $\tau_k$ until a sufficient decrease in function values is achieved. We will discuss such strategies later on in \cref{sec:meta:linesearch}.
Another highly successful variant of explicit splitting applies \term{inertia} to the iterates for faster convergence; this we will discuss in \cref{sec:meta:inertia} after developing tools for the study of convergence rates.

\begin{remark}
    Forward-backward splitting for finding the root of the sum of two monotone operators was already proposed in  \cite{lionsmercier1979splitting}. It has become especially popular under the name \term[method!iterative soft-thresholding]{iterative soft-thresholding} (ISTA\index{ISTA|see{method, iterative soft-thresholding}}) in the context of \term[problem!sparse regression]{sparse regression} (i.e., regularization of linear inverse problems with $\ell^1$ penalties), see, e.g., \cite{chambolledevore1998nonlinear,daubechiesdefriesdemol2004ista,wrightnovak2009sparse}.
\end{remark}

\section{Implicit splitting: Douglas--Rachford splitting}

Even with a line search, the restriction on the step sizes $\tau_k$ in explicit splitting remains unsatisfactory. Such restrictions are not needed in implicit splitting methods. (Compare the properties of explicit vs.~implicit Euler methods for differential equations.) Here, the proximal point formulation is applied to both subdifferential inclusions in \eqref{eq:splitting:optsys}, which yields the optimality conditions
\begin{equation*}\label{eq:splitting:optsys_imp}
    \left\{\begin{aligned}
            \realoptx &= \prox_{\tau F}(\realoptx + \tau \realopt{p}),\\
            \realoptx &= \prox_{\tau G}(\realoptx - \tau \realopt{p}).
    \end{aligned}\right.
\end{equation*}
To eliminate $\realopt{p}$ from these equations, we set $\realoptz  \defeq \realoptx+\tau \realopt{p}$ and $\realoptw \defeq \realoptx - \tau \realopt{p}=2\realoptx-\realoptz$.
It remains to derive a recursion for $\realoptz$, which we obtain from the productive zero $\realoptz = \realoptz + (\realoptx - \realoptx)$.
Further replacing some copies of $\realoptx$ by a new variable $\realopty$ leads to the overall fixed point system
\begin{equation*}
    \left\{\begin{aligned}
            \realoptx &= \prox_{\tau F}(\realoptz),\\
            \realopty &= \prox_{\tau G}(2\realoptx - \realoptz),\\
            \realoptz &= \realoptz + \realoptx - \realopty.
    \end{aligned}\right.
\end{equation*}
The corresponding fixed-point iteration leads to the \term[method!Douglas--Rachford splitting]{Douglas--Rachford splitting} (DRS\index{DRS|see{method, Douglas--Rachford splitting}}) method
\begin{algeqbox}
    \begin{equation}\label{eq:splitting:dr}
        \left\{\begin{aligned}
                x^{k+1} &\defeq \prox_{\tau F}(z^k),\\
                y^{k+1} &\defeq \prox_{\tau G}(2x^{k+1} - z^k),\\
                z^{k+1} &\defeq z^k + y^{k+1} - x^{k+1}.
        \end{aligned}\right.
    \end{equation}
\end{algeqbox}
Of course, the algorithm and its derivation generalize to arbitrary monotone operators $A, B: X \setto X$:
\begin{equation}\label{eq:splitting:dr-monotone}
    \left\{\begin{aligned}
            x^{k+1} &= \calR_{\tau B}(z^k),\\
            y^{k+1} &= \calR_{\tau A}(2x^{k+1} - z^k),\\
            z^{k+1} &= z^k + y^{k+1} - x^{k+1}.
    \end{aligned}\right.
\end{equation}

We can also write the DRS method in more implicit form. Indeed, inverting the resolvents in \eqref{eq:splitting:dr-monotone} and using the last update to change variables in the first two yields
\begin{equation*}
    \left\{\begin{aligned}
            0 & \in \tau B(\nextx) + \nexty - \nextz,\\
            0 & \in \tau A(\nexty) + \nextz - \nextx,\\
            0 & =  \nextx - \nexty + (\nextz - \thisz).
    \end{aligned}\right.
\end{equation*}
Therefore, with $u \defeq (x, y, z) \in X^3$, and the operators\footnote{Here and in the following, we identify $x\in X$ with the singleton set $\{x\}\subset X$ whenever there is no danger of confusion.}
\begin{equation}
    \label{eq:proximal:drs-form}
    \Hany(x, y, z) \defeq
    \begin{pmatrix}
        \tau B(x) + y - z \\
        \tau A(y) + z - x \\
        x - y
    \end{pmatrix}
    \quad\text{and}\quad
    \Precond \defeq
    \begin{pmatrix}
        0 & 0 & 0 \\
        0 & 0 & 0 \\
        0 & 0 & \Id
    \end{pmatrix},
\end{equation}
we can write the DRS method as the \term[method!proximal point!preconditioned]{preconditioned proximal point method}
\begin{equation}
    \label{eq:proximal:preconditioned}
    0 \in \Hany(\nextu) + \Precond(\nextu-\thisu).
\end{equation}
Indeed, the basic proximal point in implicit form \eqref{eq:proximal:ppa-monotone} is just \eqref{eq:proximal:preconditioned} with the \term{preconditioner} $\Precond=\inv\tau\Id$.
It is furthermore straightforward to verify that $0 \in \Hany(\realoptu)$ is equivalent to $0 \in A(\realoptx) + B(\realoptx)$.

The formulation \eqref{eq:proximal:preconditioned} will in the following chapter form the basis for proving the convergence of the method.
Recalling the discussion on convergence in \cref{sec:proximal}, it seems beneficial for $\Hany$ to be maximally monotone, as then (although this is not immediate from \cref{lem:proximal:firmly-nonexpansive}) it is reasonable to expect
the nonexpansivity of the iterates with respect to the semi-norm $u \mapsto \norm{u}_M \defeq \sqrt{\iprod{Mu}{u}}$ on $X^3$ induced by the self-adjoint operator $M$, i.e., that
\[
    \norm{\nextu-\realoptu}_M \le \norm{\thisu-\realoptu}_M.
\]
While it is straightforward to verify that $\Hany$ is monotone if $A$ and $B$ are, the question of maximal monotonicity is more involved and will be addressed in \cref{chap:convergence}. There, we will also show that the expected nonexpansivity holds in a slightly stronger sense and that this will yield the convergence of the method.

\begin{remark}
    The Douglas--Rachford splitting was first introduced in \cite{douglas1956numerical};
    the relationship to the proximal point method was discovered in \cite{Eckstein:1992}.
    The DRS is the \emph{unique} 2-operator splitting method that needs to propagate only one variable from each iteration to the next one, $\nexxt z$ \cite{Ruy:2019}.
    An extension of the DRS with a forward step with respect to a third operator is studied \cite{Davis:2017}.
    It is also possible to devise acceleration schemes under strong monotonicity \cite[see, e.g.,][]{bredies2016accelerated}.
\end{remark}

\section{Primal-dual proximal splitting}\label{sec:proximal:pd}

We now consider problems of the form
\begin{equation}
    \label{eq:pdps-problem}
    \min_{x\in X} F(x) + G(Kx)
\end{equation}
for $F:X\to\Rbar$ and $G:Y\to\Rbar$ proper, convex, and lower semicontinuous, and $K\in \linear(X;Y)$. Applying \cref{thm:convex:fenchel,lem:convex:fenchel-young} to such a problem yields the Fenchel extremality conditions
\begin{equation}\label{eq:splitting:fenchel}
    \left\{    \begin{aligned}
            -K^*\bar{y} &\in \partial F(\bar x),\\
            \bar{y} &\in \partial G(K\bar x) ,
    \end{aligned}\right.
    \quad\equivalent \quad
    \left\{    \begin{aligned}
            -K^*\bar{y} &\in \partial F(\bar x),\\
            K\bar x &\in \partial G^*(\bar y).
    \end{aligned}\right.
\end{equation}
With the general notation $u \defeq (x, y)$, this can be written as $0 \in \Hsaddle(\realoptu)$
for
\begin{equation}
    \label{eq:saddle-operator}
    \Hsaddle(u) \defeq \begin{pmatrix} \partial F(x) + K^*y \\ \partial G^*(y) - Kx \end{pmatrix}.
\end{equation}

It is again not difficult to see that $H$ is monotone.
This suggests that we might be able to apply the proximal point method to find a root of $\Hsaddle$. In practice we however need to work a little bit more, as the resolvent of $H$ can rarely be given an explicit, easily solvable form. If, however, the resolvents of $G^*$ and $F$ can individually be computed explicitly, it makes sense to try to \enquote{decouple} the primal and dual variables. This is what we will do.

To do so, we reformulate for arbitrary $\sigma,\tau>0$ the extremality conditions \cref{eq:splitting:fenchel} using \cref{lem:proximal:subdiff} as
\begin{equation*}
    \left\{\begin{aligned}
            \realoptx &= \prox_{\tau F}(\realoptx - \tau K^*\realopty),\\
            \realopty &= \prox_{\sigma G^*}(\realopty + \sigma K\realoptx).
    \end{aligned}        \right.
\end{equation*}
This suggests the fixed-point iterations
\begin{equation}\label{eq:splitting:pdhg}
    \left\{\begin{aligned}
            x^{k+1} &= \prox_{\tau F}(x^k - \tau K^*y^k),\\
            y^{k+1} &= \prox_{\sigma G^*}(y^{k} + \sigma K x^{k+1}).
    \end{aligned}        \right.
\end{equation}
In the first equation, we now use $\prox_{\tau F} = (\Id +\tau\partial F)^{-1}$ to obtain that
\begin{equation}
    \label{eq:splitting:pdhg-first-step-rewrite}
    \begin{aligned}[t]
        x^{k+1}  = \prox_{\tau F}(x^k - \tau K^*y^k) &\equivalent x^k - \tau K^* y^k \in  x^{k+1} + \tau \partial F(x^{k+1}) \\
        &\equivalent 0 \in \tau^{-1}(\nextx-\thisx) - K^* (\nexty-\thisy)
        \\ & \qquad\quad
        + [\partial F(\nextx) + K^*\nexty].
    \end{aligned}
\end{equation}
Similarly, the second equation of \eqref{eq:splitting:pdhg} gives
\begin{equation}
    \label{eq:splitting:pdhg-second-step-rewrite}
    \begin{aligned}[t]
        y^{k+1}  = \prox_{\sigma G^*}(y^{k} + \sigma K x^{k+1}) & \equivalent
        \sigma^{-1}y^{k} \in \sigma^{-1} y^{k+1} + \partial G^*(y^{k+1}) - Kx^{k+1}
        \\ &\equivalent 0 \in \sigma^{-1}(\nexty-\thisy) + [\subdiff G^*(\nexty) -K\nextx].
    \end{aligned}
\end{equation}
With the help of \eqref{eq:splitting:pdhg-first-step-rewrite}, \eqref{eq:splitting:pdhg-second-step-rewrite}, and the operator
\begin{equation*}
    \tilde M \defeq \begin{pmatrix} \tau^{-1}\Id &-K^* \\ 0 & \sigma^{-1}\Id \end{pmatrix},
\end{equation*}
we can then rearrange \eqref{eq:splitting:pdhg} as the preconditioned proximal point method \eqref{eq:proximal:preconditioned}. Furthermore, provided the step lengths are such that $M=\tilde M$ is invertible, this can be written as
\begin{equation}
    \label{eq:cp-prox}
    0 \in \Hsaddle(\nextu) + M(\nextu-\thisu)
    \equivalent \nextu = \calR_{M^{-1}\Hsaddle}\thisu.
\end{equation}

However, considering the remarks on convergence in the previous sections, there is a problem: $M$ is not self-adjoint and therefore does not induce a \mbox{(semi-)}norm on $X \times Y$.
We therefore change our algorithm and take
\begin{equation}
    \label{eq:pdps-basic-m}
    M \defeq \begin{pmatrix} \tau^{-1}\Id &-K^* \\ -K & \sigma^{-1}\Id \end{pmatrix}.
\end{equation}
Correspondingly, replacing \eqref{eq:splitting:pdhg-second-step-rewrite} by
\begin{equation*}
    \begin{aligned}
        y^{k+1}  = \prox_{\sigma G^*}(y^{k} + \sigma K(2x^{k+1}-x^k)) & \equivalent
        \sigma^{-1}y^{k} - Kx^{k} \in \sigma^{-1} y^{k+1} + \partial G^*(y^{k+1}) - 2Kx^{k+1}
        \\ & \equivalent 0 \in \sigma^{-1}(\nexty-\thisy) -K(\nextx-\thisx)
        \\ & \qquad\quad + [\subdiff G^*(\nexty) -K\nextx],
    \end{aligned}
\end{equation*}
we then obtain from \eqref{eq:cp-prox} the \term[method!primal-dual proximal splitting]{Primal-Dual Proximal Splitting} (\index{PDPS|see{method, primal-dual proximal splitting}}PDPS) method
\begin{algeqbox}
    \begin{equation}\label{eq:splitting:pd}
        \left\{\begin{aligned}
                x^{k+1} &\defeq \prox_{\tau F}(x^k - \tau K^*y^k),\\
                \overnextx &\defeq 2x^{k+1}-x^k,\\
                y^{k+1} &\defeq \prox_{\sigma G^*}(y^{k} + \sigma K\bar x^{k+1}).
        \end{aligned}        \right.
    \end{equation}
\end{algeqbox}
The middle over-relaxation step is a consequence of our choice of the bottom-left corner of $M$ defined in \eqref{eq:pdps-basic-m}. This itself was forced to have its current form through the self-adjointness requirement on $M$ and the choice of the top-right corner of $M$. As mentioned above, the role of the latter is to \enquote{decouple} the primal update from the dual update by shifting $K^*\nexty$ within $\Hsaddle$ to $K^*\thisy$ so that the primal iterate $\nextx$ can be computed without knowing $\nexty$. (Alternatively, we could zero out the off-diagonal of $M$ and still have a self-adjoint operator, but then we would generally not be able to compute $\nextx$ independent of $\nexty$.)

In the following chapters, we will demonstrate that the PDPS method converges if the step sizes are chosen to ensure $\sigma\tau\norm{K}_{\linear(X;Y)}^2<1$, and that in fact it has particularly good convergence properties.
Note that although the iteration \eqref{eq:splitting:pd} is implicit in $F$ and $G$, it is still explicit in $K$; it is therefore not surprising that step size restrictions based on $K$ remain. Applying, for example, the PDPS method with $\tilde G(x) \defeq G(Kx)$ (i.e., applying only the sum rule but not the chain rule) would lead to a fully implicit method. This would, however, require computing $\inv K$ in the primal proximal step involving $\prox_{\sigma\tilde G^*}$. It is precisely the point of the primal-dual proximal splitting to avoid having to invert $K$, which is often prohibitively expensive if not impossible (e.g., if $K$ does not have closed range as in many inverse problems).

\begin{remark}
    \label{rem:splitting:pdps}
    The primal-dual proximal splitting was first introduced in \cite{pock2009mumford} for specific image segmentation problems, and later more generally in \cite{Pock_PD_2010}. For this reason, it is frequently referred to as the \term[method!Chambolle--Pock]{Chambolle--Pock method}. The relation to proximal point methods was first pointed out in \cite{He:2012}.
    In \cite{esser2010general} it was classified as the \term[method!primal-dual hybrid gradient, modified]{Primal-Dual Hybrid Gradient method, Modified} or \index{PDHGM|see{method, primal-dual hybrid gradient, modified}}PDHGM after the method \eqref{eq:splitting:pdhg}, which is called the PDHG. The latter is due to \cite{zhu2008efficient}.

    Banach space generalizations of the PDPS method, based on a so-called \term[divergence, Bregman]{Bregman divergence} in place of $u \mapsto \frac{1}{2}\norm{u}^2$, were introduced in \cite{hohage2014generalization}. We will discuss Bregman divergences in further detail in \cref{sec:gap:ergodic:bregman}.

    The PDPS method has been also generalized to different types of nonconvex problems in \cite{Valkonen:2014,mollenhoff2014primal}.
    Stochastic generalizations are considered in \cite{tuomov-blockcp,chambolle2017stochastic}.
\end{remark}

\section{Primal-dual explicit splitting}
\label{sec:proximal:gist}

The PDPS method is useful for dealing with the sum of functionals where one summand includes a linear operator. However, if this is the case for \emph{both} operators, i.e.,
\[
    \min_{x\in X} F(Ax) + G(Kx)
\]
for $F:Z\to\Rbar$, $G:Y\to\Rbar$, $A\in\linear(X;Z)$ and $K\in\linear(X;Y)$, we again have the problem of dealing with a complicated proximal mapping.
One workaround is the following ``lifting trick'': we introduce
\begin{equation}
    \label{eq:proximal:dualisation-trick}
    \tilde F(x) \defeq 0,
    \quad
    \tilde G(y, z) \defeq G(y) + F(z)
    \quad\text{and}\quad
    \tilde Kx \defeq (Kx, Ax),
\end{equation}
and then apply the PDPS method to the reformulated problem $\min_x \tilde F(x) + \tilde G(\tilde K x)$.
According to \cref{lem:proximal:calculus}\,\ref{lem:proximal:calculus:iii}, the dual step of the PDPS method will then split into separate proximal steps with respect to $G^*$ and $F^*$, while the proximal point mapping in the primal step will be trivial.
However, an additional dual variable will have been introduced through the introduction of $z$ above, which can be costly.

An alternative approach is the following.
Analogously to \eqref{eq:splitting:pdhg}, but only using \cref{lem:proximal:subdiff} on the second relation of \eqref{eq:splitting:fenchel} together with the chain rule (\cref{thm:convex:chain}), we can reformulate the latter as
\begin{equation}
    \label{eq:splitting:gist:fixed-point:0}
    \left\{
        \begin{aligned}
            \realoptx & \in \realoptx - \tau[A^*\subdiff F(A\realoptx) + K^*\realopty],\\
            \realopty &= \prox_{\sigma G^*}(\realopty + \sigma K\realoptx).
        \end{aligned}
    \right.
\end{equation}
(For $K=\Id$, we can alternatively obtain \eqref{eq:splitting:gist:fixed-point:0} from the derivation of explicit splitting  by using Moreau's identity, \cref{thm:proximal:moreau}, in the second relation of \eqref{eq:splitting:optsys_exp}.)

If $F$ is Gateaux differentiable (and taking $A=\Id$ for the sake of presentation), inserting the first relation in the second relation, \eqref{eq:splitting:gist:fixed-point:0} can be further rewritten as
\begin{algeqbox*}
    \begin{equation*}
        \label{eq:splitting:gist:fixed-point}
        \left\{
            \begin{aligned}
                \realoptx & \defeq \realoptx - \tau[\grad F(\realoptx) + K^*\realopty],\\
                \realopty &\defeq \prox_{\sigma G^*}(\realopty + \sigma K\realoptx- \sigma\tau K[\grad F(\realoptx) + K^*\realopty]).
            \end{aligned}
        \right.
    \end{equation*}
\end{algeqbox*}
Reordering the lines and fixing $\tau=\sigma=1$, the corresponding fixed-point iteration leads to the \term[method!primal-dual explicit splitting]{primal-dual explicit splitting} (PDES\index{PDES|see{method, primal-dual explicit splitting}}) method
\begin{equation}
    \label{eq:splitting:gist}
    \left\{\begin{aligned}
            \nexty &= \prox_{G^*}((\Id-KK^*)\thisy + K(\thisx- \grad F(\thisx))),\\
            \nextx & = \thisx - \grad F(\thisx) - K^*\nexty.\\
    \end{aligned}        \right.
\end{equation}

Again, we can write \eqref{eq:splitting:gist} in more implicit form as
\begin{equation*}
    \left\{\begin{aligned}
            0 & \in \subdiff G^*(\nexty) - K(\thisx- \grad F(\thisx)-K^*\thisy)
            +(\nexty-\thisy),
            \\
            0 & = \grad F(\thisx) + K^*\nexty + (\nextx-\thisx). \\
    \end{aligned}        \right.
\end{equation*}
Inserting the second relation in the first, this is
\begin{equation*}
    \left\{\begin{aligned}
            0 & \in \subdiff G^*(\nexty) - K\nextx + (\Id-KK^*)(\nexty-\thisy),
            \\
            0 & = \grad F(\thisx) + K^*\nexty + (\nextx-\thisx). \\
    \end{aligned}        \right.
\end{equation*}
If we now introduce the preconditioning operator
\begin{equation}
    \label{eq:splitting:gist-m}
    M \defeq \begin{pmatrix} \Id & 0 \\ 0 & \Id-KK^* \end{pmatrix},
\end{equation}
then in terms of the monotone operator $\Hsaddle$ introduced in \eqref{eq:saddle-operator} for the PDPS method and $u=(x, y)$, the PDES method \eqref{eq:splitting:gist} can be written in implicit form as
\begin{equation}
    \label{eq:splitting:gist-implicit}
    0 \in H(\nextu) + \begin{pmatrix} \grad F(\thisx)-\grad F(\nextx) \\ 0 \end{pmatrix}
    + \Precond(\nextu-\thisu).
\end{equation}
The middle term switches the step with respect to $F$ to be explicit.
Note that \eqref{eq:splitting:fb:implicit} could have also been written with a similar middle term; we can therefore think of the PDES method as a \term[method!explicit splitting!preconditioned]{preconditioned explicit splitting} method.

The preconditioning operator $\Precond$ is self-adjoint as well as positive semi-definite if $\norm{K}_{\linear(X;Y)} \le 1$. It does not have the off-diagonal decoupling terms that the preconditioner for the PDPS method has. Instead, through the special structure of the problem, the term $\Id-KK^*$ decouples $\nexty$ from $\nextx$, allowing $\nexty$ to be computed first.

We will in \cref{sec:convergence:fb-general} see that the iterates of the PDES method converge weakly when $\grad F$ is Lipschitz with factor strictly less than $2$.

\begin{remark}
    The primal-dual explicit splitting was introduced in \cite{loris2011generalization} as \term[method!iterative soft-thresholding!generalized]{Generalized Iterative Soft Thresholding} (\index{GIST|see{method, iterative soft-thresholding, generalized}}{GIST}) for $F(x)=\frac{1}{2}\norm{b-x}^2$.
    The general case has later been called the \term[method!primal-dual fixed point]{primal-dual fixed point method} (\index{PDFP|see{method, primal-dual fixed point}}{PDFP}) in \cite{chen2013pdfp} and the \term[method!proximal alternating predictor-corrector]{proximal alternating predictor-corrector} (\index{PAPC|see{method, proximal alternating predictor-corrector}}{PAPC}) in \cite{drori2015simple}.
\end{remark}

\section{Augmented Lagrangian and alternating directions method of multipliers}
\label{sec:proximal:admm}

Let $F:X \to \Rbar$ and $G: Z \to \Rbar$ be convex, proper, and lower semicontinuous.
Also let $A \in \linear(X; Y)$, and $B \in \linear(Z; Y)$, and consider for some $c \in Y$ the problem
\begin{equation}
    \label{eq:admm-problem}
    \min_{x,z}~ F(x) + G(z)  \quad \text{s.t.}\quad Ax+Bz=c.
\end{equation}
A traditional way to handle this kind of constraint problems is by means of the {augmented Lagrangian}. We start by introducing the \term{Lagrangian}
\[
    \mathcal{L}(x, z; \lambda) \defeq  F(x) + G(z) + \iprod{Ax+Bz-c}{\lambda}_Y.
\]
Then \eqref{eq:admm-problem} has the same solutions as the saddle-point problem
\begin{equation}
    \label{eq:proximal:admm:minmax}
    \min_{x \in X, z \in Z} \max_{\lambda \in Y}~ \mathcal{L}(x, z; \lambda).
\end{equation}
We may then \enquote{augment} the Lagrangian by a squared penalty on the violation of the constraint, hence obtaining the equivalent problem
\begin{equation}
    \label{eq:admm-problem-alm}
    \min_{x \in X, z \in Z} \max_{\lambda \in Y}~ \mathcal{L}_\tau(x, z; \lambda) \defeq  F(x) + G(z) + \iprod{Ax+Bz-c}{\lambda}_Y +\frac{\tau}{2}\norm{Ax+Bz-c}_Y^2,
\end{equation}
where $\mathcal{L}_\tau$ is the \term[Lagrangian!augmented]{augmented Lagrangian}.

A classical approach for the solution of \eqref{eq:admm-problem-alm} is by alternately solving for one variable while keeping the others fixed. If we take a proximal step for the dual variable or \term[multiplier, Lagrange]{Lagrange multiplier} $\lambda$, this yields the \term[method!alternating direction]{Alternating Directions Method of Multipliers} (ADMM\index{ADMM|see{method, alternating direction}})
\begin{algeqbox}
    \begin{equation}
        \label{eq:admm-basic}
        \left\{\begin{aligned}
                \nextx & \defeq \argmin_{x \in X}~ \mathcal{L}_\tau(x, \thisz; \this{\lambda}), \\
                \quad
                \nextz & \defeq \argmin_{z \in Z}~ \mathcal{L}_\tau(\nextx, z; \this{\lambda}), \\
                \quad
                \nexxt{\lambda} & \defeq \argmax_{\lambda \in Y}~ \mathcal{L}_\tau(\nextx, \nextz; \lambda) - \frac{1}{2\tau}\norm{\lambda-\this{\lambda}}^2_Y.
        \end{aligned}\right.
    \end{equation}
\end{algeqbox}
This can be rewritten as
\begin{equation}
    \label{eq:admm}
    \left\{\begin{aligned}
            \nextx & \in \inv{(A^*A+\inv\tau \subdiff F)}(A^*(c-B\this{z}-\inv\tau\this{\lambda})), \\
            \nexxt{z} & \in \inv{(B^*B+\inv\tau \subdiff G)}(B^*(c-A\nextx-\inv\tau\this{\lambda})), \\
            \nexxt{\lambda} & \defeq \this{\lambda} + \tau(A\nextx+B\nexxt{z}-c).
    \end{aligned}\right.
\end{equation}
As can be observed, the ADMM requires inverting relatively complicated set-valued operators in place of simple proximal point operations. This is why the basic ADMM is seldom practically implementable without the application of a further optimization method to solve the $x$ and $z$ updates.

In the literature, there have been various remedies to the nonimplementability of the ADMM.
In particular, one can modify the ADMM iterations by adding to \eqref{eq:admm-basic} additional proximal terms. Introducing for some $Q_x \in \linear(X; X)$ and $Q_z \in \linear(Z; Z)$ the weighted norms $\norm{x}_{Q_x}\defeq \sqrt{\iprod{Q_xx}{x}_X}$ and $\norm{z}_{Q_z}\defeq \sqrt{\iprod{Q_zz}{z}_Z}$, this leads to the iteration
\begin{algeqbox}
    \begin{equation}
        \label{eq:admm-preconditioned}
        \left\{\begin{aligned}
                \nextx & \defeq \argmin_{x \in X}~ \mathcal{L}_\tau(x, \thisz; \this{\lambda}) + \frac{1}{2}\norm{x-\thisx}_{Q_x}^2, \\
                \quad
                \nextz & \defeq \argmin_{z \in Z}~ \mathcal{L}_\tau(\nextx, z; \this{\lambda}) + \frac{1}{2}\norm{z-\thisz}_{Q_z}^2, \\
                \quad
                \nexxt{\lambda} & \defeq \argmax_{\lambda \in Y}~ \mathcal{L}_\tau(\nextx, \nextz; \lambda) - \frac{1}{2\tau}\norm{\lambda-\this{\lambda}}_Y^2.
        \end{aligned}\right.
    \end{equation}
\end{algeqbox}
If we specifically take $Q_x \defeq \inv\sigma \Id - \tau A^*A$ and $Q_z \defeq \inv\theta \Id - \tau B^*B$  for some $\sigma,\theta>0$ with $\sigma\tau\norm{A}<1$ and $\theta\tau \norm{B}<1$, then we can expand
\[
    \begin{aligned}
        \mathcal{L}_\tau(x, z; \lambda) + \frac{1}{2}\norm{x-\thisx}_{Q_x}^2
        &
        =
        F(x) + G(z) + \iprod{Ax+Bz-c}{\lambda}_Y
        \\
        \MoveEqLeft[-1]
        +\tau\iprod{x}{A^*(Bz-c)}_X + \frac{\tau}{2}\norm{Bz-c}_Y^2
        \\
        \MoveEqLeft[-1]
        + \frac{1}{2\sigma}\norm{x-\thisx}_X^2
        +\tau \iprod{\nextx}{A^*A \thisx}_X
        -\frac{\tau}{2} \norm{A\thisx}_Y^2,
    \end{aligned}
\]
which has the \enquote{partial} subdifferential $\partial_x$ with respect to $x$ (keeping $z,\lambda$ fixed)
\[
    \subdiff_x \mathcal{L}_\tau(x, z; \lambda)
    =
    \subdiff F(x) + A^*\lambda
    +\tau A^*(Bz-c)
    +\inv\sigma(x-\thisx)
    +\tau A^*A\thisx.
\]
Similarly computing the partial subdifferential $\partial_z$ with respect to $z$, \eqref{eq:admm-preconditioned} can thus be written as the \term[method!alternating direction!preconditioned]{preconditioned ADMM}\index{ADMM!preconditioned|see{method, alternating direction, preconditioned}}
\begin{algeqbox}
    \begin{equation}
        \label{eq:admm0-preconditioned}
        \left\{\begin{aligned}
                \nextx & \defeq \prox_{\sigma F}((\Id-\sigma\tau)A^*A\thisx+\sigma A^*(\tau(c-B\this{z})-\this{\lambda})), \\
                \nexxt{z} & \defeq \prox_{\theta G}((\Id-\theta\tau)B^*B\thisz+\theta B^*(\tau(c-A\nextx)-\this{\lambda})), \\
                \nexxt{\lambda} & \defeq \this{\lambda} + \tau(A\nextx+B\nexxt{z}-c).
        \end{aligned}\right.
    \end{equation}
\end{algeqbox}
We will see in the next section that this method is just the PDPS method with the primal and dual variables exchanged.

\begin{remark}
    \label{remark:admm}
    The ADMM was introduced in \cite{gabay,arrow1958strudies} as an alternating approach to the classical Augmented Lagrangian method.
    The preconditioned ADMM is due to \cite{zhang2011unified}.
\end{remark}

\section{Connections}
\label{sec:proximal:connections}

In \cref{sec:proximal:gist} we have seen the importance and interplay of problem formulation and algorithm choice for problems with a specific structure. We will now see that many of the algorithms we have presented are actually equivalent when applied to different formulations of the problem. Hence, if one algorithm is efficient on one formulation of the problem, another algorithm may work equally well on a different formulation.

We start by considering the ADMM problem \eqref{eq:admm-problem}, which we can reformulate as
\begin{equation*}
    \min_{x,z}~ F(x) + G(z)  + \delta_{\{c\}}(Ax+Bz).
\end{equation*}
Applying the PDPS method \eqref{eq:splitting:pd} to this formulation yields the algorithm
\begin{algeqbox}
    \begin{equation}
        \label{eq:pdps-admm-problem}
        \left\{\begin{aligned}
                \nextx & \defeq \prox_{\tau F}(\thisx - \tau A^*\this{\lambda}),\\
                \nextz & \defeq \prox_{\tau G}(\thisz - \tau B^*\this{\lambda}),\\
                \overnextx & \defeq 2\nextx-\thisx,\\
                \overnext{z} & \defeq 2\nextz-\thisz,\\
                \nexxt{\lambda} & \defeq \this{\lambda} + \sigma (A\overnextx+B\overnext{z}-c).
        \end{aligned}\right.
    \end{equation}
\end{algeqbox}
Note that both the ADMM \eqref{eq:admm} and the preconditioned ADMM \eqref{eq:admm0-preconditioned} have a very similar form to this iteration.
We will now demonstrate that if $A=\Id$ and so $X=Y$, i.e., if we want to solve the (primal) problem
\begin{equation}
    \label{eq:admm-equiv-primal}
    \min_{z \in Z} F(c-Bz) + G(z),
\end{equation}
then the ADMM is equivalent to the PDPS method \eqref{eq:splitting:pd} applied to the (dual) problem
\begin{equation}
    \label{eq:pdps-admm-equiv-problem-primal}
    \min_{y \in Y}~ [G^*(B^*y) - \iprod{c}{y}_Y] + F^*(y),
\end{equation}
where the dual step will be performed with respect to $F^*$.

To make the exact way the PDPS method is applied in each instance clearer, and to highlight the primal-dual nature of the PDPS method, it will be more convenient to write the problem to which the PDPS method is applied in \term[problem!saddle-point]{saddle-point form}. Specifically, keeping in mind \eqref{eq:convex:saddle} together with the discussion following \cref{thm:convex:fenchel}, the problem $\min_x F(x)+G(Kx)$ can be written as the saddle-point problem
\[
    \min_{x \in X}~\max_{y \in Y} F(x) + \iprod{Kx}{y}_Y - G^*(y).
\]
This formulation also shows the dual variable directly in the problem formulation.
Applied to \eqref{eq:pdps-admm-equiv-problem-primal}, we then obtain the problem
\begin{equation}
    \label{eq:pdps-admm-equiv-problem}
    \min_{y \in Y} \max_{x \in X}~ [G^*(B^*y) - \iprod{c}{y}_Y] + \iprod{x}{y}_Y - F(x).
\end{equation}
Our claim is that the PDPS method applied to this saddle-point formulation is equivalent to the ADMM in case of $A=\Id$. The iterates of the two algorithms will be different, as the variables solved for will be different aside from the shared $x$. However, all the variables will be related by affine transformations.

We will also demonstrate that the preconditioned ADMM is equivalent to the PDPS method when $B=\Id$.
In fact, we will demonstrate a chain of relationships from ADMM or preconditioned ADMM (primal problem) via the PDPS (saddle-point problem) method to the DRS method (dual problem); the equivalence between the ADMM and the DRS method even holds generally.

To demonstrate the idea, we start with $A=B=\Id$. Then \eqref{eq:admm} reads
\begin{algeqbox}
    \begin{equation}
        \label{eq:admm0-identity}
        \left\{\begin{aligned}
                \nextx & \defeq \prox_{\inv\tau F}(c-\this{z}-\inv\tau\this{\lambda}), \\
                \nextz & \defeq \prox_{\inv\tau G}(c-\nextx-\inv\tau\this{\lambda}), \\
                \nexxt{\lambda} & \defeq \this{\lambda} + \tau(\nextx+\nexxt{z}-c).
        \end{aligned}\right.
    \end{equation}
\end{algeqbox}
Using the third step for the previous iteration to obtain an expression for $\thisz$, we can rewrite the first step as
\[
    \nextx \defeq \prox_{\inv\tau F}(\thisx-\inv\tau(2\this{\lambda}-\prev{\lambda})).
\]
If we use  \cref{lem:proximal:calculus}\,\ref{item:proximal:calculus:conjugate},
the second step reads
\begin{equation*}
    \nexxt{z} \defeq (c-\nextx-\inv\tau\this{\lambda}) - \inv\tau \prox_{\tau G^*}(\tau(c-\nextx)-\this{\lambda}).
\end{equation*}
Keeping in mind the third step of \eqref{eq:admm0-identity}, this yields
$\nexxt{\lambda} = - \prox_{\tau G^*}(\tau(c-\nextx)-\this{\lambda})$.
Replacing $\nexxt{\lambda}$ by $\nexty \defeq - \nexxt{\lambda}$, moving $c$ into the proximal part, and reordering the steps such that $\nextx$ becomes $\thisx$, transforms \eqref{eq:admm0-identity} into
\begin{equation}
    \label{eq:admm0-identity-transformed1}
    \left\{\begin{aligned}
            \nexty & \defeq \prox_{\tau (G^*-\iprod{c}{\freevar})}(\thisy-\tau \thisx), \\
            \nextx & \defeq \prox_{\inv\tau F}(\thisx+\inv\tau(2\nexty-\thisy)). \\
    \end{aligned}\right.
\end{equation}
This is the PDPS method applied to \eqref{eq:pdps-admm-equiv-problem} with $B=\Id$. However, the step lengths $\tau$ and $\sigma=\inv\tau$ do not satisfy $\tau\sigma\norm{K}^2<1$, which would be needed to deduce convergence of the ADMM from that of the PDPS method.
But we will see in \cref{chap:gap} that these step lengths at least lead to convergence of a certain \enquote{Lagrangian duality gap}, and for the ADMM we can in general only prove such gap estimates.

To show the relation of ADMM to implicit splitting, we further use \cref{lem:proximal:calculus}\,\ref{item:proximal:calculus:conjugate} in the second step of \eqref{eq:admm0-identity-transformed1} to obtain
\[
    \nextx = \inv\tau(2\nexty-\thisy)+\thisx - \inv\tau \prox_{\tau F^*}(2\nexty-\thisy+\tau\thisx).
\]
Introducing $\nexxt{w} \defeq \nexty-\tau\nextx$ and changing variables, we thus transform  \eqref{eq:admm0-identity-transformed1} into
\begin{algeqbox*}
    \begin{equation*}
        \label{eq:admm0-identity-transformed2}
        \left\{\begin{aligned}
                \nexty & \defeq \prox_{\tau (G^*-\iprod{c}{\freevar})}(\this{w}), \\
                \nexxt{w} & \defeq \this{w}-\nexty + \prox_{\tau F^*}(2\nexty-\this{w}). \\
        \end{aligned}\right.
    \end{equation*}
\end{algeqbox*}
But this is the DRS method \eqref{eq:splitting:dr} applied to
\[
    \min_{x \in X}~ F^*(x) + [G^*(x)-\iprod{c}{x}_X].
\]
Recall now from \cref{lem:convex:fenchel_calc}\,\ref{lem:convex:fenchel_calc:ii} that $[G(c-\freevar)]^*=G^*(- \freevar) + \iprod{c}{\freevar}_Y$. \Cref{thm:convex:fenchel} thus shows that this is the dual problem of \eqref{eq:admm-equiv-primal}, so we can at least deduce from \cref{cor:convergence:drs} the convergence of $y^k$ to a solution of the dual problem.

We can make the correspondence more general with the help of the following generalization of Moreau's identity (\cref{thm:proximal:moreau}).
\begin{lemma}
    \label{lemma:connections:moreau-generalisation}
    Let $S=G \circ K$ for convex, proper, and lower semicontinuous $G: Y \to \Rbar$ and $K \in \linear(X; Y)$. If there exists an $x_0\in \dom S$ such that $Kx_0\in \interior (\dom G)$, then for all $x \in X$ and $\gamma > 0$,
    \[
        x=\prox_{\gamma S}(x)+\gamma K^*\inv{(KK^* + \inv\gamma \subdiff G^*)}(\inv\gamma K x).
    \]
    In particular,
    \[
        \prox_{S^*}(x)=K^*\inv{(KK^* + \subdiff G^*)}(K x).
    \]
\end{lemma}
\begin{proof}
    By \cref{thm:convex:fenchel}, $w=\prox_{\gamma S}(x)$ if and only if for some $y^* \in Y^*$ it holds that
    \begin{equation*}
        \left\{
            \begin{aligned}
                -K^*y^* &\in w-x,\\
                y^* &\in \gamma \partial G(K w).
            \end{aligned}
        \right.
    \end{equation*}
    In other words, by \cref{lem:convex:fenchel-young},
    \begin{equation*}
        \left\{
            \begin{aligned}
                -K^* y^* & =  w-x, \\
                Kw &\in \partial G^*(\inv \gamma y^*).
            \end{aligned}
        \right.
    \end{equation*}
    Applying $K$ to the first relation, inserting the second, and multiplying by $\gamma^{-1}$ yields
    \begin{equation*}
        KK^*\inv\gamma y^* + \inv\gamma \subdiff G^*(\inv\gamma y^*)\ni\inv\gamma K x,
    \end{equation*}
    i.e., $\inv\gamma y^* \in \inv{(KK^* + \inv\gamma \subdiff G^*)}(\inv\gamma Kx)$. Combined with $-K^* y^* =  w-x$, this yields the first claim. The second claim then follows from \cref{thm:moreau:conjugate} together with the first claim for $\gamma=1$.
\end{proof}

\begin{theorem}
    \label{thm:connections:equivalence}
    Let $F:X \to \Rbar$ and $G: Z \to \Rbar$ be convex, proper, and lower semicontinuous.
    Also let $A \in \linear(X; Y)$, and $B \in \linear(Z; Y)$, and $c \in Y$.
    Assume the existence of a point $(x_0, z_0) \in \dom F \times \dom G$ with $Ax_0+Bz_0=c$.
    Then the iterates of the following algorithms can be transformed to one another with affine transformations and (to obtain the ADMM) the addition of elements of $\kernel A$ and $\kernel B$:
    \begin{enumerate}
        \item
            The ADMM applied to the (primal) problem
            \begin{equation}
                \label{eq:connections:equivalence:problem:admm}
                \min_{x \in X, z \in Z}~ F(x) + G(z)  \quad \text{s.t.}\quad Ax+Bz=c.
            \end{equation}
        \item
            The DRS method applied to the (dual) problem
            \begin{equation}
                \label{eq:connections:equivalence:problem:drs}
                \min_{y \in Y}~ F^*(A^* y) + [G^*(B^* y) - \iprod{c}{y}_Y].
            \end{equation}
        \item
            If $A=\Id$, $X=Y$, and $\sigma=\inv\tau$, the PDPS method applied to the (saddle-point) problem
            \begin{equation}
                \label{eq:connections:equivalence:problem:pdps}
                \min_{y \in Y} \max_{x \in X}~ [G^*(B^* y) - \iprod{c}{y}_Y] + \iprod{x}{y}_Y - F(x).
            \end{equation}
    \end{enumerate}
\end{theorem}
\begin{proof}
    We first show that the ADMM updates for \eqref{eq:connections:equivalence:problem:admm} can be transformed, via affine transformations alone, to the DRS updates for \eqref{eq:connections:equivalence:problem:drs}, and the PDPS updates for \eqref{eq:connections:equivalence:problem:pdps}.
    Observe that the assumption on the existence of $(x_0, z_0)$ ensures that the infimum in \eqref{eq:connections:equivalence:problem:admm} is finite.
    Thus, multiplying the first and second updates of \eqref{eq:admm} by $A$ and $B$, and changing variables $\nextx$ and $\nexxt{z}$ to $\nexxt{\tilde x} \defeq A\nextx$ and $\nexxt{\tilde z}\defeq B \nextz$, we obtain
    \begin{equation}
        \label{eq:admm0-multiplied}
        \left\{\begin{aligned}
                \nexxt{\tilde x} & \in A\inv{(A^*A+\inv\tau \subdiff F)}(A^*(c-\this{\tilde z}-\inv\tau\this{\lambda})), \\
                \nexxt{\tilde z} & \in B\inv{(B^*B+\inv\tau \subdiff G)}(B^*(c-\nexxt{\tilde x}-\inv\tau\this{\lambda})), \\
                \nexxt{\lambda} & \defeq \this{\lambda} + \tau(A\nextx+B\nexxt{z}-c).
        \end{aligned}\right.
    \end{equation}
    Using \cref{lemma:connections:moreau-generalisation} with $\thisy \defeq - \this{\lambda}$ and
    $-\nexty = - \thisy + \tau(\nexxt{\tilde x}+\nexxt{z}-c)$, we transform this as above to
    \begin{equation}
        \label{eq:admm0-multiplied-pdps}
        \left\{\begin{aligned}
                \nexty & \in \prox_{\tau G^* \circ B^*}(\tau(c-\this{\tilde x})+\thisy), \\
                \nexxt{\tilde x} & \in A\inv{(A^*A+\inv\tau \subdiff F)}(A^*(2\nexty-\thisy+\this{\tilde x})).
        \end{aligned}\right.
    \end{equation}
    If $A=\Id$, this is the PDPS method for \eqref{eq:connections:equivalence:problem:pdps} with the iterate equivalence $\nexxt{\tilde x}=\nextx$.
    We continue with \cref{lemma:connections:moreau-generalisation} and  $\nexxt{w} \defeq \nexty-\tau\nexxt{\tilde x}$ to transform \eqref{eq:admm0-multiplied-pdps} further into
    \begin{equation}
        \label{eq:admm0-multiplied-drs}
        \left\{\begin{aligned}
                \nexty & \defeq \prox_{\tau (G^* \circ B^* - \iprod{c}{\freevar})}(\this{w}), \\
                \nexxt{w} & \defeq \this{w}-\nexty + \prox_{\tau F^* \circ A^*}(2\nexty-\this{w}). \\
        \end{aligned}\right.
    \end{equation}
    This is the DRS method for \eqref{eq:connections:equivalence:problem:drs}.

    In the other direction, it is clear from the derivation above that the DRS \eqref{eq:admm0-multiplied-drs} generates \eqref{eq:admm0-multiplied-pdps} and, via affine transformations, its iterates.
    Likewise \eqref{eq:admm0-multiplied-pdps} can be transformed back into \eqref{eq:admm0-multiplied-drs} by reversing the steps.
    The passage from the iterates of \eqref{eq:admm0-multiplied} back to the iterates of the ADMM \eqref{eq:admm} cannot be achieved with affine transformations alone, unless $A$ and $B$ are injective. However, when there exist $\nexxt{\tilde x}$ and $\nexxt{\tilde z}$ solving \eqref{eq:admm0-multiplied}, there must exist some $\nextx$ and $\nextz$ with $\nexxt{\tilde x} \defeq A\nextx$ and $\nexxt{\tilde z}\defeq B \nextz$ that satisfy \eqref{eq:admm}.

    We still need to establish the claimed duality relationship between the problems \eqref{eq:connections:equivalence:problem:admm}, \eqref{eq:connections:equivalence:problem:drs}, and \eqref{eq:connections:equivalence:problem:pdps}.
    To pass from \eqref{eq:connections:equivalence:problem:admm} to \eqref{eq:connections:equivalence:problem:drs}, we would like to apply \cref{thm:convex:fenchel} to $\tilde F(x, z) \defeq F(x)+G(z)$, $\tilde G(y) \defeq \delta_{\{y=c\}}(y)$, and $\tilde K \defeq  (A, B)$. However, $\dom G=\{c\}$ has empty interior, so condition \ref{thm:convex:fenchel:ii} of the theorem does not hold. Recalling \cref{rem:convex:attouch-brezis,rem:fenchel:attouch-brezis}, we can however replace the interior with the relative interior $\ri \dom \tilde G =\{c\}$.
    Thus the condition reduces to the existence of $y_0 \in \dom \tilde F$ with $Ky_0=c$, which is satisfied by $y_0=(x_0, z_0)$.

    Finally, the relationship to \eqref{eq:connections:equivalence:problem:pdps} when $A=\Id$ is immediate from \eqref{eq:connections:equivalence:problem:drs} and the definition of the conjugate function $F^*$.
    The existence of a saddle point follows from the proof of \cref{thm:convex:fenchel}.
\end{proof}

The methods in the proof of \cref{thm:connections:equivalence} are rarely computationally feasible or efficient unless $A=B=\Id$, due to the difficult proximal mappings for compositions of functionals with operators or the set-valued operator inversions required.
On the other hand, the PDPS method \eqref{eq:pdps-admm-problem} only requires that we can compute the proximal mappings of $G$ and $F$. This demonstrates the importance of problem formulation.

\bigskip

Similar connections hold for the preconditioned ADMM \eqref{eq:admm0-preconditioned}.
With the help of the third step of \eqref{eq:admm0-preconditioned}, the first step can be rewritten as
\[
    \nextx \defeq \prox_{\sigma F}(\thisx-\sigma A^*(2\this{\lambda}-\prev{\lambda})).
\]
If $\theta\tau=1$ and $B=\Id$, the second step reads
\[
    \nexxt{z} \defeq \prox_{\inv\tau G}((c-A\nextx)-\inv\tau\this{\lambda}).
\]
We transform this with \cref{lem:proximal:calculus}\,\ref{item:proximal:calculus:conjugate} into
\[
    \nextz
    =(c-A\nextx)-\inv\tau\this{\lambda} - \inv\tau \prox_{\tau G^*}(\tau(c-A\nextx)-\this{\lambda}).
\]
Using the third step of \eqref{eq:admm0-preconditioned}, this is equivalent to
\begin{equation*}
    -\nexxt{\lambda} = \prox_{\tau G^*}(\tau(c-A\nextx)-\this{\lambda}).
\end{equation*}
Introducing $\nexty \defeq -\nexxt{\lambda}$ and changing the order of the first and second step, we therefore transform \eqref{eq:admm0-preconditioned} into the PDPS method
\begin{algeqbox}
    \begin{equation}
        \label{eq:admm0-preconditioned-pdps}
        \left\{\begin{aligned}
                \nexty & \defeq \prox_{\tau G^*}(\thisy-\tau A\thisx), \\
                \nextx & \defeq \prox_{\sigma F}(\thisx+\sigma A^*(2\nexty-\thisy)).
        \end{aligned}\right.
    \end{equation}
\end{algeqbox}

We therefore have obtained the following result.
\begin{theorem}
    Let $F:X \to \Rbar$ and $G: Y \to \Rbar$ be convex, proper, and lower semicontinuous.
    Also let $A \in \linear(X; Y)$ and $c \in Y$.
    Assume the existence of a  point $(x_0, z_0) \in \dom F \times \dom G$ with $Ax_0+z_0=c$.
    Take $\theta=\inv\tau$. Then subject to affine transformations to obtain iterates not explicitly generated in each case, the following are equivalent:
    \begin{enumerate}
        \item The preconditioned ADMM \eqref{eq:admm0-preconditioned} applied to the (primal) problem
            \[
                \min_{x \in X, z \in Y} F(x) + G(z)  \quad\text{s.t.}\quad Ax+z=c.
            \]
        \item The PDPS method applied to the (saddle point) problem
            \begin{equation*}
                \min_{y \in Y} \max_{x \in X}~ [G^*(y) - \iprod{c}{y}_Y] + \iprod{Ax}{y}_Y - F(x).
            \end{equation*}
        \item If $A=\Id$, $X=Y$, and $\sigma=\inv\tau$, the Douglas--Rachford splitting method applied to the (dual) problem
            \[
                \min_{y \in X}~ F^*(y) + [G^*(y)- \iprod{c}{y}_Y].
            \]
    \end{enumerate}
\end{theorem}
\begin{proof}
    We have already proved the equivalence of the preconditioned ADMM and the PDPS method.
    For equivalence to the DRS method, we observe that under the additional assumptions of this theorem, \eqref{eq:admm0-preconditioned-pdps} reduces to \eqref{eq:admm0-identity-transformed1}.
\end{proof}

\chapter{Splitting methods: weak convergence}\label{chap:convergence}

Now that we have in the previous chapter derived several iterative procedures through the manipulation of fixed-point equations, we have to show that they indeed converge to a fixed point (which by construction is then the solution of an optimization problem, making these procedures optimization algorithms). We start with weak convergence, as this is the most that can generally be expected.

The classical approach to proving weak convergence is by introducing suitable contractive (or at least firmly nonexpansive) operators related to the algorithm and then applying classical fixed-point theorems (see \cref{remark:convergence:browder} below).
We will instead introduce a very direct approach that will then extend in the following chapters to be also capable of proving convergence rates. The three main ingredients of all convergence proofs will be
\begin{enumerate}
    \item The three-point identity \eqref{eq:hilbert:three-point-identity}, which we recall here as
    \begin{equation}
        \label{eq:convergence:three-point-identity}
        \iprod{x-y}{x-z}_X
        = \frac{1}{2}\norm{x-y}_X^2
        - \frac{1}{2}\norm{y-z}_X^2
        + \frac{1}{2}\norm{x-z}_X^2 \quad\text{for all }x,y,z\in X.
    \end{equation}

\item The monotonicity of the operator $H$ whose roots we seek to find (which in the simplest case equals $\subdiff F$ for the functional $F$ we want to minimize).

    \item The nonnegativity of the preconditioning operators $M$ defining the implicit forms of the algorithms we presented in \cref{chap:proximal}.
\end{enumerate}
In the later chapters, stronger versions of the last two ingredients will be required to obtain convergence rates and the convergence of function value differences $F(\nextx)-F(\realoptx)$ or of more general gap functionals.

\section{Opial's lemma and Fej\'er monotonicity}
\label{sec:convergence:opial}

The next lemma forms the basis of all our weak convergence proofs. It is a generalized subsequence argument, showing that if all weak accumulation points of a sequence lie in a set and if the sequence does not diverge (in the strong sense) away from this set, the full sequence converges weakly.
We recall that $\optx \in X$ is a weak(-$*$) accumulation point of the sequence $\{\thisx\}_{k \in \N}$ if there exists a subsequence such that $x^{k_\ell} \weakto \optx$ weakly(-$*$) in $X$.

\begin{lemma}[Opial]\index{lemma!Opial}
    \label{lemma:opial}
    Let $X$ be a Hilbert space and $\hat X \subset X$ be a nonempty subset. If the sequence $\{\thisx\}_{k \in \N} \subset X$ satisfies
    \begin{enumerate}
        \item\label{item:opial-nonincreasing} $\norm{\nextx-\optx}_X \le \norm{\thisx-\optx}_X$ for all $\optx \in \hat X$ and $k \in \N$;
        \item\label{item:opial-limit} all weak accumulation points of $\{\thisx\}_{k \in \N}$ belong to $\hat X$;
    \end{enumerate}
    then $\thisx \weakto \realoptx$ in $X$ for some $\realoptx \in \hat X$.
\end{lemma}
\begin{proof}
    First, the assumption \ref{item:opial-nonincreasing} implies that the sequence $\{\thisx\}_{k\in\N}$ is bounded and hence by \cref{thm:ebsmul} contains a weakly convergent subsequence.
    Let now $\bar x$ and $\hat x$ be weak accumulation points.
    The assumption \ref{item:opial-nonincreasing} then implies that both $\{\norm{x^k-\bar x}_X\}_{k\in\N}$ and $\{\norm{x^k-\hat x}_X\}_{k\in\N}$ are decreasing and bounded from below and therefore convergent.
    This yields that
    \begin{equation*}
        \iprod{x^k}{\bar x - \hat x}_X = \frac12\left(\norm{x^k-\hat x}_X^2 - \norm{x^k-\bar x}_X^2 + \norm{\bar x}_X^2 - \norm{\hat x}_X^2\right)\to c\in \R.
    \end{equation*}
    Since $\bar x$ is a weak accumulation point, there exists a subsequence $\{x^{k_n}\}_{n\in \N}$ with $x^{k_n}\weakto \bar x$; similarly, there exists a subsequence $\{x^{k_m}\}_{m\in \N}$ with $x^{k_m}\weakto \hat x$. Hence,
    \begin{equation*}
        \iprod{\bar x}{\bar x - \hat x}_X =  \lim_{n\to \infty} \iprod{x^{k_n}}{\bar x - \hat x}_X = c =  \lim_{m\to \infty} \iprod{x^{k_m}}{\bar x - \hat x}_X=\iprod{\hat x}{\bar x - \hat x}_X,
    \end{equation*}
    and therefore
    \begin{equation*}
        0 = \iprod{\bar x - \hat x}{\bar x - \hat x}_X=\norm{\bar x -\hat x}_X^2,
    \end{equation*}
    i.e., $\bar x =\hat x$. Every convergent subsequence thus has the same weak limit (which lies in $\hat X$ by assumption \ref{item:opial-limit}).
    The claim now follows from a standard subsequence--subsequence argument: Assume to the contrary that there exists a subsequence of $\{x^k\}_{k\in\N}$ that does not converge to $\hat x$. Then we can apply the above argument to obtain a further subsequence converging to $\hat x$, which is a contradiction to the fact that any subsequence of a convergent sequence converges to the same limit.
\end{proof}
A sequence satisfying the condition \ref{item:opial-nonincreasing} is called \term[sequence!Fejér monotone]{Fejér monotone} (with respect to $\hat X$); this is a crucial property of iterates generated by any fixed-point algorithm.

\begin{remark}
    \Cref{lemma:opial} first appeared in the proof of \cite[Theorem 1]{opial1967weak}.
    (There $\hat X$ is assumed to be closed and convex, but we do not require this since condition \ref{item:opial-limit} is already sufficient to show the claim.)

    The concept of Fejér monotone sequences first appears in \cite{fejer1922}, where it was observed that for every point outside the convex hull of a subset of the Euclidean plane, it is always possible to construct a point that is closer to each point in the subset than the original point (and that this property in fact characterizes the convex hull). The term \emph{Fejér monotone} itself appears in \cite{motzkin_schoenberg_1954}, where this construction is used to show convergence of an iterative scheme for the projection onto a convex polytope.
\end{remark}

\section{The fundamental methods: proximal point and explicit splitting}
\label{sec:convergence:fundamental}

Using Opial's \cref{lemma:opial}, we can fairly directly show weak convergence of the proximal point and forward-backward splitting methods. As in the last chapter, we assume throughout the following that $X$ is a Hilbert space.

\subsection*{Proximal point method}\label{ssec:convergence:ppm}

We recall our most fundamental nonsmooth optimization algorithm, the proximal point method.
For later use, we treat the general version of \eqref{eq:proximal:ppa} for an arbitrary set-valued operator $H:X\setto X$, i.e.,
\begin{algeqbox}
    \begin{equation}\label{eq:convergence:ppm-general}
        x^{k+1} \defeq \calR_{\tau_k H}(x^k).
    \end{equation}
\end{algeqbox}
We will need the next lemma to allow a very general choice of the step lengths $\{\tau_k\}_{k \in \N}$. (If we assume $\tau_k \ge \varepsilon > 0$, in particular if we keep $\tau_k \equiv \tau$ constant, it will not be needed.) For the statement, note that by the definition of the resolvent, \eqref{eq:convergence:ppm-general} is equivalent to $\tau_k^{-1}(\thisx-\nextx)\in H(\nextx)$.

\begin{lemma}
    \label{lemma:convergence:tau-sequence-w}
    Let $\{\tau_k\}_{k\in\N}\subset(0,\infty)$ with $\sum_{k=0}^\infty \tau_k^2 = \infty$, and let $\Hany: X \setto X$ be monotone.
    Suppose $\{\thisx\}_{k \in \N}$ and $\nexxt{w} \defeq -\inv\tau_k(\nextx-\thisx)$ satisfies
    \begin{enumerate}
        \item\label{it:convergence:tau-sequence-w:inc} $\displaystyle 0\neq \nexxt{w} \in \Hany(\nextx)$ and
        \item\label{it:convergence:tau-sequence-w:bnd} $\displaystyle\sum_{k=0}^{\infty} \tau_k^2 \norm{\this{w}}_X^2 < \infty$.
    \end{enumerate}
    Then $\norm{\this{w}}_X \to 0$.
\end{lemma}

\begin{proof}
    Since $w^k \in \Hany(x^k)$ and $\Hany$ is monotone, we have from the definition of $w^k$ that
    \[
        \begin{aligned}
            0
            &
            \le
            \iprod{\nexxt{w}-\this{w}}{\nextx-\thisx}_X
            =
            \tau_k \iprod{\this{w}-\nexxt{w}}{\nexxt{w}}_X
            \le
            \tau_k \norm{\nexxt{w}}_X(\norm{\this{w}}_X-\norm{\nexxt{w}}_X).
        \end{aligned}
    \]
    Thus the nonnegative sequence $\{\norm{w^k}_X\}_{k \in \N}$ is decreasing and hence converges to some $M\geq 0$.
    If $M>0$, we can deduce from the choice of $\{\tau_k\}_{k\in \N}$, assumption \ref{it:convergence:tau-sequence-w:bnd}, and the decreasing property that
    \begin{equation*}
        \infty > \frac1{M^2}\sum_{k=0}^\infty \tau_k^2 \norm{w^k}_X^2 \geq \frac1{M^2}\sum_{k=0} \tau_k^2 M^2
        = \sum_{k=0}^\infty \tau_k^2 = \infty,
    \end{equation*}
    a contradiction. Hence $M=0$ and therefore $\norm{w^{k}}_X \to 0$ as claimed.
\end{proof}
This shows that the \enquote{generalized residual} $w^k$ in the inclusion $w^k\in H(x^k)$ converges (strongly) to zero. As usual, this does not (yet) imply that $\{x^k\}_{k\in\N}$ itself converges; but if it does, we expect the limit to be a root of $H$. This is what we prove next, using the three fundamental ingredients we introduced in the beginning of the chapter.

\begin{theorem}\label{thm:convergence:prox}
    Let $\Hany:X\setto X$ be monotone and weak-to-strong outer semicontinuous with $\Hany^{-1}(0)\neq \emptyset$. Furthermore, let $\{\tau_k\}_{k\in\N}\subset(0,\infty)$ with $\sum_{k=0}^\infty \tau_k^2 = \infty$.
    If $\{x^k\}_{k\in\N}\subset X$ is given by the iteration \eqref{eq:convergence:ppm-general}
    for any initial iterate $x^0 \in X$,
    then $x^k\weakto \realoptx$ for some root $\realoptx \in \inv \Hany(0)$.
\end{theorem}

\begin{proof}
    We recall that the proximal point iteration can be written in implicit form as
    \begin{equation}
        \label{eq:convergence:prox}
        0 \in \tau_k\Hany(\nextx)+(\nextx-\thisx).
    \end{equation}
    We ``test'' \eqref{eq:convergence:prox} by the application of $\iprod{\freevar}{\nextx-\realoptx}_X$ for an arbitrary $\realoptx \in \inv \Hany(0)$. Thus we obtain
    \begin{equation}
        \label{eq:convergence:prox-tested}
        0 \in \iprod{\tau_k\Hany(\nextx)+(\nextx-\thisx)}{\nextx-\realoptx}_X,
    \end{equation}
    where the right-hand side should be understood as the set of all possible inner products involving elements of $\Hany(\nextx)$.
    By the monotonicity of $\Hany$, since $0 \in H(\realoptx)$, we have
    \begin{equation*}
        \iprod{\Hany(\nextx)}{\nextx-\realoptx}_X \ge 0,
    \end{equation*}
    which again should be understood to hold for any $w\in \Hany(\nextx)$. (We will frequently make use of this notation and the one from \eqref{eq:convergence:prox-tested} throughout this and the following chapters to keep the presentation concise.)
    Thus \eqref{eq:convergence:prox-tested} yields
    \begin{equation*}
        \iprod{\nextx-\thisx}{\nextx-\realoptx}_X \le 0.
    \end{equation*}
    Applying now the three-point identity \eqref{eq:convergence:three-point-identity} for $x=\nextx$, $y=\thisx$, and $z=\realoptx$, yields
    \begin{equation}
        \label{eq:convergence:prox-est2}
        \frac{1}{2}\norm{\nextx-\realoptx}_X^2
        +\frac{1}{2}\norm{\nextx-\thisx}_X^2
        \le \frac{1}{2}\norm{\thisx-\realoptx}_X^2.
    \end{equation}
    This shows the Fejér monotonicity of $\{\thisx\}_{k \in \N}$ with respect to $\hat X=\inv \Hany(0)$.

    Furthermore, summing \eqref{eq:convergence:prox-est2} over $k=0,\ldots,N-1$ gives
    \begin{equation}
        \label{eq:convergence:prox-est3}
        \frac{1}{2}\norm{x^N-\realoptx}_X^2
        + \sum_{k=0}^{N-1}
        \frac{1}{2}\norm{\nextx-\thisx}_X^2
        \le \frac{1}{2}\norm{x^0-\realoptx}_X^2 =: C_0.
    \end{equation}
    Writing $w^{k+1} \defeq -\inv\tau_{k}(x^{k+1}-x^{k})$, the implicit iteration \eqref{eq:convergence:prox} shows that $\nexxt{w} \in \Hany(\nextx)$, and \eqref{eq:convergence:prox-est3} implies that
    \[
        \sum_{k=0}^{N-1} \tau_k^2 \norm{\nexxt{w}}_X^2 \le 2C_0.
    \]
    If $w^k=0$ for some $k\in\N$, then $0\in \Hany(x^k)$ and hence we are done. Otherwise, we can let $N\to \infty$ and apply \cref{lemma:convergence:tau-sequence-w} to deduce that $\norm{\nexxt{w}}_X \to 0$.

    Let finally $\optx$ be any weak accumulation point of $\{\thisx\}_{k \in \N}$, that is $x^{k_i} \weakto \optx$ for a subsequence $\{k_i\}_{i \in \N} \subset \N$.
    Recall that $w^{k_i} \in \Hany(x^{k_i})$. The weak-to-strong outer semicontinuity of $H$ now immediately yields $0 \in \Hany(\optx)$.
    We then finish by applying Opial's \cref{lemma:opial} for the set $\hat X = \inv \Hany(0)$.
\end{proof}

Note that the conditions of \cref{thm:convergence:prox} are in particular satisfied if $\Hany$ is either maximally monotone (\cref{cor:monoton:closed}) or monotone and BCP outer semicontinuous (\cref{lemma:bpr-strong-dual-corollary}). In particular, applying \cref{thm:convergence:prox} to $H=\partial J$ yields the convergence of the proximal point method \eqref{eq:proximal:ppa} for any proper, convex, and lower semicontinuous functional $J:X\to\Rbar$.

\begin{remark}\label{remark:convergence:browder}
    A conventional way of proving the convergence of the proximal point method is with Browder's fixed-point theorem \cite{browder1965nonexpansive},
    which shows the existence of fixed points of firmly nonexpansive or, more generally, $\alpha$-averaged mappings.
    (We have already shown in \cref{lem:proximal:firmly-nonexpansive} the firm nonexpansivity of the proximal map.)
    On the other hand, to prove Browder's fixed-point theorem itself, we can  use similar arguments as \cref{thm:convergence:prox}, see \cref{thm:convergence:browder} below.
\end{remark}

\subsection*{Explicit splitting}\label{ssec:convergence:fbs}

The convergence of the forward-backward splitting method
\begin{algeqbox}
    \begin{equation}
        \label{eq:convergence:fb}
        x^{k+1} \defeq \prox_{\tau_k G}(x^k - \tau_k \nabla F(x^k))
    \end{equation}
\end{algeqbox}
can be shown analogously.
To do so, we need to assume the Lipschitz continuity of the gradient of $F$ (since we are not using a proximal point mapping for $F$ which is always firmly nonexpansive and hence Lipschitz continuous).

\begin{theorem}\label{thm:convergence:fb}
    Let $F:X\to\R$ and $G:X\to\Rbar$ be proper, convex, and lower semicontinuous. Suppose $\inv{(\subdiff(F+G))}(0) \ne \emptyset$, i.e., that $J \defeq F + G$ has a minimizer. Furthermore, let $F$ be Gateaux differentiable with $L$-Lipschitz gradient. If $0<\tau_{\min}\leq \tau_k \leq \tau_{\max} < 2L^{-1}$, then for any initial iterate $x^0 \in X$ the sequence generated by \eqref{eq:convergence:fb} converges weakly to a root $\realoptx \in \inv{(\subdiff(F+G))}(0)$.
\end{theorem}

\begin{proof}
    We again start by writing \eqref{eq:convergence:fb} in implicit form as
    \begin{equation}
        \label{eq:convergence:fb-prox}
        0 \in \tau_k[\subdiff G(\nextx)+\grad F(\thisx)]+(\nextx-\thisx).
    \end{equation}
    By the monotonicity of $\subdiff G$ and the three-point monotonicity \eqref{eq:smoothness:three-point:monotonicity} of $F$ from \cref{cor:smoothness:three-point}, we first deduce for any $\realoptx\in \hat X \defeq \inv{(\subdiff(F+G))}(0)$ that
    \begin{equation*}
        \iprod{\subdiff G(\nextx)+\grad F(\thisx)}{\nextx-\realoptx}_X \ge -\frac{L}{4}\norm{\nextx-\thisx}_X^2.
    \end{equation*}
    Thus, again testing \eqref{eq:convergence:fb-prox} with $\iprod{\freevar}{\nextx-\realoptx}_X$ yields
    \begin{equation*}
        \iprod{\nextx-\thisx}{\nextx-\realoptx}_X \le \frac{L\tau_k}{4}\norm{\nextx-\thisx}_X^2.
    \end{equation*}
    The three-point identity \eqref{eq:convergence:three-point-identity} now implies that
    \begin{equation}
        \label{eq:convergence:fb-est}
        \frac{1}{2}\norm{\nextx-\realoptx}_X^2
        +\frac{1-\tau_k L/2}{2}\norm{\nextx-\thisx}_X^2
        \le \frac{1}{2}\norm{\thisx-\realoptx}_X^2.
    \end{equation}
    The assumption $2>\tau_k L$ then establishes the Fejér monotonicity of $\{\thisx\}_{k \in \N}$ with respect to $\hat X$.
    Let now $\optx$ be a weak accumulation point of $\{\thisx\}_{k \in \N}$, i.e., $x^{k_i} \weakto \optx$ for a subsequence $\{k_i\}_{i \in \N} \subset \N$. Since \eqref{eq:convergence:fb-est} implies $\nextx-\thisx \to 0$ by the assumption on the step lengths, we have $\grad F(x^{k_i+1})-\grad F(x^{k_i}) \to 0$ by the Lipschitz continuity of $\grad F$. Consequently, using again the subdifferential sum rule \cref{thm:subdiff:sum},  $\subdiff(G+F)(x^{k_i+1}) \ni w^{k_i+1} + \grad F(x^{k_i+1}) - \grad F(x^{k_i}) \to 0$. By the weak-to-strong outer semicontinuity of $\subdiff(G+F)$ from \cref{cor:monoton:closed,thm:monoton:subdiff}, it follows that $0 \in \subdiff (G+F)(\optx)$. We finish by applying Opial's \cref{lemma:opial} with $\hat X=\inv{(\subdiff(F+G))}(0)$.
\end{proof}

\begin{remark}
    The $L$-Lipschitz requirement on $\grad F$ can in some cases be restrictive, or upper estimates of $L$ difficult to obtain. In the latter case, line search can be used, as we will discuss in \cref{sec:meta:linesearch}.
    However, a slight modification of the forward-backward iteration can avoid the requirement entirely. Indeed, \cite{MalitskyTam:2020} prove the convergence of the  \term[method!forward-reflected-backward splitting]{forward-reflected-backward} iteration $\nextx \defeq \prox_{\tau G}(x^k - 2 \tau \grad F(x^k) + \tau \grad F(x^{k-1}))$, merely requiring $F$ to be $\tilde L$-Lipschitz, and the fixed step length $0 < \tau < 1/(2\tilde L)$.
\end{remark}

\section{Preconditioned proximal point methods: DRS and PDPS}
\label{sec:convergence:general-prox}

We now extend the analysis of the previous section to the \term[method!proximal point!preconditioned]{preconditioned proximal point method} \eqref{eq:proximal:preconditioned}, which we recall can be written in implicit form as
\begin{equation}
    \label{eq:convergence:preconditioned}
    0 \in \Hany(\nextx) + \Precond(\nextx-\thisx)
\end{equation}
for some preconditioning operator $M \in \linear(X; X)$ and includes the Douglas--Rachford splitting (DRS) and the primal-dual proximal splitting (PDPS) methods as special cases.
To deal with $M$, we need to improve \cref{thm:convergence:fb} slightly.
First, we introduce the preconditioned norm $\norm{x}_M \defeq \sqrt{\iprod{Mx}{x}}$, which satisfies the \term[identity!three-point!preconditioned]{preconditioned three-point identity}
\begin{equation}
    \label{eq:convergence:three-point-identity-cp}
    \iprod{M(x-y)}{x-z}
    = \frac{1}{2}\norm{x-y}_{M}^2
    - \frac{1}{2}\norm{y-z}_{M}^2
    + \frac{1}{2}\norm{x-z}_{M}^2\quad\text{for all }x,y,z\in X.
\end{equation}

The boundedness assumption in the statement of the next theorem holds in particular for $M=\Id$ and $\Hany$ maximally monotone by \cref{lem:proximal:lipschitz}.
\begin{theorem}
    \label{thm:convergence:precond}
    Suppose $\Hany: X \setto X$ is monotone and weak-to-strong outer semicontinuous with $\inv H(0) \ne \emptyset$, that $\Precond \in \linear(X; X)$ is self-adjoint and positive semi-definite, and that either $\Precond$ has a bounded inverse, or $\inv{(\Hany + \Precond)} \circ \Precond^{1/2}$ is bounded on bounded sets.
    Let the initial iterate $x^0 \in X$ be arbitrary, and assume that \eqref{eq:convergence:preconditioned} has a unique solution $\nextx$ for all $k\in\N$. Then the iterates $\{\thisx\}_{k \in \N}$ of \eqref{eq:convergence:preconditioned} are bounded and satisfy $0 \in \limsup_{k \to \infty} H(\thisx)$ and $M^{1/2}(\thisx-\realoptx) \weakto 0$ for some $\realoptx \in \inv H(0)$.
\end{theorem}

\begin{proof}
    Let $\realoptx \in \inv \Hany(0)$ be arbitrary.
    By the monotonicity of $\Hany$, we then have as before
    \begin{equation*}
        \iprod{\Hany(\nextx)}{\nextx-\realoptx}_X \ge 0,
    \end{equation*}
    which together with \eqref{eq:convergence:preconditioned} yields
    \begin{equation}
        \label{eq:convergence:precond-est1}
        \iprod{M(\nextx-\thisx)}{\nextx-\realoptx}_X \le 0.
    \end{equation}
    Applying the preconditioned three-point identity \eqref{eq:convergence:three-point-identity-cp} for $x=\nextx$, $y=\thisx$, and $z=\realoptx$ in \eqref{eq:convergence:precond-est1} shows that
    \begin{equation}
        \label{eq:convergence:precond-est2}
        \frac{1}{2}\norm{\nextx-\realoptx}_{M}^2
        +\frac{1}{2}\norm{\nextx-\thisx}_{M}^2
        \le \frac{1}{2}\norm{\thisx-\realoptx}_{M}^2,
    \end{equation}
    and summing \eqref{eq:convergence:precond-est2} over $k=0,\ldots,N-1$ yields
    \begin{equation}
        \label{eq:convergence:precond-est3}
        \frac{1}{2}\norm{x^N-\realoptx}_{M}^2
        + \sum_{k=0}^{N-1}
        \frac{1}{2}\norm{\nextx-\thisx}_{M}^2
        \le \frac{1}{2}\norm{x^0-\realoptx}_{M}^2.
    \end{equation}

    Let now $\thisz \defeq M^{1/2} \thisx$. Our objective is then to show $\thisz \weakto \optz$ for some $\optz \in \hat Z \defeq M^{1/2}\inv \Hany(0)$, which we do by using Opial's \cref{lemma:opial}.
    From \eqref{eq:convergence:precond-est2}, we obtain the necessary Fejér monotonicity of $\{\thisz\}_{k \in \N}$ with respect to the set $\hat Z$.
    It remains to verify that $\hat Z$ contains all weak accumulation points of $\{\thisz\}_{k \in \N}$.

    Let therefore $\optz$ be such an accumulation point, i.e., $z^{k_i} \weakto \optz$ for a subsequence $\{k_i\}_{i \in \N}$.
    We want to show that $\optz=M^{1/2}\optx$ for a weak accumulation point $\optx$ of $\{\thisx\}_{k \in \N}$.
    We proceed by first showing in two cases the boundedness of $\{\thisx\}_{k \in \N}$:
    \begin{enumerate}[label=(\roman*)]
        \item If $\Precond$ has a bounded inverse, then $\Precond \ge \theta I$ for some $\theta>0$, and thus the sequence $\{\thisx\}_{k \in \N}$ is bounded by \eqref{eq:convergence:precond-est3}.
        \item Otherwise, $\inv{(\Hany + \Precond)} \circ \Precond^{1/2}$ is bounded on bounded sets.
        Now \eqref{eq:convergence:precond-est3} only gives boundedness of $\{\thisz\}_{k \in \N}$. However, $\nextx \in \inv{(\Hany + \Precond)}(\Precond \thisx)=\inv{(\Hany + \Precond)}(\Precond^{1/2} \thisz)$, and $\{\thisz\}_{k \in \N}$ is bounded by \eqref{eq:convergence:precond-est3}, so we obtain the boundedness of $\{\thisx\}_{k \in \N}$.
    \end{enumerate}
    Thus there exists a further subsequence of $\{x^{k_i}\}_{i \in \N}$, weakly converging to some $\optx \in X$. Since $\thisz = M^{1/2} \thisx$, it follows that $\optz=M^{1/2}\optx$.
    To show that $\optz \in \hat Z$, if therefore suffices to show that the weak accumulation points of $\{\thisx\}_{k \in \N}$ belong to $\inv H(0)$.

    Let thus $\optx$ be any weak accumulation point of $\{\thisx\}_{k \in \N}$, i.e., $x^{k_i} \weakto \optx$ for some subsequence $\{k_i\}_{k \in \N} \subset \N$.
    From \eqref{eq:convergence:precond-est3}, we obtain first that $M^{1/2}(x^{k+1}-x^{k}) \to 0$ and hence that $w^{k+1} \defeq -M(x^{k+1}-x^{k}) \to 0$. From \eqref{eq:cp-prox}, we also know that $w^{k+1} \in \Hany(x^{k+1})$.
    It follows that $0 = \lim_{k \to \infty} \nexxt w \in \limsup_{k \to \infty} H(\nextx)$.
    The weak-to-strong outer semicontinuity now immediately yields $0 \in \Hany(\optx)$. Hence, $\hat Z$ contains all weak accumulation points of $\{\thisz\}_{k \in \N}$.

    The claim now follows from \cref{lemma:opial}.
\end{proof}

In the following, we verify that the DRS and PDPS methods satisfy the assumptions of this theorem.

\subsection*{Douglas--Rachford splitting}\label{ssec:convergence:drs}

Recall that the DRS method \cref{eq:splitting:dr}, i.e.,
\begin{algeqbox}
    \begin{equation}
        \label{eq:convergence:dr}
        \left\{\begin{aligned}
                x^{k+1} &\defeq \prox_{\tau F}(z^k),\\
                y^{k+1} &\defeq \prox_{\tau G}(2x^{k+1} - z^k),\\
                z^{k+1} &\defeq z^k + y^{k+1} - x^{k+1},
        \end{aligned}\right.
    \end{equation}
\end{algeqbox}
can be written as the preconditioned proximal point method \eqref{eq:proximal:preconditioned} in terms of $u=(x, y, z) \in \Space \defeq X^3$ and the operators
\begin{equation}
    \label{eq:convergence:drs-form}
    \Hany(x, y, z) \defeq
    \begin{pmatrix}
        \tau B(x) + y - z \\
        \tau A(y) + z - x \\
        x - y
    \end{pmatrix}
    \quad\text{and}\quad
    \Precond \defeq
    \begin{pmatrix}
        0 & 0 & 0 \\
        0 & 0 & 0 \\
        0 & 0 & I
    \end{pmatrix}
\end{equation}
for $B=\partial F$ and $A=\partial G$.
We are now interested in the properties of $H$ in terms of those of $A$ and $B$. For this, we can make use of the generic structure of $H$, which will reappear several times in the following.
\begin{lemma}
    \label{lemma:bcp-skew-adjoint}
    If $A: X \setto X$ is maximally monotone and $\Xi \in \linear(X; X)$ is skew-adjoint (i.e., $\Xi^*=-\Xi$), then $\Hsaddle\defeq A+\Xi$ is maximally monotone.
    In particular, any skew-adjoint operator $\Xi$ is maximally monotone.
\end{lemma}
\begin{proof}
    Let $x, z^*\in X$ be given such that
    \begin{equation*}
        \iprod{z^*-\tilde z^*}{x-\tilde x}_X \geq 0 \quad \text{for all }\tilde x\in X, \tilde z^*\in H(\tilde x).
    \end{equation*}
    Recalling \eqref{eq:monoton:max_char}, we need to show that $z^* \in H(x)$.
    By the definition of $H$, for any $\tilde z^*\in H(\tilde x)$ there exists a $\tilde x^*\in A(\tilde x)$ with $\tilde z^* = \tilde x^* + \Xi \tilde x$. On the other hand, setting $x^* \defeq z^* - \Xi x$, we have $z^* = x^* + \Xi x$. We are thus done if we can show that $x^*\in A(x)$. But using the skew-adjointness of $H$ and the symmetry of the inner product, we can write
    \begin{equation}
        \label{eq:bcp-skew-adjoint-iprod-expansion}
        \begin{aligned}[t]
            0&\leq \iprod{z^*-\tilde z^*}{x-\tilde x}_X\\
            &= \iprod{x^*-\tilde x^*}{x-\tilde x}_X + \iprod{\Xi(x-\tilde x)}{x-\tilde x}_X\\
            &= \iprod{x^*-\tilde x^*}{x-\tilde x}_X + \frac12 \iprod{\Xi(x-\tilde x)}{x-\tilde x}_X-\frac12\iprod{x-\tilde x}{\Xi(x-\tilde x)}_X\\
            &= \iprod{x^*-\tilde x^*}{x-\tilde x}_X,
        \end{aligned}
    \end{equation}
    and $x^*\in A(x)$ follows from the maximal monotonicity of $A$.

    To prove the final claim about skew-adjoint operators being maximally monotone, we take $A=\{0\}=\subdiff S$ for the constant functional $S\equiv 0$, which is maximally monotone by \cref{thm:monoton:subdiff}.
\end{proof}

\begin{corollary}\label{lemma:convergence:drs-h-bcp}
    Let $A$ and $B$ be maximally monotone.
    Then the operator $H$ defined in \eqref{eq:convergence:drs-form} is maximally monotone.
\end{corollary}
\begin{proof}
    Let
    \[
        \tilde A(u) \defeq
        \begin{pmatrix}
            \tau B(x) \\
            \tau A(y) \\
            0
        \end{pmatrix}
        \quad\text{and}\quad
        \Xi \defeq \begin{pmatrix} 0 & \Id & -\Id \\ -\Id & 0 & \Id \\ \Id & -\Id & 0 \end{pmatrix}.
    \]
    From the definition of the inner product on the product space $X^3$ together with \cref{lem:monoton:scalar_maxmon}, we have that $\tilde A$ is maximally monotone, while $\Xi$ is clearly skew-adjoint. The claim now follows from \cref{lemma:bcp-skew-adjoint}.
\end{proof}

We can now show convergence of the DRS method.
\begin{corollary}
    \label{cor:convergence:drs}
    Let $A,B:X\setto X$ be maximally monotone, and suppose $\inv{(A+B)}(0) \ne \emptyset$.
    Pick a step length $\tau>0$ and an initial iterate $z^0 \in X$. Then the iterates $\{(\thisx,\thisy,\thisz)\}_{k \in \N}$ of the DRS method \eqref{eq:convergence:dr} converge weakly to $(\realoptx, \realopty, \realopt{z})\in H^{-1}(0)$ satisfying $\realoptx = \realopty \in \inv{(A+B)}(0)$. Moreover, $\thisx-\thisy \to 0$.
\end{corollary}

\begin{proof}
    Since $A$ and $B$ are maximally monotone, \cref{lem:proximal:lipschitz} shows that the DRS iteration is always solvable for $\nextu$.
    Regarding convergence, we start by proving that the sequence $\{\thisu=(\thisx, \thisy, \thisz)\}_{k \in \N}$ is bounded, $\thisz \weakto \realoptz$ for some $\realoptz$, and $0 \in \limsup_{k \to \infty} H(\thisu)$. Note that the latter implies as claimed that $\thisx-\thisy \to 0$ strongly.
    We do this using \cref{thm:convergence:precond} whose conditions we have to verify.
    By \cref{lemma:convergence:drs-h-bcp}, $H$ is maximally monotone and hence weak-to-strong outer semicontinuous by \cref{cor:monoton:closed}.
    Since $\Precond$ is noninvertible, we also have to verify that $\inv{(\Hany+\Precond)} \circ \Precond^{1/2}$ is bounded on bounded sets.
    But since $\nextu \in (\Hany+\Precond)^{-1}(\Precond\thisu)=(\Hany+\Precond)^{-1}(\Precond^{1/2}\thisu)$ is an equivalent formulation of the iteration \eqref{eq:convergence:dr}, this follows from the Lipschitz continuity of the resolvent (\cref{lem:proximal:lipschitz}).
    Hence, we can apply \cref{thm:convergence:precond} to deduce $0 \in \limsup_{k \to \infty} H(\thisu)$ as well as $M^{1/2}(\thisu - \realoptu) \weakto 0$ for some $\realoptu=(\realoptx,\realopty,\realoptz)$ with $0 \in H(\realoptu)$.
    By the definition of $M$, this gives $\thisz\weakto \realoptz$.
    Moreover, the third line in the definition of $H$ implies that $\realoptx=\realopty$. Adding the first two lines in the same definition, we then obtain $0\in A(\realoptx)+ B(\realoptx)$.

    It remains to show weak convergence of the other variables. Since $\{\thisu\}_{k\in\N}$ is bounded, it contains a subsequence converging weakly to some $\tilde u = (\tilde x,\tilde y,\tilde z)$ which satisfies $0\in \Hany(\tilde x,\tilde y,\tilde z)$ such that $\tilde x = \tilde y$.
    Since $z^k\weakto \realoptz$, we have $\tilde z = \realoptz$. The first relation of the inclusion then can be rearranged to $\tilde y = \tilde x = \calR_{\tau B}(\tilde z) = \calR_{\tau B}(\realoptz)$ by the single-valuedness of the resolvent (\cref{lem:proximal:lipschitz}). The limit is thus independent of the subsequence, and hence a subsequence--subsequence argument shows that the full sequence converges.
\end{proof}

In particular, this convergence result applies to the special case of $B=\subdiff F$ and $A=\subdiff G$ for proper, convex, lower semicontinuous $F,G:X\to\Rbar$. However, the fixed point provided by the DRS method is related to a solution of the problem $\min_{x \in X}~F(x)+G(x)$ only if the subdifferential sum rule (\cref{thm:subdiff:sum}) holds with equality.

\subsection*{Primal-dual proximal splitting}\label{ssec:convergence:pdps}

To study the PDPS method, we recall from \eqref{eq:saddle-operator} and \eqref{eq:pdps-basic-m} the operators
\begin{equation}
    \label{eq:saddle-operator-repeat}
    \Hsaddle(u) \defeq \begin{pmatrix} \partial F(x) + K^*y \\ \partial G^*(y) - Kx \end{pmatrix},
    \quad\text{and}\quad
    M \defeq \begin{pmatrix} \tau^{-1}\Id &-K^* \\ -K & \sigma^{-1}\Id \end{pmatrix}
\end{equation}
for $u=(x,y)\in X\times Y=: U$.
With these we have already shown in \cref{sec:proximal:pd} that the PDPS method
\begin{algeqbox}
    \begin{equation}
        \label{eq:convergence:pd}
        \left\{\begin{aligned}
                x^{k+1} &\defeq \prox_{\tau F}(x^k - \tau K^*y^k),\\
                \overnextx &\defeq 2x^{k+1}-x^k,\\
                y^{k+1} &\defeq \prox_{\sigma G^*}(y^{k} + \sigma K\overnextx).
        \end{aligned}        \right.
    \end{equation}
\end{algeqbox}
has the form \eqref{eq:convergence:preconditioned} of the preconditioned proximal point method.
To show convergence, we first have to establish some basic properties of both $H$ and $M$.

\begin{lemma}\label{lemma:convergence:pd_spd}
    The operator $M:U\to U$ defined in \eqref{eq:saddle-operator-repeat} is bounded and self-adjoint. If $\sigma\tau\norm{K}_{\linear(X;Y)}^2 < 1$, then $M$ is positive definite.
\end{lemma}
\begin{proof}
    The definition of $M$ directly implies boundedness (since $K\in \linear(X;Y)$ is bounded) and self-adjointness. Let now $u=(x,y)\in \Space$ be given. Then
    \begin{equation}
        \label{eq:convergence:pd-m-estim}
        \begin{aligned}[t]
            \iprod{Mu}{u}_U &= \iprod{\tau^{-1}x - K^*y}{x}_X + \iprod{\sigma^{-1}y-Kx}{y}_Y\\
            &= \tau^{-1}\norm{x}_X^2 - 2 \iprod{x}{K^*y}_X + \sigma^{-1}\norm{y}_Y^2\\
            &\geq \tau^{-1}\norm{x}_X^2 - 2 \norm{K}_{\linear(X;Y)}\norm{x}_X\norm{y}_Y + \sigma^{-1}\norm{y}_Y^2\\
            &\geq \tau^{-1}\norm{x}_X^2 - \norm{K}_{\linear(X;Y)}\sqrt{\sigma\tau}(\tau^{-1}\norm{x}_X^2+\sigma^{-1}\norm{y}_Y^2)+ \sigma^{-1}\norm{y}_Y^2\\
            &= (1-\norm{K}_{\linear(X;Y)}\sqrt{\sigma\tau})(\tau^{-1}\norm{x}_X^2 + \sigma^{-1}\norm{y}_Y^2) \\
            &\geq C(\norm{x}_X^2 + \norm{y}_Y^2)
        \end{aligned}
    \end{equation}
    for $C\defeq(1-\norm{K}_{\linear(X;Y)}\sqrt{\sigma\tau})\min\{\tau^{-1},\sigma^{-1}\}>0$.
    Hence, $\iprod{Mu}{u}_U\geq C\norm{u}^2_U$ for all $u\in U$, and therefore $M$ is positive definite.
\end{proof}

\begin{lemma}\label{lemma:convergence:saddle-h-bcp}
    The operator $H:U\setto U$ defined in \eqref{eq:saddle-operator-repeat} is maximally monotone.
\end{lemma}
\begin{proof}
    Let $A(u) \defeq \begin{psmallmatrix} \subdiff F(x) \\ \subdiff G^*(y)\end{psmallmatrix}$ and $\Xi \defeq \begin{psmallmatrix} 0 & K^* \\ -K & 0 \end{psmallmatrix}$.
    Then $\Xi$ is skew-adjoint, and $A$ is maximally monotone by the definition of the inner product on $U=X\times Y$ and \cref{thm:monoton:subdiff}.
    The claim now follows from \cref{lemma:bcp-skew-adjoint}.
\end{proof}

With this, we can deduce the convergence of the PDPS method.

\begin{corollary}\label{thm:convergence:pd_conv}
    Let the convex, proper, and lower semicontinuous functions $F:X\to\Rbar$, $G:Y\to\Rbar$, and the linear operator $K\in \linear(X;Y)$ satisfy the assumptions of \cref{thm:convex:fenchel}. If, moreover, $\sigma\tau \norm{K}_{\linear(X;Y)}^2 <1$, then the sequence $\{u^k \defeq (x^k,y^k)\}_{k\in\N}$ generated by the PDPS method \eqref{eq:convergence:pd} for any initial iterate $u^0 \in X \times Y$ converges weakly in $\Space$ to a pair $\realoptu \defeq (\realoptx,\realopty) \in \inv \Hsaddle(0)$, i.e., satisfying \eqref{eq:splitting:fenchel}.
\end{corollary}

\begin{proof}
    By \cref{lemma:convergence:pd_spd}, $\Precond$ is self-adjoint and positive definite and thus has a bounded inverse. Keeping in mind \cref{lemma:convergence:saddle-h-bcp}, we can therefore apply \cref{thm:convergence:precond} to show that $(\thisu-\realoptu) \weakto 0$ for some $\realoptu \in \inv \Hsaddle(0)$ with respect to the inner product $\iprod{M\freevar}{\freevar}_U$.
    Since $\Precond$ has a bounded inverse, this implies that
    \begin{equation*}
        \iprod{\thisu}{Mw}_U = \iprod{M\thisu}{w}_U \to \iprod{M\realoptu}{w}_U =\iprod{\realoptu}{Mw}_U \quad\text{for all } w\in U
    \end{equation*}
    and hence $\thisu\weakto \realoptu$ in $U$ since $\range{\Precond}=U$ due to the invertibility of $\Precond$.
\end{proof}

\begin{remark}
    Through a general approach to degenerately preconditioned proximal point methods, i.e., singular $M$, the work \cite{Bredies:2022} proves the weak convergence of PDPS in the degenerate case $\tau\sigma\norm{K}^2=1$.
    Like our \cref{cor:convergence:drs}, their approach also readily establishes the convergence of $\{(\thisx, \thisy, \thisz)\}_{k \in \N}$ for the DRS; classical proofs only show the convergence of $\{\thisz\}_{k \in \N}$.
\end{remark}

\section{Preconditioned explicit splitting methods: PDES and more}
\label{sec:convergence:fb-general}

Let $A, B: X \setto X$ be monotone operators and consider the iterative scheme
\begin{equation}
    \label{eq:convergence:general:alg-monotone}
    0 \in A(\nextx)+B(\thisx) + M(\nextx-\thisx),
\end{equation}
which is implicit in $A$ but explicit in $B$.
We obviously intend to use this method to find some $\realoptx \in \inv{(A+B)}(0)$.

As we have seen, the proximal point, PDPS, and DRS methods are all of the form \eqref{eq:convergence:general:alg-monotone} with $B=0$.
The basic explicit splitting method is also of this form with $A=\subdiff G$, $B=\grad F$, and $M=\inv\tau \Id$. It is moreover not difficult to see from \eqref{eq:splitting:gist-implicit} that the primal-dual explicit splitting (PDES) method is also of the form \eqref{eq:convergence:general:alg-monotone} with nonzero $B$. So to prove the convergence of this algorithm, we want to improve \cref{thm:convergence:fb} to be able to deal with the preconditioning operator $M$ and the general monotone operators $A$ and $B$ in place of subdifferentials and gradients.

To proceed, we need a suitable notion of smoothness for $B$ to be able to deal with the explicit step.
In \cref{thm:convergence:fb} we only used the Lipschitz continuity of $\grad F$ in two places: first, to establish the three-point monotonicity using \cref{cor:smoothness:three-point}, and second, at the end of the proof for a continuity argument. To simplify dealing with $B$ that may only act on a subspace, as in the case of the primal-dual explicit splitting in \cref{sec:proximal:gist}, we now make this three-point monotonicity with respect to an operator $\Lambda$ our main assumption.

Specifically, we say that $B: X \setto X$ is \term[mapping!monotone!three-point]{three-point monotone} at $\realoptx \in X$ with respect to $\Lambda \in \linear(X; X)$ if
\begin{equation}
    \label{eq:convergence:3monotone}
    \iprod{B(z)-B(\realoptx)}{x-\realoptx}_X \ge -\frac{1}{4}\norm{z-x}^2_\Lambda \quad \text{for all }x,z\in X.
\end{equation}
If this holds for every $\realoptx$, we say that $B$ is \emph{three-point monotone with respect to $\Lambda$}.
From \cref{cor:smoothness:three-point}, it is clear that if $\grad F$ is Lipschitz continuous with constant $L$, then $B=\grad F$ is three-point monotone with respect to $\Lambda=L\, \Id$.

We again start with a lemma exploiting the structural properties of the saddle-point operator to show a \enquote{shifted outer semicontinuity}.

\begin{lemma}
    \label{lemma:convergence:bcp-lipschitz-perturbation}
    Let $H=A+B: X \setto X$ be weak-to-strong outer semicontinuous with $B$ single-valued and Lipschitz continuous.
    If $\nextw \in A(\nextx)+B(\thisz)$ for $k \in \N$ with $\thisw \to \optw$ and $\nextx-\thisz \to 0$ strongly in $X$ and $\thisx \weakto \optx$ weakly in $X$, then $\optw \in H(\optx)$.
\end{lemma}

\begin{proof}
    We have $\nextw \in A(\nextx)+B(\thisz)$  so that
    \[
        \nexxt{\tilde w}
        \defeq
        \nextw - B(\thisz)+B(\nextx)
        \in H(\nextx).
    \]
    Since $\nextw \to \optw$ and $\nextx-\thisz \to 0$ and $B$ is Lipschitz continuous, we have $\nexxt{\tilde w} \to \optw$ as well. The weak-to-strong outer semicontinuity of $H$ then immediately yields $\optw\in H(\optx)$.
\end{proof}

\begin{theorem}
    \label{theorem:convergence:abxi}
    Let $H=A+B$ with $\inv H(0) \ne \emptyset$ for $A, B: X \setto X$ with $A$ monotone and $B$ single-valued Lipschitz continuous and three-point monotone with respect to some $\Lambda \in \linear(X; X)$.
    Furthermore, let $\Precond \in \linear(X; X)$ be self-adjoint, positive definite, with a bounded inverse, and satisfy $(2-\epsilon)M \ge \Lambda$ for some $\epsilon>0$.
    Suppose $H$ is weak-to-strong outer semicontinuous.
    Let the starting point $x^0 \in X$ be arbitrary and assume that \eqref{eq:convergence:general:alg-monotone} has a unique solution $\nextx$ for every $k\in \N$. Then the iterates $\{\thisx\}_{k \in \N}$ of \eqref{eq:convergence:general:alg-monotone} satisfy $M^{1/2}(\thisx-\realoptx) \weakto 0$ for some $\realoptx \in \inv H(0)$.
\end{theorem}
\begin{proof}
    The proof follows along the same lines as that of \cref{thm:convergence:precond} with minor modifications.
    First, since $0 \in H(\realoptx)$, the monotonicity of $A$ and the three-point monotonicity \eqref{eq:convergence:3monotone} of $B$ yields
    \begin{equation*}
        \iprod{A(\nextx)+B(\thisx)}{\nextx-\realoptx}_X
        \ge -\frac{1}{4}\norm{\nextx-\thisx}_\Lambda^2,
    \end{equation*}
    which together with \eqref{eq:convergence:general:alg-monotone} leads to
    \begin{equation*}
        \iprod{M(\nextx-\thisx)}{\nextx-\realoptx}_X \le \frac{1}{4}\norm{\nextx-\thisx}_\Lambda^2.
    \end{equation*}
    From the preconditioned three-point identity \eqref{eq:convergence:three-point-identity-cp} we then obtain
    \begin{equation}
        \label{eq:convergence:precond-est2-abxi}
        \frac{1}{2}\norm{\nextx-\realoptx}_{M}^2
        +\frac{1}{2}\norm{\nextx-\thisx}_{M-\Lambda/2}^2
        \le \frac{1}{2}\norm{\thisx-\realoptx}_{M}^2.
    \end{equation}
    Our assumption that $(2-\epsilon)M \ge \Lambda$ implies that $M - \Lambda/2 \ge \epsilon M/2$. By definition, we can therefore bound the second norm on the left-hand side from below to obtain \eqref{eq:convergence:precond-est2} with an additional constant depending on $\epsilon$. We may thus proceed as in the proof of \cref{thm:convergence:precond} to establish $w^{k+1} \defeq -M(x^{k+1}-x^{k}) \to 0$. We now have $\nextw \in A(\nextx)+B(\thisx)$ and therefore use \cref{lemma:convergence:bcp-lipschitz-perturbation} with $\thisz=\thisx$ and $\optw=0$ to establish $0 \in H(\optx)$. The rest of the proof again proceeds as for \cref{thm:convergence:precond} with the application of Opial's \cref{lemma:opial}.
\end{proof}

We again apply this result to show the convergence of specific splitting methods containing an explicit step.

\subsection*{Primal-dual explicit splitting}

We now return to algorithms for problems of the form
\[
    \min_{x \in X} F(x)+G(Kx)
\]
for Gateaux differentiable $F$ and linear $K$. Recall from \eqref{eq:splitting:gist} the \index{method!primal-dual explicit splitting}primal-dual explicit splitting (PDES) method
\begin{algeqbox*}
    \begin{equation*}
        \label{eq:convergence:gist}
        \left\{\begin{aligned}
                \nexty &\defeq \prox_{G^*}((\Id-KK^*)\thisy + K(\thisx- \grad F(\thisx))),\\
                \nextx & \defeq \thisx - \grad F(\thisx) - K^*\nexty,\\
        \end{aligned}        \right.
    \end{equation*}
\end{algeqbox*}
which can be written in implicit form as
\begin{equation}
    \label{eq:convergence:gist-implicit}
    0 \in H(\nextu) + \begin{pmatrix} \grad F(\thisx)-\grad F(\nextx) \\ 0 \end{pmatrix}
    + \Precond(\nextu-\thisu)
\end{equation}
with
\begin{equation}
    \label{eq:convergence:gist-h-m}
    \Hsaddle(u) \defeq \begin{pmatrix} \partial F(x) + K^*y \\ \partial G^*(y) - Kx \end{pmatrix}
    \quad\text{and}\quad
    M \defeq \begin{pmatrix} \Id & 0 \\ 0 & \Id-KK^* \end{pmatrix},
\end{equation}
for $u=(x, y)\in X\times Y=:U$.

\begin{corollary}
    \label{thm:gist}
    Let $F:X\to\R$ and $G:Y\to\Rbar$ be proper, convex, and lower semicontinuous, and $K \in \linear(X; Y)$.
    Suppose $F$ is Gateaux differentiable with $L$-Lipschitz gradient for $L<2$, that $\norm{K}_{\linear(X; Y)} < 1$, and that the assumptions of \cref{thm:convex:fenchel} are satisfied.
    Then for any initial iterate $u^0 \in X \times Y$ the iterates $\{\thisu=(\thisx,\thisy)\}_{k \in \N}$ of the \eqref{eq:splitting:gist} converge weakly to some $\realoptu \in \inv H(0)$ with $H$ given by \eqref{eq:saddle-operator}.
\end{corollary}

\begin{proof}
    We recall that \cref{thm:convex:fenchel} guarantees that $\inv H(0) \ne \emptyset$.
    To apply \cref{theorem:convergence:abxi}, we write $H=A+B$ for
    \[
        A(u) \defeq \begin{pmatrix} 0 \\ \subdiff G^*(y) \end{pmatrix} + \Xi u,
        \quad
        B(u) \defeq \begin{pmatrix} \grad F(x) \\ 0 \end{pmatrix},
        \quad
        \Xi \defeq \begin{pmatrix} 0 & K^* \\ -K & 0 \end{pmatrix}.
    \]
    We first note that $M$ as given in \eqref{eq:convergence:gist-h-m} is self-adjoint and positive definite under our assumption $\norm{K}_{\linear(X;Y)} < 1$.
    By \cref{cor:smoothness:three-point}, the three-point monotonicity \eqref{eq:convergence:3monotone} holds for
    $\Lambda \defeq \begin{psmallmatrix} L & 0 \\ 0 & 0 \end{psmallmatrix}$.
    Since $L<2$, there furthermore exists an $\epsilon>0$ sufficiently small such that $(2-\epsilon)M \ge \Lambda$.
    Finally, \cref{lemma:convergence:saddle-h-bcp} shows that $H$ is maximally monotone and hence weak-to-strong outer semicontinuous by \cref{cor:monoton:closed}.
    The claim now follows from \cref{theorem:convergence:abxi}.
\end{proof}

\begin{remark}
    It is possible to improve the result to $\norm{K}_{\linear(X; Y)} \le 1$ if we increase the complexity of \cref{theorem:convergence:abxi} slightly to allow for $M\geq 0$. However, in this case it is only possible to show the convergence of the partial iterates $\{\thisx\}_{k \in \N}$.
\end{remark}

\subsection*{Primal-dual proximal splitting with an additional forward step}

Using a similar switching term as in the implicit formulation \eqref{eq:convergence:gist-implicit} of the PDES method, it is possible to incorporate additional forward steps in the PDPS method.
For $F=F_0+E$ with $F_0,E$ convex and $E$ Gateaux differentiable, we therefore consider
\begin{equation*}
    \label{eq:convergence:pdps:problem}
    \min_{x \in X} F_0(x)+E(x)+G(Kx).
\end{equation*}
With $u=(x, y)$ and following \cref{sec:proximal:pd}, any minimizer $\realoptx\in X$ satisfies $0 \in \Hsaddle(\realopt u)$ for
\begin{equation}
    \label{eq:convergence:pdps:h}
    \Hsaddle(u) \defeq \begin{pmatrix} \partial F_0(x)+\grad E(x) + K^*y \\ \partial G^*(y) - Kx \end{pmatrix}.
\end{equation}
Similarly, following the arguments in \cref{sec:proximal:pd}, we can show that the iteration
\begin{algeqbox}
    \begin{equation}\label{eq:convergence:pdps:forward}
        \left\{\begin{aligned}
                x^{k+1} &\defeq \prox_{\tau F_0}(x^k - \tau \grad E(\thisx) - \tau K^*y^k),\\
                \overnextx &\defeq 2x^{k+1}-x^k,\\
                y^{k+1} &\defeq \prox_{\sigma G^*}(y^{k} + \sigma K\bar x^{k+1}),
        \end{aligned}        \right.
    \end{equation}
\end{algeqbox}
is equivalent to the implicit formulation
\begin{equation*}
    0 \in
    \begin{pmatrix}
        \subdiff F_0(\nextx) + \grad E(\thisx) + K^*\nexty \\
        \subdiff G(\nexty) - K\nextx
    \end{pmatrix}
    +
    M(\nextu-\thisu)
\end{equation*}
with the preconditioner $M$ defined as in \eqref{eq:saddle-operator-repeat}. The convergence can thus be shown as for the PDES method.

\begin{corollary}
    \label{cor:convergence:pdps:forward}
    Let $E:X\to\R$, $F_0:X\to\Rbar$, and $G: Y\to\Rbar$ be proper, convex, and lower semicontinuous, and $K \in \linear(X; Y)$. Suppose $E$ is Gateaux differentiable with an $L$-Lipschitz gradient, and that the assumptions of \cref{thm:convex:fenchel} are satisfied with $F \defeq F_0+E$.
    Assume, moreover, that $\tau,\sigma>0$ satisfy
    \begin{equation}
        \label{eq:convergence:pdps:forward-stepsize}
        1 > \norm{K}_{\linear(X; Y)}^2\tau\sigma + \tau \frac{L}{2}.
    \end{equation}
    Then for any initial iterate $u^0 \in X \times Y$ the iterates $\{\thisu\}_{k \in \N}$ of \eqref{eq:convergence:pdps:forward} converge weakly to some $\realoptu \in \inv H(0)$ for $H$ given by \eqref{eq:convergence:pdps:h}.
\end{corollary}

\begin{proof}
    As before, \cref{thm:convex:fenchel} guarantees that $\inv H(0) \ne \emptyset$.
    We apply \cref{theorem:convergence:abxi} to
    \[
        A(u) \defeq \begin{pmatrix} \subdiff F_0(x) \\ \subdiff G^*(y) \end{pmatrix}+\Xi u,
        \quad
        B(u) \defeq \begin{pmatrix} \grad E(x) \\ 0 \end{pmatrix},
        \quad
        \Xi \defeq \begin{pmatrix} 0 & K^* \\ -K & 0 \end{pmatrix},
    \]
    and $M$ given by \eqref{eq:saddle-operator-repeat}.
    By \cref{cor:smoothness:three-point}, the three-point monotonicity \eqref{eq:convergence:3monotone} holds with
    $\Lambda \defeq \begin{psmallmatrix} L & 0 \\ 0 & 0 \end{psmallmatrix}$.
    We have already shown in \cref{lemma:convergence:pd_spd} that $M$ is self-adjoint and positive definite.
    Furthermore, from \eqref{eq:convergence:pd-m-estim} in the proof of \cref{lemma:convergence:pd_spd}, we have
    \begin{equation*}
        \iprod{Mu}{u}_U
        \ge (1-\norm{K}_{\linear(X;Y)}\sqrt{\sigma\tau})(\tau^{-1}\norm{x}_X^2 + \sigma^{-1}\norm{y}_Y^2).
    \end{equation*}
    Thus \eqref{eq:convergence:pdps:forward-stepsize} implies that $M$ is positive definite. Arguing similarly to \eqref{eq:convergence:pd-m-estim}, we also estimate
    \begin{equation*}
        \iprod{Mu}{u}_U
        \ge
        \tau^{-1}\norm{x}_X^2 - 2 \norm{K}_{\linear(X;Y)}\norm{x}_X\norm{y}_Y + \sigma^{-1}\norm{y}_Y^2
        \ge (1 - \norm{K}_{\linear(X;Y)}^2\sigma\tau)\tau^{-1}\norm{x}_X^2.
    \end{equation*}
    By the strict inequality in \eqref{eq:convergence:pdps:forward-stepsize}, we thus deduce $(2-\epsilon)M \ge \Lambda$ for some $\epsilon>0$.

    Now by \cref{lemma:convergence:saddle-h-bcp}, $H$ is again maximally monotone and therefore weak-to-strong outer semicontinuous by \cref{cor:monoton:closed}, and the claim follows from \cref{theorem:convergence:abxi}.
\end{proof}

\begin{remark}
    The forward step was introduced to the basic PDPS method in \cite{condat2013primaldual,vu2013splitting}, see also \cite{chambolle2014ergodic}. These papers also introduced an additional over-relaxation step that we will discuss in \cref{chap:meta}.
\end{remark}

\section{Fixed-point theorems}

Based on our generic approach, we now prove the classical \emph{Browder fixed-point theorem}, which can itself be used to prove the convergence of optimization methods and other fixed-point iterations (see \cref{remark:convergence:browder}).
We begin with a useful lemma.
\begin{lemma}\label{lem:convergence:averaged-fixed-points}
    Let $X$ be a Hilbert space and let $T: X \to X$ be $\alpha$-averaged for some $\alpha \in (0, 1)$.
    Assume that there exists a fixed point $\realoptx\in X$ of $T$.
    Then the set of fixed points is convex and closed.
\end{lemma}
\begin{proof}
    Let $T=(1-\alpha)\Id + \alpha J$ for some nonexpansive operator $J:X\to X$. Then $x\in X$ is a fixed point of $T$ if and only if $x$ is a fixed point of $J$. Hence
    \begin{equation*}
        \{\optx \mid \optx=T(\optx)\}=
        \{\optx \mid \optx=J(\optx)\}=\inv{(\Id-J)}(0)
    \end{equation*}
    so it suffices to show that the set on the right-hand side is convex and closed. But this follows from \cref{lemma:monotone:convex,cor:monoton:closed} since if $J$ is nonexpansive, then $\Id -J$ is maximally monotone by \cref{lem:monotone:nonexpansive} and hence so is $\inv{(\Id -J)}$ by \cref{lemma:monotone:inverse}.
\end{proof}

We recall from \cref{lemma:proximal:averaged} that firmly nonexpansive maps are $(1/2)$-averaged, so the following result applies by \cref{lem:proximal:firmly-nonexpansive} to the resolvents of maximally monotone maps in particular -- hence proving the convergence of the proximal point method.
\begin{theorem}[Browder]\index{theorem!Browder fixed-point}
    \label{thm:convergence:browder}
    Let $X$ be a Hilbert space and let $T: X \to X$ be $\alpha$-averaged for some $\alpha \in (0, 1)$.
    Assume that there exists a fixed point $\realoptx\in X$ of $T$.
    Let $\nextx \defeq T(\thisx)$ for $k\in \N$ and $x^0\in X$.
    Then $\thisx \weakto \optx$ weakly in $X$ for some fixed point $\optx$ of $T$.
\end{theorem}

\begin{proof}
    Finding a fixed point of $T$ is equivalent to finding a root of $H(x) \defeq T(x)-x$.
    Similarly, we can rewrite the fixed-point iteration as solving for $\nextx$ the inclusion
    \begin{equation}
        \label{eq:convergence:browder:implicit}
        0 = \thisx-T(\thisx) + (\nextx-\thisx).
    \end{equation}
    Proceeding as in the previous sections, we test this by the application of $\iprod{\freevar}{\nextx-\realoptx}_X$. After application of the three-point identity \eqref{eq:convergence:three-point-identity}, we then obtain
    \begin{equation}
        \label{eq:convergence:browder:est0}
        \frac{1}{2}\norm{\nextx-\realoptx}_X^2
        +\frac{1}{2}\norm{\nextx-\thisx}_X^2
        +\iprod{\thisx-T(\thisx)}{\nextx-\realoptx}_X
        \le
        \frac{1}{2}\norm{\thisx-\realoptx}_X^2.
    \end{equation}

    Since $\nextx = T(\thisx)$, $\realoptx$ is a fixed point of $T$, and by assumption $T=(1-\alpha)\Id + \alpha J$ for some nonexpansive operator $J: X \to X$, we have
    \[
        \begin{aligned}
            \iprod{\thisx-T(\thisx)}{\nextx-\realoptx}_X
            &
            =
            \iprod{\thisx-\realoptx-(T(\thisx)-T(\realoptx))}{T(\thisx)-T(\realoptx)}_X
            \\
            &
            =
            \alpha\iprod{\thisx-\realoptx-(J(\thisx)-J(\realoptx))}{(1-\alpha)(\thisx-\realoptx)+\alpha(J(\thisx)-J(\realoptx))}_X
            \\
            &
            =(\alpha-\alpha^2)\norm{\thisx-\realoptx}_X^2
            -\alpha^2\norm{J(\thisx)-J(\realoptx)}_X^2
            \\ \MoveEqLeft[-1]
            +(2\alpha^2-\alpha)\iprod{\thisx-\realoptx}{J(\thisx)-J(\realoptx)}_X
        \end{aligned}
    \]
    as well as
    \[
        \begin{aligned}[t]
        \frac{1}{2}\norm{\nextx-\thisx}_X^2
        &
        =
        \frac{1}{2}\norm{T(\thisx)-\thisx}_X^2
        =
        \frac{\alpha^2}{2}\norm{J(\thisx)-\thisx}_X^2
        =
        \frac{\alpha^2}{2}\norm{J(\thisx)-J(\realoptx)-(\thisx-\realoptx)}_X^2
        \\
        &
        =\frac{\alpha^2}{2}\norm{\thisx-\realoptx}_X^2
        +\frac{\alpha^2}{2}\norm{J(\thisx)-J(\realoptx)}_X^2
        -\alpha^2\iprod{\thisx-\realoptx}{J(\thisx)-J(\realoptx)}_X.
        \end{aligned}
    \]
    Thus, for any $\delta>0$,
    \begin{multline*}
        \frac{1-\delta}{2}\norm{\nextx-\thisx}_X^2
        +\iprod{\thisx-T(\thisx)}{\nextx-\realoptx}_X
        =
        ((1+\delta)\alpha^2-\alpha)\iprod{\thisx-\realoptx}{J(\thisx)-J(\realoptx)}_X
        \\
        +\frac{2\alpha-(1+\delta)\alpha^2}{2}\norm{\thisx-\realoptx}_X^2
        -\frac{(1+\delta)\alpha^2}{2}\norm{J(\thisx)-J(\realoptx)}_X^2.
    \end{multline*}
    Taking $\delta=\tfrac{1}{\alpha}-1$, we have $\delta>0$ and
    $\alpha = (1+\delta)\alpha^2$.
    Thus the factor in front of the inner product term is positive, and hence we obtain by the nonexpansivity of $J$
    \[
        \frac{1-\delta}{2}\norm{\nextx-\thisx}_X^2
        +\iprod{\thisx-T(\thisx)}{\nextx-\realoptx}_X
        =
        \frac{\alpha}{2}\norm{\thisx-\realoptx}_X^2
        -\frac{\alpha}{2}\norm{J(\thisx)-J(\realoptx)}_X^2
        \ge
        0.
    \]

    From \eqref{eq:convergence:browder:est0}, it now follows that
    \begin{equation*}
        \frac{1}{2}\norm{\nextx-\realoptx}_X^2
        +\frac{\delta}{2}\norm{\nextx-\thisx}_X^2
        \le
        \frac{1}{2}\norm{\thisx-\realoptx}_X^2.
    \end{equation*}
    As before, this implies Fejér monotonicity of $\{\thisx\}_{k\in\N}$ and that $\norm{\nextx-\thisx}_X \to 0$. The latter implies $\norm{T(\thisx)-\thisx}_X \to 0$ via \eqref{eq:convergence:browder:implicit}.
    Let $\optx$ be any weak accumulation point of $\{\thisx\}_{k \in \N}$.
    Denote by $N \subset \N$ the indices of the corresponding subsequence.
    We show that $\optx$ is a fixed point of $T$.
    Since by \cref{lem:convergence:averaged-fixed-points} the set of fixed points is convex and closed, the claim then follows from Opial's \cref{lemma:opial}.

    To show that $\optx$ is a fixed point of $T$, first, we expand
    \[
        \frac{1}{2}\norm{\thisx-T(\optx)}_X^2
        =\frac{1}{2}\norm{\thisx-\optx}_X^2
        +\frac{1}{2}\norm{\optx-T(\optx)}_X^2
        +\iprod{\thisx-\optx}{\optx-T(\optx)}_X.
    \]
    Since $\thisx \weakto \optx$, this gives
    \[
        \limsup_{N \ni k \to \infty} \frac{1}{2}\norm{\thisx-T(\optx)}_X^2
        \ge
        \limsup_{N \ni k \to \infty} \frac{1}{2}\norm{\thisx-\optx}_X^2
        + \frac12\norm{\optx-T(\optx)}_X^2.
    \]
    On the other hand, by the nonexpansivity of $T$ and $T(\thisx)-\thisx \to 0$, we have
    \[
        \begin{aligned}
            \limsup_{N \ni k \to \infty} \norm{\thisx-T(\optx)}_X^2
            &\le
            \limsup_{N \ni k \to \infty}\left( \norm{T(\thisx)-T(\optx)}_X + \norm{\thisx-T(\thisx)}_X\right)^2\\
            &\le \limsup_{N \ni k \to \infty} \norm{\thisx-\optx}_X^2.
        \end{aligned}
    \]
    Together these two inequalities show that $\norm{T(\optx)-\optx}_X=0$, from which the claim follows.
\end{proof}

\begin{remark}
    \Cref{thm:convergence:browder} in its modern form (stated for firmly nonexpansive or more generally $\alpha$-averaged maps) can be first found in \cite{browder1967convergence}.
    However, similar results for what are now called \term[method!Krasnoselskii--Mann]{Krasnoselskii--Mann iterations} -- which are closely related to $\alpha$-averaged maps -- were stated in more limited settings in \cite{Mann1953mean,schaefer1957,petryshyn1966construction,krasnoselskii1955remarks,opial1967weak}.
    Our overall approach in this book, based on \cite{tuomov-proxtest}, is an “implicit” counterpart to the more classical fixed point theorems. Instead of considering explicit iterations $\nextx \defeq T(\thisx)$, the theory is based on $\nextx$ defined implicitly through equations $0 = H(\nextx) + (\nextx-\thisx)$.
\end{remark}

\chapter{Splitting methods: rates of convergence}
\label{chap:testing}

As we have seen, minimizers of convex problems in a Hilbert space $X$ -- which again is the standard setting in this chapter -- can generally be characterized by the inclusion
\begin{equation*}
    0 \in H(\realoptx)
\end{equation*}
for the unknown $\realoptx \in X$ and a suitable monotone operator $H: X \setto X$. This inclusion in turn can be solved using a (preconditioned) proximal point iteration that converges weakly under suitable assumptions.
In the present chapter, we want to improve this analysis to obtain \emph{convergence rates}, i.e., estimates of the distance $\norm{\thisx-\realoptx}_X$ of iterates to $\realoptx$ in terms of the iteration number $k$. Our general approach will be to consider this distance multiplied by an iteration-dependent \emph{testing parameter} $\tauTest_k$ (or, for structured algorithms, consider the norm relative to a \emph{testing operator}) and to show by roughly the same arguments as in \cref{chap:convergence} that this product stays bounded: $\tauTest_k \norm{\thisx-\realoptx}_X \le C$. If we can then show that this testing parameter grows at a certain rate, the distance must decay at the reciprocal rate. Consequently, we can now avoid the complications of dealing with weak convergence; in fact, this chapter will consist of simple algebraic manipulations. However, for this to work we need to assume additional properties of $H$, namely strong monotonicity. Recall from \cref{thm:smoothness:strong-convexity} that $H$ is called strongly monotone with factor $\gamma>0$ if
\begin{equation}\label{eq:testing:def-strong-monotone}
    \iprod{H(\tilde x)-H(x)}{\tilde x-x}_X\geq \gamma\norm{\tilde x-x}_X^2 \quad (\tilde x,x\in X),
\end{equation}
where, in a slight abuse of notation, the left-hand side is understood to stand for any choice of elements from $H(\tilde x)$ and $H(x)$.

Before we turn to the actual estimates, we first define various notions of convergence rates.
Consider a function $r: \N \to [0, \infty)$ (e.g., $r(k)=\norm{\thisx-\realoptx}_X$ or $r(k)=G(\thisx)-G(\realoptx)$ for $\realoptx$ a minimizer of $G$).
\begin{enumerate}
    \item We say that $r(k)$ \term[convergence!rate]{converges (to zero as $k \to \infty$) at the rate} $O(f(k))$ if $r(k) \le C f(k)$ for some constant $C > 0$ for all $k \in \N$ and a decreasing function $f: \N \to [0, \infty)$  with $\lim_{k \to \infty} f(k)=0$ (e.g., $f(k) = 1/k$ or $f(k) = 1/k^2$).

    \item Analogously, we say that a function $R: \N \to [0, \infty)$ \term[rate of growth]{grows at the rate $\Omega(F(k))$} if $R(k) \ge c F(k)$ for all $k \in \N$ for some constant $c>0$ and an increasing function $F: \N \to [0, \infty)$ with $\lim_{k \to \infty} F(k)=\infty$.
\end{enumerate}
Clearly $r=1/R$ converges to zero at the rate $f=1/F$ if and only if $R$ grows at the rate $F$.
The most common cases are $F(k) = k$ or $F(k) = k^2$.

We can alternatively characterize \term[convergence!order]{orders} of convergence via
\[
    \mu \defeq \lim_{k \to \infty} \frac{r(k+1)}{r(k)}.
\]
\begin{enumerate}
    \item If $\mu = 1$, we say that $r(k)$ converges (to zero as $k \to \infty$) \term[convergence!sublinear]{sublinearly}.

    \item If $\mu \in (0, 1)$, then this convergence is \term[convergence!linear]{linear}. This is equivalent to a convergence at the rate $O(\tilde\mu^k)$ for any $\tilde\mu \in (\mu, 1)$.

    \item If $\mu=0$, then the convergence is \term[convergence!superlinear]{superlinear}.
\end{enumerate}
Different rates of superlinear convergence can also be studied. We say that $r(k)$ converges (to zero as $k \to \infty$) \term[convergence!superlinear with order]{superlinearly with order $q>1$} if
\[
    \lim_{k \to \infty} \frac{r(k+1)}{r(k)^q} < \infty.
\]
The most common case is $q=2$, which is also known as \term[convergence!quadratic]{quadratic convergence}. (This is not to be confused with the -- much slower -- convergence at the rate $O(1/k^2)$; similarly, linear convergence is different from -- and much faster than -- convergence at the rate $O(1/k)$.)

\section{The fundamental methods}
\label{sec:testing:fundamental}

Before going into this abstract operator-based theory, we demonstrate the general concept of testing by studying the fundamental methods, the proximal point and explicit splitting methods. These are purely primal methods with a single step length parameter, which simplifies the testing approach since we only need a single testing parameter.
(It should be pointed out that the proofs in this section can be carried out -- and in fact shortened -- without introducing testing parameters at all. Nevertheless, we follow this approach since it provides a blueprint for the proofs for the structured primal-dual methods where these are required.)

\subsection*{Proximal point method}

We start with the basic proximal point method for solving $0 \in H(\realoptx)$ for a monotone operator $H:X\setto X$, which we recall can be written in implicit form as
\begin{equation}
    \label{eq:testing:prox-repeat}
    0 \in \tau_k\Hany(\nextx)+(\nextx-\thisx).
\end{equation}

\begin{theorem}[proximal point method iterate rates]
    \label{thm:testing:prox}
    Suppose $H: X \setto X$ is strongly monotone with $\inv H(0)\neq \emptyset$.
    Let $\nextx \defeq \calR_{\tau_k H}(\thisx)$ for some $\{\tau_k\}_{k\in\N}\subset (0,\infty)$ and $x^0\in X$ be arbitrary. Then the following hold for the iterates $\{\thisx\}_{k \in \N}$ and the unique point $\realoptx \in \inv H(0)$:
    \begin{enumerate}
        \item\label{item:testing:prox:linear}
            If $\tau_k\equiv \tau$ is constant, then $\norm{x^k-\realoptx}_X \to 0$ linearly.
        \item\label{item:testing:prox:superlinear}
            If $\tau_k\upto \infty$, then $\norm{x^k-\realoptx}_X \to 0$ superlinearly.
    \end{enumerate}
\end{theorem}

\begin{proof}
    Let $\realoptx \in \inv \Hany(0)$; by assumption, such a point exists and is unique due to the assumed strong monotonicity of $H$ (since inserting any two roots $\hat x,\tilde x\in X$ of $H$ in \eqref{eq:testing:def-strong-monotone} yields $\norm{\hat x-\tilde x}_X\leq 0$).
    For each iteration $k \in \N$, pick a \term[parameter!testing]{testing parameter} $\tauTest_k>0$ and apply the \term[functional!testing]{test} $\tauTest_k\iprod{\freevar}{\nextx-\realoptx}_X$ to \eqref{eq:testing:prox-repeat} to obtain (using the same notation from \cref{thm:convergence:prox})
    \begin{equation}
        \label{eq:testing:prox:tested}
        0 \in \tauTest_k\tau_k\iprod{\Hany(\nextx)}{\nextx-\realoptx}_X+\tauTest_k\iprod{\nextx-\thisx}{\nextx-\realoptx}_X.
    \end{equation}
    By the strong monotonicity of $\Hany$, and the fact that $0 \in \Hany(\realoptx)$, for some $\gamma>0$,
    \begin{equation*}
        \iprod{\Hany(\nextx)}{\nextx-\realoptx}_X \ge \gamma\norm{\nextx-\realoptx}_X^2
    \end{equation*}
    Multiplying this inequality with $\phi_k\tau_k$ and using \eqref{eq:testing:prox:tested}, we obtain
    \begin{equation*}
        \tauTest_k\tau_k\gamma\norm{\nextx-\realoptx}_X^2 + \tauTest_k \iprod{\nextx-\thisx}{\nextx-\realoptx}_X \le 0.
    \end{equation*}
    An application of the three-point identity \eqref{eq:convergence:three-point-identity} then yields
    \begin{equation}
        \label{eq:testing:prox:est2-repeat}
        \frac{\tauTest_k(1+2\tau_k\gamma)}{2}\norm{\nextx-\realoptx}_X^2
        +\frac{\tauTest_k}{2}\norm{\nextx-\thisx}_X^2
        \le \frac{\tauTest_k}{2}\norm{\thisx-\realoptx}_X^2.
    \end{equation}
    Let us now force on the testing parameters the recursion
    \begin{equation}
        \label{eq:testing:prox:metric-update}
        \tauTest_0 = 1,\qquad \tauTest_{k+1} = \tauTest_k(1+2\tau_k\gamma).
    \end{equation}
    Then \eqref{eq:testing:prox:est2-repeat} yields
    \begin{equation}
        \label{eq:testing:prox:test-before-sum}
        \frac{\tauTest_{k+1}}{2}\norm{\nextx-\realoptx}_X^2
        +\frac{\tauTest_k}{2}\norm{\nextx-\thisx}_X^2
        \le \frac{\tauTest_k}{2}\norm{\thisx-\realoptx}_X^2.
    \end{equation}

    We now distinguish the two cases for the step sizes $\tau_k$.
    \begin{enumerate}
        \item
            Summing \eqref{eq:testing:prox:test-before-sum} for $k=0,\ldots,N-1$ gives
            \begin{equation*}
                \frac{\tauTest_N}{2}\norm{x^N-\realoptx}_X^2
                + \sum_{k=0}^{N-1}
                \frac{\tauTest_k}{2}\norm{\nextx-\thisx}_X^2
                \le \frac{\tauTest_0}{2}\norm{x^0-\realoptx}_X^2.
            \end{equation*}
            In particular, $\phi_0=1$ implies that
            \begin{equation*}
                \norm{x^N-\realoptx}_X^2 \le \tauTest_N^{-1} \norm{x^0-\realoptx}_X^2.
            \end{equation*}
            Since $\tau_k\equiv \tau$, \eqref{eq:testing:prox:metric-update} implies that $\tauTest_N = (1+2\tau\gamma)^N$. Setting $\tilde \mu \defeq (1+2\tau\gamma)^{-1/2}<1$ now gives convergence at the rate $O(\tilde \mu^{N})$ and therefore the claimed linear rate.

        \item  From \eqref{eq:testing:prox:test-before-sum} combined with \eqref{eq:testing:prox:metric-update} it follows directly that
            \begin{equation*}
                \frac{\norm{x^{k+1}-\realoptx}_X^2}{\norm{x^k-\realoptx}_X^2}
                \le
                \frac{\tauTest_k}{\tauTest_{k+1}}
                =
                (1+2\tau_k\gamma)^{-1} \to 0
            \end{equation*}
            since $\tau_k \to \infty$, which implies the claimed superlinear convergence of $\norm{x^k-\realoptx}_X$. (A similar argument can be used to directly show linear convergence for constant step sizes.)
            \qedhere
    \end{enumerate}
\end{proof}

\subsection*{Explicit splitting}

We now return to problems of the form
\begin{equation}
    \label{eq:testing:fb:problem}
    \min_{x \in X} F(x)+G(x)
\end{equation}
for Gateaux differentiable $F$ and study the convergence rates of the explicit (or forward--backward) splitting method
\begin{algeqbox}
    \begin{equation}
        \label{eq:testing:fb:fb}
        \nextx \defeq \prox_{\tau G}(x^k - \tau \grad F(x^k)),
    \end{equation}
\end{algeqbox}
which we recall can be written in implicit form as
\begin{equation}
    \label{eq:testing:fb:implicit}
    0 \in \tau[\subdiff G(\nextx)+\grad F(\thisx)]+(\nextx-\thisx).
\end{equation}

\begin{theorem}[explicit splitting iterate rates]
    \label{thm:testing:fb}
    Let $F:X\to\R$ and $G: X \to \Rbar$ be convex, proper, and lower semicontinuous. Suppose further that $F$ is Gateaux differentiable, $\grad F$ is Lipschitz continuous with constant $L>0$, and $G$ is $\gamma$-strongly convex for some $\gamma>0$.
    If $\inv{[\subdiff(F+G)]}(0) \ne \emptyset$ and the step length parameter $\tau>0$ satisfies $\tau L \le 2$, then for any initial iterate $x^0 \in X$ the iterates $\{\thisx\}_{k \in \N}$ generated by the explicit splitting method \eqref{eq:testing:fb:fb} converge linearly to the unique minimizer of \eqref{eq:testing:fb:problem}.
\end{theorem}

\begin{proof}
    Let $\realoptx \in \inv{[\subdiff(F+G)]}(0)$; by assumption, such a point exists and is unique due to the strong and therefore strict convexity of $G$. As in the proof of \cref{thm:testing:prox}, for each iteration $k \in \N$, pick a \term[parameter!testing]{testing parameter} $\tauTest_k>0$ and apply the \term[functional!testing]{test} $\tauTest_k\iprod{\freevar}{\nextx-\realoptx}_X$ to \eqref{eq:testing:fb:implicit} to obtain
    \begin{equation}
        \label{eq:testing:fb:tested}
        0 \in \tauTest_k\tau\iprod{\subdiff G(\nextx)+ \grad F(\thisx)}{\nextx-\realoptx}_X+\tauTest_k\iprod{\nextx-\thisx}{\nextx-\realoptx}_X.
    \end{equation}
    Since $G$ is strongly convex, it follows from \eqref{eq:testing:fb:implicit} and \cref{thm:smoothness:strong-convexity}\,\ref{item:strong-monotonicity} that
    \[
        \iprod{\partial G(\nextx)-\partial G(\realoptx)}{\nextx-\realoptx}_X \ge \gamma\norm{\nextx-\realoptx}_X^2.
    \]
    Similarly, since $\grad F$ is Lipschitz continuous, it follows from   \cref{cor:smoothness:three-point} that
    \[
        \iprod{\grad F(\thisx)- \grad F(\realoptx)}{\nextx-\realoptx}_X \ge -\frac{L}{4}\norm{\nextx-\thisx}_X^2.
    \]
    Combining the last two inequalities with $0 \in \subdiff G(\realoptx)+\grad F(\realoptx)$, we obtain
    \begin{equation}
        \label{eq:testing:fb:est1}
        \iprod{\partial G(\nextx)+\grad F(\thisx)}{\nextx-\realoptx}_X \ge \gamma\norm{\nextx-\realoptx}_X^2-\frac{L}{4}\norm{\nextx-\thisx}_X^2.
    \end{equation}
    Inserting this into \eqref{eq:testing:fb:tested} and using the three-point identity, we now obtain as in the proof of \cref{thm:testing:prox} that
    \begin{equation}
        \label{eq:testing:fb:est2}
        \frac{\tauTest_k(1+2\tau\gamma)}{2}\norm{\nextx-\realoptx}_X^2
        +\frac{\tauTest_k(1-\tau L/2)}{2}\norm{\nextx-\thisx}_X^2
        \le \frac{\tauTest_k}{2}\norm{\thisx-\realoptx}_X^2.
    \end{equation}
    Since $1-\tau L/2 \ge 0$, summing over $k=0,\ldots,N-1$, we arrive at
    \begin{equation*}
        \frac{\tauTest_N}{2}\norm{x^N-\realoptx}_X^2
        \le \frac{\tauTest_0}{2}\norm{x^0-\realoptx}_X^2.
    \end{equation*}
    As in \cref{thm:testing:prox}, the definition of $\tauTest_k$ shows that $\norm{\thisx-\realoptx}_X^2\to 0$ linearly.
\end{proof}

Observe that it is not possible to obtain superlinear convergence in this case since the assumption $\tau_k \le 2\inv L$ forces the step lengths to remain bounded.

\section{Structured algorithms and acceleration}
\label{sec:testing:structured}

We now extend the analysis above to the structured case where $H=A+B$, since we have already seen that most common first-order algorithms can be written as calculating in each step the next iterate $\nextx$ from a specific instance of the general preconditioned implicit--explicit splitting method
\begin{equation}
    \label{eq:testing:structured:alg}
    0 \in A(\nextx)+B(\thisx) + M(\nextx-\thisx).
\end{equation}
In the proofs of convergence of the proximal point and explicit splitting methods (e.g., in \cref{thm:testing:prox,thm:testing:fb} as well as in \cref{chap:convergence}), we had the step length $\tau_k$ in front of $H$ or $\grad F+\partial G$. On the other hand, in \cref{sec:convergence:general-prox} on structured algorithms, we incorporated the step length parameters into the preconditioning operator $M$. To transfer the testing approach from these fundamental methods to the structured methods, we will now split them out from $M$ and move them in front of $H$ as well by introducing a \term[operator!step length]{step length operator} $\Step_{k+1}$. We will also allow the preconditioner $\Precond_{k+1}$ to vary by iteration; as we will see below, this is required for accelerated versions of the PDPS method. Correspondingly, we consider the scheme
\begin{equation}
    \label{eq:testing:structured:alg-varying}
    0 \in \Step_{k+1}[A(\nextx)+B(\thisx)] + M_{k+1}(\nextx-\thisx).
\end{equation}
Since we now have a step length operator instead of a single scalar step length, we will also have to consider instead of a scalar testing parameter an iteration-dependent \term[operator!testing]{testing operator} $Z_{k+1} \in \linear(X; X)$.
The rough idea is that $Z_{k+1}M$ -- or, as needed for accelerated algorithms, $Z_{k+1}M_{k+1}$ -- will form a \enquote{local norm} that measures the rate of convergence in a nonuniform way; and rather than testing the (scalar) three-point identity \eqref{eq:testing:prox:est2-repeat}, we will build the testing already into the initial strong monotonicity inequality.
We therefore require an operator-level version of strong monotonicity, which we introduce next.

Let $A: X \setto X$ and let $\Test,\Gamma \in \linear(X; X)$ be such that $\Test\Gamma$ is positive semi-definite.
Then we say that $A$ is \term[mapping!monotone!$\Gamma$-strongly]{$\Gamma$-strongly monotone} at $\realoptx \in X$ with respect to $Z$ if
\begin{equation}
    \label{eq:testing:structured:strongmono}
    \iprod{A(x)-A(\realoptx)}{x-\realoptx}_{Z}
    \ge
    \norm{x-\realoptx}_{\Test\Gamma}^2
    \quad (x \in X).
\end{equation}
If this holds for all $\realoptx \in X$, we say that $A$ is $\Gamma$-strongly monotone with respect to $Z$.

It is clear that strongly monotone operators with parameter $\gamma>0$ are $\gamma \cdot \Id$-strongly monotone with respect to $Z=\Id$. More generally, operators with a separable block-structure, $A(x)=(A_1(x_1), \ldots, A_n(x_n))$ for $x=(x_1,\ldots,x_n)$ satisfy the property, as illustrated in more detail in the next example for the two-block case.

\begin{example}
    \label{ex:testing:structured:strongmono-block}
    Let $A(x)=(A_1(x_1), A_2(x_2))$ for $x=(x_1, x_2) \in X_1 \times X_2$ and the monotone operators $A_1: X_1 \setto X_1$ and $A_2: X_2 \setto X_2$.
    Suppose $A_1$ and $A_2$ are, respectively $\gamma_1$- and $\gamma_2$-(strongly) monotone for $\gamma_1,\gamma_2 \ge 0$.
    Then \eqref{eq:testing:structured:strongmono} holds for any $\phi_1,\phi_2 > 0$ for
    \[
        \Gamma \defeq \begin{pmatrix}\gamma_1 \Id & 0 \\ 0 & \gamma_2 \Id \end{pmatrix}
        \quad\text{and}\quad \Test \defeq \begin{pmatrix}\phi_1 \Id & 0 \\ 0 & \phi_2 \Id \end{pmatrix}
    \]
\end{example}

We do not impose $\Test\Gamma$ to be self-adjoint in \eqref{eq:testing:structured:strongmono}, although we use the norm notation. Forgoing with self-adjointness allows $\Gamma$ to have skew-adjoint parts $\Xi=-\Xi^*$, cf.~\cref{lemma:bcp-skew-adjoint}. Indeed, for the operator $H$ for the PDPS method from \eqref{eq:saddle-operator}, we can for $Z=\Id$ always choose $\Gamma=\begin{psmallmatrix} 0 & K^* \\ -K & 0 \end{psmallmatrix}$ skew-adjoint. With either $F$ or $G^*$ strongly convex, $\Gamma$ will also have corresponding components as in \cref{ex:testing:structured:strongmono-block}.

Let further $B: X \setto X$ and let $Z,\Lambda \in \linear(X; X)$ be such that $Z\Lambda$ is positive semi-definite. Then we say that $B$ is \term[mapping!monotone!three-point]{three-point monotone} at $\realoptx \in X$ with respect to $Z$ and $\Lambda$ if
\begin{equation}
    \label{eq:testing:structured:3monotone}
    \iprod{B(z)-B(\realoptx)}{x-\realoptx}_{Z}
    \ge
    -\frac{1}{4}\norm{x-z}_{\Test\Lambda}^2
    \quad (x, z \in X).
\end{equation}
If this holds for all $\realoptx \in X$, we say that $B$ is three-point monotone with respect to $Z$ and~$\Lambda$.

\begin{example}
    \label{ex:testing:structured:3monotone-block}
    Let $B(x)=(\grad E_1(x_1), \grad E_2(x_2))$ for $x=(x_1, x_2) \in X_1 \times X_2$ and the  respectively $L_1$- and $L_2$-smooth convex functions $E_1: X_1 \to \R$ and $E_2: X_2 \to \R$.
    Then a referral to \cref{cor:smoothness:three-point} shows \eqref{eq:testing:structured:3monotone} to hold for any $\phi_1,\phi_2 > 0$ for
    \[
        \Lambda \defeq \begin{pmatrix}L_1 \Id & 0 \\ 0 & L_2 \Id \end{pmatrix} \quad\text{and}\quad \Test \defeq \begin{pmatrix}\phi_1 \Id & 0 \\ 0 & \phi_2 \Id \end{pmatrix}
    \]
    More generally, we can take $B(x)=(B_1(x_1), B_2(x_2))$ for $B_1: X_1 \to X_1$ and $B_2: X_2 \to X_2$ three-point monotone as defined in \eqref{eq:smoothness:three-point:monotonicity}.
\end{example}

Clearly \cref{ex:testing:structured:3monotone-block}, like \cref{ex:testing:structured:strongmono-block}, generalizes to a large number of blocks, and both generalize to operators acting separably on more general direct sums of orthogonal subspaces.

We are now ready to forge our hammer for producing convergence rates for structured algorithms.
In the following, for any $M, N \in \linear(X; X)$, we write $M \succeq N$ to mean that $M-N$ is positive semi-definite: $\norm{x}_M^2 \ge \norm{x}_N^2$ for all $x \in X$.

\begin{theorem}
    \label{thm:testing:structured:convergence}
    Let $A, B: X \setto X$ and $H \defeq A+B$.
    For each $k \in \N$, let further $\Test_{k+1},\Step_{k+1},\Precond_{k+1} \in \linear(X; X)$ be such that $\Test_{k+1}\Precond_{k+1}$ is self-adjoint and positive semi-definite. Assume that there exists a $\realoptx \in \inv H(0)$.
    For each $k \in \N$, suppose for some $\Gamma,\Lambda \in \linear(X; X)$ that $A$ is $\Gamma$-strongly monotone at $\realoptx$ with respect to $\Test_{k+1}\Step_{k+1}$ and that $B$ is three-point monotone at $\realoptx$ with respect to $\Test_{k+1}\Step_{k+1}$ and $\Lambda$.
    Let the initial iterate $x^0 \in X$ be arbitrary, and suppose $\{\nextx\}_{k \in \N}$ are generated by \eqref{eq:testing:structured:alg-varying}.
    If for every $k \in \N$ both
    \begin{align}
        \label{eq:testing:structured:metric-update}
        \Test_{k+1}(\Precond_{k+1}+2\Step_{k+1}\Gamma) & \succeq \Test_{k+2}\Precond_{k+2}
        \quad\text{and}
        \\
        \label{eq:testing:structured:hypo-bound}
        \Test_{k+1}\Precond_{k+1} \succeq \Test_{k+1}\Step_{k+1}\Lambda/2.
    \end{align}
    hold, then
    \begin{equation}
        \label{eq:testing:structured:main-estimate}
        \frac{1}{2}\norm{x^N-\realoptx}_{\Test_{N+1}\Precond_{N+1}}^2
        \le \frac{1}{2}\norm{x^0-\realoptx}_{\Test_1\Precond_1}^2.
    \end{equation}
\end{theorem}
\begin{proof}
    For brevity, write $\Happrox_{k+1}(\nextx) \defeq \Step_{k+1}[A(\nextx)+B(\thisx)]$.
    First, from \eqref{eq:testing:structured:strongmono} and \eqref{eq:testing:structured:3monotone} we have that
    \begin{equation}
        \label{eq:testing:structured:tested}
        \iprod{\Happrox_{k+1}(\nextx)}{\nextx-\realoptx}_{\Test_{k+1}}
        \ge
        \norm{\nextx-\realoptx}_{\Test_{k+1}\Step_{k+1}\Gamma}^2 -\frac{1}{4}\norm{\thisx-\nextx}^2_{\Test_{k+1}\Step_{k+1}\Lambda}.
    \end{equation}
    Multiplying \eqref{eq:testing:structured:alg-varying} with $\Test_{k+1}$ and rearranging, we obtain
    \begin{equation*}
        -\Test_{k+1} M_{k+1}(\nextx-\thisx) \in \Test_{k+1}\Happrox_{k+1}(\nextx).
    \end{equation*}
    Inserting this into \eqref{eq:testing:structured:tested} and applying the preconditioned three-point formula \eqref{eq:convergence:three-point-identity-cp} for $M=\Test_{k+1}\Precond_{k+1}$ yields
    \begin{equation*}
        \frac{1}{2}\norm{\nextx-\realoptx}_{\Test_{k+1}(\Precond_{k+1}+2\Step_{k+1}\Gamma)}^2
        +\frac{1}{2}\norm{\nextx-\thisx}_{\Test_{k+1}(\Precond_{k+1}-\Step_{k+1}\Lambda/2)}^2
        \le \frac{1}{2}\norm{\thisx-\realoptx}_{\Test_{k+1}\Precond_{k+1}}^2.
    \end{equation*}
    Using \eqref{eq:testing:structured:metric-update} and \eqref{eq:testing:structured:hypo-bound}, this implies that
    \begin{equation}
        \label{eq:testing:structured:quantitative-fejer}
        \frac{1}{2}\norm{\nextx-\realoptx}_{\Test_{k+2}\Precond_{k+2}}^2
        \le \frac{1}{2}\norm{\thisx-\realoptx}_{\Test_{k+1}\Precond_{k+1}}^2.
    \end{equation}
    Summing over $k=0,\ldots,N-1$ now yields the claim.
\end{proof}
The inequality \eqref{eq:testing:structured:quantitative-fejer} is a quantitative or \term[sequence!Fejér monotone!variable metric]{variable metric} version of the Fejér monotonicity of \cref{lemma:opial}\,\ref{item:opial-nonincreasing} with respect to $\hat X=\{\realoptx\}$.

If \cref{thm:testing:structured:convergence} is applicable, we immediately obtain the convergence rate result.
\begin{corollary}[convergence with a rate]
    \label{cor:testing:structured:rate}
    If \eqref{eq:testing:structured:main-estimate} holds and $\Test_{N+1} \Precond_{N+1} \succeq \mu(N) I$ for some $\mu:\N\to\R$, then $\norm{x^N - \realoptu}^2 \to 0$ at the rate $O(1/\mu(N))$.
\end{corollary}

\subsection*{Primal-dual proximal splitting methods}

We now apply this operator-testing technique to primal-dual splitting methods for the solution of
\begin{equation}
    \label{eq:testing:pdps:problem}
    \min_{x \in X} F_0(x)+E(x)+G(Kx)
\end{equation}
with $F_0: X \to \Rbar$, $E: X \to \R$, and $G: Y \to \Rbar$ convex, proper, and lower semicontinuous and $K \in \linear(X; Y)$.
We will also write $F\defeq F_0+E$. The methods include in particular the PDPS method with a forward step \eqref{eq:convergence:pdps:forward}. Now allowing varying step lengths and an over-relaxation parameter $\omega_k$, this can be written
\begin{algeqbox}
    \begin{equation}%
        \label{eq:testing:pdps:forward}%
        \left\{
            \begin{aligned}
                \nextx & \defeq (I+\tau_k \subdiff F_0)^{-1}(\thisx - \tau_k K^* \thisy - \tau_k \grad E(\thisx)),\\
                \overnextx & \defeq \omega_k (\nextx-\thisx)+\nextx, \\
                \nexty & \defeq (I+\sigma_{k+1} \subdiff G^*)^{-1}(\thisy + \sigma_{k+1} K \overnextx).
            \end{aligned}
        \right.
    \end{equation}
\end{algeqbox}
For the basic version of the algorithm with $\omega_k=1$, $\tau_k \equiv \tau_0 > 0$, and $\sigma_k \equiv \sigma_0 > 0$, we have seen in \cref{cor:convergence:pdps:forward} that the iterates converge weakly if the step length parameters satisfy
\begin{equation}
    \label{eq:testing:pdps:init-cond}
    L\tau_0/2 + \tau_0\sigma_0\norm{K}_{\linear(X; Y)}^2 < 1,
\end{equation}
where $L$ is the Lipschitz constant of $\grad E$.
We will now show that under strong convexity of $F_0$, we can choose these parameters to \emph{accelerate} the algorithm to yield convergence at a rate $O(1/N^2)$. If both $F_0$ and $G^*$ are strongly convex, we can even obtain linear convergence. Throughout the following, $\realoptu=(\realoptx,\realopty)$ denotes a root of
\[
    H(u) \defeq
    \begin{pmatrix}
        \subdiff F_0(x)+\grad E(x)+K^*y \\
        \subdiff G^*(y) - Kx
    \end{pmatrix},
\]
which we assume exists. From \cref{thm:convex:fenchel}, this is the case if an interior point condition is satisfied for $G\circ K$ and \eqref{eq:testing:pdps:problem} admits a solution.

We will also require the following technical lemma in place of the simpler growth argument for the choice \eqref{eq:testing:prox:metric-update}.
\begin{lemma}
    \label{lemma:testing:pdps:quadratic}
    Pick $\tauTest_0>0$ arbitrarily, and define iteratively $\tauTest_{k+1} \defeq \tauTest_k\bigl(1+2\gamma\tauTest_k^{-1/2}\bigr)$ for some $\gamma>0$. Then there exists a constant $c>0$ such that $\tauTest_k \ge \bigl(c k+\tauTest_0^{1/2}\bigr)^2$  for all $k \in \N$.
\end{lemma}
\begin{proof}
    Replacing $\tauTest_k$ by $\tauTest_k' \defeq \gamma^{-2} \tauTest_k$, we may assume without loss of generality that $\gamma=1$.
    We claim that $\tauTest_k^{1/2} \ge c k+\tauTest_0^{1/2}$ for some $c>0$. We proceed by induction. The case $k=0$ is clear. If the claim holds for $k=0,\ldots,N-1$, we can unroll the recursion to obtain the estimate
    \[
        \tauTest_N - \tauTest_0 = \sum_{k=0}^{N-1} 2\tauTest_k^{1/2}
        \ge 2 \sum_{k=0}^{N-1} c k + 2\tauTest_0^{1/2}N
        = c N(N-1)+2\tauTest_0^{1/2}N
        = c N^2+(2\tauTest_0^{1/2}-c) N.
    \]
    Expanding $(c N+\tauTest_0^{1/2})^2=c^2 N^2 + 2c\tauTest_0^{1/2} N + \tauTest_0$, we see that the claim for $\tauTest_N$ holds if $c \ge c^2$ and $2\tauTest_0^{1/2}-c \ge 2c\tauTest_0^{1/2}$. Taking the latter with equality and solving for $c$ yields  $c = 2 \tauTest_0^{1/2}/(1+2\tauTest_0^{1/2})<1$ and hence also the former. Since this choice of $c$ does not depend on $N$, the claim follows.
\end{proof}

\begin{theorem}[accelerated and linearly convergent PDPS]
    \label{thm:testing:pdps:accel}
    Let $F_0: X \to \Rbar$, $E: X \to \R$ and $G: Y \to \Rbar$ be convex, proper, and lower semicontinuous  with $\grad E$ Lipschitz continuous with constant $L>0$. Also let $K \in \linear(X; Y)$, and suppose the assumptions of \cref{thm:convex:fenchel} are satisfied with $F \defeq F_0+E$.
    Pick initial step lengths $\tau_0, \sigma_0>0$ subject to \eqref{eq:testing:pdps:init-cond}.
    For any initial iterate $u^0 \in X \times Y$, suppose $\{\nextu = (\nextx, \nexty) \}_{k \in \N}$ are generated by \eqref{eq:testing:pdps:forward}.
    \begin{enumerate}[label=(\roman*)]
        \item If $F_0$ is strongly convex with factor $\gamma>0$, and we take
            \begin{equation}
                \label{eq:testing:pdps:accel}
                \omega_k \defeq 1/\sqrt{1+2\gamma\tau_k},
                \quad
                \tau_{k+1} \defeq \tau_k\omega_k,
                \quad\text{and}\quad
                \sigma_{k+1} \defeq \sigma_k/\omega_k,
            \end{equation}
            then $\norm{x^N-\realoptx}_X^2 \to 0$ at the rate $O(1/N^2)$.

        \item If both $F_0$ and $G^*$ are strongly convex with factor $\gamma>0$ and $\rho>0$, respectively, and we take
            \begin{equation}
                \label{eq:testing:pdps:linear}
                \omega_k \defeq 1/(1+2\theta),
                \quad
                \theta \defeq \min\{\rho\sigma_0,\gamma\tau_0\},
                \quad
                \tau_k \defeq \tau_0
                \quad\text{and}\quad
                \sigma_k \defeq \sigma_0,
            \end{equation}
            then $\norm{x^N-\realoptx}_X^2 + \norm{y^N-\realopty}_Y^2 \to 0$ linearly.
    \end{enumerate}
\end{theorem}
\begin{proof}
    Recalling \cref{cor:convergence:pdps:forward}, we write \eqref{eq:testing:pdps:forward} in the form  \eqref{eq:testing:structured:alg-varying} by taking
    \begin{gather*}
        A(u) \defeq \begin{pmatrix} \subdiff F_0(x) \\ \subdiff G^*(y) \end{pmatrix}+\Xi u,
        \quad
        B(u) \defeq \begin{pmatrix} \grad E(x) \\ 0 \end{pmatrix},
        \quad
        \Xi \defeq \begin{pmatrix} 0 & K^* \\ -K & 0 \end{pmatrix},
        \\
        \label{eq:testing:pdps:forward:step-precond}
        \Step_{k+1} \defeq \begin{pmatrix} \tau_k \Id & 0 \\ 0 & \sigma_{k+1} \Id \end{pmatrix},
        \quad\text{and}\quad
        \Precond_{k+1} \defeq \begin{pmatrix} \Id & -\tau_k K^* \\  - \omega_k\sigma_{k+1} K & \Id \end{pmatrix}.
    \end{gather*}
    As before, \cref{thm:convex:fenchel} guarantees that $\inv H(0) \ne \emptyset$.
    For some primal and dual testing parameters $\tauTest_k,\sigmaTest_{k+1} > 0$, we also take as our testing operator
    \begin{equation}
        \label{eq:testing:pdps:forward:test}
        \Test_{k+1} \defeq \begin{pmatrix} \tauTest_k \Id & 0 \\ 0 & \sigmaTest_{k+1} \Id \end{pmatrix}.
    \end{equation}
    By \cref{ex:testing:structured:strongmono-block,ex:testing:structured:3monotone-block}, $A$ is then $\Gamma$-strongly monotone with respect to $\Test_{k+1}\Step_{k+1}$ and $B$ is three-point monotone with respect to  $\Test_{k+1}\Step_{k+1}$ and $\Lambda$ for
    \[
        \Gamma \defeq \Xi+\begin{pmatrix} \gamma  \Id & 0 \\ 0 & \rho  \Id \end{pmatrix},
        \quad\text{and}\quad
        \Lambda \defeq \begin{pmatrix} L \, \Id & 0 \\ 0 & 0 \end{pmatrix},
    \]
    where $\rho=0$ if $G^*$ is not strongly convex.

    We will apply \cref{thm:testing:structured:convergence}.
    Taking $\omega_k \defeq \inv\sigma_{k+1}\inv\sigmaTest_{k+1}\tauTest_k\tau_k$, we expand
    \begin{equation}
        \label{eq:testing:pdps:proof-zm}
        \Test_{k+1}\Precond_{k+1} =
        \begin{pmatrix}
            \tauTest_k \Id & -\tauTest_k\tau_k K^* \\
            -\tauTest_k\tau_k K & \sigmaTest_{k+1} \Id
        \end{pmatrix}.
    \end{equation}
    Thus $\Test_{k+1}\Precond_{k+1}$ is self-adjoint as required.
    We still need to show that it is nonnegative and indeed grows at a rate that gives our claims. We also need to verify \eqref{eq:testing:structured:metric-update} and \eqref{eq:testing:structured:hypo-bound}, which expand as
    \begin{align}
        \label{eq:testing:pdps:metric-update}
        \begin{pmatrix}
            (\tauTest_k(1+2\gamma\tau_k) - \tauTest_{k+1}) \Id & (\tauTest_k\tau_k + \tauTest_{k+1}\tau_{k+1})K^* \\
            (\tauTest_{k+1}\tau_{k+1}-2\sigmaTest_{k+1}\sigma_{k+1}-\tauTest_k\tau_k) K & (\sigmaTest_{k+1}(1+2\rho\sigma_{k+1}) -\sigmaTest_{k+2}) \Id
        \end{pmatrix}
        &\succeq 0,
        \quad\text{and}
        \\
        \label{eq:testing:pdps:hypo-bound}
        \begin{pmatrix}
            \tauTest_k(1-\tau_k L/2) \Id & -\tauTest_k\tau_k K^* \\
            -\tauTest_k\tau_k K & \sigmaTest_{k+1} \Id
        \end{pmatrix} & \succeq 0.
    \end{align}
    We now proceed backward by deriving the step length rules as sufficient conditions for these two inequalities. First, clearly \eqref{eq:testing:pdps:metric-update} holds if for all $k \in \N$ we can guarantee that
    \begin{equation}
        \label{eq:testing:pdps:par0}
        \tauTest_{k+1} \leq \tauTest_k(1+2\gamma\tau_k),
        \quad
        \sigmaTest_{k+1} \leq \sigmaTest_k(1+2\rho\sigma_k),
        \quad
        \text{and}
        \quad
        \tauTest_k\tau_k=\sigmaTest_k\sigma_k.
    \end{equation}
    We deal with \eqref{eq:testing:pdps:hypo-bound} and the lower bounds on $\Test_{k+1}\Precond_{k+1}$ in one go. By Young's inequality, we have for any $\delta \in (0, 1)$ that
    \[
        2\tauTest_k\tau_k\iprod{Kx}{y}
        \le (1-\delta)\tauTest_k\norm{x}^2
        + \tauTest_k\tau_k^2\inv{(1-\delta)}\norm{K^*y}^2
        \quad (x \in X,\, y \in Y),
    \]
    hence recalling \eqref{eq:testing:pdps:proof-zm} also
    \begin{equation}
        \label{eq:testing:pdps:zimi-est}
        \Test_{k+1}\Precond_{k+1} \succeq
        \begin{pmatrix}
            \delta\tauTest_k \Id & 0 \\
            0 & \sigmaTest_{k+1} \Id - \tauTest_k\tau_k^2\inv{(1-\delta)} KK^*
        \end{pmatrix}.
    \end{equation}
    Similarly, for the operator from \eqref{eq:testing:pdps:hypo-bound}, we have
    \[
        \begin{pmatrix}
            \tauTest_k(1-\tau_k L/2) \Id & -\tauTest_k\tau_k K^* \\
            -\tauTest_k\tau_k K & \sigmaTest_{k+1} \Id
        \end{pmatrix}
        \succeq
        \begin{pmatrix}
            \tauTest_k(\delta-\tau_k L/2) \Id & 0 \\
            0 & \sigmaTest_{k+1} \Id - \tauTest_k\tau_k^2\inv{(1-\delta)} K^*K
        \end{pmatrix}.
    \]
    The condition \eqref{eq:testing:pdps:hypo-bound} is therefore satisfied and $\Test_{k+1}\Precond_{k+1} \ge \epsilon\Test_{k+1}$ for some $\epsilon>0$ if \eqref{eq:testing:pdps:par0} holds and both
    \begin{equation}
        \label{eq:testing:pdps:hypo-bound:2}
        \delta\tauTest_k \ge \epsilon\tauTest_k + \tauTest_k \tau_k L/2
        \quad\text{and}\quad
        \sigmaTest_{k+1} \ge \epsilon\sigmaTest_{k+1} + \tauTest_k\tau_k^2\inv{(1-\delta)} \norm{K}^2.
    \end{equation}
    By \eqref{eq:testing:pdps:par0} we have $\sigmaTest_{k+1} \ge \sigmaTest_k$, and hence using in addition the last part of \eqref{eq:testing:pdps:par0} shows that \eqref{eq:testing:pdps:hypo-bound:2} holds if
    \begin{equation}
        \label{eq:testing:pdps:step-conds-L-0}
        \delta-\epsilon \ge \tau_k L/2
        \quad\text{and}\quad
        (1-\delta)(1-\epsilon) \ge \tau_k\sigma_k \norm{K}^2.
    \end{equation}
    If we choose $\tau_k$ and $\sigma_k$ such that their product stays constant (i.e., $\tau_k\sigma_k=\sigma_0\tau_0$), then the second equality holds for $\delta = 1-\sigma_0\tau_0\norm{K}^2/(1-\epsilon)$, which has to be positive. Inserting this into the first part of \eqref{eq:testing:pdps:step-conds-L-0}, we see that it holds if
    $1 \ge \sigma_0\tau_0\norm{K}^2/(1-\epsilon) + \epsilon + \tau_k L/2$.
    This holds for some $\epsilon>0$ due to the assumed \eqref{eq:testing:pdps:init-cond}, i.e., $\tau_k L/2 + \sigma_0\tau_0\norm{K}^2 < 1$. Since $\{\tau_k\}_{k \in \N}$ is nonincreasing, we see that \eqref{eq:testing:pdps:step-conds-L-0} and hence \eqref{eq:testing:pdps:hypo-bound} is satisfied when the initialization condition \eqref{eq:testing:pdps:init-cond} holds.

    To apply \cref{thm:testing:structured:convergence}, all that remains is to verify \eqref{eq:testing:pdps:par0} and that $\tau_k\sigma_k=\tau_0\sigma_0$. To obtain convergence rates, we need to further study the rate of increase of $\tauTest_k$ and $\sigmaTest_{k+1}$, which we recall that we wish to make as high as possible.
    \begin{enumerate}
        \item If $\gamma > 0$ and $\rho=0$, the best possible choice allowed by \eqref{eq:testing:pdps:par0} is $\sigmaTest_k \equiv \sigmaTest_0$ and $\tauTest_{k+1} = \tauTest_k(1+2\gamma\tau_k)$ with $\sigma_k=\tauTest_k\tau_k/\sigmaTest_0$. Together with the condition $\tau_k\sigma_k=\sigma_0\tau_0$, this forces $\sigma_0\tau_0=\tauTest_k\tau_k^2/\sigmaTest_0$.
            If we take $\sigmaTest_0=1/(\sigma_0\tau_0)$, we thus need $\tau_k=\tauTest_k^{-1/2}$.
            Since $\sigma_{k+1}=\sigma_0\tau_0/\tau_{k+1}=1/(\sigmaTest_0\tau_{k+1})$, we obtain the relations
            \[
                \omega_k=\frac{\tauTest_k\tau_k}{\sigma_{k+1}\sigmaTest_{k+1}}
                =\frac{\tauTest_k^{1/2}}{\tauTest_{k+1}^{1/2}}
                =\frac{1}{\sqrt{1+2\gamma\tau_k}},
            \]
            which are satisfied for the choices of $\omega_k$, $\tau_{k+1}$, and $\sigma_{k+1}$ in \eqref{eq:testing:pdps:accel}.

            We now use \cref{thm:testing:structured:convergence,cor:testing:structured:rate} and \eqref{eq:testing:pdps:zimi-est} to obtain
            \begin{equation*}
                \frac{\delta\tauTest_N}{2}\norm{x^N-\realoptx}_X^2
                \le
                \frac{1}{2}\norm{u^N-\realoptu}_{\Test_{N+1}\Precond_{N+1}}^2
                \le C_0 \defeq \frac{1}{2}\norm{u^0-\realoptu}_{\Test_1\Precond_1}^2.
            \end{equation*}
            Although this does not tell us anything about the convergence of the dual iterates $\{y^N\}_{N \in \N}$ as $\sigmaTest_N \equiv \sigmaTest$ stays constant, \cref{lemma:testing:pdps:quadratic} shows that the primal test $\tauTest_{N}$ grows at the rate $\Omega(N^2)$. Hence we obtain the claimed convergence of the primal iterates at the rate $O(1/N^2)$.

        \item If $\gamma>0$ and $\rho>0$ and we take $\tau_k \equiv \tau_0$ and $\sigma_k \equiv \sigma_0$, the last condition of \eqref{eq:testing:pdps:par0} forces $\sigmaTest_k=\tauTest_k\tau_0/\sigma_0$.
            Inserting this into the second condition yields $\tauTest_{k+1} \leq \tauTest_k(1+2\rho\sigma_0)$.
            Together with the first condition, we therefore at best can take $\tauTest_{k+1}=\tauTest_k(1+2\theta)$ for $\theta \defeq \min\{\rho\sigma_0,\gamma\tau_0\}$.
            Reversing the roles of $\sigmaTest$ and $\tauTest$, we see that we can at best take $\sigmaTest_{k+1}=\sigmaTest_k(1+2\theta)$.
            This leads to the relations
            \[
                \omega_k=\frac{\tauTest_k\tau_0}{\sigma_0\sigmaTest_{k+1}}
                =\frac{\tauTest_k}{\tauTest_{k+1}}=\frac{1}{1+2\theta},
            \]
            which are again satisfied by the respective choices in \eqref{eq:testing:pdps:linear}.

            We finish the proof with  \cref{thm:testing:structured:convergence,cor:testing:structured:rate}, observing now from \eqref{eq:testing:pdps:zimi-est} that $\Test_N\Precond_N \ge C (1+2\theta)^N \Id$ for some $C>0$.
            \qedhere
    \end{enumerate}
\end{proof}
Note that if $\gamma=0$ and $\rho=0$, \eqref{eq:testing:pdps:par0} forces $\tauTest_k \equiv \tauTest_0$ as well as $\sigmaTest_k \equiv \sigmaTest_0$.
If we take $\tauTest_k \equiv 1$, then we also have to take $\tau_k=\sigma_k\sigmaTest_0$.
We can use this to define $\sigmaTest_0$ if we also fix $\tau_k \equiv \tau_0$ and $\sigma_k \equiv \sigma_0$. This also forces $\omega_k \equiv 1$.
We thus again arrive at \eqref{eq:testing:pdps:par0} as well as $\tau_k\sigma_k=\sigma_0\tau_0$.
However, we cannot obtain from this convergence rates for the iterates, merely boundedness and hence weak convergence as in \cref{sec:convergence:fb-general}.

\chapter{Splitting methods: gaps and ergodic results}
\label{chap:gap}

We continue with the testing framework introduced in \cref{chap:testing} for proving rates of convergence of iterates of optimization methods. This generally required strong convexity, which is not always available.
In this chapter, we use the testing idea to derive convergence rates of objective function values and other, more general, \emph{gap functionals} that indicate algorithm convergence more indirectly than iterate convergence. This can be useful in cases where we can only obtain weak convergence of iterates, but can obtain rates of convergence of such a gap functional. Nevertheless, this gap convergence often will only be \term[convergence!ergodic]{ergodic}, i.e., the estimates only apply to a weighted sum of the history of iterates instead of the most recent iterate. In fact, we will first derive ergodic estimates for all algorithms. If we can additionally show that the algorithm is {monotonic} with respect to this gap, we can improve the ergodic estimate to the nonergodic ones as in the previous chapters.

\section{Gap functionals}
\label{sec:gap:gap}

We recall that one of the three fundamental ingredients in the convergence proofs of \cref{chap:convergence} was the monotonicity of $H$ (with one of the points fixed to a root $\realoptx$). We now modify this requirement to be able to prove estimates on the convergence of function values when $H=\subdiff G$ for some proper, convex, and lower semicontinuous $G: X\to \Rbar$. In this case, by the definition of the convex subdifferential,
\begin{equation}
    \label{eq:gap:gap:convex}
    \iprod{\subdiff G(\nextx)}{\nextx-\optx}_X \ge G(\nextx)-G(\optx)
    \quad (\optx \in X).
\end{equation}
On the other hand, for an $L$-smooth functional $F: X\to \R$, we can use the three-point estimates of \cref{cor:smoothness:three-point} to obtain
\begin{equation}
    \label{eq:gap:gap:smooth}
    \iprod{\grad F(\thisx)}{\nextx-\optx}_X \ge F(\nextx)-F(\optx)-\frac{1}{2L}\norm{\nextx-\thisx}_X^2
    \quad (\optx \in X).
\end{equation}
These two inequalities are enough to obtain function value estimates for the more general case $H=\subdiff G + \grad F$ including a forward step with respect to $F$.
We will produce such estimates in \cref{sec:gap:function}.

\subsection*{Generic gap functionals}

More generally, when $H$ does not directly arise from subdifferentials or gradients but has a more complicated structure, we introduce several \term[functional!gap]{gap functionals}.
We identified in \cref{chap:convergence} that for some lifted functionals $\tilde F$ and $\tilde G$ and a skew-adjoint operator $\Xi=-\Xi^*$, the unaccelerated PDPS, PDES, and DRS consist in taking $H=\subdiff \tilde G + \grad \tilde F + \Xi$ and iterating
\begin{equation}
    \label{eq:gap:ergodic:general-alg}
    0 \in \subdiff \tilde G(\nextx)+\grad \tilde F(\thisx) + \Xi\nextx + M(\nextx-\thisx),
\end{equation}
where the skew-adjoint operator $\Xi$ does not arise as a subdifferential of any function. Working with this requires extra effort, especially when we later study accelerated methods.

Note that by the skew-adjointness of $\Xi$, we have $\iprod{\Xi\realoptx}{\realoptx}_X=0$. Using this and the estimates  \eqref{eq:gap:gap:convex} and \eqref{eq:gap:gap:smooth} on $\tilde F$ and $\tilde G$, we obtain for the basic unaccelerated scheme \eqref{eq:gap:ergodic:general-alg} the estimate
\[
    \iprod{\subdiff \tilde G(\nextx)+\grad \tilde F(\thisx) + \Xi\nextx}{\nextx-\realoptx}_X
    \ge
    \tilde\gap(x; \realoptx)
    -\frac{1}{2L}\norm{\nextx-\thisx}_X^2
\]
with the \term[functional!gap!generic]{generic gap functional}
\begin{equation}
    \label{eq:gap:ergodic:generic-gap}
    \tilde\gap(x; \optx)
    \defeq
    (\tilde G+\tilde F)(x) - (\tilde G+\tilde F)(\optx)
    +\iprod{\Xi\optx}{x}_X.
\end{equation}
In the next lemma, we collect some elementary properties of this functional.
Note that $\tilde\gap(x,z)=0$ is possible even for $x\neq z$.

\begin{lemma}
    \label{lem:gap:ergodic:functional}
    Let $H \defeq \subdiff \tilde G + \grad \tilde F + \Xi$, where $\Xi \in \linear(X; X)$ is skew-adjoint and $\tilde G: X \to \Rbar$ and $\tilde F: X \to \R$ are convex, proper, and lower semicontinuous.
    If $\realoptx \in \inv H(0)$, then $\tilde\gap(\freevar; \realoptx) \ge 0$ and $\tilde\gap(\realoptx; \realoptx)=0$.
\end{lemma}
\begin{proof}
    We first note that $\realoptx \in \inv H(0)$ is equivalent to $-\Xi\realoptx \in \subdiff (\tilde F+\tilde G)(\realoptx)$. Hence using the definition of the convex subdifferential and the fact that $\iprod{\Xi\realoptx}{\realoptx}_X=0$ due to the skew-adjointness of $\Xi$, we deduce for arbitrary $x \in X$ that
    \[
        (\tilde F+\tilde G)(x)
        -(\tilde F+\tilde G)(\realoptx)
        \ge
        \iprod{-\Xi\realoptx}{x-\realoptx}_X
        =
        \iprod{-\Xi\realoptx}{x}_X,
    \]
    i.e., $\tilde\gap(x, \realoptx) \ge 0$.
    The fact that $\tilde\gap(\realoptx, \realoptx)=0$ follows immediately from the skew-adjointness of $\Xi$.
\end{proof}

The function value estimates \eqref{eq:gap:gap:convex} and \eqref{eq:gap:gap:smooth} -- unlike simple monotonicity-based nonnegativity estimates -- do not depend on $\optx$ being a root of $H$.
Therefore, taking any \emph{bounded} set $B \subset X$ such that $\inv H(0) \isect B \ne \emptyset$, we see that the \term[gap!partial]{partial gap}
\[
    \tilde \gap(x; B) \defeq \sup_{\optx \in B} \tilde\gap(x; \optx)
\]
also satisfies $\tilde \gap(\freevar; B) \ge 0$.

\subsection*{The Lagrangian duality gap}

Let us now return to the problem
\begin{equation}
    \label{eq:gap:pd:problem}
    \min_{x \in X}~ F(x) + G(Kx),
\end{equation}
where we split $F=F_0+E$ assuming $E$ to have a Lipschitz-continuous gradient.
With the notation $u=(x, y)$, we recall that \cref{thm:convex:fenchel} guarantees the existence of a primal-dual solution $\realoptu$ whenever its conditions are satisfied.
This, we further recall, can be written as $0 \in H(\realoptu)$ for
\begin{subequations}%
\label{eq:gap:pd:setup}
\begin{equation}
    \label{eq:gap:pd:setup:h}
    H(u) \defeq
    \begin{pmatrix}
        \subdiff F(x)+K^*y \\
        \subdiff G^*(y) - Kx
    \end{pmatrix}.
\end{equation}
As we have already seen in, e.g., \cref{thm:testing:pdps:accel}, we can express this choice of $H$ in the framework of \eqref{eq:gap:ergodic:general-alg} with
\begin{gather}
    \label{eq:gap:pd:setup:g}
    \tilde G(u) \defeq F_0(x)+G^*(y),
    \quad
    \tilde F(u) \defeq E(x),
    \quad\text{and}\quad
    \Xi \defeq \begin{pmatrix} 0 & K^* \\ -K & 0 \end{pmatrix}.
\end{gather}%
\end{subequations}%
With this, the generic gap functional $\tilde \gap$ from \eqref{eq:gap:ergodic:generic-gap} becomes the \term[gap!duality!Lagrangian]{Lagrangian duality gap}
\begin{equation}
    \label{eq:gap:pd:gap}
    \gap_L(u; \optu) \defeq
    \bigl(F(x) + \iprod{\opty}{K x}_Y  - G^*(\opty)\bigr)
    -\bigl(F(\optx) + \iprod{y}{K \optx}_Y -  G^*(y)\bigr) \leq \bar\gap(u),
\end{equation}
where
\[
    \bar\gap(u) \defeq F(x)+G(Kx)+F^*(-K\opty)+G^*(\opty)
\]
is the real \term[gap!duality]{duality gap}, cf.~\eqref{eq:convex:fenchel:duality-gap}.
As \cref{eq:convex:fenchel:lagrangian-duality-gap-bound} shows, the Lagrangian duality gap is nonnegative when $\optu=\realoptu \in \inv H(0)$.

Since \eqref{eq:gap:gap:convex} and \eqref{eq:gap:gap:smooth} do not depend on $\optx$ being a root of $H$, convergence results for the Lagrangian duality gap can sometimes be improved slightly by taking any bounded set $B \subset X \times Y$ such that $B \isect \inv H(0) \ne \emptyset$ and defining the \term[gap!duality!partial]{partial duality gap}
\begin{equation}
    \label{eq:gap:pd:partial}
    \gap(u; B) \defeq \sup_{\optu \in B} \gap_L(u; \optu).
\end{equation}
This satisfies $0 \le \gap(u; B) \le \bar\gap(u)$. Moreover, by the definition of $F^*$ and $G^{**}=G$, we have $\gap(u; X \times Y)=\bar\gap(u)$, which explains both the importance of partial duality gaps and the term ``partial gap''.

\subsection*{Bregman divergences and gap functionals}
\label{sec:gap:ergodic:bregman}

Although we will not need this in the following, we briefly discuss a possible extension to Banach spaces.
Let $X$ be a Banach space and let $J: X \to \Rbar$ be convex. Then for $x \in \dom J$ and $p \in \subdiff J(x)$, one can define the asymmetric \term[divergence, Bregman]{Bregman divergence} (or \term[distance, Bregman|see{divergence, Bregman}]{distance})
\[
    B_J^{p}(z, x) \defeq J(z)-J(x)-\dualprod{p}{z-x}_{X}
    \quad (x \in X).
\]
Due to the definition of the convex subdifferential, this is nonnegative. It is also possible to symmetrize the distance by considering $\tilde B_J(x, z) \defeq B_J^{q}(x, z) + B_J^p(z, x)$ with $q \in \subdiff J(z)$ and $z \in \dom J$, but even the symmetrized divergence is not generally a true distance as it can happen that $B_J(x, z)=0$ even if $x\neq z$.

The Bregman divergence satisfies a three-point identity for any $\optx \in \dom J$: We have
\[
    \begin{aligned}
        B_J^{p}(\optx, x)  - B_J^{p}(\optx, z) + B_J^q(x,z)
        &
        =
        [J(\optx)-J(x)-\dualprod{p}{\optx-x}_{X}]
        -
        [J(\realoptx)-J(z)-\dualprod{q}{\realoptx-z}_{X}]
        \\
        \MoveEqLeft[-1]
        +
        [J(x)-J(z)-\dualprod{q}{x-z}_{X}],
    \end{aligned}
\]
which immediately gives the three-point identity
\begin{equation}
    \label{eq:gap:bregman:three-point}
    \dualprod{p-q}{x-\optx}_{X}
    =
    B_J^{p}(\optx, x)  - B_J^{q}(\optx, z) + B_J^q(x,z)
    \quad (\optx,x,z \in X,\, p \in \subdiff J(z),\, q \in \subdiff J(x)).
\end{equation}

If $X$ is a Hilbert space, we can take $J(x)=\frac{1}{2}\norm{x}^2$ to obtain $B_J^{x-z}(z, x)=\tilde B_J(z, x)=\frac{1}{2}\norm{z-x}_X^2$.
Therefore this three-point identity generalizes the classical three-point identity \eqref{eq:convergence:three-point-identity} in Hilbert spaces.
This could be used to generalize our convergence proofs to Banach spaces to treat methods of the general form
\[
    0 \in H(\nextx) + \subdiff_1 B_J^{\this q}(\nextx, \thisx),
\]
where $\subdiff_1$ denotes taking a subdifferential with respect to the first variable.
To see how \eqref{eq:gap:bregman:three-point} applies, observe that
\[
    \subdiff_1 B_J^{\this q}(\nextx, \thisx)=\subdiff J(\nextx)-\this q
    =\{\nexxt p - \this q \mid \nexxt q \in \subdiff J(\nextx)\}.
\]
This would, however, not provide convergence in norm but with respect to $B_J$.
For a general approach to primal-dual methods based on Bregman divergences, see \cite{tuomov-firstorder}.

Returning to our generic gap functional $\tilde \gap$ defined in \eqref{eq:gap:ergodic:generic-gap}, we have already observed in the proof of \cref{lem:gap:ergodic:functional} that $-\Xi \realoptx \in \subdiff (\tilde F + \tilde G)(\realoptx)$. Since due to the skew-adjointness of $\Xi$ we also have $\iprod{\Xi\realoptx}{x}_X=\iprod{\Xi\realoptx}{x-\realoptx}_X$ for a solution $\realoptx \in \inv H(0)$, this means that
\[
    \tilde \gap(x, \realoptx)=B_{\tilde G+\tilde F}^{-\Xi \realoptx}(x, \realoptx).
\]
In other words, the gap based at a solution $\realoptx \in \inv H(0)$ is also a Bregman divergence. In general, as we have already remarked, it can be zero for $x \ne \realoptx$.

\section{Convergence of function values}
\label{sec:gap:function}

We start with the fundamental algorithms: the proximal point method and explicit splitting. Throughout the following, we assume that $X$ is a Hilbert space and write $\minval{G} \defeq \min_{x \in X} G(x)$ whenever the minimum exists.

\begin{theorem}[proximal point method ergodic function value]
    \label{thm:gap:prox:ergodic-value}
    Let $G$ be proper, lower semicontinuous, and (strongly) convex with factor $\gamma\ge0$.
    Suppose $\inv{[\subdiff G]}(0) \ne \emptyset$.
    Pick an arbitrary $x^0 \in X$.
    Let $\tauTest_{k+1} \defeq \tauTest_k(1+\gamma\tau_k)$, and $\tauTest_0 \defeq 1$. For the iterates $\nextx \defeq \prox_{\tau_k G}(\thisx)$ of the proximal point method, define the \term[sequence!ergodic]{ergodic sequence}
    \begin{equation}
        \label{eq:gap:prox:ergodic-sequence}
        \tilde x^N \defeq \frac{1}{\zeta_N}\sum_{k=0}^{N-1} \tau_k\tauTest_k \nextx
        \quad\text{for}\quad
        \zeta_N \defeq
        \sum_{k=0}^{N-1} \tau_k\tauTest_k
        \quad (N \ge 1).
    \end{equation}
    \begin{enumerate}
        \item If $\tau_k \equiv \tau > 0$ and $G$ is not strongly convex ($\gamma=0$), then $G(\tilde x^N) \to \minval{G}$ at the rate $O(1/N)$.
        \item If $\tau_k \equiv \tau > 0$ and $G$ is strongly convex ($\gamma>0$), then $G(\tilde x^N) \to \minval{G}$ linearly.
        \item If $\tau_k \upto \infty$ and $G$ is strongly convex, then $G(\tilde x^N) \to \minval{G}$ superlinearly.
    \end{enumerate}
\end{theorem}

\begin{proof}
    Let the root $\realoptx \in \inv{[\subdiff G]}(0)$ be arbitrary; by assumption at least one exists. Then $\minval{G}=G(\realoptx)$ by \cref{thm:convex:fermat}.
    We recall that the proximal point iteration for minimizing $G$ can be written as
    \begin{equation}
        \label{eq:gap:prox}
        0 \in \tau_k\subdiff G(\nextx)+(\nextx-\thisx).
    \end{equation}
    As in the proof of \cref{thm:convergence:prox}, we test \eqref{eq:gap:prox} by the application of $\tauTest_k\iprod{\freevar}{\nextx-\realoptx}_X$ for the given testing parameter $\tauTest_k>0$ to obtain
    \begin{equation}
        \label{eq:gap:prox:tested}
        0 \in \tauTest_k\tau_k\iprod{\subdiff G(\nextx)}{\nextx-\realoptx}_X+\tauTest_k\iprod{\nextx-\thisx}{\nextx-\realoptx}_X.
    \end{equation}
    The next step will differ from the proof of \cref{thm:convergence:prox}, as we want an estimate for the function values.
    Indeed, by the subdifferential characterization of strong convexity, \cref{thm:smoothness:strong-convexity}\,\ref{item:strong-subdifferentiability},
    \begin{equation*}
        \iprod{\subdiff G(\nextx)}{\nextx-\realoptx}_X \ge G(\nextx)-G(\realoptx)+\frac{\gamma}{2}\norm{\nextx-\realoptx}_X^2.
    \end{equation*}
    Using this and the three-point-identity \eqref{eq:convergence:three-point-identity} in \eqref{eq:gap:prox:tested}, we obtain similarly to the proof of \cref{thm:testing:prox} the estimate
    \begin{equation*}
        \label{eq:gap:prox:est2-value}
        \frac{\tauTest_k(1+\tau_k\gamma)}{2}\norm{\nextx-\realoptx}_X^2
        +\tauTest_k\tau_k[G(\nextx)-G(\realoptx)]
        +\frac{\tauTest_k}{2}\norm{\nextx-\thisx}_X^2
        \le \frac{\tauTest_k}{2}\norm{\thisx-\realoptx}_X^2.
    \end{equation*}
    We now use the assumed recursion
    \begin{equation}
        \label{eq:gap:prox:metric-update}
        \tauTest_k(1+\tau_k\gamma) = \tauTest_{k+1}.
    \end{equation}
    (Observe the factor-of-two difference compared to \eqref{eq:testing:prox:metric-update}.)
    Thus
    \begin{equation*}
        \label{eq:gap:prox:test-before-sum}
        \frac{\tauTest_{k+1}}{2}\norm{\nextx-\realoptx}_X^2
        +\tauTest_k\tau_k[G(\nextx)-G(\realoptx)]
        +\frac{\tauTest_k}{2}\norm{\nextx-\thisx}_X^2
        \le \frac{\tauTest_k}{2}\norm{\thisx-\realoptx}_X^2.
    \end{equation*}
    Summing over $k=0,\ldots,N-1$ then yields
    \begin{equation}
        \label{eq:gap:prox:value-main}
        \frac{\tauTest_N}{2}\norm{x^N-\realoptx}_X^2
        +\sum_{k=0}^{N-1}
            \tauTest_k\tau_k[G(\nextx)-G(\realoptx)]
        + \sum_{k=0}^{N-1}
        \frac{\tauTest_k}{2}\norm{\nextx-\thisx}_X^2
        \le \frac{\tauTest_0}{2}\norm{x^0-\realoptx}_X^2 =: C_0.
    \end{equation}
    Using Jensen's inequality, it follows for the ergodic sequence defined in \eqref{eq:gap:prox:ergodic-sequence} that
    \begin{equation*}
        \zeta_N[G(\tilde x^N)-G(\realoptx)]
        \le C_0.
    \end{equation*}

    If $\tauTest_k \equiv \tauTest_0$ and $\gamma=0$, we therefore have that $\zeta_N=N\tauTest_0\tau$ and thus obtain $O(1/N)$ convergence of function values for the ergodic variable $\tilde x^N$.

    If $\tauTest_k \equiv \tauTest_0$ and $\gamma>0$, we deduce from \eqref{eq:gap:prox:metric-update} that $\zeta_N=\sum_{k=0}^{N-1}(1+\gamma\tau_k)^k\tau_k\tauTest_0$. This grows exponentially and hence we obtain the claimed linear convergence.

    Finally, if $\tau_k \to \infty$, we would similarly to \cref{thm:testing:prox}\,\ref{item:testing:prox:superlinear} obtain superlinear convergence if $\zeta_{N}/\zeta_{N+1} \to 0$ were to hold. To show this, we can write
    \[
        \frac{\zeta_{N}}{\zeta_{N+1}}
        =\frac{\sum_{k=0}^{N-1} \tauTest_k\tau_k}{\sum_{k=0}^N \tauTest_k\tau_k}
        =\frac{\sum_{k=0}^{N-1} \frac{\tauTest_k\tau_k}{\tauTest_N\tau_N}}{1+\sum_{k=0}^{N-1} \frac{\tauTest_k\tau_k}{\tauTest_N\tau_N}}
    \]
    So it suffices to show that $c_N \defeq \sum_{k=k_0}^{N-1} \frac{\tauTest_k\tau_k}{\tauTest_N\tau_N} \to 0$ as $N \to \infty$.
    This we obtain by estimating
    \[
        \begin{aligned}
        c_N
        &
        =
        \sum_{k=0}^{N-1} \frac{\tau_k/\tau_N}{\prod_{j=k}^{N-1}(1+\gamma\tau_j)}
        \le \sum_{k=0}^{N-1} \frac{(1+\gamma\tau_k)/(1+\gamma\tau_N)}{\prod_{j=k}^{N-1}(1+\gamma\tau_j)}
        \\
        &
        = \sum_{k=0}^{N-1} \frac{1}{\prod_{j=k+1}^{N}(1+\gamma \tau_j)}
        \le \sum_{k=0}^{N-1} (1+\gamma\tau_{k+1})^{-(N-k)}.
        \end{aligned}
    \]
    In the first and last step we have used that $\{\tau_k\}_{k \in \N}$ is increasing.
    Now we pick $a>1$ and find $k_0 \in \N$ such that $1+\gamma\tau_k \ge a$ for $k \ge k_0$. Then for $N > k_0$,
    \[
        c_N \le
        \sum_{k=0}^{k_0-1} (1+\gamma\tau_{k+1})^{-(N-k)}
        +
        \sum_{k=k_0}^{N-1} a^{-(N-k)}
        =
        \sum_{k=0}^{k_0-1} (1+\gamma\tau_{k+1})^{-(N-k)}
        +
        \sum_{j=1}^{N-k_0} a^{-j}.
    \]
    The first term goes to zero as $N \to \infty$ while the second term, as a geometric series, converges to $\inv a/(1-\inv a)$. We therefore deduce that $\lim_{N \to \infty} c_N \le \inv a/(1-\inv a)$. Letting $a \to \infty$, we see that $c_N \downto 0$.
\end{proof}

It is possible to improve the result to be nonergodic by showing that the proximal point method is in fact monotonic.

\begin{corollary}[proximal point method function value]
    \label{cor:gap:prox:value}
    The proximal point method is a descent method, i.e., $G(\nextx) \le G(\thisx)$ for all $k \in \N$.
    Therefore the convergence rates of \cref{thm:gap:prox:ergodic-value} also hold for $G(x^N) \to \minval{G}$.
\end{corollary}

\begin{proof}
    We know from \eqref{eq:gap:prox} that
    \[
        0 \le
        \inv\tau_k \norm{\nextx-\thisx}_X^2
        =
        \iprod{\subdiff G(\nextx)}{\thisx-\nextx}_X
        \le G(\thisx)-G(\nextx).
    \]
    This proves monotonicity. Now \eqref{eq:gap:prox:value-main} gives
    \[
        \zeta_N [G(x^N)-G(\realoptx)] \le C_0.
    \]
    Now we proceed using the growth estimates for $\zeta_N$ in the proof of \cref{thm:gap:prox:ergodic-value}.
\end{proof}

These results can be extended to the explicit splitting method,
\begin{algeqbox*}
    \begin{equation*}
        \nextx \defeq \prox_{\tau G}(\thisx - \tau \grad F(\thisx)),
    \end{equation*}
\end{algeqbox*}
in a straightforward manner.
In the next theorem, observe in comparison to \cref{thm:testing:fb} that $\tau L \le 1$ instead of $\tau L \le 2$. This kind of factor-of-two stricter step length or Lipschitz factor bound is a general feature of function value estimates of methods involving an explicit step, as well as of the gap estimates in the following sections. It stems from the corresponding difference between the value estimate \eqref{eq:smoothness:three-point:smoothness} and the non-value estimate \eqref{eq:smoothness:three-point:monotonicity} in \cref{cor:smoothness:three-point}.

\begin{theorem}[explicit splitting function value]
    \label{thm:gap:fb:value}
    Let $J \defeq F+G$ where $G: X \to \Rbar$ and $F: X \to \R$ are convex, proper, and lower semicontinuous, with $F$ moreover $L$-smooth.
    Suppose $\inv{[\subdiff J]}(0) \ne \emptyset$.
    If $\tau L \le 1$, the explicit splitting method satisfies $J(\tilde x^N) \to \minval{J}$ at the rate $O(1/N)$. If $G$ is strongly convex, then this convergence is linear.
\end{theorem}

\begin{proof}
    With $\tau_k \defeq \tau$, as usual, we write the method as
    \begin{equation}
        \label{eq:gap:fb:implicit}
        0 \in \tau_k[\subdiff G(\nextx) + \grad F(\thisx)] + (\nextx-\thisx).
    \end{equation}
    We then take arbitrary $\realoptx \in \inv{[\subdiff(F+G)]}(0)$ and use the three-point smoothness of $F$ proved in \cref{cor:smoothness:three-point}, and the subdifferential characterization of strong convexity of $G$, \cref{thm:smoothness:strong-convexity}\,\ref{item:strong-subdifferentiability}, to obtain
    \begin{equation*}
        \iprod{\subdiff G(\nextx)+\grad F(\thisx)}{\nextx-\realoptx}_X \ge J(\nextx)-J(\realoptx)+\frac{\gamma}{2}\norm{\nextx-\realoptx}_X^2 - \frac{L}{4}\norm{\nextx-\thisx}_X^2.
    \end{equation*}
    As in the proof of \cref{thm:gap:prox:ergodic-value}, after testing  \eqref{eq:gap:fb:implicit} by the application of $\tauTest_k\iprod{\freevar}{\nextx-\realoptx}_X$, we now obtain
    \begin{equation}
        \label{eq:gap:fb:value-ineq}
        \frac{\tauTest_{k+1}}{2}\norm{\nextx-\realoptx}_X^2
        +\tauTest_k\tau_k[J(\nextx)-J(\realoptx)]
        +\frac{\tauTest_k(1-\tau_k L)}{2}\norm{\nextx-\thisx}_X^2
        \le \frac{\tauTest_k}{2}\norm{\thisx-\realoptx}_X^2.
    \end{equation}
    Since $\tau_k L \le 1$, we may proceed as in \cref{thm:gap:prox:ergodic-value} to prove the ergodic convergences.
\end{proof}

Again, we can show nonergodic convergence due to the monotonicity of the iteration.

\begin{corollary}
    \label{thm:gap:fb:value:nonergodic}
    The convergence rates of \cref{thm:gap:fb:value} also hold for $J(x^N)\to J_{\min}$.
\end{corollary}

\begin{proof}
    We obtain from \eqref{eq:gap:fb:implicit} and the smoothness of $F$ (see \eqref{eq:smoothness:smoothness}) that
    \[
        \inv\tau_k \norm{\nextx-\thisx}_X^2
        =
        \iprod{\subdiff G(\nextx)+\grad F(\thisx)}{\thisx-\nextx}_X
        \le J(\thisx)-J(\nextx) + \frac{L}{2}\norm{\nextx-\thisx}_X^2.
    \]
    Since $L\tau_k \le 1 < 2$, we obtain monotonicity.
    The rest now follows as in \cref{thm:gap:prox:ergodic-value,cor:gap:prox:value}.
\end{proof}

\begin{remark}
    Based on \cref{cor:smoothness:three-point:sc}, any strong convexity of $F$ can also be used to obtain linear convergence by adapting the steps of the proof of \cref{thm:gap:fb:value}.
\end{remark}

\section{Ergodic gap estimates}
\label{sec:gap:ergodic}

We now study the convergence of gap functionals for general unaccelerated schemes of the form \eqref{eq:gap:ergodic:general-alg}. Since $\tilde G$ may in general not have the same factor $L$ of smoothness on all subspaces, we introduce the condition \eqref{eq:gap:ergodic:3smooth} of the next result. It is simply a version of the standard result of \cref{cor:smoothness:three-point} that allows a block-separable structure through the operator $\Lambda$ in place of the factor $L$; compare \cref{ex:testing:structured:3monotone-block}.

\begin{theorem}
    \label{thm:gap:ergodic:general}
    Let $H \defeq \subdiff \tilde G+\grad \tilde F+\Xi$, where $\Xi \in \linear(X; X)$ is skew-adjoint and $\tilde G: X \to \Rbar$ and $\tilde F: X \to \R$ are convex, proper, and lower semicontinuous. Suppose $\tilde F$ satisfies for some $\Lambda \in \linear(X; X)$ the three-point smoothness condition
    \begin{equation}
        \label{eq:gap:ergodic:3smooth}
        \iprod{\grad \tilde F(z)}{x-\optx}_X
        \ge
        \tilde F(x)-\tilde F(\optx)
        -\frac{1}{2}\norm{z-x}_\Lambda^2
        \quad (\optx, x,z \in X).
    \end{equation}
    Also let $\Precond \in \linear(X;X)$ be positive semi-definite and self-adjoint. Pick $x^0 \in X$, and let the sequence $\{\nextx\}_{k \in \N}$ be generated through the iterative solution of \eqref{eq:gap:ergodic:general-alg}.
    Then for every $\optx \in X$,
    \begin{equation}
        \label{eq:gap:ergodic:general-estimate}
        \frac{1}{2}\norm{x^N-\optx}_{\Precond}^2
        +\sum_{k=0}^{N-1}\left( \tilde\gap(\nextx; \optx) + \frac{1}{2}\norm{\nextx-\thisx}_{\Precond-\Lambda}^2\right)
        \le
        \frac{1}{2}\norm{x^0-\optx}_{\Precond}^2.
    \end{equation}
\end{theorem}

\begin{proof}
    By the convexity of $\tilde G$ we have
    \begin{equation}
        \label{eq:gap:ergodic:tested}
        \iprod{\subdiff \tilde G(x)}{x-\optx}_X
        \ge
        \tilde G(x)-\tilde G(\optx)
        \quad (x \in X).
    \end{equation}
    Using \eqref{eq:gap:ergodic:3smooth} and \eqref{eq:gap:ergodic:tested}, we obtain
    \begin{equation}
        \label{eq:gap:ergodic:first-estimate}
        \begin{aligned}[t]
        \iprod{\subdiff \tilde G(\nextx)&+\grad \tilde F(\thisx)+\Xi\nextx}{\nextx-\optx}_X
        \\
        &
        \ge
            (\tilde G+\tilde F)(\nextx)-(\tilde G + \tilde F)(\optx)
            +\iprod{\Xi\nextx}{\nextx-\optx}_X
            -\frac{1}{2}\norm{\thisx-\nextx}^2_{\Lambda}
        \\
        &
        =
        \tilde \gap(\nextx; \optx) -\frac{1}{2}\norm{\thisx-\nextx}^2_{\Lambda}.
        \end{aligned}
    \end{equation}
    In the final step we have also referred to the definition of $\tilde\gap$ in \eqref{eq:gap:ergodic:generic-gap} and the skew-adjointness of $\Xi$.

    From here on, our arguments are already standard:
    We \index{functional!testing}test \eqref{eq:gap:ergodic:general-alg} through the application of $\iprod{\freevar}{\nextx-\optx}_X$, obtaining
    \begin{equation*}
        0 \in \iprod{\subdiff \tilde G(\nextx)+\grad \tilde F(\thisx)+\Xi\nextx + \Precond(\nextx-\thisx)}{\nextx-\optx}_X.
    \end{equation*}
    Then we insert \eqref{eq:gap:ergodic:first-estimate}, which gives
    \begin{equation*}
        \frac{1}{2}\norm{\nextx-\optx}_{\Precond}^2
        +\tilde \gap(\nextx; \optx)
        +\frac{1}{2}\norm{\nextx-\thisx}_{\Precond-\Lambda}^2
        \le \frac{1}{2}\norm{\thisx-\optx}_{\Precond}^2.
    \end{equation*}
    Summing over $k=0,\ldots,N-1$ yields \eqref{eq:gap:ergodic:general-estimate}.
\end{proof}

In particular, we obtain the following corollary that shows that $\tilde \gap(\tilde x^N; \realoptx) \to  \tilde \gap(\realoptx; \realoptx) = 0$ at the rate $O(1/N)$ for any $\realoptx \in \inv H(0)$. Even further, taking any \emph{bounded} set $B \subset X$ such that $\inv H(0) \isect B \ne \emptyset$, we see that also the \term[gap!partial]{partial gap} $\tilde \gap(\tilde x^N; B) \to \tilde \gap(\realoptx; B) = 0$.

\begin{corollary}
    \label{cor:gap:ergodic:general}
    In \cref{thm:gap:ergodic:general}, suppose in addition that $\Precond \ge \Lambda$ and define the ergodic sequence
    \[
        \tilde x^N \defeq \frac{1}{N} \sum_{k=0}^{N-1} \nextx.
    \]
    Then
    \[
        \tilde\gap(\tilde x^N; \optx)
        \le
        \frac{1}{2N}\norm{x^0-\optx}_{\Precond}^2.
    \]
\end{corollary}

\begin{proof}
    This follows immediately from using $\Precond \ge \Lambda$ to eliminate the term $\frac{1}{2}\norm{\nextx-\optx}_{\Precond-\Lambda}^2$ from \eqref{eq:gap:ergodic:general-estimate} and then using Jensen's inequality on the gap.
\end{proof}

Due to the presence of $\Xi$, we cannot in general prove monotonicity of the abstract proximal point method and thus get rid of the ergodicity of the estimates.

\subsection*{Implicit splitting}

We now consider the solution of
\[
    \min_{x \in X}~ F(x)+G(x).
\]
Setting $B=\subdiff F$ and $A=\subdiff G$, \eqref{eq:convergence:drs-form}, the Douglas--Rachford or implicit splitting method can be written in the general form \eqref{eq:gap:ergodic:general-alg} with $u=(x,y,z)$,
\begin{align*}
    \tilde G(u) & \defeq \tau G(y)+ \tau F(x),
    &
    \tilde F & \equiv 0,
    \\
    \Xi & \defeq \begin{pmatrix} 0 & \Id & -\Id \\ -\Id & 0 & \Id \\ \Id & -\Id & 0 \end{pmatrix},
    \quad\text{and}
    &
    \Precond & \defeq
    \begin{pmatrix}
        0 & 0 & 0 \\
        0 & 0 & 0 \\
        0 & 0 & I
    \end{pmatrix}.
\end{align*}
Moreover,
\begin{equation}
    \label{eq:gap:drs:h}
    H(u) \defeq \subdiff \tilde G(u)+\Xi u.
\end{equation}
We then have the following ergodic estimate for
\[
    \gap_{\text{DRS}}(u; \realoptu)
    =[G(y)+F(x)] - [G(\realoptx)+F(\realoptx)] + \iprod{\realoptx-\realoptz}{x-y}_X \ge 0.
\]

\begin{theorem}
    \label{thm:gap:drs}
    Let $F:X\to\Rbar$ and $G:X\to\Rbar$ be proper, convex, and lower semicontinuous.
    Let $\realoptu \in \inv H(0)$ for $H$ given by \eqref{eq:gap:drs:h}.
    Then for any initial iterate $u^0=(x^0, y^0, z^0) \in X^3$, the iterates $\{\thisu\}_{k \in N}$ of the implicit splitting method \eqref{eq:splitting:dr} satisfy
    \[
        \gap_{\text{DRS}}(\tilde u^N; \realoptu)
        \le
        \frac{1}{2N\tau}\norm{u^1-\realoptu}_{\Precond}^2,
        \quad\text{where}\quad
        \tilde u^N \defeq \frac{1}{N} \sum_{k=0}^{N-1} \nextu.
    \]
\end{theorem}

\begin{proof}
    Clearly $\Precond$ is self-adjoint and positive semi-definite, and $\Precond \ge \Lambda \defeq 0$.
    The rest is clear from \cref{cor:gap:ergodic:general} by moving $\tau$ from $\tilde\gap$ on the right-hand side, and using that $\realoptx=\realopty$.
\end{proof}

Clearly, following the discussion in \cref{sec:gap:gap}, we can define a partial version of $\gap_{\text{DRS}}$ and obtain its convergence from \cref{thm:gap:drs}.

\subsection*{Primal-dual explicit splitting}

We recall that the PDES method \eqref{eq:splitting:gist} corresponds to \eqref{eq:gap:ergodic:general-alg} with \eqref{eq:gap:pd:setup:g} for the choice $F_0=0$ and $E=F$, while the preconditioning operator is given by
\[
    M \defeq \begin{pmatrix} \Id & 0 \\ 0 & \Id-KK^* \end{pmatrix}
\]
With this, we obtain the following estimate for the Lagrangian duality gap defined in \eqref{eq:gap:pd:gap}.

\begin{theorem}
    \label{thm:gap:pdes}
    Let $F:X\to\R$ and $G:Y\to\Rbar$ be proper, convex, and lower semicontinuous, and $K \in \linear(X; Y)$. Suppose $F$ is Gateaux differentiable with $L$-Lipschitz gradient for $L \le 1$, and that $\norm{K}_{\linear(X; Y)} \le 1$.
    Then for any initial iterate $u^0 \in X \times Y$, the iterates $\{\thisu=(\thisx,\thisy)\}_{k \in \N}$ of \eqref{eq:splitting:gist} satisfy for all $\optu=(\optx, \opty) \in X \times Y$ the ergodic gap estimate
    \[
        \gap(\tilde u^N; \optu)
        \le
        \frac{1}{2N}\norm{u^1-\optu}_{\Precond}^2,
        \quad\text{where}\quad
        \tilde u^N \defeq \frac{1}{N} \sum_{k=0}^{N-1} \nextu.
    \]
    In particular, if $B \subset X$ is bounded and $B \isect \inv H(0) \ne \emptyset$, the partial duality gap $\gap(u^N, B) \to 0$ at the rate $O(1/N)$.
\end{theorem}

\begin{proof}
    We use \cref{cor:gap:ergodic:general}.
    Using the assumed bound $\norm{K}_{\linear(X; Y)} \le 1$, clearly $\Precond$ is self-adjoint and positive semi-definite.
    By \cref{cor:smoothness:three-point}, the three-point smoothness condition \eqref{eq:gap:ergodic:3smooth} holds with
    $
        \Lambda \defeq \begin{psmallmatrix} L & 0 \\ 0 & 0 \end{psmallmatrix},
    $
    where $L$ is the Lipschitz factor of $\grad F$.
    Since $\norm{K}_{\linear(X; Y)} \le 1$ and $L \le 1$, we also verify $\Precond \ge \Lambda$.
    The rest now follows from \cref{cor:gap:ergodic:general} as well as the nonnegativity  of the partial duality gap \eqref{eq:gap:pd:partial}.
\end{proof}

\subsection*{Primal-dual proximal splitting}

We continue with the problem \eqref{eq:gap:pd:problem} and the corresponding structure \eqref{eq:gap:pd:setup} for $H$. We recall from \cref{thm:convergence:pd_conv,cor:convergence:pdps:forward} that for the unaccelerated PDPS we take the preconditioning operator as
\begin{equation}
    \label{eq:gap:pd:precond}
    \Precond \defeq \begin{pmatrix} \inv\tau \Id & -K^* \\  - K & \inv\sigma\Id \end{pmatrix}
\end{equation}%
for some primal and dual step length parameters $\tau,\sigma>0$.
We now obtain the following result for the Lagrangian duality gap defined in  \eqref{eq:gap:pd:gap}.

\begin{theorem}
    \label{thm:gap:pdps}
    Let $F_0:X\to\Rbar$, $E:X\to\R$, and $G:Y\to\Rbar$ be proper, convex, and lower semicontinuous, and $K \in \linear(X; Y)$. Suppose $E$ is Gateaux differentiable with $L$-Lipschitz gradient.
    Take $\sigma,\tau>0$ satisfying
    \[
        L\tau + \tau\sigma\norm{K}^2 < 1.
    \]
    Then for any initial iterate $u^0 \in X \times Y$ the iterates $\{\thisu=(\thisx,\thisy)\}_{k \in \N}$ of the PDPS method \eqref{eq:convergence:pdps:forward} satisfy for any $\optu=(\optx, \opty) \in X \times Y$  the ergodic gap estimate
    \[
        \gap(\tilde u^N; \optu)
        \le
        \frac{1}{2N\tau}\norm{u^1-\optu}_{\Precond}^2,
        \quad\text{where}\quad
        \tilde u^N \defeq \frac{1}{N} \sum_{k=0}^{N-1} \nextu.
    \]
    In particular, if $B \subset X$ is bounded and $B \isect \inv H(0) \ne \emptyset$, the partial duality gap $\gap(u^N, B) \to 0$ at the rate $O(1/N)$.
\end{theorem}

\begin{proof}
    We use \cref{cor:gap:ergodic:general}.
    By \cref{cor:smoothness:three-point}, the three-point smoothness condition \eqref{eq:gap:ergodic:3smooth} holds with
    $
        \Lambda \defeq \begin{psmallmatrix} L & 0 \\ 0 & 0 \end{psmallmatrix},
    $
    where $L$ is the Lipschitz factor of $\grad E$.
    In \cref{cor:convergence:pdps:forward} we have already proved that $\Test\Precond$ is self-adjoint and positive semi-definite. Similarly to the proof of the corollary, we verify that the condition $L\tau + \tau\sigma\norm{K}^2 < 1$ guarantees $\Precond \ge \Lambda$. (The only difference to the conditions in that result is the standard gap estimate factor-of-two difference in the term containing $L$.)
    The rest is clear from \cref{cor:gap:ergodic:general} as well as the nonnegativity of the partial duality gap \eqref{eq:gap:pd:partial}.
\end{proof}

\section{The testing approach in its general form}
\label{sec:testing:general}

We now want to produce gap estimates for accelerated methods. As we have seen in \cref{sec:testing:fundamental}, as an extension of \cref{eq:gap:ergodic:general-alg} these iteratively solve
\begin{equation}
    \label{eq:gap:ppext0}
    0 \in \Step_{k+1}[\subdiff \tilde G(\nextx)+\grad \tilde F(\thisx) + \Xi\nextx] + \Precond_{k+1}(\nextx-\thisx)
\end{equation}
for iteration-dependent step length and preconditioning operators $\Step_{k+1} \in \linear(X; X)$ and $\Precond_{k+1} \in \linear(X; X)$. We also introduced testing operators $\Test_{k+1} \in \linear(X; X)$ such that $\Test_{k+1}\Precond_{k+1}$ is self-adjoint and positive semi-definite.

Unless $\Test_{k+1}\Step_{k+1}$ is a scalar multiple of the identity, we will not  be able to extract in a straightforward way any of the gap functionals of \cref{sec:gap:gap} out of \eqref{eq:gap:ppext0}. Indeed, it is not clear how to provide a completely general approach to gap functionals of accelerated or otherwise complex algorithms. We will specifically see the difficulties when performing gap realignment for the accelerated PDPS in \cref{sec:gap:accel} and when developing very specific gap functionals for the ADMM in \cref{sec:gap:admm}.

For brevity in the following sections, we however do some general preparatory work. Observe that the method \eqref{eq:gap:ppext0} can be written more abstractly as
\begin{equation}
    \label{eq:gap:ppext}
    0 \in  \Happrox_{k+1}(\nextx) + \Precond_{k+1}(\nextx-\thisx)
\end{equation}
for some iteration-dependent set-valued function $\Happrox_{k+1}: X \setto X$.
The estimate \eqref{eq:gap:ppext:fundamental-condition} in the next theorem is in essence a quantitative or variable-metric version of the three-point smoothness and strong convexity estimate \eqref{eq:smoothness:three-point:smoothness-sc}.
The proof of the following result is already standard,
where the abstract value $\GenGap_{k+1}(\realoptx)$ models a suitable gap functional for iterate $\nextx$.

\begin{theorem}
    \label{thm:gap:ppext:convergence}
    On a Hilbert space $X$, let $\Happrox_{k+1}: X \setto X$, and $\Precond_{k+1}, \Test_{k+1} \in \linear(X; X)$ for $k \in \N$.
    Suppose \eqref{eq:gap:ppext} is solvable for the iterates $\{\thisx\}_{k \in \N}$.
    If $\Test_{k+1}\Precond_{k+1}$ is self-adjoint and
    \begin{equation}
        \label{eq:gap:ppext:fundamental-condition}
        \begin{aligned}[t]
            \iprod{\Happrox_{k+1}(\nextx)}{\nextx-\realoptx}_{\Test_{k+1}}
            &
            \ge
            \GenGap_{k+1}(\realoptx)
            + \frac{1}{2}\norm{\nextx-\realoptx}_{\Test_{k+2}\Precond_{k+2}-\Test_{k+1}\Precond_{k+1}}^2
            \\ \MoveEqLeft[-1]
            - \frac{1}{2}\norm{\nextx-\thisx}_{\Test_{k+1} \Precond_{k+1}}^2
        \end{aligned}
    \end{equation}
    for all $k \in \N$ and some $\realoptx \in X$ and $\GenGap_{k+1}(\realoptx) \in \R$, then both
    \begin{equation}
        \label{eq:gap:ppext:quantitative-fejer}
        \frac{1}{2}\norm{\nextx-\realoptx}_{\Test_{k+2}\Precond_{k+2}}^2
        + \GenGap_{k+1}(\realoptx)
        \le
        \frac{1}{2}\norm{\thisx-\realoptx}_{\Test_{k+1}\Precond_{k+1}}^2
        \quad
        (k \in \N)
    \end{equation}
    and
    \begin{equation}
        \label{eq:gap:ppext:convergence}
        \frac{1}{2}\norm{x^N-\realoptx}^2_{\Test_{N+1}\Precond_{N+1}}
        +
        \sum_{k=0}^{N-1} \GenGap_{k+1}(\realoptx)
        \le
        \frac{1}{2}\norm{x^0-\realoptx}^2_{\Test_{1}\Precond_{1}}
        \quad
        (N \ge 1).
    \end{equation}
\end{theorem}

\begin{proof}%
    Inserting \eqref{eq:gap:ppext} into \eqref{eq:gap:ppext:fundamental-condition}, we obtain
    \begin{equation}
        \begin{aligned}[t]
            \label{eq:gap:ppext:fundamental-condition-transformed}
            - \iprod{\nextx-\thisx}{\nextx-\realoptx}_{\Test_{k+1}\Precond_{k+1}}
            &
            \ge
            \frac{1}{2}\norm{\nextx-\realoptx}_{\Test_{k+2}\Precond_{k+2}-\Test_{k+1}\Precond_{k+1}}^2
            \\ \MoveEqLeft[-1]
            -\frac{1}{2}\norm{\nextx-\thisx}_{\Test_{k+1} \Precond_{k+1}}^2
            + \GenGap_{k+1}(\realoptx).
        \end{aligned}
    \end{equation}
    We recall for general self-adjoint $M$ the three-point formula \eqref{eq:convergence:three-point-identity}, i.e.,
    \begin{equation*}
       \iprod{\nextx-\thisx}{\nextx-\realoptx}_{M}
       = \frac{1}{2}\norm{\nextx-\thisx}_{M}^2
           - \frac{1}{2}\norm{\thisx-\realoptx}_{M}^2
           + \frac{1}{2}\norm{\nextx-\realoptx}_{M}^2.
    \end{equation*}
    Using this  with $M=\Test_{k+1}\Precond_{k+1}$, we rewrite \eqref{eq:gap:ppext:fundamental-condition-transformed} as \eqref{eq:gap:ppext:quantitative-fejer}.
    Summing \eqref{eq:gap:ppext:quantitative-fejer} over $k=0,\ldots,N-1$, we obtain \eqref{eq:gap:ppext:convergence}.
\end{proof}

\section{Ergodic gaps for accelerated primal-dual methods}
\label{sec:gap:accel}

To derive ergodic gap estimates for the accelerated primal-dual proximal splitting of \cref{thm:testing:pdps:accel}, we need to perform significant additional work due to the fact that $\eta_k \defeq \tauTest_k\tau_k \ne \sigmaTest_{k+1}\sigma_{k+1}$. The overall idea of the proof remains the same, but we need to pay special attention to the blockwise structure of the problem and to do some realignment of the blocks to get the same factor $\eta_k$ in front of both $G$ and $F$.

\subsection*{Duality gap realignment}

We continue with the problem \eqref{eq:gap:pd:problem} and the setup \eqref{eq:gap:pd:setup}.
Working with the general scheme \eqref{eq:gap:ppext}, we write
\begin{subequations}
\label{eq:gap:accel:setup}
\begin{equation}
    \Happrox_{k+1}(u)
    \defeq \Step_{k+1}(\subdiff \tilde G(\nextu)+\grad \tilde F(\thisu)+\Xi)
\end{equation}
taking as in \cref{thm:testing:pdps:accel} the testing and step length operators
\begin{gather}
    \label{eq:gap:accel:step-test}
    \Step_{k+1} \defeq \begin{pmatrix} \tau_k \Id & 0 \\ 0 & \sigma_{k+1} \Id \end{pmatrix}
    \quad\text{and}\quad
    \Test_{k+1} \defeq \begin{pmatrix} \tauTest_k \Id & 0 \\ 0 & \sigmaTest_{k+1} \Id \end{pmatrix}
\end{gather}
for some step length and testing parameters $\tau_k,\sigma_{k+1},\tauTest_k,\sigma_{k+1}>0$.
Throughout this section we also take
\begin{equation}
    \label{eq:gap:accel:gamma-lambda}
    \Gamma \defeq \begin{pmatrix} \gamma \cdot \Id & 0 \\ 0 & \rho \cdot \Id\end{pmatrix}
    \quad\text{and}\quad
    \Lambda \defeq \begin{pmatrix} L \cdot \Id & 0 \\ 0 & 0 \end{pmatrix}.
\end{equation}
\end{subequations}

For the moment, we do not yet need to know the specific structure of $\Precond_{k+1}$; hence the following estimates apply not only to the PDPS method but also to the PDES method and its potential accelerated variants.

\begin{lemma}
    \label{lemma:gap:accel:h}
    Let us be given $K \in \linear(X; Y)$, $F=F_0+E$ with $F_0: X \to \Rbar$, $E: X \to \R$, and $G^*: Y \to \Rbar$ convex, proper, and lower semicontinuous on Hilbert spaces $X$ and $Y$.
    Suppose $F_0$ and $G^*$ are (strongly) convex for some $\gamma, \rho \ge 0$, and $E$ has $L$-Lipschitz continuous gradient.
    With the setup of \eqref{eq:gap:pd:setup} and \eqref{eq:gap:accel:setup}, for any $u, \realoptu \in X \times Y$ and any $k \in \N$ we have
    \[
         \iprod{\Happrox_{k+1}(u)}{u-\realoptu}_{\Test_{k+1}}
         \ge
         \gap_{k+1}(u; \realoptu)
         +\frac{1}{2}\norm{u-\realoptu}_{\Test_{k+1}\Step_{k+1}(2\Xi+\Gamma)}^2
         -\frac{1}{4}\norm{u-\thisu}_{\Test_{k+1}\Step_{k+1}\Lambda}^2
    \]
    for
    \[
        \begin{aligned}
        \gap_{k+1}(u; \realoptu)
        &
        \defeq
        \tauTest_k\tau_k(F(x)-F(\realoptx))
        + \sigmaTest_{k+1}\sigma_{k+1}(G^*(y)-G^*(\realopty))
        \\ \MoveEqLeft[-1]
        + \iprod{(\tauTest_k \tau_kK^*)\realopty}{x}_X
        - \iprod{(\sigmaTest_{k+1} \sigma_{k+1}K)\realoptx}{y}_Y
        - \iprod{(K\tauTest_k\tau_k-\sigmaTest_{k+1}\sigma_{k+1}K)\realoptx}{\realopty}_Y.
        \end{aligned}
    \]
\end{lemma}

\begin{proof}
    Expanding $\Happrox_{k+1}$, we have
    \begin{equation*}
        \begin{aligned}
            \iprod{\Happrox_{k+1}(u)}{u-\realoptu}_{\Test_{k+1}}
            & =
            \tauTest_k\tau_k \iprod{\subdiff F_0(x)}{x-\realoptx}_X
            \\ \MoveEqLeft[-1]
            + \tauTest_k\tau_k  \iprod{\grad E(\thisx)}{x-\realoptx}_X
            \\ \MoveEqLeft[-1]
            + \sigmaTest_{k+1}\sigma_{k+1} \iprod{\subdiff G^*(y)}{y-\realopty}_Y
            \\ \MoveEqLeft[-1]
            + \iprod{(\tauTest_k \tau_kK^*)y}{x-\realoptx}_X
            - \iprod{(\sigmaTest_{k+1} \sigma_{k+1}K)x}{y-\realopty}_Y.
        \end{aligned}
    \end{equation*}
    Observe that
    \begin{equation*}
        \begin{aligned}
        \iprod{(\tauTest_k \tau_kK^*)y}{x-\realoptx}_X &
        - \iprod{(\sigmaTest_{k+1} \sigma_{k+1}K)x}{y-\realopty}_Y
        \\
        &
        =
        \iprod{(K\tauTest_k\tau_k-\sigmaTest_{k+1}\sigma_{k+1}K)(x-\realoptx)}{y-\realopty}_Y
        \\ \MoveEqLeft[-1]
        +
        \iprod{(\tauTest_k \tau_kK^*)\realopty}{x-\realoptx}_X
        - \iprod{(\sigmaTest_{k+1} \sigma_{k+1}K)\realoptx}{y-\realopty}_Y
        \\
        &
        =
        \frac{1}{2}\norm{u-\realoptu}_{2\Test_{k+1}\Step_{k+1}\Xi}^2
        - \iprod{(K\tauTest_k\tau_k-\sigmaTest_{k+1}\sigma_{k+1}K)\realoptx}{\realopty}_Y
        \\ \MoveEqLeft[-1]
        + \iprod{(\tauTest_k \tau_kK^*)\realopty}{x}_X
        - \iprod{(\sigmaTest_{k+1} \sigma_{k+1}K)\realoptx}{y}_Y.
        \end{aligned}
    \end{equation*}
    Therefore
    \begin{equation}
        \label{eq:gap:accel:estimate0}
        \begin{aligned}[t]
            \iprod{\Happrox_{k+1}(u)}{u-\realoptu}_{\Test_{k+1}}
            & =
            \tauTest_k\tau_k \iprod{\subdiff F_0(x)}{x-\realoptx}_X
            \\ \MoveEqLeft[-1]
            + \tauTest_k\tau_k  \iprod{\grad E(\thisx)}{x-\realoptx}_X
            \\ \MoveEqLeft[-1]
            + \sigmaTest_{k+1}\sigma_{k+1} \iprod{\subdiff G^*(y)}{y-\realopty}_Y
            \\ \MoveEqLeft[-1]
            + \frac{1}{2}\norm{u-\realoptu}_{2\Test_{k+1}\Step_{k+1}\Xi}^2
            - \iprod{(K\tauTest_k\tau_k-\sigmaTest_{k+1}\sigma_{k+1}K)\realoptx}{\realopty}_Y
            \\ \MoveEqLeft[-1]
            + \iprod{(\tauTest_k \tau_kK^*)\realopty}{x}_X
            - \iprod{(\sigmaTest_{k+1} \sigma_{k+1}K)\realoptx}{y}_Y.
        \end{aligned}
    \end{equation}
    Due to the smoothness three-point corollaries, specifically \eqref{eq:smoothness:three-point:smoothness}, we have
    \begin{subequations}%
    \label{eq:gap:accel:convexity-est}%
    \begin{equation}
        \iprod{\grad E(\thisx)}{x-\realoptx}_X
        \ge
        E(x)-E(\realoptx)
        -\frac{L}{2}\norm{x-\thisx}_X^2.
    \end{equation}
    Also, by the (strong) convexity of $F_0$, we have
    \begin{equation}
        \iprod{\subdiff F_0(x)}{x-\realoptx}_X
        \ge F_0(x) - F_0(\realoptx) + \frac{\gamma}{2}\norm{x-\realoptx}_X^2,
    \end{equation}%
    as well as by the (strong) convexity of $G^*$
    \begin{equation}
        \iprod{\subdiff G^*(y)}{y-\realopty}_Y
        \ge G^*(y) - G^*(\realopty) + \frac{\rho}{2}\norm{y-­\realopty}_Y^2.
    \end{equation}%
    \end{subequations}%
    Applying these estimates in  \eqref{eq:gap:accel:estimate0}, and using the structure \eqref{eq:gap:accel:step-test} and \eqref{eq:gap:accel:gamma-lambda} of the involved operators, we obtain the claim.
\end{proof}

If $\tauTest_k\tau_k=\sigmaTest_{k+1}\sigma_{k+1}$, clearly $\gap_{k+1}(\nextu; \realoptu) \ge \tauTest_k\tau_k \gap(\nextu)$. This is the case in the unaccelerated case already considered in \cref{thm:gap:pdes,thm:gap:pdps}. Some specific stochastic accelerated algorithms also satisfy this \cite[see][]{tuomov-blockcp}. Applying the techniques of \cref{sec:gap:ergodic}, we could then use Jensen's inequality to estimate $\sum_{k=0}^{n-1} \gap_{k+1}(\nextu; \realoptu) \ge \sum_{k=0}^{N-1} \tauTest_k\tau_k \gap(\nextu)$ further from below to obtain a gap on suitable ergodic sequences. However, in our primary accelerated algorithm of interest, the PDPS method, instead $\tauTest_k\tau_k=\sigmaTest_k\sigma_k$. We will therefore have to do some rearrangements.

\begin{lemma}
    \label{lemma:gap:accel:htwo}
    Let $K \in \linear(X; Y)$, $F=F_0+E$ with $F_0: X \to \Rbar$, $E: X \to \R$, and $G^*: Y \to \Rbar$ convex, proper, and lower semicontinuous on Hilbert spaces $X$ and $Y$.
    Suppose $F_0$ and $G^*$ are (strongly) convex for some $\gamma, \rho \ge 0$, and $E$ has $L$-Lipschitz gradient.
    With the setup of \eqref{eq:gap:pd:setup} and \eqref{eq:gap:accel:setup}, suppose $\tauTest_k\tau_k=\sigmaTest_k\sigma_k$.
    If $\realoptu \in \inv H(0)$, then for any $N\in\N$ it holds that
    \begin{equation}
        \label{eq:gap:accel:htwo0}
        \begin{aligned}[t]
         \iprod{\Happrox_{k+1}(\nextu)}{\nextu-\realoptu}_{\Test_{k+1}}
         &
         \ge
         \gap_{*,k+1}(\nextx, \thisy; \realoptu)
         + \frac{1}{2}\norm{\nextu-\realoptu}_{\Test_{k+1}\Step_{k+1}(2\Xi+\Gamma)}^2
         \\ \MoveEqLeft[-1]
         -\frac{1}{2}\norm{\nextu-\thisu}_{\Test_{k+1}\Step_{k+1}\Lambda}^2
         \quad (k=0,\dots,N-1)
        \end{aligned}
    \end{equation}
    for some $\gap_{*, k+1}(\nextx, \thisy; \realoptu)$ satisfying with $\gap$ given by \eqref{eq:gap:pd:gap} the estimate
    \begin{equation}
       \label{eq:gap:accel:htwo-estimate}
        \sum_{k=0}^{N-1} \gap_{*, k+1}(\nextx, \thisy; \realoptu)
        \ge \sum_{k=1}^{N-1} \tauTest_k\tau_k \gap(\nextx, \thisy; \realoptu).
    \end{equation}
\end{lemma}

\begin{proof}
    First, note that \eqref{eq:gap:accel:htwo0} holds for
    \[
        \begin{aligned}
        \gap_{*,k+1}(\nextx, \thisy; \realoptu)
        &
        \defeq
        \inf_{w^{k+1}\in \Happrox_{k+1}(\nextu)}
        \iprod{w^{k+1}}{\nextu-\realoptu}_{\Test_{k+1}}
        \\
        \MoveEqLeft[-1]
        -
        \frac{1}{2}\norm{\nextu-\realoptu}_{\Test_{k+1}\Step_{k+1}(2\Xi+\Gamma)}^2
        +\frac{1}{2}\norm{\nextu-\thisu}_{\Test_{k+1}\Step_{k+1}\Lambda}^2.
        \end{aligned}
    \]
    It remains to prove the estimate \eqref{eq:gap:accel:htwo-estimate} for this choice.

    With $N \ge 1$, let us define the set
    \begin{equation*}
        \begin{aligned}
            S_N & \defeq \sum_{k=0}^{N-1}
            \Bigl(\iprod{\Happrox_{k+1}(\nextu)}{\nextu-\realoptu}_{\Test_{k+1}} -\frac{1}{2}\norm{\nextu-\realoptu}_{\Test_{k+1}\Step_{k+1}(2\Xi+\Gamma)}^2+\frac{1}{2}\norm{\nextu-\thisu}_{\Test_{k+1}\Step_{k+1}\Lambda}^2\Bigr)
            \\
            &
            =
            \sum_{k=0}^{N-1}
            \Biggl(\tauTest_k\tau_k \Bigl(\iprod{\subdiff F_0(\nextx)+\grad E(\thisx)}{\nextx-\realoptx}_X
            -\frac{\gamma}{2}\norm{\nextx-\realoptx}_X^2
            +\frac{L}{2}\norm{\nextx-\thisx}_X^2\Bigr)
            \\ \MoveEqLeft[-1] \qquad\quad
            + \sigmaTest_{k+1}\sigma_{k+1}\Bigl(\iprod{\subdiff G^*(\nexty)}{\nexty-\realopty}_Y
            -\frac{\rho}{2}\norm{\nexty-\realopty}_Y^2\Bigr)\Biggr).
        \end{aligned}
    \end{equation*}
    Observe that in the second expression, $\Test_{k+1}\Step_{k+1}\Xi$ has canceled the corresponding component of $\Happrox_{k+1}$.
    Then it is enough to prove that $S_N \ge \sum_{k=1}^{N-1} \tauTest_k\tau_k \gap(\nextx, \thisy; \realoptu)$.
    To do this, we need to shift $\nexty$ to $\thisy$. With $N \ge 2$, we therefore rearrange terms to obtain
    \[
        S_N=A_N+B_N
    \]
    for
    \begin{align*}
        A_N
        &
        =
        \tauTest_0\tau_0 \Bigl( \iprod{\subdiff F_0(x^1)+\grad E(x^0)}{x^1-\realoptx}_X-\frac{\gamma}{2}\norm{x^1-\realoptx}_X^2+\frac{L}{2}\norm{x^1-x^0}_X^2\Bigr)
        \\ \MoveEqLeft[-1]
        + \sigmaTest_{N}\sigma_{N}\Bigl(\iprod{\subdiff G^*(y^N)}{y^N-\realopty}_Y-\frac{\rho}{2}\norm{y^N-\realopty}_Y^2\Bigr)
        \\ \MoveEqLeft[-1]
        - \iprod{(K\tauTest_0\tau_0-\sigmaTest_{N}\sigma_{N}K)\realoptx}{\realopty}_Y
        + \iprod{(\tauTest_0 \tau_0 K^*)\realopty}{x^1}_X
        - \iprod{(\sigmaTest_{N} \sigma_{N}K)\realoptx}{y^N}_Y
    \shortintertext{and}
        B_N & \defeq
        \sum_{k=1}^{N-1}\Biggl(
        \tauTest_k\tau_k \Bigl( \iprod{\subdiff F_0(\nextx)+\grad E(\thisx)}{\nextx-\realoptx}_X
        -\frac{\gamma}{2}\norm{\nextx-\realoptx}_X^2 + \frac{L}{2}\norm{\nextx-\thisx}_X^2\Bigr)
        \\ \MoveEqLeft[-1] \qquad\quad
        + \sigmaTest_{k}\sigma_{k}\Bigl(
            \iprod{\subdiff G^*(\thisy)}{\thisy-\realopty}_Y
        -\frac{\rho}{2}\norm{\nexty-\realopty}_Y^2
        \Bigr)
        \\ \MoveEqLeft[-1] \qquad\quad
        + \iprod{(\tauTest_k \tau_kK^*)\realopty}{\nextx}_X
        - \iprod{(\sigmaTest_{k} \sigma_{k}K)\realoptx}{\thisy}_Y\Biggr).
    \end{align*}
    Observe that we only sum over $k=1,\ldots,N-1$ instead of $k=0,\ldots,N-1$.

    We can now use \eqref{eq:gap:accel:convexity-est} and our assumption $\tauTest_k\tau_k=\sigmaTest_k\sigma_k$ to estimate
    \begin{equation}
        \label{eq:gap:accel:bn-est}
        B_N \ge \sum_{k=1}^{N-1}
            \tauTest_k\tau_k \gap(\nextx, \thisy).
    \end{equation}
    By \cref{cor:smoothness:three-point}, $E$ satisfies the three-point monotonicity estimate \cref{eq:smoothness:three-point:monotonicity}; in particular,
    \[
        \iprod{\grad E(x^0)-\grad E(\realoptx)}{x^1-\realoptx}_X \ge -\frac{L}{2}\norm{x^1-x^0}_X^2.
    \]
    Since $K^*\realoptx \in \subdiff G^*(\realopty)$, and $-K\realopty \in \subdiff F_0(\realoptx)+\grad E(\realoptx)$, and $\subdiff F_0$ and $\subdiff G$ are strongly monotone, we also obtain
    \begin{align*}
        \iprod{\subdiff F_0(x^1)+\grad E(\realoptx)+K^*\realopty}{x^1-\realoptx}_X
        -\frac{\gamma}{2}\norm{x^1-\realoptx}_X^2
        & \ge
        0
        \quad\text{and}\quad
        \\
        \iprod{\subdiff G^*(y^N)-K\realoptx}{y^N-\realopty}_Y
        -\frac{\rho}{2}\norm{y^N-\realopty}_Y^2
        & \ge
        0.
    \end{align*}
    Rearranging and using these estimates we obtain
    \begin{equation}
        \label{eq:gap:accel:an-est}
        \begin{aligned}[t]
        A_N & =
        \tauTest_0\tau_0\Bigl(
            \iprod{\subdiff F_0(x^1)+\grad E(x^0)+K^*\realopty}{x^1-\realoptx}_X
            -\frac{\gamma}{2}\norm{x^1-\realoptx}_X^2+\frac{L}{2}\norm{x^1-x^0}_X^2
        \Bigr)
        \\ \MoveEqLeft[-1]
        + \sigmaTest_{N}\sigma_{N}\Bigl(
            \iprod{\subdiff G^*(y^N)-K\realoptx}{y^N-\realopty}_Y
            -\frac{\gamma}{2}\norm{y^N-\realopty}_Y^2
        \Bigr)
        \ge 0.
        \end{aligned}
    \end{equation}
    The estimates \eqref{eq:gap:accel:bn-est} and \eqref{eq:gap:accel:an-est} finally give $S_N \ge \sum_{k=1}^{N-1} \tauTest_k\tau_k \gap(\nextx, \thisy; \realoptu)$ as we set out to prove.
\end{proof}

In the proof of \cref{lemma:gap:accel:htwo}, we required $\realoptu \in \inv H(0)$ to show that $A_N \ge 0$. Therefore, as the estimate \eqref{eq:gap:accel:htwo0} will not hold for an arbitrary base point $\optu$ in place of $\realoptu$, we will not be able to obtain for accelerated methods the convergence of the \term[gap!duality!partial]{partial duality gap} \cref{eq:gap:pd:partial} that converges for unaccelerated methods.

The next theorem is our main result regarding ergodic gaps for general accelerated methods.
As $\gamma$ and $\rho$ feature as acceleration parameters in algorithms, the conditions of this theorem imply that gap estimates require slower acceleration.

\begin{theorem}
    \label{thm:gap:accel:convergence}
    Let $K \in \linear(X; Y)$, $F=F_0+E$ with $F_0: X \to \Rbar$, $E: X \to \R$, and $G^*: Y \to \Rbar$ convex, proper, and lower semicontinuous on Hilbert spaces $X$ and $Y$.
    Suppose $F_0$ and $G^*$ are (strongly) convex for some $\gamma, \rho \ge 0$, and $E$ has $L$-Lipschitz gradient.
    Assume the setup \eqref{eq:gap:pd:setup} and \eqref{eq:gap:accel:setup}. For each $k \in \N$, also take $\Precond_{k+1} \in \linear(X \times Y; X \times Y)$ such that $\Test_{k+1}\Precond_{k+1}$ is self-adjoint.
    Pick an initial iterate $u^0 \in X \times Y$ and suppose $\{\nextu=(\nextx,\nexty)\}_{k \in \N}$ are generated by \eqref{eq:gap:ppext}. Let $\realoptu=(\realoptx, \realopty) \in \inv H(0)$.
    If $\tauTest_k\tau_k=\sigmaTest_k\sigma_k$, and
    \begin{equation}
        \label{eq:gap:accel:fundamental-condition}
        \frac{1}{2}\norm{\nextu-\thisu}_{\Test_{k+1}(\Precond_{k+1}-\Step_{k+1}\Lambda)}^2
        + \frac{1}{2}\norm{\nextu-\realoptu}_{\Test_{k+1}(\Precond_{k+1}+\Step_{k+1}(2\Xi+\Gamma))-\Test_{k+2}\Precond_{k+2}}^2
        \ge
        0,
    \end{equation}
    then
    \begin{equation}
        \label{eq:gap:accel:convergence}
        \frac{1}{2}\norm{u^N-\realoptu}^2_{\Test_{N+1}\Precond_{N+1}}
        +  \zeta_{*,N} \gap(\tilde x_{*,N}, \tilde y_{*,N}; \realoptu)
        \le
        \norm{u^0-\realoptu}_{\Test_1 \Precond_1}^2
        \quad (N \ge 2)
    \end{equation}
    for $\gap$ given by \eqref{eq:gap:pd:gap} and the ergodic sequences
    \begin{equation*}
        \label{eq:gap:accel:ergodicseq}
        \tilde x_{*,N} \defeq \inv\zeta_{*,N}\sum_{k=1}^{N-1} \tau_k\tauTest_k \nextx
        \quad\text{and}\quad
        \tilde y_{*,N} \defeq \inv\zeta_{*,N}\sum_{k=1}^{N-1} \sigma_k\sigmaTest_k \thisy
        \quad\text{for}\quad
         \zeta_{*,N} \defeq \sum_{k=1}^{N-1} \eta_k.
    \end{equation*}
\end{theorem}

\begin{proof}
    Using \eqref{eq:gap:accel:htwo0} in \eqref{eq:gap:accel:fundamental-condition}, we obtain \eqref{eq:gap:ppext:fundamental-condition} for $\GenGap_{k+1}(\realoptu) \defeq \gap_{*, k+1}(\nextx, \thisy; \realoptu)$.
    By Jensen's inequality,
    \[
        \sum_{k=0}^{N-1} \gap_{*, k+1}(\nextx, \thisy; \realoptu)
        \ge  \zeta_{*,N} \gap(\tilde x_{*,N}, \tilde y_{*,N}; \realoptu).
    \]
    We therefore obtain \eqref{eq:gap:accel:convergence} from \eqref{eq:gap:ppext:convergence} in \cref{thm:gap:ppext:convergence}.
\end{proof}

\subsection*{Accelerated primal-dual proximal splitting}

We now obtain gap estimates for the accelerated PDPS method.
Observe the factor-of-two differences in the definitions of $\omega_k$ and in the initial conditions for the step lengths in the following theorem compared to \cref{thm:testing:pdps:accel}. Because strong convexity with factor $\gamma$ implies strong convexity with the factor $\gamma/2$, the conditions and step length rules of this theorem imply the iterate convergence results of \cref{cor:convergence:pdps:forward,thm:testing:pdps:accel} as well.

\begin{theorem}[gap estimates for PDPS]
    \label{thm:gap:accel:pdps}
    Let $F_0: X \to \Rbar$, $E: X \to \R$ and $G: Y \to \Rbar$ be convex, proper, and lower semicontinuous on Hilbert spaces $X$ and $Y$ with $\grad E$ $L$-Lipschitz. Also let $K \in \linear(X; Y)$ and let $\realoptu=(\realoptx, \realopty)$ be a primal-dual solution to the problem \eqref{eq:gap:pd:problem}.
    Pick initial step lengths $\tau_0, \sigma_0>0$ subject to $L\tau_0 + \tau_0\sigma_0\norm{K}_{\linear(X; Y)}^2<1$.
    For any initial iterate $u^0 \in X \times Y$, suppose $\{\nextu\}_{k \in \N}$ are generated by the (accelerated) PDPS method \cref{eq:testing:pdps:forward}.
    Let the Lagrangian duality gap functional $\gap$ be given by \eqref{eq:gap:pd:gap}, and the ergodic iterates $\tilde x_{*,N}$ and $\tilde y_{*,N}$ by \eqref{eq:gap:accel:ergodicseq}.

    \begin{enumerate}[label=(\roman*)]
        \item If we take $\tau_k \equiv \tau_0$ and $\sigma_k \equiv \sigma_0$, then the ergodic gap $\gap(\tilde x_{*,N}, \tilde y_{*,N}; \realoptu) \to 0$ at the rate $O(1/N)$.

        \item If $F_0$ is strongly convex with factor $\gamma>0$, and we take
        \begin{equation*}
            \omega_k \defeq 1/\sqrt{1+\gamma\tau_k},
            \quad
            \tau_{k+1} \defeq \tau_k\omega_k,
            \quad\text{and}\quad
            \sigma_{k+1} \defeq \sigma_k/\omega_k,
        \end{equation*}
        then $\gap(\tilde x_{*,N}, \tilde y_{*,N}; \realoptu) \to 0$ at the rate $O(1/N^2)$.

        \item If both $F_0$ and $G^*$ are strongly convex with respective factors $\gamma>0$ and $\rho>0$, and we take
        \begin{equation*}
            \omega_k \defeq 1/\sqrt{1+\theta},
            \quad
            \theta \defeq \min\{\rho\sigma_0,\gamma\tau_0\},
            \quad
            \tau_k \defeq \tau_0
            \quad\text{and}\quad
            \sigma_k \defeq \sigma_0,
        \end{equation*}
        then $\gap(\tilde x_{*,N}, \tilde y_{*,N}; \realoptu) \to 0$ linearly.
    \end{enumerate}
\end{theorem}

\begin{proof}
    We use \cref{thm:gap:accel:convergence} in place of \cref{thm:testing:structured:convergence} in the proof of \cref{thm:testing:pdps:accel}.
    We recall that the latter consists of showing $\Test_{k+1}\Precond_{k+1}$ to be self-adjoint and \eqref{eq:testing:structured:metric-update} and $\Precond \ge \Lambda$  to hold, i.e.,
    \begin{align*}
        \Test_{k+1}(\Precond_{k+1}+2\Step_{k+1}\Gamma) & \succeq \Test_{k+2}\Precond_{k+2},
        \quad\text{and}&
        \Test_{k+1}(\Precond_{k+1}-\Step_{k+1}\Lambda/2) & \succeq 0,
    \intertext{Now, to prove \eqref{eq:gap:accel:fundamental-condition}, we instead prove the self-adjointness as well as}
        \Test_{k+1}(\Precond_{k+1}+\Step_{k+1}\Gamma) & \succeq \Test_{k+2}\Precond_{k+2},
        \quad\text{and}&
        \Test_{k+1}(\Precond_{k+1}-\Step_{k+1}\Lambda) & \succeq 0.
    \end{align*}
    These all follows from the proof of \cref{thm:testing:pdps:accel} with the factor-of-two differences in the formulas for $\omega_k$ and the initialization condition apparent from the statements of these two theorems. The proof of \cref{thm:testing:pdps:accel} also verifies that $\tauTest_k\tau_k=\sigmaTest_k\sigma_k$.

    All the conditions \cref{thm:gap:accel:convergence} are therefore satisfied, so \eqref{eq:gap:accel:convergence} holds; in particular, $\zeta_{*,N} \gap(\tilde x_{*,N}, \tilde y_{*,N}; \realoptu) \le C_0 \defeq  \norm{u^0-\realoptu}_{\Test_1 \Precond_1}^2$ for all $N \ge 2$.
    It remains to study the convergence rate of the gap from this estimate.
    We have $\zeta_{*,N}=\sum_{k=1}^{N-1} \tauTest_k^{1/2}$.
    In the unaccelerated case ($\gamma=0$), we get $\zeta_{*,N} = N \tauTest_0^{1/2}$.
    This gives the claimed $O(1/N)$ rate.
    In the accelerated case, $\tauTest_k$ is of the order $\Omega(k^2)$ by the proof of \cref{thm:testing:pdps:accel}. Therefore also $\zeta_{*,N}$ is of the order $\Theta(N^2)$, so we get the claimed $O(1/N^2)$ convergence.
    In the linear convergence case, likewise, $\tauTest_k$ is exponential. Therefore so is $\zeta_{*,N}$.
\end{proof}

\begin{remark}[spatially adaptive and stochastic methods]
    \label{thm:gap:accel:stochastic}
    Recalling the block-separability \cref{ex:testing:structured:strongmono-block}, consider the spaces $X=X_1 \times \cdots \times X_m$ and $Y=Y_1 \times \cdots \times Y_n$. Suppose $F(x)=\sum_{j=1}^m F_j(x_j)$ and $G^*(y)=\sum_{\ell=1}^n G_\ell^*(y_\ell)$ for $x=(x_1, \ldots, x_m) \in X$ and $y=(y_1, \ldots, y_n) \in Y$.
    Take $\Test_{k+1} \defeq \begin{psmallmatrix} \TauTest_k & 0 \\ 0 & \SigmaTest_{k+1} \end{psmallmatrix}$ as well as $\Step_{k+1} \defeq \begin{psmallmatrix} \Tau_k & 0 \\ 0 & \Sigma_{k+1} \end{psmallmatrix}$ for $\Tau_k \defeq \sum_{j=1}^n \tau_{k,j} P_j$, and similar expressions for $\TauTest_k,\Sigma_{k+1}$, and $\Sigma_{k+1}$, where $P_j x \defeq x_j$ projects into $X_j$.
    Instead of $\tauTest_k\tau_k=\sigmaTest_k\sigma_k$ that we required in \eqref{thm:testing:pdps:accel}, imposing $\E[\TauTest_k\Tau_k]=\E[\SigmaTest_k\Sigma_k]=\eta_k I$ for some scalar $\eta_k$, we may then start following through the proof of \cref{thm:testing:pdps:accel} to derive stochastic block-coordinate methods that randomly update only some of the blocks on each iteration, as well as methods that adapt the blockwise step lengths to the spatial or blockwise structure of the problem. With somewhat more effort, we can also follow through the proofs of the present \cref{sec:gap:accel}. Specifically, if we replace our ergodic sequences by
    \begin{equation*}
        \tilde x_{*,N} \defeq \inv\zeta_{*,N}\E\Biggl[\sum_{k=1}^{N-1} \Tau_k^*\TauTest_k^* \nextx\Biggr]
        \quad\text{and}\quad
        \tilde y_{*,N} \defeq \inv\zeta_{*,N}\E\Biggl[\sum_{k=1}^{N-1} \Sigma_k^*\SigmaTest_k^* \thisy\Biggr]
        \quad\text{for}\quad
         \zeta_{*,N} \defeq \sum_{k=1}^{N-1} \eta_k,
    \end{equation*}
    we then obtain in place of \eqref{eq:gap:accel:convergence} the estimate
    \begin{equation*}
        \E\left[\frac{1}{2}\norm{u^N-\realoptu}^2_{\Test_{N+1}\Precond_{N+1}}\right]
        +  \zeta_{*,N} \gap(\tilde x_{*,N}, \tilde y_{*,N})
        + \sum_{k=0}^{N-1} \E\left[\GenGap_{k+1}(\realoptu)\right]
        \le
        \norm{u^0-\realoptu}_{\Test_1 \Precond_1}^2.
    \end{equation*}
    If instead  $\E[\TauTest_k\Tau_k]=\E[\SigmaTest_{k+1}\Sigma_{k+1}]=\eta_k I$, we get the result for the ergodic sequences
    \begin{equation*}
        \tilde x_{N} \defeq \inv\zeta_{N}\E\Biggl[\sum_{k=0}^{N-1} \Tau_k^*\TauTest_k^* \nextx\Biggr]
        \quad\text{and}\quad
        \tilde y_{N} \defeq \inv\zeta_{N}\E\Biggl[\sum_{k=0}^{N-1} \Sigma_{k+1}^*\SigmaTest_{k+1}^* \nexty\Biggr]
        \quad\text{where}\quad
        \zeta_N \defeq \sum_{k=0}^{N-1} \eta_k.
    \end{equation*}
    In either case, if we do not or cannot, due to lack of strong convexity of some of the $F_\ell$, accelerate all of the blockwise step lengths $\tau_{k+1,j}$ with the same factor $\gamma=\gamma_j$, it will generally be the case that $\E\left[\GenGap_{k+1}(\realoptu)\right]<0$. This quantity will have such an order of magnitude that we get mixed $O(1/N^2)+O(1/N)$ convergence rates. We refer to \cite{tuomov-blockcp} for details on such spatially adaptive and stochastic primal-dual methods, and \cite{wright2015coordinate} for an introduction to the idea of stochastic coordinate descent.
\end{remark}

\section{Convergence of the ADMM}
\label{sec:gap:admm}

Let $G: X \to \Rbar$, $F: Z \to \Rbar$ be convex, proper, and lower semicontinuous, $A \in \linear(X; Y)$,  $B \in \linear(Z; Y)$, and $c \in Y$. Recall the problem
\begin{equation}
    \label{eq:gap:admm:problem}
    \min_{x, z}~ J(x, z) \defeq G(x) + F(z) + \delta_C(x, z),
\end{equation}
where
\[
    C \defeq \{(x, z) \in X \times Z\mid Ax+Bz=c\}.
\]
We now show an ergodic convergence result for the ADMM applied to this problem, which we recall from \eqref{eq:admm} to read
\begin{algeqbox}
    \begin{equation}
        \label{eq:gap:admm}
        \left\{\begin{aligned}
                \nextx & \in \inv{(A^*A+\inv\tau \subdiff F)}(A^*(c-B\this{z}-\inv\tau\this{\lambda})), \\
                \nexxt{z} & \in \inv{(B^*B+\inv\tau \subdiff G)}(B^*(c-A\nextx-\inv\tau\this{\lambda})), \\
                \nexxt{\lambda} & \defeq \this{\lambda} + \tau(A\nextx+B\nexxt{z}-c).
        \end{aligned}\right.
    \end{equation}
\end{algeqbox}
The general structure of the convergence proof is very similar to all the other algorithms we have studied. However, now the forward-step component does not arise as a gradient $\grad \tilde E$ but is a special non-self-adjoint preconditioner $\tilde\Precond_{i+1}$. Moreover, in the first stage of the proof we obtain a convergence estimate for a duality gap that we then refine at the end of the proof to separate function value and constraint satisfaction estimates.

\begin{theorem}
    Let $G: X \to \Rbar$ and $F: Z \to \Rbar$ be convex, proper, and lower semicontinuous, $A \in \linear(X; Y)$,  $B \in \linear(Z; Y)$, and $c \in Y$. Let $J$ be defined as in \eqref{eq:gap:admm:problem}, which we assume to admit a solution $(\realoptx, \realoptz) \in X \times Z$.
    For arbitrary initial iterates $(x^0, y^0, \lambda^0)$, let $\{(\nextx,\nextz,\nexxt\lambda)\}_{k \in \N} \subset X \times Z \times Y$ be generated by the ADMM \eqref{eq:gap:admm} for \eqref{eq:gap:admm:problem}. Define the ergodic sequences $\tilde x^N \defeq \frac{1}{N}\sum_{k=0}^{N-1} \nextx$ and $\tilde z^N \defeq \frac{1}{N}\sum_{k=0}^{N-1} \nextz$.
    Then both $(G+F)(\tilde x^N, \tilde z^N) \to \min_{(x,z)\in X\times Z} J(x,z)$ and $\norm{A\tilde x^N+B\tilde z^N - c}_Y \to 0$ at the rate $O(1/N)$.
\end{theorem}

\begin{proof}
    We consider the augmented problem
    \[
        \min_{(x,z) \in X \times Z }~ J_\tau(x, z) \defeq G(x) + F(z) + \delta_C(x, z) + \frac{\tau}{2}\norm{Ax+Bz-c}_Y^2,
    \]
    which has the same solutions as \eqref{eq:gap:admm:problem}.
    As the normal cone to the constraint set $C$ at any point $(x, z) \in C$ is given by $N_C(x, z) =\{(A^*\lambda, B^*\lambda) \mid \lambda \in Y\}$,
    setting $u=(x, z, \lambda)$ and
    \[
        H(u) \defeq
        \begin{pmatrix}
            \subdiff G(x) + A^*\lambda + \tau A^*(Ax+Bz-c) \\
            \subdiff F(z) + B^*\lambda + \tau B^*(Ax+Bz-c) \\
            -(Ax+Bz-c)
        \end{pmatrix},
    \]
    the optimality conditions for this problem can be written as $0 \in H(u)$.
    In particular, there exists $\realopt\lambda \in Y$ such that $(\realoptx,\realoptz,\realopt\lambda) \in \inv H(0)$. However, we will not be needing this, and take $\realopt\lambda$ arbitrary.

    We could rewrite the algorithm \eqref{eq:gap:admm} as \eqref{eq:gap:ppext} with
    \begin{equation*}
        \Happrox_{k+1}(u) =
        H(u)
        \quad\text{and}\quad
        \Precond_{k+1} =
        \begin{pmatrix}
            0 & -\tau A^*B & -A^*\\
            0 & 0 & -B^* \\
            0 & 0 & \inv\tau I
        \end{pmatrix}.
    \end{equation*}
    However, $\Precond_{k+1}$ is nonsymmetric, and any symmetrizing $\Test_{k+1}$ would make $\Test_{k+1}\Happrox_{k+1}$ difficult to analyze.
    We therefore take instead
    \[
        \Happrox_{k+1}(u) \defeq
        H(u)
        +\tilde\Precond_{k+1}(u-\thisu)
        \quad\text{with}\quad
        \tilde\Precond_{k+1} \defeq
        \begin{pmatrix}
            0 & -\tau A^*B & -A^*\\
            0 & -\tau B^*B & -B^* \\
            0 & 0 & 0
        \end{pmatrix},
    \]
    as well as
    \[
        \Precond_{k+1} \defeq
        \begin{pmatrix}
            0 & 0 & 0 \\
            0 & \tau B^*B & 0 \\
            0 & 0 & \inv\tau I
        \end{pmatrix},
        \quad\text{and}\quad
        \Test_{k+1} \defeq I.
    \]
    Clearly $\Test_{k+1}\Precond_{k+1}$ is self-adjoint.

    Let us set
    \[
        \Gamma \defeq \tau \begin{pmatrix} A^*A & A^*B & 0 \\ B^* A & B^* B & 0 \\ 0 & 0 & 0 \end{pmatrix}
        \quad\text{and}\quad
        \Xi \defeq
        \begin{pmatrix}
            0 & 0 & A^* \\
            0 & 0 & B^* \\
            -A & -B & 0
        \end{pmatrix}.
    \]
    Using the fact that $A\realoptx+B\realoptx=c$, observe that we can split $H=\subdiff \tilde F+\Xi$, where
    \[
        \begin{aligned}
        \tilde F(u)
        &
        \defeq
        G(x)+F(z)+\frac{\tau}{2}\norm{Ax+Bz-c}_Y^2 + \iprod{c}{\lambda}_Y
        \\
        &
        =
        G(x)+F(z)+\frac{1}{2}\norm{u-\realoptu}_\Gamma^2 + \iprod{c}{\lambda}_Y.
        \end{aligned}
    \]
    It follows that
    \[
        \begin{aligned}
        \iprod{\Hany(\nextu)}{\nextu-\realoptu}_{\Test_{k+1}}
        &
        \ge \tilde F(\nextu)-\tilde F(\realoptu) + \frac{1}{2}\norm{\nextu-\realoptu}^2_\Gamma
        +\iprod{\realoptu}{\nextu}_\Xi
        \\
        &
        = [F(\nextx)+G(\nextz)]-[F(\realoptx)+G(\realoptx)] + \iprod{c}{\nexxt\lambda-\realopt\lambda}_Y
        \\ \MoveEqLeft[-1]
        +\norm{\nextu-\realoptu}^2_\Gamma
        +\iprod{\realoptu}{\nextu}_\Xi.
        \end{aligned}
    \]
    Again using $A\realoptx+B\realoptx=c$, we expand
    \[
        \begin{aligned}
        \iprod{\realoptu}{\nextu}_\Xi
        &
        =\iprod{\realopt\lambda}{A\nextx+B\nextz}_Y
        -\iprod{A\realoptx+B\realoptz}{\nexxt\lambda}_Y
        \\
        &
        =\iprod{\realopt\lambda}{A\nextx+B\nextz-c}_Y
        -\iprod{c}{\nexxt\lambda-\realopt\lambda}_Y.
        \end{aligned}
    \]
    Thus
    \begin{equation}
        \label{eq:gap:admm:hany-est}
        \begin{aligned}[t]
        \iprod{\Hany(\nextu)}{\nextu-\realoptu}_{\Test_{k+1}}
        &
        \ge
        [F(\nextx)+G(\nextz)]-[F(\realoptx)+G(\realoptx)]
        \\ \MoveEqLeft[-1]
        +\norm{\nextu-\realoptu}^2_\Gamma
        +\iprod{\realopt\lambda}{A\nextx+B\nextz-c}_Y
        \\
        &
        = \bar F(\nextu; \realopt\lambda) - \bar F(\realoptu; \realopt\lambda)
        +\norm{\nextu-\realoptu}^2_\Gamma
        \end{aligned}
    \end{equation}
    for
    \begin{equation}
        \label{eq:gap:admm:barf}
        \bar F(u; \realopt\lambda)
        \defeq
        F(x)+G(z) + \iprod{\realopt\lambda}{Ax+Bz-c}_Y.
    \end{equation}

    On the other hand,
    \[
        \begin{aligned}
            \iprod{\nextu-\thisu}{\nextu-\realoptu}_{\Test_{k+1}\tilde\Precond_{k+1}}
            &=
            \iprod{-\tau B(\nexxt{z}-\this{z})-(\nexxt{\lambda}-\this{\lambda})}{A(\nextx-\realoptx)}_Y
            \\ \MoveEqLeft[-1]
            +\iprod{-\tau B(\nexxt{z}-\this{z})-(\nexxt{\lambda}-\this{\lambda})}{B(\nexxt{z}-\realopt{z})}_Y
            \\
            & =
            \iprod{-\tau B(\nexxt{z}-\this{z})-(\nexxt{\lambda}-\this{\lambda})}{A(\nextx-\realoptx)+B(\nextz-\realoptz)}_Y.
        \end{aligned}
    \]
    From \eqref{eq:gap:admm} we recall
    \[
        \nexxt{\lambda} -\this{\lambda} = \tau(A\nextx+B\nexxt{z}-c)
        =\tau[A(\nextx-\realoptx)+B(\nextz-\realoptz)].
    \]
    Hence
    \begin{equation}
        \label{eq:gap:admm:tildem-est}
        \begin{aligned}[t]
        \iprod{\nextu-\thisu}{\nextu-\realoptu}_{\Test_{k+1}\tilde\Precond_{k+1}}
        &
        =-\norm{\nextu-\realoptu}^2_\Gamma-\iprod{B(\nextz-\thisz)}{\nexxt{\lambda}-\this{\lambda}}_Y
        \\
        &
        \ge -\norm{\nextu-\realoptu}^2_\Gamma-\frac{1}{2}\norm{\nextu-\thisu}_{\Test_{i+1}\Precond_{i+1}}^2.
        \end{aligned}
    \end{equation}
    Combining \eqref{eq:gap:admm:hany-est} and \eqref{eq:gap:admm:tildem-est} it follows that
    \[
        \iprod{\Happrox_{k+1}(\nextu)}{\nextu-\realoptu}_{\Test_{k+1}}
        \ge \bar F(\nextu; \realopt\lambda) - \bar F(\realoptu; \realopt\lambda)
        -\frac{1}{2}\norm{\nextu-\thisu}_{\Test_{i+1}\Precond_{i+1}}^2.
    \]
    By \cref{thm:gap:ppext:convergence} now
    \begin{equation*}
        \frac{1}{2}\norm{u^N-\realoptu}^2_{\Test_{N+1}\Precond_{N+1}}
        +
        \sum_{k=0}^{N-1} \left(\bar F(\nextu; \realopt\lambda) - \bar F(\realoptu; \realopt\lambda)\right)
        \le
        \frac{1}{2}\norm{u^0-\realoptu}^2_{\Test_{1}\Precond_{1}}
        \quad
        (N \ge 1).
    \end{equation*}
    Writing $\tilde u^N =(\tilde x^N, \tilde y^N, \tilde\lambda^N) \defeq \frac{1}{N}\sum_{k=0}^{N-1} \nextu$, Jensen's inequality now shows that
    \begin{equation}
        \label{eq:gap:admm:gap0}
        \bar F(\tilde u^N; \realopt\lambda) - \bar F(\realoptu; \realopt\lambda) \le \frac{1}{2N}\norm{u^0-\realoptu}^2_{\Test_{1}\Precond_{1}}
        \quad (N \ge 1).
    \end{equation}
    Since $A\realoptx+B\realoptz=c$, observe that $\bar F(\freevar; \realopt\lambda)-\bar F(\realoptu; \realopt\lambda)$ is the Lagrangian duality gap \eqref{eq:gap:pd:gap} for the saddle-point formulation \eqref{eq:proximal:admm:minmax} of \eqref{eq:gap:admm:problem}, hence nonnegative when $\realoptu \in \inv H(0)$. So \eqref{eq:gap:admm:gap0} shows the convergence of the duality gap.
    However, we can improve the result somewhat since $\realopt\lambda$ was taken as arbitrary. Expanding $\bar F$ using \eqref{eq:gap:admm:barf} and taking the supremum over $\realopt\lambda \in \B(0, \kappa)$ in \eqref{eq:gap:admm:gap0}, we thus obtain for any $\kappa>0$ the estimate
    \begin{equation*}
        \begin{aligned}
            0 &\le [F(\tilde x^N)+G(\tilde z^N)]-[F(\realoptx)+G(\realoptx)]+
            \kappa \norm{A\tilde x^N+B\tilde z^N-c}_Y
            \\
            &=
            \sup_{\realopt\lambda \in \B(0, \kappa)}\left( F(\tilde u^N; \realopt\lambda) - \bar F(\realoptu; \realopt\lambda)\right)
            \le \sup_{\realopt\lambda \in \B(0, \kappa)} \frac{1}{2N}\norm{u^0-\realoptu}^2_{\Test_{1}\Precond_{1}}.
        \end{aligned}
    \end{equation*}
    This gives the claim.
\end{proof}

\chapter{Meta-algorithms}
\label{chap:meta}

In this chapter, we consider several \term[algorithm!meta-]{``meta-algorithms''} for accelerating minimization algorithms such as the ones derived in the previous chapters. These include \emph{inertia} and \emph{over-relaxation}, as well as \emph{line searches}.
These schemes differ from the strong convexity based acceleration of \cref{chap:testing} in that no additional assumptions are made on $F$ and $G$. Rather, through the use of an additional extrapolated or interpolated point, the first two schemes attempt to obtain a second-order approximation of the function. Line search, on the other hand, can be used to find optimal parameters or to estimate unknown parameters.
Throughout the chapter, we base our work on the abstract algorithm \eqref{eq:gap:ppext}, i.e.,
\begin{equation}
    \label{eq:meta:ppext}
    0 \in  \Happrox_{k+1}(\nextx) + \Precond_{k+1}(\nextx-\thisx),
\end{equation}
where the iteration-dependent set-valued operator $\Happrox_{k+1}: X \setto X$ in suitable sense approximates a (monotone) operator $H: X \setto X$, whose root we intend to find, and $\Precond_{k+1} \in \linear(X; X)$ is a linear preconditioner.

\section{Over-relaxation}
\label{sec:meta:over-relaxation}

We start with \term{over-relaxation}.
Essentially, this amounts to taking \eqref{eq:meta:ppext} and replacing $\thisx$ in the preconditioner by an over-relaxed point $\thisz$ defined for some parameters $\lambda_{k}>0$ through the recurrence
\begin{equation}
    \label{eq:meta:overelax-z}
    \nextz \defeq \inv\lambda_{k}\nextx+(1-\inv\lambda_{k})\thisz.
\end{equation}
We thus seek to solve
\begin{equation}
    \label{eq:meta:ppext-overrelax}
    0 \in \Happrox_{k+1}(\nextx) + \Precond_{k+1}(\nextx-\thisz).
\end{equation}
Since $\nextz-\thisz=\inv\lambda_{k}(\nextx-\thisz)$, we can write \eqref{eq:meta:ppext} as
\begin{equation}
    \label{eq:meta:overelax-zmethod}
    0 \in \Happrox_{k+1}(\nextx) + \lambda_{k}\Precond_{k+1}(\nextz-\thisz).
\end{equation}
We can therefore lift the overall algorithm into the form \eqref{eq:meta:ppext} as
\begin{equation}
    \label{eq:meta:ppext-lifted}
    0 \in \hat H_{k+1}(\nexxt{q})+\hat M_{k+1}(\nexxt{q}-\this{q})
\end{equation}
by taking $q\defeq(x, z)$ with
\begin{equation}
    \label{eq:lifting-overrelax}
    \hat H_{k+1}(q) \defeq
    \begin{pmatrix}
        \Happrox_{k+1}(x) \\
        \inv\lambda_{k}(z - x) \\
    \end{pmatrix}
    \quad\text{and}\quad
    \hat \Precond_{k+1} \defeq
    \begin{pmatrix}
        0 & \lambda_{k}\Precond_{k+1} \\
        0 & (I-\inv\lambda_{k})I \\
    \end{pmatrix}.
\end{equation}
To be able to use our previous estimate on $\iprod{\Happrox_{k+1}(\nextx)}{\nextx-\realoptx}_{\Test_{k+1}}$, we would like to test with
\[
    \hat\Test_{k+1} \defeq
    \begin{pmatrix}
        \lambda_{k}\Test_{k+1} & 0 \\
        0 & 0\\
    \end{pmatrix}.
\]
Unfortunately, $\Test_{k+1}\Precond_{k+1}$ is not self-adjoint, so \cref{thm:gap:ppext:convergence} does not apply. However, observing from \eqref{eq:meta:overelax-z} that
\begin{equation}
    \label{eq:meta:overelax-z-alt}
    \nextz-\nextx=(1-\lambda_{k})(\nextz-\thisz),
\end{equation}
we are able to proceed along the same lines of proof.

\begin{theorem}
    \label{thm:meta:overrelax:convergence}
    On a Hilbert space $X$, let $\Happrox_{k+1}: X \setto X$, and $\Precond_{k+1}, \Test_{k+1} \in \linear(X; X)$ for $k \in \N$.
    Suppose \eqref{eq:meta:ppext-overrelax} is solvable for the iterates $\{\thisx\}_{k \in \N}$.
    If $\Test_{k+1}\Precond_{k+1}$ is self-adjoint,
    \begin{equation}
        \label{eq:meta:overrelax:recursion}
        \lambda_{k}^2\Test_{k+1}\Precond_{k+1} \succeq \lambda_{k+1}^2\Test_{k+2}\Precond_{k+2},
    \end{equation}
    and
    \begin{equation}
        \label{eq:meta:overrelax:fundamental-condition}
        \iprod{\Happrox_{k+1}(\nextx)}{\nextx-\realoptx}_{\Test_{k+1}} \ge \GenGap_{k+1}(\realoptx) - \frac{1}{2}\norm{\nextx-\thisz}^2_{\Test_{k+1}Q_{k+1}}
    \end{equation}
    for some $Q_{k+1} \in \linear(X; X)$, for all $k \in \N$ and some $\realoptx \in X$ and $\GenGap_{k+1}(\realoptx) \in \R$, then
    \begin{multline}
        \label{eq:meta:overrelax:quantitative-fejer}
        \frac{\lambda_{k+1}^2}{2}\norm{\nextz-\realoptx}_{\Test_{k+2}\Precond_{k+2}}^2
        +\lambda_{k}\GenGap_{k+1}(\realoptx)
        + \frac{\lambda_k}{2}\norm{\nextz-\thisz}^2_{\lambda_k(2\lambda_k-1)\Test_{k+1}\Precond_{k+1}-\Test_{k+1}Q_{k+1}}
        \\
        \le
        \frac{\lambda_{k}^2}{2}\norm{\thisz-\realoptx}_{\Test_{k+1}\Precond_{k+1}}^2
        \quad (k \in \N).
    \end{multline}
\end{theorem}

\begin{proof}
    Taking $\realopt{q} \defeq (\realoptx, \realoptx)$, we apply $\iprod{\freevar}{\nexxt{q}-\realopt{q}}_{\hat\Test_{k+1}}$ to \eqref{eq:meta:ppext-overrelax}.
    Thus
    \[
        0 \in \iprod{\hat H_{k+1}(\nexxt{q})+\hat\Precond_{k+1}(\nexxt{q}-\this{q})}{\nexxt{q}-\realopt{q}­}_{\hat\Test_{k+1}}.
    \]
    Observe that
    \[
        \hat\Test_{k+1}\hat\Precond_{k+1}
        =
        \begin{pmatrix}
            0 & \lambda_{k}^2\Test_{k+1}\Precond_{k+1} \\
            0 & 0
        \end{pmatrix}.
    \]
    Thus
    \[
        0 \in \iprod{\Happrox_{k+1}(\nextx)}{\nextx-\realoptx}_{\lambda_{k}\Test_{k+1}}
        +\lambda_k^2\iprod{\nextz-\thisz}{\nextx-\realoptx}_{\Test_{k+1}\Precond_{k+1}}.
    \]
    Using \eqref{eq:meta:overelax-z-alt} we then get
    \[
        \begin{aligned}
            0 \in \iprod{\Happrox_{k+1}(\nextx)}{\nextx-\realoptx}_{\lambda_{k}\Test_{k+1}}
            &-
            \lambda_{k}^2(1-\lambda_{k})\norm{\nextz-\thisz}_{\Test_{k+1}\Precond_{k+1}}^2\\
            &+\lambda_{k}^2\iprod{\nextz-\thisz}{\nextz-\realoptx}_{\Test_{k+1}\Precond_{k+1}}.
        \end{aligned}
    \]
    Using the three-point-identity \eqref{eq:convergence:three-point-identity}, we rearrange this into
    \[
        \begin{aligned}
        0 \in \iprod{\Happrox_{k+1}(\nextx)}{\nextx-\realoptx}_{\lambda_{k}\Test_{k+1}}
        &+
        \frac{\lambda_{k}^2-2\lambda_{k}^2(1-\lambda_{k})}{2}\norm{\nextz-\thisz}_{\Test_{k+1}\Precond_{k+1}}^2
        \\
        &+\frac{\lambda_{k}^2}{2}\norm{\nextz-\realoptx}_{\Test_{k+1}\Precond_{k+1}}^2
        -\frac{\lambda_{k}^2}{2}\norm{\thisz-\realoptx}_{\Test_{k+1}\Precond_{k+1}}^2.
        \end{aligned}
    \]
    Observe that $\lambda_{k}^2-2\lambda_{k}^2(1-\lambda_{k})=\lambda_k^2(2\lambda_k-1)$.
    Using \eqref{eq:meta:overelax-z}, \eqref{eq:meta:overrelax:fundamental-condition}, and \eqref{eq:meta:overrelax:recursion}, this gives \eqref{eq:meta:overrelax:quantitative-fejer}.
\end{proof}

Clearly we should try to ensure $\lambda_k(2\lambda_k-1)\Test_{k+1}\Precond_{k+1} \ge \Test_{k+1}Q_{k+1}$. If $\Test_{k+1}\Precond_{k+1}=\Test_0\Precond_0$ is constant and $Q_{k+1}=0$, this holds if $\{\lambda_k\}_{k \in N}$ is nonincreasing and satisfies $\lambda_k\ge 1/2$.
Therefore, we cannot get any convergence rates from the iterates in this case. It is, however, possible to obtain convergence of a gap, and it would be possible to obtain weak convergence.

The next result is a variant of \cref{cor:gap:ergodic:general} for over-relaxed methods.

\begin{corollary}
    \label{cor:meta:overrelax}
    Let $H \defeq \subdiff \tilde F+\grad \tilde G+\Xi$, where $\Xi \in \linear(X; X)$ is skew-adjoint, and $\tilde G: X \to \Rbar$ and $\tilde F: X \to \R$ convex, proper, and lower semicontinuous. Suppose $\tilde F$ satisfies for some $\Lambda \in \linear(X; X)$ the three-point smoothness condition \eqref{eq:gap:ergodic:3smooth}.
    Also let $\Precond \in \linear(X;X)$ be positive semi-definite and self-adjoint. Pick $x^0=z^0 \in X$, and define the sequence $\{(\nextx, \nextz)\}_{k \in \N}$ through
    \begin{equation}
        \label{eq:meta:overrelax:fgxi}
        \left\{
            \begin{aligned}
                0 & \in [\subdiff \tilde G(\nextx)+\subdiff \tilde F(\thisz)+\Xi\nextx] + M(\nextx-\thisz),
                \\
                \nextz & \defeq \inv\lambda_k\nextx-(\inv\lambda_k-1)\thisz.
            \end{aligned}
        \right.
    \end{equation}
    Suppose $\{\lambda_k\}_{k \in \in \N}$ is nonincreasing and
    \begin{equation}
        \label{eq:meta:overrelax:fgxi:lambda-cond}
        \lambda_k(2\lambda_k-1)\Precond
        \succeq \Lambda
        \quad (k \in \N).
    \end{equation}
    Then for every $\realoptx \in \inv H(0)$ and the gap functional $\tilde\gap$ defined in \eqref{eq:gap:ergodic:generic-gap},
    \begin{equation}
        \label{eq:meta:overrelax:general-gap-estimate}
        \tilde\gap(\tilde x^N; \realoptx)
        \le
        \frac{\lambda_0^2}{2\sum_{k=0}^{N-1} \lambda_k}\norm{z^0-\realoptx}_{\Precond}^2,
        \quad\text{where}\quad
        \tilde x^N \defeq \frac{1}{\sum_{k=0}^{N-1} \lambda_k} \sum_{k=0}^{N-1} \lambda_k \nextx.
    \end{equation}
\end{corollary}

\begin{proof}
    The method \eqref{eq:meta:overrelax:fgxi} is \eqref{eq:meta:ppext-overrelax} with $\tilde H_{k+1}(x) \defeq \subdiff \tilde G(x) + \subdiff \tilde F(\thisz)+\Xi x$ as well as $\Precond_{k+1} \equiv \Precond $ and $\Test_{k+1} \equiv \Id$.
    Using \eqref{eq:gap:ergodic:3smooth} for $\tilde F$, the convexity of $\tilde G$, and the assumption $ZW=\eta \Id$, we obtain as in the proof of \eqref{thm:gap:ergodic:general} the estimate
    \begin{equation*}
        \iprod{\Happrox_{k+1}(\nextx)}{\nextx-\realoptx}_X
        \ge
        \tilde \gap(\nextx; \realoptx) -\frac{1}{2}\norm{\thisz-\nextx}^2_{\Lambda}
    \end{equation*}
    This provides \eqref{eq:meta:overrelax:fundamental-condition} while \eqref{eq:meta:overrelax:fgxi:lambda-cond} and the constant choice of the testing and preconditioning operators guarantee that $\lambda_k(2\lambda_k-1)\Test_{k+1}\Precond_{k+1} \succeq \Test_{k+1}Q_{k+1}$ for $Q_{k+1} \equiv \Lambda$.
    By \cref{thm:meta:overrelax:convergence}, we now obtain
    \begin{equation}
        \label{eq:meta:overrelax:mainest}
        \frac{\lambda_{k+1}^2}{2}\norm{\nextz-\realoptx}^2_{\Precond}
        +\lambda_k \tilde \gap(\nextx; \realoptx)
        \le
        \frac{\lambda_k^2}{2}\norm{\thisz-\realoptx}^2_{\Precond}.
    \end{equation}
    Summing over $k=0,\ldots,N-1$ and an application of Jensen's inequality finishes the proof.
\end{proof}

\subsection*{Over-relaxed proximal point method}

We apply the above results to the \term[method!proximal point!over-relaxed]{over-relaxed proximal point method}
\begin{algeqbox}
    \begin{equation}
        \label{eq:meta:overrelax:prox}
        \left\{\begin{aligned}
                \nextx & \defeq \prox_{\tau G}(\thisz), \\
                \nextz & \defeq \inv\lambda_k\nextx-(\inv\lambda_k-1)\thisz.
        \end{aligned}\right.
    \end{equation}
\end{algeqbox}

\begin{theorem}
    \label{thm:meta:overrelax:prox}
    Let $G: X \to \Rbar$ be convex, proper, and lower semicontinuous with $\inv{[\subdiff G]}(0) \ne \emptyset$.
    Pick an initial iterate $x^0=z^0 \in X$.
    If $\{\lambda_k\}_{k \in \N} \ge 1/2$ is nonincreasing, the ergodic sequence $\{\tilde x^N\}_{N \in \N}$ defined in \eqref{eq:meta:overrelax:general-gap-estimate} and generated from the iterates $\{\thisx\}_{k \in \N}$ of the over-relaxed proximal point method \eqref{eq:meta:overrelax:prox} satisfies $G(\tilde x^N) \to \minval G\defeq \min_{x\in X}G(x)$ at the rate $O(1/N)$.
\end{theorem}

\begin{proof}
    We apply \cref{cor:meta:overrelax} with $\tilde G=G$, $\tilde F=0$, $\Precond=\inv\tau\Id$. Clearly $\tilde F$ satisfies \eqref{eq:gap:ergodic:3smooth} with $\Lambda=0$.
    Then \eqref{eq:meta:overrelax:fgxi:lambda-cond} holds if $2\lambda_k \ge 1$, that is to say $\lambda_k \ge 1/2$.
    For $\realoptx \in \argmin G$, we have $\tilde \gap(x; \realoptx)=G(x)-G(\realoptx)=G(x)-\minval G$.
    Therefore \cref{cor:meta:overrelax} gives
    \begin{equation}
        G(\tilde x^N)
        \le
        \minval G
        +
        \frac{\lambda_0^2}{2\tau\sum_{k=0}^{N-1} \lambda_k}\norm{z^0-\realoptx}_X^2
    \end{equation}
    Since $\sum_{k=0}^{N-1} \lambda_k \ge N/2$ by the lower bound on $\lambda_k$, we get the claimed $O(1/N)$ convergence rate of the function values for the ergodic sequence.
\end{proof}

\subsection*{Over-relaxed explicit splitting}

For a smooth function $F$, the \term[method!explicit splitting!over-relaxed]{over-relaxed explicit splitting method} iterates
\begin{algeqbox}
    \begin{equation}
        \label{eq:meta:overrelax:fb}
        \left\{\begin{aligned}
                \nextx & \defeq \prox_{\tau G}(\thisz-\tau \grad F(\thisz)), \\
                \nextz & \defeq \inv\lambda_k\nextx-(\inv\lambda_k-1)\thisz.
        \end{aligned}\right.
    \end{equation}
\end{algeqbox}

\begin{theorem}
    \label{thm:meta:overrelax:fb}
    Let $J \defeq G + F$ for $G: X \to \Rbar$ and $F: X \to \R$ be convex, proper, and lower semicontinuous with $\grad F$ $L$-Lipschitz.
    Suppose $\inv{[\subdiff J]}(0) \ne \emptyset$.
    Pick an initial iterate $x^0=z^0 \in X$.
    If $\{\lambda_k\}_{k \in \N}$ is nonincreasing and satisfies
    \begin{equation}
        \label{eq:meta:over-forward-backward:bound}
        \lambda_k \ge \frac{1}{4}(1 + \sqrt{1+8 L\tau}),
    \end{equation}
    then the ergodic sequence $\{\tilde x^N\}_{N \in \N}$ defined in \eqref{eq:meta:overrelax:general-gap-estimate} and generated from the iterates $\{\thisx\}_{k \in \N}$ of the over-relaxed explicit splitting method \eqref{eq:meta:overrelax:fb} satisfies $J(\tilde x^N) \to \minval J\defeq \min_{x\in X}J(x)$ at the rate $O(1/N)$.
\end{theorem}

\begin{proof}
    We apply \cref{cor:meta:overrelax} with $\tilde G=G$, $\tilde F=F$, and $\Precond=\inv\tau\Id$.
    By \cref{cor:smoothness:three-point}, $\tilde F$ satisfies the three-point smoothness condition \eqref{eq:gap:ergodic:3smooth} with $\Lambda=L\,\Id$.
    The condition \eqref{eq:meta:overrelax:fgxi:lambda-cond} consequently holds if $\lambda_k(2\lambda_k-1)>L\tau$, which holds under the assumption \eqref{eq:meta:over-forward-backward:bound}.
    The rest follows as in the proof of \cref{thm:meta:overrelax:prox}.
\end{proof}

\subsection*{Over-relaxed PDPS}

With $F=F_0+E: X \to \Rbar$, $G^*: Y \to \Rbar$, and $K \in \linear(X; Y)$, take $H: X \times Y \setto X \times Y$ as well as $\tilde F, \tilde G$, and $\Xi$ as in \eqref{eq:gap:pd:setup}, and the preconditioner $\Precond$ as in \eqref{eq:gap:pd:precond}
 for fixed step length parameters $\tau,\sigma>0$.
 Writing $\thisz=(\this\xi, \this\upsilon)$, and, as usual $\thisu=(\thisx,\thisy)$, the method \eqref{eq:meta:overelax-zmethod} then becomes the \term[method!primal-dual proximal splitting!over-relaxed]{over-relaxed primal-dual proximal splitting (PDPS) method} with a forward step, also known as the \term[method!V\~u--Condat]{V\~u--Condat method}:
 \begin{algeqbox}
     \begin{equation}
         \label{eq:meta:overrelax:pdps}
         \left\{\begin{aligned}
                 \nextx & \defeq (I+\tau \subdiff F_0)^{-1}(\this{\xi} - \tau K^* \thisy - \tau \grad E(\this{\xi})),\\
                 \overnextx & \defeq (\nextx-\this{\xi})+\nextx, \\
                 \nexty & \defeq (I+\sigma \subdiff G^*)^{-1}(\this{\upsilon} + \sigma K \overnextx), \\
                 \nexxt{\xi} & \defeq \inv\lambda_k\nextx-(\inv\lambda_k-1)\this{\xi}, \\
                 \nexxt{\upsilon} & \defeq \inv\lambda_k\nexty-(\inv\lambda_k-1)\this{\upsilon}.
         \end{aligned}\right.
     \end{equation}
 \end{algeqbox}
For the statement of the next result, we recall that for the primal-dual saddle-point operator $\Hany$ from \eqref{eq:gap:pd:setup}, the generic gap functional $\tilde\gap$ becomes the primal-dual gap $\gap$ given in \eqref{eq:gap:pd:gap}.

\begin{theorem}
    Suppose $F_0: X \to \Rbar$, $E: X \to \R$ and $G: Y \to \Rbar$ are convex, proper, and lower semicontinuous on Hilbert spaces $X$ and $Y$ with $\grad E$ $L$-Lipschitz.
    Let also $K \in \linear(X; Y)$.
    With $F=F_0+E$, suppose the assumptions of \cref{thm:convex:fenchel} are satisfied.
    Pick an initial iterate $u^0=z^0 \in X \times Y$.
    If the sequence $\{\lambda_k\}_{k \in \N}$  is nonincreasing and satisfies
    \begin{equation}
        \label{eq:meta:over-pdps:bound}
        \lambda_k \geq
        \frac{1}{4}(1+\sqrt{1+8L\tau/(1-\tau\sigma\norm{K}^2)})
        \quad\text{and}\quad
        \tau\sigma\norm{K}^2 < 1,
    \end{equation}
    then the ergodic sequence $\{\tilde u^N=(\tilde x^N, \tilde y^N)\}_{N \in \N}$ defined as in \eqref{eq:meta:overrelax:general-gap-estimate} and generated from the iterates $\{\thisu=(\thisx, \thisy)\}_{k \in \N}$ of the over-relaxed PDPS method \eqref{eq:meta:overrelax:pdps} satisfies $\gap(\tilde x_{*,N}, \tilde y_{*,N}) \to 0$  at the rate $O(1/N)$.
\end{theorem}

\begin{proof}
    We recall that $\inv H(0) \ne \emptyset$ under the assumptions of \cref{thm:convex:fenchel}. Clearly $\Precond$ is self-adjoint.
    The condition \eqref{eq:meta:overrelax:fgxi:lambda-cond} can with \eqref{eq:testing:pdps:zimi-est} be reduced to
    \[
        \begin{pmatrix}
            \lambda_k(2\lambda_k-1)\delta\tau \Id & 0 \\
            0 & \inv\sigma I- \tau\inv{(1-\delta)} KK^*
        \end{pmatrix}
        \succeq
        \begin{pmatrix}
            L & 0 \\ 0 & 0
        \end{pmatrix}
    \]
    for some $\delta \in (0, 1)$.
    As in \eqref{eq:testing:pdps:step-conds-L-0} in the proof of \cref{thm:testing:pdps:accel}, these conditions reduce to
    \begin{equation}
        \label{eq:meta:overrelax:cpock-step-conds-L-0}
        \lambda_k(2\lambda_k-1)\delta \ge \tau L
        \quad\text{and}\quad
        1-\delta \ge \tau\sigma \norm{K}^2.
    \end{equation}
    The first inequality holds if $\lambda_k \ge \frac{1}{4}(1+\sqrt{1+8L\tau\inv\delta})$.
    Solving the second inequality as an equality for $\delta$ yields the condition
    \[
        \lambda_k \ge \frac{1}{4}(1+\sqrt{1+8L\tau\inv{[(1-\tau\sigma\norm{K}^2)]}}),
    \]
    i.e., \eqref{eq:meta:over-pdps:bound}.
    Now we obtain the gap convergence from \cref{cor:meta:overrelax}.
\end{proof}

\begin{remark}
    The method \eqref{eq:meta:overrelax:pdps} is due to \cite{vu2013splitting,condat2013primaldual}. The convergence of the ergodic gap was observed in \cite{chambolle2014ergodic}.
\end{remark}

\section{Inertia}
\label{sec:meta:inertia}

Our next \term[inertia]{inertial} meta-algorithm will likewise not yield convergence of the main iterates, but through a special arrangement of variables combined with intricate unrolling arguments, is able to do away with the word \emph{ergodic} in the gap estimates.
In essence, the meta-algorithm replaces the previous iterate $\thisx$ in the linear preconditioner of \eqref{eq:meta:ppext} by an inertial point
\begin{equation}
    \label{eq:meta:inertia-baru}
    \this{\bar x} \defeq (1+\alpha_{k})\thisx-\alpha_{k}\prevx
    \quad\text{for}\quad
    \alpha_{k} \defeq \lambda_{k}(\inv\lambda_{k-1}-1)
\end{equation}
for some \term[parameter!inertial]{inertial parameter} sequence $\{\lambda_k\}_{k \in \N}$. We thus solve
\begin{equation}
    \label{eq:meta:ppext-inertia}
    0 \in \Happrox_{k+1}(\nextx) + \Precond_{k+1}(\nextx-\this{\bar x}).
\end{equation}
We can relate this to over-relaxation as follows: we simply replace $\thisz$ in the definition \eqref{eq:meta:overelax-z} of $\nextz$ by $\thisx$, i.e., we take
\begin{equation}
    \label{eq:meta:inertia-z}
    \nextz \defeq \inv\lambda_{k}\nextx-(\inv\lambda_{k}-1)\thisx.
\end{equation}
Since
\begin{equation}
    \label{eq:meta:inertia:z-bar-u-relationship}
    \begin{aligned}[t]
    \lambda_{k}(\nextz-\thisz)
    &
    =\nextx-(1-\lambda_{k})\thisx-\lambda_{k}[\inv\lambda_{k-1}\thisx-(\inv\lambda_{k-1}-1)\prevx]
    \\
    &
    =\nextx-[1-\lambda_{k}+\lambda_{k}\inv\lambda_{k-1}]\thisx+\lambda_{k}(\inv\lambda_{k-1}-1)\prevx
    \\
    &=\nextx-[(1+\alpha_{k})\thisx-\alpha_{k}\prevx]
    = \nextx - \this{\bar x},
    \end{aligned}
\end{equation}
we obtain the method \eqref{eq:meta:overelax-zmethod}, with the differing update \eqref{eq:meta:inertia-z} of $\nextz$.
Again we can also lift the overall algorithm into the form \eqref{eq:meta:ppext}, specifically \eqref{eq:meta:ppext-lifted}, by taking $q\defeq(x, z)$ with
\[
    \hat H_{k+1}(q) \defeq
    \begin{pmatrix}
        \Happrox_{k+1}(x) \\
        z - x \\
    \end{pmatrix},
    \quad\text{and}\quad
    \hat \Precond_{k+1} \defeq
    \begin{pmatrix}
        0 & \lambda_{k}\Precond_{k+1} \\
        (I-\inv\lambda_{k})I & 0 \\
    \end{pmatrix}.
\]

Now comes the trick with inertial methods: Unlike with over-relaxed methods, where we wanted to avoid having to estimate $\iprod{\Happrox_{k+1}(\nextx)}{\nextz-\realoptz}_{\Test_{k+1}}$, with inertial methods we are brave enough to do this. Indeed, our specific choice  \eqref{eq:meta:inertia-z} makes this possible, as we shall see below. We therefore test with
\[
    \hat\Test_{k+1} \defeq
    \begin{pmatrix}
        0 & 0 \\
        \lambda_{k}\Test_{k+1} & 0\\
    \end{pmatrix}
\]
to obtain a self-adjoint and positive semi-definite
\begin{equation}
    \label{eq:meta:inertia-zm}
    \hat\Test_{k+1}\Precond_{k+1}
    =
    \begin{pmatrix}
        0 & 0 \\
        0 & \lambda_{k}^2\Test_{k+1}\Precond_{k+1}
    \end{pmatrix}.
\end{equation}
Therefore \cref{thm:gap:ppext:convergence} applies, and we obtain the following:

\begin{theorem}
    \label{thm:meta:inertia:convergence}
    Let $X$ be a Hilbert space, $\Happrox_{k+1}: X \setto X$, and $\Precond_{k+1}, \Test_{k+1} \in \linear(X; X)$ for $k \in \N$.
    Suppose \eqref{eq:meta:ppext-inertia} is solvable for the iterates $\{\thisx\}_{k \in \N}$ and inertial parameters $\{\lambda_k\}_{k \in \N} \subset (0, \infty)$.
    If $\Test_{k+1}\Precond_{k+1}$ is self-adjoint, and
    \begin{equation}
        \label{eq:meta:inertia:fundamental-condition}
        \begin{aligned}[t]
            \lambda_{k}\iprod{\Happrox_{k+1}(\nextx)}{\nextz-\realoptz}_{\Test_{k+1}}
            &
            \ge
            \GenGap_{k+1}(\realoptx)
            +
            \frac{1}{2}\norm{\nextz-\realoptz}_{\lambda_{k+1}^2\Test_{k+2}\Precond_{k+2}-\lambda_{k}^2\Test_{k+1}\Precond_{k+1}}^2
            \\ \MoveEqLeft[-1]
            - \frac{\lambda_{k}^2}{2}\norm{\nextz-\thisz}_{\Test_{k+1} \Precond_{k+1}}^2
        \end{aligned}
    \end{equation}
    for all $k \in \N$ and some $\realoptx \in X$ and $\GenGap_{k+1}(\realoptx) \in \R$, then
    \begin{equation*}
        \frac{\lambda_{N}^2}{2}\norm{z^N-\realoptx}^2_{\Test_{N+1}\Precond_{N+1}}
        +
        \sum_{k=0}^{N-1} \GenGap_{k+1}(\realoptx)
        \le
        \frac{\lambda_0^2}{2}\norm{z^0-\realoptx}^2_{\Test_{1}\Precond_{1}}
        \quad
        (N \ge 1).
    \end{equation*}
\end{theorem}

\begin{proof}
    This follows directly from \cref{thm:gap:ppext:convergence} and the expansion \eqref{eq:meta:inertia-zm}.
\end{proof}

We now provide examples of how to apply this result to the proximal point method and explicit splitting. As we recall, in these algorithms we take $\Test_{k+1}=\tauTest_k I$ and $\Step_{k+1}=\tau_k I$. To proceed, we will need a few further general-purpose technical lemmas. The first one is the fundamental lemma for inertia, which provides inertial function value unrolling.

\begin{lemma}
    \label{lemma:meta:inertia:proxest}
    Let $G: X \to \Rbar$ be convex, proper, and lower semicontinuous.
    Suppose $\lambda_k \in [0, 1]$ and $\tauTest_k,\tau­_k > 0$ for $k \in \N$ with
    \begin{equation}
        \label{eq:meta:inertia:proxest-lambda-recurrence}
        \tauTest_{k+1}\tau_{k+1}(1-\lambda_{k+1}) \le \tauTest_k\tau_k
        \quad (k \ge 0).
    \end{equation}
    Assume $\nexxt{q} \in \subdiff G(\nextx)$ for $k=0,\ldots,N-1$, and $0 \in \subdiff G(\realoptx)$.
    Then
    \begin{equation}
        \label{eq:meta:inertia:proxest}
        \begin{aligned}[t]
        s_{G,N} & \defeq \sum_{k=0}^{N-1}
            \tauTest_k\tau_k\lambda_k \iprod{\nexxt{q}}{\nextz - \realoptx}_X
        \\
        &
        \ge
        \tauTest_{N-1}\tau_{N-1}(G(x^N) - G(\realoptx))
        -\tauTest_0\tau_0(1-\lambda_0)(G(x^0)-G(\realoptx)).
        \end{aligned}
    \end{equation}
\end{lemma}

\begin{proof}
    Using \cref{eq:meta:inertia-z}, observe that
    \begin{equation}
        \label{eq:meta:inertia-nextz-realoptu-useful-for-unrolling}
        \begin{aligned}[t]
        \lambda_k(\nextz-\realoptx)
        &
        =\lambda_k[\inv\lambda_{k+1}\nextx-(\inv\lambda_k-1)\thisx-\realoptx]
        \\
        &
        =\lambda_k(\nextx-\realoptx)+(1-\lambda_k)(\nextx-\thisx).
        \end{aligned}
    \end{equation}
    Recalling from \eqref{eq:meta:inertia-nextz-realoptu-useful-for-unrolling} that $\lambda_k(\nextz-\realoptx)=\lambda_k(\nextx-\realoptx)+(1-\lambda_k)(\nextx-\thisx)$ and using the convexity of $G$, we can estimate
    \begin{equation}
        \label{eq:meta:inertia:proxest-first-ineq}
        \begin{aligned}[t]
        s_{G,N}
        &
        =
        \sum_{k=0}^{N-1} \tauTest_k\tau_k\Bigl[
            \lambda_k\iprod{\nexxt{q}}{\nextx - \realoptx}_X
            +(1-\lambda_k)\iprod{\nexxt{q}}{\nextx - \thisx}_X
        \Bigr]
        \\
        &
        \ge
        \sum_{k=0}^{N-1} \tauTest_k\tau_k\Bigl[
            \lambda_k(G(\nextx) - G(\realoptx))
            +(1-\lambda_k)(G(\nextx)-G(\thisx))
        \Bigr]
        \\
        &
        =
        \sum_{k=0}^{N-1} \left[
            \tauTest_k\tau_k(G(\nextx) - G(\realoptx))
            -\tauTest_k\tau_k(1-\lambda_k)(G(\thisx)-G(\realoptx))
        \right].
        \end{aligned}
    \end{equation}
    Since $G(\thisx)\ge G(\realoptx)$, the recurrence inequality \cref{eq:meta:inertia:proxest-lambda-recurrence} together with a telescoping argument now gives
    \[
        s_{G,N}
        \ge
        \tauTest_{N-1}\tau_{N-1}(G(x^N)-G(\realoptx))
        -\tauTest_0\tau_0(1-\lambda_0)(G(x^0)-G(\realoptx))
    \]
    as claimed.
\end{proof}

\begin{lemma}
    \label{lemma:meta:inertia:recurrence}
    Suppose $\lambda_0=1$ and $\lambda_{k}^{-2}=\lambda_{k+1}^{-2}-\inv\lambda_{k+1}$ for $k= 0,\dots,N-1$, $N\in\N$.
    Then
    \begin{equation}
        \label{eq:meta:inertia:recurrence}
        \lambda_{k+1}=\frac{2}{1+\sqrt{1+4\lambda_k^{-2}}} \qquad(k=0,\dots,N-1)
    \end{equation}
and $\inv \lambda_N \ge \frac12(N+1)$.
\end{lemma}

\begin{proof}
    First, it is straightforward to verify that the recursion \eqref{eq:meta:inertia:recurrence} defines a sequence that satisfies the assumed quadratic relation. We show the lower bound by induction on $N$. For $N=1$, the recursion gives $\inv\lambda_1=\frac12(1+\sqrt{5})\geq 1$. Assume now that $N\in\N$ is arbitrary and $\inv\lambda_N\geq \frac12(N+1)$. Then it follows from \eqref{eq:meta:inertia:recurrence} and the induction assumption that
    \begin{equation*}
        \inv\lambda_{N+1}
            \geq \frac12\left(1+\sqrt{1+4\left(\frac{N+1}{2}\right)^2}\right)
            \geq \frac12\left(1+\sqrt{4\left(\frac{N+1}{2}\right)^2}\right)
            = \frac12(N+2)
    \end{equation*}
    and hence the claim.
\end{proof}

\subsection*{Inertial proximal point method}

Let $H=\subdiff G$ and $\tilde H_{k+1}=\tau \subdiff G$ for a convex, proper, lower semicontinuous function $G$.
Take $\tau>0$ and $\lambda_{k+1}$ by \eqref{eq:meta:inertia:recurrence} for $\lambda_0=1$.
Then \eqref{eq:meta:ppext-inertia} becomes the \term[method!proximal point!inertial]{inertial proximal point method}
\begin{algeqbox}
    \begin{equation}
        \label{eq:meta:inertia:prox}
        \left\{
            \begin{aligned}
                \nextx & \defeq \prox_{\tau G}(\this{\bar x}),
                \\
                \alpha_{k+1} &\defeq \lambda_{k+1}(\inv\lambda_k-1),\\
                \nexxt{\bar x} & \defeq (1+\alpha_{k+1})\nextx-\alpha_{k+1}\thisx.
            \end{aligned}
        \right.
    \end{equation}
\end{algeqbox}
Note that $x^0$ is never needed since $\alpha_1=0$. The real initial iterate, which can be freely chosen, is $\bar x^0$.

\begin{theorem}
    \label{thm:meta:inertia:prox}
    Let $G: X \to \Rbar$ be convex, proper, and lower semicontinuous.
    Suppose $\inv{[\subdiff G]}(0) \ne \emptyset$.
    Take $\tau>0$ and $\lambda_0=1$, and pick an initial iterate $\bar x^0 \in X$. Then the inertial proximal point method \eqref{eq:meta:inertia:prox} satisfies $G(x^N) \to \minval G$ at the rate $O(1/N^2)$.
\end{theorem}

\begin{proof}
    If we take $\tau_k=\tau$ as stated and $\tauTest_k=\lambda_k^{-2}$, then \eqref{lemma:meta:inertia:recurrence} verifies \eqref{eq:meta:inertia:proxest-lambda-recurrence}.
    Since now $\lambda_{k+1}^2\tauTest_{k+1} = \lambda_k^2\tauTest_k$, \eqref{eq:meta:inertia:fundamental-condition} holds if
    \begin{equation}
        \label{eq:meta:inertia:fundamental-condition-prox}
        \lambda_k\tauTest_k\tau_k\iprod{\subdiff G(\nextx)}{\nextz-\realoptz}_X
        \ge
        \GenGap_{k+1}(\realoptx)
        - \frac{\lambda_{k}^2\tauTest_k}{2}\norm{\nextz-\thisz}_X^2
    \end{equation}
    for some $\GenGap_{k+1}(\realoptx) \in \R$.
    This is verified by \cref{lemma:meta:inertia:proxest} for some $\GenGap_{k+1}(\realoptx)$ such that
    \[
        \sum_{k=0}^{N-1} \GenGap_{k+1}(\realoptx)
        \ge  \tauTest_{N-1}\tau_{N-1}(G(x^N) - G(\realoptx))
        -\tauTest_0\tau_0(1-\lambda_0)(G(x^0)-G(\realoptx)).
    \]
    Since $\lambda_0=1$, \cref{thm:meta:inertia:convergence} gives the estimate
    \[
        \frac{\tauTest_N\lambda_N^2}{2}\norm{x^N-\realoptx}_X^2
        +\tauTest_{N-1}\tau_{N-1}(G(x^N) - G(\realoptx))
        \le
        \frac{\tauTest_0\lambda_0^2}{2}\norm{x^0-\realoptx}_X^2.
    \]
    \Cref{lemma:meta:inertia:recurrence} now yields $\tauTest_{N-1}\tau_{N-1}=\lambda_{N-1}^{-2}\tau \ge \tau \frac14N^2$, and hence we obtain the claimed convergence rate.
\end{proof}

\subsection*{Inertial explicit splitting}

Let $H=\subdiff G+\grad F$ and $\tilde H_{k+1}(x)=\tau(\subdiff G(x)+\grad F(\this{\bar x}))$ for convex, proper, lower semicontinuous functions $G$ and $F$ with $F$ smooth.
Take $\tau>0$ and $\lambda_{k+1}$ by \eqref{eq:meta:inertia:recurrence} for $\lambda_0=1$.
Then \eqref{eq:meta:ppext-inertia} becomes the \term[method!explicit splitting!inertial]{inertial explicit splitting method}
\begin{algeqbox}
    \begin{equation}
        \label{eq:meta:inertia:fista}
        \left\{
            \begin{aligned}
                \nextx & \defeq \prox_{\tau G}(\this{\bar x}-\tau \grad F(\this{\bar x})),
                \\
                \alpha_{k+1} &\defeq \lambda_{k+1}(\inv\lambda_k-1),\\
                \nexxt{\bar x} & \defeq (1+\alpha_{k+1})\nextx-\alpha_{k+1}\thisx.
            \end{aligned}
        \right.
    \end{equation}
\end{algeqbox}
Again, $x^0$ is never needed, as $\alpha_1=0$. The actual initial iterate to be chosen is $\bar x^0$.

To prove the convergence of this method, we need to incorporate the forward step into \cref{lemma:meta:inertia:proxest}.

\begin{lemma}
    \label{lemma:meta:inertia:splitest}
    Let $J \defeq F+G$ for $G: X \to \Rbar$ and $F: X \to \R$ be convex, proper, and lower semicontinuous.
    Suppose $F$ has $L$-Lipschitz gradient and that $\lambda_k \in [0, 1]$ and $\tauTest_k,\tau_k >0$ satisfy the recurrence inequality \eqref{eq:meta:inertia:proxest-lambda-recurrence} for $k \in \N$.
    Assume $\nexxt{w} \in \subdiff G(\nextx)$ for all $k=0,\ldots,N-1$, and that $0 \in \subdiff J(\realoptx)$.
    Then
    \begin{equation}
        \label{eq:meta:inertia:splitest}
        \begin{aligned}[t]
        s_N &\defeq \sum_{k=0}^{N-1}\left(
            \tauTest_k\tau_k\lambda_k \iprod{\nexxt{w} + \grad F(\this{\bar x})}{\nextz - \realoptx}_X
            +\frac{\tauTest_k\tau_k \lambda_k^2 L}{2}\norm{\nextz-\thisz}_X^2
            \right)
        \\
        &
        \ge
        \tauTest_{N-1}\tau_{N-1}(J(x^N) - J(\realoptx))
        -\tauTest_{0}\tau_{0}(1-\lambda_0)(J(x^0) - J(\realoptx)).
        \end{aligned}
    \end{equation}
\end{lemma}

\begin{proof}
    We recall from \cref{eq:meta:inertia:z-bar-u-relationship} that
    $
        \frac{\lambda_k^2}{2}\norm{\nextz-\thisz}_X^2
        =
        \frac{1}{2}\norm{\nextx-\this{\bar x}}_X^2.
    $
    We therefore estimate using \cref{cor:smoothness:three-point} that
    \begin{equation}
        \label{eq:meta:inertia:forward-step-unrolling}
        \begin{aligned}[t]
        s_{F,N}
        & \defeq \sum_{k=0}^{N-1}\left(
            \tauTest_k\tau_k\lambda_k \iprod{\grad F(\this{\bar x})}{\nextz - \realoptx}_X
            +\frac{\tauTest_k\tau_k \lambda_k^2 L}{2}\norm{\nextz-\thisz}_X^2
            \right)
        \\
        &
        =
        \sum_{k=0}^{N-1} \tauTest_k\tau_k\left[
            \lambda_k\iprod{\grad F(\this{\bar x})}{\nextx - \realoptx}_X
            +(1-\lambda_k)\iprod{\grad F(\this{\bar x})}{\nextx - \thisx}_X
            + \frac{L}{2}\norm{\nextx-\this{\bar x}}_X^2
        \right]
        \\
        &
        \ge
        \sum_{k=0}^{N-1} \tauTest_k\tau_k\left[
            \lambda_k(F(\nextx) - F(\realoptx))
            +(1-\lambda_k)(F(\nextx)-F(\thisx))
        \right]
        \\
        &
        =
        \sum_{k=0}^{N-1} \left[
            \tauTest_k\tau_k(F(\nextx) - F(\realoptx))
            -\tauTest_k\tau_k(1-\lambda_k)(F(\thisx)-F(\realoptx))
        \right].
        \end{aligned}
    \end{equation}
    Summing with the estimate \eqref{eq:meta:inertia:proxest-first-ineq} for $G$, we deduce
    \[
        s_N
        \ge
        \sum_{k=0}^{N-1} \left[
            \tauTest_k\tau_k((F+G)(\nextx) - (F+G)(\realoptx))
            -\tauTest_k\tau_k(1-\lambda_k)((F+G)(\thisx)-(F+G)(\realoptx))
        \right].
    \]
    Since $(F+G)(\thisx) \ge (F+G)(\realoptx)$, the recurrence inequality \cref{eq:meta:inertia:proxest-lambda-recurrence} together with a telescoping argument now gives the claim.
\end{proof}

\begin{theorem}
    \label{thm:meta:inertia:fb}
    Let $J \defeq G + F$ for $G: X \to \Rbar$ and $F: X \to \R$ be convex, proper, and lower semicontinuous with $\grad F$ Lipschitz.
    Suppose $\inv{[\subdiff J]}(0) \ne \emptyset$.
    Take $\tau>0$ with $\tau L \le 1$, and $\lambda_0=1$, and pick an initial iterate $\bar x^0 \in X$. Then the inertial explicit splitting \eqref{eq:meta:inertia:fista} satisfies $J(x^N) \to \min_{x\in X} J(x)$ at the rate $O(1/N^2)$.
\end{theorem}

\begin{proof}
    The proof follows that of \cref{thm:meta:inertia:prox}: in place of \eqref{eq:meta:inertia:fundamental-condition-prox} we reduce \eqref{eq:meta:inertia:fundamental-condition} to the condition
    \[
        \lambda_k\tauTest_k\tau_k\iprod{\subdiff G(\nextx)+\grad F(\this{\bar x})}{\nextz-\realoptz}_X
        \ge
        \GenGap_{k+1}(\realoptx)
        - \frac{\lambda_{k}^2\tauTest_k}{2}\norm{\nextz-\thisz}_X^2.
    \]
    This is verified for some $\GenGap_{k+1}(\realoptx)$ such that
    \[
        \sum_{k=0}^{N-1} \GenGap_{k+1}(\realoptx)
        \ge  \tauTest_{N-1}\tau_{N-1}(J(x^N) - J(\realoptx))
        -\tauTest_0\tau_0(1-\lambda_0)(J(x^0) - J(\realoptx))
    \]
    by using \cref{lemma:meta:inertia:splitest} and the bound $\tau L \le 1$ in place of \cref{lemma:meta:inertia:proxest}.
\end{proof}

\begin{remark}[accelerated gradient methods, FISTA]
    The inertial scheme was first introduced by \cite{Nesterov} for the basic gradient descent method for smooth functions. The extension to explicit splitting is due to \cite{BeckTeboulle}, which proposed a \term[method!iterative soft-thresholding!fast]{fast iterative shrinkage-thresholding algorithm} (\index{FISTA|see{method, iterative soft-thresholding, fast}}{FISTA}) for the specific problem of minimizing a least-squares term plus a weighted $\ell^1$ norm.
    (Note that in most treatments of FISTA, $\inv\lambda_k$ is written as $t_k$.)
    We refer to \cite{Nesterov:2004,Beck:2017} for a further discussion of these algorithms and more general accelerated gradient methods based on combinations of a history of iterates.
\end{remark}

\begin{remark}[PDPS, Douglas--Rachford, and correctors]
    The above unrolling arguments cannot be directly applied to PDPS, Douglas--Rachford splitting, and other methods based on \eqref{eq:meta:ppext} with non-maximally monotone $H$. Following \cite{chambolle2014ergodic}, one can apply inertia to the PDPS method with the restricted choice $\alpha_k \in (0, 1/3)$.
    This prevents the use of the FISTA rule \eqref{eq:meta:inertia:recurrence} and only yields $O(1/N)$ convergence of an ergodic gap.
    Based on alternative argumentation, when one of the functions is quadratic, \cite{patrinos2014douglas} managed to employ the FISTA rule and obtain $O(1/N^2)$ rates for inertial Douglas--Rachford splitting.
    Moreover, \cite{tuomov-inertia} observed that by introducing a \term{corrector} for the non-subdifferential component of $H$, in essence $\Xi_{k+1}$, the gap unrolling arguments can be performed. This approach also allows combining inertial acceleration with strong monotonicity based acceleration.
\end{remark}

\section{Line search}
\label{sec:meta:linesearch}

Let us return to the basic results on weak convergence (\cref{thm:convergence:fb}), strong convergence with rates (\cref{thm:testing:fb}), and function value convergence (\cref{thm:gap:fb:value}) of the explicit splitting method. These results depend on the three-point inequalities of \cref{cor:smoothness:three-point} (or, for faster rates under strong convexity, \cref{cor:smoothness:three-point:sc}), specifically either the non-value estimate
\begin{align}
    \label{eq:meta:linesearch:three-point:monotonicity}
    \iprod{\grad F(\thisx)-\grad F(\realoptx)}{\nextx-\realoptx}_X &\ge
    -\frac{L}{4}\norm{\nextx-\thisx}_X^2
\intertext{or the value estimate}
    \label{eq:meta:linesearch:three-point:smoothness}
    \iprod{\grad F(\thisx)}{\nextx-\realoptx}_X
    &\ge
    F(\nextx)-F(\realoptx) -  \frac{L}{2}\norm{\nextx-\thisx}_X^2.
\end{align}
Recall that for weak convergence of iterates, we required the step length parameters $\{\tau_k\}_{k \in \N}$ to satisfy on each iteration the bound $\tau_k L < 2$.
Under a strong convexity assumption, the bound $\tau_k L \le 2$ was sufficient for strong convergence of iterates. Function value convergence was finally shown under the bound $\tau_k L \le 1$. All cases thus hold for $\tau_k L \le 1$, which we assume in the following for simplicity.

In this section, we address the following question: What if we do not know the Lipschitz factor $L$?
A basic idea is to take $L$ large enough. But what is large enough? Finding such a large enough $L$ is the same as taking $\tau_k$ small enough and $L=1/\tau_k$. This leads us to the following rough \term{line search} rule: for some $\tau>0$ and line search parameter $\theta \in (0, 1)$, start with $\tau_k \defeq \tau$, and iterate $\tau_k \mapsto \theta \tau_k$ until \eqref{eq:meta:linesearch:three-point:smoothness} (or \eqref{eq:meta:linesearch:three-point:monotonicity}) is satisfied with $L=1/\tau_k$. Note that on each update of $\tau_k$, we need to recalculate $\nextx \defeq \prox_{\tau_k G}(\thisx-\tau_k \grad F(\thisx)))$.

Performing this line search still appears to depend on knowing $\realoptx$ through \eqref{eq:meta:linesearch:three-point:smoothness}. However, going back to the proof of \cref{cor:smoothness:three-point}, we see that what is really needed is to satisfy the smoothness (or descent) inequality \eqref{eq:smoothness:smoothness} which was used to derive \eqref{eq:meta:linesearch:three-point:smoothness}. We are therefore lead to the following practical line search method to guarantee the inequality
\begin{equation}
    \label{eq:meta:linesearch:smoothness}
    \iprod{\grad F(\thisx)}{\nextx-\thisx}_X
    \ge
    F(\nextx)-F(\thisx) -  \frac{1}{2\tau_k}\norm{\nextx-\thisx}_X^2
\end{equation}
at every iteration:
\begin{algenumbox}
    \begin{enumerate}[label=\arabic*., start=0]
        \item Pick $\theta \in (0,1)$, $\tau >0$, $\lambda_0\defeq 1$, $x^0\in X$; set $k=0$.
        \item\label{item:meta:linesearch:init} Set $\tau_k=\tau$.
        \item\label{item:meta:linesearch:prox} Calculate $\nextx \defeq \prox_{\tau_k G}(\thisx-\tau_k \grad F(\thisx))$.
        \item If \eqref{eq:meta:linesearch:smoothness}
            does not hold, update $\tau_k \defeq \theta \tau_k$, and go back to step \ref{item:meta:linesearch:prox}
        \item Set $k \defeq k+1$, and continue from step \ref{item:meta:linesearch:init}
    \end{enumerate}
\end{algenumbox}

\begin{theorem}[explicit splitting line search]
    Let $J \defeq F+G$ where $G: X \to \Rbar$ and $F: X \to \R$ are convex, proper, and lower semicontinuous, with $\grad F$ moreover Lipschitz.
    Suppose $\inv{[\subdiff J]}(0) \ne \emptyset$.
    Then the above line search method satisfies $J(x^N) \to \min_{x\in X}J(x)$ at the rate $O(1/N)$. If $G$ is strongly convex, then this convergence is linear.
\end{theorem}

\begin{proof}
    Since $\grad F$ is $\tilde L$-smooth for some unknown $\tilde L>0$, eventually the line search procedure satisfies $1/\tau_k \ge \tilde L$. Hence \eqref{eq:meta:linesearch:smoothness} is satisfied, and $\tau_k \ge \epsilon > 0$ for some $\epsilon>0$. We can therefore follow through the proof of \cref{thm:gap:fb:value} with $L=1/\tau_k$.
\end{proof}

We can also combine the line search method with the inertial explicit splitting   \eqref{eq:meta:inertia:fista}. If in place of \eqref{eq:meta:linesearch:smoothness} we seek to satisfy
\begin{equation}
    \label{eq:meta:linesearch:smoothness:inertia}
    \iprod{\grad F(\this{\bar x})}{\nextx-\thisx}_X
    \ge
    F(\nextx)-F(\thisx) - \frac{1}{2\tau_k}\norm{\nextx-\this{\bar x}}_X^2,
\end{equation}
then also
\begin{equation*}
    \iprod{\grad F(\this{\bar x})}{\nextx-\realoptx}_X
    \ge
    F(\nextx)-F(\realoptx) -  \frac{1}{2\tau_k}\norm{\nextx-\this{\bar x}}_X^2.
\end{equation*}
This allows the inequality of \eqref{eq:meta:inertia:forward-step-unrolling} to be shown.

We are therefore led to the following practical backtracking inertial explicit splitting:
\begin{algenumbox}
    \begin{enumerate}[label=\arabic*., start=0]
        \item Pick $\theta \in (0,1)$, $\tau >0$, $\lambda_0\defeq 1$, $\bar x^0=x^0\in X$; set $k=0$.
        \item\label{item:meta:linesearch:init:inertia} Set $\tau_k=\tau$.
        \item\label{item:meta:linesearch:prox:inertia} Calculate $\nextx \defeq \prox_{\tau_k G}(\this{\bar x}-\tau_k \grad F(\this{\bar x})))$.
        \item If \eqref{eq:meta:linesearch:smoothness:inertia} does not hold, update $\tau_k \defeq \theta \tau_k$, and go back to step \ref{item:meta:linesearch:prox:inertia}
        \item Set $\nexxt{\bar x} \defeq (1+\alpha_{k+1})\nextx-\alpha_{k+1}\thisx$ for $\alpha_{k+1} \defeq \lambda_{k+1}(\inv\lambda_k-1)$.
        \item Set $k \defeq k+1$, and continue from step \ref{item:meta:linesearch:init:inertia}
    \end{enumerate}
\end{algenumbox}

The proof of the following is immediate:
\begin{theorem}
    Let $J \defeq G + F$ for $G: X \to \Rbar$ and $F: X \to \R$ be convex, proper, and lower semicontinuous with $\grad F$ Lipschitz.
    Suppose $\inv{[\subdiff J]}(0) \ne \emptyset$.
    Take $\tau>0$ and $\lambda_0=1$, and pick an initial iterate $\bar x^0 \in X$. Then the above backtracking inertial explicit splitting satisfies $J(x^N) \to \min_{x\in X} J(x)$ at the rate $O(1/N^2)$.
\end{theorem}

The reader may now work out how to use line search to satisfy the nonnegativity of the metric $\Test_{k+1}\Precond_{k+1}$ in the PDPS method when $\norm{K}$ is not known, or how to satisfy the condition $L\tau_0 + \tau_0\sigma_0\norm{K}^2 < 1$ when the Lipschitz factor $L$ of the forward step component $E$ is not known.

\begin{remark}[adaptive inertial parameters, quasi-Newton methods, and primal-dual proximal line searches]
    Regarding our statement in the beginning of the chapter about inertia methods attempting to construct a second-order approximation of the function, \cite{ochs2017adaptive} show that an adaptive inertial explicit splitting, performing an optimal line search on $\lambda_k$ instead of $\tau_k$, is equivalent to a proximal quasi-Newton method.
    Such a method is a further development of variants \cite[see][]{beck2009fgp} of the method that attempt to restore the monotonicity of explicit splitting that is lost by inertia.
    Indeed, if $J(\nexxt{\bar x}) \le J(\this{\bar x})$ does not hold for $\lambda_k<1$, we can revert to $\lambda_k=1$ to ensure descent as the step reduces to basic explicit splitting, which we know to be monotone by \cref{thm:gap:fb:value}.
    Finally, a line search for the PDPS method is studied in \cite{malitsky2018linesearch}.
\end{remark}

\part{Nonconvex analysis}\label{part:nonconvex}

\chapter{Clarke subdifferentials}\label{chap:clarke}

We now turn to a concept of generalized derivatives that covers, among others, both Fréchet derivatives and convex subdifferentials. Again, we start with the general class of functionals that admit such a derivative. It is clear that we need to require some continuity properties, since otherwise there would be no relation between functional values at neighboring points and thus no hope of characterizing optimality through pointwise properties. In \cref{part:convex}, we used lower semicontinuity for this purpose, which together with convexity yielded the required properties. In this part, we want to drop the latter, global, assumption; in turn we need to strengthen the local continuity assumption. We thus consider now locally Lipschitz continuous functionals. Recall that $F:X\to\R$ is locally Lipschitz continuous near $x\in X$ if there exist a $\delta >0$ and an $L>0$ (which in the following will always denote the local Lipschitz constant of $F$) such that
\begin{equation*}
    |F(x_1)-F(x_2)|\leq L \norm{x_1-x_2}_X \qquad\text{for all }x_1,x_2\in \OB(x,\delta).
\end{equation*}
We will refer to the $\OB(x,\delta)$ from the definition as the \term[neighborhood, Lipschitz]{Lipschitz neighborhood} of $x$.
Note that for this we have to require that $F$ is (locally) finite-valued (but see \cref{rem:clarke:extended-real} below).
Throughout this chapter, we will assume that $X$ is a Banach space unless stated otherwise.

\section{Definition and basic properties}

We proceed as for the convex subdifferential and first define for $F:X\to\R$ the
\term[derivative!directional!generalized]{generalized directional derivative} in $x\in X$ in direction $h\in X$ as
\begin{equation}\label{eq:clarke:dir}
    F^\circ(x;h) \defeq \limsup_{\substack{y\to x\\t\downto 0}} \frac{F(y+th)-F(y)}{t}.
\end{equation}
Note the difference to the classical directional derivative: We no longer require the existence of a limit but merely of accumulation points.
We will need the following properties.
\begin{lemma}\label{lem:clarke:dir}
    Let $F:X\to\R$ be locally Lipschitz continuous near $x\in X$ with factor $L>0$. Then the mapping $h\mapsto F^\circ (x;h)$ is
    \begin{enumerate}
        \item \label{lem:clarke:dir:i}
            Lipschitz continuous with constant $L$ and satisfies $|F^\circ(x;h)|\leq L\norm{h}_X<\infty$;
        \item \label{lem:clarke:dir:ii}
            \term[functional!subadditive]{subadditive}, i.e., $F^\circ(x;h+g)\leq F^\circ(x;h)+F^\circ(x;g)$ for all $h,g\in X$;
        \item \label{lem:clarke:dir:iii}
            \term[functional!positively homogeneous]{positively homogeneous}, i.e., $F^\circ(x;\alpha h) = (\alpha F)^\circ(x;h)$ for all $\alpha > 0$ and $h\in X$;
        \item \label{lem:clarke:dir:iv}
            \term[functional!reflective]{reflective}, i.e., $F^\circ(x;-h) = (-F)^\circ(x;h)$ for all $h\in X$.
    \end{enumerate}
\end{lemma}
\begin{proof}
    \emph{\ref{lem:clarke:dir:i}:} Let $h,g\in X$ be arbitrary. The local Lipschitz continuity of $F$ implies that
    \begin{equation*}
        F(y+th)-F(y) \leq F(y+tg) - F(y) + t L\norm{h-g}_X
    \end{equation*}
    for all $y$ sufficiently close to $x$ and $t$ sufficiently small. Dividing by $t>0$ and taking the $\limsup$ then yields that
    \begin{equation*}
        F^\circ(x;h) \leq F^\circ(x;g) + L\norm{h-g}_X.
    \end{equation*}
    Exchanging the roles of $h$ and $g$ shows the Lipschitz continuity of $F^\circ(x;\cdot)$, which also yields the claimed boundedness since $F^\circ(x;g) = 0$ for $g=0$ from the definition.

    \emph{\ref{lem:clarke:dir:ii}:} Since $t\downto 0$ and $g\in X$ is fixed, $y\to x$ if and only if $y+tg\to x$. The definition of the $\limsup$ and the productive zero thus immediately yield
    \begin{equation*}
        \begin{aligned}
            F^\circ(x;h+g) &= \limsup_{\substack{y\to x\\t\downto 0}} \frac{F(y+th+tg)-F(y)}{t} \\
            &\leq
            \limsup_{\substack{y\to x\\t\downto 0}} \frac{F(y+th+tg)-F(y+tg)}{t}+ \limsup_{\substack{y\to x\\t\downto 0}} \frac{F(y+tg)-F(y)}{t}\\
            &= F^\circ(x;h)+ F^\circ(x;g).
        \end{aligned}
    \end{equation*}

    \emph{\ref{lem:clarke:dir:iii}:} The claim is clear for $\alpha=0$. For $\alpha>0$, we obtain again from the definition that
    \begin{equation*}
        \begin{aligned}
            F^\circ(x;\alpha h) &= \limsup_{\substack{y\to x\\t\downto 0}} \frac{F(y+t(\alpha h))-F(y)}{t}\\
            &= \limsup_{\substack{y\to x\\\alpha t\downto 0}} \alpha \frac{F(y+(\alpha t)h)-F(y)}{\alpha t} =  (\alpha F)^\circ(x;h).
        \end{aligned}
    \end{equation*}

    \emph{\ref{lem:clarke:dir:iv}:} Similarly, since $t\downto 0$ and $h\in X$ is fixed, $y\to x$ if and only if $w\defeq y-th \to x$. We thus have that
    \begin{equation*}
        \begin{aligned}[b]
            F^\circ(x;-h) &= \limsup_{\substack{y\to x\\t\downto 0}} \frac{F(y-th)-F(y)}{t}\\
            &= \limsup_{\substack{w\to x\\t\downto 0}} \frac{-F(w+th)-(-F(w))}{t} =  (-F)^\circ(x;h).
        \end{aligned}
        \qedhere
    \end{equation*}
\end{proof}
In particular, \cref{lem:clarke:dir}\,\ref{lem:clarke:dir:i}--\ref{lem:clarke:dir:iii} imply that the mapping  $h\mapsto F^\circ(x;h)$ is proper, convex, and lower semicontinuous.

We now define for a locally Lipschitz continuous functional $F:X\to \R$ the \term[subdifferential!Clarke]{Clarke subdifferential} in $x\in X$ as
\begin{equation}\label{eq:clarke:def}
    \partial_C F(x) \defeq \setof{x^*\in X^*}{\dual{x^*,h}_X \leq F^\circ (x;h)\text{ for all }h\in X}.
\end{equation}
The definition together with \cref{lem:clarke:dir}\,\ref{lem:clarke:dir:i} directly implies the following properties.
\begin{lemma}\label{lem:clarke:properties}
    Let $F:X\to\R$ be locally Lipschitz continuous and $x\in X$. Then $\partial_C F(x)$ is convex, weakly-$*$ closed, and bounded. Specifically, if $F$ is Lipschitz near $x$ with constant $L$, then $\partial _C F(x)\subset \B(0,L)$.
\end{lemma}
Furthermore, we have the following useful continuity property.
\begin{lemma}\label{lem:clarke:closed}
    Let $F:X\to\R$. Then $\partial_C F(x)$ is strong-to-weak-$*$ outer semicontinuous, i.e., if $x_n\to x$ and if $ \partial_C F(x_n)\ni x_n^*\weaktostar x^*$, then $x^*\in \partial_C F(x)$.
\end{lemma}
\begin{proof}
    Let $h\in X$ be arbitrary.
    By assumption, we then have that $\dual{x_n^*,h}_X \leq F^\circ(x_n;h)$ for all $n\in \N$.
    The weak-$*$ convergence of $\{x_n^*\}_{n\in\N}$ then implies that
    \begin{equation*}
        \dual{x^*,h}_X = \lim_{n\to\infty} \dual{x_n^*,h}_X \leq \limsup_{n\to\infty} F^\circ(x_n;h).
    \end{equation*}
    Hence we are finished if we can show that $\limsup_{n\to\infty} F^\circ(x_n;h)\leq F^\circ(x;h)$ (since then $x^*\in \partial_C F(x)$ by definition).

    For this, we use that by definition of $F^\circ(x_n;h)$, there exist sequences $\{y_{n,m}\}_{m\in\N}$ and $\{t_{n,m}\}_{m\in\N}$ with $y_{n,m}\to x_n$ and $t_{n,m}\downto 0$ for $m\to \infty$ realizing the $\limsup$ for each $x_n$. Hence, for all $n\in \N$ we can find a $y_n\defeq y_{n,m(n)}$ and a $t_n\defeq t_{n,m(n)}$ such that $\norm{y_n-x_n}_X+t_n < n^{-1}$ (and hence in particular $y_n\to x$ and $t_n\downto 0$) as well as
    \begin{equation*}
        F^\circ(x_n;h) - \tfrac1n \leq \frac{F(y_n+t_nh)-F(y_n)}{t_n}
    \end{equation*}
    for $n$ sufficiently large. Taking the $\limsup$ for $n\to\infty$ on both sides yields the desired inequality.
\end{proof}

Again, the construction immediately yields a Fermat principle.\footnote{Similarly to \cref{thm:convex:fermat}, we do not need to require Lipschitz continuity of $F$ -- the Fermat principle for the Clarke subdifferential characterizes (among others) \emph{any} local minimizer. However, if we want to use this principle to verify that a given $\bar x\in X$ is indeed a (candidate for) a minimizer, we need a suitable characterization of the subdifferential -- and this is only possible for (certain) locally Lipschitz continuous functionals.}
\begin{theorem}[Fermat principle]\label{thm:clarke:fermat}\index{principle!Fermat!Clarke}
    If $F:X\to\R$ has a local minimum in $\bar x$, then $0\in \partial_C F(\bar x)$.
\end{theorem}
\begin{proof}
    If $\bar x\in X$ is a local minimizer of $F$, then $F(\bar x) \leq F(\bar x + th)$ for all $h\in X$ and $t>0$ sufficiently small (since the topological interior is always included in the algebraic interior). But this implies that
    \begin{equation*}
        \dual{0,h}_X = 0  \leq \liminf_{t\downto 0} \frac{F(\bar x+th)-F(\bar x)}{t}\leq \limsup_{t\downto 0} \frac{F(\bar x+th)-F(\bar x)}{t}\leq F^\circ(x;h)
    \end{equation*}
    and hence $0\in \partial_C F(\bar x)$ by definition.
\end{proof}
Note that $F$ is not assumed to be convex, and hence the condition is in general not sufficient (consider, e.g., $f(t) = -|t|$).

\section{Fundamental examples}

Next, we show that the Clarke subdifferential is indeed a generalization of the derivative concepts we've studied so far.
\begin{theorem}\label{thm:clarke:frechet}
    Let $F:X\to\R$ be continuously differentiable in a neighborhood $U$ of $x\in X$. Then $\partial_C F(x) = \{F'(x)\}$.
\end{theorem}
\begin{proof}
    First, we note that $F$ is locally Lipschitz continuous near $x$ by \cref{lem:variation:c1-lipschitz}.
    We now show that $F^\circ(x;h)=F'(x)h$ $(=F'(x;h))$ for all $h\in X$. Take again sequences $\{y_n\}_{n\in\N}$ and $\{t_n\}_{n\in\N}$ with $y_n\to x$ and $t_n\downto 0$ realizing the $\limsup$ in \eqref{eq:clarke:dir}. Applying the mean value \cref{thm:frechet:mean} and using the continuity of $F'$ yields for any $h\in X$ that
    \begin{equation*}
        \begin{aligned}
            F^\circ(x;h) &= \lim_{n\to\infty} \frac{F(y_n+t_nh)-F(y_n)}{t_n}\\
            &=\lim_{n\to\infty} \int_{0}^1 \frac1{t_n}\dual{F'(y_n+s(t_nh)),t_nh}_X\,ds\\
            &=\dual{F'(x),h}_X
        \end{aligned}
    \end{equation*}
    since the integrand converges uniformly in $s\in [0,1]$ to $\dual{F'(x),h}_X$.
    Hence by definition, $x^*\in \partial_C F(x)$ if and only if $\dual{x^*,h}_X \leq \dual{F'(x),h}_X$ for all $h\in X$, which is only possible for $x^*=F'(x)$.
\end{proof}

The following example shows that \cref{thm:clarke:frechet} does \emph{not} hold if $F$ is merely Fréchet differentiable.
\begin{example}\label{ex:clarke:frechet}
    Let $F:\R\to\R$, $F(x)=x^2\sin(x^{-1})$. Then it is straightforward (if tedious) to show that $F$ is differentiable on $\R$ with
    \begin{equation*}
        F'(x) =
        \begin{cases}
            2x\sin(x^{-1}) - \cos(x^{-1}) & \text{if } x\neq 0,\\
            0 & \text{if } x=0.
        \end{cases}
    \end{equation*}
    In particular, $F$ is not continuously differentiable at $x=0$. But a similar limit argument shows that for all $h\in \R$,
    \begin{equation*}
        F^\circ(0; h) = |h|
    \end{equation*}
    and hence that
    \begin{equation*}
        \partial_C F(0) = [-1,1] \supsetneq \{0\} = \{F'(0)\}.
    \end{equation*}
    (The first equality also follows more directly from \cref{thm:clarke:gradient} below.)
\end{example}

As the example suggests, we always have the following weaker relation.
\begin{lemma}\label{lem:clarke:gateaux}
    Let $F:X\to\R$ be Lipschitz continuous near $x\in X$ and Gateaux differentiable at $x$. Then $DF(x)\in \partial_C F(x)$.
\end{lemma}
\begin{proof}
    Let $h\in X$ be arbitrary. First, note that we always have that
    \begin{equation}\label{eq:clarke:dir-gen-dir-inequality}
        F'(x;h) = \lim_{t\downto 0} \frac{F(x+th)-F(x)}{t} \leq  \limsup_{\substack{y\to x\\t\downto 0}} \frac{F(y+th)-F(y)}{t}
        = F^\circ(x;h).
    \end{equation}
    Since $F$ is Gateaux differentiable, it follows that
    \begin{equation*}
        \dual{DF(x),h}_X = F'(x; h) \leq F^\circ (x; h) \qquad\text{for all }h\in X,
    \end{equation*}
    and thus $DF(x)\in \partial_C F(x)$ by definition.
\end{proof}
Similarly, the Clarke subdifferential reduces to the convex subdifferential in some situations.
\begin{theorem}\label{thm:clarke:convex}
    Let $F:X\to\Rbar$ be convex and lower semicontinuous. Then $\partial_C F(x) = \partial F(x)$ for all $x\in \interior(\dom F)$.
\end{theorem}
\begin{proof}
    By \cref{thm:convex:cont}, $F$ is locally Lipschitz continuous near $x\in\interior(\dom F)$.
    We now show that $F^\circ(x;h)=F'(x;h)$ for all $h\in X$, which together with \cref{lem:convex:equiv} yields the claim. By \eqref{eq:clarke:dir-gen-dir-inequality}, we always have that $F'(x; h) \leq F^\circ(x; h)$.
    To show the reverse inequality, let $\delta>0$ be arbitrary.
    Since the difference quotient of convex functionals is increasing by \cref{lem:convex:direct}\,\ref{lem:convex:direct:i}, we obtain that
    \begin{equation*}
        \begin{aligned}
            F^\circ(x;h) &= \lim_{\eps\downto 0} \sup_{y\in \B(x,\delta\eps)}\sup_{0<t<\eps} \frac{F(y+th)-F(y)}{t}\\
            &\leq \lim_{\eps\downto 0} \sup_{y\in \B(x,\delta\eps)} \frac{F(y+\eps h)-F(y)}{\eps}\\
            &\leq \lim_{\eps\downto 0} \frac{F(x+\eps h)-F(x)}{\eps} + 2L \delta\\
            &= F'(x;h) + 2L \delta,
        \end{aligned}
    \end{equation*}
    where the last inequality follows by adding two productive zeros and using the local Lipschitz continuity in $x$. Since $\delta>0$ was arbitrary, this implies that $F^\circ(x;h) \leq F'(x;h)$, and the claim follows.
\end{proof}
A locally Lipschitz continuous functional $F:X\to\R$ with $F^\circ(x;h)=F'(x;h)$ for all $h\in X$ is called \term[functional!regular!Clarke]{regular (in the sense of Clarke)}\index{regularity!Clarke} in $x\in X$. We have just shown that every continuously differentiable and every convex and lower semicontinuous functional is regular; intuitively, a function is thus regular at any points in which it is either differentiable or has at most a \enquote{convex kink}.

\bigskip

Finally, similarly to \cref{thm:lebesgue:subdiff} one can show the following pointwise characterization of the Clarke subdifferential of integral functionals with Lipschitz continuous integrands. We again assume that $\Omega\subset\R^d$ is open and bounded.
\begin{theorem}\label{thm:clarke:pointwise}
    Let $f:\R\to\R$ be Lipschitz continuous and $F:L^p(\Omega)\to\Rbar$ with $1\leq p <\infty$ as in \cref{lem:lebesgue:lsc}. Then we have for all $u\in L^p(\Omega)$ with $q=\frac{p}{p-1}$ (where $q=\infty$ for $p=1$) that
    \begin{equation*}
        \partial_C F(u)  \subset \setof{u^*\in L^q(\Omega)}{u^*(x) \in\partial_C f(u(x))\text{ for almost every } x\in\Omega}.
    \end{equation*}
    If $f$ is regular at $u(x)$ for almost every $x\in \Omega$, then $F$ is regular at $u$, and equality holds.
\end{theorem}
\begin{proof}
    First, by the properties of the Lebesgue integral and the Lipschitz continuity of $f$, we have for any $u,v\in L^p(\Omega)$ that
    \begin{equation*}
        |F(u)-F(v)| \leq \int_\Omega |f(u(x))-f(v(x))|\,dx \leq L \int_\Omega |u(x)-v(x)|\,dx \leq L C_p \norm{u-v}_{L^p},
    \end{equation*}
    where $L$ is the Lipschitz constant of $f$ and $C_p$ the constant from the continuous embedding $L^p(\Omega)\hookrightarrow L^1(\Omega)$ for $1\leq p\leq \infty$. Hence $F:L^p(\Omega)\to \R$ is Lipschitz continuous and therefore finite-valued as well.

    Let now $\xi\in \partial_C F(u)\subset L^p(\Omega)^*$ be given and $h\in L^p(\Omega)$ be arbitrary. By definition, we thus have
    \begin{equation}
        \label{eq:clarke:pointwise:integral-ineq}
        \begin{aligned}[t]
            \dual{\xi,h}_{L^p} \leq F^\circ(u; h) &=\limsup_{\substack{v\to u\\t\downto 0}} \frac{F(v+th)-F(v)}{t}\\
            &\leq \int_\Omega \limsup_{\substack{v\to u\\t\downto 0}} \frac{f(v(x)+th(x))-f(v(x))}{t} \,dx\\
            &\leq \int_\Omega \limsup_{\substack{v_x\to u(x)\\t_x\downto 0}} \frac{f(v_x+t_xh(x))-f(v_x)}{t_x} \,dx\\
            &= \int_\Omega f^\circ(u(x); h(x))\,dx,
        \end{aligned}
    \end{equation}
    where we were able to use the reverse Fatou lemma to exchange the $\limsup$ with the integral in the first inequality since the integrand is bounded from above by the integrable function $L |h|$ due to \cref{lem:clarke:dir}\,\ref{lem:clarke:dir:i}; the second inequality follows by bounding for almost every $x\in \Omega$ the (pointwise) limit over the sequences realizing the $\limsup$ in the second line by the $\limsup$ over all admissible sequences.

    In order to interpret \eqref{eq:clarke:pointwise:integral-ineq} pointwise, we use that \cref{lem:clarke:dir}\,\ref{lem:clarke:dir:i} together with the (global) Lipschitz continuity of $f$ implies that the function $x\mapsto f^\circ(u(x);t)$ is integrable for any $t\in\R$.
    We can thus argue exactly as in the proof of \cref{thm:lebesgue:subdiff}: Let $t\in \R$ be arbitrary and $A\subset \Omega$ be an arbitrary measurable subset. Setting
    \begin{equation*}
        h(x) = \begin{cases} t & \text{if }x\in A, \\ 0 & \text{if } x\notin A, \end{cases}
    \end{equation*}
    (so that $h\in L^\infty(\Omega)\subset L^p(\Omega)$)
    and using $f^\circ(u(x);0)=0$, we obtain from \eqref{eq:clarke:pointwise:integral-ineq} together with the representation of $\xi \in L^p(\Omega)^*$ via some $u^*\in L^q(\Omega)$ that
    \begin{equation*}
        \int_A u^*(x)t\,dx = \dual{\xi,h}_{L^p} \leq \int_\Omega f^\circ(u(x); h(x))\,dx = \int_A f^\circ(u(x); t)\,dx.
    \end{equation*}
    Since $A$ was arbitrary, this implies that
    \begin{equation*}
        u^*(x)t \leq f^\circ(u(x); t) \qquad\text{for almost every }x \in \Omega.
    \end{equation*}
    Since $t\in \R$ was arbitrary, we obtain $u^*(x) \in \partial_C f(u(x))$ almost everywhere.

    It remains to show the remaining assertions when $f$ is regular. In this case, it follows from \eqref{eq:clarke:pointwise:integral-ineq} that for any $h\in L^p(\Omega)$,
    \begin{equation}
        \label{eq:clarke:pointwise:regular-est}
        \begin{aligned}[t]
            F^\circ(u; h) & \leq \int_\Omega f^\circ(u(x);h(x))\,dx = \int_\Omega f'(u(x); h(x))\,dx
            \\
            &
            \leq \lim_{t\downto 0} \frac {F(u+th)-F(u)}{t} = F'(u; h) \leq F^\circ(u; h),
        \end{aligned}
    \end{equation}
    where the second inequality is obtained by applying Fatou's lemma, this time appealing to the integrable lower bound $- L|h(x)|$. This shows that $F'(u; h)= F^\circ(u; h)$ and hence that $F$ is regular. We further obtain for any $u^* \in L^q(\Omega)$ with $u^*(x) \in \partial_C f(u(x))$ almost everywhere and any $h\in L^p(\Omega)$, that
    \begin{equation*}
        \dual{u^*,h}_{L^p} = \int_\Omega u^*(x) h(x) \, dx \leq \int_\Omega f^\circ(u(x); h(x))\,dx
        \leq F^\circ(u, h),
    \end{equation*}
    where we have used \eqref{eq:clarke:pointwise:regular-est} in the last inequality.
    Since $h\in L^p(\Omega)$ was arbitrary, this implies that $u^*\in \partial_C F(u)$.
\end{proof}
Under additional assumptions similar to those of \cref{thm:superpos:continuous} and with more technical arguments, this result can be extended to spatially varying integrands $f:\Omega\times \R\to \R$; see, e.g., \cite[Theorem 2.7.5]{Clarke:1990a}.

\section{Calculus rules}

We now turn to calculus rules.
The first one follows directly from the definition of the Clarke subdifferential.
\begin{theorem}\label{thm:clarke:scalar}
    Let $F:X\to\R$ be locally Lipschitz continuous near $x\in X$ and $\alpha\in\R$. Then,
    \begin{equation*}
        \partial_C(\alpha F)(x) = \alpha \partial_C(F)(x).
    \end{equation*}
\end{theorem}
\begin{proof}
    First, $\alpha F$ is clearly locally Lipschitz continuous near $x$ for any $\alpha\in\R$. If $\alpha = 0$, both sides of the claimed equality are zero (which is easiest seen from \cref{thm:clarke:frechet}).
    If $\alpha > 0$, we have that $(\alpha F)^\circ(x;h) = \alpha F^\circ (x;h)$ for all $h\in X$ from the definition. Hence,
    \begin{equation*}
        \begin{aligned}
            \alpha \partial_C F(x)
            &= \setof{\alpha x^*\in X^*}{\dual{x^*,h}_X \leq F^\circ(x;h)\quad\text{for all }h\in X}\\
            &= \setof{\alpha x^*\in X^*}{\dual{\alpha x^*,h}_X \leq \alpha F^\circ(x;h)\quad\text{for all }h\in X}\\
            &= \setof{y^*\in X^*}{\dual{y^*,h}_X \leq (\alpha F)^\circ(x;h)\quad\text{for all }h\in X}\\
            &= \partial_C(\alpha F)(x).
        \end{aligned}
    \end{equation*}
    To conclude the proof, it suffices to show the claim for $\alpha=-1$.
    For that, we use \cref{lem:clarke:dir}\,\ref{lem:clarke:dir:iv} to obtain that
    \begin{equation*}
        \begin{aligned}[b]
            \partial_C (-F)(x)
            &= \setof{x^*\in X^*}{\dual{x^*,h}_X \leq (-F)^\circ(x;h)\quad\text{for all }h\in X}\\
            &= \setof{x^*\in X^*}{\dual{-x^*,-h}_X \leq F^\circ(x;-h)\quad\text{for all }h\in X}\\
            &= \setof{-y^*\in X^*}{\dual{y^*,g}_X \leq F^\circ(x;g)\quad\text{for all }g\in X}\\
            &= -\partial_C F(x).
        \end{aligned}
        \qedhere
    \end{equation*}
\end{proof}
\begin{corollary}\label{lem:clarke:fermat2}
    Let $F:X\to\R$ be locally Lipschitz continuous near $\bar x\in X$. If $F$ has a local maximum in $\bar x$, then $0\in \partial_C F(\bar x)$.
\end{corollary}
\begin{proof}
    If $\bar x$ is a local maximizer of $F$, it is a local minimizer of $-F$. Hence, \cref{thm:clarke:fermat,thm:clarke:scalar} imply that
    \begin{equation*}
        0\in \partial_C(-F)(\bar x) = -\partial_C F(\bar x),
    \end{equation*}
    i.e., $0=-0\in \partial_C F(\bar x)$.
\end{proof}

\subsection*{Support functionals}
\label{sec:clarke:support}

The remaining rules are significantly more involved.
As in the previous proofs, a key step is to relate different sets of the form \eqref{eq:clarke:def}, which we will do with the help of the following lemmas due to \cite{Hoermander:1955}.
\begin{lemma}\label{lem:clarke:support1}
    Let $S:X\to\R$ be positively homogeneous, subadditive, and lower semicontinuous, and let
    \begin{equation*}
        A=\setof{x^*\in X^*}{\dual{x^*,x}_X\leq S(x)\quad\text{for all }x\in X}.
    \end{equation*}
    Then
    \begin{equation}\label{eq:clarke:support1}
        S(x) =  \sup_{x^*\in A}\, \dual{x^*,x}_X\qquad\text{for all }x\in X.
    \end{equation}
\end{lemma}
\begin{proof}
    By definition of $A$, the inequality $\dual{x^*,x}_X-S(x)\leq 0$ holds for all $x\in X$ if and only if $x^*\in A$. Thus a case distinction as in \cref{ex:convex:fenchel}\,\ref{ex:convex:fenchel:iii} using the positive homogeneity of $S$ (which in particular implies that $S(0)=0$) shows that
    \begin{equation*}
        S^*(x^*) = \sup_{x\in X}\, \dual{x^*,x}_X - S(x) = \begin{cases} 0 & x^* \in A,\\\infty & x^*\notin A,\end{cases}
    \end{equation*}
    i.e., $S^*=\delta_A$.
    Furthermore, by assumption $S$ is also subadditive and hence convex as well as lower semicontinuous; it is also proper.
    \Cref{thm:convex:moreau} thus yields
    \begin{equation}\label{eq:clarke:support_finite}
        S(x) = S^{**}(x) = (\delta_A)_*(x) = \sup_{x^*\in A} \, \dual{x^*,x}_X.
        \qedhere
    \end{equation}
\end{proof}
The right-hand side of \eqref{eq:clarke:support1} is called the \term[functional!support]{support functional} of $A\subset X^*$; see, e.g., \cite{Hiriart:2001} for their use in convex analysis (in finite dimensions).
Note that \eqref{eq:clarke:support_finite} implies that any set of the form $A$ is nonempty since the supremum over the empty set is $-\infty$ and $S$ was assumed to be real-valued.

\begin{lemma}\label{lem:clarke:support3}
    Let $A,B\subset X^*$ be nonempty, convex, and weakly-$*$ closed. Then $A\subset B$ if and only if
    \begin{equation}\label{eq:clarke:support3}
        \sup_{x^*\in A}\,  \dual{x^*,x}_X \leq \sup_{x^*\in B}\, \dual{x^*,x}_X \qquad\text{for all }x\in X.
    \end{equation}
\end{lemma}
\begin{proof}
    If $A\subset B$, then the right-hand side of \eqref{eq:clarke:support3} is obviously not less than the left-hand side.
    Conversely, assume that there exists an $x^*\in A$ with $x^*\notin B$.
    By the assumptions on $A$ and $B$, we then obtain from  \cref{thm:clarke:hb} an $x\in X$ and a $\lambda\in \R$ with
    \begin{equation*}\label{eq:clarke:hb1}
        \dual{z^*,x}_X \leq \lambda < \dual{x^*,x}_X  \qquad\text{for all }z^*\in B.
    \end{equation*}
    Taking the supremum over all $z^*\in B$ and estimating the right-hand side by the supremum over all $x^*\in A$ then yields that
    \begin{equation*}\label{eq:clarke:hb2}
        \sup_{z^*\in B}\, \dual{z^*,x}_X < \sup_{x^*\in A}\,\dual{x^*,x}_X.
    \end{equation*}
    Hence \eqref{eq:clarke:support3} is violated, and the claim follows by contraposition.
\end{proof}
\begin{corollary}\label{lem:clarke:support2}
    Let $A,B\subset X^*$ be nonempty, convex, and weakly-$*$ closed. Then $A=B$ if and only if
    \begin{equation}\label{eq:clarke:support2}
        \sup_{x^*\in A}\, \dual{x^*,x}_X = \sup_{x^*\in B}\, \dual{x^*,x}_X \qquad\text{for all }x\in X.
    \end{equation}
\end{corollary}
\begin{proof}
    Again, the claim is obvious if $A=B$. Conversely, if \eqref{eq:clarke:support2} holds, then in particular \eqref{eq:clarke:support3} holds, and we obtain from \cref{lem:clarke:support3} that $A\subset B$. Exchanging the roles of $A$ and $B$ now yields the claim.
\end{proof}

\bigskip

Since generalized directional derivatives are always real-valued, \cref{lem:clarke:support1} together with \cref{lem:clarke:dir} directly yields the following useful representation.
\begin{corollary}\label{cor:clarke:support-dir}
    Let $F:X\to\R$ be locally Lipschitz continuous and $x\in X$. Then
    \begin{equation*}
        F^\circ(x;h) = \sup_{x^*\in\partial_C F(x)}\dual{x^*,h}_X\quad\text{for all }h\in X.
    \end{equation*}
    In particular, $\partial_C F(x)$ is nonempty.
\end{corollary}
For example, this implies a converse result to \cref{thm:clarke:frechet}.
\begin{corollary}\label{cor:clarke:single-valued}
    Let $F:X\to\R$ be locally Lipschitz continuous near $x$. If $\partial_C F(x)=\{x^*\}$ for some $x^*\in X^*$, then $F$ is Gateaux differentiable at $x$ with $DF(x) = x^*$.
\end{corollary}
\begin{proof}
    Under the assumption, it follows from \cref{cor:clarke:support-dir} that
    \begin{equation*}
        F^\circ(x;h) = \sup_{\tilde x^*\in\partial F_C(x)}\dual{\tilde x^*,h}_X = \dual{x^*,h}_X
    \end{equation*}
    for all $h\in X$. In particular, $F^\circ(x;h)$ is linear (and not just reflective) in $h$. It thus follows from \cref{lem:clarke:dir}\,\cref{lem:clarke:dir:iv} that for any $h\in X$,
    \begin{equation*}
        \begin{aligned}
            \liminf_{\substack{y\to x\\t\downto 0}} \frac{F(y+th)-F(y)}{t} &=
            -\limsup_{\substack{y\to x\\t\downto 0}} \frac{-F(y+th)-(-F(y))}{t} \\
            &=
            -(-F)^\circ(x;h) = -F^\circ(x;-h) = F^\circ(x,h) \\
            &=
            \limsup_{\substack{y\to x\\t\downto 0}} \frac{F(y+th)-F(y)}{t}.
        \end{aligned}
    \end{equation*}
    Hence the $\limsup$ is a proper limit, and thus $F^\circ(x;h) = F'(x;h)$; i.e., $F$ is regular in $x$.
    This shows that $F'(x;h)$ is linear and bounded in $h$, and hence $x^*$ is by definition the Gateaux derivative.
\end{proof}
It is not hard to verify from the definitions of the Clarke subdifferential and of the Fréchet derivative using the Lipschitz continuity of $F$ that in this case, $x^*$ is in fact a Fréchet derivative.

We can also use this to show the promised nonemptiness of the convex subdifferential.
\begin{theorem}\label{cor:convex:nonempty}
    Let $X$ be a Banach space and let $F:X\to\Rbar$ be proper, convex, and lower semicontinuous, and $x\in\interior(\dom F)$. Then $\partial F(x)$ is nonempty, convex, weakly-$*$ closed, and bounded.
\end{theorem}
\begin{proof}
    Since $x\in \interior(\dom F)$, \cref{thm:clarke:convex} shows that $\partial F(x) = \partial_C F(x) \neq \emptyset$ by \cref{cor:clarke:support-dir}. The remaining properties follow similarly from \cref{lem:clarke:properties}.
\end{proof}
By a similar argument, we now obtain the promised converse of \cref{thm:convex:gateaux}; we combine both statements here for the sake of reference.
\begin{theorem}\label{thm:convex:singleton}
    Let $X$ be a Banach space and let $F:X\to\Rbar$ be convex. If $F$ is Gateaux differentiable at $x$, then $\partial F(x) = \{DF(x)\}$.
    Conversely, if $x\in \interior (\dom F)$ and $\subdiff F(x)=\{x^*\}$ is a singleton, then $F$ is Gateaux differentiable at $x$ with $DF(x)=x^*$.
\end{theorem}
\begin{proof}
    The first claim was already shown in \cref{thm:convex:gateaux}, while the second follows from \cref{cor:clarke:single-valued} together with \cref{thm:clarke:convex}.
\end{proof}

As another consequence, we can show that Moreau--Yosida regularization defined in \cref{sec:moreau-yosida} preserves (global!) Lipschitz continuity.
\begin{lemma}
    \label{lemma:moreau:lipschitz}
    Let $X$ be a Hilbert space and let $F: X \to \R$ be Lipschitz continuous with constant $L$.
    Then $F_\gamma$ is Lipschitz continuous with constant $L$ as well.
    If $F$ is in addition convex, then $F - \tfrac{\gamma L^2}{2} \le F_\gamma \le F$.
\end{lemma}
\begin{proof}
    Let $x, z \in X$. We expand
    \begin{equation*}
        F_\gamma(x)-F_\gamma(z)
        =\sup_{y_z \in X} \inf_{y_x \in X}\left( F(y_x)-F(y_z) + \frac{1}{2\gamma}\norm{y_x-x}_X^2 - \frac{1}{2\gamma}\norm{y_z-z}_X^2 \right).
    \end{equation*}
    Taking $y_x=y_z+x-z$, we estimate
    \begin{equation*}
        F_\gamma(x)-F_\gamma(z) \le \sup_{y_z \in X} \left( F(y_z+x-z) - F(y_z) \right) \le L\norm{x-z}_X.
    \end{equation*}
    Exchanging $x$ and $z$, we obtain the first claim.

    For the second claim, we first observe that by assumption $\dom F=X$. Hence by \cref{cor:convex:nonempty,lem:clarke:properties}, for every $x \in X$, there exists some $x^* \in \subdiff F(x)$ with $\norm{x^*}_{X^*} \le L$.
    Thus, using \cref{lem:convex:equiv}, for any $x^* \in \subdiff F(x)$,
    \begin{equation*}
        F_\gamma(x) = \inf_{y \in X} F(z) + \frac{1}{2\gamma}\norm{x-y}_X^2
        \ge
        F(x) + \iprod{x^*}{z-x}_X + \frac{1}{2\gamma}\norm{x-z}_X^2.
    \end{equation*}
    The Cauchy--Schwarz and generalized Young's inequality then yield $F_\gamma(x) \ge F(x) - \tfrac{\gamma}{2}\norm{x^*}_{X^*}^2 \ge F(x) - \tfrac\gamma2 L^2$.
    The second inequality follows by estimating the infimum in \eqref{eq:proximal:moreau-envelope} by $z=x$.
\end{proof}

\subsection*{Sum rule}

With the aid of these results on support functionals, we can now show a sum rule.
\begin{theorem}\label{thm:clarke:sum}
    Let $F,G:X\to\R$ be locally Lipschitz continuous near $x\in X$. Then,
    \begin{equation*}
        \partial_C (F+G)(x) \subset \partial_C F(x) + \partial_C G(x).
    \end{equation*}
    If $F$ and $G$ are regular at $x$, then $F+G$ is regular at $x$ and equality holds.
\end{theorem}
\begin{proof}
    It is clear that $F+G$ is locally Lipschitz continuous near $x$. Furthermore, from the properties of the $\limsup$ we always have for all $h\in X$ that
    \begin{equation}
        \label{eq:clarke:sum:gdirdiff-est0}
        (F+G)^\circ(x;h) \leq F^\circ (x;h) + G^\circ (x;h).
    \end{equation}
    If $F$ and $G$ are regular at $x$, the calculus of limits yields that
    \begin{equation*}
        F^\circ (x;h) + G^\circ (x;h) = F'(x;h) + G'(x;h) = (F+G)'(x;h) \leq (F+G)^\circ(x;h),
    \end{equation*}
    which implies that $(F+G)^\circ(x;h) = (F+G)'(x;h)$, i.e., $F+G$ is regular.

    By the definition of the Clarke subdifferential, it follows from \eqref{eq:clarke:sum:gdirdiff-est0}
    \begin{equation*}
        \partial_C (F+G)(x) \subset \setof{x^*\in X^*}{\dual{x^*,h}_X \leq F^\circ (x;h) + G^\circ (x;h)\text{ for all }h\in X}=:A
    \end{equation*}
    (with equality if $F$ and $G$ are regular); it thus remains to show that $A = \partial_C F(x) + \partial_C G(x)$.
    For this, we use that $\partial_C F(x)$ and $\partial_C G(x)$ are convex and weakly-$*$ closed by \cref{lem:clarke:properties} and nonempty by \cref{cor:clarke:support-dir}, and hence so is their sum since both sets are bounded.
    Furthermore, as shown in \cref{lem:clarke:dir}, generalized directional derivatives and hence their sums are real-valued, positively homogeneous, convex, and lower semicontinuous. We thus obtain from \cref{lem:clarke:support1} for all $h\in X$ that
    \begin{equation*}
        \begin{aligned}[b]
            \sup_{x^* \in \partial_C F(x) + \partial_C G(x)}\dual{x^*,h}_X
            &= \sup_{x_1^*\in \partial_C F(x)} \dual{x_1^*,h}_X + \sup_{x_2^*\in \partial_C G(x)} \dual{x_2^*,h}_X \\
            &= F^\circ(x;h) + G^\circ(x;h)
            = \sup_{x^*\in A}\,\dual{x^*,h}_X.
        \end{aligned}
    \end{equation*}
    The claimed equality of $A$ (which is nonempty, convex, and weakly-$*$ closed as well) and the sum of the subdifferentials now follows from \cref{lem:clarke:support2}.
\end{proof}
Note the differences from the convex sum rule: The generic inclusion is now in the other direction; furthermore, \emph{both} functionals have to be regular, and in exactly the point where the sum rule is applied.
By induction, one obtains from this a sum rule for an arbitrary number of functionals (which all have to be regular).

\subsection*{Chain rule}

To prove a chain rule, we need the following \enquote{nonsmooth} mean value theorem due to \cite{Lebourg:1975,Lebourg:1979}.
\begin{theorem}\label{thm:clarke:mean}
    Let $F:X\to\R$ be locally Lipschitz continuous near $x\in X$ and $\tilde x$ be in the Lipschitz neighborhood of $x$. Then there exists a $\lambda\in(0,1)$ and an $x^*\in \partial_C F(x+\lambda(\tilde x-x))$ such that
    \begin{equation*}
        F(\tilde x)-F(x) = \dual{x^*,\tilde x-x}_X.
    \end{equation*}
\end{theorem}
\begin{proof}
    Define $\psi,\phi:[0,1]\to\R$ as
    \begin{equation*}
        \psi(\lambda)\defeq F(x+\lambda(\tilde x-x)),\qquad \phi(\lambda)\defeq\psi(\lambda)+\lambda(F(x)-F(\tilde x)).
    \end{equation*}
    By the assumptions on $F$ and $\tilde x$, both $\psi$ and $\phi$ are Lipschitz continuous. In addition, $\phi(0) = F(x) = \phi(1)$, and hence $\phi$ has a local minimum or maximum in an interior point $\bar \lambda\in (0,1)$. From the Fermat principle \cref{thm:clarke:fermat} or \cref{lem:clarke:fermat2}, respectively, together with the sum rule from \cref{thm:clarke:sum} and the characterization of the subdifferential of the second term from \cref{thm:clarke:frechet}, we thus obtain that
    \begin{equation*}
        0\in \partial_C \phi(\bar \lambda) \subset \partial_C\psi(\bar \lambda) + \{F(x)-F(\tilde x)\}.
    \end{equation*}
    Hence we are finished if we can show for $x_{\bar \lambda} \defeq x +\bar \lambda(\tilde x-x)$ that
    \begin{equation}\label{eq:clarke:mean1}
        \partial_C\psi(\bar \lambda) \subset\setof{\dual{x^*,\tilde x-x}_X}{x^*\in\partial_C F(x_{\bar\lambda})} =:A.
    \end{equation}
    For this purpose, consider for arbitrary $s\in \R$ the generalized directional derivative
    \begin{equation*}
        \begin{aligned}
            \psi^\circ(\bar\lambda;s) &= \limsup_{\substack{\lambda\to \bar\lambda\\t\downto 0}} \frac{\psi(\lambda + ts)-\psi(\lambda)}{t}\\
            &=\limsup_{\substack{\lambda\to \bar\lambda\\t\downto 0}} \frac{F(x+(\lambda + ts)(\tilde x-x))-F(x+\lambda(\tilde x-x))}{t}\\
            &\leq \limsup_{\substack{z\to x_{\bar\lambda}\\t\downto 0}} \frac{F(z +ts(\tilde x-x))-F(z)}{t} = F^\circ(x_{\bar\lambda};s(\tilde x-x)),
        \end{aligned}
    \end{equation*}
    where the inequality follows from considering arbitrary sequences $z\to x_{\bar\lambda}$ (instead of special sequences of the form $z_n = x+\lambda_n(\tilde x-x)$) in the last $\limsup$. Again, the definition of the Clarke subdifferential thus implies that
    \begin{equation}\label{eq:clarke:mean2}
        \partial_C\psi(\bar \lambda) \subset\setof{t^*\in\R}{t^*s \leq F^\circ(x_{\bar\lambda};s(\tilde x-x))\text{ for all }s\in \R} =: B.
    \end{equation}
    It remains to show that the sets $A$ and $B$ from \eqref{eq:clarke:mean1} and \eqref{eq:clarke:mean2} coincide. But this follows again from \cref{lem:clarke:support1,lem:clarke:support2}, since for all $s\in\R$ we have that
    \begin{equation*}
        \sup_{t^*\in A}\,t^*s = \sup_{x^*\in \partial_C F(x_{\bar\lambda})} \dual{x^*,s(\tilde x-x)}_X = F^\circ (x_{\bar\lambda};s(\tilde x-x)) = \sup_{t^*\in B}\,t^*s.
        \qedhere
    \end{equation*}
\end{proof}
We also need the following generalization of the argument in \cref{thm:clarke:frechet}.
\begin{lemma}\label{lem:frechet:diffquot}
    Let $X,Y$ be Banach spaces and $F:X\to Y$ be continuously differentiable at $x\in X$. Let $\{x_n\}_{n\in\N}\subset X$ be a sequence with $x_n\to x$ and $\{t_n\}_{n\in\N}\subset(0,\infty)$ be a sequence with $t_n\downto 0$. Then for any $h\in X$,
    \begin{equation*}
        \lim_{n\to\infty} \frac{F(x_n+t_n h)-F(x_n)}{t_n} = F'(x)h.
    \end{equation*}
\end{lemma}
\begin{proof}
    Let $h\in X$ be arbitrary. By the Hahn--Banach extension \cref{thm:functan:hb-extension}, for every $n\in \N$ there exists a $y_n^*\in Y^*$ with $\norm{y_n^*}_{Y^*}=1$ and
    \begin{equation*}
        \norm{t_n^{-1}(F(x_n+t_nh)-F(x_n)) - F'(x)h}_Y = \dual{y_n^*,t_n^{-1}(F(x_n+t_nh)-F(x_n)) - F'(x)h}_Y.
    \end{equation*}
    Applying now the classical mean value theorem to the scalar functions
    \begin{equation*}
        f_n :[0,1]\to\R,\qquad f_n(s) = \dual{y_n^*,F(x_n+st_nh)}_Y,
    \end{equation*}
    we obtain similarly to the proof of \cref{thm:frechet:mean} for all $n\in\N$ that
    \begin{equation*}
        \begin{aligned}
            \norm{t_n^{-1}(F(x_n+t_nh)-F(x_n)) - F'(x)h}_Y &= t_n^{-1}\int_0^1\dual{y_n^*,F'(x_n+st_nh)t_nh}_Y\,ds - \dual{y_n^*,F'(x)h}_Y\\
            &= \int_0^1 \dual{y_n^*,[F'(x_n+st_nh)-F'(x)]h}_Y\,ds\\
            &\leq \int_0^1 \norm{F'(x_n+st_nh)-F'(x)}_{\linear(X; Y)}\,ds \,\norm{h}_X,
        \end{aligned}
    \end{equation*}
    where we have used \eqref{eq:functan:cs_banach} together with $\norm{y_n^*}_{Y^*}=1$ in the last step. Since $F'$ is continuous by assumption, the integrand goes to zero as $n\to\infty$ uniformly in $s\in[0,1]$, and the claim follows.
\end{proof}

We now come to the chain rule, which in contrast to the convex case does not require the inner mapping to be linear; this is one of the main advantages of the Clarke subdifferential in the context of nonsmooth optimization.
\begin{theorem}\label{thm:clarke:chain}
    Let $Y$ be a separable Banach space, $F:X\to Y$ be continuously differentiable at $x\in X$, and $G:Y\to \R$ be locally Lipschitz continuous near $F(x)$. Then,
    \begin{equation*}
        \partial_C (G\circ F)(x) \subset F'(x)^*\partial_C G(F(x)) \defeq \setof{F'(x)^*y^*}{y^*\in\partial_C G(F(x))}.
    \end{equation*}
    If $G$ is regular at $F(x)$, then $G\circ F$ is regular at $x$, and equality holds.
\end{theorem}
\begin{proof}
    The local Lipschitz continuity of $G\circ F$ follows from that of $G$ and $F$ (which in turn follows from \cref{lem:variation:c1-lipschitz}). For the claimed inclusion (respectively, equality), we argue as before using the support calculus.
    First we show that for every $h\in X$ there exists a $y^*\in \partial_C G(F(x))$ with
    \begin{equation}\label{eq:clarke:chain1}
        (G\circ F)^\circ(x;h) = \dual{y^*,F'(x)h}_Y.
    \end{equation}
    To this end, consider for given $h\in X$ sequences $\{x_n\}_{n\in\N}\subset X$ and $\{t_n\}_{n\in\N}\subset (0,\infty)$ with $x_n\to x$, $t_n \downto 0$, and
    \begin{equation*}
        (G\circ F)^\circ(x;h) = \lim_{n\to\infty} \frac{G(F(x_n+t_nh))-G(F(x_n))}{t_n}.
    \end{equation*}
    Let us now write $U_{F(x)}$ for the neighborhood of $F(x)$ where $G$ is Lipschitz with constant $L$.
    By continuity of $F$, we can then find $n_0\in\N$ such that $F(x_n), F(x_n+t_nh) \in U_{F(x)}$ for all $n\geq n_0$. \Cref{thm:clarke:mean} thus yields for all $n\geq n_0$ a $y_n^*\in \partial_C G(y_n)$ with $y_n: = F(x_n)+\lambda_n(F(x_n+t_nh)-F(x_n))$ for some $\lambda_n\in(0,1)$ such that
    \begin{equation}\label{eq:clarke:chain1a}
        \frac{G(F(x_n+t_nh))-G(F(x_n))}{t_n} = \dual{y_n^*,q_n}_Y
        \quad\text{with}\quad
        q_n \defeq  \frac{F(x_n+t_nh)-F(x_n)}{t_n}.
    \end{equation}
    Since $\lambda_n\in(0,1)$ is uniformly bounded, we also have that $y_n\to F(x)$ for $n\to \infty$.
    Hence $y_n$ is in the Lipschitz neighborhood of $F(x)$ for $n\in\N$ large enough,
    and \cref{lem:clarke:properties} yields that $y_n^*\in \partial_C G(y_n) \subset \B(0,L)$ for $n\in \N$ sufficiently large.
    This implies that $\{y_n^*\}_{n\in\N}\subset Y^*$ is bounded, and the
    Banach--Alaoglu theorem (\cref{thm:banachal}) yields a weakly-$*$ convergent subsequence with limit $y^*\in \partial_C G(F(x))$ by \cref{lem:clarke:closed}.
    Finally, since $F$ is continuously differentiable, $q_n \to F'(x)h$ strongly in $Y$ by \cref{lem:frechet:diffquot}. Hence, $\dualprod{y_n^*}{q_n}_Y \to \dualprod{y^*}{F'(x)h}$ as the duality pairing of weakly-$*$ and strongly converging sequences. Passing to the limit in \eqref{eq:clarke:chain1a} therefore yields \eqref{eq:clarke:chain1} (first along the subsequence chosen above; by convergence of the left-hand side of \eqref{eq:clarke:chain1a} and the uniqueness of limits then for the full sequence as well). By definition of the Clarke subdifferential, we thus have for $y^*\in\partial_C G(F(x))$ that
    \begin{equation}\label{eq:clarke:chain2}
        (G\circ F)^\circ(x;h) = \dual{y^*,F'(x)h}_Y\leq G^\circ(F(x);F'(x)h).
    \end{equation}

    If $G$ is now regular at $F(x)$, we have that $G^\circ(F(x);F'(x)h) = G'(F(x);F'(x)h)$ and hence by the local Lipschitz continuity of $G$ and the Fréchet differentiability of $F$ that
    \begin{multline*}
        G^\circ(F(x);F'(x)h) \\
        \begin{aligned}[t]
            &= \lim_{t\downto 0} \frac{G(F(x)+tF'(x)h) - G(F(x))}{t}\\
            &= \lim_{t\downto 0} \frac{G(F(x)+tF'(x)h) -G(F(x+th))+G(F(x+th))- G(F(x))}{t}\\
            &\leq \lim_{t\downto 0} \left(L\norm{h}_X\frac{\norm{F(x)+F'(x)th - F(x+th)}_Y}{\norm{th}_X} + \frac{G(F(x+th))- G(F(x))}{t} \right)\\
            &= (G\circ F)'(x;h) \leq (G\circ F)^\circ(x;h).
        \end{aligned}
    \end{multline*}
    Together with \eqref{eq:clarke:chain2}, this implies that $(G\circ F)'(x;h)  = (G\circ F)^\circ(x;h)$ (i.e., $G\circ F$ is regular at $x$) and that
    \begin{equation}
        \label{eq:clarke:chain3}
        (G\circ F)^\circ(x;h) = G^\circ(F(x);F'(x)h).
    \end{equation}

    As before, \cref{lem:clarke:support1} now implies for all $h\in X$ that
    \begin{equation*}
        \sup_{x^*\in F'(x)^*\partial_CG(F(x))}\dual{x^*,h}_X = \sup_{y^*\in \partial_C G(F(x))}\dual{y^*,F'(x)h}_Y = G^\circ (F(x); F'(x)h)
    \end{equation*}
    and hence by \cref{lem:clarke:support3} that
    \begin{equation*}
        F'(x)^*\partial_C G(F(x)) = \setof{x^*\in X^*}{\dual{x^*,h}_X\leq G^\circ(F(x);F'(x)h)\text{ for all }h\in X}.
    \end{equation*}
    Combined with \eqref{eq:clarke:chain2} or \eqref{eq:clarke:chain3} and the definition of the Clarke subdifferential in \eqref{eq:clarke:def}, this now yields the claimed inclusion or equality, respectively, for the Clarke subdifferential of the composition.
\end{proof}
Again, the generic inclusion is the reverse of the one in the convex chain rule.
Note that equality in the chain rule also holds if $-G$ is regular, since we can then apply \cref{thm:clarke:chain} to $-G\circ F$ and use that $\partial_C (-G)(F(x)) = -\partial_C G(F(x))$ by \cref{thm:clarke:scalar}.
Furthermore, if $G$ is not regular but $F'(x)$ is surjective, a similar proof shows that equality (but not the regularity of $G\circ F$) holds in the chain rule; see \cite[Theorem 10.19]{Clarke:2013}.

\begin{example}
    As a simple example, we consider
    \begin{equation*}
        F:\R^2\to\R,\qquad (x_1,x_2)\mapsto |x_1x_2|,
    \end{equation*}
    which is not convex. To compute the Clarke subdifferential, we write $F=g\circ T$ for
    \begin{equation*}
        g:\R\to\R,\quad t\mapsto |t|,\qquad T:\R^2\to\R,\quad (x_1,x_2)\mapsto x_1x_2,
    \end{equation*}
    where $g$ is finite-valued, convex, and Lipschitz continuous, and hence regular at any $t\in\R$, and $T$ is continuously differentiable for all $x\in \R^2$ with Fréchet derivative
    \begin{equation*}
      T'(x):\R^2 \to \R,\qquad
      T'(x)h \defeq x_2h_1+x_1h_2.
    \end{equation*}
    Its adjoint is easily verified to be given by
    \begin{equation*}
      T'(x)^*:\R\to\R^2,
      \qquad
        T'(x)^*t \defeq
        \begin{psmallmatrix}
            x_2 t \\ x_1 t
        \end{psmallmatrix}.
    \end{equation*}
    Hence, \cref{thm:clarke:chain} together with \cref{thm:clarke:convex} yields that $F$ is regular at any $x\in\R^2$ and that
    \begin{equation*}
        \partial_C F(x) = T'(x)^*\partial g(T(x)) =
        \begin{pmatrix}
            x_2\\x_1
        \end{pmatrix}
        \sign(x_1x_2),
    \end{equation*}
    for the set-valued sign function from \cref{ex:convex:subdiff_abs}.
\end{example}

\section{Characterization in finite dimensions}

A more explicit characterization of the Clarke subdifferential is possible in finite-dimensional spaces. The basis is the following theorem, which only holds in $\R^N$; a proof can be found in, e.g., \cite[Theorem 23.2]{DiBenedetto} or \cite[Theorem 3.1]{Heinonen}.
\begin{theorem}[Rademacher]\label{thm:rademacher}\index{theorem!Rademacher}
    Let $U\subset \R^N$ be open and $F:U\to\R$ be Lipschitz continuous. Then $F$ is Fréchet differentiable at almost every $x\in U$.
\end{theorem}
This result allows replacing the $\limsup$ in the definition of the Clarke subdifferential (now considered as a subset of $\R^N$, i.e., identifying the dual of $\R^N$ with $\R^N$ itself) with a proper limit.
\begin{theorem}\label{thm:clarke:gradient}
    Let $F:\R^N\to\R$ be locally Lipschitz continuous near $x\in\R^N$. Then $F$ is Fréchet differentiable on $\R^N\setminus E_F$ for a set $E_F\subset \R^N$ of Lebesgue measure $0$ and
    \begin{equation}\label{eq:clarke:gradient}
        \partial_C F(x) = \conv \setof{\lim_{n\to\infty} \nabla F(x_n)}{x_n\to x,\ x_n\notin E_F},
    \end{equation}
    where $\conv A$ denotes the convex hull of $A\subset \R^N$.
\end{theorem}
\begin{proof}
    We first note that the Rademacher theorem ensures that such a set $E_F$ exists and has Lebesgue measure $0$. Hence there indeed exist sequences $\{x_n\}_{n\in\N}\subset\R^N\setminus E_F$ with $x_n\to x$.
    Furthermore, the local Lipschitz continuity of $F$ yields that for any $x_n$ in the Lipschitz neighborhood of $x$ and any $h\in \R^N$, we have that
    \begin{equation*}
        |\dual{\nabla F(x_n),h}| = \left|\lim_{t\downto 0}\frac{F(x_n+th)-F(x_n)}{t}\right| \leq L\norm{h}
    \end{equation*}
    and hence that $\norm{\nabla F(x_n)} \leq L$ for all $n\in\N$ large enough. This implies that $\{\nabla F(x_n)\}_{n\in\N}\subset \R^N$ is bounded and thus contains a convergent subsequence. The set on the right-hand side of \eqref{eq:clarke:gradient} is therefore nonempty.

    Let now $\{x_n\}_{n\in\N}\subset \R^N\setminus E_F$ be an arbitrary sequence with $x_n\to x$ and $\{\nabla F(x_n)\}_{n\in\N}\to x^*$ for some $x^*\in\R^N$.
    Since $F$ is differentiable at every $x_n\notin E_F$ by definition, \cref{lem:clarke:gateaux} yields that $\nabla F(x_n) \in \partial_C F(x_n)$, and hence $x^*\in \partial_C F(x)$ by \cref{lem:clarke:closed}.
    The convexity of $\partial_C F(x)$ from \cref{lem:clarke:properties} now implies that any convex combination of such limits $x^*$ is contained in $\partial_C F(x)$, which shows the inclusion ``$\supset$'' in \eqref{eq:clarke:gradient}.

    For the other inclusion, we first show for all $h\in \R^N$ and $\eps>0$ that
    \begin{equation}\label{eq:clarke:gradient1a}
        F^\circ(x;h) - \eps \leq \limsup_{E_F\not\ni y\to x}\, \dual{\nabla F(y),h}=:M(h).
    \end{equation}
    Indeed, by definition of $M(h)$ and of the $\limsup$, for every $\eps>0$ there exists a $\delta>0$ such that
    \begin{equation*}
        \dual{\nabla F(y),h} \leq  M(h)+\eps \qquad\text{for all }y\in \OB(x,\delta)\setminus E_F.
    \end{equation*}
    Here, $\delta>0$ can be chosen sufficiently small for $F$ to be Lipschitz continuous on $\OB(x,\delta)$. In particular, $E_F\cap \OB(x,\delta)$ is a set of zero measure. Hence, $F$ is differentiable at $y+th$ for almost all $y\in \OB(x,\frac\delta2)$ and almost all $t\in (0,\frac\delta{2\norm{h}})$ by Fubini's theorem. The classical mean value theorem therefore yields for all such $y$ and $t$ that
    \begin{equation}\label{eq:clarke:gradient1}
        F(y+th)-F(y) = \int_0^t \dual{\nabla F(y+sh),h}\,ds \leq t(M(h)+\eps)
    \end{equation}
    since $y+sh\in \OB(x,\delta)$ for all $s\in (0,t)$ by the choice of $t$.
    The continuity of $F$ implies that the full inequality \eqref{eq:clarke:gradient1} even holds for \emph{all} $y\in \OB(x,\frac\delta2)$ and \emph{all} $t\in (0,\frac\delta{2\norm{h}})$. Dividing by $t>0$ and taking the $\limsup$ over all $y\to x$ and $t\downto 0$ now yields \eqref{eq:clarke:gradient1a}.

    Since $\eps>0$ was arbitrary, this implies that $F^\circ(x;h)\leq M(h)$ for all $h\in \R^N$ and hence that
    \begin{equation*}
        \partial_C F(x) \subset \setof{x^*\in\R^N}{\dual{x^*,h} \leq M(h)\quad\text{for all }h\in\R^N} =: B.
    \end{equation*}
    We are thus finished if we can show that $B$ is equal to the set on the right-hand side of \eqref{eq:clarke:gradient}, which we denote by $\conv A$.
    For this, we once again appeal to \cref{lem:clarke:support2}.
    First, we note that the definition of the convex hull implies for all $h\in\R^N$ that
    \begin{equation*}
        \sup_{x^*\in \conv A} \dual{x^*,h}
        = \sup_{\substack{x_i^*\in A\\\sum_i t_i = 1,t_i\geq 0}} \sum_i t_i\dual{x_i^*,h} = \sup_{\sum_i t_i = 1,t_i\geq 0} \sum_i t_i \sup_{x_i^*\in A}\dual{x_i^*,h} = \sup_{x^*\in A}\dual{x^*,h}
    \end{equation*}
    since the sum is maximal if and only if each summand is maximal. Next we have that
    \begin{equation*}
        M(h) = \limsup_{E_F\not\ni y\to x} \,\dual{\nabla F(y),h} = \sup_{E_F \not\ni x_n\to x}\dual{\lim\nolimits_{n\to\infty} \nabla F(x_n),h} = \sup_{x^*\in A}\dual{x^*,h}.
    \end{equation*}
    Finally, one can show as in \cref{lem:clarke:dir} that the mapping $h\mapsto M(h)$ is positively homogeneous, subadditive, and lower semicontinuous. From \cref{lem:clarke:support1}, we thus have that
    \begin{equation*}
        \sup_{x^*\in B} \dual{x^*,h} = M(h) = \sup_{x^*\in A}\dual{x^*,h} = \sup_{x^*\in \conv A} \dual{x^*,h}.
    \end{equation*}
    Since both sets are clearly convex and closed as well as nonempty (which we've already argued for $\conv A$ and which follows from \eqref{eq:clarke:gradient1a} for $B$), \cref{eq:clarke:support2} yields $B=\conv A$ and thus the claim.
\end{proof}

\begin{remark}\label{rem:clarke:extended-real}
    It is possible to extend the Clarke subdifferential defined here to extended-real valued functions using an equivalent, more geometrical, construction involving generalized normal cones to epigraphs; see \cite[Definition 2.4.10]{Clarke:1990a}. We will follow this approach when studying the more general subdifferentials for set-valued functionals in \cref{chap:cones,chap:graphical}.
\end{remark}

\chapter{Semismooth Newton methods}\label{chap:semismooth}

The proximal point and splitting methods in \cref{chap:proximal} are generalizations of gradient methods and in general have at most linear convergence. In this chapter, we will therefore consider second-order methods, specifically a generalization of Newton's method which admits (locally) superlinear convergence.

\section{Convergence of generalized Newton methods}

As a motivation, we first consider the most general form of a Newton-type method. Let $X$ and $Y$ be normed vector spaces and $F:X\to Y$ be given and suppose we are looking for an $\bar x \in X$ with $F(\bar x)=0$. A Newton-type method to find such an $\bar x$ then consists of repeating the following steps:
\begin{enumerate}[label=\arabic*.]
    \item choose an invertible $M_k \defeq M(x^k)\in \linear(X;Y)$;
    \item solve the \term[step, Newton]{Newton step} $M_k s^k = -F(x^k)$;
    \item update $x^{k+1} = x^k + s^k$.
\end{enumerate}
We can now ask under which conditions this method converges to $\bar x$, and in particular, when the convergence is \term[convergence!superlinear]{superlinear}, i.e.,
\begin{equation}\label{eq:newton:superlinear}
    \lim_{k\to\infty} \frac{\norm{x^{k+1}-\bar x}_X}{\norm{x^{k}-\bar x}_X} = 0.
\end{equation}
(Recall the discussion in the beginning of \cref{chap:testing}.)
For this purpose, we set $e^k \defeq x^k-\bar x$ and use the Newton step together with the fact that $F(\bar x)=0$ to obtain that
\begin{equation}
    \label{eq:newton:est1}
    \begin{aligned}[t]
        \norm{x^{k+1}-\bar x}_X &=  \norm{x^k - M(x^k)^{-1}F(x^{k}) - \bar x}_X\\
        &= \norm{M(x^k)^{-1}\left[F(x^k) - F(\bar x) - M(x^k)(x^k-\bar x)\right]}_X\\
        &= \norm{M(\bar x+e^k)^{-1}\left[F(\bar x+e^k) - F(\bar x) - M(\bar x+e^k) e^k\right]}_X\\
        &\leq \norm{M(\bar x+e^k)^{-1}}_{\linear(Y;X)}\norm{F(\bar x+e^k) - F(\bar x) - M(\bar x+e^k)e^k}_Y.
    \end{aligned}
\end{equation}
Hence, \eqref{eq:newton:superlinear} holds under
\begin{enumerate}
    \item a \term[condition!regularity]{regularity condition}: there exists a $C>0$ with
        \begin{equation*}
            \norm{M(x^k)^{-1}}_{\linear(Y;X)} \leq C \qquad\text{for all }k\in\N;
        \end{equation*}
    \item an \term[condition!approximation]{approximation condition}:
        \begin{equation*}\label{eq:newton:approx}
            \lim_{k\to\infty} \frac{\norm{F(\bar x+ e^k) - F(\bar x) - M(\bar x +e^k) e^k}_Y}{\norm{e^k}_X}= 0.
        \end{equation*}
\end{enumerate}

This motivates the following definition: We call $F:X\to Y$ \term[mapping!differentiable!Newton]{Newton differentiable} at $x\in X$ with \term[derivative!Newton]{Newton derivative} $D_NF(x)$ if there exists a neighborhood $U\subset X$ of $x$ and a mapping $D_N F: U \to \linear(X;Y)$ such that
\begin{equation}\label{eq:newton:semismooth}
    \lim_{\norm{h}_X\to 0} \frac{\norm{F(x+ h) - F(x) - D_N F(x +h)h }_Y}{\norm{h}_X}= 0.
\end{equation}
Note the differences to the Fréchet derivative: First, the Newton derivative is evaluated in $x+h$ instead of $x$. More importantly, we have not required \emph{any} connection between $D_N F$ and $F$, while the only possible candidate for the Fréchet derivative was the Gateaux derivative (which itself was linked to $F$ via the directional derivative). A function thus can only be Newton differentiable (or not) with respect to a concrete choice of $D_N F$. In particular, Newton derivatives are not unique.

If $F$ is Newton differentiable with Newton derivative $D_N F$, we can set $M(x^k) = D_N F(x^k)$ and obtain the \term[method!semismooth Newton]{semismooth Newton method}
\begin{algeqbox}
    \begin{equation}\label{eq:SSN}
        x^{k+1} \defeq x^k - D_N F(x^k)^{-1} F(x^k).
    \end{equation}
\end{algeqbox}
Its local superlinear convergence follows directly from the construction.
\begin{theorem}\label{thm:newton:superlinear}
    Let $X,Y$ be normed vector spaces and let $F:X\to Y$ be Newton differentiable near $\bar x\in X$ with $F(\bar x)=0$ with Newton derivative $D_N F(\bar x)$. Assume further that there exist $\delta>0$ and $C>0$ with $\norm{D_NF(x)^{-1}}_{\linear(Y;X)}\leq C$ for all $x\in \OB(\bar x,\delta)$. Then the semismooth Newton method \eqref{eq:SSN} converges superlinearly to $\bar x$ for all $x^0$ sufficiently close to $\bar x$.
\end{theorem}
\begin{proof}
    The proof is virtually identical to that for the classical Newton method. We have already shown that for any $x^0\in \OB(\bar x,\delta)$,
    \begin{equation}\label{eq:newton:superlinear1}
        \norm{e^1}_X \leq C \norm{F(\bar x+e^0) - F(\bar x) - D_NF(\bar x+e^0)e^0}_Y.
    \end{equation}
    Let now $\eps\in(0,1)$ be arbitrary. The Newton differentiability of $F$ then implies that there exists a $\rho>0$ such that
    \begin{equation}
        \label{eq:newton:superlinear2}
        \norm{F(\bar x+h) - F(\bar x) - D_NF(\bar x+h)h}_Y \leq \frac\eps{C}\norm{h}_X \qquad\text{for all }\norm{h}_X\leq \rho.
    \end{equation}
    Hence, if we choose $x^0$ such that $\norm{\bar x- x^0}_X\leq \min\{\delta,\rho\}$, the estimate \eqref{eq:newton:superlinear1} implies that  $\norm{\bar x- x^1}_X\leq \eps\norm{\bar x- x^0}_X$. By induction, we obtain from this that $\norm{\bar x-x^k}_X\leq \eps^k \norm{\bar x- x^0}_X \to 0$. Since $\eps\in(0,1)$ was arbitrary, we can take in each step $k$ a different $\eps_k\to 0$ to obtain that $\norm{x^{k+1}-\bar x}_X \leq \eps_k\norm{x^k-\bar x}_X$ and hence that the convergence is superlinear.
\end{proof}

Sometimes, the Newton derivatives $D_N F(x)$ are poorly conditioned, or the region of convergence impractically small. In that case, it may help to \term[dampening]{dampen} the method to
\begin{algeqbox*}
    \begin{equation*}
        x^{k+1} \defeq x^k - [D_N F(x^k) + \theta \Id]^{-1} F(x^k),
    \end{equation*}
\end{algeqbox*}
for some $\theta > 0$. As shown in the next theorem, this method still converges, but only linearly.
For this scheme, we would take $M(x) = D_N F(x) + \theta \Id$ in the theorem.
As we will learn in \cref{chap:sparse}, it is also possible to modify $D_N F(x)$ only on a subspace.

\begin{theorem}\label{thm:newton:dampened-linear}
    Let $X,Y$ be normed vector spaces and let $F:X\to Y$ be Newton differentiable near $\bar x\in X$ with $F(\bar x)=0$ with Newton derivative $D_N F(\bar x)$.
    Also assume to be given $M(x) \in \linear(X; Y)$ that satisfy $\norm{M(x)-D_N F(x)}_{\linear(X; Y)} \le \theta$ and $\norm{M(x)^{-1}}_{\linear(Y;X)}\leq C$ for all $x\in \OB(\bar x,\delta)$ for some $\theta,\delta>0$ and $0<C<\inv\theta$. Then $x^{k+1} \defeq x^k - M(x^k)^{-1} F(x^k)$ converge linearly to $\bar x$ for all $x^0$ sufficiently close to $\bar x$.
\end{theorem}

\begin{proof}
    Following \eqref{eq:newton:est1}, we have
    \begin{equation}\label{eq:newton:linear1}
        \norm{e^1}_X \leq C \norm{F(\bar x+e^0) - F(\bar x) - M(\bar x+e^0)e^0}_Y.
    \end{equation}
    Let $\epsilon>0$. Using the Newton differentiability of $F$, following \eqref{eq:newton:superlinear2}, we deduce the existence of $\rho>0$ such that whenever $\norm{h}_X \le \rho$, we have
    \begin{equation*}
        \begin{aligned}[t]
        \norm{F(\bar x+h) - F(\bar x) - M(\opt x + h)h}_Y
        &
        \le
        \norm{F(\bar x+h) - F(\bar x) - D_N F(\opt x + h)h}_Y
        \\
        \MoveEqLeft[-1]
        + \norm{[D_N F(\opt x + h) - M(\opt x + h)]h}_Y
        \\
        &
        \le
        \left(\frac\eps{C} + \theta\right)\norm{h}_X.
        \end{aligned}
    \end{equation*}
    Hence, if we choose $x^0$ such that $\norm{\bar x- x^0}_X\leq \min\{\delta,\rho\}$, the estimate \eqref{eq:newton:linear1} implies that  $\norm{\bar x- x^1}_X\leq (C\theta+\eps)\norm{\bar x- x^0}_X$. Since $0 < C\theta < 1$, taking $\eps>0$ small enough, we have $\beta \defeq C\theta+\eps \in (0, 1)$. By induction, we obtain from this that $\norm{\bar x-x^k}_X\leq \beta^k \norm{\bar x- x^0}_X \to 0$. This shows the linear convergence.
\end{proof}

\section{Newton derivatives}

The remainder of this chapter is dedicated to the construction of Newton derivatives that satisfy the approximation condition (although it should be pointed out that the verification of the regularity condition is usually the much more involved step in practice, which is usually very specific to the concrete problem).
We begin with the obvious connection with the Fréchet derivative.
\begin{theorem}\label{thm:newton:frechet}
    Let $X,Y$ be normed vector spaces. If $F:X\to Y$ is continuously differentiable at $x\in X$, then $F$ is also Newton differentiable at $x$ with Newton derivative $D_N F(x) = F'(x)$.
\end{theorem}
\begin{proof}
    We have for arbitrary $h\in X$ that
    \begin{equation*}
        \begin{aligned}
            \norm{F(x+h)-F(x)-F'(x+h)h}_Y &\leq \norm{F(x+h)-F(x)-F'(x)h}_Y\\
            \MoveEqLeft[-1]+ \norm{F'(x)-F'(x+h)}_{\linear(X;Y)}\norm{h}_X,
        \end{aligned}
    \end{equation*}
    where the first summand is $o(\norm{h}_X)$ by definition of the Fréchet derivative and the second by the continuity of $F'$.
\end{proof}
Calculus rules can be shown similarly to those for Fréchet derivatives. For the sum rule this is immediate; here we prove a chain rule by way of example.
\begin{theorem}\label{thm:newton:chain}
    Let $X$, $Y$, and $Z$ be normed vector spaces, and let $F:X\to Y$ be Newton differentiable at $x\in X$ with Newton derivative $D_N F(x)$ and $G:Y\to Z$ be Newton differentiable at $y\defeq F(x)\in Y$ with Newton derivative $D_N G(y)$.
    If $D_N F$ and $D_N G$ are uniformly bounded in a neighborhood of $x$ and $y$, respectively, then $G\circ F$ is also Newton differentiable at $x$ with Newton derivative
    \begin{equation*}
        D_N (G\circ F)(x) = D_N G(F(x))\circ D_N F(x).
    \end{equation*}
\end{theorem}
\begin{proof}
    We proceed as in the proof of \cref{thm:frechet_chain}. For $h\in X$ and  $g \defeq F(x+h)-F(x)$ we have that
    \begin{equation*}
        (G\circ F)(x+h) - (G\circ F)(x) = G(y+g)  - G(y).
    \end{equation*}
    The Newton differentiability of $G$ then implies that
    \begin{equation*}
        \norm{(G\circ F)(x+h ) - (G\circ F)(x) - D_N G(y+g)g}_Z \leq r_1(\norm{g}_Y)
    \end{equation*}
    with $r_1(t)/t \to 0$ for $t\to 0$.
    The Newton differentiability of $F$ further implies that
    \begin{equation*}
        \norm{g - D_N F(x+h)h}_Y \leq r_2(\norm{h}_X)
    \end{equation*}
    with $r_2(t)/t \to 0$ for $t\to 0$. In particular,
    \begin{equation*}
        \norm{g}_Y \leq \norm{D_N F(x+h)}_{\linear(X;Y)}\norm{h}_Y + r_2(\norm{h}_X).
    \end{equation*}
    The uniform boundedness of $D_N F$ now implies that $\norm{g}_Y\to 0$ for $\norm{h}_X\to 0$. Hence, using that $y+g=F(x+h)$, we obtain
    \begin{multline*}
        \norm{(G\circ F)(x+h) - (G\circ F)(x) -  D_N G(F(x+h))D_NF(x+h)h}_Z \\
        \begin{aligned}[t]
            &\leq \norm{G(y+g)-G(y)-D_NG(y+g)g}_Z\\
            \MoveEqLeft[-1]+ \norm{D_N G(y+g)\left[g-D_NF(x+h)h\right]}_Z\\
            &\leq r_1(\norm{g}_Y) +  \norm{D_N G(y+g)}_{\linear(Y;Z)} r_2(\norm{h}_X),
        \end{aligned}
    \end{multline*}
    and the claim thus follows from the uniform boundedness of $D_N G$.
\end{proof}

Finally, it follows directly from the definition of the product norm and Newton differentiability that Newton derivatives of vector-valued functions can be computed componentwise.
\begin{theorem}\label{thm:newton:vector}
    Let $X,Y_i$ be normed vector spaces and let $F_i:X\to Y_i$ be Newton differentiable with Newton derivative $D_N F_i$ for $1\leq i\leq m$. Then
    \begin{equation*}
        F:X\to (Y_1\times\cdots\times Y_m), \qquad x\mapsto (F_1(x),\dots,F_m(x))^T,
    \end{equation*}
    is also Newton differentiable with Newton derivative
    \begin{equation*}
        D_N F(x) = (D_N F_1(x),\dots, D_N F_m(x))^T.
    \end{equation*}
\end{theorem}

\bigskip

Since the definition of a Newton derivative is not constructive, allowing different choices, the question remains how to obtain a candidate for which the approximation condition in the definition can be verified. For two classes of functions, such an explicit construction is known.

\subsection*{Locally Lipschitz continuous functions on \texorpdfstring{$\scriptstyle\R^N$}{ℝᴺ}}

If $F:\R^N\to \R$ is locally Lipschitz continuous, candidates can be taken from the Clarke subdifferential, which has an explicit characterization by \cref{thm:clarke:gradient}. Under some additional assumptions, each candidate is indeed a Newton derivative.

A function $F:\R^N\to\R$ is called \term[mapping!differentiable!piecewise]{piecewise (continuously) differentiable} or \term[function!PC$^1$]{PC$^1$ function}, if
\begin{enumerate}
    \item $F$ is continuous on $\R^N$;
    \item for all $x\in\R^N$ there exists an open neighborhood $U_x\subset \R^N$ of $x$ and a finite set $\{F_i:U_x\to \R\}_{i\in I_x}$ of continuously differentiable functions with
        \begin{equation*}
            F(\tilde x) \in \{F_i(\tilde x)\}_{i\in I_x}\qquad\text{for all }\tilde x\in U_x.
        \end{equation*}
\end{enumerate}
In this case, we call $F$ a \term[selection!continuous]{continuous selection} of the $F_i$ in $U_x$. The set
\begin{equation*}
    I_a(x) \defeq \setof{i\in I_x}{F(x) = F_i(x)}
\end{equation*}
is called the \term[set!index!active]{active index set} at $x$. Since the $F_i$ are continuous, we have that $F(\tilde x) \neq F_j(\tilde x)$ for all $j\notin I_a(x)$ and $\tilde x$ sufficiently close to $x$. Hence, indices that are only active on sets of zero measure do not have to be considered in the following. We thus define the \term[set!index!essentially active]{essentially active index set}
\begin{equation*}
    I_e(x) \defeq \setof{i\in I_x}{x\in \closure\left(\interior\setof{\tilde x\in U_x}{F(\tilde x) = F_i(\tilde x)}\right)}\subset I_a(x).
\end{equation*}
An example of an active but not essentially active index set is the following.
\begin{example}
    Consider the function
    \begin{equation*}
        f:\R\to\R,\qquad t\mapsto \max\{0,t,t/2\},
    \end{equation*}
    i.e., $f_1(t)=0$, $f_2(t)=t$, and $f_3(t)=t/2$. Then $I_a(0)=\{1,2,3\}$ but $I_e(0)=\{1,2\}$, since $f_3$ is active only at $t=0$ and hence $\interior\setof{t\in\R}{f(t)=f_3(t)}=\emptyset=\closure\emptyset$.
\end{example}

Since any $C^1$ function $F_i:U_x\to\R$ is Lipschitz continuous with Lipschitz constant $L_i\defeq \sup_{\tilde x\in U_x}|\grad F(\tilde x)|$ by \cref{lem:variation:c1-lipschitz}, PC$^1$ functions are always locally Lipschitz continuous; see \cite[Corollary 4.1.1]{Scholtes:2012}.
\begin{theorem}
    Let $F:\R^N\to\R$ be piecewise differentiable. Then $F$ is locally Lipschitz continuous on $\R^N$ with local constant $L(x)=\max_{i\in I_a(x)} L_i$.
\end{theorem}
This yields the following explicit characterization of the Clarke subdifferential of a PC$^1$ function.
\begin{theorem}\label{thm:newton:clarke}
    Let $F:\R^N\to\R$ be piecewise differentiable and $x\in\R^N$. Then
    \begin{equation*}
        \partial_C F(x) = \conv\setof{\nabla F_i(x)}{i\in I_e(x)}.
    \end{equation*}
\end{theorem}
\begin{proof}
    Let $x\in\R^N$ be arbitrary. By \cref{thm:clarke:gradient} it suffices to show that
    \begin{equation*}
        \setof{\lim_{n\to\infty} \nabla F(x_n)}{x_n\to x,\ x_n\notin E_F} = \setof{\nabla F_i(x)}{i\in I_e(x)},
    \end{equation*}
    where $E_F$ is the set of Lebesgue measure $0$ where $F$ is not differentiable from Rademacher's theorem.
    For this, let $\{x_n\}_{n\in\N}\subset\R^N$ be a sequence with $x_n\to x$, $F$ is differentiable at $x_n$ for all $n\in \N$, and $\nabla F(x_n)\to x^* \in \R^N$. Since $F$ is differentiable at $x_n$, it must hold that $F(\tilde x) = F_{i_n}(\tilde x)$ for some $i_n\in I_a(x)$ and all $\tilde x$ sufficiently close to $x_n$, which implies that $\nabla F(x_n) = \nabla F_{i_n} (x_n)$.
    For sufficiently large $n\in \N$, we can further assume that $i_n\in I_e(x)$ (if necessary, by adding $x_n$ with $i_n\notin I_e(x)$ to $E_F$, which does not increase its Lebesgue measure). If we now consider subsequences $\{x_{n_k}\}_{k\in\N}$ with constant index $i_{n_k}=: i\in I_e(x)$ (which exist since $I_e(x)$ is finite), we obtain using the continuity of $\nabla F_i$ that
    \begin{equation*}
        x^*=\lim_{k\to\infty} \nabla F(x_{n_k}) =\lim_{k\to\infty} \nabla F_i(x_{n_k})  \in \setof{\nabla F_i(x)}{i\in I_e(x)}.
    \end{equation*}

    Conversely, for every $\nabla F_i(x)$ with $i\in I_e(x)$ there exists by definition of the essentially active indices a sequence $\{x_n\}_{n\in\N}$ with $x_n\to x$ and $F=F_i$ in a sufficiently small neighborhood of each $x_n$ for $n$ large enough. The continuous differentiability of the $F_i$ thus implies that $\nabla F(x_n) = \nabla F_i(x_n)$ for all $n\in\N$ large enough and hence that
    \begin{equation*}
        \nabla F_i(x) = \lim_{n\to \infty} \nabla F_i(x_n) = \lim_{n\to \infty} \nabla F(x_n).
        \qedhere
    \end{equation*}
\end{proof}
From this, we obtain the Newton differentiability of PC$^1$ functions.
\begin{theorem}\label{thm:newton:clarke_ndiff}
    Let $F:\R^N\to\R$ be piecewise differentiable. Then $F$ is Newton differentiable for all $x\in\R^N$, and every $D_NF(x)\in \partial_C F(x)$ is a Newton derivative.
\end{theorem}
\begin{proof}
    Let $x\in \R^N$ be arbitrary and $h\in \R^N$ with $x+h\in U_x$. By \cref{thm:newton:clarke}, every $D_N F(x+h)\in \partial_C F(x+h)$ is of the form
    \begin{equation*}
        D_N F(x+h) = \sum_{i\in I_e(x+h)} \lambda_i \nabla F_i(x+h)\qquad\text{for }\sum_{i\in I_e(x+h)} \lambda_i = 1,\, \lambda_i\geq 0.
    \end{equation*}
    Since $F$ is continuous, we have for all $h\in\R^N$ sufficiently small that $I_e(x+h)\subset I_a(x+h)\subset I_a(x)$, where the second inclusion follows from the fact that by continuity, $F(x)\neq F_i(x)$ implies that $F(x+h) \neq F_i(x+h)$. Hence, $F(x+h) = F_i(x+h)$ and $F(x) = F_i(x)$ for all $i\in I_e(x+h)$.
    \Cref{thm:newton:frechet} then yields that
    \begin{equation*}
        |F(x+h)-F(x)-D_N F(x+h)h| \leq \sum_{i\in I_e(x+h)} \lambda_i |F_i(x+h)-F_i(x) - \nabla F_i(x+h)h| = o(\norm{h}),
    \end{equation*}
    since all $F_i$ are continuously differentiable by assumption.
\end{proof}

A natural application of the above are proximal point mappings of convex and lower semicontinuous functionals.
\begin{example}\label{ex:newton:rn}~
    \begin{enumerate}
        \item\label{ex:newton:rn:box}
        We first consider the proximal mapping for the indicator function $\delta_A:\R^N\to \Rbar $ of the set $A\defeq\setof{x\in \R^N}{x_i\in [a,b]}$ for some $a<b\in \R$. Analogously to \cref{ex:proximal:reell}\,\ref{ex:proximal:reell:iii}, the corresponding proximal mapping is the componentwise projection
            \begin{equation*}
                [\proj_A(x)]_i = \proj_{[a,b]}x_i =
                \begin{cases}
                    a &\text{if } x_i<a,\\
                    x_i &\text{if } x_i\in [a,b],\\
                    b &\text{if } x_i>b,
                \end{cases}
            \end{equation*}
            which is clearly piecewise differentiable.
            \Cref{thm:newton:clarke} thus yields (also componentwise) that
            \begin{equation*}
                \partial_C [\proj_A(x)]_i =
                \begin{cases}
                    \{1\} &\text{if } x_i\in (a,b),\\
                    \{0\} &\text{if } x_i\notin [a,b],\\
                    [0,1] &\text{if } x_i\in\{a,b\}.
                \end{cases}
            \end{equation*}
            By \cref{thm:newton:clarke_ndiff,thm:newton:vector}, a possible Newton derivative is therefore given by
            \begin{equation*}
                [D_N \proj_{A}(x) h]_i = [\1_{[a,b]}(x)\odot h]_i \defeq
                \begin{cases}
                    h_i & \text{if }x_i\in [a,b], \\
                    0 & \text{if }x_i\notin[a,b],
                \end{cases}
            \end{equation*}
            where the choice of which case to include $x_i\in\{a,b\}$ in is arbitrary. (The componentwise product $[x\odot y]_i\defeq x_iy_i$ on $\R^N$ is also known as the \term[product!Hadamard]{Hadamard product}.)

        \item\label{ex:newton:rn:l1}
            Consider now the proximal mapping for $G:\R^N\to \R$, $G(x) \defeq \norm{x}_1$, whose proximal mapping for arbitrary $\gamma>0$ is given by \cref{ex:proximal:rn}\,\ref{ex:proximal:rn:ii} componentwise as
            \begin{equation*}
                [\prox_{\gamma G}(x)]_i =
                \begin{cases}
                    x_i-\gamma & \text{if }x_i>\gamma,\\
                    0 & \text{if }x_i \in[-\gamma,\gamma],\\
                    x_i+\gamma & \text{if }x_i < -\gamma.
                \end{cases}
            \end{equation*}
            Again, this  is clearly piecewise differentiable, and
            \Cref{thm:newton:clarke} thus yields (also componentwise) that
            \begin{equation*}
                \partial_C [(\prox_{\gamma G})(x)]_i =
                \begin{cases}
                    \{1\} & \text{if }|x_i|>\gamma,\\
                    \{0\} & \text{if }|x_i|<\gamma,\\
                    [0,1] & \text{if }|x_i|=\gamma.
                \end{cases}
            \end{equation*}
            By \cref{thm:newton:clarke_ndiff,thm:newton:vector}, a possible Newton derivative is therefore given by
            \begin{equation*}
                [D_N \prox_{\gamma G}(x) h]_i = [\1_{\R\setminus(-\gamma,\gamma)}(x)\odot h]_i \defeq
                \begin{cases}
                    h_i & \text{if }|x_i|\geq \gamma, \\
                    0 & \text{if }|x_i|<\gamma,
                \end{cases}
            \end{equation*}
            where again we could have taken the value $t h_i$ for any $t\in [0,1]$ for $|x_i|=\gamma$.
    \end{enumerate}
\end{example}

\subsection*{Superposition operators on \texorpdfstring{$\scriptstyle L^p(\Omega)$}{Lᵖ(Ω)}}

Rademacher's theorem does not hold in infinite-dimensional function spaces, and hence the Clarke subdifferential no longer yields an algorithmically useful candidate for a Newton derivative in general. One exception is the class of superposition operators defined by scalar Newton differentiable functions, for which the Newton derivative can be evaluated pointwise as well.

We thus consider as in \cref{sec:superposition} for an open and bounded domain $\Omega\subset \R^N$, a Carathéodory function $f:\Omega\times \R\to\R$ (i.e., $(x,z)\mapsto f(x,z)$ is measurable in $x$ and continuous in $z$), and $1\leq p,q\leq \infty$ the corresponding superposition operator
\begin{equation*}
    F:L^p(\Omega)\to L^q(\Omega),\qquad [F(u)](x) = f(x,u(x))\quad\text{for almost every }x\in\Omega.
\end{equation*}
The goal is now to similarly obtain a Newton derivative $D_N F$ for $F$ as a superposition operator defined by the Newton derivative $D_N f(x,z)$ of $z\mapsto f(x,z)$. Here, the assumption that $D_N f$ is also a Carathéodory function is too restrictive, since we want to allow discontinuous derivatives as well (see \cref{ex:newton:rn}). Luckily, for our purpose, a weaker property is sufficient: A function is called \term[function!Baire--Carathéodory]{Baire--Carathéodory function} if it can be written as a pointwise limit of Carathéodory functions, i.e., if
\begin{equation*}
    f(x,z) = \lim_{n\to \infty} f_n(x,z) \qquad\text{for almost every }x\in \Omega\text{ and all }z\in \R,
\end{equation*}
where $f_n$ is a Carathéodory function for all $n\in\N$; see \cite[Lemma 1.4]{Appell:1990}.

Under certain growth conditions on $f$ and $D_N f$,\footnote{which can be significantly relaxed; see \cite[Proposition \textsc{a}.1]{Schiela:2008a}} we can transfer the Newton differentiability of $f$ to $F$, but we again have to take a \term[discrepancy, two-norm]{two-norm discrepancy} into account.
\begin{theorem}\label{thm:newton:super}
    Let $f:\Omega\times\R\to\R$ be a Carathéodory function. Furthermore, assume that
    \begin{enumerate}
        \item $z\mapsto f(x,z)$ is uniformly Lipschitz continuous for almost every $x\in \Omega$ and $x\mapsto f(x,0)$ is bounded;
        \item $z\mapsto f(x,z)$ is Newton differentiable with Newton derivative $z\mapsto D_N f(x,z)$ for almost every $x\in \Omega$;
        \item\label{it:newton:super:baire} $D_N f$ is a Baire--Carathéodory function and uniformly bounded.
    \end{enumerate}
    Then for any $1\leq q<p< \infty$, the corresponding superposition operator $F:L^p(\Omega)\to L^q(\Omega)$ is Newton differentiable with Newton derivative
    \begin{equation*}
        D_N F:L^p(\Omega)\to \linear(L^p(\Omega);L^q(\Omega)), \qquad
        [D_N F(u)h](x) = D_N f(x,u(x))h(x)
    \end{equation*}
    for almost every $x\in \Omega$ and all $h\in L^p(\Omega)$.
\end{theorem}
\begin{proof}
    First, the uniform Lipschitz continuity together with the reverse triangle inequality yields that
    \begin{equation*}
        |f(x,z)| \leq |f(x,0)|+L|z| \leq C +L|z|^{q/q}\quad\text{for almost every }x\in\Omega\text{ and all } z\in\R,
    \end{equation*}
    and hence the growth condition \eqref{eq:superpos:growth} is satisfied for all $1\leq q\leq \infty$. Due to the continuous embedding $L^p(\Omega)\hookrightarrow L^q(\Omega)$ for all $1\leq q\leq p\leq \infty$, the superposition operator $F:L^p(\Omega)\to L^q(\Omega)$ is therefore well-defined and continuous by \cref{thm:superpos:continuous}.

    For any measurable $u:\Omega\to\R$, we have that $x\mapsto D_N f(x,u(x))$ is by assumption \ref{it:newton:super:baire} the pointwise limit of measurable functions and hence itself measurable. Furthermore, its uniform boundedness in particular implies the growth condition \eqref{eq:superpos:growth} for $p'\defeq p$ and $q'\defeq p-q>0$. As in the proof of \cref{thm:superpos:differentiable}, we deduce that the corresponding superposition operator $D_N F:L^p(\Omega)\to L^s(\Omega)$ is well-defined and continuous for $s\defeq\frac{pq}{p-q}$, and that for any $u\in L^p(\Omega)$, the mapping $h\mapsto D_NF(u)\cdot h$ defines a bounded linear operator $D_NF(u):L^p(\Omega)\to L^q(\Omega)$.
    (This time, we do not distinguish in notation between the linear operator and the function defining this operator by pointwise multiplication.)

    To show that $D_N F(u)$ is a Newton derivative for $F$ in $u\in L^p(\Omega)$, we consider the pointwise residual
    \begin{equation*}
        r:\Omega\times\R\to\R,\qquad r(x,z) \defeq
        \begin{cases}
            \frac{|f(x,z)-f(x,u(x))-D_N f(x,z)(z-u(x))|}{|z-u(x)|} & \text{if }z\neq u(x),\\
            0 & \text{if }z=u(x).
        \end{cases}
    \end{equation*}
    Since $f$ is a Carathéodory function and $D_N f$ is a Baire--Carathéodory function, the function $x\mapsto r(x,\tilde u(x))=:[R(\tilde u)](x)$ is measurable for any measurable $\tilde u:\Omega\to\R$ (since sums, products, and quotients of measurable functions are again measurable).
    Furthermore, for $\tilde u\in L^p(\Omega)$, the uniform Lipschitz continuity of $f$ and the uniform boundedness of $D_N f$ imply that for almost every $x\in \Omega$ with $\tilde u(x) \neq u(x)$,
    \begin{equation}\label{eq:newton:superpos1}
        |[R(\tilde u)](x)| = \frac{|f(x,\tilde u(x))-f(x,u(x))-D_N f(x,\tilde u(x))(\tilde u(x)-u(x))|}{|\tilde u(x)-u(x)|}
        \leq L + C
    \end{equation}
    and thus that $R(\tilde u)\in L^\infty(\Omega)$. Hence, the superposition operator $R:L^p(\Omega)\to L^{s}(\Omega)$ is well-defined.

    Let now $\{u_n\}_{n\in\N}\subset L^p(\Omega)$ be a sequence with $u_n\to u\in L^p(\Omega)$. Then there exists a subsequence, again denoted by $\{u_n\}_{n\in\N}$, with $u_n(x)\to u(x)$ for almost every $x\in \Omega$. Since $z\mapsto f(x,z)$ is Newton differentiable almost everywhere, we have by definition that $r(x,u_n(x))\to 0$ for almost every $x\in\Omega$.
    Together with the boundedness from \eqref{eq:newton:superpos1}, Lebesgue's dominated convergence theorem therefore yields that $R(u_n)\to 0$ in $L^{s}(\Omega)$ (and hence along the full sequence since the limit is unique).\footnote{This is the step that fails for $F:L^\infty(\Omega)\to L^\infty(\Omega)$, since pointwise convergence and boundedness together do not imply uniform convergence almost everywhere.}
    For any $\tilde u\in L^p(\Omega)$, the Hölder inequality with $\frac1p+\frac1{s}=\frac1q$ thus yields that
    \begin{equation*}
        \begin{aligned}
            \norm{F(\tilde u)-F(u)-D_NF(\tilde u)(\tilde u-u)}_{L^q} =
            \norm{R(\tilde u)(\tilde u-u)}_{L^q}
            \leq \norm{R(\tilde u)}_{L^{s}}\norm{\tilde u-u}_{L^p}.
        \end{aligned}
    \end{equation*}
    If we now set $\tilde u \defeq u+h$ for $h\in L^p(\Omega)$ with $\norm{h}_{L^p}\to 0$, we have that $\norm{R(u+h)}_{L^{s}}\to 0$ and hence by definition the Newton differentiability of $F$ in $u$ with Newton derivative $h\mapsto D_N F(u)h$ as claimed.
\end{proof}
\begin{example}\label{ex:semismooth:l2}~
    \begin{enumerate}
        \item\label{ex:semismooth:l2:box}
            Consider
            \begin{equation*}
                A \defeq \setof{u\in L^2(\Omega)}{a \leq u(x) \leq b\quad \text{for almost every }x\in \Omega}
            \end{equation*}
            and $\proj_A : L^p(\Omega)\to L^2(\Omega)$ for $p>2$, which by \cref{lem:lebesgue:proximal} can be written as a superposition operator of the corresponding Lipschitz continuous scalar projection $\proj_{[a,b]}$, whose Newton derivative is given in \cref{ex:newton:rn}\,\ref{ex:newton:rn:box}. Since this derivative is clearly bounded (by $1$) and the pointwise limit of continuous functions, \cref{thm:newton:super} yields the pointwise almost everywhere Newton derivative
            \begin{equation*}
                [D_N \proj_A(u)h](x) =
                [\1_{[a,b]}(u)h](x) \defeq
                \begin{cases}
                    h(x) & \text{if }u(x)\in [a,b], \\
                    0 & \text{if }u(x) \notin[a,b].
                \end{cases}
            \end{equation*}
        \item\label{ex:semismooth:l2:l1}
            Consider now
            \begin{equation*}
                G:L^2(\Omega)\to \R,\qquad G(u) = \norm{u}_{L^1} = \int_\Omega |u(x)|\,dx
            \end{equation*}
            and $\prox_{\gamma G}:L^p(\Omega)\to L^2(\Omega)$ for $p>2$ and $\gamma>0$, which by \cref{lem:lebesgue:proximal} can be written as a superposition operator of the corresponding Lipschitz continuous scalar soft shrinkage operator, whose Newton derivative is given in \cref{ex:newton:rn}\,\ref{ex:newton:rn:l1}. Since this derivative is clearly bounded (by $1$) and the pointwise limit of continuous functions, \cref{thm:newton:super} yields the pointwise almost everywhere Newton derivative
            \begin{equation*}
                [D_N \prox_{\gamma G}(u) h](x) = [\1_{\R\setminus(-\gamma,\gamma)}(u)h](x) \defeq
                \begin{cases}
                    h(x) & \text{if }|u(x)|\geq \gamma, \\
                    0 & \text{if }|u(x)|<\gamma.
                \end{cases}
            \end{equation*}
    \end{enumerate}
\end{example}

For $p=q\in[1,\infty]$, however, the claim is false in general, as can be shown by counterexamples.
\begin{example}
    We take
    \begin{equation*}
        f:\R\to\R,\qquad f(z) = \max\{0,z\} \defeq
        \begin{cases}
            0 &\text{if }z\leq 0, \\
            z &\text{if }z\geq 0.
        \end{cases}
    \end{equation*}
    This is a piecewise differentiable function, and hence by \cref{thm:newton:clarke_ndiff} we can for any $\delta\in[0,1]$ take as Newton derivative
    \begin{equation*}
        D_N f(z)h =
        \begin{cases}
            0 & \text{if }z<0,\\
            \delta h &\text{if } z=0,\\
            h & \text{if }z>0.
        \end{cases}
    \end{equation*}
    We now consider the corresponding superposition operators $F:L^p(\Omega)\to L^p(\Omega)$ and $D_N F(u)\in \linear(L^p(\Omega);L^p(\Omega))$ for any $p\in [1,\infty)$ and show that the approximation condition \eqref{eq:newton:semismooth} is violated for $\Omega=(-1,1)$, $u(x) = -|x|$, and
    \begin{equation*}
        h_n(x) =
        \begin{cases}
            \frac1n & \text{if } |x| < \frac1n,\\
            0 & \text{if } |x|\geq \frac1n.
        \end{cases}
    \end{equation*}
    First, it is straightforward to compute $\norm{h_n}_{L^p}^p = \frac2{n^{p+1}}$. Then since
    \begin{equation*}
        [F(u)](x) = \max\{0,-|x|\} = 0 \qquad\text{almost everywhere,}
    \end{equation*}
    we have that
    \begin{equation*}
        [F(u+h_n) - F(u) - D_N F(u+h_n)h_n](x) =
        \begin{cases}
            - |x| & \text{if }|x| < \frac1n,\\
            0    & \text{if }|x| > \frac1n,\\
            -\frac\delta{n} & \text{if } |x| = \frac1n,
        \end{cases}
    \end{equation*}
    and thus
    \begin{equation*}
        \norm{F(u+h_n) - F(u) - D_N F(u+h_n)h_n}^p_{L^p} = \int_{-\frac1n}^{\frac1n} |x|^p\,dx = \frac2{p+1}\left(\frac1n\right)^{p+1}.
    \end{equation*}
    This implies that
    \begin{equation*}
        \lim_{n\to\infty} \frac{\norm{F(u+h_n) - F(u) - D_N F(u+h_n)h_n}_{L^p}}{\norm{h_n}_{L^p}} = \left(\frac1{p+1}\right)^{\frac1p} \neq 0
    \end{equation*}
    and hence that $F$ is not Newton differentiable from $L^p(\Omega)$ to $L^p(\Omega)$ for any $p<\infty$.

    \bigskip

    For the case $p=q=\infty$, we take $\Omega=(0,1)$, $u(x) = x$, and
    \begin{equation*}
        h_n(x) =
        \begin{cases}
            nx-1 & \text{if }x\leq \frac1n,\\
            0 & \text{if }x\geq \frac1n,
        \end{cases}
    \end{equation*}
    such that $\norm{h_n}_{L^\infty}=1$ for all $n\in\N$. We also have that $x + h_n = (1+n)x-1\leq 0$ for $x\leq \frac{1}{n+1}\leq \frac1n$ and hence that
    \begin{equation*}
        [F(u+h_n) - F(u) - D_N F(u+h_n)h_n](x) =
        \begin{cases}
            (1+n)x-1 & \text{if }x \leq \frac1{n+1},\\
            0    & \text{if }x \geq \frac1{n+1},
        \end{cases}
    \end{equation*}
    since either $h_n=0$ or $F(u+h_n)=F(u)+D_NF(u)h_n$ in the second case. Now,
    \begin{equation*}
        \sup_{x\in(0,\frac{1}{n+1}]}|(1+n)x-1| = 1 \qquad\text{for all }n\in\N,
    \end{equation*}
    which implies that
    \begin{equation*}
        \lim_{n\to\infty} \frac{\norm{F(u+h_n) - F(u) - D_N F(u+h_n)h_n}_{L^p}}{\norm{h_n}_{L^p}} = 1 \neq 0
    \end{equation*}
    and hence that $F$ is not Newton differentiable from $L^\infty(\Omega)$ to $L^\infty(\Omega)$ either.
\end{example}

\bigskip

\begin{remark}
    Semismoothness was introduced in \cite{Mifflin:1977} for Lipschitz-continuous functionals $F:\R^N\to\R$ as a condition relating Clarke subderivatives and directional derivatives near a point. This definition was extended to functions $F:\R^N\to \R^M$ in \cite{Qi:1993,Qi:1993a} and shown to imply a uniform version of the approximation condition \eqref{eq:newton:semismooth} for all elements of the Clarke subdifferential and hence superlinear convergence of the semismooth Newton method in finite dimensions.
    A semismooth Newton method specifically for PC$^1$ functions was already considered in \cite{kojima1986newton}.
    In normed vector spaces, \cite{Kummer:1988} was the first to study an abstract class of Newton methods for nonsmooth equations based on the condition \eqref{eq:newton:semismooth}, unifying the previous results; see \cite{KlatteKummer:2002}. In all these works, the analysis was based on semismoothness as a property relating $F:X\to Y$ to a set-valued mapping $G:X\setto \linear(X,Y)$, whose elements (uniformly) satisfy \eqref{eq:newton:semismooth}. In contrast, \cite{Kummer:2000,Chen:2000a} considered -- as we do in this book -- single-valued Newton derivatives (named \emph{Newton maps} in the former and \term[function!slanting]{slanting functions} in the latter) in Banach spaces. This approach was later followed in
    \cite{Hintermuller:2002a,Kunisch:2008a} to show that for a specific choice of Newton derivative, the classical \term[method!primal-dual active set]{primal-dual active set method} for solving quadratic optimization problems under linear inequality constraints can be interpreted as a semismooth Newton method. In parallel, \cite{Ulbrich:2002a,Ulbrich:2011} showed that superposition operators defined by semismooth functions (in the sense of \cite{Qi:1993a}) are semismooth (in the sense of \cite{Kummer:1988}) between the right spaces. A similar result for single-valued Newton derivatives was shown in \cite{Schiela:2008a} using a proof that is much closer to the one for the classical differentiability of superposition operators; compare \cref{thm:newton:super,thm:superpos:differentiable}. It should, however, be mentioned that not all calculus results for semismooth functions are available in the single-valued setting; for example, the implicit function theorem from \cite{kruse2018implicit} requires set-valued Newton derivatives, since the selection of the Newton derivative of the implicit function need not correspond to the selection of the given mapping.
    Finally, we remark that the notion of semismoothness and semismooth Newton methods were very recently extended to set-valued mappings in \cite{outrata2019semismooth}.
\end{remark}

\chapter{Nonlinear primal-dual proximal splitting}
\label{chap:nlpdps}

In this chapter, our goal is to extend the primal-dual proximal splitting (PDPS) method to \emph{nonlinear} operators $K \in C^1(X; Y)$, i.e., to problems of the form
\begin{equation}
    \label{eq:nlpdps:problem}
    \min_X F(x)+G(K(x)),
\end{equation}
where we still assume $F: X \to \Rbar$ and $G: Y \to \Rbar$ to be convex, proper, and lower semicontinuous on the Hilbert spaces $X$ and $Y$. For simplicity, we will only consider linear convergence under a strong convexity assumption and refer to the literature for weak convergence and acceleration under partial strong convexity (see \cref{rem:pdps} below). As in earlier chapters, we use the same notation for the inner product as for the duality pairing in Hilbert spaces to distinguish them better from pairs of elements.

We recall the three-point program for convergence proofs of first-order methods from \cref{chap:convergence}, which remains fundamentally the same in the nonlinear setting. However, we need to make some of the concepts local. Thus the three main ingredients of our convergence proofs will be the following.
\begin{enumerate}
    \item The three-point identity \eqref{eq:hilbert:three-point-identity}.

    \item The \term[monotonicity!local]{local monotonicity} of the operator $H$ whose roots correspond to the (primal-dual) critical points of \eqref{eq:nlpdps:problem}. We fix one of the points in the definition of monotonicity in \cref{sec:monotone:monotone} to a root $\realoptx$ of $H$, and only vary the other point in a neighborhood of $\realoptx$. This is essentially a nonsmooth variant of the standard second-order sufficient (or local quadratic growth) condition $\grad^2 F(x) \succ 0$ (i.e., positive definiteness of the Hessian) for minimizing a smooth function $F:\R^N\to \R$.

    \item The nonnegativity of the preconditioning operators $M_{k+1}$ defining the implicit form of the algorithm. These will now in general depend on the current iterate, and thus we can only show the nonnegativity in a neighborhood of suitable $\realoptx$.
\end{enumerate}

\section{Nonconvex explicit splitting}
\label{sec:nlpdps:fb}

To motivate our more specific assumptions on $K$, we start by showing that forward-backward splitting can be applied to a nonconvex function for the forward step. We thus consider for the problem
\begin{equation}
    \label{eq:nlpdps:fb:problem}
    \min_{x \in X} G(x) + F(x),
\end{equation}
with $F$ smooth but possibly nonconvex, the algorithm
\begin{algeqbox}
    \begin{equation}
        \label{eq:nlpdps:fb}
        x^{k+1} \defeq \prox_{\tau G}(x^k - \tau \grad F(x^k)).
    \end{equation}
\end{algeqbox}
To show convergence of this algorithm, we extend the non-value three-point smoothness inequalities of \cref{cor:smoothness:three-point,cor:smoothness:three-point:sc} from convex smooth functions to $C^2$ functions. (It is also possible to obtain corresponding value inequalities.)
\begin{lemma}
    \label{lemma:nlpds:fb:c2-smoothness}
    Suppose $F \in C^2(X)$.
    Let $z, \realoptx \in X$, and suppose for some $L>0$ and $\gamma \ge 0$ for all $\zeta \in \B(\realoptx,\norm{z-\realoptx}_X)$ that $\gamma \cdot \Id \le \grad^2 F(\zeta) \le L \cdot \Id$.
    Then for any $\beta \in (0, 2]$ and $x \in X$ we have
    \begin{equation}
        \label{eq:nlpdps:fb:c2-smoothness}
        \iprod{\grad F(z)-\grad F(\realoptx)}{x-\realoptx}_X
        \ge \frac{\gamma(2-\beta)}{2}\norm{x-\realoptx}_X^2
        -\frac{L}{2\beta} \norm{x-z}_X^2.
    \end{equation}
\end{lemma}
\begin{proof}
    By the one-dimensional mean value theorem applied to $t \mapsto \iprod{\grad F(\realoptx+t(z-\realoptx))}{x-\realoptx}_X$, we obtain for $\zeta=\realoptx+s(z-\realoptx)$ for some $s \in [0, 1]$ that
    \[
        \iprod{\grad F(z)-\grad F(\realoptx)}{x-\realoptx}_X
        =\iprod{\grad^2 F(\zeta)(z-\realoptx)}{x-\realoptx}_X.
    \]
    Therefore, applying the generalized Young inequality for arbitrary $\beta>0$ yields
    \begin{equation}
        \label{eq:nlpdps:fb:c2-smoothness-proof}
        \begin{aligned}[t]
            \iprod{\grad F(z)-\grad F(\realoptx)}{x-\realoptx}_X
            &
            =\norm{x-\realoptx}^2_{\grad^2 F(\zeta)}
            +\iprod{\grad^2 F(\zeta)(z-x)}{x-\realoptx}_X
            \\
            &
            \ge \frac{2-\beta}{2}\norm{x-\realoptx}^2_{\grad^2 F(\zeta)}
            -\frac{1}{2\beta} \norm{x-z}^2_{\grad^2 F(\zeta)}.
        \end{aligned}
    \end{equation}
    By the definition of $\gamma$ and $L$, we obtain \eqref{eq:nlpdps:fb:c2-smoothness}.
\end{proof}

The following result is almost a carbon copy of \cref{thm:convergence:fb,thm:testing:fb} for convex smooth $F$. However, since our present problem is nonconvex, we can only expect local convergence to a critical point of $J \defeq F+G$.

\begin{theorem}\label{thm:nlpdps:fb}
    Let $F \in C^2(X)$ and let $G:X\to\Rbar$ be proper, convex, and lower semicontinuous.
    Given an initial iterate $x^0$ and a critical point $\realoptx \in \inv{[\subdiff G+\grad F]}(0)$ of $J \defeq F+G$, let $\neighx \defeq \B(\realoptx, \norm{x^0-\realoptx})$, and suppose for some $L>0$ and $\gamma \ge 0$ that
    \begin{equation}
        \label{eq:nlpdps:fb:ass}
        \gamma \cdot \Id \le \grad^2 F(\zeta) \le L \cdot \Id
        \quad (\zeta \in \neighx).
    \end{equation}
    Take $0< \tau < 2L^{-1}$.
    \begin{enumerate}
        \item If $\gamma>0$, then the sequence $\{\thisx\}_{k \in \N}$ generated by \eqref{eq:nlpdps:fb} converges linearly to $\realoptx$.
        \item If $\gamma=0$, then the sequence $\{\thisx\}_{k \in \N}$ converges weakly to a critical point of $J$.
    \end{enumerate}
\end{theorem}

Note that if $G$ is locally finite-valued, then by \cref{thm:clarke:sum} our definition of a critical point in this theorem means $\realoptx \in \inv{[\subdiff_C J]}(0)$.

\begin{proof}
    As usual, we write \eqref{eq:nlpdps:fb} as
    \begin{equation}
        \label{eq:nlpdps:fb:prox}
        0 \in \tau[\subdiff G(\nextx)+\grad F(\thisx)]+(\nextx-\thisx).
    \end{equation}
    Suppose $\thisx \in \neighx$ and let $\beta \in (L\tau, 2)$ be arbitrary (which is possible since $\tau L <2$).
    By the monotonicity of $\subdiff G$ and the local three-point monotonicity \eqref{eq:nlpdps:fb:c2-smoothness} of $F$ implied by \cref{lemma:nlpds:fb:c2-smoothness}, we obtain
    \begin{equation}
        \label{eq:nlpdps:fb:prox-hest}
        \iprod{\subdiff G(\nextx)+\grad F(\thisx)}{\nextx-\realoptx}_X \ge \frac{\gamma(2-\beta)}{2}\norm{\nextx-\realoptx}_X^2-\frac{L}{2\beta}\norm{\nextx-\thisx}_X^2.
    \end{equation}
    Observe that if we had $\nextx=\thisx$ (or $F=0$), this would show the local quadratic growth of $F+G$ at $\realoptx$. Since, in general, $\nextx \ne \thisx$, we however need to compensate for taking the forward step with respect to $F$.

    Testing \eqref{eq:nlpdps:fb:prox} by the application of $\tauTest_k\iprod{\freevar}{\nextx-\realoptx}_X$ for some testing parameter $\tauTest_k>0$ and afterwards applying \eqref{eq:nlpdps:fb:prox-hest} yields
    \begin{equation*}
        \frac{\tauTest_k\gamma\tau(2-\beta)}{2}\norm{\nextx-\realoptx}_X^2 - \frac{\tauTest_k L\tau}{2\beta}\norm{\nextx-\thisx}_X^2 + \tauTest_k \iprod{\nextx-\thisx}{\nextx-\realoptx}_X \le 0.
    \end{equation*}
    Taking
    \begin{equation}
        \label{eq:nlpdps:fb:test-update}
        \tauTest_{k+1} \defeq \tauTest_k(1+\gamma\tau(2-\beta))
        \quad\text{with}\quad
        \tauTest_0>0,
    \end{equation}
    the three-point formula \eqref{eq:convergence:three-point-identity} yields
    \begin{equation}
        \label{eq:nlpdps:fb:fejer}
        \frac{\tauTest_{k+1}}{2}\norm{\nextx-\realoptx}_X^2
        +\frac{\tauTest_k(1-\tau L/\beta)}{2}\norm{\nextx-\thisx}_X^2
        \le \frac{\tauTest_k}{2}\norm{\thisx-\realoptx}_X^2.
    \end{equation}
    Since $\beta \in (L\tau,2)$ and $\thisx\in \neighx$, this implies that $\nextx\in \neighx$. By induction, we thus obtain that
    $\{\thisx\}_{k \in \N} \subset \neighx$ under our assumption $x^0 \in \neighx$.

    If $\gamma>0$, the recursion \eqref{eq:nlpdps:fb:test-update} together with $\beta<2$ shows that $\tauTest_k$ grows exponentially. Using that $\tau L/\beta \le 1$ and telescoping \eqref{eq:nlpdps:fb:fejer} then shows the claimed linear convergence.

    Let us then consider weak convergence. With $\gamma=0$ and $\beta<2$, the recursion \eqref{eq:nlpdps:fb:test-update} reduces to $\tauTest_{k+1} \equiv \tauTest_0>0$. Since $\tau L\leq \beta$, the estimate \eqref{eq:nlpdps:fb:fejer} yields Fejér monotonicity of the iterates $\{\thisx\}_{k \in \N}$. Moreover, we establish for $\nexxt{w} \defeq - \inv\tau(\nextx-\thisx)$ that $\norm{\nexxt{w}}_X \to 0$ and $\nexxt{w} \in \subdiff G(\nextx)+ \grad F(\thisx)$ for all $k \in \N$.
    Let $\optx$ be any weak accumulation point of $\{\thisx\}_{k \in \N}$, i.e., there exists a subsequence $\{x^{k_n}\}_{n \in \N}$ with $x^{k_n} \weakto \optx \in \neighx$. Then also $x^{k_n+1} \weakto \optx \in \neighx$. Since $\grad F$ is by \eqref{eq:nlpdps:fb:ass} Lipschitz continuous in $\neighx$, we have $\grad F(x^{k_n+1})-\grad F(x^{k_n}) \to 0$. Consequently, $\subdiff G(x^{k_n+1}) + \grad F(x^{k_n+1}) \ni w^{k_n+1} + \grad F(x^{k_n+1})- \grad F(x^{k_n}) \to 0$. By the outer semicontinuity of $\subdiff G+ \grad F$, it follows that $0 \in \subdiff G(\optx) + \grad F(\optx)$ and therefore $\optx \in \inv{(\subdiff G + \grad F)}(0) \subset \neighx$.
    The claim thus follows by applying Opial's \cref{lemma:opial}.
\end{proof}

\section{Nonconvex primal-dual splitting: algorithm and assumptions}

As mentioned above, we consider the problem \eqref{eq:nlpdps:problem} with $F: X \to \Rbar$ and $G: Y \to \Rbar$ convex, proper, and lower semicontinuous, and $K \in C^1(X; Y)$. We will soon state more precise assumptions on $K$.
When either the null space of $\kgradconj{x}$ is trivial or $\dom G=Y$,
we can apply the chain rule \cref{thm:clarke:chain} for Clarke subdifferentials as well as the equivalences of \cref{thm:clarke:convex,thm:clarke:frechet} for convex and differentiable functions, respectively, to rewrite as in \cref{sec:proximal:pd} the critical point conditions for this problem as $0 \in H(\realoptu)$ for the set-valued operator $H: X \times Y \setto X \times Y$ defined for $u=(x,y)\in X\times Y$ as
\begin{equation}
    \label{eq:nlpdps:h}
    H(u) \defeq
    \begin{pmatrix}
        \subdiff F(x) + \kgradconj{x} y \\
        \subdiff G^*(y) - K(x)
    \end{pmatrix}.
\end{equation}
Throughout the rest of this chapter, we write $\realoptu=(\realoptx, \realopty) \in \inv H(0)$ for an arbitrary root of $H$ that we assume to exist.

In analogy to the basic PDPS method, the basic unaccelerated nonlinear PDPS method then iterates
\begin{algeqbox}
    \begin{equation}
        \label{eq:nlpdps:nlpdps}
        \left\{
            \begin{aligned}
                \nextx & \defeq (I + \tau \subdiff F)^{-1}(\thisx - \tau\kgradconj{\thisx}\thisy),\\
                \overnextx & \defeq  (1+\omega)\nextx-\omega\thisx,\\
                \nexty & \defeq  (I + \sigma \subdiff G^*)^{-1}(\thisy + \sigma K(\overnextx)),
            \end{aligned}
        \right.
    \end{equation}
\end{algeqbox}
for some acceleration parameter $\tilde\gamma_{G^*} \ge 0$ (later to be fixed to be less than the factor of strong convexity of $G^*$), and where we set the over-relaxation parameter
\begin{equation}
    \label{eq:nlpdps:omega0}
    \omega = \frac{1}{1+2\tilde\gamma_{G^*} \sigma}.
\end{equation}
We can write this algorithm in the general form of \cref{thm:gap:ppext:convergence} as follows. For each iteration $k \in \N$ with some primal and dual testing parameters $\tauTest_k,\sigmaTest_{k+1}>0$,  we define the step length and testing operators
\[
    \Step \defeq
    \begin{pmatrix}
        \tau \Id & 0 \\
        0 & \sigma \Id
    \end{pmatrix}
    \quad\text{and}\quad
    \Test_{k+1} \defeq
    \begin{pmatrix}
        \tauTest_k \Id & 0 \\
        0 & \sigmaTest_{k+1} \Id
    \end{pmatrix}.
\]
We also define the linear preconditioner $\Precond_{k+1}$ and the step length weighted partial linearization $\Happrox_{k+1}$ of $H$ by
\begin{align}
    \label{eq:nlpdps:precond}
    \Precond_{k+1} & \defeq
    \begin{pmatrix}
        \Id & -\tau \kgradconj{\thisx} \\
        -\omega \sigma \kgrad{\thisx} & \Id
    \end{pmatrix},
    \quad\text{and}
    \\
    \label{eq:nlpdps:tilde-h}
    \Happrox_{k+1}(u) & \defeq \Step\begin{pmatrix}
        \subdiff F(x) + \kgradconj{\thisx} y \\
        \subdiff G^*(y)-K(\overnextx)-\kgrad{\thisx}(x-\overnextx)
    \end{pmatrix}.
\end{align}
Observe that $\Happrox_{k+1}(u)$ simplifies to $\Step H(u)$ for linear $K$.
Then \eqref{eq:nlpdps:nlpdps} becomes
\begin{equation}
    \label{eq:nlpdps:ppext}
    0 \in \Happrox_{k+1}(\nextu)+\Precond_{k+1}(\nextu-\thisu).
\end{equation}

We will need $K$ to be locally Lipschitz differentiable.

\begin{assumption}[locally Lipschitz $\grad K$]
    \label{ass:nlpdps:k-lipschitz}
    The operator $K:X\to Y$ is Fréchet differentiable, and for some $L \ge 0$ and a neighborhood $\neighx_K$ of $\realoptx$,
    \begin{equation}%
        \label{eq:nlpdps:ass-k-lipschitz}
        \|\kgrad{x}-\kgrad{z}\|_{\linear(X,Y)} \le L\|x-z\|_X \quad (x,z\in\neighx_K).
    \end{equation}%
\end{assumption}

We also require a three-point assumption on $K$. This assumption combines a second-order growth condition with a three-point smoothness estimate. Note that the factor $\gamma_K$ can be negative; if it is, it will need to be offset by sufficient strong convexity of $F$.

\begin{assumption}[three-point condition on $K$]
    \label{ass:nlpdps:k-nonlinear}
    There exists a neighborhood $\neighx_K$ of $\realoptx$ and $\gamma_K \in \R$ as well as $\lambda, \theta \ge 0$ such that
    \begin{multline}
        \label{eq:nlpdps:ass-k-nonlinear}
        \iprod{[\kgrad{z}-\kgrad{\realoptx}]^*\realopty}{x-\realoptx}_X
        \\
        \ge \gamma_K\norm{x-\realoptx}_X^2 + \theta\norm{K(\realoptx)-K(x)-\kgrad{x}(\realoptx-x)}_Y-\frac{\lambda}{2}\norm{x-z}_X^2
        \quad (x,z\in\neighx_K).
    \end{multline}
\end{assumption}

We observe the following special cases of \cref{ass:nlpdps:k-nonlinear}:
\begin{enumerate}[label=(\alph*)]
    \item For linear $K$, the assumption trivially holds for any $\gamma_K \le 0$, $\theta \ge 0$ and $\lambda=0$.
    \item Let $G^*=\delta_{\{1\}}$, so that $K: X \to \R$ and the problem \eqref{eq:nlpdps:problem} reduces to \eqref{eq:nlpdps:fb:problem} with $K$ in place of $F$. Keeping in mind that in this case $\realopty=1$, \cref{lemma:nlpds:fb:c2-smoothness} with $\beta=1$ shows that \cref{ass:nlpdps:k-nonlinear} is satisfied for $\lambda=L$, any $\theta \ge 0$, and $\gamma_K \le \gamma$ with $\gamma, L \ge 0$ satisfying $\gamma \cdot \Id \le \grad^2 K(\zeta) \le L \cdot \Id$ for all $\zeta \in \neighx_K$.
\end{enumerate}
In more general settings, the verification of \cref{ass:nlpdps:k-nonlinear} can demand some effort.
We refer to \cite{tuomov-nlpdhgm-redo} for examples and to \cite{tuomov-nlpdhgm-general} for further generalizations.

\section{Nonconvex primal-dual splitting: convergence proof}

For simplicity of treatment, and to demonstrate the main ideas without excessive technicalities, we only show linear convergence under strong convexity of both $F$ and $G^*$.

We will base our proof on \cref{thm:gap:ppext:convergence} and thus have to verify its assumptions. Most of the work is in verifying the inequality \eqref{eq:gap:ppext:fundamental-condition}, which we do in several steps. First, we ensure that the operator $\Test_{k+1}\Precond_{k+1}$ giving rise to the local metric is self-adjoint.
Then we show that $\Test_{k+2}\Precond_{k+2}$ and the update $\Test_{k+1}(\Precond_{k+1}+\GammaLift{k+1})$ actually performed by the algorithm yield identical norms, where $\GammaLift{k+1}$ represents some off-diagonal components from the algorithm as well as any strong convexity provided by $F$ and $G^*$.
Finally, we estimate the local monotonicity of $\Happrox_{k+1}$.

We write $\gamma_F, \gamma_{G^*} \ge 0$ for the factors of (strong) convexity of $F$ and $G^*$, and recall the factor $\gamma_K \in \R$ from \cref{ass:nlpdps:k-nonlinear}.
Then for some ``acceleration parameters'' $\tilde\gamma_F, \tilde\gamma_{G^*} \ge 0$ and $\kappa \in [0, 1)$, we require that
\begin{subequations}
    \label{eq:nlpdps:basic-step-rules}%
    \begin{align}%
        \gamma_F + \gamma_K & \ge \tilde\gamma_F \ge 0, &
        \gamma_{G^*} & \ge \tilde\gamma_{G^*} \ge 0, \\
        \eta_k & \defeq \tauTest_k\tau = \sigmaTest_k\sigma,
        &
        1-\kappa & \le \tau\sigma \norm{\kgrad{\thisx}}^2,
        \\
        \tauTest_{k+1} & = \tauTest_k(1+2\tilde\gamma_F \tau),
        \quad\text{and}\quad
        &
        \sigmaTest_{k+1} & = \sigmaTest_k(1+2\tilde\gamma_{G^*} \sigma)
        \quad (k \in \N).
    \end{align}
\end{subequations}
With this, $\omega$ defined in \eqref{eq:nlpdps:omega0} satisfies
\begin{equation}
    \label{eq:nlpdps:omega}
    \omega = \inv\sigmaTest_{k+1}\sigmaTest_k = \inv\eta_{k+1}\eta_k
    \quad (k \in \N).
\end{equation}

The next lemma adapts \cref{lemma:convergence:pd_spd}.
\begin{lemma}
    \label{lemma:nlpdps:zimi-estim}
    Fix $k \in \N$ and suppose \eqref{eq:nlpdps:basic-step-rules} holds. Then $\Test_{k+1}\Precond_{k+1}$ is self-adjoint and satisfies
    \[
        \Test_{k+1}\Precond_{k+1}
        \succeq
        \begin{pmatrix}
            \delta \tauTest_k \cdot \Id & 0 \\
            0 & (\kappa-\delta)(1-\delta)^{-1}\sigmaTest_{k+1} \cdot \Id
        \end{pmatrix}
        \quad\text{ for any }\delta \in [0, \kappa].
    \]
\end{lemma}
\begin{proof}
    From \eqref{eq:nlpdps:basic-step-rules} and \eqref{eq:nlpdps:omega} we have
    $\tauTest_k\tau=\sigmaTest_{k+1}\omega\sigma=\eta_k$. By \eqref{eq:nlpdps:precond} then
    \begin{equation}
        \label{eq:nlpdps:test-precond-expansion}
        \Test_{k+1}\Precond_{k+1}
        =
        \begin{pmatrix}
            \tauTest_k \cdot \Id & -\eta_k \kgradconj{\thisx} \\
            -\eta_k \kgrad{\thisx} & \sigmaTest_{k+1} \cdot \Id
        \end{pmatrix}.
    \end{equation}
    This shows that $\Test_{k+1}\Precond_{k+1}$ is self-adjoint.
    Furthermore, since Young's inequality followed by \eqref{eq:nlpdps:basic-step-rules} and \eqref{eq:nlpdps:omega} shows that
    \[
        \begin{aligned}
        2\eta_k\iprod{\kgrad{\thisx}\tilde x}{\tilde y}
        &
        \le
        (1-\delta)\eta_k\inv\tau\norm{\tilde x}^2 + \frac{\eta_k\tau}{1-\delta}\norm{\kgrad{\thisx}^* \tilde y}^2
        \\
        &
        =
        (1-\delta)\tauTest_k\norm{\tilde x}^2 + \sigmaTest_{k+1}\omega\frac{\tau\sigma}{1-\delta}\norm{\kgrad{\thisx}^* \tilde y}^2
        \quad
        (\tilde x \in X,\,\tilde y \in Y),
        \end{aligned}
    \]
    we obtain from \eqref{eq:nlpdps:test-precond-expansion} that
    \begin{equation}
        \label{eq:nlpdps:test-precond-expansion-estimate}
        \Test_{k+1}\Precond_{k+1}
        \succeq
        \begin{pmatrix}
            \delta \tauTest_k \Id & 0 \\
            0 & \sigmaTest_{k+1}\left(\Id - \omega\frac{\tau\sigma}{1-\delta}\kgrad{\thisx}\kgradconj{\thisx}\right)
        \end{pmatrix}.
    \end{equation}
    The claimed estimate then follows from the assumptions \eqref{eq:nlpdps:basic-step-rules}.
\end{proof}

Our next step is to simplify the operator $\Test_{k+1}\Precond_{k+1}-\Test_{k+2}\Precond_{k+2}$ occurring in the inequality \eqref{eq:gap:ppext:fundamental-condition} we are trying to prove.

\begin{lemma}
    \label{lemma:nlpdps:local-metric-transfer}
    Fix $k \in \N$, and suppose \eqref{eq:nlpdps:basic-step-rules} holds.
    Then
    $
    \frac{1}{2}\norm{\freevar}_{\Test_{k+1}(\Precond_{k+1}+\GammaLift{k+1})-\Test_{k+2}\Precond_{k+2}}^2
    =
    0
    $ for
    \begin{equation}
        \label{eq:nlpdps:gammalift-def}
        \GammaLift{k+1}\defeq
        \begin{pmatrix}
            2 \tilde\gamma_F \tau \Id & 2\tau \kgradconj{\thisx} \\
            -2\sigma\kgrad{\nextx} & 2 \tilde\gamma_{G^*} \sigma \Id
        \end{pmatrix}.
    \end{equation}
\end{lemma}
\begin{proof}
    Using \eqref{eq:nlpdps:basic-step-rules} and \eqref{eq:nlpdps:test-precond-expansion}, we can write
    \[
        \Test_{k+1}(\Precond_{k+1}+\GammaLift{k+1})-\Test_{k+2}\Precond_{k+2}
        =
        D_{k+1}
    \]
    for the skew-symmetric operator
    \[
        D_{k+1}
        \defeq
        \begin{pmatrix}
            0 & [\eta_{k+1}\kgrad{\nextx} + \eta_k\kgrad{\thisx}]^* \\
            -[\eta_{k+1}\kgrad{\nextx} + \eta_k\kgrad{\thisx}] & 0 \\
        \end{pmatrix}.
    \]
    This yields the claim.
\end{proof}

For our convergence claim, we need to assume that the dual variables stay bounded within the ``nonlinear range'' of $K$. To this end, we introduce the (possibly empty) subspace $\Ylin$ of $Y$ in which $K$ acts linearly, i.e.,
\[
    \Ylin \defeq \{ y \in Y \mid \text{the mapping } x \mapsto \iprod{y}{K(x)} \text{ is linear} \}
    \quad\text{and}\quad
    \Ynl \defeq \Ylin^\perp.
\]
We then denote by $\Pnl$ the orthogonal projection to $\Ynl$. We also write
\[
    \B_\NL(\realopty, r) \defeq \{ y \in Y \mid \norm{y-\realopty}_{\Pnl} \le r\}
\]
for the closed cylinder in $Y$ of the radius $r$ with axis orthogonal to $\Ynl$.

With $\neighx_K$ given by \cref{ass:nlpdps:k-lipschitz}, we now define for some radius $\metricRhoY>0$ the neighborhood
\begin{equation}
    \label{eq:neighu-definition}
    \neighu (\metricRhoY) \defeq \neighx_K \times \B_{\NL}(\realopty, \metricRhoY).
\end{equation}
We will require that the iterates $\{\thisu\}_{k \in \N}$ of \eqref{eq:nlpdps:nlpdps} stay within this neighborhood for some fixed $\metricRhoY>0$.

The next lemma provides the necessary three-point inequality to estimate the linearizations performed within $\Happrox_{k+1}$.
\begin{lemma}
    \label{lemma:nlpdps:nonlinear-preconditioner-estimate}
    For a fixed $k \in \N$, suppose $\overnextx \in \neighx_K$, and let $\metricRhoY \ge 0$ be such that $\thisu, \nextu\in \neighu(\metricRhoY)$.
    Suppose $K$ satisfies \cref{ass:nlpdps:k-lipschitz,ass:nlpdps:k-nonlinear} with $\omega\theta\ge\metricRhoY$.
    If \eqref{eq:nlpdps:basic-step-rules} holds, then
    \begin{equation*}
        \iprod{\Happrox_{k+1}(\nextu)}{\nextu-\realoptu}_{\Test_{k+1}}
        \ge
        \frac{1}{2}\norm{\nextu-\realoptu}_{\Test_{k+1}\GammaLift{k+1}}^2
        -\frac{\eta_k[\lambda+3L\metricRhoY]}{2}\norm{\nextx-\thisx}_X^2.
    \end{equation*}
\end{lemma}
\begin{proof}
    From \eqref{eq:nlpdps:h}, \eqref{eq:nlpdps:tilde-h}, \eqref{eq:nlpdps:basic-step-rules}, and \eqref{eq:nlpdps:gammalift-def}, we calculate
    \begin{equation}
        \label{eq:nlpdps:nonlinear-preconditioner-estimate-d-def}
        \begin{aligned}[t]
            D&\defeq \iprod{\Happrox_{k+1}(\nextu)}{\nextu-\realoptu}_{\Test_{k+1}}
            -\frac{1}{2}\norm{\nextu-\realoptu}_{\Test_{k+1}\GammaLift{k+1}}^2\\
            &=\iprod{H(\nextu)}{\nextu-\realoptu}_{\Test_{k+1}\Step}
            -\eta_k\tilde\gamma_{F} \norm{\nextx-\realoptx}_X^2
            -\eta_{k+1}\tilde\gamma_{G^*} \norm{\nexty-\realopty}_Y^2
            \\
            &\quad
            +\eta_k\iprod{[\kgrad{\thisx}-\kgrad{\nextx}](\nextx-\realoptx)}{\nexty}_Y
            \\
            &\quad
            +\eta_{k+1}\iprod{K(\nextx)-K(\overnextx)-\kgrad{\thisx}(\nextx-\overnextx)}{\nexty-\realopty}_Y
            \\
            &\quad
            + \iprod{(\eta_{k+1}\kgrad{\nextx}-\eta_k\kgrad{\thisx})(\nextx-\realoptx)}{\nexty-\realopty}_Y.
        \end{aligned}
    \end{equation}
    Here the first of the terms involving $K$ comes from the first lines of $\Happrox_{k+1}$ and $H$, the second of the terms from the second line, and the third from $\GammaLift{k+1}$.
    Since $0 \in H(\realoptu)$, we have $q_F \defeq  - \kgradconj{\realoptx}\realopty \in \subdiff F(\realoptx)$ and $q_{G^*} \defeq  K(\realoptx) \in \subdiff G^*(\realopty)$.
    Using \eqref{eq:nlpdps:basic-step-rules}, we can therefore expand
    \begin{equation*}
        \begin{aligned}
            \iprod{H(\nextu)}{\nextu-\realoptu}_{\Test_{k+1}\Step}
            & =\eta_k\iprod{\subdiff F(\nextx)-q_F}{\nextx-\realoptx}_X
            \\
            \MoveEqLeft[-1]
            +\eta_{k+1} \iprod{\subdiff G^*(\nexty)-q_{G^*}}{\nexty-\realopty}_Y
            \\
            \MoveEqLeft[-1]
            +\eta_k\iprod{\kgradconj{\nextx} \nexty-\kgradconj{\realoptx}\realopty}{\nextx-\realoptx}_X
            \\
            \MoveEqLeft[-1]
            +\eta_{k+1}\iprod{K(\realoptx)-K(\nextx)}{\nexty-\realopty}_Y.
        \end{aligned}
    \end{equation*}
    Using the $\gamma_F$-strong monotonicity of $\subdiff F$ and the $\gamma_{G^*}$-strong monotonicity of $\subdiff G^*$, and rearranging terms, we obtain
    \begin{equation*}
        \begin{aligned}[t]
            \iprod{H(\nextu)}{\nextu-\realoptu}_{\Test_{k+1}\Step}
            & \ge \eta_k\gamma_F \norm{\nextx-\realoptx}_X^2 + \eta_{k+1}\gamma_{G^*} \norm{\nexty-\realopty}_Y^2
            \\
            \MoveEqLeft[-1]
            +\eta_k\iprod{\kgrad{\nextx}(\nextx-\realoptx)}{\nexty}_Y
            \\
            \MoveEqLeft[-1]
            -\eta_k\iprod{\kgrad{\realoptx}(\nextx-\realoptx)}{\realopty}_Y
            \\
            \MoveEqLeft[-1]
            +\eta_{k+1}\iprod{K(\realoptx)-K(\nextx)}{\nexty-\realopty}_Y.
        \end{aligned}
    \end{equation*}
    Combining this estimate with \eqref{eq:nlpdps:nonlinear-preconditioner-estimate-d-def} and rearranging terms, we obtain
    \begin{equation*}
        \begin{aligned}
            D&\ge
            \eta_k(\gamma_F-\tilde\gamma_{F})\norm{\nextx-\realoptx}_X^2
            + \eta_{k+1}(\gamma_{G^*}-\tilde\gamma_{G^*})\norm{\nexty-\realopty}_Y^2
            \\
            &\quad
            -\eta_k\iprod{\kgrad{\realoptx}(\nextx-\realoptx)}{\realopty}_Y+\eta_k\iprod{\kgrad{\thisx}(\nextx-\realoptx)}{\nexty}_Y
            \\
            &\quad
            +\eta_{k+1}\iprod{K(\realoptx)-K(\overnextx)-\kgrad{\thisx}(\nextx-\overnextx)}{\nexty-\realopty}_Y
            \\
            &\quad
            + \iprod{(\eta_{k+1}\kgrad{\nextx}-\eta_k\kgrad{\thisx})(\nextx-\realoptx)}{\nexty-\realopty}_Y.
        \end{aligned}
    \end{equation*}
    Further rearrangements and $\gamma_F + \gamma_K \ge \tilde\gamma_F$ and $\gamma_{G^*} \ge \tilde\gamma_{G^*}$ give
    \begin{equation}
        \label{eq:nlpdps:d2}
        \begin{aligned}[t]
            D &
            \ge
            -\eta_k\gamma_K \norm{\nextx-\realoptx}_X^2
            +\eta_k\iprod{[\kgrad{\thisx}-\kgrad{\realoptx}](\nextx-\realoptx)}{\realopty}_Y
            \\
            &\quad
            +\eta_{k+1}\iprod{K(\realoptx)-K(\nextx)-\kgrad{\nextx}(\realoptx-\nextx)}{\nexty-\realopty}_Y
            \\
            &\quad
            +\eta_{k+1}\iprod{K(\nextx)-K(\overnextx)+\kgrad{\nextx}(\overnextx-\nextx)}{\nexty-\realopty}_Y
            \\
            &\quad
            +\eta_{k+1}\iprod{(\kgrad{\thisx}-\kgrad{\nextx})(\overnextx-\nextx)}{\nexty-\realopty}_Y.
        \end{aligned}
    \end{equation}

    Using \cref{ass:nlpdps:k-lipschitz} and the mean value theorem in the form
    \[
        K(x')=K(x)+\kgrad{x}(x'-x)+\int_{0}^{1}(\kgrad{x+s(x'-x)}-\kgrad{x})(x'-x)ds,
    \]
    we obtain for any $x,x' \in \neighx_K$ and $y\in Y$ the inequality
    \begin{equation}
        \label{eq:nlpdps:ass-k-lipschitz-2}
        \iprod{K(x')-K(x)-\kgrad{x}(x'-x)}{y}_Y \le (L/2) \norm{x-x'}_X^2\norm{y}_{\Pnl}.
    \end{equation}
    Applying \cref{ass:nlpdps:k-lipschitz}, the inequality \eqref{eq:nlpdps:ass-k-lipschitz-2}, and $\overnextx-\nextx=\omega(\nextx-\thisx)$ to the last two terms of \eqref{eq:nlpdps:d2}, we obtain
    \begin{align*}
        \iprod{K(\nextx)-K(\overnextx)+\kgrad{\nextx}(\overnextx-\nextx)}{\nexty-\realopty}_Y
        &\ge -\frac{L\omega^2}{2}\norm{\nextx-\thisx}^2\norm{\nexty-\realopty}_{\Pnl}
        \\
        \shortintertext{and}
        \iprod{(\kgrad{\thisx}-\kgrad{\nextx})(\overnextx-\nextx)}{\nexty-\realopty}_Y
        &\ge -L\omega\norm{\nextx-\thisx}_X^2\norm{\nexty-\realopty}_{\Pnl}.
    \end{align*}
    These estimates together with \eqref{eq:nlpdps:basic-step-rules} and $\nextu\in \neighu(\metricRhoY)$ now imply that $D\ge\eta_k D^K_{k+1}$ for
    \begin{equation*}
        \begin{aligned}
            D^K_{k+1}
            & \defeq \iprod{[\kgrad{\thisx}-\kgrad{\realoptx}](\nextx-\realoptx)}{\realopty}_Y
            -\gamma_K\norm{\nextx-\realoptx}_X^2-L(1+\omega/2)\metricRhoY\norm{\nextx-\thisx}_X^2
            \\
            \MoveEqLeft[-1]
            -\inv\omega\norm{\nexty-\realopty}_{\Pnl} \norm{K(\realoptx)-K(\nextx)-\kgrad{\nextx}(\realoptx-\nextx)}_Y.
        \end{aligned}
    \end{equation*}
    Finally, we use \cref{ass:nlpdps:k-nonlinear} and Young's inequality to estimate
    \begin{equation*}
        \begin{aligned}[t]
            D^K_{k+1}
            &\ge(\theta-\inv\omega\norm{\nexty-\realopty}_{\Pnl})\norm{K(\realoptx)-K(\nextx)-\kgrad{\nextx}(\realoptx-\nextx)}_Y
            \\ \MoveEqLeft[-1]
            -\frac{\lambda+3L\metricRhoY}{2}\norm{\nextx-\thisx}_X^2.
        \end{aligned}
    \end{equation*}
    Now observe that $\theta-\inv\omega\norm{\nexty-\realopty}_{\Pnl}\ge\theta-\inv\omega\metricRhoY \ge 0$. Combining with the estimate $D\ge\eta_k D^K_{k+1}$, we therefore obtain our claim.
\end{proof}

We now have all the necessary tools at hand to prove the main estimate \eqref{eq:gap:ppext:fundamental-condition} needed for the application of \cref{thm:gap:ppext:convergence}.

\begin{theorem}
    \label{thm:nlpdps:nonneg-penalty-nlpdhgm}
    Let $F: X \to \Rbar$ and $G: Y \to \Rbar$ be convex, proper, and lower semicontinuous.
    Suppose $K: X \to Y$ satisfies \cref{ass:nlpdps:k-lipschitz,ass:nlpdps:k-nonlinear}.
    Fix $k \in \N$, and also suppose that $\overnextx \in \neighx_K$ and that $\thisu, \nextu\in \neighu(\metricRhoY)$ for some $\metricRhoY \ge 0$.
    Suppose \eqref{eq:nlpdps:basic-step-rules} holds for some $\kappa \in [0, 1)$ and
    \begin{equation}
        \label{eq:nlpdps:tau-upperbound}
        \tau < \frac{\kappa}{\lambda+3L\metricRhoY}
    \end{equation}
    as well as $\omega\theta \ge \metricRhoY$.
    Then
    \begin{equation}
        \label{eq:nlpdps:quantitative-fejer}
        \frac{1}{2}\norm{\nextu-\realoptu}_{\Test_{k+2}\Precond_{k+2}}^2
        \le
        \frac{1}{2}\norm{\thisu-\realoptu}_{\Test_{k+1}\Precond_{k+1}}^2
        \quad
        (k \in \N).
    \end{equation}
\end{theorem}

\begin{proof}
    We show that \eqref{eq:gap:ppext:fundamental-condition} holds with
    $\GenGap_{k+1} \equiv 0$,
    i.e., that
    \begin{equation}
        \label{eq:nlpdps:fundamental-condition}
        \iprod{\Happrox_{k+1}(\nextu)}{\nextu-\realoptu}_{\Test_{k+1}}
        \ge
        \frac{1}{2}\norm{\nextu-\realoptu}_{\Test_{k+2}\Precond_{k+2}-\Test_{k+1}\Precond_{k+1}}^2
        -\frac{1}{2}\norm{\nextu-\thisu}_{\Test_{k+1}\Precond_{k+1}}^2.
    \end{equation}
    The claim then follows from \cref{thm:gap:ppext:convergence} and \cref{lemma:nlpdps:zimi-estim}, the latter of which provides the necessary self-adjointness of $\Test_{k+1}\Precond_{k+1}$.

    Let thus $\delta \in (0, \kappa)$ be arbitrary, and define
    \[
        S_{k+1} \defeq
        \begin{pmatrix}
            (\delta\tauTest_k
            -\eta_k[\lambda+3L\metricRhoY]) \Id
            & 0 \\
            0 & \sigmaTest_{k+1}\left(\Id - \omega\frac{\tau\sigma}{1-\delta}\kgrad{\thisx}\kgradconj{\thisx}\right)
        \end{pmatrix}.
    \]
    Using \eqref{eq:nlpdps:test-precond-expansion-estimate} and \eqref{eq:nlpdps:test-precond-expansion} and, in the second and third step, \cref{lemma:nlpdps:nonlinear-preconditioner-estimate,lemma:nlpdps:local-metric-transfer}, we estimate
    \begin{equation*}
        \begin{aligned}
            \frac{1}{2}\norm{\nextu-\thisu}_{S_{k+1}-\Test_{k+1} \Precond_{k+1}}^2
            &
            \le
            -\frac{\eta_k[\lambda+3L\metricRhoY]}{2}\norm{\nextx-\thisx}_X^2
            \\
            &
            \le
            \iprod{\Happrox_{k+1}(\nextu)}{\nextu-\realoptu}_{\Test_{k+1}}
            -\frac{1}{2}\norm{\nextu-\realoptu}_{\Test_{k+1}\GammaLift{k+1}}^2
            \\
            &
            =
            \iprod{\Happrox_{k+1}(\nextu)}{\nextu-\realoptu}_{\Test_{k+1}}
            -\frac{1}{2}\norm{\nextu-\realoptu}_{\Test_{k+2}\Precond_{k+2}-\Test_{k+1}\Precond_{k+1}}^2.
        \end{aligned}
    \end{equation*}
    Then \eqref{eq:nlpdps:fundamental-condition} holds if $S_{k+1} \ge 0$.
    This readily follows from \eqref{eq:nlpdps:tau-upperbound} with $\delta \in (0, \kappa)$ for the primal variable and \eqref{eq:nlpdps:basic-step-rules} for the dual variable.
\end{proof}

\begin{theorem}
    \label{thm:nlpdps:linear-convergence}
    Let $F: X \to \Rbar$ and $G: Y \to \Rbar$ be strongly convex, proper, and lower semicontinuous.
    Suppose $K: X \to Y$ satisfies \cref{ass:nlpdps:k-lipschitz,ass:nlpdps:k-nonlinear}.
    Let $R_K>0$ be such that $\sup_{x\in\neighx_K}\norm{\grad K(x)}\le R_K$.
    Pick $0<\tau<1/(\lambda+3L\rho_y)$ for a given $\metricRhoY \ge 0$, and take $\sigma=\tau\tilde\gamma_F/\tilde\gamma_{G^*}$ for some $\tilde\gamma_F \in (0, \gamma_F+\gamma_K]$ and $\tilde\gamma_{G^*} \in (0, \gamma_{G^*}]$ such that  $\omega\theta \ge \metricRhoY$.
    Let the iterates $\{(\thisu, \overnextx)\}_{k \in \N}$ be generated by the NL-PDPS method \eqref{eq:nlpdps:nlpdps}. If $\overnextx \in \neighx_K$ and $\thisu \in \neighu(\metricRhoY)$ for all $k \in \N$ and some  $\realoptu \in \inv H(0)$ for $H$ given in \eqref{eq:nlpdps:h}, then $\thisu \to \realoptu$ linearly.
\end{theorem}

\begin{proof}
    Take $\tauTest_{k+1} \defeq \tauTest_k(1+2\tilde\gamma_F \tau)$ and $\sigmaTest_{k+1} \defeq \sigmaTest_k(1+2\tilde\gamma_{G^*} \sigma)$ for $\tauTest_0=1$ and $\sigmaTest_1 \defeq \tau/\sigma$.
    Then $\tauTest_k\tau=\sigmaTest_{k+1}\sigma$ if and only if $1+2\tilde\gamma_F\tau=1+2\tilde\gamma_{G^*}\sigma$, i.e., for $\sigma=\tau\tilde\gamma_F/\tilde\gamma_{G^*}$ as stated. Consequently \eqref{eq:nlpdps:basic-step-rules} is satisfied and the testing parameters $\tauTest_k$ and $\sigmaTest_{k+1}$ grow exponentially.
    Clearly \eqref{eq:nlpdps:tau-upperbound} holds for some $\kappa \in [0, 1)$.
    Combining \eqref{eq:nlpdps:quantitative-fejer} from \cref{thm:nlpdps:nonneg-penalty-nlpdhgm} with \cref{lemma:nlpdps:zimi-estim} now shows the claimed linear convergence.
\end{proof}

Besides step length bounds and structural properties of the problem, \cref{thm:nlpdps:linear-convergence} still requires us to ensure that the iterates stay close enough to the critical point $\realoptx$. This can be done if we initialize close enough to a critical point. As the proof is very technical, we merely state the following result.

\begin{theorem}[{\cite[Proposition 4.8]{tuomov-nlpdhgm-redo}}]
    Under the assumptions of \cref{thm:nlpdps:linear-convergence}, for any $\metricRhoY>0$ there exists an $\varepsilon>0$ such that  $\{\thisu\}_{k \in \N} \subset \neighu(\metricRhoY)$ for all initial iterates $u^0=(x^0,y^0)$ satisfying
    \begin{equation}\label{eq:neighborhood-start-close}
        \sqrt{2\inv\delta(\norm{x^0-\realoptx}^2+\tau\inv\sigma\norm{y^0-\realopty}^2)}
        \le\varepsilon.
    \end{equation}
\end{theorem}

\begin{remark}[weaker assumptions, weaker convergence]
    We have only demonstrated linear convergence of the method under the strong convexity of both $F$ and $G^*$. However, under similarly weaker assumptions as for the basic PDPS method familiar from \cref{part:convex}, both an accelerated $O(1/N^2)$ rate and weak convergence can be proved. We refer to \cite{tuomov-nlpdhgm-redo} for details, noting that Opial's \cref{lemma:opial} extends straightforwardly to the quantitative Fejér monotonicity \eqref{eq:testing:structured:quantitative-fejer} that is the basis of our proofs here.
    We also note that our linear convergence result differs from that in \cite{tuomov-nlpdhgm-redo} by taking the over-relaxation parameter $\omega=1$ in \eqref{eq:nlpdps:nlpdps} instead of $\omega=1/(1+2\tilde\gamma_{G^*}\sigma)<1$; compare \cref{thm:testing:pdps:accel}.
\end{remark}

\begin{remark}[historical development of the NL-PDPS]
    \label{rem:pdps}
    The NL-PDPS method was first introduced in \cite{Valkonen:2014} in finite dimensions with applications to inverse problems in magnetic resonance imaging.
    The method was later extended in \cite{tuomov-pdex2nlpdhgm} to infinite dimensions and applied to PDE-constrained optimization problems.
    In these works, only (weak) convergence of the iterates is shown, based on the metric regularity of the operator $H$. We discuss metric regularity later in \cref{chap:regularity,chap:stability}.
    Convergence rates were then first shown in \cite{tuomov-nlpdhgm-redo}.
    In that paper, alternative forms of the three-point condition \cref{ass:nlpdps:k-nonlinear} on $K$ are also discussed.

    Similarly to how we showed in \cref{sec:proximal:connections} that the preconditioned ADMM is equivalent to the PDPS method, it is possible to derive a preconditioned nonlinear ADMM that is equivalent to the NL-PDPS method; such algorithms are considered in \cite{benning2015preconditioned}.
    The NL-PDPS method has been extended in \cite{tuomov-nlpdhgm-general} by replacing $\iprod{K(x)}{y}_Y$ by a general saddle term $K(x,y)$, which can be applied to nonconvex optimization problems such as $\ell^0$-TV denoising or elliptic Nash equilibrium problems.
    Block-adapted and stochastic variants in the spirit of \cref{thm:gap:accel:stochastic} can be found in \cite{tuomov-nlpdhgm-block}. Finally, a simplified approach using the Bregman divergences\index{divergence, Bregman} of \cref{sec:gap:ergodic:bregman} is presented in \cite{tuomov-firstorder}.
\end{remark}

\chapter{Limiting subdifferentials}
\label{chap:limiting}

While the Clarke subdifferential is a suitable concept for nonsmooth but convex or nonconvex but smooth functionals, it has severe drawbacks for nonsmooth \emph{and} nonconvex functionals: As shown in \cref{lem:clarke:fermat2}, its Fermat principle cannot distinguish minimizers from maximizers. The reason is that the Clarke subdifferential is always convex, which is a direct consequence of its construction \eqref{eq:clarke:def} via polarity with respect to (generalized) directional derivatives. To obtain sharper results for such functionals, it is therefore necessary to construct \emph{nonconvex} subdifferentials directly via a \emph{dual} limiting process. On the other hand, deriving calculus rules for the previous subdifferentials crucially exploited their convexity by applying Hahn--Banach separation theorems, and calculus rules for nonconvex subdifferentials are thus significantly more difficult to obtain.
As in \cref{chap:clarke}, we will assume throughout this chapter that $X$ is a Banach space unless stated otherwise.

\section{Bouligand subdifferentials}

The first definition is motivated by \cref{thm:clarke:gradient}: We \emph{define} a subdifferential as a suitable limit of classical derivatives (without convexification). For $F:X\to \Rbar$, we first define the \term[point!Gateaux]{set of Gateaux points}
\begin{equation*}
    G_F \defeq \setof{x\in X}{F \text{ is Gateaux differentiable at }x} \subset \dom F
\end{equation*}
and then the \term[subdifferential!Bouligand]{Bouligand subdifferential} of $F$ at $x$ as
\begin{equation}\label{eq:limiting:bouligand}
    \partial_B F(x) \defeq\setof{x^*\in X^*}{DF(x_n)\weaktostar x^* \text{ for some } G_F\ni x_n\to x}.
\end{equation}
For $F:\R^N\to \R$ locally Lipschitz, it follows from \cref{thm:clarke:gradient} that $\partial_C F(x) = \conv \partial_B F(x)$.
However, unless $X$ is finite-dimensional and thus Rademacher's theorem is available, it is not clear a priori that the Bouligand subdifferential is nonempty even for $x\in \dom F$.\footnote{Versions of Rademacher's theorem are also available in so-called \term[space!Asplund|infn]{Asplund spaces}; see \cref{rem:epsilon:asplund}. In special cases, it is even possible to give a full characterization; see, e.g., \cite{CCMW:2017}.} Furthermore, the subdifferential does not admit a satisfactory calculus; not even a Fermat principle holds.

\begin{example}\label{ex:limiting:bouligand}
    Let $F:\R\to\R$, $F(x) \defeq |x|$. Then $F$ is differentiable at every $x\neq 0$ with $F'(x) = \sign(x)$. Correspondingly,
    \begin{equation*}
        0\notin \{-1,1\} = \partial_B F(0).
    \end{equation*}
\end{example}

To make this approach work therefore requires a more delicate limiting process. The remainder of this chapter is devoted to one such approach, where we only give an overview and state important results following \cite{Mordukhovich:2006}. The full theory is based on a geometric construction similar to \cref{lem:convex:subdiff_epi} making use of tangent and normal cones (corresponding to generalized directional derivatives and subgradients, respectively) that also allows for differentiation of set-valued mappings. We will develop this theory in \crefrange{chap:cones}{chap:superposition}. For an alternative, more axiomatic, approach to generalized derivatives of nonconvex functionals, we refer to \cite{Penot:2013,ioffe2017variational}.

\section{Fr\'echet subdifferentials}
\label{sec:limiting:frechet}

We begin with the following limiting construction, which combines the characterizations of both the Fréchet derivative and the convex subdifferential. Let $X$ be a Banach space and $F:X\to\Rbar$. The \term[subdifferential!Fréchet]{Fréchet subdifferential} (or \term[subdifferential!regular]{regular subdifferential} or \term[presubdifferential]{presubdifferential}) of $F$ at $x$ is then defined as\footnote{The equivalence of \eqref{eq:limiting:frechet} with the usual definition based on corresponding normal cones follows from, e.g., \cite[Theorem 1.86]{Mordukhovich:2006}.}
\begin{equation}\label{eq:limiting:frechet}
    \partial_F F(x) \defeq \setof{x^*\in X^*}{\liminf_{y\to x} \frac{F(y)-F(x)-\dual{x^*,y-x}_X}{\norm{y-x}_X}\geq 0}.
\end{equation}
Note how this \enquote{localizes} the definition of the convex subdifferential around the point of interest: the numerator need not be nonnegative for all $y$; it suffices if this holds for any $y$ sufficiently close to $x$. By a similar argument as for \cref{thm:convex:fermat}, we thus obtain a Fermat principle for \emph{local} minimizers.
\begin{theorem}\label{thm:limiting:frechet:fermat}
    Let $F:X\to\Rbar$ be proper and $\bar x \in \dom F$ be a local minimizer. Then $0\in \partial_F F(\bar x)$.
\end{theorem}
\begin{proof}
    Let $\bar x\in \dom F$ be a local minimizer. Then there exists an $\epsilon>0$ such that $F(\bar x)\leq F(y)$ for all $y\in \OB(\bar x, \epsilon)$, which is equivalent to
    \begin{equation*}
        \frac{F(y)-F(\bar x)-\dual{0,y-\bar x}_X}{\norm{y-\bar x}_X}\geq 0 \quad\text{for all }y\in \OB(\bar x, \epsilon).
    \end{equation*}
    Now for any strongly convergent sequence $y_n\to \bar x$, we have that $y_n\in \OB(\bar x, \epsilon)$ for $n$ large enough. Taking the $\liminf$ in the above inequality thus yields $0\in \partial_F F(\bar x)$.
\end{proof}

For convex functionals, of course, the numerator is always nonnegative by definition, and the Fréchet subdifferential reduces to the convex subdifferential.
\begin{theorem}\label{thm:limiting:frechet:convex}
    Let $F:X\to\Rbar$ be proper, convex, and lower semicontinuous and $x\in \dom F$. Then $\partial_F F(x) = \partial F(x)$.
\end{theorem}
\begin{proof}
    By definition of the convex subdifferential, any $x^*\in\partial F(x)$ satisfies
    \begin{equation*}
        F(y)-F(x)-\dual{x^*,y-x}_X \geq 0 \quad\text{for all } y\in X.
    \end{equation*}
    Dividing by $\norm{x-y}_X>0$ for $y\neq x$ and taking the $\liminf$ as $y\to x$ thus yields $x^*\in \partial_F F(x)$.

    Conversely, let $x^*\in \partial_F F(x)$ and $h\in X\setminus\{0\}$ be arbitrary. Then for an $\delta>0$, there exists an $\epsilon >0$ such that
    \begin{equation*}
        \frac{F(x+th)-F(x)-\dual{x^*,th}_X}{t\norm{h}_X} \geq -\delta \quad\text{for all }t\in (0,\epsilon).
    \end{equation*}
    Multiplying by $\norm{h}_X>0$ and letting $t\to 0$, we obtain from \cref{lem:convex:direct} that
    \begin{equation}
        \dual{x^*,h}_X \leq \frac{F(x+th)-F(x)}{t} + \delta\norm{h}_X \to F'(x; h) + \delta\norm{h}_X.
    \end{equation}
    Since $\delta>0$ was arbitrary, this implies by \cref{lem:convex:equiv} that $x^*\in \partial F(x)$.
\end{proof}

Similarly, for Fréchet differentiable functionals, the limit in \eqref{eq:limiting:frechet} is zero for all sequences.
\begin{theorem}\label{lem:limiting:frechet:frechet}
    Let $F:X\to\R$ be Fréchet differentiable at $x\in X$. Then $\partial_F F(x) = \{F'(x)\}$.
\end{theorem}
\begin{proof}
    The definition of the Fréchet derivative immediately yields
    \begin{equation*}
        \lim_{y\to x} \frac{F(y)-F(x)-\dual{F'(x),y-x}_X}{\norm{x-y}_X} = \lim_{\norm{h}_X\to 0} \frac{F(x+h)-F(x)-F'(x)h}{\norm{h}_X} = 0
    \end{equation*}
    and hence $F'(x)\in \partial_F F(x)$.

    Conversely, let $x^*\in \partial_F F(x)$ and let again $h\in X\setminus\{0\}$ be arbitrary. As in the proof of \cref{thm:limiting:frechet:convex}, we then obtain that
    \begin{equation}
        \dual{x^*,h}_X \leq F'(x; h) = \dual{F'(x),h}_X.
    \end{equation}
    Applying the same argument to $-h$ then yields $\dual{x^*,h}_X = \dual{F'(x),h}_X$ for all $h\in X$, i.e., $x^* = F'(x)$.
\end{proof}

For nonsmooth and nonconvex functionals, the Fréchet subdifferential can be strictly smaller than the Clarke subdifferential.
\begin{example}\label{ex:limiting:frechet}
    Consider $F:\R\to\R$, $F(x) \defeq -|x|$. For any $x\neq 0$, it follows from \cref{lem:limiting:frechet:frechet} that $\partial_F F(x) = \{-\sign (x)\}$. But for $x=0$ and arbitrary $x^*\in \R$, we have that
    \begin{equation*}
        \liminf_{y\to 0} \frac{F(y)-F(0)-\dual{x^*,y-0}}{|y-0|} = \liminf_{y\to 0} (-1-x^*\cdot \sign(y)) = -1-|x^*| < 0
    \end{equation*}
    and hence that
    \begin{equation*}
        \partial_F F(0) = \emptyset \subsetneq [-1,1] = \partial_C F(0).
    \end{equation*}
\end{example}
Note that $0\in \dom F$ in this example. Although the Fréchet subdifferential does not pick up a maximizer in contrast to the Clarke subdifferential, the fact that $\partial_F F(x)$ can be empty even for $x\in \dom F$ is a problem when trying to derive calculus rules that hold with equality. In fact, as \cref{ex:limiting:frechet} shows, the Fréchet subdifferential fails to be outer semicontinuous, which is also not desirable. This leads to the next and final definition.

\section{Mordukhovich subdifferentials}
\label{sec:limiting:mordukhovich}

Let $X$ be a reflexive Banach space and $F:X\to \Rbar$. The \term[subdifferential!Mordukhovich]{Mordukhovich subdifferential} (or \term[subdifferential!basic]{basic subdifferential} or \term[subdifferential!limiting]{limiting subdifferential}) of $F$ at $x\in \dom F$ is then defined as the strong-to-weak$^*$ outer \emph{sequential} closure of $\partial_F F(x)$, i.e.,%
\footnote{The equivalence of this definition with the original geometric definition -- which holds in \emph{reflexive} Banach spaces -- follows from \cite[Theorem 2.34]{Mordukhovich:2006}.}
\begin{equation}\label{eq:limiting:mordukhovich}
    \begin{aligned}[t]
        \partial_M F(x) &\defeq \weakstarlimsup_{y\to x}\partial_F F(y) \\
        &= \setof{x^*\in X^*}{x_n^* \weaktostar x^* \text{ for some } x_n^*\in \partial_F F(x_n) \text{ with }x_n\to x},
    \end{aligned}
\end{equation}
which can be seen as a generalization of the definition \eqref{eq:limiting:bouligand} of the Bouligand subdifferential. Note that in contrast to \eqref{eq:limiting:bouligand}, this definition includes the constant sequence $x_n^* \equiv x^*$ even at nondifferentiable points, which makes this a more useful concept in general. This also implies that $\partial_F F(x) \subset \partial_M F(x)$ for any $F$, and \cref{thm:limiting:frechet:fermat} immediately yields a Fermat principle.
\begin{corollary}\label{thm:limiting:mordukhovich}
    Let $F:X\to\Rbar$ be proper and $\bar x \in \dom F$ be a local minimizer. Then $0\in \partial_M F(\bar x)$.
\end{corollary}
As for the Fréchet subdifferential, maximizers do not satisfy the Fermat principle.
\begin{example}
    Consider again $F:\R\to\R$, $F(x) \defeq -|x|$. Using \cref{ex:limiting:frechet}, we directly obtain from \eqref{eq:limiting:mordukhovich} that $\partial_M F(0) = \{-1,1\} = \partial_B F(0)$.
\end{example}

Since the convex subdifferential is strong-to-weak$^*$ outer semicontinuous, the Mordukhovich subdifferential reduces to the convex subdifferential as well.
\begin{theorem}
    Let $X$ be a reflexive Banach space, $F:X\to\Rbar$ be proper, convex, and lower semicontinuous, and $x\in \dom F$. Then $\partial_M F(x) = \partial F(x)$.
\end{theorem}
\begin{proof}
    From \cref{thm:limiting:frechet:convex}, it follows that $\partial F(x) = \partial_F F(x) \subset \partial_M F(x)$. Let therefore $x^*\in \partial_M F(x)$ be arbitrary. Then by definition there exists a sequence $\{x_n^*\}_{n\in\N}\subset X^*$ with $x_n^*\weaktostar x^*$ and $x_n^*\in \partial_F F(x_n)=\partial F(x_n)$ for $x_n\to x$. From \cref{thm:monoton:subdiff,cor:monoton:closed}, it then follows that $x^*\in \partial F(x)$ as well.
\end{proof}
Note that this does not imply that the Mordukhovich subdifferential itself is outer semicontinuous; see \cite[Example 5.6]{mehlitzwachsmuth2019decomposable} for a counterexample.

A similar result holds for \emph{continuously} differentiable functionals.
\begin{theorem}
    Let $X$ be a reflexive Banach space and $F:X\to\Rbar$ be continuously differentiable at  $x\in X$. Then $\partial_M F(x) = \{F'(x)\}$.
\end{theorem}
\begin{proof}
    From \cref{thm:limiting:frechet:convex}, it follows that $\{F'(x)\} = \partial_F F(x) \subset \partial_M F(x)$. Let therefore $x^*\in \partial_M F(x)$ be arbitrary. Then by definition there exists a sequence $\{x_n^*\}_{n\in\N}\subset X^*$ with $x_n^*\weaktostar x^*$ and $x_n^*\in \partial_F F(x_n) =\{F'(x_n)\}$ for $x_n\to x$. The continuity of $F'$ then immediately implies that $F'(x_n)\to F'(x)$, and since strong limits are also weak-$*$ limits, we obtain $x^*=F'(x)$.
\end{proof}
The same function as in \cref{ex:clarke:frechet} shows that this equality does not hold if $F$ is merely Fréchet differentiable.

We also have the following relation to Clarke subdifferentials, which should be compared to \cref{thm:clarke:gradient}. We will give a proof in a more restricted setting in \cref{chap:graphical}, cf.~\cref{cor:graphical:clarke-weakstar-convex}.

\begin{theorem}[\protect{\cite[Theorem 3.57]{Mordukhovich:2006}}]\label{thm:limiting:clarke}
    Let $X$ be a reflexive Banach space and $F:X\to\R$ be locally Lipschitz continuous around $x\in X$. Then $\partial_C F(x) = \mathrm{cl}^* \conv \partial_M F(x)$, where $\mathrm{cl}^*A $ stands for the weak-$*$ closure of the set $A\subset X^*$.%
    \footnote{Of course, in reflexive Banach spaces the weak-$*$ closure coincides with the weak closure. The statement holds more general in \term[space!Asplund|infn]{Asplund spaces} which include some nonreflexive Banach spaces.}
\end{theorem}
The following example illustrates that the Mordukhovich subdifferential can be nonconvex.
\begin{example}
    Let $F:\R^2\to\R$, $F(x_1,x_2) = |x_1| - |x_2|$. Since $F$ is continuously differentiable for any $(x_1,x_2)$ where $x_1,x_2\neq 0$ with
    \begin{equation*}
        \nabla F(x_1,x_2)\in \{(1,1),(-1,1),(1,-1),(-1,-1)\},
    \end{equation*}
    we obtain from \eqref{eq:limiting:frechet} that
    \begin{equation*}
        \partial_F F(x_1,x_2) =
        \begin{cases}
            \{(1,-1)\} & \text{if } x_1>0,x_2>0,\\
            \{(-1,-1)\} & \text{if } x_1<0,x_2>0,\\
            \{(-1,1)\} & \text{if } x_1<0,x_2<0,\\
            \{(1,1)\} & \text{if } x_1>0,x_2<0,\\
            \setof{(t,-1)}{t\in [-1,1]} & \text{if } x_1=0,x_2>0,\\
            \setof{(t,1)}{t\in [-1,1]} & \text{if } x_1=0,x_2<0,\\
            \emptyset & \text{if } x_2=0.
        \end{cases}
    \end{equation*}
    In particular, $\partial_F F(0,0) = \emptyset$. However, from \eqref{eq:limiting:mordukhovich} it follows that
    \begin{equation*}
        \partial_M F(0,0) = \setof{(t,-1)}{t\in [-1,1]}\cup \setof{(t,1)}{t\in [-1,1]}.
    \end{equation*}
    In particular, $0\notin \partial_M F(0,0)$.
    On the other hand, \cref{thm:limiting:clarke} then yields that
    \begin{equation}
        \partial_C F(0,0) = \setof{(t,s)}{t,s\in[-1,1]} = [-1,1]^2
    \end{equation}
    and hence $0\in \partial_C F(0,0)$.
    (Note that $F$ attains neither a minimum nor a maximum on $\R^2$, while $(0,0)$ is a nonsmooth saddle-point.)
\end{example}

In contrast to the Bouligand subdifferential, the Mordukhovich subdifferential admits a satisfying calculus, although the assumptions are understandably more restrictive than in the convex setting. The first rule follows as always straight from the definition.

\begin{theorem}\label{thm:limiting:scalar}
    Let $X$ be a reflexive Banach space and $F:X\to\Rbar$. Then for any $\lambda\geq 0$ and $x\in X$,
    \begin{equation*}
        \partial_M (\lambda F)(x) = \lambda \partial_M F(x).
    \end{equation*}
\end{theorem}

Full calculus in infinite-dimensional spaces holds only for a rather small class of mappings.

\begin{theorem}[{\protect\cite[Proposition 1.107]{Mordukhovich:2006}}]\label{thm:limiting:sum}
    Let $X$ be a reflexive Banach space, $F:X\to\R$ be continuously differentiable, and $G:X\to\Rbar$ be arbitrary. Then for any $x\in \dom G$,
    \begin{equation*}
        \partial_M (F+G)(x) = \{F'(x)\} + \partial_M G(x).
    \end{equation*}
\end{theorem}
While the previous two theorems also hold for the Fréchet subdifferential (the latter even for merely Fréchet differentiable $F$), the following chain rule is only valid for the Mordukhovich subdifferential. Compared to \cref{thm:clarke:chain}, it also allows for the outer functional to be extended-real valued.
\begin{theorem}[{\protect\cite[Proposition 1.112]{Mordukhovich:2006}}]\label{thm:limiting:chain}
    Let $X$ be a reflexive Banach space, $F:X\to Y$ be continuously differentiable, and $G:Y\to \Rbar$ be arbitrary. Then for any $x\in X$ with $F(x)\in \dom G$ and $F'(x):X\to Y$ surjective,
    \begin{equation*}
        \partial_M (G\circ F)(x) = F'(x)^*\partial_M G(F(x)).
    \end{equation*}
\end{theorem}
More general calculus rules require $X$ to be a reflexive Banach\footnote{or Asplund} space as well as additional, nontrivial, assumptions on $F$ and $G$; see, e.g., \cite[Theorem 3.36, Theorem 3.41]{Mordukhovich:2006}.

We will illustrate how to prove the above calculus results and more in \cref{sec:graphical:subdiff,sec:colimiting:subdiff}, after studying the differentiation of set-valued mappings.

\chapter{\texorpdfstring{$\eps$}{ε}-subdifferentials and approximate Fermat principles}
\label{chap:epsilon}

We now study an approximate variant of the Fréchet subdifferential of \cref{sec:limiting:frechet} as well as related approximate Fermat principles; these will be needed in \cref{chap:cones} to study limiting tangent and normal cones.

\section{\texorpdfstring{$\eps$}{ε}-subdifferentials}
\label{sec:epsilon:epsilon}

Just like the $\eps$-minimizers in \cref{sec:variation:ekeland}, it can be useful to consider \enquote{relaxed} $\eps$-subdifferenti\-als. In particular, it is possible to derive \emph{exact} calculus rules for these relaxed subdifferentials, which can lead to tighter results than inclusions for the corresponding exact subdifferentials (in particular, for the Fréchet subdifferential). We will make use of this in \cref{chap:regularity}.

Similarly to the Fréchet subdifferential \eqref{eq:limiting:frechet}, we thus define for $F:X\to \Rbar$ the \term[subdifferential!$\epsilon$-]{$\eps$(-Fréchet)-subdifferential} by
\begin{equation}
    \label{eq:epsilon:frechet-epsilon}
    \subdiff_\epsilon F(x) \defeq \setof{x^*\in X^*}{\liminf_{y\to x} \frac{F(y)-F(x)-\dual{x^*,y-x}_X}{\norm{y-x}_X}\geq -\epsilon},
\end{equation}
where $\subdiff_0 F=\subdiff_F F$.
The following lemma provides further insight into the $\epsilon$-subdifferential.

\begin{lemma}
    \label{eq:epsilon:alt-formulations}
    Let $F: X \to \Rbar$ on a Banach space $X$, and $\epsilon \ge 0$. Then the following are equivalent:
    \begin{enumerate}
        \item $x^* \in \subdiff_\epsilon F(x)$;
        \item $x^* \in \subdiff_F[F+\epsilon\norm{\freevar - x}_X](x)$;
        \item $0 \in \subdiff_F[F+\epsilon\norm{\freevar - x}_X-\dualprod{x^*}{\freevar-x}](x)$.
    \end{enumerate}
\end{lemma}

\begin{proof}
    Each of the alternatives is by \eqref{eq:epsilon:frechet-epsilon} and \eqref{eq:limiting:frechet} equivalent to
    \begin{equation*}
        \liminf_{y\to x} \frac{\epsilon\norm{y-x}_X+F(y)-F(x)-\dualprod{x^*}{y-x}_X}{\norm{y-x}_X} \ge 0.
        \qedhere
    \end{equation*}
\end{proof}

We have the following \enquote{fuzzy} $\eps$-sum rule.

\begin{lemma}
    \label{lemma:epsilon:sumrule}
    Let $X$ be a Banach space, $G: X \to \Rbar$, and $F: X \to \R$ be convex with $\subdiff F(x) \subset \B(\optx^*, \epsilon)$ for some $\epsilon \ge 0$ and $\optx^* \in X^*$.
    Then for all $\delta \ge 0$,
    \begin{equation*}
        \subdiff_\delta G(x) + \subdiff F(x) \subset \subdiff_\delta[G+F](x) \subset \subdiff_{\epsilon+\delta} G(x) + \{\optx^*\}.
    \end{equation*}
    In particular, if $\optx^* \in \subdiff F(x)$, then
    \begin{equation*}
        \subdiff_\delta G(x) + \subdiff F(x) \subset \subdiff_\delta[G+F](x) \subset \subdiff_{\epsilon+\delta} G(x) + \subdiff F(x).
    \end{equation*}
\end{lemma}
\begin{proof}
    We start with the first inclusion. Let $\alt x^* \in \subdiff F(x)$ and $x^* \in \subdiff_\delta G(x)$.
    Then the definitions \eqref{eq:convex:def_subdiff} and \eqref{eq:epsilon:frechet-epsilon}, respectively, imply that
    \begin{multline*}
        \liminf_{y \to x} \frac{G(y)-G(x)+F(y)-F(x)-\dualprod{x^*+\alt x^*}{y-x}_X}{\norm{y-x}_X}
        \\
        \ge
        \liminf_{y \to x} \frac{G(y)-G(x)-\dualprod{x^*}{y-x}_X}{\norm{y-x}_X}
        \ge -\delta,
    \end{multline*}
    i.e., $x^*+\alt x^* \in \subdiff_\delta[G+F](x)$.

    To prove the second inclusion, let $x^* \in \subdiff_\delta[G+F](x)$ and $h \in X$ with $\norm{h}_X=1$.
    Then \eqref{eq:epsilon:frechet-epsilon} implies that for all $t_n \downto 0$ and $h_n \to h$,
    \begin{equation}
        \label{eq:epsilon:sumrule:2}
        \liminf_{n \to \infty}~ \frac{F(x+t_nh_n)-F(x)+G(x+t_nh_n)-G(x) - t_n \dualprod{x^*}{h_n}_X}{t_n}
        \ge -\delta.
    \end{equation}
    Since $F$ is directionally differentiable by \cref{lem:convex:direct} and locally Lipschitz around $x \in \interior (\dom F)=X$ by \cref{thm:convex:cont} with Lipschitz constant $L>0$, we have
    \begin{equation*}
        \lim_{n \to \infty} \frac{F(x+t_nh_n)-F(x)}{t_n}
        \le
        \lim_{n \to \infty}\left( \frac{F(x+t_nh)-F(x)}{t_n}+L \norm{h_n-h}_X
        \right)
        =
        F'(x; h).
    \end{equation*}
    Let now $\rho>0$ be arbitrary. Then by \cref{lem:convex:equiv,thm:clarke:convex,cor:clarke:support-dir} there exists an $x^*_{h,\rho} \in \subdiff F(x)$ such that $F'(x; h) \le \dualprod{x^*_{h,\rho}}{h}_X + \rho$.
    Therefore
    \begin{equation*}
        \begin{aligned}[t]
            \lim_{n \to \infty}
            \frac{F(x+t_nh_n)-F(x)-t_n\dualprod{\optx^*}{h_n}_X}{t_n}
            &
            \le
            F'(x; h) - \dualprod{\optx^*}{h}_X
            \\
            &
            \le \dualprod{x^*_{h,\rho}-\optx^*}{h}_X
            + \rho
            \\
            &
            \le \epsilon + \rho,
        \end{aligned}
    \end{equation*}
    where we have used that $\subdiff F(x) \subset \B(\optx^*, \epsilon)$ and $\norm{h}_X=1$ in the last inequality.
    Since $\rho>0$ was arbitrary, the characterization \eqref{eq:epsilon:sumrule:2} now implies
    \begin{equation*}
        \liminf_{n \to \infty}~ \frac{G(x+t_n h_n)-G(x) - t_n \dualprod{x^*-\opt x^*}{h_n}_X}{t_n}
        \ge -(\delta+\epsilon).
    \end{equation*}
    Since $y_n \defeq x+t_n h_n \to x$ was arbitrary, this proves $x^*-\opt x^* \in \subdiff_{\epsilon+\delta} G(x)$, i.e., $\subdiff_\delta[G+F](x) \subset \subdiff_{\epsilon+\delta} G(x) + \{\optx^*\}$.
\end{proof}

The following is now immediate from \cref{thm:convex:gateaux}, since we are allowed to take $\eps=0$ if $\partial F(x)$ is a singleton.

\begin{corollary}
    \label{cor:epsilon:sumrule:gateaux}
    Let $X$ be a Banach space, $G: X \to \Rbar$, and $F: X \to \R$ be convex and Gateaux differentiable at $x \in X$. Then for every $\delta \ge 0$,
    \begin{equation*}
        \subdiff_\delta[G+F](x) = \subdiff_\delta G(x) + \{DF(x)\}.
    \end{equation*}
    In particular,
    \begin{equation*}
        \subdiff_F[G+F](x) = \subdiff_F G(x) + \{DF(x)\}.
    \end{equation*}
\end{corollary}

\section{Smooth spaces}\label{sec:epsilon:smooth}

For the remaining results in this chapter, we need additional assumptions on the normed vector space $X$. In particular, we need to assume that the norm is Gateaux differentiable on $X \setminus \{0\}$; we call such spaces \term[space!Gateaux smooth]{Gateaux smooth}.

Recalling from \cref{chap:smoothness} the duality between differentiability and convexity, it is not surprising that this property can be related to the convexity of the dual norm.
Here we need the following property: a normed vector space $X$ is called \term[space!locally uniformly convex]{locally uniformly convex} if for any $x\in X$ with $\norm{x}_X=1$ and all $\eps\in (0,2]$ there exists a $\delta(\eps,x)>0$ such that
\begin{equation}\label{eq:smoothness:uniformly-convex}
    \norm{\tfrac12(x+y)}_X \leq 1 - \delta(\eps,x) \quad\text{for all }y\in X \text{ with }\norm{y}_X=1 \text{ and } \norm{x-y}_X\geq \eps.
\end{equation}
\begin{lemma}\label{lem:epsilon:gateaux-convex}
    Let $X$ be a Banach space and $X^*$ be locally uniformly convex. Then $X$ is Gateaux smooth.
\end{lemma}
\begin{proof}
    Let $x\in X\setminus\{0\}$ be given. Since norms are convex, it suffices by \cref{thm:convex:singleton} to show that $\partial\norm{\freevar}_X(x)$ is a singleton. Let therefore $x_1^*,x_2^*\in \partial\norm{\freevar}_X(x)$, i.e., satisfying by \cref{thm:subdifferential:norm}
    \begin{equation*}
        \norm{x_1^*}_{X^*} = \norm{x_2^*}_{X^*} = 1,\qquad
        \dual{x_1^*,x}_{X} = \dual{x_2^*,x}_{X} = \norm{x}_X.
    \end{equation*}
    This implies that
    \begin{equation*}
        2 = \frac1{\norm{x}_X}\left( \dual{x_1^*,x}_{X} + \dual{x_2^*,x}_{X} \right)
        =
        \dual{x_1^*+x_2^*,\tfrac{x}{\norm{x}_X}}_{X}\leq \norm{x_1^*+x_2^*}_{X^*}
    \end{equation*}
    by \eqref{eq:functan:cs_banach} and hence that $\norm{\frac12(x_1^*+x_2^*)}_{X^*}\geq 1$. Since $X^*$ is locally uniformly convex, this is only possible if $x_1^* = x_2^*$, as otherwise we could choose for $\eps\defeq \norm{x_1^*-x_2^*}_{X^*} \in (0, 2]$ a $\delta(\eps,x_1^*)>0$ such that $\norm{\frac12(x_1^*+x_2^*)}_{X^*}\leq 1-\delta(\eps,x_1^*)<1$.
\end{proof}
\begin{remark}%
  \label{rem:epsilon:radonriesz}
  In fact, if $X$ is additionally reflexive, the norm is even continuously (Fréchet) differentiable; see \cite[Proposition 4.7.10]{Schirotzek:2007}. We will not need this stronger property, however. In addition, locally uniformly convex spaces always have the Radon--Riesz property; see \cite[Lemma 4.7.9]{Schirotzek:2007}.
\end{remark}

\begin{example}\label{ex:epsilon:gateaux}
    The following spaces are locally uniformly convex:
    \begin{enumerate}
        \item\label{item:epsilon:gateaux:hilbert}
            $X$ a Hilbert space. This follows from the \term[identity!parallelogram]{parallelogram identity}
            \begin{equation*}
                \norm{\tfrac12(x+y)}_X^2 = \frac12\norm{x}_X^2 + \frac12\norm{y}_X^2 - \frac14 \norm{x-y}_X^2 \qquad\text{for all }x,y\in X,
            \end{equation*}
            which in fact characterizes precisely those norms that are induced by an inner product.
            This identity immediately yields for all $\eps>0$ and all $x,y\in X$ satisfying $\norm{x-y}_X\geq \eps$ that
            \begin{equation*}
                \norm{\tfrac12(x+y)}_X^2 \leq 1 -\frac{\eps^2}4 \leq \left(1 -\frac{\eps^2}8\right)^2,
            \end{equation*}
            which in particular verifies \eqref{eq:smoothness:uniformly-convex} with $\delta\defeq \frac{\eps^2}{8}$.

        \item\label{item:epsilon:gateaux:lebesgue} $X=L^p(\Omega)$ for $p\in (2,\infty)$. This follows from the algebraic inequality
            \begin{equation*}
                |a+b|^p + |a-b|^p \leq 2^{p-1}(|a|^p + |b|^p) \qquad\text{for all }a,b\in \R,
            \end{equation*}
            see \cite[Lemma II.4.1]{Cioranescu}. This implies that
            \begin{equation*}
                \norm{\tfrac12(u+v)}_{L^p(\Omega)}^p \leq \frac12\norm{u}_{L^p(\Omega)}^p + \frac12\norm{v}_{L^p(\Omega)}^p - \frac1{2^p} \norm{u-v}_{L^p(\Omega)}^p
                \quad\text{for all }u,v\in L^p(\Omega).
            \end{equation*}
            We can now argue exactly as in case (i).

        \item $X=L^p(\Omega)$ for $p\in (1, 2)$. This follows from the algebraic inequality
            \begin{equation*}
                |a+b|^p + |a-b|^p \leq 2(|a|^p + |b|^p)^{p/(p-1)} \qquad\text{for all }a,b\in \R,
            \end{equation*}
            see \cite[Lemma II.4.1]{Cioranescu}, implying a similar inequality for the $L^p(\Omega)$ norms from which the claim follows as for \ref{item:epsilon:gateaux:hilbert} and \ref{item:epsilon:gateaux:lebesgue}.
    \end{enumerate}
    Hence every Hilbert space (by identifying $X$ with $X^*$) and every $L^p(\Omega)$ for $p\in(1,\infty)$ (identifying $L^p(\Omega)$ with $L^q(\Omega)$, $q=\frac{p}{p-1}\in(1,\infty)$) is Gateaux smooth.
\end{example}

In fact, the celebrated Lindenstrauss and Trojanski \term[theorem!renorming]{renorming theorems} show that every reflexive Banach space admits an equivalent norm such that the space (with that norm) becomes locally uniformly convex; see \cite[Theorem III.2.10]{Cioranescu}. (Of course, even though that means that the dual space of the \emph{renormed space} is Gateaux smooth, this does not imply anything about the differentiability of the original norm, as the obvious example of $\R^N$ endowed with the $1$- or the $\infty$-norm shows.)
For many more details on smooth and uniformly convex spaces, see \cite{fabian2001differentiability,Schirotzek:2007,Cioranescu}.

Note that even in Gateaux smooth spaces, the norm will not be differentiable at $x=0$. But this can be addressed by considering $\norm{x}_X^p$ for $p>1$; for later use, we state this for $p=2$.
\begin{lemma}\label{lem:epsilon:squared-norm}
    Let $X$ be a Gateaux smooth Banach space and $F(x) = \norm{x}_X^2$. Then $F$ is Gateaux differentiable at any $x\in X$ with
    \begin{equation*}
        DF(x) = 2\norm{x}_X x^* \qquad\text{for any } x^*\in X^* \text{ with } \norm{x^*}_{X^*}=1 \text{ and } \dual{x^*,x}_X = \norm{x}_X.
    \end{equation*}
\end{lemma}
\begin{proof}
    Since norms are convex, we can apply \cref{thm:convex:increasing-post,thm:subdifferential:norm} to obtain that
    \begin{equation*}
        \subdiff F(x) = \{  2\norm{x}_X x^* \mid x^*\in X^* \text{ with } \norm{x^*}_{X^*}=1 \text{ and } \dual{x^*,x}_X = \norm{x}_X\}
        \quad (x \in X).
    \end{equation*}
    At any $x \ne 0$, this set is a singleton by \cref{thm:convex:gateaux} and the assumption that $X$ is Gateaux smooth. Clearly also $\subdiff F(0)=\{0\}$, and hence the claim follows from \cref{thm:convex:singleton}.
\end{proof}

\begin{remark}[Asplund spaces]
    \label{rem:epsilon:asplund}
    \term[space!Asplund]{Asplund spaces} are, by (one equivalent) definition, those Banach spaces where every continuous, convex, real-valued function is Fréchet differentiable on a dense set.
    (This is a limited version of Rademacher's \cref{thm:rademacher} in $\R^N$; see, e.g., the seminal work \cite{Preiss:1990}.)
    Importantly, reflexive Banach spaces are Asplund.
    We refer to \cite{yost1993asplund} for an introduction to Asplund spaces and to \cite{Lindenstrauss:2021} for an in-depth treatment of the Fréchet differentiability of Lipschitz continuous mappings on such spaces.

    The norm of an Asplund space is thus differentiable on a dense set $D$.
    It was shown in \cite{ekeland1976generic} that perturbed optimization problems on Asplund spaces have solutions on a dense set of perturbation parameters and that the objective function is differentiable at such a solution. If we worked in the following sections with perturbed optimization problems and applied such an existence result instead of the Ekeland or the Borwein--Preiss variational principles (\cref{thm:variation:ekeland} or \cref{thm:variation:borweinpreiss}, respectively), we would be able to extend the following results to Asplund spaces.
\end{remark}

\section{Fuzzy Fermat principles}

The following result generalizes the Fermat principle of \cref{thm:limiting:frechet:fermat} to sums of two functions in a \enquote{fuzzy} fashion. We will use it to show a fuzzy containment formula for $\epsilon$-subdifferentials. Its generalizations to more than two functions can also be used to derive more advanced fuzzy sum rules than \cref{lemma:epsilon:sumrule}. Our focus is, however, on exact calculus, so we will not be developing such generalizations.

\begin{lemma}[fuzzy Fermat principle]\index{principle!Fermat!fuzzy}
    \label{lemma:epsilon:fuzzy-fermat}
    Let $X$ be a Gateaux smooth Banach space and $F, G: X \to \Rbar$.
    If $F+G$ attains a local minimum at a point $\optx \in X$ where $F$ is lower semicontinuous and $G$ is locally Lipschitz, then for any $\delta,\mu>0$ we have
    \begin{equation*}
        0 \in \Union_{x, y \in \B(\optx, \delta)} \left( \subdiff_F F(x) + \subdiff_F G(y) \right) + \mu\B_{X^*}.
    \end{equation*}
\end{lemma}

\begin{proof}
    Let $\rho,\alpha>0$ be arbitrary.
    The idea is to separate the two nonsmooth functions $F$ and $G$, and hence be able to use the exact sum rule of \cref{cor:epsilon:sumrule:gateaux}, by locally relaxing the problem $\min_{x \in X} (F+G)$ to
    \begin{equation*}
        \inf_{x, y \in X}~ J_\alpha(x, y) \defeq F(x) + G(y) + \alpha\norm{x-y}_X^2 + \norm{x-\opt x}_X^2 + \delta_{\B(\optx,\rho)^2}(x, y).
    \end{equation*}
    We take $\rho>0$ small enough that $\optx$ minimizes $F+G$ within $\B(\optx, \rho)$, and both $F \ge F(\optx)-1$ and $G \ge G(\opty)-1$ on $\B(\optx, \rho)$.
    The first requirement is possible by the assumption of $F+G$ attaining its local minimum at $\optx$, while the latter follows from the lower semicontinuity of $F$ and the local Lipschitz continuity of $G$.
    In the following, we denote by $L$ the Lipschitz constant of $G$ on $\B(\optx,\rho)$. It follows that $J_\alpha(x,y) \ge F(\optx) + G(\optx) - 2$ for all $(x,y)\in \B(\optx,\rho)^2=\dom J_\alpha$, and hence $J_\alpha$ is bounded from below.

    We study the approximate solutions of the relaxed problem in several steps.

    \emph{Step 1: constrained infimal values converge to $J(\optx, \optx)$.}
    Let $x_\alpha, y_\alpha \in \B(\optx, \rho)$ be such that
    \begin{equation}
        \label{eq:epsilon:fuzzy-fermat:value-bound}
        J_\alpha(x_\alpha, y_\alpha) < j_\alpha + \inv\alpha
        \quad\text{where}\quad
        j_\alpha \defeq \inf_{x, y \in X} J_\alpha(x, y).
    \end{equation}
    We show that
    \begin{equation*}
        J_\alpha(\optx, \optx) < j_\alpha + \epsilon_\alpha
        \quad\text{for}\quad
        \epsilon_\alpha \defeq L\sqrt{\frac{\inv\alpha+2}{\alpha}} + \inv\alpha.
    \end{equation*}
    To start with, we have
    \begin{equation*}
        \begin{aligned}[t]
            F(\optx)+G(\optx) + \inv\alpha
            &
            =
            J(\optx, \optx) + \inv\alpha
            \\
            &
            \ge
            j_\alpha + \inv\alpha
            \\
            &
            > J_\alpha(x_\alpha, y_\alpha)
            \\
            &
            =  F(x_\alpha) + G(y_\alpha) + \alpha\norm{x_\alpha-y_\alpha}_X^2 + \norm{x_\alpha-\optx}_X^2
            \\
            &
            \ge F(\optx) + G(\optx) + \alpha\norm{x_\alpha-y_\alpha}_X^2 + \norm{x_\alpha-\optx}_X^2 - 2.
        \end{aligned}
    \end{equation*}
    This implies that $\norm{x_\alpha-y_\alpha}_X < \sqrt{\frac{\inv\alpha+2}{\alpha}}$.
    Since $\optx$ minimizes $F+G$ within $\B(\optx, \rho)$, we obtain the bound \eqref{eq:epsilon:fuzzy-fermat:value-bound} through
    \begin{equation*}
        \begin{aligned}[t]
            J_\alpha(\optx, \optx)
            &
            =
            F(\optx)+G(\optx)
            \\
            &
            \le F(x_\alpha) + G(x_\alpha)
            \\
            &
            \le F(x_\alpha) + G(y_\alpha) +  L\norm{x_\alpha-y_\alpha}_X.
            \\
            &
            \le J(x_\alpha, y_\alpha) + L\norm{x_\alpha-y_\alpha}_X.
            \\
            &
            < j_\alpha + \epsilon_\alpha.
        \end{aligned}
    \end{equation*}

    \emph{Step 2: exact unconstrained minimizers exist for a perturbed problem.}
    By \eqref{eq:epsilon:fuzzy-fermat:value-bound}, we can apply the Borwein--Preiss variational principle (\cref{thm:variation:borweinpreiss}) for any $\lambda,\alpha>0$, small enough $\rho>0$ (all to be fixed later), and $p=2$ to obtain a sequence $\{\mu_n\}_{n\geq 0}$ of nonnegative weights summing to $1$ and a sequence $\{(x_n, y_n)\}_{n \ge 0} \subset X\in X$ with $(x_0,y_0)=(\optx,\optx)$ converging strongly to some $(\realoptx_\alpha, \realopty_\alpha)\in X\times X$ (endowed with the euclidean product norm) such that
    \begin{enumerate}
        \item\label{item:epsilon:fuzzy-fermat:borweinpreiss:1}
            $\norm{x_n-\realoptx_\alpha}_X^2 + \norm{y_n-\realopty_\alpha}_X^2 \le \lambda^2$ for all $n \geq 0$ (in particular, $\norm{\optx-\realoptx_\alpha}_X \leq \lambda$);
        \item
            the function
            \begin{equation*}
                H_\alpha(x, y) \defeq J_\alpha(x, y) + \frac{\eps_\alpha}{\lambda^2}\sum_{n=0}^\infty \mu_n\left(\norm{x-x_n}^2+\norm{y-y_n}^2\right)
            \end{equation*}
            attains its global minimum at $(\realoptx_\alpha, \realopty_\alpha)$.
    \end{enumerate}
    Note that since $J_\alpha$ includes the constraint $(x,y)\in \B(\optx, \rho)^2$, we have $(\realoptx_\alpha, \realopty_\alpha) \in \B(\optx, \rho)^2$.
    In fact, by taking $\lambda \in (0, \rho)$, it follows from \cref{item:epsilon:fuzzy-fermat:borweinpreiss:1} and the convergence $(x_n, y_n) \to (\realoptx_\alpha, \realopty_\alpha)$ that the minimizer $(\realoptx_\alpha, \realopty_\alpha) \in  \B(\optx, \lambda)^2 \subset \interior \B(\optx, \rho)^2$ is unconstrained.

    \emph{Step 3: the perturbed minimizers satisfy the claim for large $\alpha$ and small $\lambda$.}
    Setting $\Psi_y(x) \defeq \norm{x-y}_X^2$, it follows from \cref{lem:epsilon:squared-norm} that $\Psi_y$ is Gateaux differentiable for any $y\in X$ with $D\Psi_y(x) \in 2\norm{x-y}_X \B_{X^*}$.
    Furthermore, since $(\realoptx_\alpha,\realopty_\alpha)\in \interior \B(\optx,\rho)^2$, we have $\subdiff \delta_{\B(\optx,\rho)^2}(\realoptx_\alpha,\realopty_\alpha)=(0,0)$.
    Hence the only nonsmooth component of $H_\alpha$ at $(\realoptx_\alpha, \realopty_\alpha)$ is $(x, y) \mapsto F(x)+G(y)$. We can thus apply \cref{thm:limiting:frechet:fermat,cor:epsilon:sumrule:gateaux} to obtain
    \begin{equation*}
        0  \in \subdiff_F H_\alpha(\realoptx_\alpha, \realopty_\alpha)
        = \begin{pmatrix}
            \subdiff_F F(\optx) + \alpha D\Psi_{\realopty_\alpha}(\realoptx_\alpha)
            + D \Psi_{\optx}(\realoptx_\alpha)
            + \frac{\eps_\alpha}{\lambda^2}\sum_{n=0}^\infty \mu_n D\Psi_{x_n}(\realoptx_\alpha) \\
            \subdiff_F G(\optx) + \alpha D \Psi_{\realoptx_\alpha}(\realopty_\alpha)
            + \frac{\eps_\alpha}{\lambda^2}\sum_{n=0}^\infty \mu_n D\Psi_{y_n}(\realopty_\alpha)
        \end{pmatrix}.
    \end{equation*}
    By \ref{item:epsilon:fuzzy-fermat:borweinpreiss:1} and $\realoptx_\alpha,\realopty_\alpha \in \B(\optx, \lambda)$ we have $\norm{\realoptx_\alpha-x_n}_X, \norm{\realopty_\alpha-y_n}_X \le \lambda$ for all $n \ge 0$. In addition, $\sum_{n=0}^\infty\mu_n = 1$, and thus $\frac{\eps_\alpha}{\lambda^2}\sum_{n=0}^\infty \mu_n D\Psi_{x_n}(\realoptx_\alpha) \in \frac{2\eps_\alpha}{\lambda} \B_{X^*}$ and likewise for $D\Psi_{y_n}$ (so that in fact we were justified in differentiating the series term-wise).
    By \ref{item:epsilon:fuzzy-fermat:borweinpreiss:1} also $\norm{\realoptx_\alpha-\optx}_X \le \lambda$, so that $D\Psi_{\optx}(\realoptx_\alpha) \in 2\lambda\B_{X^*}$.
    Finally, since $-x^*\in \partial\norm{\freevar}_X(-x)$ for any $x^*\in \partial\norm{\freevar}_X(x)$ and any $x\in X$, we have $D\Psi_y(x) = -D\Psi_x(y)$ for all $x,y\in X$.
    We thus have
    \begin{equation*}
        \left\{
            \begin{aligned}
                -\alpha D\Psi_{\realopty_\alpha}(\realoptx_\alpha) & \in \subdiff_F F(\realoptx_\alpha) + \left(2\lambda+\frac{2\eps_\alpha}{\lambda} \right)\B_{X^*},
                \\
                \alpha D\Psi_{\realopty_\alpha}(\realoptx_\alpha) & \in \subdiff_F G(\realopty_\alpha) + \frac{2\eps_\alpha}{\lambda} \B_{X^*},
            \end{aligned}
        \right.
    \end{equation*}
    which implies that
    \begin{equation*}
        0 \in \subdiff_F F(\realoptx_\alpha) + \subdiff_F G(\realopty_\alpha) + \left(2\lambda+\frac{4\eps_\alpha}{\lambda} \right)\B_{X^*}.
    \end{equation*}
    Since $(\realoptx_\alpha, \realopty_\alpha) \in  \B(\optx, \lambda)^2$, the claim now follows by taking $\lambda \in (0, \rho)$ small enough and then $\alpha>0$ large (and thus $\eps_\alpha$ small) enough.
\end{proof}

\begin{remark}[fuzzy Fermat principles and trustworthy subdifferentials]
    \Cref{lemma:epsilon:fuzzy-fermat} is due to \cite{fabian1988classes}.
    Such fuzzy Fermat principles are studied in more detail from the point of view of \term[principle!variational!fuzzy]{fuzzy variational principles} in \cite{ioffe2017variational}.
    Specifically, the claim of \cref{lemma:epsilon:fuzzy-fermat} has to hold for an arbitrary subdifferential operator $\subdiff_*$ for it to be called \term[subdifferential!trustworthy]{trustworthy}, whereas the converse inclusion $\subdiff_* G(x) + \subdiff_* F(x) \subset \subdiff_*[G+F](x)$ is required for the subdifferential to be called \term[subdifferential!elementary]{elementary}.
\end{remark}

\begin{remark}[notes on the proof of \cref{lemma:epsilon:fuzzy-fermat}]
    Note how we had to apply the Borwein--Preiss variational principle instead of Ekeland's to obtain a differentiable convex perturbation and thus to be able to apply the sum rule \cref{cor:epsilon:sumrule:gateaux}.
    In contrast, the proof in \cite{ioffe2017variational} is based on the \noemphterm[principle!variational!Deville--Godefroy--Zizler]{Deville--Godefroy--Zizler variational principle}, which makes no convexity assumption on the perturbation function and hence requires the stronger property of Fréchet smoothness (i.e., Fréchet instead of Gateaux differentiability of the norm outside the origin).

    Finally, with an additional argument showing $ J_\alpha(\realoptx_\alpha, \realopty_\alpha) \le j_\alpha + \beta_\alpha$ for a suitable $\beta_\alpha$, it would be possible to further constrain $\abs{F(x)-F(\optx)} \le \delta$ in the claim of \cref{lemma:epsilon:fuzzy-fermat}, as is done in \cite[Theorem 4.30]{ioffe2017variational}.
\end{remark}

\begin{corollary}
    \label{cor:epsilon:fuzzy-inclusion}
    Let $X$ be a Gateaux smooth Banach space, let $F: X \to \Rbar$ be lower semicontinuous near $\optx\in X$, and $\epsilon>0$.
    Then for any $\delta>0$ and $\epsilon'>\epsilon$ we have
    \begin{equation*}
        \subdiff_\epsilon F(\optx) \subset \Union_{z \in \B(\optx, \delta)} \subdiff_F F(z) + \epsilon'\B_{X^*}.
    \end{equation*}
\end{corollary}

\begin{proof}
    We may assume that $\optx \in \dom F$, in particular that there exists some $x^* \in \subdiff_\epsilon F(\optx)$, i.e., such that
    \begin{equation*}
        \liminf_{\optx \ne y\to \optx} \frac{F(y)-F(\optx)-\dualprod{x^*}{y-\optx}_X}{\norm{y-\optx}_X} \ge -\epsilon.
    \end{equation*}
    Taking any $\epsilon'>\epsilon$ and defining
    \begin{equation*}
        \bar F(x) \defeq F(x)-\dualprod{x^*}{x-\optx}_X
        \quad\text{and}\quad
        \bar G(x) \defeq \epsilon'\norm{x-\optx}_X,
    \end{equation*}
    we obtain as in \cref{eq:epsilon:alt-formulations} that
    \begin{equation*}
        \liminf_{\optx \ne y \to \optx} \frac{(\bar G+\bar F)(y)-(\bar G+\bar F)(\optx)}{\norm{y-\optx}_X} \ge (\epsilon'-\epsilon).
    \end{equation*}
    Thus $\bar F+\bar G$ achieves its local minimum at $\optx$.
    The function $\bar G$ is convex and Lipschitz while $\bar F$ lower semicontinuous. Hence \cref{lemma:epsilon:fuzzy-fermat} implies for any $\delta>0$ and $\mu'>0$ that
    \begin{equation*}
        0 \in \Union_{z, y \in \B(\optx, \delta)} \left( \subdiff_F \bar F(y) + \subdiff_F \bar G(z) \right) + \mu'\B_X.
    \end{equation*}
    Since $\subdiff_F \bar F(y)=\subdiff_F F(y)-\{x^*\}$ (by \cref{cor:epsilon:sumrule:gateaux} or directly from the definition) and $\subdiff_F \bar G(z)=\subdiff \bar G(z) \subset \epsilon'\B_{X^*}$, we obtain
    \begin{equation*}
        x^* \in  \Union_{z \in \B(x, \delta)} \subdiff_F \bar F(z) + (\mu'+\epsilon')\B_{X^*}.
    \end{equation*}
    Since $\mu'>0$ and $\epsilon'>\epsilon$ were arbitrary, the claim follows.
\end{proof}

\section{Approximate Fermat principles and projections}

We now introduce an \term[principle!Fermat!approximate]{approximate Fermat principle}, which can be invoked when we \emph{do not know} whether a minimizer exists; in particular, when $F$ fails to be \emph{weakly} lower semicontinuous so that \cref{thm:variation:existence} is not applicable.

\begin{theorem}
    \label{thm:epsilon:approximate-fermat}
    Let $X$ be a Banach space and $F: X \to \Rbar$ be proper, lower semicontinuous, and bounded from below. Then for every $\epsilon, \delta>0$ there exists an $\opt x_\eps \in X$ such that
    \begin{enumerate}[label=(\roman*)]
        \item\label{item:epsilon:approximate-fermat:fbound} $F(\opt x_\eps) \le \inf_{x\in X} F(x) + \epsilon$;
        \item\label{item:epsilon:approximate-fermat:min} $F(\opt x_\eps) < F(x)+\delta\norm{x-\optx_\eps}_X$ for all $x \ne \optx_\eps$;
        \item $0 \in \subdiff_\delta F(\opt x_\eps)$.
    \end{enumerate}
\end{theorem}

\begin{proof}
    Since $F$ is bounded from below, $\inf_{x\in X}F(x)>-\infty$. We can thus take a minimizing sequence $\{x_n\}_{n \in \N}$ with $F(x_n) \downto \inf_{x\in X} F(x)$ and find a $n(\eps)\in\N$ such that $x_\eps \defeq x_{n(\eps)}$ satisfies \cref{item:epsilon:approximate-fermat:fbound}. Ekeland's variational principle \cref{thm:variation:ekeland} thus yields for $\lambda\defeq \eps/\delta$ an $\optx_\eps \defeq \optx_{\eps,\lambda}$ such that $\norm{\optx_\eps-x_\eps}_X\leq \lambda$,
    \begin{equation*}
        F(\optx_\eps) \leq F(\optx_\eps) + \frac\eps\lambda\norm{\optx_\eps - x_\eps}_X \leq F(x_\eps),
    \end{equation*}
    as well as
    \begin{equation*}
        F(\optx_\eps) < F(x) + \frac{\eps}{\lambda}\norm{\optx_\eps-x}_X
        \quad (x \ne \optx_\eps).
    \end{equation*}
    Thus \cref{item:epsilon:approximate-fermat:fbound} as well as \cref{item:epsilon:approximate-fermat:min} hold.
    The latter implies for all $x \ne \optx_\eps$ that
    \begin{equation*}
        \frac{F(x)-F(\opt x_{\epsilon}) -\dualprod{0}{x-\opt x_{\epsilon}}_X}{\norm{x-\opt x_{\epsilon}}_X} \ge - \delta,
    \end{equation*}
    i.e., $0\in\partial_\delta F(\optx_{\epsilon})$ by definition.
\end{proof}

As an example for possible applications of approximate Fermat principles, we use \cref{thm:epsilon:approximate-fermat} to prove the following result on projections and approximate projections onto a \emph{nonconvex} set $C \subset X$.
For nonconvex sets, even the exact projection need no longer be unique; furthermore, for the reasons discussed before \cref{thm:epsilon:approximate-fermat}, the set of projections $P_C(x)$ may be empty when $C \ne \emptyset$ is closed but not \emph{weakly} closed.
We recall that by \cref{lem:convex_closed}, convex closed sets are weakly closed, as are, of course, finite-dimensional closed sets. However, more generally, weak closedness can be elusive. Hence we will need to perform \term[projection!approximate]{approximate projections} in \cref{part:setvalued}.
It is not surprising that this requires additional assumptions on the containing space to make up for this.

\begin{theorem}
    \label{thm:epsilon:projection}
    Let $X$ be a Gateaux smooth Banach space and let $C\subset X$ be nonempty and closed.
    Define the (possibly multi-valued) projection
    \begin{align*}
        P_C & : X \setto X,\qquad  P_C(x) \defeq \argmin_{\alt x \in C} \norm{\alt x-x}_X
        \intertext{and the corresponding distance function}
        d_C & : X \to \R, \qquad d_C(x)  \defeq \inf_{\alt x \in C} \norm{\alt x-x}_X.
    \end{align*}
    Then the following hold:
    \begin{enumerate}[label=(\roman*)]
        \item\label{item:epsilon:projection:frechet}
            For any $\optx \in P_C(x)$, there exists an $\optx^* \in \subdiff_F \delta_C(\optx)$ such that
            \begin{equation}
                \label{eq:epsilon:projection}
                \dualprod{\optx^*}{x-\optx}_X=\norm{x-\optx}_X,
                \qquad
                \norm{\optx^*}_{X^*} \le 1.
            \end{equation}

        \item\label{item:epsilon:projection:epsilon}
            For any $\epsilon>0$, there exists an approximate projection $\optx_\epsilon \in C$ satisfying
            \begin{equation*}
                \norm{\optx_\epsilon-x}_X \le d_C(x) + \epsilon
            \end{equation*}
            as well as \eqref{eq:epsilon:projection} for some $\optx^* \in \subdiff_\epsilon \delta_C(\optx_\epsilon)$.

        \item\label{item:epsilon:projection:hilbert}
            If $X$ is a Hilbert space, then $x-\optx \in \subdiff_\eps \delta_C(\optx)$ for all $\eps\geq 0$.
    \end{enumerate}
\end{theorem}

\begin{proof}
    \emph{\ref{item:epsilon:projection:frechet}:}
    Let $x \not \in C$, since otherwise $\optx^* \defeq 0 \in \subdiff_F(\optx)$ for $\optx=x \in C$ by the definition of the Fréchet subdifferential.
    Set $F(\tilde x)\defeq \norm{\tilde x-x}_X$ and assume that $\optx \in P_C(x)$.
    The Fermat principle \cref{thm:limiting:frechet:fermat} then yields that $0 \in \subdiff_F[\delta_C+F](\optx)$. Since $x \not \in C$ and $\optx \in C$, by assumption $F$ is differentiable at $\opt x$. Thus \cref{thm:convex:gateaux} shows that $\subdiff F(\optx)=\{DF(\optx)\}$ is a singleton.
    The sum rule of \cref{cor:epsilon:sumrule:gateaux} then yields that $\optx^* \defeq -DF(\optx) \in \subdiff_F \delta_C(\optx)$.
    The claim of \eqref{eq:epsilon:projection} now follows from \cref{thm:subdifferential:norm}.

    \emph{\ref{item:epsilon:projection:epsilon}:}
    Compared to \ref{item:epsilon:projection:frechet}, we merely invoke the approximate Fermat principle of \cref{thm:epsilon:approximate-fermat} in place of \cref{thm:limiting:frechet:fermat}, which establishes the existence of $\optx_\epsilon \in C$ satisfying $\norm{\optx_\epsilon-x}_X \le d_C(x)$ and $0 \in \subdiff_\epsilon[\delta_C+F](\optx)$. The sum rule of \cref{lemma:epsilon:sumrule} then shows that $\optx^* \defeq -DF(\optx) \in \subdiff_\epsilon \delta_C(\optx)$.

    \emph{\ref{item:epsilon:projection:hilbert}:}
    In a Hilbert space, we can identify $-DF(\optx)$ with the corresponding gradient $-\nabla F(\optx) = (x-\optx)/\norm{x-\optx}_X\in X$ for $\optx\neq 0$ (otherwise $-\nabla F(\optx) = 0 = x-\optx$).
    Since $\subdiff_\eps \delta_C(\optx)$ is a cone, this implies that $x-\optx \in \subdiff_F \delta_C(\optx)$ as well.
\end{proof}

In the next chapters, we will see that $\subdiff_F \delta_C(\optx)$ coincides with a suitable normal cone to $C$ at $\optx$. In other words, $\optx^*$ is a normal vector to the set $C$. In Hilbert spaces, this normal vector can be identified with the (normalized) vector pointing from $\optx$ to $x$.

\part{Set-valued analysis}\label{part:setvalued}

\addtocontents{toc}{\protect\enlargethispage{1cm}}
\chapter{Tangent and normal cones}
\label{chap:cones}

We now start our study of stability properties of the solutions to nonsmooth optimization problems. As we have characterized the latter via subdifferential inclusions, we need to study the sensitivity of such relations to perturbations. As in the smooth case, this can be done through derivatives of these conditions with respect to relevant parameters; however, these conditions are expressed as inclusions instead of simple equations. Hence we require notions of derivatives for set-valued mappings.

To motivate how we will develop differential calculus for set-valued mappings, recall from \cref{lem:convex:subdiff_epi} how the subdifferential of a convex function $F$ can be defined in terms of the normal cone to the epigraph of $F$. This idea forms the basis of differentiating general set-valued mappings $H: X \setto Y$, where instead of taking the normal cone at $(x, F(x))$ to $\epi F$, we do this at any point $(x, y)$ of $\graph H \defeq \{(x, y) \in X \times Y \mid y \in H(x)\}$. Since we are generally not in the nice convex setting -- even for a convex function $F$‚ the set $\graph \subdiff F$ is not convex unless $F$ is linear -- there are some complications which result in having to deal with various nonequivalent definitions.
In this chapter, we introduce the relevant graphical notions of tangent and normal cones.
In \cref{chap:pointcones}, we develop specific expressions for these cones to sets in $L^p(\Omega)$ defined as pointwise via finite-dimensional sets.
In the following \crefrange{chap:graphical}{chap:colimiting}, we then define and further develop notions of differentiation of set-valued mappings based on these cones.

\section{Definitions and examples}

\subsection*{The fundamental cones}

Our first type of tangent cone is defined using roughly the same limiting process on difference quotients as basic directional derivatives.
Let $X$ be a Banach space.
We define the \term[cone!tangent]{tangent cone} (or \term[cone!Bouligand]{Bouligand} or \term[cone!contingent]{contingent cone}) of the set $C \subset X$ at $x \in X$ as
\begin{equation}
    \label{eq:cones:def-tangent}
    \begin{aligned}[t]
        T_C(x)&\defeq \limsup_{\tau \downto 0} \frac{C-x}{\tau}\\
        &= \setof{\Delta x\in X}{\Delta x = \lim_{k \to \infty} \frac{x_k - x}{\tau_k} \text{ for some }C \ni x_k\to x,\, \tau_k \downto 0},
    \end{aligned}
\end{equation}
i.e., the tangent cone is the outer limit (in the sense of \cref{sec:monotone:basic}) of the ``blown up'' sets $(C-x)/\tau$ as $\tau \downto 0$. A $\Delta x\in T_C(x)$ is called a \term[vector!tangent]{tangent vector} to $C$ at $x$.

The tangent cone is closely related to the \term[cone!normal!Fréchet]{Fréchet normal cone}, which is based on the same limiting process as the Fréchet subdifferential in \cref{chap:limiting}:
\begin{equation}
    \label{eq:cones:def-frechetnormal}
    \frechetNormal_C(x) \defeq \setof{x^*\in X^*}{\limsup_{C \ni \alt x \to x} \frac{\dualprod{x^*}{\alt x-x}_X}{\norm{\alt x-x}_X} \le 0}.
\end{equation}
Correspondingly, any $x^*\in \frechetNormal_C(x)$ is called a \term[vector!normal]{normal vector} to $C$ at $x$.

\subsection*{Limiting cones in finite dimensions}

One difficulty with Fréchet normal cones is that they are not outer semicontinuous. By taking their outer limit (in the sense of set-valued mappings), we obtain the less ``irregular'' (\term[cone!normal!basic]{basic} or \term[cone!normal!limiting]{limiting} or \term[cone!normal!Mordukhovich]{Mordukhovich}) \term[cone!normal]{normal cone}.
This definition is somewhat more involved in infinite dimensions, so we first consider $C \subset \R^N$ at $x \in \R^N$. In this case, the limiting normal cone is defined as
\begin{equation}
    \label{eq:cones:def-limnormal-findim}
    \begin{aligned}[t]
        N_C(x) &\defeq \limsup_{C \ni \alt x \to x} \frechetNormal_C(\alt x)\\
        &= \setof{x^*\in \R^N}{x^* = \lim_{k\to \infty} x^*_k \text{ for some } x^*_k \in \frechetNormal_C(x_k),\, C\ni x_k\to x}.
    \end{aligned}
\end{equation}
Despite $N_C$ being obtained by the outer semicontinuous regularization of $\frechetNormal_C$, the \emph{latter} is sometimes in the literature called the \term[cone!normal!regular]{regular normal cone}. We stick to the convention of calling $\frechetNormal_C$ the \term[cone!normal!Fréchet]{Fréchet normal cone} and $N_C$ the \term[cone!normal!limiting]{limiting normal cone}.

The limiting variant of the tangent cone is the \term[cone!tangent!Clarke]{Clarke tangent cone} (also known as the \term[cone!tangent!Clarke]{regular tangent cone}), defined for a set $C \subset \R^N$ at $x \in \R^N$ as the inner limit
\begin{equation}
    \label{eq:cones:def-clarketangent-findim}
    \begin{aligned}[t]
        \clarkeTangent_C(x) &\defeq \liminf_{\substack{C \ni \alt x \to x, \\\tau \downto 0}}~ \frac{C-\alt x}{\tau}
        \\
        &= \setof{\Delta x\in\R^N}
        {\begin{array}{r}
                \text{for all }\tau_k\downto 0,\, C\ni x_k\to x \text{ there exists } C\ni \tilde x_k \to x \\
                \text{ with } (\tilde x_k - x_k)/\tau_k \to \Delta x
        \end{array} }.
    \end{aligned}
\end{equation}
We will later in \cref{cor:cones:clarke-liminf} see that for a closed set $C \subset \R^N$, we in fact have that $\clarkeTangent_C(x)=\liminf_{C \ni \alt x \to x} T_C(\alt x)$.

The following example as well as  \cref{fig:cones:tangent-normal} illustrate the different cones.

\begin{figure}
    \centering
    \begin{subfigure}[t]{0.24\textwidth}
        \centering
        \begin{asy}
            pair x=(0.3, 1); int t=2;
            path p=(0, 0)..(0.2, 1)..x--(0.3,1.5)..(1,1.2)..(1, .5)..cycle;
            fill(p, lightfill);
            label("$C$", x, 5*S);
            label("$x$", x, 2*E);
            pair t1=-dir(p, t, -1);
            pair t2=dir(p, t, 1);
            fill((x+0.4*t1)..(x-0.2*(t1+t2))..(x+0.4*t2)--x--cycle, darkfill);
            draw(x--(x+0.7*t1), primalline+linewidth(1.1), Arrow);
            draw(x--(x+0.7*t2), primalline+linewidth(1.1), Arrow);
            dot(x);
            draw(p);
        \end{asy}
        \caption{tangent cone $T_C(x)$}
    \end{subfigure}
    \begin{subfigure}[t]{0.24\textwidth}
        \centering
        \begin{asy}
            pair x=(0.3, 1); int t=2;
            path p=(0, 0)..(0.2, 1)..x--(0.3,1.5)..(1,1.2)..(1, .5)..cycle;
            fill(p, lightfill);
            label("$C$", x, 5*S);
            label("$x$", x, 2*E+N);
            pair t1=-dir(p, t, -1);
            pair t2=dir(p, t, 1);
            pair n1=orthog(t1, -1);
            pair n2=orthog(t2, -1);
            fill((x-0.4*t1)..(x-0.4/sqrt(2)*(t1+t2))..(x-0.4*t2)--x--cycle, darkfill);
            draw(x--(x-0.7*t1), primalline+linewidth(1.1), Arrow);
            draw(x--(x-0.7*t2), primalline+linewidth(1.1), Arrow);
            dot(x);
            draw(p);
        \end{asy}
        \caption{Clarke tangent cone $\clarkeTangent_C(x)$}
    \end{subfigure}
    \begin{subfigure}[t]{0.24\textwidth}
        \centering
        \begin{asy}
            pair x=(0.3, 1); int t=2;
            path p=(0, 0)..(0.2, 1)..x--(0.3,1.5)..(1,1.2)..(1, .5)..cycle;
            fill(p, lightfill);
            label("$C$", x, 5*S);
            label("$x$", x, 2*E);
            draw(p);
            dot(x, primalline);
        \end{asy}
        \caption{Fréchet normal cone $\frechetNormal_C(x)=\{0\}$}
        \label{fig:cones:tangent-normal-frechetnormal}
    \end{subfigure}
    \begin{subfigure}[t]{0.24\textwidth}
        \centering
        \begin{asy}
            pair x=(0.3, 1); int t=2;
            path p=(0, 0)..(0.2, 1)..x--(0.3,1.5)..(1,1.2)..(1, .5)..cycle;
            fill(p, lightfill);
            label("$C$", x, 5*S);
            label("$x$", x, 2*E);
            draw(p);
            dot(x);
            draw(x--(x+0.7*orthog(dir(p, t, -1), -1)), primalline+linewidth(1.1), Arrow);
            draw(x--(x+0.7*orthog(dir(p, t, 1), -1)), primalline+linewidth(1.1), Arrow);
        \end{asy}
        \caption{limiting normal cone $N_C(x)$}
    \end{subfigure}
    \caption{Illustration of the different normal and tangent cones at a nonregular point of a set $C$. The dot indicates the base point $x$. The thick arrows and dark filled-in areas indicate the directions included in the cones.}
    \label{fig:cones:tangent-normal}
\end{figure}

\begin{example}
    \label{ex:cones:basicex}
    We compute the different tangent and normal cones at all points $x \in C$ for different $C \subset \R^2$.
    \begin{enumerate}
        \item $C=\B(0,1)$: Clearly, if $x \in \interior C$, then
            \begin{align*}
                N_C(x)&=\frechetNormal_C(x)=\{0\},\\
                T_C(x)&=\clarkeTangent_C(x)=\R^2.
            \end{align*}
            For any $x \in \BD C$, on the other hand,
            \begin{align*}
                N_C(x)&=\frechetNormal_C(x)=[0, \infty)x\defeq \setof{tx}{t\geq 0},\\
                T_C(x)&=\clarkeTangent_C(x)=\{z \mid \iprod{z}{x} \le 0\}.
            \end{align*}

        \item $C=[0,1]^2$: For $x \in \interior C$, we again have that $N_C(x)=\frechetNormal_C(x)=\{0\}$ and $T_C(x)=\clarkeTangent_C(x)=\R^2$.
            Let then $x = (1, y)$ for $y \in (-1, 1)$, i.e., $x$ is on the right edge of $C$ excluding the corners.
            Then
            \begin{align*}
                N_C(x)&=\frechetNormal_C(x)=[0, \infty) \times \{0\},
                \\
                T_C(x)&=\clarkeTangent_C(x)=(-\infty, 0] \times \R.
            \end{align*}
            The other points $x \in \BD C\setminus \{(0,0), (0,1), (1,0),(1,1)\}$ are handled similarly.

            Of the corners, we concentrate on $x=(1,1)$, the others being analogous. Here
            \begin{align*}
                N_C(x)&=\frechetNormal_C(x)=\{(\Delta x, \Delta y) \mid \Delta x, \Delta y \ge 0\},\\
                T_C(x)&=\clarkeTangent_C(x)=\{(\Delta x, \Delta y) \mid \Delta x, \Delta y \le 0\}.
            \end{align*}

        \item $C=[0, 1]^2 \setminus [\frac12, 1]^2$: Here as well $N_C(x)=\frechetNormal_C(x)=\{0\}$ and $T_C(x)=\clarkeTangent_C(x)=\R^2$ for $x\in \interior C$.
            Other points on $\BD C$ are computed analogously to similar corners and edges of the square $[0, 1]^2$, but we have to be careful with the \enquote{interior corner} $x=(\frac12, \frac12)$.
            Here, similarly to \cref{fig:cones:tangent-normal-frechetnormal}, we see that $\frechetNormal_C(x)=\{0\}$. However, as a $\limsup$,
            \begin{equation*}
                N_C(x)=(0, 1)[0, \infty) \union (1,0)[0,\infty).
            \end{equation*}
            For the tangent cones, we then get
            \begin{align*}
                T_C(x)&=\{(\Delta x, \Delta y) \mid  \Delta x \le 0 \text{ or } \Delta y \le 0 \},\\
                \shortintertext{while, as a $\liminf$,}
                \clarkeTangent_C(x)&=\{(\Delta x, \Delta y) \mid  \Delta x \le 0 \text{ and } \Delta y \le 0 \}.
            \end{align*}
    \end{enumerate}
\end{example}

\subsection*{Limiting cones in infinite dimensions}

Let now $X$ be again a Banach space.
Although the fundamental cones -- the (basic) tangent cone and the Fréchet normal cone -- were defined based on strongly convergent sequences, in infinite-dimensional spaces weak modes of convergence better replicate various relationships between the different cones. We thus call an element $\dir x \in X$ \emph{weakly tangent} to $C$ at $x$ if
\begin{equation}
    \label{eq:cones:def-tangent-weak}
    \dir x=\weaklim_{k \to \infty} \frac{x_k - x}{\tau_k} \quad \text{for some} \quad C \ni x_k\to x,\, \tau_k \downto 0,
\end{equation}
where the $\weaklim$ of course stands for $\tau_k^{-1}(x_k -x)\weakto \dir x$.
We denote by the \term[cone!tangent!weak]{weak tangent cone} (or \term[cone!contingent!weak]{weak contingent cone}) $T^w_C(x) \subset X$ the set of all such $\dir x$.
Using the notion of outer limits of set-valued mappings from \cref{chap:monotone}, we can also write
\begin{equation}
    \label{eq:cones:tangent-weak}
    T^w_C(x)=\weaklimsup_{\tau \downto 0} \frac{C-x}{\tau}.
\end{equation}

Likewise, the limiting normal cone $N_C(x)$ to $C \subset X$ in a general infinite-dimensional Banach space $X$ is based on weak-$*$ limits. Moreover, several proofs will be easier if we slightly relax the definition. Therefore, given $\epsilon\ge 0$ we first introduce the \term[cone!normal!$\epsilon$-]{$\epsilon$-normal cone}  of $x^* \in X^*$ satisfying
\begin{equation}
    \label{eq:cones:def-epsiloncone}
    \frechetNormal_C^\epsilon(x) \defeq \setof{x^*\in X^*}{\limsup_{C \ni \alt x \to x} \frac{\dualprod{x^*}{\alt x-x}_X}{\norm{\alt x-x}_X} \le \epsilon}.
\end{equation}
The \term[cone!normal!Fréchet]{Fréchet normal cone} is then simply $\frechetNormal_C(x) \defeq \frechetNormal_C^0(x)$.

Now, the (\term[cone!normal!basic]{basic} or \term[cone!normal!limiting]{limiting} or \term[cone!normal!Mordukhovich]{Mordukhovich}) \emph{normal cone} is defined as
\begin{equation}
    \label{eq:cones:def-limnormal}
    N_C(x) \defeq \weakstarlimsup_{\substack{\alt x \to x,~ \epsilon \downto 0}} \frechetNormal_C^\epsilon(\alt x).
\end{equation}
In other words, $x^* \in N_C(x)$ if and only if there exist $C \ni x_k \to x$, $\epsilon_k \downto 0$ and $x_k^* \in N_C^{\epsilon_k}(x_k)$ such that $x_k^* \weaktostar x^*$.

In Gateaux smooth Banach spaces, we can fix $\epsilon \equiv 0$ in \eqref{eq:cones:def-limnormal}. Thus such spaces can be treated similarly to the finite-dimensional case in \eqref{eq:cones:def-limnormal-findim}.

\begin{theorem}
    \label{thm:cones:nonepsilon-limnormal}
    Let $X$ be a Gateaux smooth Banach space, $C \subset X$, and $x \in X$.
    Then
    \begin{equation}
        \label{eq:cones:nonepsilon-limnormal}
        N_C(x) = \weakstarlimsup_{\alt x \to x} \frechetNormal_C(\alt x).
    \end{equation}
\end{theorem}

\begin{proof}
    Denote by $K$ the set on the right hand side of \eqref{eq:cones:nonepsilon-limnormal}. Then by the definition \eqref{eq:cones:def-limnormal}, clearly $N_C(x) \supset K$.
    To show $N_C(x) \subset K$, let $x^* \in N_C(x)$. Then \eqref{eq:cones:def-limnormal} yields $x_k \to x$, $\epsilon_k \downto 0$, and  $x_k^* \weaktostar x^*$ with $x_k^* \in \frechetNormal_C^{\epsilon_k}(x_k)$.
    We need to show that there exist some $\alt x_k \to x$ and $\alt x_k^* \weaktostar x^*$ with $\alt x_k^* \in \frechetNormal_C(\alt x_k)$.
    Indeed, since $\frechetNormal_C^\epsilon=\subdiff_\epsilon \delta_C$, by \cref{cor:epsilon:fuzzy-inclusion} applied to $F=\delta_C$, we have for any sequence $\delta_k \downto 0$ that
    \begin{equation*}
        x_k^* \in N_C^{\epsilon_k}(x_k) \subset \Union_{\alt x \in \B(x_k, \delta_k)} N_C(\alt x) + \delta_k\B_{X^*} \qquad (k\in \N).
    \end{equation*}
    In particular, there exist $\alt x_k \in \B(x_k, \delta_k)$ and $\alt x_k^* \in N_C(\alt x_k) \isect \B(x_k^*, \delta_k)$, which implies that $\alt x_k \to x$ and $\alt x_k^* \weaktostar x^*$ as desired.
\end{proof}

\begin{remark}
    \Cref{thm:cones:nonepsilon-limnormal} can be extended to Asplund spaces -- in particular to reflexive Banach spaces. The equivalence of \eqref{eq:cones:nonepsilon-limnormal} and \eqref{eq:cones:def-limnormal} can, in fact, be used as a definition of an Asplund space. For details we refer to \cite[Theorem 2.35]{Mordukhovich:2006}.
\end{remark}

Finally, the \term[cone!tangent!Clarke]{Clarke tangent cone} is defined as in finite dimensions as
\begin{equation}
    \label{eq:cones:def-clarketangent}
    \begin{aligned}[t]
        \clarkeTangent_C(x) &\defeq \liminf_{\substack{C \ni \alt x \to x, \\\tau \downto 0}}~ \frac{C-\alt x}{\tau}
        \\
        &= \setof{\Delta x\in X}
        {\begin{array}{r}
                \text{for all }\tau_k\downto 0,\, C\ni x_k\to x \text{ there exists } C\ni \tilde x_k \to x \\
                \text{ with } (\tilde x_k - x_k)/\tau_k \to \Delta x
        \end{array} }.
    \end{aligned}
\end{equation}
In infinite-dimensional spaces, however, we in general only have the inclusion
\begin{equation*}
    \clarkeTangent_C(x)\subset\liminf_{C \ni \alt x \to x} T_C(\alt x);
\end{equation*}
see \cref{cor:cones:clarke-liminf}.

\begin{remark}[a much too brief history of various cones]
    The (Bouligand) tangent cone was already introduced for smooth sets by Peano in 1908 \cite{peano}; the term \term[cone!contingent]{contingent cone} is due to Bouligand \cite{bouligand1930}.
    The Clarke tangent cone (also called \term[cone!circatangent]{circatangent cone}) was introduced in \cite{clarke1973necessary,clarke1975generalized}; see also \cite{Clarke:1990a}. The limiting normal cone can be found in \cite{mordukhovich1976maximum}, who stressed the need of defining (nonconvex) normal cones directly rather than as (necessarily convex) polars of tangent cones. The history of the Fréchet normal cone is harder to trace, but it has appeared in the literature as the polar of the tangent cone.
    We will see that in finite dimensions, $\frechetNormal_C(x)=\polar{T_C(x)}$.
    In infinite dimensions, $\polar{T_C(x)}$ is sometimes called the \term[cone!normal!Dini]{Dini normal cone} and is in general not equal to the Fréchet normal cone.

    We do not attempt to do full justice to the muddier parts of the historical development here, and rather refer to the accounts in \cite{dolecki2011tangency,bigolin2014historical} as well as \cite[Commentary to Ch.~6] {Rockafellar:1998} and \cite[Commentary to Ch.~1]{mordukhovich2018variational}. Various further cones are also discussed in \cite{aubin1990sva}.
\end{remark}

\section{Basic relationships and properties}

As seen in \cref{ex:cones:basicex}, the limiting normal cone $N_C(x)$ can be larger than the Fréchet normal cone $\frechetNormal_C(x)$; conversely, the Clarke tangent cone $\clarkeTangent_C(x)$ is \emph{smaller} than the tangent cone $T_C(x)$; see \cref{fig:cones:tangent-normal}.
These inclusions hold in general.

\begin{theorem}
    \label{thm:cones:inclusions}
    Let $C \subset X$ and $x \in X$. Then
    \begin{enumerate}
        \item\label{item:cones:inclusions:tangent}
            $\clarkeTangent_C(x)  \subset T_C(x) \subset T_C^w(x)$;
        \item\label{item:cones:inclusions:normal}
            $\frechetNormal_C(x)  \subset N_C(x)$.
    \end{enumerate}
\end{theorem}

\begin{proof}
    If we fix the base point $\alt x$ as $x$ in the definition \eqref{eq:cones:def-clarketangent}  of $\clarkeTangent_C(x)$, the inclusion $\clarkeTangent_C(x)  \subset T_C(x)$ is clear from the definition \eqref{eq:cones:def-tangent} of $T_C(x)$ as an outer limit and of $\clarkeTangent_C(x)$ as an inner limit.
    The inclusion $T_C(x) \subset T_C^w(x)$ is likewise clear from the definition of $T^C(x)$ as a strong outer limit and of $T_C^w(x)$ as the corresponding weak outer limit.

    The normal inclusion $\frechetNormal_C(x) \subset N_C(x)$ follows from the definition \eqref{eq:cones:def-limnormal} of $N_C(x)$ as the outer limit of $\frechetNormal_C^\epsilon(\alt x)$ as $\alt x \to x$ and $\epsilon \downto 0$. (In finite dimensions, we can fix $\epsilon=0$ in this argument or refer to the equivalence of definitions shown in \cref{thm:cones:nonepsilon-limnormal}.)
\end{proof}

For a closed and convex set $C$, however, both the Fréchet and limiting normal cones coincide with the convex normal cone defined in \cref{lem:convex:normalcone} (which we here denote by $\partial\delta_C(x)$ to avoid confusion).

\begin{lemma}
    \label{lemma:cones:convex}
    Let $C \subset X$ be nonempty, closed, and convex. Then for all $x\in X$,
    \begin{enumerate}
        \item\label{item:cones:convex:frechet}
            $\frechetNormal_C(x) = \partial \delta_C(x)$;
        \item\label{item:cones:convex:limiting}
            if $X$ is Gateaux smooth (in particular, finite-dimensional), $N_C(x) = \partial \delta_C(x)$.
    \end{enumerate}
\end{lemma}

\begin{proof}
    If $x\notin C$, it follows from their definitions that all three cones are empty. We can thus assume that $x\in C$.

    \emph{\ref{item:cones:convex:frechet}:}
    If $x^*\in \partial \delta_C(x)$, we have by definition that
    \begin{equation*}
        \dualprod{x^*}{y-x}_X \leq 0\quad\text{for all }y\in C.
    \end{equation*}
    Taking in particular $y=\tilde x$ and passing to the limit $\tilde x\to x$ thus implies that $x^*\in \frechetNormal_C(x)$.

    Conversely, let $x^*\in \frechetNormal_C(x)$ and let $y\in C$ be arbitrary. Since $C$ is convex, this implies that $x_t\defeq x+t(y-x)\in C$ for any $t\in (0,1)$ as well. We also have that $x_t\to x$ for $t\to 0$. From \eqref{eq:cones:def-frechetnormal}, it then follows by inserting the definition of $x_t$ and dividing by $t>0$ that
    \begin{equation*}
        0 \geq \lim_{t\to 0} \frac{\dualprod{x^*}{x_t -x}_X}{\norm{x_t-x}_X} =
        \frac{\dualprod{x^*}{y-x}_X}{\norm{y-x}_X}.
    \end{equation*}
    and hence, since $y\in C$ was arbitrary, that $x^*\in \partial\delta_C(x)$.

    \emph{\ref{item:cones:convex:limiting}:}
    By \cref{lem:variation:indicator,thm:monoton:subdiff,cor:monoton:closed}, $\partial \delta_C$ is strong-to-weak-$*$ outer semicontinuous, which by \cref{thm:cones:inclusions} and the $\epsilon \equiv 0$ characterization of \cref{thm:cones:nonepsilon-limnormal} implies that
    \begin{equation*}
        N_C(x)
        = \weakstarlimsup_{\alt x \to x} \frechetNormal_C(\alt x)
        = \weakstarlimsup_{\alt x \to x} \subdiff\delta_C(\alt x)
        \subset \subdiff\delta_C(x)
        = \frechetNormal_C(x)
        \subset N_C(x).
    \end{equation*}
    Hence $\frechetNormal_C(x) = N_C(x)$.
\end{proof}

Note that convexity was only used for the second inclusion, and hence $\partial \delta_C(x)\subset N_C(x)$ always holds.
In general, comparing \eqref{eq:cones:def-frechetnormal} with \eqref{eq:limiting:frechet}, we have the following relation.
\begin{corollary}\label{lemma:cones:frechet-subdiff}
    Let $C\subset X$ and $x\in X$. Then $\frechetNormal_C(x) = \partial_F \delta_C(x)$.
\end{corollary}

The next theorem lists some of the most basic properties of the various tangent and normal cones.

\begin{theorem}
    \label{thm:cones:basic-prop}
    Let $C \subset X$ and $x \in X$. Then
    \begin{enumerate}
        \item\label{item:cones:basic-prop:cone} $T_C(x)$, $\clarkeTangent_C(x)$, $\frechetNormal_C(x)$, and $N_C(x)$ are cones;
        \item $T_C(x)$, $\clarkeTangent_C(x)$, and $\frechetNormal_C^\epsilon(x)$ for every $\epsilon \ge 0$ are closed;
        \item $\clarkeTangent_C(x)$ and $\frechetNormal_C^\epsilon(x)$ for every $\epsilon \ge 0$ are convex;
        \item if $X$ is finite-dimensional, then $N_C(x)$ is closed.
    \end{enumerate}
\end{theorem}

\begin{proof}
    We argue the different properties for each type of cone in turn.

    \emph{The Fréchet ($\epsilon$-)normal cone:}
    It is clear from the definition of $\frechetNormal_C(x)$ that it is a cone, i.e., that $x^*\in \frechetNormal_C(x)$ implies that $\lambda x^* \in\frechetNormal_C(x)$ for all $\lambda>0$.

    Let now $\epsilon\geq 0$ be arbitrary. Let $x_k^* \in \frechetNormal_C^\epsilon(x)$ converge to some $x^* \in X^*$. Also suppose $C \ni x_\ell \to x$. Then for any $\ell, k \in \N$, we have by the Cauchy--Schwarz inequality that
    \begin{equation*}
        \frac{\dualprod{x^*}{x_\ell-x}_X}{\norm{x_\ell-x}_X}
        \le
        \frac{\dualprod{x_k^*}{x_\ell-x}_X}{\norm{x_\ell-x}_X}
        + \norm{x_k^*-x^*}_X
    \end{equation*}
    and thus that
    \begin{equation*}
        \limsup_{\ell \to \infty} \frac{\dualprod{x^*}{x_\ell-x}_X}{\norm{x_\ell-x}_X}
        \le
        \epsilon + \norm{x_k^*-x^*}_X.
    \end{equation*}
    Since $k \in \N$ was arbitrary and $x_k^* \to x^*$, we see that $x^* \in \frechetNormal_C^\epsilon(x)$ and may conclude that $\frechetNormal_C^\epsilon(x)$ is closed.

    To show convexity, take $x_1^*,x_2^* \in \frechetNormal_C^\epsilon(x)$ and let $x^* \defeq \lambda x_1^*+(1-\lambda)x_2^*$ for some $\lambda \in (0, 1)$. We then have
    \begin{equation*}
        \frac{\dualprod{x^*}{x_\ell-x}_X}{\norm{x_\ell-x}_X}
        =
        \lambda\frac{\dualprod{x_1^*}{x_\ell-x}_X}{\norm{x_\ell-x}_X}
        +
        (1-\lambda)\frac{\dualprod{x_2^*}{x_\ell-x}_X}{\norm{x_\ell-x}_X}.
    \end{equation*}
    Taking the limit $x_\ell \to x$ now yields $x^* \in \frechetNormal_C^\epsilon(x)$ and hence the convexity.

    \emph{The limiting normal cone:}
    If $X$ is finite-dimensional, the set $N_C(x)$ is a closed cone as the strong outer limit of the (closed) cones $\frechetNormal_C(x_\ell)$ as $x_\ell \to x$; see \cref{lemma:limsup-setlimit}.

    \emph{The tangent cone:}
    By \cref{lemma:limsup-setlimit}, $T_C(x)$ is closed as the outer limit of the sets $C_\tau \defeq (C-x)/\tau$ as $\tau \downto 0$.
    To see that it is a cone, suppose $\dir x \in T_C(x)$. Then there exist by definition $\tau_k \downto 0$ and $C \ni x_k \to x$ such that $(x_k-x)/\tau_k \to \dir x$. Now, for any $\lambda>0$, taking $\tilde\tau_k \defeq \inv\lambda\tau_k$, we have $(x_k-x)/\tilde\tau_k \to \lambda \dir x$. Hence $\lambda \dir x \in T_C(x)$.

    \emph{The Clarke tangent cone:}
    Finally, $\clarkeTangent_C(x)$ is a closed set through its definition as an inner limit, cf.~\cref{lemma:limsup}, as well as a cone by analogous arguments as for $T_C(x)$. To see that it is convex, take $\dir x^1, \dir x^2 \in \clarkeTangent_C(x)$. Since $\clarkeTangent_C(x)$ is a cone, we only need to show that $\dir x \defeq \dir x^1+\dir x^2 \in \clarkeTangent_C(x)$. By the definition of $\clarkeTangent_C(x)$ as an inner limit, we therefore have to show that for any sequence $\tau_k \downto 0$ and any ``base point sequence'' $C\ni x_k \to x$, there exist $\alt x_k \in C$ such that  $(\alt x_k - x_k)/\tau_k \to \dir x$.
    We do this by using the varying base point in the definition of $\clarkeTangent_C(x)$ to ``bridge'' between the sequences that generate $\dir x^1$ and $\dir x^2$; see \cref{fig:cones:bridging}.
    First, since $\dir x^1 \in \clarkeTangent_C(x)$, by the very same definition of $\clarkeTangent_C(x)$ as an inner limit, we can find for the base point sequence $\{x_k\}_{k \in \N}$ points $C\ni x^1_k \to x$ with $(x^1_k-x_k)/\tau_k \to \dir x^1$. Continuing in the same way, since $\dir x^2 \in \clarkeTangent_C(x)$, we can now find with $\{x^1_k\}_{k \in \N}$ as the base point sequence points $x^2_k \in C$ such that $(x^2_k-x^1_k)/\tau_k \to \dir x^2$. It follows
    \begin{equation*}
        \frac{x^2_k-x_k}{\tau_k}=\frac{x^2_k-x^1_k}{\tau_k}+\frac{x^1_k-x_k}{\tau_k} \to \dir x^1 + \dir x^2 = \dir x.
    \end{equation*}
    Thus $\{\alt x_k\}_{k\in\N}=\{x^2_k\}_{k\in\N}$ is the sequence we are looking for, showing that $\dir x \in \clarkeTangent_C(x)$ and hence that the Clarke tangent cone is convex.
\end{proof}

\begin{figure}
    \centering
    \begin{asy}
        unitsize(150, 150);
        transform sh=rotate(-30);
        pair p1=(.8, 1);
        pair p2=(1, .2);
        path p=sh*((.6,.5)..p1..(1.1,1.2)..p2..controls(.2, .1)..cycle);
        fill(p, lightfill);
        draw(p);
        p1=sh*p1;
        pair x1=point(p, 0.7)+0.03*(1, -1);
        p2=sh*p2;
        pair tildez=p1+0.3*(1, -.6);
        pair barz=p1+0.8*(.4, -.8);
        pair w1=x1+0.1*(1, -.7);
        pair w2=w1+0.15*(.2, -.7);
        dot(p1);
        label("$x$", p1, N);
        dot(x1);
        label("$x_k$", x1, 2*N+W);
        draw(p1--barz, Arrow);
        label("$\dir x^2$",  barz, E);
        draw(p1--tildez, Arrow);
        label("$\dir x^1$",  tildez, E);

        dot(w1);
        label("$x^1_k$", w1, W+S);
        draw(x1--(x1+3.3*(w1-x1)), dashed, Arrow);

        dot(w2);
        label("$x^2_k$", w2, W+S);
        draw(w1--(w1+5*(w2-w1)), dashed, Arrow);

        draw(x1--(x1+5*(w2-x1)), dashdotted, Arrow);
        draw(p1--((tildez+barz)-p1), dotted, Arrow);
        label("$\dir x^1 + \dir x^2$", (tildez+barz)-p1, E);
    \end{asy}
    \caption{Illustration of the ``bridging'' argument in the proof of \cref{thm:cones:basic-prop}. As $x_k$ converges to $x$, the dashed arrows converge to the solid arrows, while the dash-dotted arrow converges to the dotted one, which depicts the point $\dir x^1_k+\dir x^2_k$ that we are trying to prove to be in $\clarkeTangent_C(x)$.}
    \label{fig:cones:bridging}
\end{figure}

One might expect $T_C^w(x)$ to be weakly closed and $N_C(x)$ to be weak-$*$ closed. However, this is not necessarily the case, since weak and weak-$*$ inner and outer limits need not be closed in the respective topologies. Consequently, $N_C$ may also not be (strong-to-weak-$*$) outer semicontinuous at a point $x$, as this would imply $N_C(x)$ to be weak-$*$ closed and hence closed.
However, in finite dimensions we do have outer semicontinuity.

\begin{corollary}
    \label{cor:cones:continuity}
    If $X$ is finite-dimensional, then the mapping $x \mapsto N_C(x)$ is outer semicontinuous.
\end{corollary}
\begin{proof}
    Let $C\ni x_k \to x$ and $x_k^* \in N_C(x_k)$ with $x_k^* \to x^*$. Then for $\delta_k \downto 0$, the definition \eqref{eq:cones:def-limnormal-findim} provides $\alt x_k \in C$ and $\alt x_k^* \in \frechetNormal_C(\alt x_k)$ with $\norm{\alt x_k^*-x_k^*}\le \delta_k$ and $\norm{\alt x_k-x_k}\le\delta_k$.
    It follows that $C\ni \alt x_k \to x$ and $\alt x_k^* \to x^*$ with $\alt x_k^* \in \frechetNormal_C(\alt x_k)$. Thus by definition, $x^* \in N_C(x)$, and hence $N_C$ is outer semicontinuous.
\end{proof}

\section{Polarity and limiting relationships}
\label{sec:cones:polarity}

The tangent and normal cones satisfy various polarity relationships.
To state these, recall from \cref{sec:functan:dual} for a general set $C \subset X$ the definition of the \term[cone!polar]{polar cone}
\begin{equation*}
    \polar C = \left\{ x^* \in X^* \mid \dualprod{x^*}{x}_X \le 0 \text{ for all } x \in C\right\}
\end{equation*}
as well as of the \term[cone!bipolar]{bipolar cone} $\bipolar C={(\polar C)}_\circ\subset X$.

\subsection*{The fundamental cones}

The relations in the following result will be crucial.
\begin{lemma}
    \label{lemma:cones:fundamental-polar}
    Let $X$ be a Banach space, $C \subset X$, and $x \in X$. Then
    \begin{enumerate}
        \item\label{item:cones:fundamental-polar:weak-incl}
            $\frechetNormal_C(x) \subset \polar{T_C^w(x)} \subset \polar{T_C(x)}$;
        \item\label{item:cones:fundamental-polar:reflexive}
            if $X$ is reflexive, then  $\frechetNormal_C(x) = \polar{T_C^w(x)}$;
        \item\label{item:cones:fundamental-polar:findim}
            if $X$ is finite-dimensional, then
            $\frechetNormal_C(x) = \polar{T_C(x)}$.
    \end{enumerate}
\end{lemma}

\begin{proof}
    \emph{\cref{item:cones:fundamental-polar:weak-incl}:}
    We take $\dir x \in T_C^w(x)$ and $x^* \in \frechetNormal_C(x)$. Then there exist $\tau_k \downto 0$ and $C\ni x_k\to x$ such that $(x_k-x)/\tau_k \weakto \dir x$ weakly in $X$.
    Thus
    \begin{equation*}
        \dualprod{x^*}{\dir x}_X
        =\limsup_{k \to \infty} \frac{\dualprod{x^*}{x_k-x}_X}{\tau_k}
        =\limsup_{k \to \infty} \frac{\dualprod{x^*}{x_k-x}_X}{\norm{x_k-x}_X}\cdot\frac{\norm{x_k-x}_X}{\tau_k}.
    \end{equation*}
    Since $x^* \in \frechetNormal_C(x)$ and $C \ni x_k \to x$, we have by definition that $\limsup_{k \to \infty} \dualprod{x^*}{x_k-x}_X/\norm{x_k-x}_X \le 0$. Moreover, $(x_k-x)/\tau_k \weakto \dir x$ implies that $\norm{x_k-x}_X/\tau_k$ is bounded. Passing to the limit, it therefore follows that $\dualprod{x^*}{\dir x}_X \le 0$.
    Since this holds for every $\dir x \in T_C^w(x)$, we see that $x^* \in \polar{T_C^w(x)}$.
    This shows that $\frechetNormal_C(x) \subset \polar{T_C^w(x)}$.
    Since $T_C(x) \subset T_C^w(x)$ by \cref{thm:cones:inclusions}, $\polar{T_C^w(x)} \subset \polar{T_C(x)}$ follows from \cref{lemma:functan:polar-inclusion}.

    \emph{\cref{item:cones:fundamental-polar:reflexive}:}
    Due to \ref{item:cones:fundamental-polar:weak-incl}, we only need to show “$\supset$”.
    Let $x^* \not\in \frechetNormal_C(x)$. Then, by definition, there exist $C \ni x_k \to x$ with
    \begin{equation}
        \label{eq:graphical:frechet-normal-noninclusion}
        \lim_{k \to \infty} \dualprod{x^*}{\dir x_k} > 0
        \quad\text{for}\quad
        \dir x_k \defeq \frac{x_k-x}{\norm{x_k-x}_X}.
    \end{equation}
    We now use the reflexivity of $X$ and the Eberlein--\u{S}mulyan theorem (\cref{thm:ebsmul}) to pass to a subsequence (not relabelled) such that $\dir x_k \weakto \dir x$ for some $\dir x \in X$ that by definition satisfies $\dir x \in T_C^w(x)$. However, passing to the limit in \eqref{eq:graphical:frechet-normal-noninclusion} now shows that $\dualprod{x^*}{\dir x}_X > 0$ and hence that $x^* \not\in \polar{T_C^w(x)}$.

    \emph{\cref{item:cones:fundamental-polar:findim}:}
    This is immediate from \cref{item:cones:fundamental-polar:reflexive} since $T_C(x)=T_C^w(x)$ in finite-dimensional spaces.
\end{proof}

\subsection*{The limiting cones: preliminary lemmas}

For a polarity relationship between the basic normal cone and the Clarke tangent cone, we need to work significantly harder. We start here with some preliminary lemmas shared between the finite-dimensional and infinite-dimensional setting, and then treat the two in that order.

\begin{lemma}
    \label{lemma:cones:clarke-liminf:preliminary}
    Let $X$ be a reflexive Banach space, $C \subset X$, and $x \in X$. Then
    \begin{equation}
        \label{eq:cones:clarke-liminf:preliminary}
        \clarkeTangent_C(x) \subset \liminf_{C \ni \alt x \to x} T_C^w(\alt x).
    \end{equation}
    If $X=\R^N$, then
    \begin{equation*}
        \clarkeTangent_C(x) \subset \liminf_{C \ni \alt x \to x} T_C(\alt x).
    \end{equation*}
\end{lemma}

\begin{proof}
    The case $X=\R^N$ trivially follows from \eqref{eq:cones:clarke-liminf:preliminary}.
    To prove \eqref{eq:cones:clarke-liminf:preliminary}, denote by $K$ the set on its right-hand side.
    If $\dir x \not \in K$, then there exist $\epsilon>0$ and a sequence $ C\ni x_k \to x$ such that
    \begin{equation}
        \label{eq:cones:clarke-liminf:preliminary.1}
        \inf_{\dir x_k \in T_C^w(x_k)} \norm{\dir x_k-\dir x}_X \ge 3\epsilon.
    \end{equation}
    Fix $k \in \N$ and suppose that for some $\tau_\ell \downto 0$ and $\tilde x_\ell \in C$,
    \begin{equation}
        \label{eq:cones:clarke-liminf:preliminary.2}
        \adaptnorm{\tfrac{\tilde x_\ell-x_k}{\tau_\ell}-\dir x}_X \le 2\epsilon\qquad ( \ell \in \N).
    \end{equation}
    Using the reflexivity of $X$ and the Eberlein--\u{S}mulyan theorem (\cref{thm:ebsmul}), we then find a further (not relabelled) subsequence of $\{(\alt x_\ell, \tau_\ell)\}_{\ell \in \N}$ such that $(\tilde x_\ell-x_k)/\tau_\ell\weakto \dir x_k$ as $\ell \to \infty$ for some $\dir x_k \in T_C^w(x_k)$ with $\norm{\dir x_k-\dir x}_X \le 2\epsilon$, in contradiction to \eqref{eq:cones:clarke-liminf:preliminary.1}.
    We thus have
    \begin{equation*}
        \lim_{\tau \downto 0} ~ \inf_{\tilde x \in C} \adaptnorm{\tfrac{\tilde x-x_k}{\tau}-\dir x}_X \ge 2\epsilon.
    \end{equation*}
    Since this holds for all $k \in \N$, we can find $\tau_k > 0$ with $\tau_k \downto 0$ satisfying the inequality
    \begin{equation*}
        \liminf_{k \to \infty} \inf_{\tilde x \in C} \adaptnorm{\tfrac{\tilde x-x_k}{\tau_k}-\dir x}_X \ge \epsilon
    \end{equation*}
    implying that $\dir x \not \in \clarkeTangent_C(x)$. Therefore \eqref{eq:cones:clarke-liminf:preliminary} holds.
\end{proof}

\begin{lemma}
    \label{lemma:cones:limiting-polar-inclusion}
    Let $X$ be a reflexive and Gateaux smooth (or finite-dimensional) Banach space, $C \subset X$, and $x \in X$. Then
    \begin{equation*}
        \clarkeTangent_C(x) \subset \polar{N_C(x)}.
    \end{equation*}
\end{lemma}

\begin{proof}
    Take $x^* \in N_C(x)$ and $\dir x \in \clarkeTangent_C(x)$. This gives by \cref{thm:cones:nonepsilon-limnormal} (or \eqref{eq:cones:def-limnormal-findim} if $X$ is finite-dimensional) sequences $x_k \to x$ and $x_k^* \weaktostar x^*$ with $x_k^* \in \frechetNormal_C(x_k)$.
    By \cref{lemma:cones:clarke-liminf:preliminary}, we can find for each $k \in \N$ a $\dir x_k \in T_C^w(x_k)$ such that $\dir x_k \to \dir x$.
    Since $\frechetNormal_C(x_k) = \polar{T_C^w(x_k)}$ by \cref{lemma:cones:fundamental-polar}\,\cref{item:cones:fundamental-polar:reflexive} when $X$ is reflexive, we have $\dualprod{x_k^*}{\dir x_k}_X \le 0$.
    Combining all these observations, we obtain
    \begin{equation*}
        \begin{aligned}
            \dualprod{x^*}{\dir x}_X
            &
            =\lim_{k \to\infty}( \dualprod{x_k^*}{\dir x_k}_X + \dualprod{x^*-x_k^*}{\dir x}_X+\dualprod{x_k^*}{\dir x-\dir x_k}_X )
            \\
            &
            =\lim_{k \to \infty} \dualprod{x_k^*}{\dir x_k}_X \le 0.
        \end{aligned}
    \end{equation*}
    Since $x^* \in N_C(x)$ was arbitrary, we deduce that $\dir x \in \polar{N_C(x)}$ and hence the claim.
\end{proof}

\subsection*{The limiting cones in finite dimensions}

\begin{figure}[t]
    \begin{subfigure}[t]{\textwidth}
        \begin{asy}
            unitsize(70, 70);
            pair x=(0, 0);
            pair xi=x+(0.1, -0.1);
            pair z=(-1, 0.3);
            real taui=2;
            real tautilde=1.1;
            real epsilon=.3;
            pair v0=(-.8, -1);
            pair zi=taui*(z+epsilon*v0/length(v0));
            real alpha=(180/pi)*acos(-epsilon/length(z));
            pair v=rotate(alpha)*(z/length(z));
            pair v2=rotate(-alpha)*(z/length(z));
            pair vprime=rotate(20+alpha)*(z/length(z));
            pair ztilde=tautilde*(z+epsilon*vprime);
            pair optx=xi+ztilde;

            path p=(zi+xi+0.3*z)..(zi+xi)..((zi+xi)*3/4+xi/4+0.45*v)..((zi+xi+xi)/2+0.15*v)..(optx)..x..(x-2*z)..((x+xi+zi)/2+4*v)..cycle;
            fill(p, lightfill);
            real y=(xi+taui*z).y;
            real y2=(xi-1.5*z).y;
            clip((-5, y)--(2, y)--(2, y2)--(-5, y2)--cycle);
            picture tmp=new picture;
            draw(tmp, p);
            clip(tmp, (-5, y)--(2, y)--(2, y2+0.09)--(-5, y2+0.09)--cycle);
            add(tmp);

            dot(xi);            label("$x$", xi, E);
            dot(xi+zi);         label("$\tilde x$", zi+xi, W+S);
            dot(optx);     label("$\opt x$", optx, E+S);
            dot(xi+taui*z);     label("$x+z$", xi+taui*z, N);
            dot(xi+tautilde*z); label("$x+\opt\theta z$", xi+tautilde*z, N+0.1E);

            draw(xi--(xi+taui*z));

            draw(shift(xi+taui*z)*scale(epsilon*taui)*unitcircle);
            tmp=new picture;
            draw(tmp, shift(xi+taui*z)*scale(epsilon*taui)*unitcircle, linewidth(1.1));
            clip(tmp, (xi+taui*(z+epsilon*v2))--(xi+tautilde*(z+epsilon*v2))--(xi+tautilde*(z+epsilon*v))--(xi+taui*(z+epsilon*v))--cycle);
            add(tmp);

            draw(shift(xi+tautilde*z)*scale(epsilon*tautilde)*unitcircle, dashed);
            tmp=new picture;
            draw(tmp, shift(xi+tautilde*z)*scale(epsilon*tautilde)*unitcircle, linewidth(1.1));
            clip(tmp, (xi+tautilde*(z+epsilon*v2))--xi--(xi+tautilde*(z+epsilon*v))--cycle);
            add(tmp);

            draw((xi+taui*(z+epsilon*v))--(xi+tautilde*(z+epsilon*v)), dashed);
            draw((xi+taui*(z+epsilon*v2))--(xi+tautilde*(z+epsilon*v2)), dashed);

            draw((xi+zi)--(xi+taui*z), dotted);
            label("$x^*$", (xi+zi)--(xi+taui*z), E);
            draw((optx)--(xi+tautilde*z), dotted);
            label("$\opt x^*$", (optx)--(xi+tautilde*z), W);
        \end{asy}
        \caption{By assumption, the interior of the ball around $x+z$ of radius $\epsilon$ does not intersect $C$ (shaded). In this example, the point $\alt x \in C$ intersects the boundary; however, it is not on the leading edge (thick lines) where the normal vector $x^*$ would satisfy $\iprod{z}{x^*} \ge \epsilon$. Reducing $\theta<1$ produces an intersecting point $\opt x$ on the leading edge.}
        \label{fig:cones:icecream:normal}
    \end{subfigure}
    \begin{subfigure}[t]{\textwidth}
        \begin{asy}
            unitsize(70, 70);
            pair x=(0, 0);
            pair xi=x+(0.1, -0.1);
            pair z=(-1, 0.3);
            real taui=2;
            real tautilde=1.1;
            real taubar=tautilde;
            real epsilon=.3;
            pair v0=(-.8, -1);
            pair zi=taui*(z+epsilon*v0/length(v0));
            real alpha=(180/pi)*acos(-epsilon/length(z));
            pair v=rotate(alpha)*(z/length(z));
            pair v2=rotate(-alpha)*(z/length(z));
            pair vprime=rotate(20+alpha)*(z/length(z));
            pair ztilde=tautilde*(z+epsilon*vprime);
            pair optx=xi+ztilde;
            pair nv=(xi+tautilde*z)-optx;
            pair tv=rotate(90)*(nv/length(nv));
            real tvl0=1.1;
            real tvl=1.2;

            path p=(zi+xi+0.3*z)..(zi+xi)..((zi+xi)*3/4+xi/4+0.45*v)..((zi+xi+xi)/2+0.15*v)..(xi+ztilde)..x..(x-2*z)..((x+xi+zi)/2+4*v)..cycle;
            fill(p, lightfill);
            real y=(xi+taui*z).y;
            real y2=(xi-1.5*z).y;
            picture tmp=new picture;
            draw(tmp, p);
            clip(tmp, (-5, y)--(2, y)--(2, y2+0.09)--(-5, y2+0.09)--cycle);
            add(tmp);
            fill((optx-tvl*tv)--(optx-tvl*tv-nv)--(optx+tvl*tv-nv)--(optx+tvl*tv)--cycle, darkfill);
            draw(optx--(optx-tvl0*tv),primalline+linewidth(1.1), Arrow);
            draw(optx--(optx+tvl0*tv),primalline+linewidth(1.1), Arrow);
            clip((-5, y)--(2, y)--(2, y2)--(-5, y2)--cycle);

            dot(optx);     label("$\opt x$", optx, E+S);
            dot(optx+taubar*z);     label("$\opt x+\opt\theta z$", optx+taubar*z, N);
            label("$T_C(\optx)$", optx-nv, 2*W);

            draw((optx)--(xi+tautilde*z), dotted);
            label("$\opt x^*$", (optx)--(xi+tautilde*z), N+E);

            draw(optx--(optx+taubar*z));
            draw(optx--(optx+taubar*(z+epsilon*v)), dashed);
            draw(optx--(optx+taubar*(z+epsilon*v2)), dashed);
            draw(shift(optx+taubar*z)*scale(epsilon*taubar)*unitcircle, dashed);
        \end{asy}
        \caption{The \enquote{ice cream cone} emanating from $\optx$ along the line $[\optx, \optx+\opt\theta z]$ with a ball of radius $\opt\epsilon\opt\theta$ does not intersect $C$ (light shading). From this it follows that the tangent cone $T_C(\optx)$ (incomplete dark shading) is at a distance $\opt\epsilon$ from $z$.}
        \label{fig:cones:icecream:tangent}
    \end{subfigure}
    \caption{Geometric illustration of the construction in the proof of \cref{lemma:cones:icecream}.}
    \label{fig:cones:icecream}
\end{figure}

We now start our development of polarity relationships between the limiting cones, as well as limiting relationships between the tangent and Clarke tangent cones.
Our main tool will be the following \enquote{ice cream cone lemma}, for which it is important that we endow $\R^N$ with the Euclidean norm.

\begin{lemma}\index{lemma!ice cream cone}
    \label{lemma:cones:icecream}
    Let $C \subset \R^N$ be closed and let $x \in C$.
    Let $z\in \R^N\setminus\{0\}$ and $\epsilon>0$ be such that
    \begin{equation}
        \label{eq:cones:icecream:assumption}
        \interior \B(x+z, \epsilon) \isect C=\emptyset.
    \end{equation}
    Then for any $\bar \epsilon \in (0,\epsilon)$, there exists an $\opt x\in C$ such that there exist
    \begin{enumerate}[label=(\roman*)]
        \item\label{item:cones:icecream:local-and-tangent} $\opt \theta \in (0, 1]$ satisfying $\norm{(\opt x+\opt \theta z)-(x+z)} \le \opt\epsilon$ and
            $\inf_{\dir x \in T_C(\opt x)} \norm{\dir x -z} \ge \bar\epsilon$;
        \item\label{item:cones:icecream:normal}
            $\opt x^* \in \frechetNormal_C(\opt x)$ satisfying $\iprod{\opt x^*}{z} \ge \opt\epsilon$ and $\norm{\opt x^*} \le 1$.
    \end{enumerate}
\end{lemma}

\begin{proof}
    We define the increasing real function $\phi(t) \defeq \sqrt{1+t^2}$ and $F, G: \R^N \times \R \to \Rbar$ by
    \begin{equation*}
        F(\alt x,\theta) \defeq \phi(\opt\epsilon) \theta + \phi(\norm{(\alt x+\theta z)-(x+z)})
        \quad\text{and}\quad
        G(\alt x,\theta) \defeq \delta_C(\alt x) + \delta_{[0, \infty)}(\theta).
    \end{equation*}
    Then $F+G$ is proper, coercive, and lower semicontinuous and hence admits a minimizer $(\opt x, \opt \theta)\in C\times [0,\infty)$ by \cref{thm:variation:existence}. (We illustrate the idea of such a minimizer geometrically in \cref{fig:cones:icecream}.)
    Let $\bar y\defeq (\opt x+\opt\theta z)-(x+z)$.

    \emph{\ref{item:cones:icecream:local-and-tangent}:}
    We first prove $\opt \theta \in (0, 1]$. Suppose $\opt \theta = 0$. Since $\opt x \in C$ we obtain using \eqref{eq:cones:icecream:assumption} that
    \begin{equation*}
        [F+G](\opt x, 0)
        = \phi(\norm{\opt x-(x+z)})
        \ge \phi(\epsilon)
        > \phi(\opt\epsilon)
        = [F+G](x, 1).
    \end{equation*}
    This is a contradiction to $(\opt x, 0)$ being a minimizer. Thus $\opt\theta \ne 0$.
    Likewise,
    \begin{equation*}
        \phi(\opt\epsilon)\opt\theta
        + \phi(\norm{\opt y})
        = [F+G](\opt x, \opt\theta)
        \le [F+G](x, 1)=\phi(\opt\epsilon),
    \end{equation*}
    where both terms on the left-hand side are nonnegative.
    Hence $\opt\theta \le 1$. By the monotonicity of $\phi$, this also verifies the claim $\norm{\opt y} \le \opt\epsilon$.

    We still need to prove the claim on the tangent cone.
    Since $(\opt x, \opt \theta)$ is a minimizer of $F+G$, for any $\alt\theta \ge 0$ and $\alt x \in C$ we have
    \begin{equation*}
        \phi(\opt\epsilon)\opt\theta
        + \phi(\norm{\opt y})
        \le
        [F+G](\optx,\opt\theta)
        \le [F+G](\alt x,\alt\theta)
        = \phi(\opt\epsilon)\alt\theta
        + \phi(\norm{y}).
    \end{equation*}
    Letting $y \defeq (\alt x+\alt\theta z)-(x+z)$ and using first this inequality and then the convexity of $\phi$ with $\phi'(t)=t/\phi(t) \le 1$ for all $t \ge 0$ yields
    \begin{equation*}
        \begin{aligned}
            \phi(\opt\epsilon)(\opt\theta-\alt\theta)
            &
            \le \phi(\norm{\opty}) - \phi(\norm{y}) \\
            &
            \le \phi'(\norm{\opty})\left(\norm{\opty}-\norm{y}\right) \\
            &
            \le \norm{\opty-y}
            = \norm{\alt x-\opt x-(\opt\theta-\alt\theta)z}.
        \end{aligned}
    \end{equation*}
    Dividing by $\tau=\opt\theta-\alt\theta$ for $\alt\theta \in [0, \opt\theta)$, we obtain that $\opt\epsilon \le \phi(\opt\epsilon) \le \adaptnorm{\tfrac{\alt x - \optx}{\tau}-z}$. Taking the infimum over $\alt x \in C$ and $\tau \in (0, \opt\theta]$ thus yields $\inf_{\dir x \in T_C(\opt x)} \norm{\dir x -z} \ge \bar\epsilon$.

    \emph{\ref{item:cones:icecream:normal}:}
    By \cref{lem:convex:func}\,\cref{lem:convex:func:iv}, $F$ is convex. Furthermore, $\interior(\dom F)=\R^{N+1}$ so that $F$ is Lipschitz near $(\opt x, \opt \theta)$ by \cref{thm:convex:cont}.
    Using \cref{thm:subdifferential:norm,thm:convex:increasing-post,thm:convex:chain} with $K(x,\theta) \defeq x+\theta z$, it follows that
    \begin{align}
        \label{eq:cones:icecream:subdiff-f}
        \subdiff F(\opt x, \opt \theta) & = \left\{
            \begin{pmatrix}\phi'(\norm{\opt y})y^*\\ \phi(\opt\epsilon) + \phi'(\norm{\opt y})\dualprod{z}{y^*}\end{pmatrix}
            \ \middle|\,
            \begin{array}{lr}
                \dualprod{y^*}{\opt y}=\norm{\opt y}, \norm{y^*}  = 1 & \text{if }\opt y\neq 0 \\
                \norm{y^*}\leq 1 & \text{if }\opt y = 0
            \end{array}
        \right\}.
    \end{align}
    Since $\R^N$ endowed with the Euclidean norm is a Hilbert space, $x\mapsto\norm{x}^2$ is Gateaux differentiable by \cref{ex:epsilon:gateaux}\,\cref{item:epsilon:gateaux:hilbert} and \cref{lem:epsilon:squared-norm}.
    Hence $\subdiff F(\opt x, \opt \theta)$ is a singleton, and therefore $F$ is Gateaux differentiable at $(\optx, \opt\theta)$ due to \cref{lem:clarke:gateaux,thm:clarke:convex}.
    We can thus apply the Fermat principle (\cref{thm:limiting:frechet:fermat}) and the Fréchet sum rule (\cref{cor:epsilon:sumrule:gateaux}) to deduce $0 \in \subdiff_F F(\opt x, \opt \theta) + \subdiff_F G(\opt x, \opt \theta)$.
    Since $\opt \theta > 0$, we have $\subdiff_F G(\opt x, \opt \theta) = \frechetNormal_C(\opt x) \times \{0\}$ by \cref{lemma:cones:frechet-subdiff}, which implies that
    \begin{equation}
        \label{eq:cones:icecream:expanded-oc}
        -\phi'(\norm{\opt y})y^* \in \frechetNormal_C(\opt x) \quad\text{and}\quad
        \phi(\opt\epsilon) + \phi'(\norm{\opt y})\dualprod{z}{y^*}=0.
    \end{equation}
    Since $\phi(\opt\epsilon)>0$, the second equation in \eqref{eq:cones:icecream:expanded-oc} yields $\phi'(\norm{\opt y})\neq 0$ as well.
    As $\phi'(t)\in(0,1)$ and $\phi(t)>t$ for all $t>0$, we can set $x^* \defeq -\phi'(\norm{\opt y})y^*$ to obtain $x^* \in \frechetNormal_C(\opt x)$ with $\norm{x^*} \le 1$ and
    $
    \dualprod{z}{x^*}
    = \frac{\phi(\opt\epsilon)}{\phi'(\norm{\opt y})}
    \ge \phi(\opt\epsilon)
    \ge \opt\epsilon.
    $
\end{proof}

The following consequence of the ice cream cone lemma will be useful for several polarity relations. We call a set $C$ \term[set!closed near a point]{closed near} $x\in C$, if there exists a $\delta>0$ such that $C\cap \B(x,\delta)$ is closed.
\begin{lemma}
    \label{lemma:cones:limiting-polar:findim}
    Let $C\subset \R^N$ be closed near $x$.
    If $z \not \in \clarkeTangent_C(x)$, then there exist $\alt\eps>0$ and a sequence $C \ni \alt x_k \to x$ such that for all $k \in \N$,
    \begin{enumerate}[label=(\roman*)]
        \item\label{item:cones:limiting-polar:findim:tangent}
            $\inf_{\dir \alt x_k \in T_C(\alt x_k)}\norm{\dir \alt x_k-z} \ge \alt\eps$;
        \item\label{item:cones:limiting-polar:findim:normal}
            there exists $\alt x^*_k \in \frechetNormal_C(\alt x_k)$ with $\norm{\alt x^*_k} \le 1$ and $\iprod{\alt x_k^*}{z} \ge \alt\eps$.
    \end{enumerate}
\end{lemma}

\begin{proof}
    First, $z \not \in \clarkeTangent_C(x)$ implies by \eqref{eq:cones:def-clarketangent} the existence of $\epsilon>0$, $C \ni x_k \to x$, and $\tau_k \downto 0$ such that
    \begin{equation*}
        \inf_{\alt x \in C} \adaptnorm{\tfrac{\alt x-x_k}{\tau_k}-z} \ge \epsilon
        \qquad (k \in \N),
    \end{equation*}
    implying that
    \begin{equation*}
        \interior \B(x_k+\tau_k z, \tau_k \epsilon) \isect C = \emptyset.
    \end{equation*}
    By taking $\tau_k$ small enough -- i.e., $k \in \N$ large enough -- we may without loss of generality assume that $C$ is closed.
    For any $\alt\epsilon \in (0, \epsilon)$ and every $k \in \N$, \cref{lemma:cones:icecream} now yields $\alt x_k \in C$ and $\alt\theta_k \in (0, 1]$ satisfying
    \begin{enumerate}[label=(\roman*$'$)]
        \item\label{item:cones:limiting-polar:findim:tangent:prime}
            $\norm{(\alt x_k+\alt\theta_k \tau_k z)-(x+\tau_k z)} \le \alt\epsilon\tau_k$ and
            $\inf_{\dir \alt x_k \in T_C(\alt x_k)}\norm{\dir \alt x_k-\tau_k z} \ge \alt\eps\tau_k$;
        \item\label{item:cones:limiting-polar:findim:normal:prime}
            there exists an $\alt x_k^* \in \frechetNormal_C(\alt x_k)$ such that $\dualprod{\alt x_k^*}{\tau_k z} \ge \tau_k \alt\epsilon$ and $\norm{\alt x_k^*} \le 1$.
    \end{enumerate}
    We readily obtain \ref{item:cones:limiting-polar:findim:tangent} from \ref{item:cones:limiting-polar:findim:tangent:prime} and \ref{item:cones:limiting-polar:findim:normal} from \ref{item:cones:limiting-polar:findim:normal:prime}.
    Since \ref{item:cones:limiting-polar:findim:tangent:prime} also shows that $\alt x_k \to x$ as $\tau_k \downto 0$, this finishes the proof.
\end{proof}

We can now show the converse inclusion of \cref{lemma:cones:limiting-polar-inclusion} when the set is closed near $x$.

\begin{theorem}
    \label{thm:cones:limiting-polar:findim}
    If $C\subset \R^N$ is closed near $x$, then
    \begin{equation*}
        \clarkeTangent_C(x) = \polar{N_C(x)}.
    \end{equation*}
\end{theorem}

\begin{proof}
    By \cref{lemma:cones:limiting-polar-inclusion}, we only need to prove $ \clarkeTangent_C(x) \supset \polar{N_C(x)}$. We argue by contraposition.  Let $z \not\in \clarkeTangent_C(x)$. Then \cref{lemma:cones:limiting-polar:findim} yields a sequence $\{x_k^*\}_{k\in\N}\subset \R^N$ such that $x_k^*\in \frechetNormal_C(x_k)$ for $C\ni x_k\to x$ and $\dualprod{x_k^*}{z} \geq \eps>0$ as well as $\norm{x_k^*}\leq 1$. Since $\{x_k^*\}_{k \in \N}$ is bounded, we can extract a subsequence that converges to some $x^* \in \R^N$. By definition of the limiting normal cone, $x^*\in N_C(x)$. Moreover, $\dualprod{x^*}{z} \geq \eps>0$. This provides, as required, that $z \not\in \polar{N_C(x)}$.
\end{proof}

\subsection*{The limiting cones in infinite dimensions}

We now repeat the arguments above in infinite dimensions, however, we need extra care and extra assumptions. Besides reflexivity (to obtain weak-$*$ compactness from the Eberlein--\u{S}mulyan theorem (\cref{thm:ebsmul})) and Gateaux smoothness (to obtain differentiability of the norm), we need to use
the approximate Fermat principle of \cref{thm:epsilon:approximate-fermat} since exact projections to general sets $C$ may not exist; compare \cref{thm:epsilon:projection}. This introduces $\epsilon$-normal cones into the proof.
The geometric ideas of the proof, however, are the same as illustrated in \cref{fig:cones:icecream}.

\begin{lemma}\index{lemma!ice cream cone}
    \label{lemma:cones:icecream:infdim}
    Let $X$ be a Banach space, $C \subset X$ be closed, and $x \in C$.
    Let $z\in X\setminus\{0\}$ and $\epsilon>0$ be such that
    \begin{equation}
        \label{eq:cones:icecream:infdim:assumption}
        \interior \B(x+z, \epsilon) \isect C=\emptyset.
    \end{equation}
    Then for any $\bar \epsilon \in (0,\epsilon)$ and $\rho>0$, there exists $\opt x\in C$ such that there exist
    \begin{enumerate}[label=(\roman*)]
        \item\label{item:cones:icecream:infdim:local-and-tangent} $\opt \theta \in (0, 1]$ such that $\norm{(\opt x+\opt \theta z)-(x+z)}_X \le \opt\epsilon$ and
            $\inf_{\dir x \in T_C(\opt x)} \norm{\dir x -z}_X \ge \bar\epsilon$;
        \item\label{item:cones:icecream:infdim:normal}
            if $X$ is Gateaux smooth, $\opt x^* \in \frechetNormal_C^\rho(\opt x)$ such that $\dualprod{\opt x^*}{z}_X \ge \opt\epsilon$ and $\norm{\opt x^*}_{X^*} \le 1$.
    \end{enumerate}
\end{lemma}

\begin{proof}
    We define the convex and increasing real function $\phi(t) \defeq \sqrt{1+t^2}$ and pick arbitrary
    \begin{equation}
        \label{eq:cones:icecream:infdim:rho-delta}
        \alt\epsilon \in (\opt\epsilon,\epsilon),
        \quad
        0 < \rho < \frac{\phi(\alt\epsilon)-\opt\epsilon}{2+\opt\epsilon},
        \quad\text{and}\quad
        0 < \delta < \phi(\epsilon)-\phi(\alt\epsilon).
    \end{equation}
    The latter upper bound on $\rho$ does not affect the generality of \ref{item:cones:icecream:infdim:normal} since $\frechetNormal_C^\rho(\opt x) \subset \frechetNormal_C^{\rho'}(\opt x)$ for $\rho' \ge \rho$.
    Then we define $F, G: X \times \R \to \Rbar$ by
    \begin{equation*}
        F(\alt x,\theta) \defeq \phi(\alt\epsilon) \theta + \phi(\norm{(\alt x+\theta z)-(x+z)}_X)
        \quad\text{and}\quad
        G(\alt x,\theta) \defeq \delta_C(\alt x) + \delta_{[0, \infty)}(\theta).
    \end{equation*}
    The function $F+G$ is proper and coercive, hence $\inf(F+G)>-\infty$. However, it may not admit a minimizer. Nevertheless, the approximate Fermat principle of \cref{thm:epsilon:approximate-fermat} produces an approximate minimizer $(\opt x, \opt \theta)\in C\times [0,\infty)$ with
    \begin{enumerate}[label=(\alph*)]
        \item\label{item:cones:icecream:infdim:ekeland:bound}
            $[F+G](\opt x, \opt \theta) \le \inf[F+G]+\delta$,
        \item\label{item:cones:icecream:infdim:ekeland:minimising}
            $[F+G](\opt x, \opt \theta) < [F+G](\alt x, \theta)+\rho\norm{\alt x-\opt x}_X+\rho\abs{\theta-\opt\theta}$ for all $(\alt x, \theta) \ne (\opt x, \opt \theta)$, and
        \item\label{item:cones:icecream:infdim:ekeland:subdiff}
            $0 \in \subdiff_{\rho}[F+G](\opt x, \opt \theta)$.
    \end{enumerate}
    Let again $\bar y\defeq (\opt x+\opt\theta z)-(x+z)$.

    \emph{\ref{item:cones:icecream:infdim:local-and-tangent}:}
    We first prove $\opt \theta \in (0, \tfrac{\phi(\epsilon)}{\phi(\alt\epsilon)}]$, which will
    in particular imply that $\opt\theta \in (0, 1+\epsilon)$.
    Suppose $\opt \theta =0$.
    Since $\opt x \in C$, using \eqref{eq:cones:icecream:infdim:assumption} and the convexity of $\phi$, we obtain
    \begin{equation*}
        [F+G](\opt x, 0) - \delta
        = \phi(\norm{\opt x-(x+z)}_X)  - \delta
        \ge \phi(\eps) - \delta
        > \phi(\alt\epsilon)
        = [F+G](x, 1)
    \end{equation*}
    in contradiction to \ref{item:cones:icecream:infdim:ekeland:bound}. Thus $\opt\theta \ne 0$.
    Likewise,
    \begin{equation*}
        \phi(\alt\epsilon)\opt\theta
        + \phi(\norm{\opt y}_X)
        = [F+G](\opt x, \opt\theta)
        \le [F+G](x, 1) + \delta = \phi(\alt\epsilon) + \delta < \phi(\epsilon).
    \end{equation*}
    where both terms on the left-hand side are nonnegative.
    Hence $\opt\theta \le \tfrac{\phi(\epsilon)}{\phi(\alt\epsilon)}$. By monotonicity of $\phi$, this also verifies the claim $\norm{\opt y}_X \le \epsilon$.

    We still need to prove the claim on the tangent cone.
    Letting $y \defeq (\alt x+\theta z)-(x+z)$, we rearrange \ref{item:cones:icecream:infdim:ekeland:minimising} as
    \begin{equation}
        \label{eq:cones:icecream:infdim:tangent-proof-ineq}
        \phi(\alt\epsilon)(\opt\theta-\theta)  - \rho\abs{\theta-\opt\theta}
        \le
        \phi(\norm{\opty}_X) - \phi(\norm{y}_X) +  \rho\norm{\alt x-\opt x}_X
    \end{equation}
    Using the convexity of $\phi$, we also have
    \begin{equation*}
        \phi(\norm{\opty}_X) - \phi(\norm{y}_X)
        \le \frac{1}{\phi(\norm{\opty}_X)}\left(\norm{\opty}_X-\norm{y}_X\right)
        \le \norm{\opty-y}_X
        = \norm{\alt x-\opt x-(\opt\theta-\theta)z}_X
    \end{equation*}
    Further estimating $\norm{\alt x-\opt x}_X \le \norm{\alt x-\opt x-(\opt\theta-\theta)z}_X+\abs{\opt\theta-\theta}$, \eqref{eq:cones:icecream:infdim:tangent-proof-ineq} now yields
    \begin{equation*}
        [\phi(\alt\epsilon)-2\rho](\opt\theta-\theta)
        \le
        (1+\rho)\norm{\alt x-\opt x-(\opt\theta-\theta)z}_X
        \qquad (\theta \in [0, \opt\theta),\, \alt x \in C).
    \end{equation*}
    Dividing by $(1+\rho)(\opt\theta-\theta)$ and using  \eqref{eq:cones:icecream:infdim:rho-delta} (for the first inequality), we obtain that
    \begin{equation*}
        \opt\epsilon \le
        \frac{\phi(\alt\epsilon)-2\rho}{1+\rho}
        \le \inf_{\alt x \in C,\, \theta \in [0, \opt\theta)} \adaptnorm{\frac{\alt x-\opt x}{\opt\theta-\theta}-z}_X.
    \end{equation*}
    This shows $\inf_{\dir x \in T_C(\opt x)} \norm{\dir x -z}_X \ge \bar\epsilon$.

    \emph{\ref{item:cones:icecream:infdim:normal}:}
    By \cref{lem:convex:func}\,\cref{lem:convex:func:iv}, $F$ is convex. Furthermore $\interior(\dom F)=X \times \R$, and hence $F$ is Lipschitz near $(\opt x, \opt \theta)$ by \cref{thm:convex:cont}.
    Using \cref{thm:subdifferential:norm,thm:convex:increasing-post,thm:convex:chain} with $K(x,\theta) \defeq x+\theta z$, it follows that
    \begin{align}
        \label{eq:cones:icecream:infdim:subdiff-f}\!\!
        \subdiff F(\opt x, \opt \theta) & = \left\{
            \begin{pmatrix}\phi'(\norm{\opt y}_X)y^*\\ \phi(\alt\epsilon) + \phi'(\norm{\opt y}_X)\dualprod{z}{y^*}_X\end{pmatrix}
            \ \middle|\,
            \begin{array}{lr}
                \dualprod{y^*}{\opt y}_X=\norm{\opt y}_X, \norm{y^*}_{X^*}  = 1 & \text{if }\opt y\neq 0 \\
                \norm{y^*}_{X^*}\leq 1 & \text{if }\opt y = 0
            \end{array}
        \right\}.\!\!\!\!
    \end{align}
    Again, $\subdiff F(\opt x, \opt \theta)$ is a singleton by \cref{lem:epsilon:squared-norm} and the assumption that $X$ is Gateaux smooth.

    We can thus apply the $\eps$-sum rule (\cref{lemma:epsilon:sumrule}) in \ref{item:cones:icecream:infdim:ekeland:subdiff} to deduce $0 \in \subdiff F(\opt x, \opt \theta) + \subdiff_{\rho} G(\opt x, \opt \theta)$.
    Since $\opt \theta > 0$, we have $\subdiff_{\rho} G(\opt x, \opt \theta) = \frechetNormal_C^{\rho}(\opt x) \times \{0\}$, which implies that
    \begin{equation}
        \label{eq:cones:icecream:infdim:expanded-oc}
        -\phi'(\norm{\opt y}_X)y^* \in \frechetNormal_C^{\rho}(\opt x) \quad\text{and}\quad
        \phi(\alt\epsilon) + \phi'(\norm{\opt y}_X)\dualprod{z}{y^*}_X=0.
    \end{equation}
    Since $\phi(\alt\epsilon)>0$, the second equation in \eqref{eq:cones:icecream:infdim:expanded-oc} yields $\phi'(\norm{\opt y}_X)\neq 0$ as well.
    As $\phi'(t)\in(0,1)$ and $\phi(t)>t$ for all $t>0$, we can set $x^* \defeq -\phi'(\norm{\opt y}_X)y^*$ to obtain $x^* \in \frechetNormal_C^{\rho}(\opt x)$ with $\norm{x^*}_{X^*} \le 1$ and
    $
    \dualprod{z}{x^*}_X
    = \frac{\phi(\alt\epsilon)}{\phi'(\norm{\opt y}_X)}
    \ge \phi(\alt\epsilon)
    \ge \opt\epsilon.
    $
\end{proof}

\begin{remark}
    If $X$ is in addition reflexive, we can use the Eberlein--\u{S}mulyan theorem (\cref{thm:ebsmul}) to pass to the limit as $\rho \downto 0$ in \cref{lemma:cones:icecream:infdim} and produce $\opt x^* \in \frechetNormal_C(\opt x)$ satisfying the other claims of the lemma.
\end{remark}

\begin{lemma}
    \label{lemma:cones:limiting-polar}
    Let $X$ be a Banach space and $C\subset X$ be closed near $x \in C$.
    If $z \not \in \clarkeTangent_C(x)$, then there exist $\alt\eps>0$ and a sequence $C \ni \alt x_k \to x$ such that for all $k \in \N$,
    \begin{enumerate}[label=(\roman*)]
        \item\label{item:cones:limiting-polar:tangent}
            $\inf_{\dir \alt x_k \in T_C(\alt x_k)}\norm{\dir \alt x_k-z}_X \ge \alt\eps$;
        \item\label{item:cones:limiting-polar:normal}
            if $X$ is Gateaux smooth, for any $\rho_k>0$ there exists $\alt x^*_k \in \frechetNormal_C^{\rho_k}(\alt x_k)$ with $\norm{\alt x^*_k}_{X^*} \le 1$ and $\dualprod{\alt x_k^*}{z}_{X^*} \ge \alt\eps$.
    \end{enumerate}
\end{lemma}

\begin{proof}
    The assumption $z \not \in \clarkeTangent_C(x)$ implies by \eqref{eq:cones:def-clarketangent} the existence of $\epsilon>0$, $C \ni x_k \to x$, and $\tau_k \downto 0$ such that
    \begin{equation*}
        \inf_{\alt x \in C} \adaptnorm{\tfrac{\alt x-x_k}{\tau_k}-z}_X \ge \epsilon
        \qquad (k \in \N).
    \end{equation*}
    This implies that
    \begin{equation*}
        \interior \B(x_k+\tau_k z, \tau_k \epsilon) \isect C = \emptyset.
    \end{equation*}
    Since the argument is local, by taking $\tau_k$ small enough -- i.e., $k \in \N$ large enough -- we may without loss of generality assume that $C$ is closed.
    For any $\alt\epsilon \in (0, \epsilon)$ and $\rho_k > 0$, \cref{lemma:cones:icecream:infdim} now produces $\alt x_k \in C$ and $\alt\theta_k \in (0, 1]$ satisfying
    \begin{enumerate}[label=(\roman*$'$)]
        \item\label{item:cones:limiting-polar:tangent:prime}
            $\norm{(\alt x_k+\alt\theta_k \tau_k z)-(x+\tau_k z)}_X \le \alt\epsilon\tau_k$ and  $\inf_{\dir \alt x_k \in T_C(\alt x_k)}\norm{\dir \alt x_k-\tau_k z}_X \ge \alt\eps\tau_k$;
        \item\label{item:cones:limiting-polar:normal:prime}
            if $X$ is Gateaux smooth, there exists $\alt x_k^* \in \frechetNormal_C^{\rho_k}(\alt x_k)$ such that $\dualprod{x_k^*}{\tau_k z}_X \ge \tau_k \alt\epsilon$ and $\norm{\alt x_k^*}_{X^*} \le 1$.
    \end{enumerate}
    We readily obtain \ref{item:cones:limiting-polar:tangent} from \ref{item:cones:limiting-polar:tangent:prime} and \ref{item:cones:limiting-polar:normal} from \ref{item:cones:limiting-polar:normal:prime}.
    Since \ref{item:cones:limiting-polar:tangent:prime} also shows that $\alt x_k \to x$ as $\tau_k \downto 0$, this finishes the proof.
\end{proof}

\begin{theorem}
    \label{thm:cones:limiting-polar}
    Let $X$ be a reflexive and Gateaux smooth Banach space and let $C\subset X$ be closed near $x \in C$. Then
    \begin{equation*}
        \clarkeTangent_C(x) = \polar{N_C(x)}.
    \end{equation*}
\end{theorem}
\begin{proof}
    By \cref{lemma:cones:limiting-polar-inclusion}, we only need to prove $ \clarkeTangent_C(x) \supset \polar{N_C(x)}$. Let $z \not\in \clarkeTangent_C(x)$ and $\rho_k \downto 0$. Then \cref{lemma:cones:limiting-polar}\,\cref{item:cones:limiting-polar:normal} yields a sequence $\{\alt x_k^*\}_{k \in \N}\subset \B_{X^*}$ such that $\alt x_k^*\in \frechetNormal_C^{\rho_k}(x_k)$ and $\dualprod{\alt x^*_k}{z}_X \geq \eps$. Since $X$ is reflexive, $X^*$ is reflexive as well, and so we can apply \cref{thm:ebsmul} to extract a subsequence of $\{\alt x_k^*\}_{k \in \N}$ that converges weakly and thus, again by reflexivity, also weakly-$*$ to some $x^* \in N_C(x)$ (by definition \eqref{eq:cones:def-limnormal} of the limiting normal cone) with $\dualprod{x^*}{z}_X\geq \eps>0$.
\end{proof}

\subsection*{The Clarke tangent cone}

We can now show the promised alternative characterization of the Clarke tangent cone $\clarkeTangent_C(x)$ as the inner limit of tangent cones.
\begin{corollary}
    \label{cor:cones:clarke-liminf}
    Let $X$ be a reflexive Banach space and let $C \subset X$ be closed near $x\in X$. Then
    \begin{equation}
        \label{eq:cones:clarke-liminf}
        \liminf_{C \ni \alt x \to x} T_C(\alt x) \subset
        \clarkeTangent_C(x) \subset \liminf_{C \ni \alt x \to x} T^w_C(\alt x).
    \end{equation}
    In particular, if $X$ is finite-dimensional, then
    \begin{equation*}
        \clarkeTangent_C(x) = \liminf_{C \ni \alt x \to x} T_C(\alt x).
    \end{equation*}
\end{corollary}
\begin{proof}
    We have already proved the second inclusion of \eqref{eq:cones:clarke-liminf} in \cref{lemma:cones:clarke-liminf:preliminary}.
    For the first inclusion, suppose $z \not\in \clarkeTangent_C(x)$. Then \cref{lemma:cones:limiting-polar}\,\cref{item:cones:limiting-polar:tangent} yields an $\alt\epsilon>0$ and a sequence $C \ni \alt x_k \to x$ such that $\inf_{\dir \alt x_k \in T_C(\alt x_k)}\norm{\dir \alt x_k-z}_X \ge \alt\eps$ for all $k$.
    This shows that $z \not \in \liminf_{C \ni \alt x \to x} T_C(\alt x)$.
\end{proof}

\begin{remark}
    \Cref{lemma:cones:limiting-polar} and thus the first inclusion of \eqref{eq:cones:clarke-liminf} do not actually require the reflexivity of $X$.
    In contrast, \cref{lemma:cones:clarke-liminf:preliminary} and thus the second inclusion of \eqref{eq:cones:clarke-liminf} do not require the local closedness assumption.
    Besides in finite-dimensional spaces, the claimed equality holds more generally if $X$ has the Radon--Riesz property and is Fréchet smooth; see \cite[Theorem 1.9]{Mordukhovich:2006} and compare \cref{rem:epsilon:radonriesz}.
\end{remark}

\section{Regularity}
\label{sec:cones:regularity}

It stands to reason that without any assumptions on the set $C\subset X$ such as convexity, there is little hope of obtaining precise characterizations or exact transformation rules for the various cones. Similarly, precise characterizations or exact calculus rules for the derivatives of set-valued mappings -- which, respectively, we will derive from the former -- require strong assumptions on these mappings. This is especially true of the limiting cones.
As befitting the introductory character of this textbook, we will therefore only develop calculus for the derivatives based on the limiting cones when they are equal to the corresponding basic cones.
This will allow deriving exact results that are nevertheless applicable to the situations we have been focusing on in the previous parts, such as problems of the form \eqref{eq:intro:prob}. These conditions can be compared to constraint qualifications in nonlinear optimization that guarantee that the tangent cone coincides with the linearization cone. However, \enquote{fuzzy} results are available under more general assumptions, for which we refer to the monographs \cite{aubin1990sva,Rockafellar:1998,mordukhovich2018variational,Mordukhovich:2006}.

Specifically, we say that $C \subset X$ is \term[set!regular!tangentially]{tangentially regular} at $x \in C$ if $T_C(x)=\clarkeTangent_C(x)$, and \term[set!regular!normally]{normally regular} at $x$ if $N_C(x)=\frechetNormal_C(x)$.
We call $C$ \term[set!regular]{regular} at $x$ if $C$ is both normally and tangentially regular.

\begin{example}
    Continuing from \cref{ex:cones:basicex}, we see that $C=\B(0,1)$ and $C=[0,1]^2$ are regular at every $x\in C$, while $C=[0,1]^2\setminus[\frac12,1]^2$ is regular everywhere except at $x=(\frac12,\frac12)$.
\end{example}

\begin{example}
    The set $C\subset\R^2$ shown in \cref{fig:cones:tangent-normal} is regular at every point except at the \enquote{kink} where the tangent and normal cones are illustrated.
    In that point, it is neither tangentially nor normally regular: neither $T_C(x)$ and $\clarkeTangent(x)$ nor $N_C(x)$ and $\frechetNormal_C(x)$ coincide.
\end{example}

In finite dimensions, the two concepts of regularity are equivalent and have various characterizations. By \cref{lemma:cones:convex}, these hold in particular for closed convex sets.
\begin{theorem}
    \label{thm:cones:regularity:findim}
    Let $C \subset \R^N$ be closed near $x$. Then the following conditions are equivalent:
    \begin{enumerate}
        \item\label{item:cones:regularity:i} $C$ is normally regular at $x$;
        \item\label{item:cones:regularity:ii} $C$ is tangentially regular at $x$;
        \item\label{item:cones:regularity:iv} $\frechetNormal_C$ is outer semicontinuous at $x$;
        \item\label{item:cones:regularity:iii} $T_C$ is inner semicontinuous at $x$ (relative to $C$).
    \end{enumerate}
    In particular, if any of these hold, $C$ is regular at $x$.
\end{theorem}

\begin{proof}
    \emph{\ref{item:cones:regularity:i} $\equivalent$ \ref{item:cones:regularity:ii}:}
    If \ref{item:cones:regularity:i} holds, then by \cref{lemma:functan:polar-inclusion,thm:cones:inclusions,lemma:cones:fundamental-polar,thm:cones:limiting-polar:findim}
    \begin{equation*}
        T_C(x) \subset \bipolar{T_C(x)}
        =\polar{\frechetNormal_C(x)}
        =\polar{N_C(x)}
        =\clarkeTangent_C(x)
        \subset T_C(x),
    \end{equation*}
    which shows \ref{item:cones:regularity:ii}.
    The other direction is completely analogous, exchanging the roles of ``$N$'' and ``$T$'' to obtain
    \begin{equation*}
        N_C(x) \subset \bipolar{N_C(x)}
        =\polar{\clarkeTangent_C(x)}
        =\polar{T_C(x)}
        =\frechetNormal_C(x)
        \subset N_C(x).
    \end{equation*}

    \emph{\ref{item:cones:regularity:i} $\equivalent$ \ref{item:cones:regularity:iv}:}
    If \ref{item:cones:regularity:i} holds, then the outer semicontinuity of $N_C$ (\cref{cor:cones:continuity}) and the inclusion $\frechetNormal_C(\alt x) \subset N_C(\alt x)$ from \cref{thm:cones:inclusions} show that $\limsup_{\alt x \to x} \frechetNormal_C(\alt x) \subset \frechetNormal_C(x)$, i.e., the outer semicontinuity of $\frechetNormal_C$.
    Conversely, the outer semicontinuity of $\frechetNormal_C$ and the definition $N_C(x)=\limsup_{\alt x \to x} \frechetNormal_C(\alt x)$ show that $N_C(x) \subset \frechetNormal_C(x)$. Combined with the inclusion $\frechetNormal_C(\alt x) \subset N_C(\alt x)$ from \cref{thm:cones:inclusions}, we obtain \ref{item:cones:regularity:i}.

    \emph{\ref{item:cones:regularity:ii} $\equivalent$ \ref{item:cones:regularity:iii}:}
    To show that \ref{item:cones:regularity:iii} implies \ref{item:cones:regularity:ii}, recall from \cref{cor:cones:clarke-liminf} that
    \begin{equation}
        \label{eq:cones:regularity:clarke-liminf}
        \clarkeTangent_C(x) = \liminf_{C \ni \alt x \to x} T_C(\alt x).
    \end{equation}
    By the assumed inner semicontinuity and the definition of the inner limit, we thus obtain that $T_C(x) = \liminf_{C \ni \alt x \to x} T_C(\alt x)=\clarkeTangent_C(x)$.
    For the other direction, we simply use $\clarkeTangent_C(x) = T_C(x)$ in \eqref{eq:cones:regularity:clarke-liminf}.
\end{proof}

Combining the previous result with \cref{lemma:cones:fundamental-polar,thm:cones:limiting-polar:findim}, we deduce the following.
\begin{corollary}
    \label{cor:cones:regularity}
    If $C \subset \R^N$ is regular at $x$ and closed near $x$, then both $T_C(x)$ and $N_C(x)$ are convex. Furthermore,
    \begin{enumerate}
        \item $N_C(x)=\polar{T_C(x)}$;
        \item $T_C(x)=\polar{N_C(x)}$.
    \end{enumerate}
\end{corollary}

In infinite dimensions, our main equivalent characterization of normal regularity is the following.
(We do not have a similar characterization of tangential regularity.)

\begin{theorem}%
    \label{thm:cones:regularity:infdim}
    Let $X$ be a reflexive and Gateaux smooth Banach space. Then $C\subset X$ is normally regular at $x\in C$ if and only if $\clarkeTangent_{C}(x) = \polar{\frechetNormal_{C}(x)}$.
\end{theorem}

\begin{proof}

    Suppose first that $\clarkeTangent_{C}(x) = \polar{\frechetNormal_{C}(x)}$. Since  $\clarkeTangent_{C}(x) \subset \polar{N_{C}(x)}$ by \cref{lemma:cones:limiting-polar-inclusion}, we have $\polar{\frechetNormal_{C}(x)} \subset \polar{N_{C}(x)}$. Furthermore, \cref{thm:cones:inclusions}\,\cref{item:cones:inclusions:normal} yields $\frechetNormal_{C}(x) \subset N_{C}(x)$ and thus  $\polar{\frechetNormal_{C}(x)} \supset \polar{N_{C}(x)}$ by \cref{lemma:functan:polar-inclusion}.
    It follows that $\polar{\frechetNormal_{C}(x)} = \polar{N_{C}(x)}$.
    We now recall from \cref{thm:cones:basic-prop} that $\frechetNormal_{C}(x)$ is closed and convex.
    Hence $\bar x^* \in N_{C}(x) \setminus \frechetNormal_{C}(x)$ implies by \cref{thm:clarke:hb} that there exist $\bar x \in X$ and $\lambda \in \R$ such that
    \begin{equation*}
        \dualprod{x^*}{\bar x}_X \le \lambda < \dualprod{\bar x^*}{\bar x}_X
        \qquad (x^* \in \frechetNormal_{C}(x)).
    \end{equation*}
    Since $\frechetNormal_{C}(x)$ is a cone, this is only possible for $\lambda \ge 0$.
    Thus the first inequality shows that $\bar x \in \polar{\frechetNormal_{C}(x)}$ and the second that $\bar x \not \in \polar{N_{C}(x)}$. This is in contradiction to $\polar{\frechetNormal_{C}(x)} = \polar{N_{C}(x)}$.
    Hence $N_{C}(x) = \frechetNormal_{C}(x)$, i.e., $C$ is normally regular at~$x$.

    Conversely, if $C$ is normally regular at $x$, we obtain using \cref{lemma:cones:limiting-polar-inclusion} that
    \begin{equation*}
        \clarkeTangent_C(x) \subset \polar{N_C(x)} = \polar{\frechetNormal_C(x)}.
    \end{equation*}
    By \cref{lemma:cones:fundamental-polar}\,\cref{item:cones:fundamental-polar:weak-incl}, \cref{thm:cones:inclusions}\,\cref{item:cones:inclusions:tangent}, and \cref{lemma:functan:polar-inclusion} using the fact that $\clarkeTangent_C(x)$ is a closed convex cone by \cref{thm:cones:basic-prop}, we also have
    \begin{equation*}
        \polar{\frechetNormal_C(x)}
        \supset
        \bipolar{T_C(x)}
        \supset
        \bipolar{\clarkeTangent_C(x)}
        =\clarkeTangent_C(x).
    \end{equation*}
    Therefore $\clarkeTangent_{C}(x) = \polar{\frechetNormal_{C}(x)}$ as claimed.
\end{proof}

In sufficiently regular spaces, normal regularity implies tangential regularity of closed sets.

\begin{lemma}
    \label{lemma:cones:normal-tangent-regular}
    Let $X$ be a reflexive and Gateaux smooth Banach space and let $C\subset X$ be closed near $x \in C$.
    If $C$ is normally regular at $x$, then $C$ is tangentially regular at $x$.
\end{lemma}

\begin{proof}
    Arguing as in the proof of \cref{thm:cones:regularity:findim}\,\ref{item:cones:regularity:i} $\equivalent$ \ref{item:cones:regularity:ii}, by  \cref{lemma:functan:polar-inclusion,thm:cones:inclusions,lemma:cones:fundamental-polar,thm:cones:limiting-polar} we have
    \begin{equation*}
        T_C(x) \subset T_C^w(x) \subset \bipolar{T_C^w(x)}
        =\polar{\frechetNormal_C(x)}
        =\polar{N_C(x)}
        =\clarkeTangent_C(x)
        \subset T_C(x).
    \end{equation*}
    This shows that $T_C(x)=\clarkeTangent_C(x)$.
\end{proof}

From \cref{lemma:cones:convex,lemma:cones:normal-tangent-regular}, we immediately obtain the following regularity result.

\begin{corollary}
    \label{cor:cones:convex-regularity}
    Let $X$ be a Gateaux smooth Banach space and let $C\subset X$ be nonempty, closed, and convex. Then $C$ is normally regular at every $x \in C$. If $X$ is additionally reflexive, then $C$ is also tangentially regular at every $x\in C$.
\end{corollary}

\chapter{Tangent and normal cones of pointwise-defined sets}
\label{chap:pointcones}

As we have seen in \cref{chap:cones}, the relationships between the different tangent and normal cones are less complete in infinite-dimensional spaces than in finite-dimensional ones. In this chapter, however, we show that certain \emph{pointwise-defined} sets on $L^p(\Omega)$ for $p \in (1, \infty)$ largely satisfy the finite-dimensional relations. We will use these results in \cref{chap:superposition} to derive expressions for generalized derivatives of pointwise-defined set-valued mappings, in particular for subdifferentials of integral functionals. As mentioned in \cref{sec:cones:regularity}, these relations are less satisfying for the limiting cones than for the basic cones. To treat the limiting cones, we will therefore assume the regularity of the underlying pointwise sets. For the basic cones, we also require an assumption, which however is weaker than (tangential) regularity.

\section{Derivability}
\label{sec:pointcones:derivability}

We start with the fundamental regularity assumption.
Let $X$ be a Banach space and $C\subset X$. We then say that a tangent vector $\dir x \in T_C(x)$ at $x \in C$ is \term[vector!tangent!derivable]{derivable} if there exists an $\eps>0$ and a Borel-measurable curve $\xi: [0, \eps] \to C$ that generates $\dir x$ at $0$, i.e.,
\begin{equation}
    \label{eq:superposition:geomderiv}
    \xi(0)=x \quad\text{and}\quad \dir x = \lim_{\tau \downto 0} \frac{\xi(\tau)-\xi(0)}{\tau} = \xi'(0).
\end{equation}
Note that we do not make any assumptions on the differentiability or continuity of $\xi$ except at $\tau=0$.
We say that $C$ is \term[set!derivable!geometrically]{geometrically derivable} at $x\in C$ if every $\dir x \in T_C(x)$ is derivable.

As the next lemma shows, the point of this definition is that derivable tangent vectors are characterized by a full limit instead of just an inner limit; this additional property will allow us to construct tangent vectors in $L^p(\Omega)$ from pointwise tangent vectors, similarly to how Clarke regularity was used to obtain equality in the pointwise characterization of Clarke subdifferentials of integral functionals in \cref{thm:clarke:pointwise}.

\begin{lemma}
    \label{lemma:superposition:derivable}
    Let $C \subset X$ and $x\in C$. Then the set $T_C^0(x)$ of derivable tangent vectors is given by
    \begin{equation}
        \label{eq:superposition:derivable-cone}
        T_C^0(x)=\liminf_{\tau \downto 0} \frac{C-x}{\tau}.
    \end{equation}
\end{lemma}
\begin{proof}
    We first recall that by definition of the inner limit, $\dir x$ is an element of the set on the right-hand side if for every sequence $\tau_k \downto 0$ there exist $x_k \in C$ such that $(x_k-x)/\tau_k \to \dir x$.
    For a derivable tangent vector $\dir x \in T_C^0(x)$ and any $\tau_k\downto 0$, we can simply take $x_k=\xi(\tau_k)$.
    For the converse inclusion, let $\dir x$ be an element of the right-hand side set and let
    $\{\tau_k\}_{k \in \N}$ be such that $\tau_k \downto 0$. Then there exist $x_k \in C$ realizing the inner limit in \eqref{eq:superposition:derivable-cone} for $\dir x$.
    For any $\tau>0$ sufficiently small, set $\xi(\tau) \defeq x_k$ for a $k \in \N$ such that $\tau_{k+1} < \tau \le \tau_k$.
    Then $\xi$ is a measurable curve with $\xi(\tau) \in C$ for all $\tau \in [0, \tau_1]$.
    Moreover, \eqref{eq:superposition:derivable-cone} implies that the limit in \eqref{eq:superposition:geomderiv} exists and equals the previously chosen $\dir x \in T_C^0(x)$.
\end{proof}

By taking $\alt x\equiv x$ constant in \eqref{eq:cones:def-clarketangent} and comparing with \eqref{eq:superposition:derivable-cone}, we immediately obtain that all Clarke tangent vectors are derivable.
\begin{corollary}
    \label{cor:superposition:clarke-derivable}
    Let $C \subset X$ and $x\in C$. Then every $\dir x\in \clarkeTangent_C(x)$ is derivable.
\end{corollary}

Clearly, if $C$ is tangentially regular at $x$, then also every tangent vector is derivable.
\begin{corollary}
    \label{cor:superposition:regular-derivable}
    If $C \subset X$ is tangentially regular at $x \in C$, then every $\dir x \in T_C(x)$ is derivable.
\end{corollary}
However, a set can be geometrically derivable without being tangentially regular.
\begin{example}
    Let $C \defeq ([0, \infty) \times \{0\}) \union (\{0\} \times [0, \infty))\subset \R^2$.
    Then we obtain directly from the definition of the tangent cone that
    \begin{equation*}
        T_C(x_1, x_2) = \begin{cases}
            C & \text{if }(x_1, x_2)=(0, 0), \\
            \{0\} \times \R &\text{if } x_1=0, x_2 > 0, \\
            \R \times \{0\} &\text{if } x_1>0, x_2=0, \\
            \emptyset & \text{otherwise}.
        \end{cases}
    \end{equation*}
    However, it follows from \cref{cor:cones:clarke-liminf} that $\clarkeTangent_C(0, 0)=\{(0, 0)\}$.
    Thus $C$ is not tangentially regular at $(0, 0)$.

    On the other hand, for any $\dir x=(t_1, 0) \in T_C(0, 0)$, $t_1\in\R$, setting $\xi(s) \defeq  (st_1, 0)$ yields $\xi(0)=(0, 0)$ and $\xi'(0)=(t_1, 0)=\dir x$. Hence $\dir x$ is derivable. Similarly, setting $\xi(s) \defeq  (0, st_2)$ shows that $\dir x=(0, t_2) \in T_C(0, 0)$ is derivable for every $t_2\in\R$.
    Thus $C$ is geometrically derivable at $(0,0)$.
\end{example}

\section{Tangent and normal cones}

As the goal is to define derivatives of set-valued mappings $F:X\setto Y$ via tangent cones to their epigraphs $\epi F\subset X\times Y$, we need to consider product spaces of $p$-integrable functions (with possibly different $p$). Let therefore $\Omega\subset \R^d$ be an open and bounded domain. For
$\vec p\defeq (p_1,\ldots,p_m) \in (1, \infty)^m$, we then define
\begin{equation*}
    L^{\vec p}(\Omega) \defeq L^{p_1}(\Omega) \times \cdots \times L^{p_m}(\Omega),
\end{equation*}
endowed with the canonical Euclidean product norm, i.e.,
\begin{equation*}
    \norm{u}_{L^{\vec p}} \defeq \sqrt{\sum_{k=1}^m \norm{u_k}_{L^{p_k}}^2}\quad (u = (u_1,\dots,u_m) \in L^{\vec p}).
\end{equation*}
We will need the case $m=2$ in \cref{chap:superposition}; on first reading of the present chapter, we recommend picturing $m=1$, i.e., $L^{\vec p}(\Omega)=L^p(\Omega)$ for some $p \in (1, \infty)$.
We further denote by $p^*$ the conjugate exponent of $p \in (1, \infty)$, defined as satisfying $1/p+1/p^*=1$, and write $\vec p^* \defeq (p_1^*,\ldots, p_m^*)$ so that $L^{\vec p}(\Omega)^* \cong L^{\vec p^*}(\Omega)$.
Note that $L^{\vec p}(\Omega)$ is reflexive and Gateaux smooth as the product of reflexive and Gateaux smooth spaces; cf.~\cref{ex:epsilon:gateaux}.
We also recall the characteristic function $\1_U$ of a set $U\subset L^{\vec p}(\Omega)$, which in the vector-valued case satisfies $\1_U(u) = (1,\dots,1)\in \R^m$ if $u\in U$ and $\1_U(u) = 0\in \R^m$ otherwise. (Note that $U$ need not have product form.)

Furthermore, we denote by $\lebesgue(B)$ the Lebesgue measure of a Borel set $B\in\mathcal{B}^d(\Omega)$, the Borel algebra on $\Omega$, and call a set-valued mapping $C:\Omega\setto \R^m$ \term[mapping!Borel-measurable]{Borel-measurable} if the preimage $C^{-1}(O):=\setof{x\in\Omega}{C(x)\cap O\neq\emptyset}\in \mathcal{B}^d(\Omega)$ for every open set $O\subset \R^m$; we refer to, e.g., \cite[Chapter 14]{Rockafellar:1998} for details on measurable finite-dimensional set-valued mappings.

We then call a set $U \subset L^{\vec p}(\Omega)$ for $\vec p\in (1,\infty)^m$ \term[set!defined, pointwise]{pointwise defined} if
\begin{equation*}
    U  \defeq \left\{ u \in L^{\vec p}(\Omega)
        \mid u(x) \in C(x) \text{ for a.e. } x \in \Omega
    \right\}
\end{equation*}
for a Borel-measurable mapping $C: \Omega \setto \R^m$.
We say that $U$ is \term[set!derivable!pointwise]{pointwise derivable} if $C(x)$ is geometrically derivable at every $\xi \in C(x)$  for almost every $x \in \Omega$.

\subsection*{The fundamental cones}

We now derive pointwise characterizations of the fundamental cones to pointwise defined sets, starting with the tangent cone.

\begin{theorem}
    \label{thm:superposition:cone-tangent}
    Let $U \subset L^{\vec p}(\Omega)$  be pointwise derivable.
    Then for every $u \in U$,
    \begin{equation}
        \label{eq:superposition:lp-tangent}
        T_U(u) = \bigl\{ \dir u \in L^{\vec p}(\Omega)
            \,\bigm|\, \dir u(x) \in T_{C(x)}(u(x)) \text{ for a.e. }
            x \in \Omega
        \bigr\}.
    \end{equation}%
\end{theorem}

\begin{proof}
    The inclusion \enquote{$\subset$} follows from \eqref{eq:cones:def-tangent} and the fact that a sequence convergent in $L^{\vec p}(\Omega)$ for $\vec p \in (1, \infty)$ converges, after possibly passing to a subsequence, pointwise almost everywhere.

    For the converse inclusion, we take for almost every $x \in \Omega$ a tangent vector $\dir u(x) \in T_{C(x)}(u(x))$ at $u(x) \in C(x)$. We only need to consider the case $\dir u \in L^{\vec p}(\Omega)$. By geometric derivability, we may find for almost every $x \in \Omega$ an $\eps(x)>0$ and a curve $\xi(\freevar, x): [0, \eps(x)] \to C(x)$ such that $\xi(0, x)=u(x)$ and $\xi'_+(0, x)=\dir u(x)$. In particular, for any given $\rho>0$, we may find $\eps_\rho(x) \in (0, \eps(x)]$ such that
    \begin{equation}
        \label{eq:superposition:tangent-xi}
        \frac{\abs{\xi(t, x)-\xi(0, x)-\dir u(x)t}_2}{t} \le \rho
        \qquad (t \in (0, \eps_\rho(x)],\, \text{a.e. } x \in \Omega).
    \end{equation}
    For $t>0$, let us set
    \begin{equation*}
        E_{\rho,t} \defeq \{ x \in \Omega \mid t \le \eps_\rho(x) \}
    \end{equation*}
    and define
    \begin{equation*}
        \tilde u^{\rho,t}(x) \defeq
        \begin{cases}
            \xi(t, x) &\text{if } x \in E_{\rho,t}, \\
            u(x) &\text{if } x \in \Omega \setminus E_{\rho,t}.
        \end{cases}
    \end{equation*}
    Then each $\tilde u^{\rho,t}$ is Borel-measurable by the piecewise construction on Borel sets and the measurability (by definition) of $\xi$.
    Writing $\xi=(\xi_1, \ldots, \xi_m)$ and $\dir u=(\dir u_1, \ldots, \dir u_m)$, we have from \eqref{eq:superposition:tangent-xi} that
    \begin{equation}
        \label{eq:superposition:tangent-xi:component}
        \frac{\abs{\xi_j(t, x)-\xi_j(0, x)-\dir u_j(x)t}}{t} \le \rho
        \quad (j=1,\ldots,m,\, t \in (0, \eps_\rho(x)] \text{ for a.e. } x \in \Omega).
    \end{equation}
    Therefore, using the elementary inequality $(a+b)^2\leq 2a^2+2b^2$, we obtain
    \begin{equation}
        \label{eq:superposition:tangent-approx}
        \begin{aligned}[t]
        \norm{\tilde u^{\rho,t}-u}_{L^{\vec p}}^2
        &
        =
        \sum_{j=1}^m \norm{[\tilde u_j^{\rho,t}-u]_j}_{L^{p_j}}^2
        \\
        &
        \le
        \sum_{j=1}^m \left(\int_\Omega t^{p_j} (\rho+\abs{\dir u_j(x)})^{p_j} \ddd x\right)^{2/p_j}
        \\
        &
        \le
        \sum_{j=1}^m  \left(t \rho \lebesgue(\Omega)^{1/p_j}+t\norm{\dir u_j}_{L^{p_j}}\right)^2
        \\
        &
        \leq
        2 t^2 \sum_{j=1}^m \left(\rho \lebesgue(\Omega)^{1/p_j}\right)^2+2 t^2\norm{\dir u}_{L^{\vec p}}^2.
        \end{aligned}
    \end{equation}
    Similarly, \eqref{eq:superposition:tangent-xi:component} and the same elementary inequality together with Minkowski's inequality in the form $(a^p + b^p)^{1/p}\leq |a|+|b|$ yield
    \begin{equation}
        \label{eq:superposition:tangent-construct}
        \begin{aligned}[t]
            \frac{\norm{\tilde u^{\rho,t}-u-t\dir u}_{L^{\vec p}}^2}{t^2}
            &
            =
            \sum_{j=1}^m \frac{1}{t^2} \biggl(
                \int_{E_{\rho,t}} \abs{\xi_j(t,x)-\xi_j(0,x)-t\dir u_j(x)}^{p_j} \, d x
                \\
                \MoveEqLeft[-5]
                +
                \int_{\Omega \setminus E_{\rho,t}} \abs{\dir u_j(x)t}^{p_j} \, d x
            \biggr)^{2/p_j}
            \\
            &
            \le
            \sum_{j=1}^m \left(
                \rho^{p_j} \lebesgue(\Omega)
                + \norm{\dir u\1_{\Omega \setminus E_{\rho,t}}}_{L^{\vec p}}^{p_j}
            \right)^{2/p_j}
            \\
            &
            \le 2\sum_{j=1}^m  \left(\rho\lebesgue(\Omega)^{1/p_j}\right)^2 + 2 \norm{\dir u\1_{\Omega \setminus E_{\rho,t}}}_{L^{\vec p}}^2.
        \end{aligned}
    \end{equation}
    Now for each $k \in \N$, we can find $t_k \downto 0$ such that $\norm{\dir u\1_{\Omega \setminus E_{1/k,t_k}}}_{L^{\vec p}} \le 1/k$. This follows from Lebesgue's dominated convergence theorem and the fact that $\lebesgue(\Omega \setminus E_{\rho,t}) \to 0$ as $t \to 0$.
    The estimates \eqref{eq:superposition:tangent-approx} and \eqref{eq:superposition:tangent-construct} with $\rho=1/k$ and $t=t_k$ thus show for $u_k \defeq \alt u^{1/k,t_k}$ that $u_k \to u$ and $(u_k-u)/t_k \to \dir u$, i.e., $\dir u \in T_U(u)$.
\end{proof}

We next consider the Fréchet normal cone.

\begin{theorem}
    \label{thm:superposition:cone-frechet}
    Let $U \subset L^{\vec p}(\Omega)$  be pointwise derivable.
    Then for every $u \in U$,
    \begin{equation}
        \label{eq:superposition:lp-frechetnormal}
        \frechetNormal_U(u) = \bigl\{ u^* \in L^{\vec p^*}(\Omega)
            \,\bigm|\, u^*(x) \in \frechetNormal_{C(x)}(u(x)) \text{ for a.e. }
            x \in \Omega
        \bigr\}.
    \end{equation}
\end{theorem}

\begin{proof}
    Recalling the definition of $\frechetNormal_U(u)$ from \eqref{eq:cones:def-epsiloncone}, we need to find all $u^* \in L^{\vec p^*}(\Omega)$ satisfying for every given sequence $U \ni u_k \to u$
    \begin{equation}
        \label{eq:superposition:lk-limsup}
        0\geq
        \limsup_{k \to \infty} \frac{\dualprod{u^*}{u_k-u}_{L^{\vec p}}}{\norm{u_k-u}_{L^{\vec p}}}
        =:
        \limsup_{k \to \infty} L_k.
    \end{equation}
    Let $\eps > 0$ be arbitrary and set $v_k\defeq u-u_k$ as well as
    \begin{subequations}
        \begin{equation}
            \label{eq:superposition:cones:normest}
            Z^1_k \defeq \{ x \in \Omega \mid \abs{v_k(x)}_2 \le \inv\eps \norm{v_k}_{L^{\vec p}}\}
            \qquad (k \in \N).
        \end{equation}
        Furthermore, let $Z^2 \subset \Omega$ be such that
        \begin{align}
            \label{eq:superposition:cones:boundedness}
            &u^* \text{ is bounded on } Z^2,
            \\
            \label{eq:superposition:cones:remainder}
            &\lebesgue(Z^1_k \setminus Z^2) \le \epsilon \qquad (k \in \N).
        \end{align}
    \end{subequations}
    Using Hölder's inequality, \cref{eq:superposition:cones:normest,eq:superposition:cones:remainder}, we then estimate for $k=1,\dots,m$
    \begin{equation*}
        \begin{aligned}[t]
            L_k
            &
            = \frac{\int_{\Omega \setminus (Z^1_k \isect Z^2)} \iprod{u^*(x)}{v_k(x)}_2 \,d x}{\norm{v_k}_{L^{\vec p}}}
            +\frac{\int_{Z^1_k \isect Z^2} \iprod{u^*(x)}{v_k(x)}_2 \,d x}{\norm{v_k}_{L^{\vec p}}}
            \\
            &
            \le \frac{\norm{\1_{\Omega \setminus (Z^1_k \isect Z^2)}u^*}_{L^{\vec p^*}} \norm{v_k}_{L^{\vec p}}}{\norm{v_k}_{L^{\vec p}}}
            + \int_{Z^1_k \isect Z^2} \frac{\iprod{u^*(x)}{v_k(x)}_2}{\abs{v_k(x)}_2}
            \cdot \frac{\abs{v_k(x)}_2}{\norm{v_k}_{L^{\vec p}}} \,d x
            \\
            &
            \le \norm{\1_{\Omega \setminus (Z^1_k \isect Z^2)}u^*}_{L^{\vec p^*}}
            + \inv\eps \int_{Z^2} \max\left\{0, \frac{\iprod{u^*(x)}{v_k(x)}_2}{\abs{v_k(x)}_2}\right\} \,d x.
        \end{aligned}
    \end{equation*}
    If now for almost every $x \in \Omega$ we have that $u^*(x) \in \frechetNormal_{C(x)}(u(x))$, then also $\iprod{u^*(x)}{v_k(x)}_2 \le 0$ for almost every $x \in \Omega$.
    It follows using \eqref{eq:superposition:cones:boundedness} and the reverse Fatou inequality in the previous estimate that
    \begin{equation}
        \label{eq:superposition:cones:first-normal-est}
        \limsup_{k \to \infty} L_k \le \limsup_{k \to \infty} \norm{\1_{\Omega \setminus (Z^1_k \isect Z^2)}u^*}_{L^{\vec p^*}}.
    \end{equation}
    Since $\abs{v_k(x)}_2 \ge \inv\eps \norm{v_k}_{L^{\vec p}}$ for $x \in \Omega \setminus Z^1_k$, we have
    \begin{equation*}
        \norm{v_k}_{L^{\vec p}}
        \ge \norm{\1_{\Omega \setminus Z^1_k}v_k}_{L^{\vec p}}
        \ge (\eps^{-p} \lebesgue(\Omega \setminus Z^1_k))^{1/p} \norm{v_k}_{L^{\vec p}}.
    \end{equation*}
    Hence $\lebesgue(\Omega \setminus Z^1_k) \le \eps^p$ and $\lebesgue(\Omega \setminus (Z^1_k \isect Z_2)) \le \lebesgue(\Omega \setminus Z^1_k) + \lebesgue(\Omega \setminus Z_2) \le C\eps$ for some constant $C>0$ and small enough $\eps>0$. It therefore follows from Egorov's theorem that $\1_{\Omega \setminus (Z^1_k \isect Z^2)}u^*$ converge to $0$ in measure as $k\to\infty$. Since $u^* \in L^{\vec p^*}(\Omega)$ and $\1_{\Omega \setminus (Z^1_k \isect Z^2)}u^* \le u^*$, it follows from Vitali's convergence theorem (see, e.g., \cite[Proposition 2.27]{fonseca2007mmc}) that $\limsup_{k \to \infty} \norm{\1_{\Omega \setminus (Z^1_k \isect Z^2)}u^*}_{L^{\vec p^*}}=0$.
    Since $\eps>0$ was arbitrary, we deduce from \eqref{eq:superposition:cones:first-normal-est} that \eqref{eq:superposition:lk-limsup} holds and, consequently,
    \begin{equation*}
        \frechetNormal_U(u) \supset \{ u^* \in L^{\vec p^*}(\Omega)
            \mid u^*(x) \in \frechetNormal_{C(x)}(u(x)) \text{ for a.e. }
            x \in \Omega
        \}.
    \end{equation*}
    This proves one direction of \eqref{eq:superposition:lp-frechetnormal}, which therefore holds even without geometric derivability.

    For the converse inclusion, let $u^* \in \frechetNormal_U(u)$. We have to show that $u^*(x) \in \frechetNormal_{C(x)}(u(x))$ for almost every $x \in \Omega$, which we do by contradiction. Assume therefore that the pointwise inclusion does not hold. By the polarity relationship $\frechetNormal_{C(x)}(u(x))=\polar{T_{C(x)}(u(x))}$ from \cref{lemma:cones:fundamental-polar}, we can find $\delta > 0$ and a Borel set $E \subset \Omega$ of finite positive Lebesgue measure such that for each $x \in E$, there exists $w(x) \in T_{C(x)}(u(x))$ with $\abs{w(x)}_2=1$ and $\iprod{u^*(x)}{w(x)}_2 \ge \delta$.
    We may without loss of generality assume that $C(x)$ is geometrically derivable at $w(x)$ for every $x \in E$, i.e., for each $x \in E$ there exists a curve $\xi(\freevar, x): [0, \eps(x)] \to C(x)$ such that $\xi'_+(0, x)=w(x)$ and $\xi(0, x)=u(x)$. Let now $c \in (0, \delta)$ be arbitrary. By replacing $E$ by a subset of positive measure, we may by Egorov's theorem assume the existence of $\eps > 0$ such that
    \begin{equation}
        \label{eq:superposition:cones:normal-xi}
        \abs{\xi(t, x)-\xi(0, x)-w(x)t}_2 \le c t
        \quad (t \in [0, \eps],\, x \in E).
    \end{equation}

    Let us define
    \begin{equation*}
        \tilde u^t(x) \defeq \begin{cases}
            \xi(t, x) &\text{if } x \in E, \\
            u(x) &\text{if } x \in \Omega \setminus E.
        \end{cases}
    \end{equation*}
    Setting $v^t \defeq \tilde u^t -u$, we have $v^t(x)=\xi(t, x)-\xi(0, x)$ for $x \in E$ and $v^t(x)=0$ for $x \in \Omega \setminus E$. Therefore, writing $v^t=(v_1^t, \ldots, v_m^t)$, $w=(w_1, \ldots, w_m)$, and $\xi=(\xi_1, \ldots \xi_m)$, we obtain using \eqref{eq:superposition:cones:normal-xi} for $t \in (0, \eps]$ and some $c'>0$ that
    \begin{equation*}
        \begin{aligned}[t]
            \norm{v^t}_{L^{\vec p}}^2
            &
            = \sum_{j=1}^m \left(\int_E \abs{\xi_j(t, x)-\xi_j(0, x)}^{p_j}\, d x\right)^{2/p_j}
            \\
            &
            \le
            \sum_{j=1}^m \left(\int_E (\abs{w_j(x)} t+ ct)^{p_j}\, d x\right)^{2/p_j}
            \le c't^2.
        \end{aligned}
    \end{equation*}
    Likewise,
    \begin{equation*}
        \iprod{u^*(x)}{v^t(x)}_2
        \ge \iprod{u^*(x)}{w(x)}_2 - \abs{u^*(x)}_2 \cdot \abs{\xi(t, x)-\xi(0, x)-wt}_2
        \ge \delta t - ct.
    \end{equation*}
    It follows that
    \begin{equation*}
        \limsup_{t \downto 0} \int_{E} \frac{\iprod{u^*(x)}{v^t(x)}_2}{\norm{v^t}_{L^{\vec p}}} \,d x
        \ge
        \limsup_{t \downto 0}
        \frac{\lebesgue(E)(\delta t - c t)}{c't}
        = \frac{\lebesgue(E)(\delta-c)}{c'} > 0.
    \end{equation*}
    Taking $u_k \defeq \tilde u^{1/k}$ for $k \in \N$, we obtain $\lim_{k \to \infty} L_k > 0$ and therefore
    $u^* \notin \frechetNormal_U(u)$. By contraposition, this shows that $u^*(x) \in \frechetNormal_{C(x)}(u(x))$ for almost every $x \in \Omega$.
\end{proof}

We can now derive a similar polarity relationships to the finite-dimensional one in \cref{lemma:cones:fundamental-polar}.

\begin{corollary}
    \label{cor:superposition:polarity}
    Let $U  \subset L^{\vec p}(\Omega)$  be pointwise derivable and $u\in U$.
    Then $\frechetNormal_U(u)=\polar{T_U(u)}$.
\end{corollary}

\begin{proof}
    By \cref{thm:superposition:cone-frechet,thm:superposition:cone-tangent,lemma:cones:fundamental-polar}, we have
    \begin{equation}
        \label{eq:superposition:pointwise-polar}
        \begin{aligned}[t]
            u^* \in \frechetNormal_U(u)
            & \equivalent
            u^*(x) \in \frechetNormal_{C(x)}(u(x)) && \text{(a.e.~$x \in \Omega$)}
            \\
            & \equivalent
            \iprod{u^*(x)}{\dir u(x)}_2 \le 0 && \text{(a.e.~$x \in \Omega$ when $\dir u(x) \in T_{C(x)}(u(x))$)}
            \\
            & \implies
            \dualprod{u^*}{\dir u}_{L^{\vec p}} \le 0 && \text{(when $\dir u \in T_U(u)$)}
            \\
            & \equivalent
            u^* \in \polar{T_U(u)}.
        \end{aligned}
    \end{equation}
    Hence $\frechetNormal_U(u) \subset \polar{T_U(u)}$.

    For the converse inclusion, we need to improve the implication in \eqref{eq:superposition:pointwise-polar} to an equivalence.
    We argue by contradiction. Assume that $u^* \in \polar{T_U(u)}$ and that there exists some $\dir \bar u \in T_U(u)$ and a subset $E \subset \Omega$ with $\lebesgue(\Omega \setminus E)>0$ and
    \begin{equation*}
        \iprod{u^*(x)}{\dir{\bar u}(x)}_2 > 0
        \qquad (x \in E).
    \end{equation*}
    Taking $\bar u^*(x) \defeq (1+t\1_E(x))u^*(x)$, we obtain for sufficiently large $t$ that
    $\dualprod{\bar u^*}{\dir{\bar u}}_{L^{\vec p}} > 0$.
    This contradicts that $u^* \in \polar{T_U(u)}$.
    Hence $\frechetNormal_U(u) \supset \polar{T_U(u)}$.
\end{proof}

\subsection*{The limiting cones}

For the limiting cones, we in general only have an inclusion of the pointwise cones.
\begin{theorem}
    \label{lemma:superposition:limiting-incl-tangent}
    Let $U \subset L^{\vec p}(\Omega)$  be pointwise derivable.
    Then for every $u \in U$,
    \begin{equation*}
        \clarkeTangent_U(u)
        \supset
         \bigl\{ \dir u \in L^{\vec p}(\Omega)
            \,\bigm|\, \dir u(x) \in \clarkeTangent_{C(x)}(u(x)) \text{ for a.e. }
            x \in \Omega
        \bigr\}.
    \end{equation*}%
\end{theorem}

\begin{proof}
    Let $\dir u \in L^{\vec p}(\Omega)$ with $\dir u(x)\in \clarkeTangent_{C(x)}(u(x))$ for almost every $x\in\Omega$ and let $u_k \to u$ in $L^{\vec p}(\Omega)$.
    In particular, we then have $u_k(x) \to u(x)$ for almost every $x \in \Omega$. Furthermore, by the inner limit characterization of $\clarkeTangent_{C(x)}(u(x))$ in \cref{cor:cones:clarke-liminf}, there exist $\dir \tilde u_k(x) \in T_{C(x)}(u_k(x))$ with $\dir \tilde u_k(x) \to \dir u(x)$. Egorov's theorem then yields for all $\ell \ge 1$ a Borel-measurable set $E_\ell \subset \Omega$ such that $\lebesgue(\Omega \setminus E_\ell) < 1/\ell$ and $\dir \tilde u_k \to \dir u$ uniformly on $E_\ell$.
    Since $T_{C(x)}(u_k(x))$ is a cone, we have $0 \in T_{C(x)}(u_k(x))$. It follows that
    \begin{equation*}
        T_{C(x)}(u_k(x)) \ni \dir u_{\ell,k}(x) \defeq \1_{E_\ell}(x) \dir \tilde u_k(x).
    \end{equation*}
    In particular, \eqref{eq:superposition:lp-tangent} shows that $\dir u_{\ell,k} \in T_U(u_k)$ with
    $\dir u_{\ell,k} \to \dir u_\ell \defeq \dir u \1_{E_\ell}$ in $L^{\vec p}(\Omega)$ as $k \to \infty$.
    By Vitali's convergence theorem (compare the proof of \cref{thm:superposition:cone-frechet}), $\dir u \1_{E_\ell} \to \dir u$ in $L^{\vec p}(\Omega)$ as $\ell \to \infty$.
    Therefore, we may extract a diagonal subsequence $\{\dir \alt u_k \defeq \dir u_{\ell_k,k}\}_{k \ge 1}$ of $\{\dir u_{\ell,k}\}_{k,\ell \ge 1}$ such that  $\dir \alt u_k \to \dir u$.
    Since $u_k \to u$ was arbitrary and $\dir \alt u_k \in T_U(u_k)$, we deduce that $\dir u \in \clarkeTangent_U(u)$.
\end{proof}

\begin{theorem}
    \label{lemma:superposition:limiting-incl-normal}
    Let $U \subset L^{\vec p}(\Omega)$  be pointwise derivable.
    Then for every $u \in U$,
    \begin{equation*}
        N_U(u)
        \supset
        \bigl\{ u^* \in L^{\vec p^*}(\Omega)
            \,\bigm|\, u^*(x) \in N_{C(x)}(u(x)) \text{ for a.e. }
            x \in \Omega
        \bigr\}.
    \end{equation*}
\end{theorem}
\begin{proof}
    Let $u^*\in L^{\vec p^*}(\Omega)$ with $u^*(x)\in N_{C(x)}(u(x))$ for almost every $x\in \Omega$. Then by definition, for almost all $x \in \Omega$ there exist $C(x) \ni {\tilde u}_k(x) \to u(x)$ as well as $\frechetNormal_{C(x)}(u(x)) \ni {\tilde u}_k^*(x) \to u^*(x)$.
    By Egorov's theorem, for every $\ell \ge 1$ there exists a Borel-measurable set $E_\ell \subset \Omega$ such that $\lebesgue(\Omega \setminus E_\ell) < 1/\ell$ and ${\tilde u}_k^* \to u^*$ as well as ${\tilde u}_k \to u$ uniformly on $E_\ell$.
    We set $u_{\ell,k} \defeq \1_{E_\ell} u_k + (1-\1_{E_\ell}) u$ and  $u_{\ell,k}^* \defeq \1_{E_\ell} {\tilde u}_k^*$.
    Then $u_{\ell,k}^*(x) \in \frechetNormal_{C(x)}(u_{\ell,k}(x))$ for almost every $x \in \Omega$.
    By Vitali's convergence theorem (compare the proof of \cref{thm:superposition:cone-frechet}), both $u_{\ell,k} \to u$ in $L^{\vec p}(\Omega)$ and $u_{\ell,k}^* \to u_\ell^*$ in $L^{\vec p^*}(\Omega)$ for $u_\ell^* \defeq \1_{E_\ell} u^*$.
    Since $u_\ell^* \to u^*$ in $L^{\vec p^*}(\Omega)$, we can extract a diagonal subsequence of $\{(u_{\ell,k}, u_{\ell,k}^*)\}_{\ell,k \ge 1}$ to deduce that $u^* \in N_U(u)$.
\end{proof}

If the pointwise sets $C(x)$ are regular, we have the following polarity between the cones to the pointwise-defined set $U$.
\begin{lemma}
    \label{cor:superposition:polarity2}
    Let $U  \subset L^{\vec p}(\Omega)$  be pointwise derivable and $u \in U$. If $C(x)$ is regular at $u(x)$ and closed near $u(x)$ for almost every $x \in \Omega$,
    then $T_U(u)=\polar{\frechetNormal_U(u)}$.
\end{lemma}

\begin{proof}
    By the regularity of $C(x)$ at $u(x)$ for almost every $x \in \Omega$ and \cref{thm:superposition:cone-frechet}, we have
    \begin{equation*}
        \frechetNormal_U(u) = \bigl\{ u^* \in L^{\vec p^*}(\Omega)
            \,\bigm|\, u^*(x) \in N_{C(x)}(u(x)) \text{ for a.e. }
            x \in \Omega
        \bigr\}.
    \end{equation*}
    By \cref{thm:cones:limiting-polar:findim}, $\polar{N_{C(x)}(u(x))}=\clarkeTangent_{C(x)}(u(x))$ for almost every $x \in \Omega$. Arguing as in the proof of \cref{cor:superposition:polarity}, we thus obtain
    \begin{equation*}
        \polar{\frechetNormal_U(u)} = \bigl\{ \dir u \in L^{\vec p}(\Omega)
            \,\bigm|\, \dir u(x) \in \clarkeTangent_{C(x)}(u(x)) \text{ for a.e. }
            x \in \Omega
        \bigr\}.
    \end{equation*}
    The regularity of $C(x)$ also implies that  $\clarkeTangent_{C(x)}(u(x))=T_{C(x)}(u(x))$ for almost every $x \in \Omega$.
    The claims now follow from \cref{thm:superposition:cone-tangent}.
\end{proof}

We can use this result to transfer the regularity of $C(x)$ to $U$.
\begin{lemma}
    \label{cor:superposition:weak-tangent}
    Let $U  \subset L^{\vec p}(\Omega)$  be pointwise derivable and $u \in U$. If $C(x)$ is regular at $u(x)$ and closed near $u(x)$ for almost every $x \in \Omega$, then $U$ is regular at $u$ and
    \begin{equation*}
        T_U^w(u)=T_U(u)=\clarkeTangent_U(u).
    \end{equation*}
\end{lemma}

\begin{proof}
    Since $L^{\vec p}(\Omega)$ is reflexive, we have $\frechetNormal_U(u)=\polar{T_U^w(u)}$ by \cref{lemma:cones:fundamental-polar}\,\cref{item:cones:fundamental-polar:reflexive}.
    This fact together with \cref{cor:superposition:polarity2,lemma:functan:polar-inclusion,thm:cones:inclusions} shows that
    \begin{equation*}
        T_U^w(u) \subset \bipolar{T_U^w(u)}
        =\polar{\frechetNormal_U(u)}
        =T_U(u)
        \subset T_U^w(u).
    \end{equation*}
    Furthermore, by the regularity and closedness assumptions, we obtain from \cref{thm:superposition:cone-tangent,lemma:superposition:limiting-incl-tangent} that $T_U(u) = \clarkeTangent_U(u)$, which also implies tangential regularity.

    Since $L^{\vec p}(\Omega)$ for $\vec p \in (1, \infty)^m$ is reflexive and Gateaux smooth, normal regularity follows from \cref{thm:cones:regularity:infdim} together with \cref{cor:superposition:polarity2}.
\end{proof}

From this, we obtain pointwise expressions with equality. For the Clarke tangent cone, we only require local closedness of the underlying sets.
\begin{theorem}
    \label{lemma:superposition:limiting-equal-tangent}
    Let $U \subset L^{\vec p}(\Omega)$  be pointwise derivable.
    If $C(x)$ is closed near $u(x)$ for almost every $x\in \Omega$ for every $u \in U$, then
    \begin{equation*}
        \clarkeTangent_U(u) = \bigl\{ \dir u \in L^{\vec p}(\Omega)
            \,\bigm|\, \dir u(x) \in \clarkeTangent_{C(x)}(u(x)) \text{ for a.e. }
            x \in \Omega
        \bigr\}.
    \end{equation*}%
\end{theorem}

\begin{proof}
    The inclusion \enquote{$\supset$} was already shown in \cref{lemma:superposition:limiting-incl-tangent}.
To prove the converse inclusion when $C(x)$ is closed near $u(x)$ for almost every $x \in \Omega$, we only need to observe from \cref{lemma:cones:limiting-polar-inclusion,lemma:superposition:limiting-incl-normal} that
    \begin{equation*}
        \begin{aligned}[b]
            \clarkeTangent_C(u) \subset \polar{N_C(u)}
            &\subset \polar{\bigl\{ u^* \in L^{\vec p^*}(\Omega)
                    \,\bigm|\, u^*(x) \in N_{C(x)}(u(x)) \text{ for a.e. }
                    x \in \Omega
            \bigr\}} \\
            &=
            \bigl\{ \dir u \in L^{\vec p}(\Omega)
                \,\bigm|\, \dir u(x) \in \clarkeTangent_{C(x)}(u(x)) \text{ for a.e. }
                x \in \Omega
            \bigr\},
        \end{aligned}
    \end{equation*}
    where the last equality again follows from \cref{thm:cones:limiting-polar:findim} together with an argument as in the proof of \cref{cor:superposition:polarity}.
\end{proof}

For the limiting normal cone, however, we \emph{do} require regularity.
\begin{theorem}
    \label{lemma:superposition:limiting-equal-normal}
    Let $U \subset L^{\vec p}(\Omega)$  be pointwise derivable.
    If $C(x)$ is regular at $u(x)$ and closed near $u(x)$ for almost every $x \in \Omega$, then for every $u \in U$,
    \begin{equation*}
        N_U(u) = \bigl\{ u^* \in L^{\vec p^*}(\Omega)
            \,\bigm|\, u^*(x) \in N_{C(x)}(u(x)) \text{ for a.e. }
            x \in \Omega
        \bigr\}.
    \end{equation*}
\end{theorem}
\begin{proof}
    The inclusion \enquote{$\supset$} was already shown in \cref{lemma:superposition:limiting-incl-normal}.
    The converse inclusion for regular and closed $C(x)$ follows from \cref{cor:superposition:weak-tangent,thm:superposition:cone-frechet}.
\end{proof}

\begin{remark}\label{rem:pointcones:limiting}
    \Cref{thm:superposition:cone-frechet,thm:superposition:cone-tangent} on the fundamental cones are based on \cite{tuomov-pdex2stability}.
    Without regularity, the characterization of the limiting normal cone of a pointwise-defined set is much more delicate. A full characterization was given in \cite{mehlitzwachsmuth2018limiting,mehlitzwachsmuth2019decomposable}, which showed that even for a closed nonconvex set, the limiting normal cone contains the convex hull of the strong limiting normal cone (where the limit is taken with respect to strong convergence instead of weak-$*$ convergence) and is dense in the  Dini normal cone $\polar\clarkeTangent_C(x)$ -- in the words of the authors, it may be \enquote{unpleasantly large}. This is due to an inherent convexifying effect of integration with respect to the Lebesgue measure.

    A characterization of specific pointwise-defined sets in Sobolev spaces was derived in \cite{harderwachsmuth2018sobolev}, with similar conclusions.
\end{remark}

\chapter{Derivatives and coderivatives of set-valued mappings}
\label{chap:graphical}

We are now ready to differentiate set-valued mappings; as already discussed, these generalized derivatives are based on the tangent and normal cones of the previous \cref{chap:cones}.
To account for the changed focus, we will slightly switch notation and use in this and the following chapters of \cref{part:setvalued} uppercase letters for set-valued mappings and lowercase letters for scalar-valued functionals such that, e.g., $F(x) = \partial f(x)$.
We focus in this chapter on examples, basic properties, and relationships between the various derivative concepts. In the following \crefrange{chap:gderiv}{chap:colimiting}, we then develop calculus rules for each of the different derivatives and coderivatives.

\section{Definitions}
\label{sec:graphical:def}

To motivate the following definitions, it is instructive to recall the geometric intuition behind the classical derivative of a smooth scalar function $f$ as the limit of a difference quotient: given an (infinitesimal) change $\dir x$ of the argument $x$, it gives the corresponding (infinitesimal) change $\dir y$ of the value $y=f(x)$ required to stay on the graph of $f$. In other words, $\Delta y = f'(x)\Delta x$, which in geometric terms means that $(\dir x, \dir y)$ is a tangent vector to $\graph f$.
If we want to derive calculus rules such as a chain rule (e.g., \cref{thm:clarke:chain}), we will also need a generalization of the adjoint derivative $f'(x)^*$, which by definition satisfies
$\dualprod{f'(x)\dir x}{y^*} = \dualprod{\dir x}{f'(x)^*y^*}$ for all incremental changes $\Delta x$ and \enquote{dual changes} $y^*$. In other words, $\dualprod{\dir x}{f'(x)^*y^*}+\dualprod{f'(x)\dir x}{-y^*} = 0$ for the tangent vector $(\dir x, f'(x)\dir x)$, which in geometric terms means that $(f'(x)^*y^*,-y^*)$ is orthogonal and thus a normal vector to $\graph f$; see \cref{fig:graphical:coderivative-negation} and \cref{ex:graphical:linear} below.

We thus distinguish
\begin{enumerate}
    \item \emph{graphical derivatives}, which generalize classical derivatives and are based on tangent cones;
    \item \emph{coderivatives}, which generalize adjoint derivatives and are based on normal cones.
\end{enumerate}
In each case, we can use either basic or limiting cones, leading to four different definitions.

Specifically, let $X,Y$ be Banach spaces and $F: X \setto Y$. Then we define
\begin{enumerate}
    \item the \term[derivative!graphical]{graphical derivative} of $F$ at $x\in X$ for $y\in Y$ as
        \begin{equation*}
            D F(x|y): X \setto Y, \qquad
            D F(x|y)(\dir x) \defeq \setof{ \dir y\in Y}{(\dir x, \dir y) \in T_{\graph F}(x, y)};
        \end{equation*}

    \item the \term[derivative!Clarke graphical]{Clarke graphical derivative} of $F$ at $x\in X$ for $y\in Y$ as
        \begin{equation*}
            \clarkeGD F(x|y): X \setto Y, \qquad
            \clarkeGD F(x|y)(\dir x) \defeq \setof{ \dir y\in Y}{(\dir x, \dir y) \in \clarkeTangent_{\graph F}(x, y)};
        \end{equation*}

    \item
        the \term[coderivative!Fréchet]{Fréchet coderivative} of $F$ at $x\in X$ for $y\in Y$ as
        \begin{equation*}
            \frechetCod F(x|y): Y^* \setto X^*,\qquad
            \frechetCod F(x|y)(y^*) \defeq \setof{ x^*\in X^* }{ (x^*, -y^*) \in \frechetNormal_{\graph F}(x, y)};
        \end{equation*}

    \item the (\term[coderivative!basic]{basic} or \term[coderivative!limiting]{limiting} or \term[coderivative!Mordukhovich]{Mordukhovich}) \term{coderivative} of $F$ at $x\in X$ for $y\in Y$ as
        \begin{equation*}
            \coderivative F(x|y): Y^* \setto X^*,\qquad
            \coderivative F(x|y)(y^*) \defeq \setof{ x^* \in X^*}{(x^*, -y^*) \in N_{\graph F}(x, y)}.
        \end{equation*}
\end{enumerate}

\begin{figure}[t]
    \centering
    \begin{subfigure}[b]{0.45\textwidth}
    \begin{asy}
        real f(real x){ return x^2/1.5; };
        real gf(real x){ return x*2/1.5; };
        path g=graph(f, -1, 1);
        draw(g, defaultpen);
        path l=pt(f, -.1)--pt(f, .7);
        label("$f$", midpoint(g), N);

        pair p=pt(f, .5);
        dot(p);
        pair t=1.1*(1, gf(.5));
        pair n=1.1*(gf(.5), -1);
        draw(p..(p-t), primalline+linewidth(1.1), Arrow);
        draw(p..(p+t), primalline+linewidth(1.1), Arrow);
        dot(p+n, invisible);
        dot(p-n, invisible);
        real dx=0.9;

        pair xx=p+t/t.x*dx;
        pair yy=p+n/n.x*dx;
        pair zz=p+(dx, 0);
        draw(p..zz, dotted, Arrow);
        label("$\dir x$", midpoint(p..(p+(dx, 0))), S);
        draw(zz..xx, dotted, Arrow);
        label("$f'(x)\dir x$", midpoint(zz..xx), E);

        label("$T_{\graph f}$", midpoint(p..(p-t)), 3*E);
    \end{asy}
    \caption{derivative and tangent cone}
    \end{subfigure}
    \begin{subfigure}[b]{0.45\textwidth}
    \begin{asy}
        real f(real x){ return x^2/1.5; };
        real gf(real x){ return x*2/1.5; };
        path g=graph(f, -1, 1);
        draw(g, defaultpen);
        path l=pt(f, -.1)--pt(f, .7);
        label("$f$", midpoint(g), N);

        pair p=pt(f, .5);
        dot(p);
        pair t=1.1*(1, gf(.5));
        pair n=1.1*(gf(.5), -1);
        draw(p..(p+n), primalline+linewidth(1.1), Arrow);
        draw(p..(p-n), primalline+linewidth(1.1), Arrow);
        real dx=0.6;

        pair yy=p+n/n.x*dx;
        pair zz=p+(dx, 0);
        draw(p..zz, dotted, Arrow);
        label("$f'(x)^*y^*$", p+(dx, 0), N+0.5*E);
        draw(zz..yy, dotted, Arrow);
        label("$-y^*$", midpoint(zz..yy), E);
        label("$N_{\graph f}$", midpoint(p..(p+n)), W);
    \end{asy}
    \caption{adjoint derivative and normal cone}
    \end{subfigure}

    \caption{Illustration of derivative and adjoint derivative in relation to tangent and normal cone for a smooth single-valued function. Note how $f'(x)^*y^*$ is orthogonal to $f'(x)\Delta x$ if $-y^*$ is orthogonal to $\Delta x$.}
    \label{fig:graphical:coderivative-negation}
\end{figure}

Observe how the coderivatives operate from $Y^*$ to $X^*$, while the derivatives operate from $X$ to $Y$.
It is crucial that these are defined directly via (possibly nonconvex) normal cones rather than via polarity from the corresponding graphical derivatives to avoid convexification. This will allow for sharper results involving these coderivatives.

We illustrate these definitions -- and verify that these capture the motivation -- for a single-valued linear operator.
\begin{example}[single-valued linear operators]
    \label{ex:graphical:linear}
    Let $F(x) \defeq \{Ax\}$ for  $A \in \linear(X; Y)$ and $u=(x,Ax) \in \graph F$.
    Note that $\graph F$ is a linear subspace of $X \times Y$.
    Since $\graph F$ is regular by \cref{cor:cones:convex-regularity}, both of the tangent cones are given by
    \begin{equation*}
        T_{\graph F}(u)=\clarkeTangent_{\graph F}(u)=\graph F=\{(\dir x, A\dir x) \in X \times Y \mid \dir x \in X\},
    \end{equation*}
    while the normal cones are given by
    \begin{equation*}
        \begin{aligned}
            N_{\graph F}(u)&=\frechetNormal_{\graph F}(u)
            =\{u^* \in X^* \times Y^* \mid u^* \perp \graph F\}
            \\
            &
            =\{(x^*,y^*) \in X^* \times Y^* \mid \iprod{x^*}{\dir x}_X + \iprod{y^*}{A\dir x}_Y = 0 \text{ for all } \dir x \in X\}
            \\
            &
            =\{(A^*y^*, -y^*)\in X^*\times Y^* \mid y^* \in Y^*\}.
        \end{aligned}
    \end{equation*}
    This immediately yields the graphical derivatives
    \begin{equation*}
        DF(x|Ax)(\dir x)
        =\clarkeGD F(x|Ax)(\dir x)
        = \{ A\dir x \}
    \end{equation*}
    as well as the coderivatives
    \begin{equation*}
        \coderivative F(x|y)(y^*)
        = \frechetCod F(x|y)(y^*)
        =\{A^*y^*\}.
    \end{equation*}
\end{example}

Using \eqref{eq:cones:def-tangent}, we can also write the graphical derivative as
\begin{equation}
    \label{eq:graphical:gderiv-altdef}
    DF(x|y)(\dir{x}) = \limsup_{t \downto 0,\, \dir{\alt{x}} \to \dir{x}} \frac{F(x+t \dir{\alt{x}})-y}{t},
\end{equation}
since
\begin{equation*}
    (\Delta x, \Delta y) \in \limsup_{\tau \downto 0} \frac{\graph F - (x, y)}{\tau}
\end{equation*}
if and only if there exist $\tau_k \downto 0$ and $x_k$ such that
\begin{equation}
    \label{eq:graphical:gderiv-explicit}
    \Delta x = \lim_{k \to \infty} \frac{x_k-x}{\tau_k}
    \quad
    \text{and}
    \quad
    \Delta y \in \limsup_{k \to \infty} \frac{F(x_k)-y}{\tau_k}.
\end{equation}
The former forces $x_k = x+\tau_k\Delta x_k$ for $\Delta x_k \to \Delta x$, so the latter gives \eqref{eq:graphical:gderiv-altdef}.

In infinite-dimensional spaces, we also have to distinguish the \term[derivative!graphical!weak]{weak graphical derivative} $D^w F(x|y)$ and the \term[coderivative!$\epsilon$-]{$\epsilon$-coderivative} $\frechetCod_\epsilon F(x|y)$, both constructed analogously from the weak tangent cone $T^w_{\graph F}(x, y)$ and the $\epsilon$-normal cone $\frechetNormal_{\graph F}^\epsilon(x, y)$, respectively. However, we will not be working directly with these and instead switch to the setting of the corresponding cones when they would be needed.

\begin{remark}[a much too brief history of various (co)derivatives]
    As for the various tangent and normal cones, the (more recent) development of derivatives and coderivatives of set-valued mappings is complicated, and we do not attempt to give a full account, instead referring to the commentaries to \cite[Chapter 8]{Rockafellar:1998}, \cite[Chapter 1.4.12]{Mordukhovich:2006}, and \cite[Chapter 1]{mordukhovich2018variational}.

    The graphical derivative goes back to Aubin \cite{aubin1981contingent}, who also introduced the Clarke graphical derivative (under the name \term[derivative!circatangent]{circatangent derivative}) in \cite{aubin1984lipschitz}. Coderivatives based on normal cones were mainly treated there for mappings whose graphs are convex, for which these cones can be defined as polars of the appropriate tangent cones. Graphical derivatives were further studied in \cite{thibault1983tangent}. In parallel, Mordukhovich introduced the (nonconvex) limiting coderivative via his limiting normal cone in \cite{Mordukhovich1980metric}, again stressing the need for a genuinely nonconvex direct construction. The term \emph{coderivative} was coined by Ioffe, who was the first to study these mappings systematically in \cite{ioffe1984approximate}.
\end{remark}

\section{Basic properties}
\label{sec:graphical:properties}

We now translate various results of \cref{chap:cones} on tangent and normal cones to the setting of graphical derivatives and coderivatives.
First of all, since derivatives and coderivatives are defined via tangent and normal cones, their positive homogeneity is a direct consequence of \cref{thm:cones:basic-prop}\,\ref{item:cones:basic-prop:cone}.
\begin{corollary}
    \label{cor:graphical:pos-hgen}
    For $F: X \setto Y$, $x \in X$, and $y \in Y$, and any $\alpha > 0$ as well as $\dir x \in X$ and $x^* \in Y^*$, we have
    \begin{enumerate}
        \item $DF(x|y)(\alpha \dir x)=\alpha DF(x|y)(\dir x)$;
        \item $\clarkeGD F(x|y)(\alpha \dir x)=\alpha \clarkeGD F(x|y)(\dir x)$;
        \item $\frechetCod F(x|y)(\alpha x^*)=\alpha \frechetCod F(x|y)(x^*)$;
        \item $D^* F(x|y)(\alpha x^*)=\alpha D^* F(x|y)(x^*)$.
    \end{enumerate}
\end{corollary}

Furthermore, from \cref{thm:cones:inclusions}, we immediately obtain
\begin{corollary}
    \label{cor:graphical:inclusions}
    For $F: X \setto Y$, $x \in X$, and $y \in Y$, we have the inclusions
    \begin{enumerate}
        \item
            $\clarkeGD F(x|y)(\dir x)  \subset DF(x|y)(\dir x) \subset D^w F(x|y)(\dir x)$ for all $\dir x \in X$;
        \item
            $\frechetCod F(x|y)(y^*) \subset \coderivative F(x|y)(y^*)$ for all $y^* \in Y^*$.
    \end{enumerate}
\end{corollary}

Similarly, we obtain from \cref{thm:cones:basic-prop} the following outer semicontinuity and convexity properties.
\begin{corollary}
    \label{cor:graphical:semicontinuity}
    For $F: X \setto Y$, $x \in X$, and $y \in Y$,
    \begin{enumerate}
        \item
            $D F(x|y)$, $\clarkeGD F(x|y)$, and $\frechetCod F(x|y)$ are closed;
        \item if $X$ and $Y$ are finite-dimensional, then $\coderivative F(x|y)$ is closed;
        \item $\clarkeGD F(x|y)$ and $\frechetCod F(x|y)$ are convex.
    \end{enumerate}
\end{corollary}

Graphical derivatives and coderivatives behave completely symmetrically with respect to inversion of a set-valued mapping (which we recall is always possible in the sense of preimages).
\begin{lemma}
    \label{lemma:graphical:inverse}
    Let $F:X \setto Y$, $x\in X$, and $y\in Y$. Then
    \begin{align*}
        \dir y &\in D F(x|y)(\dir x) & &\equivalent & \dir x &\in D \inv F(y|x)(\dir y),
        \\
        \dir y &\in \clarkeGD F(x|y)(\dir x) & &\equivalent  &\dir  x &\in \clarkeGD \inv F(y|x)(\dir y),
        \\
        x^* & \in \frechetCod F(x|y)(y^*) & &\equivalent & -y^* &\in \frechetCod \inv F(y|x)(-x^*),
        \\
        x^* &\in \coderivative F(x|y)(y^*) & &\equivalent &-y^* &\in \coderivative \inv F(y|x)(-x^*).
    \end{align*}
\end{lemma}

\begin{proof}
    We have
    \begin{equation*}
        \begin{aligned}
            \dir y \in D F(x|y)(\dir x)
            &
            \equivalent
            (\dir x, \dir y) \in T_{\graph F}(x, y)
            \\
            &
            \equivalent
            (\dir y, \dir x) \in T_{\graph \inv F}(y, x)
            \\
            & \equivalent
            \dir x \in D \inv F(y|x)(\dir y).
        \end{aligned}
    \end{equation*}
    The proof for the regular derivative and the coderivatives is completely analogous.
\end{proof}

\subsection*{Adjoints of set-valued mappings}

From the various relations between normal and tangent cones, we obtain corresponding relations between these derivatives. To state these relationships, we need to introduce the upper and lower adjoints of set-valued mappings. Let $H: X \setto Y$ be a set-valued mapping. Then the \term[adjoint!of set-valued mapping!upper]{upper adjoint} of $H$ is defined as
\begin{equation*}
    \upperadj H(y^*) \defeq \{ x^* \mid \iprod{x^*}{x}_X \le \iprod{y^*}{y}_Y \text{ for all } y \in H(x),\, x \in X \},
\end{equation*}
and the \term[adjoint!of set-valued mapping!lower]{lower adjoint} of $H$ as
\begin{equation*}
    \loweradj H(y^*) \defeq \{ x^* \mid \iprod{x^*}{x}_X \ge \iprod{y^*}{y}_Y \text{ for all } y \in H(x),\, x \in X \}.
\end{equation*}

As the next example shows, these notions generalize the definition of the adjoint of a linear operator.
\begin{example}[upper and lower adjoints of linear mappings]
    \label{ex:graphical:adjoints:linear}
    Let $H(x) \defeq \{Ax\}$ for $A \in \linear(X; Y)$. Then
    \begin{equation*}
        \begin{aligned}
            \upperadj H(y^*)
            &
            =
            \setof{x^*\in X^*}{\iprod{x^*}{x}_X \le \iprod{y^*}{y}_Y \text{ for all } y =Ax,\ x\in X}
            \\
            &
            =
            \setof{x^*\in X^*}{\iprod{x^*}{x}_X \le \iprod{y^*}{Ax}_Y \text{ for all } x \in X}
            \\
            &
            =
            \setof{x^*\in X^*}{\iprod{x^*-A^* y^*}{x}_X \le 0 \text{ for all } x\in X}
            \\
            &
            = \{A^* y^*\}.
        \end{aligned}
    \end{equation*}
    Similarly, $\loweradj H(y^*) = \{A^*y^*\}$.
\end{example}

For solution mappings of linear equations, we have the following adjoints.
\begin{example}[upper and lower adjoints of solution maps to linear equations]
    \label{ex:graphical:adjoints:linear-solution}
    Let $H(x) \defeq \{y \mid Ay=x\}$ for $A \in \linear(X; Y)$. Then
    \begin{equation*}
        \upperadj H(y^*)
        =
        \setof{ x^*}{\iprod{x^*}{x} \le \iprod{y^*}{y} \text{ for all } Ay = x,\ x\in X}.
    \end{equation*}
    If $y^* \not\in \range{A^*}$, then $\range{A^*} \perp \kernel{A} \ne \emptyset$, so for every $x^*\in X^*$ and $x\in X$ we can choose $y\in Y$ such that the above condition is not satisfied.
    Therefore $\upperadj H(y)=\emptyset$.
    Otherwise, if $y^*=A^*\tilde x^*$, we continue to calculate
    \begin{equation*}
        \upperadj H(y^*)
        =
        \setof{x^*\in X^*}{\iprod{x^*}{x}_X \le \iprod{\tilde x^*}{x}_X \text{ for all } x\in X}
        = \{\tilde x^*\}.
    \end{equation*}
    Therefore
    \begin{equation*}
        \upperadj H(y^*)=\{x^* \in X \mid A^* x^*=y^*\}.
    \end{equation*}
    A similar argument shows that $\loweradj H(y^*) = \upperadj H(y^*)$.
\end{example}

These examples and \cref{ex:graphical:linear} suggest the adjoint relationships of the next corollary.
Note that in infinite-dimensional spaces, we only have a relationship between the limiting derivatives, i.e., between the Clarke graphical derivative and the limiting coderivative.

\begin{corollary}
    \label{cor:graphical:adjoint}
    Let $X,Y$ be Banach spaces and $F: X \setto Y$.
    \begin{enumerate}
        \item\label{item:graphical:adjoint:fundamental}
            If $X$ and $Y$ are finite-dimensional, then
            \begin{equation*}
                \frechetCod F(x|y) = \upperadj{D F(x|y)}.
            \end{equation*}

        \item\label{item:graphical:adjoint:limiting}
            If $X$ and $Y$ are reflexive and Gateaux smooth (in particular, if they are finite-dimensional), and $\graph F$ is closed near $(x, y)$, then
            \begin{equation*}
                \clarkeGD F(x|y) = \loweradj{\coderivative F(x|y)}.
            \end{equation*}
    \end{enumerate}
\end{corollary}

\begin{proof}
    \emph{\ref{item:graphical:adjoint:fundamental}:}
    Identifying $X^*$ with $X$ and $Y^*$ with $Y$ in finite dimension, we have by definition that
    \begin{align*}
        D F(x|y)(\dir x) &= \setof{ \dir y\in Y}{(\dir x, \dir y) \in T_{\graph F}(x, y)}
        \shortintertext{and}
        \frechetCod F(x|y)(\dir y) &= \setof{\dir x \in X}{(\dir x, -\dir y) \in \frechetNormal_{\graph F}(x, y)}.
    \end{align*}
    Using \cref{lemma:cones:fundamental-polar}\,\cref{item:cones:fundamental-polar:findim}, we then see that
    \begin{equation*}
        \begin{aligned}
            x^* \in \upperadj{D F(x|y)}(y^*)
            &
            \equivalent
            \iprod{x^*}{\dir x}_X \le \iprod{y^*}{\dir y}_Y
            \text{ for }
            \dir y \in D F(x|y)(\dir x)
            \\
            &
            \equivalent
            \iprod{x^*}{\dir x}_X + \iprod{-y^*}{\dir y}_Y \le 0
            \text{ for }
            (\dir x, \dir y) \in T_{\graph F}(x, y)
            \\
            &
            \equivalent
            (x^*, -y^*) \in \polar{T_{\graph F}(x, y)} = \frechetNormal_{\graph F}(x, y)
            \\
            &
            \equivalent
            x^* \in \frechetCod F(x|y)(y^*).
        \end{aligned}
    \end{equation*}
    This proves the claim.

    \emph{\ref{item:graphical:adjoint:limiting}:}
    We proceed analogously to \ref{item:graphical:adjoint:fundamental} using \cref{thm:cones:limiting-polar} (or \cref{thm:cones:limiting-polar:findim} if $X$ and $Y$ are finite-dimensional):
    \begin{equation*}
        \begin{aligned}[b]
            \dir y \in \loweradj{D^* F(x|y)}(\dir x)
            &
            \equivalent
            \iprod{y^*}{\dir y}_Y \ge \iprod{x^*}{\dir x}_X
            \text{ for }
            x^* \in D^* F(x|y)(y^*)
            \\
            &
            \equivalent
            \iprod{x^*}{\dir x}_X + \iprod{-y^*}{\dir y}_Y \le 0
            \text{ for }
            (x^*, -y^*) \in N_{\graph F}(x, y)
            \\
            &
            \equivalent
            (\dir x, \dir y) \in \polar{N_{\graph F}(x, y)} = \clarkeTangent_{\graph F}(x, y)
            \\
            &
            \equivalent
            \dir y \in \clarkeGD F(x|y)(\dir x).
        \end{aligned}
        \qedhere
    \end{equation*}
\end{proof}

\subsection*{Limiting characterizations in finite dimensions}

In finite dimensions, we can characterize the limiting coderivative and the Clarke derivative directly as inner and outer limits, respectively.

\begin{corollary}
    \label{cor:graphical:limits:findim}
    Let $X$ and $Y$ be finite-dimensional and $F: X\setto Y$. Then for all $(x, y) \in X\times Y$ and all $y^* \in Y$,
    \begin{align}
        \label{eq:graphical:limits:findim:colimiting}
        \coderivative F(x|y)(y^*) &=
        \left\{ x^* \in X
            \;\middle|\;
            \begin{array}{r}
                \text{there exists} \graph F \ni (\alt x, \alt y) \to (x, y)
                \\
                \text{and } (\alt x^*, \alt y^*) \to (x^*, y^*)
                \\
                \text{with } \alt x^* \in \frechetCod F(\alt x|\alt y)(\alt y^*)
            \end{array}
        \right\}.
        \intertext{If $\graph F$ is closed near $(x, y)$, then for all $\dir x \in \R^N$ }
        \label{eq:graphical:limits:findim:gclarke}
        \clarkeGD F(x|y)(\dir x) &=
        \left\{ \dir y \in Y
            \;\middle|\;
            \begin{array}{r}
                \text{for all } \graph F \ni (\alt x, \alt y) \to (x, y)
                \\
                \text{there exists } (\dir \alt x, \dir \alt y) \to (\dir x, \dir y)
                \\
                \text{with } \dir \alt y \in DF(\alt x|\alt y)(\dir \alt x)
            \end{array}
        \right\}.
    \end{align}
\end{corollary}
\begin{proof}
    The characterization \eqref{eq:graphical:limits:findim:colimiting} of the limiting coderivative is a direct application of the definition of the limiting normal cone \eqref{eq:cones:def-limnormal-findim} as an outer limit of the Fréchet normal.
    The characterization \eqref{eq:graphical:limits:findim:gclarke} of the Clarke graphical derivative follows from the characterization of \cref{cor:cones:clarke-liminf} of the Clarke tangent cone as an inner limit of (basic) tangent cones.
\end{proof}

\subsection*{Regularity}

Based on the regularity concepts of sets from \cref{sec:cones:regularity}, we can define concepts of regularity of set-valued mappings.
We say that $F$ at $(x, y) \in \graph F$ (or at $x$ for $y \in F(x)$) is
\begin{enumerate}
    \item \term[mapping!regular!T-]{T-regular} if $D F(x|y)=\clarkeGD F(x|y)$ (i.e., if $\graph F$ has tangential regularity);
    \item \term[mapping!regular!N-]{N-regular} if $\coderivative F(x|y)=\frechetCod F(x|y)$ (i.e., if $\graph F$ has normal regularity).
\end{enumerate}
If $F$ is both T- and N-regular at $(x,y)$, we say that $F$ is \term[mapping!regular!graphically]{graphically regular}.

From \cref{thm:cones:regularity:infdim}, we immediately obtain the following characterization of N-regularity.

\begin{corollary}
    \label{cor:graphical:regularity:infdim}
    Let $X,Y$ be reflexive and Gateaux smooth Banach spaces, $F: X \setto Y$, and let $(x, y) \in \graph F$ with $\graph F$ closed near $(x, y)$.
    Then $F$ is N-regular at $(x,y)$ if and only if $\clarkeGD F(x|y)=\loweradj{[\frechetCod F(x|y)]}$.
\end{corollary}

Writing out various alternatives of \cref{thm:cones:regularity:findim} for set-valued mappings, we obtain full equivalence of the notions and alternative characterizations in finite dimensions.

\begin{corollary}
    \label{cor:graphical:regularity:findim}
    Let $X,Y$ be finite-dimensional and $F:X\setto Y$.
    If $\graph F$ is closed near $(x, y)$, then the following conditions are equivalent:
    \begin{enumerate}
        \item $F$ is N-regular at $x$ for $y$, i.e., $\coderivative F(x|y)=\frechetCod F(x|y)$;
        \item $F$ is T-regular at $x$ for $y$, i.e., $D F(x|y)=\clarkeGD F(x|y)$;
        \item
            $\displaystyle
            \frechetCod F(x|y)(y^*) \supset
            \left\{ x^* \in X
                \;\middle|\;
                \begin{array}{r}
                    \text{there exists } \graph F\ni (\alt x, \alt y) \to (x, y)\\
                    \text{and } (\alt x^*, \alt y^*) \to (x^*, y^*)
                    \\
                    \text{with } \alt x^* \in \frechetCod F(\alt x|\alt y)(\alt y^*)
                \end{array}
            \right\};
            $
        \item
            $\displaystyle
            D F(x|y)(\dir x) \subset
            \left\{ \dir y \in Y
                \;\middle|\;
                \begin{array}{r}
                    \text{for all }(\alt x, \alt y) \to (x, y)\hfill\\
                    \text{there exists }
                    \graph F\ni (\dir \alt x, \dir \alt y) \to (\dir x, \dir y)
                    \\
                    \text{with } \dir \alt y \in DF(\alt x|\alt y)(\dir \alt x)
                \end{array}
            \right\}.
            $
    \end{enumerate}
    In particular, if any of these hold, $F$ is graphically regular at $x$ for $y$.
\end{corollary}

\section{Examples}

As the following examples demonstrate, the graphical derivatives and coderivatives generalize classical (sub)differentials.

\begin{figure}
    \centering
    \begin{subfigure}[t]{0.32\textwidth}
        \begin{asy}
            fill((1,1)--(0,0)--(1,0)--cycle, lightfill);
            draw((-1, 0)--(0,0)--(1,1), primalline+linewidth(1.1));
            draw((0,0)--(1,0), primalline+linewidth(1.1));

            dot((0.0, -0.5), invisible);
            dot((0,0));
            label("$(0,0)$", (0,0), S);
        \end{asy}
        \caption{graph}\label{fig:graphical:polyhedron:graph}
    \end{subfigure}
    \begin{subfigure}[t]{0.32\textwidth}
        \begin{asy}
            fill((1,1)--(0,0)--(1,0)--cycle, lightfill);
            draw((-1, 0)--(0,0)--(1,1), primalline+linewidth(1.1));
            draw((0,0)--(1,0), primalline+linewidth(1.1));

            draw((0.5,0.5)--(0.0, 1.0), dualline+linewidth(1.1), Arrow);
            draw((0.75,0.0)--(0.75, -0.5), dualline+linewidth(1.1), Arrow);
            draw((-0.75,0.5)--(-0.75, -0.5), dualline+linewidth(1.1), Arrows);
            dot((0.8, 0.35), dualline);
            draw((0.0,0.0)--(0.0, -0.5), dualline+linewidth(1.1), Arrow);

            dot((-0.75, 0.0), primalline);
            dot((0.75, 0.0), primalline);
            dot((0.5, 0.5), primalline);
            dot((0.0, 0.0), primalline);
        \end{asy}
        \caption{Fréchet normal cones}\label{fig:graphical:polyhedron:frechet}
    \end{subfigure}
    \begin{subfigure}[t]{0.32\textwidth}
        \begin{asy}
            fill((1,1)--(0,0)--(1,0)--cycle, lightfill);
            draw((-1, 0)--(0,0)--(1,1), primalline+linewidth(1.1));
            draw((0,0)--(1,0), primalline+linewidth(1.1));

            draw((0.5,0.5)--(0.0, 1.0), dualline+linewidth(1.1), Arrow);
            draw((0.75,0.0)--(0.75, -0.5), dualline+linewidth(1.1), Arrow);
            draw((-0.75,0.5)--(-0.75, -0.5), dualline+linewidth(1.1), Arrows);
            dot((0.8, 0.35), dualline);
            draw((0.0,0.5)--(0.0, -0.5), dualline+linewidth(1.1), Arrows);
            draw((0.0,0.0)--(-0.5, 0.5), dualline+linewidth(1.1), Arrow);

            dot((-0.75, 0.0), primalline);
            dot((0.75, 0.0), primalline);
            dot((0.5, 0.5), primalline);
            dot((0.0, 0.0), primalline);
        \end{asy}
        \caption{limiting normal cones}\label{fig:graphical:polyhedron:limiting}
    \end{subfigure}
    \\
    \hspace{0.32\textwidth}
    \begin{subfigure}[t]{0.32\textwidth}
        \begin{asy}
            real eps = 0.15;
            real eps2 = 0.25;

            path semicircle = (-eps,0)
                ..controls(-eps/sqrt(2), eps/sqrt(2))
                ..(0,eps)
                ..controls(eps/sqrt(2), eps/sqrt(2))
                ..(eps,0)
                --cycle;

            fill((1,1)--(0,0)--(1,0)--cycle, lightfill);
            draw((-1, 0)--(0,0)--(1,1), primalline+linewidth(1.1));
            draw((0,0)--(1,0), primalline+linewidth(1.1));

            fill((eps2, eps2)--(0,0)--(eps2,0)--cycle, tertfill);
            draw((-eps2, 0)--(0,0)--(eps2,eps2), tertline+linewidth(1.2), Arrows);
            draw((0,0)--(eps2,0), tertline+linewidth(1.2), Arrow);

            fill(shift((0.8, 0.35))*scale(eps)*unitcircle, tertfill+linewidth(1.2));
            draw((-0.75-eps2, 0)--(-0.75+eps2, 0), tertline+linewidth(1.2), Arrows);
            fill(shift((0.5, 0.5))*rotate(-135)*semicircle, tertfill+linewidth(1.2));
            fill(shift((0.75, 0.0))*semicircle, tertfill+linewidth(1.2));

            dot((-0.75, 0.0), primalline);
            dot((0.75, 0.0), primalline);
            dot((0.5, 0.5), primalline);
            dot((0.0, 0.0), primalline);
            dot((0.8, 0.35), primalline);
        \end{asy}
        \caption{tangent cones}\label{fig:graphical:polyhedron:tangent}
    \end{subfigure}
    \begin{subfigure}[t]{0.32\textwidth}
        \begin{asy}
            real eps = 0.15;
            real eps2 = 0.25;

            path semicircle = (-eps,0)
                ..controls(-eps/sqrt(2), eps/sqrt(2))
                ..(0,eps)
                ..controls(eps/sqrt(2), eps/sqrt(2))
                ..(eps,0)
                --cycle;

            fill((1,1)--(0,0)--(1,0)--cycle, lightfill);
            draw((-1, 0)--(0,0)--(1,1), primalline+linewidth(1.1));
            draw((0,0)--(1,0), primalline+linewidth(1.1));

            draw((0,0)--(eps2,0), tertline+linewidth(1.2), Arrow);

            fill(shift((0.8, 0.35))*scale(eps)*unitcircle, tertfill+linewidth(1.2));
            draw((-0.75-eps2, 0)--(-0.75+eps2, 0), tertline+linewidth(1.2), Arrows);
            fill(shift((0.5, 0.5))*rotate(-135)*semicircle, tertfill+linewidth(1.2));
            fill(shift((0.75, 0.0))*semicircle, tertfill+linewidth(1.2));

            dot((-0.75, 0.0), primalline);
            dot((0.75, 0.0), primalline);
            dot((0.5, 0.5), primalline);
            dot((0.0, 0.0), primalline);
            dot((0.8, 0.35), primalline);
        \end{asy}
        \caption{Clarke tangent cones}\label{fig:graphical:polyhedron:clarke}
    \end{subfigure}
    \caption{The mapping of \cref{ex:graphical:polyhedron} along with the normal and tangent cones at several points. Observe the non-convexity of the limiting normal cone at the point $(0, 0)$. The circles and semicircles represent full and half-spaces.
    }
    \label{fig:graphical:polyhedron}
\end{figure}

\begin{example}
    \label{ex:graphical:polyhedron}
    Consider the set-valued mapping
    \[
        F:\R\setto\R,\qquad
        F(x) = \begin{cases}
            \{0\} & \text{if } x \le 0, \\
            [0, x] & \text{if } x > 0,
        \end{cases}
    \]
    see \cref{fig:graphical:polyhedron:graph}.
    To compute the derivatives and coderivatives, we consider several cases:
    \begin{enumerate}
        \item For any $x<0$, we have $F(x)=\{0\}$. Hence $F$ locally coincides with a linear operator, and hence \cref{ex:graphical:linear} gives
            \begin{align*}
                D^*F(x|0)(y^*) &= \frechetCod F(x|0)(y^*) = \{0\}
                &&\text{for all } y^* \in \R;
                \\
                DF(x|0)(\dir x) &= \clarkeGD F(x|0)(\dir x) = \{0\}
                &&\text{for all }\dir x \in \R.
            \end{align*}

        \item
            Let $x>0$ and consider the derivatives at $x$ for $y \in (0, x)$.
            Then $(x, y) \in \interior \graph F$.
            This implies that
            $
            N_{\graph F}(x, y) = \frechetNormal_{\graph F}(x, y)=\{(0,0)\},
            $
            see \cref{fig:graphical:polyhedron:frechet,fig:graphical:polyhedron:limiting}.
            Hence $\graph F$ is normally regular at $(x,y)$, and \cref{thm:cones:regularity:findim,cor:cones:regularity} yield that
            \[
                T_{\graph F}(x, y) = \clarkeTangent_{\graph F}(x, y)=\R^2,
            \]
            see \cref{fig:graphical:polyhedron:tangent,fig:graphical:polyhedron:clarke}.
            Consequently, for all $\dir x, y^* \in \R$,
            \begin{align*}
                DF(x|y)(\dir x) & = \clarkeGD F(x|y)(\dir x) = \R;
                \\
                D^*F(x|y)(y^*) & = \frechetCod F(x|y)(y^*) =
                \begin{cases}
                    \{0\} & \text{if } y^* = 0 , \\
                    \emptyset &  \text{otherwise}.
                \end{cases}
            \end{align*}
        \item Let $x>0$ and consider the derivatives at $x$ for $y=0$.
            Then it is straightforward to verify from the definitions that
            $
            N_{\graph F}(x, 0) = \frechetNormal_{\graph F}(x, 0) = \{0\} \times (-\infty, 0],
            $
            see \cref{fig:graphical:polyhedron:frechet,fig:graphical:polyhedron:limiting}.
            Hence $\graph F$ is normally regular at $(x,0)$, and \cref{thm:cones:regularity:findim,cor:cones:regularity} yield that
            \[
                T_{\graph F}(x, 0) = \clarkeTangent_{\graph F}(x, 0) = \R \times [0, \infty),
            \]
            see \cref{fig:graphical:polyhedron:tangent,fig:graphical:polyhedron:clarke}.
            Consequently, for all $\dir x, y^* \in \R$,
            \begin{align*}
                DF(x|0)(\dir x) &= \clarkeGD F(x|0)(\dir x) = [0, \infty);
                \\
                D^*F(x|0)(y^*) &= \frechetCod F(x|0)(y^*) =
                \begin{cases}
                    \{0\} &\text{if } y^* \ge 0 , \\
                    \emptyset & \text{otherwise}.
                \end{cases}
            \end{align*}

        \item Let $x>0$ and consider the derivatives at $x$ for $y=x$.
            Then it is again straightforward to verify that
            $
            N_{\graph F}(x, x) = \frechetNormal_{\graph F}(x, x) = \{(x^*, -x^*) \mid x^* \le 0\},
            $
            see \cref{fig:graphical:polyhedron:frechet,fig:graphical:polyhedron:limiting}.
            Hence $\graph F$ is normally regular at $(x,x)$, and \cref{thm:cones:regularity:findim,cor:cones:regularity} yield that
            \[
                T_{\graph F}(x, 0) = \clarkeTangent_{\graph F}(x, 0) = \{(\dir x, \dir y) \mid \dir y, \dir x \in \R,\, \dir x \ge \dir y\},
            \]
            see \cref{fig:graphical:polyhedron:tangent,fig:graphical:polyhedron:clarke}.
            Consequently, for all $\dir x, y^* \in \R$,
            \begin{align*}
                DF(x|x)(\dir x) &= \clarkeGD F(x|x)(\dir x) = (-\infty, \dir x];
                \\
                D^*F(x|x)(y^*) &= \frechetCod F(x|x)(y^*) =
                \begin{cases}
                    \{y^*\} &\text{if } y^* \le 0, \\
                    \emptyset & \text{otherwise}.
                \end{cases}
            \end{align*}
        \item Let $x=0$. Then $\graph F$ is locally a union of cones, so it is again straightforward to verify that
            \begin{equation*}
                T_{\graph F}(x, 0)=\graph F, \qquad\frechetNormal_{\graph F}(0, 0) = \{0\} \times (-\infty, 0],
            \end{equation*}
            see \cref{fig:graphical:polyhedron:tangent,fig:graphical:polyhedron:frechet}.
            However, by the previous cases,
            \[
                N_C(0, 0)
                = \limsup_{\graph F \ni (\tilde x, \tilde y) \to 0} \frechetNormal_C(\tilde x, \tilde y)
                =  \{(x^*, -x^*) \mid x^* \le 0\} \union (\{0\} \times \R)
            \]
            while by \cref{cor:cones:clarke-liminf},
            \[
                \begin{aligned}
                    \clarkeTangent(0, 0)
                    = \liminf_{\graph F \ni (\tilde x, \tilde y) \to 0} T_C(\tilde x, \tilde y)
                    &= (\R \times \{0\}) \isect \{(\dir x, \dir y) \mid \dir x \ge \dir y\}
                    \\
                    &= [0, \infty) \times \{0\},
                \end{aligned}
            \]
            see \cref{fig:graphical:polyhedron:limiting,fig:graphical:polyhedron:clarke}.
            It follows for all $\dir x, y^* \in \R$ that
            \begin{align*}
                \frechetCod F(0|0)(y^*) & =
                \begin{cases}
                    \{0\} &\text{if } y^* \ge 0, \\
                    \emptyset & \text{otherwise},
                \end{cases}
                &
                D^* F(0|0)(y^*) & =
                \begin{cases}
                    \{0\} &\text{if } y^* > 0, \\
                    \{0, y^*\} &\text{otherwise},
                \end{cases}
                \\
                \clarkeGD F(0|0)(\dir x) & = \begin{cases}
                    \{0\} &\text{if } \dir x \ge 0, \\
                    \emptyset & \text{otherwise},
                \end{cases}
                &
                DF(0|0)(\dir x) & = \begin{cases}
                    [0, \dir x] &\text{if } \dir x \ge 0,
                    \\
                    \{0\} &\text{otherwise}.
                \end{cases}
            \end{align*}
    \end{enumerate}
\end{example}

\subsection*{Single-valued mappings and their inverses}

For the Clarke graphical derivative and the limiting coderivatives (which are obtained as inner or outer limits), we have to require -- just as for the Clarke subdifferential in \cref{thm:clarke:frechet} -- slightly more than just Fréchet differentiability.
\begin{theorem}%
    \label{thm:graphical:single}
    Let $X,Y$ be Banach spaces and let $F: X \to Y$ be single-valued and Fréchet differentiable at $x\in X$.
    Then
    \begin{align*}
        DF(x|y)(\dir x)&=\begin{cases}
            \{F'(x)\dir x\} & \text{if } y=F(x), \\
            \emptyset & \text{otherwise},
        \end{cases}
        \\
        \shortintertext{and}
        \frechetCod F(x|y)(y^*)&=\begin{cases}
            \{F'(x)^*y^*\} & \text{if }  y=F(x), \\
            \emptyset & \text{otherwise}.
        \end{cases}
    \end{align*}

    If $F$ is continuously differentiable at $x$, then $F$ is graphically regular at $x$ for $F(x)$, and hence the corresponding expressions also hold for $\clarkeGD F(x|y)$ and $\coderivative F(x|y)$.
\end{theorem}

\begin{proof}
    \emph{The graphical derivative:}
    We have $(\dir x, \dir y) \in T_{\graph F}(x, y)$ if and only if for some $x_k \to x$, $y_k \defeq F(x_k)$, and $\tau_k \downto 0$ there holds
    \begin{subequations}
        \begin{align}
            \dir x &= \lim_{k \to \infty} \frac{x_k-x}{\tau_k} =: \lim_{k\to\infty} \dir x_k
            \shortintertext{and}
            \label{eq:graphical:single:gd1}
            \dir y &= \lim_{k \to \infty} \frac{y_k-y}{\tau_k} = \lim_{k \to \infty} \frac{F(x+\tau_k\dir x_k)-F(x)}{\tau_k}.
        \end{align}
    \end{subequations}
    If $\dir x_k=0$ for all sufficiently large $k \in \N$, clearly both $\dir x=0$ and $\dir y=0$.
    This satisfies the claimed expression.
    So we may assume that $\dir x_k \ne 0$ for all $k\in\N$.
    In this case, \eqref{eq:graphical:single:gd1} holds if and only if
    \begin{equation*}
        \lim_{k \to \infty} \frac{F(x+h_k)-F(x)-\tau_k\dir y_k}{\norm{h_k}_X}=0
    \end{equation*}
    for $h_k \defeq \tau_k\dir x_k$ and any $\dir y_k \to \dir y$. Since $F$ is Fréchet differentiable, this clearly holds with
    \begin{equation*}
        \dir y_k \defeq \inv \tau_k F'(x)h_k=F'(x)\dir x_k \to   F'(x)\dir x=:\dir y.
    \end{equation*}
    This shows that $DF(x|y)(\dir x)=\{F'(x)\dir x\}$.

    \emph{The Clarke graphical derivative:}
    To calculate $\clarkeGD F(x|y)$, we have to find all $\dir x$ and $\dir y$ such that
    for every $\tau_k \downto 0$ and $(\alt x_k, \alt y_k) \to (x, y)$ with $\alt y_k=F(\alt x_k)$, there exists $x_k\to x$ with
    \begin{equation*}
        \dir x = \lim_{k \to \infty} \frac{x_k-\alt x_k}{\tau_k}
        \qquad\text{and}\qquad
        \dir y = \lim_{k \to \infty} \frac{F(x_k)-F(\alt x_k)}{\tau_k}.
    \end{equation*}
    Setting $x_k=\alt x_k+\tau_k\dir x_k$ with $\dir x_k \to \dir x$, the second condition becomes
    \begin{equation*}
        \dir y = \lim_{k \to \infty} \frac{F(\alt x_k+\tau_k\dir x_k)-F(\alt x_k)}{\tau_k}.
    \end{equation*}
    Taking $\alt x_k=x$, arguing as for $DF$ shows that $\dir y = F'(x)\dir x$ is the only candidate. It just remains to show that any choice of $\alt x_k$ gives the same limit, i.e., that
    \begin{equation*}
        \lim_{k \to \infty} \frac{F(\alt x_k+\tau_k\dir x_k)-F(\alt x_k)-\tau_k F'(x)\dir x}{\tau_k}=0.
    \end{equation*}
    But this follows from the assumed continuous differentiability using \cref{lem:frechet:diffquot}. Thus for $y=F(x)$,
    \begin{equation*}
        \clarkeGD F(x|y)(\dir x)=\{F'(x)\dir x\}=DF(x|y)(\dir x).
    \end{equation*}
    This shows that $F$ is T-regular at $x$ for $y$.

    \emph{The Fréchet coderivative:}
    The claim follows from proving that
    \begin{equation}
        \label{eq:graphical:single:epsilon}
        \frechetCod_\epsilon F(x|y)(y^*)=\begin{cases}
            \B(F'(x)^*y^*, \epsilon)& \text{if } y=F(x), \\
            \emptyset & \text{otherwise}.
        \end{cases}
    \end{equation}
    To show this, we note that $x^* \in \frechetCod_\epsilon F(x|y)(y^*)$ if and only if for every sequence $x_k \to x$ with $F(x_k) \to F(x)$,
    \begin{equation*}
        \limsup_{k \to \infty} \frac{\iprod{x^*}{x_k-x}_X - \iprod{y^*}{F(x_k)-F(x)}_Y}{\sqrt{\norm{x_k-x}_X^2+\norm{F(x_k)-F(x)}_Y^2}} \le \epsilon.
    \end{equation*}
    Dividing both numerator and denominator by $\norm{x_k-x}_X>0$, we obtain the equivalent condition that
    \begin{equation*}
        \limsup_{k \to \infty} q_k \le \epsilon \qquad\text{for}\qquad q_k \defeq \frac{\iprod{x^*}{x_k-x}_X - \iprod{y^*}{F(x_k)-F(x)}_Y}{\norm{x_k-x}_X}.
    \end{equation*}
    If we take $x^* \in \B(F'(x)^*y^*, \epsilon)$, this condition is verified by the Fréchet differentiability of $F$ at $x$.
    Conversely, to show that this implies $x^* \in \B(F'(x)^*y^*, \epsilon)$, we take $x_k \defeq x+\tau_k h$ for some $\tau_k \downto 0$ and $h \in X$ with $\norm{h}_X=1$. Then again by the Fréchet differentiability of $F$,
    \begin{equation*}
        \epsilon \geq \lim_{k \to \infty} q_k = \iprod{x^*}{h} - \iprod{y^*}{F'(x)h}.
    \end{equation*}
    Since $h\in \B_{X}$ was arbitrary, this shows that $x^* \in \B(F'(x)^*y^*, \epsilon)$.

    \emph{The limiting coderivative:}
    By the definition \eqref{eq:cones:def-limnormal}, the formula \eqref{eq:graphical:single:epsilon} for $\epsilon$-coderivatives, and the continuous differentiability, we have
    \begin{equation*}
        \begin{aligned}
            N_{\graph F}(x, F(x))
            &
            = \weakstarlimsup_{\alt x \to x,\, \epsilon \downto 0} \frechetNormal_{\graph F}^{\epsilon}(\alt x, F(\alt x))
            \\
            &
            = \weakstarlimsup_{\alt x \to x\, \epsilon \downto 0}\,\{(y^*, F'(\alt x)^*y^*+z^*)\in Y^*\times X^* \mid y^* \in Y^*,\, z^* \in \B(0, \epsilon)\}
            \\
            &
            = \weakstarlimsup_{\alt x \to x}\,\{(y^*, F'(\alt x)^*y^*) \in Y^*\times X^*\mid y^* \in Y^*\}
            \\
            &
            = \{(y^*, F'(x)^*y^*) \in Y^*\times X^*\mid y^* \in Y^*\}.
        \end{aligned}
    \end{equation*}
    This shows the claimed formula for the limiting coderivative and hence N- and therefore graphical regularity.
\end{proof}

\begin{remark}
    In finite dimensional spaces, it would be possible to more concisely prove the expression for $\clarkeGD F(x|y)$ using \cref{cor:cones:clarke-liminf}.
    Likewise, we could use the polarity relationships of \cref{cor:graphical:adjoint} to obtain the expression for $\frechetCod F(x|y)$.
    These approaches will, however, not be possible in more general spaces.
\end{remark}

Combining \cref{thm:graphical:single} with \cref{lemma:graphical:inverse} allows us to compute the graphical derivatives and coderivatives of inverses of single-valued functions.

\begin{corollary}
    \label{thm:graphical:single:inverse}
    Let $X,Y$ be Banach spaces and let $F: X \to Y$ be single-valued and Fréchet differentiable at $x\in X$.
    Then
    \begin{align*}
        D\inv F(y|x)(\dir y)&=\begin{cases}
            \{\dir x\in X \mid F'(x)\dir x=\dir y\} & \text{if }y=F(x), \\
            \emptyset & \text{otherwise},
        \end{cases}
        \\
        \intertext{and}
        \frechetCod \inv F(y|x)(x^*)&=\begin{cases}
            \{y^* \in Y^* \mid F'(x)^*y^*=x^*\} &\text{if } y=F(x), \\
            \emptyset & \text{otherwise}.
        \end{cases}
    \end{align*}

    If $F$ is continuously differentiable at $x$, then $\inv F$ is graphically regular at $y=F(x)$ for $x$, and hence the corresponding expressions also hold for $\clarkeGD \inv F(y|x)$ and $\coderivative \inv F(y|x)$.
\end{corollary}

It is important that \cref{thm:graphical:single} concerns the strong graphical derivatives $DF$ instead of the weak graphical derivative $D^w F$. Indeed, as the next counter-example demonstrates, $D^w F$ is more of a theoretical tool (with the important property in reflexive spaces that $\frechetCod F(x|y) = \upperadj{D^w F(x|y)}$ by \cref{lemma:cones:fundamental-polar}\,\cref{item:cones:fundamental-polar:reflexive}) which does not enjoy a rich calculus consistent with conventional notions. In the following chapters, we will therefore not develop calculus rules for the weak graphical derivative.

\begin{example}[counter-example to single-valued weak graphical derivatives]
    Let $f\in C^1(\R)$, $\Omega\subset \R^d$ be open, and
    \begin{equation*}
        F: L^2(\Omega)\to\R,\qquad F(u)=\int_0^1 f(u(x)) \,d x.
    \end{equation*}
    Then by the above,
    \begin{equation*}
        DF(u|F(u))(\dir u)=\left\{\int_0^1 f'(u(x))\dir u(x) \,d x\right\}.
    \end{equation*}
    In particular, $DF(u|F(u))(0)=\{0\}$.

    However, choosing, e.g., $f(t)=\sqrt{1+t^2}$, $\Omega = (0,1)$, and
    $u_k(x) \defeq \sign \sin(2^k\pi x)$,
    we have $u_k \weakto 0$ in $L^2(\Omega)$ but $\abs{u_k(x)}=1$ for a.e.~$x \in [0, 1]$.
    Take now $\tilde u_k \defeq \alpha \tau_k u_k$ for any given $\tau_k \downto 0$ and $\alpha>0$. Then $\tilde u_k \weakto 0$ as well,  while
    \begin{equation*}
        F(\tilde u_k)-F(0) = \sqrt{1+\alpha^2 \tau_k^2}-1 \to 0.
    \end{equation*}
    Moreover, $(\tilde u_k-0)/\tau_k=\alpha u_k \weakto 0$ and $\lim_{k \to \infty} \left(\sqrt{1+\alpha^2\tau_k^2}-1\right)/\tau_k=\alpha^2$.
    As $\alpha>0$ was arbitrary, we deduce that $D^w F(u|F(u))(0) \supset [0, \infty)$.
\end{example}

\subsection*{Derivatives and coderivatives of subdifferentials}

We now apply these notions to set-valued mappings arising as subdifferentials of convex functionals.
First, we directly obtain from \cref{thm:graphical:single} an expression for the squared norm in Hilbert spaces.

\begin{corollary}
    \label{ex:graphical:norm22}
    Let $X$ be a Hilbert space and $f(x) = \frac12\norm{x}_X^2$ for $x\in X$.
    Then
    \begin{equation*}
        \clarkeGD[\subdiff f](x|y)(\dir x) =
        D[\subdiff f](x|y)(\dir x) =
        \begin{cases}
            \{\dir x\}  & \text{if }y=x, \\
            \emptyset & \text{otherwise},
        \end{cases}
    \end{equation*}
    and
    \begin{equation*}
        \coderivative [\subdiff f](x|y)(y^*) =
        \frechetCod [\subdiff f](x|y)(y^*) =
        \begin{cases}
            \{y^*\}  & \text{if }y=x, \\
            \emptyset & \text{otherwise}.
        \end{cases}
    \end{equation*}
    In particular, $\subdiff f$ is graphically regular at every $x\in X$.
\end{corollary}

Of course, we are more interested in subdifferentials of nonsmooth functionals.
We first study the indicator functional of an interval; see \cref{fig:graphical:indicator}.

\begin{figure}[htp!]
    \centering
    \begin{subfigure}[t]{0.32\textwidth}
        \centering
        \begin{asy}
            unitsize(40,40);
            real l=1.5;
            real eps=0.3;
            pair gfstart=(-1, -l);
            pair gfend=(1, l);
            path gf=gfstart--(-1, 0)--(1, 0)--gfend;
            draw(gf, dashed, Arrows);

            draw((-eps,0)--(eps, 0), primalline+linewidth(1.1), Arrows);
            dot((0, 0));

            draw((-1,-eps)--(-1,0)--(-1+eps, 0), primalline+linewidth(1.1), Arrows);
            dot((-1, 0));

            draw((1,0.75+eps)--(1, 0.75-eps), primalline+linewidth(1.1), Arrows);
            dot((1, 0.75));
        \end{asy}
        \caption{graphical derivative $D[\subdiff f]$}
        \label{fig:graphical:indicator:gderiv}
    \end{subfigure}
    \hfill
    \begin{subfigure}[t]{0.32\textwidth}
        \centering
        \begin{asy}
            unitsize(40,40);
            real l=1.5;
            real eps=0.3;
            pair gfstart=(-1, -l);
            pair gfend=(1, l);
            path gf=gfstart--(-1, 0)--(1, 0)--gfend;
            draw(gf, dashed, Arrows);

            draw((-eps, 0)--(eps, 0), primalline+linewidth(1.1), Arrows);
            dot((0, 0));

            fill((-1,-eps)--(-1,0)--(-1+eps, 0)..controls (-1+eps/sqrt(2), -eps/sqrt(2))..cycle, darkfill);
            draw((-1,-eps)--(-1,0)--(-1+eps, 0), primalline+linewidth(1.1), Arrows);
            dot((-1, 0));

            draw((1,0.75+eps)--(1, 0.75-eps), primalline+linewidth(1.1), Arrows);
            dot((1, 0.75));
        \end{asy}
        \caption{convex hull $\conv D[\subdiff f]$}
    \end{subfigure}
    \hfill
    \begin{subfigure}[t]{0.32\textwidth}
        \centering
        \begin{asy}
            unitsize(40,40);
            real l=1.5;
            real eps=0.3;
            pair gfstart=(-1, -l);
            pair gfend=(1, l);
            path gf=gfstart--(-1, 0)--(1, 0)--gfend;
            draw(gf, dashed, Arrows);

            draw((0, -eps)--(0, eps), primalline+linewidth(1.1), Arrows);
            dot((0, 0));

            fill((-1-eps,0)--(-1,0)--(-1, eps)..controls (-1-eps/sqrt(2), eps/sqrt(2))..cycle, darkfill);
            draw((-1-eps,0)--(-1,0)--(-1, eps), primalline+linewidth(1.1), Arrows);
            dot((-1, 0));

            draw((1-eps,0.75)--(1+eps, 0.75), primalline+linewidth(1.1), Arrows);
            dot((1, 0.75));
        \end{asy}
        \caption{Fréchet coderivative $\frechetCod[\subdiff f]$}
    \end{subfigure}

    \medskip

    \begin{subfigure}[t]{0.32\textwidth}
        \centering
        \begin{asy}
            unitsize(40,40);
            real l=1.5;
            real eps=0.3;
            pair gfstart=(-1, -l);
            pair gfend=(1, l);
            path gf=gfstart--(-1, 0)--(1, 0)--gfend;
            draw(gf, dashed, Arrows);

            draw((0, -eps)--(0, eps), primalline+linewidth(1.1), Arrows);
            dot((0, 0));

            fill((-1-eps,0)--(-1,0)--(-1, eps)..controls (-1-eps/sqrt(2), eps/sqrt(2))..cycle, darkfill);
            draw((-1-eps,0)--(-1+eps, 0), primalline+linewidth(1.1), Arrows);
            draw((-1, eps)--(-1, -eps), primalline+linewidth(1.1), Arrows);
            dot((-1, 0));

            draw((1-eps,0.75)--(1+eps, 0.75), primalline+linewidth(1.1), Arrows);
            dot((1, 0.75));
        \end{asy}
        \caption{limiting~coderivative $\coderivative[\subdiff f]$}
    \end{subfigure}
    \hfill
    \begin{subfigure}[t]{0.32\textwidth}
        \centering
        \begin{asy}
            unitsize(40,40);
            real l=1.5;
            real eps=0.3;
            pair gfstart=(-1, -l);
            pair gfend=(1, l);
            path gf=gfstart--(-1, 0)--(1, 0)--gfend;
            draw(gf, dashed, Arrows);

            draw((0, -eps)--(0, eps), primalline+linewidth(1.1), Arrows);
            dot((0, 0));

            fill((-1-eps,0)
            ..controls(-1-eps/sqrt(2), eps/sqrt(2))
            ..(-1,eps)
            ..controls(-1+eps/sqrt(2), eps/sqrt(2))
            ..(-1+eps,0)
            ..controls(-1+eps/sqrt(2), -eps/sqrt(2))
            ..(-1,-eps)
            ..controls(-1-eps/sqrt(2), -eps/sqrt(2))
            ..cycle, darkfill);
            draw((-1-eps,0)--(-1+eps, 0), primalline+linewidth(1.1), Arrows);
            draw((-1, eps)--(-1, -eps), primalline+linewidth(1.1), Arrows);
            dot((-1, 0));

            draw((1-eps,0.75)--(1+eps, 0.75), primalline+linewidth(1.1), Arrows);
            dot((1, 0.75));
        \end{asy}
        \caption{convex hull $\conv \coderivative[\subdiff f]$}
    \end{subfigure}
    \hfill
    \begin{subfigure}[t]{0.32\textwidth}
        \centering
        \begin{asy}
            unitsize(40,40);
            real l=1.5;
            real eps=0.3;
            pair gfstart=(-1, -l);
            pair gfend=(1, l);
            path gf=gfstart--(-1, 0)--(1, 0)--gfend;
            draw(gf, dashed, Arrows);

            draw((-eps,0)--(eps, 0), primalline+linewidth(1.1), Arrows);
            dot((0, 0));

            dot((-1, 0));

            draw((1,0.75+eps)--(1, 0.75-eps), primalline+linewidth(1.1), Arrows);
            dot((1, 0.75));
        \end{asy}
        \caption{Clarke graphical derivative $\clarkeGD[\subdiff f]$}
    \end{subfigure}
    \caption{Illustration of the different graphical derivatives and coderivatives of $\subdiff f$ for $f=\delta_{[-1,1]}$. The dashed line is $\graph \subdiff f$. The dots indicate the base points $(x, y)$ where $D[\subdiff f](x|y)$ is calculated, and the thick arrows and filled-in areas the directions of $(\Delta x, \Delta y)$ (resp.~$(\Delta x, -\Delta y)$ for the coderivatives) relative to the base point.  Observe that there is no graphical regularity at $(x, y)\in \{(-1, 0) ,(1, 0)\}$. Everywhere else, $\subdiff f$ is graphically regular. Observe also that cones in the last figures of each row are polar to the cones in the first and the second figures on the same row.}
    \label{fig:graphical:indicator}
\end{figure}

\begin{theorem}
    \label{lemma:graphical:indicator}
    Let $f(x) \defeq \delta_{[-1,1]}(x)$ for $x\in \R$.
    Then
    \begingroup
    \allowdisplaybreaks
    \begin{align}
        \label{eq:graphical:indicator:gderiv}
        D[\subdiff f]( x | y )(\dir x ) & =
        \begin{cases}
            \R
            &\text{if } \abs{ x }=1,\,  y  \in (0, \infty)  x ,\, \dir x =0, \\
            [0, \infty)  x
            &\text{if } \abs{ x }=1,\,  y =0,\, \dir x  = 0,  \\
            \{0\}
            &\text{if } \abs{ x }=1,\,  y  = 0,\,  x  \dir x  < 0, \\
            \{0\}
            &\text{if } \abs{ x } < 1,\,  y =0,\, \\
            \emptyset & \text{otherwise},
        \end{cases}
        \\[1ex]
        \frechetCod[\subdiff f]( x | y )(y^* ) & =
        \begin{cases}
            \R,
            &\text{if } \abs{ x }=1,\,  y  \in (0, \infty)  x , y^* = 0\\
            [0, \infty)  x
            &\text{if } \abs{ x }=1,\,  y =0,\, x y^* \ge 0,  \\
            \{0\}
            &\text{if } \abs{ x } < 1,\,  y =0, \\
            \emptyset & \text{otherwise},
        \end{cases}
        \\[1ex]
        \label{eq:graphical:indicator:gclarke}
        \clarkeGD[\subdiff f]( x | y )(\dir x ) & =
        \begin{cases}
            \mathrlap{\R}\phantom{[0, \infty)x}
            &\text{if } \abs{ x }=1,\,  y  \in (0, \infty)  x ,\, \dir x =0, \\
            \{0\}
            &\text{if } \abs{ x }=1,\,  y =0,\, \dir x  = 0,  \\
            \{0\} &\text{if } \abs{ x } < 1,\,  y =0,\, \\
            \emptyset
            &\text{otherwise},
        \end{cases}
        \shortintertext{and}
        \label{eq:graphical:indicator:colimiting}
        \coderivative[\subdiff f]( x | y )(y^* ) & =
        \begin{cases}
            \R
            &\text{if } \abs{ x }=1,\,  y  \in [0, \infty)  x , y^* = 0\\
            [0, \infty)  x
            &\text{if } \abs{ x }=1,\,  y =0,\, x y^* > 0,  \\
            \{0\}
            &\text{if } \abs{ x }=1,\,  y =0,\, x y^* < 0,  \\
            \{0\}
            &\text{if } \abs{ x } < 1,\,  y =0, \\
            \emptyset & \text{otherwise}.
        \end{cases}
    \end{align}
    \endgroup
    In particular, $\subdiff f$ is graphically regular at $x$ for $y \in \subdiff f(x)$ if and only if $\abs{x} <1$ or $y \ne 0$.
\end{theorem}

\begin{proof}
    We first of all recall from \cref{ex:convex:subdiff_ind} that $\graph \subdiff f$ is closed with
    \begin{equation}
        \label{eq:graphical:indicator:subdiff}
        \subdiff f(x)=
        \begin{cases}
            [0, \infty) x &\text{if } \abs{ x }=1, \\
            \{0\} &\text{if } \abs{ x } < 1, \\
            \emptyset & \text{otherwise}.
        \end{cases}
    \end{equation}

    We now verify \eqref{eq:graphical:indicator:gderiv}.
    If $ y  \in \subdiff f( x )$ and $\dir{ y } \in  D[\subdiff f]( x | y )(\dir x )$, there exist by \eqref{eq:graphical:gderiv-altdef} sequences $t_k \downto 0$, $x_k \to x$, and $ y_k \in \subdiff f( x +t_k\dir{ x_k})$ such that
    \begin{equation}
        \label{eq:graphical:indicator:limit}
        \dir x = \lim_{k \to \infty} \frac{x_k-x}{t_k}
        \qquad\text{and}\qquad
        \dir y = \lim_{k \to \infty} \frac{y_k-y}{t_k}.
    \end{equation}
    We proceed by case distinction.
    \begin{enumerate}
        \item\label{tem:graphical:indicator:case1} $\abs{ x }=1$, $\dir{ x }=0$, and $y \in (0, \infty) x$: Then choosing $x_k \equiv x$, any $\dir y \in \R$ and $k$ large enough, we can take $y_k=y + t_k\dir y \in [0, \infty) x  = \subdiff f(x)$. This yields the first case of \eqref{eq:graphical:indicator:gderiv}.

        \item
            $\abs{x}=1$, $\dir x = 0$, but $y=0$: In this case, choosing $x_k \equiv x$, we can take any $y_k \in \subdiff f(x+t_k\dir x_k)=\subdiff f(x)=[0, \infty)  x $. Picking any $\dir y \in [0, \infty)  x $ and setting $y_k \defeq y + t_k\dir y$, we deduce that $\dir y \in D[\subdiff f]( x | y )(\dir{ x })$. Thus ``$\supset$'' holds in the second case of \eqref{eq:graphical:indicator:gderiv}. Since $\dir y \in -(0, \infty) x$ is clearly not obtainable with $y_k \in [0, \infty) x$, also ``$\subset$'' holds.

        \item\label{tem:graphical:indicator:case3} $\abs{x}=1$ and $\dir{x}=0$, but $y \in -(0, \infty) x$: Then we have $y_k \in [0, \infty)x$ for $k$ large enough since in this case either $x_k=x$ or $x_k \in (-1, 1)$.
            Thus $\abs{y-y_k} \ge \abs{y}>0$, so the second limit in \eqref{eq:graphical:indicator:limit} cannot exist. Therefore the coderivative is empty, which is covered by the last case of \eqref{eq:graphical:indicator:gderiv}.

        \item $\abs{x}=1$ and $x \dir x>0$: Then the first limit in \eqref{eq:graphical:indicator:limit} requires that $x_k \not \in \dom \subdiff f$, and hence $\subdiff f(x_k)=\emptyset$ for $k$ large enough. This is again covered by the last case of \eqref{eq:graphical:indicator:gderiv}.

        \item $\abs{x}=1$ and $x \dir x < 0$
            (the case $x \dir x = 0$ being covered by \ref{tem:graphical:indicator:case1}--\ref{tem:graphical:indicator:case3}):
            Since $\dir x \ne 0$ has a different sign from $x$, it follows from the first limit in \eqref{eq:graphical:indicator:limit} that $x_k \in (-1, 1)$ for $k$ large enough.
            Consequently, $\subdiff f(x_k)=\{0\}$, i.e., $y_k=0$.
            The limit \eqref{eq:graphical:indicator:limit} in this case only exists if $y =0$, in which case also $\dir{ y }=0$. This is covered by the third case of~\eqref{eq:graphical:indicator:gderiv}, while $y \ne 0$ is covered by the last case.

        \item $\abs{x}<1$: Then $y=0$ and necessarily $y_k=0$ for $k$ large enough. Therefore also $\dir y=0$, which yields the fourth case in \eqref{eq:graphical:indicator:gderiv}.

        \item $\abs{x}>1$: Then $\subdiff f(x)=\emptyset$ and therefore  the coderivative is empty as well, yielding again the final case in \eqref{eq:graphical:indicator:gderiv}.
    \end{enumerate}

    The expression for $\frechetCod[\subdiff f](x|y)$ can be verified using \cref{cor:graphical:adjoint}\,\cref{item:graphical:adjoint:fundamental}.
    It can also be seen graphically from \cref{fig:graphical:indicator}.

    By the inner and outer limit characterizations of \cref{cor:graphical:limits:findim}, we now obtain the expressions for the Clarke graphical derivative $\clarkeGD[\subdiff f]( x | y )$ and the limiting coderivative $\coderivative[\subdiff f](x|y)$.
    Since $\graph \subdiff f$ is locally contained in an affine subspace outside of the \enquote{corner cases} $(x,y) \in \{(1,0),(-1,0)\}$, only the latter need special inspection. For the Clarke graphical derivative, we need to write $\dir y$ as the limit of $\dir y_k \in D[\subdiff f](x_k,y_k)(\dir x_k)$ for some $\dir x_k \to \dir x$ and \emph{all} $\graph \subdiff f\ni (x_k,y_k) \to (x, y)$. Consider for example $(x,y)=(-1,0)$. Trying both $(x_k, y_k)=(-1+1/k,0)$ and $(x_k, y_k)=(-1,-1/k)$, we see that this is only possible for $(\dir x,\dir y)=(\dir x_k,\dir y_k)=(0,0)$. This yields the second case of \eqref{eq:graphical:indicator:gclarke}.
    Conversely, for the limiting coderivative, it suffices to find \emph{one} such sequence from the Fréchet coderivative. Choosing for $(x,y)=(-1,0)$ again $(x_k, y_k)=(-1+1/k,0)$ and $(x_k, y_k)=(-1,-1/k)$ as well as the constant sequence $(x_k,y_k)=(-1,0)$ yields the second, third, and first case of \eqref{eq:graphical:indicator:colimiting}, respectively.

    Finally, in finite dimensions the mapping $\subdiff f$ is graphically regular if and only if $D[\subdiff f]( x | y )=\clarkeGD[\subdiff f]( x | y )$ by \cref{cor:graphical:regularity:findim}, which is the case exactly when $\abs{x} <1$ or $y \ne 0$ as claimed.
\end{proof}
In nonlinear optimization with inequality constraints, the case where $\subdiff f$ is graphically regular corresponds precisely to the case of \term[complementarity, strict]{strict complementarity} of the minimizer $\bar x$ and the Lagrange multiplier $\bar y$ for the constraint $x\in[-1,1]$.

We next study the different derivatives and graphical regularity of the subdifferential of the absolute value function; see \cref{fig:graphical:absvalue}.

\begin{figure}[htp!]
    \centering
    \begin{subfigure}[t]{0.32\textwidth}
        \centering
        \begin{asy}
            unitsize(40,40);
            real l=1.5;
            real eps=0.3;
            pair gfstart=(-l, -1);
            pair gfend=(l, 1);
            path gf=gfstart--(0, -1)--(0, 1)--gfend;
            draw(gf, dashed, Arrows);

            draw((0,-eps)--(0,eps), primalline+linewidth(1.1), Arrows);
            dot((0, 0));

            draw((-eps,-1)--(0,-1)--(0, -1+eps), primalline+linewidth(1.1), Arrows);
            dot((0, -1));

            draw((0.5*l-eps, 1)--(0.5*l+eps, 1), primalline+linewidth(1.1), Arrows);
            dot((0.5*l, 1));

            dot((0,-1-eps), invisible);
        \end{asy}
        \caption{graphical derivative $D[\subdiff f]$}
        \label{fig:graphical:absvalue:gderiv}
    \end{subfigure}
    \hfill
    \begin{subfigure}[t]{0.32\textwidth}
        \centering
        \begin{asy}
            unitsize(40,40);
            real l=1.25;
            real eps=0.3;
            pair gfstart=(-l, -1);
            pair gfend=(l, 1);
            path gf=gfstart--(0, -1)--(0, 1)--gfend;
            draw(gf, dashed, Arrows);

            draw((0, -eps)--(0, eps), primalline+linewidth(1.1), Arrows);
            dot((0, 0));

            fill((-eps, -1)--(0,-1)--(0, -1+eps)..controls (-eps/sqrt(2), -1+eps/sqrt(2))..cycle, darkfill);
            draw((-eps, -1)--(0,-1)--(0, -1+eps), primalline+linewidth(1.1), Arrows);
            dot((0, -1));

            draw((0.5*l-eps, 1)--(0.5*l+eps, 1), primalline+linewidth(1.1), Arrows);
            dot((0.5*l, 1));

            dot((0,-1-eps), invisible);
        \end{asy}
        \caption{convex hull $\conv D[\subdiff f]$}
    \end{subfigure}
    \hfill
    \begin{subfigure}[t]{0.32\textwidth}
        \centering
        \begin{asy}
            unitsize(40,40);
            real l=1.5;
            real eps=0.3;
            pair gfstart=(-l, -1);
            pair gfend=(l, 1);
            path gf=gfstart--(0, -1)--(0, 1)--gfend;
            draw(gf, dashed, Arrows);

            draw((-eps,0)--(eps, 0), primalline+linewidth(1.1), Arrows);
            dot((0, 0));

            fill((0, -1-eps)--(0,-1)--(eps, -1)..controls (eps/sqrt(2), -1-eps/sqrt(2))..cycle, darkfill);
            draw((0, -1-eps)--(0,-1)--(eps, -1), primalline+linewidth(1.1), Arrows);
            dot((0, -1));

            draw((0.5*l, 1+eps)--(0.5*l, 1-eps), primalline+linewidth(1.1), Arrows);
            dot((0.5*l, 1));
        \end{asy}
        \caption{Fréchet coderivative $\frechetCod[\subdiff f]$}
    \end{subfigure}

    \medskip

    \begin{subfigure}[t]{0.32\textwidth}
        \centering
        \begin{asy}
            unitsize(40,40);
            real l=1.5;
            real eps=0.3;
            pair gfstart=(-l, -1);
            pair gfend=(l, 1);
            path gf=gfstart--(0, -1)--(0, 1)--gfend;
            draw(gf, dashed, Arrows);

            draw((-eps,0)--(eps, 0), primalline+linewidth(1.1), Arrows);
            dot((0, 0));

            fill((0, -1-eps)--(0,-1)--(eps, -1)..controls (eps/sqrt(2), -1-eps/sqrt(2))..cycle, darkfill);
            draw((0, -1-eps)--(0,-1)--(eps, -1), primalline+linewidth(1.1), Arrows);
            draw((-eps, -1)--(0,-1)--(0, -1+eps), primalline+linewidth(1.1), Arrows);
            dot((0, -1));

            draw((0.5*l, 1+eps)--(0.5*l, 1-eps), primalline+linewidth(1.1), Arrows);
            dot((0.5*l, 1));
        \end{asy}
        \caption{limiting~coderivative $\coderivative[\subdiff f]$}
        \label{fig:graphical:absvalue:limiting}
    \end{subfigure}
    \hfill
    \begin{subfigure}[t]{0.32\textwidth}
        \centering
        \begin{asy}
            unitsize(40,40);
            real l=1.5;
            real eps=0.3;
            pair gfstart=(-l, -1);
            pair gfend=(l, 1);
            path gf=gfstart--(0, -1)--(0, 1)--gfend;
            draw(gf, dashed, Arrows);

            draw((-eps,0)--(eps, 0), primalline+linewidth(1.1), Arrows);
            dot((0, 0));
            fill((0,-1-eps)
            ..controls(eps/sqrt(2),-1-eps/sqrt(2))
            ..(eps,-1)
            ..controls(eps/sqrt(2),-1+eps/sqrt(2))
            ..(0,-1+eps)
            ..controls(-eps/sqrt(2),-1+eps/sqrt(2))
            ..(-eps,-1)
            ..controls(-eps/sqrt(2),-1-eps/sqrt(2))
            ..cycle, darkfill);
            draw((0, -1-eps)--(0,-1)--(eps, -1), primalline+linewidth(1.1), Arrows);
            draw((-eps, -1)--(0,-1)--(0, -1+eps), primalline+linewidth(1.1), Arrows);
            dot((0, -1));

            draw((0.5*l, 1+eps)--(0.5*l, 1-eps), primalline+linewidth(1.1), Arrows);
            dot((0.5*l, 1));
        \end{asy}
        \caption{convex hull $\conv \coderivative[\subdiff f]$}
    \end{subfigure}
    \hfill
    \begin{subfigure}[t]{0.32\textwidth}
        \centering
        \begin{asy}
            unitsize(40,40);
            real l=1.5;
            real eps=0.3;
            pair gfstart=(-l, -1);
            pair gfend=(l, 1);
            path gf=gfstart--(0, -1)--(0, 1)--gfend;
            draw(gf, dashed, Arrows);

            draw((0,-eps)--(0,eps), primalline+linewidth(1.1), Arrows);
            dot((0, 0));

            dot((0, -1));

            draw((0.5*l-eps, 1)--(0.5*l+eps, 1), primalline+linewidth(1.1), Arrows);
            dot((0.5*l, 1));

            dot((0,-1-eps), invisible);
        \end{asy}
        \caption{Clarke graphical derivative $\clarkeGD[\subdiff f]$}
    \end{subfigure}
    \caption{Illustration of the different graphical derivatives and coderivatives of $\subdiff f$ for $f=\abs{\freevar}$. The dashed line is $\graph \subdiff f$. The dots indicate the base points $(x, y)$ where $D[\subdiff f](x|y)$ is calculated, and the thick arrows and filled-in areas the directions of $(\Delta x, \Delta y)$ (resp.~$(\Delta x, -\Delta y)$ for the coderivatives) relative to the base point.  Observe that there is no graphical regularity at $(x, y)\in\{(0, -1), (0, 1)\}$. Everywhere else, $\subdiff f$ is graphically regular.  Observe that the cones in the last figure of each row are polar to the cones in the first and the second figures on the same row.
    }
    \label{fig:graphical:absvalue}
\end{figure}

\begin{theorem}
    \label{lemma:graphical:absvalue}
    Let  $f(x) \defeq \abs{x}$ for $x\in \R$. Then
    \begingroup
    \allowdisplaybreaks
    \begin{align}
        \label{eq:graphical:absvalue:gderiv}
        D[\subdiff f]( x | y )(\dir x ) & =
        \begin{cases}
            \{0\}
            &\text{if }  x  \ne 0,\,  y  = \sign  x , \\
            \{0\}
            &\text{if }  x =0,\,
            \dir x \ne 0,\, y = \sign \dir x,
            \\
            (-\infty,0] y
            &\text{if }  x =0,\, \dir x  = 0,\, \abs{ y } = 1, \\
            \R
            &\text{if }  x =0,\, \dir x  = 0,\, \abs{ y } < 1, \\
            \emptyset
            & \text{otherwise,}
        \end{cases}
        \\[1ex]
        \frechetCod[\subdiff f]( x | y )(y^* ) & =
        \begin{cases}
            \{0\}
            &\text{if }  x  \ne 0,\,  y  = \sign  x , \\
            (-\infty,0] y
            &\text{if }  x =0,\, y y^* \le 0,\, \abs{ y } = 1, \\
            \R
            &\text{if }  x =0,\, y^*  = 0,\, \abs{ y } < 1, \\
            \emptyset
            & \text{otherwise,}
        \end{cases}
        \\[1ex]
        \label{eq:graphical:absvalue:gclarke}
        \clarkeGD[\subdiff f]( x | y )(\dir x ) & =
        \begin{cases}
            \{0\}
            &\text{if }  x  \ne 0,\,  y  = \sign  x , \\
            \{0\}
            &\text{if }  x =0,\, \dir x  = 0,\, \abs{ y } = 1, \\
            \mathrlap{\R}\phantom{(-\infty,0]x}
            &\text{if }  x =0,\, \dir x  = 0,\, \abs{ y } < 1, \\
            \emptyset
            & \text{otherwise,}
        \end{cases}
        \shortintertext{and}
        \label{eq:graphical:absvalue:colimiting}
        \coderivative [\subdiff f]( x | y )(y^* ) & =
        \begin{cases}
            \{0\}
            &\text{if }  x  \ne 0,\,  y  = \sign  x , \\
            \{0\}
            &\text{if }  x =0,\,
            yy^* > 0,\, \abs{y}=1,
            \\
            (-\infty, 0] y
            &\text{if }  x =0,\, y y^* < 0,\, \abs{ y } = 1, \\
            \R
            &\text{if }  x =0,\, y^*  = 0,\, \abs{ y } \le 1, \\
            \emptyset
            & \text{otherwise.}
        \end{cases}
    \end{align}
    \endgroup
    In particular, $\subdiff f$ is graphically regular if and only if $ x \ne 0$ or $\abs{ y } < 1$.
\end{theorem}
\begin{proof}
    To start with proving \eqref{eq:graphical:absvalue:gderiv}, we recall from \cref{ex:convex:subdiff_abs} that
    \begin{equation}
        \label{eq:graphical:absvalue:subdiff}
        \subdiff f( x )= \sign(x) =
        \begin{cases}
            \{1\} &\text{if }  x  >0, \\
            \{-1\} &\text{if }  x  <0, \\
            [-1, 1] &\text{if }  x  = 0.
        \end{cases}
    \end{equation}
    To calculate the graphical derivative, we use that if $ y  \in \subdiff f( x )$ and $\dir{ y } \in  D[\subdiff f]( x | y )(\dir x )$, there exist by \eqref{eq:graphical:gderiv-altdef} sequences $t_k \downto 0$, $x_k \to x$, and $ y_k \in \subdiff f( x +t_k\dir{ x_k})$ such that
    \begin{equation}
        \label{eq:graphical:absvalue:limit}
        \dir x = \lim_{k \to \infty} \frac{x_k-x}{t_k}
        \qquad\text{and}\qquad
        \dir y = \lim_{k \to \infty} \frac{y_k-y}{t_k}.
    \end{equation}
    We proceed by case distinction:
    \begin{enumerate}
        \item $ x \ne 0$ and $y \ne \sign x$: Then $y \not \in \subdiff f(x)$ and therefore $D[\subdiff f]( x | y )=\emptyset$, which is covered by the last case of \eqref{eq:graphical:absvalue:gderiv}.

        \item $x \ne 0$ and $y = \sign x$: Then for any $x_k \to x$, we have  that $\subdiff f(x_k)=\subdiff f(x)=\{\sign x\}$ for $k$ large enough. Therefore, for any $\dir x \in \R$ we have that $\dir y=0$, which is the first case of \eqref{eq:graphical:absvalue:gderiv}.

        \item $x = 0$ and $\dir x \ne 0$: Then $x_k \ne 0$ and $y_k = \sign x_k = \sign\dir x$. Therefore the limits in \eqref{eq:graphical:absvalue:limit} will only exist if $\abs{y}=1$, which holds from $y =\sign \dir x$. Thus $\dir y=0$, i.e., we obtain the second case of \eqref{eq:graphical:absvalue:gderiv}.

        \item $ x =0$ and $\dir x =0$: Then taking $x_k \equiv x$, we can choose $y_k \in [-1, 1]$ arbitrarily. If $\abs{y}=1$, then $(y-y_k)\sign y \le 0$, so \eqref{eq:graphical:absvalue:limit} shows that $\dir y\sign y \le 0$, which is the third case of \eqref{eq:graphical:absvalue:gderiv}.
            If $\abs{y}<1$, we may obtain any $\dir y \in \R$ by the limit in \eqref{eq:graphical:absvalue:limit}. This is the fourth case of  \eqref{eq:graphical:absvalue:gderiv}.
    \end{enumerate}

    The expression for $\frechetCod[\subdiff f](x|y)$ can be verified using \cref{cor:graphical:adjoint}\,\cref{item:graphical:adjoint:fundamental}.
    It can also be seen graphically from \cref{fig:graphical:absvalue}.

    By the inner and outer limit characterizations of \cref{cor:graphical:limits:findim}, we now obtain the expressions for the Clarke graphical derivative $\clarkeGD[\subdiff f]( x | y )$ and the limiting coderivative $\coderivative[\subdiff f](x|y)$.
    Since $\graph \subdiff f$ is locally contained in an affine subspace outside of the \enquote{corner cases} $(x,y) \in \{(0,1),(0,-1)\}$, only the latter need special inspection. For the Clarke graphical derivative, we need to write $\dir y$ as the limit of $\dir y_k \in D[\subdiff f](x_k,y_k)(\dir x_k)$ for some $\dir x_k \to \dir x$ and \emph{all} $\graph \subdiff f\ni (x_k,y_k) \to (x, y)$. Consider for example $(x,y)=(0,-1)$. Trying both $(x_k, y_k)=(0,-1+1/k)$ and $(x_k, y_k)=(-1/k,-1)$, we see that this is only possible for $(\dir x,\dir y)=(\dir x_k,\dir y_k)=(0,0)$. This yields the third case of \eqref{eq:graphical:absvalue:gclarke}.
    Conversely, for the limiting coderivative, it suffices to find \emph{one} such sequence from the Fréchet coderivative. Choosing for $(x,y)=(0,-1)$ again $(x_k, y_k)=(0,-1+1/k)$ and $(x_k, y_k)=(-1/k,-1)$ as well as the constant sequence $(x_k,y_k)=(0,-1)$ yields the fourth, second, and third case of \eqref{eq:graphical:absvalue:colimiting}, respectively.

    Finally, in finite dimensions the mapping $\subdiff f$ is graphically regular if and only if $D[\subdiff f]( x | y )=\clarkeGD[\subdiff f]( x | y )$ by \cref{cor:graphical:regularity:findim}, which is the case exactly when $x \ne 0$ or $\abs{y} < 1$ as claimed.
\end{proof}

\section{Relation to subdifferentials}
\label{sec:graphical:subdiff}

All of the subdifferentials that we have studied in \cref{part:nonconvex} can be constructed from the corresponding normal cones to the epigraph of a functional $J:X\to\Rbar$ as in the convex case; see \cref{lem:convex:subdiff_epi}.

For the Fréchet and Mordukhovich (or limiting) subdifferentials, by defining the \term[mapping!epigraphical]{epigraphical mapping}
\[
    \epi_J:X\setto \R,\qquad \epi_J(x) \defeq \{ t \in \R \mid t \ge J(x)\},
\]
i.e., $\graph \epi_J=\epi J$, it is straightforward to obtain from the corresponding definitions the relationships
\begin{align}
    \label{eq:graphical:frechet-subdiff}
    \subdiff_F J(x)&=\{x^*\in X^* \mid (x^*, -1) \in \frechetNormal_{\epi J}(x, J(x))\} = \frechetCod[\epi_J](x|J(x))(1);
    \\
    \label{eq:graphical:mordukhovich-subdiff}
    \subdiff_M J(x)&=\{x^*\in X^* \mid (x^*, -1) \in N_{\epi J}(x, J(x))\} =  D^*[\epi_J](x|J(x))(1).
\end{align}
Thus the results of the following \cref{chap:cofrechet,chap:colimiting} can be used to derive the missing calculus rules for the Fréchet and Mordukhovich subdifferentials; see \cref{sec:colimiting:subdiff}.

For the Clarke subdifferential, however, we have to work a bit harder.
First, we define for $A \subset X$ and $x \in X$ the \term[cone!normal!Clarke]{Clarke normal cone}
\begin{equation}
    \label{eq:graphical:clarke-normal}
    N^C_A(x) \defeq \polar{\clarkeTangent_A(x)}.
\end{equation}

We can now extend the definition of the Clarke subdifferential to arbitrary functionals $J:X\to\Rbar$ on Gateaux smooth Banach spaces via the Clarke normal cone to their epigraph.
\begin{lemma}
    \label{lemma:graphical:clarke-subdiff}
    Let $X$ be a reflexive and Gateaux smooth Banach space and let $J: X \to \R$ be locally Lipschitz continuous around $x\in X$. Then
    \begin{equation*}
        \subdiff_C J(x)=\{x^*\in X^* \mid (x^*, -1) \in N^C_{\epi J}(x, J(x))\}.
    \end{equation*}
\end{lemma}

\begin{proof}
    The Clarke tangent cone to $\epi J$ by definition is
    \begin{equation*}
        \clarkeTangent_{\epi J}(x, J(x))
        =
        \left\{
            (\dir x, \dir t) \in X \times \R \,\middle|\,
            \begin{array}{r}
                \text{for all } \tau_k \downto 0,\, x_k \to x,\, J(x_k) \le t_k \to J(x)
                \\
                \text{there exist } \alt x_k \in X \text{ and } \alt t_k \ge J(\alt x_k)
                \\
                \text{with } (\alt x_k-x_k)/\tau_k \to \dir x \text{ and } (\alt t_k-t_k)/\tau_k \to \dir t
            \end{array}
        \right\}.
    \end{equation*}
    If $(\dir x, \dir t) \in \clarkeTangent_{\epi J}(x, J(x))$, then replacing $\alt t_k$ by $\alt t_k+\tau_k (\dir s-\dir t) \ge J(\alt x_k)$ shows that also $(\dir x, \dir s) \in \clarkeTangent_{\epi J}(x, J(x))$ for all $\dir s \ge \dir t$. Thus we may make the minimal choices $\alt t_k=J(\alt x_k)$ and $t_k=J(x_k)$ to see that
    \begin{equation*}
        \clarkeTangent_{\epi J}(x, J(x))
        =
        \left\{
            (\dir x, \dir t) \in X \times \R \,\middle|\,
            \begin{array}{r}
                \text{for all } \tau_k \downto 0,\, x_k \to x
                \text{ there exist } \alt x_k \in X
                \\
                \text{with } (\alt x_k-x_k)/\tau_k \to \dir x
                \\
                \text{ and } \limsup_{k \to \infty} (J(\alt x_k)-J(x_k))/\tau_k \le \dir t
            \end{array}
        \right\}.
    \end{equation*}
    Since $J$ is locally Lipschitz continuous, it suffices to take $\alt x_k=x_k+\tau_k\dir x$ to obtain
    \begin{equation*}
        \clarkeTangent_{\epi J}(x, J(x))
        =
        \{
            (\dir x, \dir t) \in X \times \R
            \mid
            x \in X,\, \dir t \ge J^\circ(x; \dir x)
        \}
        =\epi[J^\circ(x; \freevar)].
    \end{equation*}
    Hence $(x^*, -1) \in  N_{\epi J}^C(x, J(x)) = \polar{\clarkeTangent_{\epi J}(x, J(x))}$ if and only if $\dualprod{x^*}{\dir x}_X \le J^\circ(x; \dir x)$ for all $x \in X$, which by definition is equivalent to $x^*\in \subdiff_C J(x)$.
\end{proof}

We furthermore have the following relationship between the Clarke and limiting normal cones.
\begin{corollary}
    \label{cor:graphical:clarke-normal}
    Let $X$ be a reflexive and Gateaux smooth Banach space and $A\subset X$ be closed near $x \in A$. Then
    \begin{equation*}
        N^C_A(x) = \bipolar{N_A(x)} = \closure \conv^* N_A(x),
    \end{equation*}
    where $\closure \conv^*$ denotes the weak-$*$ closed convex hull.
\end{corollary}

\begin{proof}
    First, $N_A(x) \ne \emptyset$ since $x\in A$. Furthermore, $\closure\conv^* N_A(x)$ is the smallest weak-$*$-closed and convex set that contains $N_A(x)$, and therefore \cref{lemma:functan:polar-inclusion,lem:convex_closed} imply $\bipolar{N_A(x)}=\bipolar{\closure \conv^* N_A(x)}=\closure \conv^* N_A(x)$.
    The relationship $N^C_A(x)=\bipolar{N_A(x)}$ is an immediate consequence of \cref{thm:cones:limiting-polar}.
\end{proof}

Assuming that $X$ is Gateaux smooth, we now have everything at hand to give a proof of \cref{thm:limiting:clarke}, which characterizes the Clarke subdifferential as the weak-$*$ closed convex hull of the limiting subdifferential.

\begin{corollary}
    \label{cor:graphical:clarke-weakstar-convex}
    Let $X$ be a reflexive and Gateaux smooth Banach space and $J:X\to\R$ be locally Lipschitz continuous around $x\in X$. Then $\partial_C J(x) = \mathrm{cl}^* \conv \partial_M J(x)$.
\end{corollary}
\begin{proof}
    Together, \cref{lemma:graphical:clarke-subdiff,cor:graphical:clarke-normal} and \eqref{eq:graphical:mordukhovich-subdiff} directly yield
    \begin{equation*}
        \begin{aligned}[b]
            \subdiff_C J(x)&=\{x^*\in X^* \mid (x^*, -1) \in N^C_{\epi J}(x, J(x))\} \\
            &=\{x^*\in X^* \mid (x^*, -1) \in \mathrm{cl}^*\conv N_{\epi J}(x, J(x))\}\\
            &=\mathrm{cl}^*\conv \{x^*\in X^* \mid (x^*, -1) \in N_{\epi J}(x, J(x))\}\\
            &=\mathrm{cl}^*\conv \partial_M J(x).
        \end{aligned}
        \qedhere
    \end{equation*}
\end{proof}

\begin{corollary}
    \label{cor:graphical:subdiff-inclusions}
    Let $X$ be a reflexive and Gateaux smooth Banach space and $J:X\to\R$ be locally Lipschitz continuous around $x\in X$.
    Then $\subdiff_C J(x) \supset \subdiff_M J(x) \supset \subdiff_F J(x)$.
\end{corollary}

\begin{proof}
    The first inclusion is immediate from \cref{cor:graphical:clarke-weakstar-convex}.
    The second inclusion follows from \cref{eq:graphical:frechet-subdiff,eq:graphical:mordukhovich-subdiff} and the corresponding relation for the normal cones from \cref{thm:cones:inclusions}.
\end{proof}

\begin{remark}
    The Gateaux smoothness of $X$ can be relaxed to $X$ being an Asplund space following \cref{rem:epsilon:asplund}.
\end{remark}

\chapter{Derivatives and coderivatives of pointwise-defined mappings}
\label{chap:superposition}

Just as for tangent and normal cones, the relationships between the basic and limiting derivatives and coderivatives are less complete in infinite-dimensional spaces than in finite-dimensional ones. In this chapter, we apply the results of \cref{chap:pointcones} to derive pointwise characterizations analogous to \cref{thm:lebesgue:subdiff} for the basic derivatives of pointwise-defined set-valued mappings, which (only) in the case of graphically regularity transfer to their limiting variants.

\section{Proto-differentiability}
\label{sec:superposition:protodiff}

For our superposition formulas, we need some regularity from the finite-dimensional mappings.
The appropriate notion is that of \emph{proto-differentiability}, which corresponds to the \term[set!derivable!geometrically]{geometric derivability} of the underlying tangent cone.

Let $X,Y$ be Banach spaces.
We say that a set-valued mapping $F: X \setto Y$ is \term[mapping!differentiable!proto-]{proto-differentiable} at $x \in X$ for $y \in F(x)$ if
\begin{subequations}
    \label{eq:gderiv:protodiff}
    \begin{align}
        &\text{for every } \dir x\in X,\, \dir y \in DF(x|y)(\dir x),\text{ and }\tau_k \downto 0,\\
        &\text{there exist } x_k \in X
        \text{ with }
        \frac{x_k-x}{\tau_k} \to \dir x
        \quad \text{and}\quad y_k \in F(x_k)
        \text{ with }
        \frac{y_k-y}{\tau_k} \to \dir y.
    \end{align}
\end{subequations}
In other words, in addition to the basic limit \eqref{eq:graphical:gderiv-altdef} defining $DF(x|y)$, a corresponding inner limit holds in the graph space.

By application of \cref{lemma:superposition:derivable,cor:superposition:regular-derivable}, we immediately obtain the following equivalent characterization.
\begin{corollary}
    \label{cor:superposition:derivable-proto}
    Let $X,Y$ be Banach spaces and $F: X\setto Y$. Then $F$ is proto-differentiable at every $x \in X$ for every $y \in F(x)$ if and only if $\graph F$ is geometrically derivable at $(x, y)$.
    In particular, if $F$ is graphically regular at $(x,y)$, then $F$ is proto-differentiable at $x$ for $y$.
\end{corollary}

Clearly, differentiable single-valued mappings are proto-differentiable.
Another large class are maximally monotone set-valued mappings on Hilbert spaces.
\begin{lemma}
    Let $X$ be a Hilbert space and let $A: X\setto X$ be maximally monotone.
    Then $A$ is proto-differentiable at any $x\in \dom A$ for any  $x^*\in A(x)$.
\end{lemma}
\begin{proof}
    Let $\dir x^* \in D[A](x|x^*)(\dir x)$.
    By definition, there then exist $\tau_k \downto 0$ and $(x_k, x_k^*) \in \graph A$ such that $(x_k-x)/\tau_k \to \dir x$ and $(x_k^*-x^*)/\tau_k \to \dir x^*$.
    To show that $A$ is proto-differentiable, we will construct for an arbitrary sequence $\alt \tau_k \downto 0$ sequences $(\alt x_k, \alt x_k^*) \in \graph A$ such that $(\alt x_k-x)/\alt \tau_k \to \dir x$ and $(\alt x_k^*-x^*)/\alt \tau_k \to \dir x^*$.
    We will do so using resolvents. Similarly to \cref{lem:proximal:subdiff}, we have that
    \begin{equation*}
        \begin{aligned}
            x^*\in A(x)
            &\quad \Leftrightarrow \quad
            x \in A^{-1}(x^*)
            \quad \Leftrightarrow \quad
            x^* + x \in \{x^*\} + A^{-1}(x^*)\\
            &\quad \Leftrightarrow \quad
            x^* \in \calR_{A^{-1}}(x^*+x).
        \end{aligned}
    \end{equation*}
    Since $A$ is maximally monotone and $X$ is reflexive, $A^{-1}$ is maximally monotone by \cref{lemma:monotone:inverse} as well, and thus the resolvent $\calR_{A^{-1}}$ is single-valued by \cref{lem:proximal:lipschitz}.
    We therefore take
    \begin{equation*}
        \begin{aligned}
            \alt x_k &\defeq x + \frac{\alt \tau_k}{\tau_k}(x_k-x) + \frac{\alt\tau_k}{\tau_k}x_k^*+\left(1-\frac{\alt\tau_k}{\tau_k}\right)x^*-\alt x_k^*
            \quad\text{and}
            \\
            \alt x_k^*
            &\defeq \calR_{A^{-1}}\left(x+
                \frac{\alt \tau_k}{\tau_k}(x_k-x) + \frac{\alt\tau_k}{\tau_k}x_k^*+\left(1-\frac{\alt\tau_k}{\tau_k}\right)x^*
                - \alt\tau_k(\dir x + \dir x^*)
            \right)
            + \tilde\tau_k \dir x^*
            \\
            &
            = \calR_{A^{-1}}(\alt x_k^*+\alt x_k - \alt\tau_k(\dir x + \dir x^*)) + \tilde \tau_k \dir x^*.
        \end{aligned}
    \end{equation*}
    Since resolvents of maximally monotone operators are $1$-Lipschitz by \cref{lem:proximal:firmly-nonexpansive}, we have
    \begin{equation*}
        \begin{aligned}
        \lim_{k \to \infty} \frac{\norm{\alt x_k^*-x^* - \alt\tau_k \dir x^*}_X}{\alt\tau_k}
        &
        =
        \lim_{k \to \infty} \frac{\norm{\calR_{\inv A}(\alt x_k^*+\alt x_k-\alt\tau_k(\dir x + \dir x^*))-\calR_{\inv A}(x^*+x)}_X}{\alt\tau_k}
        \\
        &
        \leq
        \lim_{k \to \infty} \frac{\norm{(\alt x_k^* + \alt x_k-\alt\tau_k(\dir x + \dir x^*))-(x^*+x)}_X}{\alt\tau_k}
        \\
        &
        =
        \lim_{k \to \infty}
        \frac{\norm{(x_k-x-\tau_k \dir x)+(x_k^*-x^*-\tau_k \dir x^*)}_X}{\tau_k}
        =0.
        \end{aligned}
    \end{equation*}
    Likewise, by inserting the definition of $\alt x_k$ and using the triangle inequality, we obtain
    \begin{equation*}
        \begin{aligned}
        \lim_{k \to \infty} \frac{\norm{\alt x_k-x-\alt\tau_k \dir x}_X}{\alt\tau_k}
        &
        \le
        \lim_{k \to \infty} \frac{\norm{(x_k-x- \tau_k \dir x)+(x_k^*-x^* - \tau_k \dir x^*)}_X}{\tau_k}
        \\
        \MoveEqLeft[-1]
        + \lim_{k \to \infty} \frac{\norm{\alt x_k^*-x^*-\alt\tau_k \dir x^*}_X}{\alt\tau_k}
        \\
        &
        =0.
        \end{aligned}
    \end{equation*}
    This shows the claimed proto-differentiability.
\end{proof}

Since subdifferentials of convex and lower semicontinuous functionals on reflexive Banach spaces are maximally monotone by \cref{thm:monoton:subdiff}, we immediately obtain the following.
\begin{corollary}
    \label{lemma:superposition:convex-proto}
    Let $X$ be a Hilbert space and let $J: X \to \Rbar$ be proper, convex, and lower semicontinuous.
    Then $\subdiff J$ is proto-differentiable at any $x\in \dom J$ for any $x^*\in \subdiff J(x)$.
\end{corollary}

This corollary combined with \cref{lemma:graphical:indicator,lemma:graphical:absvalue} shows that proto-differentiability is a strictly weaker property than graphical regularity.

\section{Graphical derivatives and coderivatives}

As a corollary of the tangent and normal cone representations from \cref{thm:superposition:cone-tangent,thm:superposition:cone-frechet}, we obtain explicit characterizations of the  graphical derivative and the Fréchet coderivative of a class of pointwise-defined set-valued mappings. In the following, let $\Omega\subset \R^d$ be an open and bounded domain and write again $p^*$ for the conjugate exponent of $p\in (1,\infty)$ satisfying $1/p+1/p^*=1$.

\begin{theorem}
    \label{cor:superposition:functionals}
    Let $F: L^p(\Omega) \setto L^q(\Omega)$ for $p, q \in (1, \infty)$ have the form
    \begin{equation*}
        F(u) = \{ w \in L^q(\Omega) \mid w(x) \in f(u(x))\, \text{for a.e. }  x \in \Omega\}
    \end{equation*}
    for some pointwise almost everywhere proto-differentiable mapping $f: \R\setto\R$.
    Fix $u\in L^p(\Omega)$ and $w\in L^q(\Omega)$.
    Then for every $w^* \in L^{q^*}(\Omega)$ and $\dir u \in L^p(\Omega)$,
    \begin{subequations}%
        \label{eq:superposition:functionals:derivatives}%
        \begin{align}%
            \label{eq:superposition:functionals:derivatives:frechet}%
            \frechetCod{F}(u|w)(w^*)
            &
            =
            \left\{
                u^* \in L^{p^*}(\Omega)
                \,\middle|\,
                \begin{array}{r}
                    u^*(x) \in \frechetCod{f}(u(x)|w(x))(w^*(x)) \\
                    \text{ for a.e. } x \in \Omega
                \end{array}
            \right\},
            \\
            \label{eq:superposition:functionals:derivatives:gderiv}%
            D{F}(u|w)(\dir{u})
            &
            =
            \left\{
                \dir{w} \in L^q(\Omega)
                \,\middle|\,
                \begin{array}{r}
                    \dir{w}(x) \in D{f}(u(x)|w(x))(\dir{u}(x)) \\
                    \text{ for a.e. } x \in \Omega
                \end{array}
            \right\}.
        \end{align}%
    \end{subequations}%
    Moreover, if $f$ is graphically regular at $u(x)$ for $w(x)$ for almost every $x \in \Omega$, then $F$ is graphically regular at $u$ for $w$ and
    \begin{align*}
        \clarkeGD F(u|w)&=D^w F(u|w)=DF(u|w), \\
        \coderivative F(u|w)&=\frechetCod F(u|w).
    \end{align*}
\end{theorem}

\begin{proof}
    First, $\graph f$ is geometrically derivable by \cref{cor:superposition:derivable-proto} due to the assumed proto-differentiability of $f$. We further have
    \begin{equation*}
        \graph F = \setof{ (u, w) \in L^p(\Omega) \times L^q(\Omega) }{ (u(x), w(x)) \in \graph f \text{ for a.e. } x \in \Omega}.
    \end{equation*}
    Now \eqref{eq:superposition:functionals:derivatives:gderiv} and \eqref{eq:superposition:functionals:derivatives:frechet} follow from \cref{thm:superposition:cone-tangent,thm:superposition:cone-frechet}, respectively, for $C: x \mapsto  \graph f$ and $U=\graph F$ together with definitions of the graphical derivative in terms of the tangent cone and of the Fréchet coderivative in terms of the Fréchet normal cone.
    The remaining claims under graphical regularity follow similarly from \cref{cor:superposition:weak-tangent}.
\end{proof}

The above result directly applies to second derivatives of integral functionals.
\begin{corollary}
    \label{cor:superposition:2nd}
    Let $J:L^p(\Omega)\to\Rbar$ for $p \in (1, \infty)$ be given by
    \begin{equation*}
        J(u)=\int_\Omega j(u(x)) \,d x
    \end{equation*}
    for some proper, convex, and lower semicontinuous integrand $j: \R \to (-\infty, \infty]$.
    Then
    \begin{align*}
        \frechetCod{[\subdiff J]}(u|u^*)(\dir u)
        &=
        \left\{
            \dir u^* \in L^{p^*}(\Omega)
            \,\middle|\,
            \begin{array}{r}
                \dir u^*(x) \in \frechetCod{[\subdiff j]}(u(x)|u^*(x))(\dir u(x)) \\
                \text{ for a.e. } x \in \Omega
            \end{array}
        \right\},
        \\
        D{[\subdiff J]}(u|u^*)(\dir u)
        &=
        \left\{
            \dir u^*  \in L^{p^*}(\Omega)
            \,\middle|\,
            \begin{array}{r}
                \dir u^*(x) \in D{[\subdiff j]}(u(x)|u^*(x))(\dir u(x)) \\
                \text{ for a.e. } x \in \Omega
            \end{array}
        \right\}.
    \end{align*}

    Moreover, if $\subdiff j$ is graphically regular at $u(x)$ for $u^*(x)$ for almost every $x \in \Omega$, then $\subdiff J$ is graphically regular at $u$ for $u^*$ and
    \begin{align*}
        \clarkeGD[\subdiff J](u|u^*)&=D^w[\subdiff J](u|u^*)=D[\subdiff J](u|u^*),
        \\
        \coderivative [\subdiff J](u|u^*)&=\frechetCod [\subdiff J](u|u^*).
    \end{align*}
\end{corollary}

\begin{proof}
    By \cref{lemma:superposition:convex-proto}, $\subdiff j$ is proto-differentiable.
    Since
    \begin{equation*}
        \subdiff J(u) = \left\{ u^* \in L^{p^*}(\Omega) \mid u^*(x) \in \subdiff j(u(x)) \text{ for a.e. } x \in \Omega \right\}
    \end{equation*}
    by \cref{thm:lebesgue:subdiff} and therefore
    \begin{equation*}
        \graph[\subdiff J]=\left\{ (u, u^*) \in L^p(\Omega) \times L^{p^*}(\Omega) \mid u^*(x) \in \subdiff j(u(x)) \text{ for a.e. } x \in \Omega \right\},
    \end{equation*}
    the remaining claims follow from \cref{cor:superposition:functionals} with $F=\subdiff J$, $f=\subdiff j$, and $q=p^*$.
\end{proof}
\begin{remark}
    The case of vector-valued and spatially-varying set-valued mappings and convex integrands can be found in \cite{tuomov-pdex2stability}.
\end{remark}

We illustrate this result with the usual examples.
To keep the presentation simple, we focus on the case $p^*=p=2$ such that $L^2(\Omega)$ is a Hilbert space and we can identify $X\cong X^*$.

First, we immediately obtain from \cref{ex:graphical:norm22} together with \cref{cor:superposition:2nd}
\begin{corollary}
    Let $J:L^2(\Omega)\to\R$ be given by
    \begin{equation*}
        J(u) \defeq \int_\Omega \frac{1}{2}\abs{u(x)}^2 \,dx.
    \end{equation*}
    Then for $u^*=u$ and all $\dir u \in L^2(\Omega)$, we have
    \begin{align*}
        \clarkeGD[\subdiff J](u|u^*)(\dir u)&=D^w[\subdiff J](u|u^*)(\dir u)=D[\subdiff J](u|u^*)(\dir u)
        = \dir u,\\
        \coderivative [\subdiff J](u|u^*)(\dir u)&=\frechetCod [\subdiff J](u|u^*)(\dir u)
        = \dir u.
    \end{align*}
    If $u^* \ne u$, all the derivatives and coderivatives are empty.
\end{corollary}

From \cref{lemma:graphical:indicator}, we also obtain expressions for the basic derivatives of indicator functionals for pointwise constraints. For the limiting derivatives, we only obtain expressions at points where graphical regularity (corresponding to strict complementarity) holds; cf.~\cref{rem:pointcones:limiting}.

\begin{corollary}
    Let $J:L^2(\Omega)\to\Rbar$ be given by
    \begin{equation*}
        J(u) \defeq \int_\Omega \delta_{[-1, 1]}(u(x)) \,dx.
    \end{equation*}
    Let $u\in \dom J$ and $u^* \in \subdiff J(u)$. Then $\dir u^* \in D[\subdiff J](u|u^*)(\dir u) \subset L^2(\Omega)$ if and only if for almost every $x\in \Omega$,
    \begin{equation*}
        \dir u^*(x) \in
        \begin{cases}
            \R
            &\text{if } \abs{ u(x) }=1,\,  u^*(x)  \in (0, \infty)  u(x) ,\, \dir u(x) =0, \\
            [0, \infty)  u(x)
            &\text{if } \abs{ u(x) }=1,\,  u^*(x) =0,\, \dir u(x)  = 0,  \\
            \{0\}
            &\text{if } \abs{ u(x) }=1,\,  u^*(x)  = 0,\,  u(x)  \dir u(x)  < 0, \\
            \{0\}
            &\text{if } \abs{ u(x) } < 1,\,  u^*(x) =0,\, \\
            \emptyset & \text{otherwise}.
        \end{cases}
    \end{equation*}
    Similarly, $\dir u \in \frechetCod[\subdiff J](u|u^*)(\dir u^*) \subset L^2(\Omega)$ if and only if for almost every $x\in \Omega$,
    \begin{equation*}
        \dir u(x) \in
        \begin{cases}
            \R,
            &\text{if } \abs{ u(x) }=1,\,  u^*(x)  \in (0, \infty)  u(x) , \dir u^*(x) = 0,\\
            [0, \infty)  u(x)
            &\text{if } \abs{ u(x) }=1,\,  u^*(x) =0,\, u(x) \dir u^*(x) \ge 0,  \\
            \{0\}
            &\text{if } \abs{ u(x) } < 1,\,  u^*(x) =0, \\
            \emptyset & \text{otherwise}.
        \end{cases}
    \end{equation*}
    If either $\abs{u(x)} <1$ or $u^*(x) \ne 0$, then $\clarkeGD[\subdiff J](u|u^*) = \coderivative[\subdiff J](u|u^*)$.
    % $\dir u^* \in \clarkeGD[\subdiff J](u|u^*) (\dir u) = \coderivative[\subdiff J](u|u^*)(\dir u)$ if and only if for almost every $x\in \Omega$,
    % \begin{equation*}
    %     \dir u^*(x) \in \begin{cases}
    %         \R &\text{if } \abs{u(x)}=1,\, u^*(x) \in (0, \infty)u(x),\, \dir u(x)=0, \\
    %         \{0\} &\text{if } \abs{u(x)}<1,\, \dir u(x) \in \R, \\
    %         \emptyset & \text{otherwise}.
    %     \end{cases}
    % \end{equation*}
\end{corollary}

From \cref{lemma:graphical:absvalue}, we obtain a similar characterization for the basic derivatives of the $L^1$ norm (as a functional on $L^2(\Omega)$).

\begin{corollary}
    Let $J:L^2(\Omega)\to\R$ be given by
    \begin{equation*}
        J(u) \defeq \int_\Omega \abs{u(x)} \,dx.
    \end{equation*}
    Let $u\in \dom J$ and $u^* \in \subdiff J(u)$. Then $\dir u^* \in D[\subdiff J](u|u^*)(\dir u) \subset L^2(\Omega)$ if and only if for almost every $x\in \Omega$,
    \begin{equation*}
        \dir u^*(x) \in
        \begin{cases}
            \{0\} & \text{if } u(x)\neq 0,\, u^*(x) = \sign u(x),\\
            \{0\} & \text{if } u(x) = 0,\, \dir u(x) \neq 0,\, u^*(x) = \sign \dir u(x),\\
            (-\infty,0]u^*(x) & \text{if } u(x) = 0,\, \dir u(x) = 0,\, |u^*(x)|=1,\\
            \R & \text{if } u(x) = 0,\, \dir u(x) = 0,\, |u^*(x)|<1,\\
            \emptyset &\text{otherwise}.
        \end{cases}
    \end{equation*}
    Similarly, $\dir u \in \frechetCod[\subdiff J](u|u^*)(\dir u^*) \subset L^2(\Omega)$ if and only if for almost every $x\in \Omega$,
    \begin{equation*}
        \dir u(x) \in
        \begin{cases}
            \{0\} & \text{if } u(x)\neq 0,\, u^*(x) = \sign u(x),\\
            (-\infty,0]u^*(x) & \text{if } u(x) = 0,\, \dir u(x)\dir u^*(x) \leq 0,\, |u^*(x)|=1,\\
            \R & \text{if } u(x) = 0,\, \dir u^*(x) = 0,\, |u^*(x)|<1,\\
            \emptyset &\text{otherwise}.
        \end{cases}
    \end{equation*}
    If either  $u(x) \ne 0$ or $\abs{u^*(x)} < 1$, then
    $\clarkeGD[\subdiff J](u|u^*)  = \coderivative[\subdiff J](u|u^*)$.
    % $\dir u^* \in \clarkeGD[\subdiff J](u|u^*) (\dir u) = \coderivative[\subdiff J](u|u^*)(\dir u)$ if and only if for almost every $x\in \Omega$,
    % \begin{equation*}
    %     \dir u^*(x) \in
    %     \begin{cases}
    %         \{0\} & \text{if } u(x)\neq 0,\, u^*(x) = \sign u(x),\\
    %         \R & \text{if } u(x) = 0,\, \dir u(x) = 0,\, |u^*(x)|<1,\\
    %         \emptyset &\text{otherwise}.
    %     \end{cases}
    % \end{equation*}
\end{corollary}

Obtaining similar characterizations for derivatives of the Clarke subdifferential of integral functions with nonsmooth nonconvex integrands requires verifying proto-differentiability of the pointwise subdifferential mapping, which is challenging since the Clarke subdifferential in general does not have the nice properties of the convex subdifferential as a set-valued mapping. For problems of the form \eqref{eq:intro:prob} in the preface, it is therefore simpler to first apply the calculus rules from the following chapters (assuming they are applicable) and to then use the above results for the derivatives of the convex or smooth component mappings.

\chapter{Calculus for the graphical derivative}
\label{chap:gderiv}

We now turn to calculus such as sum and product rules. We concentrate on the situation where at least one of the mappings involved is classically differentiable, which allows for exact results and is already useful in practice.
For a much fuller picture of infinite-dimensional calculus in high generality, the reader is referred to \cite{Mordukhovich:2006}. For further finite-dimensional calculus, we refer to \cite{Rockafellar:1998,mordukhovich2018variational}.

The rules we develop for the various (co)derivatives are in each case based on linear transformation formulas of the underlying cones as well as on a fundamental composition lemma.
These fundamental lemmas, however, require further regularity assumptions that are satisfied in particular by (continuously) Fréchet differentiable single-valued mappings and their inverses.
For the sake of presentation, we treat each derivative in its own chapter, starting with the relevant regularity concept, then proving the fundamental lemmas, and finally deriving the calculus rules. We start with the (basic) graphical derivative, which we recall is defined for $F:X\setto Y$ as
\begin{equation}\label{eq:gderiv:defrecall}
    \begin{aligned}[t]
        D F(x|y): X \setto Y, \qquad
        D F(x|y)(\dir x) &\defeq \setof{ \dir y\in Y}{(\dir x, \dir y) \in T_{\graph F}(x, y)}\\
        &= \limsup_{t \downto 0,\, \dir{\alt{x}} \to \dir{x}} \frac{F(x+t \dir{\alt{x}})-y}{t}.
    \end{aligned}
\end{equation}

\section{Semi-differentiability}

Let $X,Y$ be Banach spaces and $F: X \setto Y$. We say that $F$ is \emph{semi-differentiable} at $x \in X$ for $y \in F(x)$ if
\begin{subequations}
    \label{eq:gderiv:semidiff}
    \begin{align}
        &\text{for every } \dir y \in DF(x|y)(\dir x) \quad\text{and}\quad x_k\to x,\ \tau_k \downto 0 \quad\text{with}\quad
        \frac{x_k-x}{\tau_k} \to \dir x\\
        &\text{there exist }
        y_k \in F(x_k)
        \quad\text{with}\quad
        \frac{y_k-y}{\tau_k} \to \dir y,
    \end{align}
\end{subequations}
i.e., if $DF(x|y)$ is a full limit.

\begin{lemma}
    \label{lemma:gderiv:semidiff}
    A mapping $F: X \setto Y$ is semi-differentiable at $x \in X$ for $y \in Y$ if and only if
    \begin{equation}
        \label{eq:gderiv:semidiff-equiv}
        DF(x|y)(\dir x) = \lim_{\tau \downto 0,\, \dir \alt x \to \dir x} \frac{F(x+\tau\dir \alt x)-y}{\tau}
        \qquad(\dir x \in X).
    \end{equation}
\end{lemma}

\begin{proof}
    First, note that the second expression of \eqref{eq:gderiv:defrecall} shows that $DF(x|y)(\dir x)$ is the outer limit corresponding to \eqref{eq:gderiv:semidiff-equiv}. Similarly, by \eqref{eq:gderiv:semidiff}, $F$ is semi-differentiable if $DF(x|y)$ equals to the corresponding inner limit. (For any sequence $\tau_k\downto 0$, we can relate $x_k$ in \eqref{eq:gderiv:semidiff} and $\dir \alt x=:\dir x_k$ in \eqref{eq:gderiv:semidiff-equiv} via $\dir x_k = (x_k-x)/\tau_k$.) Hence, $F$ is semi-differentiable if and only if the outer limit in \eqref{eq:gderiv:defrecall} is a full limit.
\end{proof}

Compared to the definition of proto-differentiability in \cref{sec:superposition:protodiff}, we now require that $\dir y$ can be written as the limit of a difference quotient taken from $F(x_k)$ for \emph{any} sequence $\{x_k\}_{k\in\N}$ similarly realizing $\dir x$ (while for proto-differentiability, this only has to be possible for \emph{one} such sequence).
Hence, semi-differentiability is a stronger property than proto-differentiability with the former implying the latter.

\begin{example}[proto-differentiable but not semi-differentiable]
    Let $F: \R \setto \R$ have $\graph F= \Q \times \{0\}$.
    Then $F$ is proto-differentiable at any $x \in \Q$ by the density of $\Q$ in $\R$.
    However, $F$ is not semi-differentiable, as we can take $x_k \not \in \Q$ in \eqref{eq:gderiv:semidiff}.
\end{example}

To characterize the semi-differentiability of the inverses of single-valued mappings, we require the next lemma.
We say that $A \in \linear(X; Y)$ has a \term[inverse!right-]{right-inverse} $\rinv A \in \linear(Y; Y)$ if $A\rinv A=\Id$; similarly, $\linv A \in \linear(Y;X)$ is a \term[inverse!left-]{left-inverse} for $A$ if $\linv A A=\Id$. For later use, we note that if $\linv A$ is a left-inverse for $A$, then $\linvstar A\in\linear(X^*; Y^*)$ is a right-inverse for $A^*\in\linear (Y^*; X^*)$, i.e., $A^*\linvstar A= \Id$; similarly, a right-inverse for $A$ furnishes a left-inverse for $A^*$.

\begin{lemma}
    \label{lemma:gderiv:inverse:single:left}
    Let $X,Y$ be  Banach spaces and $F: X \to Y$ be continuously differentiable at $x$ such that $F'(x) \in \linear(X; Y)$ has a right-inverse $\rinv{F'(x)} \in \linear(Y; X)$.
    For $P \defeq \Id - \rinv{F'(x)}F'(x)$, define
    \begin{equation*}
        \bar F: X \to Y \times \kernel F'(x),\qquad \bar F(\alt x) \defeq (F(\alt x), P\alt x)
        \quad\text{for all}\ \alt x \in X.
    \end{equation*}
    Then $\bar F$ is bijective in a neighborhood $U$ of $\bar F(x)$ with a continuously differentiable inverse satisfying
    $\inv{\bar F}(\tilde w) \in \inv F(\tilde y)$ for all $\tilde w=(\tilde y, \tilde q) \in U$ as well as
    \begin{equation}
        \label{eq:gderiv:inverse:single:left:diffinv}
        (\inv{\bar F})'(\bar F(x))(\dir y, \dir q) = \rinv{F'(x)}\dir y + \dir q
        \quad\text{for all}\ (\dir y, \dir q) \in Y \times \kernel F'(x).
    \end{equation}
\end{lemma}

\begin{proof}
    Let $A \defeq F'(x)$ and $\rinv A \defeq \rinv{F'(x)}$.
    Then $P = \Id - \rinv A A$ is a projection into $\ker A = \kernel F'(x)$, which implies that $AP=0$.
    We further define
    \begin{equation*}
        M: Y \times \ker A \to X,\qquad  M(\alt y, \alt q) \defeq \rinv A\alt y + \alt q
        \quad\text{for all}\ \alt y \in Y \text{ and } \alt q \in \ker A.
    \end{equation*}
    Then for all $\dir x\in X$,
    \begin{equation*}
        M\bar F'(x)\dir x
        =\rinv A A\dir x + P \dir x=\dir x.
    \end{equation*}
    Thus $M$ is a left-inverse of $\bar F'(x)$, and consequently $\kernel \bar F'(x) = \{0\}$.
    Since $\bar F(x)'\dir x=(A\dir x, P\dir x)$ for all $\dir x\in X$, we similarly have for all $(\alt y, \alt q) \in Y \times \kernel A$ that
    \[
        \bar F'(x)M(\alt y, \alt q)
        = (A\rinv A\tilde y + A\alt q, P\pinv A \tilde y + P \alt q)
        = (A \rinv A \tilde y, P \alt q)
         = (\tilde y, \alt q),
    \]
    which shows that $M$ is also the right-inverse of $\bar F'(x)$ on $Y \times \kernel F'(x)$.
    Hence $\bar F'(x)$ is bijective, $(\inv{\bar F})'(\bar F(x))=M$, and the construction of $M$ establishes \eqref{eq:gderiv:inverse:single:left:diffinv}.

    By the inverse function theorem (\cref{thm:inversefunctiontheorem}), a continuously differentiable $\inv{\bar F}$ exists in a neighborhood $U$ of $w=(y, q) \defeq \bar F(x)$ in $Y \times \kernel A$ with $(\inv{\bar F})'(w)=M$ and $\inv{\bar F}(w)=x$.
    By construction, $\inv{\bar F}(\tilde w) \in \inv F(\tilde y)$ for $\tilde w=(\tilde y, \tilde q) \in U$.
\end{proof}

We now have the following characterizations for the semi-differentiability of single-valued mappings and their inverses.

\begin{lemma}
    \label{lemma:gderiv:regularity:single}
    Let $X,Y$ be Banach spaces and $F: X \to Y$.
    \begin{enumerate}
        \item\label{item:gderiv:regularity:single:forward}
            If $F$ is Fréchet differentiable at $x$, then $F$ is semi-differentiable at $x$ for $y=F(x)$.
        \item\label{item:gderiv:regularity:single:inv}
            If $F$ is continuously differentiable at $x$ and $F'(x) \in \linear(X; Y)$ has a right-inverse $\rinv{F'(x)} \in \linear(Y; X)$, then $\inv F: Y \setto X$ is semi-differentiable at $y=F(x)$ for $x$.
    \end{enumerate}
\end{lemma}

\begin{proof}
    \emph{\cref{item:gderiv:regularity:single:forward}:} This follows directly from the definition of semi-differentiability and the Fréchet derivative.

    \emph{\cref{item:gderiv:regularity:single:inv}:} By \cref{thm:graphical:single:inverse},
    $
    D\inv F(y|x)(\dir y)=\{\dir x\in X \mid F'(x)\dir x=\dir y\}
    $
    for $y=F(x)$. Hence \eqref{eq:gderiv:semidiff} for $\inv F$ requires showing that for all $\tau_k \downto 0$ and $y_k\in Y$ with $(y_k-y)/\tau_k \to F'(x)\dir x$, there exist $x_k$ such that $y_k=F(x_k)$ such that $(x_k-x)/\tau_k \to \dir x$.
    Let $\bar F$ be given by \cref{lemma:gderiv:inverse:single:left}.
    Since $\bar F$ is invertible in a neighborhood of $\bar F(x) = (y, q) =: w$ for $q \defeq Px \in \kernel F'(x)$, let us take $x_k \defeq \inv{\bar F}(y_k, q + \tau_k \dir q)$ for $\dir q \defeq P \dir x$.
    Then, by construction, $\bar F(x_k) = (F(x_k), Px_k)=(y_k, q + \tau_k \dir q)$.
    Moreover
    \[
        \lim_{k \to \infty} \frac{x_k-x}{\tau_k}
        = \lim_{k \to \infty} \frac{\inv{\bar F}(y_k, q + \tau_k \dir q)-\inv{\bar F}(w)}{\tau_k}
        = (\inv{\bar F})'(w)(\dir y, \dir q).
    \]
    By \eqref{eq:gderiv:inverse:single:left:diffinv}, we then have
    \[
        (\inv{\bar F})'(w)(\dir y, \dir q)
        =\rinv{F'(x)}\dir y + \dir q
        =\rinv{F'(x)}F'(x)\dir x + (\Id - \rinv{F'(x)}F'(x))\dir x
        = \dir x,
    \]
    which establishes the claim.
\end{proof}

\begin{remark}
    \label{rem:graphical:regularity:single:d:findim}
    In \cref{lemma:gderiv:regularity:single}\,\ref{item:gderiv:regularity:single:inv}, if $X$ is finite-dimensional, it suffices to assume that $F$ is continuously differentiable with $\kernel F'(x)^*=\{0\}$. In this case we can take $\rinv{F'(x)} \defeq A^*\inv{(AA^*)}$ for $A \defeq F'(x)$.
\end{remark}

\section{Cone transformation formulas}
\label{sec:gderiv:cones}

At their heart, calculus rules for (co)derivatives of set-valued mappings derive from corresponding transformation formulas for the underlying cones.
To formulate these, let $C\subset Y$ and $R \in \linear (Y; X)$.
Define $RC\defeq \setof{Ry}{y\in C}$.
We then say that a point $y \in \closure C$ admits an \term[selection!inverse]{inverse selection} $R_y$ of $R$ at $x=Ry$ if there exists a neighborhood $U_y \subset RC \union \{x\}$ of $x$ and a (not necessarily linear!) mapping $\inv R_y: U_y \to C$ such that $\inv R_y(x)=y$ and $R \inv R_y (\alt x)=\alt x$ for every $\alt x\in U_y$.
We say that the inverse selection is continuous, Lipschitz, or Fréchet differentiable at or near $x$ if $\inv R_y$ satisfies the corresponding property.

\begin{example}[inverse selections]
    Let $C \defeq [-1,1]^2 \subset \R^2$, and $R(y_1, y_2) \defeq y_1$. Then $RC=[-1, 1]$.
    Given $x \in RC$, for any $y=(x, y_2) \in C$, we can take $\inv R_y \alt x \defeq (\alt x, y_2)$ in $U_y = RC = [-1, 1]$. Clearly $U_y$ is a neighborhood of $x$ in $RC$, and we have both $\inv R_y \alt x \in C$ and $R\inv R_y(\alt x)=\alt x$ for every $\alt x \in U_y$.
    Hence the $\inv R_y$ are continuous and even Fréchet differentiable inverse selections.

    If, on the other hand, $C \defeq \B \subset \R^2$ is the Euclidean ball, we still have $RC=[-1, 1]$.
    However, given $x \in RC$ and $y=(x, y_2) \in C$, taking $\inv R_y$ as above, we only have $\inv R_y \alt x \in C$ for $\alt x$ in the set $U_y = \{ \alt x \mid x^2 + y_2^2 \le 1$\}.
    If $y_2=\pm 1$, in which case $y \in C$ forces $\alt x=0$, we get $U_y=\{x\}$, which is not a proper neighborhood of $x$. Hence the $\inv R_y$ are not continuous inverse selections.

    However, $\inv R_y \alt x \defeq (\alt x, y_2\sqrt{1-\alt x^2})$ do satisfy $\inv R_y \alt x \in \B$ as well as $R\inv R_y\alt x = \alt x$ for all $\alt x \in U_y \defeq [-1,1]$ and hence are continuous inverse selections but are not Fréchet differentiable at $x=\pm 1$.
\end{example}

\begin{lemma}
    \label{lemma:gderiv:cone-linear}
    Let $X,Y$ be Banach spaces, $C\subset Y$, and $R\in \linear(Y; X)$.
    If $y \in \closure C$ admits a Fréchet differentiable inverse selection of $R$ at $x=Ry$, then
    \begin{equation*}
        T_{RC}(x) = R T_C(y).
    \end{equation*}%
\end{lemma}

\begin{proof}
    We first prove ``$\supset$''.
    Suppose $\dir y \in T_C(y)$.
    Then $\dir y=\lim_{k \to \infty} (y_k-y)/\tau_k$ for some $y_k \in C$ and $\tau_k \downto 0$. Consequently, since $R$ is bounded, $R(y_k-y)/\tau_k \to R\dir y$. But $Ry \in \closure RC$, so $R\dir y \in T_{RC}(x)$.
    On the other hand, if $y \not \in \closure C$, then $T_C(y)=\emptyset$ and thus there is nothing to show.

    To establish ``$\subset$'', we first of all note that $T_{RC}(x)=\emptyset$ if $x \not \in \closure RC$.
    So suppose $\dir x \in T_{RC}(x)$.
    Then $\dir x=\lim_{k \to \infty} (x_k-x)/\tau_k$ for some $x_k \in RC$, $x_k \to x$, and $\tau_k \downto 0$. For large enough $k$ that $x_k \in U_y$, we have $x_k=Ry_k$ for $y_k \defeq \inv R_y(x_k)$.
    Since $\inv R_y$ is Fréchet differentiable at $x$, letting $h_k \defeq x_k-x$ and using that $(h_k-\tau_k \dir x)/\tau_k=h_k/\tau_k-\dir x \to 0$ and $\norm{h_k}_X/\tau_k \to \norm{\dir x}_X$, we have
    \begin{equation*}
        \begin{aligned}
            \lim_{k \to \infty}\left(
                \frac{y_k-y}{\tau_k}-(\inv R_y)'(x)\dir x
            \right)
            &
            =\lim_{k \to \infty} \frac{\inv R_y(x_k)-\inv R_y(x)-\tau_k (\inv R_y)'(x)\dir x}{\tau_k}
            \\
            &
            =\lim_{k \to \infty} \frac{\inv R_y(x+h_k)-\inv R_y(x)-(\inv R_y)'(x)h_k}{\tau_k}
            =0.
        \end{aligned}
    \end{equation*}
    This proves that $\dir y \defeq (\inv R_y)'(x)\dir x \in T_C(y)$.
    We also have $\dir x=R\dir y$ because
    \[
        \frac{x_k-x}{\tau_k} = \frac{R(y_k-y)}{\tau_k} \to R\dir y,
    \]
    which establishes the inclusion and hence the claim.
\end{proof}

\begin{remark}[qualification conditions in finite dimensions]
    \label{rem:gderiv:qc}
    If $X$ and $Y$ are finite-dimensional, we could replace the existence of the family $\{\inv R_y\}$ of continuous selections in \cref{lemma:gderiv:cone-linear} by the more conventional \term[condition!qualification]{qualification condition}
    \begin{equation*}
        T_C(y)\cap \kernel R = \{0\}.
    \end{equation*}
    We do not employ such a condition, as the extension to Banach spaces would have to be based not on $T_C(y)$ but on the weak tangent cone $T^w_C(y)$ that is difficult to compute explicitly.
\end{remark}

We base all our calculus rules on the previous linear transformation lemma (\cref{lemma:gderiv:cone-linear}) and the following composition lemma for the tangent cone $T_C$.

\begin{lemma}[fundamental lemma on compositions]
    \label{lemma:gderiv:fundamental}
    Let $X,Y,Z$ be Banach spaces and
    \begin{equation*}
        C \defeq \{(x, y, z) \mid y \in F(x),\, z \in G(y)\}
    \end{equation*}
    for $F: X \setto Y$ and $G: Y \setto Z$.
    If $(x, y, z) \in C$ and either
    \begin{enumerate}
        \item\label{item:gderiv:fundamental:G} $G$ is semi-differentiable at $y$ for $z$, or
        \item\label{item:gderiv:fundamental:invF} $\inv F$ is semi-differentiable at $y$ for $x$,
    \end{enumerate}
    then
    \begin{equation}%
        \label{eq:gderiv:fundamental}
        T_C(x, y, z)=\{(\dir x, \dir y, \dir z) \mid \dir y \in D F(x|y)(\dir x),\, \dir z \in D G(y|z)(\dir y)\}.
    \end{equation}
\end{lemma}

\begin{proof}
    We only consider the case \ref{item:gderiv:fundamental:G}; the case \ref{item:gderiv:fundamental:invF} is shown analogously.
    By definition, we have $(\dir x, \dir y, \dir z) \in T_C(x, y, z)$ if and only if for some $(x_k, y_k, z_k) \in C$ and $\tau_k \downto 0$,
    \begin{equation*}
        \dir x=\lim_{k \to \infty} \frac{x_k-x}{\tau_k},
        \qquad
        \dir y=\lim_{k \to \infty} \frac{y_k-y}{\tau_k},
        \qquad
        \dir z=\lim_{k \to \infty} \frac{z_k-z}{\tau_k}.
    \end{equation*}
    On the other hand, we have $\dir y \in DF(x|y)(\dir x)$ if and only if the first two limits hold for some $(x_k, y_k) \in \graph F$ and $\tau_k \downto 0$.
    Likewise, we have $\dir z \in DG(y|z)(\dir y)$ if and only if the last two limits hold for some $(y_k, z_k) \in \graph G$.
    This immediately yields \enquote{$\subset$}.

    To prove \enquote{$\supset$}, take $\tau_k>0$ and $(x_k, y_k) \in \graph F$ such that the first two limits hold.
    By the semi-differentiability of $G$ at $y$ for $z$, for any $\dir z \in DG(y|z)(\dir y)$ we can find $z_k \in G(y_k)$ such that $(z_k-z)/\tau_k \to \dir z$. This shows the remaining limit.
\end{proof}

If one of the two mappings is single-valued, we can use \cref{lemma:gderiv:regularity:single} for verifying its semi-differentiability and \cref{thm:graphical:single} for the expression of its graphical derivative to obtain from \cref{lemma:gderiv:fundamental} the following two special cases.

\begin{corollary}[fundamental lemma on compositions: single-valued outer mapping]
    \label{lemma:gderiv:fundamental:single-outer}
    Let $X,Y,Z$ be Banach spaces and
    \begin{equation*}
        C \defeq \{(x, y, G(y)) \mid y \in F(x)\}
    \end{equation*}
    for $F: X \setto Y$ and $G: Y \to Z$.
    If $(x, y, z) \in C$ and $G$ is Fréchet differentiable at $y$, then
    \begin{equation*}
        T_C(x, y, z)=\{(\dir x, \dir y, G'(y)\dir y) \mid \dir y \in DF(x|y)(\dir x)\}.
    \end{equation*}
\end{corollary}

\begin{corollary}[fundamental lemma on compositions: single-valued inner mapping]
    \label{lemma:gderiv:fundamental:single-inner}
    Let $X,Y,Z$ be Banach spaces and
    \begin{equation*}
        C \defeq \{(x, y, z) \mid y=F(x),\, z \in G(y)\}
    \end{equation*}
    for $F: X \setto Y$ and $G: Y \to Z$.
    If $(x, y, z) \in C$, $F$ is continuously differentiable at $x$
    and $F'(x)$ has a right-inverse $\rinv{F'(x)} \in \linear(Y; X)$,
    then
    \begin{equation*}%
        T_C(x, y, z)=\{(\dir x, \dir y, \dir z) \mid \dir y=F'(x)\dir x,\, \dir z \in DG(y|z)(\dir y)\}.
    \end{equation*}%
\end{corollary}

\section{Calculus rules}%
\label{sec:gderiv:calculus}

Combining now the previous results, we quickly obtain various calculus rules. We begin as usual with a sum rule.

\begin{theorem}[addition of a single-valued differentiable mapping]
    \label{thm:gderiv:addition}
    Let $X,Y$ be Banach spaces, $G:X\to Y$ be Fréchet differentiable at $x \in X$, and $F:X\setto Y$.
    Then for any $y \in H(x)\defeq F(x) + G(x)$, we have
    \begin{equation*}
        D H(x|y)(\dir x) = DF(x|y-G(x))(\dir x) + \{ G'(x) \dir x \}
        \qquad (\dir x \in X).
    \end{equation*}
\end{theorem}

\begin{proof}
    We can write $\graph H=RC$ for
    \begin{equation}
        \label{eq:gderiv:addition:c-r}
        C \defeq \{(u, \alt x, G(\alt x)) \mid \alt x \in X,\, u \in F(\alt x)\}
        \quad\text{and}\quad
        R(u, \alt x, v) \defeq (\alt x, u+v).
    \end{equation}
    We have $(x, y)=Rp$ for $p \defeq (y-G(x), x, G(x)) \in C$.
    To use \cref{lemma:gderiv:cone-linear} to calculate $T_{RC}(x, y)$, define the inverse selection
    \begin{equation}
        \label{eq:gderiv:addition:ry}
        \inv R_p: RC \to C ,\qquad
        \inv R_p(\alt x, \alt y) \defeq (\alt y-G(\alt x), \alt x, G(\alt x)).
    \end{equation}
    Then $\inv R_p(x, y)=p$ and $\inv R_p(\alt x, \alt y) \in C$ for every $(\alt x, \alt y) \in RC$. Furthermore, by the assumed Fréchet differentiability of $G$ at $x$, $\inv R_p$ is a Fréchet differentiable inverse selection at $(x,y)$.

    Taking any neighborhood $U_p \subset RC$ of $(x, y)$, \cref{lemma:gderiv:cone-linear} now yields
    \begin{equation*}
        T_{\graph H}(x, y)
        = RT_C(p)
        = \{(\dir x, \dir u+\dir v) \mid (\dir u, \dir x, \dir v) \in T_C(p)\}.
    \end{equation*}
    Moreover, $C$ given in \eqref{eq:gderiv:addition:c-r} coincides with the $C$ defined in \cref{lemma:gderiv:fundamental:single-outer} with $\inv F$ in place of $F$.
    Inserting the expression from \cref{lemma:graphical:inverse} for $D\inv F$ into the result, it follows that
    \begin{equation*}
        T_C(p)=\{(\dir u, \dir x, G'(x)\dir x) \mid \dir u \in DF(x|y-G(x))(\dir x)\}.
    \end{equation*}
    Thus
    \begin{equation*}
        \begin{aligned}
            DH(x|y)(\dir x)
            &= \{\dir u+\dir v \mid (\dir u, \dir x, \dir v) \in T_C(p)\}
            \\
            &= \{\dir u+G'(x)\dir x \mid \dir u \in DF(x|y-G(x))(\dir x)\},
        \end{aligned}
    \end{equation*}
    which yields the claim.
\end{proof}

We now turn to chain rules, beginning with the case that the outer mapping is single-valued. To apply the cone transformation formula, we require a kind of \enquote{one-sided local invertibility} of the outer mapping. To this end, we call a mapping $G:Y\to Z$ \term[mapping!locally left-invertible]{locally left-invertible} near $z\in \range G$ if there exists a neighborhood $U_z\subset \range G$ of $z$ such that for all $\alt z\in U_z$, there exists a \emph{unique} $\alt y\in Y$ such that $\alt z = G(\alt y)$. In a slight abuse of notation, we call $\linv G:U_z\to Y$, $\alt z\mapsto \alt y$ a \term[inverse!left-!local]{local left-inverse} of $G$ near $z$.

\begin{theorem}[outer composition with a single-valued differentiable mapping]
    \label{thm:gderiv:outer}
    Let $X,Y$ be Banach spaces, $F:X\setto Y$, and $G:Y\to Z$.
    Let $x \in X$ and $z \in H(x) \defeq G(F(x))$ be given.
    If $G$ is Fréchet differentiable at $y \in F(x) \isect \inv G(\{z\})$ and locally left-invertible near $z$ such that the local left-inverse $\linv G$ is Fréchet differentiable at $z$, then
    \begin{equation*}
        D H(x|z)(\dir x) = G'(y)DF(x|y)(\dir x)
        \qquad (\dir x \in X).
    \end{equation*}
\end{theorem}

\begin{proof}
    Observing that $\graph H = RC$ for
    \begin{equation}\label{eq:gderiv:outer:RC}
        C \defeq \{(\alt x, \alt y, G(\alt y)) \mid \alt y \in F(\alt x)\}
        \quad\text{and}\quad
        R(\alt x, \alt y, \alt z) \defeq (\alt x, \alt z),
    \end{equation}
    we again use \cref{lemma:gderiv:cone-linear} to calculate $T_{RC}(x, z)$.
    Clearly, $(x, z)=Rp$ for $p \defeq (x,y,z) \in C$.
    Accordingly, we define the inverse selection
    \begin{equation}
        \label{eq:gderiv:outer:invR}
        \inv R_{p}: RC \to C,\qquad \inv R_{p}(\alt x, \alt z)\defeq(\alt x, \linv G(\alt z), \alt z).
    \end{equation}
    Then $\inv R_{p}(x, z)=p$.
    By construction and the assumption of left-invertibility of $G$, for any $(\alt x, \alt z) \in RC$ near $(x, z)$ there exists a unique $\alt y = \linv G(\alt z)$ such that $G(\alt y)=\alt z$ and -- since $(\alt x, \alt z) \in RC$ and $\alt y$ is unique -- $\alt y \in F(\alt x)$.
    Hence $\inv R_{p}(\alt x, \alt z) \in C$ and $\inv R_p$ is an inverse selection of $R$ at $(x, z)$.
    Moreover, since $\linv G$ is Fréchet differentiable at $z$, $\inv R_{p}$ is a Fréchet differentiable inverse selection at $(x, z)$.

    Taking any neighborhood $U_p \subset RC$ of $(x, z)$, we can therefore apply \cref{lemma:gderiv:cone-linear}, which yields
    \begin{equation*}
        T_{\graph H}(x, z)
        = RT_C(p)
        = \{(\dir x, \dir z) \mid (\dir x, \dir y, \dir z) \in T_C(p)\}.
    \end{equation*}
    Since $G$ is Fréchet differentiable at $y$, we can further use \cref{lemma:gderiv:fundamental:single-outer} to obtain
    \begin{equation*}
        \begin{aligned}
            D H(x|z)(\dir x) &= \{ \dir z \mid (\dir x, \dir y, \dir z) \in T_C(p)\}
            \\
            &
            =  \{ G'(y)\dir y \mid \dir y \in DF(x|y)(\dir x)\},
        \end{aligned}
    \end{equation*}
    which after further simplification yields the claimed expression.
\end{proof}
In particular, this result holds if $G$ is injective around $y$ as well as continuously differentiable such that $G'(y)$ is bijective, since in this case the inverse function theorem (\cref{thm:inversefunctiontheorem}) guarantees the local existence and differentiability of $\linv G$.

A useful special case is when the mapping $G$ is linear.

\begin{corollary}[outer composition with a linear operator]
    \label{cor:gderiv:outer:linear}
    Let $X,Y,Z$ be Banach spaces, $A\in \linear(Y; Z)$, and $F:X\setto Y$.
    If $A$ has a left-inverse $\linv A$,
    then for any $x\in X$ and $z\in H(x)\defeq AF(x)$,
    \begin{equation*}
        D H(x|z)(\dir x) = A DF(x|y)(\dir x)
        \qquad (\dir x \in X)
    \end{equation*}
    for the unique $y \in Y$ such that $Ay=z$.
\end{corollary}

\begin{proof}
    We apply \cref{thm:gderiv:outer} to $G(y) \defeq Ay$, which is clearly Fréchet differentiable at every $y \in F(x)$ and by assumption has the -- similarly differentiable -- local left-inverse $\linv G=\linv A$ near every $z\in Z$.
\end{proof}

The assumption of left-invertibility is in particular satisfied if $Y$ and $Z$ are Hilbert spaces and $A$ is injective and has closed range, since in this case we can take $\linv A=\pinv A \defeq \inv{(A^*A)}A^*$
(the \term[pseudoinverse, Moore--Penrose]{Moore--Penrose pseudoinverse} of $A$).
%and $\pinvstar A=(\pinv A)^*$.

We next consider chain rules where the inner mapping is single-valued.

\begin{theorem}[inner composition with a single-valued differentiable mapping]
    \label{thm:gderiv:inner}
    Let $X,Y,Z$ be Banach spaces, $F: X\to Y$ and $G:Y\setto Z$.
    Let $x \in X$ and $z \in H(x)\defeq G(F(x))$ be given.
    If $F$ is continuously differentiable at $x$ and $F'(x)$ has a right-inverse $\rinv{F'(x)} \in \linear(Y; X)$, then
    \begin{equation*}
        D H(x|z)(\dir x)  = D G(F(x)|z)(F'(x)\dir x)
        \qquad (\dir x \in X).
    \end{equation*}%
\end{theorem}

\begin{proof}
    Observing that $\graph H=RC$ for
    \begin{equation}\label{eq:gderiv:inner:RC}
        C \defeq \{(\alt x, \alt y, \alt z) \mid \alt y=F(\alt x),\, \alt z \in G(\alt y)\}
        \quad\text{and}\quad R(\alt x, \alt y, \alt z) \defeq (\alt x, \alt z),
    \end{equation}
    we again use \cref{lemma:gderiv:cone-linear} to compute $T_{RC}(x, z)$.
    We have $(x, z)=Rp$ for $p \defeq (x, F(x), z) \in C$.
    Accordingly, we define the inverse selection
    \begin{equation}
        \label{eq:gderiv:inner:invRp}
        \inv R_{p}: RC \to C,\qquad\inv R_{p}(\alt x, \alt z) \defeq (\alt x, F(\alt x), \alt z).
    \end{equation}
    Clearly $\inv R_{p}(x, z)=(x, F(x), z)$ and $\inv R_{p}(\alt x, \alt z) \in C$ for $(\alt x, \alt z) \in RC$. Moreover, by the Fréchet differentiability of $F$ at $x$, also $\inv R_{p}$ is Fréchet differentiable at $(x,z)$.

    Thus \cref{lemma:gderiv:cone-linear} yields
    \begin{equation*}
        T_{\graph H}(x, z)=\{(\dir x, \dir z) \mid (\dir x, \dir y, \dir z) \in T_C(p)\}.
    \end{equation*}
    On the other hand, due to the {continuous} differentiability of $F$ at $x$ and the right-invertibility of $F'(x)$, we can apply \cref{lemma:gderiv:fundamental:single-inner} to obtain
    \begin{equation*}
        T_C(x, y, z)=\{(\dir x, \dir y, \dir z) \mid \dir y=F'(x)\dir x,\, \dir z \in DG(y|z)(\dir y)\}.
    \end{equation*}
    Thus
    \begin{equation*}
        \begin{aligned}
            DH(x|z)(\dir x)
            &=
            \{ \dir z \mid (\dir x, \dir y, \dir z) \in T_C(p)\}
            \\
            &=
            \{ \dir z \mid \dir y=F'(x)\dir x,\, \dir z \in DG(F(x)|z)(\dir y)\},
        \end{aligned}
    \end{equation*}
    which yields the claim.
\end{proof}

Again, we can specialize this result to the case where the single-valued mapping is linear.
\begin{corollary}[inner composition with a linear operator]%
    \label{cor:gderiv:inner:linear}
    Let $X,Y,Z$ be Banach spaces, $A\in \linear(X; Y)$, and $G:Y\setto Z$.
    If $A$ has a right-inverse $\rinv A \in \linear(Y; X)$, then for all $x \in X$ and $z \in H(x) \defeq G(Ax)$,
    \begin{equation*}
        DH(x|z)(\dir x)  = DG(Ax|z)(A\dir x)
        \qquad (\dir x \in X).
    \end{equation*}%
\end{corollary}

We wish to apply these results to derive second-order chain rules from \cref{thm:convex:chain,thm:clarke:chain}. For the former, this is straight-forward based on the two corollaries so far obtained.

\begin{corollary}[second-order chain rule for convex subdifferentials]
    \label{cor:gderiv:second-convex}
    Let $X,Y$ be Banach spaces, $f:Y\to \Rbar$ be proper, convex, and lower semicontinuous, $A\in\linear(X; Y)$ be such that $A$ has a right-inverse $\rinv{A}  \in \linear(Y; X)$ and that $\range A \isect \interior \dom f \ne \emptyset$. Let $h\defeq f\circ A$.
    Then for any $x\in X$ and $x^*\in \partial h(x) = A^*\partial f(Ax)$,
    \begin{equation*}
        D[\subdiff h](x|x^*)(\dir x)  = A^* D[\subdiff f](Ax|y^*)(A\dir x)
        \qquad (\dir x \in X)
    \end{equation*}
    for the unique $y^* \in Y^*$ satisfying $A^*y^*=x^*$.
\end{corollary}

\begin{proof}
    The expression for $\partial h(x)$ follows from \cref{thm:convex:chain}, to which we apply \cref{cor:gderiv:inner:linear} as well as \cref{cor:gderiv:outer:linear} with $A^*$ in place of $A$, recalling that a right-inverse $\rinv A$ of $A$ produces a left-inverse $\linv{(A^*)} = (\rinv A)^*$ of $A^*$.
\end{proof}

To obtain a second-order chain rule from \cref{thm:clarke:chain}, we also need a product rule for a single-valued mapping $G$ and a set-valued mapping $F$.
In principle, this could be obtained as a composition of $x \mapsto (x_1, x_2)$, $(x_1, x_2) \mapsto \{G(x_1)\} \times F(x_2)$, and $(y_1, y_2) \mapsto y_1y_2$; however, the last one of these mappings does not possess the left-inverse required by \cref{thm:gderiv:outer}. We therefore take another route, which starts with the following lemma.

\begin{lemma}
    \label{lemma:gderiv:cartesian-product}
    Let $X,Y$ be Banach spaces and $F:X\setto Y$.
    Define $\bar F:X\setto X\times Y$ by $\bar F(x) \defeq \{x\} \times F(x)$.
    Then, for all $x, \dir x \in X$ and $y \in F(x)$, we have
    \begin{equation*}
        D \bar F(x|(x, y))(\dir x)
        = \{\dir x\} \times DF(x|y)(\dir x).
    \end{equation*}
\end{lemma}

\begin{proof}
    We have
    \begin{equation*}
        \graph \bar F = R \graph F
        \quad\text{for}\quad
        R(\alt x, \alt y) \defeq (\alt x, (\alt x, \alt y)).
    \end{equation*}
    Let now $y \in F(x)$.
    Observe that $(x, y) \in \graph F$ is the unique point satisfying $R(x,y)=p$ for $p=(x, (x, y)) \in \graph \bar F$.
    We define
    \[
        \inv R_p: R \graph F \to \graph F,
        \quad
        \inv R_p(\alt x, (\alt z, \alt y)) \defeq (\alt x,  \alt y).
    \]
    This is clearly a Fréchet differentiable inverse selection of $R$ at $p$ for any neighborhood $U_p \subset \graph \bar F$ of $(x, y)$.
    By \cref{lemma:gderiv:cone-linear}, we therefore have
    \begin{equation*}
        T_{R\graph F}(p)
        = \{(\dir x, (\dir x, \dir y)) \mid (\dir x, \dir y) \in T_{\graph F}(x, y)\},
    \end{equation*}
    which establishes the claim.
\end{proof}

\begin{theorem}[product rule]
    \label{thm:gderiv:product}
    Let $X,Y,Z$ be Banach spaces, $F:X\setto Y$, and let $G:X\to \linear (Y; Z)$ be Fréchet differentiable at $x\in X$.
    If $G(\alt x) \in \linear(Y; Z)$ has a left-inverse $\linv{G(\alt x)} \in \linear(Z; Y)$ for every $\alt x$ near $x$ and the mapping $\alt x \mapsto \linv{G(\alt x)}$ is Fréchet differentiable at $x$, then for all $z\in H(x)\defeq G(x)F(x) \defeq \Union_{y \in F(x)} G(x)y$,
    \begin{equation*}
        D H(x|z)(\dir x)
        =
        \{[G'(x)\dir x]y+G(x)DF(x|y)\dir x\}
        \quad (z \in H(x),\, \dir x \in X)
    \end{equation*}
    for the unique $y \in F(x)$ satisfying $G(x)y=z$.
\end{theorem}

\begin{proof}
    Let $\bar F$ be as in \cref{lemma:gderiv:cartesian-product}.
    Then $\graph H=R \graph(\bar G \circ \bar F)$ for
    \begin{equation*}
        \bar G(\alt x, \alt y)=(\alt x, G(\alt x)\alt y)
        \quad\text{and}\quad
        R(\alt x_1, \alt x_2, \alt z) \defeq (\alt x_1, \alt z),
    \end{equation*}
    where
    \[
        \begin{aligned}[t]
        \graph(\bar G \circ \bar F)
        &
        = \{(\alt x, \alt x, G(x)\alt y) \mid (\alt x, (\alt x, \alt y)) \in \graph \bar F \}
        \\
        &
        = \{(\alt x, \alt x, G(\alt x)\alt y) \mid \alt x \in X,\, \alt y \in F(\alt x) \}.
        \end{aligned}
    \]

    We now wish to apply \cref{thm:gderiv:outer} to $\bar G \circ \bar F$, for which we need to verify its assumptions.
    First, $\bar G$ is single-valued and differentiable at $(x, y)$.
    Since $G(\alt x)$ is assumed to be left-invertible for $\alt x$ near $x$, the mapping $\linv{\bar G}: (\alt x, \alt z) \mapsto (\alt x, \linv{G(\alt x)}\alt z)$ is a left-inverse of $\bar G$, which is Fréchet differentiable at $(x,z)$ since $\alt x \to \linv{G(\alt x)}$ is Fréchet differentiable at $x$.
    Finally, we also have
    \begin{equation*}
        \bar G'(x, y)(\dir x, \dir y)=(\dir x, [G'(x)\dir x]y+G(x)\dir y).
    \end{equation*}
    Thus \cref{thm:gderiv:outer,lemma:gderiv:cartesian-product} yield
    \begin{equation*}
        \begin{aligned}
            D[\bar G \circ \bar F](x|(x, z))(\dir x)
            &= \bar G'(x, y)D \bar F(x|(x, y))(\dir x)
            \\
            &=  \bar G'(x, y)(\dir x, DF(x|y)\dir x)
            \\
            &= \{\dir x\} \times \left(\{[G'(x)\dir x]y\}+G(x)DF(x|y)\dir x\right).
        \end{aligned}
    \end{equation*}
    It follows that
    \begin{equation*}
        T_{\graph(\bar G \circ \bar F)}(x, x, z)
        = \{(\dir x, \dir x, \dir z) \mid \dir z \in \{[G'(x)\dir x]y\}+G(x)DF(x|y)\dir x\}.
    \end{equation*}

    By the left-invertibility of $G(x)$, any $y$ satisfying $G(x)y=z$ is unique.
    Let $p=(x,w) \in \graph(\bar G \circ \bar F)$ and $w \defeq (x, G(x)y) = Rp$ and observe, from the uniqueness of $y$ that $p$ is the unique point in $\graph(\bar G) \circ \bar F$ satisfying $Rp=w$.
    We define
    \[
        \inv R_p: R \graph(\bar G \circ \bar F) \to \graph(\bar G \circ \bar F),
        \quad
        \inv R_p(\alt x_1, \alt z) \defeq (\alt x_1, (\alt x_1, \alt z)).
    \]
    This is clearly a Fréchet differentiable inverse selection of $R$ in any neighborhood $U_p \subset R \graph(\bar G \circ \bar F)$ of $w$.
    Another application of \cref{lemma:gderiv:cone-linear} now yields
    \begin{equation*}
        T_{\graph H}(x, z)
        = \{(\dir x, \dir z) \mid \dir z \in \{[G'(x)\dir x]y\}+G(x)DF(x|y)\dir x\},
    \end{equation*}
    and hence the claim follows.
\end{proof}

We are now ready to prove the second-order chain rule for the Clarke subdifferential promised before \cref{lemma:gderiv:cartesian-product}.
\begin{corollary}[second-order chain rule for Clarke subdifferentials]
    \label{cor:gderiv:second-clarke}
    Let $X,Y$ be Banach spaces, $f:Y\to \R$ be locally Lipschitz continuous, and $S:X\to Y$ be twice differentiable at $x \in X$. Set $h:X\to \R$, $h \defeq f \circ S$.
    If there exists a neighborhood $U$ of $x$ such that
    \begin{enumerate}
        \item\label{item:gderiv:second-clarke:i} $f$ is Clarke regular at $S(\alt x)$ for all $\alt x\in U$;
        \item\label{item:gderiv:second-clarke:ii} $S'(\alt x)$ has a right-inverse $\rinv{S'(\tilde x)} \in \linear(Y; X)$ for all $\tilde x\in U$;
        \item\label{item:gderiv:second-clarke:iii} the mapping $\alt x \mapsto \rinvstar{S'(\alt x)}$ is Fréchet differentiable at $x$;
    \end{enumerate}
    then for all $x^*\in \partial_C h(x) = S'(x)^*\partial_C f(S(x))$, we have
    \begin{equation*}
        D[\subdiff_C h](x|x^*)(\dir x) = \{y^*S''(x)\dir x\} + S'(x)^* D[\subdiff_C f](S(x)|y^*)(S'(x)\dir x)
        \quad (\dir x \in X)
    \end{equation*}
    for the unique $y^*\in \partial_C f(S(x))$ with $S'(x)^*y^* = x^*$ and $y^*S''(x)\dir x\in X^*$ defined by
    \begin{equation*}
        X\ni \Delta \alt x \mapsto \dualprod{y^*}{[S''(x)\Delta x]\Delta \alt x}_Y.
    \end{equation*}
\end{corollary}

\begin{proof}
    The expression for $\partial_C h(\alt x)$ for all $\alt x \in U$ follows from \cref{thm:clarke:chain} and the assumption~\ref{item:gderiv:second-clarke:i}.
    Let now $G:X\to \linear(Y^* ; X^*)$, $G(\alt x) \defeq S'(\alt x)^*$. Then $G$ is Fréchet differentiable at $x$ by the twice differentiability of $S$ at $x$, and, by the assumption~\cref{item:gderiv:second-clarke:ii}, has the left-inverse $\rinvstar{S'(\alt x)}$ for all $\alt x \in U$.
    The latter establishes the uniqueness of $y^*$.
    Together with assumption~\ref{item:gderiv:second-clarke:iii} we can now apply \cref{thm:gderiv:product} to $F=\subdiff_C f \circ S$ and $G$ to obtain
    \begin{equation*}
        D[\subdiff_C h](x|x^*)(\dir x) =
        \{[G'(x)\dir x]y^*\} + S'(x)^*D[(\subdiff_C f) \circ S](x|y^*)(\dir x)
        \quad
        (\dir \alt x \in X).
    \end{equation*}
    Observe that by assumption, $S$ is a fortiori continuously differentiable at $x$.
    Moreover, $S'(x)$ has a bounded right-inverse by \ref{item:gderiv:second-clarke:ii}.
    We may therefore apply \cref{thm:gderiv:inner} to $G=\subdiff_C f$ and $F=S$ to obtain
    \begin{equation*}
        D[(\subdiff_C f) \circ S](x|y^*)(\dir x)
        =D[\subdiff_C f](S(x)|y^*)(S'(x)\dir x)
        \qquad
        (\dir x \in X).
    \end{equation*}
    To show the claimed expression for the first term, let $\Delta \alt x\in X$ be arbitrary. By the linearity and continuity of the duality pairing and of $G'(x) = S'(x)^*\in \linear(Y^*; X^*)$, we can then write
    \begin{equation}
        \label{eq:gderiv:second-deriv-weak-form}
        \begin{aligned}[t]
            \dualprod{[G'(x)\Delta x]y^*}{\Delta \alt x}_X &= \dualprod{[\lim_{t\downto0}t^{-1}(G(x+t\Delta x)-G(x))]y^*}{\Delta \alt x}_X\\
            &=\lim_{t\downto 0}t^{-1}\dualprod{[S'(x+t\Delta x)^*-S'(x)^*]y^*}{\Delta \alt x}_X\\
            &=\lim_{t\downto 0}t^{-1}\dualprod{y^*}{[S'(x+t\Delta x)-S'(x)]\Delta \alt x}_Y\\
            &=\dualprod{y^*}{[S''(x)\Delta x]\Delta \alt x}_Y.
        \end{aligned}
    \end{equation}
    Since $\Delta \alt x$ was arbitrary and the expression on the right-hand side defines a bounded linear functional due to $S'':X\to \linear(X; \linear(X; Y))$, this establishes the claim.
\end{proof}

\chapter{Calculus for the Fr{\'e}chet coderivative}
\label{chap:cofrechet}

We continue with calculus rules for the Fréchet coderivative, which we recall is defined for $F:X\setto Y$ as
\begin{equation*}
    \frechetCod F(x|y): Y^* \setto X^*,\qquad
    \frechetCod F(x|y)(y^*) \defeq \setof{ x^*\in X^* }{ (x^*, -y^*) \in \frechetNormal_{\graph F}(x, y)}.
\end{equation*}
As in \cref{chap:gderiv}, we start with the relevant regularity concept, then prove the fundamental lemmas, and finally derive the calculus rules.

\section{Semi-codifferentiability}

Let $X,Y$ be Banach spaces.
We say that $F$ is \term[mapping!codifferentiable!semi-]{semi-codifferentiable} at $x \in X$ for $y \in F(x)$ if for each $y^* \in Y^*$ there exists some $x^* \in \frechetCod F(x|y)(y^*)$ satisfying
\begin{equation}
    \label{eq:cofrechet:semi-codiff}
    \lim_{\graph F \ni (x_k, y_k) \to (x, y)} \frac{\dualprod{x^*}{x_k-x}_X-\dualprod{y^*}{y_k-y}_Y}{\norm{(x_k-x, y_k-y)}_{X\times Y}}=0.
\end{equation}
Recalling \eqref{eq:cones:def-epsiloncone}, this is equivalent to requiring that $-x^* \in \frechetCod F(x|y)(-y^*)$ as well.
For single-valued mappings and their inverses, we have the following characterization.

\begin{lemma}
    \label{lemma:cofrechet:regularity:single}
    Let $X,Y$ be Banach spaces and let $F: X \to Y$ be single-valued.
    If $F$ is Fréchet differentiable at $x\in X$, then
    \begin{enumerate}
        \item\label{item:cofrechet:regularity:single}
        $F$ is semi-codifferentiable at $x$ for $y=F(x)$.
    \end{enumerate}
    If, moreover, $F'(x) \in \linear(X; Y)$ has a left-inverse $\linv{F'(x)} \in \linear(Y; X)$,
    then
    \begin{enumerate}[resume*]
        \item\label{item:cofrechet:regularity:single:inv}
        $\inv F$ is semi-codifferentiable at $y=F(x)$ for $x$.
    \end{enumerate}
\end{lemma}

\begin{proof}
    Recalling from \cref{thm:graphical:single} that $\frechetCod F(x|y)(y^*) = \{F'(x)^*y^*\}$  when $y=F(x)$, the claim \cref{item:cofrechet:regularity:single} follows immediately from the observation above that semi-codifferentiability is equivalent to the existence for all $y^*$ of $x^* \in \frechetCod F(x|y)(y^*)$ such that $-x^* \in \frechetCod F(x|y)(-y^*)$ as well.

    As for \cref{item:cofrechet:regularity:single:inv}, recalling the inverse relationships of \cref{lemma:graphical:inverse} and again using \cref{thm:graphical:single}, we have that
    $\frechetCod \inv F(y|x)(x^*) = \{y^* \mid x^* = F'(x)^*y^*\}$.
    Moreover, we recall that for a left-inverse $\linv{F'(x)}$ of $F'(x)$, the operator $\linvstar{F'(x)}$ is a right-inverse of $F'(x)^*$.
    Thus, for any $x^*$, we have $y^* \defeq \linvstar{F'(x)} x^* \in \frechetCod \inv F(y|x)(x^*)$, and, by linearity, $-y^* \in \frechetCod \inv F(y|x)(-x^*)$. Hence $F^{-1}$ is semi-codifferentiable at $y$ for $x$.
\end{proof}

\section{Cone transformation formulas}

In the following, we consider more general $\eps$-normal cones for $\eps\geq 0$, as these results will be needed later in \cref{chap:colimiting} for proving the corresponding expressions for the limiting normal cone. We refer to \cref{sec:gderiv:cones} for the definition of a continuous inverse selection.

\begin{lemma}
    \label{lemma:cofrechet:cone-linear}
    Let $X,Y$ be Banach spaces, $C\subset Y$, and $R\in \linear(Y; X)$.
    If $y \in \closure C$ admits a Lipschitz continuous inverse selection of $R$ at $x=Ry$ with factor $L_x$, then for all $\eps\geq 0$,
    \begin{equation}
        \label{eq:cofrechet:cone-linear:epsilon}
        \frechetNormal_{RC}^{\epsilon/\norm{R}_{\linear(Y; X)}}(x)
        \subset \{x^* \in X^* \mid R^* x^* \in \frechetNormal_C^\epsilon(y)\}
        \subset \frechetNormal_{RC}^{\epsilon L_x}(x).
    \end{equation}
    In particular,
    \begin{equation*}
        \frechetNormal_{RC}(x) = \{x^* \in X^* \mid R^* x^* \in \frechetNormal_C(y)\}.
    \end{equation*}%
\end{lemma}

\begin{proof}
    By the definition \eqref{eq:cones:def-epsiloncone}, $x^* \in \frechetNormal_{RC}^{\alt\epsilon}(x)$ for a given $\alt\epsilon>0$ if and only if
    \begin{equation}
        \label{eq:graphical:frechetnormal-outer-composition-linear:0.0}
        \limsup_{RC \ni Ry_k \to Ry} \frac{\dualprod{R^*x^*}{y_k-y}_Y}{\norm{R(y_k-y)}_X} \le \alt\epsilon.
    \end{equation}
    Since $R$ is continuous, $C \ni y_k \to y$ implies that $RC \ni Ry_k \to Ry$. Furthermore, the expression inside the limit in \eqref{eq:graphical:frechetnormal-outer-composition-linear:0.0} is invariant under perturbations in $\kernel R$, and hence \eqref{eq:graphical:frechetnormal-outer-composition-linear:0.0} implies that
    \[
        \limsup_{C \ni y_k \to y} \frac{\dualprod{R^*x^*}{y_k-y}_Y}{\norm{R(y_k-y)}_X} \le \alt\epsilon.
    \]
    Thus $x^* \in \frechetNormal_{RC}^{\alt\epsilon}(x)$ implies for every $\epsilon'>\alt\epsilon$ the existence of a $\delta>0$ such that
    \begin{equation}
        \label{eq:graphical:frechetnormal-outer-composition-linear:1.0}
        \dualprod{R^*x^*}{y_k-y}_Y \le \epsilon' \norm{R(y_k-y)}_X
        \quad (y_k \in \B(y, \delta) \isect C).
    \end{equation}
    Similarly, by definition, $R^*x^* \in \frechetNormal_C^{\epsilon}(y)$ if and only if
    \begin{equation}
        \label{eq:graphical:frechetnormal-outer-composition-linear:2}
        \limsup_{C \ni y_k \to y} \frac{\dualprod{R^*x^*}{y_k-y}_Y}{\norm{y_k-y}_Y} \le \epsilon.
    \end{equation}
    Now if $x^* \in \frechetNormal_{RC}^{\tilde\epsilon}(x)$, then using \eqref{eq:graphical:frechetnormal-outer-composition-linear:1.0} and estimating $\norm{R(y_k-y)}_X \le \norm{R}_{\linear(Y;X)}\norm{y_k-y}_Y$ yields \eqref{eq:graphical:frechetnormal-outer-composition-linear:2} for $\epsilon=\epsilon'\norm{R}_{\linear(Y;X)}$ for all $\epsilon'>\alt\epsilon$. Hence the first inclusion in \eqref{eq:cofrechet:cone-linear:epsilon} holds by taking $\alt\epsilon=\epsilon/\norm{R}_{\linear(Y;X)}$ and letting $\epsilon' \downto \alt\epsilon$.

    For the second inclusion, let $R^*x^* \in \frechetNormal_C^{\epsilon}(y)$.
    Then \eqref{eq:graphical:frechetnormal-outer-composition-linear:2} holds for $y_k=\inv R_y(x_k)$ for any $U_y \ni x_k \to x$.
    We also observe using $R\inv R_y(x_k)=x_k$ and $R\inv R_y(x)=x$ that $x^* \in \frechetNormal_{RC}^{\alt\epsilon}(x)$ if and only if
    \begin{equation}
        \label{eq:graphical:frechetnormal-outer-composition-linear:1}
        \limsup_{RC \ni x_k \to x} \frac{\dualprod{R^*x^*}{\inv R_y(x_k)-\inv R_y(x)}_X}{\norm{x_k-x}_X} \le \alt\epsilon.
    \end{equation}
    Since $\inv R_y$ is assumed to be Lipschitz continuous at $x$, we have
    \[
        L_x
        \ge \limsup_{k \to \infty} \frac{\norm{\inv R_y(x_k)-\inv R_y(x)}_Y}{\norm{x_k-x}_X}
        = \limsup_{k \to \infty} \frac{\norm{y_k-y}_Y}{\norm{x_k-x}_X}.
    \]
    Thus, \eqref{eq:graphical:frechetnormal-outer-composition-linear:2} yields \eqref{eq:graphical:frechetnormal-outer-composition-linear:1} for $\alt\epsilon=\epsilon L_x$.
    We conclude that $x^* \in \frechetNormal_{RC}^{\alt\epsilon}(x)$, which yields the second inclusion in \eqref{eq:cofrechet:cone-linear:epsilon}.
\end{proof}

\begin{remark}[polarity and qualification condition in finite dimensions]
    In finite dimensions, \cref{lemma:cofrechet:cone-linear} for $\epsilon=0$ could also be proved with the help of the polarity relationships $\frechetNormal_{RC}(x)=\polar{T_{RC}(x)}$ and $\frechetNormal_{C}(y)=\polar{T_{C}(y)}$ from \cref{lemma:cones:fundamental-polar}.
    Furthermore, the existence of a family of continuous selections could be replaced by a qualification condition as in \cref{rem:gderiv:qc}.
\end{remark}

We are now ready to prove the fundamental composition lemma, this time for the Fréchet normal cone.
We say that $G: Y \setto Z$ is \term[mapping!continuous!inner Lipschitz]{inner Lipschitz at $y$ for $z$} if for some $L, \delta>0$ and all $\tilde y \in B(\delta, y)$ we have
\[
    \inf_{\tilde z \in G(\tilde y)} \norm{\tilde z - z}_Z \le L\norm{\tilde y - y}_Y.
\]
For single-valued mappings, this property obviously reduces to Lipschitz continuity at $y$. We will return to further Lipschitz-like properties of set-valued mappings in \cref{chap:regularity}.
The statement here assumes that the product spaces are endowed with the Euclidean product norm, i.e., $\norm{(x,y)}_{X\times Y}^2 = \norm{x}_X^2+\norm{y}_Y^2$.
\begin{lemma}[fundamental lemma on compositions]
    \label{lemma:cofrechet:fundamental}
    Let $X,Y,Z$ be Banach spaces and
    \begin{equation*}
        C \defeq \{(x, y, z) \mid y \in F(x),\, z \in G(y), \, x\in X\}
    \end{equation*}
    for $F: X \setto Y$, and $G: Y \setto Z$. Let $(x,y,z)\in C$.
    \begin{enumerate}
        \item\label{item:cofrechet:fundamental:G}
            If  $G$ is semi-codifferentiable and inner Lipschitz at $y$ for $z$ with factor $L > 0$, then
            \[
                K_\epsilon \subset \frechetNormal_C^\epsilon(x, y, z) \subset K_{\sqrt{1+L^2}\epsilon}
            \]
            for all $\epsilon \ge 0$ and
            \[
                K_\epsilon \defeq
                \left\{
                    (x^*, y^*, z^*)
                    \,\middle|\,
                    \begin{array}{l}
                    x^* \in \frechetCod_\epsilon F(x|y)(-\tilde y^*-y^*),\, z^*\in Z^*
                    \\
                    \tilde y^* \in \frechetCod G(y|z)(z^*),\, -\tilde y^* \in \frechetCod G(y|z)(-z^*)
                    \end{array}
                \right\}.
            \]
        \item\label{item:cofrechet:fundamental:invF}
            If $\inv F$ is semi-codifferentiable and inner Lipschitz at $y$ for $x$ with factor $\ell>0$, then
            \[
                 Q_\epsilon \subset \frechetNormal_C^\epsilon(x, y, z) \subset Q_{\sqrt{1+\ell^2}\epsilon}
            \]
            for all $\epsilon \ge 0$ and
            \[
                Q_\epsilon
                \defeq
                \left\{
                    (x^*, y^*, z^*)
                    \,\middle|\,
                    \begin{array}{l}
                    x^* \in \frechetCod F(x|y)(-\tilde y^*-y^*),
                    \\
                    -x^* \in \frechetCod F(x|y)(\tilde y^*+y^*),
                    \\
                    -\tilde y^* \in \frechetCod_\epsilon G(y|z)(-z^*),\, z^*\in Z^*
                    \end{array}
                \right\}.
            \]
    \end{enumerate}
\end{lemma}

\begin{proof}
    We recall that $(x^*, y^*, z^*) \in \frechetNormal_C^\epsilon(x, y, z)$ if and only if
    \begin{equation}
        \label{eq:cofrechet:fundamental:epsilon:1}
        \limsup_{C \ni (x_k, y_k, z_k) \to (x, y, z)}
        \frac{\dualprod{x^*}{x_k-x}_X+\dualprod{y^*}{y_k-y}_Y+\dualprod{z^*}{z_k-z}_Z}{
        \norm{(x_k,y_k,z_k)-(x,y,z)}_{X\times Y \times Z}}
        \le \epsilon.
    \end{equation}

    In case \ref{item:cofrechet:fundamental:G}, the semi-codifferentiability of $G$ implies that for some $\tilde y^* \in \frechetCod G(y|z)(z^*)$ we have $-\tilde y^* \in \frechetCod G(y|z)(-z^*)$ or equivalently
    \begin{equation*}
        \lim_{\graph G \ni (y_k, z_k) \to (y, z)} \frac{\dualprod{\tilde y^*}{y_k-y}_Y-\dualprod{z^*}{z_k-z}_Z}{\norm{(y_k,z_k)-(y, z)}_{Y\times Z}}=0.
    \end{equation*}
    Thus \eqref{eq:cofrechet:fundamental:epsilon:1} holds if and only if for some $\tilde y^* \in \frechetCod G(y|z)(z^*)$ with $-\tilde y^* \in \frechetCod G(y|z)(-z^*)$ we have
    \begin{equation}
        \label{eq:cofrechet:fundamental:epsilon:2}
        \limsup_{C \ni (x_k, y_k, z_k) \to (x, y, z)}
        \frac{\dualprod{x^*}{x_k-x}_X+\dualprod{\tilde y^*+y^*}{y_k-y}_Y}{
        \norm{(x_k,y_k,z_k)-(x,y,z)}_{X\times Y\times Z}}
        \le \epsilon.
    \end{equation}
    But this follows from $x^* \in \frechetCod_\epsilon F(x|y)(-\tilde y^*-y^*)$, which yields the first inclusion in \ref{item:cofrechet:fundamental:G}.
    For the second inclusion, taking large enough $k$, we use the inner Lipschitz assumption to choose $z_k \in G(y_k)$ such that $\norm{z_k-z}_Z \le (L+1/k)\norm{y_k-y}_Y$. Then \eqref{eq:cofrechet:fundamental:epsilon:2} implies for large enough $k$ and any $\epsilon' > \epsilon$ that
    \[
        \begin{split}
        \dualprod{x^*}{x_k-x}_X+\dualprod{\tilde y^*+y^*}{y_k-y}_Y
        &
        \le \epsilon' \norm{(x_k,y_k,z_k)-(x,y,z)}_{X\times Y\times Z}
        \\
        &
        \le \epsilon' \sqrt{1+L^2} \norm{(x_k,y_k)-(x,y)}_{X\times Y}.
        \end{split}
    \]
    Keeping in mind that $\epsilon'>\epsilon$ was arbitrary, this yields
    \begin{equation*}
        \limsup_{\graph F \ni (x_k, y_k) \to (x, y)}
        \frac{\dualprod{x^*}{x_k-x}_X+\dualprod{\tilde y^*+y^*}{y_k-y}_Y}{
        \norm{(x_k,y_k)-(x,y)}_{X\times Y}}
        \le \sqrt{1+L^2} \epsilon,
    \end{equation*}
    which by definition is equivalent to $x^* \in \frechetCod_{\sqrt{1+L^2}\epsilon} F(x|y)(-\tilde y^*-y^*)$.

    In case \ref{item:cofrechet:fundamental:invF}, the semi-codifferentiability of $\inv F$ implies that there exists a $\tilde y^* \in \{-y^*\} + \frechetCod \inv F(y|x)(-x^*)$, i.e., satisfying $x^* \in \frechetCod F(x|y)(-\tilde y^*-y^*)$, such that $-x^* \in \frechetCod F(x|y)(\tilde y^*+y^*)$ as well.
    This is again equivalently written as
    \begin{equation*}
        \lim_{\graph \inv F \ni (y_k, x_k) \to (y, x)} \frac{-\dualprod{\tilde y^*+y^*}{y_k-y}_Y-\dualprod{x^*}{x_k-x}_X}{\norm{(y_k,x_k)-(y, x)}_{Y \times X}}=0.
    \end{equation*}
    Thus \eqref{eq:cofrechet:fundamental:epsilon:1} holds if and only if for some $\tilde y^*$ we have both $x^* \in \frechetCod F(x|y)(-\tilde y^*-y^*)$ and $-x^* \in \frechetCod F(x|y)(\tilde y^*+y^*)$, as well as
    \begin{equation*}
        \label{eq:cofrechet:fundamental:epsilon:3}
        \limsup_{C \ni (x_k, y_k, z_k) \to (x, y, z)}
        \frac{\dualprod{z^*}{z_k-z}_Z-\dualprod{\tilde y^*}{y_k-y}_Y}{
        \norm{(x_k,y_k,z_k)-(x,y,z)}_{X\times Y\times Z}}
        \le \epsilon.
    \end{equation*}
    But this follows from $-\tilde y^* \in \frechetCod_\epsilon G(y|z)(-z^*)$, which yields the first inclusion in \ref{item:cofrechet:fundamental:invF}. For the second inclusion, we again use the inner Lipschitz assumption as in case \ref{item:cofrechet:fundamental:G}.
\end{proof}

For the remaining results, we fix $\eps=0$.
If one of the two mappings is single-valued, \cref{lemma:cofrechet:fundamental} yields the following two special cases.

\begin{corollary}[fundamental lemma on compositions: single-valued outer mapping]
    \label{lemma:cofrechet:fundamental:single-outer}
    Let $X,Y,Z$ be Banach spaces and
    \begin{equation*}
        C \defeq \{(x, y, G(y)) \mid y \in F(x),\, x\in X\}
    \end{equation*}
    for $F: X \setto Y$ and $G: Y \to Z$.
    If $(x, y, z) \in C$ and $G$ is Fréchet differentiable at $y$, then
    \begin{equation*}
        \frechetNormal_C(x, y, z)=\{(x^*, y^*, z^*) \mid
            x^* \in \frechetCod F(x|y)(-[G'(y)]^*z^*-y^*),\, y^* \in Y^*,\, z^*\in Z^*
        \}.
    \end{equation*}
\end{corollary}

\begin{proof}
    We first use \cref{lemma:cofrechet:regularity:single}\,\cref{item:cofrechet:regularity:single} to show the semi-codifferentiability of $G$ at $y$ for $z$.
    The assumed Fréchet differentiability at $y$ implies that $G$ is Lipschitz and hence inner Lipschitz at $y$ for $z=G(y)$.
    Thus we may apply \cref{lemma:cofrechet:fundamental}\,\ref{item:cofrechet:fundamental:G} to obtain an explicit expression for $\frechetNormal_C(x, y, z)$, into which we insert the expression given by \cref{thm:graphical:single} for $\frechetCod G(y|z)(z^*)$.
\end{proof}

The corresponding result for a single-valued inner mapping is not quite as straightforward unless we assume full invertibility of $F'(x)$. We first do so, and then relax the assumption to mere right-invertibility.

\begin{corollary}[initial lemma on compositions: single-valued inner mapping]
    \label{lemma:cofrechet:fundamental:single-inner:initial}
    Let $X,Y,Z$ be Banach spaces and
    \begin{equation*}
        C \defeq \{(x, y, z) \mid y=F(x),\, z \in G(y),\, x\in X\}
    \end{equation*}
    for $F: X \to Y$ and $G: Y \setto Z$.
    If $(x, y, z) \in C$ and $F$ is continuously differentiable at $x$
    such that $F'(x)$ has an inverse $\inv{F'(x)} \in \linear(Y; X)$, then
    \begin{equation*}
        \frechetNormal_C(x, y, z)
        =\{(F'(x)^*(-\tilde y^*-y^*), y^*, z^*) \mid
            -\tilde y^* \in \frechetCod G(y|z)(-z^*),\, y^* \in Y^*,\, z^*\in Z^*
        \}.
    \end{equation*}
\end{corollary}

\begin{proof}
    Similarly to the previous proof, we apply \cref{thm:graphical:single} and \cref{lemma:cofrechet:regularity:single}\,\cref{item:cofrechet:regularity:single:inv} to $F$ to prove its semi-codifferentiability and then use \cref{lemma:cofrechet:fundamental}\,\ref{item:cofrechet:fundamental:invF}.
    To prove that $\inv F$ is inner Lipschitz at $y=F(x)$ for $x$, we apply the inverse function theorem (\cref{thm:inversefunctiontheorem}), which shows that $\inv F$ exists and is continuously differentiable. Then \cref{lem:variation:c1-lipschitz} shows that $\inv F$ is locally Lipschitz at $y$, which implies that $\inv F$ is inner Lipschitz, as required.
\end{proof}

\begin{lemma}[fundamental lemma on compositions: single-valued inner mapping]
    \label{lemma:cofrechet:fundamental:single-inner}
    Let $X,Y,Z$ be Banach spaces and
    \begin{equation*}
        C \defeq \{(x, y, z) \mid y=F(x),\, z \in G(y),\, x\in X\}
    \end{equation*}
    for $F: X \to Y$ and $G: Y \setto Z$.
    If $(x, y, z) \in C$, the mapping $F$ is continuously differentiable at $x$,
    and $F'(x)$ has a {right} inverse $\rinv{F'(x)} \in \linear(Y; X)$, then
    \begin{equation*}
        \frechetNormal_C(x, y, z)
        =\{(F'(x)^*(-\tilde y^*-y^*), y^*, z^*) \mid
            -\tilde y^* \in \frechetCod G(y|z)(-z^*),\, \tilde y^*, y^* \in Y^*,\, z^*\in Z^*
        \}.
    \end{equation*}
\end{lemma}

\begin{proof}
    Let $\bar F: X \to Y \times \kernel F'(x)$, $\bar F(x) \defeq (F(x), Px)$ for $P \defeq \Id - \rinv{F'(x)}F'(x)$.
    Also let $\bar G: Y \times \kernel F'(x) \setto Z$ be defined by $\bar G(y, q) \defeq G(y)$.
    Then, by \cref{lemma:gderiv:inverse:single:left}, $\bar F$ is invertible, and by either the proof of the lemma or by the inverse function theorem (\cref{thm:inversefunctiontheorem}), $\bar F'(x)$ has an inverse $\inv{\bar F'(x)} \in \linear(Y \times \kernel F'(x); X)$.
    Directly from the definition, we deduce that for every $z^* \in Z^*$,
    \[
        \frechetCod \bar G(y|z)(z^*) = \{ (z^*, q^*) \in Z^* \times [\kernel F'(x)]^* \mid  z^* \in  \frechetCod \bar G(y|z)(z^*),\, q^* = 0 \}.
    \]
    We apply \cref{lemma:cofrechet:fundamental:single-inner:initial} to $\bar F$ and $\bar G$ to obtain for
    \begin{equation*}
        \bar C \defeq \{(x, (y, q), z) \mid (y, q)=\bar F(x),\, z \in \bar G(y, q),\, x\in X\}
    \end{equation*}
    the expression
    \begin{equation*}
        \begin{aligned}[t]
        \frechetNormal_{\bar C}&(x, (y, q), z)
        \\
        &
        =
        \left\{
            (\bar F'(x)^*(-\tilde y^*-y^*, -\tilde q^*-q^*), (y^*, q^*), z^*)
        \,\middle|\,
            \begin{array}{l}
                -(\tilde y^*, \tilde q^*) \in \frechetCod \bar G(y|z)(-z^*);\, z^*\in Z^*\\
                \tilde y^*, y^* \in Y^*;\,
                \tilde q^*, q^* \in [\kernel F'(x)]^*
            \end{array}
        \right\}
        \\
        &
        =
        \left\{
            (F'(x)^*(-\tilde y^*-y^*) + P^*(-\tilde q^*-q^*), (y^*, q^*), z^*)
        \,\middle|\,
            \begin{array}{l}
                -\tilde y^* \in \frechetCod G(y|z)(-z^*);
                \\
                \tilde  y^*, y^* \in Y^*;\, z^*\in Z^*\\
                q^* \in [\kernel F'(x)]^*;\,
                \tilde q^*=0
            \end{array}
        \right\}.
        \end{aligned}
    \end{equation*}
    Now we write $C=R\bar C$ for $R(\tilde x, (\tilde y, \tilde q), \tilde z) \defeq (\tilde x, \tilde y, \tilde z)$, and observe that $\inv R_{(x, (y, q), z)}(\tilde x, \tilde y, \tilde z) \defeq (\tilde x, (\tilde y, P\tilde x), \tilde z)$ for $q=Px$ is a Lipschitz inverse selection of $R$ at $(x, (y, q), z)$ to $\bar C$.
    Therefore \cref{lemma:cofrechet:cone-linear} establishes
    \[
        \begin{aligned}
        \frechetNormal_{C}(x, y, z)
        &
        =
        \Isect_{(\tilde x, (\tilde y, \tilde q), \tilde z) \in \bar C, R(\tilde x, (\tilde y, \tilde q), \tilde z)=(x, y, z)}
        \{ (x^*, y^*, z^*) \mid (x^*, (y^*, 0), z^*) \in \frechetNormal_{\bar C}(\tilde x, (\tilde y, \tilde q), \tilde z)\}
        \\
        &
        =\{(F'(x)^*(-\tilde y^*-y^*), y^*, z^*) \mid
            -\tilde y^* \in \frechetCod G(y|z)(-z^*),\, y^* \in Y^*
        \},
        \end{aligned}
    \]
    as claimed.
\end{proof}

\section{Calculus rules}\label{sec:cofrechet:calculus}

Using the above lemmas, we again obtain calculus rules. The proofs are similar to those in \cref{sec:gderiv:calculus}, and we only note the differences.
Overall, in the assumptions of results based on the cone transformation lemmas (\cref{lemma:cofrechet:cone-linear} in place of \cref{lemma:gderiv:cone-linear}) and the fundamental single-valued outer composition results (\cref{lemma:cofrechet:fundamental:single-outer} in place of \cref{lemma:gderiv:fundamental:single-outer}), the Fréchet differentiability requirement at a point $x$ needs to be amended by a Lipschitz requirement at $x$. The assumptions of results that use the fundamental single-valued inner composition results (\cref{lemma:cofrechet:fundamental:single-inner} in place of \cref{lemma:gderiv:fundamental:single-inner}) are unchanged as both variants require the continuous differentiability at $x$, which is stronger than the assumptions of the cone transformation lemmas.

\begin{theorem}[addition of a single-valued differentiable mapping]
    \label{thm:cofrechet:addition}
    Let $X,Y$ be Banach spaces, $G:X\to Y$, and $F:X\setto Y$.
    Let $x \in X$ and $y\in H(x)\defeq F(x)+G(x)$.
    If $G$ is Fréchet differentiable and Lipschitz at $x$, then
    \begin{equation*}
        \frechetCod H(x|y)(y^*)  = \frechetCod F(x|y-G(x))(y^*) + \{[G'(x)]^* y^*\}
        \qquad (y^* \in Y^*).
    \end{equation*}%
\end{theorem}

\begin{proof}
    We have $\graph H=RC$ for $C$ and $R$ given by \eqref{eq:gderiv:addition:c-r} in the proof of \cref{thm:gderiv:addition}.
    Since $G$ is Lipschitz at $x$, the mapping $\inv R_p$ given by \eqref{eq:gderiv:addition:ry} for the unique $p$ with $(x,y)=Rp$ is Lipschitz at $(x,y)$ with factor $L_{(x,y)}$.
    We may therefore apply \cref{lemma:cofrechet:cone-linear} with $\epsilon=0$ in place of \cref{lemma:gderiv:cone-linear} in the proof of \cref{thm:gderiv:addition} to obtain
    \begin{equation*}
        \frechetNormal_{\graph H}(x, y)=\{(x^*, y^*) \mid (y^*, x^*, y^*) \in \frechetNormal_C(y-G(x), x, G(x))\}.
    \end{equation*}
    Moreover, $C$ given in \eqref{eq:gderiv:addition:c-r} coincides with the $C$ defined in \cref{lemma:cofrechet:fundamental:single-outer} with $\inv F$ in place of $F$.
    Since $G$ is Fréchet differentiable at $x$, we may apply the corollary and continue with the expression from \cref{lemma:graphical:inverse} for $\frechetCod\inv F$, to obtain
    \begin{equation*}
        \begin{aligned}
        \frechetNormal_C(u, x, v)
        &
        = \{(u^*, x^*, v^*) \mid
            u^* \in \frechetCod \inv F(u|x)(-[G'(x)]^*v^*-x^*),\, v^*\in Y^*,\, x^*\in X^*
        \}
        \\
        &=\{(u^*, x^*, v^*) \mid
            [G'(x)]^*v^*+x^* \in \frechetCod F(x|u)(-u^*),\, v^*\in Y^*,\, x^*\in X^*
        \}.
        \end{aligned}
    \end{equation*}
    Thus
    \begin{equation*}
        \begin{aligned}
            \frechetCod H(x|y)(y^*)
            &
            =
            \{ x^* \in X^* \mid (-y^*, x^*, -y^*) \in \frechetNormal_C(y-G(x), x, G(x))\}
            \\
            &
            =
            \{ x^* \in X^* \mid -[G'(x)]^*y^*+x^* \in \frechetCod F(x|y-G(x))(y^*)\},
        \end{aligned}
    \end{equation*}
    which yields the claim.
\end{proof}

\begin{theorem}[outer composition with a single-valued differentiable mapping]
    \label{thm:cofrechet:outer}
    Let $X,Y$ be Banach spaces, $F:X\setto Y$, and $G:Y\to Z$.
    Let $x \in X$ and $z \in H(x) \defeq G(F(x))$ be given.
    If $G$ is Fréchet differentiable at $y \in F(x) \isect \inv G(\{z\})$ and locally left-invertible near $z$ such that the local left-inverse $\linv G$ is Fréchet differentiable at $z$, then
    \begin{equation*}
        \frechetCod H(x|z)(z^*)=\frechetCod F(x|y)([G'(y)]^*z^*)
        \qquad (z^* \in Z^*).
    \end{equation*}
\end{theorem}

\begin{proof}
    We have $\graph H= RC$ for $R$ and $C$ as given by \eqref{eq:gderiv:outer:RC} in the proof of \cref{thm:gderiv:outer}.
    Since $\linv G$ is Lipschitz at $z$, the inverse selection $\inv R_p: RC \to C$ constructed in \eqref{eq:gderiv:outer:invR} for $p=(x, y, z)$ is also Lipschitz at $(x, z)$.
    Applying \cref{lemma:cofrechet:cone-linear} then yields
    \begin{equation*}
        \frechetNormal_{\graph H}(x, z)= \{(x^*, z^*) \mid (x^*, 0, z^*) \in \frechetNormal_C(x, y, z) \}.
    \end{equation*}
    \Cref{lemma:cofrechet:fundamental:single-outer}, which requires the Fréchet differentiability of $G$ at $y$, then shows that
    \begin{equation*}
        \begin{aligned}
            \frechetCod H(x|z)(z^*) &= \{x^* \mid (x^*, 0, -z^*) \in \frechetNormal_C(x, y, z) \}
            \\
            &
            =
             \{x^* \mid x^* \in \frechetCod F(x|y)([G'(y)]^*z^*) \}.
        \end{aligned}
    \end{equation*}
    After further simplification, we arrive at the claimed expression.
\end{proof}

\begin{corollary}[outer composition with a linear operator]
    \label{cor:cofrechet:outer:linear}
    Let $X,Y,Z$ be Banach spaces, $A\in \linear(Y; Z)$, and $F:X\setto Y$.
    If $A$ has a left-inverse $\linv A$, then for any $x\in X$ and $z\in H(x)\defeq AF(x)$,
    \begin{equation*}
        \frechetCod H(x|z)(z^*)=
        \frechetCod F(x|y)(A^*z^*)
        \qquad (z^* \in Z^*)
    \end{equation*}
    for the unique $y \in Y$ such that $Ay=z$.
\end{corollary}

\begin{proof}
    We only need to verify that $G(y) \defeq Ay$ satisfies the assumptions of \cref{thm:cofrechet:outer}, which can be done exactly as in the proof of \cref{cor:gderiv:outer:linear}.
\end{proof}

\begin{theorem}[inner composition with a single-valued mapping]
    \label{thm:cofrechet:inner}
    Let $X,Y,Z$ be Banach spaces, $F: X\to Y$ and $G:Y\setto Z$.
    Let $x \in X$ and $z \in H(x)\defeq G(F(x))$ be given.
    If $F$ is continuously differentiable at $x$ and $F'(x)$ has a right-inverse $\rinv{F'(x)} \in \linear(Y; X)$, then
    \begin{equation*}
        \frechetCod H(x|z)(z^*)  = [F'(x)]^* \frechetCod G(F(x)|z)(z^*)
        \qquad (z^* \in Z^*).
    \end{equation*}%
\end{theorem}

\begin{proof}
    We have $\graph H=RC$ for $C$ and $R$ as given by \eqref{eq:gderiv:inner:RC} in the proof of \cref{thm:gderiv:inner}.
    Since $F$ is continuously differentiable at $x$, it is Lipschitz at $x$.
    Consequently $\inv R_p$, given in \eqref{eq:gderiv:inner:invRp}, is Lipschitz at $(x, z)$ for $p=(x,F(x),z) \in C$ the unique point with $Rp=(x, z)$.
    We can therefore apply \cref{lemma:cofrechet:cone-linear} to obtain
    \begin{equation*}
        \frechetNormal_{\graph H}(x, z)=\{(x^*, z^*) \mid (x^*, 0, z^*) \in \frechetNormal_C(x, F(x), z) \}.
    \end{equation*}
    On the other hand, since $F$ is continuously differentiable at $x$, \cref{lemma:cofrechet:fundamental:single-inner} implies that
    \begin{equation*}
        \frechetNormal_C(x, y, z)
        =\{(F'(x)^*(-\tilde y^*-y^*), y^*, z^*) \mid
            -\tilde y^* \in \frechetCod G(y|z)(-z^*),\, y^* \in Y^*
        \}.
    \end{equation*}
    Thus
    \begin{equation*}
        \begin{aligned}
            \frechetCod H(x|z)(z^*)
            &=
            \{ x^* \mid (x^*, 0, -z^*) \in \frechetNormal_C(x, F(x), z) \}
            \\
            &=
            \{ F'(x)^* \tilde y^* \mid \tilde y^* \in \frechetCod G(y|z)(z^*)\},
        \end{aligned}
    \end{equation*}
    which yields the claim.
\end{proof}

\begin{corollary}[inner composition with a linear operator]
    \label{cor:cofrechet:inner:linear}
    Let $X,Y,Z$ be Banach spaces, $A\in \linear(X; Y)$, and $G:Y\setto Z$.
    If $A$ has a right-inverse $\rinv A \in \linear(Y; X)$, then for any $x\in X$ and $z\in H(x)\defeq G(Ax)$,
    \begin{equation*}
        \frechetCod H(x|z)(z^*)  = A^* \frechetCod G(Ax|z)(z^*)
        \qquad (z^* \in Z^*).
    \end{equation*}%
\end{corollary}

We again apply this to the chain rule for the convex subdifferential (\cref{thm:convex:chain}). Note that for $\partial f:X\setto X^*$, we have $\frechetCod [\partial f](x|x^*): X^{**}\to X^*$.
\begin{corollary}[second-order chain rule for convex subdifferentials]
    \label{cor:cofrechet:second-convex}
    Let $X,Y$ be Banach spaces, $f:Y\to \Rbar$ be proper, convex, and lower semicontinuous, and $A\in\linear(X; Y)$ be such that $A$ has a right-inverse $\rinv{A} \in \linear(Y; X)$ and that $\range A \isect \interior \dom f \ne \emptyset$. Let $h\defeq f\circ A$.
    Then for any $x\in X$ and $x^*\in \partial h(x) = A^*\partial f(Ax)$,
    \begin{equation*}
        \frechetCod[\subdiff h](x|x^*)(x^{**})  = A^* \frechetCod[\subdiff f](Ax|y^*)(A^{**}x^{**})
        \qquad (x^{**}\in X^{**})
    \end{equation*}
    for the unique $y^* \in Y^*$ satisfying $A^*y^*=x^*$.
\end{corollary}
\begin{proof}
    The expression for $\partial h(x)$ follows from \cref{thm:convex:chain}, to which we apply \cref{cor:cofrechet:inner:linear} as well as \cref{cor:cofrechet:outer:linear} with $A^*$ in place of $A$, recalling that a right-inverse $\rinv A$ for $A$ produces the left-inverse $\linv{(A^*)} = (\rinv A)^*$ for $A^*$.
\end{proof}

\begin{remark}
    Comparing \cref{cor:cofrechet:second-convex} with \cref{cor:gderiv:second-convex}, we see that if $X$ is reflexive, the coderivative coincides with the graphical derivative due to the linearity of $A$.
\end{remark}

For the corresponding result for the Clarke subdifferential, we again need a product rule.
We start with the following lemma.

\begin{lemma}
    \label{lemma:cofrechet:cartesian-product}
    Let $X,Y$ be Banach spaces and $F:X\setto Y$.
    Define $\bar F:X\setto X\times Y$ by $\bar F(x) \defeq \{x\} \times F(x)$.
    Then, for all $x\in X$, $y \in F(x)$, $x^* \in X^*$, and $y^* \in Y^*$, we have
    \begin{equation*}
        \frechetCod \bar F(x|(x,y))(x^*, y^*)
        = \{x^*\} + \frechetCod F(x|y)(y^*).
    \end{equation*}
\end{lemma}

\begin{proof}
    The proof is analogous to \cref{lemma:gderiv:cartesian-product} for the graphical derivative.
    We have
    \begin{equation*}
        \graph \bar F = R \graph F
        \quad\text{for}\quad
        R(\alt x, \alt y) \defeq (\alt x, (\alt x, \alt y)).
    \end{equation*}
    Observe that $(x, y) \in \graph F$ is the unique point satisfying $R(x,y)=p$ for $p=(x, (x, y)) \in \graph \bar F$.
    We define
    \[
        \inv R_p: R \graph F \to \graph F,
        \quad
        \inv R_p(\alt x, (\alt z, \alt y)) \defeq (\alt x,  \alt y).
    \]
    This is clearly a Lipschitz inverse selection of $R$ at $p$ for any neighborhood $U_p \subset \graph \bar F$ of $(x, y)$.
    Therefore, by \cref{lemma:cofrechet:cone-linear}, we have
    \begin{equation*}
        \frechetNormal_{R\graph F}(x, (x, y))
        = \{(x_0^*, (-x^*, -y^*)) \mid (x_0^*-x^*, -y^*) \in \frechetNormal_{\graph F}(x, y)\},
    \end{equation*}
    which establishes the claim.
\end{proof}

\begin{theorem}[product rule]
    \label{thm:cofrechet:product}
    Let $X,Y,Z$ be Banach spaces, $G:X\to \linear (Y; Z)$ be Fréchet differentiable at $x \in X$, and $F:X\setto Y$.
    Assume that $G(\alt x) \in \linear(Y; Z)$ has a left-inverse $\linv{G(\alt x)} \in \linear(Z; Y)$ for all $\alt x$ near $x$ and that the mapping $\alt x \mapsto \linv{G(\alt x)}$ is Lipschitz at $x$.
    Let $z\in H(x)\defeq G(x)F(x) \defeq \Union_{y \in F(x)} G(x)y$ and let $y \in F(x)$ be the unique element satisfying $G(x)y=z$.
    Then
    \begin{equation*}
        \frechetCod H(x|z)(z^*)
        =
        \{([G'(x)\freevar]y)^*z^*\} + \frechetCod F(x|y)(G(x)^*z^*)
        \qquad (z^* \in Z^*),
    \end{equation*}
    where $([G'(x)\freevar]y)^*z^* \in X^*$ is defined by
    \begin{equation*}
        X\ni \dir x \mapsto \dualprod{z^*}{[G'(x)\dir x]y}_Z.
    \end{equation*}
\end{theorem}

\begin{proof}
    The proof is analogous to \cref{thm:gderiv:product} for the graphical derivative.
    We again have $\graph H=R \graph(\bar G \circ \bar F)$ for $\bar F(x) \defeq \{x\} \times F(x)$ as in \cref{lemma:cofrechet:cartesian-product},
    \begin{equation*}
        \bar G(\alt x, \alt y)=(\alt x, G(\alt x)\alt y),
        \quad\text{and}\quad
        R(\alt x_1, \alt x_2, \alt z) \defeq (\alt x_1, \alt z).
    \end{equation*}
    Since $G$ is Fréchet differentiable at $x$, we have
    \[
        \bar G'(x, y)(\dir x, \dir y) = (\dir x, [G'(x)\dir x]y + G(x)\dir y)\in X\times Z
    \]
    for any $(\dir x, \dir y) \in X\times Y$. A straightforward calculation then shows that for any $(x_0^*, z^*)\in X^*\times Z^*$,
    \begin{equation*}
        \dualprod{(x_0^*,z^*)}{G'(x,y)(\dir x, \dir y)}_{X\times Z} = \dualprod{x_0^*}{\dir x}_X + \dualprod{z^*}{[G'(x)\dir x]y}_Z + \dualprod{G(x)^*z^*}{\dir y}_Y.
    \end{equation*}
    Since the right-hand side defines a bounded linear operator on $X\times Y$, this implies that $\bar G$ is Fréchet differentiable at $(x, y)$ with
    \begin{equation*}
        \bar G'(x, y)^*(x_0^*, z^*)=\left(x_0^*+([G'(x)\freevar]y)^*z^*, G(x)^* z^*\right)\in X^*\times Y^*.
    \end{equation*}
    Furthermore, since $G(\alt x)$ is assumed to be left-invertible for $\alt x$ near $x$, the mapping $Q: (\alt x, \alt z) \mapsto (\alt x, \linv{G(\alt x)}\alt z)$ is a left-inverse of $\bar G$, which is Lipschitz at $(x,z)$ since $\alt x \mapsto \linv{G(\alt x)}$ is Lipschitz at $x$.
    Thus we may apply \cref{thm:cofrechet:outer} to $\bar G$ and $\bar F$, which together with
    \cref{lemma:cofrechet:cartesian-product} yields
    \begin{equation*}
        \begin{aligned}
            \frechetCod[\bar G \circ \bar F](x|(x, z))(x_0^*, z^*)
            &=  \frechetCod \bar F(x|(x, y))(\bar G'(x, y)^*(x_0^*, z^*))
            \\
            &= \frechetCod \bar F(x|(x, y))\left(x_0^*+([G'(x)\freevar]y)^*z^*, G(x)^*z^*\right)
            \\
            &=  \{x_0^*\}+\frechetCod F(x|y)(G(x)^*z^*)+\{([G'(x)\freevar]y)^*z^*\}.
        \end{aligned}
    \end{equation*}
    It follows that
    \begin{equation*}
        \frechetNormal_{\graph(\bar G \circ \bar F)}(x, x, z)
        = \{
            (x^*, -x_0^*, -z^*)
            \mid
            x^*-x_0^* \in \frechetCod F(x|y)(G(x)^*z^*)
            +([G'(x)\freevar]y)^*z^*
        \}.
    \end{equation*}

    By the left-invertibility of $G(x)$, we find that $y$ satisfying $G(x)y=z$ is unique.
    Let $p=(x,w) \in \graph(\bar G \circ \bar F)$ and $w \defeq (x, G(x)y) = Rp$ and observe from the uniqueness of $y$ that $p$ is the unique point in $\graph(\bar G \circ \bar F)$ satisfying $Rp=w$.
    We define
    \begin{equation}
        \label{eq:cofrechet:product:invRp}
        \inv R_p: R \graph(\bar G \circ \bar F) \to \graph(\bar G \circ \bar F),
        \quad
        \inv R_p(\alt x_1, \alt z) \defeq (\alt x_1, (\alt x_1, \alt z)).
    \end{equation}
    This is clearly a Lipschitz inverse selection of $R$ in any neighborhood $U_p \subset R \graph(\bar G \circ \bar F)$ of $w$.
    Therefore, another application of \cref{lemma:cofrechet:cone-linear} yields
    \begin{equation*}
        \frechetNormal_{\graph H}(x, z)
        = \setof{(x^*, -z^*)}{x^* \in \frechetCod F(x|y)(G(x)^*z^*)+\{([G'(x)\freevar]y)^*z^*\} },
    \end{equation*}
    from which the claim follows.
\end{proof}

\begin{corollary}[second-order chain rule for Clarke subdifferentials]
    \label{cor:cofrechet:second-clarke}
    Let $X,Y$ be Banach spaces, $f:Y\to \R$ be locally Lipschitz continuous, and let $S:X\to Y$ be twice differentiable at $x$. Set $h:X\to Y$, $h(x)\defeq f(S(x))$.
    If there exists a neighborhood $U$ of $x \in X$ such that
    \begin{enumerate}
        \item\label{item:cofrechet:second-clarke:i}
        $f$ is Clarke regular at $S(\alt x)$ for all $\alt x\in U$;
        \item\label{item:cofrechet:second-clarke:ii}
        $S'(\alt x)$ has a right-inverse $\rinv{S'(\tilde x)} \in \linear(Y; X)$ for all $\tilde x\in U$;
        \item\label{item:cofrechet:second-clarke:iii}
        the mapping $\alt x \mapsto \rinvstar{S'(\alt x)}$ is continuously differentiable at $x$;
    \end{enumerate}
    then for all $x^*\in \partial_C h(x) = S'(x)^*\partial_C f(S(x))$ we have
    \begin{multline*}
        \frechetCod[\subdiff_C h](x|x^*)(x^{**})
        =
        \{x^{**}[S''(x)\freevar]^*y^*\}\\
        + S'(x)^* \frechetCod[\subdiff_Cf](S(x)|y^*)(S'(x)^{**}x^{**})
        \quad (x^{**}\in X^{**})
    \end{multline*}
    for the unique $y^*\in \partial_C f(S(x))$ with $S'(x)^*y^* = x^*$ and $x^{**}[S''(x)\freevar]^*y^* \in X^*$ defined by
    \begin{equation*}
        X\ni \dir x \mapsto \dualprod{x^{**}}{[S''(x)\dir x]^*y^*}_{X^*}.
 \end{equation*}
\end{corollary}

\begin{proof}
    The expression for $\partial_C h(\alt x)$ for all $\alt x \in U$ follows from \cref{thm:clarke:chain} and the assumption~\ref{item:cofrechet:second-clarke:i}.
    Let now $G:X\to \linear(Y^* ; X^*)$, $G(\alt x) \defeq S'(\alt x)^*$. Then $G$ is Fréchet differentiable at $x$ by the twice differentiability of $S$ at $x$ and has the left-inverse $\rinvstar{S'(\alt x)}$ for all $\alt x \in U$ by assumption~\cref{item:cofrechet:second-clarke:ii}.
    The latter establishes the uniqueness of $y^*$.
    Combining this with assumption~\ref{item:cofrechet:second-clarke:iii}, we can now apply \cref{thm:cofrechet:product} to $F=\subdiff_C f \circ S$ and $G$ to obtain
    \begin{multline*}
        \frechetCod[\subdiff_C h](x|x^*)(x^{**})
        =
        \{([G'(x)\freevar]y^*)^*x^{**}\} \\
        + \frechetCod[(\subdiff_C f) \circ S](x|y^*)(S'(x)^{**}x^{**})
        \quad
        (x^{**}\in X^{**}).
    \end{multline*}
    Observe that, by being twice differentiable at $x$, $S$ is continuously differentiable at $x$.
    Moreover, $S'(x)$ has a bounded right-inverse by \ref{item:cofrechet:second-clarke:ii}.
    We may therefore apply \cref{thm:cofrechet:inner} to $G=\subdiff_C f$ and $F=S$ to obtain
    \begin{equation*}
        \frechetCod [(\subdiff_C f) \circ S](x|y^*)(y^{**})
        =S'(x)^* \frechetCod [\subdiff_C f](S(x)|y^*)(y^{**})
        \qquad
        (y^{**} \in Y^{**}).
    \end{equation*}
    To show the claimed expression for the first term, we proceed as in the proof of \cref{cor:gderiv:second-clarke} to compute for any $\dir x \in X$ that
    \[
        \begin{aligned}
            \dualprod{x^{**}}{[G'(x)\dir x] y^*}_{X^{*}}
            &=\dualprod{x^{**}}{[\lim_{t\downto0}t^{-1}(G(x+t\dir x)-G(x))]y^*}_{X^*}\\
            &=\dualprod{x^{**}}{[\lim_{t\downto0}t^{-1}(S'(x+t\dir x)^*-S'(x)^*)]y^*}_{X^*}\\
            &=\dualprod{x^{**}}{[\lim_{t\downto0}t^{-1}(S'(x+t\dir x)-S'(x))]^*y^*}_{X^*}\\
            &=\dualprod{x^{**}}{[S''(x)\dir x]^*y^*}_{X^*},
        \end{aligned}
    \]
    from which the claim follows.
\end{proof}

\begin{remark}
    In contrast to \cref{cor:gderiv:second-clarke}, the linear functional $x^{**}[S''(x)\freevar]^*y^* \in X^*$ acts on $\alt x\in X$ via the \enquote{first variation} $h_1\in X$ of $(h_1,h_2)\mapsto [S''(x)h_1]h_2\in Y$ instead of the second.
    This is consistent with \cref{ex:graphical:linear,ex:graphical:adjoints:linear,thm:graphical:single}, which show that the coderivative of a differentiable single-valued mapping coincides with the adjoint of the Fréchet derivative.
\end{remark}

\section{Subdifferential calculus}\label{sec:cofrechet:subdiff}

The above results immediately yield calculus rules for the Fréchet subdifferential of \cref{sec:limiting:frechet}.
To see this, we recall from \eqref{eq:graphical:mordukhovich-subdiff} that for $f: X \to \Rbar$ we have that
\[
    \subdiff_F f(x) = D^*[\epi_f](x|f(x))(1),
\]
where $\epi_f:X\setto\R$ denotes the epigraphical mapping of $f$.
We also observe that if $g: X \to \Rbar$ is Fréchet differentiable, then $g'(x)^* \in \linear(\R; X^*)$ and hence
\[
    g'(x)^*z^* = z^* g'(x)
    \qquad\text{for all }
    z^* \in \R.
\]

The following direct corollary of \cref{thm:cofrechet:addition} then yields \cref{thm:limiting:sum} under an additional regularity assumption.
\begin{corollary}[addition of a differentiable mapping]
    \label{cor:cofrechet:addition:subdiff}
    Let $X$ be a Banach space, $g:X\to \R$ be Fréchet differentiable, and $f:X \to \Rbar$.
    If $g$ is Lipschitz at $x\in X$, then
    \begin{equation*}
        \subdiff_F (f+g)(x) = \subdiff_F f(x) + \{g'(x)\}.
    \end{equation*}%
\end{corollary}

We now turn to chain rules.

\begin{corollary}[outer composition with an increasing differentiable mapping]
    \label{cor:cofrechet:outer:subdiff}
    Let $X$ be a Banach space, $f:X\to \R$, and $g:\R \to \Rbar$.
    Let $x \in X$ and $z \defeq  h(x) \defeq g(f(x))$ be given.
    If $g$ is increasing, Fréchet differentiable at $y \defeq f(x)$, and left-invertible on $\range g$ near $z$ such that the left-inverse is Lipschitz at $z$, then
    \begin{equation*}
        \subdiff_F h(x) = g'(y)\subdiff_F f(x).
    \end{equation*}
\end{corollary}

\begin{proof}
    Since $g$ is assumed to be increasing, $g'(y) > 0$ and hence $g'(y)^* z^* = z^* g'(y) = g'(y) > 0$ for $z^* = 1$.
    Due to \cref{cor:graphical:pos-hgen}, we then have
    \[
        D^*f(x|y)( g'(y)) =   g'(y) D^*f(x|y)(1)=  g'(y) \subdiff_F f(x).
    \]
    The claim now follows from \cref{thm:cofrechet:outer}.
\end{proof}

From \cref{thm:cofrechet:inner}, we similarly obtain the following.

\begin{corollary}[inner composition with a differentiable mapping]
    \label{cor:cofrechet:inner:subdiff}
    Let $X,Y$ be Banach spaces, $f: X\to Y$ and $g:Y\to \Rbar$.
    Let $h \defeq g \circ f$ and $x \in \dom h$ be given.
    If $f$ is continuously differentiable near $x$ such that $f'(x)$ has a right-inverse $\rinv{f'(x)} \in \linear(Y; X)$, then
    \begin{equation*}
        \subdiff_F h(x)  = f'(x)^* \subdiff_F g(f(x)).
    \end{equation*}%
\end{corollary}

As a special case, we obtain from \cref{cor:cofrechet:inner:linear} the following linear chain rule (note the slightly different regularity assumption on the inner mapping).

\begin{corollary}[inner composition with a linear operator]
    \label{cor:cofrechet:inner:linear:subdiff}
    Let $X,Y$ be Banach spaces, $A\in \linear(X; Y)$, and $f:Y\to \Rbar$.
    Let $h \defeq f \circ A$ and $x \in \dom h$ be given.
    If $A$ has a right-inverse $\rinv A \in \linear(Y; X)$, then
    \begin{equation*}
        \subdiff_F h(x)  = A^* \subdiff_F f(Ax).
    \end{equation*}%
\end{corollary}

Finally, we obtain a product rule.

\begin{corollary}[product rule]
    \label{cor:cofrechet:product:subdiff}
    Let $X$ be a Banach space, $g:X\to \Rbar$, and $f:X \to \Rbar$.
    Set $h(x)\defeq g(x)f(x)$, and suppose $g$ is Fréchet differentiable at $x \in \dom h$.
    If $g(\alt x) > 0$ for all $\alt x$ near $x$ and the mapping $\alt x \mapsto 1/g(\alt x)$ is Lipschitz at $x$, then
    \begin{equation*}
        \subdiff_F h(x)
        =
        g(x)\subdiff_F f(x)+\{f(x)g'(x)\}.
    \end{equation*}
\end{corollary}

\begin{proof}
    We again use the positivity of $g(x)=g(x)^*z^*$ with $z^*=1$ and \cref{cor:graphical:pos-hgen} to deduce the claim from \cref{thm:cofrechet:product}.
\end{proof}

\chapter{Calculus for the Clarke graphical derivative}
\label{chap:gclarke}

We now turn to the limiting (co)derivatives. Compared to the basic (co)derivatives, calculus rules for these are much more challenging and require even more assumptions. In this chapter, we consider the Clarke graphical derivative, where in addition to strict differentiability we will for the sake of simplicity assume T-regularity of the set-valued mapping (so that the Clarke graphical derivative coincides with the graphical derivative) and show that this regularity is preserved under addition and composition with a single-valued mapping.
We again recall that for $F:X\setto Y$, the Clarke graphical derivative is defined as
\begin{equation*}
    \clarkeGD F(x|y): X \setto Y, \qquad
    \clarkeGD F(x|y)(\dir x) \defeq \setof{ \dir y\in Y}{(\dir x, \dir y) \in \clarkeTangent_{\graph F}(x, y)}.
\end{equation*}

\section{Strict differentiability}

The following concept generalizes the notion of strict differentiability for single-valued mappings (see \cref{rem:strictlydiff}) to set-valued mappings. Let $X,Y$ be Banach spaces.
We say that $F: X \setto Y$ is \term[mapping!differentiable!strictly]{strictly differentiable} at $x \in X$ for $y \in F(x)$ if $\graph F$ is closed near $(x, y)$ and
\begin{subequations}
    \label{eq:gclarke:strictdiff}
    \begin{align}
        &\begin{multlined}[t][0.8\linewidth]
            \text{for every }  \dir y \in \clarkeGD F(x|y)(\dir x),\quad\tau_k \downto 0,\quad \tilde x_k\to x\quad\text{with}\quad\frac{x_k-\tilde x_k}{\tau_k} \to \dir x,\\
            \text{and}\quad \tilde y_k \in F(\tilde x_k) \quad\text{with}\quad \tilde y_k \to y,
        \end{multlined}
        \\
        &\text{there exist } y_k\in F(x_k) \quad\text{with}\quad
        \frac{y_k-\tilde y_k}{\tau_k} \to \dir y.
    \end{align}
\end{subequations}

Compared to semi-differentiability, strict differentiability requires that the limits realizing the various directions are interchangeable with limits of the base points; in other words, that the graphical derivative is itself an inner limit, i.e.,
\begin{equation}\label{eq:gclarke:strictdiff-alt}
    \clarkeGD F(x|y)(\dir x) = \liminf_{\substack{\tau \downto 0,\,  \dir \alt x \to \dir x \\ \graph F \ni (\alt x, \alt y) \to (x,y)}} \frac{F(\alt x+\tau\dir \alt x)-\alt y}{\tau}
    \qquad(\dir x\in X).
\end{equation}
\begin{lemma}
    \label{lemma:gclarke:strictdiff}
    If $X$ and $Y$ are finite-dimensional, then $F: X \setto Y$ is strictly differentiable at $x \in X$ for $y \in F(x)$ if and only if
    \begin{equation}
        \label{eq:gclarke:strictdiff:outer-characterisation}
        \clarkeGD F(x|y)(\dir x) = \liminf_{\substack{\graph F \ni (\alt x, \alt y) \to (x, y), \\ \dir \alt x \to \dir x,\, DF(\alt x|\alt y)(\dir \alt x) \ne \emptyset}}  DF(\alt x|\alt y)(\dir \alt x)
        \qquad(\dir x\in X).
    \end{equation}
\end{lemma}

\begin{proof}
    We need to show that
    \[
        \graph \clarkeGD F(x, y) =
        K \defeq
        \left\{(\dir x, \dir y) \,\middle|\, \dir y \in \liminf_{\substack{\graph F \ni (\alt x, \alt y) \to (x, y), \\ \dir \alt x \to \dir x,\, DF(\alt x|\alt y)(\dir \alt x) \ne \emptyset}}  DF(\alt x|\alt y)(\dir \alt x)\right\}.
    \]
    We first show that $\graph \clarkeGD F(x, y) \subset K$.
    If $(\dir x, \dir y) \not \in K$, then there exist $\graph F \ni (\alt x_k, \alt y_k) \to (x, y)$ and $\dir x_k \to \dir x$ with $DF(\alt x_k|\alt y_k)(\dir x_k) \ne \emptyset$ such that for some $\epsilon>0$ and an infinite subset $N \subset \N$,
    \begin{equation*}
        \inf_{\dir y_k \in DF(\alt x_k|\alt y_k)(\dir x_k)} \norm{\dir y_k-\dir y}_Y \ge 2\epsilon
        \qquad (k \in N).
    \end{equation*}
    By the characterization \eqref{eq:graphical:gderiv-altdef} of $DF(\alt x_k|\alt y_k)$, this implies the existence of $\tau_k \downto 0$ such that
    \begin{equation*}
        \limsup_{k \to \infty} \inf_{y_k \in F(x_k+\tau_k\dir x_k)} \adaptnorm{\frac{y_k-\alt y_k}{\tau_k}-\dir y}_Y \ge \epsilon.
    \end{equation*}
    Thus $(\dir x, \dir y) \not \in \graph \clarkeGD F(x, y)$ and hence $\graph \clarkeGD F(x, y) \subset K$.

    By the definition of the inner limit, we can equivalently write
    \begin{equation*}
        K
        =
        \setof{ (\dir x, \dir y)}{%
            \begin{array}{r}
                (\alt x, \alt y, \dir \alt x) \to (x, y, \dir x)
                \implies
                \exists\, \dir \alt y \to \dir y
                \\
                \text{ with } \dir \alt y \in DF(\alt x|\alt y)(\dir \alt x)
        \end{array}},
    \end{equation*}
    and hence we obtain from the characterization \eqref{eq:graphical:limits:findim:gclarke} of $\clarkeGD F(x, y)$ the converse inclusion $K\subset \graph \clarkeGD F(x, y)$. Therefore \eqref{eq:gclarke:strictdiff:outer-characterisation} holds.
\end{proof}

In particular, single-valued continuously differentiable mappings and their inverses are strictly differentiable.

\begin{lemma}
    \label{lemma:gclarke:regularity:single}
    Let $X,Y$ be Banach spaces and let $F: X \to Y$ be single-valued.
    \begin{enumerate}
        \item
            If $F$ is continuously differentiable at $x\in X$, then $F$ is strictly differentiable at $x$ for $y=F(x)$.
        \item\label{item:gclarke:regularity:single:inv}
            If $F$ is continuously differentiable near $x\in X$ and $F'(x)$ has a right-inverse $\rinv{F'(x)} \in \linear(Y; X)$, then $\inv F$ is strictly differentiable at $y=F(x)$ for $x$.
    \end{enumerate}
\end{lemma}

\begin{proof}
    The proof is analogous to \cref{lemma:gderiv:regularity:single}, since the inverse function theorem (\cref{thm:inversefunctiontheorem}) establishes the continuous differentiability of $\inv{\bar F}$ and hence strict differentiability.
\end{proof}

\begin{remark}
    As in \cref{rem:graphical:regularity:single:d:findim}, if $X$ is finite-dimensional, it suffices in \cref{lemma:gclarke:regularity:single}\,\ref{item:gclarke:regularity:single:inv} to assume that $F$ is continuously differentiable with $\kernel F'(x)^*=\{0\}$.
\end{remark}

\section{Cone transformation formulas}%
\label{sec:gclarke:cones}

The main aim in the following lemmas is to show that tangential regularity is preserved under certain transformations. We do this by proceeding as in \cref{sec:gderiv:cones} to derive explicit expressions for the transformed cones and then comparing them with the corresponding expressions obtained there for the graphical derivative.

\begin{lemma}
    \label{lemma:gclarke:cone-linear}
    Let $X,Y$ be Banach spaces, $C\subset Y$, and $R\in \linear(Y; X)$.
    If $y \in \closure C$ admits a strictly differentiable inverse selection of $R$ at $x=Ry$, then
    \begin{equation*}
        \clarkeTangent_{RC}(x) = R \clarkeTangent_C(y).
    \end{equation*}%
    Moreover, if $C$ is tangentially regular at $y$, then $RC$ is tangentially regular at $x$.
\end{lemma}

\begin{proof}
    We first prove ``$\supset$''.
    Suppose $\dir y \in \clarkeTangent_C(y)$.
    Then for any $C \ni \alt y_k \to y$ there exist $y_k \in C$ and $\tau_k \downto 0$ such that $\dir y=\lim_{k \to \infty} (y_k-\alt y_k)/\tau_k$.
    Consequently, since $R$ is bounded, $R(y_k-\alt y_k)/\tau_k \to R \dir y$.
    To show that $R \dir y \in \clarkeTangent_{RC}(x)$, let $RC \ni \alt x_k \to x$ be given.
    Take now $\alt y_k=\inv R_y(\alt x_k)$, which satisfies $\alt y_k \to y=\inv R_y(x)$ due to $\alt x_k \to x$. Then $(Ry_k-\alt x_k)/\tau_k = R(y_k-\alt y_k)/\tau_k\to R\dir y$.

    To prove ``$\subset$'', let  $\dir x \in \clarkeTangent_{RC}(x)$.
    We will show that $\dir x \in R \clarkeTangent_C(y)$.
    Let $\tau_k \downto 0$.
    We then deduce from the definition \eqref{eq:cones:def-clarketangent} of $\hat T_{RC}(x)$ for $RC \ni x_k \defeq Ry_k \to Ry = x$ the existence of $RC \ni \alt x_k \to x$ such that $(x_k - \alt x_k)/\tau_k \to \dir x$.
    Then $\alt x_k = R\alt y_k$ for $\alt y_k \defeq \inv R_y(\alt x_k)$.
    By the strict differentiability of $\inv R_y$,
    \[
        \frac{y_k - \alt y_k}{\tau_k}
        = \frac{\inv R_y(x_k) - \inv R_y(\alt x_k)}{\tau_k}
        \to [\inv R_y(x)]' \dir x =: \dir y.
    \]
    Since this limit is independent of the sequence $\alt y_k \to y$, we deduce that $\dir y \in \clarkeTangent_C(y)$.
    We also have
    \[
        \frac{\alt x_k - x_k}{\tau_k}
        =
        \frac{R(\alt y_k - y_k)}{\tau_k} \to R\dir y.
    \]
    Since $(\alt x_k - x_k)/\tau_k \to \dir x$ by construction, $R\dir y=\dir x$.

    Finally, comparing now the expression for $\clarkeTangent_{RC}(x)=R\clarkeTangent_C(y)$ with the expression for $T_{RC}(x)$ provided by \cref{lemma:gderiv:cone-linear} and using the tangential regularity of $C$ shows the claimed tangential regularity of $RC$.
\end{proof}

\begin{lemma}[fundamental lemma on compositions]
    \label{lemma:gclarke:fundamental}
    Let $X,Y,Z$ be Banach spaces and
    \begin{equation*}
        C \defeq \{(x, y, z) \mid y \in F(x),\, z \in G(y)\}
    \end{equation*}
    for $F: X \setto Y$, and $G: Y \setto Z$.
    If $(x, y, z) \in C$ and either
    \begin{enumerate}[label=(\alph*)]
        \item\label{item:gclarke:fundamental:G} $G$ is inner semicontinuous, strictly differentiable, and T-regular at $y$ for $z$, or
        \item\label{item:gclarke:fundamental:invF} $\inv F$ is inner semicontinuous, strictly differentiable, and T-regular at $y$ for $x$,
    \end{enumerate}
    then
    \begin{equation}
        \label{eq:gclarke:fundamental}
        \clarkeTangent_C(x, y, z)=\{(\dir x, \dir y, \dir z) \mid \dir y \in \clarkeGD F(x|y)(\dir x),\, \dir z \in \clarkeGD G(y|z)(\dir y)\}.
    \end{equation}%
    Moreover, if $F$ is T-regular at $x$ for $y$ and $G$ is T-regular at $y$ for $z$, then $C$ is tangentially regular at $(x, y, z)$.
\end{lemma}

\begin{proof}
    We only consider the case \ref{item:gclarke:fundamental:G} as the case \ref{item:gclarke:fundamental:invF} is again proved similarly.
    The proof is analogous to \cref{lemma:gderiv:fundamental}, using in this case the strict differentiability of $G$ in place of semi-differentiability.
    First, we observe that $(\dir x, \dir y, \dir z) \in \clarkeTangent_C(x, y, z)$ if and only if for all $\tau_k \downto 0$ and $C \ni (\alt x_k, \alt y_k, \alt z_k) \to (x, y, z)$, there exist $(x_k, y_k, z_k) \in C$ such that
    \begin{equation}
        \label{eq:gclarke:fundamental:limits}
        \dir x=\lim_{k \to \infty} \frac{x_k-\alt x_k}{\tau_k},
        \qquad
        \dir y=\lim_{k \to \infty} \frac{y_k-\alt y_k}{\tau_k},
        \qquad
        \dir z=\lim_{k \to \infty} \frac{z_k-\alt z_k}{\tau_k}.
    \end{equation}

    Suppose $(\dir x, \dir y, \dir z) \in \clarkeTangent_C(x, y, z)$.
    Taking $(\alt x_k, \alt y_k, \alt z_k)=(x, y, z)$, it is immediate that $\dir y \in DF(x|y)(\dir x)$ and $\dir z \in DG(y|z)(\dir y)$. By the T-regularity of $G$, it follows that $\dir z \in \clarkeGD G(y|z)(\dir y)$.
    Now take any $\graph F \ni (\alt x_k, \alt y_k) \to (x, y)$.
    By the assumption that $G$ is inner semicontinuous, there exists some $G(\alt y_k) \ni \alt z_k \to z$. Thus, by the above characterization of $(\dir x, \dir y, \dir z) \in \clarkeTangent_C(x, y, z)$, there exist $(x_k, y_k) \in \graph F$ such that $(x_k-\alt x_k)/\tau_k \to \dir x$ and $(y_k-\alt y_k)/\tau_k \to \dir y$, i.e., $(\dir x, \dir y) \in \clarkeTangent_{\graph F}(x, y)$.
    This shows ``$\subset$'' in  \eqref{eq:gclarke:fundamental}.

    To prove ``$\supset$'', suppose $\dir y \in \clarkeGD F(x|y)(\dir x)$ and $\dir z \in \clarkeGD G(y|z)(\dir y)$ and take $\tau_k \downto 0$ and $C \ni (\alt x_k, \alt y_k, \alt z_k) \to (x, y, z)$.
    By definition of $\clarkeGD F(x|y)$, there then exist $(x_k, y_k) \in \graph F$ such that the first two limits in \eqref{eq:gclarke:fundamental:limits} hold.
    By the strict differentiability of $G$ at $y$ for $z$, we can also find $z_k \in G(y_k)$ such that $(z_k-\alt z_k)/\tau_k \to \dir z$.
    This shows the remaining limit.

    Finally, the tangential regularity of $C$ follows from the assumed T-regularities of $F$ and $G$ by comparing \eqref{eq:gclarke:fundamental} with the corresponding expression \eqref{eq:gderiv:fundamental}.
\end{proof}

If one of the two mappings is single-valued, we can use \cref{lemma:gclarke:regularity:single} for verifying its semi-differentiability and \cref{thm:graphical:single} for the regularity and the expression of its graphical derivative to obtain from \cref{lemma:gclarke:fundamental} the following two special cases.

\begin{corollary}[fundamental lemma on compositions: single-valued outer mapping]
    \label{lemma:gclarke:fundamental:single-outer}
    Let $X,Y,Z$ be Banach spaces and
    \begin{equation*}
        C \defeq \{(x, y, G(y)) \mid y \in F(x)\}
    \end{equation*}
    for $F: X \setto Y$ and $G: Y \to Z$.
    If $(x, y, z) \in C$ and $G$ is continuously differentiable at $y$, then
    \begin{equation*}
        \clarkeTangent_C(x, y, z)=\{(\dir x, \dir y, G'(y)\dir y) \mid \dir y \in \clarkeGD F(x|y)(\dir x)\}.
    \end{equation*}
    Moreover, if $F$ is T-regular at $(x, y)$, then $C$ is tangentially regular at $(x, y, G(y))$.
\end{corollary}

\begin{corollary}[fundamental lemma on compositions: single-valued inner mapping]
    \label{lemma:gclarke:fundamental:single-inner}
    Let $X,Y,Z$ be Banach spaces and
    \begin{equation*}
        C \defeq \{(x, y, z) \mid y=F(x),\, z \in G(y)\}
    \end{equation*}
    for $F: X \to Y$ and $G: Y \setto Z$.
    If $(x, y, z) \in C$, $F$ is continuously differentiable at $x$,
    and $F'(x)$ has a right-inverse $\rinv{F'(x)} \in \linear(Y; X)$,
    then
    \begin{equation*}
        \clarkeTangent_C(x, y, z)=\{(\dir x, \dir y, \dir z) \mid \dir y=F'(x)\dir x,\, \dir z \in \clarkeGD G(y|z)(\dir y)\}.
    \end{equation*}%
    Moreover, if $G$ is T-regular at $(y, z)$, then $C$ is tangentially regular at $(x, y, z)$.
\end{corollary}

\section{Calculus rules}\label{sec:gclarke:calculus}

Using these lemmas, we again obtain calculus rules under the assumption that the involved set-valued mapping is regular.
The proofs are again similar to those in \cref{sec:gderiv:calculus}, and we only note the differences.
Overall, due to the changes in assumptions of the fundamental single-valued composition results (\cref{lemma:gclarke:fundamental:single-outer,lemma:gclarke:fundamental:single-inner} instead of \cref{lemma:gderiv:fundamental:single-outer,lemma:gderiv:fundamental:single-inner}), assumptions of Fréchet differentiability at a point $x$ generally need to be strengthened to continuous differentiability at $x$. In results that only use the cone transformation results (\cref{lemma:gclarke:cone-linear} instead of \cref{lemma:gderiv:cone-linear}), the Fréchet differentiability needs to be strengthened to strict differentiability (see \cref{rem:strictlydiff}).
As the latter is implied by continuous differentiability, most of our results simply assume continuous differentiability.

\begin{theorem}[addition of a single-valued differentiable mapping]
    \label{thm:gclarke:addition}
    Let $X,Y$ be Banach spaces, $G:X\to Y$ be continuously differentiable at $x \in X$, and $F:X\setto Y$.
    Then for any $y \in H(x)\defeq F(x) + G(x)$, we have
    \begin{equation*}
        \clarkeGD H(x|y)(\dir x) = \clarkeGD F(x|y-G(x))(\dir x) + \{G'(x) \dir x\}
        \qquad (\dir x \in X).
    \end{equation*}
    Moreover, if $F$ is T-regular at $(x, y-G(x))$, then $H$ is T-regular at $(x, y)$.
\end{theorem}

\begin{proof}
    The proof is the same as that of \cref{thm:gderiv:addition}, using \cref{lemma:gclarke:cone-linear} (which requires continuous differentiability of $G$) in place of \cref{lemma:gderiv:cone-linear}.

    When $F$ is T-regular, the regularity of $H$ follows by comparing our claim to that of \cref{thm:gderiv:addition}.
\end{proof}

\begin{theorem}[outer composition with a single-valued differentiable mapping]
    \label{thm:gclarke:outer}
    Let $X,Y$ be Banach spaces, $F:X\setto Y$, and $G:Y\to Z$.
    Let $x \in X$ and $z \in H(x) \defeq G(F(x))$ be given.
    If $G$ is continuously differentiable at $y \in F(x) \isect \inv G(\{z\})$ and locally left-invertible near $z$ such that the local left-inverse $\linv G$ is strictly differentiable at $z$, then
    \begin{equation*}
        \clarkeGD H(x|z)(\dir x) =  G'(y) \clarkeGD F(x|y)(\dir x)
        \qquad (\dir x \in X).
    \end{equation*}
    Moreover, if $F$ is T-regular at $(x, y)$, then $H$ is T-regular at $(x, z)$.
\end{theorem}

\begin{proof}
    The proof is the same as that of \cref{thm:gderiv:outer}, using \cref{lemma:gclarke:cone-linear,lemma:gclarke:fundamental:single-outer} in place of \cref{lemma:gderiv:cone-linear,lemma:gderiv:fundamental:single-outer}.
    Here, \cref{lemma:gclarke:fundamental:single-outer} requires the continuous differentiability of $G$, while strict differentiability of $\inv G$ implies the strict differentiability of $\inv R_p$ required by \cref{lemma:gclarke:cone-linear}.

    When $F$ is T-regular, the regularity of $H$ follows by comparing the claim to that of \cref{thm:gderiv:outer}.
\end{proof}

The special case for a linear operator follows from this exactly as \cref{cor:gderiv:outer:linear} did from \cref{thm:gderiv:outer}.

\begin{corollary}[outer composition with a linear operator]
    \label{cor:gclarke:outer:linear}
    Let $X,Y,Z$ be Banach spaces, $A\in \linear(Y; Z)$, and $F:X\setto Y$.
    If $A$ has a left-inverse $\linv A$, then for any $x\in X$ and $z\in H(x)\defeq AF(x)$,
    \begin{equation*}
        \clarkeGD H(x|z)(\dir x)=
        A \clarkeGD F(x|y)(\dir x)
        \qquad (\dir x \in X)
    \end{equation*}
    for the unique $y\in Y$ with $Ay=z$.

    Moreover, if $F$ is T-regular at $(x,y)$, then $H$ is T-regular at $(x, z)$.
\end{corollary}

For inner composition, there are no changes in assumptions as we already required continuous differentiability for \cref{thm:gderiv:inner}.

\begin{theorem}[inner composition with a single-valued differentiable mapping]
    \label{thm:gclarke:inner}
    Let $X,Y,Z$ be Banach spaces, $F: X\to Y$ and $G:Y\setto Z$.
    Let $x \in X$ and $z \in H(x)\defeq G(F(x))$ be given.
    If $F$ is continuously differentiable at $x$ and $F'(x)$ has a right-inverse $\rinv{F'(x)} \in \linear(Y; X)$, then
    \begin{equation*}
        \clarkeGD H(x|z)(\dir x)  = \clarkeGD G(F(x)|z)(F'(x)\dir x)
        \qquad (\dir x \in X).
    \end{equation*}%
    Moreover, if $G$ is T-regular at $(F(x), z)$, then $H$ is T-regular at $(x, z)$.
\end{theorem}

\begin{proof}
    The proof is the same as that of \cref{thm:gderiv:inner}, using \cref{lemma:gclarke:cone-linear,lemma:gclarke:fundamental:single-inner} in place of \cref{lemma:gderiv:cone-linear,lemma:gderiv:fundamental:single-inner}.
    When $G$ is T-regular, the regularity of $H$ follows by comparing the claim to that of \cref{thm:gderiv:inner}.
\end{proof}

\begin{corollary}[inner composition with a linear operator]%
    Let $X,Y,Z$ be Banach spaces, $A\in \linear(X; Y)$, and $G:Y\setto Z$.
    Let $x \in X$ and $z \in H(x) \defeq G(Ax)$ be given.
    If $A$ has a right-inverse $\rinv A \in \linear(Y; X)$, then
    \begin{equation*}
        \clarkeGD H(x|z)(\dir x)  = \clarkeGD G(Ax|z)(A\dir x)
        \qquad (\dir x \in X).
    \end{equation*}%
    Moreover, if $G$ is T-regular at $(Ax,z)$, then $H$ is T-regular at $(x, z)$.
\end{corollary}

As in \cref{sec:gderiv:calculus}, we can apply these results to chain rules for subdifferentials, this time only at points where these subdifferentials are T-regular.
\begin{corollary}[second-order chain rule for convex subdifferentials]
    \label{cor:gclarke:second-convex}
    Let $X,Y$ be Banach spaces, $f:Y\to \Rbar$ be proper, convex, and lower semicontinuous, and $A\in\linear(X; Y)$ be such that $A$ has a right-inverse $\rinv{A}  \in \linear(Y; X)$ and that $\range A \isect \interior \dom f \ne \emptyset$. Let $h\defeq f\circ A$.
    Then for any $x\in X$ and $x^*\in \partial h(x) = A^*\partial f(Ax)$,
    \begin{equation*}
        \clarkeGD [\subdiff h](x|x^*)(\dir x)  = A^* \clarkeGD [\subdiff f](Ax|y^*)(A\dir x)
        \qquad (\dir x \in X)
    \end{equation*}
    for the unique $y^* \in Y^*$ satisfying $A^*y^*=x^*$.
    Moreover, if $\partial f$ is T-regular at $Ax$ for $y^*$, then $\partial h$ is T-regular at $x$ for $x^*$.
\end{corollary}

For the second-order chain rule for Clarke subdifferentials, we again use a product rule.
\begin{theorem}[product rule]
    \label{thm:gclarke:product}
    Let $X,Y,Z$ be Banach spaces, let $G:X\to \linear (Y; Z)$, and $F:X\setto Y$.
    If $G$ is continuously differentiable at $x\in X$, $G(\alt x) \in \linear(Y; Z)$ has a left-inverse $\linv{G(\alt x)} \in \linear(Z; Y)$ for all $\alt x$ near $x$ and the mapping $\alt x \mapsto \linv{G(\alt x)}$ is strictly differentiable at $x$, then for all $z\in H(x)\defeq G(x)F(x) \defeq \Union_{y \in F(x)} G(x)y$,
    \begin{equation*}
        \clarkeGD H(x|z)(\dir x)
        =
        \{[G'(x)\dir x]y\}+G(x)\clarkeGD F(x|y)\dir x
        \qquad (\dir x \in X)
    \end{equation*}
    for the unique $y \in F(x)$ satisfying $G(x)y=z$.

    Moreover, if $F$ is T-regular at $x$ for $y$, then $H$ is T-regular at $x$ for $z$.
\end{theorem}

\begin{proof}
    The proof is the same as that of \cref{thm:gderiv:product}, using \cref{lemma:gclarke:cone-linear,thm:gclarke:outer} in place of \cref{lemma:gderiv:cone-linear,thm:gderiv:outer}, and changing \cref{lemma:gderiv:cartesian-product} appropriately.
    Here, to use \cref{thm:gclarke:outer}, we need $\bar G$ to be continuously differentiable at $(x, y)$, which follows from the corresponding property of $G$ at $x$.
    Likewise, to use \cref{lemma:gclarke:cone-linear}, we need the left-inverse $\linv{\bar G}$ of $\bar G$ constructed in its proof to be strictly differentiable at $(x, z)$.
    This follows from the strict differentiability of $\alt x \mapsto \linv{G(\alt x)}$ at $x$.

    When $F$ is T-regular, the regularity of $H$ follows by comparing the claim to that of \cref{thm:gderiv:product}.
\end{proof}

Similarly to \cref{cor:gderiv:second-clarke}, we now obtain from \cref{thm:gclarke:product} a second-order chain rule for the Clarke subdifferential. Note that for $f:X\to\R$, we have $\partial_C f:X\setto X^*$ and hence $D[\partial_C f](x|x^*):X\setto X^*$ as well.

\begin{corollary}[second-order chain rule for Clarke subdifferentials]
    \label{cor:gclarke:second-clarke}
    Let $X,Y$ be Banach spaces, $f:Y\to \R$ be locally Lipschitz continuous, and $S:X\to Y$ be twice continuously differentiable at $x \in X$. Set $h:X\to \R$, $h \defeq f \circ S$.
    If there exists a neighborhood $U$ of $x$ such that
    \begin{enumerate}
        \item\label{item:gclarke:second-clarke:i} $f$ is Clarke regular at $S(\alt x)$ for all $\alt x\in U$;
        \item\label{item:gclarke:second-clarke:ii} $S'(\alt x)$ has a right-inverse $\rinv{S'(\tilde x)} \in \linear(Y; X)$ for all $\tilde x\in U$;
        \item\label{item:gclarke:second-clarke:iii} the mapping $\alt x \mapsto \rinvstar{S'(\alt x)}$ is strictly differentiable at $x$;
    \end{enumerate}
    then for all $x^*\in \partial_C h(x) = S'(x)^*\partial_C f(S(x))$ we have
    \begin{equation*}
        \clarkeGD[\subdiff_C h](x|x^*)(\dir x)  = y^*S''(x)\dir x + S'(x)^* \clarkeGD[\subdiff_C f](S(x)|y^*)(S'(x)\dir x)
        \quad (\dir x \in X)
    \end{equation*}
    for the unique $y^*\in \partial_C f(S(x))$ with $S'(x)^*y^* = x^*$ and $y^*S''(x)\dir x\in X^*$ defined as in \cref{cor:gderiv:second-clarke}.

    Moreover, if $\partial_C f$ is T-regular at $S(x)$ for $y^*$, then $\partial_C h$ is T-regular at $x$ for $x^*$.
\end{corollary}

\begin{proof}
    In the proof of \cref{cor:gderiv:second-clarke}, observe that the twice continuous differentiability of $S$ at $x$ guarantees the continuous (and not merely Fréchet) differentiability of $G$ at $x$.
    We have also appropriately strengthened the assumption~\cref{item:gclarke:second-clarke:iii} to impose strict differentiability.
\end{proof}

\begin{remark}
    Unlike the Fréchet subdifferential in \cref{sec:cofrechet:subdiff}, we cannot obtain calculus rules for the Clarke subdifferential of extended-real-valued functionals as direct corollaries of the above calculus rules for the Clarke graphical derivative.
    This is due to the fact the Clarke subdifferential is actually characterized by the (somewhat uncommon) Clarke \emph{normal} cone,  see \cref{lemma:graphical:clarke-subdiff}, while the Clarke graphical derivative is defined through the (more common) Clarke \emph{tangent} cone.
    (Calculus rules for the Clarke subdifferential of \emph{real-valued} functionals, however, can be derived directly using the generalized directional derivative; see \cref{chap:clarke}.)
\end{remark}

\chapter{Calculus for the limiting coderivative}
\label{chap:colimiting}

The limiting or Mordukhovich coderivative is the most challenging of all the graphical derivatives and coderivatives, and developing exact calculus rules for it requires the most assumptions. In particular, we will here assume a stronger variant of the assumptions of \cref{chap:cofrechet} for the Fréchet coderivative that also implies N-regularity of the set-valued mapping so that we can exploit the stronger properties of the Fréchet coderivative. To prove the fundamental composition lemmas, we will also need to introduce the concept of \term[compactness!sequential normal!partial]{partial sequential normal compactness} that will be used to prevent certain unit-length coderivatives from converging weakly-$*$ to zero. This concept will also be needed in \cref{chap:regularity}.
We again recall for convenience the definition of the limiting coderivative for $F:X\setto Y$ as
\begin{equation*}
    \coderivative F(x|y): Y^* \setto X^*,\qquad
    \coderivative F(x|y)(y^*) \defeq \setof{ x^* \in X^*}{(x^*, -y^*) \in N_{\graph F}(x, y)}.
\end{equation*}

\section{Strict codifferentiability}

Let $X,Y$ be Banach spaces.
We say that $F$ is \term[mapping!codifferentiable!strictly]{strictly codifferentiable} at $x \in X$ for $y \in F(x)$ if
\begin{equation}
    \label{eq:colimiting:strictcodiff}
    \coderivative F(x|y)(y^*) =
    \left\{ x^* \in X^*
        \,\middle|\,
        \begin{array}{l}
            \forall \graph F \ni (x_k, y_k) \to (x, y),\, \epsilon_k \downto 0:
            \\
            \exists (x_k^*, y_k^*) \weaktostar (x^*, y^*)
            \text{ with } x_k^* \in \frechetCod_{\epsilon_k} F(x_k|y_k)(y_k^*)
        \end{array}
    \right\},
\end{equation}
i.e., if \eqref{eq:cones:def-limnormal} is a full weak-$*$-limit.
From \cref{thm:graphical:single,thm:graphical:single:inverse}, it is clear that single-valued continuously differentiable mappings and their inverses are strictly codifferentiable.

\begin{lemma}
    \label{lemma:colimiting:regularity:single}
    Let $X,Y$ be Banach spaces, $F: X \to Y$, $x\in X$, and $y=F(x)$.
    \begin{enumerate}
        \item If $F$ is continuously differentiable at $x$, then $F$ is strictly codifferentiable at $x$ for $y$.
        \item If $F$ is continuously differentiable near $x$, then $\inv F$ is strictly codifferentiable at $y$ for $x$.
    \end{enumerate}
\end{lemma}

The next lemma and counterexample demonstrate that strict codifferentiability is a stronger assumption than N-regularity.

\begin{lemma}
    \label{lemma:colimiting:strictcodiff:gr}
    Let $X,Y$ be Banach spaces and let $F: X \setto Y$ be strictly codifferentiable at $x$ for $y$.
    Then $F$ is N-regular at $x$ for $y$.
\end{lemma}

\begin{proof}
    By \cref{thm:cones:inclusions}, strict codifferentiability, and the definition of the inner limit, respectively,
    \begin{equation*}
        \begin{aligned}
            \frechetNormal_{\graph F}(x, y)
            &
            \subset
            N_{\graph F}(x, y)
            \\
            &
            =
            \liminf_{\graph F \ni (\alt x,\alt y) \to (x, y),\,\epsilon \downto 0} \frechetNormal_{\graph F}^{\epsilon}(\alt x, \alt y)
            \\
            &
            \subset
            \frechetNormal_{\graph F}(x, y).
        \end{aligned}
    \end{equation*}
    Therefore $N_{\graph F}(x, y)=\frechetNormal_{\graph F}(x, y)$, i.e., $\graph F$ is normally regular at $(x, y)$.
\end{proof}

\begin{example}[graphical regularity does not imply strict codifferentiability]
    Consider $F(x) \defeq [\abs{x}, \infty)$, $x\in \R$.
    Then $\graph F=\epi \abs{\freevar}$ is a convex set and therefore graphically regular at all points with
    \begin{equation*}
        N_{\graph F}(x, \abs{x})=
        \begin{cases}
            (\sign x, -1)[0, \infty) & \text{if } x\neq 0,\\
            \polar{\graph F}=\{(x^*, y^*) \mid -y^* \ge \abs{x^*} \} &\text{if }x=0.
        \end{cases}
    \end{equation*}
    Hence $N_{\graph F}$ is not continuous and therefore, \emph{a fortiori}, $F$ is not strictly codifferentiable at $(0, 0)$.
\end{example}

\section{Partial sequential normal compactness}
\label{sec:colimiting:psnc}

One central difficulty in working with infinite-dimensional spaces is the need to distinguish weak-$*$ convergence and strong convergence. In particular, we need to prevent certain sequences whose norm is bounded away from zero from weak-$*$ converging to zero. As we cannot guarantee this in general, we need to add this as an assumption. In our specific setting, this is the \term[compactness!sequential normal!partial]{partial sequential normal compactness} (\term[PSNC|see{compactness, partial sequential normal}]{PSNC}) of $G: Y \setto Z$ at $y$ for $z$, which holds if
\begin{multline}
    \label{eq:colimiting:psnc}
    \epsilon_k \downto 0,\ (y_k, z_k) \to (y, z),\ y_k^* \weaktostar 0,\ \norm{z_k^*}_{Z^*} \to 0,\, \text{ and } y_k^* \in \frechetCod_{\epsilon_k} G(y_k|z_k)(z_k^*)
    \\
    \implies \norm{y_k^*}_{Y^*} \to 0.
\end{multline}
Obviously, if $Y^*$ is finite-dimensional, then every mapping $G: Y \setto Z$ is PSNC.
To prove the PSNC property of single-valued mappings and their inverses, we will need an estimate of $\epsilon$-coderivatives.

\begin{lemma}
    \label{lemma:colimiting:single-valued-epsilon}
    Let $X,Y$ be Banach spaces and let $F:X \to Y$ be Fréchet differentiable at $x \in X$. Then for any $\epsilon>0$, $L \defeq \norm{F'(x)}_{\linear(X; Y)}$, and $y=F(x)$,
    \begin{equation*}
        \frechetCod_\epsilon F(x|y)(y^*) \subset \B(F'(x)^*y^*, (L+1)\epsilon)
        \qquad (y^* \in Y^*).
    \end{equation*}
\end{lemma}

\begin{proof}
    By definition, $x^* \in \frechetCod_\epsilon F(x|y)(y^*)$ if and only if for every sequence $x_k \to x$,
    \begin{equation}
        \label{eq:colimiting:epsilon-single-valued-0}
        \limsup_{k \to \infty} \frac{\dualprod{x^*}{x_k-x}_{X} - \dualprod{y^*}{F(x_k)-F(x)}_{Y}}{\sqrt{\norm{x_k-x}_X^2+\norm{F(x_k)-F(x)}_Y^2}} \le \epsilon.
    \end{equation}
    Let $\ell > L$. Then by the Fréchet differentiability and therefore the Lipschitz continuity of $F$ \emph{at} $x$ (see \cref{rem:variation:frechet-lipschitz-at})
    we have $\norm{F(x_k)-F(x)}_Y \le \ell\norm{x_k-x}_X$ for large enough $k$ and therefore
    \begin{equation*}
        \limsup_{k \to \infty} \frac{\dualprod{x^*}{x_k-x}_{X} - \dualprod{y^*}{F(x_k)-F(x)}_{Y}}{\norm{x_k-x}_X} \le \epsilon(\ell+1).
    \end{equation*}
    Furthermore, the Fréchet differentiability of $F$ implies that
    \begin{equation*}
        \limsup_{k \to \infty} \frac{\dualprod{F'(x)^*y^*}{x_k-x}_{X} - \dualprod{y^*}{F(x_k)-F(x)}_{Y}}{\norm{x_k-x}_X} = 0
    \end{equation*}
    and hence that
    \begin{equation*}
        \limsup_{k \to \infty} \frac{\dualprod{x^*-F'(x)^*y^*}{x_k-x}_{X}}{\norm{x_k-x}_X} \le \epsilon(\ell+1).
    \end{equation*}
    Since $x_k \to x$ was arbitrary, this implies $\norm{x^*-F'(x)^*y^*}_{X^*} \le \epsilon(\ell+1)$, and since $\ell>L$ was arbitrary, the claim follows.
\end{proof}

\begin{lemma}
    \label{lemma:colimiting:psnc:single-valued}
    Let $Y,Z$ be Banach spaces and $G:Y\to Z$. If either
    \begin{enumerate}[label=(\alph*)]
        \item\label{item:colimiting:psnc:single-valued:cont}
        $G$ is continuously differentiable near $y \in Y$ or
        \item\label{item:colimiting:psnc:single-valued:finite}
        $Y^*$ is finite-dimensional,
    \end{enumerate}
    then $G$ is PSNC at $y$ for $z=G(y)$.
\end{lemma}

\begin{proof}
    The finite-dimensional case \cref{item:colimiting:psnc:single-valued:finite} is clear from the definition \eqref{eq:colimiting:psnc} of the PSNC property.

    For case \cref{item:colimiting:psnc:single-valued:cont}, we have from \cref{lemma:colimiting:single-valued-epsilon} that $\frechetCod_{\epsilon_k} G(y_k|z_k)(z_k^*) \subset \B(G'(y_k)^*z_k^*, \ell\epsilon_k)$ for any $\ell > \norm{G'(y_k)}_{\linear(Y; Z)}$. By the continuous differentiability of $G$, this will hold for $\ell > \norm{G'(y)}_{\linear(Y; Z)}$ and any $k \in \N$ large enough. Thus there exist $d_k^* \in \B(0, \ell \epsilon_k)$ such that
    \begin{equation*}
        y_k^*
        =G'(y_k)^*z_k^*+d_k^*
        =G'(y)^*z_k^*+[G'(y_k)-G'(y)]^*z_k^* + d_k^* \to 0
    \end{equation*}
    since $d_k^* \to 0$ (due to $\epsilon_k \downto 0$),  $\norm{z_k^*}_{Z^*}\to 0$, $y_k \to y$, and $G$ is continuously differentiable near $y$.
\end{proof}

\begin{remark}
    The assumption in \cref{lemma:colimiting:psnc:single-valued}\,\cref{item:colimiting:psnc:single-valued:cont} can be relaxed to $G$ being Fréchet differentiable and Lipschitz near $y$, which allows a similar relaxation of assumptions in the calculus rules of \cref{sec:colimiting:calculus}.
\end{remark}

\begin{lemma}
    \label{lemma:colimiting:psnc:single:inverse}
    Let $Y,Z$ be Banach spaces and $G:Y\to Z$. If either
    \begin{enumerate}[label=(\alph*)]
        \item\label{item:colimiting:psnc:single:inverse:cont}
        $G$ is continuously differentiable near $y \in Y$ and $G'(y) \in \linear(Y; Z)$ has a right-inverse $\rinv{G'(y)} \in \linear(Z; Y)$, or
        \item\label{item:colimiting:psnc:single:inverse:finite}
        $Z^*$ is finite-dimensional,
    \end{enumerate}
    then $\inv G$ is PSNC at $z=G(y)$ for $y$.
\end{lemma}

\begin{proof}
    The finite-dimensional case \cref{item:colimiting:psnc:single:inverse:finite} is clear from the definition \eqref{eq:colimiting:psnc} of the PSNC property.

    For case \cref{item:colimiting:psnc:single:inverse:cont}, we have from the definition of $\frechetCod_\epsilon G$ via $\frechetNormal^\epsilon_{\graph G}$ that $\dir z_k^* \in \frechetCod_\epsilon \inv G(z_k|y_k)(\dir y_k^*)$ if and only if $\dir y_k^* \in \frechetCod_\epsilon G(y_k|z_k)(\dir z_k^*)$. We thus have to show that
    \begin{multline*}
        \epsilon_k \downto 0,\ (y_k, z_k) \to (y, z),\ z_k^* \weaktostar 0,\ \norm{y_k^*}_{Y^*} \to 0,\, \text{ and } y_k^* \in \frechetCod_{\epsilon_k} G(y_k|z_k)(z_k^*)
        \\
        \implies \norm{z_k^*}_{Z^*} \to 0.
    \end{multline*}
    It follows from \cref{lemma:colimiting:single-valued-epsilon} that $\frechetCod_{\epsilon_k} G(y_k|z_k)(z_k^*) \subset \B(G'(y_k)^*z_k^*, \ell\epsilon_k)$ for any choice of $\ell > \norm{G'(y_k)}_{\linear(Y; Z)}$.
    As in \cref{lemma:colimiting:psnc:single-valued}, we now deduce that
    $y_k^*=G'(y_k)^*z_k^*+d_k^*$ for some $d_k^* \in \B(0, \ell \epsilon_k)$.
    In this case, both $y_k^*\to 0$ (by assumption) and $d_k^*\to 0$ (due to $\epsilon_k \to 0$).
    Furthermore, $\{z_k^*\}_{k \in \N}$ is bounded as a weakly-$*$ converging sequence. By $y_k\to y$ and the continuous differentiability of $G$, we thus have that
    \begin{equation*}
        G'(y)^*z_k^* = G'(y_k)^*z_k^* - [G'(y_k)-G'(y)]^*z_k^* = y_k^*-d_k^* - [G'(y_k)-G'(y)]^*z_k^* \to 0.
    \end{equation*}
    Since $G'(y)$ is assumed to have a right-inverse, $G'(y)^*$ has a left-inverse, which implies $z_k^* \to 0$ as required.
\end{proof}

We will use PSNC to obtain the following partial compactness property for the limiting coderivative, for which we need to assume reflexivity (or finite-dimensionality) of $Y$.

\begin{lemma}%
    \label{lemma:colimiting:compactness}
    Let $Y,Z$ be Banach spaces and $G:Y\setto Z$.
    Let $y\in Y$ and $z\in G(y)$ be given. Assume $y^* \in \coderivative G(y|z)(0)$ implies $y^*=0$ and
    either
    \begin{enumerate}[label=(\alph*)]
        \item $Y$ is finite-dimensional or
        \item $Y$ is reflexive and $G$ is PSNC at $y$ for $z$.
    \end{enumerate}
    If
    \begin{equation*}
        (y_k, z_k) \to (y, z),
        \quad
        z_k^* \weaktostar z^*,
        \quad
        \alt\epsilon_k \downto 0,
        \quad \text{and} \quad
        \bar y_k^* \in \frechetCod_{\alt\epsilon_k}
        G(y_k|z_k)(z_k^*),
    \end{equation*}
    then there exists a subsequence such that $\bar y_k^* \weaktostar \bar y^* \in \coderivative G(y|z)(z^*)$.
\end{lemma}
\begin{proof}
    We first show that $\{\bar y_k^*\}_{k\in\N}$ is bounded. We argue by contradiction and suppose that $\{\bar y_k^*\}_{k \in \N}$ is unbounded.
    We may then assume that $\norm{\bar y_k^*}_{Y^*} \to \infty$ by switching to a (not relabelled) subsequence.
    Since $\frechetCod_{\alt\epsilon_k} G(y_k|z_k)$ is formed from a cone, we also have
    \begin{equation*}
        \B_{Y^*}\ni \bar y_k^*/\norm{\bar y_k^*}_{Y^*} \in \frechetCod_{\alt\epsilon_k} G(y_k|z_k)(z^*_k/\norm{\bar y_k^*}_{Y^*}).
    \end{equation*}
    Observe that $\norm{z^*_k/\norm{\bar y_k^*}_{Y^*}}_{Z^*} \to 0$ because $\{z_k^*\}_{k \in \N}$ is bounded. Since $Y$ is reflexive, we can use the Eberlein--\u{S}mulyan theorem (\cref{thm:ebsmul}) to extract a subsequence such that $\bar y_k^*/\norm{\bar y_k^*}_{Y^*} \weaktostar \bar y^*$ for some $\bar y^* \in \coderivative G(y|z)(0)$.
    If $Y$ is finite-dimensional, clearly $\bar y^* \ne 0$.
    Otherwise we need to use the assumed PSNC property.
    If $\bar y^* = 0$, then \eqref{eq:colimiting:psnc} implies that $1=\norm{\bar y_k^*/\norm{\bar y_k^*}_{Y^*}}_{Y^*} \to 0$, which is a contradiction.
    Therefore $\bar y^* \ne 0$.
    However, $\bar y^* \in D^*G(y|z)(0)$ implies $\bar y^*=0$ by assumption and hence we obtain a contradiction.

    Therefore $\{\bar y_k^*\}_{k \in \N}$ is bounded, and thus we may again use \cref{thm:ebsmul} to extract a subsequence converging to some $\bar y^* \in Y^*$.
    By the definition of the limiting coderivative, this implies $\bar y^* \in \coderivative G(y|z)(z^*)$ and hence the claim.
\end{proof}

\begin{remark}
    The PSNC property, its stronger variant \term[compactness!sequential normal]{sequential normal compactness} (\term[SNC|see{compactness, sequential normal}]{SNC}),
    and their implications are studied in significant detail in \cite{Mordukhovich:2006}.
\end{remark}

\section{Cone transformation formulas}

As in \cref{sec:gclarke:cones}, we now show that normal regularity is preserved under certain transformations by deriving explicit expressions for the transformed cones and then comparing them with the corresponding expressions of the Fréchet coderivative.
\begin{lemma}
    \label{lemma:colimiting:cone-linear}
    Let $X,Y$ be Banach spaces, $C\subset Y$, and $R\in \linear(Y; X)$.
    If $y \in \closure C$ admits a Lipschitz inverse selection of $R$ at $x=Ry$ with factor $L_x$ and $C$ is normally regular at $y$, then $RC$ is normally regular at $x$ and
    \begin{equation*}
        N_{RC}(x) = \{x^* \in X^* \mid R^* x^* \in N_C(y)\}.
    \end{equation*}

    If $C$ is not normally regular, this identity also holds if for any $x^* \in X^*$ and $y_k^* \weaktostar R^*x^*$ there exist $x_k^* \in X^*$ with $\norm{R^*x_k^*-y_k^*}_{Y^*} \downto 0$ and $x_k^* \weaktostar x^*$.
\end{lemma}

\begin{proof}
    We first prove ``$\subset$''.
    Recall that $U_y$ denotes the neighborhood of $x$ where $\inv R_y$ is defined.
    Let $x^* \in N_{RC}(x)$. By definition, this holds if and only if there exist $\epsilon_k \downto 0$ as well as $x_k^* \weaktostar x^*$ and $RC \ni x_k \to x$ with $x_k^* \in \frechetNormal_{RC}^{\epsilon_k}(x_k)$.
    Defining $y_k \defeq \inv R_y(x_k)$, we have $Ry_k=x_k$ and $C \ni y_k \to y$. Thus, for large enough $k$ that $x_k \in U_y$, \cref{lemma:cofrechet:cone-linear} yields $R^*x_k^* \in \frechetNormal_C^{\epsilon_k L_x}(y_k)$.
    We have $R^*x_k^* \weaktostar R^*x^*$ by the continuity of $R$.
    By the definition of the limiting normal cone, this implies that $R^*x^* \in N_C(y)$.

    For ``$\supset$'', let $x^*\in X^*$ be such that $R^*x^* \in N_C(y)$.
    Then the assumption of normal regularity of $C$ at $y$ implies that $R^*x^* \in \frechetNormal_C(y)$. Hence we deduce from \cref{lemma:cofrechet:cone-linear} that $x^* \in \frechetNormal_{RC}(x)$. By  \cref{thm:cones:inclusions}, this implies that $x^* \in N_{RC}(x)$.
    Thus we have proved “$\supset$”.
    The normal regularity of $RC$ at $x$ is now clear from $N_C(y)=\frechetNormal_C(y)$ and comparison of the expression for $N_{RC}(x)$ to the expression for $\frechetNormal_{RC}(x)$ provided by \cref{lemma:cofrechet:cone-linear}.

    If $C$ is not normally regular, we still deduce from the definition \eqref{eq:cones:def-limnormal} of the limiting normal cone the existence of $y_k^* \weaktostar R^*x^*$, $\epsilon_k \downto 0$, and $C \ni y_k \to y$ with $y_k^* \in N_C^{\epsilon_k}(y_k)$.
    By assumption, there exist $x_k^* \in X^*$ with $\rho_k \defeq \norm{R^*x_k^*-y_k^*}_{Y^*} \downto 0$ and $x_k^* \weaktostar x^*$.
    By the definition of the Fréchet normal cone, it follows that $R^*x_k^* \in N_C^{\epsilon_k+\rho_k}(y_k)$.
    Setting $x_k \defeq Ry_k$, \cref{lemma:cofrechet:cone-linear} shows that $x_k^* \in  \frechetNormal_{RC}^{(\epsilon_k + \rho_k)L}(x_k)$.
    We have $x_k = Ry_k \to Ry = x$ as well as $(\epsilon_k + \rho_k)L \downto 0$ and $x_k^* \weaktostar x^*$. Hence $x^* \in N_{RC}(x)$.
    This proves ``$\supset$''.
\end{proof}

\begin{remark}[regularity assumptions]
    \label{rem:colimiting:regularity}
    The assumption on $R$ if $C$ is not regular holds for example when $R^*$ has a right-inverse $\starlinv{R}$ such that $R^*\starlinv{R}-\Id$ is compact: In this case, we can take $x_k^* = \starlinv{R} z_k^*$ because $z_k^* \weaktostar R^*x^*$ implies that $R^*x_k^*-z_k^*=(R^*\starlinv{R}-\Id)z_k^* \to (R^*\starlinv{R} - \Id)R^*x^*=0$.
    If $R$ is a projection in a Hilbert space with a finite-dimensional kernel, we can even take $\starlinv{R}=R$. More generally, if we can write $x=(x_1, x_2)$ with $Rx=x_1$ and $x_2$ in a finite-dimensional subspace, we can take $\starlinv{R}(x_1^*, x_2^*)=x_1^*$.

    Again, the assumption in \cref{lemma:colimiting:cone-linear} that $C$ is normally regular is also not needed if $\kernel R=\{0\}$ or, more generally, if $R$ is a continuously differentiable mapping with $\kernel \grad R(y)=\{0\}$.
\end{remark}

For the fundamental lemma for the limiting coderivative, we need to assume reflexivity of $Y$ in order to apply PSNC via \cref{lemma:colimiting:compactness}.

\begin{lemma}[fundamental lemma on compositions]
    \label{lemma:colimiting:fundamental}
    Let $X,Y,Z$ be Banach spaces with $Y$ reflexive and
    \begin{equation*}
        C \defeq \{(x, y, z) \mid y \in F(x),\, z \in G(y),\, x\in X\}
    \end{equation*}
    for $F: X \setto Y$, and $G: Y \setto Z$.
    Let $(x, y, z) \in C$.
    \begin{enumerate}
        \item\label{item:colimiting:fundamental:g}
            If $G$ is strictly codifferentiable and PSNC at $y$ for $z$, semi-codifferentiable near $(y,z)\in \graph G$, and $y^* \in \coderivative G(y|z)(0)$ implies $y^*=0$, then
            \begin{equation*}
                N_C(x, y, z)=\{(x^*, y^*, z^*) \mid
                    x^* \in \coderivative F(x|y)(-\tilde y^*-y^*), \tilde y^* \in \coderivative G(y|z)(z^*),z^*\in Z^*
                \}.
            \end{equation*}
        \item\label{item:colimiting:fundamental:invf}
            If $\inv F$ is strictly codifferentiable and PSNC at $y$ for $x$, semi-codifferentiable near $(y,x)\in \graph \inv F$, and $y^* \in \coderivative \inv F(y|x)(0)$ implies $y^*=0$, then
            \begin{equation*}
                N_C(x, y, z)=\{(x^*, y^*, z^*) \mid
                x^* \in \coderivative F(x|y)(-\tilde y^*-y^*), -\tilde y^*  \in \coderivative G(y|z)(-z^*), z^*\in Z^*\}.
            \end{equation*}
    \end{enumerate}
    Moreover, if $F$ is N-regular at $x$ for $y$ and $G$ is N-regular at $y$ for $z$, then $C$ is normally regular at $(x, y, z)$.
\end{lemma}

\begin{proof}
    We only consider the case \cref{item:colimiting:fundamental:g}; the case \cref{item:colimiting:fundamental:invf} is shown analogously.
    To show the inclusion \enquote{$\subset$}, let $(x^*, y^*, z^*) \in N_C(x, y, z)$, which by definition holds if and only if there exist $\epsilon_k \downto 0$ as well as $(x^*_k, y^*_k, z_k^*) \weaktostar (x^*, y^*, z^*)$ and $C \ni (x_k, y_k, z_k) \to (x, y, z)$ with $(x^*_k, y^*_k, z^*_k) \in \frechetNormal_C^{\epsilon_k}(x_k, y_k, z_k)$.
    Since by assumption $G$ is semi-codifferentiable at $(y_k,z_k)\in \graph G$ for $k\in\N$ sufficiently large, we can apply \cref{lemma:cofrechet:fundamental}\,\ref{item:cofrechet:fundamental:G}
    to obtain a $\tilde y_k^* \in \frechetCod G(y_k|z_k)(z^*_k)$ such that
    \begin{equation}
        \label{eq:colimiting:frechetcod-fapprox}
        x^*_k \in \frechetCod_{\epsilon_k} F(x_k|y_k)(-\tilde y_k^*-y^*_k).
    \end{equation}
    Since $z^*_k \weaktostar z^*$, $(y_k, z_k) \to (y, z)$, and $\epsilon_k \downto 0$, we deduce from \cref{lemma:colimiting:compactness} that $\tilde y_k^* \weaktostar \tilde y^*$ (if necessary, by passing to a subsequence) for some $\tilde y^* \in \coderivative G(y|z)(z^*)$. Since also $x^*_k \weaktostar x^*$ and $y^*_k \weaktostar y^*$, by \eqref{eq:colimiting:frechetcod-fapprox} and the definition of the limiting coderivative, this implies that $x^* \in \coderivative F(x|y)(-\tilde y^*-y^*)$.

    To show ``$\supset$'', let $x^* \in \coderivative F(x|y)(-\tilde y^*-y^*)$ and $\tilde y^* \in \coderivative G(y|z)(z^*)$.
    We can then by the definition of $\coderivative F(x|y)$ find $\epsilon_k \downto 0$ as well as $(x_k, y_k) \to (x, y)$ and $(x^*_k, \bar y_k^*) \weaktostar (x^*, \tilde y^*+y^*)$ with $x^*_k \in \frechetCod_{\epsilon_k} F(x_k|y_k)(-\bar y_k^*)$.
    Since $G$ is strictly codifferentiable at $y$ for $z$, taking any $z_k \to z$, we can now find $z_k^* \weaktostar z^*$ and $\tilde y_k^* \weaktostar \tilde y^*$ with $\tilde y_k^* \in \frechetCod_{\epsilon_k} G(y_k|z_k)(z_k^*)$. Letting $y_k^* \defeq \bar y_k^*-\tilde y_k^*$, this implies that $y_k^* \weaktostar y^*$ and that $x^*_k \in \frechetCod_{\epsilon_k} F(x_k|y_k)(-\tilde y_k^*-y_k^*)$.
    By \cref{lemma:cofrechet:fundamental}\,\ref{item:cofrechet:fundamental:G}, it follows that $(x_k^*, y_k^*, z_k^*) \in \frechetNormal_C^{\epsilon_k}(x_k, y_k, z_k)$. The claim now follows again from the definition of $N_C(x, y,z)$ as the corresponding outer limit.

    Finally, the normal regularity of $C$ follows from the N-regularity of $F$ and $G$ (via \cref{lemma:colimiting:strictcodiff:gr}) by comparing our identities for $N_C$ with \cref{lemma:cofrechet:fundamental} for $\epsilon=0$.
\end{proof}

If one of the two mappings is single-valued, we can use \cref{lemma:colimiting:regularity:single} for verifying its semi-differentiability and \cref{thm:graphical:single} for the expression of its graphical derivative to obtain from \cref{lemma:colimiting:fundamental} the following two special cases.

\begin{corollary}[fundamental lemma on compositions: single-valued outer mapping]
    \label{lemma:colimiting:fundamental:single-outer}
    Let $X,Y,Z$ be Banach spaces with $Y$ reflexive and
    \begin{equation*}
        C \defeq \{(x, y, G(y)) \mid y \in F(x),\, x\in X\}
    \end{equation*}
    for $F: X \setto Y$ and $G: Y \to Z$.
    If $(x, y, z) \in C$ and $G$ is continuously differentiable near $y$, then
    \begin{equation*}
        N_C(x, y, z)=\{(x^*, y^*, z^*) \mid
            x^* \in \coderivative F(x|y)(-[G'(y)]^*z^*-y^*),\, y^* \in Y^*,\, z^*\in Z^*
        \}.
    \end{equation*}
    Moreover, if $F$ is N-regular at $(x, y)$, then $C$ is normally regular at $(x, y, G(y))$.
\end{corollary}

\begin{proof}
    We apply \cref{lemma:colimiting:fundamental}, where the strict and semi-codifferentiability requirements on $G$ are verified by \cref{lemma:colimiting:regularity:single,lemma:cofrechet:regularity:single};
    the PSNC requirement follows from \cref{lemma:colimiting:psnc:single-valued};
    and the requirement of $y^* \in \coderivative G(y|z)(0)$ implying $y^*=0$ follows from the expression of \cref{thm:graphical:single} for $\coderivative G(y|z)(0)$.
    The claimed normal regularity of $C$ for N-regular $F$ follows from the N-regularity of $G$ established by \cref{thm:graphical:single}.
\end{proof}

\begin{corollary}[fundamental lemma on compositions: single-valued inner mapping]
    \label{lemma:colimiting:fundamental:single-inner}
    Let $X,Y,Z$ be Banach spaces with $Y$ reflexive and
    \begin{equation*}
        C \defeq \{(x, y, z) \mid y=F(x),\, z \in G(y)\}
    \end{equation*}
    for $F: X \to Y$ and $G: Y \setto Z$.
    If $(x, y, z) \in C$, $F$ is continuously differentiable near $x$, and either
    \begin{enumerate}[label=(\alph*)]
        \item $F'(x)$ has a right-inverse $\rinv{F'(x)} \in \linear(Y; X)$ or
        \item $Y^*$ is finite-dimensional,
    \end{enumerate}
    then
    \begin{equation*}
        N_C(x, y, z)
        =\{(F'(x)^*(-\tilde y^*-y^*), y^*, z^*) \mid
            -\tilde y^* \in \coderivative G(y|z)(-z^*),\, y^* \in Y^*,\, z^*\in Z^*
        \}.
    \end{equation*}%
    Moreover, if $G$ is N-regular at $(y, z)$, then $C$ is normally regular at $(x, y, z)$.
\end{corollary}

\begin{proof}
    We apply \cref{lemma:colimiting:fundamental}, where the strict and semi-codifferentiability requirements on $\inv F$ are verified by \cref{lemma:colimiting:regularity:single,lemma:cofrechet:regularity:single};
    the PSNC requirement follows from \cref{lemma:colimiting:psnc:single:inverse};
    and the requirement of $y^* \in \coderivative \inv F(y|x)(0)$ implying $y^*=0$ follows from the expression of \cref{thm:graphical:single:inverse} for $\coderivative \inv F(y|x)(0)$.
    The claimed normal regularity of $C$ for N-regular $G$ follows from the N-regularity of $F$ established by \cref{thm:graphical:single}.
\end{proof}

\section{Calculus rules}\label{sec:colimiting:calculus}

Using these lemmas, we proceed to derive calculus rules.
The proofs are again similar to those in \cref{sec:gderiv:calculus,sec:cofrechet:calculus}, and we only note the differences.
Overall, compared to the Fréchet coderivative, the fundamental single-valued composition \cref{lemma:colimiting:fundamental:single-outer,lemma:colimiting:fundamental:single-inner} (in place of \cref{lemma:cofrechet:fundamental:single-outer,lemma:cofrechet:fundamental:single-inner}) now always impose continuous differentiability, similarly to the Clarke graphical derivative in \cref{sec:gclarke:calculus}.
The cone transformation \cref{lemma:colimiting:cone-linear} (in place of \cref{lemma:cofrechet:cone-linear}) now imposes local (neighborhood) instead of pointwise Lipschitz assumptions.

\begin{theorem}[addition of a single-valued differentiable mapping]
    \label{thm:colimiting:addition}
    Let $X,Y$ be Banach spaces with $X$ reflexive, $G:X\to Y$, and $F:X\setto Y$.
    Let $x \in X$ and $y\in H(x)\defeq F(x)+G(x)$.
    Suppose $G$ is continuously differentiable at $x$.
    If $F$ is N-regular at $(x, y-G(x))$, then $H$ is N-regular at $(x, y)$ and
    \begin{equation*}
        \coderivative H(x|y)(y^*)  = \coderivative F(x|y-G(x))(y^*) + \{[G'(x)]^* y^*\}
        \qquad (y^* \in Y^*).
    \end{equation*}%
    If $Y$ is finite-dimensional, this identity holds without assuming N-regularity.
\end{theorem}

\begin{proof}
    We construct $\graph H=RC$ from $C$ and $R$ defined in \eqref{eq:gderiv:addition:c-r} and follow the proof of \cref{thm:cofrechet:addition}, using \cref{lemma:colimiting:cone-linear,lemma:colimiting:fundamental:single-outer} (the latter of which requires continuous differentiability of $G$ at $x$) in place of \cref{lemma:cofrechet:cone-linear,lemma:cofrechet:fundamental:single-outer}.
    It remains to verify the assumptions of \cref{lemma:colimiting:cone-linear}.

    Consider first the case where $Y$ is finite-dimensional.
    From \eqref{eq:gderiv:addition:c-r}, we have $R^*(\alt x, \alt y) = (\alt y, \alt x, \alt y)$.
    Thus $\bar R(u, \alt x, v) \defeq (\alt x, v)$ is a left-inverse of $R^*$ and
    $[R^* \bar R - \Id](u, \alt x, v) = (v-u, 0, 0)$.
    This is a compact mapping when $Y$ is finite-dimensional,
    which guarantees the condition of \cref{lemma:colimiting:cone-linear} in the non-regular case; see \cref{rem:colimiting:regularity}.

    Since $G$ is continuously differentiable at $x$, it is locally Lipschitz near $x$ by \cref{lem:variation:c1-lipschitz}, and hence the mapping $\inv R_p$ given by \eqref{eq:gderiv:addition:ry} is Lipschitz near $(x,y)$.
    We may therefore apply \cref{lemma:colimiting:cone-linear} both when $C$ is normally regular and when $Y$ is finite-dimensional.
    The normal regularity of $C$ is guaranteed by \cref{lemma:colimiting:fundamental:single-outer} when $F$ is N-regular.
    In this case, it follows that $RC$ is normally regular and hence that $H$ is N-regular.
\end{proof}

\begin{theorem}[outer composition with a single-valued differentiable mapping]
    \label{thm:colimiting:outer}
    Let $X,Y$ be Banach spaces with $Y$ reflexive, $F:X\setto Y$, and $G:Y\to Z$.
    Let $x \in X$ and $z \in H(x) \defeq G(F(x))$ be given.
    Suppose $G$ is continuously differentiable at $y \in F(x) \isect \inv G(\{z\})$ and locally left-invertible near $z$ such that the local left-inverse $\linv G$ is Lipschitz continuous near $z$.
    If $F$ is N-regular at $(x, y)$, then $H$ is N-regular at $(x, z)$ and
    \begin{equation*}
        \coderivative H(x|z)(z^*)=\coderivative F(x|y)([G'(y)]^*z^*)
        \qquad (z^* \in Z^*).
    \end{equation*}
    If $Y$ is finite-dimensional, this identity holds without assuming N-regularity.
\end{theorem}

\begin{proof}
    We construct $\graph H=RC$ from $C$ and $R$ defined in \eqref{eq:gderiv:outer:RC} and follow the proof of \cref{thm:cofrechet:outer}, using \cref{lemma:colimiting:cone-linear,lemma:colimiting:fundamental:single-outer} (the latter of which requires continuous differentiability of $G$ at $y$) in place of \cref{lemma:cofrechet:cone-linear,lemma:cofrechet:fundamental:single-outer}.
    It remains to verify the assumptions of \cref{lemma:colimiting:cone-linear}.

    Consider first the case where $Y$ is finite-dimensional.
    From \eqref{eq:gderiv:outer:RC}, we have $R^*(\alt x, \alt z) = (\alt x, 0, \alt z)$.
    Thus $R$ is a left-inverse of $R^*$ and $[R^* R - \Id](\alt x, \alt y, \alt z) = (0, \alt y, 0)$.
    This is a compact mapping when $Y$ is finite-dimensional,
    which guarantees the condition of \cref{lemma:colimiting:cone-linear} in the non-regular case; see \cref{rem:colimiting:regularity}.

    Since $\linv G$ is Lipschitz near $x$, the mapping $\inv R_p$ given by \eqref{eq:gderiv:outer:invR} is Lipschitz near $(x,z)$.
    We may therefore apply \cref{lemma:colimiting:cone-linear} both when $C$ is normally regular and when $Y$ is finite-dimensional.
    The normal regularity of $C$ is guaranteed by \cref{lemma:colimiting:fundamental:single-outer} when $F$ is N-regular.
    In this case, it follows that $RC$ is normally regular and hence that $H$ is N-regular.
\end{proof}

\begin{corollary}[outer composition with a linear operator]
    \label{cor:colimiting:outer:linear}
    Let $X,Y,Z$ be Banach spaces with $Y$ reflexive, $A\in \linear(Y; Z)$, and $F:X\setto Y$.
    Let $x \in X$ and $z \in H(x) \defeq AF(x)$ be given.
    If $A$ has a left-inverse $\linv A$ and $F$ is N-regular at $(x,y)$ for the unique $y\in Y$ with $Ay=z$, then $H$ is N-regular at $(x, z)$ and
    \begin{equation*}
        \coderivative H(x|z)(z^*)=
        \coderivative F(x|y)(A^*z^*)
        \qquad (z^* \in Z^*).
    \end{equation*}
    If $Y$ is finite-dimensional, this identity holds without assuming N-regularity.
\end{corollary}

\begin{theorem}[inner composition with a single-valued differentiable mapping]
    \label{thm:colimiting:inner}
    Let $X,Y,Z$ be Banach spaces with $Y$ reflexive, $F: X\to Y$ and $G:Y\setto Z$.
    Let $x \in X$ and $z \in H(x)\defeq G(F(x))$ be given.
    Suppose $F$ is continuously differentiable at $x$ such that $F'(x)$ has a right-inverse $\rinv{F'(x)} \in \linear(Y; X)$.
    If $G$ is N-regular at $(F(x), z)$, then $H$ is N-regular at $(x, z)$ and
    \begin{equation*}
        \coderivative H(x|z)(z^*)  = [F'(x)]^* \coderivative G(F(x)|z)(z^*)
        \qquad (z^* \in Z^*).
    \end{equation*}%
    If $Y$ is finite-dimensional, this identity holds without assuming N-regularity.
\end{theorem}

\begin{proof}
    We construct $\graph H=RC$ from $C$ and $R$ defined in \eqref{eq:gderiv:inner:RC} and follow the proof of \cref{thm:cofrechet:inner}, using \cref{lemma:colimiting:cone-linear,lemma:colimiting:fundamental:single-inner} (the latter of which requires continuous differentiability of $F$ at $x$) in place of \cref{lemma:cofrechet:cone-linear,lemma:cofrechet:fundamental:single-inner}.
    It remains to verify the assumptions of \cref{lemma:colimiting:cone-linear}.

    Consider first the case where $Y$ is finite-dimensional.
    From \eqref{eq:gderiv:inner:RC}, we have $R^*(\alt x, \alt z) = (\alt x, 0, \alt z)$.
    Thus $R$ is a left-inverse of $R^*$.
    We thus have $[R^* R - \Id](\alt x, \alt y, \alt z) = (0, \alt y, 0)$.
    This is a compact mapping when $Y$ is finite-dimensional,
    which guarantees the condition of \cref{lemma:colimiting:cone-linear} in the non-regular case; see \cref{rem:colimiting:regularity}.

    Since $F$ is continuously differentiable at $x$, it is locally Lipschitz near $x$ by \cref{lem:variation:c1-lipschitz}, and hence the mapping $\inv R_p$ given by \eqref{eq:gderiv:inner:invRp} is Lipschitz near $(x,z)$.
    We may therefore apply \cref{lemma:colimiting:cone-linear} both when $C$ is normally regular and when $Y$ is finite-dimensional.
    The normal regularity of $C$ is guaranteed by \cref{lemma:colimiting:fundamental:single-inner} when $F$ is N-regular.
    In this case, it follows that $RC$ is normally regular and hence that $H$ is N-regular.
\end{proof}

\begin{corollary}[inner composition with a linear operator]
    \label{cor:colimiting:inner:linear}
    Let $X,Y,Z$ be Banach spaces with $Y$ reflexive, $A\in \linear(X; Y)$, and $G:Y\setto Z$.
    Let $x \in X$ and $z \in H(x)\defeq G(Ax)$ be given.
    Suppose $A$ has a right-inverse $\rinv A \in \linear(Y; X)$.
    If $G$ is N-regular at $(Ax,z)$, then $H$ is N-regular at $(x, z)$ and
    \begin{equation*}
        \coderivative H(x|z)(z^*)  = A^* \coderivative G(Ax|z)(z^*)
        \qquad (z^* \in Z^*).
    \end{equation*}%
    If $Y$ is finite-dimensional, this identity holds without assuming N-regularity.
\end{corollary}

To apply these results for chain rules of subdifferentials, we now need to assume that \emph{both} spaces are reflexive in addition to N-regularity of the subdifferential.

\begin{corollary}[second-order chain rule for convex subdifferentials]
    \label{cor:colimiting:second-convex}
    Let $X,Y$ be reflexive Banach spaces, let $f:Y\to \Rbar$ be proper, convex, and lower semicontinuous, $A\in\linear(X; Y)$ be such that $A$ has a right-inverse $\rinv{A} \in \linear(Y; X)$, and $\range A \isect \interior \dom f \ne \emptyset$. Let $h\defeq f\circ A$.
    If $\partial f$ is N-regular at $(Ax,y^*)$ for $x\in X$ and $y^*\in \partial f(Ax)$, then $\partial h$ is N-regular at $(x,A^*y^*)$ and
    \begin{equation*}
        \coderivative[\subdiff h](x|x^*)(\dir x)  = A^* \coderivative[\subdiff f](Ax|y^*)(A\dir x)
        \qquad (\dir x \in X).
    \end{equation*}
    If $Y$ is finite-dimensional, this identity holds without assuming N-regularity.
\end{corollary}

\begin{proof}
    The expression for $\partial h(x)$ follows from \cref{thm:convex:chain}, to which we apply \cref{cor:colimiting:inner:linear} (which requires reflexivity of $Y$) as well as \cref{cor:colimiting:outer:linear} with $A^*$ in place of $A$ (which requires reflexivity of $X^*$, which holds if and only if $X$ is reflexive), recalling that a right-inverse $\rinv A$ for $A$ produces the left-inverse $\linv{(A^*)} = (\rinv A)^*$ for $A^*$.
\end{proof}

\begin{theorem}[product rule]
    \label{thm:colimiting:product}
    Let $X,Y,Z$ be Banach spaces with $X,Y$ reflexive, $G:X\to \linear (Y; Z)$ be continuously differentiable at $x \in X$, and $F:X\setto Y$.
    Assume that $G(\alt x) \in \linear(Y; Z)$ has a left-inverse $\linv{G(\alt x)} \in \linear(Z; Y)$ for $\alt x$ near $x\in X$ such that the mapping $\alt x \mapsto \linv{G(\alt x)}$ is Lipschitz near $x$.
    Let $z\in H(x)\defeq G(x)F(x) \defeq \Union_{y \in F(x)} G(x)y$ and let $y \in F(x)$ be the unique element satisfying $G(x)y=z$.
    If $F$ is N-regular at $x$ for $y$, then $H$ is N-regular at $x$ for $z$ and
    \begin{equation*}
        \coderivative H(x|z)(z^*)
        =
        \{([G'(x)\freevar]y)^*z^*\} + \coderivative F(x|y)(G(x)^*z^*)
        \qquad (z^* \in Z^*)
    \end{equation*}
    for $([G'(x)\freevar]y)^*z^*\in X^*$ as defined in \cref{thm:cofrechet:product}.
    If $X$ and $Y$ are finite-dimensional, this identity holds without assuming N-regularity.
\end{theorem}

\begin{proof}
    We follow the proof of \cref{thm:cofrechet:product}, using \cref{lemma:colimiting:cone-linear,thm:colimiting:outer} in place of \cref{lemma:cofrechet:cone-linear,thm:cofrechet:outer}.
    The proof also depends on a straightforward adaptation of \cref{lemma:cofrechet:cartesian-product} to the limiting derivative.
    \Cref{thm:colimiting:outer} is applied to $\bar G: X \times Y \to X \times Z$ and $\bar F: X \to X \times Y$ as constructed in the proof of \cref{thm:cofrechet:product};
    this is possible due to the assumed reflexivity of the intermediate space $X \times Y$ in the composition $\bar G \circ \bar F$, i.e., of both $X$ and $Y$.
    The assumptions on $G$ and its left-inverse guarantee the corresponding assumptions on $\bar G$.
    If $F$ is not N-regular, we also need here the finite-dimensionality of both $X$ and $Y$.
    It remains to verify the assumptions of \cref{lemma:colimiting:cone-linear}.

    Consider first the case when $X$ and $Y$ are finite-dimensional.
    Recalling the definition of $R$ from the proof of \cref{thm:cofrechet:product}, we have $R^*(\alt x_1, \alt z) = (\alt x_1, 0, \alt z)$.
    Thus $R$ is a left-inverse of $R^*$ and
    $[R^* R - \Id](\alt x_1, \alt x_2, \alt z) = (0, -\alt x_2, 0)$.
    This is a compact mapping when $X$ is finite-dimensional,
    which guarantees the condition of \cref{lemma:colimiting:cone-linear} in the non-regular case; see \cref{rem:colimiting:regularity}.

    The mapping $\inv R_p$ given by \eqref{eq:cofrechet:product:invRp} is clearly Lipschitz near $(x,z)$.
    We may therefore apply \cref{lemma:colimiting:cone-linear} both when $\graph(\bar G \circ \bar F)$ is normally regular and when $X$ and $Y$ are finite-dimensional.
    The normal regularity of $\graph(\bar G \circ \bar F)$ is guaranteed by \cref{thm:cofrechet:outer} when $F$ is N-regular.
    In this case, it follows that $R\graph(\bar G \circ \bar F)$ is normally regular and hence that $H$ is N-regular.
\end{proof}

The next result follows from \cref{thm:colimiting:product,thm:colimiting:inner} similarly to the proof of \cref{cor:cofrechet:second-clarke}.

\begin{corollary}[second-order chain rule for Clarke subdifferentials]
    \label{cor:colimiting:second-clarke}
    Let $X,Y$ be reflexive Banach spaces, let $f:Y\to \R$ be locally Lipschitz continuous, and let $S:X\to Y$ be twice continuously differentiable. Set $h:X\to \R$, $h(x)\defeq f(S(x))$.
    Suppose there exists a neighborhood $U$ of $x \in X$ such that
    \begin{enumerate}
        \item $f$ is Clarke regular at $S(\alt x)$ for all $\alt x\in X$;
        \item $S'(\alt x)$ has a right-inverse $\rinv{S'(\tilde x)} \in \linear(Y; X)$ for all $\tilde x\in U$;
        \item the mapping $\alt x \mapsto \rinvstar{S'(\alt x)}$ is Lipschitz near $x$.
    \end{enumerate}
    If $\partial_C f$ is N-regular at $S(x)$ for $y^*\in \partial_C f(S(x))$,
    then $\partial_C h$ is N-regular at $x$ for $x^* = S'(x)^*y^*$ and
    \begin{multline*}
        D^*[\subdiff_C h](x|x^*)(x^{**})
        =
        \{x^{**}[S''(x)\freevar]^*y^*\}\\
        + S'(x)^* D^*[\subdiff_Cf](S(x)|y^*)(S'(x)x^{**})
        \quad (x^{**}\in X)
    \end{multline*}
    for $x^{**}[S''(x)\freevar]^*y^*\in X^*$ defined as in \cref{cor:cofrechet:second-clarke}.
    If $X$ and $Y$ are finite-dimensional, this identity holds without assuming N-regularity.
\end{corollary}

\begin{remark}
    Even in finite dimensions, calculus rules for the sum $F+G$ of arbitrary set-valued mappings $F, G: \R^N \setto \R^M$ or the composition $F \circ H$ for $H: \R^N \setto \R^N$ are much more limited and in general only yield inclusions of the form
    \begin{align*}
        \coderivative [F+G](x|y)(y^*) &\subset \Union_{\substack{y=y_1+y_2, \\ y_1 \in F(x), \\ y_2 \in G(x)}} \coderivative F(x|y_1)(y^*) + \coderivative G(x|y_2)(y^*)
        \shortintertext{and}
        \coderivative[F \circ H](x|y)(y^*) &\subset \Union_{z \in H(x) \isect \inv F(y)} \coderivative H(x|z) \circ \coderivative F(z|y)(y^*).
    \end{align*}
    We refer to \cite{Rockafellar:1998,mordukhovich2018variational} for these and other results.
\end{remark}

\section{Subdifferential calculus}\label{sec:colimiting:subdiff}

The above results immediately yield calculus rules for the Mordukhovich or limiting subdifferential from \cref{sec:limiting:mordukhovich}.
To see this, we recall from \eqref{eq:graphical:mordukhovich-subdiff} that for $f: X \to \Rbar$ we have that
\[
    \subdiff_M f(x) = D^*[\epi_f](x|f(x))(1),
\]
where $\epi_f:X\setto\R$ denotes the epigraphical mapping of $f$.
We also observe that if $g: X \to \Rbar$ is Fréchet differentiable, then $g'(x)^* \in \linear(\R; X^*)$ and hence
\[
    g'(x)^*z^* = z^* g'(x)
    \qquad\text{for all }
    z^* \in \R.
\]
Finally, we call $f$ \term[functional!regular!epigraphically]{epigraphically regular}\index{regularity!epigraphical} at $x$ if $\epi_f$ is N-regular at $x$ for $f(x)$, i.e., if $N_{\epi f}(x, f(x))=\frechetNormal_{\epi f}(x, f(x))$.
For example, by \cref{cor:cones:convex-regularity} any functional that locally coincides with a convex or concave function is epigraphically regular.

In many cases, we can relate epigraphical regularity to equality of the Fréchet and Clarke subdifferentials or to Clarke regularity.\index{regularity!Clarke}

\begin{lemma}
    \label{lemma:colimiting:epigraphical-regularity-subdiff}
    Let $X$ be a reflexive and Gateaux smooth Banach space and $f: X \to \R$ be Lipschitz continuous near $x$.
    Then $f$ is epigraphically regular at $x$ if and only if $\subdiff_C f(x)=\subdiff_F f(x)$.
\end{lemma}

\begin{proof}
    By \cref{thm:cones:regularity:infdim}, epigraphical regularity of $f$ at $x$ is equivalent to
    \begin{equation}
        \label{eq:colimiting:epigraphical-regularity:alt}
        \clarkeTangent_{\epi f}(x, f(x))=\polar{\frechetNormal_{\epi f}(x, f(x))}.
    \end{equation}
    Since $(x, f(x)) \in \epi f$, the cone $\frechetNormal_{\epi f}(x, f(x))$ is nonempty by definition as well as closed and convex by \cref{thm:cones:basic-prop}.
    Hence \cref{lemma:functan:polar-inclusion} yields $\bipolar{\frechetNormal_{\epi f}(x, f(x))}=\frechetNormal_{\epi f}(x, f(x))$.
    Furthermore, by the definition \eqref{eq:graphical:clarke-normal}, $N^C_{\epi f}(x, f(x)) = \polar{\clarkeTangent_{\epi f}(x, f(x))}$ where $\clarkeTangent_{\epi f}(x, f(x))$ is a closed and convex cone by \cref{thm:cones:basic-prop}.
    It then follows from \cref{lemma:functan:polar-inclusion} that \eqref{eq:colimiting:epigraphical-regularity:alt} holds and hence that epigraphical regularity implies
    \[
        N^C_{\epi f}(x, f(x))=\frechetNormal_{\epi f}(x, f(x)).
    \]
    By \cref{eq:graphical:frechet-subdiff,lemma:graphical:clarke-subdiff}, this implies that $\subdiff_C f(x)=\subdiff_F f(x)$.

    To show the converse implication, we need to show that both the Fréchet and the Clarke subdifferential generate the full corresponding normal cones $\frechetNormal_{\epi f}(x, f(x))$ and $N^C_{\epi f}(x, f(x))$, respectively. when $\subdiff_C f(x)=\subdiff_F f(x)$.
    By assumption, there exists a neighborhood $U$ of $x$ where $f$ is Lipschitz continuous.
    Take any $\alt x \in U$.
    By definition, we have $(x^*, t) \in \frechetNormal_{\epi f}(\alt x, f(\alt x))$ if and only if
    \begin{equation}
        \label{eq:colimiting:frechet-limit-epi}
        \limsup_{\epi f \ni (x_k, t_k) \to (\alt x, f(\alt x))} \frac{\dualprod{x_k-\alt x}{x^*}_X + (t_k-f(\alt x))t}{\norm{x_k-\alt x}_X+\abs{t_k-f(\alt x)}} \le 0.
    \end{equation}
    By taking $x_k=\alt x$ and $t_k \downto f(\alt x)$, we see that $t \le 0$.
    If $t=0$, we can take $t_k=f(x_k)$ and use that $f$ is finite-valued and continuous at $\alt x$ to take $x_k=\alt x+h/k$ for some $\epsilon>0$ and any $h \in \B(0, \epsilon)$. This shows that $x^*=0$ as well.
    When $t<0$, the $\limsup$ in \eqref{eq:colimiting:frechet-limit-epi} is achieved by taking $t_k=f(x_k)$ and hence the definition \eqref{eq:limiting:frechet} of the Fréchet subdifferential and the fact that $f$ is Lipschitz in $U$ implies that $x^*\in \subdiff_F f(\alt x)$. Combining both cases, we find that
    \begin{equation}
        \label{eq:colimiting:frechet-epi-lipshitz}
        \frechetNormal_{\epi f}(\alt x, f(\alt x))
        =
        \{t(x^*, -1) \mid t \ge 0,\, x^* \in \subdiff_F f(\alt x)\}.
    \end{equation}
    On the other hand, \cref{cor:graphical:clarke-weakstar-convex,thm:cones:inclusions}, \eqref{eq:graphical:frechet-subdiff}, and \eqref{eq:graphical:mordukhovich-subdiff} imply that
    \begin{equation*}
        \label{eq:colimiting:clarke-frechet-subdiff-incl}
        \subdiff_C f(\alt x) = \closure^* \conv \subdiff_M f(\alt x) \supset \subdiff_F f(\alt x)
    \end{equation*}
    and hence that
    \begin{equation}
        \label{eq:colimiting:frechet-epi-lipshitz-clarke}
        \frechetNormal_{\epi f}(\alt x, f(\alt x))
        \subset
        \{t(x^*, -1) \mid t \ge 0,\, x^* \in \subdiff_C f(\alt x)\}.
    \end{equation}
    By assumption and \eqref{eq:colimiting:frechet-epi-lipshitz}, this holds as an equality at $x=\alt x$.
    For $\alt t>f(\alt x)$, we have by continuity that $(\alt x, \alt t) \in \interior \epi f$ and hence that $\frechetNormal_{\epi f}(\alt x, \alt t) = \{0\}$.
    \Cref{thm:cones:nonepsilon-limnormal,eq:colimiting:frechet-epi-lipshitz-clarke,lem:clarke:closed} thus yield that
    \[
        \begin{aligned}
            N_{\epi f}(x, f(x))
            &
            =
            \weakstarlimsup_{\alt x \to x} \frechetNormal_{\epi f}(\alt x, f(\alt x))
            \\
            &
            \subset
            \weakstarlimsup_{\alt x \to x}
            \{t(x^*, -1) \mid t \ge 0,\, x^* \in \subdiff_C f(\alt x)\}
            \\
            &
            \subset \{t(x^*, -1) \mid t \ge 0,\, x^* \in \subdiff_C f( x)\}.
        \end{aligned}
    \]
    From \cref{cor:graphical:clarke-normal,thm:cones:limiting-polar,lemma:functan:polar-inclusion} and the convexity of $\subdiff_C f(x)$ (\cref{lem:clarke:properties}), we then obtain that
    \[
        \begin{aligned}
            N_{\epi f}^C(x, f(x))
            =
            \polar{\clarkeTangent_{\epi f}(x, f(x))}
            =
            \bipolar{N_{\epi f}(x, f(x))}
            &
            \subset
            \bipolar{\{t(x^*, -1) \mid t \ge 0,\, x^* \in \subdiff_C f( x)\}}
            \\
            &
            =
            \{t(x^*, -1) \mid t \ge 0,\, x^* \in \subdiff_C f( x)\}.
        \end{aligned}
    \]
    By the continuity of $f$, $\epi f$ is closed.
    Hence, by the previous inclusion, $\subdiff_C f(x)=\subdiff_F f(x)$, \eqref{eq:colimiting:frechet-epi-lipshitz}, and \cref{cor:graphical:clarke-normal}, we have that
    \[
        \begin{aligned}
            N_{\epi f}(x, f(x))
            \subset
            N_{\epi f}^C(x, f(x))
            &
            \subset
            \{t(x^*, -1) \mid t \ge 0,\, x^* \in \subdiff_C f( x)\}
            \\
            &
            =
            \{t(x^*, -1) \mid t \ge 0,\, x^* \in \subdiff_F f( x)\}
            \\
            &
            =
            \frechetNormal_{\epi f}(x, f(x)).
        \end{aligned}
    \]
    Since the converse inclusion always holds by \cref{thm:cones:inclusions}, the claim follows.
\end{proof}

\begin{theorem}
    \label{thm:colimiting:clarke-epi-regularity}
    Let $X$ be a reflexive and Gateaux smooth Banach space and $f: X \to \R$ be Lipschitz continuous near $x\in X$.
    Then epigraphical regularity at $x$ implies Clarke regularity at $x$.
    If $X$ is finite-dimensional, then the converse holds as well.
\end{theorem}

\begin{proof}
    Assume first that $f$ is epigraphically regular at $x$,
    which implies that $\subdiff_C f(x) = \subdiff_F f(x)$
    by \cref{lemma:colimiting:epigraphical-regularity-subdiff}.
    It thus follows from the definition \eqref{eq:limiting:frechet} of the Fréchet subdifferential that every $x^* \in \subdiff_C f(x)$ satisfies
    \[
        \liminf_{t \downto 0} \frac{f(x+th)-f(x)-\dualprod{x^*}{th}_X}{t} \ge 0
        \qquad\text{for all } h \in X
    \]
    which implies that
    \[
        \dualprod{x^*}{h}_X \le f'(x; h)
        \qquad\text{for all } h \in X.
    \]
    \Cref{cor:clarke:support-dir} thus shows that
    \begin{equation*}
        f^\circ(x;h) = \sup_{x^*\in\partial_C f(x)} \dualprod{x^*}{h}_X \leq f'(x; h)
        \qquad\text{for all } h \in X.
    \end{equation*}
    Since $f'(x; h) \leq f^\circ(x; h)$ always holds (see \eqref{eq:clarke:dir-gen-dir-inequality}), $f$ is Clarke regular at $x$.

    For the converse, let $X$ be finite-dimensional and $f$ be Lipschitz continuous near $x$ and Clarke regular at $x$.
    Pick $x^* \in \subdiff_C f(x)$ as well as a sequence $0 \ne h_k \to 0$ satisfying
    \[
        A
        \defeq
        \lim_{k \to \infty} \frac{f(x+h_k) - f(x) - \dualprod{x^*}{h_k}_X}{\norm{h_k}_X}
        =
        \liminf_{y \to x} \frac{f(y) - f(x) - \dualprod{x^*}{y-x}_X}{\norm{y-x}_X}.
    \]
    Since $X$ is finite-dimensional, we can assume (if necessary after passing to a subsequence) that $h_k/\norm{h_k}_X \to h$ for some $h \in X$.
    Then by Clarke regularity and the definition of the Clarke subdifferential,
    $
    f'(x; h) = f^\circ(x;h) \ge \dualprod{x^*}{h}_X.
    $
    It follows that
    \[
        \begin{aligned}
            0
            \le
            f'(x; h) - \dualprod{x^*}{h}_X
            &
            = \lim_{t \downto 0} \frac{f(x+th)-f(x)- \dualprod{x^*}{th}_X}{t}
            \\
            &
            = \lim_{k \to \infty} \frac{f(x+\norm{h_k}_X h) - f(x) - \dualprod{x^*}{\norm{h_k}h}_X}{\norm{h_k}_X}
            \\
            &
            \le
            A
            + \limsup_{k \to \infty} \frac{f(x+\norm{h_k}_X h) - f(x+h_k) - \dualprod{x^*}{\norm{h_k}_Xh-h_k}_X}{\norm{h_k}_X}.
        \end{aligned}
    \]
    Since $f$ is Lipschitz continuous near $x$, and $(\norm{h_k}_X h-h_k)/\norm{h_k} = h - h_k/\norm{h_k}_X \to 0$, this establishes that $0 \le A$ and hence that $x^* \in \subdiff_F f(x)$.
    We have therefore shown that $ \subdiff_C f(x) \subset \subdiff_F f(x)$.
    As the converse inclusion always holds by \cref{cor:graphical:subdiff-inclusions}, the claim now follows from \cref{lemma:colimiting:epigraphical-regularity-subdiff}.
\end{proof}

\begin{remark}
    If $X$ is not finite-dimensional in the second implication of \cref{thm:colimiting:clarke-epi-regularity}, the Eberlein--\u{S}mulyan theorem (\cref{thm:ebsmul}) allows us to obtain $h$ as a weak limit of the $h_k/\norm{h_k}_X$.
    Then we can still show $0 \le A$ and go through with the proof as long as $f$ is Lipschitz continuous with respect to a semi-norm induced by a compact operator.
\end{remark}

We now proceed with calculus rules for the limiting subdifferential.
The following direct corollary of \cref{thm:colimiting:addition} yields \cref{thm:limiting:sum}.

\begin{corollary}[addition of a differentiable mapping]
    \label{cor:colimiting:addition:subdiff}
    Let $X$ be a reflexive Banach space, $g:X\to \R$, and $f:X \to \Rbar$.
    Suppose $g$ is continuously differentiable at $x\in X$.
    Then
    \begin{equation*}
        \subdiff_M (f+g)(x) = \subdiff_M f(x) + \{g'(x)\}.
    \end{equation*}%
    Moreover, if $f$ is epigraphically regular at $x$, then $f+g$ is epigraphically regular at $x$
\end{corollary}

\begin{proof}
    Since $Y=\R$ is finite-dimensional, we can apply \cref{thm:colimiting:addition} to obtain the claimed identity without having to assume epigraphical regularity.
    Assuming in addition epigraphical regularity, \cref{thm:colimiting:addition} yields the regularity of the sum.
\end{proof}

We now turn to chain rules for the Mordukhovich subdifferential.

\begin{corollary}[outer composition with an increasing differentiable mapping]
    \label{cor:colimiting:outer:subdiff}
    Let $X$ be a Banach space, $f:X\to \R$, and $g:\R \to \Rbar$.
    Let $x \in X$ and $z \defeq  h(x) \defeq g(f(x))$ be given.
    Suppose $g$ is increasing and continuously differentiable at $y \defeq f(x)$ and left-invertible on $\range g$ near $z$ such that the left-inverse is Lipschitz continuous near $z$.
    Then
    \begin{equation*}
        \subdiff_M h(x) = g'(y)\subdiff_M f(x).
    \end{equation*}
    Moreover, if $f$ is epigraphically regular at $x$, then $h$ is epigraphically regular at $x$.
\end{corollary}

\begin{proof}
    Since $g$ is assumed to be increasing, $g'(y) > 0$ and hence $g'(y)^* z^* = z^* g'(y) = g'(y) > 0$ for $z^* = 1$.
    Due to \cref{cor:graphical:pos-hgen}, we then have
    \[
        D^*f(x|y)( g'(y)) =   g'(y) D^*f(x|y)(1)=  g'(y) \subdiff_M f(x).
    \]
    Since $Y=\R$ is finite-dimensional, the claim now follows from \cref{thm:colimiting:outer}.
\end{proof}

From \cref{thm:colimiting:inner}, we similarly obtain \cref{thm:limiting:chain} under a regularity assumption.
\begin{corollary}[inner composition with a differentiable mapping]
    \label{cor:colimiting:inner:subdiff}
    Let $X$ and $Y$ be Banach spaces with $Y$ reflexive, $f: X\to Y$ and $g:Y\to \Rbar$.
    Let $h \defeq g \circ f$ and $x \in \dom h$ be given.
    Suppose $f$ is continuously differentiable at $x$ such that $f'(x)$ has a left-inverse $\linv{f'(x)} \in \linear(Y; X)$.
    If $g$ is epigraphically regular at $f(x)$, then $h$ is epigraphically regular at $x$ and
    \begin{equation*}
        \subdiff_M h(x)  = f'(x)^* \subdiff_M g(f(x)).
    \end{equation*}%
    If $Y$ is finite-dimensional, then this identity holds without assuming epigraphical regularity.
\end{corollary}

As a special case, we obtain from \cref{cor:colimiting:inner:linear} the following linear chain rule (note the slightly different regularity assumption on the inner mapping).
\begin{corollary}[inner composition with a linear operator]
    \label{cor:colimiting:inner:linear:subdiff}
    Let $X$ and $Y$ be Banach spaces with $Y$ reflexive, $A\in \linear(X; Y)$, and $f:Y\to \Rbar$.
    Let $h \defeq f \circ A$ and $x \in \dom h$ be given.
    If $A$ has a right-inverse $\rinv A \in \linear(Y; X)$ and $f$ is epigraphically regular at $Ax$ for  $x\in X$, then $h$ is epigraphically regular at $x$ and
    \begin{equation*}
        \subdiff_M h(x)  = A^* \subdiff_M f(Ax).
    \end{equation*}%
    If $Y$ is finite-dimensional, then this identity holds without assuming epigraphical regularity.
\end{corollary}

Finally, we obtain a product rule for the Mordukhovich subdifferential.
\begin{corollary}[product rule]
    \label{cor:colimiting:product:subdiff}
    Let $X$ be a reflexive Banach space, $g:X\to \Rbar$, and $f:X \to \Rbar$.
    Define $h(\alt x)\defeq g(\alt x)f(\alt x)$.
    Assume that $g$ is continuously differentiable at $x \in X$, that $g(\alt x) > 0$ for $\alt x$ near $x$, and that the mapping $\alt x \mapsto 1/g(\alt x)$ is Lipschitz continuous near $x$.
    If $f$ is epigraphically regular at $x$, then $h$ is epigraphically regular at $x$ and
    \begin{equation*}
        \subdiff_M h(x)
        =
        g(x)\subdiff_M f(x)+\{f(x)g'(x)\}.
    \end{equation*}
    If $X$ is finite-dimensional, then this identity holds without assuming epigraphical regularity.
\end{corollary}

\begin{proof}
    We again use the positivity of $g(x)=g(x)^*z^*$ with $z^*=1$ and \cref{cor:graphical:pos-hgen} to deduce the claim from \cref{thm:colimiting:product}.
\end{proof}

\chapter{Second-order optimality conditions}
\label{chap:secondorder}

We now illustrate the use of set-valued derivatives for optimization problems by showing how these can be used to derive second-order (sufficient and necessary) optimality conditions for non-smooth problems.
Again, we do not aim for the most general or sharpest possible results and focus instead on problems having the form
\begin{equation*}
    \min_{x\in C} \frac1p\norm{S(x)-z}_Y^p + \frac{\alpha}{q}\norm{x}_X^q
\end{equation*}
involving the composition of a nonsmooth convex functional with a smooth nonlinear operator. As in the previous chapters, we will also assume a regularity conditions that allows for cleaner results.

\section{Second-order derivatives}

Let $X$ be a Banach space and $f:X\to \Rbar$. In this chapter, we set
\begin{equation*}
    \subdiff_C f(x)\defeq \setof{x^*\in X^*}{(x^*, -1) \in N^C_{\epi f}(x, f(x))},
\end{equation*}
where $N^C_A \defeq \polar\clarkeTangent_A$ is the Clarke normal cone. By \cref{lemma:graphical:clarke-subdiff}, this coincides with the classical Clarke subdifferential if $f:X\to\R$ is locally Lipschitz continuous.

As in the smooth case, second-order conditions are based on a local quadratic model built from curvature information at a point. Since in the nonsmooth case, second derivatives, i.e., graphical derivatives of the subdifferential, are no longer unique, we need to consider the entire set of them when building this curvature information. We therefore need to distinguish a \term[model!lower curvature]{lower curvature model} at $x\in X$ for $x^*\in \subdiff_C f(x)$ in direction $\dir x\in X$
\begin{align*}
    Q_f(\dir x; x|x^*) &\defeq \inf_{\dir x^* \in D[\subdiff_C f](x|x^*)(\dir x)}  \dualprod{\dir x^*}{\dir x}_X
    \intertext{as well as an \term[model!upper curvature]{upper curvature model}}
    Q^f(\dir x; x|x^*) &\defeq \sup_{\dir x^* \in D[\subdiff_C f](x|x^*)(\dir x)}  \dualprod{\dir x^*}{\dir x}_X.
    \intertext{It turns out that even for $\dir x\neq 0$, we need to consider the \term[model!stationary upper]{stationary upper model}}
    Q^f_0(\dir x; x|x^*) &\defeq \sup_{\dir x^* \in D[\subdiff_C f](x|x^*)(0)}  \dualprod{\dir x^*}{\dir x}_X,
\end{align*}
which we use to define the \term[model!extended upper]{extended upper model}
\begin{equation*}
    \begin{aligned}
        \hat Q^f(\dir x; x|x^*) &\defeq \max\left\{Q^f(\dir x; x|x^*), Q^f_0(\dir x; x|x^*) \right\}\\
        &= \sup_{\dir x^* \in D[\subdiff_C f](x|x^*)(\dir x) \union D[\subdiff_C f](x|x^*)(0)}  \dualprod{\dir x^*}{\dir x}_X.
    \end{aligned}
\end{equation*}

For smooth functionals, these models coincide with the quadratic form \eqref{eq:variation:quadratic} induced by the second-order Fréchet derivative.\index{derivative!second-order}
\begin{theorem}\label{thm:secondorder:model-smooth}
    Let $X$ be a Banach space and let $f: X\to\R$ be twice continuously differentiable. Then for every $x,\dir x\in X$,
    \begin{align*}
        Q_f(\dir x; x|f'(x))&=Q^f(\dir x; x|f'(x))=\dualprod{f''(x)\dir x}{\dir x}_X
        \shortintertext{and}
        \hat Q^f(\dir x; x|f'(x))&=\max\left\{0, \dualprod{f''(x)\dir x}{\dir x}_X\right\}.
    \end{align*}
\end{theorem}
\begin{proof}
    Since $\subdiff_C f(x)=\{f'(x)\}$ by \cref{thm:clarke:frechet}, it follows from \cref{thm:graphical:single} that
    \begin{equation*}
        D[\subdiff_C f)](x|f'(x))(\dir x)=\dualprod{f''(x)\dir x}{\dir x}_X
    \end{equation*}
    and in particular $D[\subdiff_C f](x|f'(x))(0)=0$, which immediately yields the claim.
\end{proof}

We illustrate the nonsmooth case with the usual examples of the indicator functional of the unit ball and the norm on $\R$.

\begin{lemma}
    Let $f(x) = \delta_{[-1,1]}(x)$, $x\in \R$. Then for every $x^*\in \partial f(x)$ and $\dir x\in\R$,
    \begin{align*}
        Q_f(\dir x; x|x^*) & =
        \begin{cases}
            \infty & \text{if } \abs{x}=1,\, x^*=0,\, x\dir x > 0, \\
            \infty & \text{if } \abs{x}=1,\, x^* \in (0,\infty)x,\, \dir x \ne 0, \\
            0, & \text{otherwise},
        \end{cases}
        \\
        Q^f(\dir x; x|x^*) & =
        \begin{cases}
            -\infty & \text{if } \abs{x}=1,\, x^*=0,\, x\dir x > 0, \\
            -\infty & \text{if } \abs{x}=1,\, x^* \in (0,\infty)x,\, \dir x \ne 0, \\
            0, & \text{otherwise},
        \end{cases}
    \end{align*}
    and
    \begin{equation*}
        \hat Q^f(\dir x; x|x^*)=
        Q^f_0(\dir x; x|x^*) =
        \begin{cases}
            \infty & \text{if } \abs{x}=1,\,  x^*  \in (0, \infty) x, \\
            \infty & \text{if } \abs{x}=1,\, x^* = 0,\, x \dir x > 0, \\
            0 & \text{if } \abs{x}=1,\, x^* = 0,\, x \dir x \le 0, \\
            0 & \text{if } \abs{ x } < 1.
        \end{cases}
    \end{equation*}
\end{lemma}
\begin{proof}
    The claims follow directly from the expression \eqref{eq:graphical:indicator:gderiv} in \cref{lemma:graphical:indicator} with $\sup \emptyset=-\infty$ and $\inf\emptyset=\infty$.
\end{proof}

\begin{lemma}
    \label{lemma:secondorder:abs}
    Let $f(x) = |x|$, $x\in \R$. Then for every $x^*\in \partial f(x)$ and $\dir x\in\R$,
    \begin{align*}
        Q_f(\dir x; x|x^*) & =
        \begin{cases}
            \infty & \text{if } x=0,\, \dir x \ne 0,\, \sign \dir x \ne x^*, \\
            0 & \text{otherwise},
        \end{cases}
        \\
        Q^f(\dir x; x|x^*) & =
        \begin{cases}
            -\infty & \text{if } x=0,\, \dir x \ne 0,\, \sign \dir x \ne x^*, \\
            0 & \text{otherwise},
        \end{cases}
    \end{align*}
    and
    \begin{equation*}
        \hat Q^f(\dir x; x|x^*)=
        Q^f_0(\dir x; x|x^*) =
        \begin{cases}
            0 & \text{if }  x  \ne 0,\, x^* = \sign x, \\
            0 & \text{if } x=0,\, \abs{x^*}=1,\, x^*\dir x \ge 0, \\
            \infty & \text{if } x=0,\, \abs{x^*}=1,\, x^*\dir x < 0,  \\
            \infty & \text{if } x=0,\, \abs{x^*}<1. \\
        \end{cases}
    \end{equation*}
\end{lemma}
\begin{proof}
    The claims follow directly from the expression \eqref{eq:graphical:absvalue:gderiv} in \cref{lemma:graphical:absvalue} with $\sup \emptyset=-\infty$ and $\inf\emptyset=\infty$.
\end{proof}

These results can be lifted to the corresponding integral functionals on $L^p(\Omega)$ using the results of \cref{chap:superposition}. Similarly, we obtain calculus rules for the curvature functionals from the corresponding results in \cref{chap:gderiv}.
\begin{theorem}[sum rule]%
    \label{thm:secondorder:sum}
    Let $X$ be a Banach space, let $f:X\to\R$ be locally Lipschitz continuous, and let $g:X\to\R$ be twice continuously differentiable. Set $j(x)\defeq f(x)+g(x)$. Then for every $x\in X$ and $x^*\in \partial_C f(x)$,
    \begin{align*}
        Q_j(\dir x; x|x^*+g'(x)) &= Q_f(\dir x; x|x^*) + \dualprod{g''(x)\dir x}{\dir x}_X\qquad(\dir x\in X),\\
        Q^j(\dir x; x|x^*+g'(x)) &= Q^f(\dir x; x|x^*) + \dualprod{g''(x)\dir x}{\dir x}_X\qquad(\dir x\in X).
    \end{align*}
\end{theorem}
\begin{proof}
    We only show the expression for the upper model, the lower model being analogous.
    First, by \cref{thm:clarke:sum}, we have $\subdiff_C j(x)= \setof{x^* +g'(x)}{x^*\in \subdiff_C f(x)}$.
    The sum rule \cref{thm:gderiv:addition} for the graphical derivative together with \cref{thm:graphical:single} then yields
    \begin{equation*}
        D[\subdiff_C j](x|x^*+g'(x))(\dir x) =
        D[\subdiff_C f](x|x^*)(\dir x) + g''(x)\dir x
    \end{equation*}
    and therefore
    \begin{equation*}
        \begin{aligned}[b]
            Q^j(\dir x; x|x^*+g'(x))
            &=
            \sup_{\dir x^* \in D[\subdiff_C j](x|x^*+g'(x))(\dir x)}
            \dualprod{\dir x^*}{\dir x}_X \\
            &=
            \sup_{\dir x^* \in D[\subdiff_C f](x|x^*)(\dir x)}
            \dualprod{\dir x^*}{\dir x}_X + \dualprod{g''(x)\dir x}{\dir x}_X .
        \end{aligned}
        \qedhere
    \end{equation*}
\end{proof}
\begin{theorem}[chain rule]%
    \label{thm:secondorder:chain}
    Let $X,Y$ be Banach spaces, let $f:Y\to \R$ be convex, and let $S:X\to Y$ be twice continuously differentiable. Set $j(x) \defeq f(S(x))$.
    If there exists a neighborhood $U$ of $x \in X$ such that
    \begin{enumerate}
        \item $f$ is Clarke regular at $S(\alt x)$ for all $\alt x\in U$;
        \item $S'(\alt x)$ has a right-inverse $\rinv{S'(\tilde x)} \in \linear(X^*; Y^*)$ for all $\tilde x\in U$;
        \item the mapping $\alt x \mapsto \rinvstar{S'(\alt x)}$ is strictly differentiable at $x$;
    \end{enumerate}
    then for all $x^*\in \partial_C j(x) = S'(x)^*\partial_C f(S(x))$,
    \begin{align*}
        Q_{j}(\dir x; x|x^*)
        &= \dualprod{y^*}{[S''(x)\dir x]\dir x}_Y + Q^f( S'(x)\dir x; S(x)|y^*)
        \qquad (\dir x \in X),\\
        Q^{j}(\dir x; x|x^*)
        &= \dualprod{y^*}{[S''(x)\dir x]\dir x}_Y + Q^f( S'(x)\dir x; S(x)|y^*)
        \qquad (\dir x \in X),
    \end{align*}
    for the unique $y^*\in \partial_C f(S(x))$ such that $S'(x)^*y^* = x^*$.
\end{theorem}
\begin{proof}
    We again only consider the upper model $Q^j$, the lower model being analogous.
    Due to our assumptions, we can apply \cref{cor:gclarke:second-clarke} to obtain
    \begin{equation*}
        D[\subdiff_C(f \circ S)](x|x^*)(\dir x)
        =
        [S''(x)^*\dir x]y^* + S'(x)^*D[\subdiff f](S(x)|y^*)(S'(x)\dir x),
    \end{equation*}
    where $S'': X \to [X \to \linear(Y^*; X^*)]$.
    Thus every $\dir x^* \in D[\subdiff_C(f \circ S)](x|x^*)(\dir x)$ can be written for some $\dir y^* \in D[\subdiff f](S(x)|y^*)(S'(x)\dir x)$ as
    $
    \dir x^*=[S''(x)\dir x]^*y^*+S'(x)^*\dir y^*.
    $
    Inserting this into the definition of $Q^j$ yields
    \begin{equation*}
        \begin{aligned}[b]
            Q^{j}(\dir x; x|x^*)
            &=\sup_{\dir y^* \in D[\subdiff f](S(x)|y^*)(S'(x)\dir x)}
            \dualprod{[S''(x)\dir x]^*y^*+S'(x)^*\dir y^*}{\dir x}_X
            \\
            &= \dualprod{y^*}{[S''(x)\dir x]\dir x}_Y
            + \sup_{\dir y^* \in D[\subdiff f](S(x)|y^*)(S'(x)\dir x)} \dualprod{\dir y^*}{S'(x)\dir x}_Y.
        \end{aligned}
        \qedhere
    \end{equation*}
\end{proof}

\section{Subconvexity}

We say that $f: X \to \Rbar$ is \term[functional!subconvex]{subconvex near} $\opt x$ for $\opt x^* \in \subdiff_C f( x)$ if for all $\rho>0$, there exists $\epsilon>0$ such that
\begin{equation}
    \label{eq:secondorder:subconvex}
    f(\alt x)-f(x) \ge \dualprod{x^*}{\alt x - x}_X - \frac{\rho}{2}\norm{\alt x - x}_X^2
    \qquad
    (x,\alt x \in \B(\opt x, \epsilon);\, x^* \in \subdiff_C f(x) \isect \B(\opt x^*, \epsilon)).
\end{equation}
We say that $f$ is \emph{subconvex at} $\opt x$ for $\opt x^*$ if this holds with $\alt x=\opt x$ fixed.
It is clear that convex functions are subconvex near any point for any subderivative.
By extension, scalar functions such as $t \mapsto \abs{t}^q$ for $q \in (0, 1)$ that are locally minorized by $\alt x \mapsto f(\optx) +  \dualprod{x^*}{\alt x - x}_X$ at points of nonsmoothness are also subconvex.

The sum of two subconvex functions for which the subdifferential sum rule holds is clearly also subconvex.
The next result shows that smooth functions simply need to have a non-negative Hessian at the point $\optx$ to be subconvex. This is in contrast to the everywhere non-negative Hessian of convex functions.

\begin{lemma}
    \label{lem:secondorder:smooth-subconvex}
    Let $X$ be a Banach space and let $f:X\to \R$ be twice continuously differentiable.
    If $\dualprod{f''(\optx)\dir x}{\dir x}_X \ge 0$ for all $\dir x\in X$,
    then $f$ is subconvex near $\opt x\in X$ for $f'(\optx)$.
\end{lemma}

\begin{proof}
    Fix $\rho>0$.
    We apply \cref{thm:frechet:mean} first to $f$ to obtain for every $x,h \in X$ that
    \begin{equation*}
        f(x+h)-f(x)=\int_0^1 \dualprod{f'(x+th)}{h}_X \,dt.
    \end{equation*}
    Similarly, the same theorem applied to $t \mapsto \dualprod{f'(x+th)}{h}$ for any $x,h \in X$ yields
    \begin{equation*}
        \iprod{f'(x+th)}{h}_X-\iprod{f'(x)}{h}_X = \int_0^1 \dualprod{f''(x+sth)h}{h}_X \,ds.
    \end{equation*}
    Combined, these two expansions yield
    \begin{equation}
        \label{eq:secondorder:smooth-proxreg:mainest}
        f(x+h)-f(x)
        = \dualprod{f'(x)}{h}_X + \int_0^1 \int_0^1 \dualprod{f''(x+sth)h}{h}_X \, ds \, dt.
    \end{equation}
    Since $\dualprod{f''(\optx)h}{h}_X \ge 0$, we have
    \begin{equation*}
        \dualprod{f''(x+q)h}{h}_X
        \ge
        \dualprod{[f''(x+q)-f''(\optx)]h}{h}_X
        \qquad (x, q, h \in X).
    \end{equation*}
    Therefore, by the continuity of $f''$, for any $\rho>0$ we can find $\epsilon>0$ such that
    \begin{equation*}
        \dualprod{f''(x+q)h}{h}_X \ge -\frac{\rho}{2}\norm{h}_X^2
        \qquad (q \in \B(0, \eps),\, x \in \B(\optx, \eps),\, h \in X).
    \end{equation*}
    Taking $q=sth$, this and \eqref{eq:secondorder:smooth-proxreg:mainest} shows that
    \begin{equation*}
        f(x+h)-f(x) \ge \dualprod{f'(x)}{h}_X -\frac{\rho}{2}\norm{h}_X^2.
    \end{equation*}
    The claim now follows by taking $h=\alt x - x$.
\end{proof}

\begin{remark}
    Subconvexity, which to our knowledge has not previously been treated in the literature, is a stronger condition than the prox-regularity introduced in \cite{poliquin1996prox}. The latter requires \eqref{eq:secondorder:subconvex} to hold merely for a fixed $\rho>0$.
    The definition in \cite{Rockafellar:1998} is slightly broader and implies the earlier one. Their definition is itself a modification of the \term[functional!primal-lower-nice]{primal-lower-nice} functions of \cite{thibault1995integration}.
    Our notion of subconvexity is also related to those of \term[set!subsmooth]{subsmooth} sets and \term[mapping!submonotone]{submonotone} operators introduced in \cite{aussel2005subsmooth}.
    An alternative concept for functions, \term[functional!subsmooth]{subsmoothness} and \term[functional!lower-$C^k$]{lower-$C^k$}, has been introduced in \cite{rockafellar1981favorable}.
\end{remark}

\section{Sufficient and necessary conditions}

We start with sufficient conditions, which are based on the upper model.

\begin{theorem}
    \label{thm:secondorder:ssc}
    Let $X$ be a Banach space and $f: X \to \Rbar$.
    If for $\optx \in X$,
    \begin{enumerate}
        \item\label{item:secondorder:ssc:subconvex}
            $f$ is subconvex \emph{near} $\opt x$ for $\opt x^*=0$;
        \item\label{item:secondorder:ssc:firstorder}
            $0 \in \subdiff_C f(\opt x)$;
        \item\label{item:secondorder:ssc} there exists a $\mu>0$ such that
            \begin{equation*}
                \hat Q^f(\dir x; \opt x|0) \ge \mu\norm{\dir x}_X^2
                \qquad (\dir x \in X);
            \end{equation*}
    \end{enumerate}
    then $\optx$ is a strict local minimizer of $f$.
\end{theorem}

\begin{proof}
    Let $\opt x^* \defeq 0$ and $\dir x\in X$.
    By the assumed subconvexity, for every $\rho>0$ there exists $\epsilon_\rho>0$ such that for $x \in \B(\opt x, \epsilon_\rho/2)$ and $x^* \in \subdiff_C f(x) \isect \B(\opt x^*, \epsilon_\rho)$, we have for every $t>0$ with $t\norm{\dir x}_X < \tfrac{1}{2}\epsilon_\rho$ that
    \begin{equation*}
        \frac{f(x + t \dir x)-f(x)-t\dualprod{\opt x^*}{\dir x}_X}{t^2}
        \ge
        \frac{\dualprod{x^*-\opt x^*}{\dir x}_X}{t}
        -\frac{\rho}{2}\norm{\dir x}_X^2.
    \end{equation*}
    Since $\rho>0$ was arbitrary, we thus obtain for every $\dir \alt x \in X$ and $\dir x^* \in D[\subdiff f](x|x^*)(\dir\alt x)$ that
    \begin{equation*}
        \begin{aligned}
            A(\dir x, \dir{\alt x}, \dir x^*) &\defeq
            \liminf_{\substack{t \downto 0,\, (x-\opt x)/t \to \dir \alt x \\ (x^*-\optx^*)/t \to \dir x^*,\, x^* \in \subdiff_C f(x)}}~
            \frac{f(x + t \dir x)-f(x)-t\dualprod{\opt x^*}{\dir x}_X}{t^2}\\
            &\ge
            \liminf_{\substack{t \downto 0,\, (x - \opt x)/t \to \dir{\alt x} \\ (x^*-\optx^*)/t \to \dir x^*,\, x^* \in \subdiff_C f(x)}}
            \frac{\dualprod{x^*-\opt x^*}{\dir x}_X}{t}
            =
            \dualprod{\dir x^*}{\dir x}_X.
        \end{aligned}
    \end{equation*}
    This implies that
    \begin{equation*}
        \sup_{\dir x^* \in D[\subdiff_C f](\opt x|\opt x^*)(\dir\alt x)}  A(\dir x, \dir\alt x, \dir x^*)
        \ge
        \sup_{\dir x^* \in D[\subdiff_C f](\opt x|\opt x^*)(\dir\alt x)} \dualprod{\dir x^*}{\dir x}_X
        =: B(\dir x, \dir\alt x).
    \end{equation*}
    Since $\opt x^*=0$, we can fix $x=\opt x+t\dir x$ and $\dir \alt x=\dir x$ in the $\liminf$ above and use \ref{item:secondorder:ssc} to obtain
    \begin{equation}
        \label{eq:secondorder:ssc:liminf1}
        \liminf_{\substack{t \downto 0}}~
        \frac{f(\opt x + 2 t \dir x)-f(\opt x+t \dir x)}{t^2}
        \ge
        B(\dir x, \dir x).
    \end{equation}
    Similarly, fixing $x=\optx$ and $\dir{\alt x}=0$ yields
    \begin{equation}
        \label{eq:secondorder:ssc:liminf2}
        \liminf_{\substack{t \downto 0}}~
        \frac{f(\optx + t \dir x)-f(\optx)}{t^2}
        \ge
        B(\dir x, 0) \ge 0,
    \end{equation}
    where the final inequality follows from the definition of $B$ by taking $\dir x^*=0$ (which is possible since $\optx^* \in \subdiff_C f(\optx)$).
    We now make a case distinction.
    \begin{enumerate}[label=(\Roman*)]
        \item $B(\dir x, 0) \ge \mu\norm{\dir x}_X^2$. In this case, the $\liminf$ is strictly positive for $\dir x\neq 0$ and hence $f(\optx + t\dir x) > f(\optx)$ for all $t>0$ sufficiently small.
        \item $B(\dir x, 0) < \mu\norm{\dir x}_X^2$. In this case, it follows from \ref{item:secondorder:ssc} that
            \begin{equation*}
                \mu\norm{\dir x}_X^2 \leq \hat Q^f(\dir x; \opt x|0)=\max\{B(\dir x, \dir x), B(\dir x, 0)\}
            \end{equation*}
            and hence that $B(\dir x, \dir x) = \hat Q^f(\dir x; \opt x|0) \ge \mu\norm{\dir x}_X^2$. Summing \eqref{eq:secondorder:ssc:liminf1} and \eqref{eq:secondorder:ssc:liminf2} then yields
            \begin{equation*}
                \liminf_{\substack{t \downto 0}}~
                \frac{f(\opt x + 2t \dir x)-f(\opt x)}{t^2}
                \ge B(\dir x, \dir x) \ge \mu \norm{\dir x}_X^2,
            \end{equation*}
            which again implies for $\dir x\neq 0$ that $f(\optx + t\dir x) > f(\optx)$ for all $t>0$ sufficiently small.
    \end{enumerate}
    Since $\dir x\in X$ was arbitrary, $\optx$ is by definition a strict local minimizer of $f$.
\end{proof}

\begin{remark}
    The use of the stationarity curvature model $Q^f_0$ in the second-order condition is required since the upper curvature model may not provide any information about the growth of $f$ at $\optx$ in certain directions. However, since $D[\subdiff_C f](\opt x|\optx^*)(0)$ is a cone, if it contains \emph{any} element $\dir x^*$ such that $\dualprod{\dir x^*}{\dir x}_X>0$, then $B(\dir x, 0)=Q^f_0(\dir x; \optx|\optx^*)=\infty$, ensuring that the condition \ref{item:secondorder:ssc} holds in the direction $\dir x$ for any $\mu>0$.
    For example, if $f(x)=\abs{x}$, then \cref{lemma:secondorder:abs} shows that $Q^f(\dir x; 0|0)=0$ for $\dir x \ne 0$, which indeed does not provide any information about the growth of $f$ at $0$. Conversely, $Q^f_0(\dir x; 0|0)=\infty$ for any $\dir x \ne 0$, so the growth is more rapid than $Q^f$ can measure.
\end{remark}

Combining \cref{thm:secondorder:ssc} with \cref{thm:secondorder:model-smooth}, we obtain the classical sufficient second-order condition. (Recall that in infinite-dimensional spaces, positive definiteness and coercivity are no longer equivalent, and the latter, stronger, property is usually required.)
\begin{corollary}
    Let $X$ be a Banach space and let $f:X\to\R$ be twice continuously differentiable.
    If for $\optx \in X$,
    \begin{enumerate}
        \item $f'(\bar x) = 0$;
        \item there exists a $\mu>0$ such that
            \begin{equation*}
                \dualprod{f''(\optx)\dir x}{\dir x}_X \ge \mu\norm{\dir x}_X^2
                \qquad (\dir x \in X);
            \end{equation*}
    \end{enumerate}
    then $\optx$ is a local minimizer of $f$.
\end{corollary}
\begin{proof}
    To apply \cref{thm:secondorder:ssc}, it suffices to note that $\partial_C f(x) = \{f'(x)\}$ by \cref{thm:clarke:frechet} and that the second-order condition ensures subconvexity of $f$ at $\opt x$ for $\optx^*=0$ by \cref{lem:secondorder:smooth-subconvex}.
\end{proof}

For nonsmooth functionals, we merely illustrate the sufficient second-order condition with a simple but nontrivial scalar example.
\begin{corollary}\label{example:secondorder:sufficient-sum}
    Let $X=\R$ and $j\defeq f+g$ for $g:\R\to\R$ twice continuously differentiable and $f(x) = |x|$. Then the sufficient condition of \cref{thm:secondorder:ssc} holds at $\optx\in\R$ with $0\in j(\bar x)$ if and only if one of the following cases holds:
    \begin{enumerate}[label=(\alph*)]
        \item\label{item:secondorder:sufficient-sum:zero-inner}
            $\optx=0$ and $\abs{g'(\optx)}<1$;
        \item\label{item:secondorder:sufficient-sum:zero-bdry}
            $\optx=0$, $\abs{g'(\optx)}=1$, and $g''(\optx)>0$; or
        \item\label{item:secondorder:sufficient-sum:nonzero}
            $\optx \ne 0$, $g'(\optx)=-\sign \optx$, and $g''(\optx)>0$.
    \end{enumerate}
\end{corollary}
\begin{proof}
    We apply \cref{thm:secondorder:ssc}, for which we need to verify its conditions. First, note that
    \cref{item:secondorder:ssc:firstorder} is equivalent to $0 = x^* + g'(\optx)$ for some $x^* \in \subdiff f(\optx) = \sign \optx$ by \cref{thm:clarke:sum} and \cref{ex:convex:subdiff_abs}.

    We now verify the subconvexity of $j$ near $\optx$ for $\optx^* = 0$. Expanding the definition \eqref{eq:secondorder:subconvex}, this requires
    \begin{multline}
        \label{eq:secondorder:subconvex:sum}
        \abs{\alt x}-\abs{x}  + g(\alt x)-g(x) \ge \dualprod{x^*+g'(x)}{\alt x - x} - \frac{\rho}{2}\norm{\alt x - x}^2
        \\
        (x,\alt x \in \B(\opt x, \epsilon);\, x^* \in \subdiff_C \abs{\freevar}(x) \isect \B(\opt x^*-g'(x), \epsilon)).
    \end{multline}
    In cases \cref{item:secondorder:sufficient-sum:zero-bdry} and \cref{item:secondorder:sufficient-sum:nonzero}, we can apply \cref{lem:secondorder:smooth-subconvex} to deduce the subconvexity of $g$ and therefore of $j=f+g$ since $f$ is convex. For case \cref{item:secondorder:sufficient-sum:zero-inner}, we have $\optx=0$ with $\abs{g'(\optx)}<1$. Since $g'$ is continuous, we consequently have $\opt x^*-g'(x) = -g'(x) \in (-1, 1)$ when $\abs{x-\optx}=\abs{x}$ is small enough. Since $\partial f(x) \in\{-1,1\}$ for $x \ne 0$, it follows that $\subdiff_C \abs{\freevar}(x) \isect \B(\opt x^*-g'(x), \epsilon)=\emptyset$ for $x \in \B(\optx, \epsilon) \setminus \{\optx\}$ for small enough $\epsilon>0$. Therefore, for small enough $\epsilon>0$, the condition \eqref{eq:secondorder:subconvex:sum} reduces to
    \begin{equation}
        \label{eq:secondorder:subconvex:sum:zero}
        \abs{\alt x}  + g(\alt x)-g(0) \ge \dualprod{x^*+g'(0)}{\alt x} - \frac{\rho}{2}\abs{\alt x}^2
        \quad
        (\alt x \in [-\eps,\eps],\,  \abs{x^*} \le 1,\, \abs{x^*+g'(0)} \le \eps).
    \end{equation}
    Furthermore, $|g'(0)|<1$ implies that for every $\rho>0$ and $c>0$, we can find an $\eps>0$ sufficiently small that
    \begin{equation*}
        (1-\eps-\abs{g'(0)})\abs{\alt x} \ge  \frac{c-\rho}{2}\abs{\alt x}^2
        \qquad (\alt x \in [-\eps, \eps]).
    \end{equation*}
    Since $g:\R\to\R$ is twice continuously differentiable, we can apply a Taylor expansion in $\optx=0$ to obtain for some $c>0$ and $|\alt x|$ sufficiently small that
    \begin{equation*}
        g(0) \leq g(\alt x) + \iprod{g'(0)}{-\alt x} + \frac{c}{2}|\alt x|^2.
    \end{equation*}
    Adding this to the previous inequality, we obtain for sufficiently small $\eps>0$ and $x^* \in [-1,1]$ satisfying $|x^*+g'(0)| \leq \eps$ that
    \begin{equation*}
        \begin{aligned}
            \abs{\alt x} + g(\alt x)-g(0)
            &\ge (|g'(0)| + \eps)|\alt x|+ \iprod{g'(0)}{\alt x} - \frac{\rho}{2}\abs{\alt x}^2\\
            &\ge \iprod{x^*+g'(0)}{\alt x} - \frac{\rho}{2}\abs{\alt x}^2
        \end{aligned}
    \end{equation*}
    for every $|\alt x|\leq \eps$, which is \eqref{eq:secondorder:subconvex:sum:zero}.
    Hence $j=f+g$ is subconvex near $\optx=0$ for $0=x^*+g'(0)$.

    To verify \cref{item:secondorder:ssc}, we compute the upper curvature model. Let $\dir x\in X$. Then by \cref{thm:secondorder:sum,thm:secondorder:model-smooth},
    \begin{align*}
        Q^{j}(\dir x; x|x^*+g'(x)) &= Q^f(\dir x; x|x^*) + \dualprod{g''(x)\dir x}{\dir x}, \\
        Q^{j}_0(\dir x; x|x^*+g'(x)) &= Q^f(\dir x; x|x^*),
    \end{align*}
    where $Q^f$ is given by \cref{lemma:secondorder:abs}. It follows that
    \begin{align*}
        Q^{j}(\dir x; x|x^*+g'(x)) &= \begin{cases}
            -\infty & \text{if } x=0,\, \dir x \ne 0,\, \sign \dir x \ne x^*, \\
            \dualprod{g''(x)\dir x}{\dir x} & \text{otherwise},
        \end{cases}
        \shortintertext{and}
        Q^{j}_0(\dir x; x|x^*+g'(x)) & = \begin{cases}
            0 & \text{if }  x  \ne 0, x^* = \sign x, \\
            0 & \text{if } x=0,\, \abs{x^*}=1,\, x^*\dir x \ge 0, \\
            \infty & \text{if } x=0,\, \abs{x^*}=1,\, x^*\dir x < 0, \\
            \infty & \text{if } x=0, \abs{x^*}<1. \\
        \end{cases}
    \end{align*}
    Thus
    \begin{equation*}
        \hat Q^j(\dir x; x|x^*+g'(x)) = \begin{cases}
            \max\{0, \dualprod{g''(x)\dir x}{\dir x}\} & \text{if } x \ne 0,\, x^* = \sign x, \\
            \max\{0, \dualprod{g''(x)\dir x}{\dir x}\} & \text{if } x=0,\, \abs{x^*}=1,\, x^*\dir x \ge 0, \\
            \infty & \text{if } x=0,\, \abs{x^*}=1,\, x^*\dir x < 0, \\
            \infty & \text{if } x=0,\, \abs{x^*}<1.
        \end{cases}
    \end{equation*}
    The condition \cref{item:secondorder:ssc} is thus equivalent to
    \begin{equation*}
        \max\{0, \dualprod{g''(\optx)\dir x}{\dir x}\} \ge \mu\norm{\dir x}^2
        \quad
        \text{when}
        \quad
        \begin{cases}
            \optx \ne 0 \text{ or}\\
            \optx = 0,\, \abs{g'(\optx)}=1,\, \text{ and } g'(\optx)\dir x<0.
        \end{cases}
    \end{equation*}
    The left inequality can only hold for arbitrary $\dir x\in\R$ if $\mu = g''(\optx) >0$.
    Hence \cref{item:secondorder:ssc:firstorder} and \cref{item:secondorder:ssc} hold if and only if one of the cases \crefrange{item:secondorder:sufficient-sum:zero-inner}{item:secondorder:sufficient-sum:nonzero} holds.
\end{proof}
Note that case \cref{item:secondorder:sufficient-sum:zero-inner} corresponds to the case of strict complementarity or graphical regularity of $\partial f$ in \cref{lemma:graphical:absvalue}. Conversely, cases \cref{item:secondorder:sufficient-sum:zero-bdry,item:secondorder:sufficient-sum:nonzero} imply that $g$ and therefore $j$ is locally convex, recalling from \cref{thm:convex:fermat} that for convex functionals, the first-order optimality conditions are necessary \emph{and} sufficient.

\bigskip

Now we formulate our necessary condition, which is based on the \emph{lower} curvature model.

\begin{theorem}
    \label{thm:secondorder:snc}
    Let $X$ be a Banach space and $f: X \to \Rbar$.
    If $\optx \in X$ is a local minimizer of $f$ and $f$ is locally Lipschitz continuous and subconvex at $\optx$ for $0\in X^*$, then
    \begin{equation*}
        Q_f(\dir x; \opt x| 0) \ge 0
        \qquad (\dir x \in X).
    \end{equation*}
\end{theorem}

\begin{proof}
    We have from \cref{thm:clarke:fermat} that $\opt x^* \defeq 0 \in \subdiff_C f(\opt x)$.
    By the assumed subconvexity, for every $\rho>0$ there exists $\epsilon>0$ such that for $x \in \B(\opt x, \epsilon/2)$ and  $x_t^* \in \subdiff_C f(\opt x+t\dir \alt x) \isect \B(\opt x^*, \epsilon)$, we have for every $t>0$ with $t\norm{\dir x}_X < \epsilon/2$ that
    \begin{equation*}
        \frac{f(\opt x + t \dir \alt x)-f(\opt x)-t\dualprod{\opt x^*}{\dir \alt x}_X}{t^2}
        \le
        \frac{\dualprod{x_t^*-\opt x^*}{\dir \alt x}_X}{t}
        +\frac{\rho}{2}\norm{\dir \alt x}_X^2.
    \end{equation*}
    For every $\dir x^* \in D[\subdiff_C f](\opt x|\opt x^*)(\dir x)$, by definition there exist $\dir \alt x \to \dir x$ and, for small enough $t>0$, $x_t^* \in \subdiff_C f(x+t\dir \alt x) \isect \B(\opt x^*, \epsilon)$ such that $(x_t^*-\opt x^*)/t \to \dir x^* \in X^*$.
    Since $\rho>0$ was arbitrary and $\optx^*=0$, it follows that
    \begin{equation*}
        \begin{aligned}
            \liminf_{\substack{\dir \alt x \to \dir x\\ t \downto 0}}~
            \frac{f(\opt x + t \dir \alt x)-f(\opt x)}{t^2}
            &
            \le \liminf_{\substack{\dir \alt x \to \dir x\\ t \downto 0}}\left(
                \frac{\dualprod{x_t^*-\opt x^*}{\dir x}_X}{t}
                +
                \frac{\dualprod{x_t^*-\opt x^*}{\dir \alt x-\dir x}_X}{t}
            \right)
            \\
            &
            =
            \liminf_{t \downto 0}~
            \frac{\dualprod{x_t^*-\opt x^*}{\dir x}_X}{t}
            \\
            &
            \le \inf_{\dir x^* \in D[\subdiff_C f](\opt x|\opt x^*)(\dir x)} \dualprod{\dir x^*}{\dir x}_X
            \\
            &
            = Q_f(\dir x; \opt x|\optx^*)
            = Q_f(\dir x; \opt x|0).
        \end{aligned}
    \end{equation*}
    Since $\opt x$ is a local minimizer, we have $f(\optx)\leq f(\optx +t \dir\alt x)$ for $t>0$ sufficiently small and $\dir\alt x$ sufficiently close to $\dir x$. Rearranging and passing to the limit thus yields the claimed nonnegativity of $Q_f(\dir x;\optx |0)$.
\end{proof}

\begin{remark}
    Compared to the sufficient condition of \cref{thm:secondorder:ssc}, the necessary condition does not involve a \term[model!stationary lower]{stationary lower model}
    \begin{equation*}
        Q_{f,0}(\dir x; \optx|0)\defeq  \inf_{\dir x^* \in D[\subdiff_C f](\opt x|0)(0)} \dualprod{\dir x^*}{\dir x}_X.
    \end{equation*}
    In fact, $Q_{f,0}(\dir x; \optx|0)\geq 0$ is \emph{not} a necessary optimality condition:
    let $f(x)=\abs{x}$, $x\in \R$, and $\optx=0$. Then by \cref{lemma:graphical:absvalue}, $D[\subdiff f](0|0)(0)=\R$ and hence $Q_{f,0}(\dir x; 0|0) = -\infty$ for all $\dir x\ne 0$.
\end{remark}

For smooth functions, we recover the usual second-order necessary condition from \cref{thm:secondorder:model-smooth}.
\begin{corollary}
    \label{thm:secondorder:snc-smooth}
    Let $X$ be a Banach space and let $f: X \to \R$ be twice continuously differentiable.
    If $\optx \in X$ is a local minimizer of $f$, then
    \begin{equation*}
        \dualprod{f''(\optx)\dir x}{\dir x}_X \ge 0
        \qquad (\dir x \in X).
    \end{equation*}
\end{corollary}

We again illustrate the nonsmooth case with a scalar example.
\begin{corollary}\label{example:secondorder:necessary-sum}
    Let $X=\R$ and $j\defeq f+g$ for $g:\R\to\R$ twice continuously differentiable and $f(x) = |x|$. Then the necessary condition of \cref{thm:secondorder:snc} holds at a minimizer $\optx\in \R$ of $j$ if and only if $g''(\optx)\geq 0$.
\end{corollary}
\begin{proof}
    We apply \cref{thm:secondorder:snc}, for which we need to verify its conditions.
    Both $f$ and $g$ are locally Lipschitz continuous by \cref{thm:convex:cont} and \cref{lem:variation:c1-lipschitz}, respectively, and hence so is $j$. We have already verified the subconvexity of $j$ in \cref{example:secondorder:sufficient-sum}.

    By \cref{thm:clarke:fermat,thm:clarke:sum,ex:convex:subdiff_abs}, we again have $0= x^* + g'(\optx)$ for some $x^* \in \subdiff f(\optx) = \sign \optx$.
    It remains to compute the lower curvature model. Let $\dir x\in X$.
    By \cref{thm:secondorder:sum,thm:secondorder:model-smooth},
    \begin{equation*}
        Q_{j}(\dir x; x|x^*+g'(x)) = Q_f(\dir x; x|x^*) + \dualprod{g''(x)\dir x}{\dir x},
    \end{equation*}
    where $Q_f$ is given by \cref{lemma:secondorder:abs}. It follows that
    \begin{equation*}
        Q_{j}(\dir x; x|x^*+g'(x)) =
        \begin{cases}
            \infty & \text{if }x=0,\, \dir x \ne 0,\, \sign \dir x \ne x^*, \\
            \dualprod{g''(x)\dir x}{\dir x} & \text{otherwise}.
        \end{cases}
    \end{equation*}
    Hence the condition $Q_{j}(\dir x; \optx|0) \ge 0$ for all $\dir x \in X$ reduces to $g''(\optx) \ge 0$.
\end{proof}

\begin{remark}
    Second-order optimality conditions can also be based on \term[derivative!epigraphical]{epigraphical derivatives}, which were introduced in \cite{rockafellar1985maximal,rockafellar1988first};
    we refer to \cite{Rockafellar:1998} for a detailed discussion.
    A related approach based on second-order directional curvature functionals was used in \cite{christof2018nogap} for deriving necessary and sufficient second-order optimality conditions for smooth optimization problems subject to nonsmooth and possibly nonconvex constraints.
\end{remark}

\chapter{Lipschitz-like properties}
\label{chap:regularity}

A related issue to second-order conditions is that of \term{stability} of the solution to optimization problems under perturbation. To motivate the following, let $f:X\to\Rbar$ and suppose we wish to find $\optx\in X$ such that $0\in \partial f(\optx)$ for a suitable subdifferential. Suppose further that we are given some $\alt x\in X$ with $w\in \partial f(\alt x)$ with $\norm{w}_{X^*}\leq \eps$; say, from one of the algorithms in \cref{chap:proximal}. A natural question is then for an error estimate $\norm{\optx - \alt x}_X$ in terms of $\epsilon$. Clearly, if $\partial f$ has a single-valued and Lipschitz continuous inverse, this is the case since then
\begin{equation*}
    \norm{\optx - \alt x}_X = \norm{(\partial f)^{-1}(0) - (\partial f)^{-1}(w)}_{X} \leq L \norm{w}_{X^*}.
\end{equation*}
Of course, the situation is much more complicated in the set-valued case. To treat this, we first have to define suitable notions of Lipschitz-like behavior of set-valued mappings, which we then characterize using coderivatives (generalizing the characterization of the Lipschitz constant of a differentiable single-valued mapping through the norm of its derivative). We return to the question of stability of minimizers in the more general context of perturbations of parametrized solution mappings in \cref{chap:stability}.

\section{Lipschitz-like properties of set-valued mappings}%
\label{sec:regularity:lipschitz}

To set up the definition of Lipschitz-like properties for set-valued mappings, it is helpful to recall from \cref{sec:functan:normed} for single-valued functions the distinction between (point-based) local Lipschitz continuity \emph{at} a point and (neighborhood-based) local Lipschitz continuity \emph{near} a point. (\Cref{fig:regularity:aubin-structure} below shows a function that is locally Lipschitz at but not near the given point.) Similarly, we will have to distinguish for set-valued mappings the corresponding notions of the \term[property!Aubin]{Aubin property} (which is point-based) and \term{calmness} (which is neighborhood-based). If these properties hold for the inverse of a mapping, we will call the mapping itself \term[mapping!regular!metrically]{metrically regular} and \term[mapping!subregular, metrically]{metrically subregular}, respectively. These four properties are illustrated in \cref{fig:regularity:cone}.

\begin{figure}[t]
    \centering
    \begin{subfigure}[t]{.45\textwidth}
        \centering
        \begin{asy}
            unitsize(80, 80);

            real ell=0.3;
            real d=0.3;
            real f(real x) { return sin(x)-x/2-sin(3*x)/3+((x+0.5)/1.4)^8; };

            pair z1=(0, 0);
            pair z2=z1;

            fill((z1+(d, ell))--z1--(z1+(-d, ell))--(z2+(-d, -ell))--z2--(z2+(d, -ell))--cycle, lightfill);

            path fg=graph(f, -0.9, 0.75, n=1000);
            draw(fg, primalline + linewidth(1.5));

            dot(z2);
            label("$(\optx, f(\optx))$", z2, 3*S);
        \end{asy}
        \caption{locally Lipschitz $f$}
        \label{fig:regularity:cone:lipschitz}
    \end{subfigure}
    \begin{subfigure}[t]{.45\textwidth}
        \centering
        \begin{asy}
            unitsize(80, 80);

            real ell=0.5;
            real d=0.2;
            real f(real x) { return 2*x^3+x; };

            pair z1=(0, 0);
            pair z2=z1;

            fill((z1+(ell, d))--z1--(z1+(ell, -d))--(z2+(-ell, -d))--z2--(z2+(-ell, d))--cycle, lightfill);

            path fg=graph(f, -0.3, 0.3, n=1000);
            draw(fg, primalline + linewidth(1.5));

            dot(z2);
            label("$(\opty, f^{-1}(\opty))$", z2, 3*E);
        \end{asy}
        \caption{locally Lipschitz $\inv f$}
         \label{fig:regularity:cone:invlipschitz}
    \end{subfigure}
    \begin{subfigure}[t]{.45\textwidth}
        \centering
        \begin{asy}
            unitsize(80, 80);

            real ell=0.15;
            real d=0.2;

            pair z1=(-0.6, -0.4);
            pair z2=z1;

            fill((z1+(d, ell))--z1--(z1+(-d, ell))--(z2+(-d, -ell))--z2--(z2+(d, -ell))--cycle, lightfill);

            dot(z2);
            label("$(\optx, \opty)$", z2, 3*S);

            z1=(0, 0.4);
            z2=z1;
            d=0.2;
            ell=0.15;

            fill((z1+(d, ell))--z1--(z1+(-d, ell))--(z2+(-d, -ell))--z2--(z2+(d, -ell))--cycle, violation);

            dot(z2);
            label("$(\check x, \check y)$", z2, 3*N);

            draw((-1, -0.4)--(0, -0.4)--(0, 0.4)--(1, 0.4), primalline + linewidth(1.5));
        \end{asy}
        \caption{Aubin property of $\subdiff\abs{\freevar}$}
        \label{fig:regularity:cone:abs:aubin}
    \end{subfigure}
    \begin{subfigure}[t]{.45\textwidth}
        \centering
        \begin{asy}
            unitsize(80, 80);

            real ell=0.15;
            real d=0.2;

            pair z1=(-0.6, -0.4);
            pair z2=z1;

            fill((z1+(ell, d))--z1--(z1+(ell, -d))--(z2+(-ell, -d))--z2--(z2+(-ell, d))--cycle, violation);

            dot(z2);
            label("$(\check x, \check y)$", z2, 2*S+3*E);

            z1=(0, 0);
            z2=z1;
            d=0.2;
            ell=0.15;

            fill((z1+(ell, d))--z1--(z1+(ell, -d))--(z2+(-ell, -d))--z2--(z2+(-ell, d))--cycle, lightfill);

            dot(z2);
            label("$(\optx, \opty)$", z2, 3*E);

            draw((-1, -0.4)--(0, -0.4)--(0, 0.4)--(1, 0.4), primalline + linewidth(1.5));
        \end{asy}
        \caption{metric regularity of $\subdiff\abs{\freevar}$}
         \label{fig:regularity:cone:abs:metric}
    \end{subfigure}
    \begin{subfigure}[t]{.45\textwidth}
        \centering
        \begin{asy}
            unitsize(80, 80);

            real ell=0.15;
            real d=0.2;

            pair z1=(-0.6, -0.4);
            pair z2=z1;

            fill((z1+(d, ell))--z1--(z1+(-d, ell))--(z2+(-d, -ell))--z2--(z2+(d, -ell))--cycle, lightfill);

            dot(z2);
            label("$(\optx, \opty)$", z2, 3*S);

            z2=(0, 0.4);
            z1=(0, -0.4);
            d=0.15;
            ell=0.15;

            fill((z2+(d, ell))--z2--(z2+(-d, ell))--(z1+(-d, -ell))--z1--(z1+(d, -ell))--cycle, lightfill);

            dot(z2);
            label("$(\check x, \check y)$", z2, 3*N);

            draw((-1, -0.4)--(0, -0.4)--(0, 0.4)--(1, 0.4), primalline + linewidth(1.5));
        \end{asy}
        \caption{calmness of $\subdiff\abs{\freevar}$}
         \label{fig:regularity:cone:abs:calm}
    \end{subfigure}
    \begin{subfigure}[t]{.45\textwidth}
        \centering
        \begin{asy}
            unitsize(80, 80);

            real ell=0.15;
            real d=0.15;

            pair z2=(-0.6, -0.4);
            pair z1=(0, -0.4);
            pair z3=(-1, -0.4);

            fill((z1+(ell, d))--z1--(z1+(ell, -d))--(z3+(0, -d))--z3--(z3+(0, d))--cycle, lightfill);

            dot(z2);
            label("$(\check x, \check y)$", z2, 2*S);

            z1=(0, 0);
            z2=z1;
            d=0.15;
            ell=0.15;

            fill((z1+(ell, d))--z1--(z1+(ell, -d))--(z2+(-ell, -d))--z2--(z2+(-ell, d))--cycle, lightfill);

            dot(z2);
            label("$(\optx, \opty)$", z2, 3*E);

            draw((-0.98, -0.4)--(0, -0.4)--(0, 0.4)--(1, 0.4), primalline + linewidth(1.5));
        \end{asy}
        \caption{metric subregularity of $\subdiff\abs{\freevar}$}
         \label{fig:regularity:cone:abs:sub}
    \end{subfigure}
    \caption{Illustration of Lipschitz-like properties using cones. The thick lines are the graph of the function; if this graph is locally contained in a light blue cone, the property holds, while a green cone indicates that the property is violated.}
    \label{fig:regularity:cone}
\end{figure}

To make these definitions precise, recall from \cref{thm:epsilon:projection} the definition of the distance of a point $x \in X$ to a set $A \subset X$, which we here write for the sake of convenience as
\begin{equation*}
    \dist(A, x) \defeq \dist(x, A) \defeq d_A(x) = \inf_{\alt x \in A} \norm{x-\alt x}_X.
\end{equation*}
We then say that $F: X \setto Y$ has the \term[property!Aubin]{Aubin} or \term[property!pseudo-Lipschitz]{pseudo-Lipschitz property} \emph{at $\optx$ for $\opty$} if $\graph F$ is closed near $(\optx, \opty)$ and there exist $\delta,\kappa>0$ such that
\begin{equation}
    \label{eq:regularity:aubin-distance}
    \dist(y, F(x)) \le \kappa \dist(\inv F(y), x)
    \quad (x \in \B(\optx, \delta),\, y \in \B(\opty, \delta)).
\end{equation}
We call the infimum of all $\kappa>0$ for which \eqref{eq:regularity:aubin-distance} holds for some $\delta>0$ the \term[modulus!graphical]{graphical modulus} of $F$ at $\optx$ for $\opty$, written $\lip F(\optx|\opty)$.

When we are interested in the stability of the optimality condition $0 \in F(\optx)$, it is typically more beneficial to study the Aubin property of the inverse $\inv F$. This is called the \term[regularity!metric]{metric regularity} of $F$ at a point $(\optx, \opty) \in \graph F$, which holds if there exist $\kappa, \delta>0$ such that
\begin{equation}
    \label{eq:regularity:metric}
    \dist(x, \inv F(y)) \le \kappa \dist(y, F(x))
    \quad (x \in \B(\optx, \delta),\, y \in \B(\opty, \delta)).
\end{equation}
We call the infimum of all $\kappa>0$ for which \eqref{eq:regularity:metric} holds for some $\delta>0$ the \term[modulus!of metric regularity]{modulus of metric regularity} of $F$ at $\optx$ for $\opty$, written $\reg F(\optx|\opty)$.

Metric regularity and the Aubin property are too strong to be satisfied in many applications. A weaker notion is provided by \term[subregularity, metric]{(metric) subregularity}
at $(\optx,\opty)\in\graph F$, which holds if there exist $\kappa,\delta>0$ such that
\begin{equation}
    \label{eq:regularity:subregularity}
    \dist(x, \inv F(\opty)) \le \kappa \dist(\opty, F(x))
    \quad (x \in \B(\optx, \delta)).
\end{equation}
Compared to metric regularity, this allows much more leeway for $F$ by fixing $y=\opty\in F(\optx)$ (while still allowing $x$ to vary).
We call the infimum of all $\kappa>0$ for which \eqref{eq:regularity:subregularity} holds for some $\delta>0$ the \term[modulus!of metric subregularity]{modulus of (metric) subregularity} of $F$ at $\optx$ for $\opty$, written $\subreg F(\optx|\opty)$.

The counterpart of metric subregularity that relaxes the Aubin property is known as \term{calmness}. We say that $F: X \setto Y$ is calm at $\optx$ for $\opty$ if there exist $\kappa,\delta>0$ such that
\begin{equation}
    \label{eq:regularity:calmness-distance}
    \dist(y, F(\optx)) \le \kappa \dist(\optx, \inv  F(y))
    \quad (y \in \B(\opty,\delta)).
\end{equation}
We call the infimum of all $\kappa>0$ for which \eqref{eq:regularity:calmness-distance} holds for some $\delta>0$ the \term[modulus!of calmness]{modulus of calmness} of $F$ at $\optx$ for $\opty$, written $\calm F(\optx|\opty)$.
Clearly the Aubin property implies calmness, while metric regularity implies metric subregularity.

\begin{example}
    \label{ex:regularity:cone}
\crefrange{fig:regularity:cone:abs:aubin}{fig:regularity:cone:abs:sub} illustrate (violation of) the Aubin property, metric regularity, calmness, and subregularity of the subdifferential mapping $F(x)=\subdiff\abs{\freevar}(x)$ of the absolute value function.
    Let us verify that these properties are indeed satisfied (with arbitrary modulus $\kappa > 0$) or violated as claimed.
    \begin{enumerate}
        \item
        The Aubin property (\cref{fig:regularity:cone:abs:aubin}) holds at the point $\optx<0$ for $\opty=-1$ because locally $F(x)=\{-1\}$ is a constant single-valued function: the left-hand side of \eqref{eq:regularity:aubin-distance} is zero if we take $\delta>0$ small enough.

        At the point $\check x=0$ for $\check y=1$, we can see that the Aubin property is violated by taking $y=1$ and $x=-\delta$ since in this case $\dist(y, F(x))=2$ but $\dist(\inv F(y), x)=\dist([0,\infty), -\delta)=\delta$.
        Hence \eqref{eq:regularity:aubin-distance} cannot hold for any $\kappa>0$.

    \item
        $F$ is calm  (\cref{fig:regularity:cone:abs:calm}) at the point $\optx<0$ for $\opty=-1$ because the Aubin property holds at this point.

        $F$ is also calm at $\check x=0$ for $\check y=1$: For any $y\in \range F = [-1,1]=F(\check x)$, the left-hand side of \eqref{eq:regularity:calmness-distance} is zero, while for $y\notin\range F$ the right-hand side is infinite.

    \item
        Metric regularity (\cref{fig:regularity:cone:abs:metric}) holds at the point $\optx=0$  for $\opty \in (-1, 1)$.
        Indeed, take $x \in \B(\optx, \delta) \setminus \{\optx\}$ and $y \in \B(\opty, \delta)$.
        Then, for small enough $\delta>0$, we have $\inv F(y)=\{0\}$.
        Hence $\dist(x, \inv F(y))=\abs{x} \le \delta$.
        Clearly then \eqref{eq:regularity:metric} holds if $x=\optx$.
        On the other hand, if $x \ne 0$,
        \[
            \dist(y, F(x))\ge\dist(\opty, F(x))-\abs{y-\opty} \ge
            \begin{cases}
                \abs{\opty+1} - \delta, & x < 0, \\
                \abs{\opty-1} - \delta, & x > 0. \\
            \end{cases}
        \]
        Therefore, \eqref{eq:regularity:metric} holds if
        \[
            \begin{cases}
                {1+\opty} \ge 2\delta & \text{when } x < 0, \\
                {1-\opty} \ge 2\delta & \text{when } x > 0.
            \end{cases}
        \]
        This can be satisfied for small enough $\delta>0$.

        However, metric regularity is violated at the point $\check x<0$ for $\check y=-1$:
        For small enough $\delta>0$ and $x \in \B(\check x, \delta)$, we have $F(x)=\{-1\}$ while taking any $1 > y > -1=\check y$ yields $\dist(\inv F(y), x)=\dist(\{0\}, x)=\abs{x}$.
        Hence \eqref{eq:regularity:metric} requires $\abs{x} \le \kappa\abs{y+1}$, which cannot hold for any $\kappa>0$ as $y \downto \check y=-1$.

    \item
        Since metric regularity implies subregularity  (\cref{fig:regularity:cone:abs:sub}), the latter holds at the point $\optx=0$ for $\opty \in (-1, 1)$.

        Metric subregularity also holds at $\check x<0$ for $\check y=-1$:
        Since $\inv F(\check y)=(-\infty, 0]$, we have for $\delta>0$ small enough that $\B(\check x, \delta) \subset \inv F(\check y)$, which implies that the left-hand side of \eqref{eq:regularity:subregularity} is zero.
    \end{enumerate}
\end{example}

\begin{remark}\label{rem:strong-metric-subregularity}
    The Aubin property is due to \cite{aubin1984lipschitz}, whereas metric subregularity is due to \cite{ioffe1979regular}, first given the modern name in \cite{dontchev2004regularity}.
    Calmness was introduced in \cite{robinson1981continuity} as the \term[property!upper Lipschitz]{upper Lipschitz property}.
    Metric regularity is equivalent to \term[openness, at a linear rate]{openness at a linear rate} near $(\realoptu, \realoptw)$ and holds for smooth maps by the classical Lyusternik--Graves theorem.
    We refer in particular to \cite{dontchev2014implicit,ioffe2017variational} for further information on these and other related properties.

    In particular, related to metric subregularity is the stronger concept of \term[subregularity, metric!strong]{strong metric subregularity}, which was introduced in \cite{rockafellar1989proto} and requires the existence of $\kappa,\delta>0$ such that
    \begin{equation*}
        \norm{x-\optx}_X \le \kappa \dist(\opty, F(x))
        \quad (x \in \B(\optx, \delta)),
    \end{equation*}
    i.e., a bound on the norm distance to $\optx$ rather than the closest preimage of $\opty$.
    Its many properties are studied in \cite{cibulka2018strong}, which also introduced $q$-exponent versions. Particularly worth noting is that strong metric subregularity is invariant with respect to perturbations by smooth functions, while metric subregularity is not.

    Weaker and ``partial'' concepts of regularity have also been considered in the literature. Of particular note is the directional metric subregularity of \cite{gfrerer2013directional}. The idea here is to study necessary optimality conditions by requiring metric regularity or subregularity only along critical directions instead of all directions.
    In \cite{tuomov-subreg}, by contrast, the norms in the definition of subregularity are made operator-relative to study the partial subregularity on subspaces; compare the testing of algorithms for structured problems in \cref{sec:testing:structured}.
\end{remark}

\begin{figure}
    \centering
    \begin{subfigure}[t]{.4\textwidth}
        \centering
        \begin{asy}
            unitsize(80, 80);

            real ell=0.5;
            real f0(real x) { if(x==0) { return 0; } else { return sin(1/x); } };
            real f(real x) { return f0(x/5)*ell*abs(x); };

            path fg=graph(f, -1, 1, n=1000);

            draw(fg, primalline);

            pair z=(0, f(0));

            draw(z--(z+(1, ell)), linewidth(1.1));
            draw(z--(z+(-1, ell)), linewidth(1.1));
            draw(z--(z+(1, -ell)), linewidth(1.1));
            draw(z--(z+(-1, -ell)), linewidth(1.1));

            dot(z); label("$\optx$", z, N);
        \end{asy}
        \caption{oscillating single-valued function}
        \label{fig:regularity:twopoint-lipschitz}
    \end{subfigure}
    \hfil
    \begin{subfigure}[t]{.49\textwidth}
        \centering
        \begin{asy}
            unitsize(80, 80);

            real ell=0.3;

            pair z1=(0, .2);
            pair z2=(0, -.2);

            fill((z1+1.1*(1, ell))--z1--(z1+1.1*(-1, ell))--(z2+1.1*(-1, -ell))--z2--(z2+1.1*(1, -ell))--cycle, violation);

            draw(z1--(z1+(1, ell)), linewidth(1.1));
            draw(z1--(z1+(-1, ell)), linewidth(1.1));
            draw(z2--(z2+(1, -ell)), linewidth(1.1));
            draw(z2--(z2+(-1, -ell)), linewidth(1.1));

            draw((z1+.5*(1, ell))--(z2+.5*(1, -ell)), dashed, Bars);
            label("$\scriptstyle F(x)+\B(0, \ell\norm{\alt x-x}_X)$", (.5, 0), E);

            label("$\scriptstyle F(x)$", (0, 0), W);

            label("$x$", (0, -.6), N);
            label("$\tilde x$", (.5, -.6), N);

            draw(z1--z2, dotted);
        \end{asy}
        \caption{graph of $x \mapsto F(x)+\B(0, \ell\norm{\alt x-x}_X)$}
        \label{fig:regularity:aubin-structure}
    \end{subfigure}
    \caption{The oscillating example in (\subref{fig:regularity:twopoint-lipschitz}) illustrates a function $f$ that is locally Lipschitz (or calm) \emph{at} $\optx$, but not locally Lipschitz (or does not have the Aubin property) \emph{near} the same point: the graph of the function stays in the cone formed by the thick lines and based at $(\optx, f(\optx)) \in \graph f$. If, however, we move the cone locally along the graph, even increasing its width, the graph will not be contained in the cone.
        In (\subref{fig:regularity:aubin-structure}) we illustrate the ``fat cone'' structure $\graph(\alt x \mapsto F(x)+\B(0, \ell\norm{\alt x-x}_X))$ appearing on the right-hand side in \cref{thm:regularity:aubin-equiv}\,\cref{item:regularity:aubin-inclusion}, and varying with the second base point $x$ around $\optx$. This is to be contrasted with the leaner cone $\graph(\alt x \mapsto f(\optx)+\B(0, \ell\norm{\alt x-\optx}_X))$ bounding the function in (\subref{fig:regularity:twopoint-lipschitz}).
    }
\end{figure}

Unfortunately, the direct calculation of the different moduli is often infeasible in practice. Much of the rest of this chapter concentrates on calculating the graphical modulus and the modulus of metric regularity in special cases. We will further consider metric subregularity (as well as a related, weaker, notion of strong submonotonicity) in \cref{sec:faster:subregularity}.

\bigskip

We start by providing alternative characterizations of the Aubin property and of calmness. These extend to metric regularity and subregularity, respectively, by application to the inverse.

The right-hand side of the set-inclusion characterization \cref{item:regularity:aubin-inclusion} in the next theorem forms a ``fat cone'' that we illustrate in \cref{fig:regularity:aubin-structure}. It should locally at each base point $x$ around $\optx$ bound $F$ for the Aubin property to be satisfied. Based on the formulation \cref{item:regularity:aubin-inclusion} below, we illustrate in \cref{fig:regularity:aubin} the satisfaction and dissatisfaction of the Aubin property.
The other two new characterizations show that we do not need to restrict $\alt x$ to a tiny neighborhood of $\optx$ in neither \cref{item:regularity:aubin-inclusion} nor the original characterization \eqref{eq:regularity:aubin-distance}.

\begin{theorem}
    \label{thm:regularity:aubin-equiv}
    Let $X,Y$ be Banach spaces and $F: X \setto Y$.
    Then the following are equivalent for $\optx\in X$ and $\opty\in F(\optx)$:
    \begin{enumerate}
        \item\label{item:regularity:aubin-inclusion} There exists $\kappa,\delta>0$ such that
            \begin{equation*}
                F(\alt x) \isect \B(\opty, \delta) \subset F(x) +  \B(0, \kappa \norm{\alt x-x}_X) \quad (\alt x, x \in \B(\optx, \delta)).
            \end{equation*}
        \item\label{item:regularity:aubin-inclusion-ext} There exists $\kappa,\delta>0$ such that
            \begin{equation*}
                F(\alt x) \isect \B(\opty, \delta) \subset F(x) +  \B(0, \kappa \norm{\alt x-x}_X) \quad (x \in \B(\optx, \delta);\, \alt x \in X).
            \end{equation*}
        \item\label{item:regularity:aubin-distance} The Aubin property \eqref{eq:regularity:aubin-distance}.
        \item\label{item:regularity:aubin-distance-restr} There exists $\kappa,\delta>0$ such that
            \begin{equation*}
                \dist(y, F(x)) \le \kappa \dist(\inv F(y) \isect \B(\optx, \delta), x)
                \quad (x \in \B(\optx, \delta),\, y \in \B(\opty, \delta)).
            \end{equation*}
    \end{enumerate}
    The infimum of $\kappa>0$ for which each of these characterizations holds is equal to the graphical modulus $\lip F(\optx|\opty)$.
    (The radius of validity $\delta>0$ for any given $\kappa>0$ may be distinct in each of the characterizations, however.)
\end{theorem}
\begin{figure}[t]
    \centering
    \begin{subfigure}[t]{0.45\textwidth}
        \centering
        \begin{asy}
            unitsize(100, 100);
            real fup0(real x){ return (-x+0.2)^3*0.45+x/3; };
            real fdown0(real x) { return fup0(x)-.6-(x-0.05)^2*1.1-x*1.5; };
            real x0=-.8;
            real mul(real x) { return (x>0 ? 2.1*x : max(-.1, 2.1*x+x*x*10)); };
            real fup(real x) { return 0.1+fup0(x0+mul(x-x0))/2.5; };
            real fdown(real x) { return fdown0(x0+mul(x-x0))/2.5; };
            real delta=.45;
            real deltaprime=.15;
            real deltaell=0.05;
            real ell=max(abs(-fup(x0+deltaell)+fup(x0)),
            abs(fdown(x0+deltaell)-fdown(x0)),
            abs(-fup(x0-deltaell)+fup(x0)),
            abs(fdown(x0-deltaell)-fdown(x0)))/deltaell;

            real epsilon=1*ell*(delta-2*deltaprime);

            path pup=graph(fup, -1, 1);
            path pdown=graph(fdown, -1, 1);
            path pupext=graph(fup, -1.1, 1.15);
            path pdownext=graph(fdown, -1.15, 1.15);

            fill(pupext--reverse(pdownext)--cycle, lightfill);
            draw(pup);
            draw(pdown);

            real maxf=fup(x0)+.3;
            real minf=fdown(x0)-.6;
            clip((-2, maxf)--(2, maxf)--(2, minf)--(-2, minf)--cycle);

            pair up0=pt(fup, x0);
            pair down0=pt(fdown, x0);

            pair xbar=(x0, 1/3*fup(x0)+2/3*fdown(x0));
            dot(xbar);
            label("$\scriptstyle (\optx, \opty)$", xbar, W);

            draw((x0-delta, -.5)--(x0+delta, -.5), dotted, Bars);
            label("$\scriptstyle \B(\optx, \delta)$", (x0, -.5), S);

            draw((x0-deltaprime, -.3)--(x0+deltaprime, -.3), dotted, Bars);
            label("$\scriptstyle \B(\optx, \delta')$", (x0, -.3), 1.5*S);

            real dd=delta*3/2;
            real dd0=deltaprime*2/3;

            real tmp=x0+dd0;
            pair x=(tmp, 1/5*fup(tmp)+4/5*fdown(tmp));
            pair up=(tmp, fup(tmp));
            pair down=(tmp, fdown(tmp));
            pair up2=up+.7*(1, ell);
            pair down2=down+.7*(1, -ell);
            pair up3=up+(dd-dd0)*(1, ell);
            pair down3=down+(dd-dd0)*(1, -ell);
            pair up4=up+.4*(-1, ell);
            pair down4=down+.4*(-1, -ell);

            fill(up--up2--down2--down--down4--up4--cycle, violation);
            draw(up--up2, linewidth(1.1));
            draw(down--down2, linewidth(1.1));
            draw(up--up4, linewidth(1.1));
            draw(down--down4, linewidth(1.1));

            draw(up--down, dashed);
            label("$x$", down, 2*S);

            draw(up3--down3, dotted);
            label("$\scriptstyle F(x)+\kappa\B(0, \norm{\alt x-x}_X)$", 1/4*up3+3/4*down3, E);

            real tmp2=x0+dd;
            pair xprime=(tmp2, xbar.y);
            label("$\scriptstyle \{\alt x\} \times \B(\opty,\epsilon')$", xprime, E);
            draw((xprime+(0,-epsilon))--(xprime+(0,epsilon)), linewidth(2)+gray(.3));
            draw(shift(xprime)*((-1, -epsilon)--(+.1, -epsilon)), dashed);
            draw(shift(xprime)*((-1, epsilon)--(+.1, epsilon)), dashed);
        \end{asy}
        \caption{illustration of technique}
        \label{fig:regularity:aubin-inclusion-ext-1}
    \end{subfigure}
    \hfil
    \begin{subfigure}[t]{0.45\textwidth}
        \centering
        \begin{asy}
            unitsize(100, 100);
            real fup0(real x){ return (-x+0.2)^3*0.45+x/3; };
            real fdown0(real x) { return fup0(x)-.6-(x-0.05)^2*1.1-x*1.5; };
            real x0=-.8;
            real mul(real x) { return (x>0 ? 2.1*x : max(-.1, 2.1*x+x*x*10)); };
            real fup(real x) { return 0.1+fup0(x0+mul(x-x0))/2.5; };
            real fdown(real x) { return fdown0(x0+mul(x-x0))/2.5; };
            real delta=.45;
            real deltaprime=.15;
            real deltaell=0.05;
            real ell=max(abs(-fup(x0+deltaell)+fup(x0)),
            abs(fdown(x0+deltaell)-fdown(x0)),
            abs(-fup(x0-deltaell)+fup(x0)),
            abs(fdown(x0-deltaell)-fdown(x0)))/deltaell;

            real epsilon=1.7*ell*(delta-2*deltaprime);

            path pup=graph(fup, -1, 1);
            path pdown=graph(fdown, -1, 1);
            path pupext=graph(fup, -1.1, 1.15);
            path pdownext=graph(fdown, -1.15, 1.15);

            fill(pupext--reverse(pdownext)--cycle, lightfill);
            draw(pup);
            draw(pdown);

            real maxf=fup(x0)+.3;
            real minf=fdown(x0)-.6;
            clip((-2, maxf)--(2, maxf)--(2, minf)--(-2, minf)--cycle);

            pair up0=pt(fup, x0);
            pair down0=pt(fdown, x0);

            pair xbar=(x0, 1/3*fup(x0)+2/3*fdown(x0));

            real dd=delta*3/2;
            real dd0=deltaprime*2/3;

            real tmp=x0+dd0;
            pair x=(tmp, 1/5*fup(tmp)+4/5*fdown(tmp));
            pair up=(tmp, fup(tmp));
            pair down=(tmp, fdown(tmp));
            pair up2=up+.7*(1, ell);
            pair down2=down+.7*(1, -ell);
            pair up3=up+(dd-dd0)*(1, ell);
            pair down3=down+(dd-dd0)*(1, -ell);
            pair up4=up+.7*(-1, ell);
            pair down4=down+.7*(-1, -ell);

            fill(up--up2--down2--down--down4--up4--cycle, violation);
            draw(up--up2, linewidth(1.1));
            draw(down--down2, linewidth(1.1));
            draw(up--up4, linewidth(1.1));
            draw(down--down4, linewidth(1.1));

            real tmp2=x0+dd;
            pair xprime=(tmp2, xbar.y);
            path el1=shift(xprime)*((-1, -epsilon)--(+.1, -epsilon));
            path el2=shift(xprime)*((-1, epsilon)--(+.1, epsilon));
            draw(el1, dashed);
            draw(el2, dashed);

            picture tmp=new picture;
            fill(tmp, up4--up--up2--cycle, darkfill);
            fill(tmp, down4--down--down2--cycle, darkfill);
            clip(tmp, el1--reverse(el2)--cycle);
            add(tmp);

        \end{asy}
        \caption{critical areas}
        \label{fig:regularity:aubin-inclusion-ext-2}
    \end{subfigure}
    \caption{
        Figure (\subref{fig:regularity:aubin-inclusion-ext-1}) illustrates the technique in \cref{thm:regularity:aubin-equiv} to prove the equivalence of the two set inclusion formulations of the Aubin property:
        For $\alt x$ outside the ball $\B(\optx, \delta)$, the set $\B(\opty,\epsilon')$ indicated by the thick dark gray line, is completely contained in the fat-cone structure $F(x)+\kappa\B(0, \norm{\alt x-x}_X)$ of \cref{fig:regularity:aubin-structure}, indicated by the thick black and dotted lines.
        Closer to $x$, within $\B(\optx, \delta)$, this is not the case, although $F(\alt x) \isect \B(\opty,\epsilon')$ itself is still contained in the structure.
        Figure (\subref{fig:regularity:aubin-inclusion-ext-2}) highlights in darker color the areas that are critical for the Aubin property to hold.
    }
    \label{fig:regularity:aubin-inclusion-ext}
\end{figure}
\begin{proof}
    \emph{\ref{item:regularity:aubin-inclusion} $\equivalent$ \ref{item:regularity:aubin-inclusion-ext}:}
    Clearly \cref{item:regularity:aubin-inclusion-ext} implies \cref{item:regularity:aubin-inclusion} with the same $\kappa,\delta>0$.
    To show the implication in the other direction, we start by applying \cref{item:regularity:aubin-inclusion} with $\alt x=\optx$, which yields
    \begin{equation*}
        F(\optx) \isect \B(\opty, \delta) \subset F(x) +  \B(0, \kappa \norm{\optx-x}_X) \quad  (x \in \B(\optx, \delta)).
    \end{equation*}
    Taking $x \in \B(\optx, \delta')$ for some $\delta' \in (0, \delta]$, we thus deduce that
    \begin{equation*}
        \opty \in F(x) +  \B(0, \kappa \norm{\optx-x}_X) \subset F(x) + \B(0, \kappa \delta').
    \end{equation*}
    In particular, for any $\epsilon'>0$, we have
    \begin{equation}
        \label{eq:regularity:aubin-inclusion-barx}
        \B(\opty, \epsilon') \subset F(x) + \B(0, \kappa\delta'+\epsilon').
    \end{equation}
    For $\alt x \in \B(\optx, \delta)$, \cref{item:regularity:aubin-inclusion-ext} is immediate from \cref{item:regularity:aubin-inclusion}, so we may concentrate on $\alt x \in X \setminus \B(\optx, \delta)$. Then
    \begin{equation*}
        \norm{\alt x-x}_X
        \ge \norm{\alt x-\optx}_X-\norm{\optx - x}_X
        \ge \delta-\delta'.
    \end{equation*}
    If we pick $\epsilon',\delta'>0$ such that $\kappa\delta'+\epsilon' \le \kappa(\delta-\delta')$, it follows
    \begin{equation*}
        \kappa\delta'+\epsilon' \le \kappa\norm{\alt x-x}_X.
    \end{equation*}
    Thus \eqref{eq:regularity:aubin-inclusion-barx} gives, as illustrated in \cref{fig:regularity:aubin-inclusion-ext},
    \begin{equation*}
        F(\alt x) \isect \B(\opty, \epsilon')
        \subset \B(\opty, \epsilon')
        \subset F(x) + \B(0, \kappa\delta'+\epsilon')
        \subset F(x) + \kappa \B(0, \norm{\alt x-x}_X),
    \end{equation*}
    which is \cref{item:regularity:aubin-inclusion-ext}.

    \emph{\ref{item:regularity:aubin-inclusion-ext} $\equivalent$ \ref{item:regularity:aubin-distance}:}
    We expand \cref{item:regularity:aubin-inclusion-ext} as
    \begin{equation*}
        \{\alt y\} \isect \B(\opty, \delta) \subset F(x) + \B(0, \kappa\norm{\alt x-x}_X)
        \quad
        (\alt y \in F(\alt x); x \in \B(\optx, \delta);\, \alt x \in X).
    \end{equation*}
    By rearranging and taking the infimum over all $y\in F(x)$, this yields
    \begin{equation*}
        \inf_{y \in F(x)} \norm{\alt y-y}_Y \le \kappa\norm{\alt x-x}_X
        \quad
        (\alt y \in F(\alt x) \isect \B(\opty, \delta);\, x \in \B(\optx, \delta);\, \alt x \in X).
    \end{equation*}
    This may further be rewritten as
    \begin{equation*}
        \inf_{y \in F(x)} \norm{\alt y-y}_Y \le \inf_{\alt x \in \inv F(\alt y)} \kappa\norm{\alt x-x}_X
        \quad
        (x \in \B(\optx, \delta); \alt y \in \B(\opty, \delta)).
    \end{equation*}
    Thus \ref{item:regularity:aubin-distance} is equivalent to \ref{item:regularity:aubin-inclusion}.

    \emph{\ref{item:regularity:aubin-distance} $\implies$ \ref{item:regularity:aubin-distance-restr}:}
    This is immediate from the definition of $\dist$, which yields
    \begin{equation*}
        \dist(\inv F(y), x)
        \le
        \dist(\inv F(y) \isect \B(\optx, \delta), x).
    \end{equation*}

    \emph{\ref{item:regularity:aubin-distance-restr} $\implies$ \ref{item:regularity:aubin-inclusion}:}
    We express \cref{item:regularity:aubin-distance-restr} as
    \begin{equation*}
        \inf_{\alt y \in F(x)} \norm{y-\alt y}_Y \le \kappa \norm{\alt x-x}_X
        \quad (x \in \B(\optx, \delta),\, y \in \B(\opty, \delta),\, \alt x \in \inv F(y) \isect \B(\optx, \delta)).
    \end{equation*}
    This can be rearranged to imply that
    \begin{equation*}
        \{y\} \subset F(x) + \B(0, \kappa \norm{\alt x-x}_X)
        \quad (x \in \B(\optx, \delta),\, y \in \B(\opty, \delta) \isect F(\alt x),\, \alt x \in \B(\optx, \delta)),
    \end{equation*}
    which can be further rewritten as
    \begin{equation*}
        F(\alt x) \isect \B(\opty, \delta) \subset F(x) + \B(0, \kappa \norm{\alt x-x}_X)
        \quad (x, \alt x \in \B(\optx, \delta)),
    \end{equation*}
    yielding \cref{item:regularity:aubin-inclusion}.
\end{proof}

\begin{figure}
    \centering
    \begin{subfigure}[t]{.45\textwidth}
        \centering
        \begin{asy}
            real fup(real x){ return x^2*1.5; };
            real fdown(real x) { return cos(10*x)/20-.5-0.5*sqrt(0.7*sin(max(-x,0))); };
            real ell=.8;

            path pup=graph(fup, -1, 1);
            path pdown=graph(fdown, -1, 1);
            path pupext=graph(fup, -1.1, 1.15);
            path pdownext=graph(fdown, -1.15, 1.15);

            fill(pupext--reverse(pdownext)--cycle, lightfill);
            draw(pup);
            draw(pdown);

            real maxf=fup(0)+.5;
            real minf=fdown(0)-.5;
            clip((-2, maxf)--(2, maxf)--(2, minf)--(-2, minf)--cycle);

            pair up0=pt(fup, 0);
            pair down0=pt(fdown, 0);

            label("$F$", (-1, fdown(-1)), E+2*N);

            fill((up0+.5*(-1, ell))--up0--(up0+.5*(1, ell))--(down0+.5*(1, -ell))--down0--(down0+.5*(-1, -ell))--cycle, violation);
            draw(up0--(up0+.5*(1, ell)), linewidth(1.1));
            draw(up0--(up0+.5*(-1, ell)), linewidth(1.1));
            draw(down0--(down0+.5*(1, -ell)), linewidth(1.1));
            draw(down0--(down0+.5*(-1, -ell)), linewidth(1.1));

            real bup=fup(0)+.2;
            real bdown=fup(0)-.2;
            draw((-1.2, bup)--(1.2, bup), dashed);
            draw((-1.2, bdown)--(1.2, bdown), dashed);

            dot((0, fup(0)));

            real xprime=-.3;

            label("$x$", (0, -1.1), N+E);
            draw((0, -.8)--(0, .5), dotted);

            label("$\tilde x$", (xprime, -1.1), N);
            draw((xprime, -.8)--(xprime, .5), dotted);

            draw((xprime, bdown)--(xprime, fup(xprime)), linewidth(2)+gray(.3));

            label("$\B(\opty, \rho)$", (1.2, (bup+bdown)/2), E);

            draw((-.4, -1.2)--(.4, -1.2), Bars);
            label("$\B(\optx, \delta)$", (0, -1.3), S);
        \end{asy}
        \caption{Aubin property is satisfied}
        \label{fig:regularity:aubin:satisfied}
    \end{subfigure}
    \hfil
    \begin{subfigure}[t]{.45\textwidth}
        \centering
        \begin{asy}
            real fup(real x){ return x>0 ? sqrt(x) : 0; };
            real fdown(real x) { return x^3/2-.5; };
            real ell=1;

            path pup=graph(fup, -1, 1);
            path pdown=graph(fdown, -1, 1);
            path pupext=graph(fup, -1.1, 1.15);
            path pdownext=graph(fdown, -1.15, 1.15);

            fill(pupext--reverse(pdownext)--cycle, lightfill);
            draw(pup);
            draw(pdown);

            real maxf=fup(0)+.5;
            real minf=fdown(0)-.5;
            clip((-2, maxf)--(2, maxf)--(2, minf)--(-2, minf)--cycle);

            pair up0=pt(fup, 0);
            pair down0=pt(fdown, 0);

            label("$F$", (-1, fdown(-1)), E+4*N);

            fill((up0+.5*(-1, ell))--up0--(up0+.5*(1, ell))--(down0+.5*(1, -ell))--down0--(down0+.5*(-1, -ell))--cycle, violation);
            draw(up0--(up0+.5*(1, ell)), linewidth(1.1));
            draw(up0--(up0+.5*(-1, ell)), linewidth(1.1));
            draw(down0--(down0+.5*(1, -ell)), linewidth(1.1));
            draw(down0--(down0+.5*(-1, -ell)), linewidth(1.1));

            real bup=fup(0)+.3;
            real bdown=fup(0)-.3;
            draw((-1.2, bup)--(1.2, bup), dashed);
            draw((-1.2, bdown)--(1.2, bdown), dashed);

            dot((0, fup(0)));

            real xprime=.2;

            label("$x$", (0, -1.1), N+.5*W);
            draw((0, -.8)--(0, .5), dotted);

            label("$\tilde x$", (xprime, -1.1), N+.5*E);
            draw((xprime, -.8)--(xprime, .5), dotted);

            draw((xprime, bdown)--(xprime, bup), linewidth(2)+gray(.3));

            label("$\B(\opty, \rho)$", (1.2, (bup+bdown)/2), E);

            draw((-.4, -1.2)--(.4, -1.2), Bars);
            label("$\B(\optx, \delta)$", (0, -1.3), S);
        \end{asy}
        \caption{Aubin property is violated}
        \label{fig:regularity:aubin:dissatisfied}
    \end{subfigure}
    \caption{Illustration of \cref{ex:regularity:aubin} on verification of the Aubin property for $x=\optx$ based on \cref{thm:regularity:aubin-equiv}\,\cref{item:regularity:aubin-inclusion}.
        The dashed lines indicate $\B(\opty, \rho)$, and the dot marks $(\optx, \opty)$, while the thick vertical lines indicate $F(\alt x) \isect \B(\opty, \rho)$.
        The green ``fat'' cone shaded in green and bounded by the thick diagonal lines, indicates $F(x)+\B(0, \kappa\norm{\alt x-x}_X)$.
    }
    \label{fig:regularity:aubin}
\end{figure}

We next illustrate how the equivalent conditions of  \cref{thm:regularity:aubin-equiv} can be used to verify the Aubin property.

\begin{example}
    \label{ex:regularity:aubin}
    First, consider the function $F$ illustrated in \cref{fig:regularity:aubin:satisfied}.
    Due to \cref{thm:regularity:aubin-equiv}\,\cref{item:regularity:aubin-inclusion}, the set $F(\alt x) \isect \B(\opty, \rho)$ indicated by the thick vertical line should be contained in the ``fat'' cone $F(x)+\B(0, \kappa\norm{\alt x-x}_X)$ shaded in green.
    From the figure, this is clearly the case for $x$ equal to the base point $\optx$.
    Also from the figure, we can see that small variations of $x$ -- which will move the fat cone -- and of $\alt x$ -- which will move the thick line -- will not destroy the inclusion.
    The violation of the bounds at the bottom does not matter because we are only interested in the area between the dashed lines.
    Hence the Aubin property is satisfied for $F$ at $\optx$ for $\opty$.

    Consider now the function $F$ illustrated in \cref{fig:regularity:aubin:dissatisfied}.
    Unlike in the previous example, the set $F(\alt x) \isect \B(\opty, \rho)$ indicated by the thick vertical line is not contained in the ``fat'' cone $F(x)+\B(0, \kappa\norm{\alt x-x}_X)$ shaded in green.
    We could reduce $\rho$ to satisfy the property at this specific $\alt x$, but it would still fail as $\alt x \to x=\optx$, due to the high rate of increase of the upper bound of the graph of $F$ at $\optx$.
    Hence the Aubin property is violated for $F$ at $\optx$ for $\opty$.
\end{example}

We have similar characterizations of calmness. The proof is analogous to that of \cref{thm:regularity:aubin-equiv}, simply fixing $x=\optx$.

\begin{corollary}
    \label{cor:regularity:calmness-equiv}
    Let $X,Y$ be Banach spaces and $F: X \setto Y$.
    Then the following are equivalent for $\optx\in X$ and $\opty\in F(\optx)$:
    \begin{enumerate}
        \item\label{item:regularity:calmness-inclusion} There exists $\kappa,\delta>0$ such that
            \begin{equation*}
                F(\alt x) \isect \B(\opty, \delta) \subset F(\optx) +  \B(0, \kappa \norm{\alt x-\optx}_X) \quad (\alt x \in \B(\optx, \delta)).
            \end{equation*}
        \item\label{item:regularity:calmness-inclusion-ext} There exists $\kappa,\delta>0$ such that
            \begin{equation*}
                F(\alt x) \isect \B(\opty, \delta) \subset F(\optx) +  \B(0, \kappa \norm{\alt x-\optx}_X) \quad (\alt x \in X).
            \end{equation*}
        \item\label{item:regularity:calmness-distance} The calmness property \eqref{eq:regularity:calmness-distance} holds.
        \item\label{item:regularity:calmness-distance-restr} There exists $\kappa,\delta>0$ such that
            \begin{equation*}
                \dist(y, F(\optx)) \le \kappa \dist(\inv F(y) \isect \B(\optx, \delta), \optx)
                \quad ( y \in \B(\opty, \delta)).
            \end{equation*}
    \end{enumerate}
    The infimum of $\kappa>0$ for which each of these characterizations holds is equal to the modulus of calmness $\calm F(\optx|\opty)$.
    (The radius of validity $\delta>0$ for any given $\kappa>0$ may be distinct in each of the characterizations, however.)
\end{corollary}

\begin{example}
    Continuing from \cref{ex:regularity:aubin}, since $x=\optx$ was fixed there, the function of \cref{fig:regularity:aubin:dissatisfied} is not calm at $\optx$ for $\opty$. The function of \cref{fig:regularity:aubin:satisfied} is, of course, calm, as this follows from its Aubin property.
\end{example}

\section{Neighborhood-based coderivative criteria}
\label{sec:regularity:coderivative-neighborhood}

 Our goal is now to relate the Aubin property to ``outer norms'' of limiting coderivatives, just as the Lipschitz property of differentiable single-valued functions can be related to norms of their derivatives. Before embarking on this in the next section, as a preparatory step we relate in this section the Aubin property to neighborhood-based criteria on Fréchet coderivatives.
To this end, we define for a set-valued mapping $F: X \setto Y$, $(\optx, \opty) \in \graph F$, and $\delta,\epsilon>0$
\begin{equation}
    \label{eq:regularity:kappa-delta-epsilon}
    \kappa_\delta^\epsilon(\optx|\opty) \defeq
    \sup\left\{
        \norm{x^*}_{X^*}
        \middle|
        \begin{array}{r}
            x^* \in \frechetCod_\epsilon F(x|y)(y^*),\,
            \norm{y^*}_{Y^*} \le 1,
            \\
            x \in \B(\optx, \delta),\,
            y \in F(x) \isect \B(\opty, \delta)
        \end{array}
    \right\},
\end{equation}
which measures locally the opening of the cones $\frechetNormal^\epsilon_{\graph F}(x|y)$ around $(\opt x, \opt y)$; for smooth functions and $\epsilon=0$, it coincides with the local supremum of $\norm{DF(x)}_{\linear(X;Y)}$ around $(\optx,F(\optx))$ (cf.~\cref{thm:graphical:single}).
The next lemma bounds these openings in terms of the graphical modulus.

\begin{lemma}
    \label{lemma:regularity:morduk-nbd:main}
    Let $X,Y$ be Banach spaces and $F: X \setto Y$.
    If $\graph F$ is closed near $(\optx, \opty)$, then for any choice of $\epsilon(\delta) \ge 0$ satisfying $\epsilon(\delta) \downto 0$ as $\delta \downto 0$ it holds that
    \begin{equation*}
        \inf_{\delta>0} \kappa_\delta^{0}(\optx|\opty)
        \le
        \inf_{\delta>0} \kappa_\delta^{\epsilon(\delta)}(\optx|\opty)
        \le
        \lip F(\optx|\opty).
    \end{equation*}
\end{lemma}

\begin{proof}
    Since $\frechetCod F(x|y)(y^*) \subset \frechetCod_\epsilon F(x|y)(y^*)$ for any $\epsilon>0$, we always have $\kappa_\delta^{\epsilon(\delta)}(\optx|\opty) \ge \kappa_\delta^0(\optx|\opty)$. It hence suffices to prove that
    \begin{equation*}
        \kappa \defeq \inf_{\delta>0} \kappa_\delta^{\epsilon(\delta)}(\optx|\opty) \le \lip F(\optx|\opty).
    \end{equation*}
    We may assume that $\lip F(\optx|\opty) < \infty$, since otherwise there is nothing to prove.
    This implies in particular that the Aubin property holds, so the definition \eqref{eq:regularity:aubin-distance} yields for any $\kappa' > \lip F(x|y)$ a $\delta'>0$ such that
    \begin{equation}
        \label{eq:regularity:morduk-nbd:proof-ineq1}
        \inf_{\alt y \in F(\alt x)} \norm{\alt y-y}_Y
        \le \kappa' \norm{\alt x-x}_X
        \quad (y \in F(x) \isect \B(\opty, \delta'),\, \alt x \in \B(\optx, \delta')).
    \end{equation}
    Pick $\alt\kappa \in (0, \kappa)$ and $\delta \in (0, \delta')$.
    By the definition of $\kappa_{\delta}^{\epsilon(\delta)}(\optx|\opty)$, there exist $x \in \B(\optx, \delta)$, $y \in F(x) \isect \B(\opty, \delta)$, and $(x^*, -y^*) \in \frechetNormal_{\graph F}^{\epsilon(\delta)}(x, y)$ such that $\norm{x^*}_{X^*} \ge \alt\kappa$ and $\norm{y^*}_{Y^*} \le 1$.
    \Cref{thm:functan:hb-extension} then yields a $\dir x \in X$ such that
    \begin{equation}
        \label{eq:regularity:morduk-nbd:dirx}
        \dualprod{x^*}{\dir x}_{X}=\norm{x^*}_{X^*}
        \quad\text{and}\quad
        \norm{\dir x}_X=1.
    \end{equation}
    Let $\tau_k \downto 0$ with $\tau_k \le \delta$ and set $x_k \defeq x + \tau_k \dir x$.
    Then taking $\alt x=x_k$ in \eqref{eq:regularity:morduk-nbd:proof-ineq1}, we can find $y_k \in F(x_k)$ such that
    \begin{equation}
        \label{eq:regularity:morduk-nbd:diry-constr}
        \liminf_{k \to \infty} \inv\tau_k \norm{y_k-y}_Y \le \kappa'\norm{\dir x}_X = \kappa'.
    \end{equation}
    In particular, after passing to a subsequence if necessary, we may assume that $y_k \to y$ strongly in $Y$.
    Using \eqref{eq:regularity:morduk-nbd:dirx},  $\norm{x^*}_{X^*} \ge \alt\kappa$, and $\norm{y^*}_{Y^*} \le 1$, this leads to
    \begin{equation}
        \label{eq:regularity:morduk-nbd:dirx-diry-bound}
        \begin{aligned}[t]
            \limsup_{k \to \infty}~ & \inv\tau_k\left(
                \dualprod{x^*}{x_k-x}_{X} - \dualprod{y^*}{y_k-y}_{Y}
            \right)
            \\
            &
            =
            \limsup_{k \to \infty} \left(
                \dualprod{x^*}{\dir x}_{X} - \inv\tau_k \dualprod{y^*}{y_k-y}_{Y}
            \right)
            \\
            &
            \ge \norm{x^*}_{X^*}-\kappa' \ge \alt\kappa-\kappa'.
        \end{aligned}
    \end{equation}
    By \eqref{eq:regularity:morduk-nbd:diry-constr} (for the chosen subsequence) and the construction of $x_k$, we have
    \begin{equation}
        \label{eq:regularity:morduk-nbd:dirx-diry-norm-bound}
        \limsup_{k \to \infty} \inv\tau_k \norm{(x_k,y_k)-(x, y)}_{X \times Y}
        \le
        (1+\kappa')\norm{\dir x}_X
        =
        1+\kappa'.
    \end{equation}
    Since  $(x^*, -y^*) \in \frechetNormal_{\graph F}^{\epsilon(\delta)}(x, y)$, from the defining equation \eqref{eq:cones:def-epsiloncone} of $\frechetNormal_{\graph F}^{\epsilon(\delta)}(x, y)$ we have
    \begin{equation}
        \label{eq:regularity:morduk-nbd:eps-normal}
        \limsup_{k \to \infty} \frac{\dualprod{x^*}{x_k-x}_{X} - \dualprod{y^*}{y_k-y}_{Y}}{\norm{(x_k,y_k)-(x, y)}_{X \times Y}} \le \epsilon(\delta).
    \end{equation}
    Therefore, \eqref{eq:regularity:morduk-nbd:dirx-diry-bound}, \eqref{eq:regularity:morduk-nbd:dirx-diry-norm-bound}, and \eqref{eq:regularity:morduk-nbd:eps-normal} together yield
    \begin{equation*}
        (1+\kappa')\epsilon(\delta) \ge \alt\kappa-\kappa'.
    \end{equation*}
    Since this holds for any $\delta \in (0, \delta')$ and by assumption $\epsilon(\delta) \downto 0$ as $\delta \downto 0$, it follows that $\alt\kappa \ge \kappa'$.
    Since $\kappa' > \lip F(\optx|\opty)$ and $\alt\kappa < \kappa$ were arbitrary, we obtain $\kappa \le \lip F(\optx|\opty)$ as desired.
\end{proof}

For the next theorem, recall the definition of Gateaux smooth spaces from \cref{sec:epsilon:smooth}.

\begin{theorem}
    \label{thm:regularity:morduk-nbd:gateaux}
    Let $X,Y$ be Gateaux smooth Banach spaces and let $F: X \setto Y$ be such that $\graph F$ is closed near $(\optx, \opty)\in X\times Y$. Then $F$ has the Aubin property at $\optx$ for $\opty$ if and only if $\kappa_\delta^\delta(\optx|\opty)<\infty$ or $\kappa_\delta^0(\optx|\opty)<\infty$  for some $\delta > 0$. Furthermore, in this case
    \begin{equation*}
        \label{eq:regularity:morduk-nbd:gateaux}
        \inf_{\delta>0} \kappa_\delta^\delta(\optx|\opty) = \lip F(\optx|\opty) = \inf_{\delta>0} \kappa_\delta^0(\optx|\opty).
    \end{equation*}
\end{theorem}

\begin{proof}
    By \cref{lemma:regularity:morduk-nbd:main}, it suffices to show for any $\epsilon(\delta) \downto 0$ as $\delta \downto 0$ that
    \begin{equation*}
        \kappa \defeq \inf_{\delta>0} \kappa_\delta^{\epsilon(\delta)}(\optx|\opty)  \ge \lip F(\optx|\opty)
    \end{equation*}
    as choosing then in turn $\epsilon(\delta) = \delta$ and $\epsilon(\delta) = 0$ gives the two claimed equalities.
    We may assume that $\lip F(\optx|\opty)>0$ as otherwise there is nothing to show.
    Our plan is now to take arbitrary $0 < \alt\kappa < \lip F(\optx|\opty)$ and show that $\kappa \ge \alt \kappa$, which implies $\kappa \ge \lip F(\optx|\opty)$ and hence the claim.

    To do so, it suffices to show that $\kappa_\delta^{\epsilon(\delta)}(\optx|\opty) \ge \alt \kappa$ for all $\delta > 0$.
    To that end, we will select families $\epsilon_t \downto 0$ and $(x_t, y_t) \to (\optx, \opty)$ as $t\downto 0$
    and construct $\epsilon_t$-normals $(x_t^*, -y_t^*) \in \frechetNormal_{\graph F}^{\epsilon_t}(x_t, y_t)$ that satisfy $\liminf_{t \to 0} \norm{x_t^*}_{X^*} \ge \alt\kappa$ and $\limsup_{t \to 0} \norm{y_t^*}_{Y^*} \le 1$.
    By the definition of $\kappa_\delta^{\epsilon(\delta)}(\optx|\opty)$ in \eqref{eq:regularity:kappa-delta-epsilon}, taking for each $\delta>0$ the index $t>0$ such that both $\epsilon_t \le \epsilon(\delta)$ and $\max\{\norm{x_t-\optx}_X,\norm{y_t-\opty}_Y\} \le \delta$, this will then show the desired inequality $\kappa_\delta^{\epsilon(\delta)}(\optx|\opty) \ge \alt \kappa$.
    The rough idea is to construct the $\epsilon_t$-normals by projecting points not in $\graph F$ back onto this set. There are, however, some technical difficulties along our way.
    We divide the construction into three steps.

    \emph{Step 1: setting up the projection problem.}
    Let $0 < \alt\kappa < \lip F(\optx|\opty)$.
    Since then the Aubin property does not hold for $\alt\kappa$, by the characterization of \cref{thm:regularity:aubin-equiv}\,\ref{item:regularity:aubin-distance-restr} there exist
    \begin{equation}
        \label{eq:regularity:morduk-nbd:t-nbd}
        \alt y_t \in F(\alt x_t) \isect \B(\opty, t)
        \quad\text{and}\quad
        \alt x_t, x_t \in \B(\optx, t)
        \quad\text{for all } t > 0
    \end{equation}
    such that
    \begin{equation}
        \label{eq:regularity:morduk-nbd:fail}
        \inf_{y_t \in F(x_t)} \norm{y_t-\alt y_t}_Y > \alt\kappa \norm{x_t-\alt x_t}_X.
    \end{equation}
    Since $\inf_{y_t \in F(\alt x_t)} \norm{y_t-\alt y_t}_Y=0$, this implies that $x_t \ne \alt x_t$ and $(x_t, \alt y_t) \not \in \graph F$. We want to locally project $(x_t, \alt y_t)$ back onto $\graph F$.
    However, the nondifferentiability of the distance function $\norm{\freevar-\alt y_t}_Y$ at $\alt y_t$ would cause difficulties, so (similarly to the proof of \cref{lemma:cones:icecream}) we modify the projection by composing  the norm with the ``smoothing function''
    \begin{equation}
        \label{eq:regularity:morduk-nbd:phi-epsilon}
        \phi_\mu(r) \defeq \sqrt{\mu^2+r^2}-\mu.
    \end{equation}
    By \cref{thm:subdifferential:norm,thm:convex:increasing-post,thm:convex:gateaux} and the assumed differentiability of $\norm{\freevar}_Y$ away from the origin, $\phi_\mu(\norm{\freevar}_Y)$ is convex and has a single-valued subdifferential mapping with elements of norm less than one.
    Hence this smoothed distance function is Gateaux differentiable by \cref{lem:clarke:gateaux}.
    Due to \eqref{eq:regularity:morduk-nbd:phi-epsilon}, for every $t>0$ and $\mu_t>0$, we further have
    \begin{equation}
        \label{eq:regularity:morduk-nbd:norm-approx-bounds}
        \norm{y-\alt y_t}_Y - \mu_t \le \phi_{\mu_t}(\norm{y-\alt y_t}_Y) \le \norm{y-\alt y_t}_Y
        \quad (y \in Y).
    \end{equation}
    To locally project $(x_t, \alt y_t)$ onto $\graph F$, we thus seek to minimize the function
    \begin{equation}
        \label{eq:regularity:morduk-nbd:psi-def}
        \psi_t(x, y) \defeq \delta_{C_t}(x, y) + \alt\kappa \norm{x-x_t}_X + \phi_{\mu_t}(\norm{y-\alt y_t}_Y)
    \end{equation}
    for
    \begin{equation*}
        C_t \defeq [\B(\optx, t+2\alt\kappa) \times \B(\opty, t+2\alt\kappa)] \isect \graph F.
    \end{equation*}
    Clearly, $\psi_t$ is bounded from  below by $-\mu_t$ as well as coercive since $C_t$ is bounded. If $t$ is small enough, then $C_t$ is closed by the local closedness of $\graph F$. Therefore $\psi_t$ is lower semicontinuous (but not weakly lower semicontinuous since $\graph F$ need not be convex).

    \emph{Step 2: finding approximate minimizers.}
    We would like to find a minimizer of $\psi_t$, but the lack of weak lower semicontinuity prevents the use of Tonelli's direct method of \cref{thm:variation:existence}.
    We therefore use Ekeland's variational principle (\cref{thm:variation:ekeland}) to find an approximate minimizer. Towards this end, choose  for every $t>0$
    \begin{equation}
        \label{eq:regularity:morduk-nbd:epsilon-lambda}
        \mu_t \defeq t^{-1/2} \alt\kappa \norm{\alt x_t-x_t}_X^2 \le \alt\kappa t^{1/2}\norm{\alt x_t-x_t}_X
        \quad\text{and}\quad
        \lambda_t \defeq \norm{\alt x_t-x_t}_X + t^{1/2} \le t+t^{1/2},
    \end{equation}
    where the inequalities hold due to \eqref{eq:regularity:morduk-nbd:t-nbd}.
    Then
    \begin{equation}
        \label{eq:regularity:morduk-nbd:psi-altxy-bound}
        \psi_t(\alt x_t, \alt y_t)
        =
        \alt\kappa\norm{\alt x_t-x_t}_X
        \le
        (\alt\kappa\norm{\alt x_t-x_t}_X + \mu_t) + \inf \psi_t.
    \end{equation}
    Therefore, applying \cref{thm:variation:ekeland} for  $\lambda=\lambda_t$ and
    \begin{equation*}
        \epsilon=\alt\kappa\norm{\alt x_t-x_t}_X+\mu_t=\alt\kappa\norm{\alt x_t-x_t}_Xt^{-1/2}\lambda_t=\frac{\mu_t\lambda_t}{\norm{\alt x_t-x_t}_X},
    \end{equation*}
    we obtain for each $t>0$ a strict minimizer $(\opt x_t, \opt y_t)$ of
    \begin{subequations}%
        \begin{gather}%
            \label{eq:regularity:morduk-nbd:ekeland-min}
            \alt \psi_t(x, y)
            \defeq \psi_t(x, y) + \frac{\mu_t}{\norm{\alt x_t-x_t}_X} \left(\norm{x-\opt x_t}_X+\norm{y-\opt y_t}_Y\right)
            \shortintertext{with}
            \label{eq:regularity:morduk-nbd:ekeland-less}
            \psi_t(\opt x_t, \opt y_t) +  \frac{\mu_t}{\norm{\alt x_t-x_t}_X} \left(\norm{\alt x_t-\opt x_t}_X+\norm{\alt y_t-\opt y_t}_Y\right)
            \le \psi_t(\alt x_t, \alt y_t)
            = \kappa\norm{\alt x_t-x_t}_X
            \shortintertext{and}
            \label{eq:regularity:morduk-nbd:ekeland-bound}
            \norm{\opt x_t-\alt x_t}_X+\norm{\opt y_t-\alt y_t}_Y \le \lambda_t.
        \end{gather}
    \end{subequations}

    We claim that $\opt x_t \ne x_t$, which we show by contradiction. Assume therefore that $\optx_t=x_t$. Then $\opty_t \in F(x_t)$, and \eqref{eq:regularity:morduk-nbd:norm-approx-bounds} yields
    \begin{equation*}
        \psi_t(x_t, \opty_t) = \phi_{\mu_t}(\norm{\opt y_t-\alt y_t}_Y)
        \ge
        \norm{\opt y_t-\alt y_t}_Y - \mu_t.
    \end{equation*}
    Thus by \eqref{eq:regularity:morduk-nbd:psi-altxy-bound} and \eqref{eq:regularity:morduk-nbd:ekeland-less},
    \begin{equation*}
        \norm{\opt y_t-\alt y_t}_Y
        \le
        \norm{\opt y_t-\alt y_t}_Y - \mu_t
        +
        \frac{\mu_t}{\norm{\alt x_t-x_t}_X} \left(\norm{\alt x_t-x_t}_X+\norm{\alt y_t-\opt y_t}_Y\right)
        \le
        \alt\kappa\norm{\alt x_t-x_t}_X.
    \end{equation*}
    But this contradicts \eqref{eq:regularity:morduk-nbd:fail} as $\opt y_t \in F(x_t)$.

    \emph{Step 3: constructing $\epsilon$-normals.}
    We are now ready to construct the desired $\epsilon$-normals. We write
    \begin{equation}
        \label{eq:regularity:morduk-nbd:alt-psi-rearrange}
        \alt \psi_t(x, y)=\delta_{C_t}(x, y) + \Phi(x, y)
    \end{equation}
    for the convex and Lipschitz continuous function
    \begin{equation*}
        \Phi(x, y) \defeq \alt\kappa \norm{x-x_t}_X + \phi_{\mu_t}(\norm{y-\alt y_t}_Y)
        +\frac{\mu_t}{\norm{\alt x_t-x_t}_X} \left(\norm{x-\alt x_t}_X+\norm{y-\alt y_t}_Y\right).
    \end{equation*}
    Since we assume $X$ to be Gateaux smooth, $x \mapsto \alt\kappa \norm{x-x_t}_Y$ is Gateaux differentiable at $\alt x_t \ne x_t$. Furthermore, $y \mapsto \phi_{\mu_t}(\norm{y-\alt y_t}_Y)$ is by construction Gateaux differentiable for all $y$.
    By \eqref{eq:regularity:morduk-nbd:epsilon-lambda}, we have
    $
    \tfrac{\mu_t}{\norm{\alt x_t-x_t}_X}
    \le t^{1/2}\alt\kappa.
    $
    Since $x_t \ne \optx_t$, \cref{thm:convex:increasing-post,thm:subdifferential:norm,thm:subdiff:sum} now yield
    \begin{equation}
        \label{eq:regularity:morduk-nbd:f-subdiff}
        \subdiff \Phi(\optx_t, \opty_t) \subset \B((-x_t^*, y_t^*), t^{1/2}\alt\kappa)
        \quad\text{for}\quad
        \left\{
        \begin{aligned}
            -x_t^* &= \alt\kappa D [\norm{\freevar-x_t}_X](\opt x_t),\\
            y_t^* &= D[\phi_{\mu_t}(\norm{\freevar -\alt y_t}_Y)](\opt y_t).
        \end{aligned}
        \right.
    \end{equation}

    Since $\opt x_t \ne x_t$, we have $\norm{x_t^*}_{X^*} = \alt\kappa$ by \cref{thm:subdifferential:norm}. Moreover, $\norm{y_t^*}_{Y^*} \le 1$ as observed in Step 1.
    \Cref{thm:limiting:frechet:fermat} further yields $0 \in \subdiff_F \alt\psi_t(\optx_t, \opty_t)$.

    Due to \eqref{eq:regularity:morduk-nbd:alt-psi-rearrange} and \eqref{eq:regularity:morduk-nbd:f-subdiff}, \cref{lemma:epsilon:sumrule} now shows that
    \begin{equation*}
        (x_t^*, -y_t^*) \in \frechetNormal_{C_t}^{\epsilon_t}(\opt x_t, \opt y_t),
        \quad\text{i.e.,}\quad
        x_t^* \in \frechetCod_{\epsilon_t} F(\opt x_t|\opt y_t)(y_t^*)
        \quad\text{for}\quad
        \epsilon_t \defeq t^{1/2}\alt\kappa.
    \end{equation*}
    We illustrate this construction in \cref{fig:regularity:morduk-nbd-part2}.
    Since $\lambda_t \le t + t^{1/2}$ by \eqref{eq:regularity:morduk-nbd:epsilon-lambda}, it follows from \eqref{eq:regularity:morduk-nbd:ekeland-bound} that
    $
    \norm{\optx_t-\optx}_X, \norm{\opty_t-\opty}_Y \le 2t+t^{1/2}
    $
    and hence that $(\optx_t, \opty_t) \to (\optx, \opty)$ as $t \downto 0$. We also have both $\liminf_{t \downto 0} \norm{x_t^*}_{X^*} \ge \alt\kappa$ and  $\limsup_{t \downto 0} \norm{y_t^*}_{Y^*} \le 1$.
    Thus we have constructed the desired sequence of $\epsilon_t$-normals.
\end{proof}

\begin{figure}
    \centering
    \begin{asy}
        unitsize(120, 120);

        import markers;

        real xbar=0;
        real slope1=-1;
        real slope2=-0.1;

        real f(real x){
            if(x<xbar){
                return slope1*x;
                }else{
                return slope2*x;
            }
        }

        pair normalize(pair p){
            return p/sqrt(p.x^2+p.y^2);
        }

        pair ta=normalize((1, slope1));
        pair no=normalize((slope1, -1));

        pair project(pair p){
            return dot(p, ta)*ta;
        }

        pair graphpoint(real f(real), real x){
            return (x, f(x));
        }

        real x0=-0.9;
        real x1=0.5;
        real y1=1.05;
        path fpath=graph(f, x0, x1);
        fill(fpath--(x1, y1)--(x0, y1)--cycle, lightfill);
        draw(fpath);

        label("$F$", (x1, y1), 2*S+2*W);

        pair popt=graphpoint(f, 0);
        pair ptilde=graphpoint(f, -0.1);
        pair pprime=graphpoint(f, -0.8);
        pair poutside=(pprime.x, ptilde.y);
        pair pproj=project(poutside);

        dot(popt);
        label("$(\opt x, \opt y)$", popt, S+E);

        dot(ptilde);
        label("$(\alt x_t, \alt y_t)$", ptilde, 2*E+N);

        dot(pprime);
        label("$(x_t, y_t)$", pprime, S+W);

        dot(poutside);
        label("$(x_t, \alt y_t)$", poutside, S+W);

        dot(pproj);
        label("$(\bar x_t, \bar y_t)=\mbox{proj}(x_t, \alt y_t)$", pproj, 2.5*W+0.5*N);

        draw(pproj--poutside, primalline + linewidth(1.2), Arrow);

        draw((pprime.x, pprime.y+0.18)--pprime, dotted, Arrow);

        pair phelper=(ptilde.x, pprime.y);
        draw(ptilde--phelper--pprime, dotted);
        markangle(A=pprime, O=ptilde, B=phelper, radius=-20, dotted);

        label("$d$", pprime--phelper, N);
        label("$\ge \alt\kappa d$", phelper--ptilde, E);

        pair phelper2=(pproj.x, poutside.y);
        draw(poutside--phelper2--pproj, dotted);
        markangle(A=phelper2, O=pproj, B=poutside, radius=-15, dotted);

        label("$\ge \alt\kappa n_y$", poutside--phelper2, S+0.5*E);
        label("$n_y$", phelper2--pproj, 0.5*W+0.5*S);

        draw((pproj+0.1*ta+0.05*no)--(popt-0.05*ta+0.05*no), dashed, Arrow);

    \end{asy}
    \caption{
        The construction in the final part of the proof of \cref{thm:regularity:morduk-nbd:gateaux}. The dotted arrow indicates how $y_t$ minimizes the distance to $\tilde y_t$ within $F(x_t)$, which ensures that $\norm{y_t-\alt y_t}_Y \ge \alt \kappa d$ for $d \defeq \norm{x_t-\alt x_t}_X$. The point $(x_t, \alt y_t)$ is outside $\graph F$; when projected back as $(\opt x_t, \opt y_t)$, the normal vector to $\graph F$ indicated by the solid arrow has $x$-component larger than the $y$-component $n_y$ by the factor $\alt\kappa$.
        The dashed arrow indicates the convergence of the other points to $(\optx, \opty)$ as $t \downto 0$.
    }
    \label{fig:regularity:morduk-nbd-part2}
\end{figure}

\begin{remark}
    Our proof of \cref{thm:regularity:morduk-nbd:gateaux} differs from those in \cite{mordukhovich2018variational,Mordukhovich:2006} by the specific construction of the point $(x_t, \alt y_t) \not \in \graph F$ and the use of the smoothed distance $\phi_{\epsilon_t}(\norm{\freevar}_X)$. In contrast, the earlier proofs first translate the Aubin property (or metric regularity) into a \term[property!covering]{covering} or \term[property!linear openness]{linear openness} property to construct the point outside $\graph F$ that is to be projected back onto this set. In finite dimensions, \cite{mordukhovich2018variational} develops calculus for the limiting subdifferential of \cref{sec:limiting:mordukhovich} to avoid the lack of calculus for the Fréchet subdifferential; we instead apply the fuzzy calculus of \cref{lemma:epsilon:sumrule} to the smoothed distance function $\phi_{\mu_t}(\norm{\freevar}_X)$.
    A further alternative in finite dimensions involves the proximal subdifferentials used in \cite{Rockafellar:1998}.
    In infinite dimensions, \cite{Mordukhovich:2006} develops advanced extremal principles to work with the Fréchet subdifferential.
\end{remark}

\begin{remark}[relaxation of Gateaux smoothness]
    \label{rem:regularity:morduk-nbd:non-gateaux}
    The assumption that $Y$ (or, with somewhat more work, $X$) is Gateaux smooth in \cref{thm:regularity:morduk-nbd:gateaux} may be replaced with the assumption of the existence of a family $\{\theta_\mu: Y \to \R\}_{\mu>0}$ of Gateaux differentiable norm approximations satisfying
    \begin{equation*}
        \norm{y}_Y - \mu \le \theta_\mu(y) \le \norm{y}_Y
        \quad (y \in Y).
    \end{equation*}
    Then \eqref{eq:regularity:morduk-nbd:norm-approx-bounds} holds with $\theta_{\mu_t}(y-\alt y_t)$ in place of $\phi_{\mu_t}(\norm{y-\alt y_t}_Y)$.
    For example, with $\phi_\mu$ as in \eqref{eq:regularity:morduk-nbd:phi-epsilon}, in $L^p(\Omega)$ we can set
    \begin{equation*}
        \theta_\mu(y) \defeq  \norm{ \phi_\mu(\abs{y(\xi)}) }_{L^p(\Omega)}
        \quad (y \in L^1(\Omega)).
    \end{equation*}
    With somewhat more effort, the Gateaux smoothness of $X$ can be similarly relaxed.
\end{remark}

\section{Point-based coderivative criteria}
\label{sec:regularity:coderivative-point}

We will now convert the neighborhood-based criterion of \cref{lemma:regularity:morduk-nbd:main,thm:regularity:morduk-nbd:gateaux} into a simpler point-based criterion.
For the statement, we need to introduce a new smaller coderivative of $F: X \setto Y$ at $x$ for $y$,
the \term[coderivative!mixed]{mixed (limiting) coderivative} $\coderivative_M F(x|y): Y^* \setto X^*$,
\begin{equation}
    \label{eq:regularity:mixed-coderivative}
    \coderivative_M F(x|y)(y^*) \defeq \weakstarlimsup_{\substack{(\alt x, \alt y) \to (x, y) \\ \alt y^* \to  y^*,\, \epsilon \downto 0}} \frechetCod_\epsilon F(\alt x|\alt y)(\alt y^*),
\end{equation}
which differs from the ``normal'' coderivative
\begin{equation}
    \label{eq:regularity:normal-coderivative}
    \coderivative F(x|y)(y^*) = \weakstarlimsup_{\substack{(\alt x, \alt y) \to (x, y) \\ \alt y^* \weaktostar  y^*,\, \epsilon \downto 0}} \frechetCod_\epsilon F(\alt x|\alt y)(\alt y^*),
\end{equation}
by the use of weak-$*$ convergence in $X^*$ and strong convergence in $Y^*$ instead of weak-$*$ convergence in both.
(The mixed coderivative is not obtained directly from any of the usual normal cones, although one can naturally define corresponding mixed normal cones on product spaces.)

We further define for any $H: W \setto Z$ the \term[norm!outer]{outer norm}
\begin{equation*}
|H|^+ \defeq \sup \setof{ \norm{z}_Z}{z \in H(w),\, \norm{w}_W \le 1}.
\end{equation*}
We illustrate the outer norm by two examples in \cref{fig:regularity:outernorm}.
\begin{figure}
    \centering
    \begin{subfigure}[t]{.4\textwidth}
        \centering
        \begin{asy}
            real fup(real x){ return x^2*1.5; };
            real fdown(real x) { return cos(10*x)/20-.5-0.5*sqrt(0.7*sin(max(-x,0))); };
            real ell=.8;
            real one=.55;
            real zero=fup(0)-.1;

            path pup=graph(fup, -1, 1);
            path pdown=graph(fdown, -1, 1);
            path pupext=graph(fup, -1.1, 1.15);
            path pdownext=graph(fdown, -1.15, 1.15);

            fill(pupext--reverse(pdownext)--cycle, lightfill);
            draw(pup);
            draw(pdown);

            real maxf=fup(0)+.5;
            real minf=fdown(0)-.5;
            clip((-2, maxf)--(2, maxf)--(2, minf)--(-2, minf)--cycle);

            pair up0=pt(fup, .7);

            label("$H$", (up0+(.7, zero))/2);

            draw((-1.2, zero)--(1.2, zero), dashed);

            draw((-one, -1.2)--(one, -1.2), dashed, Bars);
            label("$[-1,1]$", (0, -1.2), S);

            real xmax=-.37;
            real vmax=fdown(xmax);
            draw((xmax, zero)--(xmax, vmax), linewidth(1.1));
            dot((xmax, vmax));
            label("$(w, z)$", (xmax, vmax), 2*E+0.5*S);

            dot((0, zero));
            label("$(0, 0)$", (0, zero), 2.5*N);
        \end{asy}
        \caption{general set-valued mapping $H$}
        \label{fig:regularity:outernorm-general}
    \end{subfigure}
    \hfil
    \begin{subfigure}[t]{.4\textwidth}
        \centering
        \begin{asy}
            real fup(real x){ return 1.6*x; };
            real fdown(real x) { return 0.2*x; };
            real ell=.8;
            real one=.55;
            real zero=0;

            path pup=graph(fup, 0, 1);
            path pdown=graph(fdown, 0, 1);
            path pupext=graph(fup, 0, 1.15);
            path pdownext=graph(fdown, 0, 1.15);

            fill(pupext--reverse(pdownext)--cycle, lightfill);
            draw(pup);
            draw(pdown);

            pair up0=pt(fup, .9);
            pair down0=pt(fdown, .9);

            label("$H$", (up0*2+down0*5)/7);

            draw((-1.2, zero)--(1.2, zero), dashed);

            draw((-one, -0.25)--(one, -0.25), dashed, Bars);
            label("$[-1,1]$", (0, -0.25), S);

            real xmax=one;
            real vmax=fup(xmax);
            draw((xmax, zero)--(xmax, vmax), linewidth(1.1));
            dot((xmax, vmax));
            label("$(w, z)$", (xmax, vmax), 2*E);

            dot((0, zero));
            label("$(0, 0)$", (0, zero), 2*N+W);
        \end{asy}
        \caption{mapping whose $\graph H$ is a cone}
        \label{fig:regularity:outernorm-cone}
    \end{subfigure}
    \caption{Points $(w, z)$ achieving the supremum in the expression of the outer norm~$\abs{H}^+$.}
    \label{fig:regularity:outernorm}
\end{figure}
We are mainly interested in the outer norms of coderivatives, in particular of
\begin{equation}
    \label{eq:regularity:morduk:kappa-expansion}
    |\coderivative_M F(\optx|\opty)|^+
= \sup \setof{ \norm{\opt x^*}_{X^*}}{\opt x^* \in \coderivative_M F(\optx|\opty)(\opt y^*),\, \norm{\opt y^*}_{Y^*} \le 1}.
\end{equation}
Recalling \cref{thm:cones:inclusions}, we have
\begin{equation}
    \label{eq:regularity:coderivative-inclusions}
    \frechetCod F(x|y)(y^*) \subset \coderivative_M F(x|y)(y^*) \subset \coderivative F(x|y)(y^*),
\end{equation}
so the outer norms satisfy
\begin{equation*}
    |\coderivative_M F(\optx|\opty)|^+ \le |\coderivative F(\optx|\opty)|^+.
\end{equation*}
We say that $F$ is \term[mapping!coderivatively normal]{coderivatively normal} at $\optx$ for $\opty$ if $|\coderivative_M F(\optx|\opty)|^+= |\coderivative F(\optx|\opty)|^+$.
Of course, if $Y$ is finite-dimensional, then $\coderivative_M F(x|y)=\coderivative F(x|y)$ and thus $F$ is always coderivatively normal.
Note that $|D^*F(\opt x|\opt y)|^+$ can be directly related to the neighborhood-based $\kappa_\delta^\delta$ defined in \eqref{eq:regularity:kappa-delta-epsilon}.
In particular, it measures the opening of the cone $N_{\graph F}(\optx, \opt y)$; compare \cref{fig:regularity:outernorm-cone}.

As the central result of this chapter, we now use this connection to derive a characterization of the Aubin property and the graphical modulus (and hence also of metric regularity and the modulus of metric regularity) through the outer norm of the mixed limiting coderivative. This \term[criterion, Mordukhovich]{Mordukhovich criterion} generalizes the classical relation between the Lipschitz constant of a $C^1$ function and the norm of its derivative.

\begin{lemma}[Mordukhovich criterion in general Banach spaces]
    \label{lemma:regularity:morduk}
    Let $X,Y$ be Banach spaces and let $F: X \setto Y$ be such that $\graph F$ is closed near $(\optx, \opty)\in X\times Y$.
    If $F$ has the Aubin property at $\optx$ for $\opty$, then
    \begin{equation}
        \label{eq:regularity:morduk:qc}
        \coderivative_M F(\optx|\opty)(0) = \{0\}
    \end{equation}
    and
    \begin{equation}
        \label{eq:regularity:morduk:lip}
        |\coderivative_M F(\optx|\opty)|^+
        \le
        \lip F(\optx|\opty).
    \end{equation}
\end{lemma}

\begin{proof}
    As the first step, we show that the Aubin property implies \eqref{eq:regularity:morduk:lip} and hence that $\kappa \defeq |\coderivative_M F(\optx|\opty)|^+ < \infty$.
    Let $\rho > 0$.
    By the definition of $\coderivative_M F(\optx|\opty)$ in \eqref{eq:regularity:mixed-coderivative}, there then exist $\delta \in (0, \rho)$, $x \in \B(\optx,\rho)$, and $y \in F(x) \isect \B(\opty, \rho)$ as well as $y^* \in Y^*$ and $x^* \in \frechetCod_\delta F(x|y)(y^*)$ such that $\norm{y^*}_{Y^*} \le 1+\rho$ and $\norm{x^*}_{X^*} \ge \kappa(1-\rho)^2$. (The upper bound on $\norm{y^*}_{Y^*}$ is why we need the \emph{mixed} coderivative, since $\norm{\freevar}_{Y^*}$ is continuous only in the strong topology. For the lower bound on $\norm{x^*}_{X^*}$, in contrast, the weak-$*$ lower semicontinuity of $\norm{\freevar}_{X^*}$ is sufficient.)
    Since $\frechetCod_\delta F(x|y)$ is formed from a cone, we may divide $x^*$ and $y^*$ by $1+\rho$ and thus assume that $\norm{y^*}_{Y^*} \le 1$ and $\norm{x^*}_{X^*} \ge \kappa(1-\rho)$.
    Consequently
    \begin{equation*}
        \kappa(1-\rho) \le \kappa_\delta^\delta(\optx|\opty)
        =
        \sup\left\{
            \norm{x^*}_{X^*}
            \middle|
            \begin{array}{r}
                x^* \in \frechetCod_\delta F(x|y)(y^*),\,
                \norm{y^*}_{Y^*} \le 1,
                \\
                x \in \B(\optx, \delta),\,
                y \in F(x) \isect \B(\opty, \delta)
            \end{array}
        \right\}.
    \end{equation*}
    Taking the infimum over $\delta>0$ and letting $\rho \downto 0$ thus shows
    \begin{equation*}
        \label{eq:regularity:kappa-normal}
        \kappa \le \inf_{\delta>0} \kappa_\delta^\delta(\optx|\opty).
    \end{equation*}
    It now follows from \cref{lemma:regularity:morduk-nbd:main} that $\kappa \le \lip F(\optx|\opty)$, which yields \eqref{eq:regularity:morduk:lip}.

    As the second step, we prove that the Aubin property implies \eqref{eq:regularity:morduk:qc}. We argue by contraposition.
    First, note that since $\graph \coderivative_M F(x|y)$ is a cone, $0 \in \coderivative_M F(x|y)(0)$.
    Hence if \eqref{eq:regularity:morduk:qc} does not hold, there exists $x^* \in X^*\setminus\{0\}$ such that
    \begin{equation*}
        x^*[0, \infty) \subset \coderivative_M F(x|y)(0).
    \end{equation*}
    By \eqref{eq:regularity:morduk:kappa-expansion} and the first step, this implies that $\infty= \kappa\leq \lip F(\optx|\opty)$ and hence that the Aubin property of $F$ at $\optx$ for $\opty$ is violated.
\end{proof}

Applied to $\inv F$, we obtain a corresponding result for metric regularity.
\begin{corollary}[Mordukhovich criterion for metric regularity in general Banach spaces]
    \label{cor:regularity:morduk-metric}
    Let $X,Y$ be Banach spaces and let $F: X \setto Y$ be such that $\graph F$ is closed near $(\optx, \opty)\in X\times Y$.
    If $F$ is metrically regular at $(\optx,\opty)$, then
    \begin{equation}
        \label{eq:regularity:morduk-metric:qc}
        0 \in \coderivative_M F(\optx|\opty)(y^*) \implies y^*= 0
    \end{equation}
    and
    \begin{equation}
        |\coderivative_M \inv F(\opty|\optx)|^+
        \le
        \reg F(\optx|\opty).
    \end{equation}
\end{corollary}

\begin{proof}
    We apply \cref{lemma:regularity:morduk} to $\inv F$, observing that \eqref{eq:regularity:morduk:qc} applied to $\inv F$ is \eqref{eq:regularity:morduk-metric:qc}.
\end{proof}

Under stronger assumptions on the spaces and the set-valued mapping, we obtain equivalence. For the following theorem, recall the definition of partial sequential normal compactness (PSNC) from \cref{sec:colimiting:psnc}.
\begin{theorem}[Mordukhovich criterion in smooth Banach spaces]
    \label{thm:regularity:morduk:iff}
    Let $X,Y$ be Gateaux smooth Banach spaces with $X$ reflexive and let $F: X \setto Y$ be such that $\graph F$ is closed near $(\optx, \opty)\in X\times Y$.
    If $F$ is PSNC at $\optx$ for $\opty$, then the following are equivalent:
    \begin{enumerate}
      \item\label{item:regularity:morduk:aubin} the Aubin property of $F$ at $\optx$ for $\opty$;
      \item\label{item:regularity:morduk:qc} the implication \eqref{eq:regularity:morduk:qc};
      \item\label{item:regularity:morduk:norm} $|\coderivative_M F(\optx|\opty)|^+<\infty$.
    \end{enumerate}
\end{theorem}

\begin{proof}
    Due to \cref{lemma:regularity:morduk}, it suffices to show that \cref{item:regularity:morduk:norm}\,$\implies$\,\cref{item:regularity:morduk:qc}\,$\implies$\,\cref{item:regularity:morduk:aubin}.
    We start with the second implication.
    Since $X$ and $Y$ are Gateaux smooth, \cref{thm:regularity:morduk-nbd:gateaux} yields
    \begin{equation}
        \label{eq:regularity:morduk:altkappa-lip}
        \lip F(\optx|\opty)
        =
        \alt\kappa \defeq
        \inf_{\delta>0}
        \sup\left\{
            \norm{x^*}_{X^*}
            \middle|
            \begin{array}{r}
                x^* \in \frechetCod_\delta F(x|y)(y^*),\,
                \norm{y^*}_{Y^*} \le 1,
                \\
                x \in \B(\optx, \delta),\,
                y \in F(x) \isect \B(\opty, \delta)
            \end{array}
        \right\}
    \end{equation}
    and that the Aubin property holds if $\alt\kappa<\infty$.
    We now argue by contradiction. Assume that the Aubin property does not hold. Then $\alt\kappa=\infty$ and hence we can find $(x_k, y_k) \to (\optx,\opty)$, $\epsilon_k \downto 0$, and $x_k^* \in \frechetCod_{\epsilon_k} F(x_k|y_k)(y_k^*)$ with $\norm{y_k^*}_{Y^*} \le 1$ and $\norm{x_k^*}_{X^*} \to \infty$. In particular, $y_k^*/\norm{x_k^*}_{X^*} \to 0$.
    Since $X$ is reflexive, we can apply the Eberlein--\u{S}muylan theorem (\cref{thm:ebsmul}) to extract a subsequence (not relabelled) such that $x_k^*/\norm{x_k^*}_{X^*} \weaktostar x^*$ for some $x^* \in X^*$. Since $\graph\frechetCod_{\epsilon_k}F(x_k|y_k)$ is a cone, we also have
    \begin{equation*}
        x_k^*/\norm{x_k^*}_{X^*} \in \frechetCod_{\epsilon_k} F(x_k|y_k)(y_k^*/\norm{x_k^*}_{X^*}).
    \end{equation*}
    By the definition \eqref{eq:regularity:mixed-coderivative} of the mixed coderivative, we deduce that $x^* \in \coderivative_M F(\optx|\opty)(0)$.
    We now make a case distinction: If $x^*\neq 0$, then this contradicts the qualification condition \eqref{eq:regularity:morduk:qc}.
    On the other hand, if $x^*=0$, the PSNC of $F$ at $\optx$ for $\opty$, implies that $1=\norm{x_k^*/\norm{x_k^*}_{X^*}}_{X^*} \to 0$, which is also a contradiction.
    Therefore \eqref{eq:regularity:morduk:qc} implies the Aubin property.

    It remains to show that \cref{item:regularity:morduk:norm}\,$\implies$\,\cref{item:regularity:morduk:qc}.
    First, since $\graph \coderivative_M F(\optx|\opty)$ is a cone, $\coderivative_M F(\optx|\opty)(0)$ is a cone as well. Hence by \eqref{eq:regularity:morduk:kappa-expansion}, $|\coderivative_M F(\optx|\opty)|^+<\infty$ implies that $\coderivative_M F(\optx|\opty)(0)=\{0\}$, which is \eqref{eq:regularity:morduk:qc}.
\end{proof}

Again, applying \cref{thm:regularity:morduk:iff} to $\inv F$ yields a characterization of metric regularity.

\begin{corollary}[Mordukhovich criterion for metric regularity in smooth Banach spaces]
    \label{cor:regularity:morduk-metric:iff}
    Let $X,Y$ be Gateaux smooth Banach spaces with $X$ reflexive and let $F: X \setto Y$ be such that $\graph F$ is closed near $(\optx, \opty)\in X\times Y$.
    If $\inv F$ is PSNC at $\opty$ for $\optx$, then the following are equivalent:
    \begin{enumerate}
        \item the metric regularity of $F$ at $(\optx,\opty)$;
        \item the implication \eqref{eq:regularity:morduk-metric:qc};
        \item $|\coderivative_M \inv F(\opty|\optx)|^+<\infty$.
    \end{enumerate}
\end{corollary}

\begin{remark}[separable and Asplund spaces]
    The reflexivity of $X$ (resp.~$Y$) was used to obtain the weak-$*$ compactness of the unit ball in $X^*$ via the Eberlein--\u{S}mulyan theorem (\cref{thm:ebsmul}) applied to $X^*$. Alternatively, this can be obtained by assuming separability of $X$ and using the Banach--Alaoglu theorem (\cref{thm:banachal}).
    More generally, dual spaces of Asplund spaces have weak-$*$-compact unit balls; we refer to \cite{Mordukhovich:2006} for the full theory in Asplund spaces.
\end{remark}

In finite dimensions, we have a full characterization of the graphical modulus via the outer norm of the limiting coderivative (which here coincides with the mixed coderivative).

\begin{corollary}[Mordukhovich criterion for the graphical modulus in finite dimensions]
    \label{cor:regularity:morduk:finite}
    Let $X,Y$ be finite-dimensional Gateaux smooth Banach spaces and let $F: X \setto Y$ be such that $\graph F$ is closed near $(\optx, \opty)\in X\times Y$.
    Then
    \begin{equation*}
        \lip F(\optx|\opty)
        =
        |\coderivative F(\optx|\opty)|^+.
    \end{equation*}
\end{corollary}

\begin{proof}
    Due to \cref{lemma:regularity:morduk}, we only have to show that
    \begin{equation}
        \label{eq:regularity:morduk:lip-converse}
        \lip F(\optx|\opty)
        \le
        |\coderivative F(\optx|\opty)|^+.
    \end{equation}
    As in the proof of \cref{thm:regularity:morduk:iff}, the smoothness of $X$ and $Y$ allows applying \cref{thm:regularity:morduk-nbd:gateaux} to obtain that $\lip F(\optx|\opty) = \alt\kappa$ given by
    \eqref{eq:regularity:morduk:altkappa-lip}. It therefore suffices to show that $\alt\kappa \le |\coderivative F(\optx|\opty)|^+$.
    Let $\kappa' < \alt\kappa$ be arbitrary. By \eqref{eq:regularity:morduk:altkappa-lip}, we can then find $(x_k, y_k) \to (\optx, \opty)$ and $\epsilon_k \downto 0$ as well as  $x_k^* \in \frechetCod_{\epsilon_k} F(x_k|y_k)(y_k^*)$ with $\norm{y_k^*}_{Y^*} \le 1$, and $\alt\kappa \ge \norm{x_k^*} \ge \kappa'$.
    Since $X$ and $Y$ are finite-dimensional, we can apply the Heine--Borel theorem to extract \emph{strongly} converging subsequences (not relabelled) such that $x_k^*\to x^*$ with $\norm{x^*}_{X^*} \ge \kappa'$ and $y_k^*\to y^*$ with $\norm{y^*}_{Y^*}\leq 1$.
    Since strongly converging sequences also converge weakly-$*$, the expression \eqref{eq:regularity:normal-coderivative} for the normal coderivative implies that $x^* \in \coderivative F(\optx|\opty)(y^*)$ and that $|\coderivative F(\optx|\opty)|^+ \ge \norm{x^*}_{X^*}\ge\kappa'$. Since $\kappa' < \alt\kappa$ was arbitrary, we obtain \eqref{eq:regularity:morduk:lip-converse}.
\end{proof}

This relation is illustrated in \cref{fig:regularity:morduk}, using that by definition of the outer norm and of the coderivative,
\begin{equation*}
    \begin{aligned}
        \bigl|\coderivative [\subdiff f](\optx|\opty)\bigr|^+
        &
        =
        \sup\setof{\norm{x^*}_{X^*}}{x^* \in \coderivative[\subdiff f](\optx|\opty)(y^*), \norm{y^*}_{Y^*} \le 1}
        \\
        &
        =
        \sup\setof{\norm{x^*}_{X^*}}{(x^*, -y^*) \in N_{\graph F}(x, y),\, \norm{y^*}_{Y^*} \le 1}.
    \end{aligned}
\end{equation*}
\begin{figure}
    \centering
    \begin{subfigure}[t]{.4\textwidth}
        \centering
        \begin{asy}
            real fup(real x){ return x^2*1.5; };
            real fdown(real x) { return cos(10*x)/20-.5-0.5*sqrt(0.7*sin(max(-x,0))); };
            real ell=.8;

            path pup=graph(fup, -1, 1);
            path pdown=graph(fdown, -1, 1);
            path pupext=graph(fup, -1.1, 1.15);
            path pdownext=graph(fdown, -1.15, 1.15);

            fill(pupext--reverse(pdownext)--cycle, lightfill);
            draw(pup);
            draw(pdown);

            real maxf=fup(0)+.5;
            real minf=fdown(0)-.5;
            clip((-2, maxf)--(2, maxf)--(2, minf)--(-2, minf)--cycle);

            pair up0=pt(fup, 0);
            pair down0=pt(fdown, 0);

            label("$F$", (up0+down0)/2);

            dot((0, fup(0)));

            draw((0, fup(0))--(0, fup(0)+.5), primalline+linewidth(1.1), Arrow);
        \end{asy}
        \caption{property is satisfied}
        \label{fig:regularity:morduk-yes}
    \end{subfigure}
    \hfil
    \begin{subfigure}[t]{.4\textwidth}
        \centering
        \begin{asy}
            real fup(real x){ return x>0 ? sqrt(x) : 0; };
            real fdown(real x) { return x^3/2-.5; };
            real ell=1;

            path pup=graph(fup, -1, 1);
            path pdown=graph(fdown, -1, 1);
            path pupext=graph(fup, -1.1, 1.15);
            path pdownext=graph(fdown, -1.15, 1.15);

            fill(pupext--reverse(pdownext)--cycle, lightfill);
            draw(pup);
            draw(pdown);

            real maxf=fup(0)+.5;
            real minf=fdown(0)-.5;
            clip((-2, maxf)--(2, maxf)--(2, minf)--(-2, minf)--cycle);

            pair up0=pt(fup, 0);
            pair down0=pt(fdown, 0);

            label("$F$", (up0+down0)/2, .5*S);

            dot((0, fup(0)));
            draw((0, fup(0))--(0, fup(0)+.5), primalline+linewidth(1.1), Arrow);
            draw((0, fup(0))--(-.5, fup(0)), primalline+linewidth(1.1), Arrow);
        \end{asy}
        \caption{property is not satisfied}
        \label{fig:regularity:morduk-no}
    \end{subfigure}
    \caption{Illustration of \cref{cor:regularity:morduk:finite}, where the arrows denote the directions contained in the normal cone.
        In (\subref{fig:regularity:morduk-yes}), $-y^* \in [0, \infty)$ but $x^*=0$, hence $|\coderivative F(x|y)|^+=0$ and the Aubin property is satisfied. In (\subref{fig:regularity:morduk-no}), we can take for $y^*=0$ any $x^* \in (-\infty, 0]$, hence $|\coderivative F(x|y)|^+=\infty$ and the Aubin property is violated.
    }
    \label{fig:regularity:morduk}
\end{figure}

We further illustrate the computation of the graphical modulus by returning to the subdifferential mapping of the absolute value function.
\begin{example}
    \label{ex:regularity:abs-subdiff:aubin}
    Recall from \cref{fig:graphical:absvalue:limiting} and \eqref{eq:graphical:absvalue:colimiting} in \cref{lemma:graphical:absvalue} that the limiting coderivative of the subdifferential of the absolute value function $f=\abs{\freevar}$ is given by
    \begin{equation}
        \label{eq:regularity:abs-subdiff:aubin:coderivative-repeat}
        \coderivative [\subdiff f]( x | y )(y^* ) =
        \begin{cases}
            \{0\}
            &\text{if }  x  \ne 0,\,  y  = \sign  x , \\
            \{0\}
            &\text{if }  x =0,\,
            yy^* > 0,\, \abs{y}=1,
            \\
            (-\infty, 0] y
            &\text{if }  x =0,\, y y^* < 0,\, \abs{ y } = 1, \\
            \R
            &\text{if }  x =0,\, y^*  = 0,\, \abs{ y } \le 1, \\
            \emptyset
            & \text{otherwise.}
        \end{cases}
    \end{equation}
    To study the graphical modulus, we follow \cref{cor:regularity:morduk:finite} and compute for $\opty \in \subdiff f(\optx)$ that
    \[
        \begin{aligned}
            \lip[\subdiff f](\optx|\opty)
            =
            \bigl|\coderivative [\subdiff f](\optx|\opty)\bigr|^+
            &
            =
            \sup\setof{\abs{x^*}}{x^* \in D^*[\subdiff f](x|y)(y^*),\, \abs{y^*} \le 1}
            \\
            &
            =
            \begin{cases}
                0
                &\text{if } \optx  \ne 0,\,  \opty  = \sign  \optx, \\
                \infty
                &\text{if } \optx =0, \opty \in [-1, 1]. \\
            \end{cases}
        \end{aligned}
    \]
    Thus $\subdiff f$ does not have the Aubin property at $\optx=0$ for any $\opty \in \subdiff f(0)$, but does have it away from zero.
\end{example}

The absolute value function -- and, indeed, all convex functions on $\R$ due to the convexity and monotonicity of the subdifferentials (\cref{lemma:monotone:convex,thm:monoton:subdiff}, respectively) -- is therefore \enquote{too nonsmooth} at zero for its subdifferential to have the Aubin property there.
We will return in \cref{ex:stability:implicit-morduk:affine} to a properly set-valued mapping that has the Aubin property, and turn now to the inverse property of metric regularity.

\begin{corollary}[Mordukhovich criterion for the modulus of metric regularity in finite dimensions]
    \label{cor:regularity:morduk-metric:finite}
    Let $X,Y$ be finite-dimensional Gateaux smooth Banach spaces and let $F: X \setto Y$ be such that $\graph F$ is closed near $(\optx, \opty)\in X\times Y$.
    Then
    \begin{equation*}
        \reg F(\optx|\opty)
        =
        |\inv{\coderivative F(\optx|\opty)}|^+.
    \end{equation*}
\end{corollary}

\begin{proof}
    By \cref{lemma:graphical:inverse}, we have
    \begin{equation*}
        \begin{aligned}
            |\coderivative \inv F(\opty|\optx)|^+
            &
            =\sup \setof{ \norm{y^*}_{Y^*}}{-y^* \in \coderivative \inv F(\opty|\optx)(-x^*),\  \norm{x^*}_{X^*} \le 1}
            \\
            &
            =\sup \setof{ \norm{y^*}_{Y^*}}{x^* \in \coderivative F(\optx|\opty)(y^*),\  \norm{x^*}_{X^*} \le 1}
            \\
            &=|\inv{[\coderivative F(\optx|\opty)]}|^+.
        \end{aligned}
    \end{equation*}
    The claim now follows by applying \cref{cor:regularity:morduk:finite} to $\inv F$ together with \cref{cor:regularity:morduk-metric}.
\end{proof}

We again illustrate this result for the subdifferential mapping of the absolute value function.
\begin{example}
    We continue from \cref{ex:regularity:abs-subdiff:aubin}.
    To study the modulus of metric regularity, we follow \cref{cor:regularity:morduk-metric:finite}.
    We first invert \eqref{eq:regularity:abs-subdiff:aubin:coderivative-repeat} to obtain
    \[
        \coderivative [\subdiff f](x|y)^{-1}(x^*)
        =
        \begin{cases}
            \R
            &\text{if }  x  \ne 0,\,  y  = \sign  x ,\, x^*=0,\\
            (-\infty, 0)\sign y
            &\text{if }  x =0,\, \abs{ y } = 1,\, yx^* \le 0, \\
            \{0\}
            &\text{if }  x =0,\, \abs{ y } \le 1, \\
            \emptyset
            & \text{otherwise.}
        \end{cases}
    \]
    We then compute for $\opty \in \subdiff f(\optx)$ that
    \[
        \begin{aligned}
            \reg[\subdiff f](\optx|\opty)
            =
            \bigl|\coderivative [\subdiff f](\optx|\opty)^{-1}\bigr|^+
            &
            =
            \sup\setof{\abs{y^*}}{y^* \in \coderivative[\subdiff f](\optx|\opty)^{-1}(x^*), \abs{x^*} \le 1}
            \\
            &
            =
            \begin{cases}
                \infty & \text{if } \abs{\opty}=1, \\
                0 & \text{if } \abs{\opty} < 1.
            \end{cases}
        \end{aligned}
    \]
    Thus $\subdiff f$ is not metrically regular away from zero, and is metrically regular at zero only for $\abs{\opty}<1$, i.e., $\opty$ is in the interior of the subdifferential.
\end{example}

The problem in the previous example is not the nonsmoothness but the lack of sufficient growth, as the next example shows.

\begin{example}
    Let us consider $f(x)=\tfrac{1}{2}x^2 + \abs{x}$.
    We then have by \cref{thm:subdiff:sum} and \cref{ex:convex:subdiff_abs} that
    \[
        \subdiff f(x) = \{x\} + \sign(x) = \begin{cases}
            \{x + 1\}   & \text{if }x > 0, \\
            [-1, 1] & \text{if }x=0, \\
            \{x - 1\}   & \text{if }x < 0.
        \end{cases}
    \]
    \Cref{ex:proximal:reell}\,\ref{ex:proximal:reell:ii} then immediately yields
    \[
        [\subdiff f]^{-1}(y) = [\Id + \subdiff |\freevar|]^{-1}(y) = \begin{cases}
            \{y - 1\} & \text{if }y > 1, \\
            \{0\} & \text{if }y \in [-1, 1], \\
            \{y + 1\} & \text{if }y < -1.
        \end{cases}
    \]
    This function is clearly Lipschitz continuous with constant $L=1$ and therefore has the Aubin property at every point.
    Since metric regularity is the Aubin property of the inverse mapping, $\subdiff f$ is metrically regular everywhere on its graph.
\end{example}

\begin{remark}
    Derivative-based characterizations of calmness and metric subregularity are significantly more involved than those of the Aubin property and metric regularity discussed above. We refer to \cite{henrion2002calmness,zheng2010metric,gfrerer2011first,gfrerer2016lipschitzian} for a few characterizations in special cases.
\end{remark}

To close this section, we relate the Mordukhovich criterion to the classical inverse function theorem (\cref{thm:inversefunctiontheorem}).

\begin{corollary}[inverse function theorem]\index{theorem!inverse function}
    Let $X,Y$ be reflexive and Gateaux smooth Banach spaces and let $F: X \to Y$ be continuously differentiable around $\optx \in X$.
    If $F'(\optx)^* \in \linear(Y^*; X^*)$ has a left-inverse $\starlinv{F'(\optx)} \in \linear(X^*; Y^*)$,
    then there exist $\kappa>0$ and $\delta>0$ such that for all $y \in \B(F(\optx), \delta)$ there exists a single-valued selection $J(y) \in \inv F(y)$ with
    \begin{equation*}
        \norm{\optx-J(y)}_X \le \kappa \norm{F(\optx)-y}_Y.
    \end{equation*}
\end{corollary}

\begin{proof}
    Let $\opty \defeq F(\optx)$.
    By \cref{thm:graphical:single} and the reflexivity of $X$ and $Y$,
    \begin{equation}
        \label{eq:regularity:inversef:single-valued-diff}
        \coderivative F(\optx|\opty)=\frechetCod F(\optx|\opty)=
            \{F'(\optx)^*\}.
    \end{equation}
    We have both $\coderivative \inv F(\opty|\optx)=\inv{[\coderivative F(\optx|\opty)]}$ and  $\frechetCod \inv F(\opty|\optx)=\inv{[\frechetCod F(\optx|\opty)]}$ by \cref{lemma:graphical:inverse}.
    Due to \eqref{eq:regularity:coderivative-inclusions}, this then implies that $\coderivative_M \inv F(\opty|\optx)\subset\coderivative \inv F(\opty|\optx)=\inv{[\coderivative F(\optx|\opty)]}$.
    The existence of a left-inverse implies that $F'(\optx)^*$ is injective, which together with \eqref{eq:regularity:inversef:single-valued-diff} yields \eqref{eq:regularity:morduk-metric:qc}.

    By the continuity of $F$, $\graph \inv F$ is closed near $(\opty, \optx)$.
    By \cref{lemma:colimiting:psnc:single:inverse}, $\inv F$ is PSNC at $\opty$ for $\optx$.
    Consequently, \cref{cor:regularity:morduk-metric:finite} shows that $F$ is metrically regular at $\optx$ for $\opty$.
    By the definition \eqref{eq:regularity:metric} of metrical regularity, there thus exists for any $\alt\kappa > \reg F(\optx|\opty)$ a $\delta>0$ such that
    \begin{equation*}
        \inf_{\alt x \in \inv F(y)} \norm{x-\alt x}_X \le \alt\kappa \norm{F(x)-y}_Y
        \quad
        (x \in \B(\optx, \delta), y \in \B(\opty, \delta)).
    \end{equation*}
    Taking in particular $x=\optx$ yields
    \begin{equation*}
        \inf_{\alt x \in \inv F(y)} \norm{\optx-\alt x}_X \le \alt\kappa \norm{F(\optx)-y}_Y
        \quad
        (y \in \B(F(\optx), \delta)).
    \end{equation*}
    Although the infimum might not be attained, this implies that we can take arbitrary $\kappa > \alt\kappa$ to obtain for any $y \in \B(F(\optx), \delta)$ the existence of some $J(y) \defeq \alt x \in \inv F(y)$ satisfying
    $
    \norm{\optx-\alt x}_X \le \kappa \norm{F(\optx)-y}_Y,
    $
    which is the claim.
\end{proof}

\addtocontents{toc}{\protect\enlargethispage{1cm}}

\chapter{Stability with respect to perturbations}
\label{chap:stability}

We now apply the Lipschitz-like properties of \cref{chap:regularity} to study the stability of optimization problems under perturbations. As a motivating problem, we recall the introductory problem \eqref{eq:intro:prob} and consider the mapping
\begin{equation*}
    j(x; y, \alpha) \defeq \frac{1}{2}\norm{Ax-y}_Y^2 + \alpha g(x).
\end{equation*}
Assuming that a minimizer $\optx = x(y,\alpha)$ of $x\mapsto j(x;y,\alpha)$ exists, we can ask further questions about \emph{stability}, i.e., the dependence of $\optx$ on $y$ and $\alpha$, in particular whether $\optx$ depends (Lipschitz-)continuously on these parameters. This is of particular relevance in \term[problem!inverse]{inverse problems}, which study the solution of ill-posed operator equations $Ax=y$ via families of approximate \term[problem!well-posed]{well-posed} problems. The central question of \term[regularization!theory of]{regularization theory} is whether $x(y,\alpha)$ converges to a solution $\hat x$ of the operator equation $A\hat x =\hat y$ as $y\to \hat y$ and $\alpha\to 0$.

We study the question of stability in \cref{sec:stability:perturbations}.
After deriving in \cref{sec:stability:subdifferentials} a convenient characterization of the metric subregularity of convex subdifferentials, we prove the convergence of minimizers in the sense of regularization theory.

\section{Stability with respect to perturbations}
\label{sec:stability:perturbations}

Let $X,P$ be Banach spaces and $f:X\times P\to \Rbar$. We then consider for some parameter $\opt p\in P$ the \term[problem!optimization!parametric]{parametric optimization problem}
\begin{equation*}
    \min_{x \in X} f(x; \opt p)
\end{equation*}
and study how a minimizer (or critical point) $\optx \in X$ behaves under perturbations of $\opt p$.
For this purpose, we introduce the set-valued \term[mapping!solution]{solution mapping} (or, if $x\mapsto f(x; p)$ is not convex, \term[mapping!critical point]{critical point mapping})
\begin{equation}
    \label{eq:stability:solution-map}
    S: P \setto X, \quad S(p) \defeq \{ x \in X \mid 0 \in \subdiff f(x; p) \},
\end{equation}
where $\subdiff$ is a suitable (convex or Clarke) subdifferential with respect to $x$ for fixed $p$.
We apply the concepts from \cref{sec:regularity:lipschitz} to this problem.
Specifically, if $S$ has the Aubin property at $\opt p$ for $\optx$, then \eqref{eq:regularity:aubin-distance} yields
\begin{equation*}
    \inf_{x \in S(p)} \norm{\optx - x}_X \le \kappa \norm{p-\opt p}_P
    \quad (p \in \B(\opt p, \delta))
\end{equation*}
for some $\delta,\kappa>0$.
In other words, the Aubin property of the solution mapping $S$ at $\opt p$ for $\optx$ implies the local Lipschitz stability of solutions $x=S(p)$ under perturbations $p$ around the parameter $\opt p$.
This of course begs the question when a solution mapping has the Aubin property.

\bigskip

We start with a simple special case. Returning to the motivation at the beginning of this chapter, $w\in \partial f(\alt x)$ is of course equivalent to $0\in \partial f(\alt x) - \{w\} = \partial (f-\dualprod{w}{\freevar}_X)(\alt x)$ since continuous linear mappings are differentiable. Such a perturbation of $f$ is called a \term[perturbation, tilt]{tilt perturbation}, with $w\in X^*$ called \term[parameter!tilt]{tilt parameter}.

To make this more precise, let $g: X \to \R$ be locally Lipschitz.
For a tilt parameter $p \in X^*$, we then define
\begin{equation}\label{eq:stability:tiltfunction}
    f(x; p)=g(x)-\dualprod{p}{x}_X
\end{equation}
and refer to the stability of minimizers (or critical points) of $f$ with respect to $p$ as \term[stability!tilt]{tilt stability}.
By \cref{thm:clarke:fermat,thm:clarke:sum}, the solution mapping for $f$ is
\begin{equation*}
    S(p)=\{x \in X \mid p \in \subdiff_C g(x)\}
    = \inv{(\subdiff_C g)}(p),
\end{equation*}
which thus has the Aubin property -- and $f$ is tilt-stable -- if and only if $\subdiff_C g$ is metrically regular at $\opt x$ for $0$, i.e., by \eqref{eq:regularity:metric} that there exist $\kappa,\delta>0$ such that
\begin{equation}
    \label{eq:stability:tilt-example}
    \dist(x, \inv{(\subdiff_C g)}(x^*)) \le \kappa \dist(\subdiff_C g(x), x^*)
    \quad (x^* \in \B(0, \delta);\, x \in \B(\opt x, \delta)).
\end{equation}

We illustrate this with two examples. The first concerns data stability of least squares fitting, which in Hilbert spaces can be formulated as tilt stability.
\begin{example}[data stability of least squares fitting]
    \label{example:stability:lsq}
    Let $X,Y$ be Hilbert spaces and $g(x)=\frac{1}{2}\norm{Ax-y}_Y^2$ for some $A \in \linear(X; Y)$ and $y \in Y$.
    Taking $p=A^*\dir y$ for some $\dir y \in Y$, we can write this in the form of \eqref{eq:stability:tiltfunction} via
    \begin{equation*}
        f(x; p)=g(x)-\iprod{A^*\dir y}{x}_X
        = \frac{1}{2}\norm{Ax-(y+\dir y)}_Y^2 - \iprod{y}{\dir y}_Y - \frac{1}{2}\norm{\dir y}_Y^2.
    \end{equation*}
    Data stability thus follows from the metric regularity of $\subdiff g$ at a minimizer $\opt x$ of the convex functional $g$. We have $\subdiff_C g(x)=\{A^*(Ax-y)\}$, so
    \begin{equation*}
        \inv{(\subdiff_C g)}(x^*)=\{ \alt x\in X \mid A^*A\alt x =A^*y+x^* \}.
    \end{equation*}
    Therefore \eqref{eq:stability:tilt-example} is equivalent to
    \begin{multline*}
        \inf_{\alt x\in X} \{ \norm{\alt x-x}_X \mid A^*A\alt x =A^*y+x^* \}
        \le \kappa \norm{A^*Ax-(A^*y+x^*)}_X
        \\
        (x^* \in \B(0, \delta);\, x \in \B(\opt x, \delta)).
    \end{multline*}
    If $A^*A$ has a bounded inverse $\inv{(A^*A)}\in \linear(X; X)$, then we can take $\kappa=\norm{\inv{(A^*A)}}_{\linear(X;X)}$ for any $\delta>0$.
    On the other hand, if $A^*A$ is not surjective, then there cannot be metric regularity (simply take an appropriate choice of $x^*\notin\range A^*A$).
\end{example}

For a genuinely nonsmooth example, we consider the (academic) problem of minimizing the (non-squared) norm on a Hilbert space.
\begin{example}[tilt stability of least norm fitting]
    Let $X$ be a Hilbert space and $g(x)=\norm{x-z}_X$ for some $z\in X$.
    To show tilt stability, we have to verify \eqref{eq:stability:tilt-example} for some $\kappa,\delta>0$.
    For $x \ne z$, we have $\subdiff g(x)=\{(x-z)/\norm{x-z}_X\}$, and for $x=z$, we have $\subdiff g(x)=\B(0, 1)$.
    Thus \eqref{eq:stability:tilt-example} reads
    \begin{gather}
        \label{eq:stability:tilt-stability-l1}
        \dist(x, \inv{(\subdiff g)}(x^*))
        \le \kappa \begin{cases}
            \Big\|\frac{x-z}{\norm{x-z}_X}-x^*\Big\|_X & \text{if } x \ne z, \\
            \dist(x^*, \B(0, 1)) & \text{if } x=z,
        \end{cases}
        \shortintertext{for all $x^* \in \B(0, \delta)$ and $x \in \B(\opt x, \delta)$ where}
        \nonumber
        \dist(x, \inv{(\subdiff g)}(x^*))
        =
        \begin{cases}
            \dist(x-z, x^*[0, \infty)) &\text{if } \norm{x^*}_{X}=1, \\
            \norm{x-z}_X &\text{if }\norm{x^*}_{X} < 1, \\
            \infty &\text{if }\norm{x^*}_{X}>1.
        \end{cases}
    \end{gather}
    As the inequality cannot hold if $\norm{x^*}_{X}>1$, we take $\delta \in (0, 1]$ to ensure that this does not happen. If $x=z$, then \eqref{eq:stability:tilt-stability-l1} trivially holds for any $\kappa>0$, both sides being zero. For $x^* \in \B(0, \delta)$ and $x \in \B(\opt x, \delta) \setminus \{z\}$, the inequality \eqref{eq:stability:tilt-stability-l1} reads
    \begin{equation*}
        \kappa \adaptnorm{\frac{x-z}{\norm{x-z}_X}-x^*}_X
        \ge
        \begin{cases}
            \dist(x-z, x^*[0, \infty)) &\text{if } \norm{x^*}_{X}=1, \\
            \norm{x-z}_X &\text{if } \norm{x^*}_{X} < 1.
        \end{cases}
    \end{equation*}
    Choosing $x^*=\lambda(x-z)/\norm{x-z}_X$, and letting $\lambda \upto 1$, we see that the inequality cannot hold unless $\delta \in (0, 1)$ (which prevents $\lambda \upto 1$).
    Thus, taking the infimum of the left-hand side over $\norm{x^*}_{X}\le\delta<1$ and the supremum of the right-hand side over $x \in \B(\opt x, \delta)$, the inequality holds if $\kappa(1-\delta) \ge \delta$.
    This can be satisfied for any $\kappa>0$ for sufficiently small $\delta \in (0, 1)$.

    Since $x^*\in X$ is comparable to the tilt parameter $p\in X$, this says that we can only stably ``tilt'' $g$ by an amount $\norm{p}_X<1$. If we tilt with $\norm{p}_X>1$, the tilted function has no minimizer, while for $\norm{p}_X=1$, every $x=z+tp$ for $t\ge 0$ is a minimizer.
\end{example}

\bigskip

We now return to the general solution mapping \eqref{eq:stability:solution-map}.
The following result applied to $F(x, p)\defeq\subdiff f(x; p)$ provides a general tool for our analysis.

\begin{theorem}
    \label{thm:stability:implicit-morduk}
    Let $P$, $X$, and $Y$ be reflexive and Gateaux smooth Banach spaces.
    For $F: X \times P \setto Y$, let
    \begin{equation*}
        S(p) \defeq \{x \in X \mid 0 \in F(x, p)\}.
    \end{equation*}
    Then $S$ has the Aubin property at $\opt p$ for $\optx \in S(\opt p)$ if
    \begin{equation}
        \label{eq:stability:implicit-morduk-cond}
        (0, p^*) \in \coderivative F(\optx, \opt p|0)(y^*) \implies y^* =0,\, p^*=0  \quad(y\in Y)
    \end{equation}
    and
    \begin{equation*}
        Q(y, p) \defeq \{x \in X \mid y \in F(x, p)\}
    \end{equation*}
    is PSNC at $(0, \opt p)$ for $\optx$.
\end{theorem}

\begin{proof}
    We have $S(p)=Q(0, p)$. Hence if we can show that $Q$ has the Aubin property at $(0, \opt p)$ for $\optx$, this will imply the Aubin property of $S$ at $\opt p$ for $\optx$ by simple restriction of the free variables in \cref{thm:regularity:aubin-equiv}\,\cref{item:regularity:aubin-inclusion} to the subspace $\{0\} \times P$.

    We do this by applying \cref{thm:regularity:morduk:iff} to $Q$, which holds if we can show that
    \begin{equation*}
        \coderivative_M Q(0, \opt p|\optx)(0) = \{0\}.
    \end{equation*}
    By \eqref{eq:regularity:coderivative-inclusions}, a sufficient assumption for this is that
    \begin{equation*}
        \coderivative Q(0, \opt p|\optx)(0) = \{0\},
    \end{equation*}
    which can equivalently be expressed as
    \begin{equation}
        \label{eq:stability:tmorduk-cond-cone}
        (y^* , p^*, 0) \in N_{\graph Q}(0, \opt p, \optx)
        \implies y^* = 0,\, p^* = 0.
    \end{equation}
    Now
    \begin{equation*}
        \graph Q=\{(y, p, x) \mid y \in F(x, p) \}
        =\pi \graph F
    \end{equation*}
    for the permutation $\pi(x, p, y) \defeq (y, p, x)$ (which applied to a set should be understood as applied to every element of that set).
    We thus also have
    \begin{equation*}
        N_{\graph Q}(y, p, x)=\pi N_{\graph F}(\pi(y, p, x)).
    \end{equation*}
    In particular, \eqref{eq:stability:tmorduk-cond-cone} becomes
    \begin{equation*}
        (0, p^*, y^*) \in N_{\graph F}(\optx, \opt p, 0)
        \implies y^* = 0,\, p^* = 0.
    \end{equation*}
    But this is equivalent to \eqref{eq:stability:implicit-morduk-cond}.
\end{proof}

\begin{remark}
    \label{rem:stability:implicit-morduk}
    \cref{thm:stability:implicit-morduk} is related to the classical implicit function theorem. If $F$ is graphically regular at $(\optx, \opt p, 0)$, it is also possible to derive explicit characterizations of $DS$ such as
    \begin{equation*}
        DS(\opt p|\optx)(\dir p)=\{ \dir x \in X \mid DF(\optx,\opt p|0)(\dir x,\dir p) \ni 0 \}.
    \end{equation*}
    For details in finite dimensions, we refer to \cite[Theorem 9.56, Proposition 8.41]{Rockafellar:1998}.
\end{remark}

We next consider an example of a simple solution mapping that has the Aubin property.

\begin{example}
    \label{ex:stability:implicit-morduk:affine}
    For fixed $c \in \R$, let
    \[
        S: \R^N \setto \R^N,
        \quad
        p \mapsto \{ x \in \R^N \mid \iprod{x}{p} = c \},
    \]
    i.e., the solution mapping to $0=F(x, p) \defeq \iprod{x}{p} - c$.
    Since we are in finite dimensions, \cref{lemma:colimiting:psnc:single-valued} yields that the mapping $Q$ defined in \cref{thm:stability:implicit-morduk} is PSNC.
    Let now $\opt p\in\R^N$ and $\opt x \in S(\opt p)$ be given.
    We then have by \cref{thm:graphical:single} that
    \[
        D^*F(\opt x, \opt p|0)(y^*)
        = (\opt p y^*, \opt x y^*) \qquad \text{for all }y^*\in \R.
    \]
    Hence the condition \eqref{eq:stability:implicit-morduk-cond} becomes
    \[
        0 = \opt p y^* \implies y^* = 0 \text{ and } \opt x y^* = 0,
    \]
    which is satisfied if $\opt p \ne 0$. On the other hand, we have $S(0)=\emptyset$ if $c \ne 0$.
    Hence $S$ has the Aubin property at every $\opt p$ for every $\optx \in S(\opt p)$ as long as $c\neq 0$.
\end{example}

We close this section by illustrating the requirements of \cref{thm:stability:implicit-morduk} for the stability of specific problems of the form \eqref{eq:intro:prob} with respect to the penalty parameter $\alpha$. (Naturally, these can be relaxed or made further explicit in more concrete situations.)
We consider for $\alpha > 0$ and $h,g:X\to\Rbar$ the problem
\[
    \min_{x\in X} h(x)+\alpha g(x).
\]
For this problem, we define the Clarke-critical point mapping
\begin{equation}\label{eq:stability:clarke-critical-mapping}
    S(\alpha) \defeq \{x \in X \mid 0 \in \subdiff_C(h+\alpha g)(x) \}.
\end{equation}
When the problem is convex, this coincides with the solution mapping.
Subject to a non-degeneracy condition, the next theorem yields a stability estimate for convex $g$ and smooth~$h$.
\begin{theorem}
    \label{thm:stability:regularization}
    Let $X$ be a finite-dimensional and Gateaux smooth Banach space
    and let $h:X\to\R $ be twice continuously differentiable and $g:X\to\Rbar$ be convex, proper, and lower semicontinuous.
    Suppose
    \begin{equation}
        \label{eq:stability:regularization-stability}
        0 \in  h''(\opt x)^*y + \opt\alpha D^*[\subdiff g](\opt x|-\inv{\opt\alpha} h'(\optx))(y)
        \implies y=0 \quad (y\in Y).
    \end{equation}
    Then $S$ has the Aubin property at $\opt\alpha$ for any $\opt x \in S(\opt\alpha)$.
\end{theorem}

\begin{proof}
    By \cref{thm:clarke:fermat,thm:clarke:frechet,thm:clarke:sum}, we can expand
    \begin{equation*}
        S(\alpha) = \{x \in X \mid 0 \in F(x; \alpha)\}\quad\text{for}\quad
        F(x; \alpha) \defeq h'(x) + \alpha \subdiff g(x).
    \end{equation*}
    To apply \cref{thm:stability:implicit-morduk} to prove the Aubin property, we need to verify its assumptions.
    First, by \cref{thm:colimiting:addition,thm:colimiting:product}, we have
    \begin{equation*}
        D^*F(\opt x; \opt\alpha|0)(y) =
        \begin{pmatrix}
            h''(\opt x)^*y + \opt\alpha D^*[\subdiff g](\opt x|-\inv{\opt\alpha} h'(\optx))(y) \\
            -\dualprod{h'(\optx)}{y}_X
        \end{pmatrix}.
    \end{equation*}
    Thus \eqref{eq:stability:implicit-morduk-cond} holds by \eqref{eq:stability:regularization-stability}.
    Furthermore, since $X^*\times \R$ is finite-dimensional, the PSNC holds at every $(y,\alpha)$ with $y\in F(\optx,\opt\alpha)$ and $\alpha>0$ by \cref{lemma:colimiting:psnc:single-valued}. Hence \cref{thm:stability:implicit-morduk} is indeed applicable and implies that $S$ has the Aubin property at $\opt\alpha$.
\end{proof}
\begin{corollary}
    Under the assumptions of \cref{thm:stability:regularization},
    \begin{equation*}
        \inf_{x\in S(\alpha)} \norm{\optx -x}_X \leq \kappa |\opt\alpha-\alpha|
    \end{equation*}
    for some $\kappa>0$ and all $\alpha$ sufficiently close to $\opt\alpha$.
\end{corollary}
\begin{proof}
    The claim follows directly from the definition \eqref{eq:regularity:aubin-distance} of the Aubin property for $S$ given by \eqref{eq:stability:clarke-critical-mapping} in $y=\optx\in S(\opt\alpha)$, which yields
    \begin{equation*}
        \inf_{x\in S(\alpha)} \norm{\optx -x}_X  = \dist(\optx, S(\alpha))  \leq \kappa \dist(S^{-1}(\opt x), \alpha) = \kappa |\opt\alpha-\alpha|.
        \qedhere
    \end{equation*}
\end{proof}
\section{Metric subregularity of convex subdifferentials}
\label{sec:stability:subdifferentials}

We recall from \eqref{eq:regularity:subregularity} that a set-valued mapping $H:X\setto X^*$ is metrically subregular at $\realoptx\in X$ for $\realoptw \in X^*$ if there exist $\delta>0$ and $\kappa>0$ such that
\begin{equation*}
    \dist(x, \inv H(\realoptw)) \le \kappa \dist(\realoptw, H(x))
    \quad (x \in \B(\realoptx, \delta)).
\end{equation*}
We also recall that the infimum of all $\kappa>0$ for which this inequality holds for some $\delta>0$ is denoted by $\subreg H(\realoptx|\realoptw)$, the modulus of (metric) subregularity of $H$ at $\realoptx$ for $\realoptw$. In the following, we will also make use of the \emph{squared} distance of $x\in X$ to a set $A\subset X$,
\begin{equation*}
    \dist^2(x,A) \defeq \inf_{\alt x\in A} \norm{x-\alt x}_X^2.
\end{equation*}
We then have the following characterization of metric subregularity of convex functionals.
\begin{theorem}
    \label{thm:stability:subregularity:convex}
    Let $g: X \to \Rbar$ be convex, proper, and lower semicontinuous and let $\realoptx \in X$ with $0 \in \subdiff g(\realoptx)$.
    If there exist $\gamma>0$ and $\delta>0$ such that
    \begin{equation}
        \label{eq:stability:subregularity:growth}
        g(x) \ge g(\realoptx) + \gamma\dist^2(x, \inv{[\subdiff g]}(0))
        \quad
        (x \in \B_X(\realoptx, \delta)),
    \end{equation}
    then $\subdiff g$ is metrically subregular at $\realoptx$ for $0$ with $\kappa=\inv\gamma$ and the same $\delta$.

    Conversely, if $\subdiff g$ is metrically subregular at $\realoptx$ for $0$ with some $\kappa,\delta>0$, then \eqref{eq:stability:subregularity:growth} holds for any $\gamma \in (0, 1/(4\kappa))$.
\end{theorem}

\begin{proof}
    Let first \eqref{eq:stability:subregularity:growth} hold for $\gamma,\delta>0$. We need to show that
    \begin{equation}
        \label{eq:stability:subregularity:convex-subreg}
        \gamma\dist(x, \inv{[\subdiff g]}(0)) \le \dist(0, \subdiff g(x))
        \quad (x \in \B(\realoptx, \delta)).
    \end{equation}
    To that end, let $x \in \B(\realoptx, \delta)$. Clearly, if $\subdiff g(x)=\emptyset$, there is nothing to prove. So assume that there exists an $x^* \in \subdiff g(x)$.
    Then $x \in \dom g$, so that \eqref{eq:stability:subregularity:growth} shows that $\dist^2(x, \inv{[\subdiff g]}(0)) < \infty$.
    Consequently $\inv{[\subdiff g]}(0) \ne \emptyset$.
    For each $\epsilon>0$, by the definition of the set-distance, we can therefore find $x_\epsilon \in \inv{[\subdiff g]}(0)$ such that
    \begin{equation}
        \label{eq:stability:subregularity:epsilon}
        \norm{x-x_\epsilon}_X \le \dist(x, \inv{[\subdiff g]}(0)) + \epsilon.
    \end{equation}
    By the definition of the convex subdifferential and $\realoptx, x_\epsilon \in \argmin g$, we have
    \begin{equation*}
        \dualprod{x^*}{x-x_\epsilon}_X \ge g(x)-g(x_\epsilon) = g(x)-g(\realoptx).
    \end{equation*}
    Combined with \eqref{eq:stability:subregularity:growth} and  \eqref{eq:stability:subregularity:epsilon}, this yields
    \begin{equation*}
        \begin{aligned}
            \gamma\dist^2(x, \inv{[\subdiff g]}(0))
            &
            \le
            \dualprod{x^*}{x-x_\epsilon}_{X}
            \\
            &
            \le
            \norm{x^*}_{X^*}\norm{x-x_\epsilon}_X
            \le \norm{x^*}_{X^*}(\dist(x, \inv{[\subdiff g]}(0)) + \epsilon).
        \end{aligned}
    \end{equation*}
    Since $\epsilon>0$ was arbitrary and $\norm{x^*}_{X^*} \le \dist(0, \subdiff g(x))$,
    we obtain \eqref{eq:stability:subregularity:convex-subreg}.

    Conversely, let $\subdiff g$ be metrically subregular at $\realoptx$ for $0$ for some parameters $\kappa,\delta>0$. Take any $\gamma \in (0, 1/(4\kappa))$.
    We argue by contradiction. Assume that \eqref{eq:stability:subregularity:growth} does not hold.
    Then we can find some $\alt x \in \B(\realoptx, 2\delta/3)$ such that
    \begin{equation}
        \label{eq:stability:subregularity:contradiction}
        g(\alt x) < g(\realoptx) + \gamma\dist^2(\alt x, \inv{[\subdiff g]}(0)).
    \end{equation}
    However, $\realoptx$ is a minimizer of $g$, so necessarily $\gamma\dist^2(\alt x, \inv{[\subdiff g]}(0))>0$.
    By Ekeland's variational principle (\cref{thm:variation:ekeland}), we can thus find $y \in X$ satisfying
    \begin{equation}
        \label{eq:stability:subregularity:ekeland1}
        \norm{y-\alt x}_X \le \frac{1}{2}\dist(\alt x, \inv{[\subdiff g]}(0))
    \end{equation}
    and for all $x \in X$ that
    \begin{equation*}
        g(x) \ge g(y) - \frac{\gamma\dist^2(\alt x, \inv{[\subdiff g]}(0))}{\frac{1}{2}\dist(\alt x, \inv{[\subdiff g]}(0))}\norm{x-y}_X
        =
        g(y) - 2\gamma \dist(\alt x, \inv{[\subdiff g]}(0)) \norm{x-y}_X.
    \end{equation*}
    It follows that $y$ minimizes $g +  2\gamma \dist(\alt x, \inv{[\subdiff g]}(0)) \norm{\freevar-y}_X$, which by \cref{thm:convex:fermat,thm:subdifferential:norm,thm:subdiff:sum} is equivalent to
    $
    0 \in \subdiff g(y) + 2\gamma\dist(\alt x, \inv{[\subdiff g]}(0)) \B_{X^*}.
    $
    Hence we can find some $y^* \in \subdiff g(y)$ satisfying
    $
    \norm{y^*}_{X^*} \le 2\gamma\dist(\alt x, \inv{[\subdiff g]}(0)).
    $
    Using \eqref{eq:stability:subregularity:ekeland1}, we now obtain
    \begin{equation*}
        \begin{aligned}
            2\kappa\dist(0, \subdiff g(y))
            &
            <
            \inv{(2\gamma)} \dist(0, \subdiff g(y))
            \\
            &
            \le
            \inv{(2\gamma)} \norm{y^*}_{X^*}
            \le
            \dist(\alt x, \inv{[\subdiff g]}(0))
            \\
            &
            =
            2\dist(\alt x, \inv{[\subdiff g]}(0))
            -\dist(\alt x, \inv{[\subdiff g]}(0))
            \\
            &
            \le 2\norm{y-\alt x}_X
            + 2\dist(y, \inv{[\subdiff g]}(0))
            - \dist(\alt x, \inv{[\subdiff g]}(0))
            \\
            &
            \le 2\dist(y, \inv{[\subdiff g]}(0)).
        \end{aligned}
    \end{equation*}
    By \eqref{eq:stability:subregularity:ekeland1} and our choice of $\alt x \in \B(\realoptx, 2\delta/3)$,
    \begin{equation*}
        \norm{y-\realoptx}_X \le \norm{y-\alt x}_X + \norm{\alt x-\realoptx}_X \le \frac{3}{2}\norm{\alt x-\realoptx}_X \le \delta.
    \end{equation*}
    Therefore $y \in \B(\realoptx, \delta)$ violates the assumed metric subregularity \eqref{eq:stability:subregularity:convex-subreg} with the factor $\tilde\gamma$, and hence \eqref{eq:stability:subregularity:growth} holds.
\end{proof}

Applying \cref{thm:stability:subregularity:convex} to $x \mapsto g(x) + \dualprod{\realoptx^*}{x}_X$ now yields the following characterization due to \cite{artacho2013metric}.
\begin{corollary}
    \label{cor:stability:subregularity:convex}
    Let $g: X \to \Rbar$ be convex, proper, and lower semicontinuous and let $\realoptx\in X$ and $\realoptx^* \in \subdiff g(\realoptx)$.
    If there exist $\gamma>0$ and $\delta>0$ such that
    \begin{equation}
        \label{eq:stability:subregularity:growth:full}
        g(x) \ge g(\realoptx) +\dualprod{\realoptx^*}{x-\realoptx}_X + \gamma\dist^2(x, \inv{[\subdiff g]}(\realoptx^*))
        \quad
        (x \in \B_X(\realoptx, \delta)),
    \end{equation}
    then $\subdiff g$ is metrically subregular at $\realoptx$ for $\realoptx^*$ with $\kappa=\inv\gamma$ and the same $\delta$.

    Conversely, if $\subdiff g$ is metrically subregular at $\realoptx$ for $\realoptx^*$ with some $\kappa,\delta>0$, then \eqref{eq:stability:subregularity:growth:full} holds for any $\gamma \in (0, 1/(4\kappa))$.
\end{corollary}

If we denote by $\hat\gamma(\realoptx|\realoptx^*)$ the supremum of $\gamma>0$ for which \eqref{eq:stability:subregularity:growth:full} holds for some $\delta>0$, then we obtain the following estimate involving the modulus of subregularity.

\begin{corollary}
    Let $g: X \to \Rbar$ be convex, proper, and lower semicontinuous and let $\realoptx\in X$ and $\realoptx^* \in \subdiff g(\realoptx)$.
    Then
    \[
        \subreg \subdiff g(\realoptx|\realoptx^*)
        \le
        \inv{\hat\gamma(\realoptx|\realoptx^*)}
        \le
        4\subreg \subdiff g(\realoptx|\realoptx^*).
    \]
\end{corollary}

\begin{remark}[strong metric subregularity]\index{subregularity, metric!strong}
    As in \cref{rem:strong-metric-subregularity}, we can also characterize \emph{strong} metric subregularity using a strong notion of local subdifferentiability. In the setting of \cref{cor:stability:subregularity:convex}, it was shown in \cite{artacho2013metric} that
    strong metric subregularity of $\subdiff g$ at $\realoptx$ for $\realoptx^*$ is equivalent to
    \begin{equation}
        \label{eq:stability:subregularity:strong}
        g(x) \ge g(\realoptx) +\dualprod{\realoptx^*}{x-\realoptx}_X + \gamma\norm{x-\realoptx}_X^2
        \quad
        (x \in \B_X(\realoptx, \delta)),
    \end{equation}
    i.e., a local form of strong subdifferentiability.
    Compared to the characterization of metric subregularity in \eqref{eq:stability:subregularity:growth:full}, intuitively the strong version does not \enquote{squeeze} $\inv{[\subdiff g]}(\realoptx^*)$ into a single point.

    Strong metric subregularity may almost trivially be used in the convergence proofs of \cref{part:convex,chap:nlpdps} as a relaxation of strong convexity; compare \cite{tuomov-nlpdhgm-general}.
    Also observe that \eqref{eq:stability:subregularity:strong} can be expressed in terms of the \term[divergence, Bregman]{Bregman divergence} (see \cref{sec:gap:ergodic:bregman}) as
    \begin{equation*}
        B_g^{\realoptx^*}(x, \realoptx) \ge \gamma\norm{x-\realoptx}_X^2
        \quad
        (x \in \B_X(\realoptx, \delta)),
    \end{equation*}
    i.e., that $B_g^{\realoptx^*}$ is \term[divergence, Bregman!elliptic]{elliptic} at $\realoptx$ in the sense of \cite{tuomov-firstorder}.
    In optimization methods based on preconditioning by Bregman divergences\index{divergence, Bregman} instead of the linear preconditioner $\Precond$ as discussed in \cref{sec:gap:ergodic:bregman}, this generalizes the positive definiteness requirement on $\Precond$.
\end{remark}

\section{Tikhonov-type regularization of inverse problems}
\label{sec:stability:tikhonov}

Let now the data $y^\delta$ depend on a \term{noise level} $\delta>0$, and consider for a corresponding parameter $\alpha_\delta>0$ the problem
\begin{equation}
    \label{eq:stability:regtheory:tikhonov}
    \min_{x \in X}~ \frac{1}{2}\norm{Ax-y^\delta}_Y^2 + \alpha_\delta g(x),
\end{equation}
where $A \in \linear(X; Y)$ between a Banach space $X$ and a Hilbert space $Y$.
This problem is called a \term[regularization!Tikhonov-type]{Tikhonov-type regularization} of the inverse problem $Ax=y^\delta$.
If $g(x)=\frac{1}{2}\norm{\freevar}_X^2$ with $X$ a Hilbert space, we talk simply of \term[regularization!Tikhonov]{Tikhonov regularization}.

We assume for some true data $\hat y$ that
\begin{equation}
    \label{eq:stability:regtheory:noise}
    \norm{y-\hat y}_Y \le \delta.
\end{equation}
Suppose there exists a solution $\hat x$ to the problem
\begin{equation}
    \label{eq:stability:regtheory:ground-truth}
    \min_{x \in C} g(x)
    \quad\text{where}\quad
    C \defeq \{x \in X \mid Ax=\hat y\}.
\end{equation}
Denote by $\hat X$ the set of solutions to \eqref{eq:stability:regtheory:ground-truth}.
In inverse problems, the question whether solutions $x_\delta$ to the Tikhonov-type problem \eqref{eq:stability:regtheory:tikhonov} converge to some $\hat x \in \hat X$ is a topic of \term[regularization!theory of]{regularization theory}.
The condition \eqref{eq:stability:regtheory:source-condition} of the next lemma is known as a \term[condition!source]{source condition} in that context.

\begin{lemma}
    Suppose $A \in \linear(X; Y)$ and that $g: X \to \Rbar$ is convex, proper, and lower semicontinuous with $\interior \dom g \isect C \ne \emptyset$.
    We have $\hat x \in \hat X$ if and only if there exists $\hat w \in Y$ such that
    \begin{equation}
        \label{eq:stability:regtheory:source-condition}
        A \hat x = \hat y
        \quad\text{and}\quad
        - A^* \hat w \in \subdiff g(\hat x).
    \end{equation}
\end{lemma}

\begin{proof}
    The condition $\interior \dom g \isect C \ne \emptyset$ guarantees that the sum rule \cref{thm:subdiff:sum} holds as an equality for $\delta_C + g$.
    Writing $\delta_C(x) = \delta_{\{\hat y\}}(Ax)$, and using the chain rule (\cref{thm:convex:chain}) and the fact that
    \[
        \subdiff\delta_{\{\hat y\}}(y) = \begin{cases}
            Y & \text{if }y = \hat y, \\
            \emptyset & \text{otherwise},
        \end{cases}
    \]
    we therefore obtain
    \[
        \subdiff[\delta_C + g](x)
        = A^* Y + \subdiff g(x)
        \quad\text{whenever}\quad
        Ax=\hat y.
    \]
    Thus $0 \in \subdiff[\delta_C + g](x)$ whenever \eqref{eq:stability:regtheory:source-condition} holds.
    Now the Fermat principle of \cref{thm:convex:fermat} establishes the claim.
\end{proof}

The next result characterizes convergence.
For brevity we write
\[
    f_\delta(x) \defeq \frac{1}{2}\norm{Ax-y^\delta}_Y^2.
\]
The condition \eqref{eq:stability:regtheory:approxmin} in the next theorem is satisfied in particular if $x_\delta$ is an $e_\delta$-minimizer of $f_\delta + \alpha_\delta g$.
Observe that by the definition of $\hat X$ after \eqref{eq:stability:regtheory:ground-truth}, the value on the right-hand side of \eqref{eq:stability:regtheory:approxmin} is independent of the specific choice of $\hat x \in \hat X$.
We directly assume the characterization \eqref{eq:stability:subregularity:growth:full} of metric subregularity to be able to use an optimal modulus $\gamma$ for which the characterization holds.

\begin{theorem}
    \label{thm:stability:regtheory:regularization-theory}
    Let $A \in \linear(X; Y)$ and $g: X \to \Rbar$ be convex, proper, and lower semicontinuous with $\interior \dom g \isect C \ne \emptyset$.
    Suppose \eqref{eq:stability:regtheory:noise} holds, $\hat x \in X$ satisfies \eqref{eq:stability:regtheory:source-condition}, and that for every $\delta>0$, for some $e_\delta>0$ there exists $x_\delta \in X$ satisfying
    \begin{equation}
        \label{eq:stability:regtheory:approxmin}
        [f_\delta + \alpha_\delta g](x_\delta)
        \le
        [f_\delta + \alpha_\delta g](\hat x)
        + e_\delta.
    \end{equation}
    Suppose for some $\tilde\delta>0$ and $\tilde X \subset \hat X$ that for all $\delta \in (0, \tilde\delta)$ and $\tilde x \in \tilde X$, the subdifferential mapping $\subdiff[f_\delta + \alpha_\delta g]$ satisfies \eqref{eq:stability:subregularity:growth:full} at $\tilde x$ for $f_\delta'(\tilde x)-\alpha_\delta A^*\hat w$ in the neighborhood $U_{\tilde x}$ (independent of $\delta$) with the factor $\gamma>0$ (independent of both $\delta$ and $\tilde x$) with respect to the norm
    \[
        \norm{x}_\delta \defeq \sqrt{\norm{Ax}_Y^2 + \alpha_\delta \norm{x}_X^2} \quad (x \in X).
    \]
    Assume further for some $\rho>0$ that
    \begin{equation}
        \label{eq:stability:regtheory:urho-containment}
        \Union_{\tilde x \in \tilde X} U_{\tilde x} \supset U_\rho \defeq \{x \in X \mid \norm{A(x-\hat x)} \le \rho,\, R(x) \le R(\hat x) + \rho\}.
    \end{equation}
    Then there exists $\bar\delta>0$ such that
    \begin{equation*}
        \dist^2(x_\delta, \hat X)
        \le
        \frac{e_\delta}{\gamma\alpha_\delta}
        +\frac{\delta^2}{2\gamma^2\alpha_\delta}
        +\frac{\alpha_\delta}{2\gamma^2}\norm{\hat w}_Y^2
        \quad(\delta \in (0, \bar\delta)).
    \end{equation*}
\end{theorem}

\begin{proof}
    Since $A\hat x=\hat y$, using Young's inequality, \eqref{eq:stability:regtheory:approxmin}, and \eqref{eq:stability:regtheory:noise}, we have
    \[
        \begin{aligned}
            \frac{1}{2}\norm{A(x_\delta-\hat x)}_Y^2
            + 2\alpha_\delta g(x_\delta)
            &
            \le
            \norm{Ax_\delta-y^\delta}_Y^2
            + 2\alpha_\delta g(x_\delta)
            +\norm{y^\delta-\hat y}_Y^2
            \\
            &
            \le
            2e_\delta + 2\norm{y^\delta-\hat y}_Y^2 + 2\alpha_\delta g(\hat x)
            \\
            &
            \le 2(e_\delta + \delta^2 + \alpha_\delta g(\hat x)).
        \end{aligned}
    \]
    Thus both
    \[
        \norm{A(x_\delta-\hat x)}_Y^2 \le 4(e_\delta + \delta^2 + \alpha_\delta g(\hat x))
        \quad\text{and}\quad
        g(x_\delta) \le g(\hat x) + \frac{e_\delta+\delta^{2}}{\alpha_\delta}.
    \]
    This implies the existence of $\bar\delta \in (0, \tilde\delta]$ such that $x_\delta \in U_\rho$ for $\delta \in (0, \bar\delta)$.
    Consequently \eqref{eq:stability:regtheory:urho-containment} establishes for every such $\delta$ an element $\hat x_\delta \in \hat X$ such that $x_\delta \in U_{\hat x_\delta}$.
    By $f_\delta+\alpha_\delta g$ satisfying \eqref{eq:stability:subregularity:growth:full} at $\hat x_\delta$ for $f_\delta'(\tilde x)-\alpha_\delta A^*\hat w$ for such $\delta$,
    therefore
    \begin{equation}
        \label{eq:stability:regtheory:strong-local-subdiff}
        [f_\delta+\alpha_\delta g](x_\delta)-[f_\delta+\alpha_\delta g](\hat x_\delta)
        \ge \dualprod{f_\delta'(\hat x_\delta) - \alpha_\delta  A^*\hat w}{x_\delta-\hat x_\delta}_X
        + \gamma\dist_\delta^2(x_\delta, \hat X),
    \end{equation}
    where $\dist_\delta$ denotes the distance-to-set function with respect to $\norm{\freevar}_\delta$.

    We next expand
    \[
        f_\delta'(\hat x_\delta) - \alpha_\delta A^*\hat w
        =A^*(A\hat x_\delta-y^\delta-\alpha_\delta \hat w)
        =A^*(\hat y-y^\delta-\alpha_\delta \hat w).
    \]
    Hence \eqref{eq:stability:regtheory:source-condition} and \eqref{eq:stability:regtheory:strong-local-subdiff} establish
    \begin{equation*}
        \begin{aligned}[t]
            e_\delta
            &
            \ge \dualprod{f_\delta'(\hat x_\delta) - \alpha_\delta  A^*\hat w}{x_\delta-\hat x}_X
            + \gamma\dist_\delta^2(x_\delta, \hat X)
            \\
            &
            = \iprod{\hat y-y^\delta-\alpha_\delta \hat w}{A(x_\delta-\hat x_\delta)}_Y
            + \gamma \inf_{\bar x \in \hat X} \left(
                \norm{A(x_\delta-\bar x)}_Y^2
                + \alpha_\delta\norm{x_\delta-\bar x}_X^2
            \right).
        \end{aligned}
    \end{equation*}
    Since $A\bar x=A\hat x_\delta$ due to $\hat X \subset C$, distributing the $\inf$ over the entire right-hand side and using Young's inequality establishes
    \begin{equation*}
        \begin{aligned}[t]
            e_\delta
            &
            \ge
            \inf_{\bar x \in \hat X}\left(
                \iprod{\hat y-y^\delta-\alpha_\delta \hat w}{A(x_\delta-\bar x)}_Y
                + \gamma\norm{A(x_\delta-\bar x)}_Y^2
                + \gamma\alpha_\delta\norm{x_\delta-\bar x}_X^2
            \right)
            \\
            &
            \ge
            \inf_{\bar x \in \hat X}\left(
                -\frac{1}{4\gamma}\norm{\hat y-y^\delta-\alpha_\delta \hat w}_Y^2
                +\gamma\alpha_\delta\norm{x_\delta-\bar x}_X^2
            \right).
        \end{aligned}
    \end{equation*}
    Thus, again using Young's inequality and \eqref{eq:stability:regtheory:noise}, we obtain
    \begin{equation*}
        \dist^2(x_\delta, \hat X)
        \le
        \frac{e_\delta}{\gamma\alpha_\delta}
        +\frac{1}{4\gamma^2\alpha_\delta}\norm{\hat y-y^\delta-\alpha_\delta \hat w}_Y^2
        \le
        \frac{e_\delta}{\gamma\alpha_\delta}
        +\frac{\delta^2}{2\gamma^2\alpha_\delta}
        +\frac{\alpha_\delta}{2\gamma^2}\norm{\hat w}_Y^2.
    \end{equation*}
    This is the claim.
\end{proof}

Immediately we obtain the following characterization of convergence of regularized solutions.

\begin{corollary}
    Under the assumptions of \cref{thm:stability:regtheory:regularization-theory}, if
    \[
        \lim_{\delta \downto 0} \left(
            \alpha_\delta, \frac{\delta^2}{\alpha_\delta}, \frac{e_\delta}{\alpha_\delta}
        \right)=0,
    \]
    then
    \[
        \lim_{\delta \downto 0} \dist(x_\delta, \hat X) = 0.
    \]
\end{corollary}

\begin{remark}
    For an introduction to \index{problem!inverse}inverse problems, we refer to \cite{Hanke:2017,Mueller:2012,ClasonIP:2020};
    a classical treatise on regularization theory is \cite{Engl} with Banach spaces and other advanced aspects covered in \cite{Kaltenbacher:2008,Schuster:2012,ItoJin}; see also \crefrange{chap:sparse}{chap:tv} and the remarks therein.
    Our specialized account is based on \cite{tuomov-regtheory},
    which also shows that using \emph{strong} metric subregularity in \cref{thm:stability:regtheory:regularization-theory} in place of metric subregularity yields convergence to a specific $\hat x \in \hat X$ instead of the set $\hat X$.
    Those results also relax the requirement $\Union_{\tilde x \in \tilde X} U_{\tilde x} \supset U_\rho$ through assumptions of weak(-$*$) closedness and openness.
\end{remark}

\chapter{Splitting methods: faster convergence from regularity}\label{chap:fasterconvergence}

As we have seen in \cref{chap:testing}, proximal point and splitting methods can be accelerated if at least one of the involved functionals is strongly convex. However, this can be a too strong requirement, and we will show in this chapter how faster convergence (even without acceleration) can be shown under the weaker requirements of metric subregularity or strong submonotonicity. We begin in \cref{sec:faster:subregularity} by introducing the latter notion before illustrating in \cref{sec:faster:acceleration} the effect of the two properties on splitting methods by showing local linear convergence of forward-backward splitting.

\section{Submonotonicity of convex subdifferentials}
\label{sec:faster:subregularity}

Throughout this section, let $X$ be a Banach space and $G:X\to\Rbar$ be convex, proper, and lower semicontinuous. Our goal is now to give conditions for metric subregularity and strong submonotonicity of $\subdiff G:X\setto X^*$ at a critical point $\realoptx\in X$ with $0\in \subdiff G(\realoptx)$.

Recall the characterization of metric subregularity of a convex subdifferential shown in \cref{sec:stability:subdifferentials}. As a weaker alternative to that result, we now relax the strong monotonicity assumption of \cref{chap:testing} more directly.
We say that a set-valued mapping $H: X \setto X^*$ is $(\gamma, \theta)$-\term[mapping!submonotone!strongly]{strongly submonotone} \emph{at $\realoptx$ for $\realoptx^* \in H(\realoptx)$} with $\theta \ge \gamma > 0$ if there exists $\delta>0$ such that for all $x \in \B_X(\realoptx, \delta)$ and $x^* \in H(x) \isect \B_{X^*}(\realoptx^*, \delta)$,
\begin{equation}
    \label{eq:faster:submonotonicity}
    \inf_{\optx \in \inv H(\realoptx^*)}
    \left(
        \dualprod{x^*-\realoptx^*}{x-\optx}_{X}
        + (\theta-\gamma)\norm{x-\optx}_X^2
    \right)
    \ge \theta\dist^2(x, \inv H(\realoptx^*)).
\end{equation}
If this only holds for $\theta\ge \gamma=0$, then we call $H$ \term[mapping!submonotone]{submonotone} at $\realoptx$ for $\realoptx^*$.

Clearly, (strong) monotonicity (see \cref{thm:smoothness:strong-convexity}) implies (strong) submonotonicity at any $\realoptx \in X$ and $\realoptx^* \in H(\realoptx)$. However, subdifferentials of convex functionals need not be strongly monotone. The next theorem shows that local second-order growth away from the set of minimizers implies strong submonotonicity of such subdifferentials at any minimizer $\realoptx$ for $\realoptx^*=0$, which is the monotonicity-based analogue of the characterization of metric subregularity in \cref{thm:stability:subregularity:convex}.

\begin{theorem}
    \label{theorem:faster:submonotonicity:convex}
    Let $G: X \to \Rbar$ be convex, proper, and lower semicontinuous and let $\realoptx\in X$ with $0 \in \subdiff G(\realoptx)$.
    If there exists $\delta>0$ such that
    \begin{equation}
        \label{eq:faster:submonotonicity:growth}
        G(x) \ge G(\realoptx) + \gamma\dist^2(x, \inv{[\subdiff G]}(0))
        \quad
        (x \in \B_X(\realoptx, \delta)),
    \end{equation}
    then $\subdiff G$ is $(\gamma,\theta)$-strongly submonotone at $\realoptx$ for $0$ for any $\theta\geq \gamma$.
\end{theorem}

\begin{proof}
    Since $\theta\geq \gamma$, \eqref{eq:faster:submonotonicity:growth} is equivalent to
    \begin{equation}
        \label{eq:faster:submonotonicity:growth-rearranged}
        \inf_{\optx \in \inv{[\subdiff G]}(0)}
        \left(
            G(x)-G(\realoptx)
            +
            (\theta-\gamma)\norm{x-\optx}_X^2
        \right)
        \ge \theta\dist^2(x, \inv{[\subdiff G]}(0))
    \end{equation}
    for all $x \in \B_X(\realoptx, \delta)$.
    By the definition of the convex subdifferential, we have for all $\optx \in \inv{[\subdiff G]}(0)$ and $\realoptx^*=0$ that
    \begin{equation*}
        \dualprod{x^*-\realoptx^*}{x-\optx}_{X} \ge G(x)-G(\optx)=G(x)-G(\realoptx).
    \end{equation*}
    Inserting this into \eqref{eq:faster:submonotonicity:growth-rearranged} yields the definition \eqref{eq:faster:submonotonicity} of strong submonotonicity for $H=\subdiff G$.
\end{proof}
Together with \cref{thm:stability:subregularity:convex}, this shows that for convex subdifferentials, metric subregularity implies strong submonotonicity, which is thus a weaker property.

We conclude this section by showing that the subdifferentials of the indicator functional of the finite-dimensional unit ball and of the absolute value function are both subregular and strongly submonotone. Note that neither of these subdifferentials is strongly monotone in the conventional sense. Here we restrict ourselves to showing $(\gamma,\gamma)$-strong submonotonicity for some $(\realoptx, \realoptx^*) \in \graph \subdiff G$, i.e., that there exists $\delta>0$ such that
\begin{equation}
    \label{eq:faster:convex-submono}
    \iprod{x^*-\realoptx^*}{x-\realoptx}_X \ge
    \gamma \dist^2(x, \inv{[\subdiff G]}(\realoptx^*))
    \quad (x \in \B(\realoptx,\delta),\, x^* \in \subdiff G(x)).
\end{equation}

\begin{lemma}
    \label{lemma:subreg-ball-indicator}
    Let $G \defeq \delta_{ \B(0, \alpha)}$ on $(\R^N,\norm{\cdot}_2)$ and $(\realoptx, \realoptx^*) \in \graph \subdiff G$. Then $\subdiff G$ is
    \begin{enumerate}[label=(\roman*)]
        \item\label{item:faster:ball-indicator:subreg}
            metrically subregular at $\realoptx$ for $\realoptx^*$ for any $\delta \in (0, \alpha]$ and
            \begin{equation*}
                \kappa \ge
                \begin{cases}
                    2\alpha/\norm{\realoptx^*}_2 &\text{if } \realoptx^* \ne 0, \\
                    0 & \text{if }\realoptx^*=0;
                \end{cases}
            \end{equation*}
        \item\label{item:faster:ball-indicator:submono}
            $(\gamma, \gamma)$-strongly submonotone at $\realoptx$ for $\realoptx^*$ for any $\delta>0$ and
            \begin{equation*}
                \gamma \le
                \begin{cases}
                    \norm{\realoptx^*}_2/(2\alpha) & \text{if }\realoptx^* \ne 0, \\
                    \infty & \text{if }\realoptx^*=0.
                \end{cases}
            \end{equation*}
    \end{enumerate}
\end{lemma}

\begin{proof}
    We first verify \eqref{eq:stability:subregularity:growth:full} for $\delta=\alpha$ and $\gamma=\inv\kappa$ as stated.
    To that end, let $x \in \B(0, \alpha)$.
    If $\realoptx^*=0$, then \eqref{eq:stability:subregularity:growth:full} trivially holds by the subdifferentiability of $G$ and $\dist^2(x, \inv{[\subdiff G]}(\realoptx^*))=\dist^2(x, \B(0, \alpha))=0$.
    Let therefore $\realoptx^* \ne 0$. Then $\inv{[\subdiff G]}(\realoptx^*)=\{\realoptx\}$ as well as $\norm{\realoptx}_2=\alpha$ and $\realoptx^*=\beta \realoptx$ for $\beta=\norm{\realoptx^*}_2/\norm{\realoptx}_2$.
    Since $\gamma \le \norm{\realoptx^*}_2/(2\alpha)$, we have $\beta \ge 2\gamma$.
    Then $\norm{x}_2 \le \alpha$ yields
    \begin{equation*}
        \begin{aligned}[t]
            \gamma \dist^2(x, \inv{[\subdiff G]}(\realoptx^*))
            &
            =
            \gamma \norm{x-\realoptx}_2^2
            \\
            &
            \le \beta \iprod{\realoptx}{\realoptx-x}_2 - \frac{\beta}{2} \norm{\realoptx}_2^2 + \frac{\beta}{2}\norm{x}_2^2
            \\
            &
            \le  \beta \iprod{\realoptx}{\realoptx-x}_2
            \\
            &
            = \iprod{\realoptx^*}{\realoptx-x}_2
            \\
            &
            \le
            \iprod{\realoptx^*}{\realoptx-x}_2 + G(x)-G(\realoptx).
        \end{aligned}
    \end{equation*}
    Since $\dom G=\B(0,\alpha)$, this shows that \eqref{eq:stability:subregularity:growth:full} holds for any $\delta>0$.

    \Cref{cor:stability:subregularity:convex} now yields \cref{item:faster:ball-indicator:subreg}. Adding
    \begin{equation*}
        G(\realoptx)-G(x) \ge \iprod{x^*}{\realoptx-x}_2
        \quad (x^* \in \subdiff G(x))
    \end{equation*}
    to \eqref{eq:stability:subregularity:growth:full}, we also obtain \eqref{eq:faster:convex-submono} and thus \cref{item:faster:ball-indicator:submono}.
\end{proof}

\begin{lemma}
    \label{lemma:subreg-1norm}
    Let $G \defeq \abs{\freevar}$ on $\R$ and $(\realoptx, \realoptx^*) \in \graph \subdiff G$. Then $\subdiff G$ is
    \begin{enumerate}[label=(\roman*)]
        \item\label{item:faster:1norm:subreg} metrically subregular at $\realoptx$ for $\realoptx^*$ for any $\kappa>0$ and
            \begin{equation*}
                \delta \le
                \begin{cases}
                    2\kappa & \text{if }\realoptx^* \in \{1, -1\}, \\
                    \kappa & \text{if }\abs{\realoptx^*} < 1;
                \end{cases}
            \end{equation*}
        \item\label{item:faster:1norm:submono} $(\gamma, \gamma)$-strongly submonotone at $\realoptx$ for $\realoptx^*$ for any $\gamma>0$ and
            \begin{equation*}
                \delta \le
                \begin{cases}
                    2\inv \gamma & \text{if }\realoptx^* \in \{1, -1\}, \\
                    \inv\gamma & \text{if } \abs{\realoptx^*} < 1.
                \end{cases}
            \end{equation*}
    \end{enumerate}
\end{lemma}

\begin{proof}
    We first verify \eqref{eq:stability:subregularity:growth:full} for any $\delta>0$ and $\gamma=\inv\kappa$ as stated.
    Suppose first that $\realoptx^*=1$ so that $\realoptx \in \inv{[\subdiff G]}(\realoptx^*)=[0, \infty)$.
    This implies that $\realoptx=\abs{\realoptx}$, and hence \eqref{eq:stability:subregularity:growth:full} becomes
    \begin{equation*}
        \abs{x} \ge x + \gamma \inf_{\optx \ge 0}(x-\optx)^2
        \quad (\abs{x-\realoptx} \le \delta).
    \end{equation*}
    If $x \ge 0$, this trivially holds by taking $\optx=x$. If $x \le 0$, the right-hand side is minimized by $\optx=0$, and thus the inequality holds for $x \ge -2\inv\gamma$. Since $\realoptx \ge 0$, this is guaranteed by our bound on $\delta$.
    The case $\realoptx^*=-1$ is analogous.

    If $\abs{\realoptx^*} < 1$, then $\realoptx \in \inv{[\subdiff G]}(\realoptx^*)=\{0\}$, and hence  \eqref{eq:stability:subregularity:growth:full} becomes
    \begin{equation*}
        \abs{x} \ge \gamma \abs{x}^2
        \quad (\abs{x} \le \delta).
    \end{equation*}
    This again holds by our choice of $\delta$.

    \Cref{cor:stability:subregularity:convex} now yields \cref{item:faster:1norm:subreg}. Adding
    \begin{equation*}
        G(\realoptx)-G(x) \ge \iprod{x^*}{\realoptx-x}
        \quad (x^* \in \subdiff G(x))
    \end{equation*}
    to \eqref{eq:stability:subregularity:growth:full}, we also obtain \eqref{eq:faster:convex-submono} and thus \cref{item:faster:1norm:submono}.
\end{proof}

\begin{remark}
    If we allow in the definition of subregularity or submonotonicity an arbitrary neighborhood of $\realoptx$ instead of a ball, then \cref{lemma:subreg-1norm} holds in a much larger neighborhood.
\end{remark}

\section{Local linear convergence of explicit splitting}
\label{sec:faster:acceleration}

Returning to the notation used in \crefrange{chap:proximal}{chap:meta}, we now assume throughout that $X$ is a Hilbert space, $F,G: X \to \Rbar$ are convex, proper, and lower semicontinuous, and that $F$ is Fréchet differentiable and has a Lipschitz continuous gradient $\grad F$ with Lipschitz constant $L\geq 0$.
Let further an initial iterate $x^0\in X$ and a step size $\tau >0$ be given and let the sequence $\{x^k\}_{k\in\N}$ be generated by the forward-backward splitting method (or basic proximal point method if $F=0$), i.e., by solving for $\nextx$ in
\begin{equation}
    \label{eq:faster:fb}
    0 \in \tau[\subdiff G(\nextx)+\grad F(\thisx)] + (\nextx-\thisx).
\end{equation}
We also write $H \defeq \subdiff G + \grad F:X\setto X$.
Finally, it is worth recalling the approach of \cref{chap:testing} for encoding the convergence rate into ``testing'' parameters $\tauTest_k>0$.

We start our analysis by adapting the proofs of \cref{thm:testing:fb,thm:gap:fb:value} to employ the squared distance function $x \mapsto \dist^2(x; \realopt X)$ to the entire solution set $\realopt X = \inv H(0)$ in place of the squared distance function $x \mapsto \norm{x-\realoptx}_X^2$ to a fixed $\realoptx \in \inv H(0)$.

\begin{lemma}
    \label{lemma:faster:fb}
    Let $\realopt X\subset X$. If
    for all $k \in \N$ and $\nexxt{w} \defeq -\grad F(\thisx)-\inv\tau(\nextx-\thisx) \in \subdiff G(\nextx)$,
    \begin{multline}
        \label{eq:faster:convergence-condition-sub}
        \inf_{\optx \in \realopt X}\left(
            \frac{\tauTest_k}{2}\norm{\nextx-\optx}_X^2
            + \tauTest_k\tau\iprod{\nexxt{w}+\grad F(\thisx)}{\nextx-\optx}_X
        \right)
        \\
        \ge
        \frac{\tauTest_{k+1}}{2}\dist^2(\nextx, \realopt X)
        - \frac{\tauTest_k}{2}\norm{\nextx-\thisx}_X^2,
    \end{multline}
    then
    \begin{equation}
        \label{eq:faster:convergence-result-sub}
        \frac{\tauTest_N}{2}\dist^2(x^N, \realopt X)
        \le
        \frac{\tauTest_0}{2}\dist^2(x^0, \realopt X)
        \quad
        (N \ge 1).
    \end{equation}
\end{lemma}

\begin{proof}
    Inserting \eqref{eq:faster:fb} into \eqref{eq:faster:convergence-condition-sub} yields
    \begin{multline}
        \label{eq:faster:convergence-condition-sub-transformed0}
        \inf_{\optx \in \inv H(0)}\tauTest_k\left(
            \frac{1}{2}\norm{\nextx-\thisx}_X^2
            + \frac{1}{2}\norm{\nextx-\optx}_X^2
            - \iprod{\nextx-\thisx}{\nextx-\optx}_X
        \right)
        \\
        \ge
        \frac{\tauTest_{k+1}}{2}\dist^2(\nextx; \inv H(0)).
    \end{multline}
    Using the three-point formula \eqref{eq:convergence:three-point-identity}, we can then rewrite \eqref{eq:faster:convergence-condition-sub-transformed0} as
    \begin{equation*}
        \frac{\tauTest_k}{2}\dist^2(\thisx; \inv H(0))
        \ge
        \frac{\tauTest_{k+1}}{2}\dist^2(\nextx; \inv H(0)).
    \end{equation*}
    The claim now follows by a telescoping sum over $k=0,\ldots,N-1$.
\end{proof}

\subsection*{Rates from error bounds and metric subregularity}
\label{sec:rates-from-error-bounds}

Our first approach for the satisfaction of \eqref{eq:faster:convergence-condition-sub} is based on \term[bounds, error]{error bounds}, which we will prove using metric subregularity.
The essence of error bounds is to prove for some $\theta >0$ that
\begin{equation*}
    \label{eq:faster:error-bound0}
    \norm{\nextx-\thisx}_X \ge \theta \norm{\nextx-\realoptx}_X.
\end{equation*}
We slightly weaken this condition, and assume the bound to be relative to the entire solution set, i.e.,
\begin{equation}
    \label{eq:faster:error-bound}
    \norm{\nextx-\thisx}_X^2 \ge \theta \dist^2(\nextx; \inv H(0)).
\end{equation}
This bound holds under metric subregularity. We first need the following technical lemma on the iteration \eqref{eq:faster:fb}.
\begin{lemma}
    \label{lemma:error-bound-first-estimate-fb}
    If $\tau>0$, we have
    \begin{equation*}
        \frac{1}{2}\norm{\nextx-\thisx}_X^2
        \ge
        \frac{\tau^2}{4(1+L^2\tau^2)} \dist^2(0, \subdiff G(\nextx)+\grad F(\nextx)).
    \end{equation*}
\end{lemma}

\begin{proof}
    Since $-(\nextx-\thisx) \in \tau[\subdiff G(\nextx)+\grad F(\thisx)]$ by \eqref{eq:faster:fb}, we have
    \begin{equation}
        \label{eq:faster:error-bound-first-estimate}
        \frac{1}{2}\norm{\nextx-\thisx}_X^2
        =
        \frac{1}{2}\dist^2(0, \{-(\nextx-\thisx)\})
        \ge
        \frac{1}{2}\dist^2(0, \tau[\subdiff G(\nextx)+\grad F(\thisx)]).
    \end{equation}
    The generalized Young's inequality for any $\alpha \in (0, 1)$ then yields
    \begin{multline*}
        \frac{1}{2}\dist^2(0, \tau[\subdiff G(\nextx)+\grad F(\thisx)])
        \\
        \begin{aligned}
            &
            =
            \frac{\tau^2}{2}\dist^2(\grad F(\nextx)-\grad F(\thisx), \subdiff G(\nextx)+\grad F(\nextx))
            \\
            &
            =
            \inf_{q \in \subdiff G(\nextx)}
            \frac{\tau^2}{2}\norm{(\grad F(\nextx)-\grad F(\thisx))-(q+\grad F(\nextx))}_{X}^2
            \\
            &
            \ge
            \frac{\tau^2(1-\inv\alpha)}{2}\norm{\grad F(\nextx)-\grad F(\thisx)}_{X}^2+\inf_{q \in \subdiff G(\nextx)} \frac{\tau^2(1-\alpha)}{2}\norm{q+\grad F(\nextx)}_{X}^2
            \\
            &
            \ge
            \frac{\tau^2(1-\inv\alpha)L^2}{2}\norm{\nextx-\thisx}_X^2 + \frac{\tau^2(1-\alpha)}{2}\dist^2(0, \subdiff G(\nextx)+\grad F(\nextx)),
        \end{aligned}
    \end{multline*}
    where we have used in the last step that $\alpha \in (0, 1)$ and that $\grad F$ is Lipschitz continuous.
    Combining this estimate with \eqref{eq:faster:error-bound-first-estimate}, we obtain that
    \begin{equation*}
        \frac{1+\tau^2(\inv\alpha-1)L^2}{2}\norm{\nextx-\thisx}_X^2 \ge \frac{\tau^2(1-\alpha)}{2}\dist^2(0, \subdiff G(\nextx)+\grad F(\nextx)).
    \end{equation*}
    Rearranging and using that $1+\tau^2(\inv\alpha-1)L^2>0$ by assumption then yields
    \begin{equation*}
        \frac{1}{2}\norm{\nextx-\thisx}_X^2 \ge \frac{\theta}{2}\dist^2(0, \subdiff G(\nextx)+\grad F(\nextx))
    \end{equation*}
    for
    \begin{equation*}
        \theta
        \defeq \frac{\tau^2(1-\alpha)}{1+\tau^2(\inv\alpha-1)L^2},
    \end{equation*}
    which for $\alpha=1/2$ yields the claim.
\end{proof}

Metric subregularity then immediately yields the error bound \eqref{eq:faster:error-bound}.

\begin{lemma}
    \label{lemma:subregularity-to-peb}
    Let $H$ be metrically subregular at $\realoptx$ for $\realoptw=0$ for $\kappa>0$ and $\delta>0$.
    If $\tau>0$
    and $\nextx \in \B(\realoptx, \delta)$, then \eqref{eq:faster:error-bound} holds with $\theta=\frac{\tau^2}{2\kappa^2(1+L^2\tau^2)}$.
\end{lemma}

\begin{proof}
    Combining \cref{lemma:error-bound-first-estimate-fb} and the definition of metric subregularity yields
    \begin{equation*}
        \frac{1}{2}\norm{\nextx-\thisx}_X^2
        \ge
        \frac{\tau^2}{4(1+L^2\tau^2)} \dist^2(0, H(\nextx))
        \ge
        \frac{\tau^2}{4\kappa^2(1+L^2\tau^2)} \dist^2(\nextx, \inv H(0)).
        \qedhere
    \end{equation*}
\end{proof}

From this lemma, we now obtain \emph{local} linear convergence of the forward-backward splitting method when $H$ is metrically subregular at a solution.

\begin{theorem}
    \label{thm:faster:convergence-result-sub-peb}
    Let $H$ be metrically subregular at $\realoptx\in\inv H(0)$ for $\realoptw=0$ for $\kappa>0$ and $\delta>0$.
    If $0 < \tau L \le 2$ and $\thisx \in \B(\realoptx, \delta)$ for all $k\in\N$,
    then \eqref{eq:faster:convergence-result-sub} holds for $\tauTest_{k+1} \defeq \tauTest_k(1+\rho)$ and $\tauTest_0=1$ with $\rho=\theta(1-\frac{L\tau}{2})>0$ for $\theta=\frac{\tau^2}{2\kappa^2(1+L^2\tau^2)}$.
    In particular, $\dist^2(x^N; \inv H(0)) \to 0$ at a linear rate.
\end{theorem}

\begin{proof}
    Let $\optx \in \inv H(0)$ and $\nexxt{w} \in \subdiff G(\nextx)$ as in \cref{lemma:faster:fb}.
    From \cref{eq:testing:fb:est1} in the proof of \cref{thm:testing:fb}, we obtain
    \begin{equation*}
        \iprod{\nexxt{w}+\grad F(\thisx)}{\nextx-\optx}_X
        \ge -\frac{L}{4}\norm{\nextx-\thisx}_X^2.
    \end{equation*}
    \Cref{lemma:subregularity-to-peb} now yields the error bound \eqref{eq:faster:error-bound}, i.e., $\norm{\nextx-\thisx}_X^2 \ge \theta\dist^2(\nextx; \inv H(0))$.
    Hence, keeping in mind that $2>L\tau$, for all $\optx \in \inv H(0)$ we have that
    \begin{equation*}
        \begin{aligned}
            \frac{\tauTest_k(1-\tfrac{L\tau}{2})}{2}\norm{\nextx-\thisx}_X^2
        &
        + \frac{\tauTest_k}{2}\norm{\nextx-\optx}_X^2
        \\
        &
        \ge
        \frac{\tauTest_k\theta(1-\tfrac{L\tau}{2})}{2}\dist^2(\nextx; \inv H(0))
        + \frac{\tauTest_k}{2}\norm{\nextx-\optx}_X^2
        \\
        &
        \ge
        \frac{\tauTest_{k+1}}{2}\dist^2(\nextx; \inv H(0)).
        \end{aligned}
    \end{equation*}
    Adding the first estimate multiplied by $\tauTest_k\tau$ to the latter estimate yields
    \begin{multline*}
        \frac{\tauTest_k}{2}\norm{\nextx-\thisx}_X^2
        +\frac{\tauTest_k}{2}\norm{\nextx-\optx}_X^2
        + \tauTest_k\tau\iprod{\nexxt{w}+\grad F(\thisx)}{\nextx-\optx}_X
        \\
        \ge \frac{\tauTest_{k+1}}{2}\dist^2(\nextx, \inv H(0)).
    \end{multline*}
    Taking the infimum over $\optx \in \inv H(0)$, we obtain \eqref{eq:faster:convergence-condition-sub} for $\realopt X=\inv H(0)$. The claim now follows from \cref{lemma:faster:fb} and the exponential growth of $\tauTest_k$.
\end{proof}

The convergence is local due to the requirement $\thisx \in \B(\realoptx, \delta)$ for applying subregularity. In finite dimensions, the weak convergence result of \cref{thm:convergence:fb} of course guarantees that the iterates enter and remain in this neighborhood after a finite number of steps.

\subsection*{Rates from strong submonotonicity}
\label{sec:rates-from-monotonicity}

If $H$ is instead strongly submonotone, we can (locally) ensure \eqref{eq:faster:convergence-condition-sub} directly.

\begin{theorem}
    \label{thm:faster:convergence-result-submonotone}
    Let $H$ be $(\gamma/2, \theta/2)$-strongly submonotone at $\realoptx\in \inv H(0)$ for $\realoptw=0$ for $\delta>0$. If $\gamma > \theta+L^2\tau$ and $x^0 \in \B(\realoptx, \epsilon)$ for some $\epsilon>0$ sufficiently small, then \eqref{eq:faster:convergence-result-sub} holds for $\tauTest_{k+1} \defeq \tauTest_k(1+(\gamma-L^2\tau)\tau)$ and $\tauTest_0=1$.
    In particular, $\dist^2(x^N; \inv H(0)) \to 0$ at a linear rate.
\end{theorem}

\begin{proof}
    Let $\nexxt{w}\defeq-\inv\tau(\nextx-\thisx)-\grad F(\thisx) \in \subdiff G(\nextx)$ by \eqref{eq:faster:fb}.
    By \eqref{eq:convergence:fb-est} in the proof of \cref{thm:convergence:fb}, if $x^0 \in \B(\realoptx, \epsilon)$ for $\epsilon>0$ small enough, then $\norm{\nextx-\thisx}_X \le \delta/(L+\inv\tau)$ for all $k\in \N$ such that the Lipschitz continuity of $\grad F$ yields
    \begin{equation*}
        \norm{\grad F(\nextx)-\grad F(\thisx)-\inv\tau(\nextx-\thisx)}_{X} \le \delta.
    \end{equation*}
    Thus $\nexxt{w} \in \subdiff G(\nextx) \isect \B(-\grad F(\nextx), \delta)$ and $\nextx \in \B(\realoptx, \delta)$  for all $k \in \N$.
    Now, for all $\optx \in \inv H(0)$, the strong submonotonicity of $H$ at $\realoptx$ for $0$ implies that
    \begin{equation*}
        \tauTest_k\tau \iprod{\nexxt{w}+\grad F(\nextx)}{\nextx-\optx}_X
        + \frac{(\theta-\gamma)\tauTest_k\tau}{2}\norm{\nextx-\optx}_X^2
        \ge
        \frac{\theta\tauTest_k\tau}{2}\dist^2(\nextx; \inv H(0))
    \end{equation*}
    for all $k \in \N$.
    Cauchy's inequality and the Lipschitz continuity of $\grad F$ then yields
    \begin{equation*}
        \tauTest_k\tau\iprod{\grad F(\thisx)-\grad F(\nextx)}{\nextx-\optx}_X
        \ge -\frac{\tauTest_k}{2}\norm{\nextx-\thisx}_X^2
        -\frac{\tauTest_k\tau^2 L^2}{2}\norm{\nextx-\optx}_X^2.
    \end{equation*}
    We now sum the last two inequalities to obtain
    \begin{multline*}
        \frac{\tauTest_k}{2}\norm{\nextx-\thisx}_X^2
        + \tauTest_k\tau \iprod{\nexxt{w}+\grad F(\thisx)}{\nextx-\optx}_X
        \\
        \ge
        \frac{\theta\tauTest_k\tau}{2}\dist^2(\nextx; \inv H(0))
        +
        \frac{(\gamma-\theta-L^2\tau)\tauTest_k\tau}{2}\norm{\nextx-\optx}_X^2.
    \end{multline*}
    Using that $\theta-\gamma+L^2\tau<0$ and taking the infimum over all $\optx \in \inv H(0)$ then yields
    \begin{multline*}
        \inf_{\optx \in \inv H(0)}\left(
            \frac{\tauTest_k}{2}\norm{\nextx-\optx}_X^2
            + \tauTest_k\tau \iprod{\nexxt{w}+\grad F(\thisx)}{\nextx-\optx}_X
        \right)
        \\
        \ge
        \frac{(\gamma-L^2\tau)\tauTest_k\tau+\tauTest_k}{2}\dist^2(\nextx; \inv H(0))
        -\frac{\tauTest_k}{2}\norm{\nextx-\thisx}_X^2.
    \end{multline*}
    Since $\gamma-L^2\tau>0$ and $\tauTest_{k+1} = \tauTest_k(1+(\gamma-L^2\tau)\tau)$, this shows \eqref{eq:faster:convergence-condition-sub} with $\realopt X=\inv H(0)$.
    The claim now follows from \cref{lemma:faster:fb} and the exponential growth of $\tauTest_k$.
\end{proof}

\begin{remark}
    Similarly to \cref{thm:testing:prox}\,\ref{item:testing:prox:superlinear}, if $F\equiv 0$ we can let $\tau \upto \infty$ to obtain local superlinear convergence of the proximal point method under strong submonotonicity of $\partial G$ at the solution.
\end{remark}

\begin{remark}[local linear convergence]
    Local linear convergence was first derived from error bounds in \cite{luo1992error} for matrix splitting problems and was studied for other methods, including the ADMM and the proximal point method among others, in \cite{han2013local,aspelmeier2016local,leventhal2009metric,li2012holder}. An alternative approach to the proximal point method was taken in \cite{aragon2012lyusternik} based on Lyusternik--Graves-style estimates, while \cite{adly2015newton} presented an approach based on metric regularity to Newton's method for variational inclusions.
    Furthermore, \cite{zhou2017unified} proposed a unified approach to error bounds for generic smooth constrained problems.
    Finally, \cite{liu2018partial,tuomov-subreg} introduced \emph{partial} or subspace versions of error bounds and showed the fast convergence of only some variables of structured algorithms such as the ADMM or PDPS.
    The relationships between error bounds and metric subregularity are studied in more detail in \cite{gfrerer2011first,ioffe2017variational,kruger2015error,dontchev2014implicit,ngai2008error}.
    Submonotonicity was introduced in \cite{tuomov-subreg}.
\end{remark}

\part{Applications}\label{part:applications}

\chapter{Sparse regularization}
\label{chap:sparse}

In this and the following chapters, we illustrate the application of the results and methods of the previous parts to selected nonsmooth optimization problems.

We first study the application of the optimization theory and methods that we have developed to the solution of some \term[problem!inverse]{inverse problems}, including imaging problems, which we treat in finite dimensions to avoid technical difficulties unrelated to nonsmooth optimization.
In a nutshell, inverse problems consist in trying to obtain quantities of interest that are not directly accessible by combining measured (incomplete, noisy) data with a mathematical model linking the desired quantity to the predicted measurements. Such problems are usually \term[problem!ill-posed]{ill-posed} in the sense that a solution may not exist, may not be unique, or may not be stable with respect to perturbations of the data. Hence one needs to apply \term{regularization} to obtain a stable approximation. For an introduction to the regularization of inverse problems, we refer the reader to the seminal work \cite{Engl} as well as to the more recent \cite{ItoJin,ClasonIP:2020}.
One particular approach is \term[regularization!Tikhonov]{Tikhonov regularization}, which consists in solving an optimization problem that involves the sum of (a) a \term[term!data]{data term} that matches the model prediction against available data and of (b) a \term[term!regularization]{regularization term} that attempts to promote expected and desirable features in the reconstruction (and is typically required to obtain well-posedness of the regularized problem). An increasingly popular class of regularization terms promotes ``sparsity'' of the solution in the sense that it can explain the data with a minimal number of features; as we will see, such terms require nonsmooth optimization.
This class (and nonsmooth optimization in general) is particularly relevant in the context of \term[problem!image processing]{mathematical image processing}, where the quantity of interest is an image rather than an abstract physical parameter; see, e.g., \cite{Scherzer,BrediesLorenz:2018}.

In this chapter we start with perhaps the simplest nonsmooth regularization of an inverse problem: $\ell^1$-regularized data-fitting, sometimes known as the \term[problem!Lasso]{Lasso problem}.
The starting point is linear regression, but we wish to explain the data ``in simple terms'' only through its most important features.
We then move on to signal and image recovery applications in the next \cref{chap:l1fit,chap:tv}.

\section{Problem description}

Let $b_i \in \R$ be a single measurement of an unknown signal $x \in \R^M$ through the filter $a_i \in \R^M$. Without the presence of noise, $b_i=a_i^T x$ for the $i=1,\ldots,N$ measurements.
In statistical contexts, $b_i$ is known as a dependent variable and $a_i$ as a data vector.
Since each $b_i$ and $a_i$ may be noisy, and the system
\begin{equation*}
    a_i^T x=b_i \quad (i=1,\ldots,N),
\end{equation*}
may be over- or under-determined, direct solution of $x$ from this system is not in general well-posed.
Basic linear regression instead seeks the least squares solution $x$ through solution of the optimization problem
\begin{equation}
    \label{eq:sparse:regression}
    \min_{x \in \R^M}~ \frac{1}{N}\sum_{i=1}^N \frac{1}{2}(b_i - a_i^T x)^2.
\end{equation}
To explain the data $\{(a_i, b_i)\}$ through its most important features, we want $x$ to be sparse, i.e., to have many zero elements, and few nonzero elements.
For example, $a_i$ might be the attributes (genre, length, etc.) of a film, and $b_i$ its rating. A sparse vector $x$ would then contain only the most relevant attributes for the rating and their relative weighting.
To perform such \term[problem!sparse regression]{sparse regression}, let us add to the data fitting term of \eqref{eq:sparse:regression} the regularization term $g(x)=\lambda\norm{x}_1$.
Then we obtain the so-called \term[problem!Lasso]{Lasso problem}
\begin{equation}
    \label{eq:sparse:problem-motiv}
    \min_{x \in \R^M}~ \frac{1}{N}\sum_{i=1}^N \frac{1}{2}(b_i - a_i^T x)^2 + \lambda\norm{x}_1.
\end{equation}
The hope is that to explain the data, the $\ell^1$-norm regularization term will cause the minimizer to select more relevant features from the data, ignoring irrelevant ones.

In the following, we write \eqref{eq:sparse:problem-motiv} more succinctly as
\begin{equation}
    \label{eq:sparse:problem}
    \min_{x \in \R^M}~ J(x) \quad\text{for}\quad J(x) \defeq F(x) + G(x),
\end{equation}
where
\[
    A \defeq (a_1, \ldots, a_N)^T \in \R^{N \times M},
    \quad
    F(x) \defeq \frac{1}{2}\norm{Ax-b}_2^2,
    \quad\text{and}\quad
    G(y) \defeq \lambda\norm{x}_1.
\]

\section{Optimality conditions}

Our first result characterizes the solutions of \eqref{eq:sparse:problem}.
\begin{theorem}\label{thm:sparse:optimal}
    The vector $\realoptx \in \R^M$ is a solution to \eqref{eq:sparse:problem-motiv} if and only if
 there exists a $\realopt p \in \R^{M}$ such that
    \begin{equation}
        \label{eq:sparse:primal-oc}
        -A^*(A\realoptx-b) = \lambda \realopt p
        \qquad\text{ and }\qquad
        \realopt p_i \in
        \begin{cases}
            \{1\} & \text{if }\realoptx_i>0,\\
            \{{-}1\} & \text{if }\realoptx_i<0,\\
            [-1,1] & \text{if }\realoptx_i=0.
        \end{cases}
    \end{equation}
\end{theorem}
\begin{proof}
    Since $A$ is linear and $F$ and $G$ are convex, $J$ is convex as well.
    Therefore the convex Fermat principle of \cref{thm:convex:fermat} is an equivalent characterization of  solutions to \eqref{eq:sparse:problem-motiv} as those $\realoptx$ satisfying $0 \in \subdiff J(\realoptx)$.
    Since both $F$ and $G$ have full domain and are proper and lower semicontinuous, we may further use the subdifferential sum rule of \cref{thm:subdiff:sum} to deduce for all $x \in \R^M$ that $\subdiff J(x) =\subdiff F(x) + \subdiff G(x)$.
    Since $F$ is differentiable, using \cref{thm:convex:gateaux} we therefore characterize the solutions as those points $\realoptx$ satisfying
    \begin{equation}
        \label{eq:sparse:primal-oc:0}
        -\grad F(\realoptx) \in \subdiff G(\realoptx).
    \end{equation}
    Since $F$ is smooth, expanding $\grad F(\realoptx)=A^*(A\realoptx-b)$ and using \cref{ex:convex:subdiff_abs} to calculate $\subdiff G(\realoptx)$ componentwise yields \eqref{eq:sparse:primal-oc}.
\end{proof}

Note the complementarity between the primal variable $\realoptx$ and the dual variable $\realopt p$, which yields the desired sparsity: a component $\realoptx_i$ is zero if the corresponding scaled and \enquote{back-propagated} residual ${\realopt p}_i$ is smaller than $1$ in magnitude. However, $\realoptx_i$ can be zero even if $|\realopt p_i|=1$; if this case can be excluded, we say that \term[complementarity, strict]{strict complementarity} holds, i.e.,
\begin{equation}
    \label{eq:apple_inverse:lasso:sc}
    \text{either } \realoptx_i \ne 0 \text{ or } \abs{\realopt p_i} < 1
    \quad (i=1,\ldots,M).
\end{equation}
Thus strict complementarity avoids, whenever $\realoptx_i=0$, the boundary cases $\abs{\realopt p_i}=1$ that happen when $\realoptx \ne 0$.

\section{Algorithms}

The starting point for deriving implementable algorithms for the solution of \eqref{eq:sparse:problem-motiv} is the following reformulation of the optimality conditions using the proximal point mapping.
\begin{lemma}
    The vector $\realoptx \in \R^M$ is a solution to \eqref{eq:sparse:problem-motiv} if and only if
    \begin{equation}
        \label{eq:sparse:prox-oc}
        \realoptx = \prox_{\tau G}(\realoptx -\tau A^* (A\realoptx-b)).
    \end{equation}
\end{lemma}
\begin{proof}
    Applying \cref{lem:proximal:subdiff} to $G$, we may rewrite \eqref{eq:sparse:primal-oc:0} for any $\tau>0$ as
    \begin{equation*}
        \realoptx = \prox_{\tau G}(\realoptx -\tau \grad F(\realoptx)),
    \end{equation*}
    which after inserting $\grad F(\realoptx)=A^*(A\realoptx-b)$ yields \eqref{eq:sparse:prox-oc}.
\end{proof}

\subsection*{Forward-backward splitting}

The forward-backward or explicit splitting method of \eqref{eq:splitting:fb} is our first iterative method for solving \eqref{eq:sparse:problem-motiv}.
As we did in the general setting in \cref{chap:proximal}, the method can be directly developed from  the proximal-form optimality conditions \eqref{eq:sparse:prox-oc}.
First, using \cref{ex:proximal:reell}\,\cref{ex:proximal:reell:ii} we write the proximal point mapping of $G$ in terms of the soft-thresholding operator as
\begin{equation*}
    \prox_{\tau G}(x) = \left(\soft_{\lambda\tau}(x_1), \ldots, \soft_{\lambda\tau}(x_M)\right)
    \quad\text{for}\quad
    \soft_\theta(t) \defeq
    \begin{cases} t-\theta & \text{if }t>\theta,\\
        0 &\text{if } t\in [-\theta,\theta],\\
        t+\theta &\text{if } t < -\theta.
    \end{cases}
\end{equation*}
Inserting this together with $\grad F(x) = A^*(Ax+b)$ into \eqref{eq:splitting:fb} then yields the \term[method!iterative soft-thresholding]{iterative soft-thresholding algorithm} (ISTA)
\begin{algeqbox}
    \begin{equation}
        \label{eq:sparse:fb}
        \nextx
        \defeq
        \soft_{\lambda\tau}((\Id-\tau A^*A)\thisx+\tau A^*b).
    \end{equation}
\end{algeqbox}

Under mild conditions, the iterates converge.

\begin{theorem}
    \label{thm:sparse:fb-convergence}
    Suppose $\tau \norm{A}^2<2$.
    Then for any starting point $x^0 \in \R^M$, the iterates $\{\thisx\}_{k \in \N}$ generated by \eqref{eq:sparse:fb} converge to a solution $\realoptx$ of \eqref{eq:sparse:problem}.
\end{theorem}

\begin{proof}
    The Lipschitz factor of $\grad F(x)=A^*(Ax-b)$ is $\norm{A}^2$. Therefore, the claim follows from \cref{thm:convergence:fb}.
\end{proof}

Convergence of function values can be similarly deduced from \cref{thm:gap:fb:value}.

\Cref{thm:sparse:fb-convergence} provides no convergence rates as, indeed, no rates for iterates are in general known for forward-backward splitting without some sort of stronger growth assumptions. However, \cref{theorem:regularity:lasso:reg} in \cref{sec:sparse:stability} below will show that $\subdiff[F+G]$ is metrically regular at $\realoptx$ for $0$, provided that the strict complementarity condition \eqref{eq:apple_inverse:lasso:sc} holds.
Since this implies metric subregularity, \cref{thm:faster:convergence-result-sub-peb} can be used to demonstrate the local linear convergence of \eqref{eq:sparse:fb} near a strictly complementary solution.

We can also apply the inertial explicit splitting method of \eqref{eq:meta:inertia:fista} to \eqref{eq:sparse:problem}.
Based on the basic explicit splitting \eqref{eq:sparse:fb}, this method becomes the FISTA\index{method!iterative soft-thresholding!fast}
\begin{algeqbox*}
    \begin{equation*}
        \label{eq:sparse:fista}
        \left\{
            \begin{aligned}
                \nextx
                &
                \defeq
                \soft_{\lambda\tau}((\Id-\tau A^*A)\this{\bar x}+\tau A^*b),
                \\
                \alpha_{k+1} &= \lambda_{k+1}(\inv\lambda_k-1),
                \\
                \nexxt{\bar x} & = (1+\alpha_{k+1})\nextx-\alpha_{k+1}\thisx.
            \end{aligned}
        \right.
    \end{equation*}
\end{algeqbox*}
The initial inertial parameter $\lambda_0=1$, while $\bar x^0 \in \R^M$ can be chosen freely.
Since $\alpha_1=0$, $x^0$ is never used. Regarding convergence, \cref{thm:meta:inertia:fb} readily gives the following result.

\begin{theorem}
    Suppose $\tau \norm{A}^2 \le 1$.
    Then for any starting point $x^0 \in \R^M$, the iterates $\{\thisx\}_{k \in \N}$ generated by \eqref{eq:sparse:fista} satisfy $J(\thisx) \to \min J$ at the rate $O(1/k^2)$.
\end{theorem}

In fact, under strict complementarity, zeros are identified in a finite number of steps.

\begin{theorem}
    \label{thm:sparse:identification}
    Assume the solution $\realoptx \in \R^M$ to \eqref{eq:sparse:problem-motiv} is unique and satisfies strict complementarity.
    Then there exists $K \in \N$ such that the iterates $\{\thisx\}_{k \in \N}$ generated by \eqref{eq:sparse:fista} satisfy $x^k_i = \realoptx_i$ for all $k \ge K$ and $i=1,\ldots,N$ with $\realoptx_i=0$.
\end{theorem}

\begin{proof}
    Indeed, forward-backward splitting for $\min_x (F+G)$ and a step length $\tau>0$ by definition satisfies
    \begin{equation}
        \label{eq:appl_inverse:lasso:fb-indentification}
        0 \in \subdiff G(\nextx)+ \grad F(\thisx)+\tau(\nextx-\thisx).
    \end{equation}
    Following the proof of \cref{thm:convergence:fb}, we have $\norm{\nextx-\thisx} \to 0$ and $\nextx \to \realoptx$ provided the solution $\realoptx$ is unique.
    It follows that $\grad F(\thisx) \to \grad F(\realoptx)$.
    Furthermore, strict complementarity yields $-[\grad F(\realoptx)]_i \in (-\lambda, \lambda)$ for all $i$ with $\realoptx_i=0$, and hence for those same $i$ it holds that $-[\grad F(\thisx)]_i \in (-\lambda, \lambda)$ for all $k \ge K$ for some $K \in \N$.
    By \eqref{eq:appl_inverse:lasso:fb-indentification} and $\norm{\nextx-\thisx} \to 0$, it is then necessary that $[\subdiff G(\nextx)]_i$ contains a point in $(-\lambda, \lambda)$. This is only possible if $\nextx_i=0$ after a finite number of steps.
\end{proof}

\begin{remark}[unconditional linear convergence and activity identification]
    It is shown in \cite{bolte2017errorbounds} through error bounds that forward-backward splitting for the Lasso problem converges linearly without any assumptions.
    Error bounds can also be proved more generally based on piecewise polynomial properties derived in \cite{li2013global}. These are used in \cite{garrigos2017thresholding} to obtain error bounds in separable Hilbert spaces.
    This follows earlier works such as \cite{bredies2008linear} with stricter assumptions.
    In the former also a ``finite identification property'' is studied following earlier efforts in \cite{lewis2002active,liang2014local}, among others; it can be shown that the forward-backward splitting and other methods converge in a finite number of steps to a smooth submanifold. As verified by elementary analysis in \cref{thm:sparse:identification}, in the case of the Lasso problem, forward-backward splitting identifies the strictly complementary zeros in a finite number of steps.
\end{remark}

\subsection*{Semismooth Newton method}

By \cref{thm:sparse:optimal}, we know that minimizers $\bar x$ of \eqref{eq:sparse:problem} satisfy
\begin{equation*}
    \bar x - \prox_{\gamma G}(\bar x - \gamma \nabla F(\bar x))=0
\end{equation*}
for any $\gamma >0$. We therefore look for a root of
\[
    H(x) \defeq x - \prox_{\gamma G}(x - \gamma \nabla F(x)).
\]
If we can produce an invertible and well-conditioned Newton derivative $D_N H(x)$ for all $x$ in a sufficiently large neighborhood of $\bar x$, finding a root can be done with the semismooth Newton method \eqref{eq:SSN}, i.e., solving $s^k$ from $D_N H(x^k) s_k = - H(x^k)$ and updating $x^{k+1} = x^k + s^k$.

In fact, let $T(x) \defeq x - \gamma \grad F(x)$ and consider the composition $\prox_{\gamma G} \circ T$. We may use \cref{thm:newton:frechet} to obtain $D_N T(x) = \Id - \gamma \nabla^2 F(x)$, and \cref{ex:newton:rn}\,\cref{ex:newton:rn:l1} to obtain $D_N\prox_{\gamma G}$.
Both $D_N \prox_{\gamma G}$ and $D_N T$ are locally uniformly bounded (obviously from the characterization and the continuous differentiability, respectively).
Thus, we are justified in using the chain rule from \cref{thm:newton:chain} on the composition to calculate
\[
    \begin{aligned}
    D_N H(x)
    &
    = \Id - D_N\prox_{\gamma G}(T(x)) \circ D_N T(x)
    \\
    &
    = \Id - \1_{\calA(x)} [\Id - \gamma \grad^2 F(x)]
    \\
    &
    = \1_{\calI(x)} + \gamma \1_{\calA(x)} \grad^2 F(x),
    \end{aligned}
\]
where we have defined the \term[set!inactive]{inactive} and \term[set!active]{active sets}, respectively, as
\begin{equation}
    \label{eq:sparse:ssn:inactive-active}
    \calI(x) \defeq \setof{i\in\{1,\dots,N\}}{|x_i - \gamma [\nabla F(x)]_i|<\gamma},\qquad
    \calA(x) \defeq \{1,\dots,M\}\setminus \calI(x),
\end{equation}
and -- in a slight abuse of notation -- we use $\1_{A}$ to denote a diagonal matrix with $[\1_{A}]_{ii}=1$ if $i\in A\subset\{1,\dots,N\}$ and $0$ otherwise.
The matrix $D_N H(x)$ may in general not be invertible, or may be poorly conditioned \emph{on the active components}, as we will soon see in more detail. For some $\theta>0$, we therefore replace it with the \term[dampening!of active components]{active-dampened} matrix
\begin{equation}
    \label{eq:sparse:ssn:m}
    M(x)
    \defeq
    \1_{\calI(x)} + \gamma  \1_{\calA(x)}\nabla^2F(x) + \theta \1_{\calA(x)}.
\end{equation}
Write $P_{\calI(x)}$ and $P_{\calA(x)}$ for the projections to the inactive and active components, so that $\1_{\calI(x)} = P_{\calI(x)}^*P_{\calI(x)}$, and likewise for the active components.
Thus the active-dampened semismooth Newton step $s^k$ is determined by
\begin{equation}
    \label{eq:sparse:ssn:first-system}
    \left(\1_{\calI(x^k)} + \gamma  \1_{\calA(x^k)}\nabla^2F(x^k) + \theta \1_{\calA(x^k)}\right)s^k = - x^k + \prox_{\gamma G}(x^k-\gamma \nabla F(x^k)).
\end{equation}
Since the proximal point mapping of $G$ is the soft shrinkage operator, we have using the definition of the inactive set that
\[
    P_{\calI(x^k)}\prox_{\gamma G}(x^k-\gamma \nabla F(x^k))=0.
\]
Hence, multiplying \eqref{eq:sparse:ssn:first-system} from the left by $P_{\calI(x^k)}$, we deduce that $P_{\calI(x^k)}s^k = - P_{\calI(x^k)}x^k$.
Thus $s^k_i=-x^k_i$ for the inactive components $i \in \calI(x^k)$.
It follows that $x^{k+1}_i = 0$ for $i \in \calI(x^k)$.
On the other hand, writing $s^k = \1_{\calA(x^k)}s^k + \1_{\calI(x^k)}s^k = P_{\calA(x^k)}^*P_{\calA(x^k)}s^k - \1_{\calI(x^k)}x^k$ and multiplying \eqref{eq:sparse:ssn:first-system} from the left by $P_{\calA(x^k)}$ yields
\begin{multline}
    \label{eq:sparse:ssn:active-system:approx}
    [\gamma P_{\calA(x^k)}\grad^2F(x^k) P_{\calA(x^k)}^* + \theta\Id] P_{\calA(x^k)}s^k
    \\
    = P_{\calA(x^k)}( - x^k + \prox_{\gamma G}(x^k-\gamma \grad F(x^k)) + \gamma \grad^2F(x^k)\1_{\calI(x^k)}x^k).
\end{multline}
Since $P_{\calA(x^k)}\grad^2F(x^k) P_{\calA(x^k)}^* + \theta\Id$ is positive definite, we can solve this for $P_{\calA(x^k)}s^k$.
Altogether, therefore, the semismooth Newton method for \eqref{eq:sparse:problem} becomes
\begin{algenumbox}
    \begin{enumerate}[label=\arabic*.]
        \item\label{item:sparse:ssn:first-step} form the inactive and active sets $\calI(x^k)$ and $\calA(x^k)$ following \eqref{eq:sparse:ssn:inactive-active};
        \item solve $P_{\calA(x^k)}s^k$ from \eqref{eq:sparse:ssn:active-system:approx};
        \item\label{item:sparse:ssn:last-step} update $x^{k+1} \defeq \1_{\calA(x^k)}(x^k + s^k)$.
    \end{enumerate}
\end{algenumbox}
This coincides with an \term[strategy!active set]{active set strategy} similar to those used for solving quadratic subproblems in sequential programming methods with inequality constraints; cf.~\cite[Chapter 8.4]{Kunisch:2008a}.

For convergence, we need to assume that $P_{\calA(\opt x)}\grad^2F(\optx) P_{\calA(\optx)}^*$ is invertible.
Practically this means that there are more measurements than attributes that describe the measurements.
Although superlinear convergence has superficially no stricter conditions than linear convergence, the convergence radius can in practice be smaller, and hence convergence may not hold for an arbitrary initial iterate $x^0$.

To improve readability of the next theorem proving these properties, we recall the following “operator Young's inequality”.

\begin{lemma}
    \label{lemma:sparse:matrix-young}
    On Hilbert spaces $X$ and $Y$, let $A \in \linear(X; Y)$ and $B \in \linear(X; Y)$. Then for any $\beta > 0$, we have
    \[
        2 A^* B \preceq \beta A^*A + \inv\beta B^*B,
    \]
    where $A \preceq B$ means that $B-A$ is positive semi-definite.
\end{lemma}

\begin{proof}
    Take any $x \in X$. Then using the Cauchy--Schwarz and Young's inequalities yields
    \[
        2\iprod{x}{A^*B x}_X
        =2\iprod{Ax}{Bx}_Y
        \le \beta \norm{Ax}_Y^2 + \inv\beta \norm{Bx}_Y^2
        = \iprod{x}{(\beta A^*A + \inv\beta B^*B)x}_X.
    \]
    Since this holds for all $x \in X$, this means that
    $(\beta A^*A + \inv\beta B^*B) - 2A^*B$ is positive semi-definite.
\end{proof}

\begin{theorem}
    \label{thm:sparse:ssn}
    Let $\opt x$ be a (unique) minimizer of \eqref{eq:sparse:problem}.
    Let $\gamma,\theta>0$ satisfy $2(1-\theta^2) > \gamma \theta \norm{A}^2$, and suppose that $P_{\calA(\opt x)}A^*AP_{\calA(\opt x)}^*$ is positive definite.
    If $x^0$ is sufficiently close to $\opt x$, then the sequence $\{x^{k+1}\}_{k \in \N}$ generated by iterating \crefrange{item:sparse:ssn:first-step}{item:sparse:ssn:last-step} above converges linearly to $\opt x$.
    If $\theta=0$, and $\gamma>0$ is arbitrary, the convergence is superlinear.
\end{theorem}

\begin{proof}
    We first consider linear convergence.
    Let $M(x)$ be given by \eqref{eq:sparse:ssn:m}.
    Then
    \[
        \norm{M(x)-D_NH(x)}_{\linear(\R^N;\R^N)} = \norm{\theta\1_{\calA(x)}}_{\linear(\R^N;\R^N)} \le \theta,
    \]
    so the corresponding assumption of \cref{thm:newton:dampened-linear} (applied to $H$ in place of  $F$) holds.
    To apply the theorem, we still need to prove $\norm{\inv{M(x)}}_{\linear(\R^N; \R^N)} \le C$ for all $x \in U$ for some neighborhood $U$ of $\opt x$ and some $C>0$ with $C\theta < 1$. That is to say, $M(x)^*M(x) \succeq C^{-2}\,\Id$.
    We expand
    \begin{equation}
        \label{eq:sparse:ssn:mm}
        \begin{aligned}[t]
            M(x)^*M(x)
            &
            =
            \1_{\calI(x)}
            + \gamma^2 A^*A \1_{\calA(x)} A^*A
            + \theta \gamma \1_{\calA(x)} A^*A
            + \theta \gamma A^*A \1_{\calA(x)}
            + \theta^2 \1_{\calA(x)}
            \\
            &
            =
            \1_{\calI(x)}
            + \gamma^2 A^*A \1_{\calA(x)} A^*A
            + 2 \theta \gamma \1_{\calA(x)} A^*A \1_{\calA(x)}
            + \theta^2 \1_{\calA(x)}
            \\
            \MoveEqLeft[-1]
            + \theta \gamma \1_{\calA(x)} A^*A \1_{\calI(x)}
            + \theta \gamma \1_{\calI(x)} A^*A \1_{\calA(x)}.
        \end{aligned}
    \end{equation}
    Eliminating the second term by positive semi-definiteness, and applying \cref{lemma:sparse:matrix-young} to the last two terms yields for any $\beta>0$ that
    \begin{equation}
        \label{eq:sparse:ssn:mm2}
            M(x)^*M(x)
            \succeq
            \1_{\calI(x)}
            + (2-\beta)\theta\gamma \1_{\calA(x)} A^*A \1_{\calA(x)}
            + \theta^2 \1_{\calA(x)}
            - \inv\beta \theta \gamma  \1_{\calI(x)} A^*A \1_{\calI(x)}.
    \end{equation}
    By assumption, $P_{\calA(x)}A^*AP_{\calA(x)}^*$ is positive definite and $\calA(x)=\calA(\opt x)$ for all $x$ in some open neighborhood $U$ of $\opt x$.
    Therefore $\1_{\calA(x)} A^*A \1_{\calA(x)} \succeq \epsilon  \1_{\calA(x)}$ for some $\epsilon>0$ and all $x \in U$.
    Consequently, it follows from \eqref{eq:sparse:ssn:mm2} that
    \[
        M(x)^*M(x)
        \succeq
        (1 - \inv\beta \gamma \theta \norm{A}^2)\1_{\calI(x)}
        + ((2-\beta)\theta\gamma\epsilon + \theta^2) \1_{\calA(x)}
        \quad\text{for all}\quad x \in U.
    \]
    We have $M(x)^*M(x) \ge C^{-2}\,\Id$ for some $C>0$ with $C\theta < 1$ if both factors in this expression are strictly greater than $\theta^2$. For the second factor, this follows from taking \emph{any} $\beta \in (0, 2)$.
    Keeping in mind our assumption $2(1-\theta^2) > \gamma \theta \norm{A}^2$, the first factor is greater than $\theta^2$ for \emph{some} $\beta \in (0, 2)$ as well.
    The linear convergence claim now follows from \cref{thm:newton:dampened-linear}.

    To show superlinear convergence when $\theta=0$, we will use \cref{thm:newton:superlinear}, which requires us to show that $\norm{D_NH(x)^{-1}}_{\linear(\R^n; \R^n)} \le C$ for some $C>0$.
    Since now $M=D_NH$, this amounts to showing $C^{-2}\,\Id \preceq M(x)^*M(x)$.
    We expand
    \[
        \begin{aligned}[t]
        A^*A \1_{\calA(x)} A^*A
        &
        =
        (\1_{\calA(x)} A^*A \1_{\calA(x)})^2
        + \1_{\calA(x)} A^*A \1_{\calA(x)} A^*A \1_{\calI(x)}
        \\
        \MoveEqLeft[-1]
        + \1_{\calI(x)} A^*A \1_{\calA(x)} A^*A \1_{\calA(x)}
        + \1_{\calI(x)} A^*A \1_{\calA(x)} A^*A \1_{\calI(x)}.
        \end{aligned}
    \]
    We then apply \cref{lemma:sparse:matrix-young} to the middle terms and follow with $\1_{\calA(x)} A^*A \1_{\calA(x)} \ge \epsilon  \1_{\calA(x)}$ to obtain for any $\mu>0$ the bound
    \[
        \begin{aligned}[t]
        \gamma^2 A^*A \1_{\calA(x)} A^*A
        &
        \succeq
        (1-\mu)(\1_{\calA(x)} A^*A \1_{\calA(x)})^2
        - (\inv\mu-1)\1_{\calI(x)} A^*A \1_{\calA(x)} A^*A \1_{\calI(x)}
        \\
        &
        \succeq
        (1-\mu)\epsilon^2 \1_{\calA(x)}
        - (1-\inv\mu) \norm{A}^4 \1_{\calI(x)}.
        \end{aligned}
    \]
    Inserting this lower bound into the expansion $M(x)^*M(x) = \1_{\calI(x)} + \gamma^2 A^*A \1_{\calA(x)} A^*A$ from \eqref{eq:sparse:ssn:mm}, we deduce for some $\mu \in (0, 1)$ the existence of $C>0$ such that $C^{-2}\,\Id \preceq M(x)^*M(x)$. Superlinear convergence now follows from \cref{thm:newton:superlinear}.
\end{proof}

\begin{remark}
    The superlinear convergence of semismooth Newton methods for \eqref{eq:sparse:problem} was proved by \cite{Lorenz:2008}.
\end{remark}

\subsection*{Numerical illustration}

To give a practical perspective on the above algorithms, we illustrate their performance on a simple numerical example.
We take $x \in \R^{1024}$ and $A \in \linear(\R^{1024}; \R^{128})$ as convolution with a Gaussian kernel (standard deviation $\sigma=7$ on the domain $[0, 1024]$) followed by subsampling.
To generate the data $b$, we apply $A$ to the true solution depicted in \cref{fig:sparse:reco}, and apply normally distributed noise of variance $0.03$.
As regularization parameter, we take $\lambda = 0.008$.
For all algorithms, we use the initial iterate $x^0=0$.
For the first-order methods we take the step length $\tau = 0.9/L^2$, where $L$ is an estimate of $\norm{A}$.
For the SSN method, we take the proximal parameter $\gamma=100/L^2$.
Since the basic SSN method ($\theta=0$) does not exhibit convergence, we use the active-dampened variant ($\theta>0$).
Further details on the experimental setup can be found in the accompanying code \cite{nonsmoothbook-codes}.

\begin{figure}[t!]
    \begin{tikzpicture}
        \pgfplotsset{set layers}
        \begin{axis}[%
            scale only axis,
            width=0.85\linewidth,
            height=0.3\linewidth,
            xmin=1,xmax=1024,
            legend pos=north east,
            legend style={xshift=-30ex},
            ylabel={signal magnitude},
            xlabel={coordinate},
            ymax=10,ymin=-1.6,
            ]

            \addplot[color=cb1, line width=0.5pt]
                coordinates { (0,0) };
            \addlegendentry{data}

            \addplot[color=cb3, line width=1pt, mark=o, ycomb,
                     y filter/.expression={y==0 ? nan : y}]
                table[x=coord,y=orig]{sparse_orig.txt};
            \addlegendentry{original}

            \addplot[color=cb2, line width=1pt, mark=square, ycomb,
                     y filter/.expression={y==0 ? nan : y}]
                table[x=coord,y=reco]{sparse_reco_fista.txt};
            \addlegendentry{reconstruction}

        \end{axis}
        \begin{axis}[%
            scale only axis,
            width=0.85\linewidth,
            height=0.3\linewidth,
            xmin=1,xmax=1024,
            axis y line*=right,
            axis x line=none,
            ylabel={data magnitude},
            ymax=0.6,ymin=-0.1,
            ]

            \addplot [color=cb1, line width=0.5pt, y filter/.expression={y}]
                table[x=coord,y=data]{sparse_data.txt};
        \end{axis}
    \end{tikzpicture}
    \caption{Sparse reconstruction data and result.}
    \label{fig:sparse:reco}
    \begin{tikzpicture}
        \SetMinMax{sparse_fb.txt}{value}{\dataFB}
        \UpdMinMax{sparse_fista.txt}{value}{\dataFISTA}%
        \UpdMinMax{sparse_ssn.txt}{value}{\dataSSN}%
        \begin{axis}[%
            xmode=log,
            ymode=log,
            xmin=1,
            ymax=1,
            axis x line*=bottom,
            axis y line*=left,
            ylabel={function value},
            xlabel={iteration count},
            ]

            \addplot [fb] table[x=iter,y=value]{\dataFB};
            \addlegendentry{FB}

            \addplot [fista] table[x=iter,y=value]{\dataFISTA};
            \addlegendentry{FISTA}

            \addplot [ssn] table[x=iter,y=value]{\dataSSN};
            \addlegendentry{inactive-dampened SSN}
        \end{axis}
    \end{tikzpicture}
    \caption{Sparse reconstruction algorithm performance: iterations versus function value.}
    \label{fig:sparse:performance}
    \begin{tikzpicture}
        \SetMinMax{sparse_fb.txt}{value}{\dataFB}
        \UpdMinMax{sparse_fista.txt}{value}{\dataFISTA}%
        \UpdMinMax{sparse_ssn.txt}{value}{\dataSSN}%
        \begin{axis}[%
            xmode=log,
            ymode=log,
            ymax=1,
            axis x line*=bottom,
            axis y line*=left,
            ylabel={function value},
            xlabel={time [s]},
            ]

            \addplot [fb] table[x=cputime,y=value]{\dataFB};
            \addlegendentry{FB}

            \addplot [fista] table[x=cputime,y=value]{\dataFISTA};
            \addlegendentry{FISTA}

            \addplot [ssn] table[x=cputime,y=value]{\dataSSN};
            \addlegendentry{inactive-dampened SSN}
        \end{axis}
    \end{tikzpicture}
    \caption{Sparse reconstruction algorithm performance: time (in seconds) versus function value.}
    \label{fig:sparse:performance:cputime}
\end{figure}

We show the data and the reconstructions in \cref{fig:sparse:reco} and algorithm performance in \cref{fig:sparse:performance,fig:sparse:performance:cputime}. As predicted by the theory, the inertial acceleration of FISTA \eqref{eq:sparse:fista} makes it faster than the unaccelerated forward-backward splitting method \eqref{eq:sparse:fb}. Since the SSN method has to be dampened and hence converges only linearly in this ill-posed setting, it is clearly outperformed by FISTA.

\section{Stability under perturbations}
\label{sec:sparse:stability}

We now study stability of solutions to the $\ell^1$-regularized least squares problem \eqref{eq:sparse:problem}.
We add a further perturbation parameter $p\in\R^N$ to $J$, setting
\begin{equation*}
    J(x; p)
    \defeq \frac{1}{2}\norm{Ax-b-p}_2^2 + \lambda\norm{x}_1,
\end{equation*}
so that $J(x) \defeq J(x; 0)$.
Then
\begin{equation*}
    \subdiff_x J(x; p) = A^*(Ax-b-p) + \lambda \subdiff\norm{\freevar}_1(x),
\end{equation*}
so that for perturbed data the solution mapping is given by
\begin{equation*}
    S(p)
    \defeq
    \{ x \in \R^M \mid 0 \in \subdiff_x J(x; p) \}
    =
    \{ x \in \R^M \mid A^*p \in \subdiff J(x) \}
    = [\inv{(\subdiff J)} \circ A^*](p).
\end{equation*}

The next result shows that the Lasso problem is data-stable at the solution $\realoptx$ for data $b$ if (for simplicity) the solution is strictly complementary, and the matrix $A^*A$ is invertible on the subspace corresponding to the active (i.e., explaining) features.
This is the same condition as in the convergence \cref{thm:sparse:ssn} for the SSN method.
Indeed, due to optimality of $\realoptx$ and the strict complementarity of $\realoptx$ and $\realopt p$, $\calI$ is the same set as $\calI(\realoptx)$ defined for the SSN method in \eqref{eq:sparse:ssn:inactive-active}. This can seen by using the proximal characterization \eqref{eq:splitting:optsys_exp} of the optimality conditions, and the zero-projection properties of the soft-thresholding operator, \cref{ex:proximal:reell}\,\cref{ex:proximal:reell:ii}.

\begin{theorem}
    \label{theorem:regularity:lasso:reg}
    For \eqref{eq:sparse:problem}, suppose $0 \in \subdiff J(\realoptx)$ and that $\realoptx$ and $\realopt p \defeq -\inv\lambda A^*(A\realoptx-b) \in \subdiff \norm{\freevar}_1(\realoptx)$ satisfy the strict complementarity condition \eqref{eq:apple_inverse:lasso:sc}.
    Let
    \begin{equation*}
        \calI \defeq \{ i \in \{1, \ldots, M \} \mid \realoptx_i = 0\}
    \end{equation*}
    be the set of inactive indices, and set
    \begin{equation*}
        V=\setof{x\in \R^M}{x_i = 0 \text{ for }i\in \calI}.
    \end{equation*}
    Denote by $P_V$ the orthogonal projection onto $V$.
    Then $S$ has the Aubin property at $0$ for $\realoptx$ if $P_V A^*A P_V^*$ is nonsingular on $V$.
\end{theorem}

\begin{proof}
    Since the solution mapping $S$ has the Aubin property at $0$ for $\realoptx$ if $\subdiff J$ is metrically regular at $\realoptx$ for $0$,
    we use \cref{cor:regularity:morduk-metric:finite} to verify the latter. To do so, we need an expression for $\coderivative[\subdiff J]$.
    For simplicity of notation, we write $g \defeq \norm{\freevar}_1$, so that $J(x)=\frac{1}{2}\norm{Ax-b}_2^2 + \lambda g(x)$.
    Then \cref{thm:subdiff:sum,ex:convex:subdiff_abs} give
    \begin{equation*}
        \subdiff J(x)=A^T(Ax-b) + \lambda \subdiff g(x)
        \quad\text{for}\quad
        \subdiff g(\realoptx)=\prod_{i=1}^M
        \begin{cases}
            \sign \realoptx_i & \text{if }\realoptx_i \ne 0, \\
            [-1, 1] & \text{if } \realoptx_i=0.
        \end{cases}
    \end{equation*}
    To calculate $\coderivative[\subdiff J](\realoptx|0)$, we need $\subdiff g$ to be graphically regular at $\realoptx$ for $\realopt p$. By \cref{lemma:graphical:absvalue}, this is equivalent to the strict complementarity assumed in \eqref{eq:apple_inverse:lasso:sc}.

    Since the first part of $\subdiff J(x)$ is single-valued and linear, using \cref{thm:colimiting:addition,thm:colimiting:outer} together with the assumption \eqref{eq:apple_inverse:lasso:sc}, we obtain for any $p^* \in \R^M$ that
    \begin{equation}
        \label{eq:lasso-second-diff}
        D^*[\subdiff J](\realoptx|0)(p^*)
        =A^*A p^* +  \lambda \coderivative[\subdiff g](\realoptx|\realopt p)(p^*).
    \end{equation}
    In \cref{lemma:graphical:absvalue}, for the strictly complementary cases \eqref{eq:apple_inverse:lasso:sc}, we have already calculated that
    \begin{equation*}
        \coderivative [\subdiff g]( \realoptx | \realopt p )(p^*) =
        \prod_{i=1}^M
        \begin{cases}
            \{0\}
            &\text{if }  \realoptx_i  \ne 0,\,  \realopt p_i  = \sign  \realoptx_i , \\
            \R
            &\text{if }  \realoptx_i =0,\, [p^*]_i  = 0,\, \abs{ \realopt p_i } < 1, \\
            \emptyset
            & \text{otherwise.}
        \end{cases}
    \end{equation*}
    From \eqref{eq:lasso-second-diff} we now obtain
    \begin{equation*}
        D^*[\subdiff J](\realoptx|0)(p^*)=
        \begin{cases}
            A^*Ap^* + V^\perp, & p^* \in V, \\
            \emptyset, & p^* \not\in V.
        \end{cases}
    \end{equation*}
    Note how $\lambda$ disappears from the expression, as $V$ and $V^\perp$ are subspaces and thus invariant under multiplication by $\lambda$.
    We then calculate
    \begin{equation*}
        \begin{aligned}
            \abs{\coderivative[\subdiff J](\realoptx|0)^{-1}}^+
            &
            =\sup\{\norm{p^*} \mid \exists \dir p^* \in \coderivative[\subdiff g](\realoptx|0)(p^*),\, \norm{p^*} \le 1\}
            \\
            &
            =\sup\{\norm{p^*} \mid \dir x \in V,\, z \in V^\perp\, \norm{A^*A p^* + z} \le 1\}
            \\
            &
            =\sup\{\norm{p^*} \mid \dir x \in V,\, \norm{P_V A^*A P_V^* p^*} \le 1\}.
        \end{aligned}
    \end{equation*}
    Thus \cref{cor:regularity:morduk-metric:finite} shows that $\subdiff J$ is metrically regular at $\realoptx$ for $0$.
\end{proof}

We can also prove sensitivity with respect to the regularization parameter under the exact same conditions as in the previous theorem.

\begin{theorem}
    Suppose that the conditions of \cref{theorem:regularity:lasso:reg} hold and that $P_V A^*A P_V^*$ is nonsingular on $V$.
    Let
    \begin{equation*}
        Z(\alt\lambda) \defeq \{x \in \R^M \mid 0 \in \subdiff \tilde J(x; \alt\lambda)\}\quad\text{for}\quad
       \tilde J(x; \alt\lambda) \defeq \frac{1}{2}\norm{Ax-b}_2^2 + \alt\lambda \norm{x}_1.
    \end{equation*}
    Then $Z$ has the Aubin property at $\lambda$ for any $x \in Z(\lambda)$.
\end{theorem}

\begin{proof}
    In \cref{thm:stability:regularization}, take $g(x)=\norm{x}_1$ and $h(x)=\frac{1}{2}\norm{Ax-b}_2^2$.
    If we verify \eqref{eq:stability:regularization-stability}, i.e.,
    \begin{equation*}
        0 \in  A^*Ay + \lambda D^*[\subdiff g](\opt x|-\inv{\lambda} A^*(A\optx-b))(y)
        \implies y=0,
    \end{equation*}
    then \cref{thm:stability:regularization} establishes that $Z$ has the Aubin property at $\lambda$.
    The strict complementarity condition \eqref{eq:apple_inverse:lasso:sc} implies that either $\optx_i \ne 0$ or $\abs{[\inv{\lambda} A^*(A\optx-b)]_i} < 1$ for all components $i=1,\ldots,M$.
    Therefore, \cref{lemma:graphical:absvalue} shows that
    \[
        D^*[\subdiff g](\opt x|-\inv{\lambda} A^*(A\optx-b))(y) = V^\perp\neq \emptyset
    \]
    if and only if $y \in V$.
    Consequently, \eqref{eq:stability:regularization-stability} becomes
    \begin{equation*}
        0 \in A^*A y + V^\perp \ \text{and}\ y \in V \implies y=0.
    \end{equation*}
    But this follows from the assumption that $P_VA^*AP_V^*$ is invertible on $V$.
\end{proof}

\begin{remark}[regularization theory]\index{regularization!theory of}
    A proof of convergence of solutions to the sparse regularization problems \eqref{eq:sparse:problem} in the sense of \cref{sec:stability:tikhonov} as $\lambda \downto 0$ can be found in \cite{tuomov-regtheory}.
\end{remark}

\chapter{\texorpdfstring{$\ell^1$}{ℓ¹} fitting}
\label{chap:l1fit}

Nonsmooth norms are not only useful as regularization terms. In \eqref{eq:sparse:problem-motiv}, the use of the sum of squares as a data fitting term was justified by statistical arguments: For Gaussian noise, its minimizer coincides with the \term{mean} of the signal, which is the \term[estimator, maximum likelihood]{maximum likelihood estimator} under this assumption. However, that connection is lost for non-Gaussian noise, in particular if the data contains \term{outliers} (rare deviations of much larger magnitude than a normal distribution would predict). A particular such error model is \term[noise!impulsive]{impulsive noise}, which is characterized by containing \emph{only} outliers. Such errors are relevant in digital signal and image processing, where they can arise from malfunctioning pixels in camera sensors, faulty memory locations in hardware, or transmission in noisy channels.
One particular model is \term[noise!impulsive!random-valued]{random-valued impulsive noise}, which corresponds to additive errors of the form
\begin{equation*}
    \eta(x) =
    \begin{cases}
        \xi& \text{with probability $r$},\\
        0 & \text{with probability $1-r$},
    \end{cases}
\end{equation*}
where $r\in[0,1]$ is the fraction of faulty channels and the normally distributed random variable $\xi$ with mean $0$ and variance $\sigma>0$ is the (independent) noise on each affected channel. A more extreme model is \term[noise!salt-and-pepper]{salt-and-pepper noise}, where the data in each affected channel is replaced by either $0$ or $1$ (modeling, e.g., pixels in CCD sensors that are either defective or saturated by cosmic noise).

Statistically, a more robust estimator in the presence of outliers is the \term{median}, which minimizes -- instead of the sum of squares -- the sum of absolute values; see \cite{Huber:2009,GelmanCarlinSternRubin:2013}. This leads to replacing in \eqref{eq:sparse:problem} the squared $\ell^2$ norm by the (nonsquared) $\ell^1$ norm. (Another motivation for this is that at least for impulsive noise, the model output should match the data everywhere except for the outliers -- i.e., that the residual data mismatch is sparse.)
Due to their relevance in signal and image processing, such problems have attracted increasing interest in the last decade; here we only mention \cite{Karkkainen:2005a,Yang:2008a,CJK:2010} as a sample of relevant work.
To avoid additional complexity, we here consider again the regression problem from \cref{chap:sparse}, where we now assume as in \cite{Karkkainen:2005a,CJK:2010} that the noise is sparse but the solution is smooth.

\section{Problem description}

We consider for $A \in \R^{N \times M}$ and $b \in \R^N$ as in \eqref{eq:sparse:problem} the \term[problem!$\ell^1$-fitting]{$\ell^1$-fitting problem}
\begin{equation}
    \label{eq:l1fit:problem}
    \min_{x \in \R^M} \lambda\norm{Ax-b}_1 + \frac{1}{2}\norm{x}_2^2,
\end{equation}
where $\lambda>0$ is a (inverse) regularization parameter related to the noise level. (The benefit of writing the problem in this form instead of using $\alpha\defeq \lambda^{-1}$ as in \cref{chap:sparse} will become apparent in the following.)
The regularization term is smooth to indicate no sparsity requirements on the reconstructed signal, merely the desire for small values.

\section{Optimality conditions}

Using the same approach as in \cref{chap:sparse}, we obtain optimality conditions for \eqref{eq:l1fit:problem}.
We write the problem in the canonical form
\begin{equation*}
    \min_{x \in \R^M}~ J(x) \quad\text{where}\quad J(x) \defeq F(Ax) + G(x)
\end{equation*}
by taking
\[
    F(y) \defeq \lambda\norm{y-b}_1
    \quad\text{and}\quad
    G(x) \defeq \frac{1}{2}\norm{x}_2^2.
\]
We then have the following explicit optimality conditions.

\begin{theorem}
    A vector $\realoptx \in \R^M$ is a solution to \eqref{eq:l1fit:problem} if and only if there exists $\realopty \in \R^{N}$ such that
    \begin{equation}
        \label{eq:l1fit:primal-oc}
        -\realoptx = A^* \realopty
        \qquad\text{ and }\qquad
        \realopt y_i \in
        \begin{cases}
            \{\lambda\} & \text{if }[A\realoptx-b]_i>0,\\
            \{-\lambda\} & \text{if }[A\realoptx-b]_i<0,\\
            [-\lambda,\lambda] & \text{if }[A\realoptx-b]_i=0.
        \end{cases}
    \end{equation}
\end{theorem}

\begin{proof}
    Since $F$ and $G$ are convex and $A$ is linear, also $J$ is convex.
    Therefore the convex Fermat principle of \cref{thm:convex:fermat} characterizes the solution of \eqref{eq:l1fit:problem} as those $\realoptx$ satisfying $0 \in \subdiff J(\realoptx)$.
    Since both $F$ and $G$ have full domains and are proper and lower semicontinuous, we may further use the subdifferential sum rule of \cref{thm:subdiff:sum} and the chain rule of \cref{thm:convex:chain} to calculate for all $x \in \R^m$ that $\subdiff J(x) = A^* \subdiff F(Ax) + \subdiff G(x)$.
    Since $G$ is differentiable, using \cref{thm:convex:gateaux} we therefore characterize the solutions as those points $\realoptx$ satisfying
    \begin{equation}
        \label{eq:l1fit:primal-oc:0}
        -\realoptx \in A^* \subdiff F(A\realoptx).
    \end{equation}
    Using \cref{ex:convex:subdiff_abs} to calculate $\subdiff F(A\realoptx)$ componentwise, we obtain \eqref{eq:l1fit:primal-oc}.
\end{proof}

Based on the Fenchel--Rockafellar theorem (\cref{thm:convex:fenchel}), we may alternatively study optimality conditions for the dual problem
\begin{equation*}
    \min_{y \in \R^N} Q(y) \defeq F^*(y) + G^*(-A^* y).
\end{equation*}
We know from \cref{lem:convex:power-conjugate} that $G^*(x)=\frac{1}{2}\norm{x}_2^2$.
By \cref{ex:convex:fenchel}\,\ref{ex:convex:fenchel:iii} and \cref{lem:convex:fenchel_calc}  we also calculate that
\begin{equation}
    \label{eq:l1fit:fstar}
    F^*(y) = \delta_{\lambda \B_\infty}(y) + \iprod{b}{y}.
\end{equation}
Therefore, the dual problem is given by
\begin{equation}
    \label{eq:l1fit:dual-problem}
    \min_{y \in \R^N} \delta_{\lambda \B_\infty}(y) + \iprod{b}{y} + \frac{1}{2}\norm{A^* y}_2^2.
\end{equation}
For this problem, we can also derive explicit optimality conditions.
\begin{theorem}
    \label{thm:l1fit:dual-oc}
    A vector $\realopty \in \R^N$ is a solution to \eqref{eq:l1fit:dual-problem} of \eqref{eq:l1fit:problem} if and only if there exists a $\realopt p\in\R^N$ such that
    \begin{equation}
        \label{eq:l1fit:dual-oc}
        -AA^*\realopty = \realopt p
        \qquad\text{ and }\qquad
        [\realopt p - b]_i\in
        \begin{cases}
            [0,\infty) & \text{if }\realopty_i=\lambda,\\
            0 & \text{if } \realopty_i \in (-\lambda, \lambda),\\
            (-\infty, 0] & \text{if }\realopty_i=-\lambda, \\
            \emptyset & \text{otherwise}.
        \end{cases}
    \end{equation}
\end{theorem}

\begin{proof}
    Again, the Fermat principle characterizes the solutions via $0 \in \subdiff Q(\realopty)$.
    Since $G^*$ has a full domain and both $F^*$ and $G^*$ are proper and lower semicontinuous, we may further use the subdifferential sum rule of \cref{thm:subdiff:sum} and the chain rule of \cref{thm:convex:chain} to calculate for all $y \in \R^N$ that $\subdiff Q(y) = -A \subdiff G^*(-A^* y) + \subdiff F^*(x)$. By the differentiability of $G^*$, again any dual solution $\hat y$ is therefore characterized by
    \begin{equation}
        \label{eq:l1fit:dual-oc:0}
        -AA^*\realopty = A \grad G^*(-A^* \realopty) \in \subdiff F^*(\realopty).
    \end{equation}
    Using from \cref{ex:convex:subdiff_ind} the expression of the subdifferential of the indicator function of an interval, we obtain \eqref{eq:l1fit:dual-oc}.
\end{proof}

We can also characterize the primal and dual solutions through a primal-dual system.

\begin{theorem}
    \label{thm:l1fit:primal-dual-oc}
    The solutions $\realoptx \in \R^M$ and $\realopty \in \R^N$ to the primal problem \eqref{eq:l1fit:problem}  and the dual problem \eqref{eq:l1fit:dual-problem} are simultaneously characterized by \eqref{eq:l1fit:primal-oc} or, equivalently,
    \begin{equation}
        \label{eq:l1fit:primal-dual-oc2}
        -A^*\realopty  = \realoptx
        \qquad\text{and}\qquad
        [A\realoptx - b]_i \in
        \begin{cases}
            [0,\infty) & \text{if }\realopty_i=\lambda,\\
            0 & \text{if }\realopty \in (-\lambda, \lambda),\\
            (-\infty, 0] & \text{if } \realopty=-\lambda, \\
            \emptyset & \text{otherwise}.
        \end{cases}
    \end{equation}
\end{theorem}

\begin{proof}
    According to \cref{thm:convex:fenchel}, the primal and dual solutions are characterized by
    \begin{equation}
        \label{eq:sparse:primal-dual-oc0}
        \realopty \in \subdiff F(A\realoptx)
        \quad\text{and}\quad
        -A^*\realopty \in \subdiff G(\realoptx).
    \end{equation}
    This expands as \eqref{eq:l1fit:primal-oc} where $\realopty$ is indeed the dual variable.
    By the Fenchel--Young lemma (\cref{lem:convex:fenchel-young}), the conditions \eqref{eq:sparse:primal-dual-oc0} can equivalently be written as
    \begin{equation*}
        A\realoptx \in \subdiff F^*(\realopty)
        \quad\text{and}\quad
        -A^*\realopty \in \subdiff G(\realoptx).
    \end{equation*}
    Similarly to the proof of \cref{thm:l1fit:dual-oc}, this condition becomes \eqref{eq:l1fit:primal-dual-oc2}.
\end{proof}

One may note that the primal-dual condition \eqref{eq:l1fit:primal-dual-oc2} implies the dual condition  \eqref{eq:l1fit:dual-oc} with $\realopt p=A\realoptx$.

\section{Algorithms}

Once more, the starting point for implementable algorithms is the proximal point reformulation of the optimality conditions, this time for the dual problem.
\begin{lemma}
    A vector $\realopty \in \R^N$ is a solution to \eqref{eq:l1fit:dual-problem} of \eqref{eq:l1fit:problem} if and only if
    \begin{equation}
        \label{eq:l1fit:prox-oc}
        \realopty = \proj_{\lambda \B_\infty}(\realopty - \tau[A A^* \realopty + b]).
    \end{equation}
\end{lemma}
\begin{proof}
    Recalling \cref{lem:proximal:subdiff}, we may also rewrite \eqref{eq:l1fit:dual-oc:0} in terms of the proximal operator of $F^*$, for any $\tau>0$ (we just multiply  \eqref{eq:l1fit:primal-oc:0} by $\tau$) as
    \begin{equation*}
        \realopty = \prox_{\tau F^*}(\realopty + \tau A \grad G^*(-A^*\realopty)).
    \end{equation*}
    Using the expression for $F^*$ in \eqref{eq:l1fit:fstar} and the definition of the conjugate, we have for any $y$ that
    \begin{equation*}
        \prox_{\tau F^*}(y) = \prox_{\tau\delta_{\lambda \B_\infty}}(y - \tau b)
        = \proj_{\lambda \B_\infty}(y-\tau b).
    \end{equation*}
    Hence we obtain \eqref{eq:l1fit:prox-oc}.
\end{proof}
\subsection*{Dual forward-backward splitting}

Following \cref{sec:splitting:explicit}, we obtain from \eqref{eq:l1fit:prox-oc} the forward-backward splitting method
\begin{algeqbox}
    \begin{equation}
        \label{eq:l1fit:dualfb}
        \nexty \defeq \proj_{\lambda\B_\infty}(\thisy - \tau[AA^*\thisy+b]).
    \end{equation}
\end{algeqbox}
As an immediate consequence of \cref{thm:convergence:fb}, the iterates of \eqref{eq:l1fit:dualfb} converge subject to a bound on the step length parameter $\tau>0$.
\begin{theorem}
    Suppose $\tau \norm{A}^2<2$.
    Then for any starting point $y^0 \in \R^N$, the iterates $\{\thisy\}_{k \in \N}$ generated by \eqref{eq:l1fit:dualfb} converge to a solution $\realopty$ of the dual problem \eqref{eq:l1fit:dual-problem}.
\end{theorem}
By \cref{thm:l1fit:primal-dual-oc}, the primal and dual solutions $\realoptx$ and $\realopty$ to \eqref{eq:l1fit:problem} satisfy $\realoptx = -A^*\realopty$. We can therefore recover a primal approximate solution $\thisx = -A^*\thisy$ from a dual approximate solution $\thisy$.

Convergence of function values can be obtained in a similar fashion from \cref{thm:gap:fb:value:nonergodic} under the stricter condition $\tau\norm{A}^2 \le 1$.
Under this condition, we also obtain from \cref{thm:meta:inertia:fb} the $O(1/k^2)$ convergence of the inertial variant

\begin{algeqbox*}
    \begin{equation*}
        \left\{\begin{aligned}
                \nexty &\defeq \proj_{\lambda\B_\infty}(\thisy - \tau[AA^*\this{\bar y}+b]),
                \\
                \alpha_{k+1} &\defeq \lambda_{k+1}(\inv\lambda_k-1),
                \\
                \nexxt{\bar y} & \defeq (1+\alpha_{k+1})\nexty-\alpha_{k+1}\thisy.
            \end{aligned}
        \right.
    \end{equation*}
\end{algeqbox*}

Here the initial inertial parameter is initialized with $\lambda_0=1$, while $\bar y^0 \in \R^N$ can be chosen freely.
Since $\alpha_1=0$, the initial iterate $y^0$ is in fact never used.

\subsection*{Primal-dual proximal splitting}

We can also apply the PDPS method \eqref{eq:splitting:pd} to \eqref{eq:l1fit:problem} by taking
\[
    F(x)=\frac{1}{2}\norm{x}_2^2,
    \qquad
    G(z)=\frac{1}{2}\norm{z-b}_1,
    \qquad
    K=A,
\]
in the canonical problem \eqref{eq:pdps-problem}.
Using \cref{ex:convex:fenchel,lem:convex:power-conjugate,lem:convex:fenchel_calc} we see that $G^*(y) = \delta_{\lambda \B_\infty}(y) + \iprod{y}{b}$. Consequently, it is not difficult to verify that we then have
\begin{equation*}
    \prox_{\sigma G^*}(y)=\proj_{\lambda \B_\infty}(y - \sigma b).
\end{equation*}
The projection reduces to a simple componentwise \enquote{clamping} of values in the range $[-\lambda,\lambda]$.
The PDPS method of \eqref{eq:splitting:pd} then becomes
\begin{algeqbox}
    \begin{equation}
        \label{eq:l1fit:pdps}
        \left\{\begin{aligned}
                x^{k+1} & \defeq \frac{1}{1+\tau}(x^k - \tau A^*y^k),\\
                \overnextx & \defeq 2x^{k+1}-x^k,\\
                y^{k+1} & \defeq proj_{\lambda \B_\infty}(y^{k} + \sigma (A\bar x^{k+1}-b)).
        \end{aligned}        \right.
    \end{equation}
\end{algeqbox}

The method converges subject to a simple step length condition.

\begin{theorem}
    Suppose $\tau \sigma\norm{A}^2<1$.
    Then for any starting point $(x^0, y^0) \in \R^M \times \R^N$, the iterates $\{(\thisx, \thisy)\}_{k \in \N}$ generated by \eqref{eq:l1fit:pdps} converge to solutions $\realoptx$ and $\realopty$ of \eqref{eq:l1fit:problem} and \eqref{eq:l1fit:dual-problem}, respectively.
\end{theorem}

Since $F$ is strongly convex with factor $\gamma=1$, we can also apply the accelerated method of \eqref{eq:testing:pdps:forward}, updating the step length parameter according to \eqref{eq:testing:pdps:accel} in
\begin{algeqbox}
    \begin{equation}
        \label{eq:l1fit:pdps:accel}
        \left\{\begin{aligned}
                \omega_k & \defeq 1/\sqrt{1+2\tau_k},
                \quad
                \tau_{k+1} \defeq \tau_k\omega_k,
                \quad
                \sigma_{k+1} \defeq \sigma_k/\omega_k,
                \\
                x^{k+1} & \defeq \frac{1}{1+\tau_k}(x^k - \tau_k A^*y^k),\\
                \overnextx & \defeq (1+\omega_k)x^{k+1} - \omega_k x^k,\\
                y^{k+1} & \defeq \proj_{\lambda \B_\infty}(y^{k} + \sigma_{k+1} (A\bar x^{k+1}-b)).
        \end{aligned}\right.
    \end{equation}
\end{algeqbox}

\Cref{thm:testing:pdps:accel} immediately yields its convergence.

\begin{theorem}
    Suppose $\tau_0 \sigma_0 \norm{A}^2<1$.
    Then for any starting point $(x^0, y^0) \in \R^M \times \R^N$, the primal iterates $\{\thisx\}_{k \in \N}$ generated by \eqref{eq:l1fit:pdps:accel} converge to a minimizer $\realoptx$  of \eqref{eq:l1fit:problem} at the rate $O(1/k^2)$.
\end{theorem}

Convergence of the Lagrangian duality gap can be obtained from \cref{thm:gap:pdps} or, in the accelerated case, \cref{thm:gap:accel:pdps}, under the same conditions as for iterate convergence.

\subsection*{Semismooth Newton method}

Similar to sparse regularization from \cref{chap:sparse}, we apply a semismooth Newton method to the proximal point reformulation \eqref{eq:l1fit:prox-oc} by looking for a root $\realopty$ of
\begin{equation*}
    H(y) \defeq  y - \proj_{\lambda \B_\infty}(y-\tau(AA^*y+b))
\end{equation*}
with arbitrary $\tau>0$. From \cref{ex:newton:rn}\,\cref{ex:newton:rn:box} and the chain rule \cref{thm:newton:chain}, a Newton derivative in direction $h$ is given componentwise by
\begin{equation*}
    \begin{aligned}
        [D_N H(y) h]_i &= \left[h - \1_{[-\lambda,\lambda]}(y-\tau (AA^*y+b)) \odot (h-\tau AA^*h)\right]_i\\
        &= \begin{cases}
            h_i & \text{if } |y_i - \tau[AA^*y+ b]_i| > \lambda,\\
            \tau [AA^* h]_i  & \text{if } |y_i - \tau[AA^*y+ b]_i| \leq \lambda.
        \end{cases}
    \end{aligned}
\end{equation*}
Here we recall the notation $[x\odot y]_i\defeq x_iy_i$ for the componentwise or Hadamard product on $\R^N$.
We can write this concisely as
\[
    D_N H(y)
    =
    \1_{\calA(y^k)} + \tau \1_{\calI(y^k)} AA^*
\]
for the active and inactive sets
\begin{subequations}
\label{eq:l1fit:active-inactive}
\begin{align}
    \calA(y^k) & \defeq \{ i \in \{1, \ldots, N\} \mid |y_i - \tau[AA^*y+ b]_i| > \lambda \},
    \quad\text{and}
    \\
    \calI(y^k) & \defeq \{1,\dots,N\}\setminus \calA(y^k)
\end{align}
\end{subequations}
and the diagonal matrix $\1_{A}$ with $[\1_{A}]_{ii}=1$ if $i\in A\subset\{1,\dots,N\}$ and $0$ otherwise.
Thus the semismooth Newton algorithm is $\nexty \defeq \thisy + \this s$, where we solve for $\this s$ in
\[
    (\1_{\calA(y^k)} + \tau \1_{\calI(y^k)} AA^*)\this s =
    - H(\thisy).
\]
Proceeding as for \eqref{eq:sparse:ssn:first-system}, we deduce that
$
    \this s_i = - H(\thisy)_i = [\proj_{\lambda \B_\infty}(\thisy-\tau(AA^*\thisy+b)) - \thisy]_i
$
for $i \in \calA(y^k)$, hence
\begin{equation}
    \label{eq:l1fit:active}
    \nexxt y_i = [\proj_{\lambda \B_\infty}(\thisy-\tau(AA^*\thisy+b))]_i
    \quad\text{for}\quad
    i \in \calA(y^k).
\end{equation}
For $i \in \calI(y^k)$, we have $[\proj_{\lambda \B_\infty}(\thisy-\tau(AA^*\thisy+b))]_i = [\thisy-\tau(AA^*\thisy+b)]_i$. Hence, by introducing the projection $P_{\calI(y^k)}$ to the inactive set and writing
\[
    \this s
    = P_{\calI(y^k)}^* P_{\calI(y^k)} \this s + \1_{\calA(y^k)} \this s
    = P_{\calI(y^k)}^* P_{\calI(y^k)} \this s - \1_{\calA(y^k)}H(\thisy)
\]
we deduce as after \eqref{eq:sparse:ssn:first-system} that the inactive components $P_{\calI(y^k)}\this s$ are characterized by
\begin{equation}
    \label{eq:l1fit:inactive-system}
    \tau P_{\calI(y^k)} AA^* P_{\calI(y^k)}^* [P_{\calI(y^k)} \this s] =
    - P_{\calI(y^k)}(\Id - \tau AA^*\1_{\calA(y^k)})H(\thisy).
\end{equation}
Altogether, therefore, the semismooth Newton method for the dual problem \eqref{eq:l1fit:dual-problem} iterates
\begin{algenumbox}
    \begin{enumerate}[label=\arabic*.]
        \item\label{item:l1:ssn:first-step} form the active and inactive sets $\calA(y^k)$ and $\calI(y^k)$ following \eqref{eq:l1fit:active-inactive};
        \item update $y^{k+1}_i$ for $i \in \calA(y^k)$ by  \eqref{eq:l1fit:active};
        \item solve $P_{\calI(y^k)}s^k$ from \eqref{eq:l1fit:inactive-system};
        \item\label{item:l1:ssn:last-step} update $y^{k+1}_i \defeq y^k_i + s^k_i$ for $i \in \calI(y^k)$.
    \end{enumerate}
\end{algenumbox}
From the dual iterate $y^k$, an approximation of the corresponding primal solution can again be recovered via $x^k \defeq  -A^*y^k$.

Completely analogously to the proof of superlinear convergence in \cref{thm:sparse:ssn}, \cref{thm:newton:superlinear} establishes the following convergence result.

\begin{theorem}
    \label{thm:l1:ssn}
    Let $\opt y$ be a (unique) minimizer of the dual problem \eqref{eq:l1fit:dual-problem} to \eqref{eq:l1fit:problem} and $\gamma>0$.
    Suppose that $P_{\calI(\opt y)}AA^*P_{\calI(\opt y)}^*$ is positive definite.
    If $y^0$ is sufficiently close to $\opt y$, then $\{y^{k+1}\}_{k \in \N}$ generated by iterating \crefrange{item:l1:ssn:first-step}{item:l1:ssn:last-step} above converge superlinearly to $\opt y$.
\end{theorem}

\subsection*{Numerical illustration}

\begin{figure}[t!]
    \begin{tikzpicture}
        \pgfplotsset{set layers}
        \begin{axis}[%
            scale only axis,
            width=0.85\linewidth,
            height=0.3\linewidth,
            xmin=1,xmax=1024,
            axis x line*=bottom,
            axis y line*=left,
            legend pos=north west,
            legend style={xshift=2ex, yshift=3ex},
            ylabel={signal magnitude},
            xlabel={coordinate},
            ymax=1,ymin=-2,
            ]

            \addplot[color=cb3, line width=1pt]
                table[x=coord,y=orig]{l1fit_orig.txt};
            \addlegendentry{original}

            \addplot[color=cb2, line width=1pt]
                table[x=coord,y=reco]{l1fit_reco_pdps_accel.txt};
            \addlegendentry{reconstruction}

            \addplot[color=cb1, line width=0.5pt] coordinates { (0,0) };
            \addlegendentry{noiseless data}

            \addplot[color=cb3, line width=0.5pt, dashed] coordinates { (0,0) };
            \addlegendentry{noisy data}
        \end{axis}
        \begin{axis}[%
            scale only axis,
            width=0.85\linewidth,
            height=0.3\linewidth,
            xmin=1,xmax=1024,
            axis y line*=right,
            axis x line=none,
            ylabel={data magnitude},
            ymax=12,ymin=-4,
            ]

            \addplot[color=cb1, line width=0.5pt]
                table[x=coord,y=noiseless]{l1fit_data.txt};

            \addplot[color=cb2, line width=0.75pt, dashed]
                table[x=coord,y=data]{l1fit_data.txt};

        \end{axis}
    \end{tikzpicture}
    \caption{$\ell^1$ fitting data and result.}
    \label{fig:l1fit:reco}
    \begin{tikzpicture}
        \SetMinMax{l1fit_dualfb.txt}{value}{\dataDualFB}%
        \UpdMinMax{l1fit_dualfista.txt}{value}{\dataDualFISTA}
        \UpdMinMax{l1fit_pdps.txt}{value}{\dataPDPS}%
        \UpdMinMax{l1fit_pdps_accel.txt}{value}{\dataPDPSAccel}%
        \UpdMinMax{l1fit_ssn.txt}{value}{\dataSSN}%
        \begin{axis}[%
            xmode=log,
            xmin=1,
            xlabel={iteration count},
            ylabel={function value},
            fixedylog = {3}{fixed,precision=0},
            ]

            \addplot [fb] table[x=iter,y=value]{\dataDualFB};
            \addlegendentry{dual FB}

            \addplot [fista] table[x=iter,y=value]{\dataDualFISTA};
            \addlegendentry{dual FISTA}

            \addplot [pdps] table[x=iter,y=value]{\dataPDPS};
            \addlegendentry{PDPS}

            \addplot [pdps accel] table[x=iter,y=value]{\dataPDPSAccel};
            \addlegendentry{accelerated PDPS}

            \addplot [ssn] table[x=iter,y=value]{\dataSSN};
            \addlegendentry{SSN}
        \end{axis}
    \end{tikzpicture}
    \caption{$\ell^1$ fitting reconstruction algorithm performance: iteration vs. function value.}
    \label{fig:l1fit:performance}
    \begin{tikzpicture}
        \SetMinMax{l1fit_dualfb.txt}{value}{\dataDualFB}%
        \UpdMinMax{l1fit_dualfista.txt}{value}{\dataDualFISTA}
        \UpdMinMax{l1fit_pdps.txt}{value}{\dataPDPS}%
        \UpdMinMax{l1fit_pdps_accel.txt}{value}{\dataPDPSAccel}%
        \UpdMinMax{l1fit_ssn.txt}{value}{\dataSSN}%
        \begin{axis}[%
            xmode=log,
            xlabel={time [s]},
            ylabel={function value},
            fixedylog = {3}{fixed,precision=0},
            ]

            \addplot [fb] table[x=cputime,y=value]{\dataDualFB};
            \addlegendentry{dual FB}

            \addplot [fista] table[x=cputime,y=value]{\dataDualFISTA};
            \addlegendentry{dual FISTA}

            \addplot [pdps] table[x=cputime,y=value]{\dataPDPS};
            \addlegendentry{PDPS}

            \addplot [pdps accel] table[x=cputime,y=value]{\dataPDPSAccel};
            \addlegendentry{accelerated PDPS}

            \addplot [ssn] table[x=cputime,y=value]{\dataSSN};
            \addlegendentry{SSN}
        \end{axis}
    \end{tikzpicture}
    \caption{$\ell^1$ fitting reconstruction algorithm performance: time (in seconds) vs. function value.}
    \label{fig:l1fit:performance:cputime}
\end{figure}

Again we illustrate the performance of the aforementioned algorithms on a simple numerical example.
We take $x \in \R^{1024}$ and $A \in \linear(\R^{1024}; \R^{128})$ as convolution with a Gaussian kernel (standard deviation $\sigma=7$ on the domain $[0, 1024]$) followed by subsampling.
To generate the data $b$, we apply $A$ to the true signal depicted in \cref{fig:sparse:reco}, and follow with salt-and-pepper noise of magnitude $1.8$.
As the inverse regularization parameter, we take $\lambda = 6.5$.
For all algorithms, we use the initial iterate $x^0=0$.
For the forward-backward type methods, we take the step length $\tau = 0.9/L^2$, where $L$ is again an estimate of $\norm{A}$.
For the PDPS method, we take the step lengths $\tau = 0.5/L$ and $\sigma = 1.9/L$.
For the SSN method, we take the proximal parameter $\gamma=9/L^2$.
The precise experimental details can be found in the accompanying code \cite{nonsmoothbook-codes}.

We illustrate the data and the reconstruction in \cref{fig:l1fit:reco} and the convergence behavior in \cref{fig:l1fit:performance,fig:l1fit:performance:cputime}.
The latter clearly show that acceleration improves performances of the first-order methods, but the superlinearly convergent SSN method requires significantly fewer iterations than any first-order method. In fact, even though each iteration of the former is much more expensive, the total time to reach the objective is still smaller. On the other hand, first-order methods are much faster in reducing the objective value in the beginning and therefore may be the method of choice if high accuracy is not desired.

\chapter{Total variation regularization}
\label{chap:tv}

We now turn to \term[problem!image processing]{mathematical image processing}, where the unknown to be reconstructed from data is a digital image.
The most basic mathematical image processing task is \term[problem!denoising]{denoising}, i.e., removing the noise in an image (for example, a photograph taken in low light conditions), which corresponds to taking the forward operator as the identity.
More advanced image processing tasks include \term[problem!inpainting]{inpainting}, \term[problem!deblurring]{deblurring}, and \term[problem!superresolution]{superresolution}. These correspond to filling in missing parts of an image, reducing blur caused by defocussed lenses or motion, and recovering additional detail, and involve more complicated linear forward operators.
For an introduction to mathematical image processing, we refer to \cite{Scherzer,BrediesLorenz:2018}.
In true \term[problem!inverse imaging]{inverse imaging problems}, the given data is not itself an image but related to it via some mathematical model describing the physical measurements; examples are magnetic resonance imaging (MRI), involving the Fourier transform \cite{nishimura1996principles}, or positron emission tomography (PET) and computed X-ray tomography (CT), both involving the Radon transform \cite{natterer2001mathematics}.
More challenging imaging modalities such as electrical impedance tomography (EIT) and more advanced MRI techniques require the forward operator $A$ to be nonlinear. We do not treat such operators here, but point towards the primal-dual method of \cref{chap:nlpdps} as one possible solution technique. Alternative Gauss--Newton type methods are introduced by \cite{jauhiainen2019gaussnewton}.

The salient point here is the particular structure of images, which requires an adapted regularization term. The key observation here is that images contain sharp edges (representing jumps in intensity) separating mostly smooth areas. Mathematically, this can be related to requiring sparsity of the \emph{gradient} of the image, rather than the image itself; the corresponding sparse regularization of the gradient is called the \term[regularization!total variation]{total variation regularization}, which was introduced for denoising by \cite{ROF} and has become very popular for other (inverse) imaging tasks such as the ones mentioned above.

Treating such problems in an infinite-dimensional function space framework is very challenging and requires the unknown to be considered in the space $BV(\Omega)$ of \term[function!of bounded variation]{functions of bounded variation} on a domain $\Omega\subset \R^2$, which are characterized by their distributional gradient being a Radon measure $y \in \mathcal{M}(\Omega; \R^2)$. This is a nonreflexive Banach space with a complicated structure; see \cite{Ambrosio,Attouch} for the rich functional analysis and geometric measure theory in this space. As our primary focus here is on algorithms that require a Hilbert space structure, we will treat this problem in a finite-dimensional discretized setting.

\section{Problem description}

We consider the problem
\begin{equation}
    \label{eq:tv:problem}
    \min_x~\frac{1}{2}\norm{Ax-b}_2^2 + \alpha \norm{Dx}_{1,2},
\end{equation}
where $x \in \R^M$ for $M=n_1n_2$ is a vectorization of the two-dimensional image, consisting of an $n_1 \times n_2$ grid of components called \term[pixel]{pixels}; $A \in \R^{N \times M}$ for some $N$ is the linear forward operator; and $D \in \R^{2M \times M}$ is a discretization of the image gradient to be specified below.
We index $x \in \R^M$ using two coordinates $i \in \{1,\ldots,n_1\}$ and $j \in \{1,\ldots,n_2\}$, identifying $x_{ij}$ with $x_{\iota(i, j)}$ for a suitable linear index $\iota$, for example $\iota(i, j)=i+n_1(j-1)$. Likewise we index variables $y \in \R^{2M}$ with $k \in \{1,2\}$ along with $i$ and $j$, identifying $y_{kij}$ with $y_{\iota_2(k, i, j)}$ for a suitable linear index $\iota_2$, for example $\iota_2(k, i, j) = k+2(\iota(i,j)-1)$. We also write $y_{\freevar i j}\defeq(y_{1 i j}, y_{2 i j})\in\R^2$. When necessary for clarity, we insert commas between the indices.
As a discretized derivative, we take forward differences with Neumann boundary conditions, which with the above notation corresponds to setting
\begin{align*}
    [D u]_{1ij} & =
    \begin{cases}
        u_{i+1,j}-u_{i,j}, & 1 \le i < n_1,\, 1 \le j \le n_2, \\
        0, & i=n_1,\, 1 \le j \le n_2,
    \end{cases}
    \\
    [D u]_{2ij} & =
    \begin{cases}
        u_{i,j+1}-u_{i,j}, & 1 \le i \le n_1,\, 1 \le j < n_2, \\
        0, & 1 \le i \le n_2,\, j=n_2.
    \end{cases}
\end{align*}

It remains to discuss the vector-sparsity penalty
\begin{equation*}
    \norm{y}_{1,2} \defeq \sum_{i=1}^{n_1} \sum_{j=1}^{n_2} \norm{y_{\freevar i j}}_2 = \sum_{i=1}^{n_1} \sum_{j=1}^{n_2} \sqrt{y_{1ij}^2 + y_{2ij}^2}\qquad (y\in \R^{2M}).
\end{equation*}
First, it is straightforward to verify that
\begin{equation*}
    (\R^{2M},\norm{\freevar}_{1,2})^* = (\R^{2M},\norm{\freevar}_{\infty,2}),
\end{equation*}
where
\begin{equation*}
    \norm{y}_{\infty,2} \defeq \max{\substack{i=1,\dots,n_1\\j=1,\dots,n_2}} \norm{y_{\freevar ij}}_2,
\end{equation*}
using that
\begin{equation*}
    \dual{y^*,y}_{1,2} \defeq \sum_{i=1}^{n_1}\sum_{j=1}^{n_2}\sum_{k=1}^2 y^*_{kij}y_{kij} \leq \norm{y^*}_{1,2}\norm{y}_{\infty,2}.
\end{equation*}
This allows us to compute various objects by applying the convex analysis of \cref{part:convex} pixelwise, i.e., separately for each pair of pixel coordinates $(i, j)$. First, by applying \cref{thm:subdifferential:norm}, we obtain an explicit expression for the subdifferential.

\begin{lemma}
  Let $y,y^*\in \R^{2M}$. Then $y^* \in \subdiff\norm{\freevar}_{1,2}(y)$ if and only if
  \begin{equation}
    \label{eq:tv:norm12-subdiff}
    y_{\freevar i j}^*
    \in
    \begin{cases}
        \left\{\frac{y_{\freevar i j}}{\norm{y_{\freevar i j}}_2}\right\} & \text{if }y_{\freevar i j} \ne 0,\\
        \B_2 & \text{if }y_{\freevar ij}=0,
    \end{cases}
  \end{equation}
  where $\B_2$ is the Euclidean unit ball in $\R^2$.
\end{lemma}
By \cref{ex:convex:fenchel}\,\ref{ex:convex:fenchel:ii}, the Fenchel conjugate of a norm is given by the indicator functional of the dual unit ball, which in this case is
\[
    \B_{\infty,2} \defeq \{ y \in \R^{2M} \mid \norm{y}_{\infty,2} \le 1\} = \{ y \in \R^{2M} \mid y_{\freevar i j} \in \B_2 \text{ for each } i, j \}.
\]
By \cref{lem:convex:fenchel_calc}\,\ref{lem:convex:fenchel_calc:i}, we thus have
\begin{equation*}
    (\alpha\norm{\freevar}_{1,2})^* = \delta_{\alpha\B_{\infty,2}}.
\end{equation*}
A case distinction similar to \cref{ex:convex:subdiff_ind} then yields the following characterization of the subdifferential.
\begin{lemma}
    Let $y,y^*\in \R^{2M}$ and $\alpha>0$. Then $y^* \in \subdiff\delta_{\alpha\B_{\infty,2}}(y)$ if and only if
  \begin{equation}
    \label{eq:tv:b12-subdiff}
    y_{\freevar i j}^* \in
    \begin{cases}
        [0, \infty) y_{\freevar ij} & \text{if }\norm{y_{\freevar i j}}_2 = \alpha, \\
        0 & \text{if }\norm{y_{\freevar i j}}_2 < \alpha, \\
        \emptyset & \text{otherwise}.
    \end{cases}
\end{equation}
\end{lemma}
Finally, similarly to \cref{lem:proximal:calculus}\,\ref{lem:proximal:calculus:iii} we can show that the corresponding proximal point mapping for $\gamma>0$ is
given pixelwise by
\begin{equation}
    \label{eq:tv:projection}
    [\proj_{\alpha \B_{\infty,2}}(y)]_{\freevar i j}
    = \proj_{\alpha \B_2}(y_{\freevar i j})
    =
    y_{\freevar i j}
    \begin{cases}
        \frac{\alpha}{\norm{y_{\freevar i j}}_2} & \text{if }\norm{y_{\freevar i j}}_2 > \alpha, \\
        1 & \text{if }\norm{y_{\freevar i j}}_2 \le \alpha.
    \end{cases}
\end{equation}

\section{Optimality conditions}

Our derivation of optimality conditions for \eqref{eq:tv:problem} follows that for sparse regularization in \cref{chap:sparse}.
Setting
\[
    F(x) \defeq \frac{1}{2}\norm{Ax-b}_2^2
    \quad\text{and}\quad
    G(y) \defeq \alpha\norm{y}_{1,2},
\]
we can write \eqref{eq:tv:problem} in canonical form as
\begin{equation}
    \label{eq:tv:problem:canonical}
    \min_{x \in \R^M}~ J(x) \quad\text{where}\quad J(x) \defeq F(x) + G(Dx).
\end{equation}
The following result characterizes the solutions of this convex problem.

\begin{theorem}
    \label{thm:tv:primal-oc}
    Let $\realoptx \in \R^M$ be a solution to \eqref{eq:tv:problem}. Then there exists a $\realopty \in \R^{2M}$ such that
    \begin{equation}
        \label{eq:tv:primal-oc}
        -A^*(A\realoptx-b) = D^* \realopty
        \qquad\text{ and }\qquad
        \realopt y_{\freevar i j} \in
        \begin{cases}
            \alpha \left\{\frac{[D\realoptx]_{\freevar i j}}{\norm{[D\realoptx]_{\freevar i j}}_2}\right\} & \text{if }[D\realoptx]_{\freevar i j} \ne 0,\\
            \alpha \B_2 & \text{if }[D\realoptx]_{\freevar i j}=0.
        \end{cases}
    \end{equation}
\end{theorem}

\begin{proof}
    Since $F$ and $G$ are convex, and $D$ is linear, $J$ is convex as well.
    Therefore the Fermat principle of \cref{thm:convex:fermat} characterizes the solution of \eqref{eq:tv:problem} as those $\realoptx$ satisfying $0 \in \subdiff J(\realoptx)$.
    Since both $F$ and $G$ have full domains and are proper and lower semicontinuous, we may further use the subdifferential sum rule of \cref{thm:subdiff:sum} and the chain rule of \cref{thm:convex:chain} to calculate for all $x \in \R^M$ that $\subdiff J(x) = \subdiff F(x) + D^*\subdiff G(Dx)$.
    Since $F$ is differentiable, we can use \cref{thm:convex:gateaux} to characterize the solutions as those points $\realoptx$ satisfying
    \begin{equation*}
        -\grad F(\realoptx) \in D^*\subdiff G(D \realoptx).
    \end{equation*}
    Together with \eqref{eq:tv:norm12-subdiff}, this yields \eqref{eq:tv:primal-oc}.
\end{proof}

The expression for $\realopty$ in \eqref{eq:tv:primal-oc} is difficult to work with in practice, in particular for deriving algorithms.
With the help of the Fenchel--Rockafellar theorem (\cref{thm:convex:fenchel}), we may alternatively study optimality conditions for the dual problem
\begin{equation}
    \label{eq:tv:dual-problem}
    \min_{y \in \R^{2M}} Q(y) \defeq F^*(-D^*y) + G^*(y),
\end{equation}
where $G^* = \delta_{\alpha \B_{\infty,2}}$.
If $A=\Id$, we also obtain a simple expression for $F^*$, which yields the following result.

\begin{theorem}
    \label{thm:tv:dual-oc}
    For $A=\Id$, the solutions $\realopty \in \R^{2M}$ to the dual problem \eqref{eq:tv:dual-problem} of \eqref{eq:tv:problem} are characterized by
    \begin{equation}
        \label{eq:tv:dual-oc}
        - D(D^*\realopty-b) = \realopt p
        \qquad\text{ and }\qquad
        \realopt p_{\freevar i j} \in
        \begin{cases}
            [0,\infty) \realopty_{\freevar i j} & \text{if } \norm{\realopty_{\freevar i j}}_2 = \alpha,\\
            \{0\} & \text{if } \norm{\realopty_{\freevar i j}}_2 < \alpha, \\
            \emptyset & \text{otherwise}.
        \end{cases}
    \end{equation}
\end{theorem}

\begin{proof}
    Again we can apply the Fermat principle.
    We calculate using \cref{lem:convex:power-conjugate} and \ref{lem:convex:fenchel_calc}\,\ref{lem:convex:fenchel_calc:ii} for $K=\Id$ that $F^*(y)=\frac{1}{2}\norm{y}_2^2 + \iprod{b}{y}$.
    Since $F^*$ has a full domain and both $G^*$ and $F^*$ are proper and lower semicontinuous, we may further use the subdifferential sum rule of \cref{thm:subdiff:sum} and the chain rule of \cref{thm:convex:chain} to calculate for all $y \in \R^{2M}$ that
    \begin{equation*}
        \subdiff Q(y) = -D\subdiff F^*(-D^*y) + \subdiff G^*(y).
    \end{equation*}
    By the differentiability of $F^*$, the dual solutions $\realopt y$ are therefore characterized by
    \begin{equation}
        \label{eq:tv:dual-oc:0}
        D\grad F^*(-D^*\realopty) \in \subdiff G^*(\realopty).
    \end{equation}
    Together with \eqref{eq:tv:b12-subdiff}, this yields \eqref{eq:tv:dual-oc}.
\end{proof}

\Cref{thm:convex:fenchel} also gives a primal-dual characterization of optimality. In contrast to the primal result \cref{thm:tv:primal-oc} and the dual result \cref{thm:tv:dual-oc}, it has simple expressions for all variables even for $A\neq\Id$.

\begin{theorem}
    The solutions $\realoptx \in \R^M$ and $\realopty \in \R^{2M}$ to the  primal problem \eqref{eq:tv:problem}  and the dual problem \eqref{eq:tv:dual-problem} are simultaneously characterized by \eqref{eq:tv:primal-oc} or, equivalently,
    \begin{equation}
        \label{eq:tv:primal-dual-oc}
        -D^* \realopty = A^*(A\realoptx-b)
        \qquad\text{ and }\qquad
        [D\realoptx]_{\freevar i j} \in
        \begin{cases}
            [0,\infty) \realopty_{\freevar i j} & \text{if } \norm{\realopty_{\freevar i j}}_2 = \alpha,\\
            \{0\} & \text{if } \norm{\realopty_{\freevar i j}}_2 < \alpha, \\
            \emptyset & \text{otherwise}.
        \end{cases}
    \end{equation}
\end{theorem}

\begin{proof}
    According to \cref{thm:convex:fenchel}, the primal and dual solutions are characterized by
    \begin{equation}
        \label{eq:tv:primal-dual-oc0}
        \realopty \in \subdiff G(D\realoptx)
        \quad\text{and}\quad
        -D^*\realopty = \grad F(\realoptx).
    \end{equation}
    This is simply \eqref{eq:tv:primal-oc}, where $\realopty$ is indeed the dual variable.
    By the Fenchel--Young lemma (\cref{lem:convex:fenchel-young}), the conditions \eqref{eq:tv:primal-dual-oc0} can equivalently be written
    \begin{equation*}
        D\realoptx \in \subdiff G^*(\realopty)
        \quad\text{and}\quad
        -D^*\realopty = \grad F(\realoptx).
    \end{equation*}
    Together with \eqref{eq:tv:b12-subdiff} for an expression of $\subdiff G^*$, this yields \eqref{eq:tv:primal-dual-oc}.
\end{proof}

\section{Algorithms}

Following the approach established in the previous chapters, we now derive some algorithms for \eqref{eq:tv:problem} based on either the dual optimality conditions \eqref{eq:tv:dual-oc} or the primal-dual optimality conditions \eqref{eq:tv:primal-dual-oc}. We start with the former and the corresponding forward-backward type methods.
We then move onto primal-dual splitting methods.
As the discretized gradient $D$ has a nontrivial kernel, semismooth Newton methods cannot be applied directly without dampening as in \cref{chap:sparse}, which would negate the performance advantage over splitting methods. We will therefore focus here on splitting methods, but refer to \cite{Hintermuller:2006a} for a modified semismooth Newton method that retains superlinear convergence.

\subsection*{Dual forward-backward splitting for denoising}

For $A=\Id$, we can directly apply the forward-backward splitting method \eqref{eq:splitting:fb} to the dual problem \eqref{eq:tv:dual-problem} by rewriting the optimality condition \eqref{eq:tv:dual-oc:0} using \cref{lem:proximal:subdiff} for any $\tau>0$ as
\begin{equation*}
    0 = \proj_{\alpha \B_{\infty,2}}(\realopty - \tau D (D^*\realopty-b)).
\end{equation*}
We recall from the proof of \cref{thm:tv:dual-oc} with $A=\Id$ that $F^*(y)=\frac{1}{2}\norm{y}_2^2 + \iprod{b}{y}$ and that the proximal point mapping for the indicator function is given by the metric projection.
Therefore, we obtain the iteration
\begin{algeqbox}
    \begin{equation}
        \label{eq:tv:dualfb}
        \nexty \defeq
        \proj_{\alpha \B_{\infty,2}}(\thisy - \tau D(D^* \thisy - b)),
    \end{equation}
\end{algeqbox}
where the projection operator is given by \eqref{eq:tv:projection}.

By \eqref{eq:tv:primal-dual-oc0}, the primal and dual solutions $\realoptx$ and $\realopty$ satisfy $-D^*\realopty \in \grad F(\realoptx)=\{\realoptx - b\}$, which allows us to recover a primal solution from a dual solution $\realopty$ via $\realoptx = b - D^*\realopty$.

Again, the method converges subject to a simple step length bound.
\begin{theorem}\label{thm:tv:dualfb-convergence}
    Suppose $\tau \norm{D}^2<2$.
    Then for any starting point $y^0 \in \R^{2M}$, the iterates $\{\thisy\}_{k \in \N}$ generated by \eqref{eq:tv:dualfb} converge to a solution $\realopty$ of the dual problem \eqref{eq:tv:dual-problem}.
\end{theorem}

\begin{proof}
    The Lipschitz factor of $\grad[F^* \circ D]$ is $\norm{D}^2$, and hence the claim follows from \cref{thm:convergence:fb}.
\end{proof}

Convergence of function values for the dual objective \eqref{eq:l1fit:dual-problem} can be obtained in a similar fashion from \cref{thm:gap:fb:value:nonergodic} under the stricter condition $\tau\norm{D}^2 \le 1$.

\subsection*{Primal-dual proximal splitting for unitary-simple forward operators}

Dual forward-backward splitting requires that we are able to compute $\grad F^*$, which can be difficult and numerically expensive for general $A\neq \Id$; compare \cref{lem:convex:fenchel_calc}\,\cref{lem:convex:fenchel_calc:iii}. Furthermore, $F^*$ may not even be a smooth function when $A$ is not invertible.
Similarly, the primal-dual proximal splitting \eqref{eq:splitting:pd} for \eqref{eq:tv:problem:canonical} is given by
\begin{algeqbox}
    \begin{equation}
        \label{eq:tv:pdps-unitary}
        \left\{\begin{aligned}
                x^{k+1} & \defeq \prox_{\tau F}(x^k - \tau D^*y^k),\\
                \overnextx &\defeq 2x^{k+1}-x^k,\\
                y^{k+1} & \defeq \proj_{\alpha \B_{\infty,2}}(y^{k} + \sigma D\bar x^{k+1}),
        \end{aligned}        \right.
    \end{equation}
\end{algeqbox}
which still requires computing the proximal mapping of $F$, which can in general be difficult.

However, suppose that $A=SU$ for an unitary operator $U$ and $S$ such that $\tau S^*S+\Id$ has a simple inverse. For example, $U$ can be the Fourier transform and $S$ can be a sub-sampling operator, in which case $\tau S^*S+\Id$ is diagonal; such type of problems appear in magnetic resonance imaging. 
In this case, we can write $x = \prox_{\tau F}(z)$ as
\begin{equation*}
    0 = \tau U^*S^*(SUx-b) + x-z.
\end{equation*}
Multiplying by $U$ yields
\begin{equation*}
    \tau S^*b+Uz=(\tau S^*S+\Id) Ux.
\end{equation*}
By assumption, we can solve this for
\begin{equation*}
    x=U^*\inv{(\tau S^*S+\Id)}(\tau S^*b+Uz).
\end{equation*}
This shows that
\begin{equation*}
    \prox_{\tau F}(z)=U^*\inv{(\tau S^*S+\Id)}(S^*b+Uz).
\end{equation*}
In this case, \eqref{eq:tv:pdps-unitary} is practical to implement. In particular, for $U=S=\Id$, i.e., for image denoising, we have
\begin{equation}
    \label{eq:tv:quadratic-prox}
    \prox_{\tau F}(z)=\frac{1}{1+\tau}(b+z).
\end{equation}
From \cref{thm:convergence:pd_conv}, we directly obtain the following convergence result which even holds for general $A$.

\begin{theorem}
    Suppose $\tau \sigma\norm{D}^2<1$.
    Then for any starting point $(x^0, y^0) \in \R^M \times \R^{2M}$, the iterates $\{(\thisx, \thisy)\}_{k \in \N}$ generated by \eqref{eq:tv:pdps-unitary} converge to solutions $\realoptx$ and $\realopty$ of \eqref{eq:tv:problem} and \eqref{eq:tv:dual-problem}.
\end{theorem}

If $F$ is $\gamma$-strongly convex (in particular if $U=S=\Id$, where $\gamma=1$), we can apply the accelerated variant \eqref{eq:testing:pdps:forward} to obtain the iteration
\begin{algeqbox}
    \begin{equation}
        \label{eq:tv:pdps-unitary:accel}
        \left\{\begin{aligned}
                \omega_k & \defeq 1/\sqrt{1+2\gamma\tau_k},
                \quad
                \tau_{k+1} \defeq \tau_k\omega_k,
                \quad
                \sigma_{k+1} \defeq \sigma_k/\omega_k,\\
                x^{k+1} & \defeq \prox_{\tau_k F}(x^k - \tau_k D^*y^k),\\
                \overnextx &\defeq (1+\omega_k)x^{k+1}-\omega_k x^k,\\
                y^{k+1} & \defeq \proj_{\alpha \B_{\infty,2}}(y^{k} + \sigma_{k+1} D\bar x^{k+1}).
        \end{aligned}\right.
    \end{equation}
\end{algeqbox}
From \cref{thm:testing:pdps:accel}, we then obtain convergence at the faster rate $O(1/k^2)$.
\begin{theorem}
    Suppose $\tau_0 \sigma_0 \norm{D}^2<1$ and that $F$ is $\gamma$-strongly convex.
    Then for any starting point $(x^0, y^0) \in \R^M \times \R^{2M}$, the primal iterates $\{\thisx\}_{k \in \N}$ generated by \eqref{eq:tv:pdps-unitary:accel} converge to a minimizer $\realoptx$ of \eqref{eq:tv:problem} at the rate $O(1/k^2)$.
\end{theorem}

Convergence of the Lagrangian duality gap can be obtained from \cref{thm:gap:pdps} or in the strongly convex case from \cref{thm:gap:accel:pdps}.

\subsection*{Primal-dual proximal splitting for general forward operators}

If $A$ is a more complex operator, $\prox_{\frac{\lambda}{2}\norm{A\freevar-b}_2^2}$ in general cannot be computed efficiently.
To overcome this, we will split the problem in two different ways.
First, as we observed in \cref{sec:proximal:gist}, we can equivalently write \eqref{eq:tv:problem} as
\begin{equation}
    \label{eq:tv:expanded-problem}
    \min_{x \in \R^N} \tilde G(x)+\tilde F(Kx)
\end{equation}
for
\[
    \tilde G \equiv 0,
    \quad
    \tilde F(y,z) \defeq \alpha\norm{y}_1 + \frac{1}{2}\norm{z-b}_2^2,
    \quad\text{and}\quad
    Kx \defeq (Dx, Ax).
\]
We write for brevity $F_0(z) \defeq \frac{1}{2}\norm{z-b}_2^2$.
By \cref{lem:proximal:calculus}\,\cref{lem:proximal:calculus:i} and \cref{ex:proximal:rn} -- or from \eqref{eq:tv:quadratic-prox} -- for any $\gamma>0$, we have
\begin{equation*}
    \prox_{\gamma F_0}(z)
    = \prox_{\gamma\frac{1}{2}\norm{\freevar}_2^2}(z-b)+b
    = \frac{1}{1+\gamma}(z-b)+b = \frac{1}{1+\gamma}(z+\gamma b).
\end{equation*}
Hence by \cref{lem:proximal:calculus}\,\cref{lem:proximal:calculus:ii},
\begin{equation*}
    \begin{aligned}
        \prox_{\sigma F_0^*}(z)
        &= z - \sigma\,\prox_{\inv\sigma \alt F}(\inv\sigma z)
        \\
        &= z - \sigma\frac{1}{1+\inv\sigma}(\inv\sigma z + \inv\sigma b)
        \\
        &= \frac{1}{1+\sigma}(z - \sigma b).
    \end{aligned}
\end{equation*}
By \cref{lem:proximal:calculus}\,\ref{lem:proximal:calculus:iii} we thus obtain
\[
    \prox_{\sigma F_0^*}(y, z)
    =
    (\proj_{\alpha\B_{\infty,2}}(y), \prox_{\sigma F_0^*}(z))
    =
    (\proj_{\alpha\B_{\infty,2}}(y), \tfrac{1}{1+\sigma}(z - \sigma b)).
\]
Therefore the primal-dual proximal splitting method \eqref{eq:splitting:pd} for \eqref{eq:tv:expanded-problem} is given by
\begin{algeqbox}
    \begin{equation}
        \label{eq:tv:expanded-pdps}
        \left\{\begin{aligned}
                x^{k+1} & \defeq x^k - \tau[D^*y^k+A^*z^k],\\
                \overnextx &\defeq 2x^{k+1}-x^k,\\
                y^{k+1} & \defeq \proj_{\alpha \B_{\infty,2}}(y^{k} + \sigma D\bar x^{k+1}), \\
                z^{k+1} & \defeq \frac{1}{1+\sigma}(z^{k} + \sigma[A\bar x^{k+1}-b]).
        \end{aligned}        \right.
    \end{equation}
\end{algeqbox}
As before, we can apply the general convergence result from \cref{thm:convergence:pd_conv} to show that the iterates converge to a solution of the problem \eqref{eq:tv:problem}.

\begin{theorem}
    Suppose $\tau\sigma(\norm{D}^2 +\norm{A}^2) < 1$.
    For any starting point $(x^0, y^0, z^0) \in \R^{M + 2M + N}$, let the iterates $\{(\thisx, \thisy, \thisz)\}_{k \in \N}$ be generated by \eqref{eq:tv:expanded-pdps}.
    Then the primal iterates $\{\thisx\}_{k \in \N}$ converge to a minimizer of \eqref{eq:tv:problem}.
\end{theorem}

The convergence of a Lagrangian duality gap corresponding to the formulation \eqref{eq:tv:expanded-problem} can be obtained from \cref{thm:gap:pdps}.

\subsection*{Primal-dual proximal splitting with a forward step}

The dualization trick of the expanded PDPS method does not require the data term $F$ to be differentiable; we could have derived \eqref{eq:tv:expanded-pdps} for an $F(x)=F_0(Ax)$ for an arbitrary convex, possibly nonsmooth $F_0$.
It does, however, require introducing the additional variable $z$, which may come at the cost of performance.
This can be avoided for smooth $F$ by using the variant of the PDPS method with a forward step introduced in \eqref{eq:convergence:pdps:forward}. To apply it, we write \eqref{eq:tv:problem} as
\[
    \min_{x \in X} F_0(x)+E(x)+G(Kx)
\]
for
\[
    F_0 \equiv 0,
    \quad
    E(x)=\frac{1}{2}\norm{Ax-b}_2^2,
    \quad
    G(x)=\norm{\freevar}_{2,1},
    \quad\text{and}\quad
    K=D.
\]
Thus $G$ and $K$ are as in \eqref{eq:tv:pdps-unitary}; however, the primal update becomes
\begin{equation*}
    \nextx \defeq \prox_{\tau F_0}(\thisx - \tau[\grad E(\thisx) + D^*\thisy])
\end{equation*}
We thus obtain from \eqref{eq:convergence:pdps:forward} the algorithm
\begin{algeqbox}
    \begin{equation}
        \label{eq:tv:forward-pdps}
        \left\{\begin{aligned}
                \nextx &\defeq \thisx - \tau[A^*(A\thisx-b)+D^*\thisy],\\
                \overnextx &\defeq 2\nextx-\thisx,\\
                \nexty & \defeq \proj_{\alpha \B_{\infty,2}}(\thisy + \sigma D\overnextx).
        \end{aligned}        \right.
    \end{equation}
\end{algeqbox}

We have the following convergence result.

\begin{theorem}
    Suppose $1 > \norm{D}^2\tau\sigma + \frac{\tau}{2}\norm{A}^2$.
    Then for any starting point $(x^0, y^0) \in \R^{M + 2M}$, the iterates $\{(\thisx, \thisy)\}_{k \in \N}$ generated by \eqref{eq:tv:forward-pdps} converge a solution $(\realoptx, \realopty)$ of the primal-dual optimality conditions \eqref{eq:tv:primal-dual-oc}.
    In particular, the primal iterates $\{\thisx\}_{k \in \N}$ converge to a minimizer of \eqref{eq:tv:problem}.
\end{theorem}

\begin{proof}
    Since $\grad E$ is Lipschitz with constant $L=\norm{A}^2$, the claim is a direct consequence of \cref{cor:convergence:pdps:forward}.
\end{proof}

Again, convergence of the Lagrangian duality gap can be obtained from \cref{thm:gap:pdps}.
We can also apply acceleration similarly to \eqref{eq:tv:pdps-unitary:accel}, for which convergence rates can be obtained from \cref{thm:testing:pdps:accel,thm:gap:accel:pdps}.

\subsection*{Primal-dual explicit splitting}

Just like the PDPS method with a forward step, the PDES method of \eqref{eq:splitting:gist} avoids the need to introduce an additional variable.
To apply the latter, we write the problem \eqref{eq:tv:problem} in the form $\min_x F(x) + G(Kx)$ for $F(x)=\frac{1}{2}\norm{Ax-b}_2^2$ and $G(y)=\alpha\norm{y}_{1,2}$.
However, since the convergence result from \cref{thm:gist} has the restriction $\norm{K} < 1$, we rescale by  taking $G=\alpha\lambda\norm{\freevar}_{1,2}$ and $K=\lambda^{-1}D$ for some $\lambda > \norm{D}$.
Then the PDES method \eqref{eq:splitting:gist} becomes
\begin{algeqbox}
    \begin{equation}
        \label{eq:inverse:gist}
        \left\{\begin{aligned}
                \nexty &\defeq \proj_{\lambda \alpha \B_{\infty,2}}((\Id-\lambda^{-2} DD^*)\thisy + K(\thisx- A^*(A\thisx-b))),\\
                \nextx &\defeq \thisx - A^*(A\thisx-b) - \inv\lambda D^*\nexty.
        \end{aligned}        \right.
    \end{equation}
\end{algeqbox}

We have the following convergence result.

\begin{corollary}
    For any initial iterate $(x^0,y^0) \in \R^{M\times 2M}$, the sequence $\{\thisx,\inv\lambda\thisy\}_{k \in \N}$ constructed by \eqref{eq:inverse:gist} converges to a solution of the primal-dual optimality conditions \eqref{eq:tv:primal-dual-oc}.
    In particular, the primal iterates $\{\thisx\}_{k \in \N}$ converge to a minimizer of \eqref{eq:tv:problem}.
\end{corollary}

\begin{proof}
    Since $\grad F$ is Lipschitz with constant $L=\norm{A}^2$ and $\norm{K} < 1$, \cref{thm:gist} immediately yields the convergence of $\{\thisx,\inv\lambda\thisy\}_{k \in \N}$ to some $(\optx, \opty)$ satisfying $-K^* \opt y = \grad F(\opt x)$ and $K\optx \in \subdiff G^*(\opt y)$, i.e., $\opty \in \subdiff G(K\optx)$.
    Since $\norm{\freevar}_{2,1}$ is positively homogeneous,
    \begin{equation*}
        \subdiff G(K\optx)=\alpha\lambda\subdiff \norm{\freevar}_{2,1}(\inv\lambda D\optx)=\alpha\subdiff \norm{\freevar}_{2,1}(D\optx).
    \end{equation*}
    Inserting the definition of $K=\lambda^{-1}D$ and dividing by $\lambda>0$, respectively, we thus obtain that $-D^* (\lambda^{-1}\opty) = A^*(A\optx-b)$ and $(\lambda^{-1}\opty) \in \alpha\subdiff \norm{\freevar}_{2,1}(D\optx)$. Hence $(\realoptx,\realopty)\defeq (\optx,\lambda^{-1} \opty)$ satisfies \eqref{eq:tv:primal-dual-oc}.
\end{proof}

Convergence of the Lagrangian duality gap can be obtained from \cref{thm:gap:pdes}.

\subsection*{Numerical illustration}

We illustrate the performance of the various variants of the forward-backward splitting, PDPS, and PDES methods on total variation denoising and superresolution.

\begin{figure}[t!]
    \centering
    \begin{subfigure}[t]{0.325\textwidth}
        \includegraphics[width=\textwidth]{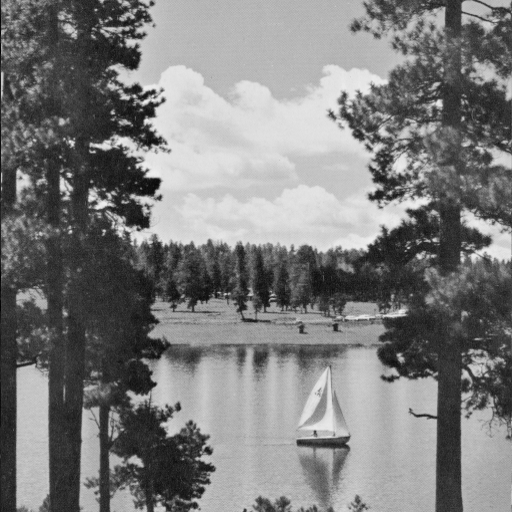}
        \caption{original}
        \label{fig:tv:denoising:reco-orig}
    \end{subfigure}
    \hfil
    \begin{subfigure}[t]{0.325\textwidth}
        \includegraphics[width=\textwidth]{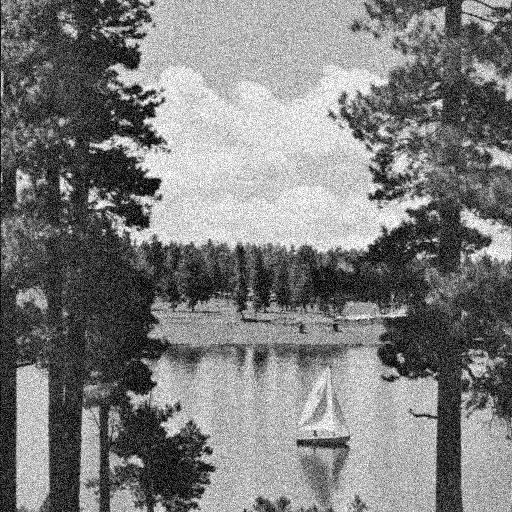}
        \caption{noisy}
        \label{fig:tv:denoising:reco-noisy}
    \end{subfigure}
    \hfil
    \begin{subfigure}[t]{0.325\textwidth}
        \includegraphics[width=\textwidth]{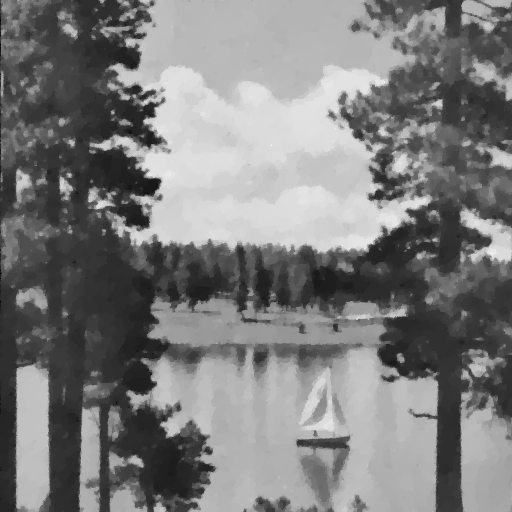}
        \caption{reconstruction}
        \label{fig:tv:denoising:reco-denoise}
    \end{subfigure}
    \caption{TV denoising data and result.}
    \label{fig:tv:denoising:reco}
\end{figure}

\begin{figure}
    \centering
    \begin{tikzpicture}
        \SetMinMax{tv_denoising_dualfb.txt}{value}{\dataDualFB}%
        \UpdMinMax{tv_denoising_dualfista.txt}{value}{\dataDualFISTA}
        \UpdMinMax{tv_denoising_pdps.txt}{value}{\dataPDPS}%
        \UpdMinMax{tv_denoising_pdps_accel.txt}{value}{\dataPDPSAccel}%
        \UpdMinMax{tv_denoising_pdps_forward.txt}{value}{\dataPDPSForward}%
        \UpdMinMax{tv_denoising_pdps_forward_accel.txt}{value}{\dataPDPSForwardAccel}%
        \begin{axis}[%
            xmode=log,
            xmin=1,
            xlabel={iteration count},
            ylabel={primal value},
            legend columns = 2,
            fixedylog = {3}{fixed,precision=0},
            ]

            \addplot [fb] table[x=iter,y=value]{\dataDualFB};
            \addlegendentry{dual FB}

            \addplot [fista] table[x=iter,y=value]{\dataDualFISTA};
            \addlegendentry{dual FISTA}

            \addplot [pdps] table[x=iter,y=value]{\dataPDPS};
            \addlegendentry{PDPS}

            \addplot [pdps accel] table[x=iter,y=value]{\dataPDPSAccel};
            \addlegendentry{accelerated PDPS}

            \addplot [pdps forward] table[x=iter,y=value]{\dataPDPSForward};
            \addlegendentry{forward PDPS}

            \addplot [pdps forward accel] table[x=iter,y=value]{\dataPDPSForwardAccel};
            \addlegendentry{accelerated forward PDPS}
        \end{axis}
    \end{tikzpicture}
    \begin{tikzpicture}
        \SetMinMax{tv_denoising_dualfb.txt}{value}{\dataDualFB}%
        \UpdMinMax{tv_denoising_dualfista.txt}{value}{\dataDualFISTA}
        \UpdMinMax{tv_denoising_pdps.txt}{value}{\dataPDPS}%
        \UpdMinMax{tv_denoising_pdps_accel.txt}{value}{\dataPDPSAccel}%
        \UpdMinMax{tv_denoising_pdps_forward.txt}{value}{\dataPDPSForward}%
        \UpdMinMax{tv_denoising_pdps_forward_accel.txt}{value}{\dataPDPSForwardAccel}%
        \begin{axis}[%
            xmode=log,
            xlabel={CPU time [s]},
            ylabel={primal value},
            legend columns = 2,
            fixedylog = {3}{fixed,precision=0},
            ]

            \addplot [fb] table[x=cputime,y=value]{\dataDualFB};
            \addlegendentry{dual FB}

            \addplot [fista] table[x=cputime,y=value]{\dataDualFISTA};
            \addlegendentry{dual FISTA}

            \addplot [pdps] table[x=cputime,y=value]{\dataPDPS};
            \addlegendentry{PDPS}

            \addplot [pdps accel] table[x=cputime,y=value]{\dataPDPSAccel};
            \addlegendentry{accelerated PDPS}

            \addplot [pdps forward] table[x=cputime,y=value]{\dataPDPSForward};
            \addlegendentry{forward PDPS}

            \addplot [pdps forward accel] table[x=cputime,y=value]{\dataPDPSForwardAccel};
            \addlegendentry{accelerated forward PDPS}
        \end{axis}
    \end{tikzpicture}
    \caption{TV denoising algorithm performance: primal function value.}
    \label{fig:tv:denoising:performance-primal}
\end{figure}
\begin{figure}
    \centering
    \begin{tikzpicture}
        \SetMinMax{tv_denoising_dualfb.txt}{dual_value}{\dataDualFB}%
        \UpdMinMax{tv_denoising_dualfista.txt}{dual_value}{\dataDualFISTA}
        \UpdMinMax{tv_denoising_pdps.txt}{dual_value}{\dataPDPS}%
        \UpdMinMax{tv_denoising_pdps_accel.txt}{dual_value}{\dataPDPSAccel}%
        \UpdMinMax{tv_denoising_pdps_forward.txt}{dual_value}{\dataPDPSForward}%
        \UpdMinMax{tv_denoising_pdps_forward_accel.txt}{dual_value}{\dataPDPSForwardAccel}%
        \begin{axis}[%
            xmode=log,
            xmin=1,
            xlabel={iteration count},
            ylabel={dual value},
            legend columns = 2,
            legend pos = north east,
            legend style={xshift=2.5ex, yshift=2.5ex},
            fixedylog = {3}{fixed,precision=0},
            ]

            \addplot [fb] table[x=iter,y=dual_value]{\dataDualFB};
            \addlegendentry{dual FB}

            \addplot [fista] table[x=iter,y=dual_value]{\dataDualFISTA};
            \addlegendentry{dual FISTA}

            \addplot [pdps] table[x=iter,y=dual_value]{\dataPDPS};
            \addlegendentry{PDPS}

            \addplot [pdps accel] table[x=iter,y=dual_value]{\dataPDPSAccel};
            \addlegendentry{accelerated PDPS}

            \addplot [pdps forward] table[x=iter,y=dual_value]{\dataPDPSForward};
            \addlegendentry{forward PDPS}

            \addplot [pdps forward accel] table[x=iter,y=dual_value]{\dataPDPSForwardAccel};
            \addlegendentry{accelerated forward PDPS}
        \end{axis}
    \end{tikzpicture}
    \begin{tikzpicture}
        \SetMinMax{tv_denoising_dualfb.txt}{dual_value}{\dataDualFB}%
        \UpdMinMax{tv_denoising_dualfista.txt}{dual_value}{\dataDualFISTA}
        \UpdMinMax{tv_denoising_pdps.txt}{dual_value}{\dataPDPS}%
        \UpdMinMax{tv_denoising_pdps_accel.txt}{dual_value}{\dataPDPSAccel}%
        \UpdMinMax{tv_denoising_pdps_forward.txt}{dual_value}{\dataPDPSForward}%
        \UpdMinMax{tv_denoising_pdps_forward_accel.txt}{dual_value}{\dataPDPSForwardAccel}%
        \begin{axis}[%
            xmode=log,
            xlabel={CPU time [s]},
            ylabel={dual value},
            legend columns = 2,
            legend pos = north east,
            legend style={xshift=2.5ex, yshift=2.5ex},
            fixedylog = {3}{fixed,precision=0},
            ]

            \addplot [fb] table[x=cputime,y=dual_value]{\dataDualFB};
            \addlegendentry{dual FB}

            \addplot [fista] table[x=cputime,y=dual_value]{\dataDualFISTA};
            \addlegendentry{dual FISTA}

            \addplot [pdps] table[x=cputime,y=dual_value]{\dataPDPS};
            \addlegendentry{PDPS}

            \addplot [pdps accel] table[x=cputime,y=dual_value]{\dataPDPSAccel};
            \addlegendentry{accelerated PDPS}

            \addplot [pdps forward] table[x=cputime,y=dual_value]{\dataPDPSForward};
            \addlegendentry{forward PDPS}

            \addplot [pdps forward accel] table[x=cputime,y=dual_value]{\dataPDPSForwardAccel};
            \addlegendentry{accelerated forward PDPS}
        \end{axis}
    \end{tikzpicture}
    \caption{TV denoising algorithm performance: dual function value.}
    \label{fig:tv:denoising:performance-dual}
\end{figure}

We start with denoising. We include in our experiments the dual forward-backward splitting \eqref{eq:tv:dualfb}, the PDPS method \eqref{eq:tv:pdps-unitary}, the forward PDPS method \eqref{eq:tv:forward-pdps}, and their accelerated variants.
We use as $b$ the noisy image shown in \cref{fig:tv:denoising:reco-noisy}, which was obtained from the original (\enquote{ground-truth}) image in \cref{fig:tv:denoising:reco-orig} by applying normally-distributed noise with mean $0$ and standard deviation $0.1$.
As the regularization parameter, we take $\alpha=0.1$; the corresponding denoised image is shown in \cref{fig:tv:denoising:reco-denoise}.
For forward-backward splitting and its inertial variant, we take $\tau=0.99/M^2$, where $M$ is an estimate of $\norm{D}$.
For the basic PDPS method and its accelerated variant we take $\tau=1.99/M$ and $\sigma=0.5/M$ to satisfy $\tau\sigma M^2 < 1$.
For the forward PDPS method and its accelerated variant we take $\tau=0.35 \cdot 2/L$ and $\sigma=0.95(1-\tau L/2)/(\tau M^2)$ to satisfy \eqref{eq:convergence:pdps:forward-stepsize}, where $L=1$ is the Lipschitz factor of $\grad F$.
Further experimental details can be found in the accompanying code \cite{nonsmoothbook-codes}.

We plot the convergence behavior in \cref{fig:tv:denoising:performance-primal} with respect to the primal functional value \eqref{eq:tv:problem}. For the primal-dual methods, we use the iterates $\this x$ to directly calculate the primal function values. For the forward-backward methods, which do not directly generate primal variables, we use the first part of the optimality conditions \eqref{eq:tv:primal-dual-oc} (with $A=\Id$) to generate $\this x$ from $\this y$. Although initially the accelerated variants seem to be slower, they eventually outperform the unaccelerated variants, in line with their better \emph{asymptotic} convergence rates. The same phenomenon can be observed in relation to the different base algorithms: Asymptotically, all algorithms converge at the same rate even though in the beginning, the simpler dual forward-backward splitting outperforms the forward PDPS method which outperforms the PDPS method.

The picture is clearer when considering convergence of the dual function values, which can be directly calculated from all iterates, and for which \cref{thm:tv:dualfb-convergence} ensures convergence for dual forward-backward splitting. Note that since all algorithms involve a dual projection step, the dual iterates are feasible, so the dual functional reduces to the strongly convex $F^*(y) = \frac12\norm{y}_2^2 + \iprod{b}{y}$.
Here, \cref{fig:tv:denoising:performance-dual} shows the expected behavior of the algorithms, with the PDPS method  outperforming the dual forward-backward splitting method and the forward PDPS method (albeit at the same asymptotic rate), and the accelerated variants clearly outperforming the base algorithms (at a higher asymptotic rate).

\begin{figure}[t!]
    \centering
    \begin{subfigure}[t]{0.325\textwidth}
        \includegraphics[width=\textwidth]{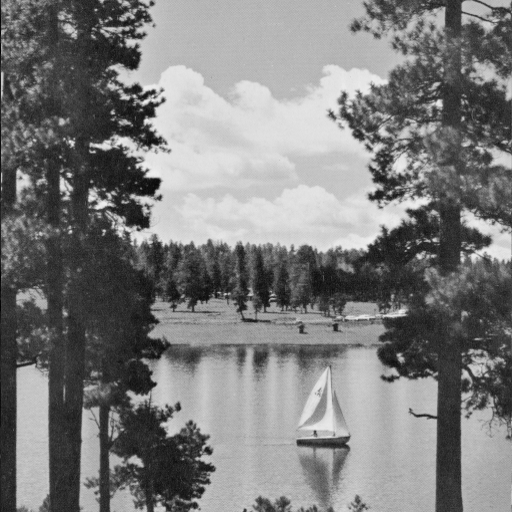}
        \caption{original}
    \end{subfigure}
    \hfil
    \begin{subfigure}[t]{0.325\textwidth}
        \includegraphics[width=\textwidth]{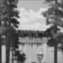}
        \caption{low-resolution data}
    \end{subfigure}
    \hfil
    \begin{subfigure}[t]{0.325\textwidth}
        \includegraphics[width=\textwidth]{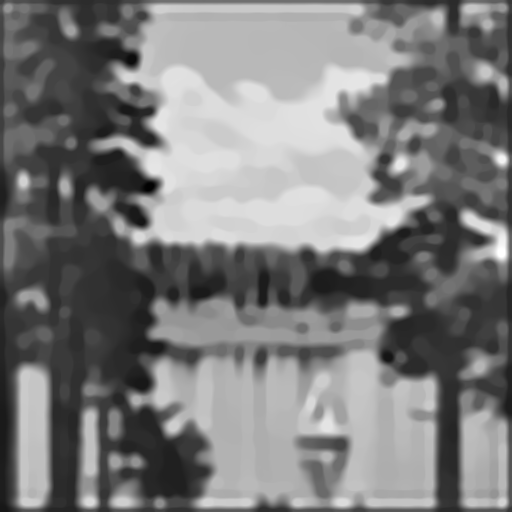}
        \caption{reconstruction}
    \end{subfigure}
    \caption{TV superresolution data and result.}
    \label{fig:tv:superresolution:reco}
    \begin{tikzpicture}
        \begin{axis}[%
            xmode=log,
            xmin=1,
            ymode=log,
            axis x line*=bottom,
            axis y line*=left,
            legend pos = south west,
            xlabel={iteration count},
            ylabel={function value},
            ]

            \addplot [pdps forward] table[x=iter,y=value]{tv_superresolution_pdps_forward.txt};
            \addlegendentry{PDPS with forward step}

            \addplot [pdps] table[x=iter,y=value]{tv_superresolution_pdps_expanded.txt};
            \addlegendentry{PDPS expanded}

            \addplot [pdes] table[x=iter,y=value]{tv_superresolution_pdes.txt};
            \addlegendentry{PDES}
        \end{axis}
    \end{tikzpicture}
    \begin{tikzpicture}
        \begin{axis}[%
            xmode=log,
            ymode=log,
            axis x line*=bottom,
            axis y line*=left,
            legend pos = south west,
            xlabel={time [s]},
            ylabel={function value},
            ]

            \addplot [pdps forward] table[x=cputime,y=value]{tv_superresolution_pdps_forward.txt};
            \addlegendentry{PDPS with forward step}

            \addplot [pdps] table[x=cputime,y=value]{tv_superresolution_pdps_expanded.txt};
            \addlegendentry{PDPS expanded}

            \addplot [pdes] table[x=cputime,y=value]{tv_superresolution_pdes.txt};
            \addlegendentry{PDES}
        \end{axis}
    \end{tikzpicture}
    \caption{TV superresolution algorithm performance.}
    \label{fig:tv:superresolution:performance}
\end{figure}

For the superresolution demonstration, we consider the forward PDPS method \eqref{eq:tv:forward-pdps}, the expanded PDPS method \eqref{eq:tv:expanded-pdps}, and the PDES method \eqref{eq:inverse:gist}.
In this experiment, the operator $A \in \linear(\R^{512^2}; \R^{64^2})$ performs convolution with a Gaussian kernel (standard deviation $\sigma=5$ on the domain $\Omega=[0, 512]^2$) followed by subsampling by factor of $8$.
We illustrate the data and the reconstruction in \cref{fig:tv:superresolution:reco}.
The low-resolution data is obtained from the original image by applying $A$ and adding normally-distributed noise of mean $0$ and standard deviation $0.001$.
As the regularization parameter, we take $\alpha=0.0001$.
For the forward PDPS method, we take $\tau=0.95 \cdot 2/L$ and $\sigma=0.95(1-\tau L/2)/(\tau M^2)$ to satisfy \eqref{eq:convergence:pdps:forward-stepsize}, where $L$ is an upper estimate of the Lipschitz factor of $\grad F$, i.e., of $\norm{A}^2$.
For the expanded PDPS method, we take $\tau=1.9/\sqrt{M^2+L}$ and $\sigma=0.5/\sqrt{M^2+L}$ to satisfy $\tau\sigma(\norm{D}^2 +\norm{A}^2) < 1$ via $\tau\sigma(M^2 + L) < 1$.
The PDES method has no step length parameters.

\enlargethispage{1cm}
We illustrate the convergence behavior in \cref{fig:tv:superresolution:performance}.
As we can see, the expanded variant of the PDPS method is slower than the other two algorithms that do not introduce additional variables. Moreover, the accelerated algorithms eventually outperform all the unaccelerated variants. The PDPS method with forward step is somewhat faster than the PDES method.

\chapter{Optimal control with constraints}\label{chap:control}

We now illustrate the applications of the theory of \cref{part:convex,part:nonconvex} in infinite-dimensional spaces, in particular function spaces.
Specifically, we consider \term[problem!optimal control]{optimal control problems}, where the solution of a (partial) differential equation -- the \term{state} -- is sought to be brought as close as possible to a desired state by adjusting a relevant \term{control}.
Typically this control is the right-hand side, boundary conditions, or coefficients of the differential equation.
Optimal control problems occur in a wide variety of applications such as autonomous vehicles, process engineering, and optimal design; they are also closely related to inverse problems for partial differential equations. Typically, this involves minimizing a weighted sum of a \term[term!tracking]{tracking term} involving the state and a \term[term!control cost]{control cost} involving the control; these are linked through the differential equation as an equality constraint, and hence this is also known as \term[problem!optimization!PDE-constrained]{PDE-constrained optimization}. Using the implicit function theorem, one can use this constraint to define a \term[mapping!control-to-state]{control-to-state mapping}; much of optimal control theory is concerned with analyzing the properties (in particular regarding differentiability) of this mapping, especially for (systems of) time-dependent and/or nonlinear equations or controls appearing as the coefficients. On these and other issues, we refer the reader to the seminal monograph \cite{Lions:1971}, to the standard textbook \cite{Troeltzsch}, as well as to \cite{Hinze2009,DelosReyes} in particular regarding applications and numerical methods.

Here we focus on dealing with optimal control problems where either the tracking term or the control costs are nonsmooth, which allows imposing additional structure on the optimal state or control. Correspondingly, such problems have received increasing attention in recent years. To avoid unnecessary technical difficulties, we restrict ourselves to the simplest possible partial differential equation: the Poisson equation with homogeneous boundary conditions and the control appearing as a right-hand side. We briefly introduce the required notation and refer to, e.g., \cite{Troeltzsch} for details and proofs of the claimed properties. Let $\Omega\subset\R^d$ be a bounded domain with Lipschitz boundary. We then introduce for $k\in\N$ and $1<p<\infty$ the \term[space!Sobolev]{Sobolev space}
\begin{equation*}
    W^{k,p}(\Omega)\defeq \setof{v\in L^p(\Omega)}{D^\alpha v\in L^p(\Omega) \text{ for all }|\alpha|\leq k},
\end{equation*}
where $D^\alpha v$ is the \term[derivative!weak]{weak derivative} of $v$ of order $|\alpha|$. These are Banach spaces with the natural norm; for $p=2$, $H^k(\Omega)\defeq W^{k,p}(\Omega)$ is a Hilbert space. Under the assumptions on the domain $\Omega$, we have the continuous embeddings
\begin{align*}
    &W^{k,p}(\Omega) \hookrightarrow L^q(\Omega)  && \text{ for } 1\leq q\leq \frac{dp}{d-kp} \quad(\defeq\infty \text{ if $kp\geq d$}),\\
    &W^{k,p}(\Omega) \hookrightarrow C(\overline\Omega)  && \text{ for } kp>d;
\end{align*}
see, e.g., \cite[Theorem~7.1]{Troeltzsch}. Furthermore, the embedding $W^{k,p}(\Omega)\hookrightarrow L^p(\Omega)$ is compact for every $k\in \N$ and $1<p<\infty$; see, e.g., \cite[Theorem~7.4]{Troeltzsch}. In particular, weakly convergent sequences in $W^{k,p}(\Omega)$ for $k \ge 2$ converge strongly in $L^p(\Omega)$. Finally, we denote by $W^{k,p}_0(\Omega)$ the closure of $C^\infty_0(\overline\Omega)$ with respect to the $W^{k,p}$-norm, whose elements have vanishing trace on the boundary of $\Omega$.

We now consider for given $u\in L^2(\Omega)$ the \term[formulation, weak]{weak formulation} of the Poisson equation $-\Delta y = u$ with homogeneous boundary condition, i.e., we look for $y\in H^1_0(\Omega)$ satisfying
\begin{equation}\label{eq:control:poisson}
    \int_\Omega \nabla y(x) \cdot \nabla v(x) \,dx = \int_\Omega u(x)v(x)\,dx \qquad\text{for all }v\in H^1_0(\Omega).
\end{equation}
Under the assumptions on $\Omega$, this equation admits a unique solution $y\in H^1_0(\Omega)$ which depends continuously on $u$; see, e.g., \cite[Theorem~2.4]{Troeltzsch}.
This allows defining a linear bounded control-to-state mapping $S:L^2(\Omega)\to L^2(\Omega)$ (which is even compact since the range $\range S\subset H^1_0(\Omega)$ embeds compactly into $L^p(\Omega)$ for any $1<p<\infty$).
If $d\leq 3$ and $\Omega\subset\R^d$ is convex, we even have $y\in H^2(\Omega)\hookrightarrow C(\overline\Omega)$; see \cite[Theorem~3.2.1.2]{Grisvard:2011}.

We will also need the adjoint $S^*:L^2(\Omega)\to L^2(\Omega)$ of $S$. Using either the implicit function theorem or formal Lagrange multiplier calculus, we can characterize $p\defeq S^*h\in L^2(\Omega)$ for given $h\in L^2(\Omega)$ as the unique solution to the \term[equation, adjoint]{adjoint equation}
\begin{equation}\label{eq:control:adjoint}
    \int_\Omega \nabla w(x) \cdot \nabla p(x) \,dx = \int_\Omega w(x)h(x)\,dx \qquad\text{for all }w\in H^1_0(\Omega);
\end{equation}
see, e.g., \cite[Lemma~2.24, Chapter~2.10]{Troeltzsch} or \cite[Chapter~1.6]{Hinze2009}. This implies that $\range S^*\subset H^1_0(\Omega)\hookrightarrow L^p(\Omega)$ for any $1<p<\infty$ as well.

\section{Control constraints}\label{chap:control:cconstraints}

We start with the simplest nonsmooth optimal control problems: quadratic control problems with pointwise constraints on the control or state. Although these problems can be treated by well-known standard methods of constrained smooth optimization (cf., e.g., \cite[Chapters~2 and 6.2]{Troeltzsch}), they serve well to illustrate the application of the abstract results of \cref{part:convex}.

\subsection*{Problem description}

Let $y^d\in L^2(\Omega)$ be a desired state and $\alpha>0$ as well as $a,b\in\R$ with $a>b$ be given. We then consider the \enquote{mother problem}
\begin{equation*}
    \begin{aligned}
        &\min_{u\in L^2(\Omega),y\in H^1_0(\Omega)} \frac12\norm{y-y^d}_{L^2(\Omega)}^2 + \frac\alpha2 \norm{u}_{L^2(\Omega}^2\\
        &\text{subject to \eqref{eq:control:poisson}} \qquad\text{and}\qquad a\leq u(x) \leq b \quad\text{for almost every }x\in \Omega.
    \end{aligned}
\end{equation*}
Introducing the \term[set!admissible]{admissible set}
\begin{equation*}
    \Uad \defeq \setof{u\in L^2(\Omega)}{a \leq u(x) \leq b\quad \text{for almost every }x\in \Omega}
\end{equation*}
and using the control-to-state-mapping $S:L^2(\Omega)\to L^2(\Omega)$, $u\mapsto y$ solving \eqref{eq:control:poisson}, introduced above, we can write this problem in \term[problem!reduced form]{reduced form} as
\begin{equation*}
    \min_{u\in \Uad} \frac12\norm{Su-y^d}_{L^2(\Omega)}^2 + \frac\alpha2\norm{u}_{L^2(\Omega)}^2.
\end{equation*}
To apply the general theory of the previous parts, we write this as $\min_{u\in L^2(\Omega)} J(u)$ for $J = F+G$ with
\begin{align*}
    F(u) &\defeq  \frac12\norm{Su-y^d}_{L^2(\Omega)}^2 + \frac\alpha2\norm{u}_{L^2(\Omega)}^2,\\
    G(u) &\defeq \delta_{\Uad}(u).
\end{align*}

\subsection*{Existence}
Since $S$ is linear and bounded (and hence weakly continuous) and the norm is weakly lower semicontinuous by \cref{cor:variation:norm} and convex, it follows from \cref{lem:variation:wlsc,lem:convex:func} that $F$ is weakly lower semicontinuous and convex; it is even strictly convex due to the control costs. Furthermore, $\dom F = L^2(\Omega)$ since $S$ is well-defined on this space.
Similarly, it can be shown that $\Uad\subset L^2(\Omega)$ is nonempty, closed, convex, and bounded and thus $G$ is proper, lower semicontinuous, convex, and coercive by \cref{lem:variation:indicator}. We thus immediately obtain from \cref{thm:convex:existence} the existence of a unique optimal control $\bar u\in \Uad$ as well as a corresponding optimal state $\bar y \defeq  S\bar u\in H^1_0(\Omega)$.

\subsection*{Optimality conditions}
To derive optimality conditions, we apply the Fermat principle as well as the calculus rules from \cref{chap:subdiff}. Although $\dom G = \Uad\subset L^2(\Omega)$ does \emph{not} contain any interior points, we have $\dom F=L^2(\Omega)$ and hence we can still apply the sum rule from \cref{thm:subdiff:sum}. In fact, since the squared norm in the Hilbert space $L^2(\Omega)$ (which we always identify with its dual via the Fréchet--Riesz \cref{thm:frechetriesz}) is Fréchet differentiable, we obtain using the chain rule from \cref{thm:frechet_chain} that
\begin{equation*}
    \nabla F(u) = S^*(Su - y^d) + \alpha u.
\end{equation*}
Using \cref{thm:subdiff:sum,thm:convex:gateaux,lem:convex:normalcone} and introducing the adjoint state $\bar p\in H^1_0(\Omega)$, we thus arrive at the primal-dual optimality conditions\footnote{If the control-to-state mapping $S$ is nonlinear but continuously differentiable, we can proceed in exactly the same fashion by using \cref{thm:clarke:sum,thm:clarke:frechet,thm:clarke:convex} instead to arrive at \eqref{eq:control:constraint:varineq} with $S'(\bar u)^*$ in place of~$S^*$.}
\begin{equation}\label{eq:control:constraint:varineq}
    \left\{
        \begin{aligned}
            &\bar p  = S^*(S\bar u -y^d),\\
            &\iprod{\bar p + \alpha \bar u}{u-\bar u}_{L^2(\Omega)} \geq 0\quad\text{for all }u\in \Uad,
        \end{aligned}
    \right.
\end{equation}
where the second relation is often called a \term[inequality!variational]{variational inequality} for the optimal control; cf.~\cite[Theorem~2.25]{Troeltzsch}.
This relation, which is the explicit form of $-\bar p - \alpha \bar u\in \partial\delta_{\Uad}(\bar u)$, can by \cref{lem:proximal:subdiff,ex:proximal:hilbert}\,\ref{ex:proximal:hilbert:iii} be written equivalently for any $\gamma>0$ as
\begin{equation}\label{eq:control:prox}
    \bar u = \prox_{\gamma\delta_{\Uad}} \left(\bar u + \gamma (-\bar p - \alpha \bar u)\right).
\end{equation}
Using the special choice $\gamma = \alpha^{-1}$ in the first expression as well as the pointwise characterization of proximal mappings on $L^2(\Omega)$ from \cref{lem:lebesgue:proximal} together with \cref{ex:proximal:reell}\,\ref{ex:proximal:reell:iii}, we obtain the well-known \term[formula, projection]{projection formula}
\begin{equation}\label{eq:control:constraint:proj}
    \bar u(x) = \proj_{[a,b]}\left(-\frac1\alpha\bar p(x)\right) =
    \begin{cases}
        a & \text{if }-\tfrac1\alpha\bar p(x)<a,\\
        -\tfrac1\alpha p(x) & \text{if }-\tfrac1\alpha\bar p(x)\in[a,b],\\
        b & \text{if }-\tfrac1\alpha\bar p(x)>b;
    \end{cases}
\end{equation}
cf. \cite[Theorem~2.28]{Troeltzsch}.

\begin{remark}
    The relation \eqref{eq:control:constraint:proj} could also have been obtained by recognizing that $G(u) + \frac\alpha2\norm{u}_{L^2(\Omega)}^2 = (G^*_\alpha)^*$ by \cref{thm:moreau:conjugate}, where $G^*_\alpha$ is the Moreau envelope of $G^*$. We therefore obtain via \cref{thm:moreau-yosida}
    \begin{equation*}
        \left\{
            \begin{aligned}
                \bar p &= S^*(S\bar u-y^d),\\
                \bar u &= (\partial G^*)_\alpha (- \bar p),
            \end{aligned}
        \right.
    \end{equation*}
    where $(\partial G^*)_\alpha$ is the Yosida approximation of $\partial G^*$.
    Using its definition \eqref{eq:proximal:yosida} together with \cref{lem:proximal:calculus}\,\ref{lem:proximal:calculus:ii}, it is straightforward to verify that the second relation is in fact equivalent to \eqref{eq:control:constraint:proj}.
\end{remark}

\subsection*{Explicit splitting methods}

Since $F$ and $G$ are proper, convex, and lower semicontinuous, and $F$ is Fréchet differentiable with Lipschitz continuous gradient (since $\nabla F(u)$ is affine, it is globally Lipschitz with constant $L\defeq \norm{S^*S+\alpha \Id}_{\linear(L^2(\Omega);L^2(\Omega))} = \norm{S}_{\linear(L^2(\Omega);L^2(\Omega))}^2+\alpha$), the optimal control $\bar u$ can be computed using the explicit splitting method \eqref{eq:convergence:fb}. In our specific instance, this becomes the \term[method!projected gradient]{projected gradient method}: Choose $u^0\in L^2(\Omega)$ (e.g., $u^0=0$) and $\tau < 2L^{-1}$ and compute for $k=0,\dots$
\begin{algeqbox}
    \begin{equation}\label{eq:control:projgradient}
        \left\{
            \begin{aligned}
                y^{k+1} &\defeq Su^k &&\text{by solving \eqref{eq:control:poisson}},\\
                p^{k+1} &\defeq S^*(y^{k+1}-y^d) &&\text{by solving \eqref{eq:control:adjoint} for $h=y^{k+1}-y^d$},\\
                u^{k+1} &\defeq \proj_{[a,b]}\left((1-\tau\alpha)u^k - \tau p^{k+1}\right) &&\text{almost everywhere}.
            \end{aligned}
        \right.
    \end{equation}
\end{algeqbox}
By \cref{thm:convergence:fb}, we then have $u^k\weakto \bar u$ in $L^2(\Omega)$. (Since $G$ is not strongly convex, we do not obtain any rates.)

We can also apply the acceleration strategies from \cref{chap:meta}.
Specifically, the \term[method!projected gradient!inertial]{inertial projected gradient method} for $z^0=u^0\in L^2(\Omega)$, $\tau > 0$, and $\lambda_0 = 1$ consists in computing for $k=0,\dots$
\begin{algeqbox}
    \begin{equation}\label{eq:control:inertial}
        \left\{
            \begin{aligned}
                y^{k+1} &\defeq Sz^k &&\text{by solving \eqref{eq:control:poisson}},\\
                p^{k+1} &\defeq S^*(y^{k+1}-y^d) &&\text{by solving \eqref{eq:control:adjoint} for $h=y^{k+1}-y^d$},\\
                u^{k+1} &\defeq \proj_{[a,b]}\left((1-\tau\alpha)u^k - \tau p^{k+1}\right) &&\text{almost everywhere},\\
                \lambda_{k+1} &\defeq 2\left(1+\sqrt{1+4\lambda_k^{-2}}\right), && \beta_{k+1} = \lambda_{k+1}(\lambda_k^{-1}-1), \\
                z^{k+1} &\defeq (1+\beta_{k+1}) u^{k+1} - \beta_{k+1} u^k.
            \end{aligned}
        \right.
    \end{equation}
\end{algeqbox}
By \cref{thm:meta:inertia:fb}, we obtain the convergence of the function values $J(\tilde u^k)\to J(\bar u)$ at the rate $O(1/k^2)$ as $k\to \infty$ (for the \emph{nonergodic} sequence).

Similarly, we could also derive the \term[method!projected gradient!over-relaxed]{over-relaxed projected gradient method} from \eqref{eq:meta:overrelax:fb}; however, since this method does not show any benefit over the projected gradient method for this problem, this is left as an exercise to the reader.

Instead, we will consider an alternative splitting. Since the maximal step length is constrained by the Lipschitz constant of $F$, it is beneficial to include as many parts of the functional as possible in the proximal point mapping. We thus turn to the splitting
\begin{align*}
    F(u) &\defeq  \frac12\norm{Su-y^d}_{L^2(\Omega)}^2,\\
    G_\alpha(u) &\defeq \delta_{\Uad}(u) + \frac\alpha2\norm{u}_{L^2(\Omega)}^2.
\end{align*}
To compute $\prox_{\gamma G_\alpha}$, we first observe that completing the square yields the scalar equality
\begin{equation*}
    \frac{1}{2\gamma}(z-t)^2 + \frac\alpha2 z^2 = \frac{1+\alpha\gamma}{\gamma}\left(z-\frac{1}{1+\alpha\gamma}t\right)^2 + \frac{\gamma}{1+\alpha\gamma}t^2.
\end{equation*}
By ignoring the constant term, we hence have pointwise almost everywhere that for all $\gamma>0$ and $v\in L^2(\Omega)$,
\begin{equation*}
    \begin{aligned}
        [\prox_{\gamma G_\alpha}(v)](x) &= \argmin_{z\in [a,b]} \frac{1}{2\gamma}(z-v(x))^2 + \frac\alpha2 z^2 \\
        &= \argmin_{z\in [a,b]} \frac{1+\alpha\gamma}{\gamma}\left(z-\frac{1}{1+\alpha\gamma}v(x)\right)^2 \\
        &= \proj_{[a,b]}\left(\frac{1}{1+\alpha\gamma}v(x)\right).
    \end{aligned}
\end{equation*}
In place of \eqref{eq:control:prox}, we thus have the equivalent optimality conditions
\begin{equation*}
    \bar u = \prox_{\gamma\delta_{\Uad}}\left(\frac{1}{1+\alpha\gamma}\left(\bar u-\gamma\bar p\right)\right).
\end{equation*}
From this, we obtain the corresponding (inertial) explicit splitting method for $G_\alpha$ by replacing the update for $u^{k+1}$ in \eqref{eq:control:projgradient} (or \eqref{eq:control:inertial}) by
\begin{algeqbox*}
    \begin{equation*}
        u^{k+1} \defeq \proj_{[a,b]}\left(\frac{1}{1+\tau\alpha}\left(u^k - \tau p^{k+1}\right)\right) \qquad\text{almost everywhere},
    \end{equation*}
\end{algeqbox*}
where $\tau$ now is only constrained by the smaller Lipschitz constant $L = \norm{S}_{\linear(L^2(\Omega);L^2(\Omega))}^2$, allowing larger steps.

In addition, since $G_\alpha$ is now strongly convex, we even get from \cref{thm:testing:fb} strong convergence of $u^{k}$ at a linear rate.

\subsection*{Semismooth Newton method}

Using again the specific choice $\gamma=\alpha^{-1}$ and the definition of the adjoint state $\bar p$, we can write
\eqref{eq:control:prox} as the nonsmooth equation $H(\bar u) = 0$ for
\begin{equation}\label{eq:control:nseq}
    H:L^2(\Omega)\to L^2(\Omega),\qquad   H(u) = u - \proj_{\Uad}\left(-\frac1\alpha S^*(Su-y^d)\right).
\end{equation}
Since $\range S^*\subset H^1_0(\Omega)\hookrightarrow L^p(\Omega)$ for any $p>2$ and $y^d\in L^2(\Omega)$, it follows from \cref{ex:semismooth:l2}\,\ref{ex:semismooth:l2:box} together with the chain rule \cref{thm:newton:chain} (since both $D_N\proj_{\Uad}$ and $S^*S$ are clearly uniformly bounded) that $H$ is Newton differentiable with a Newton derivative whose application to any $\delta u\in L^2(\Omega)$ is given pointwise almost everywhere by
\begin{equation*}
    [D_N H(u) \delta u](x) = \delta u(x) + \frac1\alpha \1_{[a,b]}\left(-\frac1\alpha[S^*Su](x)\right) [S^*S\delta u](x),
\end{equation*}
where $\1_{[a,b]}(t) = 1$ for $t\in [a,b]$ and $0$ else.
Under the usual regularity assumption, \cref{thm:newton:superlinear} thus guarantees that for any $u^0\in L^2(\Omega)$, the semismooth Newton iteration
\begin{algeqbox}
    \begin{equation}\label{eq:control:semismooth}
        u^{k+1} \defeq u^k - D_N H(u^k)^{-1}H(u^k)
    \end{equation}
\end{algeqbox}
is locally superlinearly convergent. The properties of $D_N H(u)$ also imply that the Newton step \eqref{eq:control:semismooth} can be solved efficiently using a suitable matrix-free Krylov space method (where for each Krylov iteration, one needs to solve two partial differential equations to apply $S$ and $S^*$, followed by setting the result to zero almost everywhere where $u^k(x)\notin[a,b]$).

\begin{figure}[t!]
    \begin{subfigure}[t]{0.49\textwidth}
        \includegraphics[width=\textwidth]{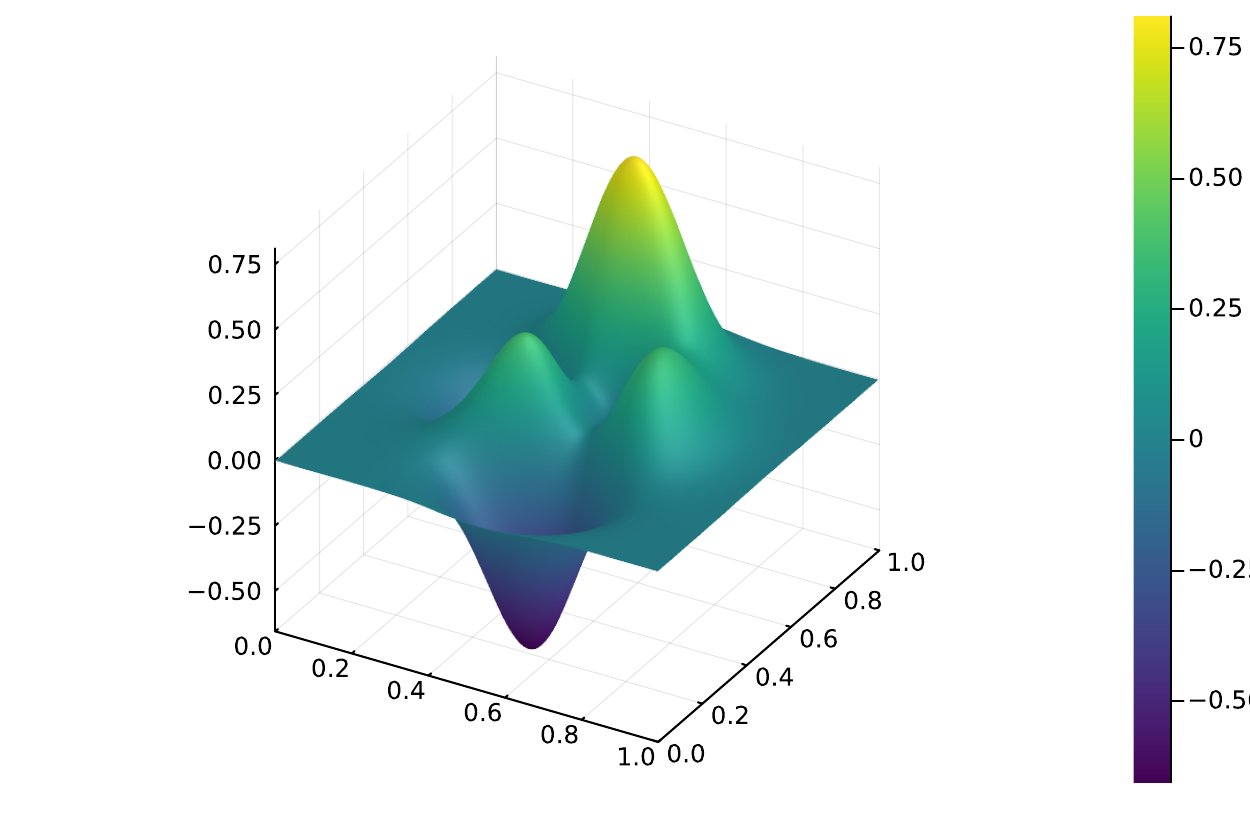}
        \caption{target $y^d$}
    \end{subfigure}
    \\
    \begin{subfigure}[t]{0.49\textwidth}
        \includegraphics[width=\textwidth]{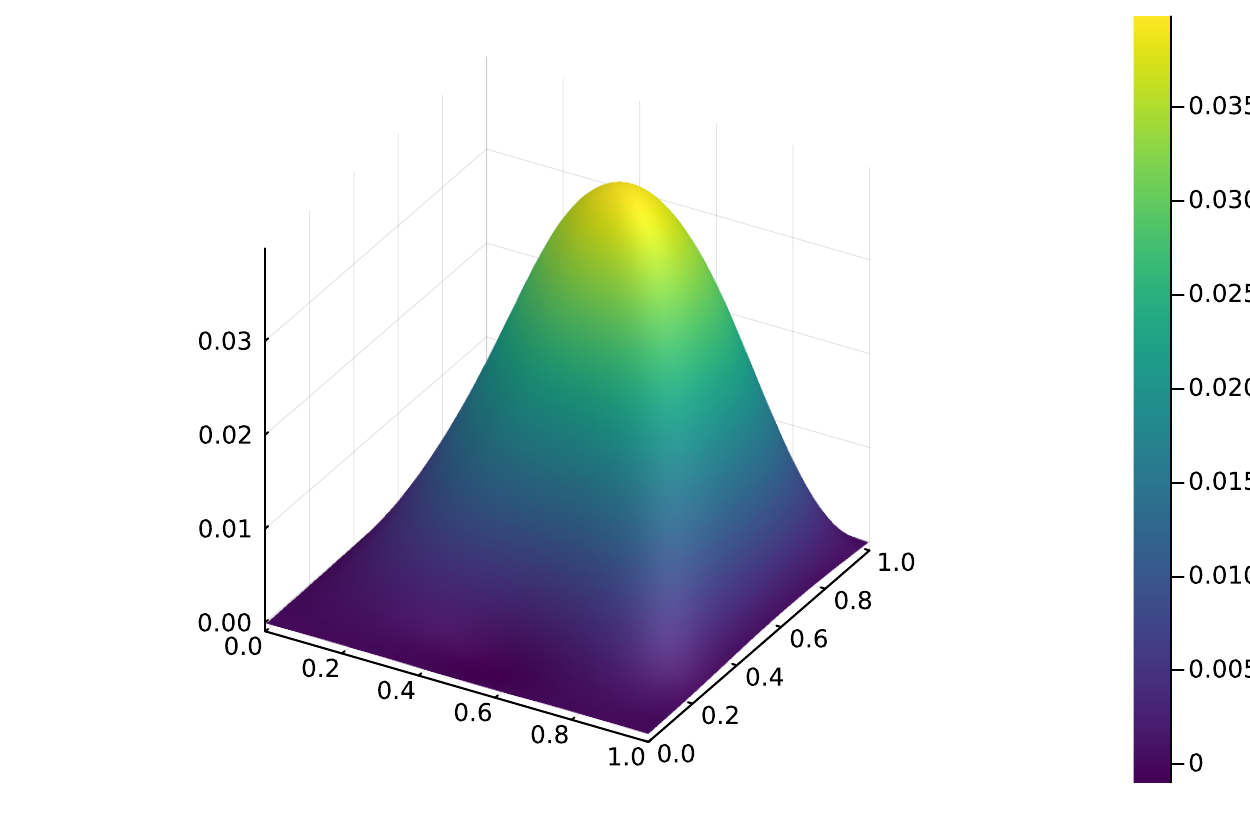}
        \caption{optimal state $\bar y$}
    \end{subfigure}
    \hfil
    \begin{subfigure}[t]{0.49\textwidth}
        \includegraphics[width=\textwidth]{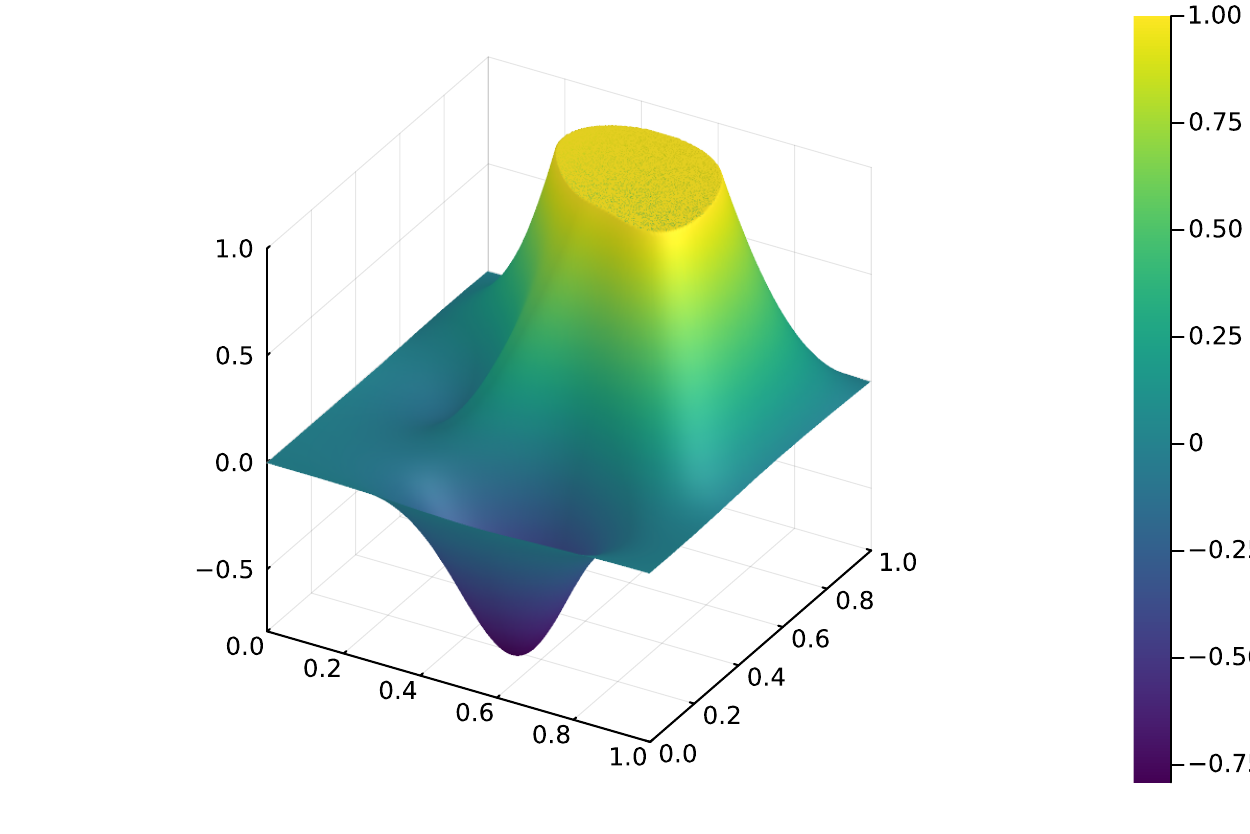}
        \caption{optimal control $\bar u$}
    \end{subfigure}
    \caption{Control constraints: target and optimal control and state.}
    \label{fig:control:control-constraints-data}
\end{figure}

We indicate the performance of the projected gradient method, the explicit splitting method with $G_\alpha$, its inertial variant, and the semismooth Newton (SSN) method for the control constraints problem with 
admissible set $[a, b]=[-1,1]$, control cost parameter $\alpha=0.005$, and target
\begin{multline}
    y^d(x_1,x_2) = \frac3{10}(4-6x_1)^2 e^{-(6x_1-3)^2 - (6x_2-2)^2} \\
   - \left(\frac{1}{5}(6x_1-3) - (6x_1-3)^3 - (6x_2-3)^5\right)e^{-(6x_1-3)^2-(6x_2-3)^2} \\
   - \frac1{30}e^{-(6x_1-2)^2 - (6x_2-3)^2};
\end{multline}
see \cref{fig:control:control-constraints-data}, which also shows the corresponding computed optimal control and state. Here and in the following, variables are discretized to a $N \times N$ grid for $N=256$.
For the splitting methods we take $\tau = 0.9/L^2$, where $L$ is an estimate of $\norm{S}_{\linear(L^2(\Omega);L^2(\Omega))}$.
More details can again be found in the accompanying code \cite{nonsmoothbook-codes}.

\begin{figure}[t!]
    \begin{tikzpicture}
        \begin{axis}[%
            axis x line*=bottom,
            axis y line*=left,
            xlabel={iteration count},
            ylabel={residual},
            ymode=log,
            legend pos = north east,
            legend style={yshift=-3ex,xshift=6ex},
            ]

            \addplot [fb] table[x=iter,y=residual]{control_constraints_fb.txt};
            \addlegendentry{explicit splitting ($G_\alpha$)}

            \addplot [projgrad] table[x=iter,y=residual]{control_constraints_projgrad.txt};
            \addlegendentry{projected gradient}

            \addplot [fista] table[x=iter,y=residual]{control_constraints_fb_inertia.txt};
            \addlegendentry{inertial explicit splitting ($G_\alpha$)}

            \addplot [ssn] table[x=iter,y=residual]{control_constraints_ssn_cg.txt};
            \addlegendentry{SSN}

        \end{axis}
    \end{tikzpicture}
    \begin{tikzpicture}
        \begin{axis}[%
            axis x line*=bottom,
            axis y line*=left,
            xlabel={time [s]},
            ylabel={residual},
            ymode=log,
            xmode=log,
            legend pos = south west,
            ]

            \addplot [fb] table[x=cputime,y=residual]{control_constraints_fb.txt};
            \addlegendentry{explicit splitting ($G_\alpha$)}

            \addplot [projgrad] table[x=cputime,y=residual]{control_constraints_projgrad.txt};
            \addlegendentry{projected gradient}

            \addplot [fista] table[x=cputime,y=residual]{control_constraints_fb_inertia.txt};
            \addlegendentry{inertial explicit splitting ($G_\alpha$)}

            \addplot [ssn] table[x=cputime,y=residual]{control_constraints_ssn_cg.txt};
            \addlegendentry{SSN}

        \end{axis}
    \end{tikzpicture}
    \caption{Algorithm performance for the control constraints example.
    We plot the residual $\norm{H(u^k)}$ for $H$ given by \eqref{eq:control:nseq}.}
    \label{fig:control:control-constraints-performance}
\end{figure}
\begin{figure}[t!]
    \begin{tikzpicture}
        \begin{axis}[%
            xmin=0,
            xmax=5,
            xtick = {0,1,...,5},
            axis x line*=bottom,
            axis y line*=left,
            xlabel={iteration count (SSN)},
            ylabel={residual},
            legend columns = 2,
            ymode=log,
            cycle list name=primaldual,
            ]

            \pgfplotsinvokeforeach{32, 64, 128, 256, 512}{
                \addplot+[line width=1pt] table[
                    x = iter,
                    y = residual
                ]{control_constraints_dims_ssn_cg_N#1.txt};
                \addlegendentry{$N={#1}$}
            }

        \end{axis}
    \end{tikzpicture}
    \begin{tikzpicture}
        \begin{axis}[%
            xmin=0,
            xmax=100,
            axis x line*=bottom,
            axis y line*=left,
            xlabel={iteration count (explicit splitting)},
            ylabel={residual},
            legend columns = 2,
            ymode=log,
            cycle list name=primaldual,
            ]

            \pgfplotsinvokeforeach{32, 64, 128, 256, 512}{
                \addplot+[line width=1pt] table[
                    x = iter,
                    y = residual
                ]{control_constraints_dims_fb_N#1.txt};
                \addlegendentry{$N={#1}$}
            }

        \end{axis}
    \end{tikzpicture}
    \caption{Control constraints: SSN performance versus dimension $N$.
        We plot the norm of residual $H(u^k)$.}
    \label{fig:controlconstr:N}
\end{figure}
As function values are not meaningful in this problem (since they can be infinite for infeasible controls), we compare the residual norm $\norm{H(u^k)}_{\linear(L^2(\Omega); L^2(\Omega))}$ for $H$ given by \eqref{eq:control:nseq}; these are shown in \cref{fig:control:control-constraints-performance}.
Inertial explicit splitting turns out to be the slowest algorithm; but this is not surprising since it only has a $O(1/k^2)$ rate of convergence, while for this strongly convex problem, the explicit splitting method has linear convergence, and the SSN method has superlinear convergence. The projected gradient method shows linear convergence as well, albeit with a smaller constant than the explicit splitting method with $G_\alpha$.
Theoretically, indeed, the iterates of both methods converge linearly due to the strongly convex regularization term $\frac{\alpha}{2}\norm{u}_{L^2(\Omega)}^2$. For the explicit splitting method with $G_\alpha$, this is a direct consequence of \cref{thm:testing:fb}. For the projected gradient method, where the strongly convex term is in $F$, we would need to adapt the proof to use \cref{cor:smoothness:three-point:sc} in place of \cref{cor:smoothness:three-point}.
That the projected gradient method is slightly slower than the explicit splitting method with $G_\alpha$ can be attributed to the fact that when the regularization term is included in $F$, the Lipschitz constant $L$ of $F$ is higher, and consequently the step length parameter $\tau$ smaller. When the proximal step can be easily calculated, it is often more efficient to do more in the proximal step and less in the gradient step, as the former does not constrain the step length parameter. Indeed, taking iteration-dependent step length parameters $\tau_k \upto \infty$, the plain proximal point method converges by \cref{thm:testing:prox} superlinearly for strongly convex objectives. Of course, its steps can be very expensive.

As in the $\ell^1$ fitting example, even though each SSN iteration involves solving a large indefinite system, the total computational time for getting the residual to machine precision is still lower than for the first-order splitting methods. Conversely, the (non-accelerated) splitting methods are faster in achieving a higher tolerance of about $10^{-6}$ and hence again may be preferable if high accuracy is not required.

Finally, as the convergence of these methods was shown on the infinite-dimensional level, it can be expected that the number of iterations required to solve optimality conditions for discretizations of the problem is independent of the fineness of the discretization. This beneficial property is referred to as \term[independence, mesh]{mesh independence}; see, e.g., \cite{Hintermuller:2004b} for its proof for a semismooth Newton method.
We numerically indicate the dimension independence of both the SSN and explicit splitting methods in \cref{fig:controlconstr:N}.

\section{State constraints}

\subsection*{Problem description}

There are also occasions when one wishes to put pointwise bounds on the state, for example when looking for optimal heat sources to achieve on average a comfortable temperature in a room without risking hot spots of a dangerous temperature.
Staying in the current setting otherwise, we thus want to solve for a given upper bound $y_{\max}>0$ (for simplicity) the \term[problem!optimal control!state-constrained]{state-constrained optimal control problem}
\begin{equation*}
    \begin{aligned}
        &\min_{u\in L^2(\Omega),y\in H^1_0(\Omega)} \frac12\norm{y-y^d}_{L^2(\Omega)}^2 + \frac\alpha2 \norm{u}_{L^2(\Omega)}^2\\
        &\text{subject to \eqref{eq:control:poisson}} \qquad\text{and}\qquad y(x) \leq y_{\max} \quad\text{for almost every }x\in \Omega.
    \end{aligned}
\end{equation*}
This has a similar structure as the control-constrained problem, and we will follow the same general approach. However, this is more delicate here, since we now have to apply the chain rule for the subdifferential of the indicator functional, which requires a nonempty interior of the corresponding set -- which does not hold in $L^2(\Omega)$, and the dual space of $L^\infty(\Omega)$ is very difficult to characterize. We thus instead assume that $\Omega\subset\R^d$ is convex for $d\leq 3$ so that the solutions to the state equation are continuous and we can impose the state constraints \emph{everywhere}. We then define the admissible set
\begin{equation*}
    \Yad \defeq \setof{w\in C(\overline{\Omega})}{w(x)\leq y_{\max}\text{ for all } x\in \overline \Omega}
\end{equation*}
as well as
\begin{equation*}
    \tilde S:L^2(\Omega) \to C(\overline{\Omega}),\qquad u \mapsto y \text{ solving }\eqref{eq:control:poisson},
\end{equation*}
which is well-defined and continuous under our assumptions on $\Omega$. The problem in reduced form is then
\begin{equation}\label{eq:control:stateconstraint}
    \min_{u\in L^2(\Omega)} \frac12\norm{Su-y^d}_{L^2(\Omega)}^2 + \frac\alpha2\norm{u}_{L^2(\Omega)}^2 + \delta_{\Yad}(\tilde Su),
\end{equation}
which has the general form $J = F+\tilde G$ with $F:L^2(\Omega)\to\R$ as above and $\tilde G=\delta_{\Yad}\circ \tilde S:L^2(\Omega)\to C(\overline\Omega)$.

\subsection*{Existence}

Since $\tilde S$ is continuous, $\Yad$ clearly is nonempty, convex, and closed, and $F$ is coercive on $L^2(\Omega)$ due to the control costs while $\tilde G$ is nonnegative, we immediately obtain the existence of an optimal control $\bar u\in L^2(\Omega)$ and an admissible optimal state $\bar y\in \Yad \cap H^1_0(\Omega)$ by \cref{thm:variation:existence}. Since $F$ is strictly convex, this control is again unique.

\subsection*{Optimality conditions}

Setting $G = \delta_{\Yad}:C(\overline\Omega)\to \Rbar$, the problem \eqref{eq:control:stateconstraint} has the form $\min_{u} F(u) + G(\tilde Su)$. Since we are working with continuous functions here and the state equation is linear, we have for $u_0=0\in L^2(\Omega)$ that $y_0\defeq \tilde S u_0 = 0 < y_{\max}$ and hence that $y_0\in\interior \Yad$. We can thus apply the Fenchel--Rockafellar \cref{thm:convex:fenchel} to obtain the primal-dual optimality condition
\begin{equation}\label{eq:control:fenchel}
    \left\{
        \begin{aligned}
            \bar\mu &\in \partial \delta_{\Yad}(\tilde S\bar u),\\
            -\tilde S^*\bar \mu &= S^*(S\bar u-y^d) + \alpha \bar u,
        \end{aligned}
    \right.
\end{equation}
where we have again used the fact that $F$ is Fréchet differentiable with the given gradient. Since $\bar \mu \in \partial \delta_{\Yad}\subset C(\overline\Omega)^* \cong \mathcal{M}(\Omega)$ is a Radon measure (cf.~\cref{ex:functan:dual}\,\ref{ex:functan:dual:iv}), a more explicit, \enquote{pointwise}, interpretation analogous to \eqref{eq:control:constraint:varineq} and \eqref{eq:control:prox} is more involved and involves results from measure theory.

First, by \cref{lem:convex:normalcone}, $\bar\mu\in\mathcal{M}(\Omega)$ and $\bar y\in C(\overline\Omega)$ satisfy
\begin{equation*}
    \int_\Omega (\tilde y(x) - \bar y(x))\,d\bar\mu (x) \leq 0 \qquad\text{for all }\tilde y \leq y_{\max}.
\end{equation*}
By a pointwise argument similar to \cref{ex:convex:subdiff_ind}, it follows that
\begin{equation}\label{eq:control:complementarity}
    \bar\mu\geq 0\qquad\text{and}\qquad \int_\Omega (\bar y (x) - y_{\max})\,d\bar\mu(x) = 0,
\end{equation}
i.e., that $\bar\mu$ is a nonnegative Radon measure whose support is contained in the active set $\setof{x\in\Omega}{\bar y(x) = y_{\max}}$.

Second, using the continuous (and dense) embedding of $W^{1,p}(\Omega)\hookrightarrow C(\overline\Omega)$ for $p$ sufficiently large, it is possible to show that any $\mu\in \mathcal{M}(\Omega)$ satisfies $S^*\mu\in W^{1,q}(\Omega)$ for some sufficiently small $q>1$ and can therefore be characterized as the unique solution $\tilde p = S^*\mu$ to
\begin{equation}\label{eq:control:sconstraint:adjoint}
    \int_\Omega \nabla w(x) \cdot \nabla \tilde p(x) \,dx = \int_\Omega w(x)\,d\mu(x)\qquad \text{for all }w\in W^{1,p}_0(\Omega),
\end{equation}
see, e.g., \cite{Schiela:2010,ClasonSchiela:2015} and the references therein.

Combining \eqref{eq:control:poisson}, \eqref{eq:control:sconstraint:adjoint} added to \eqref{eq:control:adjoint}, and \eqref{eq:control:complementarity}, we obtain from \eqref{eq:control:fenchel} the (suitably interpreted) necessary and sufficient optimality conditions
\begin{equation}\label{eq:control:state_opt}
    \left\{
        \begin{aligned}
            \alpha \bar u + \bar p &= 0,\\
            -\Delta \bar y &= \bar u, \\
            -\Delta \bar p &= \bar y - y^d + \bar \mu,\\
            \bar y &\leq y_{\max},\qquad \bar \mu \geq 0, \qquad \int_\Omega (\bar y (x) - y_{\max})\,d\bar\mu(x) = 0,
        \end{aligned}
    \right.
\end{equation}
compare \cite[Theorem~6.5]{Troeltzsch}.
(The last line corresponds again to the classical complementarity conditions from nonlinear optimization.)

\subsection*{Semismooth Newton method}\label{chap:control:sconstraints}

Since \eqref{eq:control:state_opt} cannot fully be expressed pointwise, a numerical solution is difficult. We thus instead apply the Moreau--Yosida regularization from \cref{sec:moreau-yosida} to $G$, which entails replacing $\partial G:C(\overline\Omega)\setto \mathcal{M}(\Omega)$ in \eqref{eq:control:fenchel} by its Yosida approximation $(\partial G)_\gamma:L^2(\Omega)\to L^2(\Omega)$ for $\gamma>0$ (and, as a consequence, $\tilde S$ by $S$). Following the computation in \cref{ex:moreau}\,\ref{it:moreau:indicator} and using \cref{lem:lebesgue:proximal}, we obtain the pointwise almost everywhere expression
\begin{equation*}
    [H_\gamma(y)](x) \defeq  [(\partial G)_\gamma(y)](x) = \frac1\gamma (y(x)-y_{\max})^+ \defeq \frac1\gamma\max\{0,y(x)-y_{\max}\}
\end{equation*}
and hence the regularized optimality conditions for $(u_\gamma,y_\gamma,p_\gamma)$
\begin{equation}\label{eq:control:state_my_opt}
    \left\{
        \begin{aligned}
            \alpha u_\gamma + p_\gamma &= 0,\\
            -\Delta y_\gamma &= u_\gamma, \\
            -\Delta p_\gamma &= y_\gamma - y^d + \frac1\gamma (y_\gamma - y_{\max})^+,\\
        \end{aligned}
    \right.
\end{equation}
where we have used the single-valued regularized relation $\mu_\gamma = H_\gamma(y_\gamma)$ to eliminate $\mu_\gamma$ in the last line.
By \cref{thm:moreau-yosida} and the computation in \cref{ex:moreau}\,\cref{it:moreau:indicator}, $u_\gamma\in L^2(\Omega)$ is the (unique) minimizer of
\begin{equation*}
    \min_{u\in L^2(\Omega)} \frac12\norm{Su-y^d}_{L^2(\Omega)}^2 + \frac\alpha2\norm{u}_{L^2(\Omega)}^2 + \frac1{2\gamma}\norm{(Su-y_{\max})^+}_{L^2(\Omega)}^2,
\end{equation*}
which guarantees the existence of a (unique) solution $(u_\gamma,y_\gamma,p_\gamma)\in L^2(\Omega)\times H^1_0(\Omega)\times H^1_0(\Omega)$. Of course, we cannot expect $y_\gamma\in \Yad$ in general; but since $S$ is compact, \cref{thm:moreau:convergence} shows that $u_\gamma\weakto \bar u$ as $\gamma\to 0$ (up to subsequences). In fact, using the special structure of the functional to be minimized, a lower semicontinuity argument as in \cref{thm:variation:existence} even shows that $u_\gamma \to \bar u$ and $y_\gamma\to \bar y\in \Yad$ strongly in $L^2(\Omega)$ as $\gamma\to 0$; see \cite[Theorem 6.5]{DelosReyes}

To apply a semismooth Newton method, we first eliminate $u_\gamma=-\frac1\alpha p_\gamma$ from the first relation of \eqref{eq:control:state_my_opt} in the second relation to obtain the \term[condition!optimality, reduced]{reduced optimality system}
\begin{equation}\label{eq:control:state_my_opt_red}
    \left\{
        \begin{aligned}
            -\Delta y_\gamma + \frac1\alpha p_\gamma &= 0, \\
            -\Delta p_\gamma - y_\gamma - \frac1\gamma (y_\gamma - y_{\max})^+ + y^d &=0,\\
        \end{aligned}
    \right.
\end{equation}
which is a nonsmooth system of equations for $(p_\gamma,y_\gamma)$.
Since $y_\gamma\in H^1_0(\Omega)\hookrightarrow L^r(\Omega)$ for some $r>2$, the superposition operator $H_\gamma: y\mapsto \frac1\gamma(y-y_{\max})^+$ is semismooth from $L^r(\Omega)\to L^2(\Omega)$ by \cref{ex:semismooth:l2}\,\ref{ex:semismooth:l2:box} with Newton derivative given pointwise almost everywhere by
\begin{equation*}
    [D_NH_\gamma(y)](x) = \frac1\gamma\1_{(y_{\max},\infty)}(y(x)),
\end{equation*}
where the right-hand side is now to be understood as the linear operator acting by pointwise multiplication with the given function in $L^\infty(\Omega)$.
A semismooth Newton step thus consists in solving for $(\delta p,\delta y)\in H^1_0(\Omega)$ in
\begin{equation}\label{eq:control:state_ssn}
    \begin{pmatrix}
        \frac1\alpha \Id & -\Delta \\
        -\Delta & -\Id - \frac1\gamma \1_{(y_{\max},\infty)}(y^k)
    \end{pmatrix}
    \begin{pmatrix}
        \delta p \\ \delta y
    \end{pmatrix}
    =
    -
    \begin{pmatrix}
            -\Delta y^k + \frac1\alpha p^k\\
            -\Delta p^k - y^k - \frac1\gamma (y^k - y_{\max})^+ + y^d
    \end{pmatrix}
\end{equation}
and then setting $p^{k+1}\defeq  p^k+\delta p$, $y^{k+1}\defeq y^k+\delta y$.
The block operator on the left-hand side of \eqref{eq:control:state_ssn} is a self-adjoint block operator that can be shown to be boundedly invertible (by using the fact that $-\Delta$ is a self-adjoint and positive definite operator) for any $y\in H^1_0(\Omega)$. Hence this semismooth Newton method converges locally superlinearly according to \cref{thm:newton:superlinear}.

Since the PDE constraint is linear, we can further rewrite the Newton step to avoid applying differential operators when evaluating the right-hand side.
Using $\delta p = p^{k+1}-p^{k}$ and similarly for $\delta y$, and writing $(y-y_{\max})^+ = \1_{(y_{\max},\infty)}(y)(y-y_{\max})$, we can rearrange the Newton step as
\begin{algeqbox}
    \begin{equation}\label{eq:control:state_ssn_red}
        \begin{pmatrix}
            \frac1\alpha \Id & -\Delta \\
            -\Delta & -\Id - \frac1\gamma \1_{(y_{\max},\infty)}(y^k)
        \end{pmatrix}
        \begin{pmatrix}
            p^{k+1} \\ y^{k+1}
        \end{pmatrix}
        =
        \begin{pmatrix}
            0 \\
            -y^d - \frac1\gamma  \1_{(y_{\max},\infty)}(y^k) y_{\max}
        \end{pmatrix}.
    \end{equation}
\end{algeqbox}
This is closely related to the \term[method!primal-dual active set]{primal-dual active set method} for quadratic optimization problems with box constraints; see \cite{Hintermuller:2002a,Kunisch:2008a}.
Furthermore, if
\begin{equation*}
    \1_{(y_{\max},\infty)}(y^{k+1}) = \1_{(y_{\max},\infty)}(y^{k})
\end{equation*}
almost everywhere, it is straightforward to verify that \eqref{eq:control:state_ssn_red} coincides with the reduced optimality conditions \eqref{eq:control:state_my_opt_red}, which implies that $u^{k+1} \defeq  -\frac1\alpha p^{k+1}=-\Delta y^{k+1}$ is the desired optimal control. (This \term[property!finite termination]{finite termination property} of semismooth Newton methods for quadratic optimization problems is one reason for their efficiency for such problems.)

In practice, the radius of convergence for the semismooth Newton method applied to such a Moreau--Yosida regularization shrinks with $\gamma\to0$. A possible way of dealing with this is the following \term[strategy!continuation]{continuation strategy}: Starting with a sufficiently large value of $\gamma$, solve a sequence of problems with decreasing $\gamma$ (e.g., $\gamma^k = \gamma^0/2^k$), taking the solution of the previous problem as the starting point for the next (which is hopefully close enough to the solution to lie within the convergence region; otherwise the continuation has to be terminated or the reduction strategy for $\gamma$ adapted).

\begin{figure}
    \begin{subfigure}[t]{0.49\textwidth}
        \includegraphics[width=\textwidth]{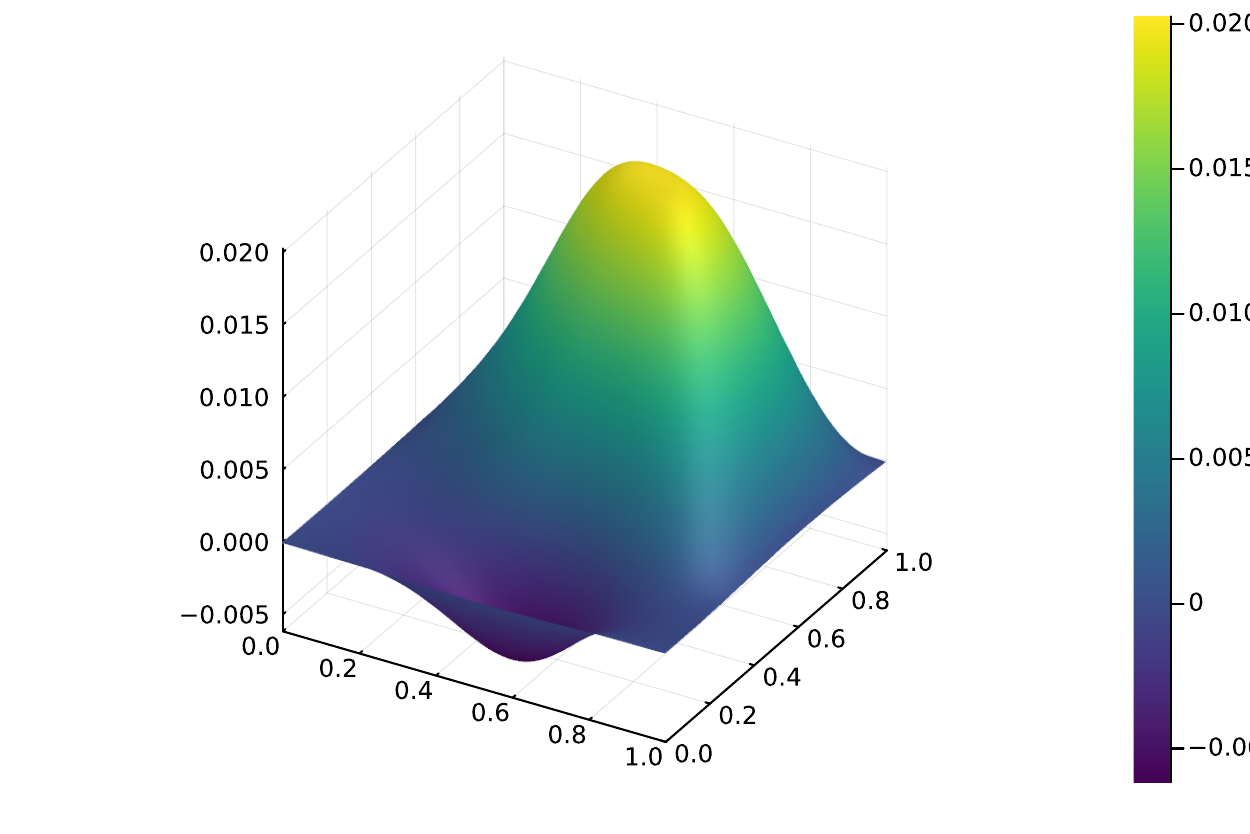}
        \caption{optimal state $y_\gamma$}
    \end{subfigure}
    \hfil
    \begin{subfigure}[t]{0.49\textwidth}
        \includegraphics[width=\textwidth]{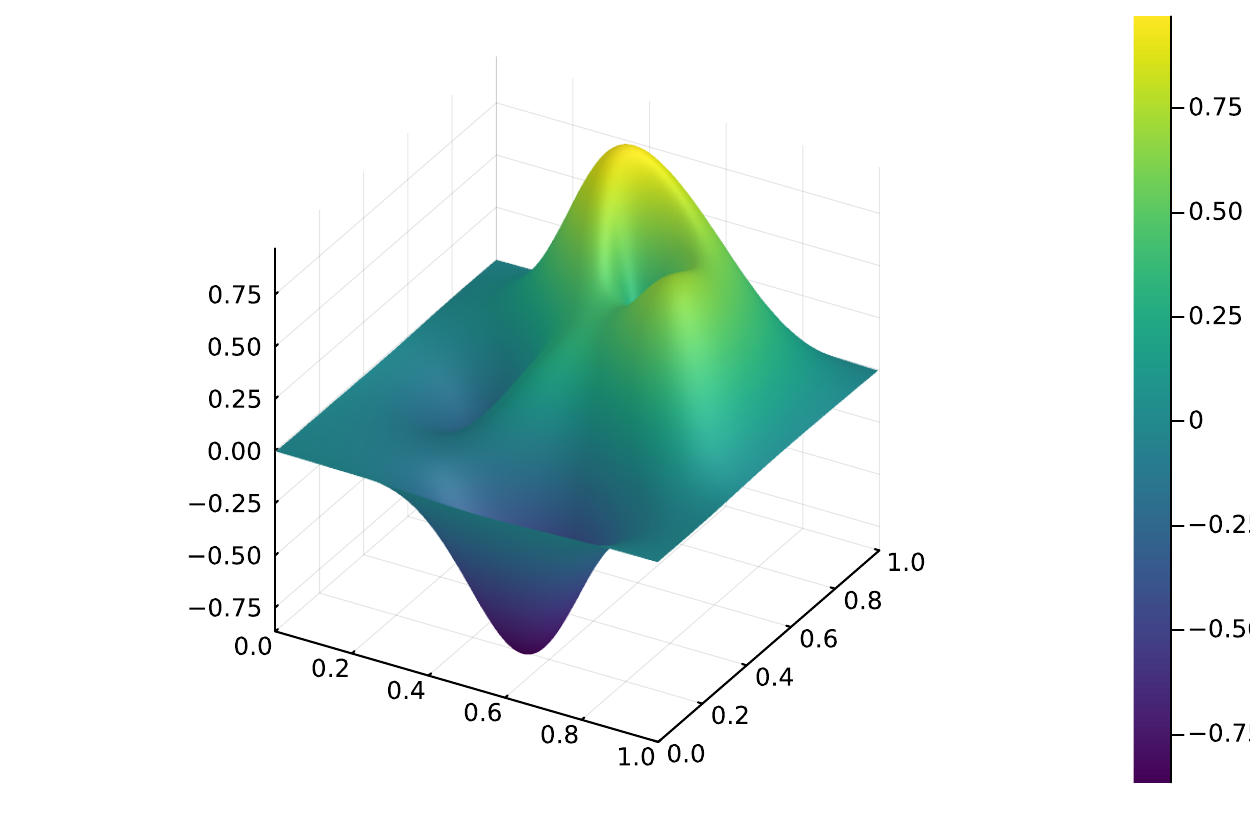 }
        \caption{optimal control $u_\gamma$}
    \end{subfigure}
    \caption{State constraints: optimal control and state for $\gamma=10^{-4}$.}
    \label{fig:control:state-constraints-data}
\end{figure}
\begin{figure}
    \begin{tikzpicture}
        \begin{axis}[%
            xmin=0,
            xmax=20,
            axis x line*=bottom,
            axis y line*=left,
            xlabel={iteration count},
            ylabel={residual},
            ymode=log,
            legend pos = north east,
            cycle list name=primaldual,
            ]

            \pgfplotsinvokeforeach{1.0, 0.1, 0.01, 0.001, 0.0001}{
                \addplot+[line width=1pt] table[
                    x=iter,
                    y = residual
                ]{state_constraints_gamma_ssn_gamma#1.txt};
                \addlegendentry{$\gamma={#1}$}
            }

        \end{axis}
    \end{tikzpicture}
    \caption{State constraints: SSN performance versus Moreau--Yosida parameter $\gamma$.
        We plot the norm of residual, which corresponds to the left hand side of \eqref{eq:control:state_my_opt_red} instantiated at $(p^k, y^k)$.}
    \label{fig:stateconstr:gamma}
\end{figure}

Our target $y^d$, dimension $N=256$, and control cost parameter $\alpha=0.005$ are exactly the same as for control constraints in the previous section. We take $y_{\max}=0.02$; all other details are specified in \cite{nonsmoothbook-codes}.
The optimal state, control, and adjoint state for $\gamma=10^{-4}$ are exemplarily shown in \cref{fig:control:state-constraints-data}.
We illustrate in \cref{fig:stateconstr:gamma} the dependence of the convergence speed of the SSN method on the Moreau--Yosida parameter $\gamma$.

\chapter{Discrete-valued optimal control}\label{chap:discretecontrol}

The final example illustrates the application to a challenging class of \term[problem!optimization!mixed-integer PDE-constrained]{mixed-integer PDE-constrained optimization problems}, where the desired controls are functions that should only take values from a specified discrete set. Such problems arise in, e.g., topology optimization, material parameter identification with a priori information, and joint image reconstruction and segmentation. The purpose of this example is to demonstrate how nonsmooth optimization can be used to impose strong, non-trivial, structural properties on the solution.

Specifically, for a given set of values $u_1< u_2 <\dots < u_m\in \R$, we consider the \term[set!admissible]{admissible set}
\begin{equation*}
    \Uad \defeq \setof{u\in L^2(\Omega)}{u(x)\in\{u_1,\dots,u_m\} \quad \text{for almost every }x\in \Omega}.
\end{equation*}
This set is nonconvex and not weakly closed, which makes the standard theory inapplicable. The usual approach of replacing $\Uad$ with its closed convex hull
\begin{equation*}
    \overline\conv\ \Uad = \setof{u\in L^2(\Omega)}{u(x)\in [u_1,u_m] \quad \text{for almost every }x\in \Omega}
\end{equation*}
however is insufficient as it loses information about the interior values $u_2,\dots,u_{m-1}$. We therefore proceed differently by first adding a \emph{pointwise} quadratic penalty that promotes discrete values of lower magnitude (assuming that lower magnitude is preferable, all other things being equal), i.e., we consider instead of $\delta_{\{u_1,\dots,u_m\}}$ the weighted indicator function
\begin{equation*}
    \hat g(t) \defeq \frac12|t|^2 + \delta_{\{u_1,\dots,u_m\}}(t)
\end{equation*}
whose convex envelope is readily seen by graphical arguments to be
\begin{equation}\label{eq:discretecontrol:multibang}
    g(t) =
    \begin{cases}
        \frac12 \left((u_{i}+u_{i+1})t - u_iu_{i+1}\right) & \text{if }t \in [u_i,u_{i+1}], \quad1\leq i < m,\\
        \infty & \text{else};
    \end{cases}
\end{equation}
see \cref{fig:discretecontrol:g}. (This will be rigorously verified in \cref{rem:discretecontrol:conjugate} below.)

\begin{figure}
    \centering
    \begin{asy}
        import graph;
        unitsize(40,40);
        draw((-2,0)..(2,0),linewidth(0.5),Arrow);
        draw((0,-0.2)..(0,2.25),linewidth(0.5),Arrow);
        int m = 5;
        real a = -1.0;
        real b = 1.0;
        real sqa = -1.5;
        real sqb = 1.5;
        real u = a;
        real high = sqa^2;
        real step = (b-a) / (m-1);
        guide g = (u, high)--(u, u^2);
        for(int i=1; i<m; ++i){
            real u_next = u + step;
            g = g--(u_next, u_next^2);
            u = u_next;
        }
        path p = g--(u, high);
        draw(p, primalline + linewidth(1.5));
        label("$g$", (a, high*0.75), E);

        real f(real x) { return x^2; }
        draw(graph(f, sqa, sqb), dashed + linewidth(0.5));
        label("$x^2$", (sqa, high*0.7));
    \end{asy}
    \caption{Plot of $g$ given by \eqref{eq:discretecontrol:multibang} for $u_1,\ldots,u_5=-1,-0.5,0,0.5,1$. The graph of $x \mapsto x^2$ is also drawn with a dashed line.}
    \label{fig:discretecontrol:g}
\end{figure}

\section{Problem description}

We now consider for given $u_1<\dots<u_m$ and $y^d\in L^2(\Omega)$ the model \term[problem!optimal control!discrete-valued]{discrete-valued control problem}
\begin{equation}\label{eq:discretecontrol:mother}
    \min_{u\in L^2(\Omega)} \frac12\norm{Su-y^d}_{L^2(\Omega)}^2 + \alpha G(u),
\end{equation}
where $S:L^2(\Omega)\to L^2(\Omega)$ is again the control-to-state mapping for \eqref{eq:control:poisson} introduced at the beginning of \cref{chap:control}, $\alpha>0$, and
\begin{equation*}
    G:L^2(\Omega)\to\Rbar,\qquad G(u) \defeq \int_\Omega g(u(x))\,dx.
\end{equation*}
Since $g$, given by \eqref{eq:discretecontrol:multibang}, is proper, convex, and lower semicontinuous, so is $G$ by \cref{lem:lebesgue:lsc}. We thus again obtain from \cref{thm:convex:existence} the existence of an optimal control $\bar u\in L^2(\Omega)$ as well as a corresponding optimal state $\bar y \defeq S\bar u\in H^1_0(\Omega)$. (Since $G$ is convex but not strictly convex, we cannot directly conclude uniqueness; however, since $F$ is strictly convex, the optimal state $\bar y = S\bar u$ must be unique, which then yields uniqueness of $\bar u$ by the continuous invertibility of $S$.)

\begin{remark}
    In the special case that $m=3$ and $u_1 = -M \ll u_2 = 0 \ll u_3 = M$, the convex penalty \eqref{eq:discretecontrol:multibang} simplifies to
    \begin{equation*}
        g(v) = \begin{cases}
            \frac{M}2 |t| & \text{if } |t|\leq M,\\
            \infty &\text{else}.
        \end{cases}
    \end{equation*}
    In other words -- after rescaling $\alpha \mapsto \frac{2}M\alpha$ -- \eqref{eq:discretecontrol:mother} becomes the \term[problem!optimal control!sparse]{sparse control problem}
    \begin{equation*}
        \min_{u\in L^2(\Omega)} \frac12\norm{Su-y^d}_{L^2(\Omega)}^2 + \alpha \norm{u}_{L^1} + \delta_{\{[-M,M]\}}(u),
    \end{equation*}
    which seeks to find a (bounded) optimal control that is zero on as large a part of the domain $\Omega$ as possible; see \cite{Stadler:2007a,Vossen:2006a}. Hence all results in this chapter can be specialized to this problem as well.

    However, in the absence of the control constraints $-M \leq u(x) \leq M$ almost everywhere (or additional $L^2(\Omega)$ regularization), the problem is no longer coercive in $L^1(\Omega)$, and an optimal control must be sought in the space $\mathcal{M}(\Omega)$ of Radon measures \cite{Bidaut:1975}. In this case, it is still possible to exploit similar arguments using the \enquote{preduality}\index{duality!pre-} of $C_0(\Omega)$ and $\mathcal{M}(\Omega)$; see \cref{rem:fenchel:predual} and \cite{ClasonKunisch:2013,CasasClasonKunisch:2012,ClasonSchiela:2015}.
\end{remark}

\section{Optimality conditions}

Again we can derive optimality conditions from the Fermat principle together with calculus rules. As in \cref{chap:control}, \cref{thm:convex:gateaux,thm:subdiff:sum} yield the primal-dual optimality conditions
\begin{equation*}
    \left\{
        \begin{aligned}
            -\bar p  &= S^*(S\bar u -y^d),\\
            \bar p &\in \partial\left(\alpha G(\bar u)\right),
        \end{aligned}
    \right.
\end{equation*}
for the adjoint state $\bar p\in H^1_0(\Omega)$. The last relation implies by \cref{lem:convex:subdiff_calc}\,\ref{lem:convex:subdiff_calc:i} that $\bar p = \alpha\bar q$ for some $\bar q\in \partial G(\bar u)$, i.e., $\frac1\alpha \bar p\in\partial G(\bar u)$. Further applying the \enquote{convex inverse function} \cref{lem:convex:fenchel-young} leads to the equivalent optimality conditions
\begin{equation}\label{eq:discretecontrol:opt_abstract}
    \left\{
        \begin{aligned}
            -\bar p  &= S^*(S\bar u -y^d),\\
            \bar u &\in \partial G^*\left(\frac1\alpha\bar p\right).
        \end{aligned}
    \right.
\end{equation}
To derive from this system some information on the structure of optimal controls, we need to obtain an explicit representation for $\partial G^*$, which we can do pointwise via \cref{thm:lebesgue:fenchel,thm:lebesgue:subdiff}.

We first compute the subdifferential $\partial g(v)$ at a point $v\in [u_1,u_m]$. To that end, we write
\begin{equation*}
    g(t) = g_1 (t) + \delta_{[u_1,u_m]}
\end{equation*}
for the real-valued extension
\begin{equation*}
    g_1:\R\to\R, \qquad g_1(t) = \begin{cases}
        u_1 t-\frac12u_1^2 &\text{if } t \leq u_1,\\
        \frac12 \left((u_{i}+u_{i+1})t - u_iu_{i+1}\right) &\text{if } t \in
        [u_i,u_{i+1}],\quad 1\leq i <m,\\
        u_m t - \frac12 u_m^2 &\text{if } t \geq u_m.
    \end{cases}
\end{equation*}
This is a convex $PC^1$ function, hence by \cref{thm:newton:clarke,thm:clarke:convex} we have that
\begin{equation*}
    \partial g_1(t) =
    \begin{cases}
        \{u_1\} &\text{if } t < u_1,\\
        [u_1, \frac12(u_1+u_2)] &\text{if } t = u_1,\\
        \{\frac12(u_i+u_{i+1})\} &\text{if } t \in (u_i,u_{i+1}), \quad\ 1\leq i<m,\\
        \left[\tfrac12(u_{i-1}+u_{i}),\tfrac12(u_{i}+u_{i+1})\right] &\text{if } t = u_{i},\qquad\qquad 1\leq i<m,\\
        [ \frac12(u_{m-1}+u_m),u_m] &\text{if } t = u_m,\\
        \{u_m\} &\text{if } t > u_m.
    \end{cases}
\end{equation*}
Furthermore, since $g_1$ is continuous in any $t\in\R$, we can apply the sum rule \cref{thm:subdiff:sum} together with the characterization of the subdifferential of the indicator function as a normal cone analogous to \cref{ex:convex:subdiff_ind} to obtain
\begin{equation}\label{eq:discretecontrol:subdiff}
    \partial g(t) =
    \begin{cases}
        (-\infty, \frac12(u_1+u_2)] &\text{if } t = u_1,\\
        \{\frac12(u_i+u_{i+1})\} &\text{if } t \in (u_i,u_{i+1}), \quad\ 1\leq i<m,\\
        \left[\tfrac12(u_{i-1}+u_{i}),\tfrac12(u_{i}+u_{i+1})\right] &\text{if } t = u_{i},\qquad\qquad 1\leq i<m,\\
        [ \frac12(u_{m-1}+u_m),\infty) &\text{if } t = u_m,\\
        \emptyset & \text{else.}
    \end{cases}
\end{equation}
We can now simply appeal to \cref{lem:convex:fenchel-young} (keeping in mind that subdifferentials are always closed) to obtain
\begin{equation}
    \label{eq:discretecontrol:subdiffstar}
    \partial g^*(q) \in \begin{cases}
        \{u_1\} & \text{if }q \in \left(-\infty,\tfrac12(u_1+u_2)\right),\\
        [u_i,u_{i+1}] &\text{if } q = \tfrac12(u_{i}+u_{i+1}),\qquad\qquad\qquad\quad 1\leq i<m,\\
        \{u_i\} &\text{if } q \in \left(\tfrac12(u_{i-1}+u_{i}),\tfrac12(u_{i}+u_{i+1})\right), \quad 1< i<m,\\
        \{u_d\} &\text{if } q \in \left(\tfrac12(u_{m-1}+u_m),\infty\right),\\
        \emptyset &\text{else.}
    \end{cases}
\end{equation}
We illustrate $\subdiff g^*$ in \cref{fig:discretecontrol:subdiffg}.

\begin{figure}
    \centering
    \begin{asy}
        import graph;
        unitsize(40,40);
        int m = 5;
        real a = -1.0;
        real b = 1.0;
        real a_inf = -2.5;
        real b_inf = 2.5;
        real gamma = 0.5;
        real t_prev = a_inf;
        real t = a;
        real u = a;
        real step = (b-a) / (m-1);
        guide g = (t_prev, u);
        for(int i=1; i<=m; ++i){
            real u_next = (i==m ? u : u + step);
            real t_next = t + step;
            real t2 = (t + t_next) / 2;
            real v2 = u;
            real t3 = (i==m ? b_inf : (t + t_next) / 2);
            real v3 = u_next;
            g = g--(t2, v2)--(t3, v3);

            ytick((a_inf-0.2, u), L = "$u_" + (string) i + "$");
            xtick((u, a-0.2),  L = "$u_" + (string) i + "$");

            t_prev = t;
            u = u_next;
            t = t_next;
        }
        draw(g, primalline + linewidth(1.5));
    \end{asy}
    \caption{Plot of $\partial g^*$ given by \eqref{eq:discretecontrol:subdiffstar} for $u_1,\ldots,u_5=-1,-0.5,0,0.5,1$.}
    \label{fig:discretecontrol:subdiffg}
\end{figure}

Applying \eqref{eq:discretecontrol:subdiffstar} and \cref{thm:lebesgue:fenchel,thm:lebesgue:subdiff} in \cref{eq:discretecontrol:opt_abstract}, we now obtain the explicit primal-dual optimality conditions
\begin{equation}\label{eq:discretecontrol:opt_explicit}
    \left\{
        \begin{aligned}
            -\bar p  &= S^*(S\bar u -y^d),\\
            \bar u(x) &\in
            \begin{cases}
                \{u_i\} & \text{if } \bar p(x) \in Q_i,\\
                [u_i,u_{i+1}]& \text{if } \bar p(x) \in Q_{i,i+1},\\
            \end{cases}
        \end{aligned}
    \right.
\end{equation}
for the sets
\begin{align*}
    Q_i &=\setof{q}{\tfrac\alpha2(u_{i-1}+u_i) < q < \tfrac\alpha2(u_{i}+u_{i+1})},\quad  1\leq i\leq m,\\
    Q_{i,i+i}&= \setof{q}{q = \tfrac\alpha2(u_{i}+u_{i+1})}, \qquad\qquad\qquad\quad 1\leq i < m,
\end{align*}
where we have set $u_{0}=-\infty$ and $u_{m+1}=\infty$ to avoid the need for further case distinctions.
This immediately implies that even after convex relaxation, the optimal control will take on almost everywhere one of the prescribed discrete values except where the adjoint state happens to attain one of the critical values $\frac\alpha2(u_i+u_{i+1})$, $i=1,\dots,m$. If this attainment can be excluded -- as in our case, where $\bar p$ is harmonic as the solution of a Poisson equation and thus cannot be constant on a set of positive measure unless it vanishes everywhere -- the relaxed control will still be admissible for the original nonconvex problem and thus locally optimal for the (weighted) discrete problem.
We also see the effect of $\alpha$ on the control: the larger $\alpha$, the more likely that $\bar p(x) \in Q_i$ corresponding to an $u_i$ of lower magnitude.

\begin{remark}\label{rem:discretecontrol:conjugate}
    We point out that it was not necessary to derive the explicit form of the conjugate itself in order to obtain explicit primal-dual optimality conditions. Nevertheless, this can be useful for verifying that $g$ is indeed the convex envelope of $\hat g$.

    First, we have by definition that
    \begin{equation*}
        \hat g^*(q) \defeq \sup_{t\in\{u_1,\dots,u_m\}} q\cdot t - \frac12|t|^2 = u_i q - \frac12|u_i|^2
    \end{equation*}
    for some $1\leq i\leq m$. Since the $u_i$ are assumed to be ordered by increasing magnitude, it therefore suffices to check for given $q\in \R$ whether
    \begin{equation*}
        u_i q - \frac12|u_i|^2 \leq u_{i+1} q - \frac12 |u_{i+1}|^2
    \end{equation*}
    or, equivalently, whether
    \begin{equation*}
        q(u_{i+1}-u_i) \leq \frac12(u_{i+1}^2 - u_i^2).
    \end{equation*}
    Since by assumption $u_{i+1}-u_i>0$, this in turn is equivalent to
    \begin{equation*}
        q\leq \frac12 (u_{i+1} + u_i).
    \end{equation*}
    Hence
    \begin{equation*}
        \hat g^*(q) =
        \begin{cases}
            qu_1 - \frac{1}{2} u_1^2  &\text{if } q\leq\frac12(u_1+u_{2}),\\
            qu_i - \frac{1}{2} u_i^2 &\text{if }\frac12(u_{i-1}+u_i)\leq q \leq \frac12(u_i+u_{i+1}), 1<i<m,\\
            qu_m - \frac{1}{2} u_m^2 &\text{if }\frac12(u_{m}+u_{m-1}) \leq q.
        \end{cases}
    \end{equation*}
    A similar -- albeit more tedious -- calculation using the piecewise differentiability of $g$ shows that
    \begin{equation*}
        g^*(q) = \hat g^*(q).
    \end{equation*}
    By \cref{thm:convex:moreau} and the convexity of $g$, we thus have
    \begin{equation*}
        \hat g^\Gamma = \hat g^{**} = (\hat g^*)^* = (g^*)^* = g.
    \end{equation*}
\end{remark}

\section{Algorithms}
\subsection{Proximal gradient methods}

As in \cref{chap:control:cconstraints}, we can compute a solution to \eqref{eq:discretecontrol:mother} via an explicit splitting method, for which we only need an explicit characterization of the proximal point mapping $\prox_{\gamma(\alpha G)}$. By \cref{lem:lebesgue:proximal}, this is given pointwise almost everywhere by the proximal point mapping for $\alpha g$, which we can derive analogously to \cref{ex:proximal:reell}\,\ref{ex:proximal:reell:ii}. For the sake of presentation, we fix $\alpha = 1$ for now.

By the definition of the proximal point mapping, $w=\prox_{\gamma g}(t) = (\Id + \gamma \partial g)^{-1}(t)$ holds for any $t\in \R$ if and only if $t \in \{w\} + \gamma \partial g(w)$.
Using \eqref{eq:discretecontrol:subdiff}, we thus distinguish the following cases for $w$:
\begin{enumerate}
    \item $w=u_1$: In this case,
        \begin{equation*}
            t \in \{w\} + \gamma \left(-\infty,\tfrac12(u_1+u_2)\right] = \left(-\infty,(1+\tfrac{\gamma}2) u_1 + \tfrac{\gamma}2 u_2\right].
        \end{equation*}
    \item $w\in (u_i,u_{i+1})$ for $1\leq i< m$: In this case,
        \begin{equation*}
            t \in \{w\} + \gamma \{\tfrac12(u_i + u_{i+1})\},
        \end{equation*}
        which first can be solved for $w$ to yield
        \begin{equation*}
            w = t - \tfrac{\gamma}2 (u_i+u_{i+1});
        \end{equation*}
        inserting this into $w\in (u_i,u_{i+1})$ and simplifying then gives
        \begin{equation*}
            t \in \left((1+\tfrac{\gamma}2)u_i + \tfrac{\gamma}2 u_{i+1}, \tfrac{\gamma}2 u_i + (1+\tfrac{\gamma}2)u_{i+1}\right).
        \end{equation*}
    \item $w = u_i$, $1<i<m$: Proceeding as in the first case, we obtain
        \begin{equation*}
            t \in \left[\tfrac{\gamma}2 u_{i-1} + (1+\tfrac{\gamma}2)u_i,(1+\tfrac{\gamma}2) u_i + \tfrac{\gamma}2 u_{i+1}\right].
        \end{equation*}
    \item $w=u_m$: Similarly, this implies that
        \begin{equation*}
            t \in \left[\tfrac{\gamma}2 u_{m-1} + (1+\tfrac{\gamma}2)u_m,\infty\right).
        \end{equation*}
\end{enumerate}
Since this is a complete and disjoint case distinction for $t\in \R$, we obtain that
\begin{equation}\label{eq:discretecontrol:prox-g}
    \prox_{\gamma g}(t) =
    \begin{cases}
        u_i &\text{if } t\in \left[\left(1+\tfrac{\gamma}2\right)u_{i}+\tfrac{\gamma}2 u_{i-1},\left(1+\tfrac{\gamma}2\right)u_{i}+\tfrac{\gamma}2 u_{i+1}\right],\\
        t-\frac{\gamma}2 (u_{i} + u_{i-1}) & \text{if } t \in \left(\left(1+\tfrac{\gamma}2\right)u_{i-1}+\tfrac{\gamma}2 u_{i},\left(1+\tfrac{\gamma}2\right)u_{i}+\tfrac{\gamma}2 u_{i-1}\right),
    \end{cases}
\end{equation}
again with the convention that $u_{0}=-\infty$ and $u_{m+1}=\infty$.
The proximal point mapping therefore has the form of a \term[operator!soft-shrinkage!generalized]{generalized soft-shrinkage operator}.
We illustrate this mapping in \cref{fig:discretecontrol:proxg}.

\begin{figure}
    \centering
    \begin{asy}
        import graph;
        unitsize(40,40);
        int m = 5;
        real a = -1.0;
        real b = 1.0;
        real a_inf = -2.5;
        real b_inf = 2.5;
        real gamma = 0.5;
        real u_prev = a;
        real u = a;
        real step = (b-a) / (m-1);
        guide g = (a_inf, u);
        for(int i=1; i<=m; ++i){
            real u_next = (i==m ? b_inf : u + step);
            real t_prev = u + (u + u_prev)*gamma/2;
            real t = u + (u + u_next)*gamma/2;
             g = g--(t_prev, u)--(t, u);

            ytick((a_inf-0.2, u), L = "$u_" + (string) i + "$");
            xtick((u, a-0.2),  L = "$u_" + (string) i + "$");

            u_prev = u;
            u = u_next;
        }
        g = g--(b_inf, u_prev);
        draw(g, primalline + linewidth(1.5));
    \end{asy}
    \caption{Plot of $\prox_{\gamma g}$ given by \eqref{eq:discretecontrol:prox-g} for $u_1,\ldots,u_5=-1,-0.5,0,0.5,1$ and $\gamma=0.5$.}
    \label{fig:discretecontrol:proxg}
\end{figure}

\begin{remark}
    In the special case of sparse control ($m=3$ and $u_1 = -M \ll u_2 = 0 \ll u_3 = M$), the proximal point mapping reduces to a projection of the well-known \term[operator!soft-shrinkage]{soft-shrinkage operator} from \cref{ex:proximal:reell}.
\end{remark}

\bigskip

Choosing $\tau < 2L^{-1}$ for $L=\norm{S}_{\linear(L^2(\Omega); L^2(\Omega))}^2$ and $u^0 = u_i$ for some $1\leq i\leq m$, we can thus apply the \term[method!proximal gradient]{proximal gradient method}
\begin{algeqbox*}
    \begin{equation*}
        \left\{
            \begin{aligned}
                y^{k+1} &\defeq Su^k &&\text{by solving \eqref{eq:control:poisson}},\\
                p^{k+1} &\defeq S^*(y^d-y^{k+1}) &&\text{by solving \eqref{eq:control:adjoint} for $h=y^d-y^{k+1}$},\\
                u^{k+1}(x) &\defeq \prox_{(\tau\alpha)g}\left(u^k(x) + \tau p^{k+1}(x)\right) &&\text{almost everywhere}.
            \end{aligned}
        \right.
    \end{equation*}
\end{algeqbox*}
By \cref{thm:convergence:fb}, we then have $u^k\weakto \bar u$ in $L^2(\Omega)$. (Since $G$ is not strongly convex, we do not obtain any rates.)

Similarly, we can apply the acceleration strategies from \cref{chap:meta}: The \term[method!proximal gradient!over-relaxed]{over-relaxed proximal gradient method} for $z^0=u^0\in L^2(\Omega)$, $\tau > 0$, and $\lambda = \frac14(1+\sqrt{1+8L\tau})$ consists in computing for $k=0,\dots$
\begin{algeqbox*}
    \begin{equation*}
        \left\{
            \begin{aligned}
                y^{k+1} &\defeq Sz^k &&\text{by solving \eqref{eq:control:poisson}},\\
                p^{k+1} &\defeq S^*(y^d-y^{k+1}) &&\text{by solving \eqref{eq:control:adjoint} for $h=y^d-y^{k+1}$},\\
                u^{k+1}(x) &\defeq \prox_{(\tau\alpha)g}\left(u^k(x) + \tau p^{k+1}(x)\right) &&\text{almost everywhere},\\
                z^{k+1} &\defeq \lambda^{-1} u^{k+1} - (\lambda^{-1}-1) z^k.
            \end{aligned}
        \right.
    \end{equation*}
\end{algeqbox*}
By \cref{thm:meta:overrelax:fb}, we obtain the convergence of the function values $J(\tilde u^N)\to J(\bar u)$ at the rate $O(1/N)$ as $N\to \infty$ for the ergodic sequence $u^N \defeq \frac1N\sum_{k=0}^N u^{k+1}$.

The \term[method!proximal gradient!inertial]{inertial proximal gradient method} for $z^0=u^0\in L^2(\Omega)$, $\tau > 0$, and $\lambda_0 = 1$ consists in computing for $k=0,\dots$
\begin{algeqbox*}
    \begin{equation*}
        \left\{
            \begin{aligned}
                y^{k+1} &\defeq Sz^k &&\text{by solving \eqref{eq:control:poisson}},\\
                p^{k+1} &\defeq S^*(y^d-y^{k+1}) &&\text{by solving \eqref{eq:control:adjoint} for $h=y^d-y^{k+1}$},\\
                u^{k+1}(x) &\defeq \prox_{(\tau\alpha)g}\left(u^k(x) + \tau p^{k+1}(x)\right) &&\text{almost everywhere},\\
                z^{k+1} &\defeq (1+\beta_{k+1}) u^{k+1} - \beta_{k+1} u^k.
            \end{aligned}
        \right.
    \end{equation*}
\end{algeqbox*}
By \cref{thm:meta:inertia:fb}, we obtain the convergence of the function values $J(\tilde u^k)\to J(\bar u)$ at the rate $O(1/k^2)$ as $k\to \infty$ (for the nonergodic sequence).

Note that in all these algorithms, the number $m$ of desired values only enters (linearly!) through the case distinction in \eqref{eq:discretecontrol:prox-g} for the proximal point mapping. In particular, the cost of each step -- which in practice is dominated by computing the solutions $y^{k+1}$ and $p^{k+1}$ of the state and adjoint equation, respectively -- is only mildly affected by $m$. The convex relaxation thus avoids the combinatorial complexity of classical (e.g., branch-and-bound) approaches to mixed-integer optimization.

\subsection{Semismooth Newton method}

The starting point for applying a semismooth Newton method is the dual Moreau--Yosida regularization of \eqref{eq:discretecontrol:opt_abstract}, i.e., replacing the set-valued subdifferential $\partial G^*$ by its single-valued Yosida approximation
\begin{equation*}
    (\partial G^*)_\gamma = \frac1\gamma\left(\Id - \prox_{\gamma G^*}\right)
\end{equation*}
for some $\gamma>0$. Again, we can exploit \cref{lem:lebesgue:proximal} for carrying out the computation pointwise. By \cref{lem:proximal:calculus}\,\ref{lem:proximal:calculus:ii}, we have that
\begin{equation*}
    \begin{aligned}
        \prox_{\gamma g^*}(t) &= t - \gamma \prox_{\gamma^{-1} g}(\tfrac1\gamma t)\\
        &= \begin{cases}
            t-\gamma u_i &\text{in case (i)},\\
            t - \gamma\left(\tfrac1\gamma t - \tfrac{1}{2\gamma}(u_i+u_{i+1})\right) = \frac12(u_i+u_{i+1}) &\text{in case (ii)},
        \end{cases}
    \end{aligned}
\end{equation*}
where case (i) corresponds to
\begin{equation*}
    \tfrac1\gamma t \in \left[1+\tfrac1{2\gamma}u_{i} + \tfrac1{2\gamma} u_{i-1}, (1+\tfrac1{2\gamma})u_i + \tfrac1{2\gamma}u_{i+1}\right],
\end{equation*}
i.e.,
\begin{equation*}
    t \in \left[\gamma u_{i} + \frac12(u_{i-1}+u_i), \gamma u_{i} + \frac12(u_{i}+u_{i+1})\right];
\end{equation*}
and case (ii) corresponds to
\begin{equation*}
    \tfrac1\gamma t \in \left(1+\tfrac1{2\gamma}u_{i-1} + \tfrac1{2\gamma} u_i, (1+\tfrac1{2\gamma})u_i + \tfrac1{2\gamma}u_{i-1}\right),
\end{equation*}
i.e.,
\begin{equation*}
    t \in \left(\gamma u_{i-1} + \frac12(u_{i-1}+u_i), \gamma u_{i} + \frac12(u_{i-1}+u_i)\right);
\end{equation*}
again with the convention $u_0=-\infty$ and $u_{m+1} = \infty$. Hence
\begin{equation*}
    H_\gamma(p) \defeq \left(\partial G^*\right)_\gamma\left(\tfrac1\alpha p\right)
\end{equation*}
is given pointwise almost everywhere by
\begin{equation}
    \label{eq:discretecontrol:h}
    [H_\gamma(p)](x) = h_\gamma(p(x))\defeq
    \begin{cases}
        u_i&\text{if }p(x) \in Q^\gamma_{i},\\
        \frac1{\alpha\gamma}\left(p(x) - \frac{\alpha}{2}(u_{i-1}+u_i)\right) &\text{if }p(x) \in Q^\gamma_{i,i+1},
    \end{cases}
\end{equation}
for
\begin{align*}
    Q^\gamma_{i} &\defeq \left[\alpha\gamma u_{i} + \frac\alpha2(u_{i-1}+u_i), \alpha\gamma u_{i} + \frac\alpha2(u_{i}+u_{i+1})\right],\\
    Q^\gamma_{i,i+1} &\defeq \left(\alpha\gamma u_{i} + \frac\alpha2(u_{i}+u_{i+1}), \alpha\gamma u_{i+1} + \frac\alpha2(u_{i}+u_{i+1})\right).
\end{align*}
We illustrate $h_\gamma$ in \cref{fig:discretecontrol:h}.
Replacing $\subdiff G^*(\tfrac1\alpha\freevar)$ by $H_\gamma$ in \eqref{eq:discretecontrol:opt_abstract} leads to the regularized optimality conditions
\begin{equation}\label{eq:discretecontrol:opt_my}
    \left\{
        \begin{aligned}
            - p_\gamma  &= S^*(Su_\gamma -y^d),\\
            u_\gamma &= H_\gamma(p_\gamma).
        \end{aligned}
    \right.
\end{equation}
Comparing this system with the expansion \eqref{eq:discretecontrol:opt_explicit} of \eqref{eq:discretecontrol:opt_abstract}, we see that the general structure -- in particular, the fact that $u_\gamma(x) = [H_\gamma(p_\gamma)](x)\in \{u_1,\dots,u_m\}$ in the first case -- is conserved; the main difference is that the set-valued second case at a point has been replaced by an affine function (with slope $\frac1\gamma$) in an interval, for which the case distinctions have been adjusted to make room.
(This relates to the fact that by \cref{thm:moreau:conjugate}, the Moreau--Yosida regularization \eqref{eq:discretecontrol:opt_my} is equivalent to replacing $G$ in \eqref{eq:discretecontrol:mother} by $G+\frac\gamma2\norm{\cdot}_{L^2}^2$, i.e., the regularized problem still has the original nonsmooth structure and has merely been made \emph{strongly} convex.)
Comparing \eqref{eq:discretecontrol:h} and \eqref{eq:discretecontrol:opt_explicit}, it is straightforward to verify that a solution satisfying $u_\gamma(x) \in \{u_1,\dots,u_m\}$ for almost every $x\in \Omega$ also satisfies the unregularized optimality conditions \eqref{eq:discretecontrol:opt_explicit} and is therefore optimal for \eqref{eq:discretecontrol:mother} as well; in this sense, the Moreau--Yosida regularization is an \term[penalization, exact]{exact (dual) penalization}.

\begin{figure}
    \centering
    \begin{asy}
        import graph;
        unitsize(40,40);
        int m = 5;
        real a = -1.0;
        real b = 1.0;
        real a_inf = -2.5;
        real b_inf = 2.5;
        real gamma = 0.5;
        real alpha = 1.0;
        real u_prev = a_inf;
        real p1 = u_prev;
        real u = a;
        real step = (b-a) / (m-1);
        guide g = (a_inf, u);
        for(int i=1; i<=m; ++i){
            real u_next = (i==m ? b_inf : u + step);
            real p2 = alpha * gamma * u + (u + u_prev) * alpha/2;
            real v2 = u;
            real p3 = i==m ? b_inf : alpha * gamma * u + (u + u_next) * alpha/2;
            real v3 = u;
            g = g--(p2, v2)--(p3, v3);

            ytick((a_inf-0.2, u), L = "$u_" + (string) i + "$");
            if(i>1){
                xtick((p2, a-0.2));
                label("$\scriptstyle Q^\gamma_{" + (string) (i-1) + "," + (string) i + "}$",
                      ((p1+p2)/2, a-0.4), S);
            }
            if(i<m){
                xtick((p3, a-0.2));
            }
            label("$\scriptstyle Q^\gamma_{" + (string) i + "}$",
                    ((p2+p3)/2, a-0.2), S);

            p1 = p3;
            u_prev = u;
            u = u_next;
        }
        draw(g, primalline + linewidth(1.5));
    \end{asy}
    \caption{Plot of $h_\gamma$ given by \eqref{eq:discretecontrol:h} for $u_1,\ldots,u_5=-1,-0.5,0,0.5,1$ and $\gamma=0.5$ and $\alpha=1$.}
    \label{fig:discretecontrol:h}
\end{figure}

\bigskip

We now derive the semismooth Newton iteration for solving \eqref{eq:discretecontrol:opt_my}. First, it is again advantageous to reformulate the system using the definition of $S$ and $S^*$ as well as the second equation of \eqref{eq:discretecontrol:opt_my} as
\begin{equation}\label{eq:discretecontrol:opt_my_red}
    \left\{
        \begin{aligned}
            -\Delta p_\gamma + y_\gamma - y^d &=0,\\
            -\Delta y_\gamma - H_\gamma(p_\gamma) &= 0,
        \end{aligned}
    \right.
\end{equation}
cf. \eqref{eq:control:state_my_opt_red}, which we can consider as a nonlinear equation $T(y,p) = 0$ for $T:H^1_0(\Omega)\times H^1_0(\Omega)\to H^1_0(\Omega)^*\times H^1_0(\Omega)^*$. (The corresponding optimal control can be recovered from its solution via $u_\gamma = H_\gamma(p_\gamma)$, which is a simple pointwise evaluation.)

To obtain a Newton derivative $D_N T(y,p)$, we clearly only need to compute one for $H_\gamma$, which we again do pointwise. First, it is straightforward to verify that $h_\gamma$ is continuous and piecewise linear, so that by \cref{thm:newton:clarke,thm:newton:clarke_ndiff} we have that
\begin{equation}\label{eq:discretecontrol:dh}
    D_N h_\gamma(t) \defeq
    \begin{cases}
        \frac1{\alpha\gamma} &\text{if } t \in Q^\gamma_{i,i+1},\\
        0 & \text{else},
    \end{cases}
\end{equation}
is a Newton derivative for $h_\gamma$ at $t$. Clearly, this function is uniformly bounded by $\frac1{\alpha\gamma}$. For fixed $\gamma>0$, the intervals $Q^\gamma_{i,i+1}$ are also separated, and hence $D_N h_\gamma$ is a Baire--Carathéodory function. Since $p_\gamma\in H^1_0(\Omega) \hookrightarrow L^r(\Omega)$ for some $r>2$, it thus follows from \cref{thm:newton:super} that a Newton derivative of $H_\gamma$ at $p$ in direction $\delta p\in L^r(\Omega)$ is given pointwise almost everywhere by
\begin{equation*}
    [D_NH_\gamma(p)\delta p](x) =
    \begin{cases}
        \frac1{\alpha\gamma}\delta p(x) &\text{if } p(x) \in Q^\gamma_{i,i+1},\\
        0 & \text{else}.
    \end{cases}
\end{equation*}
Setting $Q^\gamma \defeq \bigcup_{i=1}^m Q^\gamma_{i,i+1}$, we thus obtain as a Newton derivative for $T$ at $(y,p)\in H^1_0(\Omega)\times H^1_0(\Omega)$ the block operator
\begin{equation*}
    D_N T(y,p) = \begin{pmatrix}
        \Id & -\Delta \\
        -\Delta & -\frac1{\alpha\gamma}\1_{Q^\gamma}(p)
    \end{pmatrix},
\end{equation*}
where again $[\1_{Q^\gamma}(p)](x) = 1$ if $p(x)\in Q^\gamma$ and $0$ otherwise, and the bottom-right block is to be understood as the linear operator acting by pointwise multiplication with this function in $L^\infty(\Omega)$.
This is a self-adjoint block operator that can be shown to be uniformly (with respect to $p$) boundedly invertible; see \cite[Proposition 4.3]{ClasonKunisch:2013}. Hence by \cref{thm:newton:superlinear}, the following semismooth Newton method converges locally superlinearly to a solution to \eqref{eq:discretecontrol:opt_my_red}: Given $(p^k,y^k)\in H^1_0(\Omega)\times H^1_0(\Omega)$,
\begin{algenumbox}
    \begin{enumerate}[label=\arabic*.]
        \item solve for $(\delta p,\delta y)\in H^1_0(\Omega)\times H^1_0(\Omega)$ the coupled linear system
            \begin{align*}
                -\Delta \delta p + \delta y &= y^d - y^k + \Delta p^k,\\
                -\Delta \delta y - \frac1{\alpha\gamma}\1_{Q^\gamma}(p^k) \delta p &= H_\gamma(p^k) + \Delta y^k,
            \end{align*}
        \item set
            \begin{equation*}
                y^{k+1} = y^k + \delta y,\qquad p^{k+1} = p^k + \delta p.
            \end{equation*}
    \end{enumerate}
\end{algenumbox}
Using the linearity of the state equation and comparing \eqref{eq:discretecontrol:h} with \eqref{eq:discretecontrol:dh}, this can again be reformulated as a linear system for $(p^{k+1},y^{k+1})$.
Similarly to the proximal gradient methods, the number $m$ of desired states only enters linearly via the case distinction in $Q^\gamma$. In particular, the computation of the Newton step itself is independent of the value of $m$, hence avoiding combinatorial complexity.
As in \cref{chap:control:sconstraints}, this will in practice be embedded in a continuation strategy for $\gamma\to 0$.
\begin{figure}[t!]
    \begin{subfigure}[t]{0.49\textwidth}
        \includegraphics[width=\textwidth]{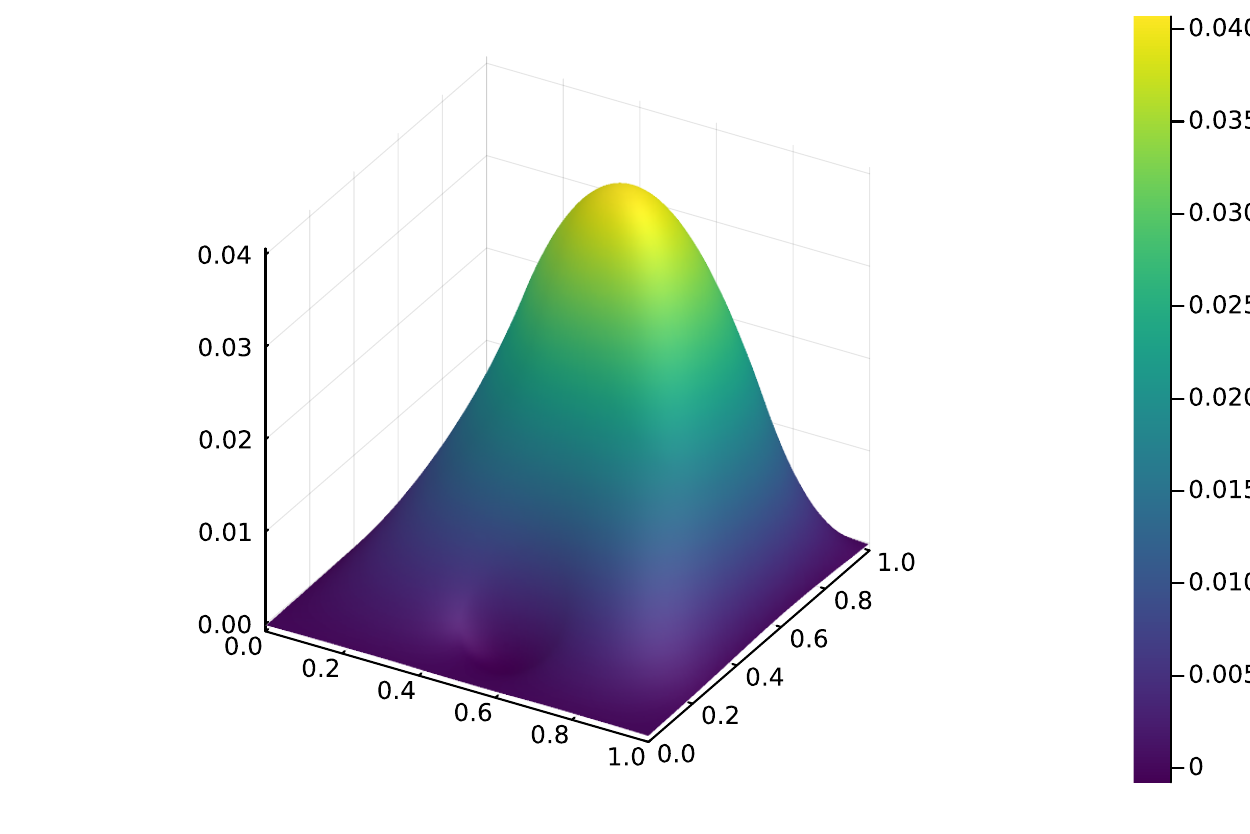}
        \caption{optimal state $y$, $m=3$}
    \end{subfigure}
    \hfil
    \begin{subfigure}[t]{0.49\textwidth}
        \includegraphics[width=\textwidth]{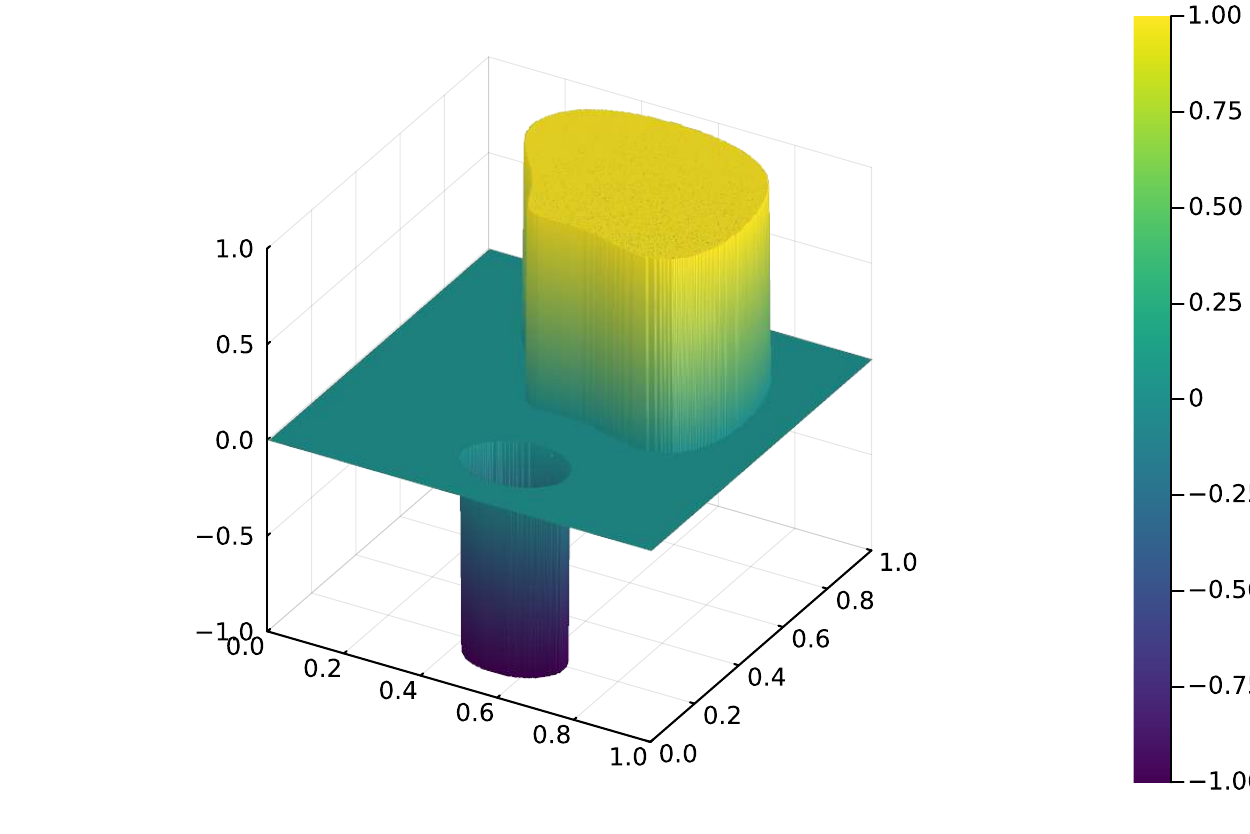}
        \caption{optimal control $u$, $m=3$}
    \end{subfigure}
    \\
    \begin{subfigure}[t]{0.49\textwidth}
        \includegraphics[width=\textwidth]{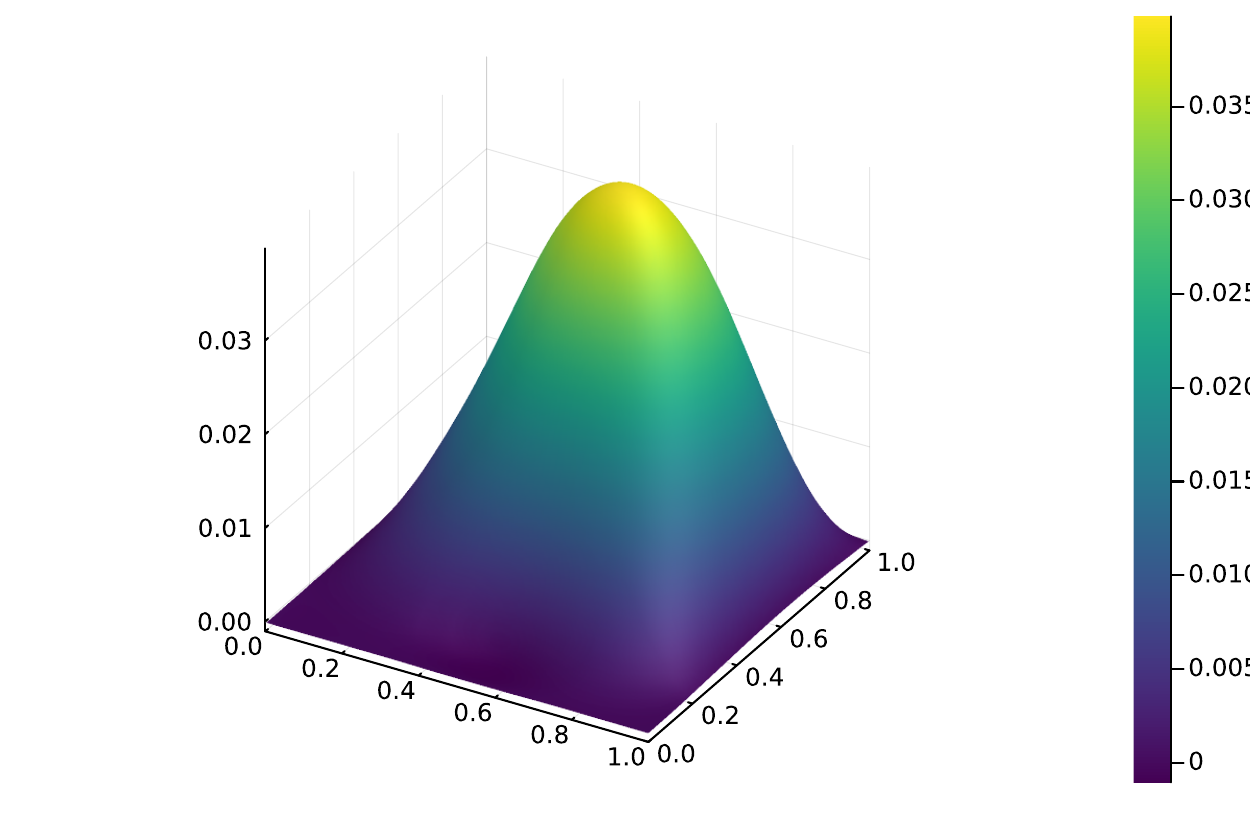}
        \caption{optimal state $y$, $m=10$}
    \end{subfigure}
    \hfil
    \begin{subfigure}[t]{0.49\textwidth}
        \includegraphics[width=\textwidth]{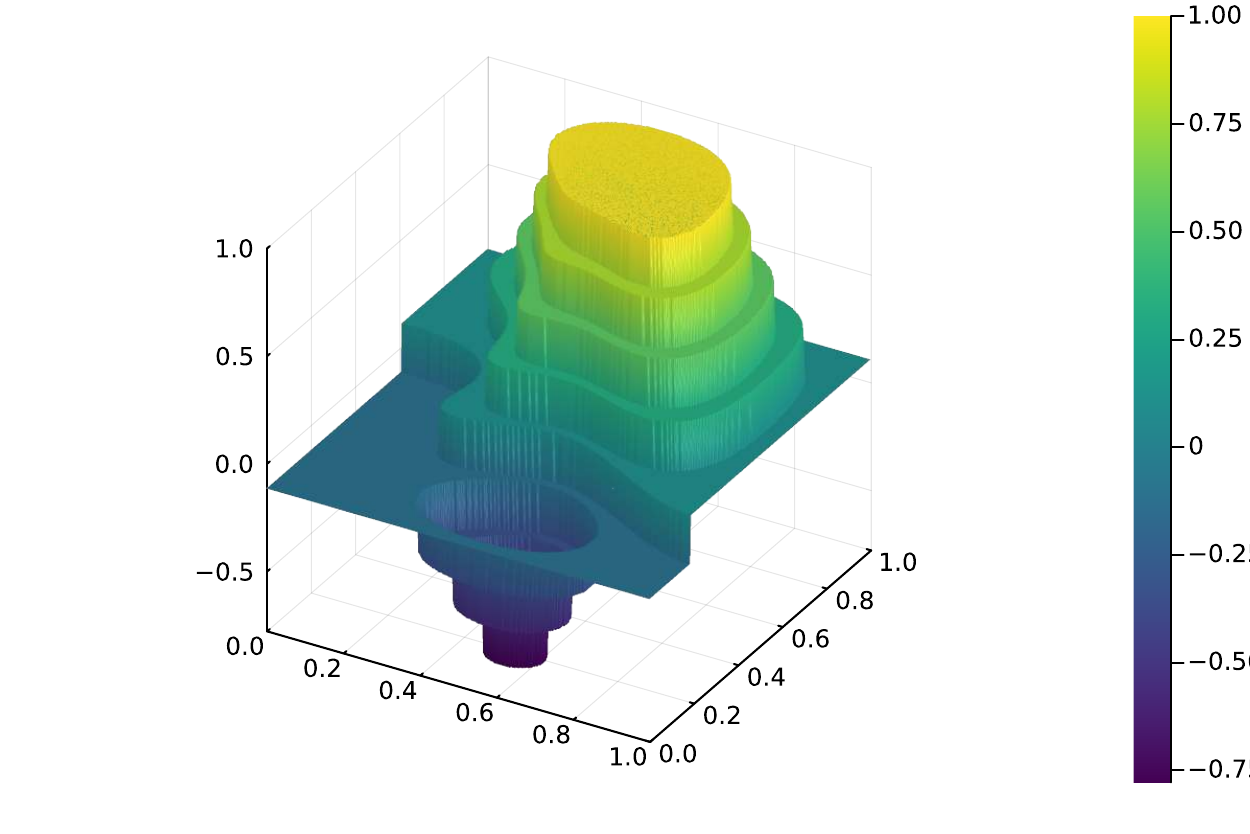}
        \caption{optimal control $u$, $m=10$}
    \end{subfigure}
    \\
    \begin{subfigure}[t]{0.49\textwidth}
        \includegraphics[width=\textwidth]{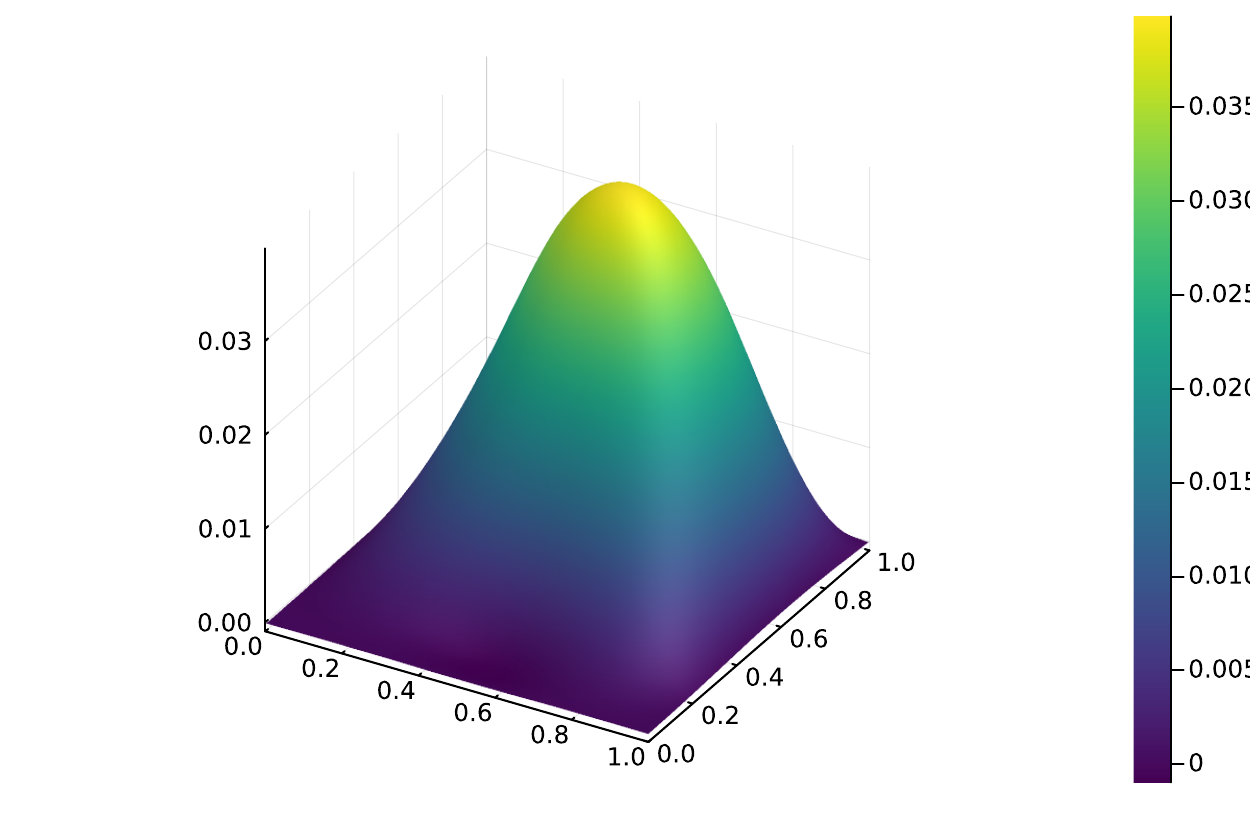}
        \caption{optimal state $y$, $m=20$}
    \end{subfigure}
    \hfil
    \begin{subfigure}[t]{0.49\textwidth}
        \includegraphics[width=\textwidth]{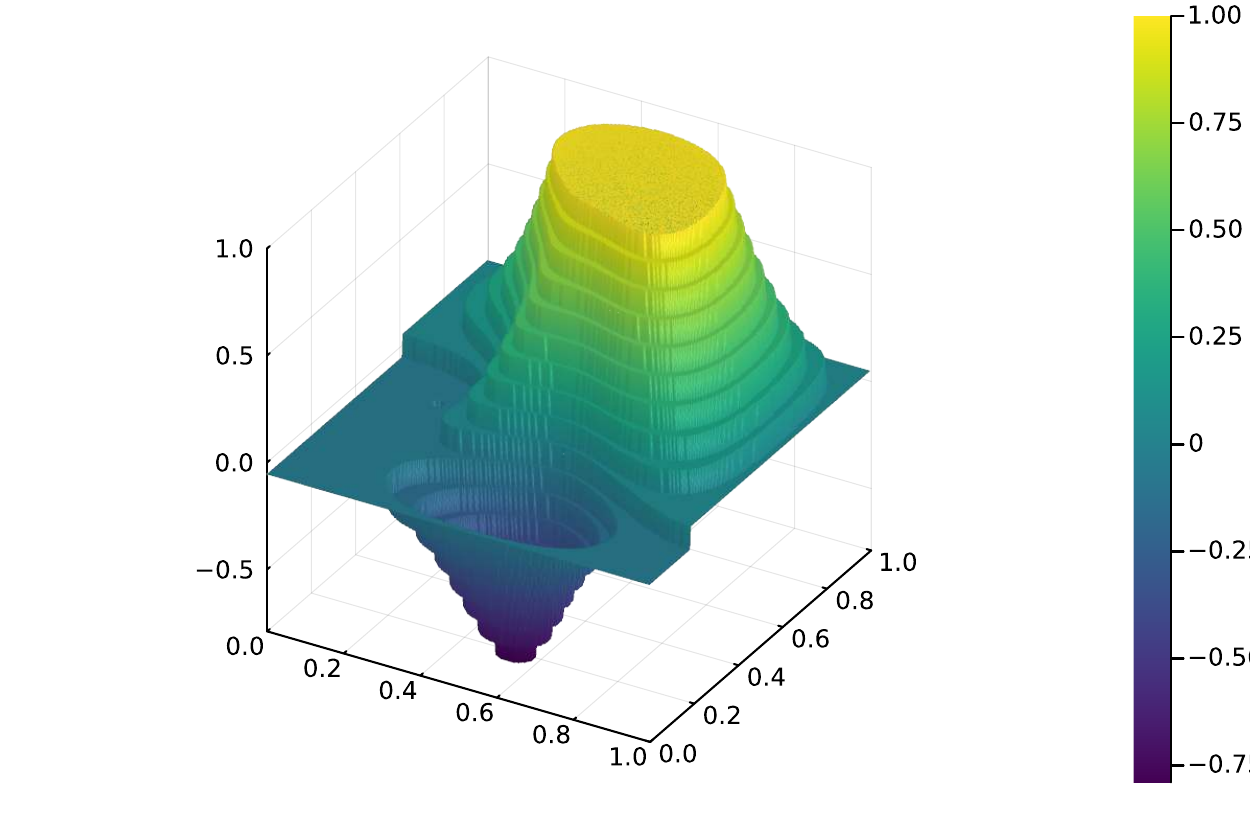}
        \caption{optimal control $u$, $m=20$}
    \end{subfigure}
    \caption{Discrete control example control and state for $\gamma=10^{-6}$ and $m=3,10,20$.
        The target $y^d$ is the same as for control constraints in \cref{fig:control:control-constraints-data} and state constraints in \cref{fig:control:state-constraints-data}.}
    \label{fig:discretecontrol:reco}
\end{figure}

We indicate the dependence of solutions to the model discrete-valued control problem \eqref{eq:discretecontrol:mother} on the set of allowed values $u_1,\ldots,u_m$, by taking $m=3,10,20$ equally spaced controls on $[a,b]=[-1,1]$.
Other than this restriction on the values of the control $u$, the experimental setup is the same as for control constraints in \cref{chap:control:cconstraints}.
For the first-order methods, we take the step length parameter $\tau=0.9/L^2$, where $L$ is an estimate of $\norm{S}$. For the SSN method, the Moreau--Yosida regularization parameter is set to $\gamma=10^{-6}$.
The corresponding optimal controls $u_\gamma$ for this value of $\gamma$ are verified to only take on admissible values almost everywhere and thus are also optimal for \eqref{eq:discretecontrol:mother}; see \cref{fig:discretecontrol:reco}. Note also how the discrete-valued controls better approximate the solution shown in \cref{fig:control:control-constraints-data} as $m$ increases.
Regarding performance, the SSN method converges to near machine precision within $15$ iterations. The forward-backward splitting methods are significantly slower both in number of iterations and actual runtime; see \cref{fig:discretecontrol:performance}. In fact, for the shown parameters, all three splitting methods will yield a control taking only allowed values only after $5000$ iterations, at which point the residual norm drops to machine precision similarly to the SSN method.

\begin{figure}[t!]
    \begin{tikzpicture}
        \begin{axis}[%
            axis x line*=bottom,
            axis y line*=left,
            xlabel={iteration count},
            ylabel={residual},
            ymode=log,
            xmin = 0,
            xmax = 50,
            legend pos = south east,
            ]

            \addplot [fb] table[x=iter,y=residual]{discrete_control_fb_m10.txt};
            \addlegendentry{FB}

            \addplot [fista] table[x=iter,y=residual]{discrete_control_fb_inertia_m10.txt};
            \addlegendentry{FISTA}

            \addplot [fbrelax] table[x=iter,y=residual]{discrete_control_fb_overrelax_m10.txt};
            \addlegendentry{over-relaxed FB}

            \addplot [ssn] table[x=iter,y=residual]{discrete_control_ssn_m10.txt};
            \addlegendentry{SSN}
        \end{axis}
    \end{tikzpicture}
    \begin{tikzpicture}
        \begin{axis}[%
            axis x line*=bottom,
            axis y line*=left,
            xlabel={time [s]},
            ylabel={residual},
            ymode=log,
            xmode=log,
            legend pos = south west,
            ]

            \addplot [fb] table[x=cputime,y=residual]{discrete_control_fb_m10.txt};
            \addlegendentry{FB}

            \addplot [fista] table[x=cputime,y=residual]{discrete_control_fb_inertia_m10.txt};
            \addlegendentry{FISTA}

            \addplot [fbrelax] table[x=cputime,y=residual]{discrete_control_fb_overrelax_m10.txt};
            \addlegendentry{over-relaxed FB}

            \addplot [ssn] table[x=cputime,y=residual]{discrete_control_ssn_m10.txt};
            \addlegendentry{SSN}
        \end{axis}
    \end{tikzpicture}
    \caption{Algorithm performance for the discrete control problem, $m=10$, $\gamma=10^{-6}$. For the SSN method, we plot the residual $\norm{H_\gamma(p^k)}_2$, while for the first-order methods, with $\gamma=0$, the residual is similarly given by the violation of \eqref{eq:discretecontrol:opt_explicit}.}
    \label{fig:discretecontrol:performance}
\end{figure}

\begin{remark}
    The convex relaxation described in this chapter was first proposed in \cite{ClasonKunisch:2013} (corresponding to the formal limit $\beta\to \infty$ there) and later applied to topology optimization \cite{ClasonKunisch:2016,ClasonKunischTrautmann:2019} and parameter identification \cite{ClasonDo:2018} problems. Vector-valued problems were considered in \cite{ClasonTamelingWirth:2021} and \cite{ClasonItoKunisch:2016}, the latter treating the related problem of \enquote{switching controls}, where at most one of a pair $(u,v)$ of distributed controls should be active at any point, i.e., $u(x)v(x)=0$ should hold pointwise almost everywhere.
    The presentation here is condensed from \cite{ClasonKunisch:2016,ClasonDo:2018,ClasonKunischTrautmann:2019,ClasonItoKunisch:2016}.
\end{remark}

\backmatter

\printbibliography[heading=bibintoc]

@string{advneural = "Advances in Neural Information Processing Systems"}

@string{amo = "Applied Mathematics \& Optimization"}

@string{cocv = "ESAIM: Control, Optimisation and Calculus of Variations"}

@string{computoptappl = "Computational Optimization and Applications"}

@string{convexanal = "Journal of Convex Analysis"}

@string{cpam = "Communications on Pure and Applied Mathematics"}

@string{doklady = "Doklady Akademii Nauk SSSR"}

@string{etna = "Electronic Transactions on Numerical Analysis"}

@string{ieeeimg = "IEEE Transactions on Image Processing"}

@string{ieeesig = "IEEE Transactions on Signal Processing"}

@string{invprob = "Inverse Problems"}

@string{jmaa = "Journal of Mathematical Analysis and Applications"}

@string{jmiv = "Journal of Mathematical Imaging and Vision"}

@string{jota = "Journal of Optimization Theory and Applications"}

@string{jscicomp = "Journal of Scientific Computing"}

@string{mathprog = "Mathematical Programming"}

@string{procams = "Proceedings of the American Mathematical Society"}

@string{siamjapl = "SIAM Journal on Applied Mathematics"}

@string{siamjcon = "SIAM Journal on Control and Optimization"}

@string{siamjimg = "SIAM Journal on Imaging Sciences"}

@string{siamjnum = "SIAM Journal on Numerical Analysis"}

@string{siamjopt = "SIAM Journal on Optimization"}

@string{siamjscicomp = "SIAM Journal on Scientific Computation"}

@string{svaa = "Set-Valued and Variational Analysis"}

@string{transams = "Transactions of the American Mathematical Society"}

@article{adly2015newton,
  author = {Adly, Samir and Cibulka, Radek and Ngai, Huynh Van},
  date = {2015},
  doi = {10.1137/130926730},
  journaltitle = siamjopt,
  number = {1},
  pages = {159--184},
  title = {{N}ewton's method for solving inclusions using set-valued approximations},
  volume = {25},
}

@article{Alber:1998,
  author = {Alber, Ya. I. and Iusem, A. N. and Solodov, M. V.},
  date = {1998},
  doi = {10.1007/BF01584842},
  journaltitle = mathprog,
  number = {1},
  pages = {23--35},
  title = {On the projected subgradient method for nonsmooth convex optimization in a {H}ilbert space},
  volume = {81},
}

@book{Alt:2016,
  author = {{Alt}, Hans Wilhelm},
  location = {London},
  publisher = {Springer},
  date = {2016},
  doi = {10.1007/978-1-4471-7280-2},
  series = {Universitext},
  subtitle = {An Application-Oriented Introduction},
  title = {Linear Functional Analysis},
}

@book{Ambrosio,
  author = {Ambrosio, Luigi and Fusco, Nicola and Pallara, Diego},
  publisher = {The Clarendon Press, Oxford University Press, New York},
  date = {2000},
  doi = {10.1093/oso/9780198502456.001.0001},
  series = {Oxford Mathematical Monographs},
  title = {Functions of Bounded Variation and Free Discontinuity Problems},
}

@book{Appell:1990,
  author = {Appell, Jürgen and Zabrejko, Petr P.},
  location = {New York},
  publisher = {Cambridge University Press},
  date = {1990},
  doi = {10.1017/cbo9780511897450},
  title = {Nonlinear Superposition Operators},
}

@article{aragon2012lyusternik,
  author = {Aragón Artacho, Francisco J. and Gaydu, Michaël},
  date = {2012},
  doi = {10.1007/s10589-011-9439-6},
  journaltitle = computoptappl,
  number = {3},
  pages = {785--803},
  title = {A {L}yusternik--{G}raves theorem for the proximal point method},
  volume = {52},
}

@article{artacho2013metric,
  author = {Aragón Artacho, Francisco J. and Geoffroy, Michel H.},
  date = {2014},
  eprint = {1303.3654},
  eprinttype = {arXiv},
  journaltitle = {Journal of Nonlinear and Convex Analysis},
  number = {1},
  pages = {35--47},
  title = {Metric subregularity of the convex subdifferential in {B}anach spaces},
  volume = {15},
}

@book{arrow1958strudies,
  author = {Arrow, K. J. and Hurwicz, L. and Uzawa, H.},
  location = {Stanford},
  publisher = {Stanford University Press},
  date = {1958},
  title = {Studies in Linear and Non-Linear Programming},
}

@article{aspelmeier2016local,
  author = {Aspelmeier, Timo and Charitha, C. and Luke, D. Russell},
  date = {2016},
  doi = {10.1137/15M103580X},
  journaltitle = siamjimg,
  number = {2},
  pages = {842--868},
  title = {Local linear convergence of the {ADMM}/ {D}ouglas--{R}achford algorithms without strong convexity and application to statistical imaging},
  volume = {9},
}

@book{Attouch:1984,
  author = {Attouch, H.},
  publisher = {Pitman (Advanced Publishing Program), Boston, MA},
  date = {1984},
  series = {Applicable Mathematics Series},
  title = {Variational Convergence for Functions and Operators},
}

@incollection{AttouchBrezis,
  author = {Attouch, Hédy and Brezis, Haı̈m},
  publisher = {North-Holland, Amsterdam},
  booktitle = {Aspects of Mathematics and its Applications},
  date = {1986},
  doi = {10.1016/S0924-6509(09)70252-1},
  pages = {125--133},
  series = {North-Holland Math. Library},
  title = {Duality for the sum of convex functions in general {B}anach spaces},
  volume = {34},
}

@book{Attouch,
  author = {Attouch, Hedy and Buttazzo, Giuseppe and Michaille, Gérard},
  location = {Philadelphia},
  publisher = {Society for Industrial and Applied Mathematics},
  date = {2014},
  doi = {10.1137/1.9781611973488},
  edition = {2},
  series = {MOS-SIAM Series on Optimization},
  title = {Variational Analysis in {S}obolev and {BV} Spaces},
  volume = {6},
}

@book{aubin1990sva,
  author = {Aubin, J.P. and Frankowska, H.},
  publisher = {Birkhäuser Basel},
  date = {1990},
  doi = {10.1007/978-0-8176-4848-0},
  title = {Set-Valued Analysis},
}

@incollection{aubin1981contingent,
  author = {Aubin, Jean-Pierre},
  publisher = {Academic Press, New York-London},
  booktitle = {Mathematical Analysis and Applications, {P}art {A}},
  date = {1981},
  pages = {159--229},
  series = {Adv. in Math. Suppl. Stud.},
  title = {Contingent derivatives of set-valued maps and existence of solutions to nonlinear inclusions and differential inclusions},
  volume = {7},
}

@article{aubin1984lipschitz,
  author = {Aubin, Jean-Pierre},
  date = {1984},
  doi = {10.1287/moor.9.1.87},
  journaltitle = {Mathematics of Operations Research},
  number = {1},
  pages = {87--111},
  title = {Lipschitz behavior of solutions to convex minimization problems},
  volume = {9},
}

@article{aussel2005subsmooth,
  author = {Aussel, D and Daniilidis, Aris and Thibault, L},
  date = {2005},
  doi = {10.1090/S0002-9947-04-03718-3},
  journaltitle = transams,
  number = {4},
  pages = {1275--1301},
  title = {Subsmooth sets: functional characterizations and related concepts},
  volume = {357},
}

@article{Aze:1995,
  author = {Azé, Dominique and Penot, Jean-Paul},
  url = {https://eudml.org/doc/73364},
  date = {1995},
  journaltitle = {Annales de la Faculté des sciences de Toulouse: Mathématiques},
  number = {4},
  pages = {705--730},
  title = {Uniformly convex and uniformly smooth convex functions},
  volume = {4},
}

@article{Bacak:2017,
  author = {Bačák, Miroslav and Kohlenbach, Ulrich},
  url = {https://www.heldermann.de/JCA/JCA25/JCA254/jca25079.htm},
  date = {2018},
  eprint = {1709.04700},
  eprinttype = {arXiv},
  journaltitle = convexanal,
  number = {4},
  pages = {1291--1318},
  title = {On proximal mappings with {Young} functions in uniformly convex {Banach} spaces},
  volume = {25},
}

@book{bagirov2014nonsmooth,
  author = {Bagirov, Adil and Karmitsa, Napsu and Mäkelä, Marko M.},
  publisher = {Springer, Cham},
  date = {2014},
  doi = {10.1007/978-3-319-08114-4},
  note = {Theory, practice and software},
  title = {Introduction to Nonsmooth Optimization},
}

@book{Barbu:2012,
  author = {Barbu, Viorel and Precupanu, Teodor},
  publisher = {Springer, Dordrecht},
  date = {2012},
  doi = {10.1007/978-94-007-2247-7},
  edition = {4},
  series = {Springer Monographs in Mathematics},
  title = {Convexity and Optimization in {B}anach Spaces},
}

@book{Bauschke,
  author = {Bauschke, Heinz H. and Combettes, Patrick L.},
  location = {New York},
  publisher = {Springer},
  date = {2017},
  doi = {10.1007/978-3-319-48311-5},
  edition = {2},
  series = {CMS Books in Mathematics/Ouvrages de Mathématiques de la SMC},
  title = {Convex Analysis and Monotone Operator Theory in {H}ilbert Spaces},
}

@book{BauschkeMoursi:2023,
  author = {Bauschke, Heinz H. and Moursi, Walaa M.},
  location = {Philadelphia, PA},
  publisher = {Society for Industrial and Applied Mathematics},
  date = {2023},
  doi = {10.1137/1.9781611977806},
  title = {An Introduction to Convexity, Optimization, and Algorithms},
}

@book{Beck:2017,
  author = {Beck, Amir},
  location = {Philadelphia},
  publisher = {Society for Industrial and Applied Mathematics},
  date = {2017},
  doi = {10.1137/1.9781611974997},
  title = {First-Order Methods in Optimization},
}

@article{BeckTeboulle,
  author = {Beck, Amir and Teboulle, Marc},
  date = {2009},
  doi = {10.1137/080716542},
  journaltitle = siamjimg,
  number = {1},
  pages = {183--202},
  title = {A fast iterative shrinkage-thresholding algorithm for linear inverse problems},
  volume = {2},
}

@article{beck2009fgp,
  author = {Beck, Amir and Teboulle, Marc},
  date = {2009},
  doi = {10.1109/TIP.2009.2028250},
  journaltitle = ieeeimg,
  number = {11},
  pages = {2419--2434},
  title = {Fast gradient-based algorithms for constrained total variation image denoising and deblurring problems},
  volume = {18},
}

@book{Beck:2024,
  author = {Beck, Lisa and Schmidt, Bernd},
  publisher = {Cham: Birkhäuser},
  date = {2025},
  doi = {10.1007/978-3-031-59138-9},
  series = {Mathematik Kompakt},
  title = {Variationsrechnung},
}

@inproceedings{benning2015preconditioned,
  author = {Benning, Martin and Knoll, Florian and Schönlieb, Carola-Bibiane and Valkonen, Tuomo},
  editor = {Bociu, Lorena and Désidéri, Jean-Antoine and Habbal, Abderrahmane},
  location = {Cham},
  publisher = {Springer International Publishing},
  url = {https://tuomov.iki.fi/m/nonlinearADMM.pdf},
  booktitle = {System Modeling and Optimization: 27th IFIP TC 7 Conference, CSMO 2015, Sophia Antipolis, France, June 29--July 3, 2015, Revised Selected Papers},
  date = {2016},
  doi = {10.1007/978-3-319-55795-3_10},
  eprint = {1511.00425},
  eprinttype = {arXiv},
  pages = {117--126},
  title = {Preconditioned {ADMM} with nonlinear operator constraint},
}

@article{Bidaut:1975,
  author = {Bidaut, Marie-Françoise},
  date = {1975},
  journaltitle = {C. R. Acad. Sci. Paris Sér. A},
  number = {9},
  pages = {A273--A276},
  title = {Un problème de contrôle optimal à fonction coût en norme {$L^{1}$}},
  volume = {281},
}

@article{bigolin2014historical,
  author = {Bigolin, Francesco and Golo, Sebastiano Nicolussi},
  date = {2014},
  doi = {10.1016/j.jmaa.2013.10.035},
  journaltitle = jmaa,
  number = {1},
  pages = {63--76},
  title = {A historical account on characterizations of {$C^1$}-manifolds in {E}uclidean spaces by tangent cones},
  volume = {412},
}

@article{bolte2017errorbounds,
  author = {Bolte, Jérôme and Nguyen, Trong Phong and Peypouquet, Juan and Suter, Bruce W.},
  date = {2017},
  doi = {10.1007/s10107-016-1091-6},
  journaltitle = mathprog,
  number = {2},
  pages = {471--507},
  title = {From error bounds to the complexity of first-order descent methods for convex functions},
  volume = {165},
}

@article{borwein1987smooth,
  author = {Borwein, Jonathan M. and Preiss, David},
  date = {1987},
  doi = {10.2307/2000681},
  journaltitle = transams,
  pages = {517--527},
  title = {A smooth variational principle with applications to subdifferentiability and to differentiability of convex functions},
  volume = {303},
}

@book{BorweinZhu:2005,
  author = {Borwein, Jonathan M. and Zhu, Qiji J.},
  publisher = {Springer-Verlag, New York},
  date = {2005},
  doi = {10.1007/0-387-28271-8},
  series = {CMS Books in Mathematics/Ouvrages de Mathématiques de la SMC},
  title = {Techniques of Variational Analysis},
  volume = {20},
}

@article{bouligand1930,
  author = {Bouligand, Georges},
  date = {1930},
  journaltitle = {Rev. Gén. des Sciences},
  pages = {39--43},
  title = {Sur quelques points de méthodologie géométrique},
  volume = {41},
}

@book{Boyd:2004,
  author = {Boyd, Stephen and Vandenberghe, Lieven},
  publisher = {Cambridge University Press, Cambridge},
  date = {2004},
  doi = {10.1017/CBO9780511804441},
  title = {Convex Optimization},
}

@book{Braides:2002,
  author = {Braides, Andrea},
  publisher = {Oxford University Press, Oxford},
  date = {2002},
  doi = {10.1093/acprof:oso/9780198507840.001.0001},
  series = {Oxford Lecture Series in Mathematics and its Applications},
  title = {{$\Gamma$}-Convergence for Beginners},
  volume = {22},
}

@article{Bredies:2022,
  author = {Bredies, Kristian and Chenchene, Enis and Lorenz, Dirk A. and Naldi, Emanuele},
  date = {2022},
  doi = {10.1137/21M1448112},
  journaltitle = siamjopt,
  number = {3},
  pages = {2376--2401},
  title = {Degenerate preconditioned proximal point algorithms},
  volume = {32},
}

@book{BrediesLorenz:2018,
  author = {Bredies, Kristian and Lorenz, Dirk},
  publisher = {Birkhäuser/Springer, Cham},
  date = {2018},
  doi = {10.1007/978-3-030-01458-2},
  series = {Applied and Numerical Harmonic Analysis},
  title = {Mathematical Image Processing},
}

@article{bredies2008linear,
  author = {Bredies, Kristian and Lorenz, Dirk A.},
  date = {2008},
  doi = {10.1007/s00041-008-9041-1},
  journaltitle = {Journal of Fourier Analysis and Applications},
  number = {5},
  pages = {813--837},
  title = {Linear convergence of iterative soft-thresholding},
  volume = {14},
}

@misc{bredies2016accelerated,
  author = {Bredies, Kristian and Sun, Hongpeng},
  date = {2016},
  eprint = {1604.06282},
  eprinttype = {arXiv},
  note = {Preprint},
  title = {Accelerated {D}ouglas--{R}achford methods for the solution of convex-concave saddle-point problems},
}

@article{BrezisCrandallPazy,
  author = {Brezis, Haı̈ and Crandall, M. G. and Pazy, A.},
  date = {1970},
  doi = {10.1002/cpa.3160230107},
  journaltitle = cpam,
  pages = {123--144},
  title = {Perturbations of nonlinear maximal monotone sets in {B}anach space},
  volume = {23},
}

@book{Brezis:2010a,
  author = {Brezis, Haı̈m},
  location = {New York},
  publisher = {Springer},
  date = {2010},
  doi = {10.1007/978-0-387-70914-7},
  title = {Functional Analysis, Sobolev Spaces and Partial Differential Equations},
}

@booklet{Brokate,
  author = {Brokate, Martin},
  url = {https://mediatum.ub.tum.de/doc/1701116/1701116.pdf},
  date = {2014},
  howpublished = {Zentrum Mathematik, TU München},
  title = {Konvexe Analysis und Evolutionsprobleme},
  type = {Lecture notes},
}

@article{browder1965nonexpansive,
  author = {Browder, Felix E},
  publisher = {National Academy of Sciences},
  date = {1965},
  doi = {10.1073/pnas.54.4.1041},
  journaltitle = {Proceedings of the National Academy of Sciences of the United States of America},
  number = {4},
  pages = {1041},
  title = {Nonexpansive nonlinear operators in a Banach space},
  volume = {54},
}

@article{browder1967convergence,
  author = {Browder, Felix E.},
  date = {1967},
  doi = {10.1007/BF01109805},
  journaltitle = {Mathematische Zeitschrift},
  number = {3},
  pages = {201--225},
  title = {Convergence theorems for sequences of nonlinear operators in {Banach} spaces},
  volume = {100},
}

@article{CasasClasonKunisch:2012,
  author = {Casas, Eduardo and Clason, Christian and Kunisch, Karl},
  date = {2012},
  doi = {10.1137/110843216},
  journaltitle = siamjcon,
  number = {4},
  pages = {1735--1752},
  title = {Approximation of elliptic control problems in measure spaces with sparse solutions},
  volume = {50},
}

@book{Cegielski:2012,
  author = {Cegielski, Andrzej},
  publisher = {Springer, Heidelberg},
  date = {2012},
  doi = {10.1007/978-3-642-30901-4},
  series = {Lecture Notes in Mathematics},
  title = {Iterative Methods for Fixed Point Problems in {H}ilbert Spaces},
  volume = {2057},
}

@article{chambolledevore1998nonlinear,
  author = {Chambolle, Antonin and DeVore, Ronald A. and Lee, Nam-yong and Lucier, Bradley J.},
  date = {1998},
  doi = {10.1109/83.661182},
  journaltitle = ieeeimg,
  number = {3},
  pages = {319--335},
  title = {Nonlinear wavelet image processing: variational problems, compression, and noise removal through wavelet shrinkage},
  volume = {7},
}

@article{chambolle2017stochastic,
  author = {Chambolle, Antonin and Ehrhardt, M. and Richtárik, P. and Schönlieb, C.},
  date = {2018},
  doi = {10.1137/17M1134834},
  journaltitle = siamjopt,
  number = {4},
  pages = {2783--2808},
  title = {Stochastic primal-dual hybrid gradient algorithm with arbitrary sampling and imaging applications},
  volume = {28},
}

@article{Pock_PD_2010,
  author = {Chambolle, Antonin and Pock, Thomas},
  date = {2011},
  doi = {10.1007/s10851-010-0251-1},
  journaltitle = jmiv,
  number = {1},
  pages = {120--145},
  title = {A first-order primal-dual algorithm for convex problems with applications to imaging},
  volume = {40},
}

@article{chambolle2014ergodic,
  author = {Chambolle, Antonin and Pock, Thomas},
  date = {2015},
  doi = {10.1007/s10107-015-0957-3},
  journaltitle = mathprog,
  pages = {1--35},
  title = {On the ergodic convergence rates of a first-order primal--dual algorithm},
}

@article{chen2013pdfp,
  author = {Chen, Peijun and Huang, Jianguo and Zhang, Xiaoqun},
  date = {2013},
  doi = {10.1088/0266-5611/29/2/025011},
  journaltitle = invprob,
  number = {2},
  pages = {025011},
  title = {A primal-dual fixed point algorithm for convex separable minimization with applications to image restoration},
  volume = {29},
}

@article{Chen:2000a,
  author = {Chen, Xiaojun and Nashed, Zuhair and Qi, Liqun},
  date = {2000},
  doi = {10.1137/s0036142999356719},
  journaltitle = siamjnum,
  number = {4},
  pages = {1200--1216},
  title = {Smoothing methods and semismooth methods for nondifferentiable operator equations},
  volume = {38},
}

@article{CCMW:2017,
  author = {Christof, Constantin and Clason, Christian and Meyer, Christian and Walter, Stefan},
  date = {2018},
  doi = {10.3934/mcrf.2018011},
  journaltitle = {Mathematical Control and Related Fields},
  number = {1},
  pages = {247--276},
  title = {Optimal control of a non-smooth semilinear elliptic equation},
  volume = {8},
}

@article{christof2018nogap,
  author = {Christof, Constantin and Wachsmuth, Gerd},
  date = {2018},
  doi = {10.1137/17M1140418},
  journaltitle = siamjopt,
  number = {3},
  pages = {2097--2130},
  title = {No-gap second-order conditions via a directional curvature functional},
  volume = {28},
}

@article{cibulka2018strong,
  author = {Cibulka, R. and Dontchev, A.L. and Kruger, A.Y.},
  date = {2018},
  doi = {10.1016/j.jmaa.2016.11.045},
  journaltitle = jmaa,
  note = {Special Issue on Convex Analysis and Optimization: New Trends in Theory and Applications},
  number = {2},
  pages = {1247--1282},
  title = {Strong metric subregularity of mappings in variational analysis and optimization},
  volume = {457},
}

@book{Cioranescu,
  author = {Cioranescu, Ioana},
  location = {Dordrecht},
  publisher = {Springer},
  date = {1990},
  doi = {10.1007/978-94-009-2121-4},
  series = {Mathematics and Its Applications},
  title = {Geometry of {Banach} Spaces, Duality Mappings and Nonlinear Problems},
  volume = {62},
}

@book{Clarke:2013,
  author = {Clarke, Francis},
  location = {London},
  publisher = {Springer},
  date = {2013},
  doi = {10.1007/978-1-4471-4820-3},
  title = {Functional Analysis, Calculus of Variations and Optimal Control},
}

@thesis{clarke1973necessary,
  author = {Clarke, Frank H},
  institution = {University of Washington},
  date = {1973},
  title = {Necessary Conditions for Nonsmooth Problems in Optimal Control and the Calculus of Variations},
  type = {PhD thesis},
}

@article{clarke1975generalized,
  author = {Clarke, Frank H},
  date = {1975},
  doi = {10.1090/s0002-9947-1975-0367131-6},
  journaltitle = transams,
  pages = {247--262},
  title = {Generalized gradients and applications},
  volume = {205},
}

@book{Clarke:1990a,
  author = {Clarke, Frank H.},
  location = {Philadelphia},
  publisher = {Society for Industrial and Applied Mathematics},
  date = {1990},
  doi = {10.1137/1.9781611971309},
  series = {Classics Appl. Math.},
  title = {Optimization and {N}onsmooth {A}nalysis},
  volume = {5},
}

@book{Clason,
  author = {Clason, Christian},
  location = {Basel},
  publisher = {Springer International Publishing},
  date = {2020},
  doi = {10.1007/978-3-030-52784-6},
  series = {Compact Textbooks in Mathematics},
  title = {Introduction to Functional Analysis},
}

@booklet{ClasonIP:2020,
  author = {Clason, Christian},
  date = {2020},
  eprint = {2001.00617},
  eprinttype = {arxiv},
  title = {{Regularization of Inverse Problems}},
  type = {Lecture notes},
}

@incollection{ClasonDo:2018,
  author = {Clason, Christian and Do, Thi Bich Tram},
  editor = {Hofmann, Bernd and Leitão, Antonio and Zubelli, Jorge},
  publisher = {Springer},
  booktitle = {New Trends in Parameter Identification for Mathematical Models},
  date = {2018},
  doi = {10.1007/978-3-319-70824-9_2},
  eprint = {1707.01041},
  eprinttype = {arxiv},
  pages = {31--51},
  series = {Trends in Mathematics},
  title = {Convex regularization of discrete-valued inverse problems},
}

@article{ClasonItoKunisch:2016,
  author = {Clason, Christian and Ito, Kazufumi and Kunisch, Karl},
  date = {2016},
  doi = {10.1051/cocv/2015017},
  eprint = {1702.07540},
  eprinttype = {arxiv},
  journaltitle = cocv,
  number = {2},
  pages = {581--609},
  title = {A convex analysis approach to optimal controls with switching structure for partial differential equations},
  volume = {22},
}

@article{CJK:2010,
  author = {Clason, Christian and Jin, Bangti and Kunisch, Karl},
  date = {2010},
  doi = {10.1137/090758003},
  journaltitle = siamjimg,
  number = {2},
  pages = {199--231},
  title = {A semismooth {N}ewton method for {L$^1$} data fitting with automatic choice of regularization parameters and noise calibration},
  volume = {3},
}

@article{Clason:2010a,
  author = {Clason, Christian and Kunisch, Karl},
  date = {2011},
  doi = {10.1051/cocv/2010003},
  journaltitle = cocv,
  number = {1},
  pages = {243--266},
  title = {A duality-based approach to elliptic control problems in non-reflexive {B}anach spaces},
  volume = {17},
}

@article{ClasonKunisch:2013,
  author = {Clason, Christian and Kunisch, Karl},
  date = {2014},
  doi = {10.1016/j.anihpc.2013.08.005},
  journaltitle = {Annales de l'Institut Henri Poincaré (C) Analyse Non Linéaire},
  number = {6},
  pages = {1109--1130},
  title = {Multi-bang control of elliptic systems},
  volume = {31},
}

@article{ClasonKunisch:2016,
  author = {Clason, Christian and Kunisch, Karl},
  date = {2016},
  doi = {10.1051/m2an/2016012},
  eprint = {1702.07525},
  eprinttype = {arxiv},
  journaltitle = {ESAIM: Mathematical Modelling and Numerical Analysis},
  number = {6},
  pages = {1917--1936},
  title = {A convex analysis approach to multi-material topology optimization},
  volume = {50},
}

@article{ClasonKunischTrautmann:2019,
  author = {Clason, Christian and Kunisch, Karl and Trautmann, Philip},
  date = {2021},
  doi = {10.1007/s00245-020-09733-9},
  eprint = {1912.08672},
  eprinttype = {arxiv},
  journaltitle = amo,
  number = {3},
  pages = {2889--2921},
  title = {Optimal control of the principal coefficient in a scalar wave equation},
  volume = {84},
}

@article{tuomov-nlpdhgm-redo,
  author = {Clason, Christian and Mazurenko, Stanislav and Valkonen, Tuomo},
  date = {2019},
  doi = {10.1137/18M1170194},
  eprint = {1802.03347},
  eprinttype = {arXiv},
  issue = {1},
  journaltitle = siamjopt,
  origdate = {2018},
  pages = {933--963},
  title = {Acceleration and global convergence of a first-order primal--dual method for nonconvex problems},
  volume = {29},
}

@article{tuomov-nlpdhgm-general,
  author = {Clason, Christian and Mazurenko, Stanislav and Valkonen, Tuomo},
  url = {https://tuomov.iki.fi/m/nlpdhgm_general.pdf},
  date = {2020},
  doi = {10.1007/s00245-020-09676-1},
  eprint = {1901.02746},
  eprinttype = {arXiv},
  journaltitle = amo,
  number = {2},
  pages = {1239--1284},
  title = {Primal–dual proximal splitting and generalized conjugation in non-smooth non-convex optimization},
  volume = {84},
}

@article{ClasonSchiela:2015,
  author = {Clason, Christian and Schiela, Anton},
  date = {2017},
  doi = {10.1051/cocv/2015046},
  eprint = {1702.07528},
  eprinttype = {arxiv},
  journaltitle = cocv,
  number = {1},
  pages = {217--240},
  title = {Optimal control of elliptic equations with positive measures},
  volume = {23},
}

@article{ClasonTamelingWirth:2021,
  author = {Clason, Christian and Tameling, Carla and Wirth, Benedikt},
  date = {2021},
  doi = {10.1137/21M1426237},
  eprint = {2108.10077},
  eprinttype = {arxiv},
  journaltitle = {SIAM Review},
  number = {4},
  pages = {783--821},
  title = {Convex relaxation of discrete vector-valued optimization problems},
  volume = {63},
}

@article{tuomov-pdex2nlpdhgm,
  author = {Clason, Christian and Valkonen, Tuomo},
  date = {2017},
  doi = {10.1137/16M1080859},
  eprint = {1606.06219},
  eprinttype = {arXiv},
  journaltitle = siamjopt,
  number = {3},
  pages = {1314--1339},
  title = {Primal-dual extragradient methods for nonlinear nonsmooth {PDE}-constrained optimization},
  volume = {27},
}

@article{tuomov-pdex2stability,
  author = {Clason, Christian and Valkonen, Tuomo},
  url = {https://tuomov.iki.fi/m/pdex2_stability.pdf},
  date = {2017},
  doi = {10.1007/s11228-016-0366-7},
  eprint = {1509.06582},
  eprinttype = {arXiv},
  issue = {1},
  journaltitle = svaa,
  pages = {69--112},
  title = {Stability of saddle points via explicit coderivatives of pointwise subdifferentials},
  volume = {25},
}

@misc{nonsmoothbook-codes,
  author = {Clason, Christian and Valkonen, Tuomo},
  date = {2025},
  doi = {10.5281/zenodo.13128264},
  howpublished = {Software on Zenodo},
  title = {Code accompanying ``{Introduction to Nonsmooth Analysis and Optimization}''},
  version = {1.0.1},
}

@article{Combettes2013,
  author = {Combettes, Patrick L. and Reyes, Noli N.},
  date = {2013},
  doi = {10.1007/s10107-013-0663-y},
  journaltitle = mathprog,
  number = {1},
  pages = {103--114},
  title = {Moreau's decomposition in {Banach} spaces},
  volume = {139},
}

@article{condat2013primaldual,
  author = {Condat, Laurent},
  date = {2013},
  doi = {10.1007/s10957-012-0245-9},
  journaltitle = jota,
  number = {2},
  pages = {460--479},
  title = {A primal--dual splitting method for convex optimization involving {Lipschitzian}, proximable and linear composite terms},
  volume = {158},
}

@book{CuiPang:2021,
  author = {Cui, Ying and Pang, Jong-Shi},
  location = {Philadelphia, PA},
  publisher = {Society for Industrial and Applied Mathematics},
  date = {2021},
  doi = {10.1137/1.9781611976748},
  title = {Modern Nonconvex Nondifferentiable Optimization},
}

@book{DalMaso:1993,
  author = {Dal Maso, Gianni},
  publisher = {Birkhäuser Boston, Inc., Boston, MA},
  date = {1993},
  doi = {10.1007/978-1-4612-0327-8},
  series = {Progress in Nonlinear Differential Equations and Their Applications},
  title = {An Introduction to {$\Gamma$}-Convergence},
  volume = {8},
}

@article{daubechiesdefriesdemol2004ista,
  author = {Daubechies, I. and Defrise, M. and De Mol, C.},
  date = {2004},
  doi = {10.1002/cpa.20042},
  journaltitle = cpam,
  number = {11},
  pages = {1413--1457},
  title = {An iterative thresholding algorithm for linear inverse problems with a sparsity constraint},
  volume = {57},
}

@article{Davis:2017,
  author = {Davis, Damek and Yin, Wotao},
  publisher = {Springer Science and Business Media LLC},
  date = {2017},
  doi = {10.1007/s11228-017-0421-z},
  issn = {1877-0541},
  journaltitle = svaa,
  number = {4},
  pages = {829--858},
  title = {A three-operator splitting scheme and its optimization applications},
  volume = {25},
}

@book{DelosReyes,
  author = {De los Reyes, Juan Carlos},
  publisher = {Springer},
  date = {2015},
  doi = {10.1007/978-3-319-13395-9},
  title = {{Numerical PDE-Constrained Optimization}},
}

@book{DiBenedetto,
  author = {DiBenedetto, Emmanuele},
  publisher = {Birkhäuser Boston, Inc., Boston, MA},
  date = {2002},
  doi = {10.1007/978-1-4612-0117-5},
  title = {Real Analysis},
}

@article{dolecki2011tangency,
  author = {Dolecki, S. and Greco, G. H.},
  url = {https://www.heldermann.de/JCA/JCA18/JCA182/jca18017.htm},
  date = {2011},
  eprint = {1003.1332},
  eprinttype = {arXiv},
  journaltitle = convexanal,
  number = {2},
  pages = {301--339},
  title = {Tangency vis-à-vis differentiability by {P}eano, {S}everi and {G}uareschi},
  volume = {18},
}

@book{dontchev2014implicit,
  author = {Dontchev, A. L. and Rockafellar, R. T.},
  publisher = {Springer New York},
  date = {2014},
  doi = {10.1007/978-1-4939-1037-3},
  edition = {2},
  series = {Springer Series in Operations Research and Financial Engineering},
  subtitle = {A View from Variational Analysis},
  title = {Implicit Functions and Solution Mappings},
}

@book{Dontchev:2021,
  author = {Dontchev, Asen L.},
  publisher = {Cham: Springer},
  date = {2021},
  doi = {10.1007/978-3-030-79911-3},
  series = {Applied Mathematical Sciences},
  title = {Lectures on Variational Analysis},
  volume = {205},
}

@article{dontchev2004regularity,
  author = {Dontchev, Asen L. and Rockafellar, R. Tyrrell},
  date = {2004},
  doi = {10.1023/B:SVAN.0000023394.19482.30},
  journaltitle = svaa,
  number = {1-2},
  pages = {79--109},
  title = {Regularity and conditioning of solution mappings in variational analysis},
  volume = {12},
}

@article{douglas1956numerical,
  author = {Douglas, Jim, Jr. and Rachford, H. H., Jr.},
  publisher = {American Mathematical Society},
  date = {1956},
  doi = {10.2307/1993056},
  journaltitle = transams,
  number = {2},
  pages = {421--439},
  title = {On the numerical solution of heat conduction problems in two and three space variables},
  volume = {82},
}

@article{drori2015simple,
  author = {Drori, Yoel and Sabach, Shoham and Teboulle, Marc},
  date = {2015},
  doi = {10.1016/j.orl.2015.02.001},
  journaltitle = {Operations Research Letters},
  number = {2},
  pages = {209--214},
  title = {A simple algorithm for a class of nonsmooth convex--concave saddle-point problems},
  volume = {43},
}

@article{Eckstein:1992,
  author = {Eckstein, Jonathan and Bertsekas, Dimitri P.},
  date = {1992},
  doi = {10.1007/bf01581204},
  journaltitle = mathprog,
  number = {1-3},
  pages = {293--318},
  title = {On the {D}ouglas--{R}achford splitting method and the proximal point algorithm for maximal monotone operators},
  volume = {55},
}

@article{ekeland1976generic,
  author = {Ekeland, Ivar and Lebourg, Gérard},
  date = {1976},
  doi = {10.1090/S0002-9947-1976-0431253-2},
  journaltitle = transams,
  number = {2},
  pages = {193--216},
  title = {Generic {Fréchet}-differentiability and perturbed optimization problems in {Banach} spaces},
  volume = {224},
}

@book{Ekeland:1999a,
  author = {Ekeland, Ivar and Témam, Roger},
  location = {Philadelphia},
  publisher = {Society for Industrial and Applied Mathematics},
  date = {1999},
  doi = {10.1137/1.9781611971088},
  series = {Classics Appl. Math.},
  title = {Convex {A}nalysis and {V}ariational {P}roblems},
  volume = {28},
}

@book{Engl,
  author = {Engl, Heinz W. and Hanke, Martin and Neubauer, Andreas},
  publisher = {Springer Netherlands},
  date = {1996},
  series = {Mathematics and its Applications},
  title = {Regularization of Inverse Problems},
  volume = {375},
}

@article{esser2010general,
  author = {Esser, Ernie and Zhang, Xiaoqun and Chan, Tony F},
  date = {2010},
  doi = {10.1137/09076934X},
  journaltitle = siamjimg,
  number = {4},
  pages = {1015--1046},
  title = {A general framework for a class of first order primal-dual algorithms for convex optimization in imaging science},
  volume = {3},
}

@article{fabian1988classes,
  author = {Fabian, Marián},
  date = {1988},
  doi = {10.1016/0362-546X(88)90013-2},
  journaltitle = {Nonlinear Analysis: Theory, Methods \& Applications},
  number = {1},
  pages = {63--74},
  title = {On classes of subdifferentiability spaces of {Ioffe}},
  volume = {12},
}

@inbook{fabian2001differentiability,
  author = {Fabian, Marián and Habala, Petr and Hájek, Petr and Santalucía, Vicente Montesinos and Pelant, Jan and Zizler, Václav},
  location = {New York, NY},
  publisher = {Springer New York},
  booktitle = {Functional Analysis and Infinite-Dimensional Geometry},
  date = {2001},
  doi = {10.1007/978-1-4757-3480-5_8},
  pages = {241--284},
  title = {Differentiability of norms},
}

@book{facchinei2003finite-1,
  author = {Facchinei, Francisco and Pang, Jong-Shi},
  publisher = {Springer-Verlag, New York},
  date = {2003},
  doi = {10.1007/b97543},
  series = {Springer Series in Operations Research},
  title = {Finite-dimensional Variational Inequalities and Complementarity Problems. {V}ol. {I}},
}

@book{facchinei2003finite-2,
  author = {Facchinei, Francisco and Pang, Jong-Shi},
  publisher = {Springer-Verlag, New York},
  date = {2003},
  doi = {10.1007/b97544},
  series = {Springer Series in Operations Research},
  title = {Finite-dimensional Variational Inequalities and Complementarity Problems. {V}ol. {II}},
}

@article{fejer1922,
  author = {Fejér, Leopold},
  date = {1922},
  doi = {10.1007/BF01449600},
  journaltitle = {Mathematische Annalen},
  number = {1},
  pages = {41--48},
  title = {{Über die Lage der Nullstellen von Polynomen, die aus Minimumforderungen gewisser Art entspringen}},
  volume = {85},
}

@book{fonseca2007mmc,
  author = {Fonseca, I. and Leoni, G.},
  publisher = {Springer},
  date = {2007},
  doi = {10.1007/978-0-387-69006-3},
  title = {Modern Methods in the Calculus of Variations: {$L^p$} Spaces},
}

@incollection{gabay,
  author = {Gabay, D.},
  editor = {Fortin, M. and Glowinski, R.},
  location = {Amsterdam},
  publisher = {North-Holland},
  booktitle = {Augmented {Lagrangian} Methods: Applications to the Numerical Solution of Boundary-Value Problems},
  date = {1983},
  doi = {10.1016/S0168-2024(08)70034-1},
  pages = {299--331},
  series = {Studies in Mathematics and its Applications},
  title = {Applications of the method of multipliers to variational inequalities},
  volume = {15},
}

@article{garrigos2017thresholding,
  author = {Garrigos, Guillaume and Rosasco, Lorenzo and Villa, Silvia},
  date = {2020},
  doi = {10.1051/cocv/2019011},
  journaltitle = cocv,
  pages = {28},
  title = {Thresholding gradient methods in {Hilbert} spaces: support identification and linear convergence},
  volume = {26},
}

@book{GelmanCarlinSternRubin:2013,
  author = {Gelman, Andrew and Carlin, John B. and Stern, Hal S. and Rubin, Donald B.},
  location = {Boca Raton},
  publisher = {Chapman \& Hall/CRC},
  url = {https://www.stat.columbia.edu/~gelman/book/},
  date = {2013},
  doi = {10.1201/9780429258480},
  edition = {3},
  title = {Bayesian Data Analysis},
}

@online{Gerd:2022,
  author = {gerw},
  organization = {MathOverflow},
  url = {https://mathoverflow.net/q/414752},
  date = {2022},
  title = {Subgradient in a predual under weak* continuity},
  urldate = {2022-01-26},
}

@article{outrata2019semismooth,
  author = {Gferer, Helmut and Outrata, Jiři},
  date = {2021},
  doi = {10.1137/19M1257408},
  journaltitle = siamjopt,
  number = {1},
  pages = {489--517},
  title = {On a semismooth$^*$ {Newton} method for solving generalized equations},
  volume = {31},
}

@article{gfrerer2016lipschitzian,
  author = {Gfrerer, H. and Outrata, J.},
  date = {2016},
  doi = {10.1137/15M1052299},
  journaltitle = siamjopt,
  number = {4},
  pages = {2160--2189},
  title = {On {Lipschitzian} properties of implicit multifunctions},
  volume = {26},
}

@article{gfrerer2011first,
  author = {Gfrerer, Helmut},
  date = {2011},
  doi = {10.1137/100813415},
  journaltitle = siamjopt,
  number = {4},
  pages = {1439--1474},
  title = {First order and second order characterizations of metric subregularity and calmness of constraint set mappings},
  volume = {21},
}

@article{gfrerer2013directional,
  author = {Gfrerer, Helmut},
  date = {2013},
  doi = {10.1007/s11228-012-0220-5},
  journaltitle = svaa,
  number = {2},
  pages = {151--176},
  title = {On directional metric regularity, subregularity and optimality conditions for nonsmooth mathematical programs},
  volume = {21},
}

@article{gribonval2019characterization,
  author = {Gribonval, R. and Nikolova, M.},
  date = {2020},
  doi = {10.1007/s10851-020-00951-y},
  journaltitle = jmiv,
  pages = {773--789},
  title = {A characterization of proximity operators},
  volume = {62},
}

@article{Lorenz:2008,
  author = {Griesse, R. and Lorenz, D. A.},
  publisher = {IOP Publishing},
  date = {2008},
  doi = {10.1088/0266-5611/24/3/035007},
  issn = {1361--6420},
  journaltitle = invprob,
  number = {3},
  pages = {035007},
  title = {A semismooth Newton method for {T}ikhonov functionals with sparsity constraints},
  volume = {24},
}

@book{Grisvard:2011,
  author = {Grisvard, Pierre},
  location = {Philadelphia},
  publisher = {Society for Industrial and Applied Mathematics},
  date = {2011},
  doi = {10.1137/1.9781611972030},
  title = {Elliptic Problems in Nonsmooth Domains},
}

@article{han2013local,
  author = {Han, Deren and Yuan, Xiaoming},
  date = {2013},
  doi = {10.1137/120886753},
  journaltitle = siamjnum,
  number = {6},
  pages = {3446--3457},
  title = {Local linear convergence of the Alternating Direction Method of Multipliers for quadratic programs},
  volume = {51},
}

@book{Hanke:2017,
  author = {Hanke, M.},
  publisher = {Society for Industrial and Applied Mathematics},
  date = {2017},
  doi = {10.1137/1.9781611974942},
  isbn = {9781611974942},
  series = {Other Titles in Applied Mathematics},
  title = {A Taste of Inverse Problems: Basic Theory and Examples},
}

@article{harderwachsmuth2018sobolev,
  author = {Harder, Felix and Wachsmuth, Gerd},
  date = {2018},
  doi = {10.1080/02331934.2018.1484467},
  journaltitle = {Optimization},
  number = {5},
  pages = {1579--1603},
  title = {The limiting normal cone of a complementarity set in {Sobolev} spaces},
  volume = {67},
}

@article{He:2012,
  author = {He, B. and Yuan, X.},
  date = {2012},
  doi = {10.1137/100814494},
  journaltitle = siamjimg,
  number = {1},
  pages = {119--149},
  title = {Convergence analysis of primal-dual algorithms for a saddle-point problem: from contraction perspective},
  volume = {5},
}

@booklet{Heinonen,
  author = {Heinonen, Juha},
  publisher = {University of Jyväskylä},
  url = {https://urn.fi/URN:ISBN:951-39-2318-5},
  date = {2005},
  series = {Rep. Univ. Jyväskylä Dept. Math. Stat.},
  title = {Lectures on {L}ipschitz Analysis},
  volume = {100},
}

@article{henrion2002calmness,
  author = {Henrion, R. and Jourani, A. and Outrata, J.},
  date = {2002},
  doi = {10.1137/S1052623401395553},
  journaltitle = siamjopt,
  number = {2},
  pages = {603--618},
  title = {On the calmness of a class of multifunctions},
  volume = {13},
}

@article{Hintermuller:2002a,
  author = {Hintermüller, M. and Ito, K. and Kunisch, K.},
  date = {2002},
  doi = {10.1137/s1052623401383558},
  journaltitle = siamjopt,
  number = {3},
  pages = {865--888 (2003)},
  title = {The primal-dual active set strategy as a semismooth {N}ewton method},
  volume = {13},
}

@article{Hintermuller:2004a,
  author = {Hintermüller, M. and Kunisch, K.},
  date = {2004},
  doi = {10.1137/S0036139903422784},
  journaltitle = siamjapl,
  number = {4},
  pages = {1311--1333},
  title = {Total bounded variation regularization as a bilaterally constrained optimization problem},
  volume = {64},
}

@article{Hintermuller:2006a,
  author = {Hintermüller, Michael and Stadler, Georg},
  date = {2006},
  doi = {10.1137/040613263},
  journaltitle = {SIAM J. Sci. Comput.},
  number = {1},
  pages = {1--23 (electronic)},
  title = {An infeasible primal-dual algorithm for total bounded variation-based inf-convolution-type image restoration},
  volume = {28},
}

@article{Hintermuller:2004b,
  author = {Hintermüller, Michael and Ulbrich, Michael},
  date = {2004},
  doi = {10.1007/s10107-004-0540-9},
  journaltitle = {Math. Program.},
  number = {1, Ser. B},
  pages = {151--184},
  title = {A mesh-independence result for semismooth {N}ewton methods},
  volume = {101},
}

@book{Hinze2009,
  author = {Hinze, M. and Pinnau, R. and Ulbrich, M. and Ulbrich, S.},
  location = {New York},
  publisher = {Springer},
  date = {2009},
  doi = {10.1007/978-1-4020-8839-1},
  title = {Optimization with {PDE} Constraints},
}

@book{lemarechal1993-1,
  author = {Hiriart-Urruty, Jean-Baptiste and Lemaréchal, Claude},
  publisher = {Springer-Verlag, Berlin},
  date = {1993},
  doi = {10.1007/978-3-662-02796-7},
  series = {Grundlehren der Mathematischen Wissenschaften},
  subtitle = {Fundamentals},
  title = {Convex Analysis and Minimization Algorithms {I}},
  volume = {305},
}

@book{lemarechal1993-2,
  author = {Hiriart-Urruty, Jean-Baptiste and Lemaréchal, Claude},
  publisher = {Springer-Verlag, Berlin},
  date = {1993},
  doi = {10.1007/978-3-662-06409-2},
  series = {Grundlehren der Mathematischen Wissenschaften},
  subtitle = {Advanced Theory and Bundle Methods},
  title = {Convex Analysis and Minimization Algorithms {II}},
  volume = {306},
}

@book{Hiriart:2001,
  author = {Hiriart-Urruty, Jean-Baptiste and Lemaréchal, Claude},
  location = {Berlin},
  publisher = {Springer-Verlag},
  date = {2001},
  doi = {10.1007/978-3-642-56468-0},
  title = {Fundamentals of Convex Analysis},
}

@misc{hohage2014generalization,
  author = {Hohage, Thorsten and Homann, Carolin},
  date = {2014},
  eprint = {1412.0126},
  eprinttype = {arXiv},
  note = {Preprint},
  title = {A generalization of the {C}hambolle--{P}ock algorithm to {B}anach spaces with applications to inverse problems},
}

@article{Hoermander:1955,
  author = {Hörmander, Lars},
  date = {1955},
  doi = {10.1007/BF02589354},
  journaltitle = {Arkiv för Matematik},
  pages = {181--186},
  title = {Sur la fonction d'appui des ensembles convexes dans un espace localement convexe},
  volume = {3},
}

@book{Huber:2009,
  author = {Huber, Peter J.},
  location = {New York},
  publisher = {John Wiley \& Sons Inc.},
  date = {2009},
  doi = {10.1002/9780470434697},
  edition = {2},
  title = {Robust Statistics},
}

@article{ioffe1984approximate,
  author = {Ioffe, A. D.},
  date = {1984},
  doi = {10.2307/1999541},
  journaltitle = transams,
  number = {1},
  pages = {389--416},
  title = {Approximate subdifferentials and applications. {I}. {T}he finite-dimensional theory},
  volume = {281},
}

@book{ioffe2017variational,
  author = {Ioffe, A. D.},
  publisher = {Springer International Publishing},
  date = {2017},
  doi = {10.1007/978-3-319-64277-2},
  series = {Springer Monographs in Mathematics},
  title = {Variational Analysis of Regular Mappings: Theory and Applications},
}

@article{ioffe1979regular,
  author = {Ioffe, Alexander D},
  date = {1979},
  doi = {10.1090/S0002-9947-1979-0531969-6},
  journaltitle = transams,
  pages = {61--69},
  title = {Regular points of {Lipschitz} functions},
  volume = {251},
}

@book{ItoJin,
  author = {Ito, Kazufumi and Jin, Bangti},
  location = {Singapore},
  publisher = {World Scientific},
  date = {2014},
  doi = {10.1142/9120},
  series = {Series on Applied Mathematics},
  subtitle = {Tikhonov Theory and Algorithms},
  title = invprob,
  volume = {22},
}

@book{Kunisch:2008a,
  author = {Ito, Kazufumi and Kunisch, Karl},
  location = {Philadelphia},
  publisher = {Society for Industrial and Applied Mathematics},
  date = {2008},
  doi = {10.1137/1.9780898718614},
  series = {Advances in Design and Control},
  title = {Lagrange {M}ultiplier {A}pproach to {V}ariational {P}roblems and {A}pplications},
  volume = {15},
}

@article{jauhiainen2019gaussnewton,
  author = {Jauhiainen, Jyrki and Kuusela, Petri and Seppänen, Aku and Valkonen, Tuomo},
  date = {2020},
  doi = {10.1137/20M1321711},
  eprint = {2002.08044},
  eprinttype = {arxiv},
  journaltitle = siamjimg,
  number = {3},
  pages = {1415--1445},
  title = {Relaxed {G}auss--{N}ewton methods with applications to electrical impedance tomography},
  volume = {13},
}

@book{Kaltenbacher:2008,
  author = {Kaltenbacher, B. and Neubauer, A. and Scherzer, O.},
  publisher = {De Gruyter},
  date = {2008},
  isbn = {9783110208276},
  number = {6},
  series = {Radon Series on Computational and Applied Mathematics},
  title = {Iterative Regularization Methods for Nonlinear Ill-Posed Problems},
}

@article{Karkkainen:2005a,
  author = {Kärkkäinen, T. and Kunisch, K. and Majava, K.},
  date = {2005},
  doi = {10.1007/s00607-004-0097-8},
  journaltitle = {Computing},
  number = {4},
  pages = {353--376},
  title = {Denoising of smooth images using {$L^1$}-fitting},
  volume = {74},
}

@book{KlatteKummer:2002,
  author = {Klatte, Diethard and Kummer, Bernd},
  publisher = {Kluwer Academic Publishers, Dordrecht},
  date = {2002},
  doi = {10.1007/b130810},
  note = {Regularity, calculus, methods and applications},
  series = {Nonconvex Optimization and its Applications},
  title = {Nonsmooth Equations in Optimization},
  volume = {60},
}

@article{kojima1986newton,
  author = {Kojima, Masakazu and Shindo, Susumu},
  date = {1986},
  doi = {10.15807/jorsj.29.352},
  journaltitle = {Journal of the Operations Research Society of Japan},
  number = {4},
  pages = {352--375},
  title = {Extension of {N}ewton and quasi-{N}ewton methods to systems of {PC{$^1$}} equations},
  volume = {29},
}

@article{krasnoselskii1955remarks,
  author = {Krasnosel'ski\weaki, M. A.},
  date = {1955},
  journaltitle = {Uspekhi Matematicheskikh Nauk},
  number = {1(63)},
  pages = {123--127},
  title = {Two remarks on the method of successive approximations},
  volume = {10},
}

@article{kruger2015error,
  author = {Kruger, Alexander Y.},
  publisher = {Taylor \& Francis},
  date = {2015},
  doi = {10.1080/02331934.2014.938074},
  journaltitle = {Optimization},
  number = {1},
  pages = {49--79},
  title = {Error bounds and metric subregularity},
  volume = {64},
}

@article{kruse2018implicit,
  author = {Kruse, Florian},
  url = {https://www.heldermann.de/JCA/JCA25/JCA252/jca25037.htm},
  date = {2018},
  journaltitle = convexanal,
  number = {2},
  pages = {595--622},
  title = {Semismooth implicit functions},
  volume = {28},
}

@article{Kummer:1988,
  author = {Kummer, Bernd},
  publisher = {Akademie-Verlag},
  date = {1988},
  journaltitle = {Mathematical Research},
  pages = {114--125},
  title = {{N}ewton's method for non-differentiable functions},
  volume = {45},
}

@article{Kummer:2000,
  author = {Kummer, Bernd},
  date = {2000},
  doi = {10.7151/dmdico.1013},
  journaltitle = {Discussiones Mathematicae. Differential Inclusions, Control and Optimization},
  number = {2},
  pages = {209--244},
  title = {Generalized {N}ewton and {NCP}-methods: convergence, regularity, actions},
  volume = {20},
}

@book{Lang:1993,
  author = {Lang, Serge},
  publisher = {Springer-Verlag, New York},
  date = {1993},
  doi = {10.1007/978-1-4612-0897-6},
  edition = {3},
  series = {Graduate Texts in Mathematics},
  title = {Real and Functional Analysis},
  volume = {142},
}

@article{Lebourg:1975,
  author = {Lebourg, Gérard},
  date = {1975},
  journaltitle = {C. R. Acad. Sci. Paris Sér. A},
  number = {19},
  pages = {A795--A797},
  title = {Valeur moyenne pour gradient généralisé},
  volume = {281},
}

@article{Lebourg:1979,
  author = {Lebourg, Gérard},
  date = {1979},
  doi = {10.2307/1998104},
  journaltitle = transams,
  pages = {125--144},
  title = {Generic differentiability of {L}ipschitzian functions},
  volume = {256},
}

@article{leventhal2009metric,
  author = {Leventhal, D.},
  date = {2009},
  doi = {10.1016/j.jmaa.2009.07.012},
  journaltitle = jmaa,
  number = {2},
  pages = {681--688},
  title = {Metric subregularity and the proximal point method},
  volume = {360},
}

@article{lewis2002active,
  author = {Lewis, A. S.},
  date = {2002},
  doi = {10.1137/S1052623401387623},
  journaltitle = siamjopt,
  number = {3},
  pages = {702--725},
  title = {Active sets, nonsmoothness, and sensitivity},
  volume = {13},
}

@article{li2013global,
  author = {Li, Guoyin},
  date = {2013},
  doi = {10.1007/s10107-011-0481-z},
  journaltitle = mathprog,
  number = {1},
  pages = {37--64},
  title = {Global error bounds for piecewise convex polynomials},
  volume = {137},
}

@article{li2012holder,
  author = {Li, Guoyin and Mordukhovich, Boris S.},
  date = {2012},
  doi = {10.1137/120864660},
  journaltitle = siamjopt,
  number = {4},
  pages = {1655--1684},
  title = {Hölder metric subregularity with applications to proximal point method},
  volume = {22},
}

@article{liang2014local,
  author = {Liang, Jingwei and Fadili, Jalal and Peyré, Gabriel},
  editor = {Ghahramani, Z. and Welling, M. and Cortes, C. and Lawrence, N. D. and Weinberger, K. Q.},
  url = {https://papers.nips.cc/paper/5260-local-linear-convergence-of-forward-backward-under-partial-smoothness.pdf},
  date = {2014},
  journaltitle = advneural,
  pages = {1970--1978},
  title = {Local linear convergence of forward--backward under partial smoothness},
  volume = {27},
}

@book{Lindenstrauss:2021,
  author = {Lindenstrauss, Joram and Preiss, David and Tišer, Jaroslav},
  publisher = {Princeton University Press, Princeton, NJ},
  date = {2012},
  series = {Annals of Mathematics Studies},
  title = {Fréchet Differentiability of {L}ipschitz Functions and Porous Sets in {B}anach Spaces},
  volume = {179},
}

@book{Lions:1971,
  author = {Lions, J.-L.},
  publisher = {Springer-Verlag, New York-Berlin},
  date = {1971},
  series = {Die Grundlehren der mathematischen Wissenschaften},
  title = {Optimal Control of Systems Governed by Partial Differential Equations},
  volume = {170},
}

@article{lionsmercier1979splitting,
  author = {Lions, P. and Mercier, B.},
  date = {1979},
  doi = {10.1137/0716071},
  journaltitle = siamjnum,
  number = {6},
  pages = {964--979},
  title = {Splitting algorithms for the sum of two nonlinear operators},
  volume = {16},
}

@article{liu2018partial,
  author = {Liu, Y. and Yuan, X. and Zeng, S. and Zhang, J.},
  date = {2018},
  doi = {10.1137/17M1144623},
  journaltitle = siamjnum,
  number = {4},
  pages = {2095--2123},
  title = {Partial error bound conditions and the linear convergence rate of the Alternating Direction Method of Multipliers},
  volume = {56},
}

@article{loris2011generalization,
  author = {Loris, Ignace and Verhoeven, Caroline},
  date = {2011},
  doi = {10.1088/0266-5611/27/12/125007},
  journaltitle = invprob,
  number = {12},
  pages = {125007},
  title = {On a generalization of the iterative soft-thresholding algorithm for the case of non-separable penalty},
  volume = {27},
}

@article{luo1992error,
  author = {Luo, Zhi-Quan and Tseng, Paul},
  date = {1992},
  doi = {10.1137/0802004},
  journaltitle = siamjopt,
  number = {1},
  pages = {43--54},
  title = {Error bound and convergence analysis of matrix splitting algorithms for the affine variational inequality problem},
  volume = {2},
}

@book{makela1992nonsmooth,
  author = {Mäkelä, Marko M. and Neittaanmäki, Pekka},
  publisher = {World Scientific Publishing Co., Inc., River Edge, NJ},
  date = {1992},
  doi = {10.1142/1493},
  note = {Analysis and algorithms with applications to optimal control},
  title = {Nonsmooth Optimization},
}

@article{malitsky2018linesearch,
  author = {Malitsky, Y. and Pock, T.},
  date = {2018},
  doi = {10.1137/16M1092015},
  journaltitle = siamjopt,
  number = {1},
  pages = {411--432},
  title = {A first-order primal-dual algorithm with linesearch},
  volume = {28},
}

@article{MalitskyTam:2020,
  author = {Malitsky, Yura and Tam, Matthew K.},
  date = {2020},
  doi = {10.1137/18M1207260},
  journaltitle = siamjopt,
  number = {2},
  pages = {1451--1472},
  title = {A forward-backward splitting method for monotone inclusions without cocoercivity},
  volume = {30},
}

@article{Mann1953mean,
  author = {Mann, W. Robert},
  date = {1953},
  doi = {10.2307/2032162},
  journaltitle = procams,
  pages = {506--510},
  title = {Mean value methods in iteration},
  volume = {4},
}

@article{martinet1970regularisation,
  author = {Martinet, B.},
  date = {1970},
  journaltitle = {Rev. Française Informat. Recherche Opérationnelle},
  number = {Sér. {R}-3},
  pages = {154--158},
  title = {Régularisation d'inéquations variationnelles par approximations successives},
  volume = {4},
}

@article{tuomov-nlpdhgm-block,
  author = {Mazurenko, Stanislav and Jauhiainen, Jyrki and Valkonen, Tuomo},
  date = {2020},
  doi = {10.1553/etna_vol52s509},
  eprint = {1911.06284},
  eprinttype = {arxiv},
  journaltitle = etna,
  origdate = {2019},
  pages = {509--552},
  title = {Primal-dual block-proximal splitting for a class of non-convex problems},
  volume = {52},
}

@article{mehlitzwachsmuth2018limiting,
  author = {Mehlitz, Patrick and Wachsmuth, Gerd},
  date = {2018},
  doi = {10.1007/s11228-016-0393-4},
  journaltitle = svaa,
  number = {3},
  pages = {449--467},
  title = {The limiting normal cone to pointwise defined sets in {L}ebesgue spaces},
  volume = {26},
}

@article{mehlitzwachsmuth2019decomposable,
  author = {Mehlitz, Patrick and Wachsmuth, Gerd},
  date = {2019},
  doi = {10.1007/s11228-017-0464-1},
  journaltitle = svaa,
  number = {1},
  pages = {265--294},
  title = {The weak sequential closure of decomposable sets in {L}ebesgue spaces and its application to variational geometry},
  volume = {27},
}

@article{Schiela:2010,
  author = {Meyer, C. and Panizzi, L. and Schiela, A.},
  date = {2011},
  doi = {10.1080/01630563.2011.587074},
  journaltitle = {Numerical Functional Analysis and Optimization},
  number = {9},
  pages = {983--1007},
  title = {Uniqueness criteria for solutions of the adjoint equation in state-constrained optimal control},
  volume = {32},
}

@article{Mifflin:1977,
  author = {Mifflin, Robert},
  date = {1977},
  doi = {10.1137/0315061},
  journaltitle = siamjcon,
  number = {6},
  pages = {959--972},
  title = {Semismooth and semiconvex functions in constrained optimization},
  volume = {15},
}

@article{mollenhoff2014primal,
  author = {Möllenhoff, Thomas and Strekalovskiy, Evgeny and Moeller, Michael and Cremers, Daniel},
  date = {2015},
  doi = {10.1137/140976601},
  journaltitle = siamjimg,
  number = {2},
  pages = {827--857},
  title = {The primal-dual hybrid gradient method for semiconvex splittings},
  volume = {8},
}

@article{Mordukhovich1980metric,
  author = {Morduhovič, B. Š.},
  date = {1980},
  journaltitle = doklady,
  number = {5},
  pages = {1072--1076},
  title = {Metric approximations and necessary conditions for optimality for general classes of nonsmooth extremal problems},
  volume = {254},
}

@book{Mordukhovich:2006,
  author = {Mordukhovich, Boris S.},
  publisher = {Springer},
  date = {2006},
  doi = {10.1007/3-540-31247-1},
  series = {Grundlehren der mathematischen Wissenschaften},
  subtitle = {Basic Theory},
  title = {Variational Analysis and Generalized Differentiation {I}},
  volume = {330},
}

@book{mordukhovich2018variational,
  author = {Mordukhovich, Boris S.},
  publisher = {Springer International Publishing},
  date = {2018},
  doi = {10.1007/978-3-319-92775-6},
  series = {Springer Monographs in Mathematics},
  title = {Variational Analysis and Applications},
}

@book{Mordukhovich:2024,
  author = {Mordukhovich, Boris S.},
  publisher = {Springer, Cham},
  date = {2024},
  doi = {10.1007/978-3-031-53476-8},
  series = {Springer Series in Operations Research and Financial Engineering},
  subtitle = {Theory, Algorithms, Applications},
  title = {Second-order Variational Analysis in Optimization, Variational Stability, and Control},
}

@article{mordukhovich1976maximum,
  author = {Mordukhovich, Boris Shlomovich},
  date = {1976},
  doi = {10.1016/0021-8928(76)90136-2},
  journaltitle = jmaa,
  number = {6},
  pages = {960--969},
  title = {Maximum principle in the problem of time optimal response with nonsmooth constraints},
  volume = {40},
}

@article{moreau1965proximite,
  author = {Moreau, Jean-Jacques},
  date = {1965},
  doi = {10.24033/bsmf.1625},
  journaltitle = {Bulletin de la Société mathématique de France},
  pages = {273--299},
  title = {Proximité et dualité dans un espace hilbertien},
  volume = {93},
}

@article{motzkin_schoenberg_1954,
  author = {Motzkin, T. S. and Schoenberg, I. J.},
  publisher = {Cambridge University Press},
  date = {1954},
  doi = {10.4153/CJM-1954-038-x},
  journaltitle = {Canadian Journal of Mathematics},
  pages = {393--404},
  title = {The relaxation method for linear inequalities},
  volume = {6},
}

@book{Mueller:2012,
  author = {Mueller, Jennifer L. and Siltanen, Samuli},
  publisher = {Society for Industrial and Applied Mathematics},
  date = {2012},
  doi = {10.1137/1.9781611972344},
  title = {Linear and Nonlinear Inverse Problems with Practical Applications},
}

@book{natterer2001mathematics,
  author = {Natterer, F.},
  publisher = {Society for Industrial and Applied Mathematics},
  date = {2001},
  doi = {10.1137/1.9780898719284},
  title = {The Mathematics of Computerized Tomography},
}

@article{Nesterov,
  author = {Nesterov, Yu. E.},
  date = {1983},
  journaltitle = doklady,
  number = {2},
  pages = {372--376},
  title = {A method for solving the convex programming problem with convergence rate {$O(1/k^{2})$}},
  volume = {27},
}

@book{Nesterov:2004,
  author = {Nesterov, Yurii},
  publisher = {Kluwer Academic Publishers, Boston, MA},
  date = {2004},
  doi = {10.1007/978-1-4419-8853-9},
  series = {Applied Optimization},
  title = {Introductory Lectures on Convex Optimization},
  volume = {87},
}

@book{Nesterov:2018,
  author = {Nesterov, Yurii},
  publisher = {Springer, Cham},
  date = {2018},
  doi = {10.1007/978-3-319-91578-4},
  series = {Springer Optimization and Its Applications},
  title = {Lectures on Convex Optimization},
  volume = {137},
}

@article{ngai2008error,
  author = {Ngai, Huynh Van and Théra, Michel},
  date = {2008},
  doi = {10.1137/060675721},
  journaltitle = siamjopt,
  number = {1},
  pages = {1--20},
  title = {Error bounds in metric spaces and application to the perturbation stability of metric regularity},
  volume = {19},
}

@book{nishimura1996principles,
  author = {Nishimura, D.G.},
  publisher = {Stanford University},
  date = {1996},
  title = {Principles of Magnetic Resonance Imaging},
}

@article{ochs2017adaptive,
  author = {Ochs, Peter and Pock, Thomas},
  date = {2019},
  doi = {10.1137/17M1156678},
  eprint = {1711.04343},
  eprinttype = {arXiv},
  journaltitle = siamjopt,
  number = {4},
  pages = {2482--2503},
  title = {Adaptive {FISTA} for non-convex optimization},
  volume = {29},
}

@article{opial1967weak,
  author = {Opial, Zdzis{ł}aw},
  date = {1967},
  doi = {10.1090/S0002-9904-1967-11761-0},
  journaltitle = {Bulletin of the American Mathematical Society},
  number = {4},
  pages = {591--597},
  title = {Weak convergence of the sequence of successive approximations for nonexpansive mappings},
  volume = {73},
}

@book{outrata1998nonsmooth,
  author = {Outrata, Jiři and Kočvara, Michal and Zowe, Jochem},
  publisher = {Kluwer Academic Publishers, Dordrecht},
  date = {1998},
  doi = {10.1007/978-1-4757-2825-5},
  note = {Theory, applications and numerical results},
  series = {Nonconvex Optimization and its Applications},
  title = {Nonsmooth Approach to Optimization Problems with Equilibrium Constraints},
  volume = {28},
}

@article{Boyd:2014,
  author = {Parikh, Neil and Boyd, Stephen},
  date = {2014},
  doi = {10.1561/2400000003},
  journaltitle = {Foundations and Trends in Optimization},
  number = {3},
  pages = {123--231},
  title = {Proximal algorithms},
  volume = {1},
}

@inproceedings{patrinos2014douglas,
  author = {Patrinos, P. and Stella, L. and Bemporad, A.},
  booktitle = {53rd IEEE Conference on Decision and Control},
  date = {2014},
  doi = {10.1109/CDC.2014.7040049},
  pages = {4234--4239},
  title = {{D}ouglas--{R}achford splitting: complexity estimates and accelerated variants},
}

@book{peano,
  author = {Peano, Giuseppe},
  publisher = {Fratelli Boca, Torino},
  date = {1908},
  title = {Formulario Mathematico},
}

@book{Penot:2013,
  author = {Penot, Jean-Paul},
  publisher = {Springer, New York},
  date = {2013},
  doi = {10.1007/978-1-4614-4538-8},
  series = {Graduate Texts in Mathematics},
  title = {Calculus Without Derivatives},
  volume = {266},
}

@article{petryshyn1966construction,
  author = {Petryshyn, W. V.},
  date = {1966},
  doi = {10.1016/0022-247X(66)90027-8},
  journaltitle = jmaa,
  number = {2},
  pages = {276--284},
  title = {Construction of fixed points of demicompact mappings in {H}ilbert space},
  volume = {14},
}

@book{Peypouquet:2015,
  author = {Peypouquet, Juan},
  publisher = {Springer, Cham},
  date = {2015},
  doi = {10.1007/978-3-319-13710-0},
  series = {SpringerBriefs in Optimization},
  subtitle = {Theory, Methods and Examples},
  title = {Convex Optimization in Normed Spaces},
}

@inproceedings{pock2009mumford,
  author = {Pock, T. and Cremers, D. and Bischof, H. and Chambolle, A.},
  booktitle = {12th IEEE Conference on Computer Vision},
  date = {2009},
  doi = {10.1109/ICCV.2009.5459348},
  pages = {1133--1140},
  title = {An algorithm for minimizing the {Mumford--Shah} functional},
}

@article{poliquin1996prox,
  author = {Poliquin, René and Rockafellar, R. T.},
  date = {1996},
  doi = {10.1090/S0002-9947-96-01544-9},
  journaltitle = transams,
  number = {5},
  pages = {1805--1838},
  title = {Prox-regular functions in variational analysis},
  volume = {348},
}

@article{Preiss:1990,
  author = {Preiss, D.},
  date = {1990},
  doi = {10.1016/0022-1236(90)90147-D},
  journaltitle = {Journal of Functional Analysis},
  number = {2},
  pages = {312--345},
  title = {Differentiability of {L}ipschitz functions on {B}anach spaces},
  volume = {91},
}

@article{Qi:1993,
  author = {Qi, Liqun},
  date = {1993},
  doi = {10.1287/moor.18.1.227},
  journaltitle = {Mathematics of Operations Research},
  number = {1},
  pages = {227--244},
  title = {Convergence analysis of some algorithms for solving nonsmooth equations},
  volume = {18},
}

@article{Qi:1993a,
  author = {Qi, Liqun and Sun, Jie},
  date = {1993},
  doi = {10.1007/bf01581275},
  journaltitle = mathprog,
  number = {3, Ser. A},
  pages = {353--367},
  title = {A nonsmooth version of {N}ewton's method},
  volume = {58},
}

@book{Rindler:2018,
  author = {Rindler, Filip},
  publisher = {Cham: Springer},
  date = {2018},
  doi = {10.1007/978-3-319-77637-8},
  series = {Universitext},
  title = {Calculus of Variations},
}

@incollection{robinson1981continuity,
  author = {Robinson, Stephen M.},
  editor = {König, H. and Korte, B. and Ritter, K.},
  location = {Berlin, Heidelberg},
  publisher = {Springer Berlin Heidelberg},
  booktitle = {Mathematical Programming at Oberwolfach},
  date = {1981},
  doi = {10.1007/BFb0120929},
  pages = {206--214},
  title = {Some continuity properties of polyhedral multifunctions},
}

@article{Rockafellar:1970,
  author = {Rockafellar, R. Tyrrell},
  date = {1970},
  doi = {10.2140/pjm.1970.33.209},
  journaltitle = {Pacific Journal of Mathematics},
  pages = {209--216},
  title = {On the maximal monotonicity of subdifferential mappings},
  volume = {33},
}

@incollection{Rockafellar:1976,
  author = {Rockafellar, R. Tyrrell},
  location = {Berlin},
  publisher = {Springer},
  booktitle = {Nonlinear Operators and the Calculus of Variations ({S}ummer {S}chool, {U}niv. {L}ibre {B}ruxelles, {B}russels, 1975)},
  date = {1976},
  doi = {10.1007/bfb0079944},
  pages = {157--207},
  series = {Lecture Notes in Math.},
  title = {Integral functionals, normal integrands and measurable selections},
  volume = {543},
}

@article{rockafellar1976proximal,
  author = {Rockafellar, R. Tyrrell},
  date = {1976},
  doi = {10.1137/0314056},
  journaltitle = siamjcon,
  number = {5},
  pages = {877--898},
  title = {Monotone operators and the proximal point algorithm},
  volume = {14},
}

@inproceedings{rockafellar1981favorable,
  author = {Rockafellar, R. Tyrrell},
  location = {IIASA, Laxenburg, Austria},
  url = {https://pure.iiasa.ac.at/id/eprint/1760/},
  booktitle = {Progress in Nondifferentiable Optimization},
  date = {1981},
  pages = {125--143},
  title = {Favorable classes of {L}ipschitz continuous functions in subgradient optimization},
}

@article{rockafellar1985maximal,
  author = {Rockafellar, R. Tyrrell},
  date = {1985},
  doi = {10.1016/S0294-1449(16)30401-2},
  journaltitle = {Annales de l'Institut Henri Poincare (C) Non Linear Analysis},
  number = {3},
  pages = {167--184},
  title = {Maximal monotone relations and the second derivatives of nonsmooth functions},
  volume = {2},
}

@article{rockafellar1988first,
  author = {Rockafellar, R. Tyrrell},
  date = {1988},
  doi = {10.1090/S0002-9947-1988-0936806-9},
  journaltitle = transams,
  number = {1},
  pages = {75--108},
  title = {First-and second-order epi-differentiability in nonlinear programming},
  volume = {307},
}

@article{rockafellar1989proto,
  author = {Rockafellar, R. Tyrrell},
  date = {1989},
  doi = {10.1016/S0294-1449(17)30034-3},
  issn = {0294-1449},
  journaltitle = {Annales de l'Institut Henri Poincare (C) Non Linear Analysis},
  number = {S6},
  pages = {449--482},
  title = {Proto-differentiability of set-valued mappings and its applications in optimization},
  volume = {6},
}

@book{Rockafellar:1998,
  author = {Rockafellar, R. Tyrrell and Wets, Roger J.-B.},
  publisher = {Springer},
  date = {1998},
  doi = {10.1007/978-3-642-02431-3},
  series = {Grundlehren der mathematischen Wissenschaften},
  title = {Variational Analysis},
  volume = {317},
}

@book{RoysetWets:2022,
  author = {Royset, Johannes O. and Wets, Roger J.-B.},
  location = {Cham},
  publisher = {Springer},
  date = {2021},
  doi = {10.1007/978-3-030-76275-9},
  title = {An Optimization Primer},
}

@article{ROF,
  author = {Rudin, Leonid I and Osher, Stanley and Fatemi, Emad},
  date = {1992},
  doi = {10.1016/0167-2789(92)90242-f},
  journaltitle = {Physica D: Nonlinear Phenomena},
  number = {1},
  pages = {259--268},
  title = {Nonlinear total variation based noise removal algorithms},
  volume = {60},
}

@book{Rudin:1991,
  author = {Rudin, Walter},
  location = {New York},
  publisher = {De Gruyter},
  date = {2021},
  doi = {10.1515/9783110750430},
  edition = {2},
  title = {Analysis},
}

@book{Ruszczynski:2006a,
  author = {Ruszczyǹski, Andrzej},
  location = {Princeton, NJ},
  publisher = {Princeton University Press},
  date = {2006},
  doi = {10.2307/j.ctvcm4hcj},
  title = {{N}onlinear {O}ptimization},
}

@book{rynne2008functional,
  author = {Rynne, Bryan P. and Youngson, Martin A.},
  publisher = {Springer, London},
  date = {2008},
  doi = {10.1007/978-1-84800-005-6},
  edition = {2},
  series = {Springer Undergraduate Mathematics Series},
  title = {Linear Functional Analysis},
}

@article{Ruy:2019,
  author = {Ryu, Ernest K.},
  publisher = {Springer Science and Business Media LLC},
  date = {2019},
  doi = {10.1007/s10107-019-01403-1},
  issn = {1436-4646},
  journaltitle = mathprog,
  number = {1–2},
  pages = {233--273},
  title = {Uniqueness of DRS as the 2 operator resolvent-splitting and impossibility of 3 operator resolvent-splitting},
  volume = {182},
}

@article{schaefer1957,
  author = {Schaefer, Helmut},
  url = {https://eudml.org/doc/146424},
  date = {1957},
  journaltitle = {Jahresbericht der {D}eutschen {M}athematiker-{V}ereinigung},
  pages = {131--140},
  title = {Über die Methode sukzessiver Approximationen},
  volume = {59},
}

@book{Scherzer,
  author = {Scherzer, Otmar and Grasmair, Markus and Grossauer, Harald and Haltmeier, Markus and Lenzen, Frank},
  publisher = {Springer, New York},
  date = {2009},
  doi = {10.1007/978-0-387-69277-7},
  title = {Variational Methods in Imaging},
}

@article{Schiela:2008a,
  author = {Schiela, Anton},
  date = {2008},
  doi = {10.1137/060674375},
  journaltitle = siamjopt,
  number = {3},
  pages = {1417--1432},
  title = {A simplified approach to semismooth Newton methods in function space},
  volume = {19},
}

@book{Schirotzek:2007,
  author = {Schirotzek, Winfried},
  location = {Berlin},
  publisher = {Springer},
  date = {2007},
  doi = {10.1007/978-3-540-71333-3},
  series = {Universitext},
  title = {Nonsmooth Analysis},
}

@book{Scholtes:2012,
  author = {Scholtes, Stefan},
  publisher = {Springer, New York},
  date = {2012},
  doi = {10.1007/978-1-4614-4340-7},
  series = {Springer Briefs in Optimization},
  title = {Introduction to Piecewise Differentiable Equations},
}

@book{Schuster:2012,
  author = {Schuster, T. and Kaltenbacher, B. and Hofmann, B. and Kazimierski, K.S.},
  publisher = {De Gruyter},
  date = {2012},
  doi = {10.1515/9783110255720},
  isbn = {9783110255720},
  series = {Radon Series on Computational and Applied Mathematics},
  title = {Regularization Methods in Banach Spaces},
}

@article{Simons:2009,
  author = {Simons, S.},
  url = {https://www.heldermann.de/JCA/JCA16/JCA161/jca16008.htm},
  date = {2009},
  journaltitle = convexanal,
  number = {1},
  pages = {165--168},
  title = {A new proof of the maximal monotonicity of subdifferentials},
  volume = {16},
}

@article{Stadler:2007a,
  author = {Stadler, Georg},
  date = {2009},
  doi = {10.1007/s10589-007-9150-9},
  journaltitle = computoptappl,
  number = {2},
  pages = {159--181},
  title = {Elliptic optimal control problems with {L$^1$}-control cost and applications for the placement of control devices},
  volume = {44},
}

@article{thibault1995integration,
  author = {Thibault, L. and Zagrodny, D.},
  date = {1995},
  doi = {10.1006/jmaa.1995.1003},
  journaltitle = jmaa,
  number = {1},
  pages = {33--58},
  title = {Integration of subdifferentials of lower semicontinuous functions on {Banach} spaces},
  volume = {189},
}

@article{thibault1983tangent,
  author = {Thibault, Lionel},
  date = {1983},
  doi = {10.2307/1999227},
  journaltitle = transams,
  number = {2},
  pages = {601--621},
  title = {Tangent cones and quasi-interiorly tangent cones to multifunctions},
  volume = {277},
}

@book{Troeltzsch,
  author = {Tröltzsch, Fredi},
  publisher = {American Mathematical Society, Providence, RI},
  date = {2010},
  doi = {10.1090/gsm/112},
  note = {Translated from the 2005 German original by Jürgen Sprekels},
  series = {Graduate Studies in Mathematics},
  subtitle = {Theory, methods and applications},
  title = {Optimal Control of Partial Differential Equations},
  volume = {112},
}

@article{Ulbrich:2002a,
  author = {Ulbrich, Michael},
  date = {2002},
  doi = {10.1137/s1052623400371569},
  journaltitle = siamjopt,
  number = {3},
  pages = {805--841},
  title = {Semismooth {N}ewton methods for operator equations in function spaces},
  volume = {13},
}

@book{Ulbrich:2011,
  author = {Ulbrich, Michael},
  location = {Philadelphia},
  publisher = {Society for Industrial and Applied Mathematics},
  date = {2011},
  doi = {10.1137/1.9781611970692},
  series = {MOS-SIAM Series on Optimization},
  title = {Semismooth {N}ewton Methods for Variational Inequalities and Constrained Optimization Problems in Function Spaces},
  volume = {11},
}

@article{Valkonen:2014,
  author = {Valkonen, Tuomo},
  date = {2014},
  doi = {10.1088/0266-5611/30/5/055012},
  journaltitle = invprob,
  number = {5},
  pages = {055012},
  title = {A primal-dual hybrid gradient method for nonlinear operators with applications to {MRI}},
  volume = {30},
}

@article{tuomov-blockcp,
  author = {Valkonen, Tuomo},
  url = {https://tuomov.iki.fi/m/blockcp.pdf},
  date = {2019},
  doi = {10.1553/etna_vol51s15},
  eprint = {1609.07373},
  eprinttype = {arXiv},
  journaltitle = etna,
  origdate = {2016},
  pages = {15--49},
  title = {Block-proximal methods with spatially adapted acceleration},
  volume = {51},
}

@article{tuomov-inertia,
  author = {Valkonen, Tuomo},
  date = {2020},
  doi = {10.1137/18M1182851},
  eprint = {1804.08736},
  eprinttype = {arxiv},
  journaltitle = siamjopt,
  number = {2},
  origdate = {2018},
  pages = {1391--1420},
  title = {Inertial, corrected, primal-dual proximal splitting},
  volume = {30},
}

@article{tuomov-proxtest,
  author = {Valkonen, Tuomo},
  url = {https://tuomov.iki.fi/m/proxtest.pdf},
  date = {2020},
  doi = {10.1007/s00245-018-9541-6},
  eprint = {1703.05705},
  eprinttype = {arxiv},
  journaltitle = amo,
  number = {2},
  pages = {591--636},
  title = {Testing and non-linear preconditioning of the proximal point method},
  volume = {82},
}

@incollection{tuomov-firstorder,
  author = {Valkonen, Tuomo},
  editor = {Chen, Ke and Schönlieb, Carola-Bibiane and Tai, Xue-Cheng and Younes, Laurent},
  publisher = {Springer, Cham},
  booktitle = {Handbook of Mathematical Models and Algorithms in Computer Vision and Imaging},
  date = {2021},
  doi = {10.1007/978-3-030-03009-4_93-1},
  eprint = {1910.00115},
  eprinttype = {arxiv},
  isbn = {978-3-030-03009-4},
  origdate = {2019},
  title = {First-order primal-dual methods for nonsmooth nonconvex optimisation},
}

@article{tuomov-subreg,
  author = {Valkonen, Tuomo},
  url = {https://tuomov.iki.fi/m/proxtest.pdf},
  date = {2021},
  eprint = {1711.05123},
  eprinttype = {arxiv},
  journaltitle = convexanal,
  number = {1},
  origdate = {2017},
  pages = {251--278},
  title = {Preconditioned proximal point methods and notions of partial subregularity},
  volume = {28},
}

@article{tuomov-regtheory,
  author = {Valkonen, Tuomo},
  date = {2021},
  doi = {10.1088/1361-6420/abe4aa},
  eprint = {2011.07575},
  eprinttype = {arxiv},
  journaltitle = invprob,
  number = {4},
  pages = {045010},
  title = {Regularisation, optimisation, subregularity},
  volume = {37},
}

@article{Vossen:2006a,
  author = {Vossen, G. and Maurer, H.},
  date = {2006},
  doi = {10.1002/oca.781},
  journaltitle = {Optimal Control Appl. Methods},
  number = {6},
  pages = {301--321},
  title = {On {L$\sp 1$}-minimization in optimal control and applications to robotics},
  volume = {27},
}

@article{vu2013splitting,
  author = {Vũ, Bằng Công},
  date = {2013},
  doi = {10.1007/s10444-011-9254-8},
  journaltitle = {Advances in Computational Mathematics},
  number = {3},
  pages = {667--681},
  title = {A splitting algorithm for dual monotone inclusions involving cocoercive operators},
  volume = {38},
}

@article{wright2015coordinate,
  author = {Wright, Stephen J.},
  date = {2015},
  doi = {10.1007/s10107-015-0892-3},
  journaltitle = mathprog,
  number = {1},
  pages = {3--34},
  title = {Coordinate descent algorithms},
  volume = {151},
}

@article{wrightnovak2009sparse,
  author = {Wright, Stephen J. and Nowak, Robert D. and Figueiredo, Mário A. T.},
  date = {2009},
  doi = {10.1109/TSP.2009.2016892},
  journaltitle = ieeesig,
  number = {7},
  pages = {2479--2493},
  title = {Sparse reconstruction by separable approximation},
  volume = {57},
}

@article{Yang:2008a,
  author = {Yang, Junfeng and Zhang, Yin and Yin, Wotao},
  date = {2009},
  doi = {10.1137/080732894},
  journaltitle = siamjscicomp,
  number = {4},
  pages = {2842--2865},
  title = {An efficient {TVL1} algorithm for deblurring multichannel images corrupted by impulsive noise},
  volume = {31},
}

@article{yost1993asplund,
  author = {Yost, David},
  publisher = {Charles University in Prague},
  url = {https://dml.cz/dmlcz/702006},
  date = {1993},
  journaltitle = {Acta Universitatis Carolinae. Mathematica et Physica},
  number = {2},
  pages = {159--177},
  title = {Asplund spaces for beginners},
  volume = {34},
}

@book{Zalinescu:2002,
  author = {Zălinescu, C.},
  publisher = {World Scientific Publishing Co., Inc., River Edge, NJ},
  date = {2002},
  doi = {10.1142/9789812777096},
  title = {Convex Analysis in General Vector Spaces},
}

@article{zhang2011unified,
  author = {Zhang, Xiaoqun and Burger, Martin and Osher, Stanley},
  date = {2011},
  doi = {10.1007/s10915-010-9408-8},
  journaltitle = jscicomp,
  number = {1},
  pages = {20--46},
  title = {A unified primal-dual algorithm framework based on {B}regman iteration},
  volume = {46},
}

@article{zheng2010metric,
  author = {Zheng, Xi Yin and Ng, Kung Fu},
  date = {2010},
  doi = {10.1137/090772174},
  journaltitle = siamjopt,
  number = {5},
  pages = {2119--2136},
  title = {Metric subregularity and calmness for nonconvex generalized equations in {Banach} spaces},
  volume = {20},
}

@article{zhou2017unified,
  author = {Zhou, Zirui and So, Anthony Man-Cho},
  date = {2017},
  doi = {10.1007/s10107-016-1100-9},
  journaltitle = mathprog,
  number = {2},
  pages = {689--728},
  title = {A unified approach to error bounds for structured convex optimization problems},
  volume = {165},
}

@report{zhu2008efficient,
  author = {Zhu, Mingqiang and Chan, Tony},
  institution = {UCLA},
  url = {https://ww3.math.ucla.edu/camreport/cam08-34.pdf},
  date = {2008},
  number = {08-34},
  title = {An efficient primal-dual hybrid gradient algorithm for total variation image restoration},
  type = {CAM Report},
}

\printindex

\end{document}